\documentclass[twoside,11pt,dvipsnames, openleft]{book}

\usepackage{varioref}
\usepackage{jmlr2e}
\usepackage{bm}
\usepackage{dsfont}

\renewcommand{\arraystretch}{1.15}

\usepackage[english]{babel}

\usepackage{imakeidx}
\makeindex[columns=2, title=Alphabetical Index]

\usepackage{stackengine}

\usepackage{cancel}
\usepackage{dsfont}

\usepackage{extarrows}

\usepackage[framemethod=tikz]{mdframed}

\usepackage{amsmath}
\usepackage{arydshln}
\numberwithin{equation}{chapter}

\usepackage{amsmath}
\usepackage{blkarray}

\usepackage[absolute,overlay]{textpos}

\usepackage[final]{pdfpages}

\usepackage[noframe]{showframe}
\usepackage{framed}
\usepackage{lipsum}

\usepackage{xcolor}

\usepackage{tcolorbox}
\usepackage{graphicx}    
\usepackage[hang]{subfigure}

\usepackage{algorithm}
\usepackage{algpseudocode}
\usepackage{amsmath}
\usepackage{graphics}
\usepackage{epsfig}
\usepackage{blkarray}
\usepackage{listings}
\usepackage{booktabs}  

\usepackage{tikz}
\usetikzlibrary{mindmap,trees,backgrounds}
\usetikzlibrary{arrows.meta}
\usetikzlibrary{patterns}

\usetikzlibrary{decorations.text}

\usetikzlibrary{calc,backgrounds}

\newcommand*\circled[1]{\tikz[baseline=(char.base)]{
		\node[shape=circle,draw,inner sep=2pt] (char) {#1};}}

\usepackage{changepage}                 

\usepackage{hhline}

\usetikzlibrary{positioning,shapes,shadows,arrows}
\tikzstyle{condition}=[rectangle, draw=black, rounded corners, fill=colorqr, drop shadow,
text centered, anchor=north, text=black, text width=3cm]
\tikzstyle{abstract}=[rectangle, draw=black, rounded corners, fill=blue!30, drop shadow,
text centered, anchor=north, text=black, text width=3cm]
\tikzstyle{comment}=[rectangle, draw=black, rounded corners, fill=color1, drop shadow,
text centered, anchor=north, text=black, text width=3cm]
\tikzstyle{myarrow}=[->, >=open triangle 90, thick]
\tikzstyle{line}=[-, thick]

\usepackage{tikz-3dplot}
\tikzset{>=latex}
\usetikzlibrary{matrix}

\usepackage{sidecap}
\sidecaptionvpos{figure}{t}

\usepackage{verbatimbox}

\usepackage[mathscr]{euscript}

\usepackage{stmaryrd}

\newlength{\offsetpage}
	{\end{adjustwidth}}
	{\end{adjustwidth}}

\usepackage{minitoc}   
\let\cleardoublepage\clearpage
\usepackage[titletoc]{appendix}

\definecolor{titlepagecolor}{cmyk}{75,68,67,90}
\definecolor{titlepagecolor2}{rgb}{1.0, 0.08, 0.58}
\definecolor{emerald}{rgb}{0.31, 0.78, 0.47}
\definecolor{deeppink}{HTML}{D14064}
\definecolor{lowpink}{HTML}{ffe6ec}
\newcommand{\partcolor}{gray!65} 
\definecolor{lowblue}{HTML}{E1EBFE}

\usepackage{afterpage}

\makeatletter \renewcommand*\cleardoublepage{
	\clearpage
	\if@twoside   
	\ifodd\c@page 
	\hbox{}\newpage
	\if@twocolumn\hbox{}   
	\newpage
	\fi
	\fi
	\fi
} \makeatother
\let\originalpart=\part
\def\part#1{\cleardoublepage\clearpage \pagecolor{\partcolor} \originalpart{#1}\nopagecolor }

\usepackage{blkarray}

\renewcommand\thesection{\thechapter.\arabic{section}}

\newenvironment{bmatrixscript}
{\scriptsize\begin{bmatrix}}
	{\end{bmatrix}\normalsize}

\newenvironment{bmatrixfoot}
{\footnotesize\begin{bmatrix}}
	{\end{bmatrix}\normalsize}

\newcommand{\longeq}{\mathrel{\rule[0.5ex]{2em}{0.5pt}\hspace{-2em}\rule[-0.5ex]{2em}{0.5pt}}}

\newcommand{\domain}{\mathrm{dom}}
\newcommand{\cond}{\text{cond} }

\newcommand{\conv}{\text{conv}}
\newcommand{\cone}{\text{cone}}

\newcommand\mathopmax[1]{\mathop{\max}_{#1}}
\newcommand\mathopmin[1]{\mathop{\min}_{#1}}

\newcommand{\holders}{{H\"older's }}

\newcommand{\projectS}{\mathcal{P}_{\mathbb{S}}}
\newcommand{\projectT}{\mathcal{P}_{\mathbb{T}}}
\newcommand{\project}{\mathcal{P}}
\newcommand{\bprox}{\mathcal{B}\mathrm{rox}}

\newcommand{\bproxfphi}{\mathcal{B}\mathrm{rox}_{f,\phi}}

\newcommand{\prox}{\mathcal{P}\mathrm{rox}}
\newcommand{\proxf}{\mathcal{P}\mathrm{rox}_f}

\newcommand{\indicatorS}{\delta_{\mathbb{S}}}
\newcommand{\indicatorG}{\delta}

\newcommand{\dom}{\mathrm{dom}}
\newcommand{\epi}{\mathrm{epi}}
\newcommand{\lev}{\mathrm{Lev}}
\newcommand{\closure}{\text{cl}}
\newcommand{\interior}{\text{int}}
\newcommand{\relint}{\text{relint}}
\newcommand{\aff}{\text{aff}}

\newcommand{\vol}{\mathrm{Vol}}
\newcommand{\spark}{\mathrm{spark}}
\newcommand{\supp}{\mathrm{supp}}
\newcommand\comple[1]{#1^c}

\newcommand{\hadaprod}{\circ}

{\endMakeFramed}

\makeatletter
\renewcommand{\l@section}{\@dottedtocline{1}{1.5em}{2.2em}}
\renewcommand{\l@subsection}{\@dottedtocline{2}{4.0em}{3.2em}}
\renewcommand{\l@subsubsection}{\@dottedtocline{3}{7.1em}{4.3em}}
\makeatother

\usepackage[labelfont=bf]{caption}
\newcommand\myhrulefill[1]{\leavevmode\leaders\hrule height#1\hfill\kern0pt}
\DeclareCaptionFormat{myformat}{
	{\color[RGB]{0,0,0}\myhrulefill{0.12em}}
	\\#1#2#3
}
\def\algoalign#1{\parbox[t]{\dimexpr\linewidth-\algorithmicindent}{#1}}

\usepackage{yfonts,color,lettrine}
\definecolor{caligraphcolor}{HTML}{74AECB}

\usepackage{setspace}

\AtBeginDocument{\setlength{\DefaultFindent}{0.5em}}
\usepackage[shortlabels]{enumitem}  
\usepackage{enumitem}

\newcommand*{\eitemi}{\tikz \draw [baseline, ball color=structurecolor,draw=none] circle (2pt);}
\newcommand*{\eitemii}{\tikz \draw [baseline, fill=structurecolor,draw=none,circular drop shadow] circle (2pt);}
\newcommand*{\eitemiii}{\tikz \draw [baseline, fill=structurecolor,draw=none] circle (2pt);}
\setlist[enumerate,1]{label=\color{black}\arabic*.,itemsep=0pt,partopsep=0pt,parsep=\parskip,topsep=5pt}
\setlist[enumerate,2]{label=\color{black}(\alph*).,itemsep=0pt,partopsep=0pt,parsep=\parskip,topsep=5pt}
\setlist[enumerate,3]{label=\color{black}\Roman*.,itemsep=0pt,partopsep=0pt,parsep=\parskip,topsep=5pt}
\setlist[enumerate,4]{label=\color{black}\Alph*.,itemsep=0pt,partopsep=0pt,parsep=\parskip,topsep=5pt}
\setlist[itemize,1]{label={\eitemi},itemsep=0pt,partopsep=0pt,parsep=\parskip,topsep=5pt}
\setlist[itemize,2]{label={\eitemii},itemsep=0pt,partopsep=0pt,parsep=\parskip,topsep=5pt}
\setlist[itemize,3]{label={\eitemiii},itemsep=0pt,partopsep=0pt,parsep=\parskip,topsep=5pt}

\usepackage{xcolor}
\usepackage[Conny]{fncychap}

\makeatletter  

\def\mghrulefill#1{\color{black}\leavevmode\leaders\hrule\@height #1\hfill\kern\z@}
\makeatother

\definecolor{chaptertitle}{RGB}{0,0,128}
\definecolor{chapternum}{RGB}{255,0,0}

\usepackage[colorlinks]{hyperref}
\usepackage{hyperref}
\hypersetup{
	colorlinks=true, 
	linktoc=all,     
	linkcolor=mydarkblue,  
	anchorcolor=blue,
	citecolor=mydarkgreen,
}

\usepackage{fancyhdr, blindtext}

\newcommand{\changefonts}{%
	\fontsize{9}{11}\selectfont
}

\fancypagestyle{plain}{%

	\fancyhf{}  
}

\usepackage{amsthm}
\newtheoremstyle{normalfontstyle} 
{3pt}                           
{3pt}                           
{\normalfont}                   
{}                              
{\bfseries}                     
{}                             
{ }                             
{}                              

\usepackage{thmtools}
\declaretheoremstyle[
spaceabove=3pt,
spacebelow=3pt,
headfont=\bfseries,
notefont=\bfseries, 
notebraces={(}{)}, 
bodyfont=\normalfont,
postheadspace=1em, 
]{normalfontboldhead}

\theoremstyle{normalfontstyle}
\newcommand{\BlackBox}{\rule{1.5ex}{1.5ex}}  
\renewenvironment{proof}{\par\noindent{\bf Proof\ }}{\hfill\BlackBox\\[2mm]}

\declaretheorem[style=normalfontboldhead, name=Definition, numberlike=theo]{definitionT}
\newmdenv[skipabove=7pt,
skipbelow=7pt,
rightline=false,
leftline=true,
topline=false,
bottomline=false,
linecolor=mydarkblue,
innerleftmargin=5pt,
innerrightmargin=5pt,
innertopmargin=0pt,
leftmargin=2cm,
rightmargin=0cm,
linewidth=4pt,
innerbottommargin=0pt]{dBox}
\newenvironment{definition}{\begin{dBox}\begin{definitionT}}{\end{definitionT}\end{dBox}}

\declaretheorem[style=normalfontboldhead, name=Exercise, numberlike=theo]{exerciseC}
\newmdenv[skipabove=7pt,
skipbelow=7pt,
rightline=false,
leftline=true,
topline=false,
bottomline=false,
linecolor=mydarkgreen,
innerleftmargin=5pt,
innerrightmargin=5pt,
innertopmargin=0pt,
leftmargin=2cm,
rightmargin=0cm,
linewidth=4pt,
innerbottommargin=0pt]{eBox}
\newenvironment{exercise}{\begin{eBox}\begin{exerciseC}}{\end{exerciseC}\end{eBox}}

\declaretheorem[style=normalfontboldhead, name=Remark, numberlike=theo]{remarekC}
\newmdenv[skipabove=7pt,
skipbelow=7pt,
rightline=false,
leftline=true,
topline=false,
bottomline=false,
linecolor=mydarkpurple,
innerleftmargin=5pt,
innerrightmargin=5pt,
innertopmargin=0pt,
leftmargin=2cm,
rightmargin=0cm,
linewidth=4pt,
innerbottommargin=0pt]{rBox}
\newenvironment{remark}{\begin{rBox}\begin{remarekC}}{\end{remarekC}\end{rBox}}

\declaretheorem[style=normalfontboldhead, name=Assumption, numberlike=theo]{assumptionC}
\newmdenv[skipabove=7pt,
skipbelow=7pt,
rightline=false,
leftline=true,
topline=false,
bottomline=false,
linecolor=mydarkpurple,
innerleftmargin=5pt,
innerrightmargin=5pt,
innertopmargin=0pt,
leftmargin=2cm,
rightmargin=0cm,
linewidth=4pt,
innerbottommargin=0pt]{asBox}
\newenvironment{assumption}{\begin{asBox}\begin{assumptionC}}{\end{assumptionC}\end{asBox}}

\declaretheorem[style=normalfontboldhead, name=Condition, numberlike=theo]{conditionC}
\newmdenv[skipabove=7pt,
skipbelow=7pt,
rightline=false,
leftline=true,
topline=false,
bottomline=false,
linecolor=mydarkpurple,
innerleftmargin=5pt,
innerrightmargin=5pt,
innertopmargin=0pt,
leftmargin=2cm,
rightmargin=0cm,
linewidth=4pt,
innerbottommargin=0pt]{cdBox}
\newenvironment{condition}{\begin{cdBox}\begin{conditionC}}{\end{conditionC}\end{cdBox}}

\declaretheorem[style=normalfontboldhead, name=Example, numberlike=theo]{exampleC}
\newmdenv[skipabove=7pt,
skipbelow=7pt,
rightline=false,
leftline=false,
topline=false,
bottomline=false,
linecolor=mydarkgreen,
innerleftmargin=1pt,
innerrightmargin=5pt,
innertopmargin=0pt,
leftmargin=2cm,
rightmargin=0cm,
linewidth=4pt,
innerbottommargin=0pt]{xBox}
\newenvironment{example}{\begin{xBox}\begin{exampleC}}{\exampbar\end{exampleC}\end{xBox}}

\usepackage{adforn}  
\newcommand{\xchaptertitle}{Chapter~\thechapter~}
\newcommand{\problemname}{Problems}
\newenvironment{problemset}[1][\xchaptertitle~\problemname]{
	\vspace*{10pt}
	\begin{center}
		\phantomsection\addcontentsline{toc}{section}{\texorpdfstring{\xchaptertitle~\problemname}{\problemname}}
		\markright{#1}
		\textcolor{structurecolor}{\Large\bfseries\adftripleflourishleft~#1~\adftripleflourishright}
	\end{center}
	\begin{enumerate}[ref=\thechapter.\theenumi]}{
\end{enumerate}}

\usepackage{color}   

\definecolor{winestain}{rgb}{0.5,0,0}
\definecolor{colorGreenOcre}{RGB}{51,102,0} 
\definecolor{colorBlue2}{RGB}{200,207,248}
\definecolor{mydarkblue}{rgb}{0,0.08,0.45}

\newcommand{\mdframecolor}{gray!10}
\newcommand{\mdframehideline}{true}

\definecolor{mylightbluetitle}{RGB}{60,113,183}
\definecolor{mylightbluetext}{rgb}{0,0.08,0.45}
\definecolor{structurecolorblue}{RGB}{60,113,183}
\definecolor{structurecolorgreen}{RGB}{63,145,182}

\colorlet{structurecolor}{structurecolorblue}
\definecolor{structurecolorelegant}{RGB}{60,113,183}
\definecolor{structurecolorlt}{RGB}{31,119,185}

\definecolor{structurecolorHighTheoremBlue}{RGB}{220,227,248}
\definecolor{structurecolorHighTheoremGreen}{RGB}{188,222,231}
\colorlet{structurecolorHighTheorem}{structurecolorHighTheoremBlue}
\definecolor{mdframecolorRemark}{RGB}{186,94,103}

\definecolor{mydarkblue}{rgb}{0,0.08,0.45}
\definecolor{mydarkred}{rgb}{0.70,0.00,0.00}
\definecolor{mydarkgreen}{rgb}{0.00,0.30,0.00}
\definecolor{mydarkyellow}{RGB}{197,151,13}
\definecolor{mydarkpurple}{RGB}{90,35,140}
\definecolor{mydarkgray}{RGB}{64,64,64}

\definecolor{color0}  {RGB}{174,225,254} 
\definecolor{color1}  {RGB}{220,227,248} 
\definecolor{color2}  {RGB}{28,130,185} 
\definecolor{color3}  {RGB}{255,253,250} 
\definecolor{colormiddleright}  {RGB}{245,253,250} 
\definecolor{colorbottomleft}  {RGB}{255,243,250} 
\definecolor{coloruppermiddle}  {RGB}{255,253,230} 
\definecolor{colormiddleleft}  {RGB}{255,244,237}
\definecolor{colorcr}  {RGB}{249,253,232} 
\definecolor{colorreduction}  {RGB}{255,235,254} 
\definecolor{colorqr}  {RGB}{254,221,199} 
\definecolor{colorbiconjugate}  {RGB}{251,149,161} 
\definecolor{colorsvd}  {RGB}{215,247,235} 
\definecolor{colorupperright}  {RGB}{239,246,251} 
\definecolor{colorspectral}  {RGB}{206,226,243} 
\definecolor{colorbottomright}  {RGB}{220,224,236} 
\definecolor{coloreigenvalue}  {RGB}{197,203,224} 
\definecolor{colorcp} {RGB}{217, 234, 186} 
\definecolor{colorcpborder} {RGB}{233, 243, 216} 
\definecolor{colorupperleft}  {RGB}{235,243,240} 
\definecolor{colorsemidefinite}  {RGB}{217,232,226} 
\definecolor{colormiddle} {RGB}{235, 240,255}
\definecolor{colorlu}  {RGB}{220,227,255} 
\definecolor{colorals}  {RGB}{240,230,255} 
\definecolor{coloralsbkg}  {RGB}{248,243,255} 
\definecolor{canaryyellow}{rgb}{1.0, 0.75, 0.0}
\definecolor{bluepigment}{rgb}{0.0, 0.0, 1.0}
\definecolor{canarypurple}{RGB}{208, 13, 241}
\definecolor{colorGreenOcre}{RGB}{51,102,0} 
\definecolor{colorBlue1}  {RGB}{220,227,248}
\definecolor{colorBlue2}{RGB}{200,207,248}
\definecolor{shadecolor}{gray}{0.75}

\definecolor{color0}  {RGB}{174,225,254} 
\definecolor{color1}  {RGB}{220,227,248} 
\definecolor{color2}  {RGB}{28,130,185} 
\definecolor{color3}  {RGB}{255,253,250} 
\definecolor{color0}  {RGB}{174,225,254} 
\definecolor{color1}  {RGB}{220,227,248} 
\definecolor{color2}  {RGB}{28,130,185} 
\definecolor{color3}  {RGB}{255,253,250} 

\definecolor{colormiddleright}  {RGB}{245,253,250} 
\definecolor{colorbottomleft}  {RGB}{255,243,250} 
\definecolor{coloruppermiddle}  {RGB}{255,253,230} 
\definecolor{colormiddleleft}  {RGB}{255,244,237}

\definecolor{colorcr}  {RGB}{249,253,232} 
\definecolor{colorreduction}  {RGB}{255,235,254} 
\definecolor{colorqr}  {RGB}{254,221,199} 
\definecolor{colorbiconjugate}  {RGB}{251,149,161} 

\definecolor{colorsvd}  {RGB}{215,247,235} 

\definecolor{colorupperright}  {RGB}{239,246,251} 
\definecolor{colorspectral}  {RGB}{206,226,243} 

\definecolor{colorbottomright}  {RGB}{220,224,236} 
\definecolor{coloreigenvalue}  {RGB}{197,203,224} 

\definecolor{colorcp} {RGB}{217, 234, 186} 
\definecolor{colorcpborder} {RGB}{233, 243, 216} 

\definecolor{colorupperleft}  {RGB}{235,243,240} 
\definecolor{colorsemidefinite}  {RGB}{217,232,226} 

\definecolor{colormiddle} {RGB}{235, 240,255}
\definecolor{colorlu}  {RGB}{220,227,255} 

\definecolor{colorals}  {RGB}{240,230,255} 
\definecolor{coloralsbkg}  {RGB}{248,243,255} 

\definecolor{canaryyellow}{rgb}{1.0, 0.75, 0.0}
\definecolor{bluepigment}{rgb}{0.0, 0.0, 1.0}

\definecolor{canarypurple}{RGB}{208, 13, 241}
\definecolor{brightlavender}{rgb}{0.44, 0.16, 0.39}

\newcommand{\uniformdist}{\text{Uniform}}
\newcommand{\exponential}{\mathcal{E}}

\newcommand{\gammadist}{\mathcal{G}}

\newcommand{\subnormal}{\mathcal{SG}}

\newcommand{\strsubnormal}{\mathcal{SSG}}
\newcommand{\normal}{\mathcal{N}}

\newcommand{\chisquared}{\chi^2}

\newcommand{\real}{\mathbb{R}}
\newcommand{\prob}{\Pr}
\newcommand{\leadto}{\qquad\underrightarrow{ \text{leads to} }\qquad}

\mathchardef\mhyphen="2D
\newcommand{\integer}{\mathbb{Z}}
\newcommand{\gap}{\,\,\,\,\,\,\,\,}

\newcommand{\diag}{\mathrm{diag}}

\newcommand{\indicator}{\mathds{1}}

\newcommand{\bernoullidist}{\mathrm{Bern}}

\newcommand{\smu}{\mu}
\newcommand{\ssigma}{\sigma}

\newcommand{\tr}{\mathrm{tr}}

\usepackage{amsmath}

\newcommand{\exampbar}{\hfill $\square$\par}

\newcommand{\cspace}{\mathcal{C}}
\newcommand{\nspace}{\mathcal{N}}

\newcommand{\Corr}{\mathbb{C}\mathrm{orr}}
\newcommand{\Cov}{\mathbb{C}\mathrm{ov}}

\newcommand{\Exp}{\mathbb{E}}
\newcommand{\Mom}{\mathbb{M}}
\newcommand{\Cum}{\mathbb{K}}
\newcommand{\Var}{\mathbb{V}\mathrm{ar}}

\newcommand{\argmax}{\operatorname*{\text{arg max}}}
\newcommand{\argmin}{\operatorname*{\text{arg min}}} 

\newcommand\abs[1]{\left\lvert#1\right\rvert}
\newcommand\absbig[1]{\big\lvert#1\big\rvert}
\newcommand\absBig[1]{\Big\lvert#1\Big\rvert}
\newcommand\norm[1]{\left\lVert#1\right\rVert}

\newcommand\normbig[1]{\big\lVert#1\big\rVert}
\newcommand\normzero[1]{\left\lVert#1\right\rVert_0}
\newcommand\normzerobig[1]{\big\lVert#1\big\rVert_0}
\newcommand\normone[1]{\left\lVert#1\right\rVert_1}
\newcommand\normonebig[1]{\big\lVert#1\big\rVert_1}

\newcommand\normtwo[1]{\left\lVert#1\right\rVert_2}
\newcommand\normtwobig[1]{\big\lVert#1\big\rVert_2}

\newcommand\norms[1]{\left\lVert#1\right\rVert_s}
\newcommand\normt[1]{\left\lVert#1\right\rVert_t}
\newcommand\norma[1]{\left\lVert#1\right\rVert_a}
\newcommand\normb[1]{\left\lVert#1\right\rVert_b}
\newcommand\normsbig[1]{\big\lVert#1\big\rVert_s}

\newcommand\normf[1]{\left\lVert#1\right\rVert_F}
\newcommand\normfbig[1]{\big\lVert#1\big\rVert_F}
\newcommand\norminf[1]{\left\lVert#1\right\rVert_{\infty}}
\newcommand\norminfbig[1]{\big\lVert#1\big\rVert_{\infty}}

\newcommand\normA[1]{\left\lVert#1\right\rVert_{\bA}}

\newcommand\normC[1]{\left\lVert#1\right\rVert_{\bC}}

\newcommand\normH[1]{\left\lVert#1\right\rVert_{\bH}}

\newcommand\normAbig[1]{\big\lVert#1\big\rVert_{\bA}}

\newcommand\normCbig[1]{\big\lVert#1\big\rVert_{\bC}}

\newcommand\innerproduct[1]{\left\langle#1\right\rangle}
\newcommand\innerproductbig[1]{\big\langle#1\big\rangle}

\newcommand{\sgn}{\mathrm{sign}}
\DeclareMathOperator{\sign}{sign}

\newcommand{\rank}{\mathrm{rank}}

\newcommand{\trace}{\mathrm{tr}}

\newcommand{\spn}{\mathrm{span}}

\mathchardef\mhyphen="2D
\newcommand{\complex}{\mathbb{C}}
\newcommand{\naturalset}{\mathbb{N}}
\newcommand{\integerset}{\mathbb{Z}}

\def\1{\bm{1}}

\newcommand{\R}{\mathbb{R}}

\newcommand{\entropy}{\mathrm{H}}
\newcommand{\KL}{D_{\mathrm{KL}}}

\let\ab\allowbreak

\newcommand{\topone}{{(1)}}
\newcommand{\toptwo}{{(2)}}
\newcommand{\topzero}{{(0)}}
\newcommand{\toptminus}{{(t-1)}}
\newcommand{\toptminusTWO}{{(t-2)}}
\newcommand{\toptminusTOP}{{(t-1)\top}}
\newcommand{\toptzero}{{(t)}}
\newcommand{\toptzeroTOP}{{(t)\top}}

\newcommand{\toptone}{{(t+1)}}

\newcommand{\toplzero}{{(\ell)}}

\newcommand{\bzero}{\boldsymbol{0}}

\newcommand{\balpha}{{\boldsymbol\alpha}}
\newcommand{\bbeta}{{\boldsymbol\beta}}

\newcommand{\bdelta}{{\boldsymbol\delta}}

\newcommand{\bepsilon}{{\boldsymbol\epsilon}}

\newcommand{\bgamma}{{\boldsymbol\gamma}}

\newcommand{\blambda}{{\boldsymbol\lambda}}
\newcommand{\bmu}{{\boldsymbol\mu}}
\newcommand{\bnu}{{\boldsymbol\nu}}

\newcommand{\bsigma}{{\boldsymbol\sigma}}

\newcommand{\btheta}{{\boldsymbol\theta}}

\newcommand{\bxi}{{\boldsymbol\xi}}
\newcommand{\bzeta}{{\boldsymbol\zeta}}

\newcommand{\bLambda}{{\boldsymbol\Lambda}}
\newcommand{\bOmega}{{\boldsymbol\Omega}}
\newcommand{\bPhi}{{\boldsymbol\Phi}}

\newcommand{\bPsi}{{\boldsymbol\Psi}}
\newcommand{\bSigma}{{\boldsymbol\Sigma}}

\newcommand{\widehatbeta}{{\widehat\beta}}

\newcommand{\widehatgamma}{{\widehat\gamma}}

\newcommand{\widehatlambda}{{\widehat\lambda}}

\newcommand{\widehatnu}{{\widehat\nu}}

\newcommand{\widetildebeta}{{\widetilde\beta}}

\newcommand{\widetildeeta}{{\widetilde\eta}}
\newcommand{\widetildegamma}{{\widetilde\gamma}}

\newcommand{\widehatbalpha}{{\widehat\balpha}}
\newcommand{\widehatbbeta}{{\widehat\bbeta}}

\newcommand{\widehatbgamma}{{\widehat\bgamma}}

\newcommand{\widehatblambda}{{\widehat\blambda}}
\newcommand{\widehatbmu}{{\widehat\bmu}}
\newcommand{\widehatbnu}{{\widehat\bnu}}

\newcommand{\widehatbtheta}{{\widehat\btheta}}

\newcommand{\widetildebalpha}{{\widetilde\balpha}}
\newcommand{\widetildebbeta}{{\widetilde\bbeta}}

\newcommand{\widetildeblambda}{{\widetilde\blambda}}

\newcommand{\widetildebnu}{{\widetilde\bnu}}

\newcommand{\widetildebtheta}{{\widetilde\btheta}}

\newcommand{\widebarbbeta}{{\overline\bbeta}}

\newcommand{\widebarblambda}{{\overline\blambda}}

\newcommand{\widebarbQ}{\overline{\bm{Q}}}

\newcommand{\widebarbg}{\overline{\bm{g}}}

\newcommand{\widebarby}{\overline{\bm{y}}}

\newcommand{\widehatC}{\widehat{C}}

\newcommand{\widehatT}{\widehat{T}}

\newcommand{\widehatc}{\widehat{c}}
\newcommand{\widehatd}{\widehat{d}}

\newcommand{\widehatf}{\widehat{f}}

\newcommand{\widehatw}{\widehat{w}}

\newcommand{\widehatbX}{\widehat{\bm{X}}}

\newcommand{\widehatbc}{\widehat{\bm{c}}}
\newcommand{\widehatbd}{\widehat{\bm{d}}}
\newcommand{\widehatbe}{\widehat{\bm{e}}}

\newcommand{\widehatbq}{\widehat{\bm{q}}}
\newcommand{\widehatbr}{\widehat{\bm{r}}}

\newcommand{\widehatbu}{\widehat{\bm{u}}}

\newcommand{\widehatbw}{\widehat{\bm{w}}}

\newcommand{\widehatby}{\widehat{\bm{y}}}

\newcommand{\mathcalA}{\mathcal{A}}
\newcommand{\mathcalB}{\mathcal{B}}
\newcommand{\mathcalC}{\mathcal{C}}
\newcommand{\mathcalD}{\mathcal{D}}
\newcommand{\mathcalE}{\mathcal{E}}
\newcommand{\mathcalF}{\mathcal{F}}
\newcommand{\mathcalG}{\mathcal{G}}
\newcommand{\mathcalH}{\mathcal{H}}

\newcommand{\mathcalL}{\mathcal{L}}

\newcommand{\mathcalO}{\mathcal{O}}
\newcommand{\mathcalP}{\mathcal{P}}

\newcommand{\mathcalT}{\mathcal{T}}
\newcommand{\mathcalU}{\mathcal{U}}
\newcommand{\mathcalV}{\mathcal{V}}

\newcommand{\mathcalX}{\mathcal{X}}

\newcommand{\widetildebA}{\widetilde{\bm{A}}}
\newcommand{\widetildebB}{\widetilde{\bm{B}}}
\newcommand{\widetildebC}{\widetilde{\bm{C}}}
\newcommand{\widetildebD}{\widetilde{\bm{D}}}

\newcommand{\widetildebS}{\widetilde{\bm{S}}}

\newcommand{\widetildebU}{\widetilde{\bm{U}}}

\newcommand{\widetildebX}{\widetilde{\bm{X}}}

\newcommand{\widetildebZ}{\widetilde{\bm{Z}}}
\newcommand{\widetildeba}{\widetilde{\bm{a}}}

\newcommand{\widetildebg}{\widetilde{\bm{g}}}

\newcommand{\widetildebn}{\widetilde{\bm{n}}}

\newcommand{\widetildebu}{\widetilde{\bm{u}}}

\newcommand{\widetildeby}{\widetilde{\bm{y}}}

\newcommand{\widetildeC}{\widetilde{C}}

\newcommand{\widetildeF}{\widetilde{F}}

\newcommand{\widetildeT}{\widetilde{T}}

\newcommand{\widetildek}{\widetilde{k}}

\newcommand{\widetildep}{\widetilde{p}}

\newcommand{\bone}{{\bm{1}}}
\newcommand{\ba}{{\bm{a}}}
\newcommand{\bA}{{\bm{A}}}
\newcommand{\bb}{{\bm{b}}}
\newcommand{\bB}{{\bm{B}}}
\newcommand{\bc}{{\bm{c}}}
\newcommand{\bC}{{\bm{C}}}
\newcommand{\bd}{{\bm{d}}}
\newcommand{\bD}{{\bm{D}}}
\newcommand{\be}{{\bm{e}}}
\newcommand{\bE}{{\bm{E}}}
\newcommand{\bff}{{\bm{f}}}
\newcommand{\bF}{{\bm{F}}}
\newcommand{\bg}{{\bm{g}}}
\newcommand{\bG}{{\bm{G}}}
\newcommand{\bh}{{\bm{h}}}
\newcommand{\bH}{{\bm{H}}}

\newcommand{\bI}{{\bm{I}}}

\newcommand{\bJ}{{\bm{J}}}

\newcommand{\bL}{{\bm{L}}}

\newcommand{\bM}{{\bm{M}}}
\newcommand{\bn}{{\bm{n}}}

\newcommand{\bP}{{\bm{P}}}
\newcommand{\bq}{{\bm{q}}}
\newcommand{\bQ}{{\bm{Q}}}
\newcommand{\br}{{\bm{r}}}

\newcommand{\bs}{{\bm{s}}}
\newcommand{\bS}{{\bm{S}}}
\newcommand{\bt}{{\bm{t}}}

\newcommand{\bu}{{\bm{u}}}
\newcommand{\bU}{{\bm{U}}}
\newcommand{\bv}{{\bm{v}}}
\newcommand{\bV}{{\bm{V}}}
\newcommand{\bw}{{\bm{w}}}
\newcommand{\bW}{{\bm{W}}}
\newcommand{\bx}{{\bm{x}}}
\newcommand{\bX}{{\bm{X}}}
\newcommand{\by}{{\bm{y}}}
\newcommand{\bY}{{\bm{Y}}}
\newcommand{\bz}{{\bm{z}}}
\newcommand{\bZ}{{\bm{Z}}}

\def\vmu{{\bm{\mu}}}
\def\vtheta{{\bm{\theta}}}
\def\va{{\bm{a}}}

\def\ve{{\bm{e}}}

\def\vx{{\bm{x}}}

\def\mH{{\bm{H}}}

\def\mJ{{\bm{J}}}

\def\mX{{\bm{X}}}

\def\mSigma{{\bm{\Sigma}}}

\DeclareMathAlphabet{\mathsfit}{\encodingdefault}{\sfdefault}{m}{sl}
\SetMathAlphabet{\mathsfit}{bold}{\encodingdefault}{\sfdefault}{bx}{n}

\def\sA{{\mathbb{A}}}
\def\sB{{\mathbb{B}}}
\def\sC{{\mathbb{C}}}

\def\sE{{\mathbb{E}}}
\def\sF{{\mathbb{F}}}
\def\sG{{\mathbb{G}}}
\def\sH{{\mathbb{H}}}
\def\sI{{\mathbb{I}}}
\def\sJ{{\mathbb{J}}}

\def\sL{{\mathbb{L}}}
\def\sM{{\mathbb{M}}}

\def\sP{{\mathbb{P}}}
\def\sQ{{\mathbb{Q}}}

\def\sS{{\mathbb{S}}}
\def\sT{{\mathbb{T}}}
\def\sU{{\mathbb{U}}}

\def\sX{{\mathbb{X}}}

\def\ra{{\textnormal{a}}}
\def\rb{{\textnormal{b}}}
\def\rc{{\textnormal{c}}}

\def\rs{{\textnormal{s}}}

\def\ru{{\textnormal{u}}}

\def\rx{{\textnormal{x}}}
\def\ry{{\textnormal{y}}}
\def\rz{{\textnormal{z}}}

\def\rva{{\mathbf{a}}}
\def\rvb{{\mathbf{b}}}

\def\rvu{{\mathbf{i}}}

\def\rvu{{\mathbf{u}}}
\def\rvv{{\mathbf{v}}}

\def\rvx{{\mathbf{x}}}
\def\rvy{{\mathbf{y}}}
\def\rvz{{\mathbf{z}}}

\def\rmA{{\mathbf{A}}}
\def\rmB{{\mathbf{B}}}

\def\rmD{{\mathbf{D}}}

\def\rmH{{\mathbf{H}}}

\def\rmU{{\mathbf{U}}}
\def\rmV{{\mathbf{V}}}

\def\rmX{{\mathbf{X}}}
\def\rmY{{\mathbf{Y}}}
\def\rmZ{{\mathbf{Z}}}

\firstpageno{1}

\newcommand{\mytitle}{A First Course in Sparse Optimization}
\begin{document}

\newpage
\thispagestyle{empty}  
\title{\mytitle}

\author{
\begin{center}
\name Jun Lu \\ 
\email jun.lu.locky@gmail.com
\end{center}
}

\frontmatter

\newpage 
\maketitle

\chapter*{\centering \begin{normalsize}Preface\end{normalsize}}

In an age defined by data---its abundance, its velocity, and its complexity---we face a paradox: the more information we gather, the harder it becomes to extract meaning from it. High-dimensional datasets now permeate every corner of science, engineering, and industry, from medical imaging and genomics to wireless communications and machine learning. Yet, beneath this deluge of numbers lies a powerful insight: much of the data we collect is not as complex as it first appears. In fact, many real-world signals are sparse---they can be accurately described using only a few essential components when viewed through the right lens.

This observation has given rise to a transformative field at the intersection of applied mathematics, signal processing, statistics, and optimization: sparse optimization. At its core, sparse optimization seeks to recover or construct solutions that are as simple as possible---solutions in which most entries are zero or negligible---while still faithfully representing the underlying phenomenon. Whether the goal is to reconstruct an image from a handful of measurements, identify the key genetic markers linked to a disease, or build a predictive model with interpretable features, sparsity offers both computational efficiency and conceptual clarity.

Two major paradigms have emerged within this field. The first, rooted in compressed sensing, reimagines how we acquire data altogether. By recognizing that many signals are sparse in some domain, compressed sensing shows that we can bypass traditional sampling limits and reconstruct signals from far fewer measurements than previously thought possible. The second paradigm, exemplified by the LASSO and related regularization techniques, embeds sparsity directly into the modeling process---encouraging models that are not only accurate but also parsimonious, selecting only the most relevant variables from a sea of candidates.

Though distinct in motivation and application, these approaches share deep mathematical foundations and algorithmic strategies. This book is an invitation to explore those connections. It is designed to serve both as an accessible introduction for students and researchers entering the field, and as a rigorous reference for practitioners seeking to apply sparse optimization methods in their work. We begin with the essential mathematical and algorithmic building blocks, then systematically develop the theory and practice of sparse signal recovery and sparse regularization. Along the way, we emphasize intuition without sacrificing precision, and theory without neglecting application.

The book will take us through the geometry of high-dimensional spaces, the design of measurement systems that preserve sparsity, and the algorithms that make large-scale sparse recovery feasible. We will see how randomness can be harnessed to guarantee recovery, how convex relaxation turns intractable problems into solvable ones, and how modern extensions---such as structured sparsity and non-convex penalties---push the boundaries of what's possible.

\paragraph{Keywords.}
Sparse optimization, LASSO, Generalized LASSO, Group LASSO,  Sparse recovery, Convex optimization, Convex relaxation, Design matrices, Nullspace property, Restricted eigenvalue property, 
Restricted isometry property, Exact recovery, Gaussian random matrices, Sub-Gaussian random matrices, Primal-dual algorithm, 
Gradient descent.

\newpage
\begingroup
\hypersetup{
linkcolor=structurecolor,
linktoc=page,  
}
\dominitoc
\pdfbookmark{\contentsname}{toc} 
\tableofcontents 

\endgroup

%

\chapter*{Notation}\label{notation}


This section provides a concise reference describing notation used throughout this
book.
If you are unfamiliar with any of the corresponding mathematical concepts,
the book describes most of these ideas in Chapter~\ref{chapter_lsintroduction} (p.~\pageref{chapter_lsintroduction}).

\vspace{0.4in}
\begin{minipage}{\textwidth}
\centerline{\bf Numbers and Arrays}
\bgroup
\def\arraystretch{1.5}
\begin{tabular}{cp{4.25in}}
$\displaystyle a$   & A scalar (integer or real, italics font)\\
$\displaystyle \ba$ & A vector (italics font)\\
$\displaystyle \bA$ & A matrix (italics font)\\
$\displaystyle \bI_n$ & Identity matrix with $n$ rows and $n$ columns\\
$\displaystyle \bI$   & Identity matrix with dimensionality implied by context\\
$\displaystyle \ve_i$ & Standard basis vector $[0,\dots,0,1,0,\dots,0]$ with a 1 at position $i$\\
$\displaystyle \text{diag}(\va)$ & A square, diagonal matrix with diagonal entries given by $\va$\\
$\displaystyle \ra$   & A scalar random variable  (normal fonts)\\
$\displaystyle \rva$  & A vector-valued random variable  (normal fonts)\\
$\displaystyle \rmA$  & A matrix-valued random variable (normal fonts)\\
\end{tabular}
\egroup
\index{Scalar}
\index{Vector}
\index{Matrix}
\end{minipage}

\index{Sets}
\vspace{0.2in}
\begin{minipage}{\textwidth}
\centerline{\bf Sets}
\bgroup
\def\arraystretch{1.5}
\begin{tabular}{cp{4.25in}}
$\displaystyle \sA$ & A set\\
$\displaystyle \varnothing$ & The null set \\
$\displaystyle \real, \complex, \sF\equiv \{\real \text{ or }\complex\}$ & The set of real, complex, either real or complex numbers\\
$\displaystyle \naturalset, \naturalset_0$ & The set of natural numbers $\naturalset=\{1,2,\ldots\}$, $\naturalset_0=\{0\}\cup\naturalset$  \\
$\displaystyle \{0, 1\}$ & The set containing 0 and 1 \\
$\displaystyle \{0, 1, \dots, n \}$ & The set of all integers between $0$ and $n$\\
$\displaystyle [a, b]$ & The real interval including $a$ and $b$\\
$\displaystyle (a, b]$ & The real interval excluding $a$ but including $b$\\
$\displaystyle \sA \backslash \sB$ & Set subtraction, i.e., the set containing the elements of $\sA$ that are not in $\sB$\\
$\displaystyle \sB_s(\bc, r)$  & Open balls: $\sB_s(\bc, r) = \{\bx\mid\norms{\bc-\bx}< r \}$\\
$\displaystyle \sB_s[\bc, r]$  & Closed balls: $\sB_s[\bc, r] = \{\bx\mid\norms{\bc-\bx}\leq r \}$\\
$\displaystyle \sB(\bc, r)$  & $\ell_2$ open balls: $\sB(\bc, r) = \{\bx\mid\normtwo{\bc-\bx}< r \}$\\
$\displaystyle \sB[\bc, r]$  & $\ell_2$ closed balls: $\sB[\bc, r] = \{\bx\mid\normtwo{\bc-\bx}\leq r \}$\\
\end{tabular}
\egroup
\index{Scalar}
\index{Vector}
\index{Matrix}
\index{Set}
\end{minipage}

\vspace{0.2in}
\begin{minipage}{\textwidth}
\centerline{\bf Indexing}
\bgroup
\def\arraystretch{1.5}
\begin{tabular}{cp{4.25in}}
$\displaystyle x_i$ & Element $i$ of vector $\bx$, with indexing starting at 1 \\
$\displaystyle \bx_{-i}$ & All elements of vector $\bx$ except for element $i$ \\
$\displaystyle  x_{ij}$ & Element $i, j$ of matrix $\bX$ \\
$\displaystyle \bX_{i, :}=\bX[i,:],\, \bx^{(i)}$ & Row $i$ of matrix $\bX$ \\
$\displaystyle \bX_{:, i}=\bX[:, i],\, \bx_i$ & Column $i$ of matrix $\bX$ \\
$\displaystyle \bX_\sS=\bX[:, \sS], \bX(\sS)$ & $\bX_\sS \in\real^{n\times \abs{\sS}}$ and $\bX(\sS)\in\real^{n\times p}$ if $\bX\in\real^{n\times p}$ \\
\end{tabular}
\egroup
\end{minipage}

\vspace{0.2in}
\begin{minipage}{\textwidth}
\centerline{\bf Linear Algebra Operations}
\bgroup
\def\arraystretch{1.5}
\begin{tabular}{cp{4.25in}}
$\displaystyle \bX^\top$ & Transpose of matrix $\bX$ \\
$\displaystyle \bX^+$ & Moore-Penrose pseudo-inverse of $\bX$\\
$\displaystyle \bX \hadaprod \bY $ & Element-wise (Hadamard) product of $\bX$ and $\bY$ \\
$\displaystyle \mathrm{det}(\bX)$ & Determinant of $\bX$ \\
$\displaystyle \cspace(\bX), \nspace(\bX), \mathcalV$ & Column, null  space of $\bX$, and a general space \\
$\displaystyle \rank(\bX)$ & Rank of $\bX$ \\
$\displaystyle \trace(\bX)$ & Trace of $\bX$ \\
$\displaystyle \lambda_{\max}(\bX), \lambda_{\min}(\bX)$ & Largest and smallest eigenvalue of $\bX$ \\
$\displaystyle \sigma_{\max}(\bX), \sigma_{\min}(\bX)$ & Largest and smallest singular value of $\bX$ \\
$\displaystyle  \diag(\bbeta)$ & Diagonal matrix with entries of $\bbeta$
\end{tabular}
\egroup
\index{Transpose}
\index{Element-wise product, Hadamard product}
\index{Hadamard product}
\end{minipage}

\vspace{0.4in}
\begin{minipage}{\textwidth}
\centerline{\bf Calculus}
\bgroup
\def\arraystretch{1.5}
\begin{tabular}{cp{4.25in}}
$\displaystyle\frac{d y} {d x}$ & Derivative of $y$ with respect to $x$\\ [2ex]
$\displaystyle \frac{\partial y} {\partial x} $ & Partial derivative of $y$ with respect to $x$ \\
$\displaystyle \nabla_{\bx} y $ & Gradient of $y$ with respect to $\bx$ \\
$\displaystyle \nabla_{\bX} y $ & Matrix derivatives of $y$ with respect to $\bX$ \\
$\displaystyle \frac{\partial f}{\partial \vx} $ & Jacobian matrix $\mJ \in \R^{m\times n}$ of $f: \R^n \rightarrow \R^m$\\
$\displaystyle \nabla_\vx^2 f(\vx)\text{ or }\mH( f)(\vx)$ & The Hessian matrix of $f$ at input point $\vx$\\
$\displaystyle \int f(\vx) d\vx $ & Definite integral over the entire domain of $\vx$ \\
$\displaystyle \int_\sS f(\vx) d\vx$ & Definite integral with respect to $\vx$ over the set $\sS$ \\
\end{tabular}
\egroup
\index{Integral}
\index{Hessian matrix}
\end{minipage}

\vspace{0.4in}
\begin{minipage}{\textwidth}
\centerline{\bf Probability and Information Theory}
\bgroup
\def\arraystretch{1.5}
\begin{tabular}{cp{4.25in}}
$\displaystyle \ra \bot \rb$ & The random variables $\ra$ and $\rb$ are independent\\
$\displaystyle \ra \bot \rb \mid \rc $ & They are conditionally independent given $\rc$\\
$\displaystyle \Pr(\bx)$ & A probability distribution over a discrete variable\\
$\displaystyle p(\bx), p_{\rvx}(\bx), f(\bx), f_{\rvx}(\bx)$ & A probability distribution over a continuous variable, or over
a variable whose type has not been specified\\
$\displaystyle \ra \sim P$ & Random variable $\ra$ has distribution $P$\\
$\displaystyle  \Exp_{\rx\sim P} [ f(x) ]\text{ or } \Exp [f(x)]$ & Expectation of $f(x)$ with respect to $P(\rx)$ \\
$\displaystyle \Var[f(x)] $ &  Variance of $f(x)$ under $P(\rx)$ \\
$\displaystyle \Cov[f(x),g(x)] $ & Covariance of $f(x)$ and $g(x)$ under $P(\rx)$\\
$\displaystyle \Corr[f(x),g(x)] $ & Correlation of $f(x)$ and $g(x)$ under $P(\rx)$\\
$\displaystyle \normal ( \mu , \sigma^2)$ &  Gaussian distribution %
 with mean $\mu$ and variance $\sigma^2$ \\
$\displaystyle \normal ( \vmu , \mSigma)$ & Multivariate Gaussian distribution %
 with mean $\vmu$ and covariance $\mSigma$ \\
$\displaystyle \bernoullidist(p) $ & Bernoulli distribution with mean $p$ \\
$\displaystyle \exponential(\lambda) $ & Exponential distribution with scale $\lambda$ \\
$\displaystyle \chisquared_{(s)} $ & Chi-squared distribution with  $s$ degrees of freedom (df)\\
$\subnormal(\sigma^2)$ & Sub-Gaussian distribution\\
$\strsubnormal(\sigma^2)$ & Strictly sub-Gaussian distribution
\end{tabular}
\egroup
\index{Independence}
\index{Conditional independence}
\index{Variance}
\index{Covariance}
\index{Kullback-Leibler divergence}
\index{Shannon entropy}
\end{minipage}

\vspace{0.4in}
\begin{minipage}{\textwidth}
\centerline{\bf Functions}
\bgroup
\def\arraystretch{1.5}
\begin{tabular}{cp{4.25in}}
$\displaystyle f: \sA \rightarrow \sB$ & The function $f$ with domain $\sA$ and range $\sB$\\
$\displaystyle f \circ g $ & Composition of the functions $f$ and $g$ \\
  $\displaystyle f(\vx ; \vtheta) $ & A function of $\vx$ parametrized by $\vtheta$.
  (Sometimes we write $f(\vx)$ and omit the argument $\vtheta$ to lighten notation) \\
$\displaystyle \ln(x),\log(x)$ & Natural logarithm of $x$ \\
$\displaystyle \sgn(x)$ & Sign $x$, taking a value among $+1$, $-1$, and 0 \\
$\displaystyle \sigma(x),\, \text{Sigmoid}(x)$ & Logistic sigmoid, i.e., $\displaystyle \frac{1} {1 + \exp\{-x\}}$ \\
$\displaystyle \text{logit}(\pi)$ & Logit function, i.e., $\text{logit}(\pi) = \ln(\pi/(1-\pi))$, where $\pi\in(0,1)$ \\
$\displaystyle \zeta(x)$ & Softplus, $\log(1 + \exp\{x\})$ \\
$\displaystyle \norm{\bx}_p, \norm{\bx}_s $ & $\ell_p$-norm of $\vx$ \\
$\displaystyle \norm{\bx}=\normtwo{\bx} $ & $\ell_2$-norm of $\vx$ \\
$\displaystyle \norm{\bx}=\normone{\bx} $ & $\ell_1$-norm of $\vx$ \\
$\displaystyle \norm{\bx}=\norminf{\bx} $ & $\ell_\infty$ norm of $\vx$ \\
$\displaystyle [x]_+$ & Positive part of $x$, i.e., $\max(0,x)$\\
$\displaystyle u(x)$ & Step function with value 1 when $x\geq0$ and value 0 otherwise\\
$\displaystyle \indicator\{\mathrm{condition}\}$ & is 1 if the condition is true, 0 otherwise\\
$\displaystyle \indicatorS(\balpha)$ & is 0 if $\balpha\in\sS$, $+\infty$ otherwise\\
$\displaystyle \Phi(x), \Phi^{-1}(\pi)$ & Standard Gaussian cdf, and the probit function, where $\pi\in(0,1)$\\
$\displaystyle \text{Negative binomial}(\alpha, x)$ & $\eta = \ln(x/(x+1/\alpha))$\\
$\displaystyle \mathcalT_\lambda(\bbeta)$ & Soft-thresholding operator\\
$\sigma_k(\bbeta)_s$ & $\ell_s$-error of best $k$-sparse approximation: $\min_{\widehat{\bbeta} \in \sB_0[k]} \normsbig{\bbeta - \widehat{\bbeta}}$ 
\end{tabular}
\egroup
\index{Sigmoid function}
\index{Softplus}
\index{Norm}
\end{minipage}
Sometimes we use a function $f$ whose argument is a scalar but apply
it to a vector, matrix: $f(\vx)$, $f(\mX)$.
This denotes the application of $f$ to the
array element-wise. For example, if $\bC = \sigma(\bX)$, then $c_{ij} = \sigma(x_{ij})$
for all valid values of $i$ and  $j$.

\vspace{0.4in}
\begin{minipage}{\textwidth}
\centerline{\bf Other General Notastions}
\bgroup
\def\arraystretch{1.5}
\begin{tabular}{cp{4.25in}}
$\displaystyle \triangleq$ & Equals by definition\\
$\displaystyle :=, \leftarrow $ & Equals by assignment \\
$\displaystyle \equiv $    & Equals by equivalence \\
$\displaystyle \pi $       & A probability value or 3.141592....\\
$\displaystyle e, \exp $   & 2.71828...
\end{tabular}
\egroup
\end{minipage}

\vspace{0.4in}
\begin{minipage}{\textwidth}
\centerline{\bf Abbreviations}
\bgroup
\def\arraystretch{1.5}
\begin{tabular}{cp{4.25in}}
PD & Positive definite  \\
PSD & Positive semidefinite \\
i.i.d. & Independently and identically distributed \\
p.d.f., PDF & Probability density function \\
p.m.f., PMF & Probability mass function \\
MGF & Moment generating function\\
MVT & Mean value theorem\\
LSC & Lower semicontinuity \\
LP & Linear programming\\
LS, OLS & Ordinary least squares\\
GD & Gradient descent\\
SGD & Stochastic gradient descent \\
MSE & Mean squared error\\
CLT & Central limit theorem \\
SVD & Singular value decomposition\\
BP & Basis pursuit\\
BPDN & Basis pursuit denoising\\
LASSO & Least absolute shrinkage and selection operator\\
DPP & Dual polytope projections\\
SAFE & Safe feature elimination\\
RIP & Restricted isometry property\\
RIC & Restricted isometry constant\\
ROP & Restricted orthogonality property\\
ROC & Restricted orthogonality constant\\
NSP & Nullspace property\\
JL lemma &  Johnson--Lindenstrauss lemma\\
SVM & Support vector machine\\
FW & Frank--Wolfe  algorithm\\

\end{tabular}
\egroup
\end{minipage}

\vspace{1.4in}
\begin{minipage}{\textwidth}
\bgroup
\def\arraystretch{1.5}
\begin{tabular}{cp{4.25in}}
CS & Compressed sensing \\
LASSO & Least absolute shrinkage and selection operator\\
LARS & Least angle regression\\
IHT & Iterative hard-thresholding \\
ADMM & Alternating direction method of multipliers \\
BB rule & Barzilai--Borwein rule \\
BCD & Block coordinate descent\\
BCQP & Bound-constrained quadratic program \\
PGD & Projected gradient method\\
PGM & Proximal gradient method\\
SG & Sub-Gaussian \\
SSG & Strictly sub-Gaussian\\
ISTA & Iterative shrinkage-thresholding algorithm \\
FISTA & Fast proximal gradient method\\
MP & Matching pursuit \\
ROMP & Regularized orthogonal matching pursuit \\
OMP & Orthogonal matching pursuit \\
StOMP & Stagewise orthogonal matching pursuit \\
CoSaMP & Compressive sampling matching pursuit\\
SP & Subspace basis pursuit\\
BT & Basic thresholding\\
HTP & Hard thresholding pursuit\\
IRLS & Iteratively reweighted least squares \\
ALF & Augmented Lagrangian function\\
FAL & First-order augmented Lagrangian\\

\end{tabular}
\egroup
\end{minipage}


\clearpage


\mainmatter

\newpage 
\chapter{Introduction}\label{chapter_lsintroduction}
\begingroup
\hypersetup{
linkcolor=structurecolor,
linktoc=page,  
}
\minitoc \newpage
\endgroup
\section{Introduction and Background}
\lettrine{\color{caligraphcolor}I}
In today's data-driven world, we are encountering a profound shift in how we collect, store, and analyze vast amounts of information. 
As sensing technologies and computational methods continue to advance, a key challenge has emerged: how to efficiently manage high-dimensional data that is inherently sparse or compressible. 
This challenge lies at the heart of \textit{sparse optimization}, a field that focuses on both the recovery of sparse signals and the formulation of optimization problems that promote sparsity in solutions.


Thanks to rapid progress in digital sensing technologies, we can now capture data with unprecedented resolution. 
However, these systems often generate enormous volumes of information, posing significant challenges for storage, processing, and transmission. Fortunately, many real-world signals---ranging from images and audio recordings to genomic sequences and medical scans---are \textit{sparse} or \textit{compressible} in some transform domain. 
In other words, although such signals may reside in high-dimensional spaces, they can typically be represented accurately using only a small number of nonzero coefficients. Recognizing and leveraging this inherent sparsity is essential for developing efficient algorithms that reduce dimensionality without sacrificing critical information.

\paragraph{Sparse signal recovery: the role of compressed sensing.}
One of the primary applications of sparse optimization techniques is \textit{sparse signal recovery}. This task involves recovering a signal from a reduced number of measurements, taking advantage of the fact that the signal has a sparse representation in some basis or transform domain.

The conventional  approach to signal acquisition, governed by the \textit{Nyquist--Shannon sampling theorem}, requires sampling a signal at twice the highest frequency present in the signal to ensure perfect reconstruction \citep{nyquist2002certain, shannon2006communication}. 
While this theory works well for band-limited signals, it is highly inefficient when the signal is sparse. \textit{Compressed sensing (CS, or compressive sensing)}, a field that emerged in the early 2000s, fundamentally rethinks this paradigm. 
Instead of first acquiring signals at a high rate and then compressing them (as in JPEG, MP3, or MPEG), CS advocates for acquiring the data in a compressed form, i.e., directly sampling at a lower rate. 
By exploiting sparsity, it becomes possible to reconstruct a signal from a dramatically reduced number of nonadaptive, linear measurements---far fewer than its ambient dimension would suggest. This is not mere data compression; it is compression at the point of acquisition.
This reduces both the number of measurements needed and the computational cost of storing and processing the data \citep{candes2006near, candes2011compressed, davenport2012introduction}.

At the heart of CS  lies the insight  that if a signal is sparse in a known basis (e.g., wavelets, Fourier, or learned dictionaries), we can reconstruct it from far fewer measurements than traditional methods would suggest. The recovery process in CS is typically framed as an \textit{optimization problem} that involves finding the sparsest signal consistent with the available measurements. Formally, given a measurement vector $\by$, the goal is to recover the signal $ \bbeta $ that solves
$$
\min_{\bbeta} \normzero{\bbeta} \quad \text{s.t.} \quad \by = \bX\bbeta + \bepsilon,
$$
where $ \normzero{\bbeta}  $ denotes the \textit{$ \ell_0 $-norm} (i.e.,  the number of nonzero coefficients in $ \bbeta $), $ \bX $ is the measurement matrix, and $ \bepsilon $ represents noise.

In practice, minimizing the $ \ell_0 $-norm is computationally intractable because the problem is non-convex and combinatorial. 
To overcome this, CS commonly replaces the  $ \ell_0 $-norm with the convex \textit{$ \ell_1 $-norm},
leading to the well-known \textit{basis pursuit} formulation:
$$
\min_{\bbeta} \normone{\bbeta} \quad \text{s.t.} \quad \by  = \bX\bbeta + \bepsilon.
$$
This relaxation enables efficient and robust recovery of sparse signals using convex optimization algorithms---a cornerstone of modern compressed sensing theory and practice.

\paragraph{Sparse optimization with LASSO: minimizing loss with sparsity constraints.}
While sparse signal recovery focuses on reconstructing signals from compressed measurements, \textit{LASSO (least absolute shrinkage and selection operator)}  addresses a related but distinct problem in optimization: it seeks to minimize a loss function while enforcing \textit{sparsity constraints} or \textit{sparse regularization} on the solution \citep{tibshirani1996regression}. 
LASSO is widely used in machine learning and statistics for tasks such as regression, where the goal is to identify a small subset of features that best explain the observed data.

In a standard linear regression setting, the objective is to minimize the squared error between the observed response vector  $\by$ and the model's predicted values $ \bX\bbeta $. LASSO modifies this objective by adding an $ \ell_1 $-norm regularization term, which encourages sparsity in the estimated coefficients. 
The LASSO optimization problem can be written as:
$$
\min_{\bbeta}  \frac{1}{2} \normtwo{\bX\bbeta - \by}^2 + \lambda \normone{\bbeta}.
$$
Here, $\frac{1}{2}\normtwo{\bX\bbeta - \by}^2 $ represents the \textit{least squares loss}, and $ \normone{\bbeta} $ is the \textit{$\ell_1$ penalty} or \textit{$\ell_1$ regularization}, which promotes sparsity by penalizing the sum of the absolute values of the coefficients. 
The hyperparameter $ \lambda $ controls the trade-off between fitting the data and enforcing sparsity. A higher value of $ \lambda $ leads to more coefficients being shrunk to zero, promoting greater sparsity.

LASSO is a versatile and computationally tractable method that has found broad application in areas such as \textit{regression}, \textit{classification}, and \textit{signal denoising}, where sparsity is a desirable property for interpretability, generalization, or efficiency. 
While LASSO is not specifically focused on signal recovery from limited measurements, it provides a powerful and general framework for incorporating sparsity directly into optimization-based modeling.

\paragraph{Sparse optimization.}
The key distinction  between \textit{sparse signal recovery} (as in compressed sensing) and \textit{sparse regularization} or \textit{sparse constraint} (as in LASSO) lies in their objectives and the nature of the underlying optimization problems:
\begin{itemize}
\item \textit{Sparse signal recovery (CS).} The objective is to recover a sparse (or a approximately sparse) signal from a limited number of measurements. The problem is typically posed as finding the sparsest solution consistent with the available data. 
CS often leverages optimization techniques like \textit{basis pursuit} (using the $ \ell_1 $-norm) to recover signals from compressed data.
\item \textit{Sparse regularization (LASSO).} In sparse optimization problems like LASSO, the goal is to minimize a data-fitting loss function (e.g., least squares) while encouraging sparsity in the model parameters. 
This is achieved by adding an $\ell_1$-norm penalty to the objective, which shrinks less important coefficients toward zero and effectively performs feature selection.
\end{itemize}
Although closely related, these two paradigms serve different purposes and require distinct algorithmic approaches. Compressed sensing is primarily concerned with signal acquisition and reconstruction from limited measurements, whereas LASSO focuses on building predictive models with interpretable, sparse parameter sets---typically in statistical learning and regression contexts.

The field of sparse optimization has profoundly influenced numerous domains, especially those involving high-dimensional data. In medical imaging, compressed sensing enables faster MRI scans with reduced radiation exposure; in communications, it facilitates more efficient data transmission \citep{lustig2006rapid, lustig2007sparse}. Meanwhile, LASSO and its extensions have become standard tools in machine learning, genomics, and econometrics, where sparse models are crucial for managing large-scale datasets and producing interpretable results \citep{tibshirani1996regression, hastle2015statistical, puig2011multidimensional}.

This book aims to provide a comprehensive overview of \textit{sparse optimization}, with a focus on both \textit{sparse signal recovery} and \textit{sparse regularization} techniques. We will begin by exploring the foundations of sparse optimization, delving into the mathematical tools and models that underpin sparse signal recovery and LASSO. We will then discuss key algorithms for both \textit{sparse recovery} (e.g., basis pursuit, matching pursuit) and \textit{sparse regularization} (e.g., LASSO, elastic net), along with their applications in real-world problems. Throughout the text, we balance intuitive explanations with rigorous mathematical formulations to provide a comprehensive resource for both newcomers and experts in the field.
Our aim is twofold: to provide a self-contained entry point for students and researchers new to the field, and to offer a rigorous reference for practitioners seeking to apply sparse optimization in science and engineering. 

The book is organized as follows. 
Chapters~\ref{chapter:mathback} and \ref{chapter:sparse_opt_cond} lay the groundwork by reviewing essential concepts from (convex) optimization and their relevance to sparse problems.
Chapters~\ref{chapter:design}$\sim$\ref{chapter:ensur_rips} discuss the design of measurement matrices that guarantee exact or stable recovery of sparse signals, including constructions based on random matrix theory and the restricted isometry property (RIP).
Chapters~\ref{chapter:algouni}$\sim$\ref{chapter:spar_recov} present algorithms for solving various sparse optimization problems, ranging from greedy methods to convex relaxations.

In the remainder of this chapter, we briefly review basic mathematical notation and foundational concepts. Additional definitions will be introduced as needed throughout the text to ensure clarity.

\section{Linear Algebra}

In all cases, scalars will be denoted in a non-bold font possibly with subscripts (e.g., $a$, $\alpha$, $\alpha_i$). 
Vectors will be represented by \textbf{boldface} lowercase letters, possibly with subscripts  (e.g., $\bmu$, $\bx$, $\bx_n$, $\bz$), and
matrices by \textbf{boldface} uppercase letters, also possibly with subscripts  (e.g., $\bX$, $\bL_j$). 
The $i$-th element of a vector $\bz$ will be written as $z_i$ in  a non-bold font.
To distinguish deterministic quantities from random ones, we adopt the following convention:
\begin{itemize}
\item Random variables are denoted using \textit{normal (upright) font} for scalars  (e.g., $\textnormal{a}$ and $\textnormal{b}_1$ are random variables, whereas italics $a$ and $b_1$ denote deterministic scalars).
\item Random vectors are indicated by \textit{normal-font bold lowercase letters}, possibly with subscripts (e.g., $\rva$ and $\rvb_1$ are random vectors, while italics bold $\ba$ and $\bb_1$ represent deterministic vectors) .
\item Random matrices are denoted by \textit{normal-font bold uppercase letters}, possibly with subscripts (e.g., $\rmA$ and $\rmB_1$ are random matrices, whereas italics bold $\bA$ and $\bB_1$ denote deterministic matrices).
\end{itemize}
Subarrays are formed by fixing a subset of indices of a matrix.
The element located in the $i$-th row and $j$-th column of a matrix $\bX$ (i.e., the $(i,j)$ entry) is denoted by $x_{ij}$. 
Consequently, a matrix $\bX\in\real^{n\times p}$ may be expressed as $\bX=\{x_{ij}\}_{i,j=1}^{n,p}=[x_{ij}]$.
For notational convenience, we adopt \textbf{Matlab-style indexing}: 
\begin{itemize}
\item The submatrix consisting of rows $i$ through $j$ and columns $k$ through $m$ of $\bX$ is denoted by  $\bX_{i:j,k:m} \equiv \bX[i:j,k:m]$. 
\item A colon is used to indicate all elements along a dimension. For example, $\bX_{:,k:m} \equiv\bX[:,k:m]$ denotes the submatrix formed by columns $k$ through $m$; and $\bX_{:,k}\equiv\bX[:,k]$ denotes the $k$-th column of $\bX$. 
\item   Alternatively, the $k$-th column of $\bX$ may be denoted more compactly by $\bx_k$; and the $k$-th row of $\bX$ can be denoted as $\bx^{(k)}$.
\end{itemize}
When the selected indices are not contiguous, let $\sI$ and $\sJ$ be ordered index sets. Then: 
\begin{itemize}
\item $\bX[\sI, \sJ]$ denotes the submatrix of $\bX$ obtained by extracting the rows and columns of $\bX$ indexed by $\sI$ and $\sJ$, respectively.
\item $\bX[:, \sJ]$ denotes the submatrix of $\bX$ obtained by extracting the columns of $\bX$ indexed by $\sJ$, where the colon again signifies ``all rows."
\end{itemize}
\index{Matlab notation}
\begin{definition}[Matlab notation]\label{definition:matlabnotation}
Let $\bX\in \real^{n\times p}$, and let $\sI=\{i_1, i_2, \ldots, i_k\}$ and $\sJ=\{j_1, j_2, \ldots, j_l\}$ be two index sets. 
Then $\bX[\sI,\sJ]$ denotes the $k\times l$ submatrix
$$
\bX[\sI,\sJ]=
\begin{bmatrix}
x_{i_1,j_1} & x_{i_1,j_2} &\ldots & x_{i_1,j_l}\\
x_{i_2,j_1} & x_{i_2,j_2} &\ldots & x_{i_2,j_l}\\
\vdots & \vdots&\ddots & \vdots\\
x_{i_k,j_1} & x_{i_k,j_2} &\ldots & x_{i_k,j_l}\\
\end{bmatrix}.
$$
Similarly, $\bX[\sI,:]$ denotes the $k\times p$ submatrix, and $\bX[:,\sJ]$ denotes the $n\times l$ submatrix analogously.
We should also notice that the range of the index satisfies:
$$
\left\{
\begin{aligned}
0&\leq \min(\sI) \leq \max(\sI)\leq n;\\
0&\leq \min(\sJ) \leq \max(\sJ)\leq p.
\end{aligned}
\right.
$$
More compactly, for $\bX\in\real^{n\times p}$ and $\bbeta\in\real^p$, we let 
\begin{equation}
\bX_{\sI} \triangleq \bX[:,\sI]  \in\real^{n\times k}
\qquad 
\text{and}
\qquad 
\bbeta_{\sI} \triangleq \bbeta[\sI]\in\real^k.
\end{equation}
Note that $\bX_\sI^\top$ denotes the transpose of the submatrix $\bX_{\sI}$, not a submatrix of $\bX^\top$.

On the other hand, when an index set $\sI$ appears  as an argument of a function---written as $\bX(\sI)$ or $\bbeta(\sI)$---it denotes a projection operator that zeroes out all entries not indexed by $\sI$.
Specifically:
\begin{itemize}
\item For $\bbeta\in\real^p$, the vector $\bbeta(\sI)\in\real^p$ equals $\bbeta$ on the entries in $\sI$ and is zero elsewhere; that is, $\bbeta(\sI)_i = \beta_i$ if $i\in\sI$, and 0 otherwise.
\item For $\bX\in\real^{n\times p}$, the matrix $\bX(\sI)\in\real^{n\times p}$ retains the columns of $\bX$ indexed by $\sI$ and sets all other columns to zero.
\end{itemize}
With this notation, the matrix-vector product decomposes as
$$
\bX\bbeta 
= \bX(\sI)\bbeta(\sI) + \bX(\comple{\sI})\bbeta(\comple{\sI}) 
= \bX[:,\sI]\bbeta[\sI] + \bX[:,\comple{\sI}]\bbeta[\comple{\sI}] 
= \bX_{\sI}\bbeta_{\sI} + \bX_{\comple{\sI}}\bbeta_{\comple{\sI}},
$$
where $\comple{\sI}$ denotes the complement set of $\sI$.
\end{definition}

Throughout this book, all vectors are assumed to be column vectors unless explicitly transposed. A row vector is denoted as the transpose of a column vector, e.g., $\bx^\top$. 
Concrete vectors are written using MATLAB-style syntax.
A column vector is written with semicolons $``;"$ separating entries,  e.g., 
$$\bx=[1;2;3] \qquad \text{(column vector)}
$$ 
is a column vector in $\real^3$. Similarly, A row vector uses commas $``,"$, e.g., 
$$\by=[1,2,3]\qquad \text{(row vector)}
$$ 
is a row vector with 3 values. 
Equivalently, a column vector may be expressed as the transpose of a row vector: $\by=[1,2,3]^\top$ is a column vector.

The transpose of a matrix $\bX$ is denoted by $\bX^\top$,
and its inverse (when it exists)  by $\bX^{-1}$. 
The $p \times p$ identity matrix is denoted by $\bI_p$ or simply by $\bI$. 
A vector or matrix of all zeros is denoted by the \textbf{boldface} symbol $\bzero$; its dimensions are inferred from context, or explicitly indicated as $\bzero_p$ for a $p$-dimensional zero vector.
Similarly, a vector or matrix of all ones is  denoted by a \textbf{boldface} one $\bone$, whose sizes are clear from  context; or we denote $\bone_p$ to be the vector of all ones with $p$ entries.
Subscripts are often omitted when the dimensions are evident from  context.

\index{Eigenvalue}
\index{Eigenvector}
\begin{definition}[Eigenvalue, eigenvector]
Let $\bX\in\complex^{p\times p}$. A scalar $\lambda \in \complex$ is called a \textit{(right) eigenvalue} (also known as a  \textit{proper value}, or \textit{characteristic value}) of $\bX$ if there exists  a nonzero vector $\bu \in \complex^p$ such that
\begin{equation*}
\bX \bu = \lambda \bu.
\end{equation*}
In this case, $\bu$ is called a \textit{(right) eigenvector} of $\bX$ associated with $\lambda$.
\end{definition}
For simplicity, we restrict our attention to real-valued matrices unless otherwise stated. Unless explicitly noted, all eigenvalues discussed are assumed to be real as well.  

Intuitively, an eigenvector $\bu$ of a matrix $\bX$ represents a direction in $\real^p$ that is invariant under the linear transformation defined by $\bX$: applying $\bX$ to $\bu$ does not change its direction---it only scales it by the corresponding eigenvalue $\lambda$.
Note that while real matrices can have complex eigenvalues in general, symmetric real matrices always have real eigenvalues (see Theorem~\ref{theorem:spectral_theorem}).

The pair $(\lambda, \bu)$  is commonly  referred to as an \textit{eigenpair}. 
Importantly, eigenvectors are not unique:  if  $\bu$ is an eigenvector, then so is any nonzero scalar multiple $\eta\bu$ (with $\eta\in\real\setminus \{\bzero\}$).
To resolve this ambiguity, eigenvectors are typically normalized---for example, by requiring $\normtwo{\bu}=1$ (unit $\ell_2$-norm; see Section~\ref{section:vec_norms}). 
Furthermore, since both $\bu$ and $-\bu$ correspond to the same eigenvalue and represent the same direction up to sign, it is common to fix the sign by convention (e.g., requiring the first nonzero entry of $\bu$ to be positive).

An important result concerning the eigenvalues of a complex matrix is the \textit{Gershgorin circle theorem} (also known as the \textit{Gershgorin disk theorem}), introduced by the Soviet mathematician Semen Aronovich Gershgorin in 1931 \citep{gershgorin1931uber}.
\begin{theoremHigh}[Gershgorin circle theorem \citep{gershgorin1931uber}\index{Gershgorin circle theorem}]\label{theorem:eign_disc} 
Let $ \bX = [x_{ij}] $ be a  $ p \times p $ complex matrix.
For each row $i$, define the {Gershgorin radius} $ R_i $ as
$$
R_i \triangleq  \sum_{\substack{j=1, j \neq i}}^p \abs{x_{ij}}.
$$
Then every eigenvalue $ \lambda $ of   $ \bX $ lies within at least one of the \textit{Gershgorin disks} centered at the diagonal entries $ x_{ii} $ with radius $ R_i $. 
In other words, for every eigenvalue $ \lambda $, there exists an index $ i $ such that:
$$
\abs{\lambda - x_{ii}} \leq R_i.
$$
\end{theoremHigh}

This means that all eigenvalues of  $ \bX $ are located within the union of these $p$ disks in the complex plane. Each disk is centered at the diagonal element $ x_{ii} $and its radius equals the sum of the absolute values of the off-diagonal entries in the corresponding row.
The Gershgorin circle theorem provides  a simple and effective way to estimate the location of eigenvalues without computing them explicitly.
It is especially useful in analyzing matrix properties, stability, and determining whether a matrix is invertible. 
For instance, if none of the Gershgorin disks contains the origin, then 0 cannot be an eigenvalue of $\bX$, and hence $\bX$ is invertible.

In linear algebra, every vector space has a basis, and every vector in that space  can be expressed uniquely as a linear combination of the basis vectors. 
Building on this idea, we define the span and dimension of a subspace using the concept of a basis.

\index{Subspace}
\index{Span}
\begin{definition}[Subspace and span]
A nonempty subset $\mathcalV$ of $\real^n$ is called a subspace if it is closed under linear combinations; that is, $x\ba+y\bb\in \mathcalV$ for every $\ba,\bb\in \mathcalV$ and every $x,y\in \real$.
Moreover, if every vector $\bv$ in a subspace $\mathcalV$ can be written as a linear combination of the vectors $\{\bx_1, \bx_2, \ldots,$ $\bx_n\}$, then we say that these vectors \textit{span} the subspace $\mathcalV$.
\end{definition}

\index{Linearly independent}
The concept of linear independence of a set of vectors is central to linear algebra. Two equivalent definitions are given below.
\begin{definition}[Linearly independent]
A set of vectors $\{\bx_1, \bx_2, \ldots, \bx_n\}$ is said to be \textit{linearly independent} if there is no combination can get $a_1\bx_1+a_2\bx_2+\ldots+a_n\bx_n=0$ except all $a_i$'s are zero. An equivalent definition is that $\bx_1\neq \bzero$, and for every $k>1$, the vector $\bx_k$ does not lie in the span of $\{\bx_1, \bx_2, \ldots, \bx_{k-1}\}$.
\end{definition}

\index{Basis}
\index{Dimension}
\begin{definition}[Basis and dimension]
A set of vectors $\{\bx_1, \bx_2, \ldots, \bx_n\}$ is called a \textit{basis} of a subspace $\mathcalV$ if the vectors are linearly independent and span $\mathcalV$. 
Every basis of a given subspace contains the same number of vectors. 
This number is called the \textit{dimension} of the subspace $\mathcalV$. 
By convention, the trivial subspace $\{\bzero\}$ has dimension zero.
Furthermore, any subspace of positive dimension admits an \textit{orthogonal basis}---that is, one can always choose a basis consisting of mutually orthogonal vectors (see Definition~\ref{definition:orthogn_mat}).
\end{definition}

\index{Column space}
\begin{definition}[Column space (range) and row space]
Let $\bX$ be an $n \times p$ real matrix. The \textit{column space (or range)} of $\bX$ is defined as  the set of all vectors that can be expressed as a linear combination of its columns:
\begin{equation*}
\cspace (\bX) = \{ \bv\in \real^n\mid \exists\, \bu \in \real^p, \, \bv = \bX \bu \}.
\end{equation*}
Similarly, the \textit{row space} of $\bX$ is the set of all vectors spanned by its rows. Equivalently, it is the column space of $\bX^\top$:
\begin{equation*}
\cspace (\bX^\top) = \{ \bu\in \real^p\mid  \exists\, \bv \in \real^n, \, \bu = \bX^\top \bv \}.
\end{equation*}
\end{definition}

\index{Null space}
\index{Left null space}
\index{Right null space}
\begin{definition}[Null space (nullspace, kernel)]\label{definition:nullspace}
Let $\bX$ be an $n \times p$ real matrix. The \textit{null space (or kernel, or nullspace)} of $\bX$ is defined as  the set:
\begin{equation*}
\nspace (\bX) = \{\bv \in \real^p\mid  \, \bX \bv = \bzero \}.
\end{equation*}
In some cases, the null space of $\bX$ is also referred to as the \textit{right null space} of $\bX$.
And the null space of $\bX^\top$ is given by 
\begin{equation*}
\nspace (\bX^\top) = \{\bu \in \real^n\mid  \, \bX^\top \bu = \bzero \},
\end{equation*}
and is commonly called the \textit{left null space} of $\bX$.
\end{definition}

\index{Rank}
\index{Dimension}
\begin{definition}[Rank]\label{definition:rank}
The $rank$ of a matrix $\bX\in \real^{n\times p}$ is defined as the dimension of its column space. That is, the rank of $\bX$ is equal to the maximum number of linearly independent columns of $\bX$, which is also equal to  the maximum number of linearly independent rows of $\bX$. The matrix $\bX$ and its transpose $\bX^\top$ always have the same rank. 
We say that $\bX$ has \textit{full rank} if its rank  equals $\min\{n,p\}$.  
Specifically, given a vector $\bu \in \real^n$ and a vector $\bv \in \real^p$, then the $n\times p$ matrix $\bu\bv^\top$ obtained by their outer product of vectors is of rank 1. 
In summary, the rank of a matrix can be characterized as:
\begin{itemize}
\item the number of linearly independent columns;
\item the number of linearly independent rows;
\item and remarkably, these two numbers are always equal (see Lemma~\ref{lemma:equal-dimension-rank}).
\end{itemize}
\end{definition}

Both the column space of $\bX$ and the null space of $\bX^\top$ are subspaces of $\real^n$. In fact, every vector in $\nspace(\bX^\top)$ is orthogonal  to every vector in $\cspace(\bX)$, and vice versa. 
Similarly, every vector in $\nspace(\bX)$ is also orthogonal  to every vector in $\cspace(\bX^\top)$, and vice versa.
This leads to the concept of \textit{orthogonal complementary subspaces}.

\index{Orthogonal complement}
\begin{definition}[Orthogonal complement in general]
The \textit{orthogonal complement} $\mathcalV^\perp\subseteq\real^n$ of a subspace $\mathcalV\subseteq\real^n$ consists of all vectors that are perpendicular to every vector in $\mathcalV$. Formally,
$$
\mathcalV^\perp = \{\bv\in\real^n\mid  \bv^\top\bu=0, \ \forall\, \bu\in \mathcalV  \}.
$$
These two subspaces are disjoint and together span the entire space $\real^n$. 
The dimensions of $\mathcalV$ and $\mathcalV^\perp$ add up to the dimension of the entire space: $\dim(\mathcalV)+\dim(\mathcalV^\perp)=n$. Furthermore, $(\mathcalV^\perp)^\perp=\mathcalV$.
\end{definition}

\begin{definition}[Orthogonal complement of column space]\label{definition:ortho_comp_col}
Let $\bX$ be an $n \times p$ real matrix. The orthogonal complement of the column space $\cspace(\bX)$, denoted  $\cspace^{\bot}(\bX)$, is the subspace defined by:
\begin{equation*}
\begin{aligned}
\cspace^{\bot}(\bX) 
&= \{\bv\in \real^n\mid \, \bv^\top \bX \bu=\bzero, \, \forall\, \bu \in \real^p \} \\
&=\{\bv\in \real^n\mid \, \bv^\top \bw = \bzero, \, \forall\, \bw \in \cspace(\bX) \}.
\end{aligned}
\end{equation*}
\end{definition}

For any matrix $\bX \in \real^{n \times p}$ of rank $r$, we obtain the \textit{four fundamental subspaces}:
\begin{itemize}
\item  $\cspace(\bX)$: Column space of $\bX$, i.e., linear combinations of columns with dimension $r$.
\item  $\nspace(\bX)$: (Right) null space of $\bX$, i.e., all $\bu$ satisfying $\bX\bu=\bzero$ with dimension $p-r$.
\item  $\cspace(\bX^\top)$: Row space of $\bX$, i.e., linear combinations of rows with dimension $r$.
\item  $\nspace(\bX^\top)$: Left null space of $\bX$, i.e., all $\bv$ satisfying $\bX^\top \bv=\bzero$ with dimension $n-r$. 
\end{itemize}
Furthermore, $\nspace(\bX)$ is the orthogonal complement of $\cspace(\bX^\top)$, and $\cspace(\bX)$ is the orthogonal complement of $\nspace(\bX^\top)$. This is known as the \textit{fundamental theorem of linear algebra}.
See \citet{lu2022matrix} for a proof. 

\index{Rank}
\index{Dimension}
We now  establish the equivalence of the dimensions of column  and row spaces stated in Definition~\ref{definition:rank}.
\begin{lemma}[Dimension of column space and row space]\label{lemma:equal-dimension-rank}
The dimension of the column space of a matrix $\bX\in \real^{n\times p}$ equals the dimension of its
row space, i.e., the row rank and the column rank of a matrix $\bX$ are identical.
\end{lemma}
\begin{proof}[of Lemma~\ref{lemma:equal-dimension-rank}]
We first notice that the null space of $\bX$ is orthogonal complementary to the row space of $\bX$: $\nspace(\bX) \bot \cspace(\bX^\top)$ (where the row space of $\bX$ is equivalent to the column space of $\bX^\top$). 
That is, vectors in the null space of $\bX$ are orthogonal to vectors in the row space of $\bX$. To see this, suppose $\bX=[\bx_1^\top; \bx_2^\top; \ldots; \bx_n^\top]$ is the row partition of $\bX$. For any vector $\bbeta\in \nspace(\bX)$, we have $\bX\bbeta = \bzero$; that is, $[\bx_1^\top\bbeta; \bx_2^\top\bbeta; \ldots; \bx_n^\top\bbeta]=\bzero$. 
And since the row space of $\bX$ is spanned by $\{\bx_1^\top, \bx_2^\top, \ldots, \bx_n^\top\}$, thus, $\bbeta$ is perpendicular to any vectors from $\cspace(\bX^\top)$. This indicates $\nspace(\bX) \bot \cspace(\bX^\top)$.

Now, assuming the dimension of the row space of $\bX$ is $r$,  \textcolor{mylightbluetext}{let $\br_1, \br_2, \ldots, \br_r$ be a set of vectors in $\real^p$ and form a basis for the row space}. 
Consequently, the $r$ vectors $\bX\br_1, \bX\br_2, \ldots, \bX\br_r$ are in the column space of $\bX$; furthermore, they are linearly independent. To see this, suppose we have a linear combination of the $r$ vectors: $\beta_1\bX\br_1 + \beta_2\bX\br_2+ \ldots+ \beta_r\bX\br_r=0$, that is, $\bX(\beta_1\br_1 + \beta_2\br_2+ \ldots+ \beta_r\br_r)=0$, and the vector $\bv=\beta_1\br_1 + \beta_2\br_2+ \ldots+ \beta_r\br_r$ is in null space of $\bX$. But since $\{\br_1, \br_2, \ldots, \br_r\}$ is a basis for the row space of $\bX$, $\bv$ is thus also in the row space of $\bX$. We have shown that vectors from null space of $\bX$ is perpendicular to vectors from row space of $\bX$, thus $\bv^\top\bv=0$ and $\beta_1=\beta_2=\ldots=\beta_r=0$. Then, \textcolor{mylightbluetext}{$\bX\br_1, \bX\br_2, \ldots, \bX\br_r$ are in the column space of $\bX$ and they are independent}. This means that the dimension of the column space of $\bX$ is larger than $r$. This result shows that \textbf{row rank of $\bX\leq $ column rank of $\bX$}. 

Applying the same argument to $\bX^\top$ yields the reverse inequality: \textbf{column rank of $\bX\leq $ row rank of $\bX$}. Hence, the two ranks are equal, completing the proof.
\end{proof}

From the previous proof, we can also conclude that if $\{\br_1, \br_2, \ldots, \br_r\}$ is a basis for the row space of $\bX$, then \textcolor{black}{$\{\bX\br_1, \bX\br_2, \ldots, \bX\br_r\}$ forms a basis for the column space of $\bX$}. 
We formalize this observation in the following lemma.

\index{Column basis}
\index{Row basis}
\begin{lemma}[Column basis from row basis]\label{lemma:column-basis-from-row-basis}
For any matrix $\bX\in \real^{n\times p}$, let $\{\br_1, \br_2, \ldots, \br_r\}$ be a set of vectors in $\real^p$, which forms a basis for the row space of $\bX$. Then, the set $\{\bX\br_1, \bX\br_2, \ldots, \bX\br_r\}$ forms a basis for the column space of $\bX$.
\end{lemma}


\index{Permutation matrix}
\begin{definition}[Permutation matrix]\label{definition:permutation-matrix}
A \textit{permutation matrix} $\bP\in \real^{p\times p}$ is a square binary matrix that has exactly one entry of 1 in each row and each column, and 0's elsewhere. 
\paragraph{Row point.} That is, the permutation matrix $\bP$ has the rows of the identity $\bI$ in any order, and the order decides the sequence of the row permutation. If we want to permute the rows of a matrix $\bX$, we multiply on the left $\bP\bX$. 
\paragraph{Column point.} Or, equivalently, the permutation matrix $\bP$ has the columns of the identity $\bI$ in any order, and the order decides the sequence of the column permutation. To apply a column permutation to $\bX$, we multiply on the right $\bX\bP$.
\end{definition}

The permutation matrix $\bP$ can be represented more compactly by a index set $\sJ \in \integer_+^p$ such that $\bP = \bI[:, \sJ]$, where $\bI$ is the $p\times p$ identity matrix. 
 Notably, the entries of $\sJ$ are a permutation of $\{1, 2, \dots, p\}$, so their sum is $1+2+\ldots+p= \frac{p^2+p}{2}$.


\index{Determinant}
\begin{definition}[Determinant: Laplace expansion by minors]\label{definition:determinant}
Let $\bX\in\real^{p\times p}$ be a square matrix, and let $\bX_{ij}\in\real^{(p-1)\times (p-1)}$ denote the submatrix of $\bX$ obtained by deleting the $i$-th row and $j$-th column of $\bX$. 
The \textit{determinant} of $\bX$ can be computed recursively using the following formulas:
\begin{equation}\label{equation:def_det}
\det(\bX)=\sum_{k=1}^{p} (-1)^{i+k} x_{ik}\det(\bX_{ik})=\sum_{k=1}^{p}(-1)^{k+j} x_{kj}\det(\bX_{kj}),
\end{equation}
where the first equation is the \textit{Laplace expansion by minors along row $i$}, and the second equation is the \textit{Laplace expansion by minors along column $j$}.
Equivalently, for any integer $r$ with $1 \leq r \leq p$, let $\sJ\subseteq\{1,2,\ldots,p\}$ be  an index set of cardinality $r$ ($\abs{\sJ}=r$), and let $\comple{\sJ}=\{1,2,\ldots,n\}\backslash \sJ$  denote its complementary set.
Then
$$
\begin{aligned}
\det(\bX) 
=\sum_{\sI} (-1)^{\gamma} \det(\bX[\sI,\sJ])\det(\bX[\comple{\sI}, \comple{\sJ}])
=\sum_{\sI} (-1)^{\gamma} \det(\bX[\sJ,\sI])\det(\bX[\comple{\sJ}, \comple{\sI}]),
\end{aligned}
$$
where $\gamma=\sum_{i\in \sI} i +\sum_{j\in \sJ}j$, and the sum is taken over all  index sets $\sI\subseteq\{1,2,\ldots,p\}$ with cardinality $r$.
When $r=1$, this reduces to \eqref{equation:def_det}.
\end{definition}

The determinant assigns to each square matrix a scalar value that encodes geometric and algebraic properties of the associated linear transformation.
For a $2 \times 2$ matrix, the absolute value of the determinant equals the area of the parallelogram spanned by its column vectors. The sign indicates orientation: it is positive if the columns form a counterclockwise basis, and negative if clockwise.
For a $3 \times 3$ matrix, the determinant gives the signed volume of the parallelepiped defined by its three column vectors, with the sign again reflecting orientation.
In general, for a linear transformation represented by a $p \times p$ matrix, the absolute value of the determinant measures the factor by which $p$-dimensional volume is scaled under the transformation. A positive determinant means the orientation (or ``handedness") of the space is preserved; a negative determinant means it is reversed.
A square matrix is invertible if and only if its determinant is nonzero. Geometrically, this means the transformation does not collapse the space into a lower-dimensional subspace or a single point---both of which occur precisely when the determinant is zero.
Algebraically, the determinant of a matrix equals the product of its eigenvalues (counting multiplicities). Thus, it reflects the cumulative scaling effect of the transformation along its eigenvector directions.

We now summarize key properties of the determinant.
\begin{lemma}[Properties of determinant]\label{lemma:determinant-intermezzo}
The determinant satisfies the following properties:
\begin{itemize}
\item  The determinant of the product of two matrices is $\det(\bX\bY)=\det (\bX)\det(\bY)$;

\item The determinant of the transpose is $\det(\bX^\top) = \det(\bX)$;

\item If $\lambda$ is an eigenvalue of $\bX$, then $\det(\bX-\lambda\bI) =0$;

\item Determinant of any identity matrix is $1$;

\item For any orthogonal matrix $\bQ$ (i.e., $\bQ\bQ^\top=\bQ^\top\bQ=\bI$), we have 
$$
\det(\bQ) = \det(\bQ^\top) = \pm 1, \qquad \text{since  } \det(\bQ^\top)\det(\bQ)=\det(\bQ^\top\bQ)=\det(\bI)=1;
$$

\item For any square matrix $\bX$ and orthogonal matrix $\bQ$, we have 
$$
\det(\bX) = \det(\bQ^\top) \det(\bX)\det(\bQ) =\det(\bQ^\top\bX\bQ);
$$
\item Suppose $\bX\in\real^{p\times p}$, then $\det(-\bX) = (-1)^p \det(\bX)$.
\end{itemize}
\end{lemma}

\index{Positive definite}\index{Positive semidefinite}
Positive definiteness and positive semidefiniteness are among the most important properties a matrix can possess.

\begin{definition}[Positive definite and positive semidefinite]\label{definition:psd-pd-defini}
A symmetric matrix $\bX\in \real^{p\times p}$ is said to be \textit{positive definite (PD)} if $\bbeta^\top\bX\bbeta>0$ for all nonzero $\bbeta\in \real^p$, denoted by $\bX\succ \bzero$.
It is called \textit{positive semidefinite (PSD)} if $\bbeta^\top\bX\bbeta \geq 0$ for all $\bbeta\in \real^p$, denoted by $\bX\succeq \bzero$. 
\footnote{
In this book, positive (semi)definite matrices are always assumed to be symmetric. The notions of positive definiteness and semidefiniteness are meaningful only for symmetric matrices.
}
\footnote{A symmetric matrix $\bX\in\real^{p\times p}$ is called \textit{negative definite} (ND) if $\bbeta^\top\bX\bbeta<0$ for all nonzero $\bbeta\in\real^p$; 
a symmetric matrix $\bX\in\real^{p\times p}$ is called \textit{negative semidefinite} (NSD) if $\bbeta^\top\bX\bbeta\leq 0$ for all $\bbeta\in\real^p$;
and a symmetric matrix $\bX\in\real^{p\times p}$ is called \textit{indefinite} (ID) if there exist $\bbeta$ and $\balpha\in\real^p$ such that $\bbeta^\top\bX\bbeta<0$ and $\balpha^\top\bX\balpha>0$.
}
\end{definition}

A fundamental result in linear algebra states that a symmetric matrix $\bX$ is positive definite if and only if all of its eigenvalues are strictly positive. Similarly, $\bX$ is positive semidefinite if and only if all of its eigenvalues are nonnegative.
This leads to the following theorem:
\begin{theoremHigh}[Eigenvalue characterization theorem]\label{theorem:eigen_charac}
A matrix $\bX\in\real^{p\times p}$ is positive definite if and only if all of its eigenvalues are \textit{positive}. Similarly, a matrix $\bX$ is positive semidefinite if and only if all of its eigenvalues  \textit{nonnegative}.~\footnote{The trace, determinant, and principal minors of a positive (semi)definite matrix is discussed in Problem~\ref{prob:tr_de_pd}.}
Moreover, we have the following implications:
\begin{itemize}
\item $\bX-\gamma\bI\succeq \bzero$ if and only if $\lambda_{\min}(\bX) \geq \gamma$;
\item  $\bX-\gamma\bI\succ \bzero$ if and only if $\lambda_{\min}(\bX) > \gamma$;
\item $\bX-\gamma\bI\preceq \bzero$ if and only if $\lambda_{\max}(\bX) \leq \gamma$;
\item $\bX-\gamma\bI\prec \bzero$ if and only if $\lambda_{\max}(\bX) < \gamma$;
\item $\lambda_{\min}(\bX)\bI\preceq \bX \preceq \lambda_{\max}(\bX)\bI$,
\end{itemize}
\noindent where $\lambda_{\min}(\bX)$ and $\lambda_{\max}(\bX)$ represent the minimum and maximum eigenvalues of $\bX$, respectively, and $\bB \prec \bC$ means $\bC-\bB$ is PSD.
\end{theoremHigh}
This theorem provides an alternative characterization  of positive definiteness and positive semidefiniteness in terms of  eigenvalues.
Given an eigenpair $(\lambda, \bbeta)$ of $\bX$, the forward implication can be verified that $\bbeta^\top\bX\bbeta=\lambda\bbeta^\top\bbeta>0$ such that $ \lambda=(\bbeta^\top\bX\bbeta)/(\bbeta^\top\bbeta)>0$ (resp., $\geq 0$) if $\bX$ is PD (resp., PSD).
The full equivalence can be established using the spectral theorem (Theorem~\ref{theorem:spectral_theorem}).

A notable subclass of positive semidefinite matrices consists of certain diagonally dominant matrices, which we define rigorously below.
\begin{definition}[Diagonally dominant matrices\index{Diagonally dominant}]
Let $\bX\in\real^{p\times p}$ be a symmetric matrix. 
Then $\bX$ is called \textit{diagonally dominant} if 
$$
\abs{x_{ii}} \geq \sum_{j\neq i} \abs{x_{ij}}, \qquad \forall\, i\in\{1,2,\ldots, p\},
$$
and \textit{strictly diagonally dominant} if 
$$
\abs{x_{ii}} > \sum_{j\neq i} \abs{x_{ij}}, \qquad \forall\, i\in\{1,2,\ldots, p\}.
$$
\end{definition}

We  now show that a \textit{diagonally dominant matrix} with nonnegative diagonal entries is positive semidefinite,  and that a \textit{strictly diagonally dominant matrix} with positive diagonal entries is positive definite. 

\begin{theoremHigh}[Positive definiteness of diagonally dominant matrices]\label{theorem:pd_diag_domi}
Let  $\bX\in\real^{p\times p}$ be a   symmetric matrix.
Then:
\begin{enumerate}[(i)]
\item If $\bX$ is diagonally dominant with nonnegative diagonals, then $\bX$ is positive semidefinite.
\item If $\bX$ is strictly diagonally dominant with positive diagonals, then $\bX$ is positive definite.
\end{enumerate}
\end{theoremHigh}
This result follows directly from the Gershgorin circle theorem (Theorem~\ref{theorem:eign_disc}). However, we provide an elementary proof below based on eigenvalue analysis.
\begin{proof}[of Theorem~\ref{theorem:pd_diag_domi}]
\textbf{(i).} Suppose, for contradiction, that $\bX$ is not positive semidefinite with a negative eigenvalue $\lambda$ associated with an eigenvector $\bv$ such that $\bX\bv=\lambda\bv$. 
Let $i$ be an index such that $\abs{v_i}=\max_{1\leq j\leq p}\abs{v_j}$.
Consider the $i$-th component of $(\bX-\lambda\bI)\bv=\bzero$, we have
$$
\abs{x_{ii} - \lambda}\cdot \abs{v_i}
=
\abs{\sum_{j\neq i} x_{ij} v_j}
\leq 
\left( \sum_{ j\neq i} \abs{x_{ij}} \right) \abs{v_i}
\leq \abs{x_{ii}} \cdot \abs{v_i}.
$$
This implies $\abs{x_{ii} - \lambda}  \leq \abs{x_{ii}}$. 
Since the diagonal entries are nonnegative, we have  $x_{ii}< {x_{ii} - \lambda}  \leq {x_{ii}}$, and this leads to a contradiction.

\paragraph{(ii).} From part (i), we know that $\bX$ is positive semidefinite. Suppose, for contradiction, that  $\bX$ is not positive definite with a zero eigenvalue $0$ associated with a nonzero eigenvector $\bv$ such that $\bX\bv=\bzero$.   Similarly, we have 
$$
\abs{x_{ii}}\cdot \abs{v_i}
=
\abs{\sum_{j\neq i} x_{ij} v_j}
\leq 
\left( \sum_{ j\neq i} \abs{x_{ii}} \right) \abs{v_i}
< \abs{x_{ii}} \cdot \abs{v_i},
$$
which is impossible and the result follows.
\end{proof}

\index{Positive semidefinite}
\begin{exercise}[Power of PSD]
Let $\bX$ be PSD. Show that  $\bX^k$ is also PSD for $k=1,2,\ldots$.
\end{exercise}

\index{Nonsingularity}
\index{Nonsingular matrix}
\index{Invertible matrix}
From an introductory linear algebra course, we recall the following well-known equivalences characterizing nonsingular (i.e., invertible) square matrices.	
\begin{remark}[List of equivalence of nonsingularity for a matrix]\label{remark:equiva_nonsingular}
Let $\bX\in \real^{p\times p}$ be a square matrix. 
The following statements are equivalent:
\begin{itemize}
\item $\bX$ is nonsingular;~\footnote{The source of the name  is a result of the singular value decomposition (SVD; see Section~\ref{section:SVD}).}
\item $\bX$ is invertible, i.e., $\bX^{-1}$ exists;
\item $\bX\bu=\by$ has a unique solution $\bu = \bX^{-1}\by$;
\item $\bX\bu = \bzero$ has a unique, trivial solution: $\bu=\bzero$;
\item Columns of $\bX$ are linearly independent;
\item Rows of $\bX$ are linearly independent;
\item $\det(\bX) \neq 0$; 
\item $\dim(\nspace(\bX))=0$;
\item $\nspace(\bX) = \{\bzero\}$, i.e., the null space is trivial;
\item $\cspace(\bX)=\cspace(\bX^\top) = \real^p$, i.e., the column space or row space span the entire $\real^p$;
\item $\bX$ has full rank $r=p$;
\item $\bX^\top\bX$ is symmetric positive definite (PD);
\item $\bX$ has $p$ nonzero (positive) singular values;
\item All eigenvalues are nonzero.
\end{itemize}
\end{remark}
It is essential to keep these equivalences in mind, as they appear frequently in both theoretical and applied contexts.
Similarly, the following remark collects equivalent characterizations of singular matrices---particularly in the context of eigenvalues.

\index{Singular matrix}
\begin{remark}[List of equivalence of singularity for a matrix]\label{remark:equiva_singular}
Let $\bX\in \real^{p\times p}$ be a square matrix, and let  $(\lambda, \bu)$ be an eigenpair of $\bX$. 
The following statements are equivalent:
\begin{itemize}
\item $(\bX-\lambda\bI)$ is singular;
\item $(\bX-\lambda\bI)$ is not invertible;
\item $(\bX-\lambda\bI)\bv = \bzero$ has nonzero $\bv\neq \bzero$ solutions, and $\bv=\bu$ is one of such solutions;
\item $(\bX-\lambda\bI)$ has linearly dependent columns;
\item $\det(\bX-\lambda\bI) = 0$; 
\item $\dim(\nspace(\bX-\lambda\bI))>0$;
\item Null space of $(\bX-\lambda\bI)$ is nontrivial;
\item Columns of $(\bX-\lambda\bI)$ are linearly dependent;
\item Rows of $(\bX-\lambda\bI)$ are linearly dependent;
\item $(\bX-\lambda\bI)$ has rank $r<p$;
\item Dimension of column space = dimension of row space = $r<p$;
\item $(\bX-\lambda\bI)^\top(\bX-\lambda\bI)$ is symmetric semidefinite;
\item $(\bX-\lambda\bI)$ has $r<p$ nonzero (positive) singular values;
\item Zero is an eigenvalue of $(\bX-\lambda\bI)$.
\end{itemize}
\end{remark}

\index{Vector norm}
\index{Norm}
\subsection{Vector Norms, Dual Norm, Balls}\label{section:vec_norms}

The concept of a \textit{norm} is fundamental for measuring the magnitude of vectors and, consequently, enables the definition of metrics on normed linear spaces. Norms provide a quantitative measure of the ``size" of a vector or matrix, which is essential in many applications---such as computing the length of a vector in Euclidean space or assessing the scale of a matrix in multidimensional settings.

Moreover, norms allow us to define distances between vectors or matrices. Specifically, the distance between two vectors $\bu$ and $\bv$ can be computed as the norm of their difference, $\norm{\bu-\bv}$. This notion is crucial in tasks that rely on proximity measures, including clustering algorithms in machine learning.\footnote{We restrict our discussion to norms (and inner products) on real vector and matrix spaces. Most results extend directly to the complex case.}

For any vector or matrix, a norm must satisfy the following three properties; see, for example, \citet{lu2021numerical} for a more detailed introduction.
\begin{definition}[Vector norm and matrix nrom\index{Matrix norm}\index{Vector norm}]\label{definition:matrix-norm}
Let $\norm{\cdot}$ denote a norm defined on either vectors or matrices. 
Then, for any matrix $\bX \in \real^{n\times p}$ and any vector $\bbeta \in \real^p$, the following properties must hold:
\begin{itemize}
\item \textit{Nonnegativity}. $\norm{\bX} \geq 0$ and $\norm{\bbeta}\geq 0$, with equality if and only if $\bX=\bzero $ or $\bbeta=\bzero$, respectively. 
\item \textit{Positive homogeneity}. $\norm{\lambda \bX} = \abs{\lambda} \cdot \norm{\bX}$ and  $\norm{\lambda \bbeta} = \abs{\lambda} \cdot \norm{\bbeta}$ for any scale $\lambda \in \real$.
\item \textit{Triangle inequality}. $\norm{\bX+\bY} \leq \norm{\bX}+\norm{\bY}$ for all matrices $\bX, \bY\in \real^{n\times p}$, and $\norm{\balpha+\bbeta} \leq \norm{\balpha}+\norm{\bbeta}$ for all vectors $\balpha,\bbeta\in \real^p$.
\end{itemize}
\end{definition}

Building on this definition, we now introduce specific types of vector norms: the $\ell_1$, $\ell_2$, and $\ell_\infty$-norms for a vector.
\begin{definition}[Vector $\ell_1, \ell_2, \ell_\infty$, $\ell_s$-norms]\label{definition:vec_l2_norm}
For a vector $\bbeta\in\real^p$, the \textit{$\ell_2$ vector norm} is defined as $\normtwo{\bbeta} = \sqrt{\beta_1^2+\beta_2^2+\ldots+\beta_p^2}$.
Similarly, the \textit{$\ell_1$-norm} can be obtained by 
$
\norm{\bbeta}_1 = \sum_{i=1}^{p} \abs{\beta_i} .
$
And the \textit{$\ell_\infty$-norm} can be obtained by 
$
\norm{\bbeta}_\infty = \mathop{\max}_{i=1,2,\ldots,p} \abs{\beta_i} .
$
More generally, the $\ell_s$-norm is defined as $\norms{\bbeta}=\sqrt[s]{ \sum_{i=1}^{p}\abs{\beta_i}^s  }$ for $s\geq 1$.~\footnote{In this book, we use the term ``$\ell_s$-norm" instead of ``$\ell_p$-norm" to avoid confusion with the dimension $p$.}
\end{definition}

\begin{exercise}[Standard bounds on vector norms\index{Standard bounds on vector norms}]\label{exercise:cauch_sc_l1l2}
Let $\bbeta\in\real^p$. Use the Cauchy--Schwarz inequality (see Theorem~\ref{theorem:cs_matvec}), prove the following relationships:
$$
\begin{aligned}
\norminf{\bbeta} &\leq \normone{\bbeta} \leq p\norminf{\bbeta}; \\
\norminf{\bbeta} &\leq \normtwo{\bbeta} \leq \sqrt{p}\norminf{\bbeta}; \\
\normtwo{\bbeta} &\leq \normone{\bbeta} \leq \sqrt{p}\normtwo{\bbeta}. \\
\end{aligned}
$$
\end{exercise}
This above exercise  shows the equivalence of the $\ell_1, \ell_2$, and $\ell_\infty$ vector norms  on finite-dimensional vector spaces.
Given positive constants $c,C>0$, for two norms $\norma{\cdot}$ and $\normb{\cdot}$, the inequality $c\cdot \norma{\bbeta}\leq \normb{\bbeta} \leq C\cdot \norma{\bbeta}$ for all vectors $\bbeta\in\real^p$
 means that the norms induce the same notion of ``closeness" or topology---they define the same convergent sequences, open sets, continuity, etc.

\subsection*{Dual Norm}
\index{Dual norm}
\index{H\"older's inequality}

Consider the $\ell_s$ vector norm. By \textit{\holders inequality} (see Theorem~\ref{theorem:holder-inequality}), we have
$
\bbeta^\top\balpha \leq \norm{\bbeta}_s \norm{\balpha}_t,
$
where $s,t>1$ satisfy $\frac{1}{s}+\frac{1}{t}=1$, and $\balpha,\bbeta\in \real^p$. Equality holds if the two sequences $\{\abs{\beta_i}^s\}$ and $\{\abs{\alpha_i}^t\}$ are linearly dependent. This implies
\begin{equation}\label{equation:dual_norm_equa}
\mathop{\max}_{\norm{\balpha}_t=1} \bbeta^\top\balpha = \norm{\bbeta}_s.
\end{equation}
For this reason, $\norm{\cdot}_t$ is called the \textit{dual norm} of $\norm{\cdot}_s$.
On the other hand, for each $\bbeta\in \real^p$ with $\norm{\bbeta}_s=1$, there exists a vector $\balpha\in \real^p$ such that $\norm{\balpha}_t=1$ and $\bbeta^\top\balpha=1$.
Notably, the $\ell_2$-norm is \textit{self-dual}, while the $\ell_1$- and $\ell_\infty$-norms are dual to each other.

\begin{definition}[Set of primal counterparts]\label{definition:set_primal}
Let $ \norm{\cdot} $ be any norm on $ \real^p $. Then the \textit{set of primal counterparts of $\ba\in\real^p$} is defined as 
\begin{equation}
\Lambda_{\ba} \triangleq  \argmax_{\bu \in \real^p} \left\{ \innerproduct{\ba, \bu} \mid  \norm{\bu} \leq 1 \right\}.
\end{equation}
That is, $\innerproduct{\ba, \ba^\dagger} = \norm{\ba}_*$ for any $\ba^\dagger\in \Lambda_{\ba}$, where $\norm{\cdot}_{*}$ denotes the dual norm of $\norm{\cdot}$.
It follows that 
\begin{enumerate}[(i)]
\item If $\ba \neq \bzero$, then $\norm{\ba^\dagger} = 1$ for any $\ba^\dagger \in \Lambda_{\ba}$.
\item If $\ba = \bzero$, then $\Lambda_{\ba} = \{\bbeta\in\real^p\mid  \norm{\bbeta} \leq 1\}$.
\end{enumerate}
\end{definition}

\begin{example}[Set of primal counterparts]\label{example:set_primal_count}
Below are a few illustrative examples of sets of primal counterparts:
\begin{itemize}
\item If the norm is the  $\ell_2$-norm, then for any $\ba \neq \bzero$,
$
\Lambda_{\ba} = \left\{ {\ba}/{\normtwo{\ba}} \right\}.
$

\item If the norm is the  $\ell_1$-norm, then for any $\ba \neq \bzero$,
$$
\Lambda_{\ba} = \left\{ \sum_{i \in \sI(\ba)} \lambda_i \sign(a_i) \be_i \mid \sum_{i \in \sI(\ba)} \lambda_i = 1, \lambda_j \geq 0, j \in \sI(\ba) \right\},
$$
where $\sI(\ba) \triangleq \argmax_{i=1,2,\ldots,p} \abs{a_i}$.

\item If the norm is the  $\ell_\infty$-norm, then for any $\ba \neq \bzero$,
$$
\Lambda_{\ba} = \left\{ \bbeta \in \real^p \mid \beta_i = \sign(a_i), i \in \sI_{\neq}(\ba), \abs{\beta_j} \leq 1, j \in \sI_0(\ba) \right\},
$$
where
$
\sI_{\neq}(\ba) \triangleq \left\{ i \in \{1, 2, \ldots, p\} \mid a_i \neq 0 \right\} $ and $ \sI_0(\ba) \triangleq \left\{ i \in \{1, 2, \ldots, p\} \mid a_i = 0 \right\}.
$
\end{itemize}
These examples play a crucial role in the development of non-Euclidean gradient descent methods \citet{lu2025practical}.
\end{example}

\index{$k$-sparse vectors}
\subsection*{Balls}
Given a specific norm, we introduce the concepts of an open ball and a closed ball as follows:
\begin{definition}[Open ball, closed ball\index{$\ell_s$-ball}]\label{definition:open_closed_ball}
Let $\norm{\cdot}_s: \real^p\rightarrow \real_+$ denote the $\ell_s$-norm. 
The \textit{open ball} centered at $\bc\in\real^p$ with radius $r>0$  is defined as
$$
\sB_s(\bc, r) \triangleq \{\bbeta\in\real^p\mid  \norm{\bbeta-\bc}_s <r\}.
$$
Similarly, the \textit{closed ball} centered at $\bc\in\real^p$ with radius $r>0$  is defined as 
$$
\sB_s[\bc,r] \triangleq \{\bbeta\in\real^p\mid \norm{\bbeta-\bc}_s \leq r\}.
$$
\end{definition}
For example, $\sB_2[\bzero,1]$ represents  the \textit{unit closed ball} w.r.t. to  the $\ell_2$-norm.
To simplify notation, we omit the subscript ``2" for  $\ell_2$-norms and the center\ $\bzero$ for balls centered at the origin, e.g., 
\begin{equation}
\sB[1] \triangleq \sB_2[\bzero,1] .
\end{equation}
As a special case, the notation 
\begin{equation}
\sB_0[k] \triangleq\sB_0[\bzero, k]
\end{equation}
is used to  denote the set of \textit{$k$-sparse vectors}, i.e., containing vectors in $\real^p$ that have at most $k$  nonzero entries.
More generally, let $\norm{\cdot}$ be any norm, the corresponding open and closed balls are denoted by
\begin{equation}
\sB_{\norm{\cdot}} (\bc, r)
\qquad \text{and}\qquad 
\sB_{\norm{\cdot}} [\bc, r].
\end{equation}

\index{Matrix norm}
\index{Submultiplicativity}
\subsection{Matrix Norms: Frobenius and Spectral Norms}
While vector norms measure the ``size" or ``length" of a vector, matrix norms extend this concept to matrices---but with the crucial distinction that they must also account for how a matrix acts as a linear operator, often capturing its effect on vectors through properties like scaling or stretching.
For a matrix $\bX\in\real^{n\times p}$, we  define the (matrix) Frobenius norm as follows.
\begin{definition}[Matrix Frobenius norm\index{Frobenius norm}]\label{definition:frobernius-in-svd}
The \textit{Frobenius norm} of a matrix $\bX\in \real^{n\times p}$ is defined as 
$$
\normf{\bX} = \sqrt{\sum_{i=1,j=1}^{n,p} x_{ij}^2}=\sqrt{\trace(\bX\bX^\top)}=\sqrt{\trace(\bX^\top\bX)} = \sqrt{\sigma_1^2+\sigma_2^2+\ldots+\sigma_r^2}, 
$$
where $\sigma_1, \sigma_2, \ldots, \sigma_r$ are nonzero singular values of $\bX$, and  $\trace(\bX^\top\bX)$ denotes the trace of $\bX^\top\bX$, i.e., sum of diagonal elements of $\bX^\top\bX$.
\end{definition}

The spectral norm is defined as follows.
\begin{definition}[Matrix spectral norm]\label{definition:spectral_norm}
The \textit{spectral norm} of a matrix $\bX\in \real^{n\times p}$ is defined as 
\begin{equation}\label{equation:spectral_norm_eq1}
\normtwo{\bX} = \mathop{\max}_{\bbeta\neq\bzero} \frac{\normtwo{\bX\bbeta}}{\normtwo{\bbeta}}  =\mathop{\max}_{\bu\in \real^p: \norm{\bu}_2=1}  \normtwo{\bX\bu} .
\end{equation}
~\footnote{Note that in this book we will not use the sup/inf notation but rather use only the min/max notation, where the usage of this notation does not imply that the maximum or minimum is actually attained.}
This quantity equals the largest singular value of  $\bX$, i.e., $\normtwo{\bX} = \sigma_{\max}(\bX)$; see Section~\ref{section:SVD}.
\end{definition}
If $\bX$ is \textbf{symmetric}, then its spectral norm also equals the maximum absolute eigenvalue:
\begin{equation}\label{equation:spectral_norm_eq2}
\normtwo{\bX}  = \max_{\bu\in \real^p: \norm{\bu}_2=1}\abs{\bu^\top\bX\bu}, \quad \text{if $\bX=\bX^\top$}.
\end{equation}

\begin{remark}[Standard bounds on matrix norms]\label{remark:stad_matnorm}
Analogous  to  standard bounds on vector norms (Exercise~\ref{exercise:cauch_sc_l1l2}), the spectral norm is bounded above by the Frobenius norm
\begin{equation}
\normtwo{\bX} \leq \normf{\bX}.
\end{equation}
This follows simply by the definition $\normtwo{\bX} = \sigma_{\max}(\bX)$ and $\normf{\bX}=\sqrt{\sigma_1^2+\sigma_2^2+\ldots+\sigma_r^2}$.
Alternatively, this can be proved using
$$
\normtwo{\bX\bbeta}^2 
= \sum_{i=1}^{n} \left(\sum_{j=1}^{p}  x_{ij}\beta_j\right)^2
\leq \sum_{i=1}^{n} \left(\sum_{j=1}^{p} x_{ij}^2\right) \left(\sum_{j=1}^{p} \beta_j^2\right)
=\normf{\bX}\normtwo{\bbeta},
$$
where the inequality follows from the Cauchy--Schwarz inequality; see Theorem~\ref{theorem:cs_matvec}.
\end{remark}


We note that the Frobenius norm serves as the matrix counterpart of the vector $\ell_2$-norm.
For simplicity, we e omit the full subscript for the vector $\ell_2$-norm or Frobenius norm when the context makes it clear which norm is intended: 
$$
\norm{\bX}=\normf{\bX}
\qquad \text{and}\qquad
\norm{\bx}=\normtwo{\bx}.
$$
However, for the spectral norm, the subscript in $\normtwo{\bX}$ should \textbf{not} be omitted.

The vector space $\real^p$, equipped with a given norm $\norm{\cdot}$, is called a \textit{normed vector space}.
On the other hand, another way to define norms for matrices is by treating a matrix $\bX\in\real^{n\times p}$ as a vector in $\real^{np}$, e.g., via vectorization.
What distinguishes a matrix norm from a generic vector norm on matrices is a key property called \textit{submultiplicativity}: $\norm{\bX\bY}\leq \norm{\bX}\norm{\bY}$ if $\norm{\cdot}$ is a submultiplicative matrix norm (see discussions below). 

In some texts, a matrix norm that is not \textit{submultiplicative} is referred to as a \textit{vector norm on matrices} or a \textit{generalized matrix norm}.
Although the definition of a matrix norm applies to both square and rectangular matrices, submultiplicativity is especially important in the analysis of square matrices.
For a submultiplicative matrix norm $\norm{\cdot}$ satisfying $\norm{\bX\bY}\leq \norm{\bX}\norm{\bY}$, and for a square matrix $\bX\in\real^{n\times n}$, it follows by induction that
\begin{equation}\label{equation:power_subm}
\norm{\bX^2} \leq \norm{\bX}^2
\quad\implies\quad
\norm{\bX^k} \leq \norm{\bX}^k, \ \forall\, k\in\{1,2,\ldots,\}.
\end{equation}
Therefore, if $\bX$ is idempotent (i.e., $\bX^2=\bX$), then $\norm{\bX}\geq 1$.
In particular, this implies
\begin{equation}\label{equation:power_subm2}
\norm{\bI}\geq 1,
\quad\text{if}\quad
\norm{\cdot} \text{ is submultiplicative}.
\end{equation}
On the other hand, if $\bX$ is nonsingular, then using submultiplicativity we obtain:
$$
1\leq \norm{\bI}=\norm{\bX\bX^{-1}}\leq \norm{\bX}\norm{\bX^{-1}}.
$$  
That is, a submultiplicative norm satisfies  $\norm{\bI}\geq 1$, and  it is said to be \textit{normalized} if and only if $\norm{\bI}=1$.

\index{Submultiplicativity}
\index{Orthogonally invariance}
\begin{theoremHigh}[Submultiplicativity and orthogonally invariance]\label{theorem:submul_ortho_matnorm}
The Frobnenius and spectral norms are submultiplicative. That is, 
$$
\norm{\bX\bY}_F\leq \norm{\bX}_F\norm{\bY}_F
\qquad \text{and}\qquad 
\norm{\bX\bY}_2\leq \norm{\bX}_2\norm{\bY}_2.
$$
Moreover, both norms are \textit{orthogonally invariant}. That is, for any orthogonal matrices $\bU\in \real^{n\times n}$ and $\bV\in \real^{p\times p}$, and any $\bX\in\real^{n\times p}$,
$$
\norm{\bX }_F =  \norm{\bU\bX\bV }_F
\qquad \text{and}\qquad 
\norm{\bX }_2 =  \norm{\bU\bX\bV }_2.
$$
\end{theoremHigh}
\begin{proof}
See \citet{lu2021numerical}.
\end{proof}

\subsection*{Induced Matrix Norm: Generalizing Spectral Norm}
More generally, many matrix norms can be derived using the concept of induced norms.

\index{Induced norm}
\begin{definition}[Induced matrix norm: General matrix norm]\label{definition:induced_norm_app}
Let $\bX\in \real^{n\times p}$ be a matrix, and let $\norma{\cdot}$ and $\normb{\cdot}$ be vector norms on $\real^p$ and $\real^n$ respectively.
The \textit{induced matrix norm} $\norm{\bX}_{a,b}$ is defined as
$$
\norm{\bX}_{a,b} =  \mathop{\max}_{\norma{\bbeta} \neq  0 } \frac{\normb{\bX\bbeta}}{\norma{\bbeta}}  
=  \mathop{\max}_{\norma{\bbeta} = 1} \normb{\bX\bbeta} .~\footnote{Although the norms  $\norma{\cdot}$ and $\normb{\cdot}$ can be chosen from a variety of options, $\ell_s$-norms are most commonly discussed in this context. 
In such cases, the induced norm is sometimes referred to as the \textit{\holders induced matrix norm}.}
$$
From this definition, the \textit{matrix-vector product inequality} follows immediately:
\begin{equation}\label{equation:induced_ineqy_intern}
\normb{\bX\bbeta} \leq \norm{\bX}_{a,b} \cdot \norma{\bbeta}.
\end{equation}
The induced matrix norm can also be referred to as the $(a,b)$-norm. 
When $a=b$, it is simply called the \textit{$a$-norm}, and we write  $\norm{\bX}_{a}$ instead of $\norm{\bX}_{a,a}$. 
In particular, when $a=b=2$, we have the classic spectral norm. 
\end{definition}
Intuitively, the induced norm measures the maximum factor by which the matrix $\bX$ can stretch a nonzero vector: the input vector $\bbeta$ is measured using the  norm $\norma{\cdot}$, while the output vector $\bX\bbeta$ is measured using norm $\normb{\cdot}$.

From the definition above, it follows that the spectral norm is a special case of an induced norm with $a=b=2$, i.e.,  $\normtwo{\bX} = \norm{\bX}_{2,2}$. Similarly, the \textit{matrix 1-norm} is obtained as
\begin{equation}\label{equation:mat_one_norm}
\textbf{Matrix 1-norm:}
\quad\normone{\bX} = \mathop{\max}_{j=1,2,\ldots,p} \sum_{i=1}^{n}\abs{x_{ij}}
=\mathop{\max}_{j=1,2,\ldots,p} \normone{\bx_j}
,
\end{equation}
which is also known as  the \textit{maximum absolute column sum norm}. 
To verify this, consider $\bX\in\real^{n\times p}$. 
For any vector $\bbeta\in\real^p$, we have 
$$
\normone{\bX\bbeta}=\normone{\sum_{i=1}^{p} \beta_i\bx_i}
\leq 
\sum_{i=1}^{p} \abs{\beta_i}\normone{\bx_i}
\leq 
\sum_{i=1}^{p} \abs{\beta_i} \left(\mathop{\max}_{j=1,2,\ldots,p}\normone{\bx_j} \right)
=
\normone{\bbeta} \normone{\bX}.
$$
Therefore, $ \mathop{\max}_{\normone{\bbeta} = 1} \normone{\bX\bbeta} \leq \normone{\bX}$.
On the other hand, consider $\mathop{\max}_{\normone{\bbeta} = 1}\normone{\bX\bbeta}$ for $\bbeta=\be_j,\,\forall\, j\in\{1,2,\ldots,p\}$, we have
$
\mathop{\max}_{\normone{\bbeta} = 1}\normone{\bX\bbeta}\geq \mathopmax{j=1,2,\ldots,p}\normone{\bx_j}=\normone{\bX}.
$
By ``sandwiching," we have $\mathop{\max}_{\normone{\bbeta} = 1}\normone{\bX\bbeta}=\normone{\bX}$.
Likewise, the \textit{matrix $\infty$-norm} is given by
\begin{equation}\label{equation:mat_inf_norm}
\textbf{Matrix $\infty$-norm:}
\quad \norm{\bX}_\infty = \mathop{\max}_{i=1,2,\ldots,n} \sum_{j=1}^{p}\abs{x_{ij}} ,
\end{equation}
which is also called the \textit{maximum absolute row sum norm}.

\begin{exercise}[Maximum absolute row sum norm]
Prove that $\norminf{\bX}$ is the matrix norm induced by the vector $\ell_\infty$-norm.
\end{exercise}
\index{Equi-induced matrix norm}
\index{Normalized matrix norm}
\begin{exercise}[Equi-induced and normalized norm]
Let $a=b$ and let $\bX\in\real^{p\times p}$ be a square matrix in Definition~\ref{definition:induced_norm_app}. 
In this case, the induced  norm is  called an \textit{equi-induced norm}. Show that the equi-induced norm is normalized: $\norma{\bI_p}=1$ and $\norma{\bX}\geq 1$.
\end{exercise}

We have already seen that the spectral norm---an important example of an induced norm---is submultiplicative. In fact, all induced norms are submultiplicative.
\index{Submultiplicativity}
\begin{theoremHigh}[Submultiplicativity of induced norms]\label{theorem:Submultiplicativity_induced}
For $a,b,c\geq 1$ and matrices $\bX,\bY$ (provided the matrix product $\bX\bY$ is defined), we have 
$$
\norm{\bX\bY}_{a,b} \leq \norm{\bX}_{c,b}\norm{\bY}_{a,c}.
$$
That is, the norm of a
product is bounded by the product of the norms.
\end{theoremHigh}
\begin{proof}[of Theorem~\ref{theorem:Submultiplicativity_induced}]
From the definition of the induced norm, we have
$$
\begin{aligned}
\norm{\bX\bY}_{a,b} &= \mathop{\max}_{\norm{\bbeta}_a=1} \norm{\bX\bY\bbeta}_b
\leq \mathop{\max}_{\norm{\bbeta}_a=1}  \norm{\bX}_{c, b} \norm{\bY\bbeta}_{c} 
\leq \mathop{\max}_{\norm{\bbeta}_a=1}\norm{\bX}_{c,b} \norm{\bY}_{a,c} \norm{\bbeta}_a .
\end{aligned}
$$
This completes the proof.
\end{proof}

\begin{remark}[Condition number\index{Condition number}]\label{remark:cond_number}
For a nonsingular matrix $\bX\in\real^{p\times p}$,
we note that $1 = \norm{\bI}_{b,b} \leq \norm{\bX}_{a,b} \norm{\bX^{-1}}_{b,a}$. The value $\kappa_{a,b}=\norm{\bX}_{a,b} \norm{\bX^{-1}}_{b,a}$ is normally known as the \textit{condition number} with respect to the induced vector norms $\norm{\cdot}_a$ and $\norm{\cdot}_b$. 
The common case is 
$$
\kappa (\bX) = {\normtwo{\bX}}{\normtwo{\bX^{-1}}} = \frac{\sigma_{\max}(\bX)}{\sigma_{\min}(\bX)},
$$
where $\sigma_{\max}$ and $\sigma_{\min}$ denote the largest and smallest singular values of $\bX$, respectively (Section~\ref{section:SVD}).
See \citet{lu2021numerical} for more details.
\end{remark}

\subsection{Orthogonal and Projection Matrices}\label{section:ortho_proj_mat}

We often use the letters $\bQ, \bU, \bV \in \real^{n\times p}$ to denote matrices whose columns are orthonormal (called \textit{semi-orthogonal matrices} when $p\neq n$). 
If such a matrix is square (i.e., $p=n$), it is called an orthogonal matrix. 
Orthogonal matrices represent linear transformations that preserve the length (magnitude) of any vector---they may change its direction, but not its norm.

\begin{definition}[Orthogonal matrix, semi-orthogonal matrix\index{Orthogonal matrix}\index{Semi-orthogonal matrix}]\label{definition:orthogn_mat}
A real square matrix $\bQ\in \real^{n\times n}$ is an \textit{orthogonal} matrix if its inverse  equals its transpose; 
that is $\bQ^{-1}=\bQ^\top$ and $\bQ\bQ^\top = \bQ^\top\bQ = \bI$. 
In other words, suppose $\bQ=[\bq_1, \bq_2, \ldots, \bq_n]$, where $\bq_i \in \real^n$ for all $i \in \{1, 2, \ldots, n\}$, then $\bq_i^\top \bq_j = \delta(i,j)$ with $\delta(i,j)$ being the Kronecker delta function. 
For any vector $\bx$, an orthogonal matrix  preserves  the Euclidean length: $\normtwo{\bQ\bx} = \normtwo{\bx}$, where $\normtwo{\cdot}$ denotes the $\ell_2$-norm (Definition~\ref{definition:vec_l2_norm}).
Note that because the columns of an orthogonal matrix are unit vectors and mutually orthogonal, they form an orthonormal set. However, by historical convention, the term \textit{orthogonormal matrix} is \textbf{not} used; instead, such matrices are always referred to as \textit{orthogonal matrices}, even though their columns are orthonormal.

On the other hand, if $\bQ\in\real^{n\times \gamma}$ contains only $\gamma<n$ of these columns, then $\bQ^\top\bQ = \bI_\gamma$ stills holds, where $\bI_\gamma$ is the $\gamma\times \gamma$ identity matrix. But $\bQ\bQ^\top=\bI$ will not be true. 
In this case, $\bQ$ referred to as a \textit{semi-orthogonal matrix}.
\end{definition}

\index{Orthonormal basis}
\index{Orthogonal matrix}
\index{Orthogonal vs orthonormal}
The vectors $\bq_1, \bq_2, \ldots, \bq_\gamma\in \real^n$ are \textit{mutually orthogonal} if their dot products $\bq_i^\top\bq_j$ satisfies $\bq_i^\top\bq_j=0$  whenever $i \neq j$. 
If, in addition, each vector has unit length (i.e., $\normtwo{\bq_i}=1$ for all i), then the vectors are mutually orthonormal. We typically collect such orthonormal vectors as the columns of a matrix $\bQ$.
\begin{itemize}
\item When $n\neq \gamma$: the matrix $\bQ$ is easy to work with because $\bQ^\top\bQ=\bI \in \real^{\gamma\times \gamma}$.

\item When $n= \gamma$: the matrix $\bQ$ is square, $\bQ^\top\bQ=\bI$ means that $\bQ^\top=\bQ^{-1}$, i.e., the transpose of $\bQ$ is the inverse of $\bQ$. Then we also have $\bQ\bQ^\top=\bI$, i.e., $\bQ^\top$ is the two-sided inverse of $\bQ$. 
Such a matrix $\bQ$ is called an \textit{orthogonal matrix}. 
\end{itemize}

By definition, the columns $\bq_1, \bq_2, \ldots, \bq_n$ of any $n \times n$ orthogonal matrix have unit norm and are orthogonal to each other, so they form an orthonormal basis. Applying $\bQ^\top$ to a vector $\bx \in \real^n$ is equivalent to computing the coefficients of its representation in the basis formed by the columns of $\bQ$. Applying $\bQ$ to $\bQ^\top\bx$ recovers $\bx$ by scaling each basis vector with the corresponding coefficient:
$$
\bx = \bQ \bQ^\top \bx = \sum_{i=1}^n \langle \bq_i, \bx \rangle \bq_i.
$$

Since orthogonal matrices represent pure rotations (and possibly reflections), it is intuitive that the product of two orthogonal matrices is also orthogonal.
Indeed, if $\bQ, \bU \in \real^{n\times n}$ are orthogonal matrices, then their product $\bQ\bU$ is also an orthogonal matrix.

Specifically, the following lemma shows that orthogonal matrices preserve the $\ell_2$-norm of any vector.
\begin{lemma}\label{lemma:or_pref_len}
Let $\bQ \in \real^{n\times n}$ be an orthogonal matrix. 
Then, for any vector $\bx \in \real^n$,
$$
\normtwo{\bQ\bx} =\normtwo{\bx}.
$$
\end{lemma}

\begin{proof}[of Lemma~\ref{lemma:or_pref_len}]
By the definition of an orthogonal matrix, we have 
$\normtwo{\bQ\bx}^2 
= \bx^\top \bQ^\top \bQ \bx 
= \bx^\top \bx 
= \normtwo{\bx}^2$.
Taking square roots on both sides yields the desired result.
\end{proof}

Semi-orthogonal matrices are closely related to the concept of \textit{orthogonal projection matrices}.
Formally, we define the {projection matrix} as follows:
\begin{definition}[Projection matrix\index{Projection matrix}]\label{definition:projection-matrix}
A matrix $\bH\in \real^{n\times n}$ is called a \textit{projection matrix} or \textit{projector} onto a subspace $\mathcalV \in \real^n$ if and only if $\bH$ satisfies the following properties:
\begin{itemize}
\item (P1). $\bH\by \in \mathcalV$ for all $\by \in \real^n$: Any vector can be projected onto the subspace $\mathcalV$.
\item (P2). $\bH\by =\by$ for all $\by \in \mathcalV$: Vectors already in $\mathcalV$ remain unchanged under the projection.
\item (P3). $\bH^2 = \bH$: Applying the projection twice has the same effect as applying it once. This property is known as {idempotence}.
\end{itemize}
\end{definition}

Since we project vectors from $\real^n$ onto a subspace of $\real^n$, any projection matrix must be square. 
If $\bH$ were $k\times n$ with $k\neq n$, the image would lie in $\real^k$, not in a subspace of $\real^n$.

Observe that  $\bH\by$ always lies in the column space of $\bH$, denoted $\cspace(\bH)$. 
One may then ask: what is the relationship between the target subspace $\mathcalV$ and the column space $\cspace(\bH)$. 
In fact, they are identical (see Lemma~\ref{lemma:projection-from-matrix} for a hint).

Suppose $\mathcalV=\cspace(\bH)$, and let $\by\in \mathcalV$. 
Then there exists some vector $\balpha$ such that $\by=\bH\balpha$. 
Using only condition (P3), we obtain
$$
\bH\by = \bH\bH\balpha = \bH\balpha = \by.
$$
Thus, condition (P3) implies both conditions (P1) and (P2). 
Consequently, a matrix $\bH$ is a projection matrix if and only if it is idempotent ($\bH^2=\bH$).

Intuitively, we often require more than just idempotence: we want the projected vector $\widehatby=\bH\by$ to be the closest point in $\mathcalV$ to the original vector $\by$, in the sense of Euclidean distance. This occurs precisely when the residual $\by-\widehatby$ is orthogonal to  $\widehatby$ (and, more generally, to the entire subspace $\mathcalV$). Such a projection is called an \textit{orthogonal projection}, and it aligns with the principle of least-squares error minimization.

\begin{figure}[h!]
\centering  
\vspace{-0.35cm} 
\subfigtopskip=2pt  
\subfigbottomskip=2pt  
\subfigcapskip=-5pt  
\subfigure[Orthogonal projection: project $\by$ to $\widehatby$.]{\label{fig:ls-geometric1}
\includegraphics[width=0.47\linewidth]{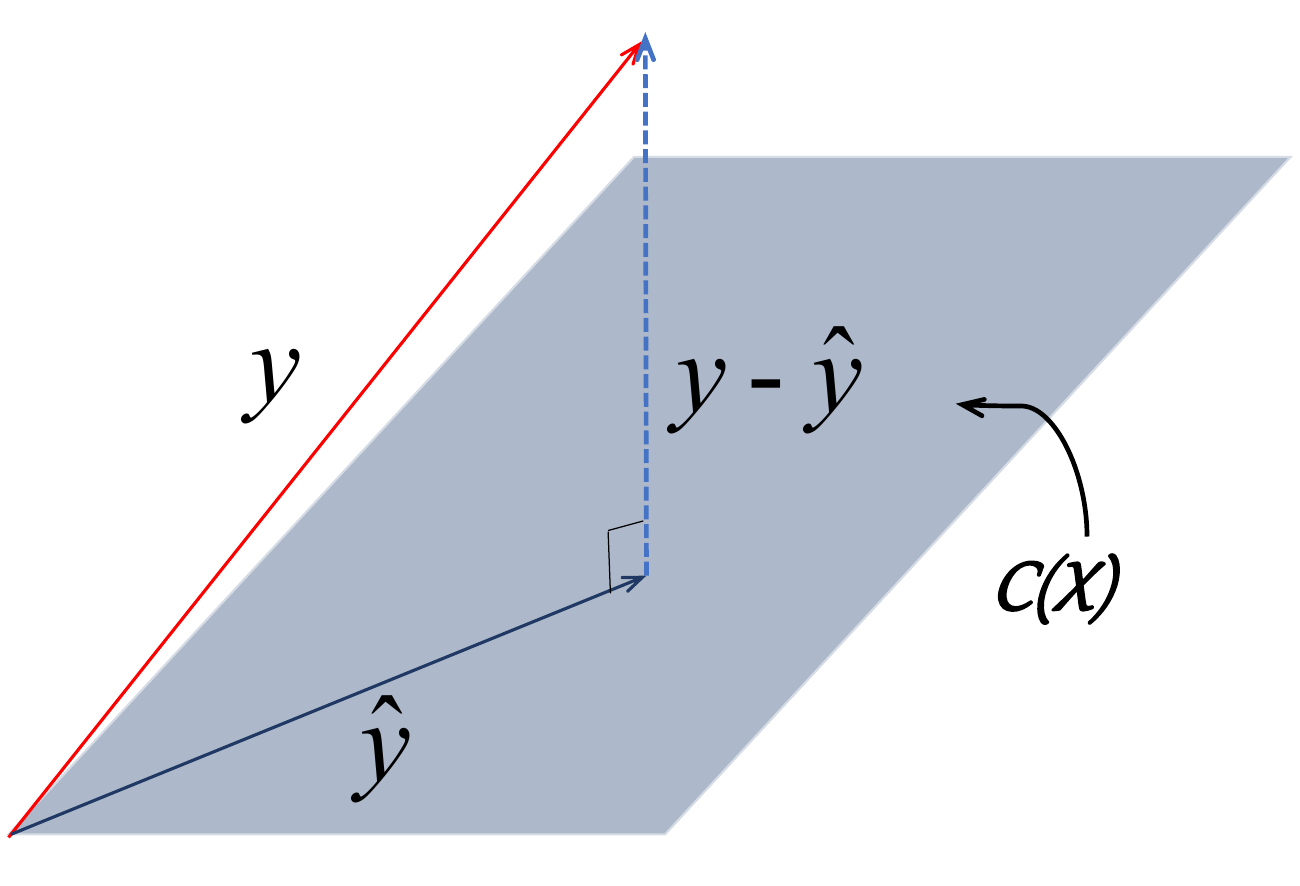}}
\quad 
\subfigure[Oblique projection: project $\by$ to $\widehatby_1$ or $\widehatby_2$.]{\label{fig:ls-geometric1-oblique}
\includegraphics[width=0.47\linewidth]{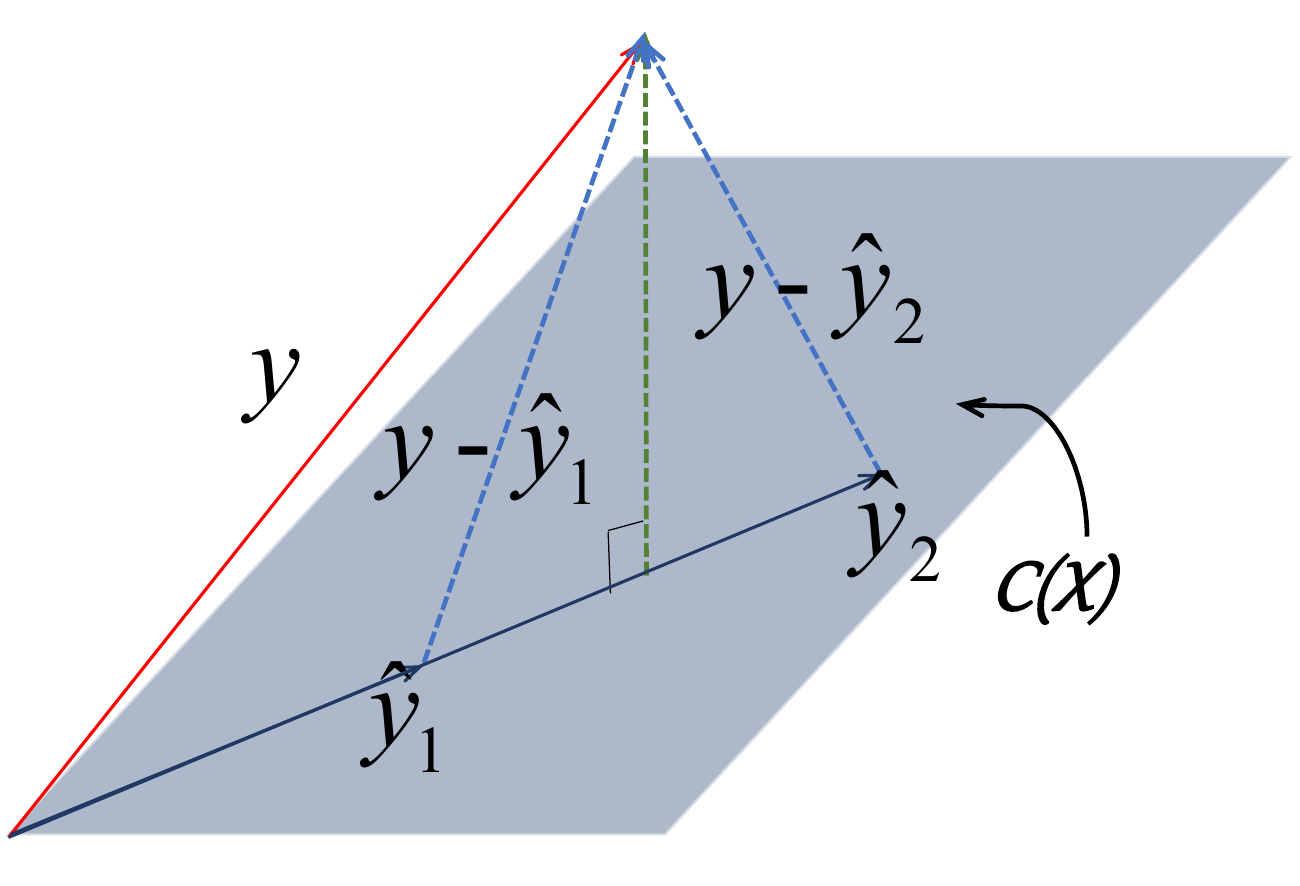}}
\caption{Projection onto the hyperplane  $\cspace(\bX)$, i.e., the column space of $\bX$.}
\label{fig:ls-geometric1-compare}
\end{figure}

\begin{definition}[Orthogonal and oblique projection matrix\index{Orthogonal projection}\index{Oblique projection}]\label{definition:orthogonal-projection-matrix}
A matrix $\bH$ is called an \textit{orthogonal projection matrix} or an \textit{orthogonal projector}  onto a subspace $\mathcalV \subseteq \real^n$ if and only if $\bH$ is a projection matrix, and the projection $\widehatby$ of any vector $\by\in \real^n$ is orthogonal to $\by - \widehatby$, i.e., $\bH$ projects onto $\mathcalV$ and along $\mathcalV^\perp$, the orthogonal complement of $\mathcalV$.

Otherwise, if $\widehatby$ is not orthogonal to $\by - \widehatby$, then the projection matrix is called an \textit{oblique projection matrix} or an \textit{oblique projector}. 
A geometric comparison between orthogonal and oblique projections is illustrated in Figure~\ref{fig:ls-geometric1-compare}.
\end{definition}

\begin{exercise}\label{exercise:sym_proj_mat} 
Show that a projection matrix $\bH$ is an orthogonal projection matrix if and only if it is symmetric.
\end{exercise}

Intuitively, let $\mathcalV$ be a subspace of $\real^n$, and let $\bH\in\real^{n\times n}$ be the orthogonal projection matrix onto $\mathcalV$. 
As shown in  Figure~\ref{fig:ls-geometric1-compare}, for any vector $\by\in\real^n$, 
the orthogonal projection minimizes the Euclidean distance to $\mathcalV$:
\begin{equation}
\normtwo{\by - \bH\by}^2 \leq \normtwo{\by - \bv}^2, \qquad \forall\, \bv \in \mathcalV.
\end{equation}

It is important to note that the term orthogonal projection does not imply that the projection matrix $\bH$ itself is an orthogonal matrix (in the sense of Definition~\ref{definition:orthogn_mat}, where $\bQ\bQ^\top=\bI$).
Rather, ``orthogonal" refers to the geometric orthogonality between the projected vector $\widehatby$ and the residua $\by - \widehatby$. 
Unless otherwise stated, the term ``projection matrix" in the remainder of this book will refer to an orthogonal projection matrix.

The following result shows that the orthogonal projection onto the column space of a set of linearly independent vectors admits a specific closed-form expression, which is intimately connected to the ordinary least squares (OLS) solution (see Section~\ref{section:optcd_ls}). 
\begin{lemma}[Projection matrix from a set of vectors]\label{lemma:projection-from-matrix}
Let $\bx_1, \bx_2, \ldots, \bx_p\in\real^n$ be linearly independent vectors such that $\cspace([\bx_1, \bx_2, \ldots, \bx_p]) = \mathcalV$, and assume $n\geq p$. 
Then the orthogonal projection onto the subspace $\mathcalV$ is given by
$$
\bH \triangleq \bX(\bX^\top\bX)^{-1}\bX^\top,
$$
where $\bX\in \real^{n\times p}$ is the matrix whose columns are  $\bx_1, \bx_2, \ldots, \bx_p$.
\end{lemma}
\begin{proof}[of Lemma~\ref{lemma:projection-from-matrix}]
It is straightforward  to verify that $\bH$ is symmetric and idempotent, which characterizes it as an orthogonal projection matrix. 
Considering the SVD of $\bX=\bU\bSigma\bV^\top$ (see Section~\ref{section:SVD}), we have $\bH = \bX(\bX^\top\bX)^{-1}\bX^\top = \bU\bSigma(\bSigma^\top\bSigma)^{-1} \bSigma^\top \bU^\top$. 
Let $\bU=[\bu_1, \bu_2, \ldots, \bu_n]$ be the column partition of $\bU$. From Theorem~\ref{theorem:svd-four-orthonormal-Basis}, $\{\bu_1, \bu_2, \ldots, \bu_p\}$ is an orthonormal basis of $\cspace(\bX)$. And $\bSigma(\bSigma^\top\bSigma)^{-1} \bSigma^\top$ in $\bH$ is an $n\times n$ matrix, where the upper-left part is a $p\times p$ identity matrix and the other parts are zero. 
Consequently, $\bH=\sum_{i=1}^{p}\bu_i\bu_i^\top$, which implies that $\cspace(\bH)=\spn\{\bu_1, \bu_2, \ldots, \bu_p\}=\cspace(\bX)$.
Thus, $\bH$ is the orthogonal projector onto $\cspace(\bX)=\mathcalV$, as required.
\end{proof}

Lemma~\ref{lemma:projection-from-matrix} shows that $\bH=\bX(\bX^\top\bX)^{-1}\bX^\top$ is the orthogonal projector onto $\cspace(\bX)$ if $\bX$ has full column rank. 
More generally, matrices with mutually orthonormal columns (semi-orthogonal matrices) can also be used to construct orthogonal projection matrices.
\begin{lemma}[Orthogonal projector onto general subspaces]\label{lemma:orthogo_genspa}
Let $ \mathcalV $ be  a subspace of $ \real^n $ with dimension $r$. Let $ \bQ_1\in\real^{n\times r} $ and $ \bQ_2 \in\real^{n\times (n-r)}$ be semi-orthogonal matrices (i.e., their columns are mutually orthonormal; Definition~\ref{definition:orthogn_mat}) such that $ \cspace(\bQ_1) = \mathcalV $ and $ \cspace(\bQ_2) = \mathcalV^{\perp} $, where $ \mathcalV^{\perp} $ denotes the orthogonal complement of $ \mathcalV $. Then the orthogonal projectors {onto} $ \mathcalV$ and $ \mathcalV^\perp$ are given by
\begin{equation}\label{equation:orthog_genspa1}
\bH_1 \triangleq \bQ_1  \bQ_1^\top
\qquad \text{and}\qquad
\bH_2 \triangleq \bQ_2\bQ_2^\top,
\end{equation}
respectively.
\end{lemma}
\begin{proof}[of Lemma~\ref{lemma:orthogo_genspa}]
We have $ \bH_1^2 = \bQ_1\bQ_1^\top\bQ_1\bQ_1^\top = \bQ_1\bQ_1^\top = \bH_1 $ since $\bQ_1^\top\bQ_1=\bI_r$. This shows that $ \bH_1 $ is a projector onto $ \mathcalV $. 
Since $\bH_1$ is symmetric, this completes the proof for the first part.
An identical argument applies to $\bH_2=\bQ_2\bQ_2^\top$, showing that it is the orthogonal projector onto $\mathcalV^{\perp}$.
\end{proof}

\subsection{Singular Value Decomposition (SVD)}\label{section:SVD}

We introduce singular value decomposition (SVD) of a matrix in this section.
Before presenting the general form of the SVD, we begin with the spectral decomposition of a symmetric matrix.
The spectral theorem---also known as the spectral decomposition for symmetric matrices---states that any real symmetric matrix has real eigenvalues and can be diagonalized using a real orthonormal basis.~\footnote{Note that the spectral decomposition for \textit{Hermitian matrices} similarly guarantees real eigenvalues and diagonalization via a complex orthonormal basis.}

\begin{theoremHigh}[Spectral decomposition\index{Spectral decomposition}\index{Spectral theorem}]\label{theorem:spectral_theorem}
A real matrix $\bX \in \real^{n\times n}$ is symmetric if and only if there exist an orthogonal matrix $\bQ$ and a diagonal matrix $\bLambda$ such that
\begin{equation*}
\bX = \bQ \bLambda \bQ^\top,
\end{equation*}
where the columns of $\bQ = [\bq_1, \bq_2, \ldots, \bq_n]$ are  mutually orthonormal eigenvectors of $\bX$, and the entries of $\bLambda=\diag(\lambda_1, \lambda_2, \ldots, \lambda_n)$ are the corresponding eigenvalues of $\bX$, which are real. 
In particular, the following properties hold:
\begin{enumerate}
\item All eigenvalues of a symmetric matrix are \textbf{real}.
\item The eigenvectors can be chosen to form an \textbf{orthonormal} set.
\item The rank of $\bX$ equals the number of its nonzero eigenvalues.
\item If all eigenvalues are distinct, then the corresponding eigenvectors are linearly independent.
\end{enumerate}
\end{theoremHigh}
\begin{proof}
See \citet{lu2022matrix}.
\end{proof}

Using spectral decomposition, we can factor a symmetric matrix into a diagonal form. 
However, this approach is not always applicable: when $\bX$ is not symmetric or not square, such a diagonalization is generally impossible.
The singular value decomposition (SVD) addresses this limitation. Instead of relying on a single orthogonal matrix of eigenvectors, the SVD expresses any matrix as a product of two orthogonal matrices and a diagonal (or diagonal-like) matrix of singular values. We now state the SVD theorem.

%
%
%
%

\begin{theoremHigh}[Full singular value decomposition\index{Singular value decomposition}]\label{theorem:full_svd_rectangular}\label{theorem:reduced_svd_rectangular}
Every real $n\times p$ matrix $\bX$ of rank $r$  admits a decomposition of the form
$$
\bX = \bU \bSigma \bV^\top,
$$ 
where $\bSigma\in \real^{n\times p}$ has the block structure $\bSigma=\footnotesize\begin{bmatrix}
\bSigma_r & \bzero \\
\bzero & \bzero
\end{bmatrix}$ with $\bSigma_r=\diag(\sigma_1, \sigma_2 \ldots, \sigma_r)\in \real^{r\times r}$, and the singular values satisfy $\sigma_1 \geq \sigma_2 \geq \ldots \geq \sigma_r$.
\begin{itemize}
\item The scalars  $\sigma_i$  are the nonzero \textit{singular values} of  $\bX$. Each $\sigma_i$ equals the positive square root of a nonzero eigenvalue of both $\bX^\top \bX$ and $\bX \bX^\top$.

\item $\bU\in \textcolor{black}{\real^{n\times n}}$ is an orthogonal matrix whose first $r$ columns are eigenvectors of $\bX \bX^\top$ corresponding to its $r$ nonzero eigenvalues; the remaining $n - r$ columns form an orthonormal basis for $\nspace(\bX^\top)$.

\item $\bV\in \textcolor{black}{\real^{p\times p}}$ is an orthogonal matrix whose first $r$ columns are eigenvectors of $\bX^\top \bX$ corresponding to its $r$ nonzero eigenvalues; the remaining $p - r$ columns form an orthonormal basis for $\nspace(\bX)$.

\item The columns of $\bU$ and $\bV$ are called the \textit{left and right singular vectors} of $\bX$, respectively. 
\end{itemize}
Moreover, the decomposition can be expressed as a sum of rank-one matrices: $ \bX = \bU \bSigma \bV^\top = \sum_{i=1}^r \sigma_i \bu_i \bv_i^\top$.
\end{theoremHigh}

\begin{figure}[h!]
\centering  
\vspace{-0.35cm}  
\subfigtopskip=2pt  
\subfigbottomskip=2pt  
\subfigcapskip=-5pt  
\subfigure[Reduced SVD decomposition.]{\label{fig:svdhalf}
\includegraphics[width=0.47\linewidth]{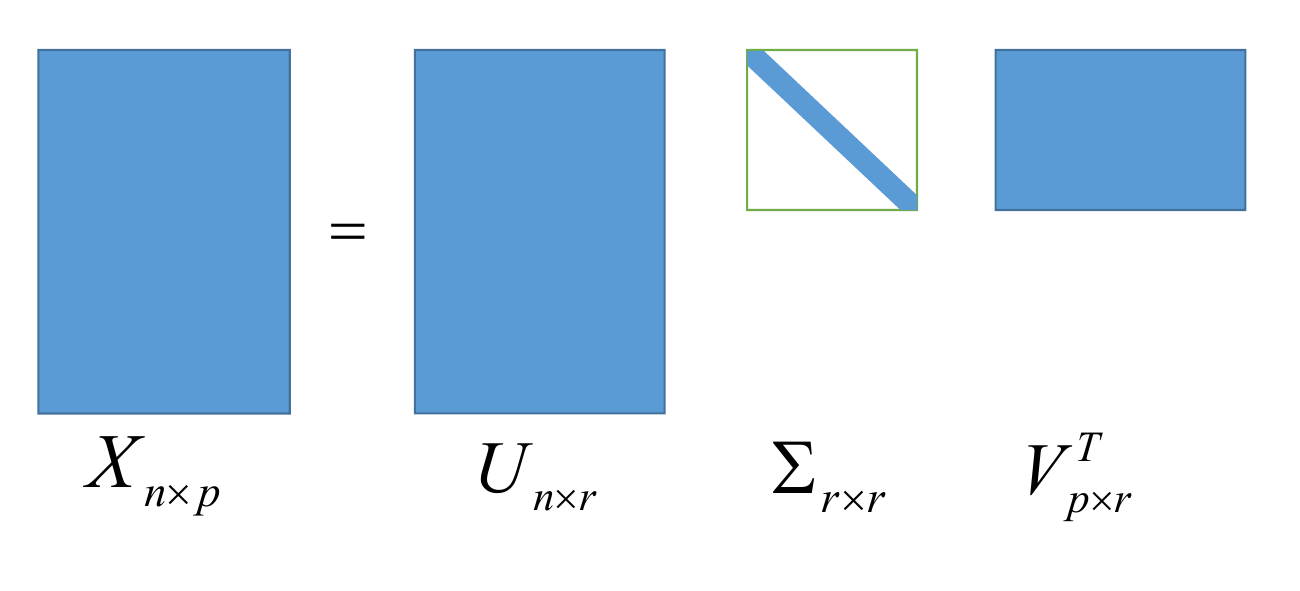}}
\quad 
\subfigure[Full SVD decomposition.]{\label{fig:svdall}
\includegraphics[width=0.47\linewidth]{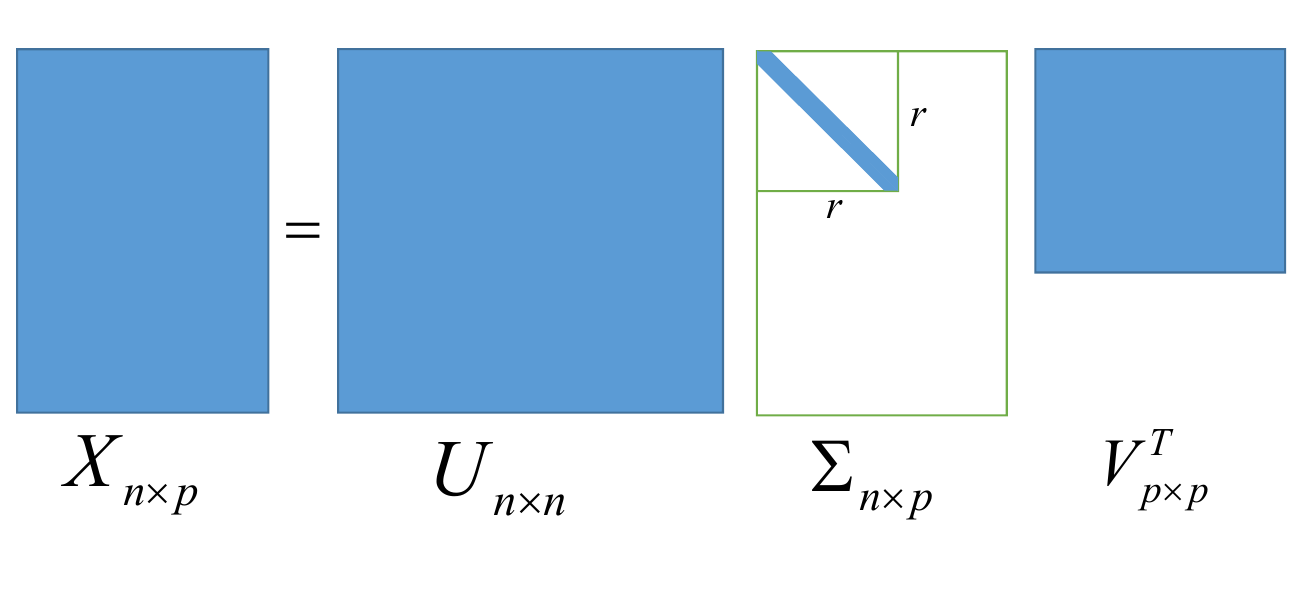}}
\caption{Comparison between the reduced and full SVD. White entries are zero, and blue entries are not necessarily zero}
\label{fig:svd-comparison}
\end{figure}

If $\bX$ has rank $r$, then its singular values satisfy $\sigma_1,  \sigma_2, \ldots, \sigma_r > 0$ while $\sigma_{r+1} = \sigma_{r+2} = \ldots = 0$. 
In many applications, it is more convenient to work with the \textit{reduced singular value decomposition (reduced SVD)}. 
For $\bX$ of rank $r$ with (full) SVD $\bX = \bU\bSigma \bV^\top$, we take the submatrices $\bU_r \in \real^{n\times r}$, $\bV_r \in \real^{p\times r}$ such that $\bU = [\bU_r, \bu_{r+1}, \ldots, \bu_n]$, $\bV = [\bV_r, \bv_{r+1}, \ldots,\bv_p]$, and $\bSigma_r = \diag([\sigma_1,\ldots,\sigma_r]) \in \real^{r\times r}$.
Therefore, we obtain the reduced singular value decomposition (reduced SVD) of $\bX$: $\bX=\bU_r\bSigma_r\bV_r^\top$.
The comparison between the reduced and full SVD is shown in Figure~\ref{fig:svd-comparison}, where white entries are zero, and blue entries are not necessarily zero.

Given $\bX\in\real^{n\times p}$ with reduced SVD $\bX = \bU_r\bSigma_r\bV_r^\top$, we observe that
\begin{align*}
\bX^\top\bX &= \bV_r\bSigma_r^\top\bU_r^\top\bU_r\bSigma_r\bV_r^\top = \bV_r\bSigma_r^2\bV_r^\top;\\
\bX\bX^\top &= \bU_r\bSigma_r\bV_r^\top\bV_r\bSigma_r^\top\bU_r^\top = \bU_r\bSigma_r^2\bU_r^\top.
\end{align*}
Thus, we recover the (reduced) spectral decompositions of symmetric positive semidefinite matrices $\bX^\top\bX$ and $\bX\bX^\top$, respectively. 
In particular, the singular values $\sigma_i = \sigma_i(\bX)$ satisfy
\begin{equation}\label{equation:sigbd_nearortho_eqp}
\sigma_i(\bX) = \sqrt{\lambda_i(\bX^\top\bX)} = \sqrt{\lambda_i(\bX\bX^\top)}, \quad i = 1,\ldots,\min\{n,p\}, 
\end{equation}
where $\lambda_1(\bX^\top\bX) \geq \lambda_2(\bX^\top\bX) \geq \ldots$ are the eigenvalues of $\bX^\top\bX$ in decreasing order. Moreover, the left and right singular vectors listed in $\bU,\bV$ can be obtained by the spectral decomposition of the positive semidefinite matrices $\bX^\top\bX$ and $\bX\bX^\top$.
Consequently, one can construct the SVD of $\bX$ directly from the spectral decompositions of $\bX^\top\bX$ and $\bX\bX^\top$---which also provides an alternative proof of the existence of the SVD.



\begin{exercise}[Proof of SVD]
	Using the discussion above and the spectral decomposition, prove the existence of the singular value decomposition of a matrix $\bX\in\real^{n\times p}$.
\end{exercise}

\subsection*{Four Orthonormal Bases in SVD}

\index{Fundamental theorem of linear algebra}
For any matrix $\bX\in\real^{n\times p}$, the following properties hold:
\begin{itemize}
\item The null space $\nspace(\bX)$ is the orthogonal complement of the row space $\cspace(\bX^\top)$ in $\real^p$: $\dim(\nspace(\bX))+\dim(\cspace(\bX^\top))=p$.

\item The left null space $\nspace(\bX^\top)$ is the orthogonal complement of the column space $\cspace(\bX)$ in $\real^n$: $\dim(\nspace(\bX^\top))+\dim(\cspace(\bX))=n$.
\end{itemize}
This result is known as  the fundamental theorem of linear algebra  (see the discussion preceding Lemma~\ref{lemma:equal-dimension-rank}).
In particular, the construction of the SVD provides orthonormal bases for all four fundamental subspaces described in this theorem. To establish this, we first require the following lemma.
\begin{lemma}[Subspace of $\bX^\top \bX$ and $\bX\bX^\top$]\label{lemma:rank-of-ttt}
Let $\bX\in \real^{n\times p}$ be given. 
Then,
\begin{itemize}
\item The column space of $\bX^\top \bX$ is identical  to the column space of $\bX^\top$ (i.e., row space of $\bX$): $\cspace(\bX^\top\bX)=\cspace(\bX^\top)$; this also shows $\nspace(\bX^\top\bX)=\nspace(\bX)$ by fundamental theorem of linear algebra (see the discussion preceding Lemma~\ref{lemma:equal-dimension-rank}).
\item The column space of $\bX\bX^\top$ is identical  to the column space of $\bX$: $\cspace(\bX\bX^\top)=\cspace(\bX)$. 
Hence, this also shows $\nspace(\bX\bX^\top)=\nspace(\bX^\top)$.
\end{itemize}
\end{lemma}
\begin{proof}[of Lemma~\ref{lemma:rank-of-ttt}]
Let $\bbeta\in \nspace(\bX)$, we have 
$
\bX\bbeta  = \bzero \implies \bX^\top\bX \bbeta =\bzero, 
$
i.e., $\bbeta\in \nspace(\bX) \implies \bbeta \in \nspace(\bX^\top \bX)$. Therefore, $\nspace(\bX) \subseteq \nspace(\bX^\top\bX)$. 
Furthermore, let $\bbeta \in \nspace(\bX^\top\bX)$, we have 
$$
\bX^\top \bX\bbeta = \bzero\implies \bbeta^\top \bX^\top \bX\bbeta = 0\implies \normtwo{\bX\bbeta}^2 = 0 \implies \bX\bbeta=\bzero, 
$$
i.e., $\bbeta\in \nspace(\bX^\top \bX) \implies \bbeta\in \nspace(\bX)$. Therefore, $\nspace(\bX^\top\bX) \subseteq\nspace(\bX) $. 
Combining both inclusions yields
$
\nspace(\bX) = \nspace(\bX^\top\bX).
$
By the fundamental theorem of linear algebra (see the discussion preceding Lemma~\ref{lemma:equal-dimension-rank}), it follows that
$$
\cspace(\bX^\top)=\cspace(\bX^\top\bX).
$$
Applying the same argument to $\bX^\top$  establishes the second claim.
\end{proof}

\index{Orthonormal basis}
\index{Fundamental theorem}
\begin{theoremHigh}[Four orthonormal bases in SVD]\label{theorem:svd-four-orthonormal-Basis}
Let $\bX = \bU \bSigma \bV^\top$ be the full SVD of $\bX\in\real^{n\times p}$, where $\bU=[\bu_1, \bu_2, \ldots,\bu_n]$ and $\bV=[\bv_1, \bv_2, \ldots, \bv_p]$ are the column partitions of $\bU$ and $\bV$, respectively. 
Then:
\begin{itemize}
\item $\{\bv_1, \bv_2, \ldots, \bv_r\} $ is an orthonormal basis for the row space  $\cspace(\bX^\top)$;

\item $\{\bv_{r+1},\bv_{r+2}, \ldots, \bv_p\}$ is an orthonormal basis for the null space $\nspace(\bX)$;

\item $\{\bu_1,\bu_2, \ldots,\bu_r\}$ is an orthonormal basis for the column space $\cspace(\bX)$;

\item $\{\bu_{r+1}, \bu_{r+2},\ldots,\bu_n\}$ is an orthonormal basis for the left null space $\nspace(\bX^\top)$. 
\end{itemize}
\end{theoremHigh}
\begin{proof}[of Theorem~\ref{theorem:svd-four-orthonormal-Basis}]
By the spectral decomposition, for the symmetric matrix $\bX^\top\bX$, its column space $\cspace(\bX^\top\bX)$ is spanned by the eigenvectors. 
Therefore, the set $\{\bv_1,\bv_2, \ldots, \bv_r\}$ forms an orthonormal basis for $\cspace(\bX^\top\bX)$.
Thus, $\{\bv_1, \bv_2,\ldots, \bv_r\}$ also serves as an orthonormal basis for $\cspace(\bX^\top)$ by Lemma~\ref{lemma:rank-of-ttt}. 

Since $\bV$ is orthogonal, the space spanned by the remaining vectors $\{\bv_{r+1}, \bv_{r+2},\ldots, \bv_n\}$ is the orthogonal complement to the space spanned by $\{\bv_1,\bv_2, \ldots, \bv_r\}$, which is precisely $\nspace(\bX)$.
Thus, $\{\bv_{r+1},\bv_{r+2}, \ldots, \bv_n\}$ constitutes an orthonormal basis for $\nspace(\bX)$. 

A similar argument applied to $\bX \bX^\top$ shows that $\{\bu_1,\bu_2, \ldots, \bu_r\}$ spans $\cspace(\bX)$ and $\{\bu_{r+1},\bu_{r+2}, \ldots, \bu_n\}$ spans $\nspace(\bX^\top)$.
Alternatively, we can see that $\{\bu_1,\bu_2, \ldots,\bu_r\}$ forms a basis for the column space of $\bX$ by Lemma~\ref{lemma:column-basis-from-row-basis}~\footnote{Recall that  for any matrix $\bX$, let $\{\br_1, \br_2, \ldots, \br_r\}$ be a set of vectors in $\real^p$, which forms a basis for the row space, then $\{\bX\br_1, \bX\br_2, \ldots, \bX\br_r\}$ is a basis for the column space of $\bX$.}, since $\bu_i = \frac{\bX\bv_i}{\sigma_i},\, \forall\, i \in\{1, 2, \ldots, r\}$ (which is a key property of SVD). 
\end{proof}

\begin{figure}[h!]
\centering
\includegraphics[width=0.8\textwidth]{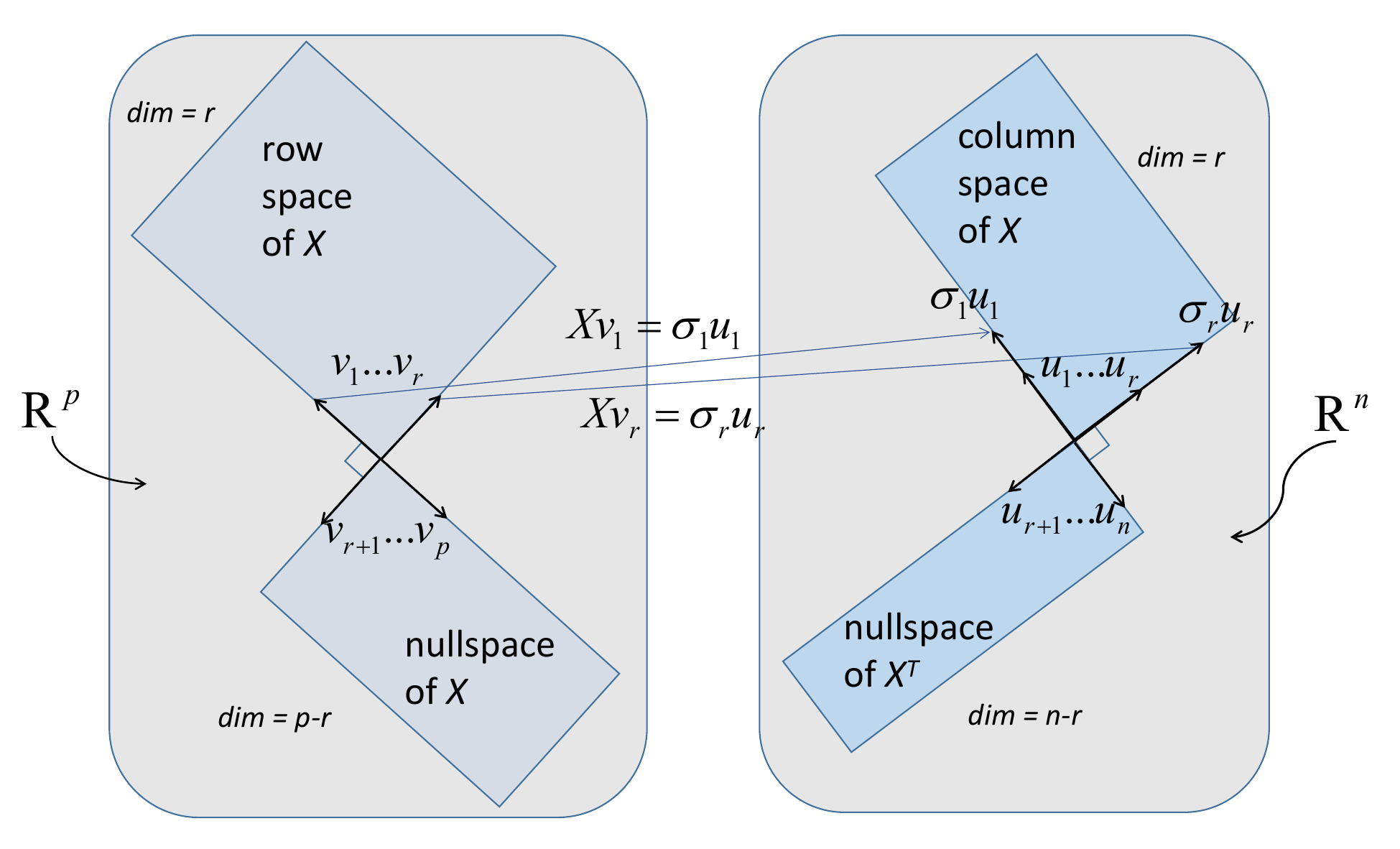}
\caption{Orthonormal bases that diagonalize $\bX$ via the SVD. The set $\{\bv_1, \bv_2, \ldots, \bv_r\} $ forms an orthonormal basis for the row space  $\cspace(\bX^\top)$, and $\{\bu_1,\bu_2, \ldots,\bu_r\}$ forms an orthonormal basis for the column space $\cspace(\bX)$.
The action of $\bX$ links these bases: for each $i \in \{1, 2, \ldots, r\}$, it transforms the row-space basis vector $\bv_i$ into the column-space basis vector $\bu_i$ scaled by the singular value $\sigma_i$, i.e., $\bX \bv_i = \sigma_i \bu_i$.
}
\label{fig:lafundamental3-SVD}
\end{figure}

The relationship among the four subspaces is illustrated  in Figure~\ref{fig:lafundamental3-SVD}.
Specifically, for each $i \in \{1, 2, \ldots, r\}$, the matrix $\bX$ maps the row-space basis vector $\bv_i$ to the column-space basis vector $\bu_i$ according to the relation $\sigma_i\bu_i=\bX\bv_i$.

\subsection*{Singular Value Bounds}
In its SVD form, a linear map can be interpreted as performing three operations on vectors: an initial rotation, followed by scaling along orthogonal directions according to the singular values, and then a final rotation. This implies that the maximum possible scaling factor equals the largest singular value, $\sigma_1$, and occurs in the direction of the corresponding right singular vector $\bv_1$.
The following theorem formalizes this intuition. It shows that if we restrict attention to the orthogonal complement of $\bv_1$, the maximum scaling is given by the second singular value, $\sigma_2$, due to the orthogonality of the singular vectors. More generally, the direction of maximum scaling that is orthogonal to the first $i-1$ right singular vectors corresponds to the $i$-th singular value, $\sigma_i$, and is aligned with the $i$-th right singular vector $\bv_i$.

\begin{theoremHigh}[Courant--Fischer-type characterization for singular values]\label{theorem:cou_fish_sing}
Let $\bX \in \real^{n\times p}$ have the singular value decomposition $\bU\bSigma\bV^\top$. 
Then the singular values satisfy
\begin{align}
\sigma_1 &= \max_{\left\{\normtwo{\bbeta}=1 \mid \bbeta \in \real^p\right\}} \normtwo{\bX\bbeta}, \label{eq:sigma1}\\
&= \max_{\left\{\normtwo{\balpha}=1 \mid \balpha \in \real^n\right\}} \normtwo{\bX^\top\balpha}, \label{eq:sigma1_alt}\\
\sigma_i &= \max_{\left\{\normtwo{\bbeta}=1 \mid \bbeta \in \real^p,\, \bbeta \perp \bv_1,\ldots,\bv_{i-1}\right\}} \normtwo{\bX\bbeta}, \label{eq:sigma_i}\\
&= \max_{\left\{\normtwo{\balpha}=1 \mid \balpha \in \real^n,\, \balpha \perp \bu_1,\ldots,\bu_{i-1}\right\}} \normtwo{\bX^\top\balpha}, \qquad 2 \leq i \leq \min\{n, p\}, \label{eq:sigma_i_alt}
\end{align}
the right singular vectors satisfy
\begin{align}
\bv_1 &= \argmax_{\left\{\normtwo{\bbeta}=1 \mid \bbeta \in \real^p\right\}} \normtwo{\bX\bbeta}, \label{eq:v1}\\
\bv_i &= \argmax_{\left\{\normtwo{\bbeta}=1 \mid \bbeta \in \real^p,\, \bbeta \perp \bv_1,\ldots,\bv_{i-1}\right\}} \normtwo{\bX\bbeta}, \qquad 2 \leq i \leq p, \label{eq:vi}
\end{align}
and the left singular vectors satisfy
\begin{align}
\bu_1 &= \argmax_{\left\{\normtwo{\balpha}=1 \mid \balpha \in \real^n\right\}} \normtwo{\bX^\top\balpha}, \label{eq:u1}\\
\bu_i &= \argmax_{\left\{\normtwo{\balpha}=1 \mid \balpha \in \real^n,\, \balpha \perp \bu_1,\ldots,\bu_{i-1}\right\}} \normtwo{\bX^\top\balpha}, \qquad 2 \leq i \leq n. \label{eq:ui}
\end{align}
\end{theoremHigh}
\begin{proof}[of Theorem~\ref{theorem:cou_fish_sing}]
Consider a unit vector $\bbeta \in \real^p$ (i.e., $\normtwo{\bbeta}=1$) that is orthogonal to $\bv_1,\ldots,\bv_{i-1}$, where $1 \leq i \leq p$ (if $i = 1$ then $\bbeta$ is arbitrary among unit vectors). 
We expand $\bbeta$ in the basis of right singular vectors of $\bX$:
$\bbeta 
= \sum_{j=i}^p \gamma_j \bv_j $, 
where $1 = \normtwo{\bbeta}^2 = \sum_{j=i}^p \gamma_j^2$. 
By the ordering of the singular values in Theorem~\ref{theorem:reduced_svd_rectangular},
\begin{align*}
\normtwo{\bX\bbeta}^2 
&= \innerproduct{\sum_{k=1}^p \sigma_k \innerproduct{\bv_k, \bbeta} \bu_k, \sum_{k=1}^p \sigma_k \innerproduct{\bv_k, \bbeta} \bu_k} 
= \sum_{k=1}^p \sigma_k^2 \innerproduct{ \bv_k, \bbeta}^2  \\
&= \sum_{k=1}^p \sigma_k^2 \innerproduct{\bv_k, \sum_{j=i}^p \gamma_j \bv_j}^2
= \sum_{j=i}^p \sigma_j^2 \gamma_j^2 
\leq \sigma_i^2 \sum_{j=i}^p \gamma_j^2 
= \sigma_i^2. 
\end{align*}
This establishes~\eqref{eq:sigma1} and~\eqref{eq:sigma_i}. To verify \eqref{eq:v1} and~\eqref{eq:vi}, observe that $\bv_i$ achieves this upper bound:
\begin{align*}
\normtwo{\bX\bv_i}^2 
&= \sum_{k=1}^p \sigma_k^2 \innerproduct{\bv_k, \bv_i}^2 
= \sigma_i^2.
\end{align*}
Applying the same reasoning to $\bX^\top$ yields~\eqref{eq:sigma1_alt},~\eqref{eq:sigma_i_alt},~\eqref{eq:u1}, and \eqref{eq:ui}.
\end{proof}

From the proof above, it also follows that the largest and smallest singular values of a matrix $\bX\in\real^{n\times p}$ (assuming $p\leq n$) satisfy
\begin{align}
\sigma_{\max}(\bX) &= \sigma_1(\bX) = \normtwo{\bX} = \max_{\normtwo{\bbeta}=1}\normtwo{\bX\bbeta}; \label{equation:svd_stre_bdP1}\\
\sigma_{\min}(\bX) &= \sigma_p(\bX) = \min_{\normtwo{\bbeta}=1} \normtwo{\bX\bbeta}. \label{equation:svd_stre_bdP2}
\end{align}
\footnote{See \citet{lu2021numerical} for a detailed proof, which relies on Rayleigh--Ritz theorem for symmetric matrices.}
This shows that 
\begin{equation}\label{equation:svd_stre_bd}
\sigma_{\min}(\bX) \normtwo{\bbeta} 
\leq \normtwo{\bX\bbeta} 
\leq \sigma_{\max}(\bX) \normtwo{\bbeta} ,
\quad \text{for any }\bbeta\in\real^p.
\end{equation}

For the purposes of this book, the following observation is especially useful.
\begin{theoremHigh}[Singular value bounds under near-orthogonality]\label{theorem:sigbd_nearortho}
Let $\bX\in\real^{n\times p}$ with $n\geq p$. 
Suppose that for some $\delta \in [0,1]$,
\begin{equation}\label{equation:sigbd_nearortho_eq1}
\normtwo{\bX^\top\bX - \bI} \leq \delta.
\end{equation}
Then the largest and smallest singular value of $\bX$ satisfy
\begin{equation}\label{equation:sigbd_nearortho_eq2}
\sigma_{\max}(\bX) \leq \sqrt{1+\delta} 
\qquad \text{and} \qquad 
\sigma_{\min}(\bX) \geq \sqrt{1-\delta}. 
\end{equation}
Conversely, if both inequalities in \eqref{equation:sigbd_nearortho_eq2} hold, then \eqref{equation:sigbd_nearortho_eq1} follows.
\end{theoremHigh}

\begin{proof}[of Theorem~\ref{theorem:sigbd_nearortho}]
Recall that by \eqref{equation:sigbd_nearortho_eqp}, the eigenvalues of $\bX^\top\bX$ are precisely the squared singular values of  $\bX$, $\lambda_i(\bX^\top\bX) = \sigma_i^2(\bX)$, $i=1,2,\ldots,p$.
Consequently, the eigenvalues of $\bX^\top\bX - \bI$ are given by $\sigma_i^2(\bX) - 1$. 
By \eqref{equation:svd_stre_bdP1} and the assumption \eqref{equation:sigbd_nearortho_eq1}, we have 
$$
\max\{\sigma_{\max}^2(\bX)-1,\, 1-\sigma_{\min}^2(\bX)\} 
= \normtwo{\bX^\top\bX - \bI} \leq \delta.
$$
This establishes the claim.
\end{proof}

The extremal singular values---both the largest and the smallest---are {1-Lipschitz functions} (see Definition~\ref{definition:lipschi_funs}) with respect to the spectral norm and the Frobenius norm. 
In other words, the difference between the extremal singular values of two matrices is bounded above by the distance between the matrices in either norm.
\begin{theoremHigh}[Lipschitz continuity of extremal singular values]\label{theorem:lipcon_extrmsing}
The smallest and largest singular values $\sigma_{\min}, \sigma_{\max}$, satisfy for all matrices $\bX,\bY$ of the same dimension,
\begin{align}
\abs{\sigma_{\max}(\bX) - \sigma_{\max}(\bY)} 
&\leq \normtwo{\bX-\bY} \leq \normf{\bX-\bY}; \label{equ:lipcon_extrmsing_eq1} \\
\abs{\sigma_{\min}(\bX) - \sigma_{\min}(\bY)} 
&\leq \normtwo{\bX-\bY} \leq \normf{\bX-\bY}. \label{equ:lipcon_extrmsing_eq2}
\end{align}
\end{theoremHigh}

\begin{proof}[of Theorem~\ref{theorem:lipcon_extrmsing}]
By the identification of the largest singular value with the operator norm we have
$$
\abs{\sigma_{\max}(\bX) - \sigma_{\max}(\bY)} 
= \absbig{\normtwo{\bX} - \normtwo{\bY}} 
\leq \normtwo{\bX-\bY}.
$$
By \eqref{equation:svd_stre_bdP2}, for the smallest singular value, we have
\begin{align*}
\sigma_{\min}(\bX) 
&= \min_{\normtwo{\bbeta}=1} \normtwo{\bX\bbeta} 
\leq \min_{\normtwo{\bbeta}=1} \left(\normtwo{\bY\bbeta} + \normtwo{(\bX-\bY)\bbeta}\right)\\
&\leq \min_{\normtwo{\bbeta}=1} \left(\normtwo{\bY\bbeta} + \normtwo{\bX-\bY}\right) = \sigma_{\min}(\bY) + \normtwo{\bX-\bY}.
\end{align*}
Therefore, $\sigma_{\min}(\bX) - \sigma_{\min}(\bY) \leq \normtwo{\bX-\bY}$ and \eqref{equ:lipcon_extrmsing_eq2} follows by symmetry. 
Finally, the inequalities involving the Frobenius norm follow from the standard relation $\normtwo{\cdot}\leq\normf{\cdot}$ (see Remark~\ref{remark:stad_matnorm}).
\end{proof}

\index{Pseudo-inverse}
\subsection{Pseudo-Inverse}\label{section:subsec_pseudo_inv}\label{section:pseudo-in-svd}
\begin{figure}[h!]
\centering
\includegraphics[width=0.98\textwidth]{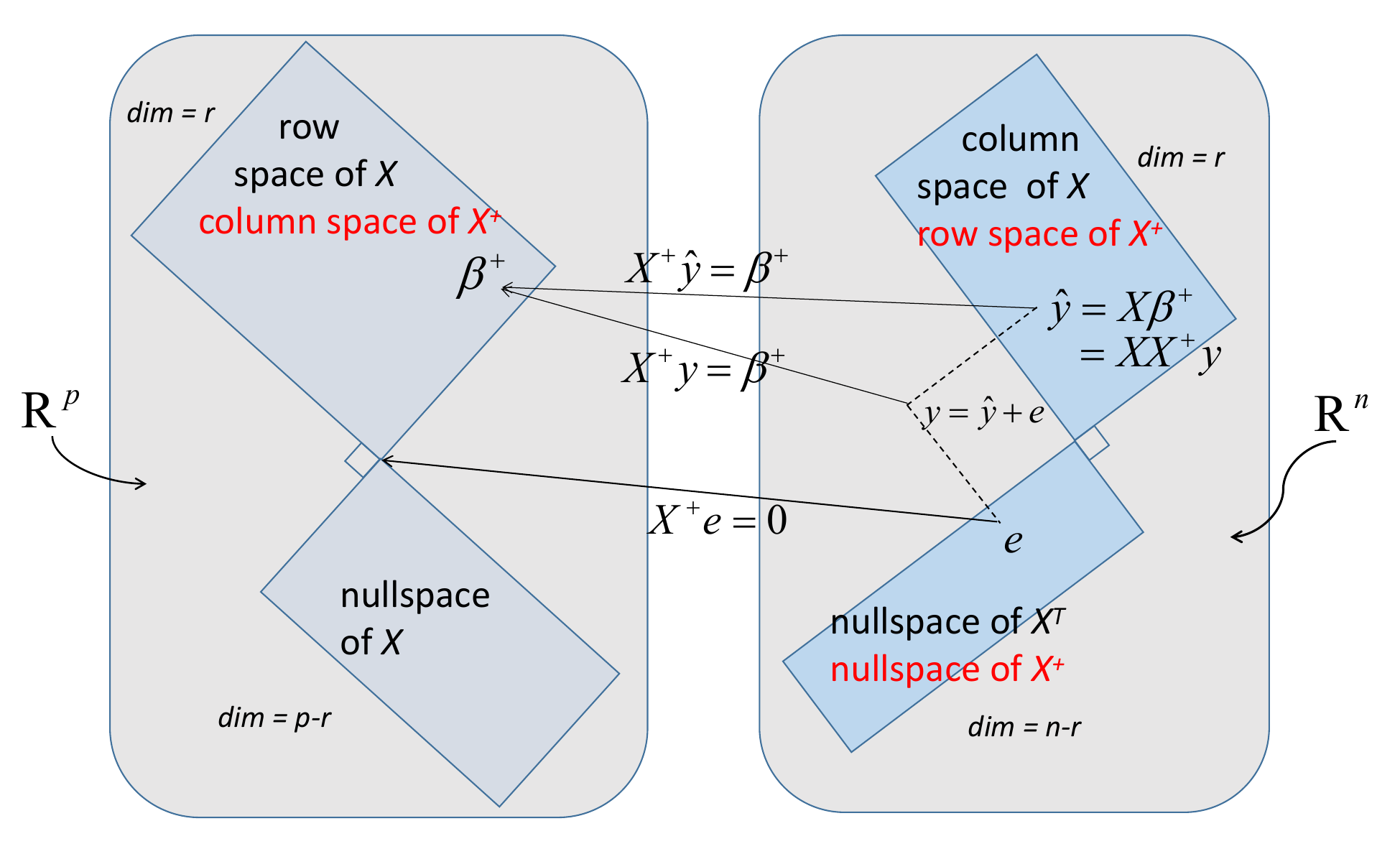}
\caption{Column space and row space of the pseudo-inverse $\bX^+$. $\bX$ transfers from the row space to the  column space. $\bX^+$ maps from the column space to the row space.}
\label{fig:lafundamental5-pseudo}
\end{figure}

For a matrix $\bX\in \real^{n\times p}$, we can define its \textit{pseudo-inverse}, denoted by $\bX^+$, which is a $p\times n$ matrix.
Intuitively, when $\bX$ multiplies a vector $\bbeta$ that lies in its row space (i.e., $\bbeta\in\cspace(\bX^\top)$), the result $\bX\bbeta$ lies  in the column space of $\bX$. 
Both of these spaces have the same dimension $r$, namely the rank of $\bX$. 
When restricted to these subspaces, $\bX$ behaves like an invertible linear transformation, and $\bX^+$ acts as its inverse.
More precisely:
\begin{itemize}
\item If $\bbeta$ belongs to the row space of $\bX$, then  $\bX^+\bX\bbeta = \bbeta$. 
\item If $\by$ belongs to the column space of $\bX$, then $\bX\bX^+\by = \by$  (see Figure~\ref{fig:lafundamental5-pseudo}).
\end{itemize}

The null space of $\bX^+$ coincides with  the null space of $\bX^\top$. 
It consists of the vectors $\by\in\real^n$  such that $\bX^\top\by = \bzero$. 
Those vectors $\by$ are orthogonal to every vector  in the column space of $\bX$, i.e., to every vector of the form $\bX\bbeta$. 
We delay the proof of this property in Theorem~\ref{theorem:pseudo-four-basis-space}.

More formally, the \textit{pseudo-inverse},  also known as the \textit{Moore-Penrose pseudo-inverse}, $\bX^+$, is defined by the \textbf{unique} (see Proposition~\ref{proposition:existence-of-pseudo-inverse}) $p\times n$ matrix satisfying the following four conditions, commonly called the \textit{Penrose conditions} \citep{penrose1955generalized}:
\begin{equation}\label{equation:pseudi-four-equations}
\boxed{
\begin{aligned}
&(C1) \qquad  \bX\bX^+\bX &=& \bX \qquad & \\
&(C2) \qquad  \bX^+\bX\bX^+ &=&\bX^+ \qquad &\\
&(C3) \qquad  (\bX\bX^+)^\top &=&\bX\bX^+ \qquad & (\text{$ \bX \bX^+ $ is symmetric})\\
&(C4) \qquad  	(\bX^+\bX)^\top &=& \bX^+\bX\qquad & (\text{$ \bX \bX^+ $ is symmetric})\\
\end{aligned}
}
\end{equation}
Conditions $(C1)$ and $(C2)$ imply that $\bX\bX^+$ and $\bX^+\bX$ are idempotent matrices, and hence they are  projection matrices (see Definition~\ref{definition:projection-matrix}). 
Conditions $(C3)$ and $(C4)$ further ensure that  these projection matrices are symmetric, which  means they are not merely oblique projections but orthogonal projection matrices (see Exercise~\ref{exercise:sym_proj_mat}).

The existence and the uniqueness of the pseudo-inverse for any matrix follows from its rank decomposition.
\begin{proposition}[Existence and uniqueness of pseudo-inverse]\label{proposition:existence-of-pseudo-inverse}
Every matrix $\bX$ has a unique pseudo-inverse. 
\end{proposition}
\begin{proof}[of Proposition~\ref{proposition:existence-of-pseudo-inverse}] 
\textbf{Existence.}
Let  $\bX=\bW\bZ\in\real^{n\times p}$ be a rank decomposition of $\bX$, where $\bW\in\real^{n\times r}$ and $\bZ\in\real^{r\times p}$ both have  full rank. 
Define  
$$
\bX^+ = \bZ^+\bW^+ = \bZ^\top (\bZ\bZ^\top)^{-1} (\bW^\top\bW)^{-1}\bW^\top,
$$
where $\bZ^+=\bZ^\top (\bZ\bZ^\top)^{-1}$ and $\bW^+ =(\bW^\top\bW)^{-1}\bW^\top$.~\footnote{It can be easily verified that $\bZ^+$ is the pseudo-inverse of $\bZ$ and $\bW^+$ is the pseudo-inverse of $\bW$.} Notably, $\bZ\bZ^\top$ and $\bW^\top\bW$ are invertible since $\bW\in \real^{n\times r}$ and $\bZ\in \real^{r\times p}$ have full rank $r$.

We now verify that $\bX^+$ satisfies the four Penrose conditions:
$$
\begin{aligned}
&(C1) \qquad \bX\bX^+\bX &=& \bW\bZ\left(\bZ^\top (\bZ\bZ^\top)^{-1} (\bW^\top\bW)^{-1}\bW^\top\right)\bW\bZ = \bW\bZ = \bX, \\
&(C2) \qquad \bX^+\bX\bX^+ &=&\left(\bZ^\top (\bZ\bZ^\top)^{-1} (\bW^\top\bW)^{-1}\bW^\top\right) \bW\bZ\left(\bZ^\top (\bZ\bZ^\top)^{-1} (\bW^\top\bW)^{-1}\bW^\top\right)\\
&&=&\bZ^\top(\bZ\bZ^\top)^{-1}(\bW^\top\bW)^{-1}\bW^\top = \bX^+,\\
&(C3) \qquad (\bX\bX^+)^\top &=& \bW(\bW^\top\bW)^{-1} (\bZ\bZ^\top)^{-1} \bZ\bZ^\top\bW^\top = \bW(\bW^\top\bW)^{-1} \bW^\top\\
&&=&\bW\bZ\bZ^\top (\bZ\bZ^\top)^{-1} (\bW^\top\bW)^{-1}\bW^\top = \bX\bX^+,\\
&(C4) \qquad (\bX^+\bX)^\top &=& \bZ^\top\bW^\top \bW(\bW^\top\bW)^{-1} (\bZ\bZ^\top)^{-1} \bZ= \bZ^\top (\bZ\bZ^\top)^{-1} \bZ\\
&&=& \bZ^\top (\bZ\bZ^\top)^{-1} (\bW^\top\bW)^{-1}\bW^\top\bW\bZ = \bX^+\bX.
\end{aligned}
$$
This shows the existence of the pseudo-inverse.

\paragraph{Uniqueness.}
Suppose $\bX_1^+$ and $\bX_2^+$ are two matrices satisfying the four Penrose conditions for $\bX$. We will show that $\bX_1^+=\bX_2^+$.
Starting from  $\bX_1^+$ and repeatedly applying the Penrose conditions, we obtain:
$$
\begin{aligned}
\bX_1^+ &= \bX_1^+\bX\bX_1^+ = \bX_1^+(\bX\bX_2^+\bX)\bX_1^+ = \bX_1^+ (\bX\bX_2^+)(\bX\bX_1^+) \qquad &(\text{by $(C2), (C1)$})\\
&=\bX_1^+ (\bX\bX_2^+)^\top(\bX\bX_1^+)^\top =  \bX_1^+\bX_2^{+\top}\bX^\top \bX_1^{+\top}\bX^\top \qquad &(\text{by $(C3)$})\\
&=\bX_1^+\bX_2^{+\top}(\bX \bX_1^{+}\bX)^\top= \bX_1^+\bX_2^{+\top}\bX^\top \qquad &(\text{by $(C1)$})\\
&= \bX_1^+(\bX \bX_2^{+})^\top =\bX_1^+ \bX \bX_2^{+} = \bX_1^+ (\bX \bX_2^+\bX) \bX_2^{+}  \qquad &(\text{by $(C3),(C1)$})\\
&= (\bX_1^+ \bX) (\bX_2^+\bX) \bX_2^{+} =  (\bX_1^+ \bX)^\top (\bX_2^+\bX)^\top \bX_2^{+}  \qquad &(\text{by $(C4)$})\\
&=(\bX \bX_1^{+} \bX)^\top \bX_2^{+\top} \bX_2^{+} = \bX^\top \bX_2^{+\top} \bX_2^{+}   \qquad &(\text{by $(C1)$})\\
&=  (\bX_2^{+}\bX)^\top \bX_2^{+} =\bX_2^{+}\bX \bX_2^{+} = \bX_2^{+}. \qquad &(\text{by $(C4),(C2)$})\\
\end{aligned}
$$
Thus, $\bX_1^+=\bX_2^+$, proving that the pseudo-inverse is unique.
\end{proof}

\begin{lemma}[Property of pseudo-inverse]\label{lemma:xtxxpllus_pseudo}
For any matrix $\bX\in\real^{n\times p}$, it holds that
$$ \bX^\top \bX \bX^+ = \bX^\top .$$
\end{lemma}
\begin{proof}[of Lemma~\ref{lemma:xtxxpllus_pseudo}]
By Condition $(C1)$ $\bX \bX^+ \bX = \bX$, we have by taking the  transpose of both sides:
$$
(\bX \bX^+ \bX)^\top = \bX^\top
\qquad \implies\qquad 
\bX^\top (\bX \bX^+)^\top = \bX^\top.
$$
Substituting the Condition $(C3)$ $ (\bX \bX^+)^\top = \bX \bX^+ $  into the equation above, we get
$
\bX^\top \bX \bX^+ = \bX^\top.
$
This completes the proof.
\end{proof}
The equality $ \bX^\top \bX \bX^+ = \bX^\top $ encodes the idea that applying $ \bX^+ $ to $ \bX $ does not distort the action of $ \bX^\top $ on vectors.

We have shown the four fundamental subspaces associated with the SVD of a matrix in Theorem~\ref{theorem:svd-four-orthonormal-Basis}. 
Similarly, we now describe the four fundamental subspaces associated with the pseudo-inverse.
\begin{theoremHigh}[Four subspaces in pseudo-inverse]\label{theorem:pseudo-four-basis-space} 
Let $\bX\in\real^{n\times p}$ be any matrix, and let  $\bX^+$ denote its  pseudo-inverse.
Then the following relationships hold:
\begin{itemize}
\item The column space of $\bX^+$ is the same as the row space of $\bX$;
\item The row space of $\bX^+$ is the same as the column space of $\bX$; 
\item The null space of $\bX^+$ is the same as the null space of $\bX^\top$; 
\item The null space of $\bX^{+\top}$ is the same as the null space of $\bX$.
\end{itemize}
\noindent These relationships among the four fundamental subspaces are illustrated in Figure~\ref{fig:lafundamental5-pseudo}.
\end{theoremHigh}
\begin{proof}[of Theorem~\ref{theorem:pseudo-four-basis-space}]
Following arguments similar to those in the proof of Lemma~\ref{lemma:rank-of-ttt}, we can show that
$$
\begin{aligned}
\cspace(\bX\bX^+) &= \cspace(\bX) \qquad \text{and} \qquad \nspace(\bX^+\bX) = \nspace(\bX)\\
\cspace(\bX^+\bX) &= \cspace(\bX^+)\qquad \text{and} \qquad \nspace(\bX\bX^+) = \nspace(\bX^+).
\end{aligned}
$$
Additionally, from Conditions $(C3)$ and $(C4)$ in the definition of pseudo-inverses, we know that:
$$
(\bX^+\bX)^\top = \bX^+\bX \qquad \text{and} \qquad (\bX\bX^+)^\top =\bX\bX^+.
$$
Using  the fundamental theorem of linear algebra (see the discussion preceding Lemma~\ref{lemma:equal-dimension-rank}), we realize that $\cspace(\bX\bX^+)$ is the orthogonal complement to $\nspace((\bX\bX^+)^\top)$, and $\cspace(\bX^+\bX)$ is the orthogonal complement to $\nspace((\bX^+\bX)^\top)$: 
$$
\begin{aligned}
\cspace(\bX\bX^+) \perp \nspace((\bX\bX^+)^\top) &\leadto \cspace(\bX\bX^+) \perp \nspace(\bX\bX^+) \\
\cspace(\bX^+\bX) \perp \nspace((\bX^+\bX)^\top) &\leadto \cspace(\bX^+\bX) \perp \nspace(\bX^+\bX).
\end{aligned}
$$
This implies 
$$
\begin{aligned}
\cspace(\bX) \perp \nspace(\bX^+) \qquad \text{and} \qquad \cspace(\bX^+)\perp \nspace(\bX).
\end{aligned}
$$
That is, $\nspace(\bX^+) = \nspace(\bX^\top)$ and $\cspace(\bX^+) = \cspace(\bX^\top)$. 
By the fundamental theorem of linear algebra again, this also implies: $\cspace(\bX^{+\top})=\cspace(\bX)$ and $\nspace(\bX^{+\top})=\nspace(\bX)$.
\end{proof}

\index{Pseudo-inverse in SVD}
\subsection*{Pseudo-Inverse in SVD}
Given the SVD of a matrix $\bX\in\real^{n\times p}$, we provide further discussion on its pseudo-inverse in different cases.
For the full SVD of  $\bX = \bU\bSigma\bV^\top$,  we consider the following scenarios:
\paragraph{Case $n>p=r$.} Since the matrix $\bX$ has independent columns, the left-pseudo-inverse can be obtained by 
\begin{equation}\label{equation:pseudo_svd_case1}
\begin{aligned}
	\text{left-pseudo-inverse} 
	&= \bX^+ = (\bX^\top\bX)^{-1}\bX^\top 
	= (\bV \bSigma^\top \bU^\top\bU\bSigma\bV^\top)^{-1}\bV \bSigma^\top \bU^\top\\
	&= \bV (\bSigma^\top\bSigma)^{-1}\bV^\top \bV \bSigma^\top \bU^\top 
	= \bV[(\bSigma^\top\bSigma)^{-1}\bSigma^\top] \bU^\top 
	= \bV \bSigma^+ \bU^\top,
\end{aligned}
\end{equation}
where the last equality holds because $\bSigma^+=(\bSigma^\top\bSigma)^{-1}\bSigma^\top$.
\paragraph{Case $p>n=r$.} Since the matrix $\bX$ has independent rows, the right-pseudo-inverse can be obtained by 
\begin{equation}\label{equation:pseudo_svd_case2}
\begin{aligned}
\text{right-pseudo-inverse} 
&= \bX^+ = \bX^\top(\bX\bX^\top)^{-1}
= (\bU\bSigma\bV^\top)^\top[(\bU\bSigma\bV^\top)(\bU\bSigma\bV^\top)^\top]^{-1}\\
&= \bV\bSigma^\top\bU^\top(\bU\bSigma\bV^\top\bV\bSigma^\top\bU^\top)^{-1} 
= \bV\bSigma^\top\bU^\top\bU^{-\top}(\bSigma\bSigma^\top)^{-1}\bU^{-1} \\
&= \bV\bSigma^\top(\bSigma\bSigma^\top)^{-1}\bU^{-1} 
=\bV \bSigma^+ \bU^\top,
\end{aligned}
\end{equation}
where the last equality follows because $\bSigma^+=\bSigma^\top(\bSigma\bSigma^\top)^{-1}$.

\paragraph{Case rank-deficient.} Even when $\bX$ is rank-deficient, the pseudo-inverse is still given by $\bX^+=\bV \bSigma^+ \bU^\top$,
where  $\bSigma^+ \in \real^{p\times n}$ has  the diagonal block  $\diag(\frac{1}{\sigma_1}, \frac{1}{\sigma_2}, \ldots, \frac{1}{\sigma_r})$ in its upper-left corner, and zeros elsewhere. 
It is straightforward to verify that this definition satisfies the four Penrose conditions in \eqref{equation:pseudi-four-equations}.

In all cases, $\bSigma^+$ is obtained by taking the reciprocal of each nonzero singular value in $\bSigma$. 
Thus, the pseudo-inverse via SVD is summarized as:
\begin{equation}\label{equation:pseudo_use_svd}
\bX^+ = \bV 
\begin{bmatrix}
\bSigma^{-1}_1 & \bzero \\
\bzero & \bzero
\end{bmatrix} 
\bU^\top,
\quad 
\text{with }
\bSigma = 
\begin{bmatrix}
\bSigma_1 & \bzero \\
\bzero & \bzero
\end{bmatrix} 
\text{ and }
\bSigma_1\in\real^{r\times r}.
\end{equation}
If $\bX$ is nonsingular, then $\bX^+ = \bX^{-1}$, so \eqref{equation:pseudo_use_svd} is a generalization of the usual inverse.
The pseudo-inverse of a scalar is
$$
\sigma^+ = 
\begin{cases}
1/\sigma, & \text{if } \sigma \neq 0; \\
0, & \text{if } \sigma = 0.
\end{cases}
$$
This shows the important fact that the pseudo-inverse $\bX^+$ is not a continuous function of $\bX$, unless we allow only perturbations that do not change the rank of $\bX$. The pseudo-inverse has the property
$$
\bX^+ = \lim_{\delta \to 0} (\bX^\top \bX + \delta \bI)^{-1} \bX^\top.
$$

Using the SVD representation of pseudo-inverses in \eqref{equation:pseudo_use_svd}, we can derive additional properties of the pseudo-inverse of a matrix $\bX$.
\begin{lemma}[Properties of pseudo-inverse using SVD]\label{lemma:prop_pseu_use_svd}
Let $\bX\in\real^{n\times p}$. The following properties of the pseudo-inverse follow from \eqref{equation:pseudo_use_svd}:
\begin{enumerate}[(i)]
\item $\bX^+ = (\bX^\top \bX)^+ \bX^\top = \bX^\top (\bX\bX^\top)^+$.
\item \label{equation:pseu_full_col_trans} When $\rank(\bX) = p$, this becomes $\bX^+ = (\bX^\top \bX)^{-1} \bX^\top, \ (\bX^\top)^+ = \bX(\bX^\top \bX)^{-1}$.
\item $(\bX^+)^+ = \bX$.
\item $(\alpha \bX)^+ = \alpha^+ \bX^+$.
\item $(\bX^+)^\top = (\bX^\top)^+$.
\item $(\bX^\top \bX)^+ = \bX^+ (\bX^\top)^+$.
\item $\bX$, $\bX^\top$, $\bX^+$, and $\bX^+ \bX$ all have the same rank, which equals $\trace(\bX^+ \bX)$.
\item If $\bU$ and $\bV$ are orthogonal, then $(\bU \bX \bV^\top)^+ = \bV \bX^+ \bU^\top$.
\item If $\bX = \sum_i \bX_i$, where $\bX_i \bX_j^\top = \bzero$, $\bX_i^\top \bX_j = \bzero$, $i \neq j$, then $\bX^+ = \sum_i \bX_i^+$.
\item If $\bX$ is normal ($\bX\bX^\top = \bX^\top \bX$), then $\bX^+ \bX = \bX\bX^+$ and $(\bX^n)^+ = (\bX^+)^n$.
\end{enumerate}
\end{lemma}

It is important to note that, unlike the standard matrix inverse, the pseudo-inverse does not generally satisfy $\bX\bX^+ = \bX^+ \bX$ or $(\bX \bY)^+ = \bY^+ \bX^+$. 
For example, let $\bX = \footnotesize [1 , 0] $ and $\bY = [1 , 1]^\top$. Then $\bX \bY = 1$, but
$ \bY^+ \bX^+ = \frac{1}{2} [1 , 1] \begin{bmatrixfoot} 1 \\ 0 \end{bmatrixfoot} =  \frac{1}{2}. $
Necessary and sufficient conditions for the identity $(\bX \bY)^+ = \bY^+ \bX^+$ to hold were established by \citet{greville1966note}. The following theorem provides an important sufficient condition under which this equality is valid.

\begin{proposition}[Sufficient condition for $(\bX \bY)^+ = \bY^+ \bX^+$]\label{proposition:suff_pro_pseudo}
If $\bX \in \real^{n \times p}$, $\bY \in \real^{p \times m}$, and $\rank(\bX) = \rank(\bY) = p$, 
then
$ (\bX \bY)^+ = \bY^+ \bX^+ = \bY^\top (\bY \bY^\top)^{-1} (\bX^\top \bX)^{-1} \bX^\top. 
$
\end{proposition}
\begin{proof}[of Proposition~\ref{proposition:suff_pro_pseudo}]
The last equality follows from Lemma~\ref{lemma:prop_pseu_use_svd}.\ref{equation:pseu_full_col_trans}. The first equality follows from the proof of Proposition~\ref{proposition:existence-of-pseudo-inverse} and  is verified by showing that the four Penrose conditions are satisfied.
\end{proof}

\index{Cauchy--Schwarz inequality}\index{Inequalities}
\subsection{Famous Inequalities}\label{section:inequalities}
In this section, we introduce several classical inequalities that will be used frequently throughout the text.
The \textit{Cauchy--Schwarz inequality}  is considered one of the most important and widely used inequalities in mathematics.

\index{Schwarz matrix inequality}
\begin{theoremHigh}[Cauchy--Schwarz matrix (vector) inequality\index{Cauchy--Schwarz inequality}]\label{theorem:cs_matvec}
For any  matrices $\bX, \bY\in\real^{n\times p}$, the following inequality holds:
$$
\normf{\bX^\top \bY} \leq \normf{\bX} \cdot \normf{\bY}.
$$
This is a direct consequence of the Cauchy--Schwarz inequality applied to the Frobenius inner product $\innerproduct{\bX,\bY} = \trace{(\bX^\top \bY)}$.

Similarly, for any vectors $\bu, \bv\in\real^n$, we have 
\begin{equation}\label{equation:vector_form_cauchyschwarz}
\abs{\bu^\top \bv} \leq \normtwo{\bu} \cdot  \normtwo{\bv}.
\end{equation}
In two dimensions, this becomes the familiar algebraic identity:
$$
(ac+bd)^2 \leq (a^2 +b^2)(c^2+d^2).
$$
\end{theoremHigh}
The vector form of the Cauchy--Schwarz inequality plays a crucial  role in various branches of modern mathematics, including Hilbert space theory and numerical analysis \citep{wu2009various}. 
For completeness, we present a simple proof of the vector case. 
Let $\bu,\bv\in \real^n$. 
Consider the nonnegative sum
$$
\begin{aligned}
0\leq \sum_{i=1}^{n}\sum_{j=1}^{n} (u_i v_j - u_j v_i)^2 &= 
\sum_{i=1}^{n}\sum_{j=1}^{n} u_i^2v_j^2 + \sum_{i=1}^{n}\sum_{j=1}^{n} v_i^2 u_j^2 - 2\sum_{i=1}^{n}\sum_{j=1}^{n} u_iu_j v_iv_j\\
&=\left(\sum_{i=1}^{n} u_i^2\right) \left(\sum_{j=1}^{n} v_j^2\right) +
\left(\sum_{i=1}^{n} v_i^2\right) \left(\sum_{j=1}^{n} u_j^2\right) - 
2\left(\sum_{i=1}^{n} u_iv_i \right)^2\\
&=2 \normtwo{\bu}^2 \cdot \normtwo{\bv}^2 -2 \normtwo{\bu^\top\bv}^2,
\end{aligned}
$$
from which the result follows.
The equality holds if and only if $\bu = \eta\bv$ for some constant $\eta\in \real$, i.e., $\bu$ and $\bv$ are linearly dependent.

\holders inequality, named after \textit{Otto  H{\"o}lder}, is another cornerstone result with broad applications in optimization, machine learning, probability, and functional analysis. It generalizes the Cauchy--Schwarz inequality (which corresponds to the special case $s=t=2$).

\begin{theoremHigh}[\holders  inequality\index{\holders  inequality}]\label{theorem:holder-inequality}
Let $s,t>1$ satisfy $\frac{1}{s}+\frac{1}{t} = 1$. Then  for any vector $\bx,\by\in \real^n$, we have
$$
\sum_{i=1}^{n}x_i y_i
\leq 
\abs{\sum_{i=1}^{n}x_i y_i}
\leq \sum_{i=1}^{n}\abs{x_i} \abs{y_i} \leq \left(\sum_{i=1}^{n}  \abs{x_i}^s\right)^{1/s}  \left(\sum_{i=1}^{n} \abs{y_i}^t\right)^{1/t}=\norm{\bx}_s\norm{\by}_t,
$$
where $\norm{\bx}_s = \left(\sum_{i=1}^{n}  \abs{x_i}^s\right)^{1/s} $ denotes the \textit{$\ell_s$-norm}  of $\bx$. 
The equality holds if the two sequences $\{\abs{x_i}^s\}$ and $\{\abs{y_i}^t\}$ are linearly dependent~\footnote{To be more concrete, the equality attains if and only if $\abs{\bx^\top\by}=\abs{\bx}^\top\abs{\by}$ and 
$$
\left\{
\begin{aligned}
&\abs{\bx}\hadaprod\abs{\by}=\norminf{\by}\abs{\bx}, &\text{if }& s=1; \\ 
&\abs{\bx}\hadaprod\abs{\by}=\norminf{\bx}\abs{\by}, &\text{if }& s=\infty;\\ 
&\norm{\by}_t^{1/s} \abs{\bx}^{\hadaprod1/t} = \norm{\bx}_s^{1/t}\abs{\by}^{\hadaprod1/s}, &\text{if }&1<s<\infty, \\ 
\end{aligned}
\right.
$$
where $\hadaprod$ denotes Hadamard product/power; see \citet{bernstein2008matrix} and the references therein.
}.
When $s=t=2$, this reduces to the vector Cauchy--Schwarz inequality \eqref{equation:vector_form_cauchyschwarz}.
\end{theoremHigh}
\begin{proof}
	See \citet{lu2021numerical} for a proof.
\end{proof}


\index{Random variable}
\section{Probability Theory and Common Probability Distributions}\label{section:stat_prob}

A \textit{random variable} is a variable that takes on different values according to some underlying randomness, used to model uncertain outcomes or events. 
We denote the random variable itself using a lowercase letter in \textit{normal (non-italic) font}, and its possible realizations (observed or potential values) using lowercase letters in \textit{italic font}.  
For instance, $y_1$ and $y_2$ are  possible values (realizations) of the random variable $\ry$. 
For vector-valued random variables, we use \textit{boldface normal} notation: the random vector is written as $\rvy$, and a specific realization is denoted by $\by$ in \textit{italic} bold---though to avoid confusion, it's common to distinguish them typographically.
In this text, we adopt the convention that the random vector is $\rvy$ (normal bold), and a realization is $\by$ (italic bold).
Similarly, for matrix-valued random variables, we denote the random matrix as $\rmY$ (upright bold) and a particular realization as $\bY$ (italic bold).

However, a random variable alone only describes the set of possible outcomes---it must be paired with a probability distribution that specifies how likely each outcome is.
In probability theory, the outcome $\rx$ of a random phenomenon is modeled as a random variable, and the phenomenon itself is described by its \textit{cumulative distribution function (c.d.f., CDF)}, defined as:
$$
F(x) \equiv \prob[\rx \leq x].
$$ 
Now suppose the phenomenon has a characteristic $\btheta$ that influences the probabilities of the outcomes of $\rx$. Such a characteristic is called a \textit{parameter}. Because the probability of the event $\{\rx \leq x\}$ depends on $\btheta$, the CDF must also depend on $\btheta$. We therefore write:
$$
F(x; \btheta) \equiv {\prob}_{\btheta}[\rx \leq x].
$$
In probability theory, if both the functional form $F(x; \btheta)$ and the true value of $\btheta$ are known, then for any possible outcome $x$, we can compute the probability: $\prob_{\btheta}[\rx \leq x] = F(x; \btheta)$.

A useful result in probability theory is the union bound, which provides an upper bound on the probability of the union of a collection of events.
\begin{theoremHigh}[Union bound\index{Union bound}]\label{theorem:union_bound_proof}
Let $\sA_1, \sA_2, \dots, \sA_n$ be a collection of events in a probability space. Then
$$
\prob\left(\cup_i \sA_i\right) 
\leq \sum_{i=1}^n \prob(\sA_i).
$$
\end{theoremHigh}
\begin{proof}[of Theorem~\ref{theorem:union_bound_proof}]
For each $i\in\{1,2,\ldots,n\}$, refine the set
$
\widetilde{\sA}_i \triangleq \sA_i \cap \bigcap_{j=1}^{i-1} \comple{\sA_j}
$
It is straightforward to show that $ \bigcup_{j=1}^k \widetilde{\sA}_j = \bigcup_{j=1}^k \sA_j $ for any $ k $ by induction, so $ \bigcup_i \sA_i = \bigcup_i \widetilde{\sA}_i $. 
By construction, the sets $ \widetilde{\sA}_1, \widetilde{\sA}_2, \dots $ are pairwise disjoint, so
\begin{align*}
\prob\left( \bigcup_i \sA_i \right) 
= \prob\left( \bigcup_i \widetilde{\sA}_i \right) = \sum_i \prob\left( \widetilde{\sA}_i \right) 
\leq \sum_i \prob(\sA_i) 
\end{align*}
where the last inequality follows because $\widetilde{\sA}_i \subseteq \sA_i$.
\end{proof}

\paragrapharrow{Discrete  random variables.}
Random variables can be either discrete or continuous. A \textit{discrete random variable} takes on a finite or countably infinite number of distinct states. 
These states need not be numerical; they may be symbolic or categorical states without numerical values  (e.g., ``red," ``blue," ``success," ``failure"). 
In contrast, a \textit{continuous random variable} takes values in an uncountable set, typically a subset of the real numbers.

A probability distribution describes how probability is assigned to the possible outcomes of one or more random variables. The way we represent this distribution depends on whether the variables are discrete or continuous.

For discrete random variables, we use a \textit{probability mass function (p.m.f., PMF)}, often denoted by a capital $\prob$ (or sometimes $f$).
The PMF maps each possible state of the random variable to the probability that the variable assumes that state. For example, the notation
$$
\prob(\ry = y)
\quad (\text{or equivalently } \prob(y), \;f_{\ry}(y), \; 
\text{Pr}_{\ry}(y)) 
$$
represents the probability that the random variable $\ry$  equals the value $y$.
By definition, probabilities lie in $[0,1]$, where 1 indicates certainty and 0 indicates impossibility.
An alternative and common convention is to first define a random variable and then specify its distribution using the ``$\sim$" symbol:
$$
\ry \sim \prob(\ry).
$$

PMFs can also be defined over multiple random variables simultaneously. This is called a \textit{joint probability mass distribution functions} or \textit{joint frequency functions}. 
For instance, $\prob(\rx = x, \ry = y)$ denotes the probability that $x = x$ and $y = y$ occur together.
We often use shorthand notations such as $\prob(x, y)$ or $\prob_{\rx, \ry}(x,y)$. 
Moreover, if the PMF depends on some known parameters $\btheta$, we may write it as $\prob(x, y \mid \btheta)$  (or $\prob_{\btheta}(x, y)$, $f(x, y; \btheta)$) for brevity.

In many applications, we are interested in the probability of one event given that another event has occurred. This is known as \textit{conditional probability}.
The conditional probability that $\rx = x$ given $\ry = y$ is written as  $\prob(\rx = x \mid  \ry = y)$, 
and is computed using the formula: 
$$
\prob(\rx = x \mid  \ry = y) = \frac{\prob(\rx = x, \ry = y)}{\prob(\ry = y)},
$$
provided that $\prob(\ry = y)>0$.
This relationship is fundamental to \textit{Bayes' theorem} \citep{bayes1958essay}.

Conversely, when the joint distribution over a set of variables is known, we may wish to find the distribution over a subset of those variables. This is called the \textit{marginal probability mass distribution} (or simply \textit{marginal distribution}).
For example, if $\rx$ and $\ry$ are discrete random variables with known joint PMF $\prob(\rx, \ry)$, the marginal distribution of $\rx$ is obtained by summing over all possible values of $\ry$:
$$
\prob(\rx=x) = \sum_{y} \prob(\rx=x, \ry=y).
$$

\paragrapharrow{Continuous random variables.}

When working with continuous random variables, we describe probability distributions using a \textit{probability density function (p.d.f., PDF)} instead of a probability mass function (PMF).
For a function $p$ to qualify as a valid PDF, it must satisfy the following properties:
\begin{itemize}
\item The domain of $p$ must be the set of all possible values (or states) of the random variable $\ry$.
\item Unlike a PMF, the PDF is not required to satisfy $p(y) \leq  1$. However, it must be nonnegative everywhere: $\forall\, y\in\ry, p(y)\geq 0$.
\item It must integrate to 1 over its entire domain: $\int p(y)dy = 1$.
\end{itemize}
Importantly, the value $p(y)$ (also denoted as $f_{\ry}(y)$ or $p_{\ry}(y)$) does not represent a probability. Instead, it describes a density. The probability that the random variable falls within an infinitesimally small interval around y of width $\delta y$ is approximately $p(y)\delta y$.
More formally, for any measurable set $\sA$, the probability that $y\in\sA$ is given by
$$
\prob(y\in\sA) = \int_{\sA}p(y)dy.
$$
If the PDF depends on known parameters $\btheta$, it is commonly written as
$$
p(x\mid \btheta),
\quad 
f_{\rx}(x; \btheta),
\quad  
\text{or simply }\; f(x; \btheta).
$$
where the semicolon or vertical bar indicates dependence on the parameter(s).

\index{Variance}
\paragrapharrow{Distribution function.}
On the other hand, any probability distribution can also be characterized by its cumulative distribution function (c.d.f., CDF).
Let $\rx$ be a random variable  with CDF defined as  $F(x)=\prob[\rx\leq x]$, where $F(x)$ is nondecreasing, right-continuous, and satisfies
$$
0 \leq F(x) \leq 1, \qquad \lim_{x\rightarrow -\infty}F(x) = 0, 
\qquad\text{and}\qquad 
\lim_{x\rightarrow \infty}F(x) = 1.
$$
\begin{subequations}
The \textit{expected value (mean)} $\mu$ and \textit{variance} $\omega^2$ of $\rx$ in the continuous case are defined as
\begin{equation}
\mu = \Exp[\rx] = \int_{-\infty}^{\infty} x dF(x), \qquad \omega^2=\Var[\rx] = \Exp[(\rx - \mu)^2] = \int_{-\infty}^{\infty} (x - \mu)^2 dF(x).
\end{equation}
Similarly, in the discrete case, where $\prob(\rx=x)$ denotes the probability mass function, they are defined as 
\begin{equation}
\mu = \Exp[\rx] = \sum_{x} x \prob(x), \qquad \omega^2 =\Var[\rx]= \Exp[(\rx - \mu)^2] = \sum_{x} (x - \mu)^2 \prob(x).
\end{equation}
It is straightforward to verify the useful identity: 
\begin{equation}
\Var[\rx] = \Exp[\rx^2] - (\Exp[\rx])^2.
\end{equation}
The \textit{median} of a random variable $\rx$ is any  number $M$ satisfying
\begin{equation}
\prob(\rx \geq M) \geq \frac{1}{2} 
\qquad \text{and} \qquad 
\prob(\rx \leq M) \geq \frac{1}{2}.
\end{equation}
\end{subequations}
\begin{exercise}[Uniform distribution\index{Uniform distribution}]\label{exercise:uniform_dist}
Let $\rx\sim \uniformdist(x\mid a,b)$ be a uniform distributed variable such that $f_{\rx}(x) = \frac{1}{b-a}$ if $a\leq x\leq b$ and $f_{\rx}(x) =0$ otherwise. Show that
$$
\Exp[\rx] = \frac{a+b}{2}
\qquad \text{and}\qquad 
\Var[\rx] =\frac{(b-a)^2}{12}.
$$
\end{exercise}

Now consider a random vector $\rvx = [\rx_1, \rx_2 \ldots, \rx_n]^\top$. 
Its \textit{joint probability  distribution function} (or simply called  the probability distribution function)---denoted as $F_{\rvx}(\bx) $, $F(x_1, x_2, \ldots, x_n)$, or $F(\bx)$---is defined as 
$$
F_{\rvx}(x_1, x_2, \ldots, x_n) = \prob[\rx_1\leq x_1, \rx_2, \leq x_2, \ldots, \rx_n\leq x_n].
$$
Depending on whether the components are discrete or continuous, we characterize the joint distribution as follows:
\begin{itemize}
\item When all the variables are discrete, as discussed above, the joint probability mass function or the joint frequency function can be characterized as 
$$
f_{\rvx}(x_1, x_2, \ldots, x_n) = \prob[\rx_1=x_1, \rx_2=x_2, \ldots, \rx_n=x_n].
$$
\item On the contrary, if the variables are continuous, the \textit{joint probability density function} is a function $f_{\rvx} :\real^n\rightarrow [0, \infty)$ (or simply $f(\bx)$, $p_{\rvx}(\bx)$, $p(\bx)$) such that 
$$
F_{\rvx}(x_1, x_2, \ldots, x_n) = \int_{-\infty}^{x_1}\ldots \int_{-\infty}^{x_n} f_{\rvx}(y_1, y_2, \ldots, y_n)d y_1 \ldots dy_n.
$$
In this case, when $f_{\rvx}$ is continuous at $\bx=[x_1,x_2, \ldots,x_n]^\top$, we have 
$$
f_{\rvx}(x_1,x_2, \ldots,x_n) =\frac{\partial^n}{\partial x_1 \partial x_2\ldots \partial x_n} F_{\rvx}(x_1,x_2, \ldots,x_n).
$$
The random variables are independent if and only if 
$$
\begin{aligned}
F_{\rvx}(x_1,x_2, \ldots,x_n) &= F_{\rx_1}(x_1) \cdot F_{\rx_2}(x_2)\cdot\ldots \cdot F_{\rx_n}(x_n);\\
\text{or} \quad f_{\rvx}(x_1,x_2, \ldots,x_n) &= f_{\rx_1}(x_1) \cdot f_{\rx_2}(x_2)\cdot\ldots \cdot f_{\rx_n}(x_n).
\end{aligned}
$$
The \textit{conditional or marginal probability density functions} for the continuous cases are analogous to those in the discrete case, except that integration is employed instead of summation. For example, $f_{\rx}(x) = \int f_{\rx,\ry}(x, y)dy$.
\end{itemize}

\index{Covariance}
\index{Correlation}
Specifically, the joint distribution of two random variables $\rx_i$ and $\rx_j$ is  is given by the joint cumulative distribution function $F(x_i, x_j) \equiv \prob[\rx_i\leq x_i, \rx_j\leq x_j]$. Then the covariance  between $\rx_i$ and $\rx_j$, denoted $\omega_{ij}$,  is defined as
$$
\omega_{ij} \triangleq \Cov[\rx_i, \rx_j] \triangleq \Exp[(\rx_i - \mu_i)(\rx_j - \mu_j)] = \int_{x_i, x_j = -\infty}^{\infty} (x_i - \mu_i)(x_j - \mu_j) dF(x_i, x_j),
$$
where $\mu_i = \Exp[\rx_i]$
Equivalently, $\omega_{ij} = \Exp[\rx_i \rx_j] - \mu_i \mu_j$. 
For a random vector $\rvx\in\real$, the covariance matrix $\bOmega \in \real^{n \times n}$  is defined as
\begin{equation}
\Cov[\rvx] \triangleq \bOmega = \Exp[(\rvx - \bmu)(\rvx - \bmu)^\top] = \Exp(\rvx\rvx^\top) - \bmu \bmu^\top,
\end{equation}
where $\bmu = \Exp[\rvx] = [\mu_1, \mu_2, \ldots, \mu_n]^\top$. 

The variance of two random variables quantifies the degree of linear dependence between the two.
A closely related measure is the \textit{correlation}, defined as
\begin{equation}
\Corr[\rx_1, \rx_2] \triangleq \frac{\Cov[\rx_1, \rx_2]}{\sqrt{\Var[\rx_1]\Var[\rx_2]}}.
\end{equation}
While covariance and correlation both capture linear dependence, correlation is scale-invariant and always lies in the interval $[-1,1]$.
It also holds that 
\begin{equation}
\abs{\Corr[\rx_1, \rx_2]} \leq \sqrt{\Var[\rx_1]\Var[\rx_2]}.
\end{equation}

\begin{lemma}[Linear transformation]\label{lemma:lin_tran_prob}
Let $\rvy = \bA\rvx +\bb$, where $\bA \in \real^{m \times n}$  and $\bb\in\real^m$ are deterministic, and let $\rvx \in \real^n$ be a random vector with  mean  $\Exp[\rvx] = \bmu\in\real^n$ and covariance  $\bOmega\in\real^{n\times n}$. 
Then
$$
\Exp[\rvy] = \bA\bmu+\bb 
\qquad\text{and}\qquad  
\Cov[\rvy] = \bA\bOmega\bA^\top.
$$
\end{lemma}
\begin{proof}[of Lemma~\ref{lemma:lin_tran_prob}]
The first identity follows from the linearity of expectation. For the second, observe that
$$
\Cov[\rvy] = \Cov[\bA\rvx] = \Exp[\bA(\rvx - \bmu)(\rvx - \bmu)^\top \bA^\top ]= \bA \Exp[(\rvx - \bmu)(\rvx - \bmu)^\top] \bA^\top = \bA\bOmega\bA^\top.  
$$
This completes the proof.
\end{proof}
In the special case when $\bA=\ba^\top$ is a row vector, $\rvy = \ba^\top\rvx$ is a linear functional of $\rvx$. Then, if $\Cov[\rvx] = \sigma^2 \bI$, $\Cov[\rvy] = \sigma^2 \ba^\top \ba$. 
\begin{lemma}[Quadratic transformation]\label{lemma:quad_tra_prob}
Let $\bA \in \real^{n \times n}$ be a symmetric matrix, and let $\rvx\in\real^n$ be a random vector with mean $\bmu=[\mu_i]\in\real^n$ and covariance  $\bOmega=[\omega_{ij}]\in\real^{n\times n}$. 
Then
$$
\Exp[\rvx^\top\bA\rvx] = \bmu^\top \bA \bmu + \trace(\bA\bOmega),
$$
where $\trace(\bA\bOmega)$ denotes the trace of $\bA\bOmega$, i.e., sum of diagonal elements of the matrix.
\end{lemma}
\begin{proof}[of Lemma~\ref{lemma:quad_tra_prob}]
Since $\rvx^\top \bA \rvx = \sum_{i=1}^n \sum_{j=1}^n a_{ij} \rx_i \rx_j$, 
it follows that  $\Exp[\rvx^\top \bA \rvx] = \sum_{i=1}^n \sum_{j=1}^n a_{ij} \Exp[\rx_i \rx_j]$.
Substitute the expectations $ \Exp[\rx_i \rx_j] = \mu_i \mu_j + \omega_{ij} $:
$$
\Exp[\rvx^\top \bA \rvx] 
= {\sum_{i=1}^n \sum_{j=1}^n a_{ij} \mu_i \mu_j} + {\sum_{i=1}^n \sum_{j=1}^n a_{ij} \omega_{ij}} = 
\bmu^\top \bA \bmu + \trace(\bA\bOmega).
$$
This completes the proof.
\end{proof}

\index{Trace trick}
\index{Quadratic expectation}

In many applications, we center the random vector so that $\bmu=\bzero$. In this case, the quadratic expectation simplifies to $\Exp[\rvx^\top\bA\rvx ] = \trace(\bA\bOmega)$.
If, further, $\bOmega=\sigma^2\bSigma$, then $\Exp[\rvx^\top\bA\rvx ] = \sigma^2\trace(\bA\bSigma)$.
To obtain an unbiased estimator of $\sigma^2$, we can choose a matrix $\bA$ such that $\trace(\bA\bSigma) = 1$.

\paragrapharrow{Other moments.}
The quantities $\Exp[\rx^s]$, $s > 0$ are called the \textit{moments} of the random variable $\rx$, while $\Exp[\abs{\rx}^s]$ are known as its \textit{absolute moments}. 
In particular, the mean $\mu=\Exp[\rx]$ and the  variance $\Exp[(\rx - \Exp\rx)^2] = \Exp[\rx^2] - (\Exp[\rx])^2$ are special cases of moments of $\rx$. 
For $1 \leq s < \infty$, $(\Exp[\abs{\rx}^s])^{1/s}$ defines a norm; in particular, the \textit{triangle inequality} holds:
\begin{equation}
(\Exp[\abs{\rx + \ry}^s])^{1/s} \leq (\Exp[\abs{\rx}^s])^{1/s} + (\Exp[\abs{\ry}^s])^{1/s}
\quad \text{for all }\rx, \ry.
\end{equation}
Similarly, \textit{H\"older's inequality} for random variables states that, for random variables $\rx, \ry$ on a common probability space, we have
\begin{subequations}
\begin{equation}
\text{H\"older's inequality:} \qquad 
\abs{\Exp[\rx \ry]} \leq \big(\Exp[\abs{\rx}^s]\big)^{1/s} \big(\Exp[\abs{\ry}^t]\big)^{1/t}, 
\quad 
\text{$s, t \geq 1$, $\frac{1}{s} + \frac{1}{t} = 1$}.
\end{equation}
The special case $s = t = 2$ yields the \textit{Cauchy--Schwarz inequality} for random variables:
\begin{equation}
\text{Cauchy--Schwarz inequality:}
\qquad 
\abs{\Exp[\rx \ry]}
 \leq \sqrt{\Exp[\abs{\rx}^2] \Exp[\abs{\ry}^2]}.
\end{equation}
This extends naturally to random matrices. 
For any random $n\times p$ matrices $\rmX$ and $\rmY$, the Frobenius norm version of  Cauchy--Schwarz reads
\begin{equation}
\text{Cauchy--Schwarz inequality:} \quad 
\Exp\left[\normf{\rmX^\top \rmY}\right] \leq \Exp\left[\normf{\rmX}^2\right]^{1/2} \Exp\left[\normf{\rmY}^2\right]^{1/2},
\end{equation}
where the inner product is defined as $\innerproduct{\rmX,\rmY} = \Exp\left[\normf{\rmX^\top \rmY}\right]$. 
Analogous matrix extensions also hold for \holders inequality.
Since the constant (deterministic) random variable $1$ has expectation $\Exp[1] = 1$, \holders inequality shows that 
\begin{equation}
\Exp[\abs{\rx}^a] = \Exp[1 \cdot \abs{\rx}^a] \leq (\Exp[\abs{\rx}^b])^{a/b}, 
\quad \text{for all $a > 0, \frac{b}{a} \geq 1$}.
\end{equation}
Therefore,
\begin{equation}
\big(\Exp[\abs{\rx}^s]\big)^{1/s} \leq \big(\Exp[\abs{\rx}^t]\big)^{1/t}, 
\quad 
\text{for all $0 < s \leq t < \infty$}.
\end{equation}
\end{subequations}

Absolute moments can be expressed in terms of tail probabilities, as shown in the following lemma.
\begin{lemma}[Absolute moments lemma]\label{lemma:mean_abs_mom}
For any random variable $\rx$ and any $s>0$,
$$
\Exp[\abs{\rx}^s] = s \int_0^\infty \prob(\abs{\rx} \geq \tau) \tau^{s-1} d\tau, \quad s > 0.
$$
\end{lemma}
\begin{proof}[of Lemma~\ref{lemma:mean_abs_mom}]
Recall that the indicator function $\indicator(\abs{\rx}^s \geq x)$ is the random variable that takes the value $1$ on the event $\{\abs{\rx}^s \geq x\}$ and $0$ otherwise. 
Using Fubini's theorem to interchange the order of integration, we compute:
\begin{align*}
\Exp[\abs{\rx}^s] &= \int_\Omega \abs{x}^s dF(x) 
= \int_\Omega \int_0^{\abs{x}^s} 1 dy dF(x) 
= \int_\Omega \int_0^\infty \indicator(\abs{x}^s \geq y) dy dF(x) \\
&= \int_0^\infty \int_\Omega \indicator(\abs{x}^s \geq y) dF(x) dy 
= \int_0^\infty \prob(\abs{\rx}^s \geq y) dy \\
&\stackrel{y\triangleq \tau^s}{\longeq} s \int_0^\infty \prob(\abs{\rx}^s \geq \tau^s) \tau^{s-1} d\tau 
= s \int_0^\infty \prob(\abs{\rx} \geq \tau) \tau^{s-1} d\tau,
\end{align*}
where the last equality follows from a change of variables.
This completes the proof.
\end{proof}

\begin{exercise}\label{exercise:exp_as_int}
Use Lemma~\ref{lemma:mean_abs_mom} to show that the expectation of a random variable $\rx$ satisfies
$$
\Exp[\rx] = \int_0^\infty \prob(\rx \geq \tau) d\tau - \int_0^\infty \prob(\rx \leq -\tau) d\tau.
$$
\textit{Hint: Express the random variable $\rx$ as $\rx = \rx \cdot \indicator\big(\rx \in [0, \infty)\big) 
+ \rx \cdot\indicator\big(-\rx \in (0, \infty )\big)$.}
\end{exercise}

\paragrapharrow{Moment generating function (MGF).}
The function $\theta \mapsto \Exp[\exp(\theta \rx)]$ is commonly referred to as the \textit{Laplace transform} or the \textit{moment generating function (MGF)} of the random variable $\rx$.
It is a powerful tool in probability theory and statistics, as it encodes information about the moments of a random variable (or random vector) and can simplify the analysis of sums of independent random variables.

\begin{definition}[Moment generating function (MGF)\index{Moment generating function}]\label{definition:momengenfunc}
Let $\rx$ be a real-valued random variable. 
The \textit{moment generating function (MGF)} of $\rx$ is defined as
\begin{subequations}
\begin{equation}
\Mom_\rx(\theta) =\Exp\left[\exp(\theta\rx)\right]: \real \to \real \cup \{\infty\}, \quad \theta \in \real
\end{equation}
Similarly, for a random vector $\rvx\in\real^n$, the MGF is defined by
\begin{equation}
\Mom_\rvx(\btheta) = \Exp\left[\exp(\btheta^\top \rvx)\right]: \real^n \to \real \cup \{\infty\}, \quad \btheta \in \real^n.
\end{equation}
On the other hand, the \textit{cumulant generating function (CGF)} of $\rx$ is the logarithm of its MGF:
\begin{equation}
\Cum_\rx(\theta) = \ln \Mom_\rx(\theta) .
\end{equation}
\end{subequations}
\end{definition}
Below are some key properties of moment generating functions:
\begin{enumerate}

\item \textit{Moments from MGF.} If the MGF exists in an open interval around zero, then all the moments of $\rx$ exist and can be computed by differentiating the MGF and evaluating at $\theta=0$:
$$
\Exp[\rx^k] = \frac{d^k}{d\theta^k} \Mom_\rx(0),
$$
where $k$ is a positive integer representing the order of the moment.

\item \textit{Uniqueness.} If two random variables have the same MGF in an open interval containing zero, then they have the same distribution. That is, if $\Mom_\rx(\theta) = \Mom_\ry(\theta)$ for all $\theta$ in an open interval around zero, then the cumulative distribution functions (CDFs) of $\rx$ and $\ry$ are identical ($F_\rx = F_\ry$).

\item \textit{Product of independent variables.} For independent random variables $\rx$ and $\ry$, the MGF of their sum $\rx + \ry$ is the product of their individual MGFs:
$$
\Mom_{\rx+\ry}(\theta) = \Mom_\rx(\theta) \Mom_\ry(\theta),
\quad \text{if $\rx$ and $\ry$ are independent.}
$$

\item On the other hand, the CGF has  \textit{additivity} for independent random variables
$$
\Cum_{\rx+\ry}(\theta) = \Cum_{\rx}(\theta)+\Cum_{\ry}(\theta), 
\quad \text{if $\rx$ and $\ry$ are independent.}
$$
This makes cumulants particularly useful in limit theorems and asymptotic analysis. 
The CGF often simplifies the analysis of sums of independent random variables. 
\end{enumerate}

These properties make the moment generating function a valuable tool for analyzing sums of independent random variables, finding moments of distributions, and proving distributional equivalences. However, it should be noted that not all distributions have moment generating functions; for example, the Cauchy distribution does not have an MGF because its moments do not exist. In such cases, other tools like characteristic functions may be used instead.

\subsection{Some Concentration Inequalities}

Markov's inequality and Chebyshev's inequality are fundamental concentration inequalities in probability theory.
They provide upper bounds on the probability that a random variable deviates from its expected value, with Markov's inequality applying to nonnegative variables and Chebyshev's inequality extending this concept by considering the variance of any random variable.

To be more specific,  the function $\tau \mapsto \prob(\abs{\rx} \geq \tau)$ is called the \textit{tail} (or \textit{tail probability}) of
the random variable $\rx$. 
Markov's inequality offers a simple yet often useful way to bound this tail using expectations.
\begin{theoremHigh}[Markov's inequality\index{Markov's inequality}]\label{theorem:markov-inequality}
Let $\rx$ be a \textbf{nonnegative} random variable. Then, given any $\tau >0$, we have 
$$
\prob(\rx \geq \tau) \leq \frac{\Exp[\rx]}{\tau}.
$$
\end{theoremHigh}

\begin{proof}[of Theorem~\ref{theorem:markov-inequality}]
Since $\rx$ is nonnegative, we have the pointwise inequality $0\leq \tau \indicator\{\rx \geq \tau\} \leq \rx$. 
This implies $\Exp[\tau \indicator\{\rx\geq \tau\}] \leq \Exp[\rx]$. We also have 
$$
\Exp\left[\tau \indicator\{\rx\geq \tau\}\right] = \tau\Exp\left[ \indicator\{\rx\geq \tau\}\right] = \tau \left(1\cdot \prob(\rx\geq \tau) + 0\cdot \prob(\rx< \tau) \right) = \tau \cdot \prob(\rx\geq \tau) \leq \Exp[\rx].
$$
This completes the proof.

Alternatively, let $f(x)$ denote the probability density function for $\rx$. Then,
$$
\Exp[\rx] = \int_{0}^{\infty} x f(x) \, dx \geq \int_{\tau}^{\infty} x f(x) \, dx \geq \int_{\tau}^{\infty} \tau f(x) \, dx = \tau \prob(\rx \geq \tau),
$$
from which the result follows.
\end{proof}

\begin{subequations}
Let $\rx$ be an arbitrary random variable (not necessarily nonnegative). 
By applying Markov's inequality to the nonnegative random variable  $\abs{\rx}$, we obtain a tail bound for $\rx$:
\begin{equation}
\prob(\abs{\rx} \geq \tau) \leq \frac{\Exp[\abs{\rx}]}{\tau}, \quad \text{for all } \tau > 0.
\end{equation}
Similarly, for $\theta > 0$, we obtain
\begin{equation}
\prob(\rx \geq \tau) 
= \prob\big(\exp(\theta \rx)  \geq \exp(\theta \tau)\big) 
\leq  \frac{\Exp[\exp(\theta \rx)]}{\exp(\theta \tau)}, 
\quad \text{for all } \tau \in \real.
\end{equation}
An important consequence is the following generalized Markov-type bound: for any $s > 0$,
\begin{equation}
\prob(\abs{\rx} \geq \tau) = \prob(\abs{\rx}^s \geq \tau^s) \leq \tau^{-s} \Exp[\abs{\rx}^s], \quad \text{for all } t > 0.
\end{equation}
The special case $s = 2$ yields  \textit{Chebyshev's inequality}, provided the variance exists.
\end{subequations}

\begin{theoremHigh}[Chebyshev's inequality\index{Chebyshev's inequality}]
Let $\rx$ be a random variable with finite mean $\Exp[\rx]< \infty$. 
Then, for any $\tau >0$, we have 
$$
\prob(\abs{\rx - \Exp[\rx]} \geq \tau) \leq \frac{\Var[\rx]}{\tau^2}.
$$
\end{theoremHigh}
This follows immediately by applying Markov's inequality to the nonnegative random variable $\ry \triangleq (\rx - \Exp[\rx])^2$.

The \textit{Chernoff method} leverages the moment generating function to derive exponentially decaying bounds on tail probabilities. The resulting inequalities are known as \textit{Chernoff bounds}.
\begin{theoremHigh}[Chernoff bound\index{Chernoff bound}]\label{theorem:chernoff_bound}
Let $\rx$ be a random variable. Then, given any $\theta>0$
$$
\prob(\rx \geq \tau) \leq \min_{\theta > 0} \left( \exp({-\theta \tau}) \Mom_\rx(\theta) \right) .
$$
\end{theoremHigh}
This can be obtained by applying  Markov's inequality  to $\exp(\theta\rx)$: $\prob(\rx \geq \tau) = \prob(\exp(\theta\rx)\geq \exp(\theta \tau)) \leq \Mom_\rx(\theta)\exp({-\theta \tau})$ and taking the infimum.

A powerful extension of this idea is \textit{Cramér's theorem}, which provides exponential tail bounds for sums of independent random variables.
\begin{theoremHigh}[Cramér's theorem \citep{cramer1994nouveau}\index{Cramér's theorem}]\label{theorem:cramer_theo}
Let $\rx_1, \rx_2, \ldots, \rx_n$ be a sequence of independent random variables, with cumulant generating functions $\Cum_{\rx_i}$, $i \in \{1,2,\ldots,n\}$. Then, for any $\tau > 0$,
$$
\prob\left(\sum_{i=1}^n \rx_i \geq \tau\right) 
\leq 
\exp\left(\min_{\theta > 0} \left\{-\theta \tau + \sum_{i=1}^n \Cum_{\rx_i}(\theta)\right\}\right).
$$
\end{theoremHigh}

\begin{proof}[of Theorem~\ref{theorem:cramer_theo}]
We use the Chernoff method.
For any $\theta > 0$, Markov's inequality (Theorem~\ref{theorem:markov-inequality}) and independence yield
\begin{align*}
\prob\left(\sum_{i=1}^n \rx_i \geq \tau\right) 
&= \prob\left(\exp\left(\theta \sum_{i=1}^n \rx_i\right) \geq \exp(\theta \tau)\right) \leq \exp({-\theta \tau}) \Exp\left[\exp\left(\theta \sum_{i=1}^n \rx_i\right)\right] \\
&= \exp({-\theta \tau}) \Exp\left[\prod_{i=1}^n \exp(\theta \rx_i)\right] = \exp({-\theta \tau}) \prod_{i=1}^n \Exp\left[\exp(\theta \rx_i)\right] \\
&= \exp({-\theta \tau}) \prod_{i=1}^n \exp(\Cum_{\rx_i}(\theta)) = \exp\left(-\theta \tau + \sum_{i=1}^n \Cum_{\rx_i}(\theta)\right).
\end{align*}
Since this holds for every $\theta>0$, taking the infimum over $\theta > 0$ yields the desired bound.
\end{proof}

Cramér's theorem provides an exponential upper bound for the tail probability of the sum of independent random variables. In contrast, Hoeffding's inequality gives a sharper bound when the summands are almost surely bounded---a common scenario in empirical process theory and learning theory.

\begin{lemma}[Hoeffding's lemma]\label{lemma:hoeffdin_lemm}

Let $\rx$ be a real-valued random variable such that  $\Exp[\rx]=0$ and $\rx\in[a,b]$ almost surely. Then for any $\theta\in\real$
$$
\Mom_\rx(\theta) \equiv \Exp\big[\exp(\theta \rx)\big]\leq \exp\!\Big(\frac{\theta^{2}(b-a)^{2}}{8}\Big).
$$
\end{lemma}
\begin{proof}[of Lemma~\ref{lemma:hoeffdin_lemm}]
Let $t \triangleq \theta(b-a) $ and $\ry\triangleq\frac{\rx-a}{\,b-a\,}\in[0,1]$.
Then $\rx=a+(b-a)\ry$ and $\Exp[\ry]=\dfrac{\Exp[\rx]-a}{b-a}=\dfrac{-a}{b-a}\triangleq \mu$ because $\Exp[\rx]=0$, where $\mu \in[0,1]$ because $a\leq \Exp[\rx]\leq b$.
The exponential function is convex, so for any $\lambda\in[0,1]$ and any real $s$,
$$
\exp(\lambda s) \leq (1-\lambda) + \lambda e^{s}.
$$
Applying this inequality with $\lambda=\ry$ and $s=t$, where $\ry=\frac{\rx-a}{\,b-a\,}\in[0,1]$ and $t=\theta(b-a)$, yields
$$
\exp(\theta \rx)
=\exp(\theta a) \exp(t\ry) 
\leq \exp(\theta a) \big(1-\ry + \ry e^{t}\big).
$$
Taking expectations on both sides gives
$$
\Exp[\exp(\theta \rx)] \leq \exp(\theta a)\big(1-\Exp[\ry]+\Exp[\ry] e^{t}\big)
= \exp(\theta a)\big(1-\mu + \mu e^{t}\big).
$$
Take logarithms and define
$
g(t)\triangleq\ln\big(1-\mu+\mu e^{t}\big).
$
Then the cumulant generating function satisfies
\begin{equation}\label{equation:hoeffdin_lemm_pv1}
\Cum_{\rx}(\theta) \equiv \ln \Exp[\exp(\theta \rx)] \leq \theta a + g(t),\qquad t=\theta(b-a).
\end{equation}
We analyze $g(t)$ via its second derivative:
$$
g''(t)=\frac{\mu(1-\mu)e^{t}}{\big(1-\mu+\mu e^{t}\big)^{2}}.
$$
To bound this, observe the elementary inequality $4uv\le(u+v)^2$. 
Setting $u=\sqrt{\mu e^{t}}$ and $v=\sqrt{1-\mu}$), we obtain
$$
\mu(1-\mu)e^{t} = (uv)^2 \leq \Big(\frac{u^2+v^2}{2}\Big)^2
= \frac{\big(\mu e^{t}+1-\mu\big)^2}{4}.
$$
Dividing both sides by $(1-\mu+\mu e^{t})^{2}$ yields
$$
g''(t)\leq \frac14.
$$
Since $g''(t)\leq \tfrac14$ for all $t$, and $g(0)=\ln(1)=0$, $g'(0)=\dfrac{\mu e^{0}}{1-\mu+\mu e^{0}}=\mu$, we can integrate twice.
Moreover, note that the derivative with respect to $\theta$ of the right-hand side in \eqref{equation:hoeffdin_lemm_pv1} at $\theta=0$ vanishes:
$$
\left.\frac{d}{d\theta}\big(\theta a+g(\theta(b-a))\big)\right|_{\theta=0}
= a + g'(0)\,(b-a) = a+\mu(b-a)=0.
$$
Therefore, the Taylor expansion of the right-hand side in \eqref{equation:hoeffdin_lemm_pv1} about $t=0$ with $g''(t)\leq 1/4$ yields
$$
\ln\Exp[\exp(\theta \rx)] 
\leq  \theta a + g(t)
\leq  \frac{1}{2}\cdot\frac{1}{4}\,t^{2}
= \frac{t^{2}}{8}
= \frac{\theta^{2}(b-a)^{2}}{8}.
$$
Exponentiating both sides gives the desired bound:
\end{proof}

Using the above lemma, we prove Hoeffding's inequality, which quantifies how tightly the sum of bounded independent random variables concentrates around its mean. 
\begin{theoremHigh}[Hoeffding's inequality \citep{hoeffding1963probability}\index{Hoeffding's inequality}]\label{theorem:hoeffding_ineq}
Let $ \rx_1, \rx_2, \dots, \rx_n $ be independent random variables such that each $ \rx_i $ is almost surely bounded in the interval $ [a_i, b_i] $, where $ a_i $ and $ b_i $ are known constants.
Let $ \rs_n \triangleq \sum_{i=1}^n \rx_i $ be the sum of these random variables, and let $ \Exp[\rs_n] $ be its mean value.
Then, for any $ \tau > 0 $, Hoeffding's inequality states:
$$
\prob(\rs_n - \Exp[\rs_n] \geq \tau) \leq \exp\left( -\frac{2\tau^2}{\sum_{i=1}^n (b_i - a_i)^2} \right).
$$
Similarly, for the lower tail:
$$
\prob(\rs_n - \Exp[\rs_n] \leq -\tau) \leq \exp\left( -\frac{2\tau^2}{\sum_{i=1}^n (b_i - a_i)^2} \right).
$$
Consequently, the two-sided version holds:
$$
\prob(\abs{\rs_n - \Exp[\rs_n]} \geq \tau) \leq 2\exp\left( -\frac{2\tau^2}{\sum_{i=1}^n (b_i - a_i)^2} \right).
$$
\end{theoremHigh}
The bound decays  exponentially in $ \tau^2 $, which implies that large deviations from the mean are extremely unlikely.
The denominator $ \sum (b_i - a_i)^2 $ captures the total variability of the random variables: the wider the individual bounds $ [a_i, b_i]$, the weaker (i.e., looser) the concentration around the mean.
\begin{proof}[of Theorem~\ref{theorem:hoeffding_ineq}]
We use the Chernoff method.
For any $ \theta > 0 $, by Markov's inequality (Theorem~\ref{theorem:markov-inequality}) and the independence of $\rx_i$:
$$
\prob(\rs_n - \Exp[\rs_n] \geq \tau) 
= \prob(e^{\theta(\rs_n - \Exp[\rs_n])} \geq \exp(\theta \tau)) 
\leq \frac{\Exp[e^{\theta(\rs_n - \Exp[\rs_n])}]}{\exp(\theta \tau)} 
=\exp({-\theta \tau}) \prod_{i=1}^n \Exp[e^{\theta(\rx_i - \Exp[\rx_i])}].
$$
To proceed, we bound the MGF of $ \rx_i - \Exp[\rx_i] $ for each $\rx_i$.
Let $ \ry_i \triangleq  \rx_i - \Exp[\rx_i] $. Then $ \ry_i \in \big[a_i - \Exp[\rx_i], b_i - \Exp[\rx_i]\big] \triangleq [c_i, d_i] $, where $ c_i \triangleq a_i - \Exp[\rx_i] $, $ d_i \triangleq b_i - \Exp[\rx_i] $,  $ d_i - c_i = b_i - a_i $, and $ \Exp[\ry_i] = 0 $. 
Applying Hoeffding's lemma (Lemma~\ref{lemma:hoeffdin_lemm}) to each $ \ry_i $:
$$
\Exp\big[\exp(\theta(\rx_i - \Exp[\rx_i]))\big] 
\leq \exp\left( \frac{\theta^2 (b_i - a_i)^2}{8} \right),
$$
whence we have
$$
\prob(\rs_n - \Exp[\rs_n] \geq \tau) 
\leq \exp(-\theta \tau) \prod_{i=1}^n \Exp\big[\exp(\theta(\rx_i - \Exp[\rx_i]))\big] 
\leq  \exp\left( -\theta \tau + \frac{\theta^2}{8}\sum_{i=1}^n (b_i - a_i)^2 \right).
$$
The optimal choice $\theta = \frac{4\tau}{ \sum_{i=1}^n (b_i - a_i)^2} $ in the above yields the upper tail. 
The lower-tail inequality follows analogously by applying the same argument to $ -\rx_i $.
Finally, the two-sided bound follows from the union bound:
$$
\prob(\abs{\rs_n - \Exp[\rs_n]} \geq \tau) \leq \prob(\rs_n - \Exp[\rs_n] \geq \tau) + \prob(\rs_n - \Exp[\rs_n] \leq -\tau) \leq 2\exp\left( -\frac{2\tau^2}{\sum_{i=1}^n (b_i - a_i)^2} \right).
$$
This completes the proof.
\end{proof}

\begin{exercise}
Let $\Exp[\rx_i]=0$ for all $i\in\{1,2,\ldots,n\}$ in Theorem~\ref{theorem:hoeffding_ineq}. Use Cramér's theorem (Theorem~\ref{theorem:cramer_theo}) to complete an alternative proof.
\end{exercise}

If the variables $ \rx_i \in \{0,1\} $ are Bernoulli random variables, then $ a_i = 0, b_i = 1 $, so $ b_i - a_i = 1 $, and $ \sum_{i=1}^n (b_i - a_i)^2 = n $.
Let $ \overline{\rx} \triangleq \frac{1}{n} \sum \rx_i $ and $ \mu \triangleq \Exp[\overline{\rx}] $. 
Then Hoeffding's inequality gives
$$
\prob(\abs{\overline{\rx} - \mu} \geq \varepsilon) \leq 2\exp(-2n\varepsilon^2).
$$
This form is widely used in learning theory---for example, in proving the uniform convergence of empirical risk minimization.

A \textit{Rademacher variable} (also known as a \textit{symmetric Bernoulli variable}) is a variant of the Bernoulli random variable: it   is a random variable $\epsilon$ that takes the values $+1$ and $-1$ with equal probability. A \textit{Rademacher sequence} $\bepsilon$ is a vector of independent Rademacher variables. 
For such sequences, we obtain the following version of Hoeffding's inequality.

\begin{corollary}[Hoeffding's inequality for Rademacher sequences\index{Rademacher sequence}]\label{corollary:hoeffd_rademac}
Let $\ba \in \real^n$ and $\bepsilon = \{\epsilon_1, \epsilon_2, \ldots, \epsilon_n\}$ be a Rademacher sequence. Then, for $\mu > 0$,
$$
\prob\left(\abs{\sum_{i=1}^n \epsilon_i a_i} \geq \normtwo{\ba} \mu \right) \leq 2 \exp(-\frac{\mu^2 }{2}).
$$
\end{corollary}
\begin{proof}[of Corollary~\ref{corollary:hoeffd_rademac}]
The random variable $ \epsilon_i a_i$ has mean zero and is bounded in absolute value by $\abs{a_i}$. 
The stated inequality therefore follows directly from Hoeffding's inequality.
\end{proof}

\textit{Bernstein's inequality} provides a useful generalization of Hoeffding's inequality  to sums of  independent random variables 
that are either bounded or even unbounded, while also incorporating information about their variance (or higher moments). We begin with the version stated below and then derive several useful variants as corollaries.

\begin{theoremHigh}[Bernstein's inequality \citep{bernstein1924modification}\index{Bernstein's inequality}]\label{theorem:bernst_ineq}
Let $\rx_1, \rx_2, \ldots, \rx_n$ be independent zero-mean random variables (i.e., $\Exp[\rx_i]=0$ for all $i$) such that, for all integers $k \geq 2$,
$$
\Exp[\abs{\rx_i}^k] 
\leq k! C^{k-2} \sigma_i^2 / 2, \quad \text{for all } i \in \{1,2,\ldots,n\} 
$$
for some constants $C > 0$ and $\sigma_i > 0$, $i \in \{1,2,\ldots,n\}$. Then, for all $\tau > 0$,
$$
\prob\left(\abs{\sum_{i=1}^n \rx_i} \geq \tau\right) \leq 2 \exp\left(-\frac{\tau^2 / 2}{\sigma^2 + C\tau}\right), 
$$
where $\sigma^2 \triangleq \sum_{i=1}^n \sigma_i^2$.
\end{theoremHigh}
\begin{proof}[of Theorem~\ref{theorem:bernst_ineq}]
Using the well-known Taylor expansion $\exp(z) = \sum_{i=0}^{\infty} \frac{z^i}{i!}$ for all $z\in\real$ around 0 and the fact that $\Exp[\rx_i] = 0$, we have
$$
\Exp[\exp(\theta \rx_i)] 
= 1 + \theta \Exp[\rx_i] + \sum_{k=2}^\infty \frac{\theta^k \Exp[\rx_i^k]}{k!}
 = 1 + \frac{\theta^2 \sigma_i^2}{2} \sum_{k=2}^\infty \frac{\theta^{k-2} \Exp[\rx_i^k]}{k! \sigma_i^2 / 2},
$$
For each $i\in\{1,2,\ldots,n\}$, define $
F_i(\theta) \triangleq \sum_{k=2}^\infty \frac{\theta^{k-2} \Exp[\rx_i^k]}{k! \sigma_i^2 / 2}$.
Then,
$$
\Exp[\exp(\theta \rx_i)] = 1 + \frac{\theta^2 \sigma_i^2 F_i(\theta)}{2}
 \leq \exp(\frac{\theta^2 \sigma_i^2 F_i(\theta)}{2}),
$$
where we used the inequality $1+x\leq \exp(x)$ for all $x\in\real$.
Now let $F(\theta) \triangleq \max_{i \in \{1,2,\ldots,n\}} F_i(\theta)$. 
Then by the fact that $\sigma^2 = \sum_{i=1}^n \sigma_i^2$ and Cramér's theorem (Theorem~\ref{theorem:cramer_theo}), we obtain
$$
\prob\left(\sum_{i=1}^n \rx_i 
\geq \tau\right) \leq \min_{\theta > 0} \exp\left(- \theta \tau + \frac{\theta^2 \sigma^2 F(\theta)}{2} \right) 
\leq \min_{0 < C\theta < 1} \exp\left(- \theta \tau + \frac{\theta^2 \sigma^2 F(\theta)}{2} \right),
$$
where we suppose $C\theta < 1$. 
Combining the inequality $\Exp[\rx_i^k] \leq \Exp[\abs{\rx_i}^k]$ and the assumption of the bound for $\Exp[\abs{\rx_i}^k]$ yields
$$
F_i(\theta) \leq \sum_{k=2}^\infty \frac{\theta^{k-2} \Exp[\abs{\rx_i}^k]}{k! \sigma_i^2 / 2} \leq \sum_{k=2}^\infty \frac{(C\theta)^{k-2}}{1 - C\theta} = \frac{1}{1 - C\theta}.
$$
Hence, $F(\theta) \leq (1 - C\theta)^{-1}$, and therefore
$$
\prob\left(\sum_{i=1}^n \rx_i \geq \tau\right) \leq \min_{0 < C\theta < 1} \exp\left(- \theta \tau + \frac{\theta^2 \sigma^2}{2(1 - C\theta)} \right).
$$
The choice $\theta = \tau / (\sigma^2 + C\tau)$ satisfies $C\theta < 1$, which yields
\begin{align*}
\prob\left(\sum_{i=1}^n \rx_i \geq \tau\right) 
\leq \exp\left(-\frac{\tau^2 / 2}{\sigma^2 + C\tau}\right).
\end{align*}
Applying the same argument to $-\rx_i$ gives an identical bound for $\prob\left(\sum_{i=1}^n \rx_i \leq -\tau\right) $.
Combining both tails via the union bound completes the proof.
\end{proof}

We present two important consequences of Bernstein's inequality.
The first is Bernstein's inequality for bounded random variables, which appears more frequently in the literature.
\begin{corollary}[Bernstein's inequality for bounded random variables]\label{corollary:berns_bound}
Let $\rx_1, \rx_2, \ldots, \rx_n$ be independent zero-mean random variables   such that $\abs{\rx_i} \leq B$ almost surely for $i \in \{1,2,\ldots,n\}$, where  $B > 0$ is a constant. 
Suppose further that  $\Exp[\rx_i^2] \leq \sigma_i^2$ for some constants $\sigma_i > 0$, $i \in \{1,2,\ldots,n\}$. Then, for all $\tau > 0$,
$$
\prob\left(\abs{\sum_{i=1}^n \rx_i} \geq \tau\right) 
\leq 2 \exp\left(-\frac{\tau^2 / 2}{\sigma^2 + B\tau / 3}\right), 
$$
where $\sigma^2 \triangleq \sum_{i=1}^n \sigma_i^2$.
\end{corollary}
\begin{proof}[of Corollary~\ref{corollary:berns_bound}]
For $k = 2$, the moment condition in Theorem~\ref{theorem:bernst_ineq} is satisfied by assumption. 
So let $k \in \naturalset$, $k \geq 3$. Since then $k! \geq 3 \cdot 2^{k-2}$, we obtain
$$
\Exp[\abs{\rx_i}^k] = \Exp[\abs{\rx_i}^{k-2} \rx_i^2] \leq B^{k-2} \sigma_i^2 \leq \frac{k! B^{k-2}}{k!} \sigma_i^2 \leq \frac{k! B^{k-2}}{2 \cdot 3^{k-2}} \sigma_i^2. 
$$
Thus, the moment condition of Theorem~\ref{theorem:bernst_ineq} is satisfied for all $k \geq 2$ with  $C = B/3$ and the same $\sigma_i$. 
The desired inequality follows directly from Theorem~\ref{theorem:bernst_ineq}.
\end{proof}

As a second consequence, we state Bernstein's inequality for sub-exponential random variables (see Definition~\ref{definition:subgau_concen_ineq} for more details).
\begin{corollary}[Bernstein's inequality for sub-exponential random variables]\label{corollary:berns_subexp}
Let $\rx_1, \rx_2, \ldots, \rx_n$ be independent zero-mean sub-exponential random variables; that is, $\prob(\abs{\rx_i} \geq \tau) \leq \Phi \exp(-\xi \tau)$ for some constants $\Phi, \xi > 0$ for all $\tau > 0$, $i \in \{1,2,\ldots,n\}$. Then, for all $\tau > 0$,
$$
\prob\left(\abs{\sum_{i=1}^n \rx_i} \geq \tau\right) 
\leq 2 \exp\left(-\frac{(\xi \tau)^2 / 2}{2 \Phi n + \xi \tau}\right). 
$$
\end{corollary}
\begin{proof}[of Corollary~\ref{corollary:berns_subexp}]
Using the absolute moments lemma (Lemma~\ref{lemma:mean_abs_mom}), for $k \in \naturalset$, $k \geq 2$,
\begin{align*}
\Exp[\abs{\rx_i}^k ]
&= k \int_0^\infty \prob(\abs{\rx_i} \geq x) x^{k-1} dx 
\leq \Phi k \int_0^\infty \exp(-\xi x) x^{k-1} dx \\
&\stackrel{u=\xi x}{\longeq} \Phi k \xi^{-k} \int_0^\infty e^{-u} u^{k-1} du = \Phi k! \xi^{-k} = k! \xi^{-(k-2)} \frac{2 \Phi \xi^{-2}}{2}, 
\end{align*}
where we have used the Gamma function  $\Gamma(x)=\int_{0}^{\infty} t^{x-1}\exp(-t)dt$ and  the fact that $\Gamma(k) = (k-1)!$ for any integer $k$. Hence, the moment condition in Theorem~\ref{theorem:bernst_ineq} holds with $C = \xi^{-1}$ and $\sigma_i^2 = 2 \Phi \xi^{-2}$. The result follows again from Theorem~\ref{theorem:bernst_ineq}. 
\end{proof}

\index{Integration by parts}
Recall that the \textit{Gamma function} is defined for $y\geq 0$ by
\begin{equation}\label{equation:gamma_func}
\Gamma(y) = \int_{0}^{\infty} t^{y-1} e^{-t} dt,  \qquad y\geq 0.
\end{equation}
Utilizing integration by parts $\int_{a}^{b} u(t) v^\prime(t) dt = u(t)v(t)|_a^b - \int_a^b u^\prime(t) v(t)dt$, where $u(t) =t^{y-1}$ and $v(t)=-e^{-t}$, we derive 
$$
\begin{aligned}
	\Gamma(y) &= -t^{y-1}e^{-t}|_0^{\infty} - \int_{0}^{\infty} (y-1)t^{y-2}(-e^{-t}) dt \\
	&= 0+ (y-1) \int_{0}^{\infty}t^{y-2}e^{-t}dt 
	= (y-1)\Gamma(y-1).
\end{aligned}
$$ 
This recurrence relation shows that when $y$ is a positive integer, $\Gamma(y) = (y-1)!$.

\subsection{Gaussian Distribution and Properties}
The \textit{Gaussian}  (or \textit{normal}) random variable is a widely preferred choice in statistical modeling and signal processing. This preference arises from a fundamental principle: the sum of independent random variables tends to follow a Gaussian distribution---a result formalized by the central limit theorem (see Theorem~\ref{thm:clt} below). Consequently, any measurement or quantity formed by the additive combination of many independent, unrelated components is likely to be approximately Gaussian. For this reason, noise in fields such as signal processing and engineering is commonly modeled using a Gaussian distribution due to its desirable mathematical properties.

\subsection*{Gaussian Random Variables}
\index{Gaussian distribution}
\begin{definition}[Gaussian or normal distribution]\label{definition:gaussian_distribution}
A random variable $\rx$ is said to follow the \textit{Gaussian distribution} (a.k.a., the \textit{normal distribution}) with mean $\mu$ and variance   $\sigma^2>0$, denoted  $\rx \sim \normal(\mu,\sigma^2)$ \footnote{Note if two random variables $\ra$ and $\rb$ have the same distribution, then we write $\ra \sim \rb$.}, if 
$$
f(x; \mu,\sigma^2)=\frac{1}{\sqrt{2\pi\sigma^2}} \exp \left\{-\frac{1}{2\sigma^2 }(x-\mu)^2 \right\}
=\sqrt{\frac{\tau}{2\pi}}\exp \left\{ -\frac{\tau}{2}(x-\mu)^2 \right\},
$$
where $\tau$ is also known as the \textit{precision} of the  distribution.
The mean,  variance, and MGF of $\rx \sim \normal( \mu,\sigma^2)$ are given by 
$$
\Exp[\rx] = \mu, \qquad \Var[\rx] =\sigma^2=\tau^{-1}, 
\qquad 
\Mom_\rx(\theta) =  \exp\{\theta\mu + \theta^2\sigma^2/2\},
$$
The cumulative distribution function (c.d.f., CDF) of $\rx$ is 
$$
F(x; \mu, \sigma^2) = \prob(\rx<x) = \frac{1}{\sqrt{2\pi\sigma^2}} \int_{-\infty}^x \exp\left\{-\frac{1}{2\sigma^2 }(z-\mu)^2 \right\}dz.
$$
In particular, we denote $\Phi(y) = \int_{-\infty}^{y} \normal(u\mid 0,1)du= \frac{1}{\sqrt{2\pi}} \int_{-\infty}^{y} \exp(-\frac{u^2}{2}) du $ as the cumulative distribution function of $\normal(0,1)$, the \textit{standard normal distribution}. 
\end{definition}

\begin{SCfigure}
\centering
\includegraphics[width=0.5\textwidth]{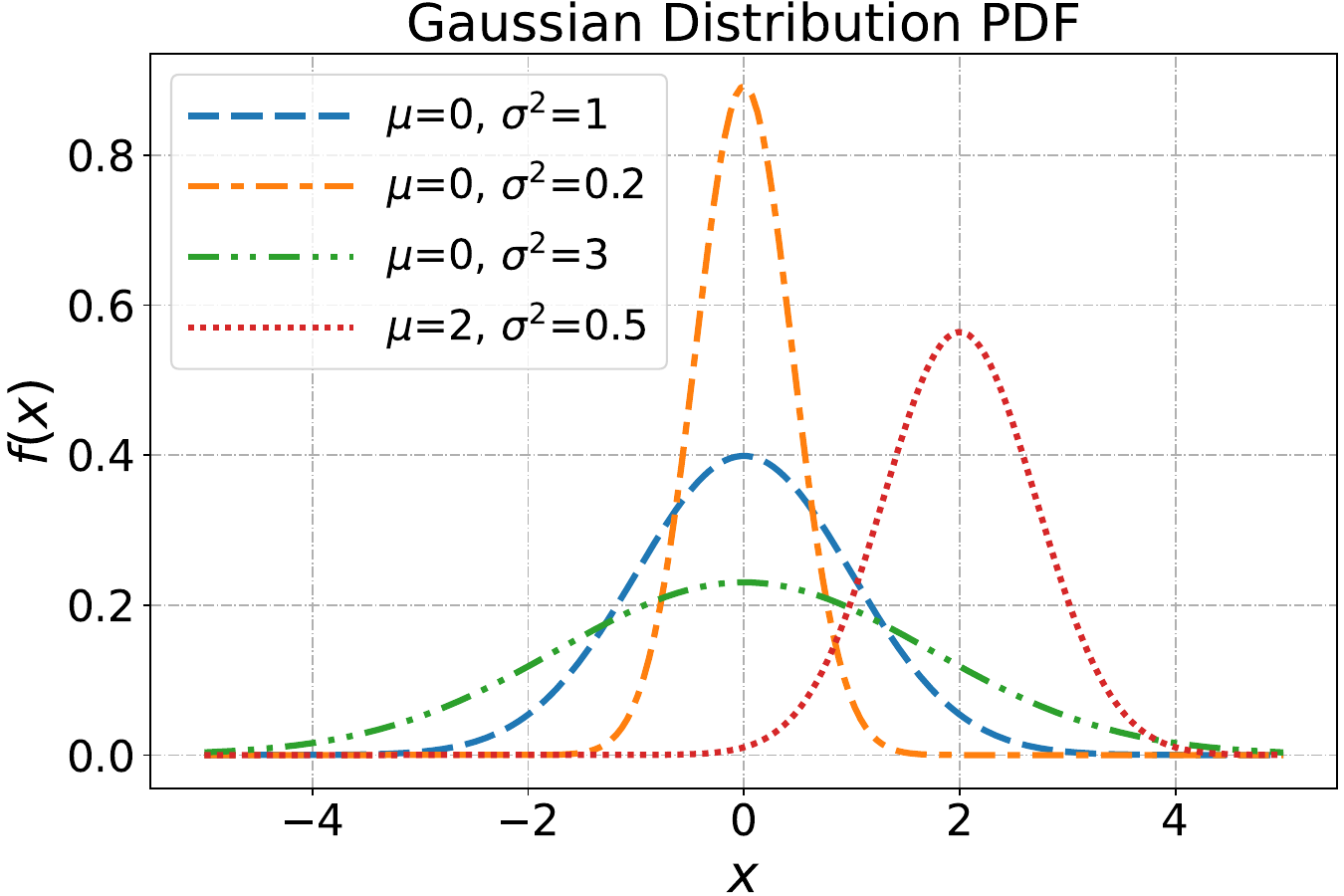}
\caption{Gaussian probability density functions for different values of the mean and variance parameters $\mu$ and $\sigma^2$.}
\label{fig:dists_gaussian}
\end{SCfigure}

Figure~\ref{fig:dists_gaussian} illustrates how different values of the parameters $\mu, \sigma^2$ affect the shape of the Gaussian distribution.

As noted above, the central limit theorem states that the distribution of the average of a large number of random variables with finite variances converges to a Gaussian distribution.
\begin{theoremHigh}[Central limit theorem (CLT)\index{Central limit theorem}]\label{thm:clt}
Let $\rx_1, \rx_2, \rx_3, \ldots$ be a sequence of i.i.d. random variables with mean $\mu$ and finite variance $\sigma^2<\infty$. We define the sequence of averages $\ry_1, \ry_2, \ry_3, \ldots$ as
$$
\ry_n \triangleq \frac{1}{n} \sum_{i=1}^{n} \rx_i.
$$
Then the normalized sequence $\rz_1, \rz_2, \rz_3, \ldots$
$$
\rz_n \triangleq \sqrt{n}(\ry_n - \mu)
$$
converges in distribution to a Gaussian random variable with mean $0$ and variance $\sigma^2$. 
That is, for any $x \in \real$,
$$
\lim_{n \to \infty} f_{\rz_n}(x) = \frac{1}{\sqrt{2\pi}\sigma} e^{-\frac{x^2}{2\sigma^2}}.
$$
\end{theoremHigh}

For sufficiently large $n$, the above theorem also implies that the sample  average $\ry_n$ is approximately Gaussian with mean $\mu$ and variance $\sigma^2/n $. 

While the product of two Gaussian variables remains analytically challenging, the sum of Gaussian variables is itself Gaussian and can be fully characterized by its mean and variance.
\begin{remark}[Sum of Gaussians]
Let $\rx\sim \normal{\mu_x, \sigma_x^2}$ and $\ry\sim \normal{\mu_y, \sigma_y^2}$ be two Gaussian random variables.
\begin{itemize}
\item If $\rx$ and $\ry$ are uncorrelated, then
$$
\rx+\ry \sim \normal(\mu_x+\mu_y, \sigma_x^2+\sigma_y^2).
$$
\item If $\rx$ and $\ry$ have correlation coefficient $\rho$, then
$$
\rx+\ry \sim \normal(\mu_x+\mu_y, \sigma_x^2+\sigma_y^2+2\rho\sigma_x\sigma_y).
$$
\end{itemize}
\end{remark}

An important property of Gaussian random variables is that affine transformations (i.e., scaling and shifting) preserve the Gaussian form.
\begin{lemma}[Linear transformation of Gaussian variables]\label{lem:lintrans_uni_gaus}
Let $\rx$ ba a Gaussian random variable with mean $\mu$ and variance $\sigma^2$. 
Then, for any $a, b \in \real$, the random variable
\begin{equation}
\ry \triangleq a\rx + b
\end{equation}
is also Gaussian  with mean $a\mu + b$ and variance $a^2\sigma^2$.
\end{lemma}
The proof follows directly from the linear transformation lemma (Lemma~\ref{lemma:lin_tran_prob}).

Gaussian distributions exhibit strong concentration properties. A fundamental result in this context is the following tail bound for a centered Gaussian random variable.
\begin{exercise}[Chernoff bound for centered Gaussian]\label{exercise:chernoff_gausvari}
For $\rx \sim \normal(0, \sigma^2)$, show that
$$
\prob(\abs{\rx} \geq \tau) \leq 2 e^{-\frac{\tau^2}{2\sigma^2}},
\quad \text{for all } \tau>0.
$$
\textit{Hint: Apply the Chernoff bound (Theorem~\ref{theorem:chernoff_bound})}
\end{exercise}

Although the Chernoff bound provide a valid upper bound on $\prob(\abs{\rx} \geq \tau)$ for a centered Gaussian, the bound is not tight.
A sharper bound can be derived as follows.
\begin{lemma}[Tighter tail bound for a centered Gaussian]\label{lemma:bd_cen_gaus}
For $\rx \sim \normal(0, \sigma^2)$, it holds that
$$
\prob(\abs{\rx} \geq \tau) 
\leq
\begin{cases}
&e^{-\frac{\tau^2}{2\sigma^2}};\\
&\frac{\sigma \sqrt{2}}{\tau \sqrt{\pi}} \exp\left( -\frac{\tau^2}{2\sigma^2} \right),
\end{cases}
\quad \text{for all } \tau>0.
$$
\end{lemma}
Note that the latter inequality is called \textit{Mill's inequality.}
\begin{proof}[of Lemma~\ref{lemma:bd_cen_gaus}]
By definition, $\prob(\abs{\rx} \geq \tau) = \frac{2}{\sqrt{2\pi\sigma^2}} \int_\tau^\infty e^{-x^2/2\sigma^2}dx $.
Applying the change of variables:
\begin{align*}
\int_{\tau}^{\infty} e^{-x^{2}/2\sigma^2} dx 
&= \int_{0}^{\infty} e^{-(x+\tau)^{2}/2\sigma^2} dx 
= e^{-\tau^{2}/2\sigma^2} \int_{0}^{\infty} e^{-x\tau/\sigma^2} e^{-x^{2}/2\sigma^2} dx\\
&\leq e^{-\tau^{2}/2\sigma^2} \int_{0}^{\infty} e^{-x^{2}/2\sigma^2} dx 
= \sqrt{\frac{\pi\sigma^2}{2}} e^{-\tau^{2}/2\sigma^2},
\end{align*}
where the inequality follows from that $e^{-x\tau/\sigma^2} \leq 1$ for $x,\tau \geq 0$.
This shows that $\prob(\abs{\rx} \geq \tau) \leq  e^{-\tau^{2}/2\sigma^2}$.

For the second part, by symmetry and union bound, $\prob(\abs{\rx} \geq \tau) 
\leq 2\prob(\rx \geq \tau)$:
\begin{align*}
2\prob(\rx \geq \tau)
&= \frac{\sqrt{2}}{\sqrt{\pi}\sigma} \int_\tau^{+\infty} \exp\left( -\frac{x^2}{2\sigma^2} \right) dx 
\leq \frac{\sigma^2 \sqrt{2}}{\sqrt{\pi}\sigma^2} \int_\tau^{+\infty} \frac{x}{\sigma^2 \tau} \exp\left( -\frac{x^2}{2\sigma^2} \right) dx \\
&= \frac{\sigma \sqrt{2}}{\tau \sqrt{\pi}} \int_\tau^{+\infty} -\frac{\partial}{\partial x} \exp\left( -\frac{x^2}{2\sigma^2} \right) dx
= \frac{\sigma \sqrt{2}}{\tau \sqrt{\pi}} \exp\left( -\frac{\tau^2}{2\sigma^2} \right).
\end{align*}
This completes the proof.
\end{proof}

The definition of a standard Gaussian random variable $\rx\sim\normal(0,1)$ shows that $\Mom_{\rx}(\theta) = \Exp[\exp(\theta \rx)] = \exp(\theta^2/2)$. A more general result holds for expectations involving quadratic functions of $\rx$.
\begin{lemma}[Moment of standard Gaussian variables]\label{lemma:mom_stand_norm}
Let $\rx\sim \normal(0,1)$. Then, for $a, b,c \in \real$ with $a < 1/2$,
$$
\Exp[\exp(a \rx^2 + b \rx + c)] 
= \frac{1}{\sqrt{1 - 2a}} \exp\left(\frac{b^2}{2(1 - 2a)} + c\right).
$$
\end{lemma}
\begin{proof}[of Lemma~\ref{lemma:mom_stand_norm}]
By definition of expectation under the standard normal density,  
$$
\Exp[\exp(a \rx^2 + b \rx +c)] 
= \frac{1}{\sqrt{2\pi}} \int_{-\infty}^{\infty} \exp(ax^2 + b x + c) \exp\left(-\frac{x^2}{2}\right) dx.
$$
Noting the identity
$$
ax^2 - \frac{x^2}{2} + b x = -\frac{1 - 2a}{2} \left(x - \frac{b}{1 - 2a}\right)^2 + \frac{b^2}{2(1 - 2a)} + c.
$$
After a change of variable, the latter integral reduces to the integral of the Gaussian probability density function with variance $1/(1-2a)$, so it equals one. This completes the proof.
\end{proof}

Note that the above lemma implies a bound on the higher-order moments of a standard Gaussian $\rx$, since
\begin{align*}
	\Exp[\exp(\theta\rx^2)]
	&= \Exp\left[\sum_{i=0}^\infty \frac{(\theta\rx^2)^i}{i!}\right]
	= \sum_{i=0}^\infty \frac{\Exp[(\theta\rx^2)^i]}{i!},
\end{align*}
where the first equality follows from the well-known Taylor expansion $\exp(z) = \sum_{i=0}^{\infty} \frac{z^i}{i!}$ for all $z\in\real$ around 0. 
Bounds that exploit the behavior of higher-order moments to control tail probabilities through the expectation of an exponential are often called Chernoff bounds or Chernoff methods; see also Theorem~\ref{theorem:chernoff_bound}.

\subsection*{Gaussian Random Vectors}\label{section:multi_gaussian_dist}

\begin{figure}[h]
\subfigure[Gaussian, $\bSigma =\begin{bmatrix}
1&0\\
0&1
\end{bmatrix}. $ ]{\includegraphics[width=0.31
\textwidth]{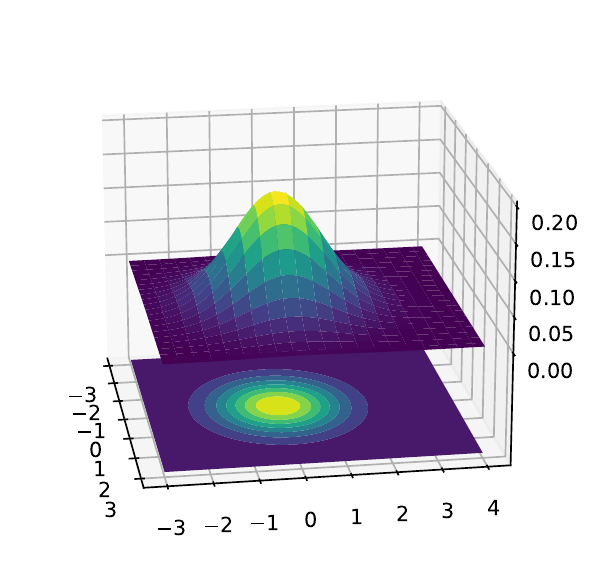} \label{fig:dists_multiGauss_sigma1}}
\subfigure[Gaussian, $\bSigma =\begin{bmatrix}
1&0\\
0&3
\end{bmatrix}.$]{\includegraphics[width=0.31
\textwidth]{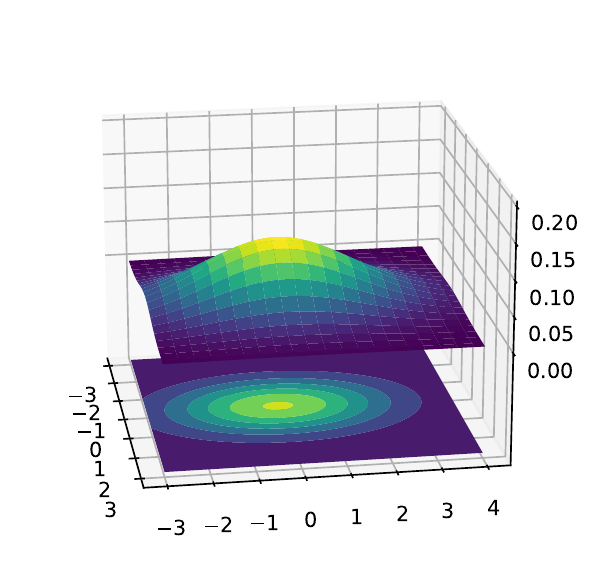} \label{fig:dists_multiGauss_sigma2}}
\subfigure[Gaussian, $\bSigma =\begin{bmatrix}
1&\textendash0.5\\
\textendash0.5&1.5
\end{bmatrix}.$]{\includegraphics[width=0.31 
\textwidth]{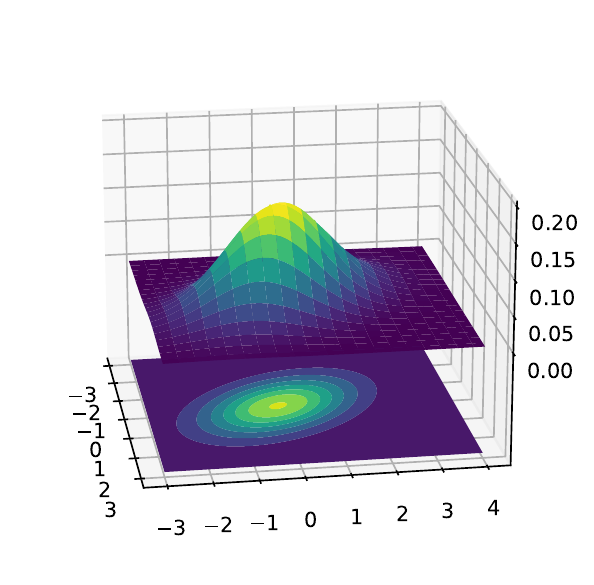} \label{fig:dists_multiGauss_sigma3}}
\subfigure[Gaussian, $\bSigma =\begin{bmatrix}
2&0\\
0&2
\end{bmatrix}. $ ]{\includegraphics[width=0.31
\textwidth]{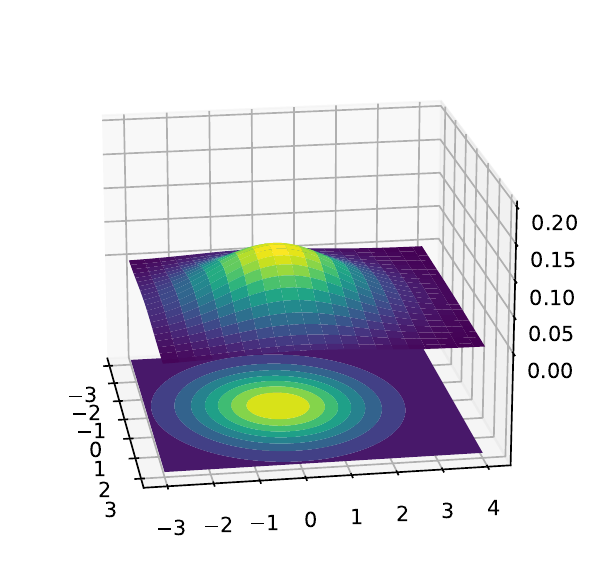} \label{fig:dists_multiGauss_sigma4}}
\subfigure[Gaussian, $\bSigma =\begin{bmatrix}
3&0\\
0&1
\end{bmatrix}.$]{\includegraphics[width=0.31
\textwidth]{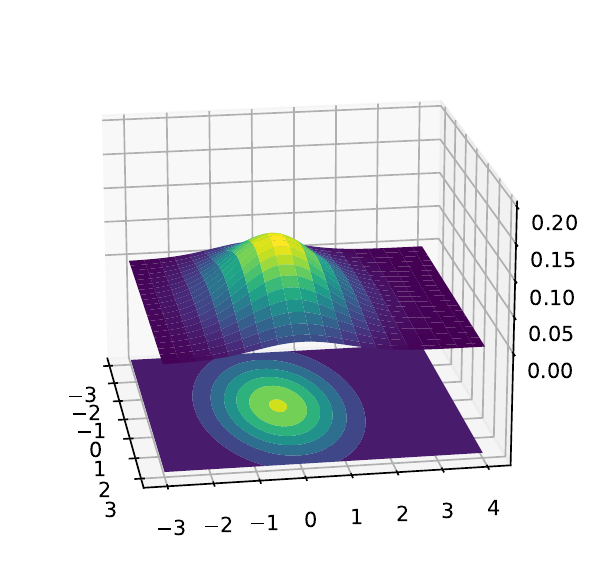} \label{fig:dists_multiGauss_sigma5}}
\subfigure[Gaussian, $\bSigma =\begin{bmatrix}
3&\textendash0.5\\
\textendash0.5&1.5
\end{bmatrix}.$]{\includegraphics[width=0.31
\textwidth]{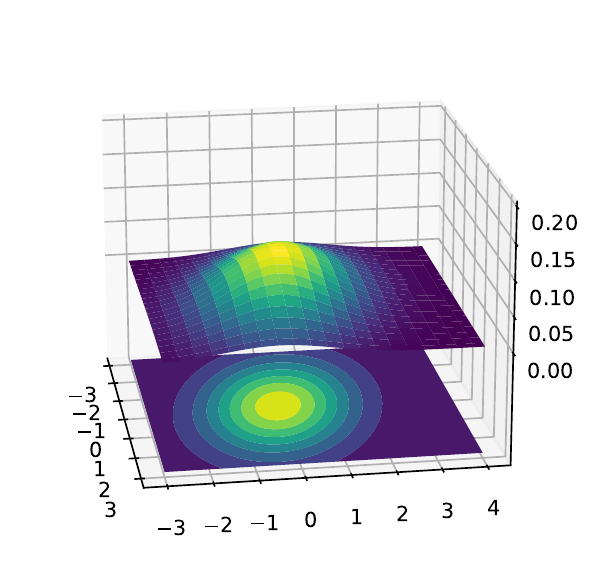} \label{fig:dists_multiGauss_sigma6}}
\centering
\caption{Density and contour plots (\textcolor{mydarkblue}{blue}=low, \textcolor{mydarkyellow}{yellow}=high) of the multivariate Gaussian distribution over the $\real^2$ space for various values of the covariance/scale matrix with a zero-mean vector.  Fig~\ref{fig:dists_multiGauss_sigma1} and \ref{fig:dists_multiGauss_sigma4}: A spherical covariance matrix has a circular shape; 
Fig~\ref{fig:dists_multiGauss_sigma2} and \ref{fig:dists_multiGauss_sigma5}:  A diagonal covariance matrix yields axis-aligned elliptical contours; 
Fig~\ref{fig:dists_multiGauss_sigma3} and \ref{fig:dists_multiGauss_sigma6}: A full covariance matrix results in rotated elliptical contours.}
\centering
\label{fig:multi_gaussian_density}
\end{figure}

A \textit{multivariate Gaussian distribution} (also known as a \textit{multivariate normal distribution} or simply a Gaussian distribution) is a continuous probability distribution that describes a set of jointly normally distributed random variables.
It is completely characterized by  its mean vector and covariance matrix. 
The covariance matrix captures the pairwise covariances between all variables, thereby encoding both their individual variances (along the diagonal) and their linear dependencies (off-diagonal entries).
Due to its mathematical tractability and flexibility, the multivariate Gaussian distribution is widely used across fields such as machine learning, statistics, and signal processing. 
We now provide the formal definition of the multivariate Gaussian distribution.
\begin{definition}[Multivariate Gaussian distribution]\label{definition:multivariate_gaussian}
A random vector $\rvx \in \real^p$ is said to follow the \textit{multivariate Gaussian distribution (multivariate normal, MVN)} with parameters $\bmu\in\real^p$ and $\bSigma\in\real^{p\times p}$, denoted  $\rvx\sim \normal(\bmu, \bSigma)$, if
$$
\begin{aligned}
f(\bx; \bmu, \bSigma)&= (2\pi)^{-p/2} \abs{\bSigma}^{-1/2}\exp\left\{-\frac{1}{2}(\bx - \bmu)^\top \bSigma^{-1}(\bx - \bmu)\right\},~\footnote{The form of which can be proved using the moment generating function of $p$ i.i.d. univariate standard Gaussian variables.}
\end{aligned}
$$
where $\bmu \in \real^p$ is called the \textit{mean vector}, and $\bSigma\in \real^{p\times p}$ is positive definite and is called the \textit{covariance matrix}. $\abs{\bSigma} = \det(\bSigma)$ denotes the determinant of  $\bSigma$.
The mean, mode, and covariance of the multivariate Gaussian distribution are given by 
\begin{equation*}
\begin{aligned}
\Exp [\rvx] &= \bmu, \qquad 
\mathrm{Mode}[\rvx] = \bmu, \qquad\text{and}\qquad 
\Cov [\rvx] = \bSigma. 
\end{aligned}
\end{equation*}
The covariance matrix can also be expressed as:
$$
\Cov[\rvx]=\Exp[(\rvx-\bmu)(\rvx-\bmu)^\top]=\Exp[\rvx\rvx^\top]-\bmu\bmu^\top.
$$
\end{definition}
Figure~\ref{fig:multi_gaussian_density} illustrates Gaussian density plots under different covariance structures.
When the covariance matrix of a Gaussian vector is the identity and its mean is zero, then its entries are i.i.d. standard Gaussians with mean zero and unit variance. We refer to such vectors as \textit{i.i.d. standard Gaussian vectors}.
Additionally, a multivariate Gaussian random vector can be generated from univariate Gaussian samples; see Problem~\ref{problem:multiGauss}.

Two important distributions related to the multivariate Gaussian distribution are the \textit{Gamma distribution} and its special case, the \textit{Chi-squared distribution}. Definitions for these distributions are provided below.
\begin{definition}[Gamma distribution]\label{definition:gamma_distri}\index{Gamma distribution}
A random variable $\rx$ is said to follow the \textit{Gamma distribution} with shape parameter $r>0$ and rate parameter $\lambda>0$, denoted by $\rx \sim \gammadist(r, \lambda)$, if its probability density function is given by 
$$ 
f(x; r, \lambda)=\left\{
\begin{aligned}
&\frac{\lambda^r}{\Gamma(r)} x^{r-1} \exp(-\lambda x) ,& \mathrm{\,\,if\,\,} x \geq 0.  \\
&0 , &\mathrm{\,\,if\,\,} x <0,
\end{aligned}
\right.
$$
where $\Gamma(x)=\int_{0}^{\infty} t^{x-1}\exp(-t)dt$ is the {Gamma function} (see \eqref{equation:gamma_func}),  which normalizes the distribution so that it integrates to 1. 
The mean and variance of $\rx \sim \gammadist(r, \lambda)$ are given by 
\begin{equation}
\Exp[\rx] = \frac{r}{\lambda}, \qquad \Var[\rx] = \frac{r}{\lambda^2}. \nonumber
\end{equation}
\end{definition}

\begin{figure}[h]
\centering
\vspace{-0.35cm}
\subfigtopskip=2pt
\subfigbottomskip=2pt
\subfigcapskip=-5pt
\subfigure[Gamma probability density
functions for different values of the parameters $r$ and $\lambda$.]{\label{fig:dists_gamma}
\includegraphics[width=0.47\linewidth]{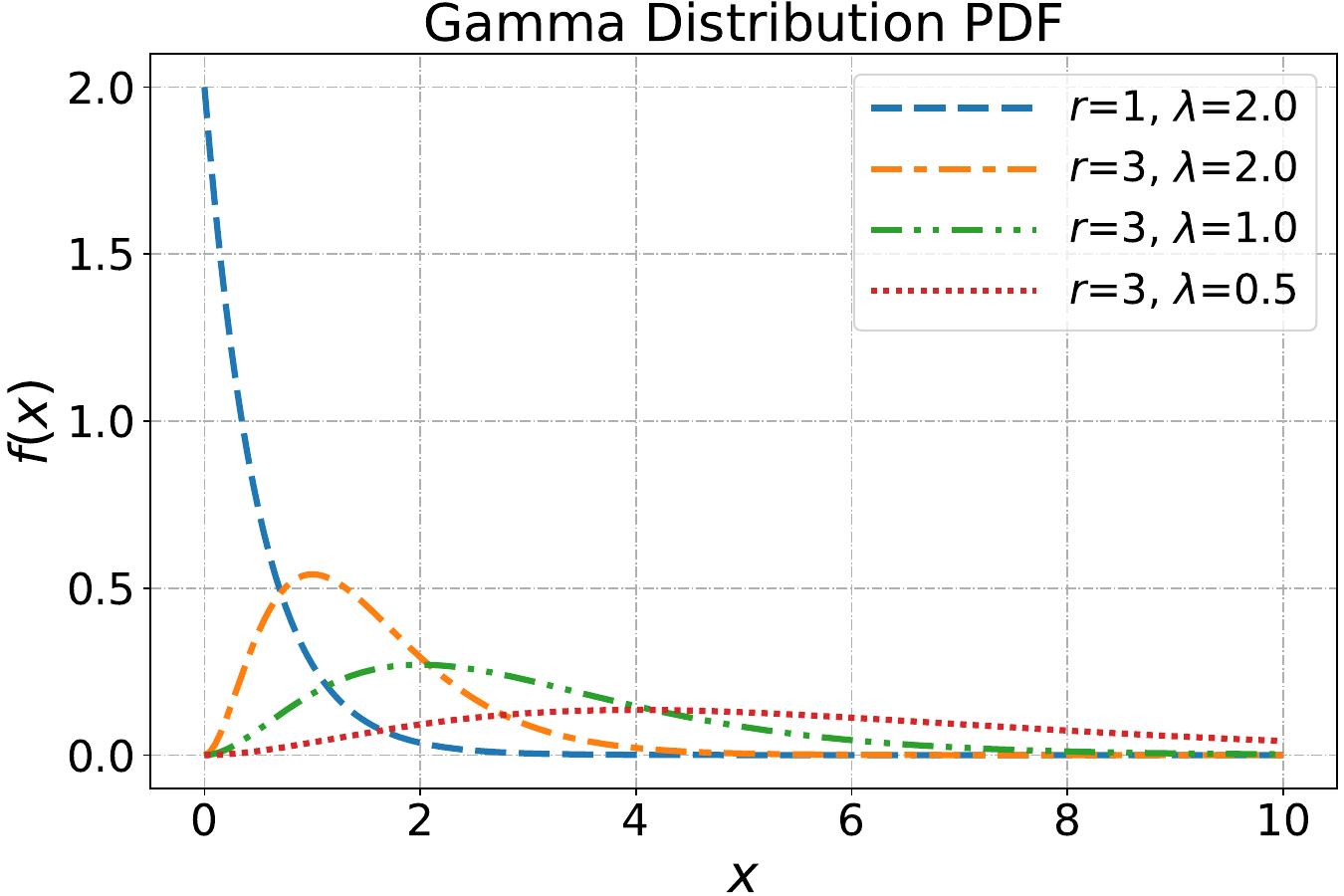}}
\quad 
\subfigure[Chi-squared probability density
functions for different values of the parameter $p$.]{\label{fig:dists_chisquare2}
\includegraphics[width=0.47\linewidth]{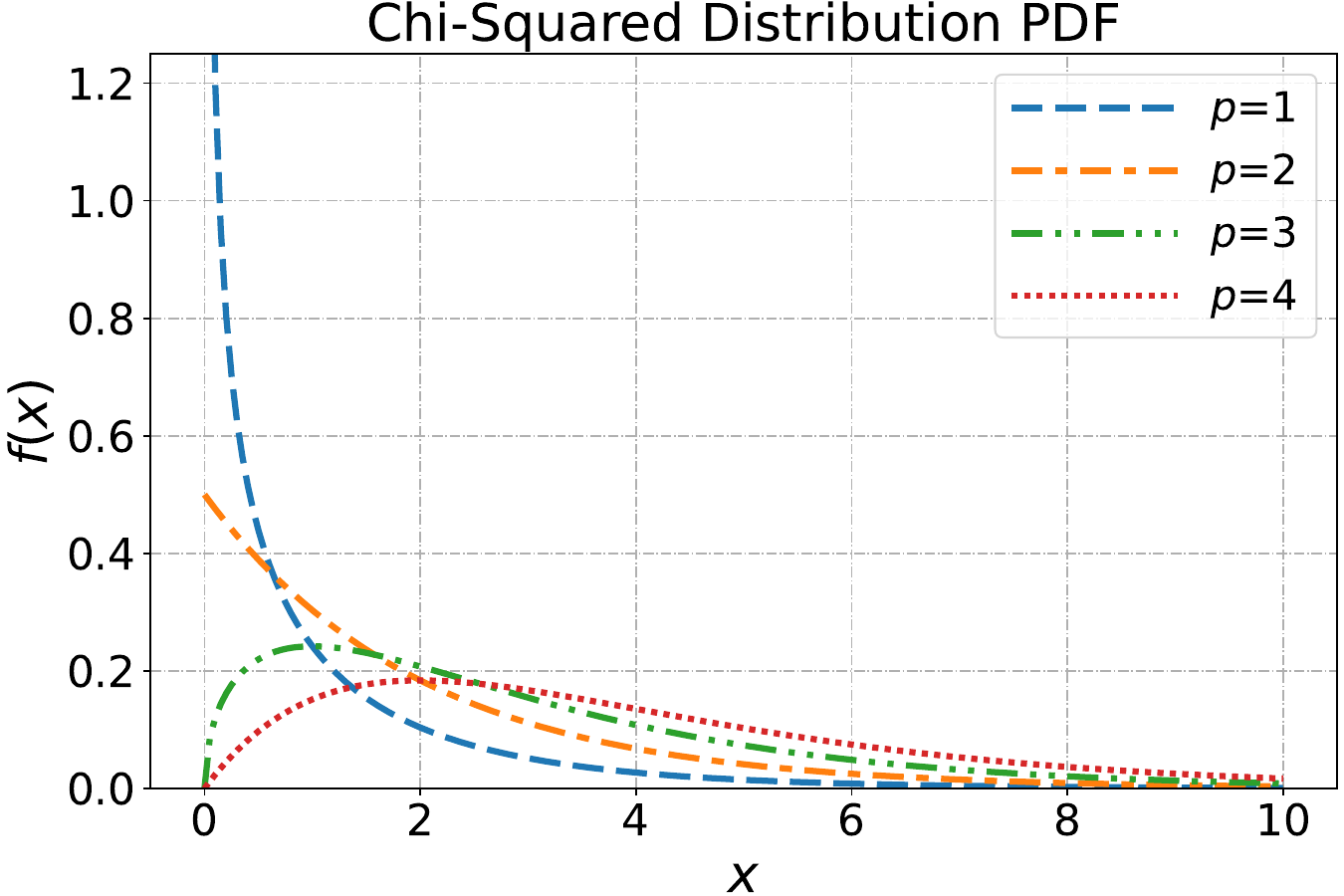}}
\caption{Comparison between the Gamma distribution and the Chi-squared distribution.}
\label{fig:gamma_chisquare_compare}
\end{figure}

Figure~\ref{fig:dists_gamma} illustrates different parameters for the Gamma distribution.
It is worth noting that the shape parameter $r$ in the Gamma distribution is not restricted to natural numbers; instead, it can take any positive value.
However, when $r$ is a positive integer, the Gamma distribution can be interpreted as the sum of $r$ independent \textit{exponential variables}, each with rate  $\lambda$ \citep{lu2023bayesian}.
This summation property extends more generally to Gamma-distributed variables with the same rate parameter. If $\rx_1$ and $\rx_2$ are random variables drawn from $\gammadist(r_1, \lambda)$ and $\gammadist(r_2, \lambda)$, respectively, then their sum $\rx_1+\rx_2$ follows a Gamma distribution $\gammadist(r_1+r_2, \lambda)$.

Although the Chi-squared distribution is a special case of the Gamma distribution, it plays a particularly important role in statistical theory.
\begin{definition}[Chi-squared distribution, $\chi^2$-distribution]\label{definition:chisquare_dist}\index{Chi-squared distribution}
Let $\ra_i\sim \normal(0, 1)$ for $i\in\{1,2,\ldots,p\}$ (equivalently, $\rva \sim \normal(\bzero, \bI_{p})$; see Definition~\ref{definition:multivariate_gaussian}). 
Then the random variable $\rx=\sum_{i=1}^p \ra_{i}^2$ follows the \textit{Chi-squared distribution} (or \textit{Chi-square distribution}, \textit{$\chi^2$-distribution}) with \textit{$p$ degrees of freedom}. 
We write $\rx \sim \chi_{(p)}^2$, which is equivalent to $\rx\sim \gammadist(p/2, 1/2)$ in Definition~\ref{definition:gamma_distri}. The probability density function of $\rx$ is given by 
$$ f(x; p)=\left\{
\begin{aligned}
&\frac{1}{2^{p/2}\Gamma(\frac{p}{2})} x^{\frac{p}{2}-1} \exp(-\frac{x}{2}) ,& \mathrm{\,\,if\,\,} x \geq 0;  \\
&0 , &\mathrm{\,\,if\,\,} x <0.
\end{aligned}
\right.
$$
The mean and variance of $\rx\sim \chi_{(p)}^2$ are given by 
$$
\Exp[\rx]=p, \qquad\Var[\rx]=2p.
$$
\end{definition}

Figure~\ref{fig:dists_chisquare2} shows the effect of varying the degrees of freedom  $p$ on the shape of the Chi-squared distribution.
The definition shows that if $\rva=[\ra_1, \ra_2, \ldots, \ra_p]^\top \sim \normal(\bzero, \bI_{p})$ (the multivariate Gaussian distribution; see Definition~\ref{definition:multivariate_gaussian}), then $\rx=\rva^\top\rva\sim \chi_{(p)}^2$. 
More generally, let   $\bH$ be an orthogonal projection matrix of rank $r<p$ (see Definition~\ref{definition:orthogonal-projection-matrix}). 
Then it follows that 
\begin{subequations}
\begin{equation}
\rva^\top \bH \rva \sim \chi_{(r)}^2, \text{ with orthogonal projector $\rank(\bH)=r<p$}.
\end{equation}
Furthermore, suppose  $ \rvx \sim \normal(\bmu, \bSigma) $,  where $\bmu\in\real^p$ and $\bSigma^{-1}$ is nonsingular. Then,
\begin{equation}
(\rvx - \bmu)^\top \bSigma^{-1} (\rvx - \bmu) \sim \chi^2_{(p)}.
\end{equation}
\end{subequations}
To see this, note that since $\bSigma$ is symmetric and positive definite, it admits the spectral decomposition $\bSigma = \bQ^\top \bLambda \bQ$, where $\bQ$ is a $p \times p$ orthogonal matrix and $\bLambda = \diag(\lambda_1,\lambda_2, \ldots, \lambda_p)$ with $\lambda_j > 0$ (see Theorem~\ref{theorem:spectral_theorem}).
Define the standardized vector $\rvy \triangleq \bLambda^{-1/2} \bQ (\rvx - \bmu)$. 
Vector $\rvy$ then is normally distributed (as a hindsight, see Lemma~\ref{lemma:affine_mult_gauss}):
$$
\rvy \sim \normal(\bzero, \bLambda^{-1/2} \bQ \bSigma \bQ^\top \bLambda^{-1/2}) = \normal(\bzero, \bLambda^{-1/2} \bQ \bQ^\top \bLambda \bLambda^{-1/2}) = \normal(\bzero, \bI_p).
$$
It follows that the components of $\rvy$ are independent (see Lemma~\ref{lemma:uncor_inde_mvn}) and that $\ry_i \sim \normal(0, 1)$. So,
$(\rvx - \bmu)^\top \bSigma^{-1} (\rvx - \bmu) = \rvy^\top\rvy \sim \chi^2_{(p)}$.

We may also be interested in the quadratic form of $\rvx^\top\bA\rvx$ where $\bA$ is symmetric. 
The following theorem provides key characterizations of when such quadratic forms follow a Chi-squared distribution.
\begin{proposition}[Quadratic of Gaussians]
We have the following results with quadratic forms of Gaussians:
\begin{itemize}
\item 
Given $\rvx\sim\normal(\bzero,\lambda\bI)$ (of length $p$) and symmetric matrix $\bA\in\real^{p\times p}$. Then, it follows that 
$$
\frac{\rvx^\top\bA\rvx}{\lambda} \sim \chisquared_{(r)},
$$
if and only if $\bA$  is idempotent ($\bA^2=\bA$) with rank $r<p$.
\item Given $\rvx\sim\normal(\bzero,\bSigma)\in\real^p$ and symmetric matrix $\bA\in\real^{p\times p}$. Then, it follows that 
$$
\rvx^\top\bA\rvx \sim \chisquared_{(r)},
$$
if and only if $\bA\bSigma$ is idempotent with rank $r<p$.
\end{itemize}
\end{proposition}

\subsection*{Properties of Multivariate Gaussian Distribution}\label{section:multi_gauss}

A fundamental property of Gaussian random vectors is that affine transformations preserve Gaussianity. 
This is a multidimensional generalization of Lemma~\ref{lem:lintrans_uni_gaus}, and it can be derived using the linear transformation rule for probability distributions (Lemma~\ref{lemma:lin_tran_prob}) or via moment generating functions of the multivariate Gaussian distribution.
 
\begin{lemma}[Affine transformation of multivariate Gaussian distribution]\label{lemma:affine_mult_gauss}
Given fixed matrices and vector, $\bA,\bB\in\real^{p\times d}$ and $\bc\in\real^p$, 
let $\rvx\sim \normal(\bmu_x, \bSigma_x)$ and $\rvy\sim \normal(\bmu_y, \bSigma_y)$ be independent random vectors (of length $d$).
Then, 
$$
\rvz=\bA\rvx+\bB\rvy +\bc \sim \normal(\bA\bmu_x+\bB\bmu_y+\bc, \bA\bSigma_x\bA^\top +\bB\bSigma_y\bB^\top).
$$
Moreover, for any fixed vector $\bd\in\real^{d}$, the linear combination $\bd^\top\rvx$ follows from a univariate Gaussian:
$$
\bd^\top\rvx \sim \normal(\bd^\top\bmu_x, \bd^\top\bSigma_x\bd).
$$
\end{lemma}
The  result also relies on the \textit{sum of independent Gaussians}:
$$
\sum_{i=1}^{n} \rvx_i \sim \normal\big(\sum_{i=1}^{n}\bmu_i, \sum_{i=1}^{n}\bSigma_i\big)
\gap 
\text{if } \rvx_i\sim\normal(\bmu_i,\bSigma_i), \ \forall\, i\in\{1,2,\ldots,n\}.
$$

As a special case, let  $\bA=\be_i^\top$ be the $i$-th standard basis vector in $\real^d$. 
Then $\rx_i = \be_i^\top\rvx \sim \normal(\mu_{x,i}, \sigma_{x,ii}^2)$, where $\mu_{x,i}$ represents the $i$-th component of $\bmu_x$ and $\sigma_{x,ii}^2$ denotes the $i$-th diagonal of $\bSigma_x$.

An immediate consequence of Lemma~\ref{lemma:affine_mult_gauss} is that any subvector of a multivariate Gaussian vector is itself Gaussian.
Another consequence of Lemma~\ref{lemma:affine_mult_gauss} is that an i.i.d. standard Gaussian vector is \textit{isotropic}, meaning its distribution is invariant under rotations. 
Formally, for any orthogonal matrix $\bQ$, if $\rvx$ is an i.i.d. standard Gaussian vector, then by Lemma~\ref{lemma:affine_mult_gauss}, $\bQ\rvx$ has the same distribution, since its mean equals $\bQ\bzero = \bzero$ and its covariance matrix equals $\bQ\bI\bQ^\top = \bQ\bQ^\top = \bI$. 
This rotational invariance is stronger than merely having equal variance in all directions (which holds for any vector with uncorrelated, identically distributed components). Isotropy implies that the entire distribution---not just second-order moments---is unchanged under orthogonal transformations.

\begin{lemma}[Rotations on  multivariate Gaussian distribution]\label{lemma:rotat_multi_gauss}
Rotations on the Gaussian distribution do not affect the distribution.
That is, for any orthogonal matrix $\bQ$ with $\bQ\bQ^\top=\bQ^\top\bQ=\bI$, if $\rvv\sim \normal(\bzero, \sigma^2\bI)$, then $\bQ\rvv\sim \normal(\bzero, \sigma^2\bI)$.
\end{lemma}

\paragrapharrow{``Standardization and decorrelation."}
The distribution $\normal(\bzero, \bI)$ is called the \textit{i.i.d. standard  Gaussian vector}. 
Given $\rvx\sim\normal(\bmu, \bSigma)$, then the \textit{decorrelation} of $\rvx$ follows that 
\begin{equation}\label{equation:std_mugau_recov}
\rvx\sim\normal(\bmu, \bSigma)
\quad\implies \quad
\rvz =\bSigma^{-1/2}(\rvx-\bmu) \sim \normal(\bzero,\bI).
\end{equation}
This also shows that if $\rvx\sim\normal(\bmu,\bSigma)$, then 
\begin{equation}
\rvx=\bmu+\bSigma^{1/2}\bepsilon,
\gap \text{where }\bepsilon\sim\normal(\bzero,\bI).
\end{equation}

Suppose $\{\bx_i,\bx_2,\ldots,\bx_n\}$ are $n$ random samples of $\normal(\bmu,\bSigma)$ and let $\overline{\bx} = \frac{1}{n} \sum_{i=1}^{n} \bx_i $. Then, it follows that 
\begin{equation}\label{equation:mean_mutigau}
\sqrt{n} (\overline{\bx}-\bmu) \sim \normal(\bzero, \bSigma).
\end{equation}

\begin{lemma}[Uncorrelation implies mutual independence for Gaussian random vectors]\label{lemma:uncor_inde_mvn}
If all  components of a Gaussian random vector $\rvx$ are uncorrelated, then they are  mutually independent.
\end{lemma}
This result follows directly by observing that the joint p.d.f. of $\rvx$ factorizes into the product of its marginal p.d.f.s when the covariance matrix is diagonal (i.e., when components are uncorrelated).

We now turn to the partitioning of Gaussian vectors.
Let $\rvx \sim \normal(\bmu, \bSigma)$ where $\bmu\in\real^p$,
and consider a partition of $\rvx$ into two subvectors:  
$$
\begin{bmatrix}
\rvx_1 \\
\rvx_2
\end{bmatrix}
\sim \normal(\bmu, \bSigma) = \normal
\left(
\begin{bmatrix}
\bmu_1 \\
\bmu_2
\end{bmatrix},
\begin{bmatrix}
\bSigma_{11} & \bSigma_{12} \\
\bSigma_{21} & \bSigma_{22}
\end{bmatrix}
\right).
$$
Then $\rvx_1$ and $\rvx_2$ are independent if and only if $\bSigma_{12} = \bzero$. 
More generally, let  $\rvx=[\rx_1, \rx_1, \ldots, \rx_p ]^\top\sim \normal(\bmu, \bSigma)$. Then, 
\begin{equation}\label{equation:iid_mulg_iff7}
\text{the $\rx_i$'s are mutually independent if and only if $\bSigma$ is diagonal.}
\end{equation}
This can be rigorously established as follows:
\begin{proof}[of Lemma~\ref{lemma:uncor_inde_mvn} and Equation~\eqref{equation:iid_mulg_iff7}]
Suppose that the $\rx_i$'s are independent. Then $\rx_i \sim \normal(\mu_i, \sigma_i^2)$ for some $\sigma_i > 0$. Thus the density of $\rvx$ is
\begin{align*}
p_{\rvx}(\bx)\ &= \prod_{i=1}^{p} p_{\rx_i}(x_i) = \prod_{i=1}^{p} \frac{1}{\sigma_i \sqrt{2\pi}} \exp \left\{ -\frac{1}{2} \frac{(x_i - \mu_i)^2}{\sigma_i^2} \right\} \\
&= \frac{1}{(2\pi)^{p/2} \abs{\diag(\sigma_1^2, \ldots, \sigma_p^2)}^{1/2}} \exp \left\{ -\frac{1}{2} (\bx - \bmu)^\top \diag(\sigma_1^{-2}, \ldots, \sigma_p^{-2}) (\bx - \bmu) \right\}.
\end{align*}
Hence $\rvx \sim \normal(\bmu, \diag(\sigma_1^2, \ldots, \sigma_p^2))$, i.e., the covariance $\bSigma$ is diagonal.

Conversely, assume $\bSigma$ is diagonal, say $\bSigma = \diag(\sigma_1^2, \ldots, \sigma_p^2)$. Then we can reverse the steps of the first part to see that the joint density $p_{\rvx}(\bx)$ can be written as a product of the marginal densities $p_{\rx_i}(x_i)$, thus proving independence.
\end{proof}

\index{Independence}
The results above immediately imply the following corollary.
\begin{corollary}[Independence of linear combinations in Gaussian Distributions]
Let $ \rvx \sim \normal(\bmu, \bSigma)$ be a random vector in $\real^{p} $ , and let $ \bA\in\real^{m \times p} $, $ \bB\in\real^{d \times p} $ be fixed matrices. 
Then,
\begin{equation}
\text{$\bA\rvx$ \text{ is independent of } $\bB\rvx$ $\quad\iff\quad$ $\bA\bSigma \bB^\top = \bzero$.}
\end{equation}
\end{corollary}
An alternative proof can be obtained using properties of the moment generating function of the multivariate Gaussian distribution; we omit the details here.

It also follows from basic probability theory that measurable functions of independent random vectors remain independent. That is, if $\rvx_1$ and $\rvx_2$ are independent, then $g_1(\rvx_1)$ and $g_2(\rvx_2)$ are independent for any measurable functions $g_1(\rvx_1)$ and $g_2(\rvx_2)$.

As an important example, suppose  variables $\rx_i$ are i.i.d. $\normal(\mu, \sigma^2)$ for $i\in\{1,2,\ldots,p\}$. Then, we can define a vector populated by $\overline{\rx}$ and $\rx_i - \overline{\rx}$:
$$
\begin{bmatrix}
\overline{\rx} \\
\rx_1 - \overline{\rx} \\
\vdots \\
\rx_p - \overline{\rx}
\end{bmatrix}
=
\begin{bmatrix}
\frac{1}{p} & \frac{1}{p} & \ldots & \frac{1}{p} \\
& \bI_p - \frac{1}{p}\bJ_p &
\end{bmatrix}
\begin{bmatrix}
\rx_1 \\
\rx_2 \\
\vdots \\
\rx_p
\end{bmatrix}
\qquad \text{where} \qquad
\bJ_p =
\begin{bmatrix}
1 & 1 & \ldots & 1 \\
1 & \ddots & & 1 \\
1 & & \ddots & \vdots \\
1 & 1 & \ldots & 1
\end{bmatrix}\in\real^{p\times p}.
$$
It then follows that 
$$
\begin{bmatrix}
\overline{\rx} \\
\rx_1 - \overline{\rx} \\
\vdots \\
\rx_p - \overline{\rx}
\end{bmatrix}
\sim
\mathcal{N}
\left(
\begin{bmatrix}
\mu \\
0 \\
\vdots \\
0
\end{bmatrix},
\sigma^2
\begin{bmatrix}
\frac{1}{p} & \bzero \\
\bzero & \bI_p - \frac{1}{p}\bJ_p
\end{bmatrix}
\right).
$$
Therefore, we find that $\overline{\rx}$ is independent of $\rx_1 - \overline{\rx},\rx_2 - \overline{\rx}, \ldots, \rx_p - \overline{\rx}$. 
In many applications, we may construct a random variable:
$$
t \triangleq  \sqrt{p} \frac{(\overline{\rx} - \mu) / \sigma}{\sqrt{\frac{1}{p-1} \sum_i (x_i - \overline{\rx})^2 / \sigma^2}},
$$
in which case, the numerator and the denominator are independent variables.


\begin{problemset}

\item \label{prob:tr_de_pd} \textbf{Trace, det of PD/PSD/ND matrices.} Let $\bX$ be positive definite (resp. positive semidefinite). Show that $\trace(\bX)$ and $\det(\bX)$  are all positive (resp. nonnegative). Moreover, show that $\trace(\bX)=0$ if and only if $\bX=\bzero$. 
Let $\bY\in\real^{p\times p}$ be negative definite. Show that $\trace(\bY)$ is negative; $\det(\bY)$ is negative for odd $p$ and is positive for even $p$.
\textit{Hint: Use Theorem~\ref{theorem:eigen_charac}.}

\item Prove Remark~\ref{remark:equiva_nonsingular} and Remark~\ref{remark:equiva_singular}.

\item Demonstrate that the vector $\ell_2$-norm, the matrix Frobenius norm, and the matrix spectral norm satisfy the three conditions outlined in Definition~\ref{definition:matrix-norm}.

\item \textbf{Induced norm properties \citep{lu2021numerical}.} 
\label{prob:induced_norm_property}
Consider the induced matrix norm in Definition~\ref{definition:induced_norm_app}.
\begin{enumerate}
\item
Given $\balpha\in \real^n$, $\bbeta\in \real^p$, and $a,b\geq 1$ with $\frac{1}{a}+\frac{1}{a^\star}=1$ and $\frac{1}{b}+\frac{1}{b^\star}=1$, show that 
$$
\norm{\balpha\bbeta^\top}_{a,b} = \normb{\balpha} \norm{\bbeta}_{a^\star}.
$$
\item 
Given $\balpha\in \real^n$ and $\bbeta\in \real^p$ with $\norma{\bbeta} =\normb{\balpha}=1$, show that there exists a matrix $\bX\in\real^{n\times p}$ such that $\norm{\bX}_{a,b}=1$ and $\bX\bbeta=\balpha$.
\item 
Particularly, show that $\norm{\bX^\top}_{a,b} = \norm{\bX}_{b^\star, a^\star}$.
\end{enumerate}

%
%

\item \textbf{Equivalence of matrix norms.} \label{prob:equiv_mat_norm}
Similarly to the vector norms in Exercise~\ref{exercise:cauch_sc_l1l2},  matrix norms also admit an equivalence statement.
Let $\norma{\cdot}$ and $\normb{\cdot}$ be two different matrix norms: $\real^{n\times p}\rightarrow \real$. Show that  there exist positive scalars $\alpha$ and $\beta$ such that for all $\bX \in \real^{n\times p}$,
$$
\alpha\norma{\bX} \leq \normb{\bX} \leq \beta\norma{\bX}.
$$
The equivalence theorem again states that if a matrix is small in one norm, it is also small in other norms, and vice versa.

\item \textbf{Construct induced norms from induced norms.}\label{prob:cons_ind}
Let $\norm{\cdot}$ be a matrix norm (resp., a submultiplicative matrix norm) on $\real^{p\times p}$, and let $\bS\in\real^{p\times p}$ be nonsingular. Show that  the following function defined for $\bX\in\real^{p\times p}$ is also a matrix norm (resp., a submultiplicative matrix norm):
$$
\norm{\bX}_{\bS} = \norm{\bS\bX\bS^{-1}}.
$$
Moreover, if $\norm{\cdot }$ is induced from the vector norm $\norm{\cdot}_v$ on $\real^p$, show that the matrix norm $\norm{\bX}_{\bS}$ is induced from the vector norm $\norm{\bbeta}_{\bS}=\norm{\bS\bbeta}_v$ on $\real^p$.


\item \label{problem:multiGauss} Suppose we can generate the univariate Gaussian variable $\normal(0, 1)$. Provide a way to generate the multivariate Gaussian variable $\normal(\bmu, \bSigma)$, where $\bSigma=\bC\bC^\top$, $\bmu\in\real^n$, and $\bSigma\in\real^{n\times n}$. \textit{Hint: if $\rx_1, \rx_2, \ldots, \rx_n$ are i.i.d. from $\normal(0, 1)$, and let $\rvx=[\rx_1, \rx_2, \ldots, \rx_n]^\top$, then it follows that $\bC\rvx +\bmu\sim \normal(\bmu, \bSigma)$}.

\item Prove Lemma~\ref{lemma:affine_mult_gauss} and Lemma~\ref{lemma:rotat_multi_gauss} rigorously.

\item \label{problem:kl_vae}
\textbf{KL of Gaussians.} Given two probability distribution function $P(\bx)$ and $Q(\bx)$, we denote  the \textit{Kullback-Leibler (KL) divergence} between $P$ and $Q$ by $\KL[P \parallel Q] = \int P(\bx) \ln \left( \frac{P(\bx)}{Q(\bx)} \right) d\bx \geq 0$, where the equality is obtained only when $P=Q$.
Given $q(\bx)=\normal( \bmu, \diag(\bsigma^2))$ and $p(\bx)=\normal( \bzero_D, \bI_D)$ where $\bmu, \bsigma^2\in\real^D$. Show that $\KL[q\parallel p]=\frac{1}{2}\sum_{i=1}^{D}(\mu_i^2+\sigma^2_i-\ln \sigma^2_i-1)$. 

\item \textbf{KL of Gaussians.} Suppose $p(x)=\normal(\mu_1, \sigma_1^2)$ and $q(x)=\normal(\mu_2,\sigma_2^2)$. Show that 
$$
\KL[p \parallel q ] = \ln\frac{\sigma_2}{\sigma_1} +\frac{\sigma_1^2+(\mu_1-\mu_2)^2}{2\sigma_2^2}-\frac{1}{2}.
$$ 
Consider the multivariate case. 
Suppose $\normal_1(\bx)=\normal(\bmu_1, \bSigma_1)$ and $\normal_2(\bx)=\normal(\bmu_2,\bSigma_2)$. Show that 
$$
\KL[\normal_1\parallel \normal_2] = 
\frac{1}{2}\ln\abs{\bSigma_2\bSigma_1^{-1}}+\frac{1}{2}\trace\bSigma_2^{-1}
\big( (\bmu_1-\bmu_2)(\bmu_1-\bmu_2)^\top +\bSigma_1-\bSigma_2 \big).
$$
More generally, consider a general distribution $p(\bx)$ and a multivariate Gaussian $\normal(\bx)=\normal(\bmu, \bSigma)$ with $\bx\in\real^D$. Show that 
$$
\KL[p\parallel \normal ]=
\int \frac{1}{2}(\bx-\bmu)^\top\bSigma^{-1}(\bx-\bmu) d\bx 
+
\frac{1}{2}\ln \abs{\bSigma}+\frac{D}{2}\ln2\pi + \int p(\bx)\ln p(\bx) d\bx.
$$


\item \label{problem:entropy_mgau} \textbf{Entropy of Gaussians.}  A close quantity related to KL divergence is the \textit{entropy}. The entropy $\entropy[p(\bx)]$ of a distribution $p(\bx)$ is defined as 
$$
\entropy[p(\bx)] \triangleq -\int p(\bx)\ln p(\bx) d\bx.
$$
Given a multivariate Gaussian distribution $\rvx\sim \normal(\bmu,\bSigma)\in\real^D$, show that 
$$
\entropy[\normal(\bmu,\bSigma)] = \frac{1}{2}\ln \abs{\bSigma}+ \frac{D}{2} \ln(2\pi e).
$$

\item \textbf{Joint Gaussian density \citep{lu2023bayesian}.} Let $\mathcalX=\{x_1, x_2, \ldots, x_N\}$ be drawn i.i.d. from a Gaussian distribution of $\normal(x\mid \mu, \sigma^2)$.
Show that the joint probability distribution of $\mathcalX$ can be written as
\begin{equation}\label{equation:uni_gaussian_likelihood}
\begin{aligned}
p(\mathcal{X} \mid \smu, \ssigma^2) 
&= (2\pi)^{-N/2}  (\ssigma^2)^{-N/2} \exp\left\{-\frac{1}{2 \ssigma^2}  \left[  N(\overline{x} - \smu)^2 +  N S_{\overline{x}} \right] \right\},
\end{aligned}
\end{equation}
where $S_{\overline{x}}\triangleq \sum_{n=1}^N(x_n - \overline{x})^2$ and $\overline{x} \triangleq \frac{1}{N} \sum_{i=1}^{N}x_i$.
Similar to the likelihood under univariate Gaussian distribution \eqref{equation:uni_gaussian_likelihood},
show that the likelihood of $N$ random observations $\mathcal{X} = \{\bx_1, \bx_2, \ldots , \bx_N \}$,  generated by a multivariate Gaussian with mean vector $\bmu$ and covariance matrix $\bSigma$, is given by 
$$
\begin{aligned}
&\gap p(\mathcal{X} \mid \bmu, \bSigma) =\prod^N_{n=1} \mathcal{N} (\bx_n\mid \bmu, \bSigma) 
= (2\pi)^{-ND/2} \abs{\bSigma}^{-N/2}\exp\left\{-\frac{1}{2} \tr( \bSigma^{-1}\bS_{\bmu} )  \right\}\\
&= (2\pi)^{-ND/2} \abs{\bSigma}^{-N/2}\exp\left\{-\frac{N}{2}(\bmu - \overline{\bx})^\top \bSigma^{-1}(\bmu - \overline{\bx})\right\}  \exp\left\{-\frac{1}{2}\tr( \bSigma^{-1}\bS_{\overline{x}} )\right\},
\end{aligned}
$$
where 
$$
\bS_{\bmu} \triangleq \sum^N_{n=1}(\bx_n - \bmu)(\bx_n - \bmu)^\top,\quad
\bS_{\overline{x}} \triangleq \sum^N_{n=1}(\bx_n - \overline{\bx})(\bx_n - \overline{\bx})^\top, \quad
\overline{\bx} \triangleq\frac{1}{N}\sum^N_{n=1}\bx_n.
$$
The matrix $\bS_{\overline{x}}$ is the \textit{matrix of sum of squares} and is also known as the \textit{scatter matrix}.

\item \label{prob:marg_mvn} \textbf{Marginal and conditional of MVN \citep{lu2023bayesian}.}
Let $\rvx$ and $\rvy$ be jointly Gaussian random vectors with 
$$
\rvz=
\begin{bmatrix}
\rvx\\
\rvy 
\end{bmatrix}
\sim 
\normal\left(
\begin{bmatrix}
\bmu_x\\
\bmu_y 
\end{bmatrix}
,
\begin{bmatrix}
\bA & \bC\\
\bC^\top & \bB
\end{bmatrix}
\right)=
\normal\left(
\begin{bmatrix}
\bmu_x\\
\bmu_y 
\end{bmatrix}
,
\begin{bmatrix}
\widetildebA & \widetildebC\\
\widetildebC^\top & \widetildebB
\end{bmatrix}^{-1}
\right).~\footnote{
Given nonsingular $\bM$ and its inverse $\bM^{-1}$; and suppose appropriate sizes for the following partitions \citep{williams2006gaussian}:
$$
\bM=
\begin{bmatrix}
\bA & \bB\\
\bC & \bD
\end{bmatrix},
\gap 
\bM^{-1}=
\begin{bmatrix}
\widetildebA & \widetildebB\\
\widetildebC & \widetildebD
\end{bmatrix}.
$$
We have 
\begin{equation}\label{equation:mt_inv}
\begin{aligned}
&\widetildebA=\bA^{-1}+\bA^{-1}\bB\widetildebD\bC\bA^{-1}&=& (\bA-\bB\bD^{-1}\bC)^{-1} , \\
&\widetildebB=-\bA^{-1}\bB\widetildebD&=& -\widetildebA\bB\bD^{-1},\\
&\widetildebC=-\widetildebD\bC\bA^{-1}&=& -\bD^{-1}\bC\widetildebA, \\
& \widetildebD=(\bD-\bC\bA^{-1}\bB)^{-1}&=& \bD^{-1}+\bD^{-1}\bC\widetildebA\bB\bD^{-1},
\end{aligned}
\end{equation}
}
$$
where $\rvx$ and $\rvy$ are \textit{independent} if and only if $\Cov[\rvx,\rvy]=\bC=\bzero$.
Show that every marginal distribution of a multivariate Gaussian distribution is itself a multivariate Gaussian distribution, and the conditional distribution $\rvx\mid \rvy$ also follows a multivariate Gaussian distribution:
\begin{equation}
\begin{aligned}
\rvx
\sim\normal(\bmu_x,\bA),
\gap 
\rvx\mid \rvy=\by 
&\sim \normal(\bmu_x+\bC\bB^{-1}(\by-\bmu_y), \bA-\bC\bB^{-1}\bC^\top)\\
&=\normal(\bmu_x-\widetildebA^{-1}\widetildebC(\by-\bmu_y), \widetildebA^{-1});\\
\rvy
\sim\normal(\bmu_y,\bB),
\gap 
\rvy\mid \rvx=\bx
&\sim \normal(\bmu_y+\bC^\top\bA^{-1}(\bx-\bmu_x), \bB-\bC^\top\bA^{-1}\bC)\\
&=\normal(\bmu_y-\widetildebB^{-1}\widetildebC^\top(\bx-\bmu_x), \widetildebB^{-1}).\\
\end{aligned}
\end{equation}

\item \textbf{Affine dependence of Gaussian variables.}
Suppose random vectors $\rvx\sim\normal(\bmu, \bSigma)$ and $\rvy\mid \rvx=\bx\sim\normal(\bA\bx+\bb, \bM)$.
Note $\rvy$ is not simply the affine transformation $\bA\rvx+\bb$, but it follows that $\rvy=\bA\bx+\bb+\bepsilon$ where $\bepsilon\sim \normal(\bzero, \bM)$.
Show that 
$$
\rvy\sim \normal(\bA\bmu+\bb, \bM+\bA\bSigma\bA^\top),
\gap 
\rvx \mid \rvy \sim \normal\big(\bL\big\{\bA^\top\bM^{-1}(\rvy-\bb)+\bSigma^{-1}\bmu  \big\}, \bL\big), 
$$
where $\bL=(\bSigma^{-1}+\bA^\top\bM^{-1}\bA)^{-1}$.
\textit{Hint: Use Problem~\ref{prob:marg_mvn}, compute the cross-covariance of $\rvx$ and $\rvy$ by $\Cov[\rvx,\rvy]=\Exp[(\rvx-\bmu_x)(\rvy-\bmu_y)^\top]=\bSigma\bA^\top$ where $\bmu_x=\bmu$ and $\bmu_y=\Exp[\rvy]=\bA\bmu+\bb$, and use Woodbury matrix identity: $(\bA+\bB\bD\bC)^{-1} = \bA^{-1} - \bA^{-1} \bB(\bD^{-1} + \bC\bA^{-1}\bB)^{-1}\bC\bA^{-1}$ for appropriate matrices $\bA,\bB,\bC$, and $\bD$; see, for example,  \citet{lu2021numerical}}.

\item \textbf{Properties of expectation.}
Let $\rx$ and $\ry$ be two random variables, and let $a, b$ be scalars. Show that $\Exp[a\rx+b\ry] = a\cdot\Exp[\rx]+b\cdot \Exp[\ry]$.
Given further a  function $h$, show that 
$$
\Exp[h(\rx)] = \sum_{x} h(x) \prob(x)
\qquad \text{and}\qquad 
\Exp[h(\rx)] = \int_{-\infty}^{\infty} h(x) d F(x),
$$
in the discrete and continuous cases, respectively.

\item \label{prob:prop_expcond}\textbf{Properties of expectation of conditionals.}
Let $\rx$ and $\ry$ be two random variables, and let $h$ be a  function.
Show that
\begin{enumerate}
\item Note that $\Exp[\rx\mid\ry]$ is a function of $\ry$. However, if $\rx$ is independent of $\ry$, then $\Exp[\rx\mid \ry] = \rx$.
\item $\Exp[c\mid \rx] =\rx$, where $c$ is a constant.
\item Linearity: $\Exp[a\rx_1 + b\rx_2\mid \ry] = a\cdot\Exp[\rx_1\mid \ry] + b\cdot\Exp[\rx_2\mid \ry]$.
\item Conditional constant: $\Exp[h(\ry)\rx\mid \ry] = h(\ry)\Exp[\rx\mid \ry]$, where $h(\ry)$ is called a \textit{conditional constant} w.r.t. $\ry$.
\item Monotonicity: if $\rx_1 \leq \rx_2$, then $\Exp[\rx_1\mid \ry] \leq \Exp[\rx_2\mid \ry]$.
\item Tower property: $\Exp\left[\Exp[\rx\mid \ry]\mid h(\ry)\right] = \Exp[\rx\mid h(\ry)]$; that is, $h(\ry)$ conveys information at most as $\ry$.
\item Unbiasedness: $\Exp\big\{\Exp[h(\rx,\ry)\mid \ry]\big\} = \Exp[h(\rx,\ry)]$; specially, $\Exp\big[\Exp[\rx\mid \ry]\big] = \Exp[\rx]$.
\item Least squares: $\Exp\big[\big(\ry - \Exp[\ry\mid \rx]\big)^2\big] \leq \Exp\big[\big(\ry - h(\rx)\big)^2\big]$ for any function $h$. This also means  $g(\rx)\triangleq  \Exp[\ry\mid \rx]$ is the best estimate in the least squares sense.
\end{enumerate}

\item \textbf{Sum of random variables by convolution.}
Let $ \rx $ and $ \ry $ be continuous random variables with probability density functions $ f_{\rx} $ and $ f_{\ry} $. Show that the density function of $ \rx + \ry $ is the convolution of $ f_{\rx} $ with $ f_{\ry} $:
$$
f_{\rx+\ry}(u) = \int_{-\infty}^{+\infty} f_{\rx}(u - v) f_\ry(v) \, dv.
$$

\item \textbf{Properties of variance and correlation.}
Let $\rx, \rx_1, \rx_2, \ry, \rvx$ be random variables or vectors, and  let $a,b$ be constants. 
Show that 
\begin{itemize}
\item Let $\bOmega$ be a real symmetric matrix. Then $\bOmega$ is positive semidefinite (Definition~\ref{definition:psd-pd-defini}) if and only if $\bOmega$ is the covariance matrix of some random vector $\rvx$.
\item $\Var[\rx] = \Exp[\rx^2] - (\Exp[\rx])^2 = \Cov[\rx,\rx]$.
\item $\Var[a\rx + b] = a^2 \Var[\rx]$.
\item $\Var\big[\sum_i \rx_i\big] = \sum_i \Var[\rx_i] + \sum_{i \neq j} \Cov[\rx_i, \rx_j]$.
\item $\Cov[\rx_1, \rx_2] = \Exp[\rx_1 \rx_2] - \Exp[\rx_1]\Exp[\rx_2]$.
\item $\Cov[a\rx_1 + b\rx_2, \ry] = a\cdot\Cov[\rx_1, \ry] + b\cdot \Cov[\rx_2, \ry]$; that is, covariance is linear in one variable.
\item If $\Exp[\rx_1^2] + \Exp[\rx_2^2] < \infty$, then the following are equivalent:
\begin{enumerate}
\item[(i)] $\Exp[\rx_1 \rx_2] = \Exp[\rx_1]\Exp[\rx_2]$;
\item[(ii)] $\Cov[\rx_1, \rx_2] = 0$;
\item[(iii)] $\Var[\rx_1 \pm \rx_2] = \Var[\rx_1] + \Var[\rx_2]$.
\end{enumerate}
Note that independence will imply these three last properties, but none of these properties imply independence.
\item Let $h$ be a nondecreasing function such that $\Exp[\rx^2]<\infty$ and $\Exp[h(\rx)^2]<\infty$. Then $\Cov[\rx, h(\rx)]>0$.
\end{itemize}

\item \textbf{Conditional variance and law of total variance.}
Let $\rx$ and $\ry$ be two random variables. The \textit{condition variance} of $\rx$ given $\ry$ is defined as 
\begin{equation}
\Var[\rx \mid \ry] \triangleq \Exp\big[ \big(\rx - \Exp[\rx \mid \ry]\big)^2 \mid \ry \big] =\Exp[\rx^2\mid \ry ] - (\Exp[\rx \mid \ry ])^2.
\end{equation}
The conditional variance tells us how much variance is left if we use $\Exp[\rx\mid \ry]$ to ``predict" $\rx$.
Prove the \textit{law of total variance}:
\begin{equation}
\Var[\rx] = \Exp\big[\Var[\rx\mid \ry]\big] + \Var\big[\Exp[\rx\mid \ry]\big].
\end{equation}
\textit{Hint: $\Var[\rx] = \Exp[\rx^2] -(\Exp[\rx])^2 = \Exp[\Exp[\rx^2\mid \ry ]] -(\Exp[\Exp[\rx\mid \ry]])^2$ by the unbiasedness property in Problem~\ref{prob:prop_expcond}.}
\end{problemset}

\part{Theory}
\newpage 
\chapter{Basics of Optimization}\label{chapter:mathback}
\begingroup
\hypersetup{
linkcolor=structurecolor,
linktoc=page,  
}
\minitoc \newpage
\endgroup

\lettrine{\color{caligraphcolor}T}
This chapter provides a concise yet self-contained introduction to the fundamental concepts of convex analysis and convex optimization. It covers essential topics such as convex sets and functions, subdifferentials, and optimality conditions. While the presentation is designed to equip the reader with the core tools needed for applications in statistics, machine learning, and signal processing, it is not intended to be exhaustive. 
For a deeper and more comprehensive treatment---including advanced theoretical results, duality theory, constrained optimization, and specialized algorithms---the reader is encouraged to consult standard references in the field, such as  \citet{bertsekas90001convex, boyd2004convex, bertsekas2009convex, foucart2013invitation, rockafellar2015convex, lu2025practical}.

\section{Differentiable Functions and Differential Calculus}\label{section:differ_calc}

Differentiability and differential calculus form the backbone of mathematical analysis, particularly in the study of functions defined on  multidimensional spaces. 
This section introduces the fundamental concepts that allow us to quantify how a function changes with respect to its input variables.

Central to this discussion is the directional derivative, which measures the instantaneous rate of change of a function $f$ at a point $\balpha$ in the direction of a given vector $\bd$.

\begin{definition}[Directional derivative, partial derivative\index{Directional derivative}\index{Partial derivative}]\label{definition:partial_deri}
Let $f$ be a real-valued function defined on a set $\sS\subseteq \real^p$, and let $\balpha\in\real^p$. 
For any nonzero vector $\bd\in\real^p$, the \textit{directional derivative} of $f$ at $\balpha\in\real^p$ in the direction $\bd$ is defined---provided the limit exists---by
$$
f^\prime(\balpha; \bd)\triangleq 
\mathop{\lim}_{\mu\rightarrow 0}
\frac{f(\balpha+\mu\bd) - f(\balpha)}{\mu}.
$$
This quantity is also sometimes referred to as the   \textit{G\^ateaux derivative}.

For each $i\in\{1,2,\ldots,p\}$, the directional derivative at $\balpha$  in the direction of the $i$-th standard basis vector  $\be_i$ (if it exists) is called the $i$-th \textit{partial derivative}. 
It is commonly denoted by $\frac{\partial f}{\partial \alpha_i} (\balpha)$, $D_{\be_i}f(\balpha)$, or $\partial_i f(\balpha)$.
\end{definition}

\begin{figure*}[h]
\centering  
\subfigtopskip=2pt 
\subfigbottomskip=9pt 
\subfigcapskip=-5pt 
\includegraphics[width=0.75\textwidth]{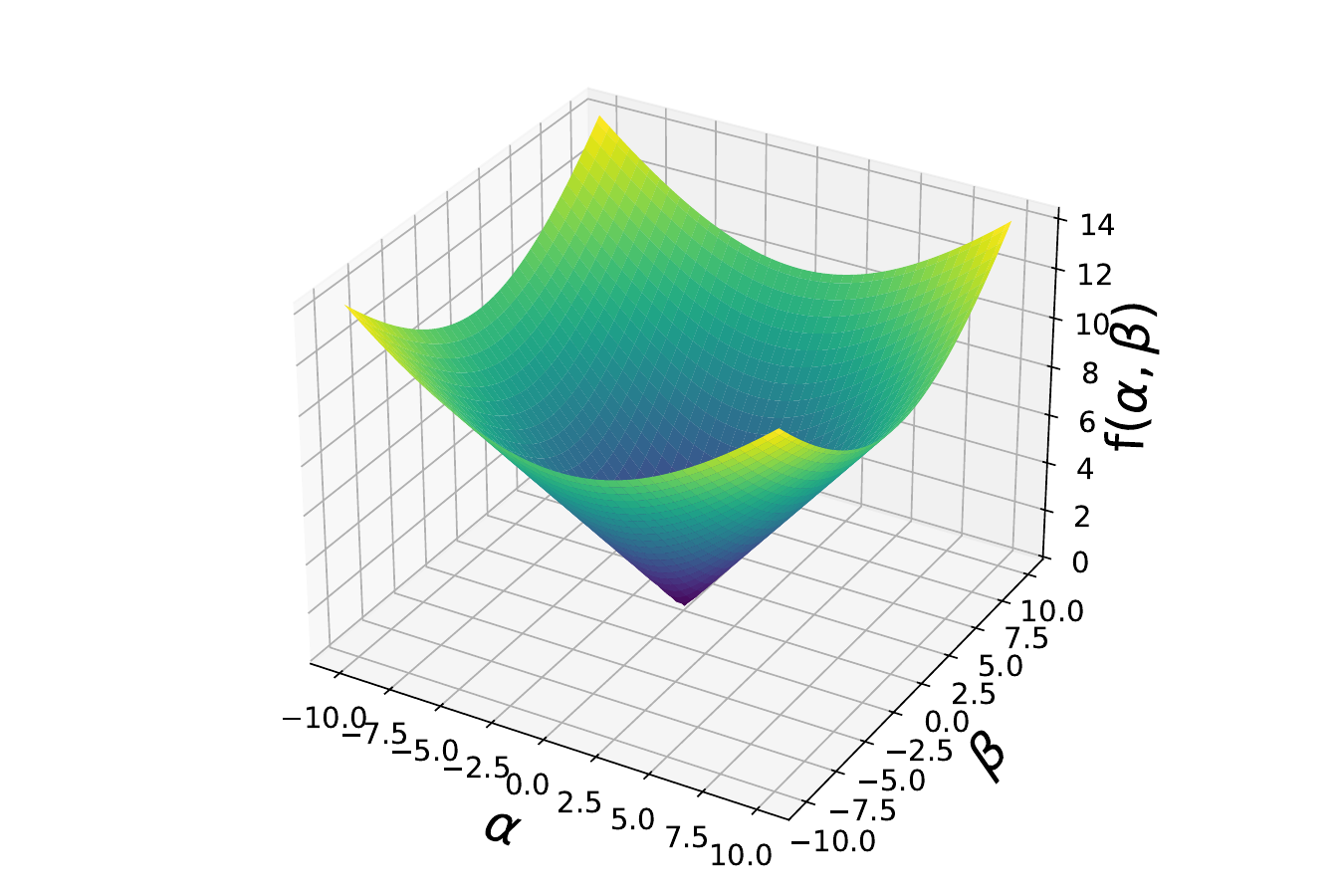}
\caption{Plot for the function $f(\alpha, \beta) = \sqrt{\alpha^2+\beta^2}$. 
The partial derivatives at $[0,0]^\top$ do not exist.}
\label{fig:direc_deri_no_partial}
\end{figure*}

It is important to note that the existence of directional derivatives in all nonzero directions at a point does not guarantee the existence of partial derivatives---or, more generally, differentiability---at that point.
A classic example is the function $f(\alpha,\beta)=\sqrt{\alpha^2+\beta^2}$ (see Figure~\ref{fig:direc_deri_no_partial}).
Consider a direction   $\bd=[a,b]^\top$. 
The partial derivatives at $[0,0]^\top$ are 
$$
\begin{aligned}
\frac{\partial f}{\partial \alpha}(0,0) &= \mathop{\lim}_{h\rightarrow 0} \frac{\sqrt{(0+h)^2+0^2}}{h} = \mathop{\lim}_{h\rightarrow 0} \frac{\abs{h}}{h},\\
\frac{\partial f}{\partial \beta}(0,0) &= \mathop{\lim}_{h\rightarrow 0} \frac{\sqrt{0^2+(0+h)^2}}{h} = \mathop{\lim}_{h\rightarrow 0} \frac{\abs{h}}{h}.
\end{aligned}
$$
When $h>0$, the partial derivatives are 1; when $h<0$, the partial derivaties are $-1$. Since the left- and right-hand limits differ, the limit does not exist.

\begin{definition}[Gradient\index{Gradient}]
If all the partial derivatives of a function $f$ exist at a point $\balpha\in\real^p$, then the \textit{gradient} of $f$ at $\balpha$ is defined as the column vector consisting of these partial derivatives:
$$
\nabla f(\balpha)=
\begin{bmatrix}
\frac{\partial f}{\partial \alpha_1} (\balpha),
\frac{\partial f}{\partial \alpha_2} (\balpha),
\ldots,
\frac{\partial f}{\partial \alpha_p} (\balpha)	
\end{bmatrix}^\top
\in \real^p.
$$

\end{definition}

\begin{exercise}
Show that a function may have all partial derivatives at a point yet fail to be continuous there (see Definition~\ref{definition:conti_funs}).
Conversely, demonstrate that a function can be continuous everywhere but still lack partial derivatives at some points.
\end{exercise}

A function $f$ defined on an open set $\sS\subseteq \real^p$ (see Definition~\ref{definition:open_close_sets}) is said to be \textit{differentiable}
at a point $\balpha\in\sS$ if it admits a linear approximation near $\balpha$.
This notion is formalized by \textit{Fr\'echet differentiability}, which is stronger than merely requiring the existence of partial derivatives.

\begin{definition}[(Fr\'echet) differentiability\index{(Fr\'echet) differentiability}]
Let $f:\real^p\rightarrow \real$. 
The function $f$ is said to be \textit{differentiable} or \textit{Fr\'echet differentiability} at $\balpha$ if there exists a vector $\bg\in\real^p$ such that 
$$
\mathop{\lim}_{\bd\rightarrow\bzero} \frac{f(\balpha+\bd)-f(\balpha)-\bg^\top\bd}{\normtwo{\bd}}=0.
$$
The  vector $\bg$ is unique and equals the gradient $\nabla f(\balpha)$.
\end{definition}

\index{Continuously differentiable}
Moreover, a function $f$ defined on an open set $\sS\subseteq \real^p$ is called \textit{continuously differentiable} on $\sS$ if all its partial derivatives exist and are also continuous on $\sS$.

Note that a function can be differentiable without being continuously differentiable. 
A classical example is the univariate function:
$$
f(\alpha) = 
\begin{cases}
\alpha^2 \sin(\frac{1}{\alpha}), & \text{if } \alpha \neq 0, \\
0, & \text{if } \alpha = 0.
\end{cases}
$$
This function is differentiable everywhere on $\real$, including at $\alpha=0$. To see this, for $\alpha \neq 0$, $f(\alpha)$ is a product of two differentiable functions ($\alpha^2$ and $\sin(1/\alpha)$), hence it is differentiable.
At $\alpha = 0$, the derivative is computed via the limit definition:
$$
f^\prime(0) = \mathop{\lim}_{\mu\rightarrow 0}
\frac{f(\mu) - f(0) }{\mu}
=\mathop{\lim}_{\mu\rightarrow 0} \frac{\mu^2 \sin(\frac{1}{\mu}) - 0 }{\mu}
=\mu \sin(\frac{1}{\mu}) 
$$
Since $\abs{\sin(1/\mu)}\leq 1$ for all $\mu \neq 0$, the limit exists and is equal to 0. Thus, $f(\alpha)$ is differentiable at $\alpha = 0$ with $f^\prime(0) = 0$.
However,  $f(\alpha)$ is not continuously differentiable at $\alpha = 0$. The derivative of $f(\alpha)$ when $\alpha\neq 0$ is 
$$
f^\prime(\alpha) = 2\alpha \sin(\frac{1}{\alpha}) - \cos(\frac{1}{\alpha}).
$$
The limit $\mathop{\lim}_{\alpha\rightarrow 0}f^\prime(\alpha)$ does not exist because the sine and cosine  functions oscillate as $\alpha$ approaches 0.

\begin{lemma}[Gradient formula for directional derivatives]\label{lemma:direc_contdiff}
If $f$ is differentiable at $\balpha\in\sS\subseteq\real^p$, then for every direction $\bd\in\real^p$, the directional derivative of $f$ at $\balpha$ in the direction $\bd$ satisfies
\begin{equation}\label{equation:direc_contdiff}
f^\prime(\balpha; \bd) = \nabla f(\balpha)^\top \bd, \gap \text{for all }\balpha\in\sS \text{ and }\bd\in\real^p.
\end{equation} 
\end{lemma}
\begin{proof}[of Lemma~\ref{lemma:direc_contdiff}]
The formula is obviously correct for $\bd = \bzero$. We then assume that $\bd \neq \bzero$. The differentiability of $f$ implies that
$$
0 = \lim_{\mu \to 0^+} \frac{f(\balpha + \mu \bd) - f(\balpha) - \innerproduct{\nabla f(\balpha), \mu \bd}}{\normtwo{\mu \bd}} 
= \lim_{\mu \to 0^+} \left[ \frac{f(\balpha + \mu \bd) - f(\balpha)}{\mu \normtwo{\bd}} - \frac{\innerproduct{\nabla f(\balpha), \bd}}{\normtwo{\bd}} \right].
$$
Therefore,
$$
\begin{aligned}
&f'(\balpha; \bd) = \lim_{\mu \to 0^+} \frac{f(\balpha + \mu \bd) - f(\balpha)}{\mu}\\
&= \lim_{\mu \to 0^+} \left\{ \normtwo{\bd} \left[ \frac{f(\balpha + \mu \bd) - f(\balpha)}{\mu \normtwo{\bd}} - \frac{\innerproduct{\nabla f(\balpha), \bd}}{\normtwo{\bd}} \right] +\innerproduct{\nabla f(\balpha), \bd} \right\}
&= \innerproduct{ \nabla f(\balpha), \bd}.
\end{aligned}
$$
This proves \eqref{equation:direc_contdiff}.
\end{proof}

Recalling the differentiability of $f$ at $\balpha$, we have 
\begin{equation}
\mathop{\lim}_{\bd\rightarrow \bzero}
\frac{f(\balpha+\bd) - f(\balpha) - \nabla f(\balpha)^\top \bd}{\normtwo{\bd}} = 0,
\quad 
\text{for all }\balpha\in\sS,
\end{equation}
or 
\begin{equation}
f(\bbeta) = f(\balpha)+\nabla f(\balpha)^\top (\bbeta-\balpha) + o(\normtwo{\bbeta-\balpha}),
\end{equation}
where the \textit{small-oh} function $o(\cdot): \real_+\rightarrow \real$ is a one-dimensional function satisfying $\frac{o(\mu)}{\mu}\rightarrow 0$ as $\mu\rightarrow 0^+$.~\footnote{Note that we also use the standard \textit{big-Oh} notation to describe the asymptotic behavior of functions.
Specifically, the notation $g(\bd) = \mathcalO(\normtwo{\bd}^p)$ means that there are positive numbers
$C_1$ and $\delta$ such that $\abs{g(\bd)} \leq  C_1 \normtwo{\bd}^p$ for all $\normtwo{\bd}\leq\delta$. In practice it is
often equivalent to $\abs{g(\bd)} \approx C_2\normtwo{\bd}^p$ for  sufficiently small $\bd$, where $C_2$ is  another positive constant.
The \textit{soft-Oh} notation is employed to hide poly-logarithmic factors i.e., $f = \widetilde{\mathcalO}(g)$ will
imply $f = \mathcalO(g \log^c(g))$ for some absolute constant $c$.}
Therefore, any differentiable function ensures that $\mathop{\lim}_{\balpha\rightarrow \ba} f(\balpha) = f(\ba)$. Hence, \textbf{any differentiable function is continuous}.

\index{Descent direction}
\paragrapharrow{Descent direction and steepest direction.}
An important consequence of the gradient-directional derivative relationship \eqref{equation:direc_contdiff} is that the gradient encodes the directions of steepest ascent and descent.
\begin{lemma}[Steepest direction]\label{lemma:steep_direction}
Let $f:\real^p \rightarrow \real$ be differentiable at $\balpha$. Among all unit directions $\bd\in\real^p$ (i.e., $\normtwo{\bd}=1$), the directional derivative $f'(\balpha;\bd)$ is maximized when  $\bd $ points in the direction of  $\nabla f(\balpha) $, and minimized (most negative) when $\bd$ points in the opposite direction.
\end{lemma}

\begin{proof}[of Lemma~\ref{lemma:steep_direction}]
For all $\bd\in\real^p$ satisfying $\normtwo{\bd}=1$,
by the Cauchy--Schwarz inequality (Theorem~\ref{theorem:cs_matvec}), we have 
\begin{align*}
\abs{f'(\balpha; \bd)} 
&= \abs{\nabla f(\balpha)^\top \bd} 
\leq \normtwo{\nabla f(\balpha)} \normtwo{\bd} 
= \normtwo{\nabla f(\balpha)} .
\quad 
\end{align*}
Equality holds if and only if $\bd$ is a scalar multiple of ${\nabla f(\balpha)}$. Since $\normtwo{\bd}=1$, the maximizing direction is $\bd = \pm \frac{\nabla f(\balpha)}{\normtwo{\nabla f(\balpha)}}$.
\end{proof}

\begin{SCfigure}
\centering
\includegraphics[width=0.45\textwidth]{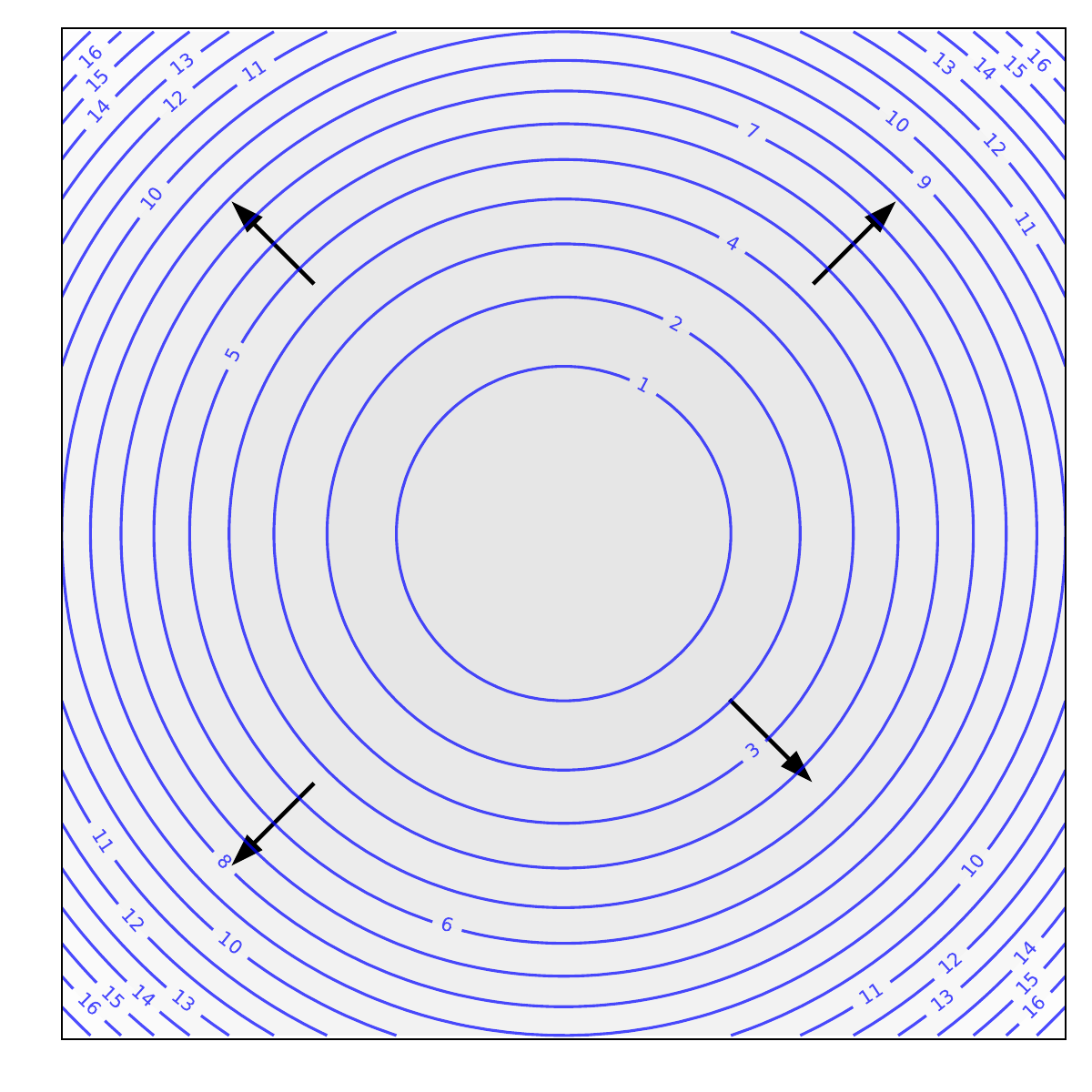}
\caption{Contour lines of a function $f : \real^2 \to \real$. The gradients at different points are represented by
	black arrows, which are orthogonal to the contour lines.}
\label{fig:gradient_directions}
\end{SCfigure}

Figure~\ref{fig:gradient_directions} illustrates the gradient field of a function $f : \real^2 \to \real$. 
The gradient vectors (shown as black arrows) are orthogonal to the contour lines. This orthogonality arises because the function is constant along each contour line, so the directional derivative in any tangent direction is zero. By \eqref{equation:direc_contdiff}, this implies $\nabla f(\balpha)^\top\bd=0$ for all tangent vectors $\bd$, meaning the gradient must be normal to the level set.

Note that $\bd=-\nabla f(\balpha)$ is called the \textit{steepest descent direction}. More generally, any direction $\bd$ satisfying 
$$
f'(\balpha;\bd)<0
$$  
is called a \textit{descent direction} at $\balpha$.
\begin{lemma}[Unconstrained and constrained descent property]\label{lemma:descent_property}
Let $ f $ be a continuously differentiable function over $\sS\subseteq\real^p$ (resp., $\real^p$), and let $\balpha \in \sS$ (resp., $\balpha\in\real^p$). Suppose that $\bd$ is a descent direction of $ f $ at $\balpha$. Then, there exists $\varepsilon > 0$ such that 
$$
f(\balpha + \eta \bd) < f(\balpha),\quad \text{for any $ \eta \in (0, \varepsilon] $.}
$$
\end{lemma}
\begin{proof}[of Lemma~\ref{lemma:descent_property}]
Since $ f'(\balpha; \bd)=\nabla f(\balpha)^\top \bd < 0 $, it follows from the definition of the directional derivative (Definition~\ref{definition:partial_deri}) that
$$
\lim_{\eta \to 0^+} \frac{f(\balpha + \eta \bd) - f(\balpha)}{\eta} = f'(\balpha; \bd) < 0.
$$
Therefore, there exists an $\varepsilon > 0$ such that
$
\frac{f(\balpha + \eta \bd) - f(\balpha)}{\eta} < 0
$
for any $ \eta \in (0, \varepsilon] $, which completes the proof.
\end{proof}

\subsection{Second-Order Partial Derivative}
The partial derivative $\frac{\partial f}{\partial \alpha_i} (\balpha)$ is itself a real-valued function of $\balpha\in\sS$, and it may be partially differentiated again. 
The $j$-th partial derivative of $\frac{\partial f}{\partial \alpha_i} (\balpha)$ is defined as 
$$
\frac{\partial^2 f}{\partial \alpha_j\partial \alpha_i} (\balpha)=
\frac{\partial \left(\frac{\partial f}{\partial \alpha_i} (\balpha)\right)}{\partial \alpha_j} (\balpha).
$$
This is called the ($j,i$)-th \textit{second-order partial derivative} of function $f$.

A function $f$ defined on an open set $\sS\subseteq$ is called \textit{twice continuously differentiable} over $\sS$ if all its second-order partial derivatives exist and are continuous on $\sS$. 
In the setting of twice continuously differentiability, 
the mixed partial derivatives are symmetric:
$$
\frac{\partial^2 f}{\partial \alpha_j\partial \alpha_i} (\balpha)=
\frac{\partial^2 f}{\partial \alpha_i\partial \alpha_j} (\balpha).
$$
The \textit{Hessian} of the function $f$ at a point $\balpha\in\sS$ 
is the $p\times p$
matrix of second-order partial derivatives:
$$
\nabla^2f(\balpha)=
\begin{bmatrix}
\frac{\partial^2 f}{\partial \alpha_1^2} (\balpha) & 
\frac{\partial^2 f}{\partial \alpha_1\partial \alpha_2} (\balpha) & \ldots &
\frac{\partial^2 f}{\partial \alpha_1\partial \alpha_p} (\balpha)\\
\frac{\partial^2 f}{\partial \alpha_2\partial \alpha_1} (\balpha) & 
\frac{\partial^2 f}{\partial \alpha_2\partial \alpha_2} (\balpha) & \ldots &
\frac{\partial^2 f}{\partial \alpha_2\partial \alpha_p} (\balpha)\\
\vdots & 
\vdots & \ddots &
\vdots\\
\frac{\partial^2 f}{\partial \alpha_p\partial \alpha_1} (\balpha) & 
\frac{\partial^2 f}{\partial \alpha_p\partial \alpha_2} (\balpha) & \ldots &
\frac{\partial^2 f}{\partial \alpha_p^2} (\balpha)
\end{bmatrix}.
$$
When $f$ is twice continuously differentiable, this matrix is symmetric due to the equality of mixed partials.

\begin{definition}[Twice differentiable]
By its definition, a  differentiable function $f:\sS \subseteq\real^p\rightarrow \real$ is \textit{twice differentiable} at a point $\balpha\in \sS$ if there exist a matrix $\nabla^2 f(\balpha)\in\real^{p\times p}$ such that 
$$
\mathop{\lim}_{\bd\rightarrow\bzero} 
\frac{\normtwo{\nabla f(\balpha+\bd)-\nabla f(\balpha)-\nabla^2 f(\balpha)\bd}}{\normtwo{\bd}}=0.
$$
The  matrix $\nabla^2 f(\balpha)\in\real^{p\times p}$ is unique and coincides with the Hessian matrix of $f$ at $\balpha$.
\end{definition}

Geometrically, the curvature of a function describes how quickly its slope changes; in one dimension, this is captured by the second derivative. In higher dimensions, the Hessian matrix encodes the curvature of the function in every direction.

Analogously to the gradient formula for directional derivatives (Lemma~\ref{lemma:direc_contdiff}), the second directional derivative is related to the Hessian as follows.

\begin{lemma}[Hessian-directional derivative identity]\label{lemma:direc_twicontdiff}
If a function $f : \sS\subseteq\real^p \to \real$ is twice differentiable, the second directional derivative $f''(\balpha; \bd)$ of $f$ at $\balpha\in \sS$ equals
\begin{equation}\label{equation:direc_twicontdiff}
f''(\balpha; \bd) = \bd^\top \nabla^2 f(\balpha) \bd,
\end{equation}
for any unit-norm vector $\bd \in \real^p$.
\end{lemma}

\subsection{Functions of Several Real Variables}
We now present several classical approximation theorems that rely on derivatives and related concepts.

\begin{theoremHigh}[Rolle's theorem\index{Rolle's theorem}]\index{Rolle's theorem}\label{theorem:rolles_theo}
Let $f(\alpha)$ be continuous on the closed interval $[a, b]$ and differentiable on $(a, b)$. 
If $f(a) = f(b)$, then there exist at least one point $\xi\in(a, b)$ such that  $f'(\xi) = 0$.
\end{theoremHigh}
\begin{proof}[of Theorem~\ref{theorem:rolles_theo}]
We simply note that either $f$ is constant, in which case $f'(\alpha) = 0$ for all $\alpha$ in $(a,b)$, or it has a local maximum or minimum at $\xi\in(a,b)$, in which case $f'(\xi) = 0$.
\end{proof}
Rolle's Theorem is primarily used as a key step in proving the ordinary mean value theorem (OMVT).
\begin{theoremHigh}[Ordinary mean value theorem (OMVT)]\label{theorem:uni_mvt}
Let $f(\alpha)$ be continuous on the closed interval $[a, b]$ and differentiable on $(a, b)$. Then, there exists a point $\xi\in(a, b)$ such that 
$$
f(b) - f(a) = f'(\xi)(b - a).
$$
\end{theoremHigh}
\begin{proof}[of Theorem~\ref{theorem:uni_mvt}]
Define the auxiliary function
$
g(\alpha) \triangleq f(\alpha) - \left(\frac{f(b) - f(a)}{b - a}\right)(\alpha - a).
$
It holds that $g(a) = g(b)$. Then by Rolle's theorem, there exists $\xi \in (a,b)$ such that $g'(\xi) = 0$, which rearranges to the desired identity:
$$
f(b) - f(a) = f'(\xi)(b - a).
$$
This completes the proof.
\end{proof}

For functions of several variables, an analogous result holds, though it takes a slightly different form due to the vector-valued nature of the gradient.
\begin{theoremHigh}[Mean value theorem (MVT)\index{Mean value theorem}]\label{theorem:mean_approx}
Let $f:\sS\rightarrow \real$ be a  continuously differentiable function on an open set $\sS\subseteq\real^p$, and given two points $\balpha, \bbeta\in\sS$. Then, there exists a point $\bxi\in[\balpha,\bbeta]$ such that 
$$
f(\bbeta) = f(\balpha)+ \nabla f(\bxi)^\top (\bbeta-\balpha).
$$

\end{theoremHigh}
\begin{proof}[of Theorem~\ref{theorem:mean_approx}]
Without loss of generality, we assume $\bxi(t) = \balpha + t(\bbeta - \balpha)$ with $t\in[0,1]$. 
Then we define $g(t) \triangleq f(\bxi(t))$. We can apply the ordinary mean value theorem to $g(t)$, to get
$$
g(1) = g(0) + g'(\gamma),
$$
for some $\gamma$ in the interval $[0, 1]$. The derivative of $g(t)$ is
$$
g'(t) = \innerproduct{\nabla f(\bxi(t)), \bbeta - \balpha},
$$
where
$
\nabla f(\bxi(t)) = \left[\frac{\partial f}{\partial \alpha_1}(\bxi(t)), 
\frac{\partial f}{\partial \alpha_2}(\bxi(t)),
\ldots, \frac{\partial f}{\partial \alpha_p}(\bxi(t))\right]^\top.
$
Therefore,
$$
g'(\gamma) = \innerproduct{\nabla f(\bxi(\gamma)), \bbeta - \balpha}.
$$
Since $\bxi(\gamma) = (1 - \gamma)\balpha + \gamma \bbeta$, this completes the proof.
\end{proof}

\begin{theoremHigh}[Linear approximation theorem\index{Linear approximation}]\label{theorem:linear_approx}
Let $f:\sS\rightarrow \real$ be a twice continuously differentiable function on an open set $\sS\subseteq\real^p$, and let $\balpha, \bbeta\in\sS$. Then, there exists a point $\bxi\in[\balpha,\bbeta]$ such that 
$$
\begin{aligned}
f(\bbeta) &= f(\balpha)+ \nabla f(\balpha)^\top (\bbeta-\balpha) + \frac{1}{2} (\bbeta-\balpha)^\top \nabla^2 f(\bxi) (\bbeta-\balpha),\\
\text{or}\quad f(\bbeta) &= f(\balpha)+\nabla f(\balpha)^\top (\bbeta-\balpha) + o(\normtwo{\bbeta-\balpha}),\\
\text{or}\quad f(\bbeta) &= f(\balpha)+\nabla f(\balpha)^\top (\bbeta-\balpha) + \mathcalO(\normtwo{\bbeta-\balpha}^2).
\end{aligned}
$$
\end{theoremHigh}
This theorem, also known as  the \textit{first-order Taylor series expansion}, suggests that the error in the linear approximation is of the order of the square of the distance between  $\balpha$ and $\bbeta$.
And the  first-order Taylor series expansion $f(\bbeta) \approx  f(\balpha)+\nabla f(\balpha)^\top (\bbeta-\balpha)$ is exact only for the neighborhood of $\bbeta$ at the point $\balpha$,

\begin{theoremHigh}[Quadratic approximation theorem]\label{theorem:quad_app_theo}
Let $f:\sS\rightarrow \real$ be a twice continuously differentiable function on an open set $\sS\subseteq\real^p$, and let $\balpha, \bbeta\in\sS$. Then it follows that 
$$
\begin{aligned}
f(\bbeta) &= f(\balpha)+ \nabla f(\balpha)^\top (\bbeta-\balpha) + \frac{1}{2} (\bbeta-\balpha)^\top \nabla^2 f(\balpha) (\bbeta-\balpha)
+
o(\normtwo{\bbeta-\balpha}^2),\\
\text{or}\quad f(\bbeta) 
&= f(\balpha)+ \nabla f(\balpha)^\top (\bbeta-\balpha) + \frac{1}{2} (\bbeta-\balpha)^\top \nabla^2 f(\balpha) (\bbeta-\balpha)
+
\mathcalO(\normtwo{\bbeta-\balpha}^3).
\end{aligned}
$$
\end{theoremHigh}
This theorem indicates that the error in the quadratic approximation is of the order of the cube of the distance between $\balpha$ and $\bbeta$, making it a more accurate approximation when $\bbeta$ is close to $\balpha$.

\subsection{Fundamental Theorem of Calculus}

We now present the \textit{undamental theorem of calculus} in a multivariate setting. This result plays a pivotal role in connecting differential and integral calculus, offering deep insight into the structure and behavior of smooth functions.
\begin{theoremHigh}[Fundamental theorem of calculus]\label{theorem:fund_theo_calculu}
Let $f:\real^p\rightarrow \real$ be a continuously differentiable function, and let $\bbeta, \balpha\in\real^p$. 
Then the difference between function values can be expressed as
\begin{subequations}
\begin{equation}\label{equation:fund_theo_calculu3}
f(\bbeta) - f(\balpha) = \int_{0}^{1} \innerproduct{\nabla f(\balpha+\mu(\bbeta-\balpha)), \bbeta-\balpha } d\mu.
\end{equation}
If, in addition, $f$ is twice continuously differentiable, then the difference between gradients satisfies
\begin{align}
\nabla f(\bbeta) -\nabla f(\balpha) &= \left( \int_{0}^{1} \nabla^2f(\balpha+\mu(\bbeta-\balpha))d\mu \right) \cdot (\bbeta-\balpha); \label{equation:fund_theo_calculu1}\\
\nabla f(\balpha+\eta\bd) - \nabla f(\balpha) &= \int_{0}^{\eta} \nabla^2 f(\balpha+\mu\bd)\bd d\mu.
\end{align}
where $ \nabla f(\bbeta) $ denotes the gradient of $ f $ evaluated at $ \bbeta $,
$\nabla^2 f(\cdot)$ the Hessian matrix, and $ \innerproduct{\cdot, \cdot} $ the standard Euclidean inner product.
Furthermore, use directional derivative for twice continuously differentiable functions, we also have
\begin{equation}
f(\bbeta) = f(\balpha) + \innerproduct{\nabla f(\balpha),  (\bbeta-\balpha)} + \int_{0}^{1} (1-\mu) \frac{\partial^2 f(\balpha+\mu(\bbeta-\balpha))}{\partial \mu^2} d\mu.
\end{equation}
\end{subequations}
\end{theoremHigh}
The first equality \eqref{equation:fund_theo_calculu3} can be derived by defining the auxiliary scalar function $g(t)\triangleq f(\balpha+t(\bbeta-\balpha))$. Then we can obtain $f(\bbeta)$ as 
$
f(\bbeta) = g(1) = g(0) + \int_{0}^{1} g'(t)dt.
$
The remaining identities can be demonstrated similarly.

\section{Convex Combinations and Convex  Sets}
This section introduces the concepts of affine, convex, and conic combinations and their associated hulls, in the context of vector sets in $\real^p$. 
We provide precise definitions and mathematical formulations for these fundamental geometric constructs.

\subsection*{Interior Points and Closed Sets}

Given a set of vectors $\balpha_1,\balpha_2,\ldots,\balpha_k\in\real^p$, 
an \textit{affine combination} of these $k$ vectors is a vector of the form  $\sum_{i=1}^{k} \lambda_i\balpha_i$, where $\sum_{i=1}^{k} \lambda_i=1$, i.e., a linear combination of these points where the coefficients (not necessarily nonnegative) sum to 1.
A \textit{convex combination} of these $k$ vectors is a vector of the form $\sum_{i=1}^{k} \lambda_i\balpha_i$, where $\lambda_i\geq 0$ for $i\in\{1,2,\ldots,k\}$ satisfying $\sum_{i=1}^{k}\lambda_i=1$ (i.e., $\{\lambda_i\}$ belongs to the unit-simplex in $\real^p$: $\blambda\in\Delta_k$); 
when only requiring $\lambda_i\geq 0$ without the constraint on their sum, the combination is referred to as a \textit{conic combination} of these vectors.

After fixing a norm on $\real^p$, we introduced open and closed balls in Definition~\ref{definition:open_closed_ball}. 
These notions lead naturally to the concepts of interior and relative interior points, which are foundational in topology and convex analysis.

\begin{definition}[Affine hulls\index{Affine hulls}]\label{definition:aff_hulls}
The \textit{affine hull} of a set $ \sS $, denoted  $ \aff(\sS) $, is the smallest affine subspace (a translation of a vector subspace) that contains $ \sS $.
Formally, given a set $ \sS \subseteq \real^p $, the affine hull of $ \sS $ is defined as the set of all affine combinations of points from $ \sS $:
$$
\aff(\sS) \triangleq \left\{ \sum_{i=1}^k \lambda_i \balpha_i \mid   \balpha_1,\balpha_2,\ldots,\balpha_k\in\sS, \lambda_i \in \real, \bone^\top\blambda=1, k\in\naturalset \right\}.
$$

\end{definition}

\index{Interior points}
\index{Relative interior points}
\begin{definition}[Interior and relative interior points]
The \textit{interior} of a set $ \sS $ in a topological space is the set of all points in $ \sS $ that admit a neighborhood entirely contained in $ \sS $. 
Formally, the interior of $ \sS $, denoted  $ \interior(\sS) $, is defined as:
$$
\interior(\sS) = \{ \balpha \in \sS \mid \exists\, \epsilon > 0 \text{ such that } \sB(\balpha, \epsilon) \subseteq \sS \}.
$$
where $ \sB(\balpha, \epsilon) $ denotes the open ball centered at $ \balpha $ with radius $ \epsilon $.

The \textit{relative interior} of a set $ \sS $ in a topological space is the interior of $ \sS $ relative to its affine hull.  
The relative interior of $ \sS $, denoted  $ \relint(\sS) $, is given by:
$$
\relint(\sS) = \{ \balpha \in \sS \mid \exists\,\epsilon > 0 \text{ such that } \sB_{\aff(\sS)}(\balpha, \epsilon) \subseteq \sS \}.
$$
where $ \sB_{\aff(\sS)}(\balpha, \epsilon) $ is the open ball centered at $ \balpha $ with radius $ \epsilon $ in the subspace topology of $ \aff(\sS) $.
\end{definition}

In other words, the \textit{interior} of a set is concerned with the open balls in the ambient space; the \textit{relative interior} is concerned with open balls in the affine hull of the set, i.e., it considers neighborhoods only within the lower-dimensional affine subspace that actually contains $\sS$.
This is particularly important in higher dimensions and in convex analysis where the set may lie in a lower-dimensional subspace.
\begin{example}[Interior and relative interior points]
\textbf{Interval in $\real$.} Let $ \sS = [0, 1] $. Then $ \interior(\sS) = (0, 1) $ and $ \relint(\sS) = (0, 1) $. Here, the interior and relative interior are the same because the set $ [0, 1] $ is a subset of $\real$ and $\aff(\sS) = \real$.

\paragraph{Line Segment in $\real^2$.} Let $ \sS = \{ (x, y) \in \real^2 \mid 0 \leq x \leq 1, y = 0 \} $ with
$\aff(\sS)=\{(x,0)\in\real^2 \mid x\in\real\}$. Then  $ \interior(\sS) = \varnothing $ and $ \relint(\sS) = \{ (x, 0) \mid 0 < x < 1 \} $.
The interior is empty because there are no open balls in $\real^2$ entirely contained within $ \sS $. However, the relative interior is the open interval (0, 1) in the context of the line segment (the affine hull is the line $ y = 0 $).

\paragraph{Triangle in $\real^2$.} Let $ \sS $ be the (filled) triangle with vertices at $(0,0)$, $(1,0)$, and $(0,1)$. Then the interior consists of all points strictly inside the triangle, excluding the boundary; 
and the  relative interior is the same as the interior in this case because the affine hull of the triangle is the entire plane $ \real^2 $.
\end{example}

\index{Open sets}
\index{Closed sets}
\index{Compact sets}
\index{Level sets}
\begin{definition}[Open, closed, compact, and level sets]\label{definition:open_close_sets}
A set $\sS_1$ is said to be \textit{open} if it consists of only interior points. That is, for every $\balpha\in\sS$, there exists a scalar $s>0$ such that $\sB(\balpha, s)\subseteq\sS_1$.

A set $\sS_2$ is said to be \textit{closed} if its complement $\comple{\sS_2}$ is open.
Alternatively, $\sS_2$ is closed if it contains all the limit points of convergent sequence of points in $\sS_2$; that is, for every sequence of points $\{\balpha_i\}_{i\geq 1}\subseteq \sS_2$ satisfying $\balpha_i\rightarrow \balpha^\star$ as $i\rightarrow \infty$, it follows that $\balpha^\star\in\sS_2$.
In the meantime, the \textit{closure} of a set $\sS$, denoted  $\closure(\sS)$, is defined as the smallest closed set containing $\sS$:
$$
\closure(\sS) = \cap\left\{ \sT\mid  \sS \subseteq \sT, \sT \text{ is closed} \right\}.
$$

A set $\sS \subseteq \real^p$ is called \textit{bounded} if there exists $M > 0$ for which $\sS \subseteq \sB(\bzero, M)$, where $\sB(\bzero, M)$ is the ball of finite radius $M$ centered at the origin. 
A set $\sS \subseteq \real^p$ is called \textit{compact} if it is closed and bounded.

Given a function $ f: \sS \rightarrow \real $ where $ \sS \subseteq \real^p $, the \textit{level set} of $ f $ at a particular value $ c \in \real $ is defined as:
$$
\lev[f, c] = \{ \balpha \in \sS \mid f(\balpha) \leq c \}.~\footnote{Note that we use square brackets instead of parentheses to indicate that the equality can be obtained; consistent with the notation for closed balls in Definition~\ref{definition:open_closed_ball}.}
$$
In other words, the level set $ \lev[f, c] $ consists of all points $ \balpha $ in the domain $ \sS $ for which the function $ f $ evaluates to smaller than or equal to  the constant $ c $.
\end{definition} 

\begin{exercise}[Interior and closure]\label{exercise:int_clos}
Let $\sS\subseteq\real^p$ be a convex set with a nonempty interior. Show that 
$$
\closure(\interior(\sS)) = \closure(\sS)
\qquad\text{and}\qquad 
\interior(\closure(\sS)) =\interior(\sS).
$$
\end{exercise}

\subsection*{Convex Sets}

For vectors $\balpha$ and $\bbeta$,  any point of the form $\lambda\balpha+(1-\lambda)\bbeta$ with $\lambda\in[0,1]$ is called a  \textit{convex combination} of the two vectors.
A set that is closed under all such convex combinations is known as a convex set.
The standard definition follows:
\begin{definition}[Convex sets\index{Convex sets}]\label{definition:convexset}
A set $\sS\subseteq \real^p$ is called \textit{convex} if, for any $\balpha,\bbeta\in\sS$ and $\lambda\in[0,1]$, the point $\lambda\balpha+(1-\lambda)\bbeta$ also belongs to $\sS$.
\end{definition}

By definition, it is straightforward to verify that a set $ \sS \subseteq \real^p $ is convex if and only if it contains every convex combination of finitely many of its points. 
That is, for any $k\in \naturalset$, any $ \balpha_1, \balpha_2, \ldots, \balpha_k \in \sS $ and $ \lambda_1, \lambda_2, \ldots, \lambda_k \geq 0 $ such that $ \sum_{i=1}^k \lambda_i = 1 $, the convex combination $ \sum_{i=1}^k \lambda_i \balpha_i $ is also contained in $ \sS $.
For an arbitrary (possibly non-convex) set $\sS\subseteq\real^p$, its \textit{convex hull} is defined as the collection of all such convex combinations:
\begin{definition}[Convex hulls]\label{definition:cvx_hull}
The \textit{convex hull} of a set $\sS\subseteq\real^p$, denoted  $\conv(\sS)$, is the set of all convex combinations of points in $\sS$: 
$$
\conv(\sS) \triangleq \left\{ \sum_{i=1}^{k} \lambda_i\balpha_i\mid \balpha_1,\balpha_2,\ldots,\balpha_k\in\sS, \blambda\in\Delta_k, k\in\naturalset \right\},
$$
where $\Delta_k = \{(\lambda_1, \lambda_2, \ldots,\lambda_k)\mid \lambda_i\geq 0, \sum_i^k \lambda_i=1\}$ is the standard $(k-1)$ unit-simplex.
\end{definition}

\begin{figure}[h]
\centering       
\vspace{-0.25cm}                 
\subfigtopskip=2pt               
\subfigbottomskip=-2pt         
\subfigcapskip=-10pt      
\includegraphics[width=0.98\textwidth]{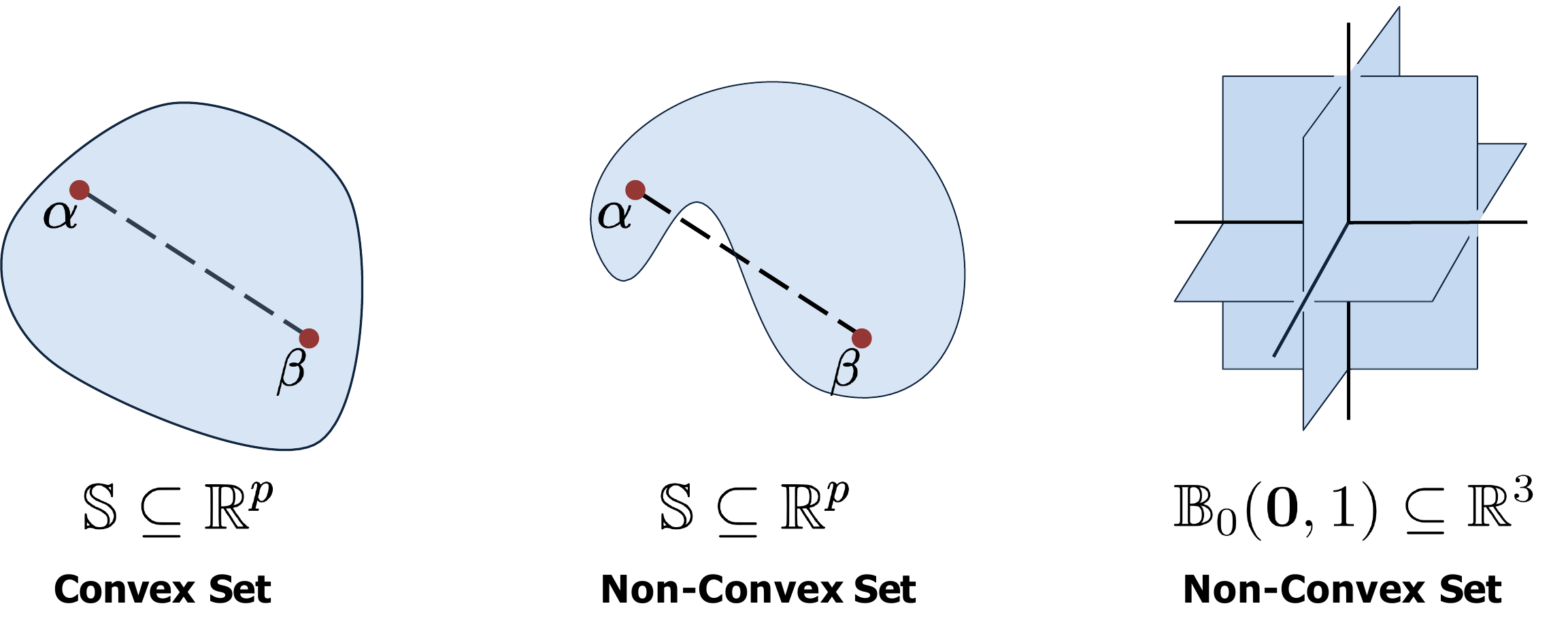}
\caption{
A set is convex if it contains all convex combinations of its points. If even a single convex combination of points in the set lies outside the set, then the set is not convex. Thus, convex sets must be ``bulge-free" and lack inward indentations. Notably, the set of sparse vectors fails this property and is therefore non-convex. The illustration is adapted from \citet{jain2017non}.
}
\label{fig:cvxset}
\end{figure}

The intersection of any collection of convex sets is itself convex. Geometrically, a convex set contains the entire line segment joining any two of its points (see Figure~\ref{fig:cvxset}). Consequently, such sets have no ``inward dents" or concave regions.

\begin{example}
Common examples of convex sets include affine subspaces, linear subspaces,  half spaces,
polygons, and norm balls of the form $\sB_{\norm{\cdot}}[\bc, r]$, where $\bc$ is the center, $r>0$ is the radius, and $\norm{\cdot}$ is a given norm (see Definition~\ref{definition:open_closed_ball}).
\end{example}

\begin{exercise}[Closure of convex]\label{exercise:clo_conv}
Let $\sS\subseteq \real^p$ be a convex set. Show that the closure $\closure(\sS)$ is also convex.
\end{exercise}

\begin{definition}[Cones and convex cones\index{Convex cone}]\label{definition:convex_cone}
A set $ \sS \subseteq \real^p $ is called a \textit{cone} if, for every $ \balpha \in \sS $ and $ \lambda \geq 0 $, the scaled vector $ \lambda \balpha $ also belongs to $\sS $.
A set $ \sS $ is a \textit{convex cone}  if it is both a cone and a convex set.
\end{definition}

The conic hull of a set $\sS\subseteq\real^p$ is the smallest \textbf{convex cone} that contains $\sS$.
\begin{definition}[Conic hulls\index{Conic hulls}]
The \textit{conic hull} of $\sS$, denoted  $\cone(\sS)$, is the set consisting of all the conic combinations of vectors from $\sS$:
$$
\cone(\sS) \triangleq \left\{\sum_{i=1}^{k} \lambda_i\balpha_i\mid \balpha_1,\balpha_2,\ldots,\balpha_k\in\sS, \blambda\in\real^k_+, k\in\naturalset  \right\}.
$$
\end{definition}

Obviously, the zero vector always belongs to any cone. 
A set $ \sS $ is a convex cone if, for all $ \balpha, \bbeta \in \sS $ and $ a, b \geq 0 $, the conic combination $ a \balpha + b \bbeta \in \sS $.

\begin{example}[Convex cone]\label{example:cvx_conv}
Common examples of convex cones include: the \textit{positive orthant} $ \real_+^p = \{ \bbeta \in \real^p \mid  \beta_i \geq 0 \text{ for all } i \in \{1,2,\ldots,p \} \}$ or the set of \textit{positive semidefinite matrices} in $ \real^{p \times p} $. 
Another important example of a convex cone is the \textit{second-order cone} (also known as the \textit{Lorentz cone}):
$
\left\{ \bbeta \in \real^{p+1} \mid \sqrt{\sum_{i=1}^p \beta_i^2} \leq \beta_{p+1} \right\}.
$
\end{example}

\begin{definition}[Dual cone and polar cone\index{Dual cone}\index{Polar cone}]\label{definition:dual_cone}
Let $ \sS \subset \real^p $ be a cone. 
Its \textit{dual cone}j, denoted $ \sS^* $, is defined as
\begin{equation}
\sS^* \triangleq \left\{ \bbeta \in \real^p \mid \innerproduct{\balpha, \bbeta} \geq 0 \text{ for all } \balpha \in \sS \right\}.
\end{equation}
The \textit{polar cone}, denoted $\sS^\circ$, is closely related to the dual cone, which is defined as
\begin{equation}
\sS^\circ \triangleq \left\{ \bbeta \in \real^p \mid\innerproduct{\balpha, \bbeta}\leq 0 \text{ for all } \balpha \in \sS \right\} = -\sS^*.
\end{equation}
\end{definition}

Since the dual cone is defined as the intersection of half spaces, $ \sS^* $ is always closed and convex, and it is straightforward to check that $ \sS^* $ is itself a cone. 
If $ \sS $ is a closed cone, then $ \sS^{**} = \sS $. Moreover, if $ \sT, \sS \subset \real^p $ are cones such that $ \sT \subset \sS $, then $ \sS^* \subset \sT^* $. 

\begin{example}[Dual cone]
The dual cone of the positive orthant $ \real_+^p $ is $ \real_+^p $ itself; in other words, $ \real_+^p $ is \textit{self-dual}. 
\end{example}

We now turn to the notion of extreme points.

\begin{definition}[Extreme points\index{Extreme points}]\label{definition:extrem_points}
Let $ \sS \subset \real^p $ be a convex set. A point $ \bz \in \sS $ is called an \textit{extreme point} of $ \sS $ if, 
whenever $ \bz = \lambda \balpha + (1 - \lambda) \bbeta $ for $ \balpha, \bbeta \in \sS $ and $ \lambda \in (0, 1) $,
it follows that $ \bz = \balpha = \bbeta $.
\end{definition}

Compact (i.e., closed and bounded) convex sets can be fully characterized by their extreme points, as formalized below.
\begin{theoremHigh}[Compact convex sets via extreme points]\label{theorem:compconv_set_extmpoint}
Every compact (i.e., closed and bounded) convex set is equal to the convex hull of its extreme points.
\end{theoremHigh}
\begin{proof}
See \citet{rockafellar2015convex}, Section~18.
\end{proof}

For instance, if $ \sS $ is a polygon in $\real^2$, its extreme points are precisely its vertices (the zero-dimensional faces), and the theorem aligns with the intuitive idea that the polygon is ``built" from its corners.

\section{Convex Functions}

A closely related concept is that of convex functions, which exhibit a specific inequality behavior under convex combinations.
Before introducing convex functions, we note that we work with \textit{extended real-valued} functions $ f : \real^p \to \widetilde{\real}\triangleq (-\infty, \infty] = \real \cup \{+\infty\} $~\footnote{While the extended real line is often defined as $\widehat{\real}\triangleq \real\cup \{-\infty, \infty\}$, we adopt the convention $\widetilde{\real}=(-\infty, \infty]$ here to simplify the treatment of proper functions.}. 
Arithmetic operations and inequalities in $ (-\infty, \infty] $ are interpreted in the natural way; for instance, $ \alpha + \infty = \infty $ for all $ \alpha \in \real $, $ \lambda \cdot \infty = \infty $ for $ \lambda > 0 $, $ \alpha < \infty $ for all $ \alpha \in \real $, $ \infty \leq \infty $. The \textit{domain} of an extended real-valued function $ f $ is defined as
$$
\domain(f) = \{ \bbeta \in \real^p\mid  f(\bbeta) \neq \infty \}.
$$
A function with $ \domain(f) \neq \varnothing $ is called \textit{proper}: that is, (i). The function never takes on the value $-\infty$; (ii). The function does not attain the value $+\infty$ everywhere. 
A function $ f : \sS \to \real $ defined on a subset $ \sS \subset \real^p $ can be extended to an extended valued function by setting $ f(\bbeta) = \infty $ for $ \bbeta \notin \sS $. 
This extension satisfies $ \domain(f) = \sS $ and is referred to as the \textit{canonical extension} of $f$.

We now recall the definition of a convex function.
\begin{definition}[Convex functions\index{Convex functions}]\label{definition:convexfuncs}
A function $f:\sS\rightarrow \real $~\footnote{
or $f:\textcolor{mylightbluetext}{\real^p}\rightarrow \real \textcolor{mylightbluetext}{\cup \{+\infty\}}$ with $\domain(f)\equiv \sS$ is convex.
} defined over a convex set $\sS\subseteq \real^p $ is called \textit{convex} if 
$$
f(\lambda\balpha+(1-\lambda)\bbeta)
\leq \lambda f(\balpha) + (1-\lambda) f(\bbeta), \text{ for any }\balpha,\bbeta\in\sS, \lambda\in[0,1].
$$
Moreover,   $f$ is called \textit{strictly convex} if 
$$
f(\lambda\balpha+(1-\lambda)\bbeta)
< \lambda f(\balpha) + (1-\lambda) f(\bbeta),\text{ for any }\balpha \neq \bbeta \in\sS, \lambda\in(0,1).
$$
\end{definition}

A function $f$ is \textit{concave} (resp., \textit{strictly concave}) if $-f$ is convex (resp., strictly convex).
Note that the domain of any convex function is necessarily a  convex set, and a function $ f : \sS \to \real $ on a convex subset $ \sS \subset \real^p $ is called convex if its canonical extension to $ \real^p $ is convex.

An equivalent and geometrically insightful characterization is that a function $ f $ is convex if and only if its \textit{epigraph}
\begin{equation}\label{equation:def_epigraph}
\epi(f) \triangleq \{ (\balpha, r) \in \sS\times \real \mid f(\balpha) \leq r \} 
\subseteq \real^p \times \real
\end{equation}
is a convex set of $\real^{p+1}$.

\begin{example}[Convex, strictly convex]
Let $f(\balpha) = \balpha^\top \bX \balpha + \by^\top \balpha + z$.
Then the function is 
\begin{itemize}
\item Convex if and only if $\bX \succeq \bzero$; strictly convex if and only if $\bX \succ \bzero$.
\item Concave if and only if $\bX \preceq \bzero$; strictly concave if and only if $\bX \prec \bzero$.
\end{itemize}
The proofs are easy if we use the second-order characterization of convexity (see Theorem~\ref{theorem:psd_hess_conv}).
\end{example}
\begin{exercise}[Convexity of directional derivatives]
Suppose the  function $f$ is a continuously differentiable convex function. 
Show that the directional derivative (Definition~\ref{definition:partial_deri}) $g(\bd)\triangleq f^\prime(\balpha;\bd)$ is also convex.
\end{exercise}

\begin{exercise}[Convexity of quadratic functions]\label{exercise:conv_quad}
Let $f(\balpha)=\frac{1}{2}\balpha^\top\bX\balpha+\by^\top\balpha+z$, where $\bX\in\real^{p\times p}$ is symmetric, $\by\in\real^p$, and $z\in\real$. Show that $f(\balpha)$ is convex (resp., strictly convex) if and only if $\bX\succeq \bzero$ (resp., $\bX\succ\bzero$).
\end{exercise}

\begin{exercise}[Convexity of indicator]\label{exercise_convex_indica}
Show that the \textit{indicator function} (also called the \textit{characterization function}) $\indicatorS (\balpha):\real^p\rightarrow \real$, which is $0$ if $\balpha$ belongs to the set $\sS$ and $+\infty$ otherwise,  is convex if and only if the set $\sS$ is  convex.
\end{exercise}

\begin{exercise}[Convexity of norm functions]
Show that the norm functions $\normtwo{\balpha}^2$  and $\normone{\balpha}$ are convex.
More generally, show that
\begin{itemize}
\item The $ \ell_s $-norm functions $ \norms{\cdot} $ are strictly convex if $ 1 < s < \infty $, and they are \textbf{not} strictly convex if $ s = 1 $ or $ s = \infty $.

\item For a nondecreasing convex function $ f : \real \to (-\infty, \infty] $ and a norm $ \norm{\cdot} $ on $ \real^p $, the function $ g(\balpha) = f(\norm{\balpha}) $ is convex. 
In particular, the function $ \balpha \mapsto \norm{\balpha}^k $ is convex provided that $ k \geq 1 $.
\end{itemize}

\end{exercise}

The exercises above show that standard norm functions are convex. However, the so-called ``$\ell_0$-norm"---defined as the number of nonzero entries in a vector---is not convex. 
Despite its name, it is not a true norm because it fails to satisfy positive homogeneity: for any nonzero $\balpha\in\real^p$, we have  $\normzero{2\balpha} = \normzero{\balpha}\neq 2\normzero{\balpha}$; see Definition~\ref{definition:matrix-norm}.

\begin{lemma}[$\ell_0$-norm is not convex]\label{lemma:ell0_noncvx}
The $\ell_0$-norm function, defined as the number of nonzero entries in a vector, is not convex.
\end{lemma}

\begin{proof}[of Lemma~\ref{lemma:ell0_noncvx}]
Consider a simple counterexample in $\real^2$. 
Let $\balpha \triangleq \begin{bmatrixfoot} 1 \\ 0 \end{bmatrixfoot}$ and $\bbeta \triangleq \begin{bmatrixfoot} 0 \\ 1 \end{bmatrixfoot}$. 
Then, for any $\lambda \in (0,1)$,
$$
\normzero{\lambda \balpha + (1 - \lambda) \bbeta} = 2 > 1 = \lambda \normzero{\balpha} + (1 - \lambda)\normzero{\bbeta}.
$$
This violates the definition of convexity. The same idea extends directly to $\real^p$, completing the proof.
\end{proof}

The rank of a matrix, when viewed as a function of its entries, is also not convex.
\begin{lemma}[The rank is not convex]\label{lemma:rk_noncvx}
The rank function $\rank(\cdot): \real^{p\times p}\rightarrow \real$ is not convex on the space of square matrices.
\end{lemma}
\begin{proof}[of Lemma~\ref{lemma:rk_noncvx}]
We provide a counterexample analogous to the one used in the proof of Lemma~\ref{lemma:ell0_noncvx}. 
Let
$$
\bX \triangleq \begin{bmatrix} 1 & 0 \\ 0 & 0 \end{bmatrix}, \qquad \bY \triangleq \begin{bmatrix} 0 & 0 \\ 0 & 1 \end{bmatrix}.
$$
For any $\lambda \in (0,1)$,
$$
\rank\left(\lambda \bX + (1 - \lambda) \bY\right) = 2 > 1 = \lambda\,\rank(\bX) + (1 - \lambda)\,\rank(\bY).
$$
Extending this counterexample into $\real^{p\times p}$ completes the proof.
\end{proof}

\subsection{Properties of Convex Functions}

A well-known inequality that follows from the definition of convex functions is stated below without proof.
\begin{theoremHigh}[Jensen's inequality\index{Jensen's inequality}]\label{theorem:jensens_ineq}
Let $f: \sS \rightarrow \real$ be a convex function defined on a convex subset $\sS \subseteq \real^p$. For any finite sequence of points $\balpha_{1}, \balpha_{2}, \ldots, \balpha_{k} \in \sS$ and any sequence of nonnegative weights $\lambda_{1}, \lambda_{2}, \lambda_{3}, \ldots, \lambda_{k}$ such that $\sum_{i=1}^{k} \lambda_{i} = 1$, Jensen's inequality states:
$$
f\left(\sum_{i=1}^{k} \lambda_{i} \balpha_{i}\right) \leq \sum_{i=1}^{k} \lambda_{i} f(\balpha_{i}).
$$
If $f$ is concave, the inequality is reversed:
$$
f\left(\sum_{i=1}^{k} \lambda_{i} \balpha_{i}\right) \geq \sum_{i=1}^{k} \lambda_{i} f(\balpha_{i}).
$$
In the context of probability theory, if $\rvx$ is a random vector taking values in $\sS$ and $f$ is a convex function, 
Jensen's inequality takes the form:
$$
f(\Exp[\rvx]) \leq \Exp[f(\rvx)],
$$
where $\Exp[\cdot]$ denotes the expectation operator over the random vector $\rvx$. For a concave function, the inequality is again reversed.
\end{theoremHigh}

Next, we consider the maximum of a convex function over a convex set.

\begin{theoremHigh}[Attainment of maximum]\label{theorem:attain_max_cvx}
Let $ \sS \subset \real^p $ be a compact  (closed and bounded) convex set, and $ f : \sS \to \real $ be a convex function. Then $ f $ attains its maximum at an extreme point of $ \sS $.
\end{theoremHigh}
\begin{proof}[of Theorem~\ref{theorem:attain_max_cvx}]
Let $ \balpha \in \sS $ such that $ f(\balpha) \geq f(\bbeta) $ for all $ \bbeta \in \sS $. Since $ \sS $ is the convex hull of its extreme points (Theorem~\ref{theorem:compconv_set_extmpoint}), we can write $ \balpha = \sum_{i=1}^k \lambda_i \balpha_i $ for some $ k $, 
where $ \lambda_i > 0 $, $ \sum_{i=1}^k \lambda_i = 1 $, and each  $ \balpha_i $ is an extreme point of $ \sS $. 
By convexity and the assumption of $\balpha$, we have 
$$
f(\balpha) = f\left( \sum_{i=1}^k \lambda_i \balpha_i \right) \leq \sum_{i=1}^k \lambda_i f(\balpha_i) \leq \sum_{i=1}^k \lambda_i f(\balpha) = f(\balpha).
$$
Hence, both inequalities are equalities, which implies $ f(\balpha_i) = f(\balpha) $ for all $ i = 1, 2, \ldots, k $. Therefore, the maximum of $ f $ is attained at an extreme point of $ \sS $.
\end{proof}

We say that a function $ f(\balpha, \bbeta) $ of two arguments $ \balpha \in \real^q $, $ \bbeta \in \real^p $ is jointly convex if it is convex as a function of the variable $ \bz = (\balpha, \bbeta) $. 
Partial minimization of a jointly convex function with respect to one block of variables preserves convexity, as stated below.
\begin{theoremHigh}[Preservation under partial infimum]\label{theorem:preser_partinf}
Let $ f : \real^q \times \real^p \to (-\infty, \infty] $ be a jointly convex function. Then the function $ g(\balpha) \triangleq \min_{\bbeta \in \real^p} f(\balpha, \bbeta) $, $ \balpha \in \real^q $, is convex.
\end{theoremHigh}
\begin{proof}[of Theorem~\ref{theorem:preser_partinf}]
For simplicity, assume that the infimum is attained for every $\balpha$.  
The general case has to be treated with an $ \varepsilon $-argument.
Given $ \balpha_1, \balpha_2 \in \real^q $, there exist $ \bbeta_1, \bbeta_2 \in \real^p $ such that
$$
f(\balpha_1, \bbeta_1) = \min_{\bbeta \in \real^p} f(\balpha_1, \bbeta), \quad f(\balpha_2, \bbeta_2) = \min_{\bbeta \in \real^p} f(\balpha_2, \bbeta).
$$
For any $ \lambda \in [0, 1] $, the minimization of $g$ and joint convexity of $f$ implies that
$$
\begin{aligned}
g(\lambda \balpha_1 + (1 - \lambda) \balpha_2) 
&\leq f(\lambda \balpha_1 + (1 - \lambda) \balpha_2, \lambda \bbeta_1 + (1 - \lambda) \bbeta_2)\\
&\leq \lambda f(\balpha_1, \bbeta_1) + (1 - \lambda) f(\balpha_2, \bbeta_2) = \lambda g(\balpha_1) + (1 - \lambda) g(\balpha_2).
\end{aligned}
$$
Thus, $g$ is convex, completing the proof.
\end{proof}

Clearly, the same reasoning shows that partial maximization of a jointly concave function yields a concave function.
The following operations also preserve convexity. These results are left as an exercise.
\begin{exercise}[Preservation of convexity]\label{exercise:pres_conv_clos}
Prove the following results:
\begin{enumerate}[(i)]
\item Let $\bX\in\real^{p\times p}$ and $\by\in\real^p$, and let $f: \real^p \rightarrow (-\infty, \infty]$ be an extended real-valued convex function. Then the function 
$
g(\balpha) = f(\bX\balpha+\by)
$
is convex.

\item Let $ f : \real \to \real $ be convex and nondecreasing, and let $ g : \real^p \to \real $ be convex. Then the function $ h(\balpha) = f(g(\balpha)) $ is convex.

\item Let $f_1, f_2, \ldots, f_k: \real^p \rightarrow (-\infty, \infty]$ be extended real-valued convex functions, and let $\sigma_1, \sigma_2, \ldots, \sigma_k \in \real_+$. Then the function $f = \sum_{i=1}^{k} \sigma_i f_i$ is convex.

\item Let $f_i: \real^p \rightarrow (-\infty, \infty], i \in \sI$ be extended real-valued convex functions, where $\sI$ is a given index set. Then the function
$
f(\balpha) = \max_{i \in \sI} f_i(\balpha)
$
is convex.
\end{enumerate}
\textit{Hint: the epigraph of $\max_{i \in \sI} f_i(\balpha)$ is the intersection of the epigraphs of functions $f_i$.}
\end{exercise}

\begin{theoremHigh}[Local is global in convex functions]\label{theorem:local_glo_cvx}
Let $ f : \sS\subseteq\real^p \to \real $ (or $f : \real^p \to \real \cup \{\infty\}$) be convex over a convex set $\sS$.
\begin{enumerate}[(i)]
\item A local minimum of $ f $ is a global minimum. If $ f $ is strictly convex, then the minimum is unique.
\item The set of minimizers  of $ f $ is convex. If $f$ is strictly convex, then the set is a \textbf{singleton}.
\end{enumerate}
\end{theoremHigh}
\begin{proof}[of Theorem~\ref{theorem:local_glo_cvx}]
\textbf{(i).} Suppose  $\widehatbtheta$ is a local minimum of $f$ over $\sS$. There exists a scalar $\tau > 0$ such that $f(\btheta) \geq f(\widehatbtheta)$ for any $\btheta \in \sS$ satisfying $\btheta \in \sB[\widehatbtheta, \tau]$. Now let $\bgamma \in \sS$ satisfy $\bgamma \neq \widehatbtheta$. It suffices to show that $f(\bgamma) \geq f(\widehatbtheta)$. Let $\lambda \in (0, 1]$ be such that $\widehatbtheta + \lambda(\bgamma - \widehatbtheta) \in \sB[\widehatbtheta, \tau]$. An example of such $\lambda$ is $\lambda = \frac{\tau}{\normtwo{\widehatbtheta - \bgamma}}$. Since $\widehatbtheta + \lambda(\bgamma - \widehatbtheta) \in \sB[\widehatbtheta, \tau] \cap \sS$, it follows that
$
f(\widehatbtheta) \leq f(\widehatbtheta + \lambda(\bgamma - \widehatbtheta)),
$
and hence by Jensen's inequality (Theorem~\ref{theorem:jensens_ineq})
$
f(\widehatbtheta) \leq f(\widehatbtheta + \lambda(\bgamma - \widehatbtheta)) \leq (1 - \lambda) f(\widehatbtheta) + \lambda f(\bgamma).
$
Therefore, we obtain $f(\widehatbtheta) \leq f(\bgamma)$.

A slight modification of the above argument  shows that any local minimum of a strictly convex function over a convex set is indeed  a strict global minimum of the function over the set.

\paragraph{(ii).}
Let $\Theta\triangleq \argmin\{f(\btheta) : \btheta \in \sS\}$.	
If $\Theta = \varnothing$, the result follows trivially. 
We then assume that $\Theta \neq \varnothing$ and denote the optimal value by $f^*$. 
Let $\btheta, \bgamma \in \Theta$ be two different minimizers and $\lambda \in [0, 1]$. Then, by Jensen's inequality $f(\lambda \btheta + (1 - \lambda) \bgamma) \leq \lambda f^* + (1 - \lambda) f^* = f^*$, and hence $\lambda \btheta + (1 - \lambda) \bgamma$ is also a minimmizer, i.e., it belongs to $\Theta$, establishing the convexity of $\Theta$. 
Suppose now that $f$ is strictly convex and $\Theta$ is nonempty; to show that $\Theta$ is a singleton, suppose in contradiction that there exist $\btheta, \bgamma \in \Theta$ such that $\btheta \neq \bgamma$. Then $\frac{1}{2}\btheta + \frac{1}{2}\bgamma \in \sS$, and by the strict convexity of $f$, we have
$$
f\left(\frac{1}{2}\btheta + \frac{1}{2}\bgamma\right) < \frac{1}{2}f(\btheta) + \frac{1}{2}f(\bgamma) = \frac{1}{2}f^* + \frac{1}{2}f^* = f^*,
$$
which leads to a contradiction to the fact that $f^*$ is the optimal value.
Therefore, $\Theta$ must be a singleton.
\end{proof}

The fact that every local minimum of a convex function is automatically a global minimum is the key reason why convex optimization problems admit efficient solution methods; see Section~\ref{section:opt_conds} for further details.

\subsection{Characterization of Convex Functions}
\paragrapharrow{First-order characterizations of convex functions.}
Convex functions need not be differentiable. However, when they are differentiable, they can be characterized by the \textit{gradient inequality}.
\begin{theoremHigh}[Gradient inequality]\label{theorem:conv_gradient_ineq}
Let $f:\sS\rightarrow \real$ be a {continuously differentiable} function defined on a convex set $\sS\subseteq\real^p$. Then, $f$ is convex over $\sS$ if and only if 
\begin{subequations}
\begin{equation}\label{equation:conv_gradient_ineq1}
f(\balpha) +\nabla f(\balpha)^\top (\bbeta-\balpha)\leq f(\bbeta), \text{ for any $\balpha,\bbeta\in\sS$}.
\end{equation}
Similarly, the function is strictly convex on $\sS$ if and only if 
\begin{equation}\label{equation:conv_gradient_ineq2}
f(\balpha) +\nabla f(\balpha)^\top (\bbeta-\balpha)< f(\bbeta), \text{ for any $\balpha\neq \bbeta\in\sS$}.
\end{equation}
\end{subequations}
\end{theoremHigh}
Geometrically, this means that the graph of a convex function lies above its tangent hyperplane at any point. For concave (or strictly concave) functions, the inequality signs are reversed.
\begin{proof}[of Theorem~\ref{theorem:conv_gradient_ineq}]
For brevity, we  prove only \eqref{equation:conv_gradient_ineq1}; the proof of~\eqref{equation:conv_gradient_ineq2} is analogous.
Suppose  that $ f $ is convex. Let $ \balpha, \bbeta \in \sS $ and $ \lambda \in [0, 1] $. If $ \balpha = \bbeta $, then \eqref{equation:conv_gradient_ineq1} trivially holds. We will therefore assume that $ \balpha \neq \bbeta $. 
By convexity
$$
f(\lambda \bbeta + (1-\lambda) \balpha) \leq \lambda f(\bbeta) + (1-\lambda) f(\balpha)
\;\implies\;
\frac{f(\balpha + \lambda (\bbeta - \balpha)) - f(\balpha)}{\lambda} \leq f(\bbeta) - f(\balpha).
$$
Taking $ \lambda \to 0^+ $, the left-hand side converges to the directional derivative of $ f $ at $ \balpha $ in the direction $ \bbeta - \balpha $ (Definition~\ref{definition:partial_deri}), whence we have 
$
f'(\balpha; \bbeta - \balpha) \leq f(\bbeta) - f(\balpha).
$
Since $ f $ is continuously differentiable, it follows that $ f'(\balpha; \bbeta - \balpha) = \nabla f(\balpha)^\top (\bbeta - \balpha) $ by \eqref{equation:direc_contdiff}, and hence \eqref{equation:conv_gradient_ineq1} follows.

Conversely, assume that the gradient inequality \eqref{equation:conv_gradient_ineq1} holds. 
Let $ \balpha, \bbeta \in \sS $, and let $ \lambda \in [0, 1] $. We will show that $ f(\lambda \balpha + (1 - \lambda) \bbeta) \leq \lambda f(\balpha) + (1 - \lambda) f(\bbeta) $. The cases for $\lambda=0$ or $\lambda=1$ holds trivially. We will therefore assume that $\lambda\in(0,1)$.
Let $ \bz \triangleq \lambda \balpha + (1 - \lambda) \bbeta \in \sS $. Then
$
-\frac{\lambda}{1 - \lambda}  (\balpha - \bz)  = (\bbeta - \bz).
$
Combining with the gradient inequality implies that 
$$
f(\bz) + \nabla f(\bz)^\top (\balpha - \bz) \leq f(\balpha)
\qquad\text{and}\qquad 
f(\bz) - \frac{\lambda}{1 - \lambda} \nabla f(\bz)^\top (\balpha - \bz) \leq f(\bbeta).
$$
Multiplying the first inequality by $ \frac{\lambda}{1 - \lambda} $ and adding it to the second one, we obtain
$
\frac{1}{1 - \lambda} f(\bz) \leq \frac{\lambda}{1 - \lambda} f(\balpha) + f(\bbeta),
$
establishing convexity of $f$.
This completes the proof.
\end{proof}

\begin{theoremHigh}[Monotonicity of gradient of convex functions]\label{theorem:monoton_convgrad}
Let $f:\sS\rightarrow \real$ be a {continuously differentiable} function defined on  a convex set $\sS\subseteq\real^p$. Then  $f$ is convex over $\sS$ if and only if 
\begin{subequations}
\begin{equation}\label{equation:monoton_convgrad}
\innerproductbig{\nabla f(\balpha) - \nabla f(\bbeta), (\balpha-\bbeta)}\geq 0, 
\quad \text{for all } \balpha, \bbeta \in \sS.
\end{equation}
Similarly,  $f$ is strictly convex over $\sS$ if and only if 
\begin{equation}\label{equation:monoton_convgrad2}
\innerproductbig{\nabla f(\balpha) - \nabla f(\bbeta), (\balpha-\bbeta)}> 0,
\quad \text{for all } \balpha \neq \bbeta \in \sS.
\end{equation}
\end{subequations}
\end{theoremHigh}
\begin{proof}[of Theorem~\ref{theorem:monoton_convgrad}]
We prove only~\eqref{equation:monoton_convgrad}; the strict case follows similarly.
Assume  that $ f $ is convex over $ \sS $. Then by the gradient inequality, for any $ \balpha, \bbeta \in \sS $,  we have 
$$
f(\balpha) \geq f(\bbeta) + \nabla f(\bbeta)^\top (\balpha - \bbeta)
\qquad\text{and}\qquad 
f(\bbeta) \geq f(\balpha) + \nabla f(\balpha)^\top (\bbeta - \balpha).
$$
Adding these two inequalities eliminates the function values and yields  \eqref{equation:monoton_convgrad}.

Conversely, assume that \eqref{equation:monoton_convgrad} holds, and let $ \balpha, \bbeta \in \sS $. Define $g(\mu) \triangleq f(\balpha + \mu (\bbeta - \balpha)),  \mu \in [0, 1]$ as a one-dimensional function.
By the {fundamental theorem of calculus (Theorem~\ref{theorem:fund_theo_calculu})}, we have
$$
\begin{aligned}
f(\bbeta) 
&= g(1) = g(0) + \int_0^1 g'(\mu) \, d\mu
= f(\balpha) + \int_0^1 (\bbeta - \balpha)^\top \nabla f\big(\balpha + \mu (\bbeta - \balpha)\big) \, d\mu\\
&= f(\balpha) + \nabla f(\balpha)^\top (\bbeta - \balpha) + \int_0^1 (\bbeta - \balpha)^\top \big(\nabla f(\balpha + \mu (\bbeta - \balpha)) - \nabla f(\balpha)\big) \, d\mu\\
&\geq f(\balpha) + \nabla f(\balpha)^\top (\bbeta - \balpha),
\end{aligned}
$$
where the last inequality follows from the monotonicity of the gradient. This shows that $f$ is convex by Definition~\ref{definition:convexfuncs}.
\end{proof}

\paragrapharrow{Second-order characterizations of convex functions.}

When the function is further assumed to be twice continuously differentiable,
its convexity can be characterized by the positive semidefiniteness (see Definition~\ref{definition:psd-pd-defini}) of its Hessian matrix
\begin{theoremHigh}[PSD Hessian of convex functions]\label{theorem:psd_hess_conv}
Let $f:\sS\subseteq\real^p\rightarrow \real$ be a twice continuously differentiable function defined on an \textbf{open} convex set $\sS$. Then $f$ is convex \textbf{if and only if} $\nabla^2 f(\balpha)\succeq \bzero $ for any $\balpha\in\sS$.~\footnote{Note that the openness of $\sS$ can be relaxed to prove the convexity of $\nabla^2 f(\balpha)\succeq \bzero $ for any $\balpha\in\interior(\sS)$, where 
\begin{equation}
\interior(\sS) = \{ \balpha \in \sS \mid \exists\, \epsilon > 0 \text{ such that } \sB(\balpha, \epsilon) \subseteq \sS \},
\end{equation}
and $\sB(\balpha,\epsilon)$ denotes the open ball centered at $\balpha$ with radius $\epsilon$ (Definition~\ref{definition:open_closed_ball}).
}
Moreover, \textbf{if} $\nabla^2 f(\balpha)\succ \bzero $ for any $\balpha\in\sS$, then the function is strictly convex over $\sS$~\footnote{Notice that the former condition is both sufficient and necessary, the latter condition is merely sufficient but not necessary. For example, $f(x)=x^6$ is strictly convex, but $f^{\prime\prime}(x)=30x^4$ is equal to zero at $x=0$.}. 
\end{theoremHigh}
\begin{proof}[of Theorem~\ref{theorem:psd_hess_conv}]
Assume that $\nabla^2 f(\balpha) \succeq \bzero$ for all $\balpha \in \sS$, and let $\balpha, \bbeta \in \sS$. 
By the linear approximation theorem (Theorem~\ref{theorem:linear_approx}), there exists $\bxi \in [\balpha, \bbeta]$ (and hence $\bxi \in \sS$) such that
\begin{equation}\label{equation:psd_hess_conv1}
f(\bbeta) = f(\balpha) + \nabla f(\balpha)^\top (\bbeta - \balpha) + \frac{1}{2} (\bbeta - \balpha)^\top \nabla^2 f(\bxi) (\bbeta - \balpha), \quad \bxi \in [\balpha, \bbeta].
\end{equation}
Given that  $\nabla^2 f(\bxi) \succeq \bzero$, it follows that $(\bbeta - \balpha)^\top \nabla^2 f(\bxi) (\bbeta - \balpha) \geq 0$, whence by \eqref{equation:psd_hess_conv1} we have  $f(\bbeta) \geq f(\balpha) + \nabla f(\balpha)^\top (\bbeta - \balpha)$, establishing the convexity by Theorem~\ref{theorem:conv_gradient_ineq}.

Conversely, assume that $f$ is convex over $\sS$. Let $\balpha \in \sS$ and let $\bbeta \in \real^p$. Since $\sS$ is open, there exists a scalar $\varepsilon>0$ such that  $\balpha + \lambda \bbeta \in \sS$ for any  $\lambda \in(0, \varepsilon)$. 
The gradient inequality (Theorem~\ref{theorem:conv_gradient_ineq}) shows that 
\begin{equation}\label{equation:psd_hess_conv2}
f(\balpha + \lambda \bbeta) \geq f(\balpha) + \lambda \nabla f(\balpha)^\top \bbeta.
\end{equation}
Additionally, by the quadratic approximation theorem (Theorem~\ref{theorem:quad_app_theo}), we have
$$
f(\balpha + \lambda \bbeta) = f(\balpha) + \lambda \nabla f(\balpha)^\top \bbeta + \frac{\lambda^2}{2} \bbeta^\top \nabla^2 f(\balpha) \bbeta + o(\lambda^2 \normtwo{\bbeta}^2),
$$
which combined with \eqref{equation:psd_hess_conv2} yields the inequality
$
\frac{\lambda^2}{2} \bbeta^\top \nabla^2 f(\balpha) \bbeta + o(\lambda^2 \normtwo{\bbeta}^2) \geq 0
$
for any $\lambda \in (0, \varepsilon)$. Dividing the latter inequality by $\lambda^2$ yields that 
$
\frac{1}{2} \bbeta^\top \nabla^2 f(\balpha) \bbeta + \frac{o(\lambda^2 \normtwo{\bbeta}^2)}{\lambda^2} \geq 0.
$
Taking $\lambda \to 0^+$, we conclude that
$
\bbeta^\top \nabla^2 f(\balpha) \bbeta \geq 0
$
for any $\bbeta \in \real^p$, implying that $\nabla^2 f(\balpha) \succeq \bzero$ for any $\balpha \in \sS$.
The strict convexity can be proved analogously.
\end{proof}

\section{Continuous and Closed Functions, and Weierstrass Theorem}

Since we will use the notion of closed functions to analyze  properties of proximal operators in Section~\ref{section:prop_proxoper}---a step that relies on the closedness condition in the Weierstrass theorem (Theorem~\ref{theorem:weierstrass_them})---and since this theorem is, in turn, essential for studying the convergence of \textit{proximal gradient methods} in Section~\ref{section:proxiGD_inClasso}, we review the relevant concepts in this section.

A continuous function is one that does not exhibit abrupt changes in value---i.e., it has no \textit{discontinuities}. More formally, continuity means that arbitrarily small changes in the input produce arbitrarily small changes in the output. Continuous functions are fundamental in calculus and mathematical analysis and enjoy several key properties.

\begin{definition}[Continuous functions]\label{definition:conti_funs}
Let $f:\sS \subseteq \real^p\rightarrow \real$. The function $f$ is said to be \textit{continuous} at a point $\bbeta \in \sS$ if for every $\epsilon > 0$, there exists a $\delta > 0$ such that for all $\balpha$ in $\sS$:
\begin{equation}
\forall\, \epsilon > 0, \exists\, \delta > 0 :\quad \forall\, \balpha \in \sS,
\normtwo{\balpha - \bbeta} < \delta 
\quad\implies \quad
\abs{f(\balpha) - f(\bbeta)} < \epsilon.
\end{equation}
In simpler terms, this means that the values $f(\balpha)$ can be made arbitrarily close to $f(\bbeta)$ by choosing $\balpha$ sufficiently close to $\bbeta$.

If $f$ is continuous at every point in its domain $\sS$, then $f$ is said to be \textit{continuous} on $\sS$, or simply continuous.
\end{definition}

\begin{remark}[Properties of continuous functions]
Continuous functions satisfy several important properties:
\begin{enumerate}
\item \textit{Intermediate value theorem.} If $f$ is continuous on a closed interval $[a, b]$ and $y$ lies between $f(a)$ and $f(b)$, then there exists some $c \in [a, b]$ such that $f(c) = y$. This is known as Rolle's theorem (Theorem~\ref{theorem:rolles_theo}).
\item \textit{Extreme value theorem.} If $f$ is continuous on a closed interval $[a, b]$, then $f$ attains both a maximum and a minimum value on that interval. This will be discussed further in the Weierstrass theorem (Theorem~\ref{theorem:weierstrass_them}).
\item \textit{Composition of continuous functions.} If $g$ is continuous at $c$ and $f$ is continuous at $g(c)$, then the composition $f \circ g$ is continuous at $c$.
\item \textit{Arithmetic operations.} If $f$ and $g$ are continuous at $c$, then so are $f + g$, $f - g$, $f g$, and $\frac{f}{g}$ (provided $g(c) \neq 0$).
\item \textit{Continuity and limits.} A function $f$ is continuous at $c$ if and only if $\lim_{x \to c} f(x) = f(c)$.
\end{enumerate}
\end{remark}

\index{Continuous}
\begin{example}[Continuous functions]
The following are examples of continuous functions:
\begin{itemize}
\item Polynomial functions are continuous everywhere.
\item Rational functions are continuous wherever they are defined (i.e., where the denominator is not zero).
\item Trigonometric functions like sine and cosine are continuous everywhere.
\item Exponential and logarithmic functions are continuous on their respective  domains.
\item Absolute value function $\abs{\alpha}$ is continuous everywhere.
\end{itemize}
In contrast, common types of discontinuities include:
\begin{itemize}
\item \textit{Removable discontinuity.} The function has a hole at a point, but it can be ``filled" to become continuous.
\item \textit{Jump discontinuity.} The left-hand limit and right-hand limit exist but are not equal.
\item \textit{Infinite discontinuity.} The function approaches  positive or negative infinity at a point.
\item \textit{Oscillating discontinuity.} The function does not approach a single value as the input approaches a point.
\end{itemize}
Continuous functions are generally easier to analyze and manipulate, while discontinuities can lead to complex behaviors that require special handling. Identifying the type of discontinuity is crucial for selecting appropriate analytical or numerical strategies in applications.
\end{example}

\index{Lower semicontinuity}
\index{Lower semicontinuous}
\begin{definition}[Lower semicontinuity (LSC)]
A function $ f: \real^p \rightarrow \real\cup \{\infty\} $ is called \textit{lower semicontinuous (LSC) at $ \balpha^* \in \real^p $} if
$$
f(\balpha^*) \leq \liminf_{t \rightarrow \infty} f(\balpha^\toptzero)
$$
for any sequence $ \{\balpha^\toptzero\}_{t \geq 1} \subseteq \real^p $ such that  $ \balpha^\toptzero \rightarrow \balpha^* $ as $ t \rightarrow \infty $.
In other words, the limit inferior of the function values at points approaching $\balpha^*$ must be greater than or equal to the function value at $\balpha^*$. 
This means the function cannot exhibit a sudden drop at $\balpha^*$.
A function $ f: \real^p \rightarrow  \real \cup \{\infty \} $ is called \textit{lower semicontinuous} if it is {lower semicontinuous} at each point in $ \real^p $.
\end{definition}

\begin{figure}[h!]
\centering                      
\vspace{-0.35cm}                
\subfigtopskip=2pt            
\subfigbottomskip=2pt           
\subfigcapskip=-5pt           
\subfigure[A lower semicontinuous funtion.]{\label{fig:aaa-1}
\includegraphics[width=0.31\linewidth]{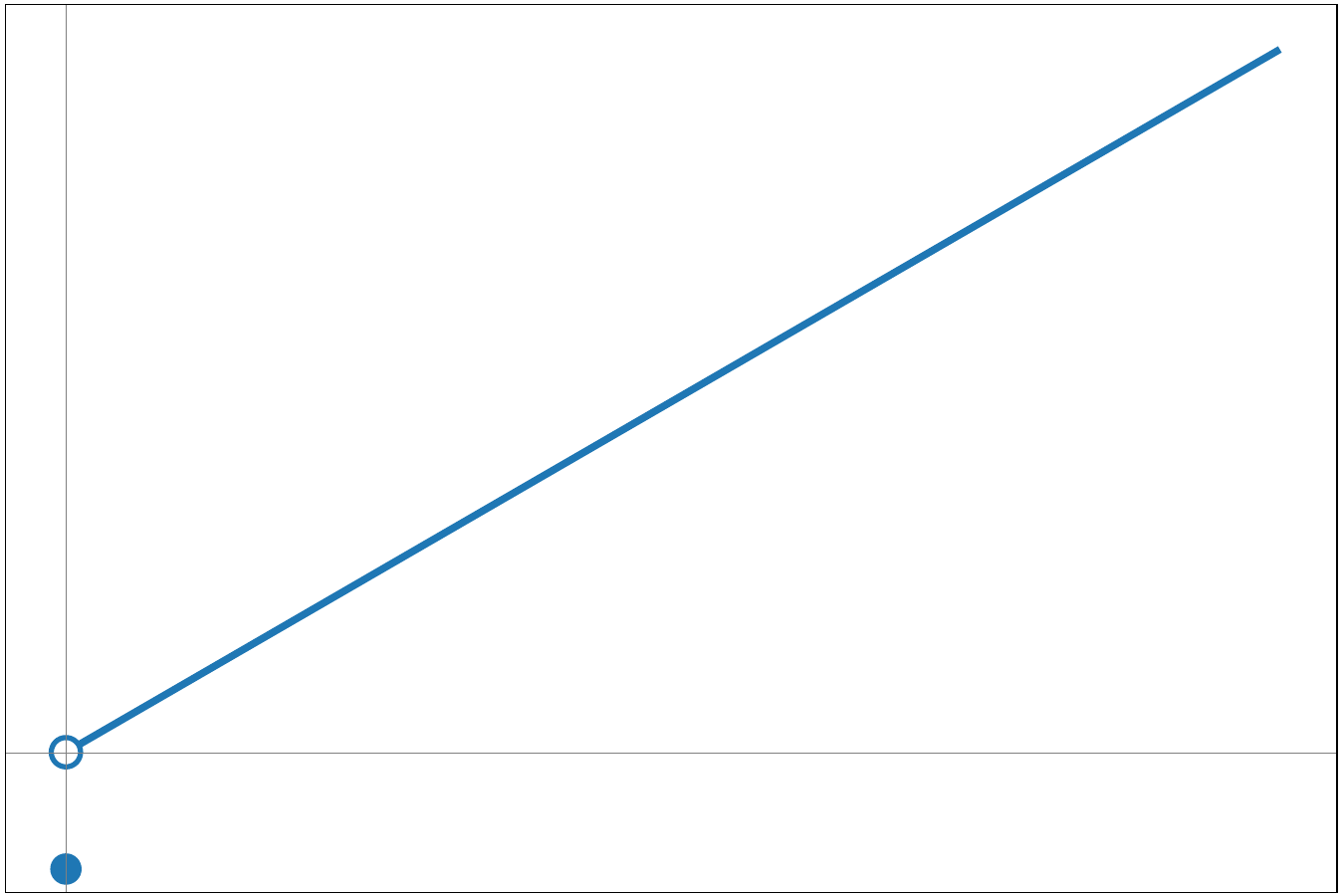}}
\subfigure[A non-lower semicontinuous funtion]{\label{fig:aaa-2}
\includegraphics[width=0.31\linewidth]{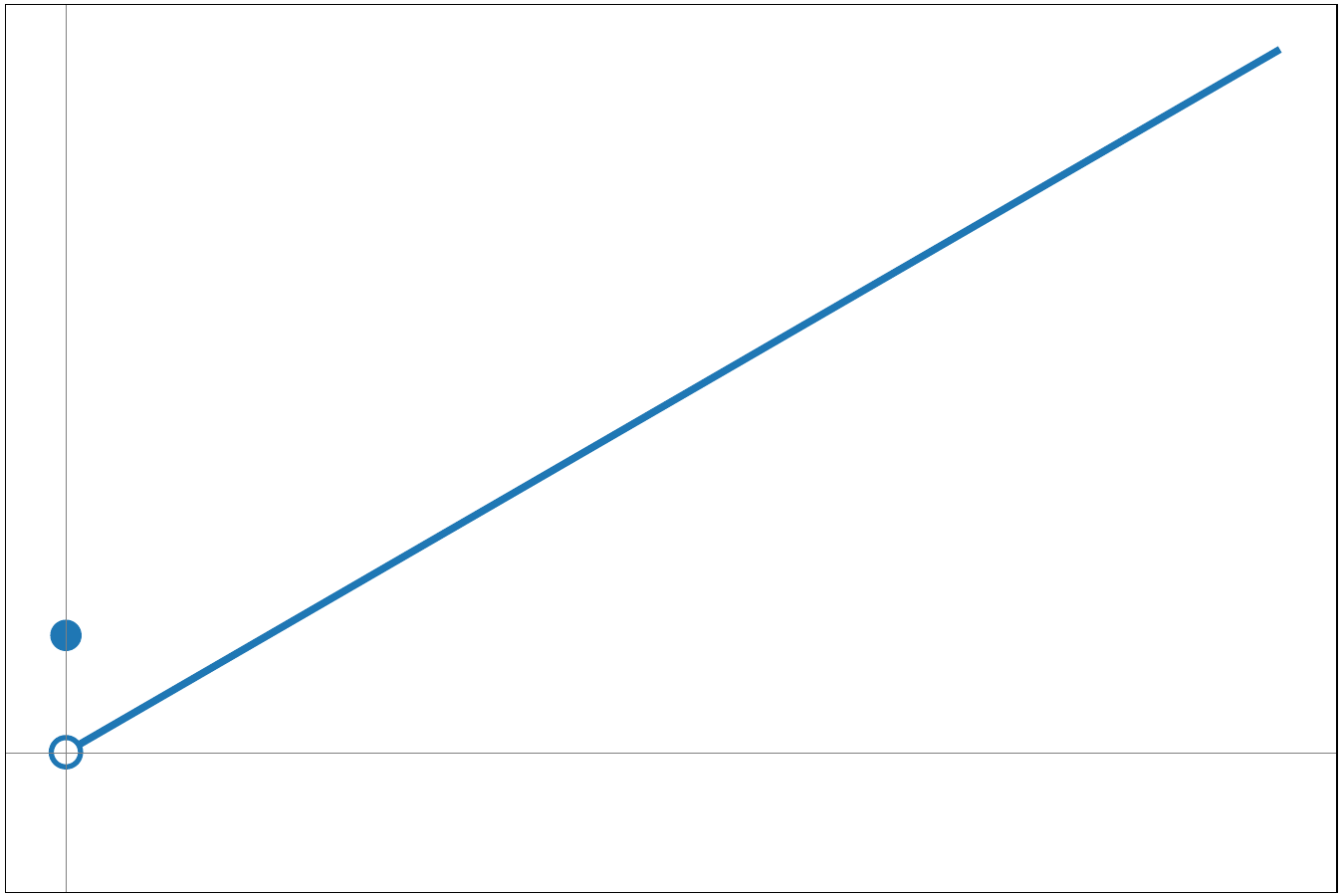}}
\subfigure[A lower semicontinuous funtion.]{\label{fig:aaa-3}
\includegraphics[width=0.31\linewidth]{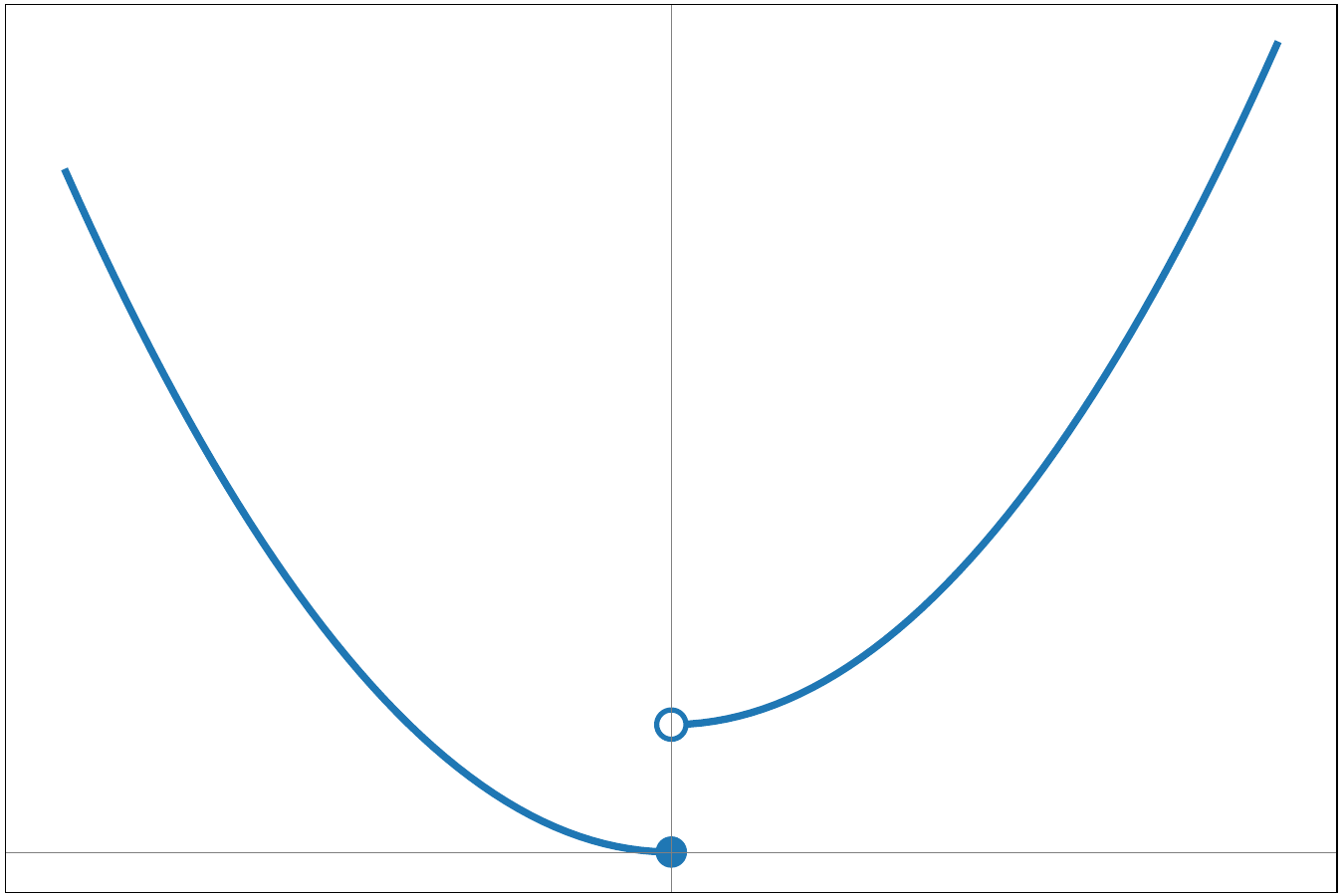}}
\caption{Lower semicontinuous and non-lower semicontinous funtions. In numerical methods and iterative algorithms, lower semicontinuity can help in proving convergence properties. For instance, if a sequence generated by an algorithm is minimizing a lower semicontinuous function, it can be shown that the limit points of the sequence are stationary points.}
\label{fig:low_nonlowersemis}
\end{figure}

Figure~\ref{fig:low_nonlowersemis} illustrates examples of lower semicontinuous and non-lower semicontinuous functions.
Note that allowing the extended real value $\{\infty\}$ in the codomain is standard when defining lower semicontinuity. While not strictly necessary for all functions, this extension provides greater generality and is essential in optimization and variational analysis. It enables us to treat functions that become unbounded above (e.g., indicator functions of constraints) while still discussing their continuity-like properties.

In optimization, lower semicontinuity plays a crucial role because it helps ensure the existence of minimizers. Specifically, if a sequence of feasible points approaches a candidate solution and the objective function is lower semicontinuous, then the function value at the limit point does not exceed the limit inferior of the sequence's values. Under suitable compactness or coercivity conditions, this guarantees that a global minimum is attained.

For instance, in the calculus of variations and optimal control, one often proves existence of minimizers by showing that a functional is both lower semicontinuous and coercive (i.e., $ f(\balpha) \to \infty $ as $ \normtwo{\balpha} \to \infty $). Combined with weak compactness arguments, this yields existence results (see Property~\ref{weier2_prop_close} in Theorem~\ref{theorem:weierstrass_them}).

In numerical optimization and iterative algorithms, lower semicontinuity is also valuable for analyzing convergence. For example, if an algorithm generates a minimizing sequence for a lower semicontinuous function, then any cluster point of the sequence is guaranteed to have a function value no larger than the limit inferior---often implying that such points are stationary or optimal. Moreover, when combined with strong convexity, lower semicontinuity (or closedness) ensures uniqueness and existence of the optimizer (Theorem~\ref{theorem:exi_close_sc}).

For these reasons, throughout this work we primarily consider lower semicontinuous (i.e., \textit{closed}) functions in our analysis.

\begin{definition}[Closed functions]\label{definition:closed_func}
A function $f:\sS\subseteq\real^p\rightarrow \real\cup\{\infty\}$ is said to be \textit{closed} if its epigraph is a closed subset in the product topology of $\sS \times \real$. In other words, for every sequence $(\balpha^{(t)}, r_t)$ in $\text{epi}(f)$ that converges to a limit $(\balpha^*, r^*)$, we have $f(\balpha^*) \leq r^*$, meaning that the limit point $(\balpha^*, r^*)$ also belongs to the epigraph.
\end{definition}

\index{Indicator function}
\begin{exercise}[Closedness of indicator]\label{exercise_closed_indica}
Show that the indicator function $\indicatorS$, which is $0$ if $\balpha$ belongs to a set $\sS$ and $+\infty$ otherwise,  is closed if and only if $\sS$ is a closed set.
\end{exercise}

The following exercise demonstrates that closedness is preserved under several common operations: affine transformations, nonnegative linear combinations, and pointwise maximization.
\begin{exercise}[Preservation of closedness]\label{exercise:pres_clos}
Prove the following results:
\begin{enumerate}[(i)]
\item Let $\bX\in\real^{p\times p}$, and let $f: \real^p \rightarrow (-\infty, \infty]$ be an extended real-valued closed function. Then the function 
$
g(\balpha) = f(\bX\balpha+\by)
$
is closed.

\item Let $f_1, f_2, \ldots, f_k: \real^p \rightarrow (-\infty, \infty]$ be extended real-valued closed functions, and let $\sigma_1, \sigma_2, \ldots, \sigma_k \in \real_+$. Then the function $f = \sum_{i=1}^{k} \sigma_i f_i$ is closed.

\item Let $f_i: \real^p \rightarrow (-\infty, \infty], i \in \sI$ be extended real-valued closed functions, where $\sI$ is a given index set. Then the function
$
f(\balpha) = \max_{i \in \sI} f_i(\balpha)
$
is closed.
\end{enumerate}
\textit{Hint: the epigraph of $\max_{i \in \sI} f_i(\balpha)$ is the intersection of the epigraphs of functions $f_i$.}
\end{exercise}

An important property of closed functions is the equivalence among three fundamental concepts: lower semicontinuity, closedness of the function (i.e., closed epigraph), and closedness of all its sublevel sets.
\begin{theoremHigh}[Equivalence of closedness, LSC, and closedness of level sets]\label{theorem:equiv_close_clkos_semicon}
Let $ f: \real^p\rightarrow  \real \cup \{-\infty, \infty \} $. Then the following three statements  are equivalent:
\begin{enumerate}
\item[(i)] $ f $ is lower semicontinuous (LSC).
\item[(ii)] $ f $ is closed.
\item[(iii)] For any $ \gamma \in \real $, the level set
$
\lev[f, \gamma] = \{ \balpha \in \real^p \mid f(\balpha) \leq \gamma \}
$
is closed.
\end{enumerate}
\end{theoremHigh}
\begin{proof}[of Theorem~\ref{theorem:equiv_close_clkos_semicon}]
\textbf{(i) $\implies$ (ii).} Assume  that $ f $ is lower semicontinuous. We will demonstrate  that $ \epi(f) $ is closed. 
Let  $ \{(\balpha^\toptzero, y^\toptzero)\}_{t \geq 1} \subseteq \epi(f) $ be a sequence such that $ (\balpha^\toptzero, y^\toptzero) \rightarrow (\balpha^*, y^*) $ as $ t \rightarrow \infty $. By definition, for any $ t \geq 1 $,
$$
f(\balpha^\toptzero) \leq y^\toptzero.
$$
Due to the lower semicontinuity of $ f $ at $ \balpha^* $, we have
$$
f(\balpha^*) \leq \liminf_{t \rightarrow \infty} f(\balpha^\toptzero) \leq \liminf_{t \rightarrow \infty} y^\toptzero = y^*,
$$
which shows that $ (\balpha^*, y^*) \in \epi(f) $, proving that $f$ is closed.

\paragraph{(ii) $\implies$ (iii).} Suppose that $ f $ is closed, i.e., the epigraph $ \epi(f) $ is closed. 
For any $ \gamma \in \real $, we aim to show that  $ \lev[f, \gamma] $ is closed. If $ \lev[f, \gamma] = \varnothing $, the claim holds trivially. 
Otherwise, take a sequence $ \{\balpha^\toptzero\}_{t \geq 1} \subseteq \lev[f, \gamma] $ that converges to $ \widehatbalpha $. Obviously, $ (\balpha^\toptzero, \gamma) \in \epi(f) $ for any $ t $ and $ (\balpha^\toptzero, \gamma) \rightarrow (\widehatbalpha, \gamma) $ as $ t \rightarrow \infty $. By the closedness of $ \epi(f) $, it follows that $ (\widehatbalpha, \gamma) \in \epi(f) $, establishing the fact that $ \widehatbalpha \in \lev[f, \gamma] $.

\paragraph{(iii) $\implies$ (i).} Assuming all level sets of $ f $ are closed, we prove that it is lower semicontinuous. Suppose for contradiction that $ f $ is not lower semicontinuous, meaning that there exists $ \balpha^* \in \real^p $ and $ \{\balpha^\toptzero\}_{t \geq 1} \subseteq \real^p $ such that $ \balpha^\toptzero \rightarrow \balpha^* $ and $ \liminf_{t \rightarrow \infty} f(\balpha^\toptzero) < f(\balpha^*) $. 
Choose  $ \gamma $ satisfying
\begin{equation}\label{equation:equiv_lowsemi}
\liminf_{t \rightarrow \infty} f(\balpha^\toptzero) < \gamma < f(\balpha^*). 
\end{equation}
Then there exists a subsequence $ \{\balpha_{n_t}\}_{t \geq 1} $ such that $ f(\balpha_{n_t}) \leq \gamma $ for all $ t \geq 1 $. By the closedness of the level set $ \lev[f, \gamma] $ and the fact that $ \balpha_{n_t} \rightarrow \balpha^* $ as $ t\rightarrow \infty $, it follows that $ f(\balpha^*) \leq \gamma $, which is a contradiction to \eqref{equation:equiv_lowsemi}, establishing that (iii) implies (i).
\end{proof}

To guarantee the existence of a \textit{global minimizer} (also called a global optimum point)---that is, a point where a function attains its smallest value over a given set---the classical Weierstrass theorem provides sufficient conditions. Below, we present the Weierstrass theorem and several of its extensions for continuous and proper closed functions.
\begin{theoremHigh}[Weierstrass theorem and variants\index{Weierstrass theorem}]\label{theorem:weierstrass_them}
We consider the Weierstrass theorem and its variants for different types of functions and sets:
\begin{enumerate}[(i)]
\item \label{weier1_continus} Let $f:\sS\rightarrow (-\infty, \infty]$ be a \textbf{continuous function} defined over a \textit{nonempty and compact (closed and bounded) set} $\sS\subseteq\real^p$.
Then, there exists a global minimum point  and a global maximum point of $f$ over $\sS$.

\item \label{weier1_continus_lev}  Let $f:\sS\rightarrow (-\infty, \infty]$ be a \textbf{continuous function} defined over a \textit{nonempty and closed set} $\sS\subseteq\real^p$. Suppose that all the level sets $\lev[f, \gamma]=\{\balpha\in\sS\mid f(\balpha)\leq \gamma\}$ are bounded. Then, $f$ has a global minimum point over $\sS$.

\item  \label{weier2_continus_coerc} Let $f:\sS\rightarrow  (-\infty, \infty]$ be a \textbf{continuous and coercive function} and let $\sS\subseteq\real^p$ be a \textit{nonempty and closed set}. Then, $f$ has a global minimum point over $\sS$. (\text{The coerciveness ensures the function is bounded over a subset.})

\item \label{weier2_prop_close} Let $f:\sS\rightarrow (-\infty, \infty]$ be a \textbf{proper closed function}, and one of the following is satisfied:
\begin{enumerate}[(a)]
\item \label{weier2_prop_close_v1}  $f$ is defined over a \textit{nonempty and compact set} $\sS$. 
\item \label{weier2_prop_close_v2}  There exists a \textit{nonempty and bounded level set} $\lev[f, \gamma]=\{\balpha\in\sS\mid f(\balpha)\leq \gamma\}$. 
\item \label{weier2_prop_close_v3} $f$ is coercive, i.e., $\mathop{\lim}_{\normtwo{\balpha}\rightarrow\infty} f(\balpha)=\infty$.
\end{enumerate}
Then, there exists a global minimum point of $f$ over $\sS$. And the set of minimizers $\{\balpha\in\sS\mid   f(\balpha)\leq  f(\bbeta), \forall\, \bbeta\in \sS\}$  of $\mathopmin{\balpha\in\sS} f(\balpha)$ is nonempty and compact.
\end{enumerate}

\end{theoremHigh}
Note the coerciveness in \ref{weier2_continus_coerc} ensures the function is bounded over a subset. The three conditions of \ref{weier2_prop_close}.\ref{weier2_prop_close_v1}$\sim$\ref{weier2_prop_close}.\ref{weier2_prop_close_v3} essentially ensure that the minimum value of $f(\balpha)$ cannot be attained at infinity.
The results in \ref{weier2_prop_close} are used more extensively in optimization analysis. To provide a counterexample, let $f=\exp(-x)$ be  defined on $\real$, which is a proper closed function.
However, the domain and some level sets of $f$ are not bounded, and $f$ is not coercive. Therefore, the attainment of the global optimal point is not guaranteed.

\section{Subgradient and Conjugate Functions}\label{section:sub_conjug}

The gradient inequality for convex functions applies specifically to continuously differentiable functions. However, this concept can be generalized using the notion of a \textit{subgradient}, which plays a crucial role in optimization, particularly  for non-differentiable or non-smooth functions.

\begin{definition}[Subgradient and subdifferential\index{Subgradient}\index{Subdifferential}]\label{definition:subgrad}
Let $f: \sS\subseteq\real^p \rightarrow \real$. A vector $\bg\in\real^p$ is called a \textit{subgradient} of $f$ at $\balpha\in\sS$ if
$$
f(\bbeta) \geq  f(\balpha) + \innerproduct{\bg, \bbeta-\balpha}, \quad \text{for all $\bbeta\in\sS$}.
$$
This inequality is known as the  \textit{subgradient inequality}.
The set of all subgradients of $f$ at $\balpha$, denoted by $\partial f(\balpha)$, is called the \textit{subdifferential} of $f$ at $\balpha$:
$$
\partial f(\balpha) \triangleq\left\{\bg\in\real^p\mid f(\bbeta) \geq  f(\balpha) + \innerproduct{\bg, \bbeta-\balpha},
\text{ for all $\bbeta\in\sS$} \right\}.
$$
Any subgradient of $f$ at $\balpha$  may be denoted by
$$
f'(\balpha) \in \partial f(\balpha).
$$
\end{definition}

It is worth noting that the concept of a subgradient applies not only to convex functions but also extends to non-convex settings.
However, for convex functions, it is guaranteed that the subdifferential at any point in the interior of the domain is nonempty (Theorem~\ref{theorem:nonemp_relint_conv}).
If $f$ is differentiable at $\balpha$, then the subdifferential contains only the gradient:
$$
\partial f(\balpha) = \{\nabla f(\balpha)\},
$$

\begin{exercise}[Subdifferential of norms]\label{exercise:sub_norms}
Let $f(\balpha)=\normtwo{\balpha}$. Show that 
$$
\partial f(\balpha)=
\begin{cases}
\left\{\frac{\balpha}{\normtwo{\balpha}}\right\}, & \balpha\neq \bzero;\\
\sB_2[\bzero, 1], & \balpha= \bzero,
\end{cases}
$$
where $\sB_2[\bzero, 1]$ denotes the closed  unit ball in the $\ell_2$-norm.
Additionally, let $g(\balpha)=\normone{\balpha}$. Show that $\partial g(\balpha) = \sum_{i=1}^{p} \partial g_i(\balpha)$, where $g_i(\balpha) = \abs{\alpha_i}$ and 
$$
\partial  g_i(\balpha)=
\begin{cases} 
\{\sign(\alpha_i) \be_i\}, & \alpha_i \neq 0; \\
[-\be_i, \be_i], & \alpha_i = 0.
\end{cases}
$$
Let $\bu\in\partial g(\balpha)$. We  have $\norminf{\bu}\leq 1$ and $\bu^\top\balpha = \normone{\balpha}$.
This also indicates that 
$$
\sign(\balpha)\in\partial g(\balpha),
$$ 
where $\sign(\balpha)$ returns the sign for each component.
This also shows the subgradient inequality for the absolute value function. Specifically,  for any $a\neq 0$ and any $b\in\real$, we have 
$$
\abs{a+b} \geq \abs{a} + \sign(a)b,
$$
with equality if and only if $ \abs{a+b}-\abs{a}$ is  linear in $b$, i.e., when 
$b$ does not cause a sign change or zero-crossing.
\end{exercise}

An established result in convex analysis states that the relative interior of a convex set is always nonempty. Moreover, for a proper convex function, the subdifferential is guaranteed to be nonempty at every point in the relative interior of its effective domain. 
The following theorem summarizes these key facts without proof.
\begin{theoremHigh}[Nonemptiness of relative interior \citep{rockafellar2015convex}]\label{theorem:nonemp_relint_conv}
Let $\sS\subseteq \real^p$ be a nonempty convex set. Then, the relative interior of $\sS$, denoted  $\relint(\sS)$, is nonempty.
Moreover, let $f : \real^p\rightarrow  (-\infty, \infty]$ be a proper convex function, 	and let $\balpha \in \relint(\dom(f))$. 
Then the subdifferential of $f$ at $\balpha$, denoted  $\partial f(\balpha)$, is nonempty.
\end{theoremHigh}

\index{Fermat's theorem}
\textit{Fermat's theorem}, also known as \textit{Fermat's theorem on stationary points}, is a fundamental result in calculus and mathematical optimization. It provides a necessary condition for a differentiable function to attain a local optimum (either a local maximum or a local minimum) at an interior point of its domain.
In the univariate case, Fermat's theorem gives the optimality condition for a point that lies in the interior of an interval---i.e., for a one-dimensional constrained optimization problem where the constraint set is an open interval.

\begin{theoremHigh}[Fermat's theorem: Necessary condition for univariate functions\index{Fermat's theorem}]\label{theorem:fermat_theorem}
Let $f: (a,b)\rightarrow \real$ be a one-dimensional differentiable function defined on an open interval ($a, b$). 
If a point $\alpha^*\in(a,b)$ (i.e., $\alpha^*\in\interior([a,b])$) is a local maximum or minimum of $f$, then $f^\prime(\alpha^*)=0$. 
\end{theoremHigh}
In other words, if a differentiable function $f$ has a local maximum or minimum at a point $\alpha^*$, then the slope of the tangent line at that point is zero; the derivative of the function at that point must be zero.
It's important to emphasize  that this condition is necessary but \textbf{not} sufficient for $\alpha^*$ to be an optimum. There are cases where the derivative is zero but the point is neither a maximum nor a minimum, such as in the case of an inflection point. Additionally, Fermat's theorem does \textbf{not} apply to boundary points of the domain of $f$ or to points where $f$ is not differentiable.

We now present the first-order necessary condition for a local minimum in multivariate optimization.
\begin{corollary}[First-order necessary condition for a  minimum point]\label{corollary:fermat_fist_opt}
Let $f: \real^p \rightarrow \real$ be a  differentiable function. If $\widehatbtheta$ is a (local) minimizer of $f$, then
$$
\nabla f(\widehatbtheta) = \bzero.
$$
Such a point is called a \textit{stationary point} of $f$.
\end{corollary}
\begin{proof}[of Corollary~\ref{corollary:fermat_fist_opt}]
Fix an index $i \in \{1, 2, \ldots, p\}$, and define the one-dimensional function $g(\mu) = f(\widehatbtheta + \mu \be_i)$. Note that $g$ is differentiable at $\mu = 0$ and that $g'(0) = \frac{\partial f}{\partial \theta_i}(\widehatbtheta)$. Since $\widehatbtheta$ is a local minimum point of $f$, it follows that $\mu = 0$ is a local minimum point of $g$, which immediately implies that $g'(0) = 0$ by {Theorem~\ref{theorem:fermat_theorem}}. This equality is exactly the same as $\frac{\partial f}{\partial \theta_i}(\widehatbtheta) = 0$. Since this holds for any $i \in \{1, 2, \ldots, p\}$, we obtain $\nabla f(\widehatbtheta) = \bzero$.
\end{proof}

Most objective functions, especially those with multiple local minima, also contain local maxima and other critical points that satisfy the necessary condition from  Fermat's theorem. 
To distinguish true local minima from these irrelevant or non-optimal critical points, we require additional criteria, such as second-order conditions or convexity assumptions \citep{lu2025practical}.
Subdifferentials are particularly useful for characterizing minimizers, especially in non-smooth or non-differentiable settings. A fundamental result in convex optimization provides a complete characterization of global minima via the subdifferential.
In essence, this is a generalization of Fermat's optimality condition at points of differentiability: $\nabla f(\balpha^*) = \bzero$.

\begin{theoremHigh}[Necessity/sufficiency  of unconstrained convex]\label{theorem:fetmat_opt}
Let $f: \real^p \rightarrow (-\infty, \infty]$ be a proper convex function. Then,
$$
\balpha^* \in \mathop{\argmin}_{\balpha \in \real^p}f(\balpha) 
$$
if and only if $\bzero \in \partial f(\balpha^*)$.
\end{theoremHigh}
\begin{proof}[of Theorem~\ref{theorem:fetmat_opt}]
By the definition of the subgradient (Definition~\ref{definition:subgrad}), it holds that $\balpha^* \in \mathop{\argmin}_{\balpha \in \real^p}f(\balpha) $ if and only if 
$$
f(\balpha) \geq f(\balpha^*) + \innerproduct{\bzero, \balpha - \balpha^*}, \quad \text{for any }  \balpha \in \domain(f),
$$
which is the same as the inclusion $\bzero \in \partial f(\balpha^*)$, where the existence of this subdifferential is stated by Theorem~\ref{theorem:fermat_theorem}.
\end{proof}


We now introduce the definition and key properties of conjugate functions.
\begin{definition}[Conjugate functions\index{Conjugate functions}]\label{definition:conjug_func}
Let $f : \real^p \rightarrow \real \cup \{ \infty \}$ be an extended real-valued function. The \textit{conjugate function}  $f^* : \real^p \rightarrow  \real \cup \{\infty \}$ of $f$ is defined by
$$
f^*(\bbeta) = \max_{\balpha \in \real^p} \{ \innerproduct{\bbeta, \balpha } - f(\balpha) \}, \quad \bbeta \in \real^p.
$$
\end{definition}

Note that $f^*(\bbeta)$ is the pointwise maximum of affine functions of $\bbeta$, each of which is convex and closed.
Consequently, the conjugate function $f^*(\bbeta)$ is always convex and closed, regardless of whether the original function $f$ possesses these properties (see Exercises~\ref{exercise:pres_conv_clos} and \ref{exercise:pres_clos}).
This leads to the following fundamental result:
\begin{lemma}[Closedness and convexity of conjugate]\label{lemma:closedconv_conj}
Let $f:\real^p\rightarrow \real \cup \{ \infty \}$ be an extended real-valued  function. Then its conjugate  $f^*$ is  closed and convex.
\end{lemma}

A direct consequence of the definition is Fenchel's inequality.
\index{Fenchel's inequality}
\begin{theoremHigh}[Fenchel's inequality]\label{theorem:fenchel_ineq}
Let $f: \real^p \rightarrow \real \cup \{ \infty \}$ be a proper   extended real-valued  function. Then for any $\balpha,\bbeta \in \real^p$,
$$
f(\balpha) + f^*(\bbeta) \geq \innerproduct{\bbeta, \balpha}.
$$
\end{theoremHigh}

The \textit{biconjugate} $f^{**}$ (i.e., the conjugate of $f^*$) plays a central role in convex analysis.
\begin{exercise}[Bounds on biconjugate \citep{rockafellar2015convex}]\label{exercise:biconjugate}
Let $f:\real^p\rightarrow [-\infty,\infty]$ be an extended real-valued function. Show that $f(\balpha)\geq f^{**}(\balpha)$ for any $\balpha\in\real^p$.
Additionally, let $g:\real^p\rightarrow \real \cup \{ \infty \}$ be a proper \textbf{closed and convex} function. Show that $g(\balpha)= g^{**}(\balpha)$
\end{exercise}

Because of this property, the biconjugate $f^{**}$ 
is often referred to as the \textit{convex envelope} (or \textit{convex closure}, \textit{convex relaxation}) of $f$---it is the largest convex lower semicontinuous function that lies below $f$.

\begin{exercise}[Properness of conjugate]\label{exercise:proper_conj}
Let $f:\real^p\rightarrow \real \cup \{ \infty \}$ be a proper convex function. Show that $f^*$ is also proper.
\end{exercise}

\begin{exercise}[Conjugate calculus]\label{exercise:conj_calc}
Let $ f : \real^p \to \real \cup \{ \infty \} $. Show that
\begin{enumerate}[(i)]
\item Let$ f_\tau(\balpha) \triangleq f(\tau \balpha) $, where  $ \tau \neq 0 $. Then $ (f_\tau)^*(\bbeta) = f^*(\bbeta/\tau) $.
\item For $ \tau > 0 $, $ (\tau f)^*(\bbeta) = \tau f^*(\bbeta/\tau) $.
\item Let $f_{(\bzeta)} \triangleq f(\balpha - \bzeta)$, where For $\bzeta \in \real^p$. Then $(f_{(\bzeta)})^*(\bbeta) = f^*(\bbeta) + \innerproduct{\bzeta, \bbeta}$.
\end{enumerate}
\end{exercise}

Let us compute the convex conjugate for some examples.

\begin{example}[Self-conjugate]\label{example:self_conj}
Let $f(\balpha) = \frac{1}{2} \normtwo{\balpha}^2$ for  $\balpha \in \real^p$. 
Then $f^*(\bbeta) = \frac{1}{2} \normtwo{\bbeta}^2 = f(\bbeta)$, $\bbeta \in \real^p$. Indeed, since
$ \innerproduct{\bbeta, \balpha} \leq \frac{1}{2} \normtwo{\balpha}^2 + \frac{1}{2} \normtwo{\bbeta}^2$,
we have
$$
f^*(\bbeta) = \max_{\balpha \in \real^p} \{ \innerproduct{\bbeta, \balpha} - f(\balpha)\} \leq \frac{1}{2} \normtwo{\bbeta}^2.
$$
For the reverse inequality, we just set $\balpha = \bbeta$ in the definition of the convex conjugate to obtain
$$
f^*(\bbeta) \geq \normtwo{\bbeta}^2 - \frac{1}{2} \normtwo{\bbeta}^2 = \frac{1}{2} \normtwo{\bbeta}^2.
$$
Hence, equality holds. This function is self-conjugate: it is the unique (up to scaling and translation in certain contexts) function on $\real^p$ satisfying $f=f^*$.
\end{example}

The central result relating subgradients and conjugates is the so-called \textit{conjugate subgradient theorem}, 
which establishes an elegant duality between a convex function and its Fenchel conjugate through their subdifferentials. In essence, it characterizes when a pair of points $(\balpha, \bbeta)$ are dual to each other---not only in the sense of achieving equality in Fenchel's inequality $f(\balpha) + f^*(\bbeta) \geq \innerproduct{\balpha, \bbeta} $, but also in the geometric sense that one is a subgradient of the function at the other. This connection lies at the heart of convex analysis and plays a crucial role in optimization, duality theory, and sensitivity analysis. Moreover, when the function is closed (i.e., lower semicontinuous), the relationship becomes symmetric: subgradients of $f$ correspond precisely to subgradients of $f^*$, reflecting the involutive nature of the conjugate operation.

\index{Conjugate subgradient theorem}
\begin{theoremHigh}[Conjugate subgradient theorem]\label{theorem:conju_subgra}
Let $ f : \real^p \to \real \cup \{ \infty \} $ be a  proper  \textbf{convex} function. The following two claims are equivalent for any $ \balpha, \bbeta \in \real^p $:
\begin{enumerate}[(i)]
\item $ \innerproduct{\balpha, \bbeta} = f(\balpha) + f^*(\bbeta) $. 
\item $ \bbeta \in \partial f(\balpha) $.
\end{enumerate}
\noindent If, in addition,   $ f $ is   \textbf{closed}, then (i) and (ii) are also equivalent to
\begin{itemize}[(iii)]
\item $ \balpha \in \partial f^*(\bbeta) $.
\end{itemize}
\noindent Alternatively, let $ f : \real^p \to \real \cup \{ \infty \} $ be a proper \textbf{closed convex} function. Then for any $\balpha, \bbeta \in \real^p$,
\begin{subequations}
\begin{align}
\partial f(\balpha) &= \argmax_{\bu \in \real^p} \left\{ \innerproduct{\balpha, \bu} - f^*(\bu) \right\};\\
\partial f^*(\bbeta) &= \argmax_{\bv \in \real^p} \left\{ \innerproduct{\bbeta, \bv} - f(\bv) \right\}.
\end{align}
\end{subequations}
\end{theoremHigh}
\begin{proof}[of Theorem~\ref{theorem:conju_subgra}]
By the subgradient inequality of $ \bbeta \in \partial f(\balpha) $, it follows that 
$$
f(\bu) \geq f(\balpha) + \innerproduct{\bbeta, \bu - \balpha} 
\quad\implies\quad 
\innerproduct{\bbeta, \balpha } - f(\balpha) \geq \innerproduct{ \bbeta, \bu} - f(\bu),
\text{ for all } \bu \in \real^p.
$$
Taking the maximum over $ \bu $, the above inequality is equivalent to 
$$
\innerproduct{\bbeta, \balpha} - f(\balpha) \geq f^*(\bbeta),
$$
which by Fenchel's inequality (Theorem~\ref{theorem:fenchel_ineq}) is the same as the equality $\innerproduct{\balpha, \bbeta} = f(\balpha) + f^*(\bbeta)$, establishing the equivalence between (i) and (ii). 
When $f$ is proper closed and convex, by Exercise~\ref{exercise:biconjugate}, $f^{**} = f$, which  implies that (i) is equivalent to
$
\innerproduct{\balpha, \bbeta}= f^*(\bbeta) + f^{**}(\balpha)
$. By the above argument, we conclude that (i) is equivalent to $\balpha \in \partial f^*(\bbeta)$.
This completes the proof.
\end{proof}

\section{Lipschitz Functions}

This section further explores the properties of Lipschitz continuity, defining the relevant concepts and examining their implications for function behavior---particularly the boundedness of subdifferential sets for convex functions.
\index{Lipschitzness}
\begin{definition}[Lipschitz functions]\label{definition:lipschi_funs}
Let $L\geq 0$, and let $f:\sS\rightarrow \real$ be a function defined on a set $\sS\subseteq  \real^p$. Then, the function $f$ is called \textit{L-Lipschitz} if, for all $\balpha,\bbeta\in\sS$, it follows that
$$
\abs{f(\balpha)-f(\bbeta)} \leq L\cdot \normtwo{\balpha-\bbeta}.
$$
If $f$ is continuously differentiable,  its gradient is said to be \textit{Lipschitz continuous} over $\sS$ (or, equivalently, $f$ is called \textit{Lipschitz continuously differentiable}) with constant $L$ if, for all $\balpha,\bbeta\in\sS$, it follows that
$$
\normtwo{\nabla f(\balpha) - \nabla f(\bbeta)} \leq L\cdot \normtwo{\balpha-\bbeta}.~\footnote{The concept of gradient Lipschitzness can be defined using dual norms: $\norm{\nabla f(\balpha) - \nabla f(\bbeta)}_* \leq L\cdot \norm{\balpha-\bbeta}$, where $\norm{\cdot}_*$ is the dual norm of $\norm{\cdot}$. Notably, the vector $\ell_2$-norm is self-dual. In our case, we only consider the vector $\ell_2$-norm.}
$$
The class of functions whose gradients are Lipschitz continuous with constant $L$ is denoted by $C_L^{1,1}(\sS)$ or simply by $C_L^{1,1}$.

Similarly, the function is called \textit{twice Lipschitz continuously differentiable} with constant $L$, denoted by $C_L^{2,2}(\sS)$, if, for all $\balpha,\bbeta\in\sS$, it follows that
$$
\normtwo{\nabla^2 f(\balpha) - \nabla^2 f(\bbeta)} \leq L\cdot \normtwo{\balpha-\bbeta}.
$$

More generally, let $\sS\subseteq\real^p$. We denote by $C_L^{k,s}(\sS)$ the class of functions satisfying the following properties:
\begin{itemize}
\item Any $f\in C_L^{k,s}(\sS)$ is $k$ times continuously differentiable on $\sS$.
\item  Its $s$-th derivative is Lipschitz continuous on $\sS$ with constant $L$:
$\normtwo{\nabla^s f(\balpha) - \nabla^s f(\bbeta)} \leq L\cdot \normtwo{\balpha-\bbeta}$
for all $\balpha,\bbeta\in\sS$. 
\end{itemize}
In this book, we primarily consider the cases $s = 0$, $s=1$, and $s = 2$. Clearly, we always have $s\leq k$. 
Moreover,  if $t\geq k$, then $C_L^{t,s}(\sS) \subseteq  C_L^{k,s}(\sS)$.
For simplicity, we also denotes $C^k(\sS)$ as the set of $k$ times continuously differentiable functions on $\sS$.
\end{definition}

\begin{example}\label{example:lipschitz_spar}
Consider the following examples:\begin{itemize}
\item Since $\absbig{\normone{\balpha} - \normone{\bbeta}} \leq \normone{\balpha-\bbeta} \leq \sqrt{p}\normtwo{\balpha-\bbeta}$ for any $\balpha,\bbeta\in\real^p$ (see Exercise~\ref{exercise:cauch_sc_l1l2}), the function $f_1(\balpha)=\normone{\balpha}$ is thus $\sqrt{p}$-Lipschitz.

\item Since $\absbig{\normtwo{\by-\bX\balpha} - \normtwo{\by-\bX\bbeta}} \leq \normtwo{(\by-\bX\balpha) - (\by-\bX\bbeta)} \leq \normtwo{\bX}\normtwo{\balpha-\bbeta}$, the function $f_2(\balpha) = \normtwo{\by-\bX\balpha}$ is $\normtwo{\bX}$-Lipschitz, where $\normtwo{\bX}$ denotes the spectral norm of $\bX$.
\end{itemize}
The function $f_3(\balpha) \triangleq f_2^2(\balpha) = \normtwo{\by-\bX\balpha}^2$ is not Lipschitz. However, it can be shown that $f_3(\balpha)$ is (gradient) Lipschitz continuously differentiable with constant $2\normtwo{\bX^\top\bX}$.
\end{example}

\begin{exercise}[Gradient Lipschitzness of quadratic functions]
Let $f(\balpha) = \frac{1}{2} \balpha^\top\bX\balpha+\by^\top\balpha+z$, where $\bX\in\real^{p\times p}$ is symmetric, $\by\in\real^p$, and $z\in\real$. Show that the gradient of $f(\balpha)$ is $\normtwo{\bX}$-Lipschitz. This also indicates that the gradient function is 0-Lipschitz when $\bX=\bzero$ (an affine function). 
\end{exercise}

\begin{remark}[Lipschitz continuity implies continuity]
Recall that  a function $ f $ is continuous (Definition~\ref{definition:conti_funs}) at a point $ \bbeta $ if, for every $ \epsilon > 0 $, there exists a $ \delta > 0 $ such that:
$$
\normtwo{\balpha - \bbeta} < \delta \implies \abs{f(\balpha) - f(\bbeta)} < \epsilon.
$$
Now, suppose $ f $ is Lipschitz continuous with constant $ L $, i.e., $\abs{f(\balpha) - f(\bbeta)} \leq L \normtwo{\balpha - \bbeta}$.
Given $ \epsilon > 0 $, choose $ \delta = \frac{\epsilon}{L} $. Then, if $ \normtwo{\balpha - \bbeta} < \delta $, we have:
$$
\abs{f(\balpha) - f(\bbeta)} \leq L \normtwo{\balpha - \bbeta} < L \cdot \frac{\epsilon}{L} = \epsilon.
$$
Thus, $ f $ is continuous at $ \bbeta $. Since $ \bbeta $ was arbitrary, $ f $ is continuous on $ \sS $.~\footnote{Moreover, the choice of $ \delta = \epsilon / L $ {does not depend on $ \bbeta $}---it works uniformly across the domain. This means $ f $ is actually {uniformly continuous}.}
The Lipschitz condition bounds the ``rate of change" of the function, preventing any infinite slopes or wild oscillations that could break continuity.
\end{remark}

\begin{exercise}
Let $ f_1 \in C^{k,s}_{L_1}(\sS) $, $ f_2 \in C^{k,s}_{L_2}(\sS) $, and let $ a_1, a_2 \in \real $. 
Define $ L_3 = \abs{a_1} \cdot L_1 + \abs{a_2} \cdot L_2 $.
Show that $ a_1 f_1 + a_2 f_2 \in C^{k,s}_{L_3}(\sS) $.
\end{exercise}

\begin{exercise}\label{exercise:cl2111_hess_bound}
Show that a function $ f:\real^p\rightarrow \real $ belongs to the class $ C^{2,1}_L(\real^p) \subset C^{1,1}_L(\real^p) $ if and only if for all $ \balpha \in \real^p $ we have
$
\normtwo{\nabla^2 f(\balpha)} \leq L.
$~\footnote{Since the spectral norm of $\nabla^2 f(\balpha)$ is the largest singular value of $\nabla^2 f(\balpha)$, where $\nabla^2 f(\balpha)$ is symmetric. This also indicates that $-L\bI\preceq \nabla^2 f(\balpha) \preceq L\bI$.}
\textit{Hint: use the {fundamental theorem of calculus (Theorem~\ref{theorem:fund_theo_calculu})}.}
\end{exercise}

We now establish the boundedness of the subdifferential under Lipschitz continuity.
\begin{theoremHigh}[Lipschitzness and boundedness of subdifferential sets]\label{theorem:lipsc_equiv}
Let $f:\sS \rightarrow \real$ be a  \textbf{convex} function, and let $\sX \subseteq \interior(\sS)$. Consider the following two claims:
\begin{enumerate}[(i)]
\item $\abs{f(\balpha) - f(\bbeta)} \leq L\normtwo{\balpha - \bbeta}$ for any $\balpha, \bbeta \in \sX$.
\item $\normtwo{\bg} \leq L$ for any $\bg \in \partial f(\balpha), \balpha \in \sX$.
\end{enumerate}
Then,
\begin{enumerate}[(a)]
\item The implication (ii) $\implies$ (i) holds,
\item If $\sX$ is open, then statement  (i) holds if and only if  statement  (ii) holds.
\end{enumerate}
\end{theoremHigh}
\begin{proof}[of Theorem~\ref{theorem:lipsc_equiv}]
\textbf{(a).} Assume that (ii) is satisfied, and let $\balpha, \bbeta \in \sX$. Let $\bg_x \in \partial f(\balpha)$ and $\bg_y \in \partial f(\bbeta)$.
By the definitions of $\bg_x, \bg_y$ and the  Cauchy--Schwarz inequality (Theorem~\ref{theorem:cs_matvec}),
\begin{align*}
f(\balpha) - f(\bbeta) &\leq \innerproduct{\bg_x, \balpha - \bbeta} \leq \normtwo{\bg_x} \normtwo{\balpha - \bbeta} \leq L\normtwo{\balpha - \bbeta}; \\
f(\bbeta) - f(\balpha) &\leq \innerproduct{\bg_y, \bbeta - \balpha} \leq \normtwo{\bg_y} \normtwo{\balpha - \bbeta} \leq L\normtwo{\balpha - \bbeta},
\end{align*}
which establishes the validity of (i).

\paragraph{(b).} The implication (ii) $\implies$ (i) was shown in (a). Conversely, assume that (i) is satisfied. Take $\balpha \in \sX$ and $\bg \in \partial f(\balpha)$. Let $\widetildebg\triangleq\frac{\bg}{\normtwo{\bg}}$ such that $\normtwo{\widetildebg} = 1 $ and $ \innerproduct{\widetildebg, \bg} = \normtwo{\bg}$. Since $\sX$ is open, we can take $\epsilon > 0$ small enough such that $\balpha + \epsilon \widetildebg \in \sX$. 
Given the subgradient inequality, it follows that
$
f(\balpha + \epsilon \widetildebg) \geq f(\balpha) + \epsilon \innerproduct{\bg, \widetildebg},
$
whence we have
$$
\epsilon \normtwo{\bg}  \leq f(\balpha + \epsilon \widetildebg) - f(\balpha) \leq L \normtwo{\balpha + \epsilon \widetildebg - \balpha} = L\epsilon
\quad\implies\quad \normtwo{\bg} \leq L.
$$
This completes the proof.
\end{proof}

\section{Strongly Convex, Strongly Smooth Functions}
This section provides rigorous definitions of \textit{strong convexity (SC)} and \textit{strong smoothness (SS, or simply smoothness)}. Although both concepts can be extended to non-differentiable convex functions using subgradients, we restrict our attention in this book to differentiable functions.
\begin{figure}[htp]
\centering       
\vspace{-0.25cm}                 
\subfigtopskip=2pt               
\subfigbottomskip=-2pt         
\subfigcapskip=-10pt      
\includegraphics[width=0.98\textwidth]{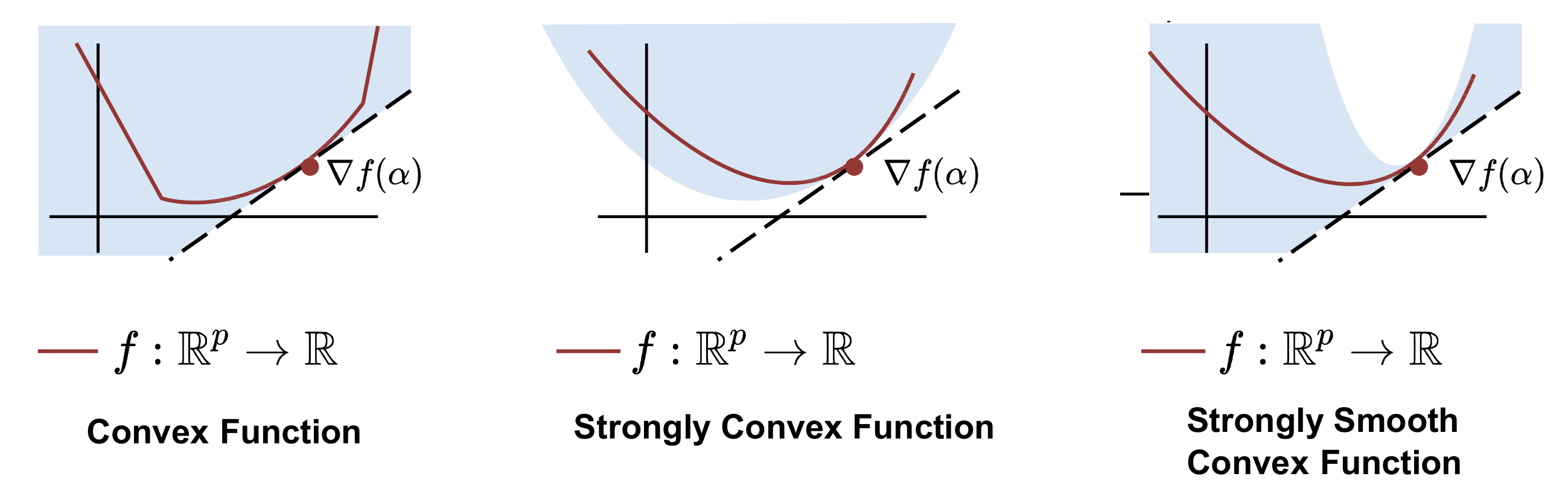}
\caption{
A convex function always lies above its tangent at any point. In contrast, a strongly convex function is not only bounded below by its tangent but also by a quadratic function, which limits how slowly it can grow. Similarly, a smooth function is bounded above by a quadratic model, preventing it from increasing too rapidly. Together, strong convexity and smoothness confine the function's behavior between two quadratic envelopes. The shaded region in each plot indicates where the graph of the function is allowed to lie. The figure is adapted from\citet{jain2017non}.}
\label{fig:cvxfunc}
\end{figure}
\index{Strongly convex}
\index{Strongly smooth}
\index{Smooth}
\begin{definition}[Strongly convex/smooth functions]\label{definition:scss_func}
Let $f:\sS\rightarrow \real$ be a  differentiable function defined on a \textbf{convex} set $\sS\subseteq  \real^p$~\footnote{Note the notions of strong convexity and smoothness are primarily meaningful for differentiable functions (at least mostly for differentiable) defined on a convex set. For example, $\sS$ can be set to $\real^p$. The reason will be clear in the sequel.}. Then, $f$ is said to be \textit{$L_a$-strongly convex (SC, or simply $L_a$-convex)} and \textit{$L_b$-strongly smooth (SS, or simply $L_b$-smooth)} if, for every $\balpha,\bbeta\in\sS$, it follows that 
\begin{equation}\label{equation:scss_func1}
\frac{L_a}{2} \normtwo{\balpha-\bbeta}^2
\leq f(\bbeta)-f(\balpha)- \innerproduct{\nabla f(\balpha), (\bbeta-\balpha)}
\leq \frac{L_b}{2} \normtwo{\balpha-\bbeta}^2.~\footnote{In many texts, the second inequality is called the \textit{descent lemma} for $L_b$-strongly smooth functions since $f(\bbeta) \leq f(\balpha)+\nabla f(\balpha)^\top (\bbeta-\balpha)+\frac{L_b}{2} \normtwo{\balpha-\bbeta}^2$, i.e., an update of the function value  from $f(\balpha)$ to $f(\bbeta)$ is upper-bounded.}
\end{equation}
When the function is twice differentiable, Theorem~\ref{theorem:psd_hess_conv} also indicates  that 
\begin{equation}\label{equation:scss_func2}
L_a\bI \preceq \nabla^2 f(\balpha) \preceq L_b\bI, \quad\text{for all } \balpha\in\interior(\sS).~\footnote{We relax to $\balpha\in\interior(\sS)$ as opposed to claiming an open set in Theorem~\ref{theorem:psd_hess_conv}.}
\end{equation}
That is, the smallest eigenvalue of $\nabla^2 f(\balpha)$ is at least $L_a$, and the largest eigenvalue of $\nabla^2 f(\balpha)$ is at most $L_b$:
\begin{equation}\label{equation:scss_func3_ra}
L_a\leq \frac{\bv^\top \nabla^2 f(\balpha) \bv}{\normtwo{\bv}^2} \leq L_b, \quad\text{for all } \balpha\in\sS, \text{and nonzero } \bv\in\real^p.~\footnote{$\frac{\bv^\top\bX\bv}{\normtwo{\bv}^2}$ is called the \textit{Rayleigh quotient} of vector $\bv$ associated with the matrix $\bX$, which lies between the largest and smallest eigenvalues of $\bX$; see, for example, \citet{lu2021numerical, lu2025practical} for more details.}
\end{equation}
It is evident that if the function is $L_b'$-strongly smooth, it is also $L_b''$-strongly smooth for any $L_b''\geq L_b'$.
An $L_a$-strongly convex function is necessarily a convex function (when $L_a=0$, it reduces to a standard convex function). 
While if the function is $L_a'$-strongly convex, it is also $L_a''$-strongly convex for any $L_a''\in(0, L_a']$.
\end{definition}

Strong convexity and smoothness are fundamental in the analysis of gradient-based optimization algorithms (see Chapter~\ref{chapter:spar}). Strong convexity ensures that the function curves upward sufficiently fast, providing a positive lower bound on its curvature and preventing gradients from vanishing too quickly near minima. Conversely, strong smoothness imposes an upper bound on the curvature, guaranteeing that the gradient does not change too abruptly. This controls the size of allowable steps in gradient descent and prevents overly aggressive updates (see Figure~\ref{fig:cvxfunc}).

Definition~\ref{definition:res_scss_func} later extends these notions to functions defined on non-convex domains.

\begin{exercise}[Sum of SC and convex, SS]\label{exercise:sum_sc_conv}
Prove the following two statements:
\begin{itemize}
\item Let $f$ be  $L_a$-SC  and $g$ be  convex. Show that $f+g$ is $L_a$-SC. 
\item Let $f$ be  $L_b'$-SS  and $g$ be  $L_b''$-SS. Show that $f+g$ is $(L_b'+L_b'')$-SS. 
\end{itemize}
\end{exercise}

\begin{theoremHigh}[SS Property-O: Equivalence between gradient Lipschitzness and smoothness]\label{theorem:equi_gradsch_smoo}
Let $f:\sS\rightarrow \real$ be a  differentiable function defined on a \textbf{convex} set $\sS\subseteq  \real^p$.
If $f$ is gradient  Lipschitz continuous with constant $L_b$, i.e., $f\in C_{L_b}^{1,1}(\sS)$, then it is also $L_b$-smooth.
\end{theoremHigh}
\begin{proof}[of Theorem~\ref{theorem:equi_gradsch_smoo}]
By the fundamental theorem of calculus (Equation~\eqref{equation:fund_theo_calculu3}), it follows that 
$
f(\bbeta) - f(\balpha) = \innerproduct{\nabla f(\balpha), \bbeta - \balpha}+ \int_0^1 \innerproduct{\nabla f(\balpha + \mu(\bbeta - \balpha)) - \nabla f(\balpha), \bbeta - \balpha} d\mu.
$
Since $\sS$ is a convex set, the function $f(\balpha + \mu(\bbeta - \balpha))$ is well-defined for $\mu\in[0,1]$.
This results in 
$$
\begin{aligned}
&\gap \abs{f(\bbeta) - f(\balpha) - \innerproduct{\nabla f(\balpha), \bbeta - \balpha}} =  \abs{\int_0^1 \innerproduct{\nabla f(\balpha + \mu(\bbeta - \balpha)) - \nabla f(\balpha), \bbeta - \balpha} d\mu } \\
&\overset{\dag}{\leq} \int_0^1 \normtwo{\nabla f(\balpha + \mu(\bbeta - \balpha)) - \nabla f(\balpha)}\cdot  \normtwo{\bbeta - \balpha} d\mu
\leq   \int_0^1 L_b \mu \normtwo{\bbeta - \balpha}^2 d\mu 
= \frac{L_b}{2} \normtwo{\bbeta - \balpha}^2,
\end{aligned}
$$
where the inequality $(\dag)$ follows from the  Cauchy--Schwarz inequality (Theorem~\ref{theorem:cs_matvec}).
Hence, $f$ is $L_b$-smooth.
\end{proof}

If $ f $ is $L_b$-smooth and admits a global minimizer $ \balpha^* $, an important consequence  is 
that the quadratic upper bound provided by smoothness can be used to relate the function suboptimality $f(\balpha) - f(\balpha^*)$ to the norm of the gradient at any point $ \balpha $. This relationship is particularly powerful when combined with strong convexity, as it yields tight two-sided bounds.

\begin{theoremHigh}[SS\&SC Property-I: Bound of $f(\balpha) - f(\balpha^*)$]\label{theorem:smoo_prop2_bound}
Let $ f:\real^p\rightarrow \real $ be  a differentiable function  defined on $ \real^p $  that possesses a global minimizer $ \balpha^* $. Suppose $ f(\balpha) $ is  $ L_b $-smooth and $L_a$-strongly convex. 
Then,
$$
\frac{1}{2L_b} \normtwo{\nabla f(\balpha)}^2 \leq f(\balpha) - f(\balpha^*)
\leq \frac{1}{2L_a}  \normtwo{\nabla f(\balpha)}^2,\quad  \text{ for any } \balpha\in\real^p.
$$
\end{theoremHigh}
\begin{proof}[of Theorem~\ref{theorem:smoo_prop2_bound}]
Since $ \balpha^* $ is a global minimuizer, for any $\bbeta\in\real^p$, applying the quadratic upper bound from the smoothness gives
$$
f(\balpha^*) \leq f(\bbeta) \leq f(\balpha) + \nabla f(\balpha)^\top (\bbeta - \balpha) + \frac{L_b}{2} \normtwo{\bbeta - \balpha}^2.
$$
Fixing $ \balpha $, note that the above inequality holds for any $ \bbeta $, thus we can take the infimum on the right side of the inequality:
$$
\begin{aligned}
f(\balpha^*) &\leq \min_{\bbeta \in \real^p} \left\{ f(\balpha) + \nabla f(\balpha)^\top (\bbeta - \balpha) + \frac{L_b}{2} \normtwo{\bbeta - \balpha}^2 \right\} 
= f(\balpha) - \frac{1}{2L_b} \normtwo{\nabla f(\balpha)}^2,
\end{aligned}
$$
where the minimum is attained at $\bbeta = \balpha-\frac{1}{L_b}\nabla f(\balpha)$ (since the domain is $\real^p$).

For the second part, since $f$ is strongly convex, we have 
$$
\begin{aligned}
f(\bbeta)&\geq f(\balpha) +\nabla f(\balpha)^\top (\bbeta-\balpha) + \frac{L_a}{2} \normtwo{\bbeta-\balpha}^2\\
&\geq f(\balpha) +\nabla f(\balpha)^\top (\widetildebbeta-\balpha) + \frac{L_a}{2} \normtwo{\widetildebbeta-\balpha}^2
=f(\balpha)-\frac{1}{2L_a} \normtwo{\nabla f(\balpha)}^2,
\end{aligned}
$$
where $\widetildebbeta \triangleq \balpha-(1/L_a) \nabla f(\balpha)$ minimizes the right-hand side of the above inequality. 
Letting $\bbeta\triangleq\balpha^*$ completes the proof.
\end{proof}

An immediate implication of the upper bound is the following practical criterion: for any tolerance $\varepsilon>0$,
\begin{equation}
\normtwo{\nabla f(\balpha)} \leq (2L_a\varepsilon)^{1/2}
\quad\implies\quad 
f(\balpha) - f(\balpha^*)\leq \varepsilon,
\end{equation}
Thus, controlling the gradient norm provides a certificate of near-optimality.
This generalizes the  classical  \textit{optimality condition for convex functions} stated  in Theorem~\ref{theorem:fetmat_opt}, which states that $\nabla f(\balpha)=\bzero$ characterizes exact minimizers in the differentiable case.
Moreover, the theorem reveals a key growth property: near the minimizer, a strongly convex and smooth function behaves like a quadratic---neither too flat nor too steep.

The next result establishes that strong convexity alone (even without smoothness) guarantees both existence and uniqueness of a minimizer, along with a quadratic growth condition around it.

\begin{theoremHigh}[SC Property-II: Existence and uniqueness of a minimizer of closed SC functions]\label{theorem:exi_close_sc}
Let $f :  \real^p \rightarrow (-\infty, \infty]$ be a differentiable~\footnote{Note that this theorem also holds for non-differentiable functions, and subgradients can be used in the proof by Theorem~\ref{theorem:nonemp_relint_conv}.}, proper closed, and $L_a$-strongly convex function ($L_a > 0$). Then,
\begin{enumerate}[(i)]
\item  $f$ has a unique minimizer.
\item  $f(\balpha) - f(\balpha^*) \geq \frac{L_a}{2} \normtwo{\balpha - \balpha^*}^2$ for all $\balpha \in \real^p$, where $\balpha^*$ is the unique minimizer of $f$.
\end{enumerate}
\end{theoremHigh}
\begin{proof}[of Theorem~\ref{theorem:exi_close_sc}]
\textbf{(i).} Fix any $\balpha_0\in\real^p$. By  strong convexity, it follows that
$$
\begin{aligned}
&f(\balpha) \geq f(\balpha_0) + \innerproduct{\nabla f(\balpha_0), \balpha - \balpha_0} + \frac{L_a}{2} \normtwo{\balpha - \balpha_0}^2 \\
&\quad\implies 
f(\balpha) \geq f(\balpha_0) - \frac{1}{2L_a} \normtwo{\nabla f(\balpha_0)}^2 + \frac{L_a}{2} \normtwo{\balpha - \left(\balpha_0 - \frac{1}{L_a} \nabla f(\balpha_0)\right)}^2 \quad \text{for any } \balpha \in \real^p.
\end{aligned}
$$
This indicates that
$$
\lev[f, f(\balpha_0)] \subseteq \sB_{2}\left[\balpha_0 - \frac{1}{L_a} \nabla f(\balpha_0), \frac{1}{L_a} \normtwo{\nabla f(\balpha_0)}\right].
$$
Since $f$ is closed, the above level set is closed (Theorem~\ref{theorem:equiv_close_clkos_semicon}); and since it is contained in a ball, it is also bounded. Therefore, $\lev[f, f(\balpha_0)]$ is compact. We can thus deduce that the optimal set for minimizing $f$ over $\dom(f)$ is identical to the optimal set for minimizing $f$ over the nonempty compact set $\lev[f, f(\balpha_0)]$. Invoking Weierstrass theorem  for proper closed functions (\ref{weier2_prop_close}.\ref{weier2_prop_close_v1} in Theorem~\ref{theorem:weierstrass_them}), it follows that a minimizer exists. To show the uniqueness, assume that $\widetildebalpha$ and $\balpha^*$ are minimizers of $f$. Then $f(\widetildebalpha)  = f(\balpha^*)$, where $f(\balpha^*)$ is the minimal value of $f$. Then by property (iv) of Theorem~\ref{theorem:charac_stronconv},
$$
f(\balpha^*) \leq f\left(\frac{1}{2}\widetildebalpha + \frac{1}{2}\balpha^*\right) \leq \frac{1}{2} f(\widetildebalpha) + \frac{1}{2} f(\balpha^*) - \frac{L_a}{8} \normtwo{\widetildebalpha - \balpha^*}^2 = f(\balpha^*) - \frac{L_a}{8} \normtwo{\widetildebalpha - \balpha^*}^2,
$$
implying that $\widetildebalpha = \balpha^*$, thus proving the uniqueness of the minimizer of $f$.

\paragraph{(ii).} Let $\balpha^*$ be the unique minimizer of $f$. Then by optimality condition for convex functions (Theorem~\ref{theorem:fetmat_opt}), $ \nabla f(\balpha^*) = \bzero$. Thus, by the definition of strong convexity,
$$
f(\balpha) - f(\balpha^*) \geq \innerproduct{\bzero, \balpha - \balpha^* } + \frac{L_a}{2} \normtwo{\balpha - \balpha^*}^2 = \frac{L_a}{2} \normtwo{\balpha - \balpha^*}^2 ,
\;\text{ for any $\balpha \in \real^p$},
$$
establishing claim (ii).
\end{proof}

The definition of smoothness does not require the function to be convex. However, when it indeed is, the function can be equivalently characterized as follows. 
\begin{theoremHigh}[Characterization theorem of SS and convexity \citep{beck2017first}]\label{theorem:charac_smoo}
Let $f: \sS\subseteq\real^p \rightarrow \real$ be a differentiable \textbf{convex} function defined on a convex set $\sS$, and let $L_b > 0$. Then the following statements  are equivalent:
\begin{enumerate}[(i)]
\item $f$ is $L_b$-smooth.
\item The function $g(\balpha)\triangleq \frac{L_b}{2}\balpha^\top\balpha - f(\balpha)$ is convex.
\item $f(\bbeta) \leq f(\balpha) + \innerproduct{\nabla f(\balpha), \bbeta - \balpha} + \frac{L_b}{2} \normtwo{\balpha - \bbeta}^2$ for all $\balpha, \bbeta \in \sS$.
\item $ f(\bbeta) \geq f(\balpha) + \innerproduct{\nabla f(\balpha), \bbeta - \balpha} + \frac{1}{2L_b} \normtwo{\nabla f(\balpha) - \nabla f(\bbeta)}^2 $ for all $\balpha, \bbeta \in \sS$.
\item $f(\lambda \balpha + (1 - \lambda)\bbeta) \geq \lambda f(\balpha) + (1 - \lambda)f(\bbeta) - \frac{L_b}{2} \lambda (1 - \lambda) \normtwo{\balpha - \bbeta}^2$ for any $\balpha, \bbeta \in \sS$ and $\lambda \in [0, 1]$.
\end{enumerate}
Reversing $\balpha$ and $\bbeta$ in (iii) and (iv), we can also obtain 
\begin{equation}\label{equ:charac_smoo_3} 
\frac{1}{L_b} \normtwo{\nabla f(\balpha) - \nabla f(\bbeta)}^2 
\leq  \innerproduct{\nabla f(\balpha) - \nabla f(\bbeta), \balpha - \bbeta} 
\leq  
L_b\normtwo{\balpha-\bbeta}^2,  
\text{ for all  $\balpha, \bbeta \in \sS$.}
\end{equation}
\end{theoremHigh}
\begin{proof}[of Theorem~\ref{theorem:charac_smoo}]
\textbf{(i)$\implies$(ii).} For any $\balpha,\bbeta\in\sS$, it follows that 
$$
\begin{aligned}
\big(\nabla g(\balpha) - \nabla g(\bbeta)\big)^\top (\balpha - \bbeta) 
&= L_b \normtwo{\balpha - \bbeta}^2 - \big(\nabla f(\balpha) - \nabla f(\bbeta)\big)^\top (\balpha - \bbeta)\\
&\geq L_b \normtwo{\balpha - \bbeta}^2 - \normtwo{\balpha - \bbeta} \normtwo{\nabla f(\balpha) - \nabla f(\bbeta)} \geq 0,
\end{aligned}
$$
where the last inequality follows from Theorem~\ref{theorem:equi_gradsch_smoo} and the definition of smoothness.
Thus, $ g(\balpha) $ is a convex function by Theorem~\ref{theorem:monoton_convgrad}.

\paragraph{(ii)$\implies$(iii).}
Since $ g(\balpha) \triangleq \frac{L_b}{2} \normtwo{\balpha}^2 - f(\balpha) $ is convex, this implies
$$
g(\bbeta) \geq g(\balpha) + \nabla g(\balpha)^\top (\bbeta - \balpha), \quad \forall\, \balpha, \bbeta \in \sS,
$$
where
$
\nabla g(\balpha) = L_b \balpha - \nabla f(\balpha).
$
Substituting the gradient into the convexity condition yields that
\begin{align}
\frac{L_b}{2} \normtwo{\bbeta}^2 - f(\bbeta) \geq \frac{L_b}{2} \normtwo{\balpha}^2 - f(\balpha) + \big(L_b \balpha - \nabla f(\balpha)\big)^\top (\bbeta - \balpha) \nonumber\\
\quad\implies\quad f(\bbeta) - f(\balpha) \leq \nabla f(\balpha)^\top (\bbeta - \balpha) + \frac{L_b}{2} \normtwo{\bbeta - \balpha}^2. \label{equation:charac_smoo11}
\end{align}

\paragraph{(iii)$\implies$(iv).} We can assume that $\nabla f(\balpha) \neq \nabla f(\bbeta)$ since otherwise the inequality (iv) is trivial by the convexity of $f$. For a fixed $\balpha \in \sS$, consider the function
$$ 
g_{\balpha}(\bbeta) \triangleq f(\bbeta) - f(\balpha) - \innerproduct{\nabla f(\balpha), \bbeta - \balpha}, \quad \bbeta \in \sS, 
$$
which is convex w.r.t. $\bbeta$.
Note that $\nabla g_{\balpha}(\bbeta) = \nabla f(\bbeta)-\nabla f(\balpha)$ for any $\bbeta\in\sS$.
By the condition of (iii), for any $\bbeta, \bu \in \sS$, we have 
\begin{equation}\label{equation:char_ss_eqq1}
\begin{aligned}
g_{\balpha}(\bu) &= f(\bu) - f(\balpha) - \innerproduct{\nabla f(\balpha), \bu - \balpha} \\
&\leq \big\{f(\bbeta) + \innerproduct{\nabla f(\bbeta), \bu - \bbeta} + \frac{L_b}{2} \normtwo{\bu - \bbeta}^2\big\} - f(\balpha) - \innerproduct{\nabla f(\balpha), \bu - \balpha} \\
&= f(\bbeta) - f(\balpha) - \innerproduct{\nabla f(\balpha), \bbeta - \balpha} + \innerproduct{\nabla f(\bbeta) - \nabla f(\balpha), \bu - \bbeta} + \frac{L_b}{2} \normtwo{\bu - \bbeta}^2 \\
&= g_{\balpha}(\bbeta) + \innerproduct{\nabla g_{\balpha}(\bbeta), \bu - \bbeta} + \frac{L_b}{2} \normtwo{\bu - \bbeta }^2,
\end{aligned}
\end{equation}
Since $\nabla g_{\balpha}(\balpha) = \bzero$, which by the convexity of $g_{\balpha}$ implies that $\balpha$ is a global minimizer of $g_{\balpha}$ such that
$
0=g_{\balpha}(\balpha) \leq g_{\balpha}(\bu),\, \text{for all } \bu \in \sS.  
$
Let $\bbeta \in \sS$, and let $\widetildebg \triangleq \frac{\nabla g_{\balpha}(\bbeta)}{\normtwo{\nabla g_{\balpha}(\bbeta)}}$ such that  $\normtwo{\widetildebg} = 1$ and $\innerproduct{\nabla g_{\balpha}(\bbeta), \widetildebg} = \normtwo{\nabla g_{\balpha}(\bbeta)}$. Defining
$ \bu \triangleq \bbeta - \frac{\normtwo{\nabla g_{\balpha}(\bbeta)}}{L_b} \widetildebg  $, the optimum of $g_{\balpha}(\balpha)$ shows that 
$$
0 = g_{\balpha}(\balpha) \leq g_{\balpha} \left( \bbeta - \frac{\normtwo{\nabla g_{\balpha}(\bbeta)}}{L_b} \widetildebg \right). 
$$
Combining the preceding inequality with \eqref{equation:char_ss_eqq1}, we obtain
\begin{align*}
0 &= g_{\balpha}(\balpha) 
\leq g_{\balpha}(\bbeta) - \frac{\normtwo{\nabla g_{\balpha}(\bbeta)}}{L_b} \innerproduct{\nabla g_{\balpha}(\bbeta), \widetildebg} + \frac{\normtwo{\nabla g_{\balpha}(\bbeta)}^2 \cdot \normtwo{\widetildebg}^2 }{2L_b} 
= g_{\balpha}(\bbeta) - \frac{\normtwo{\nabla g_{\balpha}(\bbeta)}^2}{2L_b}  \\
&= f(\bbeta) - f(\balpha) - \innerproduct{\nabla f(\balpha), \bbeta - \balpha} - \frac{1}{2L_b} \normtwo{\nabla f(\balpha) - \nabla f(\bbeta)}^2,
\end{align*}
establishing (iv).

\paragraph{(iii)$\implies$(v).}
Let $\balpha, \bbeta \in \sS$ and $\lambda \in [0, 1]$, and denote $\bxi \triangleq \lambda \balpha + (1 - \lambda) \bbeta$. Since $\sS$ is convex, we also have $\bxi\in\sS$.
Invoking (iii) with $\balpha,\bxi$ and $\bbeta,\bxi$, and plugging in the expression of $\bxi$ yields that 
$$
\begin{aligned}
f(\balpha) &\leq f(\bxi) + (1 - \lambda) \innerproduct{\nabla f(\bxi), \balpha - \bbeta} + \frac{L_b (1 - \lambda)^2}{2} \normtwo{\balpha - \bbeta}^2, \\
f(\bbeta) &\leq f(\bxi) + \lambda \innerproduct{\nabla f(\bxi), \bbeta - \balpha} + \frac{L_b \lambda^2}{2} \normtwo{\balpha - \bbeta}^2.
\end{aligned}
$$
Multiplying the first inequality by $\lambda$ and the second by $1 - \lambda$ and adding them yields the inequality (v).

\paragraph{(v)$\implies$(iii).} Rearranging terms in the inequality (v), we obtain that it is equivalent to
$$
f(\bbeta) \leq f(\balpha) + \frac{f\big(\balpha + (1 - \lambda)(\bbeta - \balpha)\big) - f(\balpha)}{1 - \lambda} + \frac{L_b}{2} \lambda \normtwo{\balpha - \bbeta}^2. 
$$
Taking $\lambda \to 1^-$, the preceding inequality becomes
$ f(\bbeta) \leq f(\balpha) + f'(\balpha; \bbeta - \balpha) + \frac{L_b}{2} \normtwo{\balpha - \bbeta}^2$.
Since $f$ is differentiable, \eqref{equation:direc_contdiff} shows that $f'(\balpha; \bbeta - \balpha) = \innerproduct{\nabla f(\balpha), \bbeta - \balpha}$, establishing the desired result. 
\end{proof}

\begin{theoremHigh}[Characterization theorem of SC \citep{beck2017first}]\label{theorem:charac_stronconv}
Let $ f:\sS\subseteq\real^p\rightarrow \real $ be a convex and differentiable function defined on a convex set $\sS$, and let $ L_a > 0 $. Then the following  statements are equivalent:
\begin{enumerate}[(i)]
\item $ f $ is $ L_a $-strongly convex.
\item $
f(\bbeta) \geq f(\balpha) + \innerproduct{\nabla f(\balpha) , \bbeta - \balpha} + \frac{L_a}{2} \normtwo{\bbeta - \balpha}^2
$
for any $ \balpha , \bbeta \in \sS$.
\item $
\innerproduct{\nabla f(\balpha)  - \nabla f(\bbeta) , \balpha - \bbeta} \geq L_a \normtwo{\balpha - \bbeta}^2
$
for any $ \balpha, \bbeta \in \sS $.
\item $f(\lambda \balpha+(1-\lambda)\bbeta )\leq \lambda f(\balpha)+(1-\lambda)f(\bbeta) -\frac{L_a}{2}\lambda(1-\lambda)\normtwo{\balpha-\bbeta}^2$ for any $\lambda\in[0,1]$ and $\balpha,\bbeta\in\sS$.
\end{enumerate}
\end{theoremHigh}
The proof is similar to that of Theorem~\ref{theorem:charac_smoo} and is left as an exercise.
The theorem shows that a strongly convex function is strictly convex.
Notably, combining Theorem~\ref{theorem:charac_smoo} and Theorem~\ref{theorem:charac_stronconv}, if $f:\sS\subseteq\real^p\rightarrow \real$ is differentiable, $L_a$-SC, and $L_b$-SS defined on a convex set $\sS$ with $0<L_a<L_b$, then for any $\balpha,\bbeta\in\sS$,
\begin{subequations}
\begin{align}
\frac{L_a}{2} \normtwo{\balpha-\bbeta}^2
\;&\leq f(\bbeta) -f(\balpha) - \innerproduct{\nabla f(\balpha), \bbeta - \balpha}
&\leq&  \frac{L_b}{2} \normtwo{\balpha - \bbeta}^2; \label{equation:sssc_eq1} \\
L_a \normtwo{\balpha - \bbeta}^2 
\;&\leq \innerproduct{\nabla f(\balpha)  - \nabla f(\bbeta) , \balpha - \bbeta} 
&\leq& L_b\normtwo{\balpha-\bbeta}^2;\label{equation:sssc_eq2} \\
\frac{L_a\zeta}{2}\normtwo{\balpha-\bbeta}^2 
\;&\leq \big\{\lambda f(\balpha)+(1-\lambda)f(\bbeta)\} - f\big(\lambda \balpha+(1-\lambda)\bbeta \big) 
&\leq& \frac{L_b\zeta}{2}  \normtwo{\balpha - \bbeta}^2, \label{equation:sssc_eq3}
\end{align}
\end{subequations}
where $\zeta\triangleq \lambda(1-\lambda)$ for any $\lambda\in[0,1]$, and \eqref{equation:sssc_eq1} is equivalent to Definition~\ref{definition:scss_func}.
In some texts, the definitions of SC and SS are defined using \eqref{equation:sssc_eq3} instead since it does not require the function to be differentiable.
Again, we focus on differentiable functions in this book. 
When the function is not differentiable, the gradient in Theorem~\ref{theorem:charac_smoo} or  Theorem~\ref{theorem:charac_stronconv} can be replaced with subgradients due to the nonemptyness of subdifferential by Theorem~\ref{theorem:nonemp_relint_conv}.
Additionally, we have the following characterization theorem for SC and SS functions.
\begin{theoremHigh}[Characterization theorem of SC and SS]\label{theorem:charac_smoo_n_stronconv}
Let $f:\sS\subseteq\real^p\rightarrow \real$ be a differentiable, $L_a$-strongly convex, and $L_b$-smooth function defined on a convex set $\sS$ with $L_a\leq L_b$. Then, for any $\balpha,\bbeta\in\sS$, it follows that 
\begin{equation}\label{equation:baillon_hadded}
\innerproduct{\nabla f(\balpha) - \nabla f(\bbeta), \balpha -\bbeta}
\geq
\frac{L_a L_b}{L_a+L_b} \normtwo{\balpha-\bbeta}^2 + \frac{1}{L_a+L_b} \normtwo{\nabla f(\balpha) - \nabla f(\bbeta)}^2.
\end{equation}
When $L_a=0$, this reduces to the implication in \eqref{equ:charac_smoo_3} of Theorem~\ref{theorem:charac_smoo}.
\end{theoremHigh}
\begin{proof}[of Theorem~\ref{theorem:charac_smoo_n_stronconv}]
Let $g(\balpha) \triangleq f(\balpha) -\frac{1}{2} L_a\normtwo{\balpha}^2$. Then $\nabla g(\balpha) = \nabla f(\balpha) -L_a\balpha$.
Using (iii) in Theorem~\ref{theorem:charac_stronconv} and \eqref{equ:charac_smoo_3} in Theorem~\ref{theorem:charac_smoo}, this implies $g(\balpha)$ is $(L_b-L_a)$-smooth and 
$$
\frac{1}{L_b-L_a} \normtwo{\nabla g(\balpha) - \nabla g(\bbeta)}^2 
\leq  \innerproduct{\nabla g(\balpha) - \nabla g(\bbeta), \balpha - \bbeta},
$$
which is equivalent to the result.
\end{proof}

\begin{figure}
\centering       
\vspace{-0.25cm}                 
\subfigtopskip=0pt               
\subfigbottomskip=0pt         
\subfigcapskip=0pt      
\includegraphics[width=0.98\textwidth]{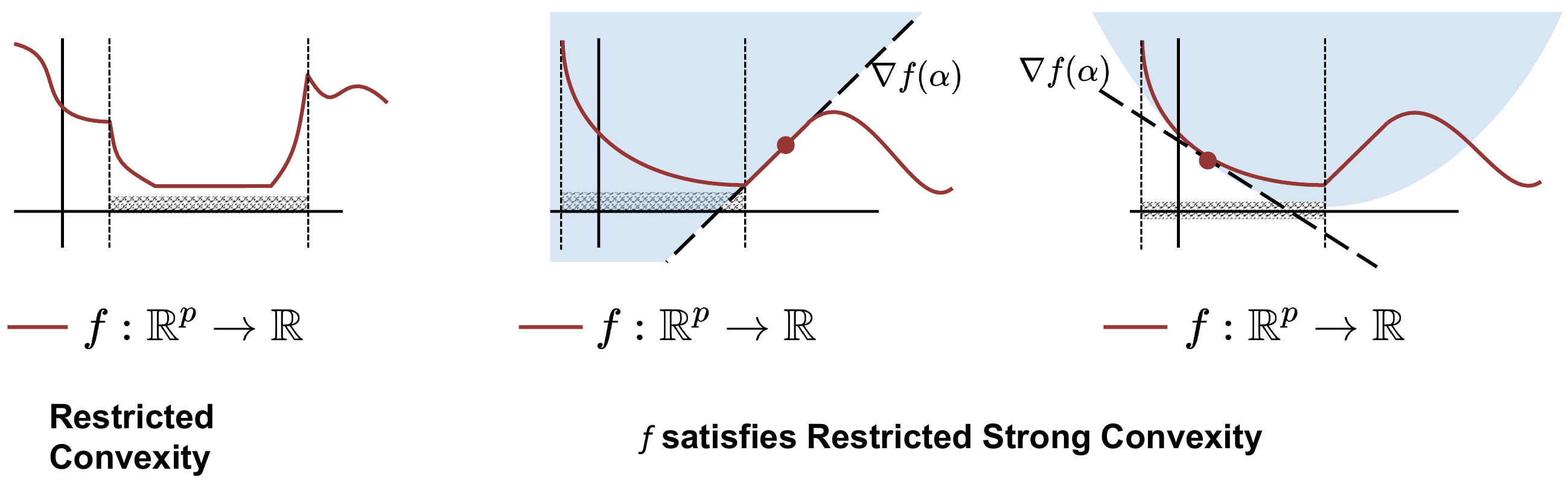}
\caption{
An illustration of restricted convexity properties. The first function $f$ exhibits non-convex behavior over the entire real line, yet it shows convex characteristics within the shaded region delineated by the dotted vertical lines. The second and third functions are also non-convex overall but demonstrate restricted strong convexity. Outside the shaded area, these functions fail to maintain convexity. However, within the defined region, they exhibit strong convexity. Figure is adapted from \citet{jain2017non}.}
\label{fig:rssfunc}
\end{figure}

In non-convex optimization, we are often interested not in global convexity or smoothness, but in restricted versions of these properties that hold only over certain subsets of the domain---even when the function is non-convex overall.
\begin{definition}[Restricted strong convexity/smoothness]\label{definition:res_scss_func}
Let $f:\sS\subseteq \real^p\rightarrow \real$ be a  differentiable function. 
Then $f$ is said to satisfy \textit{restricted convexity} over a (possibly non-convex) set $\sS\subseteq  \real^p$ if, for every $\balpha,\bbeta\in\sS$, it follows that 
\begin{equation}\label{equation:res_scss_func0}
f(\bbeta) \geq f(\balpha) +\nabla f(\balpha)^\top (\bbeta-\balpha).
\end{equation}
Additionally, $f$ is said to satisfy \textit{$L_a$-restricted strong convexity (RSC)} and \textit{$L_b$-restricted strong smoothness (RSS)} over a (possibly non-convex) set $\sS\subseteq  \real^p$ if, for every $\balpha,\bbeta\in\sS$, it follows that 
\begin{equation}\label{equation:res_scss_func1}
\frac{L_a}{2} \normtwo{\balpha-\bbeta}^2
\leq f(\bbeta)-f(\balpha)-\nabla f(\balpha)^\top (\bbeta-\balpha)
\leq \frac{L_b}{2} \normtwo{\balpha-\bbeta}^2.
\end{equation}
When the function is twice differentiable, Theorem~\ref{theorem:psd_hess_conv} also indicates  that 
\begin{equation}\label{equation:res_scss_func2}
L_a\bI \preceq \nabla^2 f(\balpha) \preceq L_b\bI,
\end{equation}
for any $\balpha\in\interior(\sS)$ and there exists a vector $\bbeta\in\real^p$ such that $\balpha+\lambda\bbeta\in\sS$ for some  sufficiently small $\lambda$ (see proof of Theorem~\ref{theorem:psd_hess_conv}).
That is, 
\begin{equation}\label{equation:res_scss_func3_ra}
L_a\leq \frac{\bv^\top \nabla^2 f(\balpha) \bv}{\normtwo{\bv}^2} \leq L_b,
\end{equation}
for any $\balpha\in\interior(\sS)$ and $\bv\in\real^p$ such that $\balpha+\lambda\bv\in\sS$ for some sufficiently small $\lambda$.
\end{definition}
Figure~\ref{fig:rssfunc} illustrates these restricted properties.
We note that when $f(\balpha) = \frac{1}{2n}\normtwo{\bX\balpha-\by}^2$, where $\bX\in\real^{n\times p}$ and $\by\in\real^n$, \eqref{equation:res_scss_func3_ra} becomes $L_a\leq \frac{\bv^\top (\frac{1}{n}\bX^\top\bX) \bv}{\normtwo{\bv}^2} \leq L_b$. 
This yields the RSC and RSS properties for sparse problems (Definition~\ref{definition:res_scss_mat}).

\section{Projection, Proximal, Bregman-Proximal, and Separation}\label{section:proj_prox_sep}

In this section, we introduce the fundamental concepts of projection and proximal operators, which are essential tools in optimization theory and convex analysis. These concepts play a central role in the projected and proximal gradient descent methods discussed in Chapter~\ref{chapter:spar}, as well as in the proofs of various theoretical results.

We begin by defining a divergence measure based on a convex function.
\begin{definition}[Bregman distance \citep{bregman1965method}\index{Bregman distance}]\label{definition:breg_dist}
The \textit{Bregman distance} (also called \textit{Bregman divergence}) generalizes the concept of squared Euclidean distance and measures the difference between two points in a space.
It is defined for a  convex and differentiable function $ \phi: \sS \rightarrow \real $, where $ \sS \subseteq \real^p $ is a convex set. The Bregman distance from point $ \bbeta $ to point $ \balpha $ with respect to $ \phi $ is given by:
\begin{equation}
\mathcalD_{\phi}(\balpha,\bbeta) = \phi(\balpha) - \phi(\bbeta) - \innerproduct{\nabla \phi(\bbeta), \balpha - \bbeta }.
\end{equation}
If $\phi$ is a non-smooth convex function, let $\bg \in \partial \phi(\bbeta)$ be a subgradient vector of the function $\phi$ at the point $\bbeta$.
The {Bregman distance} between the points $\balpha$ and $\bbeta$, denoted $D_\phi^{\bg}(\balpha, \bbeta)$, is defined as
\begin{equation}
D_\phi^{\bg}(\balpha, \bbeta) = \phi(\balpha) - \phi(\bbeta) - \innerproduct{\bg, \balpha - \bbeta}.
\end{equation}
\end{definition}

For brevity, we summarize the key properties of the Bregman distance in the following remark.
\begin{remark}[Properties of Bregman distance]\label{remark:bregnan_dist}
The Bregman distance satisfies the following properties:
\begin{enumerate}
\item \textit{Nonnegativity.}
$
\mathcalD_{\phi}(\balpha, \bbeta) \geq 0, \text{for all }  \balpha, \bbeta \in \sS.
$
When $\phi$ is strictly convex, equality holds if and only if $ \balpha = \bbeta $. 

\item The Bregman distance for the same point is zero, i.e., $D_\phi(\balpha, \balpha) = 0$.

\item \textit{Convexity in the first argument.}
For a fixed $ \bbeta $, $ \mathcalD_{\phi}(\cdot, \bbeta) $ is convex in its first argument.  
This property ensures that minimization problems involving $ \mathcalD_{\phi}(\balpha, \bbeta) $ over $ \balpha $ are well-posed---a feature frequently exploited in optimization algorithms.

\item \textit{Not symmetric.}
In general, $ \mathcalD_{\phi}(\balpha, \bbeta) \neq \mathcalD_{\phi}(\bbeta, \balpha) $. Thus, Bregman distances/divergences do not satisfy the symmetry property required for a distance metric. The term ``distance" is retained for historical reasons.

\item \textit{Closedness.} The Bregman distance can measure the closeness of two points $\bbeta$ and $\balpha$ since $D_\phi(\bbeta, \balpha) \geq D_\phi(\bxi, \balpha)$ for any point $\bxi$ on the line connecting $\bbeta$ and $\balpha$.

\item \textit{Linearity in the gradient.}
If $ \phi $ is a quadratic function, then the Bregman distance reduces  to a scaled version of the squared Euclidean distance. 
If $ \phi $ is linear, the Bregman distance becomes zero.

\item \textit{Pythagorean property.}
Let  $ \sT \subseteq \sS$ be a nonempty, closed, and convex subset of $ \sS $, and let $ \projectT^{\phi}(\bbeta) $ denote the projection of $ \bbeta $ onto $ \sT $ with respect to $ \mathcalD_{\phi} $ (Definition~\ref{definition:projec_prox_opt}). Then for any $ \balpha \in \sT $,
$$
\mathcalD_{\phi}(\balpha, \bbeta) = \mathcalD_{\phi}(\balpha, \projectT^{\phi}(\bbeta)) + \mathcalD_{\phi}(\projectT^{\phi}(\bbeta), \bbeta)
$$
This identity mirrors the classical Pythagorean theorem and is instrumental in analyzing the convergence of iterative algorithms.

\item \textit{Three-point property.}
For any three points $ \balpha, \bbeta, \bxi \in \sS $, the following relation holds:
\begin{equation}\label{equation:three_point_breg}
\mathcalD_{\phi}(\balpha, \bbeta) + \mathcalD_{\phi}(\bbeta, \bxi) 
= \mathcalD_{\phi}(\balpha, \bxi) + \innerproduct{\nabla \phi(\bxi) - \nabla \phi(\bbeta), \balpha - \bbeta}
\end{equation}
This identity is often used to derive bounds and establish algorithmic guarantees.

\item \textit{Local behavior.}
When $\balpha$ is close to $ \bbeta $, $ \mathcalD_{\phi}(\balpha, \bbeta) $ approximates  the second-order Taylor expansion of $ \phi $ around $ \bbeta $ (by the quadratic approximation theorem in Theorem~\ref{theorem:quad_app_theo}). Specifically, if $ \phi $ is twice continuously differentiable, then for small $ \normtwo{\balpha - \bbeta} $,
$$
\mathcalD_{\phi}(\balpha, \bbeta) \approx \frac{1}{2} (\balpha - \bbeta)^\top \nabla^2 \phi(\bbeta) (\balpha - \bbeta)
$$
where $ \nabla^2 \phi(\bbeta) $ denotes the Hessian  of $ \phi $ at $ \bbeta $.
\end{enumerate}
\end{remark}

Two commonly used examples of strongly convex functions are given below.
\begin{example}[Bregman with Euclidean and negative entropy]\label{example:breg_examp}
Note that when $\phi(\balpha)\triangleq\frac{1}{2}\normtwo{\balpha}^2: \real^p\rightarrow \real \cup \{ \infty \}$, $\mathcalD_{\phi}(\balpha,\bbeta)$ becomes 
$$
\mathcalD_{\phi}(\balpha,\bbeta)=\frac{1}{2}\normtwo{\balpha-\bbeta}^2,
$$
which corresponds to the squared Euclidean distance.
Let $\phi(\balpha)$ be the negative entropy function, defined as
$\phi(\balpha)\triangleq
\begin{cases}
\sum_{i=1}^{p}\alpha_i\ln (\alpha_i)& \balpha\in\real^p_+;\\
\infty&\text{otherwise},
\end{cases} 
$ 
the Bregman distance $\mathcalD_{\phi}(\balpha,\bbeta)$ becomes 
$$
\begin{aligned}
\mathcalD_{\phi}(\balpha,\bbeta)
&=
\sum_{i=1}^{p} \alpha_i \ln(\alpha_i) - \sum_{i=1}^{p} \beta_i \ln(\beta_i)
-
\sum_{i=1}^{p} (\ln(\beta_i)+1)(\alpha_i-\beta_i)\\
&=\sum_{i=1}^{p} \alpha_i\ln(\alpha_i/\beta_i) - \sum_{i=1}^{p} (\alpha_i-\beta_i).
\end{aligned}
$$
If $\balpha$ and $\bbeta$ lie on the probability simplex (i.e., $\sum_{i=1}^{p} \alpha_i= \sum_{i=1}^{p}\beta_i=1$), the linear terms cancel, and we obtain the standard \textit{Kullback--Leibler divergence}.
\end{example}

Bregman distances have broad applications across optimization, machine learning, statistics, and information theory. They are especially valuable in the design of algorithms for convex optimization because they provide a geometry-aware measure of progress toward an optimal solution---one that respects the structure induced by the convex function $\phi$. 
For instance, Bregman distances are central to \textit{mirror descent methods}, a powerful generalization of gradient descent that operates in non-Euclidean geometries; see, e.g., \citet{lu2025practical}.

Next, we introduce the definitions of the projection and proximal operators, along with a generalization known as the \textit{Bregman-proximal operator} \citep{lu2025practical}.
\begin{definition}[Projection, proximal, Bregman-proximal operators\index{Projection operators}\index{Proximal operators}]\label{definition:projec_prox_opt}
Let $\sS\subseteq \real^p$ be a nonempty \textbf{closed} set, and let $\bbeta\in\real^p$. 
The \textit{(orthogonal) projection} of $\bbeta$, denoted  $\projectS(\bbeta):\real^p\rightarrow \sS$, is defined as 
\begin{equation}\label{equation:projec_def}
\textbf{(Projection)}:\quad 
\widetildebbeta\in\projectS(\bbeta)
\triangleq
\mathop{\argmin}_{\balpha\in\sS} \normtwo{\balpha-\bbeta}
=
\mathop{\argmin}_{\balpha\in\sS} \normtwo{\balpha-\bbeta}^2.
\end{equation}
Let $f:\real^p\rightarrow \real \cup \{ \infty \}$ be  a \textbf{proper} function, and let  $\bbeta\in\real^p$. 
The \textit{proximal mapping} of $f$, denoted  $\proxf(\bbeta)$, is defined as
\begin{equation}\label{equation:prox_def}
\textbf{(Proximal)}:\quad 
\widehatbbeta \in \proxf(\bbeta) \triangleq 
\mathop{\argmin}_{\balpha\in\real^p} \left( f(\balpha) + \frac{1}{2}\normtwo{\balpha-\bbeta}^2 \right).
\end{equation}
Using the indicator function $\indicatorS$  (defined in Exercise~\ref{exercise_convex_indica}),  and assuming $\sS$ is nonempty, 
one can verify that the projection and proximal operators are related by
\begin{equation}\label{equation:equi_proj_prox}
\prox_{\indicatorS}(\bbeta) = \projectS(\bbeta),\quad \text{for all }\bbeta\in\sS.
\end{equation}
Replacing the squared Euclidean distance in \eqref{equation:prox_def} with the Bregman distance, we obtain the \textit{Bregman-proximal mapping}, denoted  $\bprox_f(\bbeta)$, as 
\begin{equation}\label{equation:bregprox_def}
\textbf{(Bregman)}:\qquad 
\widebarbbeta \in \bproxfphi(\bbeta) \triangleq 
\mathop{\argmin}_{\balpha\in\real^p} \big( f(\balpha) + \mathcalD_\phi(\balpha, \bbeta) \big),
\end{equation}
where $\phi:\real^p\rightarrow \real \cup \{ \infty \}$ is a proper \textbf{convex} function. 
By Example~\ref{example:breg_examp}, 
$$
\bproxfphi(\bbeta) \equiv \proxf(\bbeta)
\quad 
\text{if $\phi(\balpha)=\frac{1}{2}\normtwo{\balpha}^2$}.
$$
\end{definition}

Unless otherwise stated, we assume that projection operators act on nonempty closed sets and that proximal operators are applied to proper functions. The following example illustrates why these assumptions are necessary.
\begin{example}[Counterexample for non-closed $\sS$ and non-proper functions]\label{example:non_clo_proj}
Consider the open set $\sS = \{\alpha \in \real \mid \alpha > 0\}$ and let $\beta = 0$. The projection problem then becomes:
$
\widehat{\beta} = \arg\min_{\alpha > 0} \alpha^2.
$
\begin{itemize}
\item The infimum of $\alpha^2$ is $0$, but this value is not attained because $\alpha = 0\notin \sS$. Thus, there is no point in $\sS$ achieving the infimum.
\item Any sequence $\{\alpha_t\}_{t>0} \subset \sS$ with $\alpha_t \rightarrow 0$ minimizes $\alpha^2$, but the ``solution" is not unique since it depends on the sequence chosen.
\end{itemize}
Hence, the projection is not well-defined for non-closed sets. A similar issue arises for proximal operators when $f$ is not proper---for instance, if $f\equiv +\infty$, the minimization problem has no feasible solution.
\end{example}

\begin{exercise}[Proximal calculus]\label{exercise:prox_calc}
Let $g : \real^p \to (-\infty, \infty]$ be proper, and let $\eta \neq 0$.  Show that
\begin{enumerate}[(i)]
\item Let  $f(\balpha) = \eta g(\balpha/\eta)$. Then 
$$
\prox_f(\bbeta) = \eta \prox_{g/\eta}(\bbeta/\eta).
$$
\item Let $\eta \neq 0$ and $\bxi \in \real^p$. Define $f(\balpha) \triangleq g(\eta \balpha + \bxi)$. Then
$$
\prox_f(\bbeta) = \frac{1}{\eta} \left(\prox_{\eta^2 g}(\eta \bbeta + \bxi) - \bxi\right).
$$
\end{enumerate}
\end{exercise}

In the following subsections, we discuss the properties of projection and proximal operators. 
The properties of Bregman-proximal operators are discussed in Problems~\ref{prob:breg_prox_propo}$\sim$\ref{prob:breg_prox_prop1}.

\subsection{Properties of Projection Operators}
We now present some key properties of projection operators.
\begin{lemma}[Projection Property-O]\label{lemma:proj_prop0}
Let $\sS\subseteq \real^p$ be any \textbf{closed set}, and let $\bbeta\in\real^p$ such that $\widetildebbeta\triangleq\projectS(\bbeta)$ is the projection of $\bbeta$ onto the set $\sS$. Then, for all $\balpha\in\sS$, we have $\normtwobig{\widetildebbeta - \bbeta}\leq \normtwobig{\balpha-\bbeta}$.
\end{lemma}
\begin{proof}[of Lemma~\ref{lemma:proj_prop0}]
This follows directly from the definition of the projection:
$\projectS(\bbeta)=\mathop{\argmin}_{\balpha\in\sS} \normtwo{\balpha-\bbeta}^2$. 
\end{proof}

Note that Projection Property-O holds for any closed set, whether convex or not.
However, the following three properties require convexity of the set.
\begin{lemma}[Projection Property-I]\label{lemma:proj_prop1}
Let $\sS\subseteq \real^p$ be a closed \textbf{convex set}, and let $\bbeta\in\real^p$ such that $\widetildebbeta\triangleq\projectS(\bbeta)$. Then, for all $\balpha\in\sS$, we have $\innerproductbig{\balpha-\widetildebbeta, \bbeta-\widetildebbeta}\leq 0$, i.e., the angle between the two vectors is greater than or equal to 90\textdegree.
\end{lemma}
\begin{proof}[of Lemma~\ref{lemma:proj_prop1}]
Assume, for contradiction, that $\innerproductbig{\balpha-\widetildebbeta, \bbeta-\widetildebbeta}> 0$ for some $\balpha\in\sS$. 
Since $\sS$ is convex and both $\widetildebbeta, \balpha \in \sS$, for any $\lambda \in [0, 1]$, we have
$
\bxi_{\lambda} \triangleq \lambda \cdot \balpha + (1 - \lambda) \cdot \widetildebbeta \in \sS.
$
We will now show that  there exists a value of  $\lambda \in [0, 1]$ such that $\normtwo{\bbeta - \bxi_{\lambda}}^2 < \normtwobig{\bbeta - \widetildebbeta}^2$, which contradicts the fact that $\widetildebbeta$ is the closest point in the convex set to $\bbeta$ and prove the lemma. 
To see this, let $\lambda_1\triangleq\frac{2\innerproductbig{\balpha - \widetildebbeta, \bbeta - \widetildebbeta}}{\normtwo{\balpha - \widetildebbeta}^2}$ for the moment. We have 
$$
\normtwo{\bbeta-\bxi_{\lambda_1}}^2
= 
\normtwo{\bbeta-\widetildebbeta}^2 - 2\lambda_1 \innerproduct{\balpha-\widetildebbeta, \bbeta-\widetildebbeta} + \lambda_1^2 \normtwo{\balpha-\widetildebbeta}^2
=\normtwo{\bbeta-\widetildebbeta}^2.
$$
Let $g(\lambda) \triangleq \normtwobig{\bbeta-\widetildebbeta}^2 - 2\lambda \innerproductbig{\balpha-\widetildebbeta, \bbeta-\widetildebbeta} + \lambda^2 \normtwobig{\balpha-\widetildebbeta}^2$. It can be shown that $\lambda=\lambda_1$ does not obtain the minimum of $g(\lambda)$ by the property of quadratic. Indeed, the minimal value is achieved at $\lambda=\frac{\lambda_1}{2}$.
Since we assumed $\innerproductbig{\balpha - \widetildebbeta, \bbeta - \widetildebbeta} > 0$,
letting $0 < \lambda < \min\left\{1, \frac{2\innerproductbig{\balpha - \widetildebbeta, \bbeta - \widetildebbeta}}{\normtwo{\balpha - \widetildebbeta}^2}\right\}$, we have $\normtwo{\bbeta - \bxi_{\lambda}}^2 < \normtwobig{\bbeta - \widetildebbeta}^2$. This contradict the property of projection operators, and completes the proof.
\end{proof}

\begin{lemma}[Projection Property-II]\label{lemma:proj_prop2}
Let $\sS\subseteq \real^p$ be a closed  \textbf{convex set}, and let $\bbeta\in\real^p$ such that $\widetildebbeta\triangleq\projectS(\bbeta)$. Then, for all $\balpha\in\sS$, we have $\normtwobig{\widetildebbeta - \balpha}^2\leq \normtwo{\bbeta-\balpha}^2 - \normtwobig{\bbeta-\widetildebbeta}^2$, which also implies $\normtwobig{\widetildebbeta - \balpha} \leq \normtwo{\bbeta-\balpha}$ (the former inequality  can be viewed as a generalized Pythagorean theorem for projections onto convex sets.).
\end{lemma}
\begin{proof}[of Lemma~\ref{lemma:proj_prop2}]
We expand the squared norm as follows:
$$
\begin{aligned}
\normtwo{\bbeta-\balpha}^{2} &= \normtwo{(\widetildebbeta-\balpha)-(\widetildebbeta-\bbeta)}^{2} 
= \normtwo{\widetildebbeta-\balpha}^{2} + \normtwo{\widetildebbeta-\bbeta}^{2} - 2\innerproduct{\widetildebbeta-\balpha, \widetildebbeta-\bbeta} \\
&\stackrel{\dag}{\geq} \normtwo{\widetildebbeta-\balpha}^{2} + \normtwo{\widetildebbeta-\bbeta}^{2}
\geq \normtwo{\widetildebbeta-\balpha}^{2},
\end{aligned}
$$
where inequality ($\dag$) follows from the Projection Property-I.
\end{proof}

\begin{lemma}[Projection Property-III]\label{lemma:proj_prop3}
Let $\sS\subseteq \real^p$ be a closed  \textbf{convex set}, and let $\bbeta\in\real^p$. Then, $\projectS(\bbeta)
\triangleq
\mathop{\argmin}_{\balpha\in\sS} \normtwo{\balpha-\bbeta}^2$ has a \textbf{unique} optimal solution.
\end{lemma}
\begin{proof}[Lemma~\ref{lemma:proj_prop3}]
Since  the  objective function in $\projectS(\bbeta)
\triangleq
\mathop{\argmin}_{\balpha\in\sS} \normtwo{\balpha-\bbeta}^2$ is a quadratic function associated with a positive definite
matrix, it follows by Problem~\ref{prob:coerci_quad} that the objective function is coercive and hence, by Weierstrass theorem \ref{weier2_prop_close} in Theorem~\ref{theorem:weierstrass_them}, that the problem has at least one optimal solution. In addition, since the
objective function is strictly convex (again, since the objective function is quadratic associated with a
positive definite matrix), it follows by Theorem~\ref{theorem:local_glo_cvx} that there exists only one optimal
solution.
\end{proof}

\begin{figure}[h]
\centering       
\vspace{-0.25cm}                 
\subfigtopskip=2pt               
\subfigbottomskip=-2pt         
\subfigcapskip=-10pt      
\includegraphics[width=0.98\textwidth]{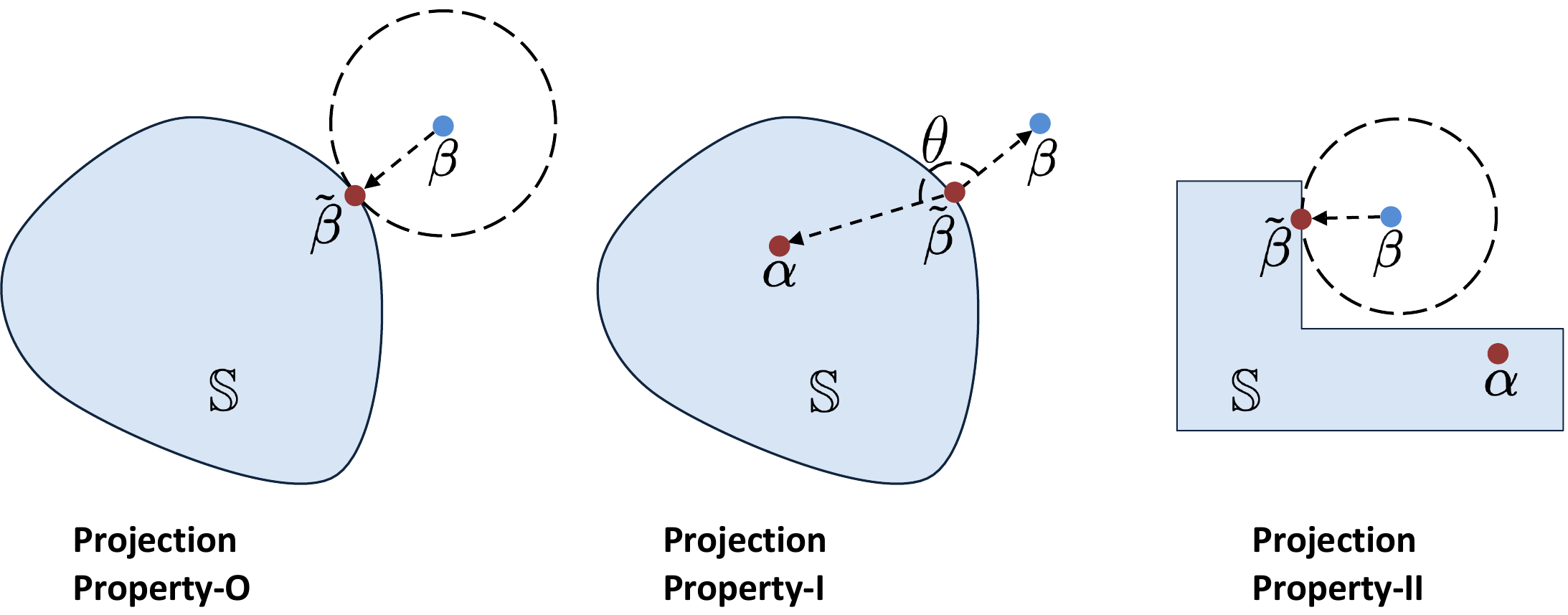}
\caption{
Projection operators map a point $\bbeta$ to the closest point $\widetildebbeta$ in a set $\sS$.
For convex sets, Projection Property-I ensures that the angle $\theta$ between the vector $\bbeta-\widetildebbeta$ (from the projection to the original point) and any vector $\balpha - \widetildebbeta$ (from the projection to another point in the set) is at least $90^\circ$.
Sets satisfying Property-I automatically satisfy Property-II, which guarantees that projecting $\bbeta$ onto $\sS$ never increases its distance to any other point $\balpha \in \sS$.
However, non-convex sets may violate  Property-II. 
}
\label{fig:proj_prop}
\end{figure}

\paragraph{Convexity matters.} 
Note that Projection Properties-I, II, and III are often referred to as  first-order properties of projections and rely crucially on the convexity of $\sS$. However, Projection Property-O, often called a zeroth-order property, holds for any nonempty closed set, convex or not; see Figure~\ref{fig:proj_prop}.

\begin{example}[Projection onto subsets of $\real^p$]\label{example:proj_subs}
Below are examples of nonempty closed and convex sets together with their orthogonal projections:
\begin{enumerate}[(i)]
\item \textbf{Nonnegative orthant $\sS_1 = \real^p_+$}: $\project_{\sS_1}(\bbeta) = \{\max\{\beta_i, 0\}\}_{i=1}^p$.
\item \textbf{Box $\sS_2 = \text{Box}[\bm{\ell}, \bu]\triangleq\{\balpha\in\real^p\mid \bm{\ell}\leq \balpha\leq \bu\}$}: $\project_{\sS_2}(\bbeta) = \left\{\min(\max(\beta_i, \ell_i), u_i)\right\}_{i=1}^p$.
\item \textbf{Closed $\ell_2$-ball $\sS_3 = \sB_2[\ba, \gamma]$}: $\project_{\sS_3}(\bbeta)=\ba + \frac{\gamma}{\max\{\normtwo{\bbeta - \ba}, \gamma\}}(\bbeta - \ba)$.
\item \textbf{Half-space $\sS_4 = \{\balpha\mid \bc^\top \balpha \leq \eta\}$}: $\project_{\sS_4}(\bbeta)=\bbeta - \frac{\max\{\bc^\top \bbeta - \eta, 0\}}{\normtwo{\bc}^2}\bc$.
\item \textbf{Affine set $\sS_5 = \{\balpha \in \real^p \mid \bX\balpha = \by\}$}: $\project_{\sS_5}(\bbeta) = \bbeta - \bX^\top(\bX\bX^\top)^{-1}(\bX\bbeta - \by)$.
\end{enumerate}
where $\bm{\ell} \in  [-\infty, \infty)^p, \bu \in (-\infty, \infty]^p$ satisfying $\bm{\ell} \leq \bu$, $\ba \in \real^p$, $\gamma > 0$, $\bc \in \real^p \setminus \{\bzero\}$,  $\eta \in \real$, $\bX \in \real^{n\times p}$ has full row rank such that $\bX\bX^\top$ is nonsingular, and $\by \in \real^n$.
The results in (i)--(iv) holds trivially. 
For (v),
it's equivalent to solving  the following optimization problem:
$$
\min_{\balpha \in \real^p} \quad \frac{1}{2} \normtwo{ \balpha - \bbeta}^2 \quad\text{s.t.}\quad\bX\balpha = \by,
$$
where $\bX \in \real^{n\times p}$, $\by \in \real^n$ and $\bbeta \in \real^p$ are given matrices and vectors. 
We assume that the matrix $\bX$ has full row rank (for matrices that do not have full row rank, we can remove linearly dependent rows to obtain a reduced constraint condition). 
For the equality constraint, we introduce the Lagrange multiplier $\blambda \in \real^n$ and construct the Lagrangian function (see Section~\ref{section:opt_conds}):
$$
L(\balpha, \blambda) = \frac{1}{2} \normtwo{\balpha - \bbeta}^2 + \blambda^\top(\bX\balpha - \by).
$$
Since the objective is convex and the constraints are affine, Slater's condition holds trivially. 
Theorem~\ref{theorem:opt_cond_sd} shows that $\balpha^*$ is a global optimum if and only if there exists $\blambda^* \in \real^n$ such that
$
\begin{cases}
\balpha^* - \bbeta + \bX^\top\blambda^* = \bzero, \\
\bX\balpha^* = \by.
\end{cases}
$
This yields that 
$$
\bX \balpha^* - \bX \bbeta + \bX\bX^\top\blambda^* = \bzero
\qquad \implies\qquad 
\blambda^* = (\bX\bX^\top)^{-1}(\bX\bbeta - \by).
$$
Plugging $\blambda^*$ back yields the claimed formula.
\end{example}

\subsection{Separation Theorem}

Using properties of projections, we can prove the following strict separation theorem, which is essential for establishing Farkas' lemma (Problem~\ref{prob:farka_lemm}). Farkas' lemma, in turn, plays a fundamental role in the derivation of the KKT conditions; see, for example, \citet{lu2025practical} for further details. Additional separation results are discussed in Problems~\ref{prob:stric_sep_theo2}$\sim$\ref{prob:supp_hyp_theo}.
\begin{theoremHigh}[Strict separation theorem]\label{theorem:stric_sep_theo}
Let $\sS \subseteq \real^p$ be a closed \textbf{convex} set, and let $\bbeta \notin \sS$. Then there exist $\bw \in \real^p \setminus \{\bzero\}$ and $\gamma \in \real$ such that
$$
\bw^\top \bbeta > \gamma
\qquad\text{and}\qquad
\bw^\top \balpha \leq \gamma, \quad \text{for all } \balpha \in \sS.
$$
\end{theoremHigh}
\begin{proof}[of Theorem~\ref{theorem:stric_sep_theo}]
Let $\widetildebbeta \triangleq \projectS(\bbeta) \in \sS$. By the Projection Property-I (Lemma~\ref{lemma:proj_prop1}), 
$$
\innerproduct{\balpha - \widetildebbeta, \bbeta - \widetildebbeta} \leq 0
\quad\implies\quad (\bbeta - \widetildebbeta)^\top \balpha \leq (\bbeta - \widetildebbeta)^\top \widetildebbeta
,\quad  \text{for all } \balpha \in \sS.
$$
Denote $\bw \triangleq \bbeta - \widetildebbeta \neq \bzero$ (since $\bbeta \notin \sS$) and $\gamma \triangleq (\bbeta - \widetildebbeta)^\top \widetildebbeta$. Then we have that $\bw^\top \balpha \leq \gamma$ for all $\balpha \in \sS$. On the other hand,
$$
\bw^\top \bbeta = (\bbeta - \widetildebbeta)^\top \bbeta = (\bbeta - \widetildebbeta)^\top (\bbeta - \widetildebbeta) + (\bbeta - \widetildebbeta)^\top \widetildebbeta =\normtwo{\bbeta - \widetildebbeta}^2 + \gamma > \gamma,
$$
establishing the desired result.
\end{proof}

\subsection{Properties of Proximal Operators}\label{section:prop_proxoper}
Analogously to projection operators, for the proximal operator, if the function is closed and satisfies a mild coerciveness condition (i.e., $ f(\balpha) \to \infty $ as $ \normtwo{\balpha} \to \infty $) then the proximal map $\proxf(\bbeta)$ is guaranteed to be nonempty for every $\bbeta\in\real^p$.
\begin{lemma}[Proximal Property-O: Nonemptiness of  prox under closedness and coerciveness]\label{lemma:prox_prop0}
Let $f: \real^p \rightarrow (-\infty, \infty]$ be a proper \textbf{closed} function, and assume that the following condition holds:
$$
\text{the function } \balpha \mapsto f(\balpha) + \frac{1}{2} \normtwo{\balpha - \bbeta}^2 \text{ is coercive for any } \bbeta \in \real^p. 
$$
Then $\proxf(\bbeta)$ is \textbf{nonempty} for any $\bbeta \in \real^p$.
\end{lemma}
\begin{proof}[of Lemma~\ref{lemma:prox_prop0}]
For any $\balpha \in \real$, the proper function $g(\balpha) \triangleq f(\balpha) + \frac{1}{2} \normtwo{\balpha - \bbeta}^2$ is closed as a sum of two closed functions. 
Therefore, by the standard existence result for minimizers of proper, closed, and coercive functions (e.g., Weierstrass theorem \ref{weier2_prop_close} in Theorem~\ref{theorem:weierstrass_them}),
the set of minimizers of $g$, which is precisely $\proxf(\bbeta)$, is nonempty.
\end{proof}

\begin{lemma}[Proximal Property-I]\label{lemma:prox_prop1}
Let $ f: \real^p \rightarrow (-\infty, \infty] $ be a proper \textbf{closed and convex} function. Then for any $ \bbeta, \widehatbbeta \in \real^p $, the following three statements  are equivalent:
\begin{enumerate}[(i)]
\item $ \widehatbbeta = \proxf(\bbeta) $.
\item $ \bbeta - \widehatbbeta \in \partial f(\widehatbbeta) $.
\item $ \innerproduct{\bbeta - \widehatbbeta, \bxi - \widehatbbeta} \leq f(\bxi) - f(\widehatbbeta) $ for any $ \bxi \in \real^p $.
\end{enumerate}
\end{lemma}
\begin{proof}[of Lemma~\ref{lemma:prox_prop1}]
By definition, $ \widehatbbeta = \proxf(\bbeta) $ if and only if $ \widehatbbeta $ is the minimizer of the problem
$
\mathopmin{\bv} \left\{ f(\bv) + \frac{1}{2} \normtwo{\bv - \bbeta}^2 \right\},
$
which, by the optimality condition (Theorem~\ref{theorem:fetmat_opt}) and the sum rule of subdifferential calculus, is equivalent to the relation
$
\bzero \in \partial f(\widehatbbeta) + \widehatbbeta - \bbeta. 
$
This establishes the equivalence between (i) and (ii). 
Finally, by the definition of the subgradient (Definition~\ref{definition:subgrad}), the membership relation of claim (ii) is equivalent to (iii).
\end{proof}

When $ f = \indicatorS $, with $ \sS $ being a nonempty closed and convex set, the equivalence between claims (i) and (iii) in the Proximal Property-I reduces exactly   to the Projection Property-I (Lemma~\ref{lemma:proj_prop1}).

For consistency with the projection properties, we also state Proximal Property-II, although it is less frequently used in practice.
\begin{lemma}[Proximal Property-II]\label{lemma:prox_prop2}
Let $ f: \real^p \rightarrow (-\infty, \infty] $ be a proper \textbf{closed and convex} function. Then for any $ \widehatbbeta \triangleq \proxf(\bbeta) $ and $\bxi\in\real^p$, it follows that
$
\normtwo{\bbeta-\bxi}^2
\geq \normtwobig{\widehatbbeta-\bxi}^2 + \normtwobig{\widehatbbeta-\bbeta}^2 + 2\left( f(\widehatbbeta) - f(\bxi) \right).
$
\end{lemma}

\begin{lemma}[Proximal Property-III]\label{lemma:prox_prop3}
Let $f: \real^p \rightarrow (-\infty, \infty]$ be a proper \textbf{closed and convex} function. Then $\proxf(\bbeta)$ is a \textbf{singleton} for any $\bbeta \in \real^p$.
\end{lemma}

\begin{proof}[of Lemma~\ref{lemma:prox_prop3}]
For any $\balpha \in \real^p$,
\begin{equation}\label{equation:prox_prop0_e1}
\proxf(\bbeta) = \arg\min_{\balpha \in \real^p} g(\bbeta;\balpha) ,
\end{equation}
where $g(\bbeta;\balpha) = f(\balpha) + \frac{1}{2}\normtwo{\balpha - \bbeta}^2$. The function $g(\bbeta;\balpha)$ is a closed and strongly convex function as a sum of the closed and strongly convex function $\frac{1}{2}\normtwo{\balpha - \bbeta}^2$ and the closed and convex function $f$ (Exercise~\ref{exercise:sum_sc_conv} and Exercise~\ref{exercise:pres_conv_clos}). The properness of $g(\bbeta;\balpha)$ immediately follows from the properness of $f$. Therefore, by the SC Property-II (Theorem~\ref{theorem:exi_close_sc}, or Theorem~\ref{theorem:local_glo_cvx} for the uniqueness under strict convexity), there exists a unique minimizer to the problem in \eqref{equation:prox_prop0_e1}.
\end{proof}

The Proximal Property-III also proves the Projection Property-III (Lemma~\ref{lemma:proj_prop3}) by noting that $\prox_{\indicatorS}(\bbeta) = \projectS(\bbeta)$ for any $\bbeta\in\sS$ if $\sS$ is convex and closed (the indicator function is closed convex due to Exercises~\ref{exercise_convex_indica} and \ref{exercise_closed_indica}).

When $f$ is proper closed and convex, the preceding lemma shows that $\proxf(\bbeta)$ is a singleton for any $\bbeta \in \real^p$. In these cases, which will constitute the vast majority of cases that will be discussed in this book, we will treat $\proxf$ as a single-valued mapping from $\real^p$ to $\real^p$, meaning that we will simply write $\proxf(\bbeta) = \widehatbbeta$ instead of $\proxf(\bbeta) = \{\widehatbbeta\}$.

\paragraph{Convexity and Closedness.} Note again  that Proximal Properties-I, II, and III are often referred to as first-order properties and can be violated if the underlying function is non-convex. However, Projection Property-O, often called a zeroth-order property, always holds, whether the underlying function is convex or not.
In all cases, we require the function to be closed to ensure the attainment due to the Weierstrass theorem (Theorem~\ref{theorem:weierstrass_them}).

\begin{theoremHigh}[Projection and Proximal Property-IV: Nonexpansiveness]\label{theorem:proj_nonexpan}
Let $\sS\subseteq \real^p$ be a closed \textbf{convex set}, and let $f:\real^p\rightarrow (-\infty, \infty]$ be a proper \textbf{closed and convex} function. Then, 
\begin{enumerate}
\item \textit{Firm nonexpansiveness.} For any $\balpha, \bbeta\in\real^p$, it follows that 
$$
\begin{aligned}
\normtwo{\projectS(\balpha)-\projectS(\bbeta)}^2 &\leq \innerproduct{\projectS(\balpha)-\projectS(\bbeta), \balpha-\bbeta};\\
\normtwo{\proxf(\balpha)-\proxf(\bbeta)}^2 &\leq \innerproduct{ \proxf(\balpha)-\proxf(\bbeta), \balpha-\bbeta}.
\end{aligned}
$$
\item \textit{Nonexpansiveness.} For any $\balpha, \bbeta\in\real^p$, it follows that 
$$
\begin{aligned}
\normtwo{\projectS(\balpha)-\projectS(\bbeta)} &\leq \normtwo{\balpha-\bbeta};\\
\normtwo{\proxf(\balpha)-\proxf(\bbeta)}&\leq \normtwo{\balpha-\bbeta}.
\end{aligned}
$$
\end{enumerate} 
\end{theoremHigh}
\begin{proof}[of Theorem~\ref{theorem:proj_nonexpan}]
\textbf{Projection operator.} By the Projection Property-I (Lemma~\ref{lemma:proj_prop1}),  for any $\bu \in \real^p$ and $\bv \in \sS$, we have that
$
\innerproduct{ \bv - \projectS(\bu), \bu - \projectS(\bu)} \leq 0. 
$
Invoking this inequality with $\bu \triangleq \balpha, \bv \triangleq \projectS(\bbeta)$ and $\bu \triangleq \bbeta, \bv \triangleq \projectS(\balpha)$:
$$
\begin{aligned}
\innerproduct{\projectS(\bbeta) - \projectS(\balpha), \balpha - \projectS(\balpha)} &\leq 0 
\qquad\text{and}\qquad
\innerproduct{\projectS(\balpha) - \projectS(\bbeta), \bbeta - \projectS(\bbeta)} &\leq 0. 
\end{aligned}
$$
Adding the two inequalities yields that
$
\innerproduct{\projectS(\bbeta) - \projectS(\balpha),  \balpha - \bbeta + \projectS(\bbeta) - \projectS(\balpha)} \leq 0,
$
showing the firm nonexpansiveness of projection operators.

To prove the nonexpansiveness, note that if $\projectS(\balpha) = \projectS(\bbeta)$, the inequality  holds trivially. We will therefore assume that $\projectS(\balpha) \neq \projectS(\bbeta)$. By the Cauchy--Schwarz inequality we have
$$
\innerproduct{\projectS(\balpha) - \projectS(\bbeta),  \balpha - \bbeta} \leq \normtwo{\projectS(\balpha) - \projectS(\bbeta)} \cdot \normtwo{\balpha - \bbeta},
$$
which combined with the firm nonexpansivess yields the desired result.

\paragraph{Proximal operator.} Denoting $ \bu \triangleq \proxf(\balpha) $, $ \bv \triangleq \proxf(\bbeta) $, by the equivalence of (i) and (ii) in the Proximal Property-I (Lemma~\ref{lemma:prox_prop1}), it follows that
$$
\balpha - \bu \in \partial f(\bu)
\qquad \text{and}\qquad  
\bbeta - \bv \in \partial f(\bv).
$$
By the subgradient inequality, we have
$$
f(\bv) \geq f(\bu) + \innerproduct{\balpha - \bu, \bv - \bu}
\qquad \text{and}\qquad  
f(\bu) \geq f(\bv) + \innerproduct{\bbeta - \bv, \bu - \bv}.
$$
Adding the two inequalities yields that
$
0 \geq \innerproduct{\bbeta - \balpha + \bu - \bv, \bu - \bv},
$
showing the firm nonexpansiveness of proximal operators.

To prove the nonexpansiveness, if $ \proxf(\balpha) = \proxf(\bbeta) $, then the inequality is trivial. Therefore, we assume that $ \proxf(\balpha) \neq \proxf(\bbeta) $. By the Cauchy--Schwarz inequality, it follows that
$$
\innerproduct{\proxf(\balpha) - \proxf(\bbeta), \balpha - \bbeta} \leq \normtwo{\proxf(\balpha) - \proxf(\bbeta)} \cdot \normtwo{\balpha - \bbeta},
$$
which combined with the firm nonexpansivess yields the desired result.
\end{proof}

Finally, note that the (firm) nonexpansiveness of the projection operator is a special case of the proximal operator result. Indeed, setting  $f\triangleq\indicatorS$ (the indicator function of $\sS$), we have $\prox_f = \projectS$, and since $\indicatorS$ is proper, closed, and convex whenever $\sS$ is nonempty, closed, and convex (Exercises~\ref{exercise_convex_indica} and \ref{exercise_closed_indica}), the proximal nonexpansiveness theorem directly implies the projection version.

\begin{SCfigure}
\centering
\includegraphics[width=0.55\textwidth]{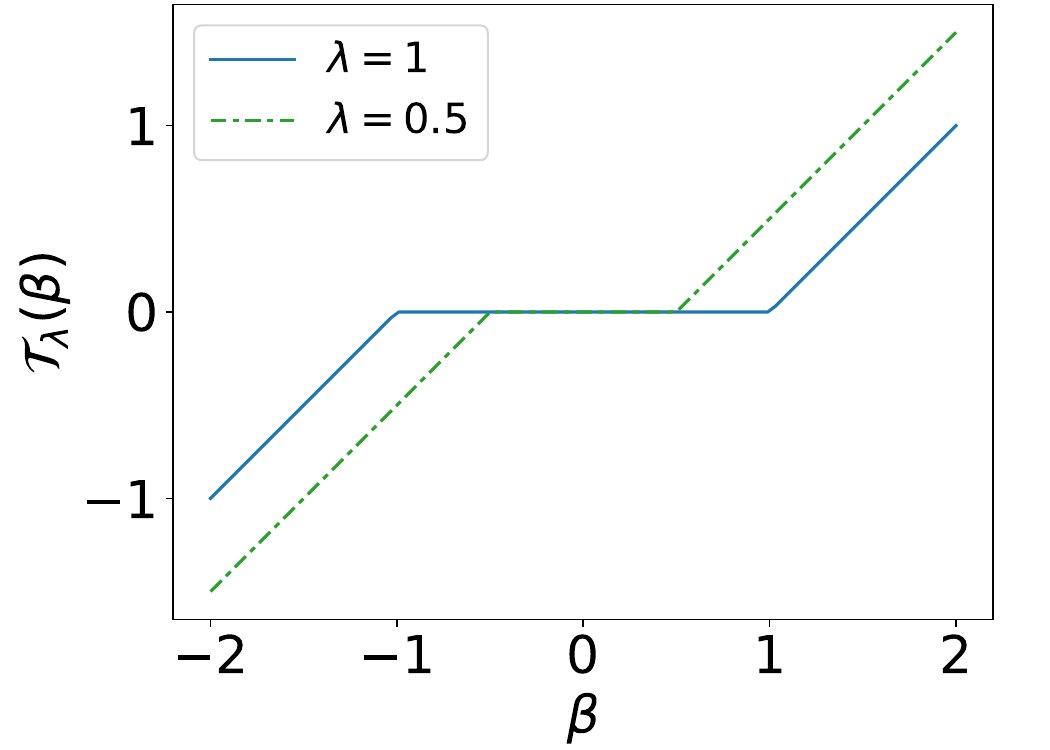} 
\caption{Illustration of the soft-thresholding function $\mathcalT_{\lambda}(\beta)$ for different values of $\lambda$.}
\label{fig:soft_threshold}
\end{SCfigure}

\begin{example}[Soft-thresholding, proximal of $\ell_1$-norms]\label{example:soft_thres}
Let  $f(\balpha) = \lambda \normone{\balpha} =\sum_{i=1}^p \psi(\alpha_i)$, where $\balpha\in\real^p$, $\lambda > 0$, and $\psi(u) = \lambda \abs{u}$. 
Then $
\proxf(\bbeta) = \big\{\prox_{\psi}(\beta_i)\big\}_{i=1}^p,
$
where $\prox_{\psi}(\beta) = \mathcalT_{\lambda}(\beta)$, and  $\mathcalT_{\lambda}$ is known as the \textit{soft-thresholding function} and is defined as
\begin{equation}
\mathcalT_{\lambda}(\beta) \triangleq
[\abs{\beta} - \lambda]_+ \cdot \sign(\beta)=
\begin{cases}
\beta - \lambda, & \beta \geq \lambda; \\
0, & \abs{\beta} < \lambda; \\
\beta + \lambda, & \beta \leq -\lambda.
\end{cases}
\end{equation}
where $[u]_+ = \max\{u, 0\}$ for any $u\in\real$, and the sign function is given by $\sign(u) = u / \abs{u}$ for $u \neq 0$, as usual. 
See Figure~\ref{fig:soft_threshold} for an illustration of this soft-thresholding function.
More compactly, the proximal operation of $f$ can be denoted as 
$$
\prox_f(\bbeta) = \mathcalT_\lambda(\bbeta) = \left\{\mathcalT_\lambda(\beta_i)\right\}_{i=1}^p
= [\abs{\bbeta} - \lambda \bone]_+ \hadaprod \sign(\bbeta),
$$
where $\hadaprod$ denotes the Hadamard product, and $\bone$ is the all-ones vector.
\end{example}

\begin{example}[Block soft-thresholding]\label{example:block_soft_thres}
The {proximal operator} of the function $ f(\balpha) = \lambda \normtwo{\balpha}$ with $\balpha\in\real^p$, the scaled $\ell_2$-norm, is known as the \textit{block soft-thresholding function}.
This is  often used in group sparsity (see Section~\ref{section:group_las}).
The proximal operator has a closed-form solution:
\begin{equation}
\prox_{\lambda \normtwo{\cdot}}(\bbeta) 
= \left[1 - \frac{\lambda}{\normtwo{\bbeta}}\right]_+ \bbeta
=
\begin{cases}
\displaystyle \left(1 - \frac{\lambda}{\normtwo{\bbeta}}\right) \bbeta, & \text{if } \normtwo{\bbeta} > \lambda; \\
\bzero, & \text{if } \normtwo{\bbeta} \leq \lambda,
\end{cases}
\end{equation}
where $ [u]_+ = \max(u, 0) $.
To see this, we observe that the objective is radially symmetric, so the minimizer $ \balpha^* $ must be colinear with $ \bbeta $: $ \balpha^* = \gamma \bbeta $ for some $ \gamma \geq 0 $.
Substitute into the objective:
$$
\lambda \normtwo{\gamma \bbeta} + \frac{1}{2} \normtwo{\gamma \bbeta - \bbeta}^2 
= \lambda \gamma \normtwo{\bbeta} + \frac{1}{2} (\gamma - 1)^2 \normtwo{\bbeta}^2.
$$
Minimize over $ \gamma \geq 0 $, we differentiate with respect to $ \gamma $:
$$
\lambda \normtwo{\bbeta} + (\gamma - 1)\normtwo{\bbeta}^2 = 0 
\quad\implies\quad 
\gamma = 1 - \frac{\lambda}{\normtwo{\bbeta}}.
$$
However, since $ \gamma \geq 0 $, this solution is valid only when $\normtwo{\bbeta}>\lambda$. 
If $ \normtwo{\bbeta} \leq \lambda $, the minimum occurs at $ \gamma = 0 $, yielding $\balpha^*=\bzero$.

As a special case, when the dimension $ p = 1 $, then $ \normtwo{\beta} = \abs{\beta} $ and
$$
\prox_{\lambda \abs{\cdot}}(\beta) = \sign(\beta) \max(\abs{\beta} - \lambda, 0),
$$
which reduces to the scalar {soft-thresholding operator} from Example~\ref{example:soft_thres}.
\end{example}

\subsection{Moreau Decomposition}
A key property of the proximal operator is the so-called \textit{Moreau decomposition} theorem,
which establishes a relationship between the proximal operator of a proper closed convex function and that of its convex conjugate (Definition~\ref{definition:conjug_func}).
\begin{theoremHigh}[Moreau decomposition]\label{theorem:moreau_iden}
Let $f: \real^p \to (-\infty, \infty]$ be a proper closed and convex function. Then, for all $\bbeta \in \real^p$,
\begin{equation}\label{equation:moreau}
\prox_f(\bbeta) + \prox_{f^*}(\bbeta) = \bbeta.
\end{equation}
\end{theoremHigh}

\begin{proof}[of Theorem~\ref{theorem:moreau_iden}]
Let $\widehatbbeta \triangleq \prox_f(\bbeta)$. Then by the equivalence between claims (i) and (ii) in the Proximal Property-I (Lemma~\ref{lemma:prox_prop1}), it follows that $\bbeta - \widehatbbeta \in \partial f(\widehatbbeta)$, which by the conjugate subgradient theorem (Theorem~\ref{theorem:conju_subgra}) is equivalent to $\widehatbbeta \in \partial f^*(\bbeta - \widehatbbeta)$. 
Using the Proximal Property-I (Lemma~\ref{lemma:prox_prop1}) again, we conclude that $\bbeta - \widehatbbeta = \prox_{f^*}(\bbeta)$. Therefore,
$$
\prox_f(\bbeta) + \prox_{f^*}(\bbeta) = \widehatbbeta + (\bbeta - \widehatbbeta) = \bbeta. 
$$
This completes the proof.
\end{proof}

The next result provides a useful generalization of the Moreau decomposition theorem.
\begin{theoremHigh}[Extended Moreau decomposition]\label{theorem:ext_moreau_iden}
Let $f: \real^p \to (-\infty, \infty]$ be a {proper closed and convex} function, and let $\eta > 0$. Then, for all $\bbeta \in \real^p$,
\begin{equation}\label{equation:extended_moreau}
\prox_{\eta f}(\bbeta) + \eta \prox_{\eta^{-1} f^*}(\bbeta/\eta) = \bbeta.
\end{equation}
\end{theoremHigh}
\begin{proof}[of Theorem~\ref{theorem:ext_moreau_iden}]
Using Moreau decomposition, for any $\bbeta \in \real^p$,
\begin{equation}\label{equation:ext_moreau_iden_pv1}
\prox_{\eta f}(\bbeta) = \bbeta - \prox_{(\eta f)^*}(\bbeta) = \bbeta - \prox_{\eta f^*(\cdot/\eta)}(\bbeta),
\end{equation}
where the second equality follows from conjugate calculus (Exercise~\ref{exercise:conj_calc}). 
By proximal calculus (Exercise~\ref{exercise:prox_calc}),
$$
\prox_{\eta f^*(\cdot/\eta)}(\bbeta) = \eta \prox_{\eta^{-1} f^*}(\bbeta/\eta),
$$
which, combined with \eqref{equation:ext_moreau_iden_pv1}, yields the desired identity.
\end{proof}

\section{Convex Optimization and Optimaliy Conditions}\label{section:opt_conds}

We have already encountered the general form of optimization problems in previous chapters.
More precisely, a \textit{convex optimization problem} (or simply a convex problem) involves minimizing a convex function over a convex set:
\begin{equation}\label{equation:convex_optim1} 
\textbf{(Convex optimization)}:\qquad
\begin{aligned}
&\min \quad f(\btheta) \quad &&\text{(convex function)}\\
&\text{s.t.}\quad \btheta\in\sS \quad &&\text{(convex set)},
\end{aligned}
\end{equation}
where the function $f: \real^p \to (-\infty, \infty]$ is called the \textit{objective function} (or \textit{loss function}), 
and $\sS$ is  the \textit{constraint set}.
A point $\btheta$ that satisfies the constraints is called \textit{feasible}.
The  problem~\eqref{equation:convex_optim1}   is said to be   feasible if at least one feasible point exists. 
A feasible point $\btheta^*$ that achieves the minimum---i.e., $f(\btheta^*) \leq f(\btheta)$ for all feasible $\btheta$---is called a \textit{minimizer}, \textit{minimum point}, or \textit{optimal point}. 
And the value $f(\btheta^*)$ is referred to as the \textit{optimal value}.

%

\begin{example}[Equality and inequality constraints]\label{example:equa_ineq_cons}
Consider the following problem:
\begin{equation}\label{equation:conv_eqineq_const}
\begin{aligned}
\min & \quad f(\btheta) \\
\text{s.t.} & \quad g_i(\btheta) \leq 0, \quad i = \{1, 2, \ldots, m\},\\
&\quad h_j(\btheta) = 0, \quad j = \{1,2, \ldots, n\},
\end{aligned}
\end{equation}
where $f, g_1, g_2, \ldots,g_m:\real^p\rightarrow \real$ are convex functions, and $h_1, h_2, \ldots,h_n: \real^p\rightarrow \real$ are affine functions.
Since the objective function is convex and the feasible set is a convex set, this problem is a convex optimization problem.
The feasible set can be expressed as
$$
\sS = \bigg( \bigcap_{i=1}^m \lev[g_i, 0] \bigg) \cap \bigg( \bigcap_{j=1}^n \{ \btheta \mid  h_j(\btheta) = 0 \} \bigg),
$$
which is a convex set since it is the intersection of level sets of convex functions and hyperplanes, both of which are also convex sets.
\end{example}

\begin{example}[Linear programming (LP)]\label{example:linear_program}
A \textit{linear programming (LP) problem}  involves minimizing a linear objective function subject to linear equality and inequality constraints:
$$
\text{(LP)}\qquad 
\begin{aligned}
\min & \quad \bc^\top \btheta +d \\
\text{s.t.}
&\quad \bG\btheta \leq \bh, \\
&\quad \bA\btheta = \bb,
\end{aligned}
$$
where $ \bG \in \real^{n \times p} $, $ \bh \in \real^n $, $ \bA \in \real^{m \times p} $, $ \bb \in \real^m $, and $ \bc \in \real^p $. (For convenience, we write $\ba \leq  \bzero$ if $a_i \leq 0$ for all $i$.)
This is a convex optimization problem because affine functions are inherently convex. The constant term $d$ in the objective does not affect the location of the minimizer and can be omitted when identifying optimal solutions.
A commonly used formulation is the ``standard form" of an LP:
$$
\text{(Standard LP)}\qquad 
\begin{aligned}
\min & \quad \bc^\top \btheta \\
\text{s.t.} & \quad \bA\btheta = \bb, \\
& \quad \btheta \geq \bzero.
\end{aligned}
$$
Any LP can be converted into this standard form. To do so, we introduce  nonnegative slack variables $\bs$ to turn inequalities into equalities: $\bG\btheta +\bs=  \bh$. 
Additionally, we   express $\btheta$ as the difference of two nonnegative variables $\btheta\triangleq\btheta^+-\btheta^-$, where $\btheta^+, \btheta^-\geq \bzero$:
$$
\text{(LP$'$)}\qquad 
\begin{aligned}
\min & \quad \bc^\top \btheta^+ - \bc^\top\btheta^- +d \\
\text{s.t.}
&\quad \bG\btheta^+ - \bG\btheta^-  +\bs=  \bh,   \\
&\quad \bA\btheta^+ - \bA\btheta^-  = \bb, \\
&\quad \bs \geq \bzero, \btheta^+ \geq \bzero, \btheta^-\geq \bzero, 
\end{aligned}
\quad\iff \quad
\begin{aligned}
\min & \quad [\bc^\top , - \bc^\top, \bzero] \widetildebtheta +d \\
\text{s.t.}
&\quad  [\bG, -\bG, \bI] \widetildebtheta =\bh,   \\
&\quad [\bA, - \bA, \bzero ] \widetildebtheta  = \bb, \\
&\quad \widetildebtheta \geq \bzero, 
\end{aligned}
$$
where 
$
\widetildebtheta\triangleq
\scriptsize
\begin{bmatrix}
\btheta^+ \\  \btheta^- \\\bs 
\end{bmatrix}
$. 
This reformulation is in standard LP form with decision variable $\widetildebtheta$. For each feasible solution $\btheta^*$ of (LP), we can set $\bs\triangleq\bh-\bG\btheta^*$, $\theta_i^+\triangleq \max\{0, \theta_i\}$, and $\theta_i^-\triangleq \max\{0, -\theta_i\}$ for all $i\in\{1,2,\ldots, p\}$, which is a feasible solution of (LP$'$). This shows the optimal value of (LP$'$) is less than or equal to the optimal value of (LP). 
Conversely, suppose $\btheta^+, \btheta^-$, and $\bs$ are feasible solutions of (LP$'$). Then, $\btheta\triangleq\btheta^+-\btheta^-$ is a feasible solution of (LP), demonstrating that the optimal value of (LP) is less than or equal to the optimal value of (LP$'$).  Therefore, the optimization problems (LP) and (LP$'$) are equivalent.
\end{example}

\begin{example}[Quadratic Programming]
The convex optimization problem \eqref{equation:conv_eqineq_const} is called a \textit{quadratic program (QP)} if the objective function is quadratic with positive semidefinite (convex), and the constraint functions  are affine. A general form can be written as:
$$
\begin{aligned}
\text{min}& \quad \btheta^\top \bP \btheta + 2 \bq^\top \btheta + r \\
\text{s.t.} & \quad \bG\btheta \leq \bh, \\
&\quad \bA\btheta=\bb,
\end{aligned}
$$
where $\bP \in \real^{p \times p}$ is positive semidefinite, $\bq \in \real^p$, $\bG\in\real^{n\times p}$, $\bh\in\real^n$, $\bA \in \real^{m \times p}$, and $\bb \in \real^m$. 
If both the objective function and inequality constraints are convex quadratic:
$$
\begin{aligned}
\text{min}& \quad \btheta^\top \bP \btheta + 2 \bq^\top \btheta + r \\
\text{s.t.} & \quad \btheta^\top\bG_i\btheta + 2\bh_i^\top\btheta+ c_i \leq 0,\quad i=\{1,2,\ldots,k\},\\
&\quad \bA\btheta=\bb,
\end{aligned}
$$
where $\bP, \bG_1, \bG_2, \ldots,\bG_k$ are positive semidefinite, the problem is called a \textit{quadratically constrained quadratic program (QCQP)}.
In a QCQP, we minimize a convex quadratic function over a feasible region defined as the intersection of ellipsoids (when each $\bG_i\succ \bzero$) and affine subspaces.
\end{example}

\subsection{Constrained Optimization over Convex Sets}\label{section:constr_convset}

We begin by restricting our attention to optimization problems defined over a \textbf{convex} set $\sS$. 
In this setting, if the objective function $f$ is also convex, the problem becomes a {convex optimization problem}, as previously discussed.
Earlier, in Corollary~\ref{corollary:fermat_fist_opt} or Theorem~\ref{theorem:fetmat_opt},  we introduced \textit{stationary points} for unconstrained problems---specifically, points $\bbeta^*$ where gradient vanishes $\nabla f(\bbeta^*)=\bzero$, or $\bzero \in\partial f(\bbeta^*)$ when the function is not smooth.
However, for constrained problems, this definition no longer suffices, and a more general notion is required.
\begin{definition}[Stationary points of constrained problems on a convex set]\label{definition:stat_point_uncons_convset}
Consider the constrained optimization problem 
\begin{equation}
\text{(P$_C$)}\quad  \min_{\bbeta}f(\bbeta) \quad \text{s.t.}\quad \bbeta\in\sS.
\end{equation}
Let $f$ be a continuously differentiable function over a closed \textbf{convex} set $\sS$. Then $\bbeta^* \in \sS$ is called a stationary point of (P$_C$) if $\nabla f(\bbeta^*)^\top (\bbeta - \bbeta^*) \geq 0$ for any $\bbeta \in \sS$.
\end{definition}


The following result shows that a local minimum point of a  constraint optimization problem must be a stationary point.

\begin{theoremHigh}[First-order necessary condition for (P$_C$) over a convex set]\label{theorem:stat_point_uncons_convset}
Consider the constrained optimization problem (P$_C$): $\min_{\bbeta}f(\bbeta) $ subject to $\bbeta\in\sS$.
Let $f$ be a continuously differentiable function on a closed \textbf{convex} set $\sS$, and let $\bbeta^*$ be a local minimum of (P$_C$). Then $\bbeta^*$ is a stationary point of (P$_C$).
\end{theoremHigh}
\begin{proof}[of Theorem~\ref{theorem:stat_point_uncons_convset}]
Let $\bbeta^*$ be a local minimum of (P$_C$), and assume in contradiction that $\bbeta^*$ is not a stationary point of (P$_C$). Then there exists some $\bbeta \in \sS$ such that $\nabla f(\bbeta^*)^\top (\bbeta - \bbeta^*) < 0$. By {Lemma~\ref{lemma:descent_property}}, it follows that there exists $\varepsilon \in (0, 1)$ such that $f(\bbeta^* + \eta \bd) < f(\bbeta^*)$ for all $\eta \in (0, \varepsilon)$, where $\bd\triangleq \bbeta-\bbeta^*$. Since $\sS$ is convex we have that $\bbeta^* + \eta \bd = (1 - \eta)\bbeta^* + \eta\bbeta \in \sS$, contradicting the assumption that $\bbeta^*$ is  a local minimum point of (P$_C$). 
Hence, $\bbeta^*$ must be a stationary point.
\end{proof}

When the projection onto the convex set $\sS$ admits a closed-form expression, the first-order condition can be equivalently stated using the projection operator. This leads to the following characterization.
\begin{corollary}[Necessity/sufficiency for (P$_C$) on a convex set under projection]\label{corollary:stat_point_uncons_convset_proj}
Consider the constrained optimization problem (P$_C$): $\min_{\bbeta}f(\bbeta) $ subject to $\bbeta\in\sS$.
Let $f$ be a continuously differentiable function on a closed \textbf{convex} set $\sS$, and let $\eta > 0$ be an arbitrary positive scalar. Then $\bbeta^*$ is a stationary point of (P$_C$)
if and only if
$$
\bbeta^* = \projectS(\bbeta^* - \eta \nabla f(\bbeta^*)).
$$
\end{corollary}
\begin{proof}[of Corollary~\ref{corollary:stat_point_uncons_convset_proj}]
By the Projection Property-I (Lemma~\ref{lemma:proj_prop1}), it follows that $\bbeta^* = \projectS(\bbeta^* - \eta \nabla f(\bbeta^*))$ if and only if
$$
\begin{aligned}
&\innerproduct{\bbeta^* - \eta \nabla f(\bbeta^*) - \bbeta^*, \bbeta - \bbeta^*} \leq 0 
\quad \implies\quad \nabla f(\bbeta^*)^\top (\bbeta - \bbeta^*) \geq 0, \text{ for any } \bbeta \in \sS.
\end{aligned}
$$
This is precisely the definition of a stationary point in Definition~\ref{definition:stat_point_uncons_convset}, completing the proof.
\end{proof}

\subsection{Equality/Inequality Constraints}
We now consider the problem presented in Example~\ref{example:equa_ineq_cons}.
Such an {optimization problem} is typically formulated as
\begin{align}
&\min_{\bbeta \in \real^p} f_0(\bbeta) \label{opt:primal_gen} \tag{P}\\
\text{s.t.}\quad &\bX\bbeta = \by, \nonumber\\
&f_i(\bbeta) \leq z_i, \quad i \in \{1,2,\ldots,q\},\nonumber
\end{align}
where the function $f_0: \real^p \to (-\infty, \infty]$ is again called the {objective function} or {loss function}, 
the functions $f_1, f_2, \ldots, f_q: \real^p \to (-\infty, \infty]$ are referred to as  \textit{(inequality) constraint functions}, and $\bX \in \real^{n\times p}$, $\by \in \real^n$ define  the \textit{equality constraint}. 
Note that if $f_0, f_1, \ldots, f_q$ are all convex functions, then  problem \eqref{opt:primal_gen} is a \textbf{convex} optimization problem.

The equality constraint can be equivalently expressed using inequality constraints. Specifically, for each row $\balpha_i^\top$ of $ \bX$ (i.e., the $i$-th row), the condition $\balpha_i^\top\bbeta=y_i$ is equivalent to the pair of inequalities
$$
\balpha_i^\top\bbeta \leq y_i
\qquad\text{and}\qquad 
-\balpha_i^\top\bbeta \leq -y_i.
$$

Note that the set of feasible points described by the constraints can be obtained as 
\begin{equation}\label{equation:cvx_obj_gen_S}
\sS = \{ \bbeta \in \real^p \mid \bX\bbeta = \by, f_i(\bbeta) \leq z_i, i \in \{1,2,\ldots,q\} \}.
\end{equation}
Thus, problem~\eqref{opt:primal_gen} is equivalent to minimizing $f_0$ over the feasible set $\sS$,
\begin{equation}
\min_{\bbeta \in \sS} f_0(\bbeta).
\end{equation}
Furthermore, by introducing the indicator function $\indicatorS$ (see Example~\ref{exercise_convex_indica}), the constrained problem can be reformulated as an unconstrained optimization problem involving the composite objective
\begin{equation}
\min_{\bbeta \in \real^p} F(\bbeta)\triangleq f_0(\bbeta) + \indicatorS(\bbeta).
\end{equation}

\subsection{Dual Problem}
The \textit{Lagrange function} (or simply the \textit{Lagrangian}) associated with an optimization problem of the form \eqref{opt:primal_gen} is defined for $\bbeta \in \real^p$, $\blambda \in \real^n$, $\bnu \in \real^q$, $\nu_i \geq 0$ , $i \in \{1,2,\ldots,q\}$ (for brevity, we write $\bnu \geq \bzero$ when all components of $\bnu$ are nonnegative), by
\begin{equation}
L(\bbeta, \blambda, \bnu) 
\triangleq f_0(\bbeta) + \blambda^\top(\bX\bbeta - \by) + \sum_{i=1}^q \nu_i (f_i(\bbeta) - z_i)
\triangleq f_0(\bbeta) + \blambda^\top(\bX\bbeta - \by) + \bnu^\top\bff,
\end{equation}
where the vectors $\blambda$ and $\bnu$ are called \textit{Lagrange multipliers}, and 
$
\bff \triangleq \{f_i(\bbeta) - z_i\}_{i=1}^q\in\real^q.
$

The \textit{Lagrange dual function} (or simply the \textit{dual function}) is defined by
\begin{equation}
D(\blambda, \bnu) \triangleq \min_{\bbeta \in \real^p} L(\bbeta, \blambda, \bnu), \quad \blambda \in \real^n, \bnu \in \real^q, \bnu \geq \bzero.
\end{equation}
If the function $\bbeta \mapsto L(\bbeta, \blambda, \bnu)$ is unbounded from below,  we set $D(\blambda, \bnu) = -\infty$. 
Importantly, the dual function is always \textbf{concave} because it is the pointwise infimum of a family of affine functions in $(\blambda, \bnu)$, even if the original problem \eqref{opt:primal_gen} is not convex (see Exercise~\ref{exercise:pres_conv_clos} by noting that the negative of a concave function is convex). 

Now, suppose $\bbeta$ is a feasible point for \eqref{opt:primal_gen}. 
Then $\bX\bbeta - \by = \bzero$ and $f_i(\bbeta) \leq z_i$, $i = 1, 2, \ldots, q$. 
Since $\bnu \geq \bzero$, it follows that
$$
\begin{aligned}
&\blambda^\top(\bX\bbeta - \by) + \bnu^\top\bff \leq 0\\
\implies &L(\bbeta, \blambda, \bnu) = f_0(\bbeta) + \blambda^\top(\bX\bbeta - \by) + \bnu^\top\bff \leq f_0(\bbeta).
\end{aligned}
$$
Taking the infimum over all $\bbeta \in \real^p$ on the left-hand side (which yields $D(\blambda, \bnu)$), and over all feasible $\bbeta$ on the right-hand side (which gives the optimal value $f_0(\bbeta^*)$) 
leads to the fundamental inequality
\begin{equation}\label{equation:dual_bound_opt}
D(\blambda, \bnu) \leq f_0(\bbeta^*), \quad \text{for all } \blambda \in \real^n, \bnu \geq \bzero.
\end{equation}
Thus, the dual function provides a lower bound on the optimal value of the primal problem. To obtain the tightest possible lower bound, we maximize this dual function, which gives rise to the \textit{dual problem} to \eqref{opt:primal_gen}:
\begin{equation}\label{opt:dual_gen}
\max_{\blambda \in \real^n, \bnu \geq \bzero} D(\blambda, \bnu). \tag{D}
\end{equation}
Consequently, the original problem~\eqref{opt:primal_gen} is often referred to as the \textit{primal problem}.

\subsection*{Strong Duality}
Since problem~\eqref{opt:dual_gen} involves maximizing a concave function, it is equivalent to the convex optimization problem of minimizing the convex function $-D$ subject to the positivity constraint $\bnu \geq \bzero$. 
A pair $(\blambda, \bnu)$ with $\blambda \in \real^n$ and $\bnu \geq \bzero$ is called \textit{dual feasible}. 
A (dual feasible) maximizer $(\blambda^*, \bnu^*)$ of \eqref{opt:dual_gen} is referred to as \textit{dual optimal} or as the \textit{optimal Lagrange multipliers}. If $\bbeta^*$ is optimal for the primal problem \eqref{opt:primal_gen}, then the triple $(\bbeta^*, \blambda^*, \bnu^*)$ is called \textit{primal-dual optimal}. 
Inequality~\eqref{equation:dual_bound_opt} implies that, for  dual optimal $(\blambda^*, \bnu^*)$ and  primal optimal $\bbeta^*$,
\begin{equation}\label{equation:dual_bound_opt2}
D(\blambda^*, \bnu^*) \leq f_0(\bbeta^*).
\end{equation}
This relationship is known as \textit{weak duality}. 
For most---but not all---convex optimization problems, \textit{strong duality} holds, meaning that equality is achieved:
\begin{equation}
D(\blambda^*, \bnu^*) = f_0(\bbeta^*).
\end{equation}
See Problem~\ref{prob:farka_lemm}$\sim$\ref{prob:lp_strongdual} for the strong duality of linear programs.

A well-known sufficient condition for strong duality in convex problems is Slater's condition, stated below in a simplified form.

\begin{theoremHigh}[Slater's condition qualification for convex problems\index{Slater's condition}]\label{theorem:slater_cond}
Assume that $f_0, f_1, f_2, \ldots, f_q$ are convex functions with $\domain(f_0) = \real^p$. If there exists $\bbeta \in \real^p$ such that $\bX\bbeta = \by$ and $f_i(\bbeta) < 0$ for all $i \in \{1,2,\ldots,q\}$, then strong duality holds for the optimization problem \eqref{opt:primal_gen}.
\end{theoremHigh}
\begin{proof}
See, for example,  Chapter 6 of \citet{bertsekas90001convex} or Chapter 2 of \citet{lu2025practical}.
\end{proof}

We also recall some well-known optimality conditions expressed in terms of the Lagrangian function in cases where strong duality holds

\begin{theoremHigh}[KKT conditions: Optimality  under strong duality]\label{theorem:opt_cond_sd}
Suppose strong duality holds; that is, the optimal values of problems~\eqref{opt:primal_gen} and~\eqref{opt:dual_gen} are finite and equal.
Then $\widehatbbeta$ and $(\widehatblambda, \widehatbnu)$ are optimal solutions to the primal and dual problems, respectively, if and only if the following conditions hold:
\begin{enumerate}[(i)]
\item \textit{Primal feasibility.} $\widehatbbeta$ is a feasible solution of \eqref{opt:primal_gen}, i.e., $\bX\widehatbbeta=\by$ and $f_i(\widehatbbeta)\leq z_i$ for all $i=1,2,\ldots,q$.
\item \textit{Dual feasibility.} $\widehatbnu \geq \bzero$.
\item \textit{Complementary slackness.} $\widehatnu_i (f_i(\widehatbbeta)-z_i) = 0$ 
for all $i=1,2,\ldots,q$.
\item \textit{Stationarity.} $\widehatbbeta \in \argmin_{\bbeta \in \real^p} L(\bbeta, \widehatblambda, \widehatbnu)$.
\end{enumerate}
\end{theoremHigh}
\begin{proof}[of Theorem~\ref{theorem:opt_cond_sd}]
\textbf{Necessity.}
Let $\widehatf$ and $\widehatd$ denote the optimal values of the primal and dual problems, respectively. 
By assumption, strong duality holds, so $\widehatf = \widehatd$. 
If $\widehatbbeta$ and $(\widehatblambda, \widehatbnu)$ are the optimal solutions of problems \eqref{opt:primal_gen} and \eqref{opt:dual_gen}, then (i) and (ii) follow immediately from feasibility. 
In addition,
\begin{align*}
\widehatf &= \widehatd = D(\widehatblambda, \widehatbnu)
= \min_{\bbeta \in \real^p} L(\bbeta, \widehatblambda, \widehatbnu) 
\leq L(\widehatbbeta, \widehatblambda, \widehatbnu) \\
&= f(\widehatbbeta) + \innerproduct{\widehatblambda, \bX\widehatbbeta-\by} 
+ \sum_{i=1}^{q} \widehatnu_i \big(f_i(\widehatbbeta)-z_i\big) 
\leq f(\widehatbbeta),
\end{align*}
where the last inequality follows by the facts that $\bX\widehatbbeta=\by$, and $f_i(\widehatbbeta) \leq z_i$. Since $\widehatf = f(\widehatbbeta)$, all inequalities above must be equalities.
This implies  that $\widehatbbeta \in \argmin_{\bbeta \in \real^p} L(\bbeta, \widehatblambda, \widehatbnu)$, meaning property (iv), 
and that $\sum_{i=1}^{q} \widehatnu_i (f_i(\widehatbbeta)-z_i) = 0$, which by the fact that $\widehatnu_i (f_i(\widehatbbeta) -z_i) \leq 0$ for any $i$, implies that $\widehatlambda_i (f_i(\widehatbbeta) -z_i) = 0$ for any $i$, showing the property (iii).

\paragraph{Sufficiency.}
Conversely, assume that properties (i)-(iv) hold. Then,
\begin{align*}
D(\widehatblambda, \widehatbnu) 
&= \min_{\bbeta \in \real^p} L(\bbeta, \widehatblambda, \widehatbnu)
= L(\widehatbbeta, \widehatblambda, \widehatbnu)  \\
&= f(\widehatbbeta) + \innerproduct{\widehatblambda, \bX\widehatbbeta-\by} + \sum_{i=1}^{q} \widehatnu_i (f_i(\widehatbbeta)-z_i) 
= f(\widehatbbeta). 
\end{align*}
where the second equality follows from property (iv), and the last equality follows from  {properties (i) and (iii)}.
By the weak duality theorem, since $\widehatbbeta$ and $(\widehatblambda, \widehatbnu)$ are primal and dual feasible solutions with equal primal and dual objective functions, it follows that they are the optimal solutions of their corresponding problems. 
\end{proof}

\begin{remark}[KKT conditions for convex optimization with equality constraints]\label{remark:kkt_nutshell_cvx}
Consider the special case of an optimization problem with only equality constraints:
$$
\min_{\bbeta \in \real^p} f_0(\bbeta) \quad \text{s.t.} \quad \bX\bbeta = \by,
$$
where $f_0$ is \textbf{convex} and differentiable.
Then Theorems~\ref{theorem:slater_cond} and \ref{theorem:opt_cond_sd} show that  $ \widehatbbeta $ is a solution if and only if $ \widehatbbeta $ is feasible and there exists $ \blambda \in \real^n $ such that:
$$
\nabla_{\bbeta} L(\widehatbbeta, \blambda) = \bzero = \nabla f_0(\widehatbbeta) - \bX^\top \blambda,
$$
where the Lagrangian is defined as $L(\bbeta, \blambda) = f_0(\bbeta) + \innerproduct{\blambda, \by-\bX\bbeta}$, and the vector $\blambda$ is the Lagrange multiplier associated with the equality constraint.
The geometric interpretation is that  $\nabla f_0(\widehatbbeta) \perp \nspace(\bX)$.

When $ f_0 $ is not differentiable,  the stationarity condition generalizes to: 
$ \widehatbbeta $ is feasible and there exists $ \blambda \in \real^n$ such that
$$
\bX^\top \blambda \text{ is a subgradient of } f_0 \text{ at } \widehatbbeta.
$$
\end{remark}

Given a primal-dual feasible triple $(\bbeta, \blambda, \bnu)$---that is, $\bbeta$ is feasible for \eqref{opt:primal_gen} and $\blambda \in \real^n$, $\bnu \geq \bzero$---the \textit{primal-dual gap} is defined as
\begin{equation}\label{equation:cvx_prim_dual_gap}
H(\bbeta, \blambda, \bnu) = F(\bbeta) - D(\blambda, \bnu), 
\quad \text{with }F(\bbeta)\triangleq f_0(\bbeta) + \indicatorS(\bbeta).
\end{equation}
This quantity serves as a measure of how close $\bbeta$ is to a primal minimizer, and how close  $(\blambda, \bnu)$ is to a dual maximizer. 
If $(\widehatbbeta, \widehatblambda, \widehatbnu)$ is a primal-dual optimal solution and strong duality holds, then the gap vanishes: $H(\widehatbbeta, \widehatblambda, \widehatbnu) = 0$. 
In practice, the primal-dual gap is commonly used as a stopping criterion in iterative optimization algorithms.

To illustrate these concepts, we now derive the dual of the $\ell_\infty$-minimization problem. (Its counterpart, the $\ell_1$-minimization problem, will be discussed in greater detail in Chapter~\ref{chapter:sparse_opt_cond}.)
\begin{example}[$\ell_\infty$-minimization problem]\label{example:linf_min_dual}
Consider the $\ell_\infty$-minimization problem:
\begin{equation}\label{equation:ellinf_exp}
\min_{\bbeta \in \real^p} \norminf{\bbeta} \quad \text{s.t.}\quad \bX\bbeta = \by.
\end{equation}
The Lagrangian associated with this problem is
$$
L(\bbeta, \blambda) = \norminf{\bbeta} + \innerproduct{\blambda, (\bX\bbeta - \by)}.
$$
The corresponding dual function is obtained by minimizing the Lagrangian over  $\bbeta$:
$$
D(\blambda) = \min_{\bbeta \in \real^p} \left\{ \norminf{\bbeta} + \innerproduct{\bX^\top \blambda, \bbeta} - \innerproduct{\blambda, \by} \right\}.
$$
By \holders inequality, we have $\innerproduct{ \bX^\top \blambda, \bbeta} \leq \normone{\bX^\top\blambda}\norminf{\bbeta}$. 
This implies that 
if $\normone{\bX^\top \blambda } > 1$, then there exists $\bbeta \in \real^p$ (by picking its sign appropriately) such that $\innerproduct{ \bX^\top \blambda, \bbeta} < -\norminf{\bbeta}$. 
Replacing $\bbeta$ by $\mu \bbeta$ and letting $\mu \to \infty$ shows that $D(\blambda) = -\infty$ in this case. 
On the other hand, if $\normone{\bX^\top \blambda} \leq 1$, then $\norminf{\bbeta} + \innerproduct{\bX^\top \blambda, \bbeta} \geq 0$. The choice $\bbeta = \bzero$ yields  the infimum, and $D(\blambda) = -\innerproduct{\blambda, \by}$. Therefore, the dual problem is equivalent to:
$$
\max_{\blambda \in \real^n} 
D(\blambda) 
= 
\begin{cases}
-\innerproduct{\blambda, \by}, & \text{if } \normone{\bX^\top \blambda} \leq 1; \\
-\infty, & \text{otherwise}.
\end{cases}
$$
Since maximizing $D(\blambda)$ is only meaningful where it is finite, the dual problem is equivalent to
\begin{equation}
\max_{\blambda \in \real^n} - \innerproduct{\blambda, \by} \quad \text{s.t.}\quad \normonebig{\bX^\top \blambda} \leq 1.
\end{equation}
Equivalently, one may write the objective as $\max_{\blambda \in \real^n} \innerproduct{\blambda, \by}$ by absorbing the minus sign into the maximization (though the constraint remains unchanged).
Finally, by Slater's condition (Theorem~\ref{theorem:slater_cond}), strong duality holds for this primal-dual pair provided the primal problem~\eqref{equation:ellinf_exp} is feasible---since the objective $\norminf{\cdot}$ is convex and the equality constraints are affine.
\end{example}

\subsection*{Saddle-Point Property}
Lagrange duality admits a natural \textit{saddle-point interpretation}. 
For simplicity, we consider the primal problem~\eqref{opt:primal_gen} without inequality constraints; however, the same reasoning extends straightforwardly to problems with inequality constraints.

Let $(\bbeta^*, \blambda^*)$ be a primal-dual optimal pair. 
By definition of the Lagrangian $L$, we have
\begin{align}
\max_{\blambda \in \real^n} L(\bbeta, \blambda) 
&= \max_{\blambda \in \real^n} f_0(\bbeta) + \blambda^\top(\bX\bbeta - \by)
= \begin{cases}
f_0(\bbeta), & \text{if } \bX\bbeta = \by, \\
\infty, & \text{otherwise}.
\end{cases}
\end{align}
In other words, this supremum is infinite whenever $\bbeta$ violates the equality constraints. 
Consequently, the (feasible) minimizer  $\bbeta^*$ of the primal problem \eqref{opt:primal_gen} satisfies 
\begin{equation}
f_0(\bbeta^*) = \min_{\bbeta \in \real^p} \max_{\blambda \in \real^n} L(\bbeta, \blambda).
\end{equation}
On the other hand, by the definition of the Lagrange dual function, a dual optimal vector $\blambda^*$ achieves
\begin{equation}
D(\blambda^*) = \max_{\blambda \in \real^n} \min_{\bbeta \in \real^p} L(\bbeta, \blambda)
\end{equation}
Weak duality then implies the general inequality
\begin{equation}
\max_{\blambda \in \real^n} \min_{\bbeta \in \real^p} L(\bbeta, \blambda) \leq \min_{\bbeta \in \real^p} \max_{\blambda \in \real^n} L(\bbeta, \blambda).
\end{equation}
When strong duality holds, equality is achieved:
\begin{equation}
\max_{\blambda \in \real^n} \min_{\bbeta \in \real^p} L(\bbeta, \blambda) = \min_{\bbeta \in \real^p} \max_{\blambda \in \real^n} L(\bbeta, \blambda).
\end{equation}
In this case, the order of minimization and maximization can be interchanged in the case of strong duality. 
This is known as the strong \textit{max-min property} or the \textit{saddle-point property}.
Indeed, under strong duality, any primal-dual optimal pair $(\bbeta^*, \blambda^*)$ forms a saddle point of the Lagrangian:
\begin{equation}\label{equation:saddl_point_gen}
	L(\bbeta^*, \blambda) \leq L(\bbeta^*, \blambda^*) \leq L(\bbeta, \blambda^*) \quad \text{for all } \bbeta \in \real^p, \blambda \in \real^n.
\end{equation}
Thus, solving the primal and dual problems jointly is equivalent to finding a saddle point of the Lagrangian, provided that strong duality holds.

\subsection{Composite Function}

Conjugate functions arise naturally in dual optimization problems, most notably in the celebrated \textit{Fenchel's duality theorem}, which we now recall. Consider the following composite optimization problem:
$$
\min_{\bbeta \in \real^p} F(\bbeta) + G(\bbeta).
$$
All the relevant optimization problems discussed in this book belong to this class; see Chapter~\ref{chapter:sparse_opt_cond}. 
\begin{example}
Many problems covered in this book take this form.
For example, setting $F(\bbeta)\triangleq \frac{1}{2}\normtwo{\bX\bbeta-\by}^2$ and $G(\bbeta) = \lambda\normone{\bbeta}$ yields the \textit{Lagrangian LASSO problem} \eqref{opt:ll}; the choice $G(\bbeta) = \normone{\bbeta}$ and $\widetildeF(\bz) \triangleq F(\bX\bbeta)$ with $\widetildeF= \delta_{\{\by\}}$ (the indicator function of the singleton $\{\by\}$) gives the $\ell_1$-minimization problem \eqref{opt:p1};
the choice $G(\bbeta) = \normone{\bbeta}$ and $\widetildeF(\bz) \triangleq F(\bX\bbeta)$ with $\widetildeF = \delta_{\sB[\by, \epsilon]}$ (the indicator function of the the $\ell_2$-ball with center $\by$ and radius $\epsilon$) leads to the quadratically-constrained  $\ell_1$-minimization problem \eqref{opt:p1_epsilon}.
\end{example}

We begin by reformulating the original problem as a constrained optimization problem:
$$
\min_{\bbeta,\balpha \in \real^p} \left\{ F(\bbeta) + G(\balpha)\right\}  
\quad \text{s.t.}\quad  \bbeta = \balpha .
$$
The associated Lagrangian is then
$$
L(\bbeta, \balpha, \blambda) = F(\bbeta) + G(\balpha) + \innerproduct{\blambda, \balpha - \bbeta}
 = -[ \innerproduct{\blambda, \bbeta} - F(\bbeta) ] - [ \innerproduct{-\blambda, \balpha} - G(\balpha) ].
$$
The dual function is obtained by minimizing the Lagrangian over the primal variables $\bbeta, \balpha$:
$$
D(\blambda) = \min_{\bbeta,\balpha} L(\bbeta, \balpha; \blambda) = -F^*(\blambda) - G^*(-\blambda).
$$
Thus, the corresponding dual problem---commonly referred to as \textit{Fenchel's dual}---is
$$
\max_{\blambda \in \real^p} \left\{ -F^*(\blambda) - G^*(-\blambda) \right\}.
$$

Fenchel's duality theorem, stated below, provides sufficient conditions under which strong duality holds between the primal problem and its Fenchel dual.

\begin{theoremHigh}[Fenchel's duality theorem; \citet{rockafellar2015convex}, Theorem 31.1\index{Fenchel's duality theorem}]\label{theorem:fenchel_dual_ori}
Let $F, G: \real^p \to (-\infty, \infty]$ be proper convex functions. If $\relint(\domain(F)) \cap \relint(\domain(G)) \neq \varnothing$, then
$$
\min_{\bbeta \in \real^p} \{ F(\bbeta) + G(\bbeta) \} = \max_{\blambda \in \real^p} \left\{ -F^*(\blambda) - G^*(-\blambda) \right\},
$$
and the maximum on the right-hand side is attained whenever it is finite.
\end{theoremHigh}

\begin{problemset}
\item  \textbf{$\ell_0$ and $\ell_1$.}
For any scalar $\beta$, define the vector
$$
\bu_\beta \triangleq [\beta,
\beta, 
\beta-1 ]^\top.
$$
Choose  $\beta$ so that $\bu_\beta$ is as sparse as possible; that is, minimize $\normzero{\bu_\beta}$. 
Show that minimizing $\normzero{\bu_\beta}$ is equivalent to minimizing $\normone{\bu_\beta}$, in the sense that both problems attain their global minimum at the same value(s) of $\beta$.

\item \textbf{Rank and nuclear norm.}
For any scalar $\beta$, consider the  matrix
$$
\bX_\beta 
\triangleq 
\begin{bmatrix}
0.5 + \beta & 1 & 1 \\
0.5 & 0.5 & \beta \\
0.5 & 1 - \beta & 0.5
\end{bmatrix}.
$$
Show that the value of $\beta$ minimizing  the rank of $\bX_\beta$ also minimizes its nuclear norm (i.e., the sum of its singular values).

\item \textbf{Necessity/sufficiency  of unconstrained convex optimization under proximal.}
Let $ f:  \real^p \rightarrow (-\infty, \infty]  $ be a proper closed and convex function. Show that $ \balpha^* $ is a minimizer of $ f $ if and only if $ \balpha^* = \proxf(\balpha^*) $.
\textit{Hint: Use Theorem~\ref{theorem:fetmat_opt} and Proximal Property-I (Lemma~\ref{lemma:prox_prop1}).}


\item \label{prob:coerci_quad} \textbf{Coerciveness of quadratic functions.}
Let $f(\bbeta)=\frac{1}{2}\bbeta^\top\bX\bbeta+\by^\top\bbeta+z$, where $\bX\in\real^{p\times p}$ is symmetric, $\by\in\real^p$, and $z\in\real$. Show that $f(\bbeta)$ is coercive (i.e., $ f(\bbeta) \to \infty $ as $ \normtwo{\bbeta} \to \infty $) if $\bX$ is positive definite.


\item Using the proof of Theorem~\ref{theorem:psd_hess_conv}, establish the equivalence among \eqref{equation:scss_func1}, \eqref{equation:scss_func2}, and \eqref{equation:scss_func3_ra} as characterizations of strong convexity and strong smoothness.

\item Prove the properties of Bregman distances stated in Remark~\ref{remark:bregnan_dist}.

\item Let $ f \in C^{2,2}_L(\real^p) $, i.e., $f$ is twice continuously differentiable and its Hessian is $L$-Lipschitz continuous. 
Show that  for all $ \balpha, \bbeta \in \real^p $,  the following inequalities hold:
$$
\begin{aligned}
&\normtwo{\nabla f(\bbeta) - \nabla f(\balpha) - \nabla^2 f(\balpha)(\bbeta - \balpha)} \leq \frac{L}{2} \normtwo{\bbeta - \balpha}^2;\\
&\abs{f(\bbeta) - f(\balpha) - \innerproduct{\nabla f(\balpha), \bbeta - \balpha} - \frac{1}{2} \innerproduct{\nabla^2 f(\balpha)(\bbeta - \balpha), \bbeta - \balpha}} \leq \frac{L}{6} \normtwo{ \bbeta - \balpha}^3.
\end{aligned}
$$
Furthermore, if $\normtwo{\balpha-\bbeta}\leq M$, show that:
$$
\nabla^2 f(\balpha) -L M \bI \preceq \nabla^2 f(\bbeta) \preceq \nabla^2 f(\balpha)+ LM  \bI.
$$
\textit{Hint: Refer to Theorem~\ref{theorem:equi_gradsch_smoo} and Exercise~\ref{exercise:cl2111_hess_bound}.}


\item \label{prob:stric_sep_theo2} 
\textbf{General separation theorem.}
Let $\sS, \sT \subseteq \real^p$ be \textbf{convex} sets with empty intersection $\sS \cap \sT = \varnothing$. 
Show that there exists a vector $\bz \in \real^p$ and a scalar $\gamma \in \real$ such that
\begin{enumerate}[(i)]
\item For all $\balpha \in \sS$, we have $\innerproduct{\bz, \balpha} \geq \gamma$.
\item For all $\balpha \in \sT$, we have $\innerproduct{\bz, \balpha} \leq \gamma$.
\end{enumerate}
Moreover, if both $\sS$ and $\sT$ are closed and at least one of them is bounded,then the inequalities can be strengthened to strict inequalities:
\begin{enumerate}[(i)]
\item For all $\balpha \in \sS$, we have $\innerproduct{\bz, \balpha} > \gamma$.
\item For all $\balpha \in \sT$, we have $\innerproduct{\bz, \balpha} < \gamma$.
\end{enumerate}


\item \label{prob:supp_hyp_theo} \textbf{Supporting hyperplane theorem.}
Let $\sS \subseteq \real^p$ be a \textbf{convex} set (not necessarily closed), and let $\bbeta \notin \sS$. 
Show that there exists $\bzero \neq \bw \in \real^p$ such that
$$
\bw^\top \balpha \leq \bw^\top \bbeta, \quad  \text{ for any } \balpha \in \sS.
$$

\item \label{prob:breg_prox_propo} \textbf{Bregman-Proximal Property.}
Assume the following:
\begin{itemize}
\item $\phi : \real^p \to (-\infty, \infty]$ is  proper \textbf{closed  convex}, differentiable over $\dom(\partial \phi)$.
\item $f : \real^p \to (-\infty, \infty]$ is  proper \textbf{closed  convex} with $\dom(f) \subseteq \dom(\phi)$.
\item $\phi + \delta_{\dom(f)}$ is $\sigma$-strongly convex ($\sigma > 0$).
\end{itemize}
\noindent
Show that the optimization problem
$$
\min_{\balpha \in \real^p} \{f(\balpha) + \phi(\balpha)\}
$$
has a unique minimizer, which lies in $\dom(f) \cap \dom(\partial \phi)$.


\item \label{prob:breg_prox_prop1} \textbf{Bregman-Proximal Property.}
Assume the following:
\begin{itemize}
\item $\phi : \real^p \to (-\infty, \infty]$ is  proper \textbf{closed convex}, differentiable over $\dom(\partial \phi)$;
\item $f : \real^p \to (-\infty, \infty]$ is  proper \textbf{closed  convex} with  $\dom(f) \subseteq \dom(\phi)$;
\item $\phi + \delta_{\dom(f)}$ is $\sigma$-strongly convex ($\sigma > 0$).
\end{itemize}
\noindent 
Show that for any $ \bbeta, \widehatbbeta \in \dom(\partial \phi) $, the following two statements  are equivalent:
\begin{enumerate}[(i)]
\item $ \widebarbbeta = \bproxfphi(\bbeta) \triangleq \mathop{\argmin}_{\balpha \in \real^p} \left\{ f(\balpha) + \mathcalD_\phi(\balpha, \bbeta) \right\}$ .
\item $ \innerproduct{\nabla \phi(\bbeta) - \nabla \phi(\widebarbbeta), \bxi - \widebarbbeta} \leq f(\bxi) - f(\widebarbbeta)$ for any $ \bxi \in \dom(f) $.
\end{enumerate}


\item \label{prob:farka_lemm}\textbf{Farkas' lemma.}
Let $\bc \in \real^n$ and $\bA \in \real^{m \times n}$. Show that \textbf{exactly one} of the following systems has a solution:
\begin{enumerate}[(i)]
\item $\bA\bx \leq \bzero, \bc^\top \bx > 0.$
\item $\bA^\top \by = \bc, \by \geq \bzero.$  (i.e., $\bc$ is a conic combination of rows of $\bA$.)
\end{enumerate}
\noindent Equivalently, this implies  the following two claims are equivalent:
\begin{enumerate}[(a)]
\item The implication $\bA\bx \leq \bzero \implies \bc^\top \bx \leq 0$ holds true.
\item There exists $\by\in \real_{+}^{m}$ such that $\bA^\top \by = \bc$.
\end{enumerate}

\item \label{prob:farka_lemm_prim}\textbf{Farkas' lemma$'$.}
Show that the  Farkas' lemma (Problem~\ref{prob:farka_lemm}) also implies the following generalized result.
Given matrices $ \bD, \bE, $ and $ \bF $ of appropriate dimensions, \textbf{exactly one} of the following two systems has a solution:
\begin{enumerate}[(1)]
\item $ \bD\bx < \bzero, \bE \bx \leq \bzero, \bF \bx = \bzero, $
\item $ \bD^\top \bu + \bE^\top \bw + \bF^\top \bv = \bzero, $ $ \bu \geq \bzero, \bw \geq \bzero, \bone^\top \bu = 1. $
\end{enumerate}

\item \label{prob:lp_strongdual} Use Farkas' lemma to prove the strong duality theorem for the linear program presented in  Example~\ref{example:linear_program}: If the primal linear program is feasible and has a finite optimal value, and the dual is also feasible, then both problems attain the same optimal objective value.

\end{problemset}

\newpage 
\chapter{Sparse Optimization Problems and Optimality Conditions}\label{chapter:sparse_opt_cond}
\begingroup
\hypersetup{
linkcolor=structurecolor,
linktoc=page,  
}
\minitoc \newpage
\endgroup

\lettrine{\color{caligraphcolor}I}
In the fields of \textit{sparse representation} and \textit{compressed sensing} (also known as \textit{compressive sensing}, \textit{compressive sampling}, or \textit{sparse sampling}), a central challenge is to recover a high-dimensional signal from a relatively small number of linear measurements. This task is mathematically formulated as solving an underdetermined linear system  
$$
\bX \boldsymbol{\beta} = \by,
$$
where $\bX \in \real^{n \times p}$ (with $n \ll p$) is a measurement or sensing matrix, $\by \in \real^n$ is the observed data vector, and $\boldsymbol{\beta} \in \real^p$ is the unknown coefficient vector assumed to be \textit{sparse}, meaning that only a small subset of its entries are nonzero. The sparsity assumption is crucial: without it, the system admits infinitely many solutions due to its underdetermined nature. However, when $\boldsymbol{\beta}$ is known (or believed) to be sparse---or approximately sparse---it becomes possible, under suitable conditions on $\bX$ (e.g., the restricted isometry property or incoherence; see Chapters~\ref{chapter:design} and \ref{chapter:recovery}), to uniquely and stably reconstruct $\boldsymbol{\beta}$ from $\by$.

This chapter focuses on how tools from optimization theory can be leveraged to solve such sparse inverse problems. Specifically, we explore two closely related but conceptually distinct classes of sparse optimization problems:
\begin{enumerate}[(i)]
\item \textit{Optimization with sparse regression (or sparse constraints).} In this setting, sparsity is imposed as a constraint or regularizer within a broader regression or estimation framework. Typical formulations include minimizing a loss function (e.g., least squares) subject to a constraint on the number of nonzero coefficients (an ``$\ell_0$-norm'' constraint) or, more tractably, adding an $\ell_1$-norm penalty (as in the \textit{LASSO}). These approaches are widely used in statistics, machine learning, and signal processing for feature selection and model interpretability.

\item \textit{Sparse recovery problems.} Here, the primary goal is exact or stable reconstruction of a sparse signal from incomplete measurements, often in the presence of noise. This perspective is common in applied mathematics and engineering, particularly in compressed sensing. Recovery is typically achieved by solving convex relaxations (e.g., \textit{basis pursuit}) or non-convex approximations of the original combinatorial problem, with theoretical guarantees provided under assumptions about the sensing matrix $\bX$.
\end{enumerate}

While both paradigms aim to exploit sparsity, they differ in emphasis: sparse regression often prioritizes predictive performance and model simplicity in data-rich contexts, whereas sparse recovery emphasizes faithful signal reconstruction from minimal measurements. In what follows, we develop the mathematical foundations, and theoretical insights underlying these two approaches, illustrating how optimization serves as a unifying framework for sparse modeling.

\section{Sparse Optimization Problems}\label{section:sparse_opt_gen}

On the one hand, \textit{compressed sensing}---a term pioneered by David Donoho, Emmanuel Candès, and Terence Tao \citep{donoho2001uncertainty, donoho2003optimally, candes2005decoding}---has emerged as a groundbreaking paradigm for sparse signal recovery, fundamentally reshaping how we tackle underdetermined linear systems. 
It exploits the fact that many natural signals admit sparse representations in some basis or frame, enabling accurate reconstruction from far fewer measurements than required by the Nyquist rate. This directly challenges the classical view that such inverse problems are inherently ill-posed. Beyond its theoretical elegance, compressed sensing has had a profound practical impact, transforming fields such as signal processing, information theory, and statistics, while also catalyzing new research directions---including low-rank matrix completion and phase retrieval \citep{baraniuk2011introduction, foucart2013invitation, jain2017non}. In doing so, it stands as a quintessential example of how a unifying mathematical insight can redefine data acquisition and analysis in today's data-intensive era.

A key advantage of compressed sensing is its ability to combine sampling and compression into a single, efficient step. Instead of collecting large volumes of raw data and compressing them afterward, compressed sensing leverages signal sparsity at the point of acquisition, dramatically reducing the number of measurements needed. This not only lowers data storage and transmission costs but also accelerates the acquisition process---making it especially valuable in time-critical applications such as magnetic resonance imaging (MRI) and wireless communications, where it enables shorter scan times or reduced bandwidth usage without compromising fidelity \citep{han2013compressive}. Moreover, its inherent robustness to noise and resilience to missing observations allow high-quality signal reconstruction even from highly incomplete data, yielding substantial gains in speed, efficiency, and overall system performance.

Although the basis pursuit formulation for sparse recovery was introduced by S. Chen, D. Donoho, and M. Saunders \citep{chen2001atomic} as a principled method to find the sparsest solution to an underdetermined linear system via $\ell_1$-norm minimization, the \textit{LASSO (least absolute shrinkage and selection operator)}, proposed earlier by R. Tibshirani \citep{tibshirani1996regression}, has become the dominant framework within the statistics and machine learning communities \citep{hastle2015statistical}. 
This divergence in popularity stems largely from differences in problem framing and disciplinary priorities. In the compressed sensing literature---rooted in applied mathematics, signal processing, and information theory---the emphasis is on exact or stable recovery of a true sparse signal from incomplete measurements, often under stringent assumptions about the sensing matrix (e.g., the restricted isometry property or mutual coherence; see Chapters~\ref{chapter:design} and \ref{chapter:recovery}). Basis pursuit and its noise-aware variants, such as \textit{basis pursuit denoising (BPDN)}, are naturally aligned with this goal, as they directly minimize the $\ell_1$-norm subject to explicit data fidelity constraints.

By contrast, the LASSO originated in the context of regression modeling, where the primary objectives are prediction accuracy, model interpretability, and variable selection in high-dimensional settings. 
Rather than enforcing hard constraints on residual error, the LASSO formulates sparsity as a trade-off: it minimizes the sum of squared residuals plus an $\ell_1$-penalty on the coefficients. This Lagrangian formulation offers greater practical flexibility---it eliminates the need to pre-specify a noise level or tolerance parameter and integrates seamlessly into cross-validation pipelines for tuning the regularization strength. 
Consequently, the LASSO resonated strongly with statisticians and data scientists seeking robust, scalable tools for feature selection and model regularization. 
Over time, the distinction between the two approaches has blurred, as they can be shown to be equivalent under appropriate choices of parameters. Nevertheless, their historical origins and methodological emphases continue to shape how researchers in different fields conceptualize and apply sparse optimization.

We now turn to the precise mathematical formulations of these two canonical problems, along with their closely related variants.

\subsection{Sparse  Regression Problems}
In this subsection, we start by the challenge of minimizing a (general) continuously differentiable objective function subject to a sparsity constraint. Specifically, we consider the following problem:
\begin{equation}\label{opt:s0_gen}
\min_{\bbeta} f(\bbeta) \quad \text{s.t.}  \quad\normzero{\bbeta} \leq R, \tag{S$_0$}
\end{equation}
where $f : \real^p \to \real$ is continuously differentiable with a Lipschitz continuous gradient of constant $L_b$ (in our case, if we  take $f(\bbeta) = \frac{1}{2}\normtwo{\by - \bX \bbeta}^2$, the least squares objective, then $L_b=\normtwo{\bX^\top\bX}$; see Example~\ref{example:lipschitz_spar}), $R > 0$ is an integer satisfying $R<p$, and $\normzero{\bbeta}$ is the so-called ``$\ell_0$-norm" of $\bbeta$, which counts the number of nonzero components in $\bbeta$.
For any positive number $s > 0$, the $\ell_s$-norm~\footnote{To avoid notational conflict, we use the term $\ell_s$-norm here, although the literature commonly uses $\ell_p$-norm for the same concept.} of a vector $\bbeta$ is defined as
$
\norm{\bbeta}_s = \left( \sum_{i} \abs{\beta_i}^s \right)^{1/s}.
$
Thus, the $\ell_0$-norm of a $p \times 1$ vector $\bbeta$ can be defined as
$$
\normzero{\bbeta} = \lim_{s\to 0} \norms{\bbeta}^s 
= \lim_{s\to 0} \sum_{i=1}^p \abs{\beta_i}^s 
= \sum_{i=1}^p \indicator(\beta_i \neq 0) 
= \abs{\{i\mid \beta_i \neq 0 , \forall\, i\in\{1,2,\ldots,p\}\}}.
$$
Thus, if $\normzero{\bbeta} \ll p$, the vector $\bbeta$ is said to be sparse.
It is important to note that the $\ell_0$-norm is not actually a norm, as it fails to satisfy the homogeneity property: for any nonzero scalar $\lambda$, we have $\normzero{\lambda \bbeta} =  \normzero{\bbeta}$ (Definition~\ref{definition:matrix-norm}).

Problem~\eqref{opt:s0_gen} is a constrained optimization problem in which we seek the best model---i.e., the one that minimizes $f$---under a hard sparsity constraint. By restricting the solution to have at most $R$ nonzero coefficients, we explicitly control model complexity, promoting interpretability and mitigating overfitting.

\subsection*{LASSO}
From the above sparse optimization analysis, 
we observe that---although convergence of algorithms like \textit{iterative hard-thresholding} (IHT; see Section~\ref{section:ell0_const_algo_IHT}) can be guaranteed under very specific conditions---the use of the $\ell_0$-norm (or hard-thresholding), which counts the number of nonzero elements in a vector, presents several significant challenges and drawbacks:
\begin{itemize}
\item \textit{Non-convexity.} The  $\ell_0$-norm is inherently non-convex. Consequently, optimization problems involving it are generally NP-hard. Finding the global minimum becomes computationally prohibitive, especially as the problem dimension grows.
\item 	\textit{Combinatorial nature.} Minimizing the $\ell_0$-norm requires evaluating all possible subsets of variables to identify the sparsest solution satisfying the constraints. This combinatorial search has exponential complexity, rendering it infeasible for high-dimensional datasets.

\item \textit{Instability.} Solutions obtained by direct $\ell_0$ minimization are highly sensitive to small perturbations in the data. This lack of robustness is particularly problematic in noisy settings or when working with real-world data containing errors or outliers.

\item \textit{Lack of continuity and differentiability.} The $\ell_0$-norm is neither continuous nor differentiable, making it incompatible with standard gradient-based optimization methods. This absence of smoothness hinders the use of iterative algorithms designed for differentiable objectives.

\item \textit{Difficulty in incorporating into optimization algorithms.} Due to its discrete and discontinuous nature, the $\ell_0$-norm is difficult to incorporate directly into most optimization frameworks. As a result, practical approaches often rely on relaxations or surrogates---most notably the convex $\ell_1$-norm---to render the problem tractable.
\end{itemize}

In contrast, the $\ell_1$-norm also promotes sparsity: it encourages solutions in which many coefficients are exactly zero. This behavior arises from its geometric structure---in high dimensions, the corners of the $\ell_1$-ball (a polytope) lie along the coordinate axes, which favors sparse solutions; see Figures~\ref{fig:l1l2_ball_loss} and \ref{fig:s_norm_2d_BETA_intersection}.
Importantly, the $\ell_1$-norm is convex. Therefore, any local minimum of an $\ell_1$-regularized problem is also a global minimum. This property dramatically simplifies optimization compared to the non-convex $\ell_0$ setting.

To be more concrete, consider the standard linear regression model in matrix-vector form:
\begin{equation}
\by = \bX \bbeta^* + \bepsilon,
\end{equation}
where $\bX \in \real^{n\times p}$ is the design (or model) matrix, $\bepsilon \in \real^n$ is a vector of noise variables, and $\bbeta^* \in \real^p$ is the unknown coefficient vector. 
The \textit{standard linear regression} or \textit{least squares problem} then corresponds to the following unconstrained optimization:
\begin{equation}\label{opt:l2}
(\textbf{Least squares}):\qquad 
\mathopmin{\bbeta\in\real^p}f(\bbeta)\triangleq \frac{1}{2}\normtwo{\by - \bX \bbeta}^2.
 \tag{L$_{2}$}
\end{equation}
By employing convex relaxation via the $\ell_1$-norm, later sections of this chapter will develop the \textit{LASSO  (least absolute selection and shrinkage operator)} framework, covering both its constrained and Lagrangian formulations. Practical algorithms and theoretical guarantees are discussed in Sections~\ref{section:alg_cons_lasso} and \ref{section:sparse_lin_reg_ana}, respectively:
\begin{subequations}\label{equation:lassos_all}
\begin{align}
(\textbf{Constrained LASSO}):\qquad &\mathopmin{\normone{\bbeta}\leq \Sigma}f(\bbeta)\triangleq \frac{1}{2}\normtwo{\by - \bX \bbeta}^2;
\label{opt:lc} \tag{L$_{C}$}\\
(\textbf{Lagrangian LASSO}):\qquad &\mathopmin{\bbeta \in \real^p}F(\bbeta)\triangleq \left\{ \frac{1}{2} \normtwo{\by - \bX \bbeta}^2 + \lambda \normone{\bbeta} \right\}. 
\label{opt:ll} \tag{L$_{L}$}
\end{align}
\end{subequations}
The constrained LASSO uses the same least squares loss as in \eqref{opt:l2}, but restricts the solution to lie within the $\ell_1$-ball $\normone{\bbeta}<\Sigma$.
The Lagrangian LASSO, also known as the \textit{$\ell_1$-regularized least squares} problem, incorporates sparsity through an additive penalty term.
Through Lagrangian duality, these two formulations are equivalent in the following sense: for every $\Sigma>0$, there exists a $\lambda \geq 0$ such that the solutions to \eqref{opt:lc} and \eqref{opt:ll} coincide, with $\lambda$ interpretable as the Lagrange multiplier associated with the constraint $\normone{\bbeta}\leq \Sigma$; see Corollary~\ref{corollary:equigv_p1_epsi_pen_lag}.

For sufficiently large $\lambda$ (or sufficiently small $\Sigma$), any solution   $\widehatbbeta$ to either \eqref{opt:lc} or \eqref{opt:ll} is typically sparse---i.e., it has only a few nonzero entries---and thus identifies the relevant features (columns of $\bX$) that are useful for predicting $\by$. 
With the gradual decrease of $\lambda$, the sparsity of the solution vector $\bbeta$ also gradually decreases. As $\lambda$ tends to 0, the solution vector $\bbeta$ becomes the vector such that $\normtwo{\by - \bX \bbeta}^2$ is minimized. That is to say, $\lambda > 0$ can balance the twin objectives by minimizing the error squared sum cost function $\frac{1}{2} \normtwo{\by - \bX \bbeta}^2$ and the $\ell_1$-norm cost function $\normone{\bbeta}$.
See Section~\ref{section:ls_regular} for further details.

Note that a scaled version of the LASSO is also commonly used, particularly in statistical learning contexts where normalization by the sample size $n$ improves interpretability:
\begin{subequations}\label{equation:lassos_all_scaled}
\begin{align}
(\textbf{Constrained LASSO-II}):\qquad &\mathopmin{\normone{\bbeta}\leq \Sigma}f(\bbeta)\triangleq \frac{1}{2n}\normtwo{\by - \bX \bbeta}^2;
\label{equation:loss_cons_lasso_scaled} 
\tag{L$_{C}'$}\\
(\textbf{Lagrangian LASSO-II}):\qquad &\mathopmin{\bbeta \in \real^p}F(\bbeta)\triangleq \left\{ \frac{1}{2n} \normtwo{\by - \bX \bbeta}^2 + \lambda \normone{\bbeta} \right\}. 
\label{opt:ll_scaled}  
\tag{L$_{L}'$}
\end{align}
\end{subequations}

The LASSO admits numerous generalizations; here are a few examples. In machine learning, sparse kernel regression is sometimes  referred to as the \textit{generalized LASSO} \citep{roth2004generalized}. 
The \textit{group LASSO} extends sparsity to predefined groups of variables \citep{puig2011multidimensional}. 
The \textit{distributed LASSO} addresses sparse linear regression in decentralized or networked settings \citep{mateos2010distributed}.
See Chapter~\ref{chapter:spar} for a more comprehensive discussion of these extensions.

\subsection{Sparse Recovery Problems}

The success of the LASSO in promoting sparsity and enabling feature selection is not an isolated phenomenon; it is deeply rooted in a broader mathematical framework known as \textit{sparse recovery}, which originated in signal processing and compressed sensing. 
As noted above, while the LASSO was developed in a statistical context---where the primary goals are prediction, estimation, and interpretability in the presence of  noise---the same $\ell_1$-regularization principle also arises in deterministic settings,
where the objective  is to exactly reconstruct a sparse signal from incomplete linear measurements. 
This leads naturally to formulations such as \textit{basis pursuit} and its noise-robust extension, \textit{basis pursuit denoising (BPDN)}, which frame sparse recovery as a constrained optimization problem rather than a penalized one. 
Although their motivations differ---statistical inference versus exact or stable signal reconstruction---the underlying geometry of the $\ell_1$-norm and the conditions that guarantee successful recovery (e.g., restricted isometry or nullspace properties; see Chapters~\ref{chapter:design} and \ref{chapter:recovery}) provide a unifying foundation for these approaches. 
This subsection introduce this complementary perspective.

There  exist a wide variety of approaches to recover a sparse  signal $\bbeta$ from a small number of linear measurements. We begin by considering a natural first approach to the problem of sparse recovery.
Given measurements $\by\in\real^n$ (e.g., generated via $\by = \bX\bbeta$,  where $\bX\in\real^{n\times p}$) and the knowledge that our original signal $\bbeta\in\real^p$ is sparse, it is natural to attempt recovery by solving an optimization problem of the form
\begin{equation}\label{equation:spar_rec_gen1}
\widehatbbeta = \argmin_{\bbeta} \normzero{\bbeta} \quad \text{s.t.} \quad \bbeta \in \mathcalB(\by),
\tag{PS$_0$}
\end{equation}
where $\mathcalB(\by)$ enforces consistency between $\bbeta$ and the observed measurements $\by$. 
Recall that $\normzero{\bbeta} = \abs{\supp(\bbeta)}$ simply counts the number of nonzero entries in $\bbeta$; 
thus, \eqref{equation:spar_rec_gen1} seeks the sparsest signal consistent with the data. 
Common choices for $\mathcalB(\by)$ include:
\begin{itemize}
\item If our measurements are \textbf{exact and noise-free}, then we can set $\mathcalB(\by) = \{\bz \mid  \bX \bz = \by\}$. 
\item If the measurements are contaminated by a small amount of \textbf{bounded noise}, we may instead use $\mathcalB(\by) = \{\bz \mid \normtwo{\bX \bz - \by} \leq \epsilon\}$. 
\end{itemize}
In both cases, \eqref{equation:spar_rec_gen1} identifies the sparsest $\bbeta$ that remains consistent with the observations. Further details on these two scenarios are provided below.

The key differences between problem  \eqref{opt:s0_gen} (p.~\pageref{opt:s0_gen}) and \eqref{equation:spar_rec_gen1} are as follows:
\begin{enumerate}[(i)]
\item \textit{Objective function.} In \eqref{opt:s0_gen}, the goal is to minimize a data-fidelity term (e.g., least squares error) subject to a sparsity constraint. In contrast, \eqref{equation:spar_rec_gen1} directly minimizes sparsity (i.e., $\normzero{\bbeta}$) subject to data consistency.

\item \textit{Constraints.} Problem \eqref{opt:s0_gen} imposes a hard constraint on the sparsity level (e.g., $\normzero{\bbeta} \leq R$), allowing flexibility in how well the model fits the data. Conversely, \eqref{equation:spar_rec_gen1} enforces strict data fidelity (exact or approximate), with sparsity as the objective.

\item \textit{Application.} Formulation \eqref{opt:s0_gen} is typically used in statistical settings where data are noisy, and a trade-off between fit and sparsity is desired (e.g., regression). In contrast, \eqref{equation:spar_rec_gen1} arises in signal processing and compressed sensing, where the priority is exact or stable reconstruction of a truly sparse signal from limited measurements.

\item \textit{Computational complexity.} Although both problems involve the nonconvex $\ell_0$-norm, \eqref{opt:s0_gen} can sometimes be tackled with iterative heuristics (e.g., IHT) under additional assumptions. However, \eqref{equation:spar_rec_gen1} is generally NP-hard due to the combinatorial nature of minimizing $\normzero{\bbeta}$, as it requires evaluating all possible support patterns of $\bbeta$ \citep{mallat1993matching, natarajan1995sparse, davis1997adaptive}.

\end{enumerate}

\paragrapharrow{Data transformation.}
Note that in \eqref{equation:spar_rec_gen1}, we implicitly assume that the signal $\bbeta$ itself is sparse. 
In many practical settings, however, sparsity holds not in the standard basis but in some transformed domain, i.e.,  $\bbeta = \bPhi \balpha$, where $\balpha\in\real^p$ is sparse and $\bPhi\in\real^{p\times p}$ is a known dictionary (e.g., a wavelet or Fourier basis). In this case, we can reformulate the recovery problem as 
\begin{equation}\label{equation:spar_rec_gen2}
\widehat{\balpha} = \argmin_{\balpha} \normzero{\balpha} 
\quad \text{s.t.} \quad 
\balpha \in \mathcalB(\by),
\tag{PS$_0'$}
\end{equation}
where $\mathcalB(\by) = \{\bz \mid \bX \bPhi \bz = \by\}$ for noise-free measurements, or $\mathcalB(\by) = \{\bz \mid \normtwo{\bX \bPhi \bz - \by} \leq \epsilon\}$ in the presence of bounded noise. 
By defining $\widetildebX = \bX \bPhi$, we see that \eqref{equation:spar_rec_gen2} is structurally identical to \eqref{equation:spar_rec_gen1}, with $\widetildebX$ replacing $\bX$. 
Moreover, in many applications, the introduction of $\bPhi$ does not substantially complicate the design of $\bX$ such that  
$\widetildebX$ satisfies desirable properties (e.g., restricted isometry).
Therefore, for most of the remainder of this book, we restrict our attention to the case $\bPhi=\bI$, i.e., when sparsity is assumed directly in the canonical basis. It is important to note, however, that this simplification imposes limitations when $\bPhi$ is a general (possibly overcomplete) dictionary rather than an orthonormal basis. In particular,
$$
\normtwobig{\bbeta - \widehatbbeta} 
= \normtwo{\bPhi \balpha - \bPhi \widehat{\balpha}} 
\neq \normtwo{\balpha - \widehat{\balpha}},
$$
so a bound on $\normtwo{\balpha - \widehat{\balpha}}$ does not directly translate into a bound on $\normtwobig{\bbeta - \widehatbbeta}$, which is often the quantity of interest in applications. Further discussion of these implications can be found in \citet{candes2011compressed}.

\paragrapharrow{Convex relaxation.\index{Convex relaxation}}
Although the performance of \eqref{equation:spar_rec_gen1} can be analyzed under suitable assumptions on $\bX$, we do not pursue this route because the objective function $\normzero{\cdot}$ is nonconvex---and consequently, the optimization problem is generally intractable. In fact, for an arbitrary measurement matrix $\bX$, even approximating the optimal solution is NP-hard.
Recall that the number of nonzero entries of a vector $\bbeta \in \real^p$ can be viewed as the limit of its   $\ell_s$-quasinorm as $s\to 0$:
$$
\sum_{i=1}^{p} \abs{\beta_i}^s 
\quad \stackrel{s \to 0}{\longrightarrow} \quad
\sum_{i=1}^{p} \indicator(\{\beta_i \neq 0\}) = \normzero{\bbeta}.
$$
This observation suggests to replace the $\ell_0$-minimization problem \eqref{equation:spar_rec_gen1} with the following family of problems:
\begin{equation}\label{opt:pss}
\widehatbbeta = \argmin_{\bbeta} \norms{\bbeta} \quad \text{s.t.} \quad \bbeta \in \mathcalB(\by), \quad 0<s<1.
\tag{PS$_s$}
\end{equation}
However, this approach has significant drawbacks.
For $s > 1$, the $\ell_s$-norm fails to promote sparsity---even 1-sparse vectors may not be recovered (see Problem~\ref{prob:1spar_ells}). For $0<s<1$, although sparsity is better encouraged, the resulting problem remains non-convex and is still NP-hard (see Problem~\ref{prob:nphard_ells}).
A more practical alternative is to replace $\normzero{\cdot}$ with its convex envelope over the unit $\ell_\infty$-ball (see Lemma~\ref{lemma:l1_ball_intersec} for further details): the $\ell_1$-norm. This leads to the convex program
\begin{equation}\label{opt:ps1}
\widehatbbeta = \argmin_{\bbeta} \normone{\bbeta} \quad \text{s.t.} \quad \bbeta \in \mathcalB(\by).
\tag{PS$_1$}
\end{equation}
Provided that $\mathcalB(\by)$ is convex, \eqref{opt:ps1} is computationally feasible. 
In the noiseless setting $\mathcalB(\by) = \{\bz\mid \bX \bz = \by\}$, it can even be cast as a linear program \citep{chen1994basis, chen2001atomic, foucart2013invitation}.

It is clear that replacing \eqref{equation:spar_rec_gen1} with \eqref{opt:ps1} transforms an intractable combinatorial problem into a tractable convex one. Less obvious, however, is whether the solution to \eqref{opt:ps1} resembles that of \eqref{equation:spar_rec_gen1}. Nevertheless, there are strong intuitive and geometric reasons to expect that $\ell_1$-minimization promotes sparsity.
For instance, as illustrated later in Example~\ref{example:ell0_1_2} and Figure~\ref{fig:l0l1l2_opt}, the solution to the $\ell_1$-minimization problem coincides exactly with the solution to $\ell_s$-minimization for any $s<1$---and is indeed sparse.

Moreover, historical precedent supports this intuition. As early as 1965, \citet{logan1965properties} showed that a bandlimited signal corrupted on a small interval can be perfectly reconstructed by minimizing the $\ell_1$-norm of the error. This method---searching for the bandlimited signal closest to the observations in $\ell_1$---can be interpreted as an early instance of robust recovery via $\ell_1$-minimization, further reinforcing the idea that the $\ell_1$-norm is well-suited to handling sparse errors.

\subsection*{Basis Pursuit and other Forms\index{Basis pursuit}}
In summary,  the core problem of sparse representation is the $\ell_0$-norm minimization
\begin{equation}\label{opt:p0}
\min_{\bbeta} \normzero{\bbeta} \quad \text{s.t.} \quad \by = \bX\bbeta, \tag{P$_0$} 
\end{equation}
where $\bX \in \real^{n\times p}, \bbeta \in \real^p$, and $ \by \in \real^n$.
Since this problem is NP-hard, we instead consider its convex relaxation known as the  \textit{basis pursuit} problem also referred to as \textit{noise-free signal recovery}, \textit{$\ell_1$-norm minimization}, or simply  \textit{$\ell_1$-minimization} problem):
\begin{equation}\label{opt:p1}
\min_{\bbeta\in\real^p} \normone{\bbeta}  \qquad\text{s.t.}  \quad \by = \bX \bbeta. \tag{P$_1$} 
\end{equation}
This is a convex optimization problem because the objective function $\normone{\bbeta}$ is convex and the constraint  $\by = \bX \bbeta$ is affine.
Thus, this basis pursuit problem can be seen as a convex relaxation of \eqref{opt:p0}.
Problem \eqref{opt:p1} is equivalent to 
\begin{equation}\label{equation:prob_p1_equiv}
\begin{aligned}
\min_{\bbeta, \bu} \sum_{i=1}^{p} u_i 
\qquad \text{s.t.} 
&\quad \by = \bX \bbeta;\\
& \quad -u_i \leq \beta_i \leq u_i, \quad \forall\, i\in\{1,2,\ldots,p\}.
\end{aligned}
\tag{P$_1'$}
\end{equation}
This is a linear program formulation (see Example~\ref{example:linear_program}).
While \eqref{opt:p1} enforces an exact equality constraint, it can also be reformulated as a penalized problem:
\begin{equation}\label{opt:p1_penalize}
\min_{\bbeta\in\real^p}   \normone{\bbeta} + \sigma \normtwo{\bX\bbeta-\by}^2,
\tag{P$_{1,\sigma}$}
\end{equation}
where $\sigma>0$ is a \textit{penalty parameter}. 
This formulation is known as the \textit{basis pursuit  denoising (BPDN)} problem \citep{chen2001atomic}. 
As $\sigma\rightarrow\infty$, the solution of \eqref{opt:p1_penalize} converges to that of \eqref{opt:p1}, making the two problems asymptotically equivalent.
\index{Penalty function}
Apparently, the penalized $\ell_1$ problem \eqref{opt:p1_penalize} is equivalent to the Lagrangian LASSO problem~\eqref{opt:ll} (p.~\pageref{opt:ll}) by setting $\sigma \triangleq \frac{1}{2\lambda}$.

\paragrapharrow{Signal recovery vs sparse regression.}
As noted previously, although the formulations of \eqref{opt:p1_penalize} and \eqref{opt:ll} are mathematically equivalent under the correspondence $\sigma \triangleq \frac{1}{2\lambda}$, their interpretations and typical application contexts differ significantly.
In the Lagrangian LASSO problem~\eqref{opt:ll}, the primary objective is to minimize the least squares loss  $\normtwo{\bX\bbeta-\by}$ subject to an $\ell_1$ regularization term $\lambda\normone{\bbeta}$. 
Here, the underlying goal is predictive modeling: we seek a good fit to the data while encouraging sparsity in the estimated coefficients to improve interpretability and prevent overfitting.
When the setting arises in overdetermined systems ($\bX\in\real^{n\times p}$ with $n\geq p$, also the LASSO can also work when $p\geq n$), where the equation $\bX\bbeta=\by$ generally has no exact solution due to noise or model mismatch, the regularization term $\lambda\normone{\bbeta}$ thus serves as a mechanism to control model complexity \citep{lu2021numerical}.

In contrast, the penalized $\ell_1$-minimization problem~\eqref{opt:p1_penalize} originates in sparse signal recovery, particularly in compressed sensing. Here, the goal is not prediction but exact or stable reconstruction of a truly sparse signal from incomplete and possibly noisy measurements. The system is usually underdetermined ($p\gg n$), meaning there are far more unknowns than observations. Sparsity is not a convenience---it is a structural assumption that makes recovery possible.

This distinction has practical implications in communication systems. When $p\gg n$, transmitting the full signal $\bbeta\in\real^p$ is inefficient. Instead, the sender can transmit only the low-dimensional measurement vector $\by=\bX\bbeta\in\real^n$, which requires far less bandwidth since $n\gg p$. If the receiver knows both $\bX$ and $\by$, it can reconstruct the original sparse signal $ \bbeta$ by solving a problem such as \eqref{opt:p1_penalize}.

\begin{figure}[h]
\centering
\includegraphics[width=0.999\textwidth]{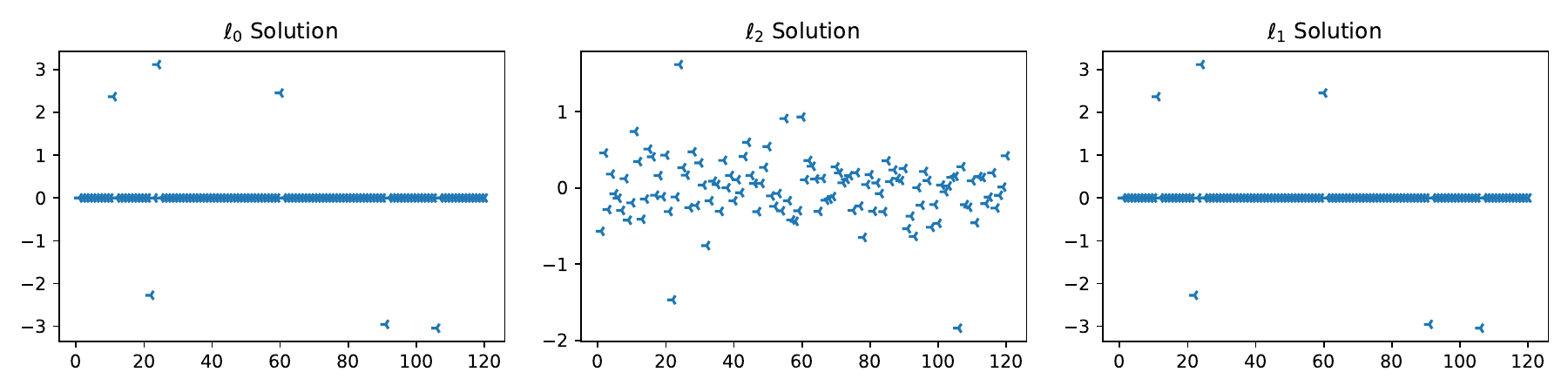}
\caption{Illustration of $\ell_0,\ell_1$, and $\ell_2$-norm minimization for sparse signal recovery.}
\label{fig:l0l1l2_opt}
\end{figure}

\begin{example}[Sparse signal recovery using $\ell_0,\ell_1$, and $\ell_2$-norm minimization]\label{example:ell0_1_2}
We compare the recovery of a sparse signal from underdetermined measurements using $\ell_0$- and $\ell_1$-norm minimization.
We also compare with the $\ell_2$-norm minimization:
\begin{equation}\label{opt:p2}
\min_{\bbeta\in\real^p} \normtwo{\bbeta} \quad \text{s.t.} \quad \by = \bX \bbeta.
\tag{P$_2$}
\end{equation}
The process begins with a random matrix $ \bX $ of size $ 60 \times 120 $ and a sparse vector $ \bbeta^* $ of length 120, where sparsity is enforced by zeroing out small values in a dense random vector. The measurements $ \by $ are then  computed as $ \by = \bX  \bbeta^* $, simulating a compressed sensing scenario where only 60 linear measurements are available for a 120-dimensional sparse  signal.

The first subplot of Figure~\ref{fig:l0l1l2_opt} the solution obtained via $\ell_0$-minimization. 
The second subplot presents the $\ell_2$-minimization solution. Because the $\ell_2$-norm favors energy spreading and smoothness, the recovered signal is dense and distributes error across all coordinates, leading to a large reconstruction error relative to $\bbeta^*$.
The third subplot shows the $\ell_1$-minimization solution. Thanks to its sparsity-promoting geometry, the $\ell_1$-norm yields a solution that closely matches the support and magnitude of the original sparse signal. As expected, the reconstruction error is significantly lower than that of the $\ell_2$ approach, illustrating the effectiveness of $\ell_1$-minimization for sparse recovery.
\end{example}

\subsection*{Constrained Sparse Recovery}
On the other hand, observed signals are typically contaminated by noise. Therefore, the equality constraint in the idealized optimization problem should be relaxed to allow for a bounded measurement error. This leads to the following $\ell_0$-norm minimization with an inequality constraint:
\begin{equation}\label{opt:p0epsilon}
\min_{\bbeta\in\real^p} \normzero{\bbeta} \quad \text{s.t.} \quad \normtwo{\bX\bbeta - \by} \leq \epsilon, 
\tag{P$_{0,\epsilon}$}
\end{equation}
where $\epsilon \geq 0$ controls the permissible deviation between the measurements $\by$ and the model prediction $\bX\bbeta$.

In realistic settings, it is also impossible to measure a signal $\bbeta\in\real^p$ with infinite precision. Consequently, the observed measurement vector $\by\in\real^n$ is only an approximation of the true linear measurement $\bX\bbeta$, satisfying
$$
\norm{\bX\bbeta - \by} \leq \epsilon
$$
for some $\epsilon \geq 0$ and  some norm $\norm{\cdot}$ on $\real^n$ (most commonly the $\ell_2$-norm).
In such cases, a good reconstruction method should produce an estimate  $\bbeta\in\real^p$ whose distance to the true signal $\bbeta$ is controlled by the measurement error $\epsilon$. 
This desirable property is known as \textit{robustness} to measurement noise.

Accordingly, the noiseless basis pursuit problem \eqref{opt:p1} is replaced by the following convex optimization problem:
\begin{equation}\label{opt:p1_epsilon}
\min_{\bbeta \in \real^p} \normone{\bbeta} \quad \text{s.t.} \quad \normtwo{\bX\bbeta - \by} \leq \epsilon,
\tag{P$_{1,\epsilon}$}
\end{equation}
which is commonly referred to as \textit{quadratically-constrained basis pursuit} (or \textit{quadratically-constrained $\ell_1$-minimization}, \textit{$\ell_1$-minimization with noise measurement}) problem and is known as the problem of \textit{signal recovery in noise}.

Another prominent variant is the\textit{Dantzig selector (DS)}\citep{candes2007dantzig}, defined as
\begin{equation}\label{opt:p1_dantzig}
\min_{\bbeta \in \real^p} \normone{\bbeta} \quad \text{s.t.} \quad \norminf{\bX^\top(\bX\bbeta - \by)} \leq D.
\tag{P$_{1,D}$}
\end{equation}
This is again a convex optimization problem. 
The intuition behind the constraint is that the residual $\bepsilon = \bX\bbeta - \by$ should exhibit only small correlations with all columns $\bx_i$ of the measurement  matrix $\bX$. 
Indeed, 
$$
\norminf{\bX^\top(\bX\bbeta - \by)} = \max_{i \in \{1,2,\ldots,p\}} \abs{\innerproduct{ \bepsilon, \bx_i}}.
$$ 
A theoretical framework analogous to the one developed later in this chapter for problems \eqref{opt:p1} and \eqref{opt:p1_epsilon} also applies to the Dantzig selector. However, we do not delve into those details here. For further insight, see Problem~\ref{prob:dantzig_dual_cert}.

\subsection{Equivalence of Some Problems}
The solution $\widehatbbeta$ of the quadratically-constrained basis pursuit problem \eqref{opt:p1_epsilon} is closely related to the solutions of the Lagrangian and constrained LASSO problems in \eqref{opt:ll} and \eqref{opt:lc} (p.~\pageref{opt:lc}). To establish this connection, we first prove a more general result.

\begin{theoremHigh}[Equivalence between constrained and regularized problems]\label{theorem:equiv_cons_reg_gen}
Let $\norm{\cdot}_a$, $\norm{\cdot}_b$ be two norms on $\real^p$. 
Given  $\bX \in \real^{n\times p}$, $\by \in \real^n$ and $\epsilon > 0$, 
consider the optimization problem
\begin{equation}\label{equation:opt1_normab}
\min_{\bbeta \in \real^p} \norm{\bbeta}_a \quad \text{s.t.}\quad \norm{\bX\bbeta - \by}_b \leq \epsilon.
\end{equation}
Assume that \eqref{equation:opt1_normab} is strictly feasible; that is, there exists some vector $\bbeta \in \real^p$ such that $\norm{\bX\bbeta - \by}_b < \epsilon$. 
Let $\widehatbbeta$ be a minimizer of  \eqref{equation:opt1_normab}. Then there exists a parameter $\sigma \geq 0$ such that $\widehatbbeta$ also minimizes the regularized problem
\begin{equation}\label{equation:opt2_normab}
\min_{\bbeta \in \real^p} \norm{\bbeta}_a + \sigma \norm{\bX\bbeta - \by}_b.
\end{equation}
Conversely, for any $\sigma > 0$, if $\widehatbbeta$ minimizes \eqref{equation:opt2_normab}, then there exists $\epsilon \geq 0$ such that $\widehatbbeta$ is a minimizer of \eqref{equation:opt1_normab}.~\footnote{Although we restrict attention to real-valued vectors, the theorem also holds in the complex case by identifying $\complex^p = \real^{2p}$.}
\end{theoremHigh}

\begin{proof}[of Theorem~\ref{theorem:equiv_cons_reg_gen}]
The Lagrange function of problem \eqref{equation:opt1_normab} is
$$
L(\bbeta, \lambda) = \norm{\bbeta}_a + \lambda (\norm{\bX\bbeta - \by}_b - \epsilon).
$$
By Slater's condition (Theorem~\ref{theorem:slater_cond}), strong duality holds for \eqref{equation:opt1_normab} under the stated strict feasibility assumption. 
Therefore, there exists an optimal dual multiplier $\widehatlambda \geq 0$. The {saddle point property \eqref{equation:saddl_point_gen}} implies that 
$$
L(\widehatbbeta, \widehatlambda) \leq L(\bbeta, \widehatlambda)
$$ 
for all $\bbeta \in \real^p$. 
Thus, $\widehatbbeta$ minimizes the function $\bbeta \mapsto L(\bbeta, \widehatlambda)$. Since the constant term $-\widehatlambda \epsilon$ does not affect the location of the minimizer, $\widehatbbeta$ also minimizes 
$$
\norm{\bbeta}_a + \widehatlambda \norm{\bX\bbeta - \by}_b,
$$
which establishes the first claim with $\sigma = \widehatlambda$.

Conversely, let $\widehatbbeta$ be the minimizer of \eqref{equation:opt2_normab} and set $\lambda = \sigma$. Choose $\epsilon = \normbig{\bX\widehatbbeta - \by}_b$. Then the Lagrange dual function $D$ satisfies 
$$
D(\lambda) = L(\widehatbbeta, \lambda) = \normbig{\widehatbbeta}_a.
$$ 
By weak duality, it holds that $D(\lambda) \leq \norm{\bbeta}_a$ for all feasible $\bbeta$. Because $\widehatbbeta$ is feasible and achieves this lower bound by the choice of $\epsilon= \normbig{\bX\widehatbbeta - \by}_b$, it must be a minimizer of \eqref{equation:opt1_normab}.
\end{proof}

This general equivalence immediately clarifies the relationships among the problems \eqref{opt:p1_epsilon}, \eqref{opt:p1_penalize}, \eqref{opt:lc}, and \eqref{opt:ll}.
\begin{corollary}[Equivalence of \eqref{opt:p1_epsilon}, \eqref{opt:p1_penalize}, \eqref{opt:lc}, and \eqref{opt:ll}]\label{corollary:equigv_p1_epsi_pen_lag}
Consider the problems \eqref{opt:p1_epsilon}, \eqref{opt:p1_penalize}, \eqref{opt:lc}, and \eqref{opt:ll}:
\begin{enumerate}[(i)]
\item As mentioned above, the penalized $\ell_1$ problem \eqref{opt:p1_penalize} is equivalent to the Lagrangian LASSO problem~\eqref{opt:ll} by setting $\sigma \triangleq \frac{1}{2\lambda}$.
\item If $\widehatbbeta$ is a minimizer of Lagrangian LASSO problem \eqref{opt:ll}, then there exists $\epsilon = \epsilon_{\widehatbbeta}$ such that $\widehatbbeta$ is a minimizer of the quadratically-constrained basis pursuit problem \eqref{opt:p1_epsilon}.

\item If $\widehatbbeta$ is a minimizer of the constrained LASSO problem \eqref{opt:lc}, then there exists  $\lambda = \lambda_{\widehatbbeta}$ such that $\widehatbbeta$ is a minimizer of the Lagrangian LASSO problem \eqref{opt:ll}.
\item If $\widehatbbeta$ is a minimizer of quadratically-constrained basis pursuit \eqref{opt:p1_epsilon}, then there is $\Sigma = \Sigma_{\widehatbbeta}$ such that $\widehatbbeta$ is a minimizer of the constrained LASSO problem \eqref{opt:lc}.
\end{enumerate}
\end{corollary}

It is important to note that these equivalences are implicit: the mapping between parameters (e.g., between $\lambda$ and $\epsilon$ or $\lambda$ and $\Sigma$) depends on the specific minimizer $\widehatbbeta$. Consequently, the transformation cannot be performed a priori---it requires knowledge of the solution itself. For this reason, while the equivalences are theoretically insightful, they are of limited practical use when selecting regularization parameters in real-world applications.

\subsection{$\ell_1$-Analysis Models\index{$\ell_1$-analysis}}
Note that the  models \eqref{opt:p1}, \eqref{opt:p1_penalize}, and \eqref{opt:p1_epsilon} are collectively  known as the \textit{$\ell_1$-synthesis models}.
We also consider the following $\ell_1$-minimization formulations:
\begin{subequations}\label{opt:q1_analysis}
\begin{align}
&\min_{\bbeta\in\real^p} \normone{\bPsi^\top \bbeta}, \quad \text{s.t.}\quad \bX \bbeta = \by, 
 \label{opt:q1} \tag{Q$_{1}$}\\
&\min_{\bbeta\in\real^p}  \normone{\bPsi^\top \bbeta}  + \sigma \normtwo{\bX \bbeta - \by}^2,  
\label{opt:q1_penalize} \tag{Q$_{1,\sigma}$} \\
&\min_{\bbeta\in\real^p} \normone{\bPsi^\top \bbeta}, \quad \text{s.t.}\quad \normtwo{\bX\bbeta-\by} \leq \epsilon, 
\label{opt:q1_epsilon} \tag{Q$_{1,\epsilon}$}
\end{align}
\end{subequations}
where $\epsilon, \sigma$ are positive parameters. 
Model \eqref{opt:q1} is appropriate in the noiseless setting, i.e., when the measurement error satisfies $\bepsilon = \bzero$. 
If $\bPsi = \bI$, the identity matrix, these models reduce exactly to the $\ell_1$-synthesis models. 
If $\bPsi \neq \bI$, they are referred to as the \textit{$\ell_1$-analysis models}  \citep{elad2007analysis, candes2011compressed, liu2012compressed, nam2013cosparse, zhang2016one}.

It is important to distinguish these models from the data transformation described in \eqref{equation:spar_rec_gen2} (p.~\pageref{equation:spar_rec_gen2}): in the analysis framework, the transformation is applied to the signal itself, not to the measurements.

In the $\ell_1$-synthesis approaches, the signal of interest is synthesized as $\bbeta^* = \bD\bb$, where $\bD$ is a certain dictionary and $\bb$ is a sparse coefficient vector. 
The $\ell_1$-analysis model, including the cosparse analysis model~\citep{nam2013cosparse} and the total variation model~\citep{rudin1992nonlinear}, has widely known examples and has  attracted a lot of attention. The underlying signal is expected to make sparse correlations with the columns (atoms) in an (overcomplete) dictionary $\bPsi\in\real^{p\times q}$ with $q \geq p$, i.e., $\bPsi^\top \bbeta^*$ is sparse.

Typically, we require the analysis operator $\bPsi$ to have full row rank such that $\norm{\cdot}_{\bPsi^\top} \triangleq \normonebig{\bPsi^\top\cdot}$ is a valid norm for $\bbeta$ (Problem~\ref{prob:qnorm}).
When $\bPsi$ has full row rank, we can invoke Theorem~\ref{theorem:equiv_cons_reg_gen} with $\norm{\cdot}_a \triangleq \norm{\cdot}_{\bPsi^\top}$ to prove the equivalence between problems \eqref{opt:q1_penalize} and \eqref{opt:q1_epsilon}.
\begin{corollary}[Equivalence of \eqref{opt:q1_penalize} and \eqref{opt:q1_epsilon}]
Suppose $\bPsi$ has full row rank and 
consider the problems \eqref{opt:q1_penalize} and \eqref{opt:q1_epsilon}:
\begin{enumerate}[(i)]
\item If $\widehatbbeta$ is a minimizer of problem \eqref{opt:q1_penalize}, then there exists $\epsilon = \epsilon_{\widehatbbeta}$ such that $\widehatbbeta$ is a minimizer of the problem \eqref{opt:q1_epsilon}.

\item If $\widehatbbeta$ is a minimizer of the problem \eqref{opt:q1_epsilon}, then there exists  $\sigma = \sigma_{\widehatbbeta}$ such that $\widehatbbeta$ is a minimizer of the problem \eqref{opt:q1_penalize}.
\end{enumerate}
\end{corollary}

\section{Optimality Condition of Least Squares Problems}\label{section:optcd_ls}

To introduce optimality conditions for sparse optimization problems, we begin with the basic least squares problem: $\min_{\bbeta}\normtwo{\bX\bbeta-\by}^2$. 
For this problem, we consider the following questions:
\begin{itemize}
	\item Q1: What is the least squares solution?
	\item Q2: When can uniqueness of the least squares solution be claimed?
	\item Q3: When can uniqueness of the least squares solution with minimum-norm be claimed?
\end{itemize}

To address these questions, we rely on the following fundamental result characterizing the optimality condition for the least squares problem:
\begin{equation}\label{equation:sca_nor_ls_eq}
\bX^\top (\by - \bX\bbeta) = \bzero.
\end{equation}
This condition arises because the objective function $f(\bbeta)=\normtwo{\bX\bbeta-\by}^2$ is convex and differentiable. Therefore, any minimizer $\widehatbbeta$ must satisfy the first-order optimality condition $\nabla f(\widehatbbeta) = 2\bX^\top(\by-\bX\bbeta)=\bzero$ (Corollary~\ref{corollary:fermat_fist_opt}).

We now present a unified perspective on the least squares problem in the following theorem. The underlying ideas will become clearer as we proceed. This theorem answers Question Q1.
\begin{theoremHigh}[A unified view of least squares problems \citep{lu2021rigorous}\index{Normal equation}]\label{theorem:unif_ls}
Let $\bX\in\real^{n\times p}$ and $\by\in\real^n$. 
Then the least squares problem $f(\bbeta)=\normtwo{\bX\bbeta-\by}^2$ admits  a minimizer $\widehatbbeta\in\real^n$ if and only if there exists a vector $\balpha\in\real^p$ such that 
\begin{equation}\label{equation:unif_ls}
\widehatbbeta=\bX^+\by+(\bI-\bX^+\bX)\balpha,
\end{equation}
where $\bX^+$ denotes the  pseudo-inverse of $\bX$ (see Section~\ref{section:subsec_pseudo_inv}):
\begin{itemize}
\item This shows that the least squares problem has a \textbf{unique} minimizer of $\widehatbbeta=\bX^+\by$ only when $\bX^+$ is a left inverse of $\bX$ (i.e., $\bX^+\bX=\bI$, and $\bX$ is left-invertible only when $\bX$ has full column rank \citep{lu2021numerical}). 

\item The optimal value is $f(\widehatbbeta)=\by^\top(\bI-\bX\bX^+)\by$.
\item If $\balpha\neq \bzero$: $\normtwo{\bX^+\by}\leq \normtwo{\bX^+\by+(\bI-\bX^+\bX)\balpha}$.
\end{itemize}
\noindent
This means that any vector $\bbeta$ that minimizes $f(\bbeta)$ must be of the form given in \eqref{equation:unif_ls}, where:
\begin{itemize}
\item $\bX^+\by\perp \nspace(\bX)$  (by Theorem~\ref{theorem:pseudo-four-basis-space}) is the {particular solution} (the minimum-norm solution).
\item $(\bI - \bX^+ \bX)\balpha \in \nspace(\bX)$  is the {homogeneous solution}, which captures the degrees of freedom in $\bbeta$ arising from the null space of $\bX$.
\end{itemize}
\end{theoremHigh}
\begin{proof}[of Theorem~\ref{theorem:unif_ls}]
From \eqref{equation:sca_nor_ls_eq}, any minimizer	 $\widehatbbeta$ must satisfy the normal equation: $\bX^\top(\by-\bX\bbeta)=\bzero$.
This equation means that the vector $\by-\bX\bbeta$ is {orthogonal to the column space of $\bX$}, i.e., the error $\by-\bX\bbeta$ lies in the {null space of $\bX^\top$}:
\begin{equation}
\by-\bX\bbeta \perp \cspace(\bX)
\qquad\text{and} \qquad 
\by-\bX\bbeta\in \nspace(\bX^\top).
\end{equation}
To solve for $\bbeta$, note that the system $\bX\bbeta =\by$ may not have an exact solution when $\by\notin \cspace(\bX)$. 
In such cases, we seek the {minimum-norm} solution that minimizes  $\normtwo{\by - \bX\bbeta}$.
The {pseudo-inverse} $\bX^+$ provides the best possible solution by giving the unique {minimum-norm least squares solution} (see the argument in the sequel):
\begin{equation}
\bbeta_{\text{particular}} = \bX^+\by.
\end{equation}
To see $\bbeta_{\text{particular}}$ satisfies \eqref{equation:sca_nor_ls_eq}, we have 
$\bX^\top\bX\bbeta_{\text{particular}} =\bX^\top\bX\bX^+\by =\bX^\top\by$, where we used the fact that $\bX^\top\bX\bX^+=\bX^\top$ (Lemma~\ref{lemma:xtxxpllus_pseudo}).

However, this is only \emph{one} solution among possibly many. To characterize the full solution set, we consider two cases.
\paragraph{Case 1: $\by\in\cspace(\bX)$.}
For this case, the solution in \eqref{equation:unif_ls} is obvious by the properties ($\bX(\bI-\bX^+\bX)\balpha=\bzero$) and uniqueness of the pseudo-inverse of a matrix (Proposition~\ref{proposition:existence-of-pseudo-inverse}).

\paragraph{Case 2: $\by\notin \cspace(\bX)$.} 
Define the residual $\be\triangleq \by-\bX\bbeta_{\text{particular}} = \by-\bX\bX^+\by$. We seek all vectors $\widetildebbeta$ such that $\normtwobig{\by-\bX\widetildebbeta}^2 = \normtwo{\be}^2$. 
Two algebraic possibilities arise:
$$
\bX\widetildebbeta-\by = \be
\qquad \text{and}\qquad 
\by-\bX\widetildebbeta = \be.
$$
The former scenario implies that $\be= \by-\bX\bbeta_{\text{particular}} = \bX\widetildebbeta-\by$, which leads to $\by = \bX\frac{\widetildebbeta+\bbeta_{\text{particular}}}{2}$, i.e., $\by\in\cspace(\bX)$, contradicting the assumption of Case 2. Hence, only the second scenario is valid.

Since $\bbeta$ is in $\real^p$, the complete set of minimizers consists of the particular solution $\bbeta_{\text{particular}} = \bX^+\by$ plus any vector in the {null space of $\bX$} (i.e., any $\bbeta_{\text{null}}$ that satisfies $\bX\bbeta_{\text{null}} = \bzero$, $\bbeta_{\text{null}}\in	 \nspace(\bX)\triangleq \{\bz \mid \bX\bz = \bzero\}$).
Thus, the {general solution} for $\bbeta$ must then take the form:
\begin{equation}
\bbeta = \bX^+\by + \bbeta_{\text{null}}, \qquad \text{with}\quad \bbeta_{\text{null}}\in\nspace(\bX).
\end{equation}
Since $\bX-\bX\bX^+\bX = \bzero$ by \eqref{equation:pseudi-four-equations}, the projection onto the {null space of $\bX$} is given by:
\begin{equation}
\bbeta_{\text{null}} = (\bI - \bX^+ \bX)\balpha, \quad \text{for some }\balpha \in \real^p.
\end{equation}
This follows because $(\bI - \bX^+ \bX)$ is an orthogonal {projection matrix} onto $\nspace(\bX)$ (see Sections~\ref{section:ortho_proj_mat} and \ref{section:subsec_pseudo_inv}).

On the other hand, since $\bbeta_{\text{null}} \in \nspace(\bX)$ and $\cspace(\bX^+) \equiv \cspace(\bX^\top)$ (see Theorem~\ref{theorem:pseudo-four-basis-space}), it follows that
$$
\normtwo{\bX^+\by+(\bI-\bX^+\bX)\balpha}^2 = \normtwo{\bX^+\by}^2 + \normtwo{(\bI-\bX^+\bX)\balpha}^2.
$$
This shows that $\normtwo{\bX^+\by}\leq \normtwo{\bX^+\by+(\bI-\bX^+\bX)\balpha}$ if $\balpha\neq \bzero$ and completes the proof.
\end{proof}

\subsection{LS in Convex Optimization}\label{section:ls_cvx}
The theorems on convex functions provide rigorous answers to Questions Q2 and Q3 posed at the beginning of this section.
Theorem~\ref{theorem:unif_ls} partially addresses Question Q2---concerning the uniqueness of the least squares solution---in the context of the large-sample setting.
Specifically, if $n>p=\rank(\bX)$, the least squares solution $\widehatbbeta=(\bX^\top\bX)^{-1}\bX^\top\by$ is unique. 
In this case, the matrix 
$$
\bH\triangleq \bX(\bX^\top\bX)^{-1}\bX^\top
$$ 
is is known as the \textit{hat matrix}, which orthogonally projects $\by$ onto the column space of $\bX$ (see Section~\ref{section:ortho_proj_mat}).
Indeed, this is the \textbf{only} scenario in which the least squares solution is guaranteed to be unique.

Note that both Q2 and Q3 can be formulated as convex optimization problems:
\begin{equation}
\begin{aligned}
\text{(P2)}: \qquad &\min_{\bbeta\in\real^p} f(\bbeta) = \normtwo{\bX\bbeta-\by}=\bbeta^\top\bX^\top\bX\bbeta - 2\by^\top\bX\bbeta + \by^\top\by;\\
\text{(P3)}:\qquad &\min_{\bbeta\in\sQ_2} g(\bbeta) = \normtwo{\bbeta}^2, \quad \sQ_2 \triangleq \{\bbeta\in \real^p\mid \normtwo{\bX\bbeta-\by}^2 =\min \}.
\end{aligned}
\end{equation}
Here, both $f(\bbeta)$ and $g(\bbeta)$ are  convex functions (Exercise~\ref{exercise:conv_quad}), and the set $\sQ_2$ is convex  (see Theorem~\ref{theorem:local_glo_cvx}). 
Therefore, both (P2) and (P3) are convex optimization problems.
Theorem~\ref{theorem:local_glo_cvx} also states that if the function is strictly convex, then the minimizer is unique.
Apparently, $g(\bbeta)$ is strictly convex (this again confirms Q3). 
And $f(\bbeta)$ is strictly convex only when $\bX^\top\bX$ is positive definite, i.e., when $\bX$ has full column rank (i.e., $n>p=\rank(\bX)$).
This provides a rigorous answer to question Q2.

Combining this analysis with Theorem~\ref{theorem:unif_ls}, we establish the uniqueness of the minimum-norm least squares solution in the following corollarys.
\begin{corollary}[Minimum-norm and unique LS solutions]\label{corollary:minimum_norm_ls}
Let $\bX \in \real^{n \times p}$ and $\by \in \real^n$. 
\begin{itemize}
\item Let $\sQ_2 \subset \real^p$ denote the set of all minimizers of $\bbeta \mapsto \normtwo{\bX\bbeta - \by}$. The optimization problem
\begin{equation}\label{equation:minimum_norm_ls_eq1}
	\min_{\bbeta \in \sQ_2} \normtwo{\bbeta}^2
\end{equation}
has the unique solution $\bbeta^+ = \bX^+ \by$.
\item If $n \geq p$ and $\bX$ has  full column rank, then the least squares problem
\begin{equation}\label{equation:minimum_norm_ls_eq2}
\min_{\bbeta \in \real^p} \normtwo{\bX\bbeta - \by}^2
\end{equation}
has the unique solution $\widehatbbeta = \bX^+ \by$.
\end{itemize}
\end{corollary}

Note that the minimizer of \eqref{equation:minimum_norm_ls_eq2} corresponds to the orthogonal projection of $\by$ onto the column space of $\bX$. 
Consequently, the matrix $\bX\bX^+$ represents the orthogonal projector onto $\bX$. 
When $\bX$ has full rank and $n \geq p$, the matrix $\bX^\top \bX$ is invertible. 
By \eqref{equation:pseudo_svd_case1}, we have  $\bX^+ = (\bX^\top \bX)^{-1} \bX^\top$. Therefore, $\widehatbbeta = \bX^+ \by$ satisfies the normal equation
\begin{equation}\label{equation:ne_firstkind}
	\bX^\top \bX \widehatbbeta = \bX^\top \by.
\end{equation}

\begin{corollary}[High-dimensional LS solution, $\ell_2$-minimization]\label{corollary:high_dim_ls}
Let $\bX \in \real^{n \times p}$ with $p \geq n$, assume $\bX$ has full rank, and let $\by \in \real^n$. 
Then the constrained optimization problem
\begin{equation}\label{equation:high_dim_ls}
\min_{\bbeta \in \real^p} \normtwo{\bbeta} \quad \text{s.t.} \quad  \bX\bbeta = \by
\end{equation}
has the unique solution $\widehatbbeta = \bX^+ \by$.
\end{corollary}

In this setting, by \eqref{equation:pseudo_svd_case2}, we have $\bX^+ = \bX^\top (\bX \bX^\top)^{-1}$ if $\bX$ has full rank and $p \geq n$; thus, the system $\bX\bbeta=\by$ is consistent.
Therefore, in the situation of the above corollary, the minimizer $\widehatbbeta$ of \eqref{equation:high_dim_ls} satisfies the \textit{normal equation of the second kind}
\begin{equation}\label{equation:ne_secondkind}
\widehatbbeta = \bX^\top \bgamma, \quad \text{where } \bX\bX^\top \bgamma = \by.
\end{equation}

Finally, although the normal equations \eqref{equation:ne_firstkind} and \eqref{equation:ne_secondkind} suggest that least squares problems can be solved by standard linear system solvers (e.g., Gaussian elimination), it is numerically advantageous to use specialized algorithms designed for least squares. These include methods based on Cholesky factorization, QR decomposition, or singular value decomposition (SVD), which offer greater stability and accuracy \citep{lu2021rigorous, bjorck2024numerical}.

\index{Weighted LS}
\subsection{Weighted LS Problems}
Corollary~\ref{corollary:high_dim_ls} introduces the unique solution of the $\ell_2$-minimization problem.
We can also consider the \textit{weighted $\ell_2$-minimization problem}:
\begin{equation}\label{opt:l2_weighted}
\min_{\balpha \in \real^p} \norm{\balpha}_{2,\bw}
\triangleq \left( \sum_{i=1}^p \alpha_i^2 w_i \right)^{1/2} 
\quad \text{s.t.}\quad \bX\balpha = \by,
\tag{L$_{2,\bw}$}
\end{equation}
where $\bw $ is a vector of positive weights ($w_i > 0$ for all $i$). 
Introduce the diagonal weight matrix $\bD_{\bw} = \diag([w_1, w_2, \ldots, w_p]) \in \real^{p \times p}$, and define the change of variables 
$$
\bbeta \triangleq \bD_{\bw}^{1/2} \balpha.
$$
Then $\norm{\balpha}_{2,\bw}^2 = \normtwo{\bbeta}^2$, and the constraint becomes $\bX\bD_{\bw}^{-1/2}\bbeta=\by$.
Thus, the minimizer $\widehatbalpha$ of \eqref{opt:l2_weighted} corresponds to the minimizer $\widehatbbeta$ of the standard minimum-norm problem:
$$
\min_{\bbeta \in \real^p} \normtwo{\bbeta} \quad \text{s.t.}\quad \bX\bD_{\bw}^{-1/2} \bbeta = \by,
$$
By Corollary~\ref{corollary:minimum_norm_ls}, the solution is
\begin{equation}\label{equation:weightls_sol1}
\widehatbalpha = \bD_{\bw}^{-1/2} \widehatbbeta = \bD_{\bw}^{-1/2} (\bX\bD_{\bw}^{-1/2})^+ \by.
\end{equation}
In particular, if $p \geq n$ and $\bX$ has full rank, then Corollary~\ref{corollary:high_dim_ls} implies that the pseudo-inverse simplifies, yielding
\begin{equation}\label{equation:weightls_sol2}
\widehatbalpha = \bD_{\bw}^{-1} \bX^\top (\bX\bD_{\bw}^{-1} \bX^\top)^{-1} \by.
\end{equation}

Theorem~\ref{theorem:unif_ls} shows that the standard least squares solution $\widehatbbeta=\bX^+\by$ is orthogonal to $\nspace(\bX)$. Similar result can be observed in the weighted LS case.
\begin{theoremHigh}[Orthogonal relation of LS-III]\label{theorem:orthogonal_weight_lsii}
A vector $\widehatbalpha \in \real^p$ is a minimizer of \eqref{opt:l2_weighted} if and only if
\begin{equation}\label{equation:orthogonal_weight_lsii_eq1}
\innerproduct{\widehatbalpha, \bn}_{\bw} = 0, \quad \text{for all } \bn \in \nspace(\bX),
\end{equation}
where $\innerproduct{\widehatbalpha, \bn}_{\bw} \triangleq \sum_{i=1}^p \beta_i n_i w_i$.
\end{theoremHigh}
\begin{proof}[of Theorem~\ref{theorem:orthogonal_weight_lsii}]
Let $\widehatbalpha$ satisfy $\bX\widehatbalpha = \by$. 
Any feasible vector  $\balpha \in \real^p$  for \eqref{opt:l2_weighted} can be written as   $\balpha = \widehatbalpha + \bn$ with $\bn \in \nspace(\bX)$. 
For any scalar $\eta \in \real$ and $\bn \in \nspace(\bX)$, consider the perturbed vector $\widehatbalpha + \eta\bn$.
Its weighted norm squared is
\begin{equation}\label{equation:orthogonal_weight_lsii_pv1}
\norm{\widehatbalpha + \eta\bn}_{2,\bw}^2 
=\norm{\widehatbalpha}_{2,\bw}^2 + \eta^2 \norm{\bn}_{2,\bw}^2 + 2\eta \innerproduct{\widehatbalpha, \bn}_{\bw}.
\end{equation}
Therefore, if $\innerproduct{\widehatbalpha, \bn}_{\bw} = 0$ for all $\bn\in\nspace(\bX)$, then $\eta = 0$ is the minimizer of $\eta \mapsto \norm{\widehatbalpha + \eta\bn}_{2,\bw}$. Hence, \eqref{equation:orthogonal_weight_lsii_eq1} implies that $\widehatbalpha$ is a minimizer of \eqref{opt:l2_weighted}. 

Conversely, suppose $\widehatbalpha$ is a minimizer of \eqref{opt:l2_weighted}. 
Then $\eta = 0$ is a minimizer of $\eta \mapsto \norm{\widehatbalpha + \eta\bn}_{2,\bw}$ for all $\bn \in \nspace(\bX)$. However, if $\innerproduct{\widehatbalpha, \bn}_{\bw}$ would be nonzero, then by \eqref{equation:orthogonal_weight_lsii_pv1}, we could find a nonzero $\eta$ sufficiently close to $0$ and of opposite sign as $\innerproduct{\widehatbalpha, \bn}_{\bw}$ such that $\norm{\widehatbalpha + \eta\bn}_{2,\bw} < \norm{\widehatbalpha}_{2,\bw}$, which leads to a contradiction to $\widehatbalpha$ being a minimizer.
\end{proof}

\subsection{Regression and other Forms}\label{section:ls_regular}
Since the least squares model is commonly used for regression problems, it typically assumes that the number of samples exceeds the number of features ($n > p$).
In this setting, if the design matrix $\bX\in\real^{n\times p}$ has full column rank, the least squares solution is given by
$$
\bbeta^{OLS}
\triangleq 
(\bX^\top\bX)^{-1}\bX^\top\by,
$$
which is widely known as the \textit{ordinary least squares (OLS)} solution.

\index{Contion number}
\index{Tikhonov regularization}
\index{$\ell_2$-regularization}
A common issue with the ordinary least squares solution arises when $\bX$ is nearly singular.
Let the SVD of $\bX$ be $\bX=\bU\bSigma\bV^\top\in\real^{n\times p}$, where $\bU\in\real^{n\times n}$ and $\bV\in\real^{p\times p}$ are orthogonal, and the main diagonal of $\bSigma\in\real^{n\times p}$ contains the singular values $\sigma_1\geq \sigma_2\geq \ldots\geq \sigma_p\geq 0$.
Consequently, $\bX^\top\bX = \bV(\bSigma^\top\bSigma)\bV^\top \triangleq \bV\bS\bV^\top$, where $\bS\triangleq \bSigma^\top\bSigma  = \diag(\sigma_1^2, \sigma_2^2, \ldots,\sigma_p^2)\in\real^{p\times p}$ contains the squared singular values of $\bX$. When $\bX$ is nearly singular, the smallest singular value satisfies $\sigma_p^2\approx 0$, which renders the inversion  $(\bX^\top\bX)^{-1} = \bV\bS^{-1}\bV^\top$ numerically unstable. 
As a result, the solution $\bbeta^{OLS} =(\bX^\top\bX)^{-1}\bX^\top\by $ may become highly sensitive to perturbations or even diverge.
To mitigate this instability, an $\ell_2$-regularization term is often added, leading to the following regularized optimization problem:
\begin{equation}
	\bbeta^{Tik} = \mathop{\argmin}_{\bbeta} \normtwo{\by-\bX\bbeta}^2 +\lambda\normtwo{\bbeta}^2.
\end{equation}
This approach is known as \textit{Tikhonov regularization} (or the \textit{ridge regression}) \citep{tikhonov1963solution}.
The gradient of the objective function is $2(\bX^\top\bX+\lambda\bI)\bbeta-2\bX^\top\by$, and setting it to zero yields the closed-form solution
$
\bbeta^{Tik} = (\bX^\top\bX+\lambda\bI)^{-1}\bX^\top\by.
$
Using the SVD of $\bX$, the inverse becomes $(\bX^\top\bX+\lambda\bI)^{-1} = \bV(\bS+\lambda\bI)^{-1}\bV^\top$, where $\widetildebS\triangleq(\bS+\lambda\bI)=\diag(\sigma_1^2+\lambda, \sigma_2^2+\lambda, \ldots,\sigma_p^2+\lambda)$. 
Thus, the OLS and Tikhonov-regularized solutions can be expressed as
\begin{subequations}
\begin{align}
\bbeta^{OLS} &= (\bX^\top\bX)^{-1}\bX^\top\by = \bV\left(\bS^{-1}\bSigma\right)\bU^\top\by;\\
\bbeta^{Tik} &= (\bX^\top\bX+\lambda\bI)^{-1}\bX^\top\by = \bV\left((\bS+\lambda\bI)^{-1}\bSigma\right)\bU^\top\by,
\end{align}
\end{subequations}\footnote{By the push-through identity (Problem~\ref{prob:push_through_ide}), the Tikhonov solution can also be written as $\bbeta^{Tik} = \bX^\top(\bX\bX^\top+\lambda\bI)^{-1}\by$.}
where the main diagonals of $\left(\bS^{-1}\bSigma\right)$ are $\diag(\frac{1}{\sigma_1}, \frac{1}{\sigma_2}, \ldots, \frac{1}{\sigma_p})$; and the main diagonals of $\left((\bS+\lambda\bI)^{-1}\bSigma\right)$ are $\diag(\frac{\sigma_1}{\sigma_1^2+\lambda}, \frac{\sigma_2}{\sigma_2^2+\lambda}, \ldots, \frac{\sigma_p}{\sigma_p^2+\lambda})$. 
The latter expression remains well-behaved even when $\sigma_p\approx 0$, provided that $\lambda$ is chosen larger than the smallest nonzero squared singular value. This regularization also improves the condition number of the system (see Remark~\ref{remark:cond_number}): when $\sigma_p$ is close to zero,
$$
\kappa(\bX^\top\bX) = {\sigma_1^2}/{\sigma_p^2}
\qquad \rightarrow \qquad
\kappa(\bX^\top\bX+\lambda\bI) = {(\lambda+\sigma_1^2)}/{(\lambda+\sigma_p^2)}.
$$
In summary, Tikhonov regularization effectively prevents divergence in the least squares estimate $\bbeta^{OLS}$ when $\bX$ is ill-conditioned or rank-deficient. This stabilization enhances the numerical convergence of least squares algorithms and their variants---such as alternating least squares \citep{lu2022matrix}---making Tikhonov regularization a widely adopted technique in practice.

\begin{exercise}
Use SVD to show that the optimal solution $\bbeta^{Tik} = (\bX^\top\bX+\lambda\bI)^{-1}\bX^\top\by$
has a norm that is non-increasing as $\lambda$ increases.
\end{exercise}

\begin{figure}[h!]
\centering  
\vspace{-0.35cm}  
\subfigtopskip=2pt  
\subfigbottomskip=2pt  
\subfigcapskip=-5pt  
\subfigure[$\ell_2$-constrained loss curve.]{\label{fig:l2_ball_loss}
\includegraphics[width=0.43\linewidth]{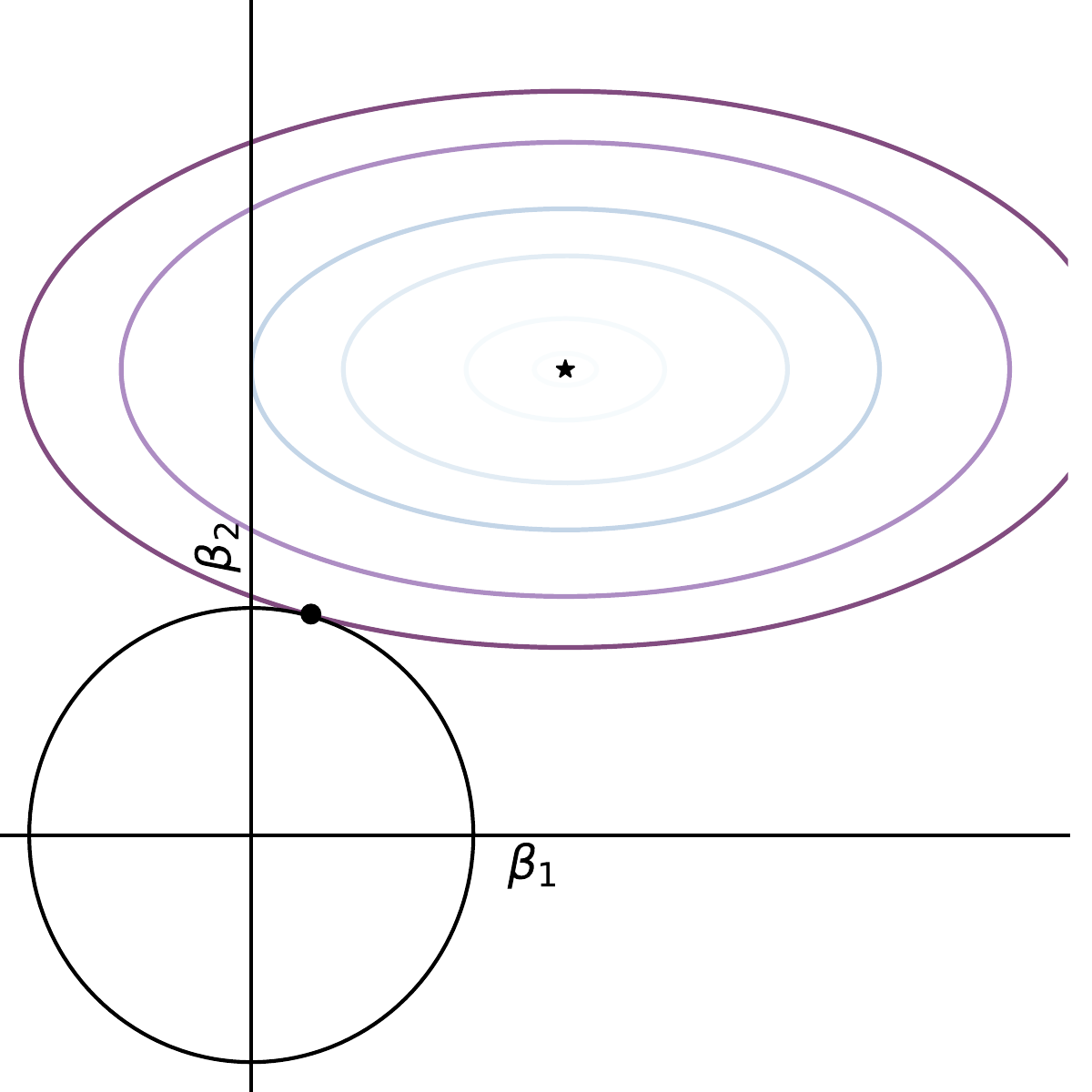}}
\quad 
\subfigure[$\ell_1$-constrained loss curve.]{\label{fig:l1_ball_loss}
\includegraphics[width=0.43\linewidth]{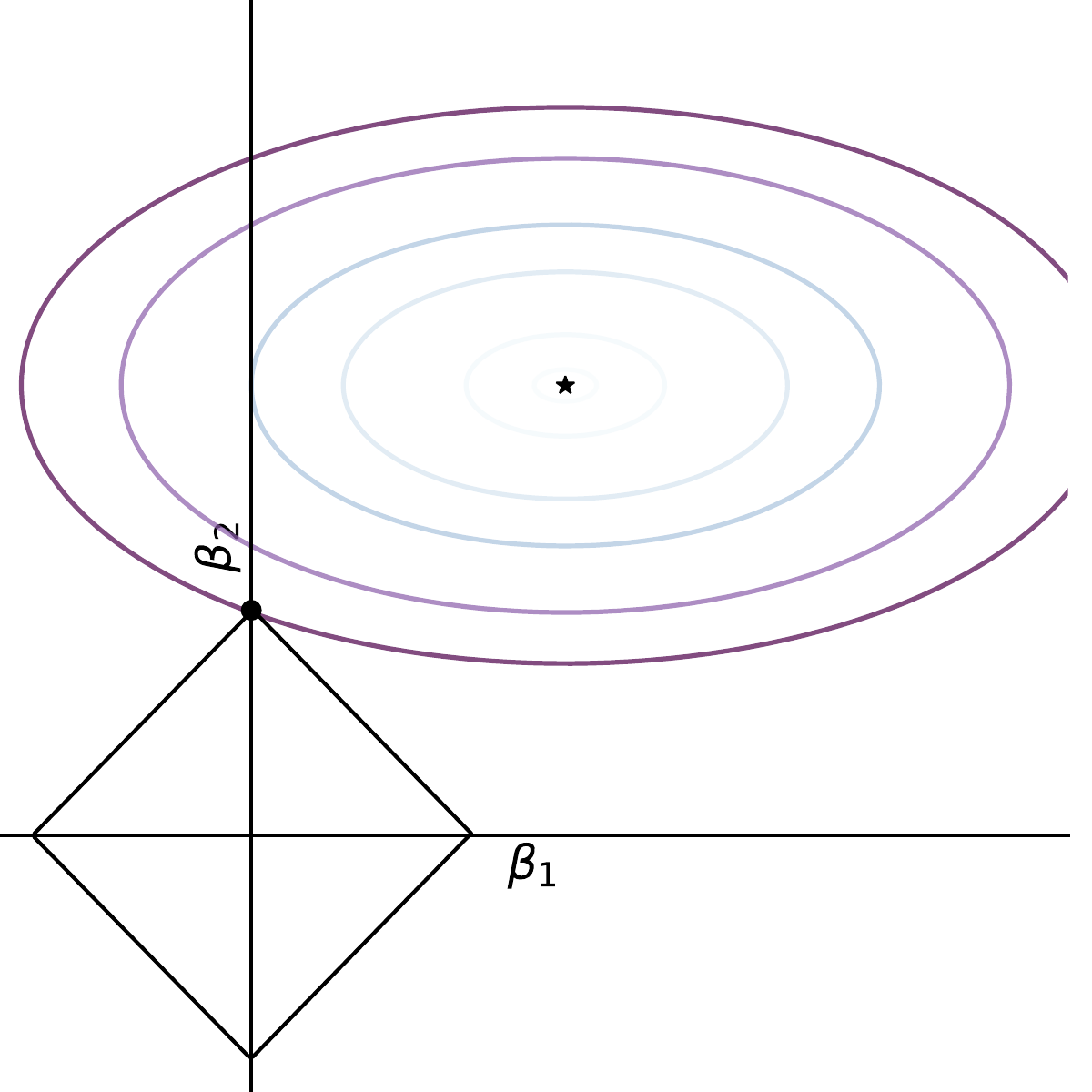}}
\caption{Comparison between the $\ell_2$-constrained and $\ell_1$-constrained least squares problems. The colored elliptical contour lines represent the objective function with respect to the parameters $\bbeta = [\beta_1,\beta_2]^\top$, assuming no constraints.}
\label{fig:l1l2_ball_loss}
\end{figure}

\paragrapharrow{$\ell_1$-regularization.}

\index{Regularization}
\index{Constraint}
\index{Sparsity}

We have shown in Corollary~\ref{corollary:equigv_p1_epsi_pen_lag} that the $\ell_2$-regularized model is equivalent to the $\ell_2$-constrained counterpart.
Figure~\ref{fig:l2_ball_loss}  illustrates the geometric interpretation  of the  $\ell_2$-constrained model for the case $p=2$. 
In this setting, the estimated coefficients are shrunk toward zero but are never exactly zero.
In the figure, the colored elliptical contour lines depict the unconstrained least squares objective as a function of  $\bbeta = [\beta_1,\beta_2]^\top$.
The global minimum of this objective lies at the center of the innermost ellipse. However, the addition of a constraint restricts the feasible region for $\bbeta$, preventing the solution from reaching this unconstrained minimum. This restriction helps mitigate overfitting.

In contrast, the $\ell_1$-regularized model~\eqref{opt:ll} (p.~\pageref{opt:ll}) encourages sparsity in the solution; see Figure~\ref{fig:l1_ball_loss} for its constrained formulation.
Among them, the $\ell_1$-regularization term restricts the parameters within a diamond-shaped region, whereas the $\ell_2$-regularization term confines the parameters within a circular region. 
As illustrated in the figure, the optimal solution under $\ell_1$-regularization often occurs at the corners (vertices) of the diamond---typically where one coefficient, such as $\beta_1$, is exactly zero.
By contrast, the circular geometry of the $\ell_2$-constraint rarely leads to exact zeros, resulting in dense (non-sparse) solutions.

\begin{SCfigure}
\includegraphics[width=0.45\textwidth]{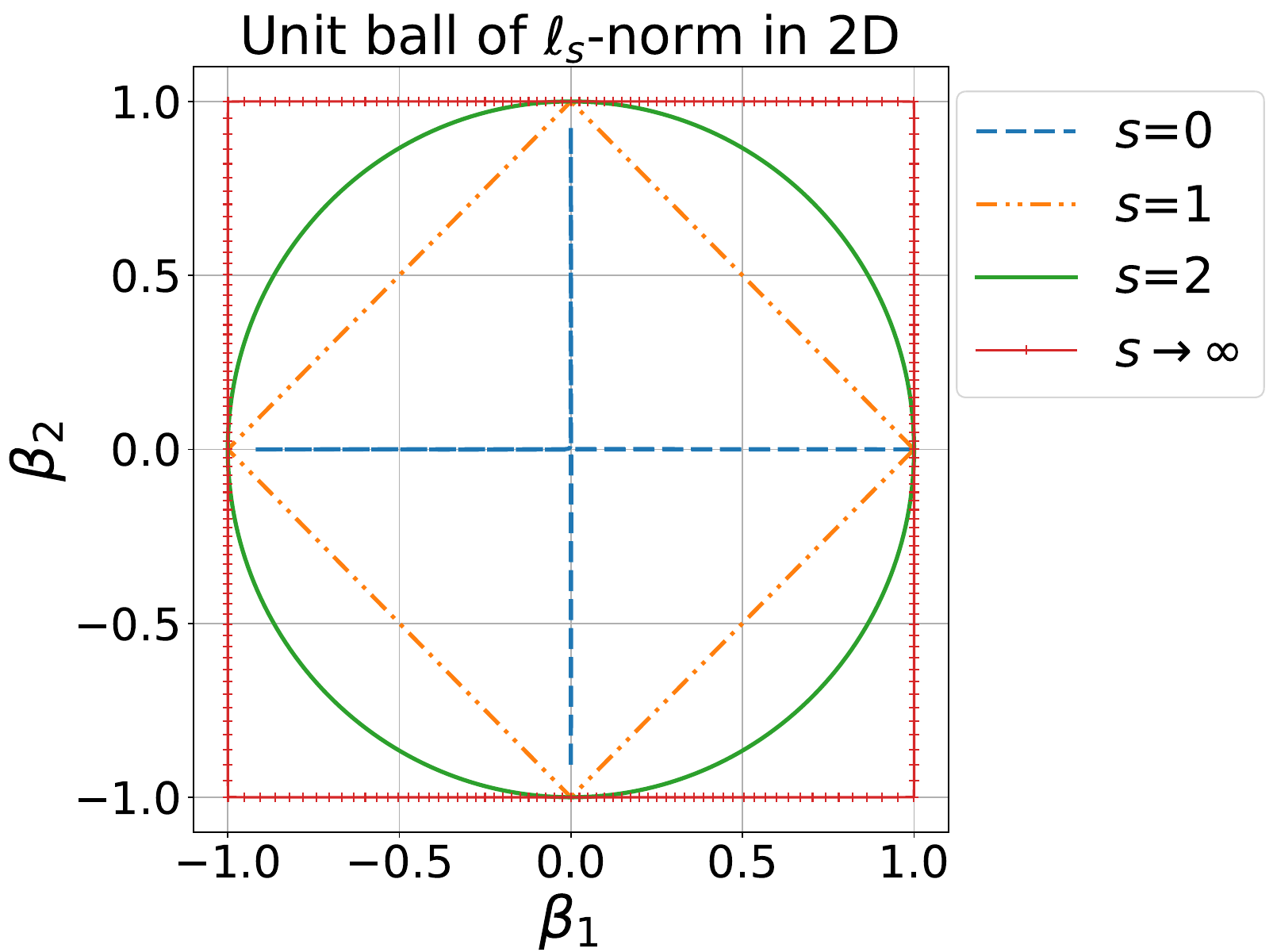}
\caption{Illustration of the intersection of $\ell_0$ and $\ell_\infty$-balls in 2D, and their convex hull forming the $\ell_1$-ball.}
\label{fig:s_norm_2d_BETA_intersection}
\end{SCfigure}

The $\ell_0$-unit ball  promotes sparsity, while it is not convex. 
And the $\ell_\infty$-unit ball includes all vectors where each entry has magnitude less than or equal to $ 1$.
The intersection of the $\ell_0$- and $\ell_\infty$-unit balls yields precisely the set of 1-sparse vectors with entries in $[-1,1]$.
This observation provides another justification for using the $\ell_1$-norm to induce sparsity:
the $\ell_1$-unit ball is the convex hull (Definition~\ref{definition:cvx_hull}) of all 1-sparse vectors with unit $\ell_\infty$-norm---that is, the convex hull of the intersection between the $\ell_0$- and $\ell_\infty$-unit balls.
This relationship is illustrated in two dimensions in Figure~\ref{fig:s_norm_2d_BETA_intersection}.

\begin{lemma}[$\ell_1$-unit ball]\label{lemma:l1_ball_intersec}
The $\ell_1$-unit ball is the convex hull of the intersection between the $\ell_0$-unit ball and the $\ell_\infty$-unit ball.
\end{lemma}
\begin{proof}[of Lemma~\ref{lemma:l1_ball_intersec}]
\textbf{$\sB_1 \subseteq \conv(\sB_0 \cap \sB_\infty)$.}
Let $\bbeta\in\sB_1\subset \real^p$ so that $\normone{\bbeta}\leq 1$. 
If we set $\theta_i \triangleq \abs{\beta_i}$ and $\theta_0 \triangleq 1 - \sum_{i=1}^p \theta_i$, we have $\sum_{i=0}^{p} \theta_i = 1$ by construction, $\theta_i = \abs{\beta_i} \geq 0$ and
\begin{align}
\theta_0 &= 1 - \sum_{i=1}^p \theta_i 
= 1 - \normone{\bbeta}
\geq 0.
\end{align}
We can express now $\bbeta$ as a convex combination of the standard basis vectors multiplied by the sign of the entries of $\bbeta$ as $\sign(\beta_1)\be_1$, $\sign(\beta_2)\be_2$, \ldots, $\sign(\beta_p)\be_p$ and the zero vector $\bzero$:
\begin{equation}
\bbeta = \sum_{i=1}^p \theta_i \sign(\beta_i)\be_i + \theta_0 \bzero, 
\end{equation}
where the vectors all belong to $\sB_0 \cap \sB_\infty$. Hence, $\sB_1 \subseteq \conv(\sB_0 \cap \sB_\infty)$.

\paragraph{$\conv(\sB_0 \cap \sB_\infty) \subseteq \sB_1$.}
Let $\bbeta\in\conv(\sB_0 \cap \sB_\infty) \subset\real^p$. 
By definition of the convex hull, we can express $\bbeta$ as
\begin{equation}
\bbeta = \sum_{i=1}^m \theta_i \balpha_i,
\end{equation}
where $m > 0$, $\balpha_1, \dots, \balpha_m \in \real^p$ have a single entry bounded by one, $\theta_i \geq 0$ for all $1 \leq i \leq m$ and $\sum_{i=1}^m \theta_i = 1$. 
By  the triangle inequality, we have
\begin{align}
\normone{\bbeta} 
&\leq \sum_{i=1}^m \theta_i \normone{\balpha_i}
= \sum_{i=1}^m \theta_i \norminf{\balpha_i}
\leq \sum_{i=1}^m \theta_i
\leq 1,
\end{align}
where the equality follows since each vector $\balpha_i$ has only one nonzero entry.
This immediately implies $\bbeta \in \sB_1$ and completes the proof.
\end{proof}

\section{Optimality Condition of $\ell_1$-Regularized Problems}\label{section:optcd_el1_reg}

We have previously demonstrated that the $\ell_1$-regularized least squares estimator promotes  sparsity in the solution: many components of the estimated coefficient vector $\bbeta$ are driven exactly to zero, leading to interpretable and compact models. This property makes $\ell_1$-regularization especially valuable in high-dimensional settings where feature selection is desired.

To analyze this behavior more formally, we now turn to the optimality conditions of the $\ell_1$-regularized least squares problem---commonly known as the Lagrangian LASSO (see \eqref{opt:ll}, p.~\pageref{opt:ll}). In its penalized (Lagrangian) form, the problem is stated as
\begin{equation}\label{equation:ll_opt}
\text{(L$_{L}$)}\quad 
\min_{\bbeta} F(\lambda, \bbeta) \triangleq  \left\{ \frac{1}{2} \normtwo{\by - \bX \bbeta}^2 + \lambda \normone{\bbeta} \right\}.
\end{equation}
The subgradient of the objective function \eqref{equation:ll_opt} with respect to $\bbeta$ is given by
\begin{equation}\label{equation:lass_grad}
\partial_{\bbeta} F(\lambda, \bbeta) = \frac{\partial F(\lambda, \bbeta)}{\partial \bbeta} 
= -\bX^\top (\by - \bX \bbeta) + \lambda \nabla_{\bbeta} \normone{\bbeta} 
\triangleq -\bc + \lambda \partial \normone{\bbeta},
\end{equation}
where $\bc \triangleq \bX^\top (\by - \bX \bbeta)$ is the \textit{residual correlation vector}, and the vector $\partial  \normone{\bbeta} = [\partial_{\beta_1} \normone{\bbeta},\partial_{\beta_2} \normone{\bbeta}, \ldots, \partial_{\beta_p} \normone{\bbeta}]^\top$ denotes the subdifferential of the $\ell_1$-norm $\normone{\bbeta}$ (Exercise~\ref{exercise:sub_norms}), with $i$-th entry
\begin{equation}\label{equation:norm1_grad_lass}
\partial \abs{\beta_i}  =
\begin{cases} 
\{+1\}, & \beta_i > 0; \\
\{-1\}, & \beta_i < 0; \\
[-1, +1], & \beta_i = 0,
\end{cases}
\qquad \forall\, i = 1,2, \ldots, p.
\end{equation}
From \eqref{equation:lass_grad}, it follows that a stationary point of the Lagrangian LASSO problem \eqref{opt:ll} must satisfy $\partial_{\bbeta} F(\lambda, \bbeta) = -\bc + \lambda \partial \normone{\bbeta} = \bzero$, i.e.,
\begin{equation}\label{equation:stat_lass}
\bc = \lambda \partial \normone{\bbeta}.
\end{equation}
Equivalently, for each $i=1,2,\ldots,p$, we have  $c_i = \lambda \partial \abs{\beta_i}$.
Because the Lagrangian LASSO problem is convex, this stationary-point condition is both necessary and sufficient for optimality (Theorem~\ref{theorem:fetmat_opt}).

The optimality condition \eqref{equation:stat_lass} can be expressed in terms of the support set $\sS \triangleq \{i \mid \beta_i\neq 0\}$ and its complement $\comple{\sS} = \{1, 2, \ldots, p\} \setminus \sS$ as
\begin{equation}\label{equation:lass_stat_q2}
\bc[\sS] = \lambda \cdot \sgn(\bbeta[\sS]) 
\qquad \text{and} \qquad 
\abs{\bc[\comple{\sS}]} \leq \lambda\bone,
\end{equation}
This reveals that, for indices in the support $\sS$, the residual correlations have magnitude exactly equal to $\lambda$ and share the same sign as the corresponding coefficients in $\bbeta$.
Equivalently, condition \eqref{equation:lass_stat_q2} can be written as
\begin{equation}\label{equation:lass_stat_q3}
\abs{c_i} = \lambda, \quad \forall\, i \in \sS
\qquad \text{and} \qquad 
\abs{c_i} \leq \lambda, \quad \forall\, i \in \comple{\sS}.
\end{equation}
In other words, the absolute value of the residual correlation equals $\lambda$ for active (nonzero) coefficients and does not exceed $\lambda$ for inactive (zero) coefficients. Consequently, the maximum absolute residual correlation satisfies $\norminf{\bc} = \max\{c_i\} = \lambda$.

\begin{remark}[KKT conditions of constrained LASSO]
The conditions in \eqref{equation:lass_stat_q3} are similar to the KKT conditions of the constrained LASSO problem~\eqref{opt:lc} (p.~\pageref{opt:lc}).
This equivalence arises from the shared optimality structure between the penalized and constrained formulations of LASSO.
At the optimal solution $\widehatbbeta$ for a given $\ell_1$-norm constraint $\Sigma$, the following KKT conditions (Theorem~\ref{theorem:opt_cond_sd}) must hold:
\begin{enumerate}[(i)]
\item \textit{Subgradient condition.} The gradient of the least squares loss ($- \bX^\top(\by - \bX\widehatbbeta)$) must be in the subdifferential of $\nu\normonebig{\widehatbbeta}$, where $\nu$ is the Lagrange multiplier associated with the constraint. 
This translates to:
\begin{itemize}
\item For any $i$ where $\widehatbeta_i \neq 0$: $\absbig{\bx_i^\top(\by - \bX\widehatbbeta)} = \nu$ (the absolute correlation between predictor $i$ and the residual is exactly $\nu$).
\item For any $i$ where $\widehatbeta_i = 0$: $\absbig{\bx_i^\top(\by - \bX\widehatbbeta)} \leq \nu$ (the absolute correlation  does not exceed $\nu$).
\end{itemize}
\item \textit{Complementary slackness.} $\nu(\normonebig{\widehatbbeta} - \Sigma) = 0$. If $\normonebig{\widehatbbeta} < \Sigma$, then $\nu = 0$ (no regularization). If $\nu > 0$, then $\normonebig{\widehatbbeta} = \Sigma$ (the constraint is active).
\end{enumerate}
\end{remark}

\subsection{LASSO Dual}
For the Lagrangian LASSO problem~\eqref{opt:ll},
introducing the residual vector $\br = \by - \bX\bbeta$, we can equivalently reformulate the primal problem as
\begin{equation}\label{equation:laglasso_equiv}
\min_{\bbeta \in \real^p, \br\in\real^n} \frac{1}{2}\normtwo{\br}^2 + \lambda\normone{\bbeta} 
\quad\text{s.t.}\quad \br = \by - \bX\bbeta.
\end{equation}
Letting $\bnu \in \real^n$  denote the vector of Lagrange multipliers associated with the equality constraint. The Lagrangian function for this problem is
\begin{equation}\label{equation:laglasso_lag}
L(\bbeta, \br, \bnu) = \frac{1}{2}\normtwo{\br}^2 + \lambda\normone{\bbeta} 
- \bnu^\top(\br - \by + \bX\bbeta).
\end{equation}
The Lagrangian dual function is obtained by minimizing $L(\bbeta, \br, \bnu)$ over $\bbeta$ and $\br$. 
In order to obtain the Lagrangian dual function $\min_{\bbeta, \br} L(\bbeta, \br, \bnu)$,
isolating those terms involving $\bbeta$, and noting that $\abs{\bnu^\top\bX\bbeta} \leq \normone{\bbeta}\norminf{\bX^\top\bnu}$ by \holders inequality (Theorem~\ref{theorem:holder-inequality}),  we find
\begin{equation}\label{equation:laglasso_dual1}
\min_{\bbeta \in \real^p} \lambda\normone{\bbeta} - \bnu^\top \bX\bbeta =
\begin{cases}
0, & \text{if } \norminf{\bX^\top\bnu} \leq \lambda; \\
-\infty, & \text{otherwise}.
\end{cases}
\end{equation}
Next, we isolate terms involving $\br$ and find
\begin{equation}\label{equation:laglasso_dual2}
\min_{\br\in\real^n} \frac{1}{2}\normtwo{\br}^2 - \bnu^\top\br = -\frac{1}{2}\bnu^\top\bnu,
\quad \text{with $\br = \bnu$}.
\end{equation}
Substituting relations~\eqref{equation:laglasso_dual1} and~\eqref{equation:laglasso_dual2} into the Lagrangian representation~\eqref{equation:laglasso_lag}, we obtain
\begin{equation}\label{equation:laglasso_dual_all}
\begin{aligned}
\text{LASSO Dual:}\qquad  
\phi(\lambda)
&\triangleq 
\max_{\bnu} \frac{1}{2}\left(\normtwo{\by}^2 - \normtwo{\by - \bnu}^2\right) 
 \triangleq \max_{\bnu} D(\bnu)\\
&\text{s.t.} \quad\norminf{\bX^\top\bnu} \leq \lambda,
\end{aligned}
\end{equation}
where $D(\bnu) \triangleq \frac{1}{2}\left(\normtwo{\by}^2 - \normtwo{\by - \bnu}^2\right) $ is the dual objective function, 
and $\phi(\lambda)$ denotes the optimal dual value corresponding to the penalty parameter $\lambda$.

\begin{SCfigure}
\centering
\includegraphics[width=0.45\textwidth]{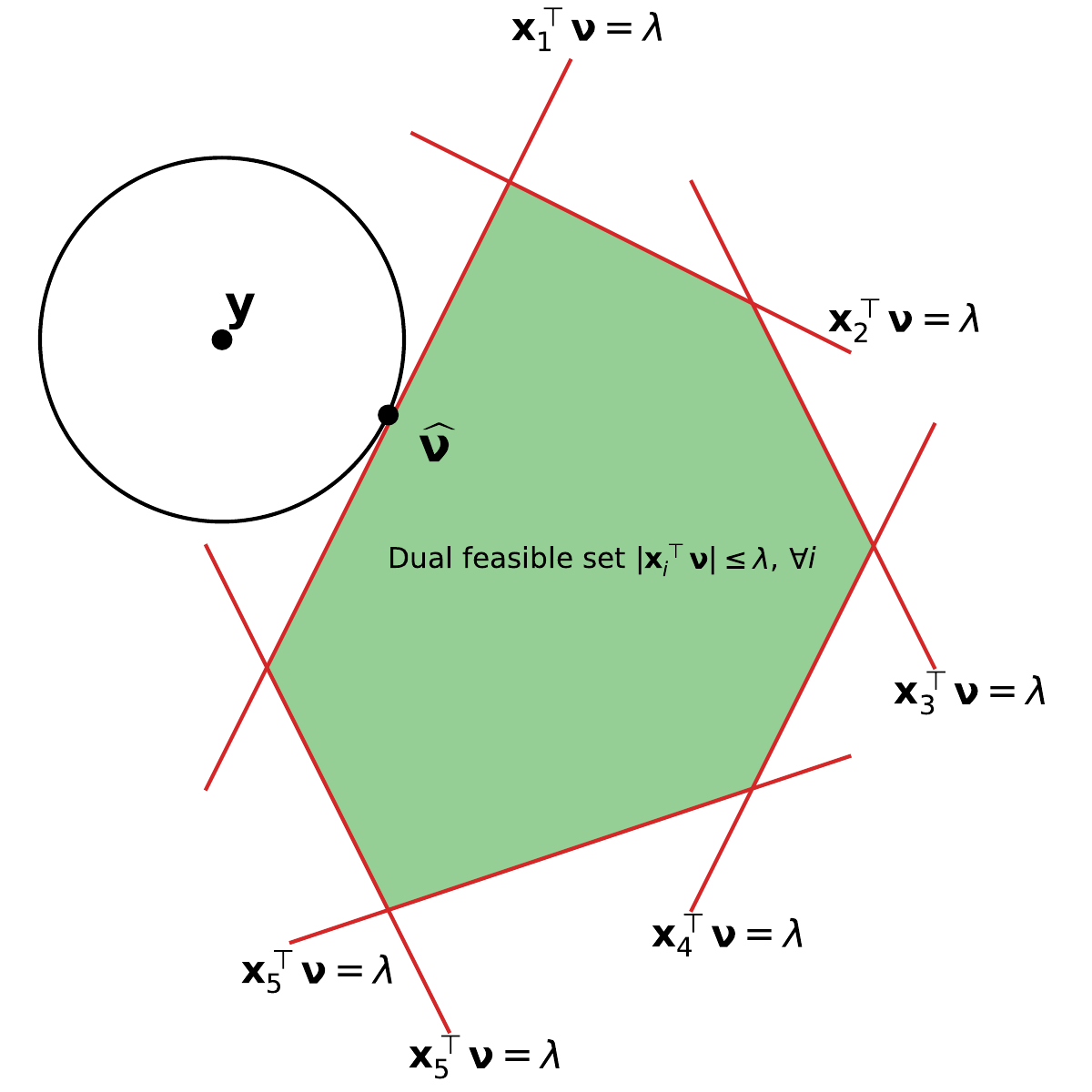}
\caption{The Lagrange dual of the LASSO, with $D(\bnu) = \frac{1}{2}(\normtwo{\by}^2 - \normtwo{\by - \bnu}^2)$. The  shaded region is the feasible set $\sS_{\lambda}$. The unconstrained dual solution is $\bnu = \by$, the null residual. 
The dual solution $\widehatbnu = \project_{{\sS_{\lambda}}}(\by)$, the projection of $\by$ onto the convex set $\sS_{\lambda}$.}
\label{fig:laglasso_dual_fea}
\end{SCfigure}

Geometrically, this dual formulation amounts to projecting the response vector $\by$ onto the feasible set $ \sS_{\lambda} \triangleq \{ \bnu \in \real^n \mid \norminf{\bX^\top\bnu} \leq \lambda \} $. 
The set $ \sS_{\lambda} $ is the intersection of the $ 2p $ half-spaces defined by $ \abs{\bx_i^\top\bnu} \leq \lambda $, each one is a convex polytope in $ \real^n $. In the language of projection and proximal operators, the solution is given by the projection or proximal map $ \widehatbnu = \project_{{\sS_{\lambda}}}(\by) \equiv\prox_{\indicatorG_{\sS_{\lambda}}}(\by) $.
In other words, $\widehatbnu$ is the projection onto the polytope $\sS_{\lambda}$. 
Figure~\ref{fig:laglasso_dual_fea} illustrates this geometric interpretation.

\paragrapharrow{Dual feasibility condition for sparsity.}
For the primal LASSO problem \eqref{opt:ll}, strong duality holds (Theorem~\ref{theorem:fenchel_dual_ori}).
The optimal primal $\widehatbbeta$ and dual $\widehatbnu$ solutions then satisfy
\begin{equation}\label{equation:laglasso_prim_dual_rel}
\by = \bX\widehatbbeta +\widehatbnu.
\end{equation}
The dual problem plays a crucial role in identifying sparsity in the primal solution: by examining the dual variable  $\widehatbnu$, we can determine which coefficients in  $\widehatbbeta$ must be exactly zero. Specifically, the optimality condition
\begin{equation}\label{equation:lasso_stationarity_ori}
\abs{\bx_i^\top \widehatbnu} 
< \lambda 
\quad\implies\quad 
\widehatbeta_i(\lambda) = 0, \quad i = 1,2,\dots,p,
\end{equation}
follows directly from \eqref{equation:lass_stat_q3} and \eqref{equation:laglasso_dual2}.
Note that the above condition is sufficient for that $\widehatbeta_i(\lambda)$ must be an inactive predictor, but it is not necessary.
Indeed, it is possible for  $\widehatbeta_i(\lambda) = 0$ even when $\abs{\bx_i^\top \widehatbnu} 
= \lambda $. In such cases, the predictor lies on the boundary of the active set and is considered marginally inactive.

\paragrapharrow{Largest $\lambda$ for a nontrivial solution.}

This form of the LASSO dual also reveals the largest value of the penalty parameter, denoted $\lambda=\lambda_{\max}$, for which the primal LASSO problem admits a nontrivial solution $\widehatbbeta_{\max} \neq \bzero$.
If $\lambda=\lambda_{\max} $ is chosen large enough, the optimal primal solution becomes identically zero: $\widehatbbeta= \widehatbbeta_{\max}=\bzero $.
To determine  such a $ \lambda_{\max} $, we substitute $ \widehatbbeta_{\max} = \bzero $ in the primal LASSO problem~\eqref{opt:ll} to obtain $ \phi(\lambda_{\max}) = \frac{1}{2}\normtwo{\by}^2 $ since $\by=\bnu$ in \eqref{equation:laglasso_dual2}. 
By strong duality, the dual problem under $\lambda=\lambda_{\max}$ achieves a value of $ \phi(\lambda_{\max}) = \frac{1}{2}\normtwo{\by}^2 $ at the unique solution 
\begin{equation}\label{equation:safe_numax}
\widehatbnu_{\max} =  \by. 
\end{equation}
The point $ \widehatbnu_{\max} $ is a dual feasible point and satisfies the constraints $ \abs{\bx_i^\top\by} \leq \lambda_{\max}  $, $ i = 1,2,\dots,p $. 
Therefore, this form of the LASSO dual shows the largest $\lambda=\lambda_{\max}$ such that the solution $\widehatbbeta_{\max}$ is nontrivial is
\begin{equation}\label{equation:safe_lambdamax}
\lambda_{\max} = \max_{1 \leq i \leq p} \abs{\bx_i^\top\by} = \norminf{\bX^\top\by}.
\end{equation}

\subsection{Dual Polytope Projections and Screening Rules}\label{section:opcond_dpp}

The quantity in \eqref{equation:safe_lambdamax} reveals an important intuition: covariates that have small inner products with the response vector $\by$ are less likely to receive nonzero coefficients in the LASSO solution compared to those with larger inner products

This observation suggests a practical strategy: we may be able to eliminate irrelevant predictors before solving the full optimization problem, thereby reducing computational cost. For instance, in genomic applications---where datasets often involve millions of single-nucleotide polymorphisms (SNPs)---we typically expect only a handful of variables to appear in the final model.

In this subsection, we discuss \textit{screening rules} that exploit this intuition. These rules can dramatically accelerate LASSO computations by efficiently identifying and removing inactive predictors---i.e., features that are guaranteed to have zero coefficients in the optimal solution---without compromising numerical accuracy. The resulting reduced problem yields the exact same solution as the original, but at a fraction of the computational cost.

We begin with a geometric screening rule based on projections onto the dual feasible set, known as the \textit{dual polytope projections (DPPs)} rule \citep{ghaoui2010safe, wang2013lasso}.
\begin{theoremHigh}[Dual polytope projections (DPP)\index{Dual polytope projection}]\label{theorem:ddp}
Suppose we are given a solution $\widehatbnu(\lambda')$ to the LASSO dual~\eqref{equation:laglasso_dual_all}  and a primal solution $\widehatbbeta(\lambda')$ at some penalty parameter $0<\lambda' \leq \lambda_{\max}$. 
Let $\lambda$ be a nonnegative value different from $\lambda'$. 
If, for a given index $i$, the following inequality holds:
\begin{equation}\label{equation:lasso_dpp_thres}
\abs{\bx_i^\top \widehatbnu(\lambda')} 
< \lambda' - \normtwo{\bx_i} \normtwo{\by} \abs{1 - \frac{\lambda'}{\lambda}},
\end{equation}
then the $i$-th coefficient in the LASSO solution at $\lambda$ is zero:  
$\widehatbeta_i(\lambda) = 0$, where $\bx_i$ denotes the $i$-th column of $\bX$.
\end{theoremHigh}
\begin{proof}[of Theorem~\ref{theorem:ddp}]
Introduce the change of variable $\bmu\triangleq \bnu/\lambda$. 
Then the LASSO dual problem \eqref{equation:laglasso_dual_all} becomes
\begin{equation}\label{equation:laglasso_dual_all_chg}
\text{LASSO Dual:}\quad  
\max_{\bmu} \frac{1}{2}\left(\normtwo{\by}^2 - \lambda^2\normtwo{\by/\lambda -  \bmu}^2\right) 
\quad \text{s.t.} \quad\norminf{\bX^\top\bmu} \leq 1.
\end{equation}
Thus, the feasible set in the scaled dual space is the fixed convex polytope $ \widetilde{\sS} \triangleq \{ \bz \in \real^n \mid \norminf{\bX^\top\bz} \leq 1 \}$, which is independent of $\lambda$.
The optimal scaled dual variable $\widehatbmu(\lambda)$   is the projection of $\by/\lambda$ into the feasible set $\widetilde{\sS}$.

By \eqref{equation:lasso_stationarity_ori}, the stationarity conditions for the LASSO problem, we have 
\begin{equation}\label{equation:lasso_stationarity_dpp}
\abs{\bx_i^\top \widehatbmu(\lambda)} 
< 1 
\quad\implies\quad \widehatbeta_i(\lambda) = 0, \quad i = 1,2,\dots,n.
\end{equation}
Since the polytopes for different penalty parameters $\lambda$ are the same, 
by the nonexpansiveness of projection operators (Theorem~\ref{theorem:proj_nonexpan}), we have
\begin{equation}\label{equation:lasso_stationarity_dpp2}
\normtwo{\widehatbmu(\lambda) - \widehatbmu(\lambda')}
=
\normtwo{\project_{\widetilde{\sS}}(\frac{\by}{\lambda}) - \project_{\widetilde{\sS}}(\frac{\by}{\lambda'})}
\leq \normtwo{ \frac{\by}{\lambda} - \frac{\by}{\lambda'}} 
= \normtwo{\by} \abs{\frac{1}{\lambda} - \frac{1}{\lambda'}}.
\end{equation}
Using triangle inequality and Cauchy--Schwarz inequality, the above inequality and the assumption~\eqref{equation:lasso_dpp_thres} yield
\begin{align*}
\abs{\bx_i^\top \widehatbmu(\lambda)}
&\leq \abs{\bx_i^\top \widehatbmu(\lambda) - \bx_i^\top \widehatbmu(\lambda')} + \abs{\bx_i^\top \widehatbmu(\lambda')} \\
&< \normtwo{\bx_i} \normtwo{\widehatbmu(\lambda) - \widehatbmu(\lambda')} + 1 - \normtwo{\bx_i} \normtwo{\by} \abs{ \frac{1}{\lambda'} - \frac{1}{\lambda} } \nonumber \\
&\leq \normtwo{\bx_i} \normtwo{\by} \abs{\frac{1}{\lambda'} - \frac{1}{\lambda}} + 1 - \normtwo{\bx_i} \normtwo{\by} \abs{\frac{1}{\lambda'} - \frac{1}{\lambda}} 
= 1.
\end{align*}
Using \eqref{equation:lasso_stationarity_dpp}  completes the proof.
\end{proof}

\begin{subequations}
The above theorem shows that if
\begin{equation}\label{equation:lasso_dpp_thres2}
	\abs{\bx_i^\top \big(\by-\bX\widehatbbeta(\lambda')\big)} 
	< \lambda' - \normtwo{\bx_i} \normtwo{\by} \abs{1 - \frac{\lambda'}{\lambda}},
\end{equation}
then  variable $i$ does not belong to the active set at penalty level $\lambda$ (i.e., $\widehatbeta_i(\lambda) = 0$).
In particular, at $\lambda'=\lambda_{\max} \equiv \norminfbig{\bX^\top\by}$ by \eqref{equation:safe_lambdamax}, the above inequality becomes
\begin{equation}\label{equation:lasso_dpp_thres3}
	\abs{\bx_i^\top \by} 
	< \lambda_{\max} - \normtwo{\bx_i} \normtwo{\by} \left(\frac{\lambda_{\max}}{\lambda} - 1\right).
\end{equation}
The \textit{dual polytope projection (DPP)} rule then suggests to discard the $i$-th predictor whenever  \eqref{equation:lasso_dpp_thres3} holds \citep{wang2013lasso}.
\end{subequations}

\paragrapharrow{Variant of DPP rule.}
On the other hand, suppose $\frac{\by}{\lambda'} \notin \widetilde{\sS} \triangleq \{ \bz \in \real^n \mid \norminf{\bX^\top\bz} \leq 1 \}$, i.e., $\lambda' \in (0, \lambda_{\max})$. 
In this case, the scaled dual solution satisfies $\frac{\by}{\lambda'} \neq \project_{\widetilde{\sS}}(\frac{\by}{\lambda'}) \equiv \widehatbmu(\lambda')$ (recall that $\widehatbmu(\lambda') = \widehatbnu(\lambda')/\lambda'$ is the scaled dual solution under the penalth parameter $\lambda=\lambda'$). 
Let $\bmu(\eta) \triangleq  \widehatbmu(\lambda') + \eta \left( \frac{\by}{\lambda'} - \widehatbmu(\lambda') \right)$ for $\eta \geq 0$, i.e., $\bmu(\eta)$ is a point lying on the ray starting from $\widehatbmu(\lambda')$ and pointing to the same direction as $\frac{\by}{\lambda'} - \widehatbmu(\lambda')$. We can observe that $\project_{\widetilde{\sS}}(\bmu(\eta)) = \widehatbmu(\lambda')$, i.e., the projection of $\bmu(\eta)$ onto the set $\widetilde{\sS}$ is $\widehatbmu(\lambda')$ as well. By the nonexpansiveness of the projection operator (Theorem~\ref{theorem:proj_nonexpan}), we have 
\begin{align}
\normtwo{\widehatbmu(\lambda) - \widehatbmu(\lambda')}
&= \normtwo{\project_{\widetilde{\sS}}\big(\frac{\by}{\lambda}\big) - \project_{\widetilde{\sS}}\big(\bmu(\eta)\big)} \nonumber \\
&\leq \normtwo{\frac{\by}{\lambda} - \bmu(\eta)} 
= \normtwo{ \eta \big(\frac{\by}{\lambda'} - \widehatbmu(\lambda')\big) - \big(\frac{\by}{\lambda} - \widehatbmu(\lambda')\big)} . 
\label{equation:vdpp_eq1}
\end{align}
Clearly, when $\eta=1$, the inequality in \eqref{equation:vdpp_eq1} reduces to the one in \eqref{equation:lasso_stationarity_dpp2}. 
Because the inequality in \eqref{equation:vdpp_eq1} holds for all $\eta \geq 0$, we may get a tighter bound by
\begin{equation}\label{equation:vdpp_eq2}
\normtwo{ \widehatbmu(\lambda) - \widehatbmu(\lambda') } 
\leq \min_{\eta \geq 0} \normtwo{\eta \bu_1 - \bu_2},
\end{equation}
where $\bu_1 \triangleq \frac{\by}{\lambda'} - \widehatbmu(\lambda')$ and 
$\bu_2 \triangleq \frac{\by}{\lambda} - \widehatbmu(\lambda')$. 
When $\lambda' = \lambda_{\max}$, we can set $\bu_1 = \sign(\bx_*^\top \by) \bx_*$ where $\bx_* \triangleq \arg\max_{\bx_i} \abs{\bx_i^\top \by}$. 
The minimization problem on the right-hand side of the inequality \eqref{equation:vdpp_eq2} has a closed-form solution:
\begin{equation}
\min_{\eta \geq 0} \normtwo{\eta \bu_1 - \bu_2} 
= \varphi(\lambda', \lambda) \triangleq 
\begin{cases}
\normtwo{\bu_2}, & \text{if } \innerproduct{\bu_1, \bu_2} < 0; \\
\normtwo{\frac{\innerproduct{\bu_1, \bu_2}}{\normtwo{\bu_1}^2} \bu_1  - \bu_2  }, & \text{otherwise}.
\end{cases}
\end{equation}
This refined bound leads to an enhanced DPP screening rule, as stated below.
\begin{corollary}[Enhanced DPP rule]\label{corollary:enhanced_dpop}

Suppose we are given a solution $\widehatbnu(\lambda')$ to the LASSO dual~\eqref{equation:laglasso_dual_all}  and a primal solution $\widehatbbeta(\lambda')$ at some penalty parameter $0<\lambda' \leq \lambda_{\max}$. 
Let $\lambda$ be a nonnegative value different from $\lambda'$. 
If, for a given index $i$, the following inequality holds:
\begin{equation}\label{equation:lasso_edpp_thres}
\abs{\bx_i^\top \widehatbnu(\lambda')} 
< \lambda' - \lambda'\normtwo{\bx_i} \varphi(\lambda', \lambda),
\end{equation}
then $\widehatbeta_i(\lambda) = 0$, where $\bx_i$ is the $i$-th column of $\bX$.
\end{corollary}

\begin{subequations}
The above corollary shows that if
\begin{equation}\label{equation:lasso_edpp_thres2}
\abs{\bx_i^\top \big(\by-\bX\widehatbbeta(\lambda')\big)} 
< \lambda' - \lambda'\normtwo{\bx_i} \varphi(\lambda', \lambda),
\end{equation}
then  variable $i$ is not part of the activate set at penalty level $\lambda$.
At $\lambda'=\lambda_{\max} \equiv \norminfbig{\bX^\top\by}$ by \eqref{equation:safe_lambdamax}, this condition becomes
\begin{equation}\label{equation:lasso_edpp_thres3}
\abs{\bx_i^\top \by} 
< \lambda_{\max} -  \lambda_{\max}\normtwo{\bx_i} \varphi(\lambda_{\max}, \lambda).
\end{equation}
Again, the enhanced DPP rule then suggests to discard the $i$-th variable if \eqref{equation:lasso_edpp_thres3} is satisfied.
\end{subequations}

\paragrapharrow{The screening rule.} 
Equation~\eqref{equation:lasso_dpp_thres2} or \eqref{equation:lasso_edpp_thres2} provides a practical screening rule for efficiently solving the Lagrangian LASSO problem~\eqref{opt:ll}.
Suppose we know the optimal dual solution $\widehatbnu(\lambda')$ at some penalty parameter $\lambda'$ (e.g., the well-known case $\widehatbnu(\lambda_{\text{max}}) = \by$). 
Using this information, we can bound the possible location of the dual solution $\widehatbnu(\lambda)$ at a nearby parameter $\lambda$. 
By bounding this possible location of $\widehatbnu(\lambda)$, we can check if $\abs{\bx_i^\top\widehatbnu(\lambda)}$ is \textit{guaranteed} to be less than $\lambda$. If so, the corresponding predictor $\bx_i$ is inactive at $\lambda$ (i.e., $\widehatbeta_i(\lambda)=0$) and can be safely removed from the optimization problem without affecting the final solution.

\subsection*{Sequential DDP}
Motivated by the dual polytope projection theorem (Theorem~\ref{theorem:ddp}) and the  screening principle above, we can construct  a sequential   DPP strategy. 
Specifically, consider a decreasing sequence of penalty parameters: $\lambda_0=\lambda_{\max}> \lambda_1 > \lambda_2 > \ldots > \lambda_m$, 
We begin by applying the DPP screening rule at $\lambda_0=\lambda_{\max}$ to eliminate predictors that are guaranteed to be inactive.
Since  $\widehatbnu(\lambda_{\max})=\by$, this initial screening is particularly simple and effective.

After discarding inactive variables, we solve the reduced LASSO problem at $\lambda_1$ using any algorithm from Section~\ref{section:algo_laglasso}, obtaining the exact primal solution  $\widehatbbeta(\lambda_1)$. From the primal-dual relationship \eqref{equation:laglasso_prim_dual_rel}, we then recover the dual solution:
$\widehatbnu
(\lambda_1)=\by-\bX 
\widehatbbeta(\lambda_1)$.
This newly computed dual solution enables us to apply the DPP rule again---now with tighter bounds---to screen out inactive predictors for the next parameter $\lambda_2$. Repeating this process iteratively yields a sequential screening framework that progressively reduces the problem size at each step, significantly accelerating computation while preserving exactness.

\begin{corollary}[Sequential DPP]\label{corollary:seq_dpp}
For the Lagrangian LASSO problem \eqref{opt:ll}, suppose we are given a sequence of parameter values $\lambda_{\max} = \lambda_0 > \lambda_1 > \ldots > \lambda_m$. Then for any integer $0 \leq t < m$, we have $\widehatbeta_i(\lambda_{t+1}) = 0$ if $\widehatbbeta(\lambda_t)$ is known, and the following holds:
$$
\abs{\bx_i^\top  \widehatbnu(\lambda_t)}
\equiv \abs{\bx_i^\top \big(\by-\bX\widehatbbeta(\lambda_t)\big)}
< \lambda_t - \normtwo{\bx_i} \normtwo{\by} \left( \frac{\lambda_t}{\lambda_{t+1}} - 1 \right), 
\quad i=1,2,\ldots,p.
$$
or 
$$
\abs{\bx_i^\top  \widehatbnu(\lambda_t)}
\equiv \abs{\bx_i^\top \big(\by-\bX\widehatbbeta(\lambda_t)\big)}
< \lambda_t - \lambda_t\normtwo{\bx_i} \varphi(\lambda_t, \lambda_{t+1}), 
\quad i=1,2,\ldots,p.
$$
\end{corollary}

\subsection{Safe Feature Elimination (SAFE)}\label{section:safe}

The \textit{SAFE (Safe FEature Elimination)} rule is another powerful technique used in sparse optimization problems, especially in LASSO-type regression, to efficiently identify and eliminate features (columns of the design matrix) that are guaranteed to have zero coefficient in the optimal solution without solving the full problem.
This significantly speeds up computation, particularly when dealing with high-dimensional data where many features are irrelevant.

The intuition behind SAFE stems from the stationarity condition in \eqref{equation:lasso_stationarity_ori}, which we recall here:
\begin{equation}\label{equation:lasso_stationarity_safe}
\abs{\bx_i^\top \widehatbnu} 
< \lambda 
\quad\implies\quad 
\widehatbeta_i(\lambda) = 0, \quad i = 1,2,\dots,p,
\end{equation}
The SAFE rule constructs  a constraint set $\Theta$ known to contain the dual optimal solution $\widehatbnu$. 
It then replaces $\widehatbnu$ in the inequality test with this entire set $\Theta$. Specifically, if
\begin{equation}
\abs{\bx_i^\top\bnu} < \lambda, \quad 
\text{for all } \bnu \in \Theta,
\end{equation}
then the condition in \eqref{equation:lasso_stationarity_safe} holds uniformly over $\Theta$. Consequently, the $i$-th component of the primal optimal solution must be zero:  $\widehatbeta_i=0$.

\subsection*{Constructing the Dual Constraint Set}

We begin by constructing a set $\Theta$ that is guaranteed to contain the dual optimal solution $\widehatbnu$ of the LASSO dual problem~\eqref{equation:laglasso_dual_all}:
\begin{equation}\label{equation:laglasso_dual_all_safe}
\begin{aligned}
\text{LASSO Dual:}\qquad  
\phi(\lambda)
&\triangleq 
\max_{\bnu}  \left\{D(\bnu)\triangleq \frac{1}{2}\left(\normtwo{\by}^2 - \normtwo{\by - \bnu}^2\right) 
=\by^\top\bnu - \frac{1}{2}\normtwo{\bnu}^2 \right\}
\\
&\text{s.t.} \quad\norminf{\bX^\top\bnu} \leq \lambda,
\end{aligned}
\end{equation} 
under the penalty parameter $\lambda$. 
To this end, we define $\Theta$ as the intersection of two sets, $\Theta_1$ and $\Theta_2$, each derived from a different optimality condition.

\paragrapharrow{Constructing $\Theta_1$.}
The first set, $\Theta_1$, is based on the fact that the dual objective $D(\bnu)$ achieves its maximum at $\widehatbnu$. That is, for any dual-feasible point $\bnu$, we have $D(\widehatbnu) \geq D(\bnu)$.  
Let $\bnu_\eta$ be a specific dual-feasible point for problem~\eqref{equation:laglasso_dual_all_safe}, and define $\gamma \triangleq D(\bnu_\eta)$. 
Then we must have $D(\widehatbnu) \geq \gamma$, so the set $\Theta_1 \triangleq \{\bnu \mid D(\bnu) \geq \gamma\}$ contains the dual optimum $\widehatbnu$, i.e., $\widehatbnu \in \Theta_1$.

One practical way to obtain such a feasible point $\bnu_\eta$ is via \textit{dual scaling}. 
Specifically, let $\widehatbnu_{0}$ be a known dual-feasible point corresponding to a larger penalty parameter $\lambda_{0}>\lambda$ (e.g., as in \eqref{equation:safe_numax}).
We then define $\bnu_\eta \triangleq \eta \widehatbnu_{0}$ where the scalar $\eta \in \real$ is chosen so that $\bnu_\eta$ remains feasible for the dual problem with penalty $\lambda$; that is, $\norminf{\bX^\top \bnu_\eta} \leq \lambda$ or $\abs{\eta} \leq \lambda / \lambda_{0}$. 
We now choose $\eta$ to maximize the dual objective over this interval:
$$
\gamma = \max_\eta \left\{ D(\eta\widehatbnu_{0}) \mid \abs{\eta} \leq \frac{\lambda}{\lambda_{0}} \right\} 
= \max_\eta \left\{ -\frac{1}{2} a\eta^2   + b \eta\mid \abs{\eta} \leq \frac{\lambda}{\lambda_{0}} \right\},
$$
where $a \triangleq \normtwo{\widehatbnu_{0} }^2> 0$
and $b \triangleq \abs{\by^\top \widehatbnu_{0}}$. 
We obtain by $\eta=\lambda/\lambda_{0}$ that
$
\gamma = -\frac{a}{2} \left(\frac{\lambda}{\lambda_{0}}\right)^2 +b \frac{\lambda}{\lambda_{0}}$, 
such that 
\begin{equation}
	\widehatbnu \in \Theta_1 
	\triangleq 
	\{\bnu \mid D(\bnu) \geq \gamma\}.
\end{equation}

\paragrapharrow{Constructing $\Theta_2$.}
The second set, $\Theta_2$, arises from a first-order optimality condition for the dual  problem~\eqref{equation:laglasso_dual_all_safe}.
Recall that for a convex maximization problem over a convex feasible set, a point $\widehatbnu_{0}$ is optimal if and only if $\bg^\top (\bnu - \widehatbnu_{0}) \geq 0$ for every dual point $\bnu$ that is feasible for the dual problem~\eqref{equation:laglasso_dual_all_safe} under the penalty parameter $\lambda_{0}$, 
where $\bg \triangleq \nabla D(\widehatbnu_{0}) =  \by - \widehatbnu_{0}$ 
(see Theorem~\ref{theorem:stat_point_uncons_convset} for the general optimality condition in convex-constrained  optimization). 
Now consider a penalty parameter  $\lambda \leq \lambda_{0}$. 
Any dual point $\bnu$ feasible for the dual problem~\eqref{equation:laglasso_dual_all_safe} under the penalty parameter $\lambda$ satisfies $\norminf{\bX^\top\bnu}\leq \lambda \leq \lambda_0$, and is therefore also dual feasible for the dual problem~\eqref{equation:laglasso_dual_all_safe} with parameter $\lambda_{0}$.  
Since $\widehatbnu$  (the dual optimum for penalty $\lambda$) is feasible under $\lambda$, it is also feasible under $\lambda_{0}$.
Therefore, we define
\begin{equation}
\widehatbnu \in \Theta_2 
\triangleq 
\left\{ \bnu \mid \bg^\top (\bnu - \widehatbnu_{0}) \geq 0 \right\}.
\end{equation}

Combining both conditions, the dual optimal solution must lie in the intersection of these two sets: 
$$
\Theta  = \Theta_1\cap \Theta_2.
$$
This set $\Theta$ provides a tight region that provably contains  $\widehatbnu$, enabling safe feature elimination via the rule described earlier.

\subsection*{Safe Feature Elimination}
The task of identifying whether the $i$-th component of the primal solution $\widehatbbeta$ is zero---and thus safely removing the $i$-th feature  (i.e., the $i$-th column $\bx_i$ of the design matrix $\bX$) from the Lagrangian LASSO problem~\eqref{opt:ll}---reduces to verifying the condition:
\begin{equation}
	\lambda > \abs{\bnu^\top \bx_i} = \max(\bnu^\top \bx_i, -\bnu^\top \bx_i) 
	\quad\text{s.t.}\quad 
	\bnu \in \Theta.
\end{equation}
Equivalently, this can be expressed as:
$$
\lambda > \max(S(\gamma, \bx_i), S(\gamma, -\bx_i)),
$$
where $S(\gamma, \bx_i)$ denotes the optimal value of the following convex optimization problem with constraints $\bnu \in \Theta_1$ and $\bnu \in \Theta_2$:
\begin{equation}\label{equation:safe_sgamm_bxi}
S(\gamma, \bx_i) \triangleq \max_{\bnu} 
\left\{\bx_i^\top \bnu 
\quad\text{s.t.}\quad 
D(\bnu) \geq \gamma,\ \bg^\top(\bnu - \widehatbnu_{0}) \geq 0
\right\}. 
\end{equation}
Remarkably, this problem admits a closed-form solution:
\begin{equation}\label{equation:sol_safe_sgamm_bxi}
S(\gamma, \bx_i) =
\begin{cases}
\widehatbnu_{0}{}^\top \bx_i + A_i B, &\text{ if } \normtwo{\bg}^2 \normtwo{\bx_i} < C \bx_i^\top \bg; \\
\by^\top \bx_i - \normtwo{\bx_i} C, & \text{ if } \normtwo{\bg}^2 \normtwo{\bx_i} \geq  C \bx_i^\top \bg,
\end{cases}
\end{equation}
where $A_i \triangleq  ( \normtwo{\bx_i}^2 - {(\bx_i^\top \bg)^2}/{\normtwo{\bg}^2} )^{1/2}$, 
$B \triangleq ( C^2 - \normtwo{\bg}^2 )^{1/2}$, and $C \triangleq ( \normtwo{\by}^2 - 2\gamma )^{1/2}$.

\begin{proof}[of \eqref{equation:sol_safe_sgamm_bxi}]
We derive the result by forming the Lagrangian dual of problem~\eqref{equation:safe_sgamm_bxi}. Introducing nonnegative dual variables $\mu_1$ and $\mu_2$ for the two constraints, we obtain:
\begin{align*}
S(\gamma, \bx_i) 
&= \max_{\mu_1,\mu_2 \geq 0} L(\mu_1, \mu_2)
\triangleq \max_{\mu_1,\mu_2 \geq 0} \min_\bnu \bx_i^\top \bnu - \mu_1 (D(\bnu) - \gamma) - \mu_2 \bg^\top (\bnu - \widehatbnu_{0}) \\
&= \max_{\mu_1,\mu_2 \geq 0} \left\{\{\mu_1 \gamma + \mu_2 \bg^\top \widehatbnu_{0}\}
+ \min_\bnu \{\bx_i^\top \bnu - \mu_1 D(\bnu) - \mu_2 \bg^\top \bnu\}\right\} \\
&= \max_{\mu_1,\mu_2 \geq 0} 
\left\{\{\mu_1 \gamma + \mu_2 \bg^\top \widehatbnu_{0}\} 
+ \mu_1 \min_\bnu \left( \frac{\bx_i^\top - \mu_1 \by^\top - \mu_2 \bg^\top}{\mu_1} \bnu + \frac{1}{2} \normtwo{\bnu}^2 \right)
\right\}
.
\end{align*}
Therefore,
$$
S(\gamma, \bx_i) = \max_{\mu_1,\mu_2 \geq 0} L(\mu_1, \mu_2),
$$
where by $\bg=\by-\widehatbnu_0$, we have 
\begin{equation}\label{equation:dual_safeL}
L(\mu_1, \mu_2) = \bx_i^\top \by - \frac{\mu_1}{2} C^2 - \frac{1}{2\mu_1} \normtwo{\bx_i}^2 - \frac{\mu_2^2}{2\mu_1} \normtwo{\bg}^2 + \frac{\mu_2}{\mu_1} \bx_i^\top \bg 
- \mu_2 \normtwo{\bg}^2,
\end{equation}
and $C \triangleq ( \normtwo{\by}^2 - 2\gamma )^{1/2}$.
To maximize $L(\mu_1, \mu_2)$, we first differentiate~\eqref{equation:dual_safeL} w.r.t. $\mu_2$ and set it to zero:
$$
-\mu_2 \normtwo{\bg}^2 + \bx_i^\top \bg - \mu_1 \normtwo{\bg}^2 = 0.
$$
Since $\mu_2\geq0$, we obtain the threshold condition: $\mu_2 = \max\{0,  {\bx_i^\top \bg}/{\normtwo{\bg}^2-\mu_1}\}$. 
This leads to two cases.
When $\mu_1 \geq {\bx_i^\top \bg}/{\normtwo{\bg}^2}$, we have $\mu_2 = 0$, $\mu_1 = {\normtwo{\bx_i}}/{C}$, and consequently, $S(\gamma, \bx_i)$ takes the value:
$$
S(\gamma, \bx_i) =  \bx_i^\top\by - \normtwo{\bx_i} C.
$$
On the other hand, when $\mu_1 \leq  {\bx_i^\top \bg}/{\normtwo{\bg}^2}$, we have $\mu_2={\bx_i^\top \bg}/{\normtwo{\bg}^2-\mu_1}$. 
Substituting this into $L(\mu_1, \mu_2)$ and differentiating with respect to $\mu_1$ leads to the optimality condition:
$$
B^2 \mu_1^2 = A_i^2,
$$
where $A_i \triangleq ( \normtwo{\bx_i}^2 - {(\bx_i^\top \bg)^2}/{\normtwo{\bg}^2} )^{1/2}$ and $B \triangleq ( C^2 - \normtwo{\bg}^2 )^{1/2}$. Substituting $\mu_1$ and $\mu_2$ in~\eqref{equation:dual_safeL}, $S(\gamma, \bx_i)$ takes the value:
$
S(\gamma, \bx_i) = \widehatbnu_{0}^\top \bx_i + A_i B.
$ 
\end{proof}

The analysis above leads directly to the SAFE rule for LASSO.
\begin{theoremHigh}[SAFE--LASSO \citep{ghaoui2010safe}]\label{theorem:safe_lasso}
Consider the LASSO problem \eqref{opt:ll} with penalty parameter $\lambda$. Let $\lambda_{0} \geq \lambda$  be a regularization parameter for which an optimal primal solution $\widehatbbeta_{0} \in \real^p$ is known. Denote by $\bx_i$ the $i$-th column of the design matrix $\bX$. Define the index set
$$
\sI = \left\{ i \mid \lambda > \max(S(\gamma, \bx_i), S(\gamma, -\bx_i)) \right\},
$$
where the function $S(\gamma, \bx_i)$ is given by
$$
S(\gamma, \bx_i) =
\begin{cases}
\widehatbnu_{0}{}^\top \bx_i + A_i B(\gamma), &\text{ if } \normtwo{\bg}^2 \normtwo{\bx_i} < C(\gamma) \bx_i^\top \bg; \\
\by^\top \bx_i - \normtwo{\bx_i} C(\gamma), & \text{ if } \normtwo{\bg}^2 \normtwo{\bx_i} \geq  C(\gamma) \bx_i^\top \bg,
\end{cases}
$$
with
\begin{align*}
\widehatbnu_{0} 
&=  \by-\bX \widehatbbeta_{0}, \quad \bg \triangleq  \by-\widehatbnu_{0} , 
\quad a_0 \triangleq \normtwo{\widehatbnu_{0}}^2, \\
b_0 
&\triangleq \abs{\by^\top \widehatbnu_{0}}, \quad
\gamma \triangleq-\frac{a_0}{2} \left(\frac{\lambda}{\lambda_{0}}\right)^2 +b_0 \frac{\lambda}{\lambda_{0}}, 
\\
C(\gamma) 
&= \left( \normtwo{\by}^2 - 2\gamma \right)^{1/2}, \quad
B(\gamma) = \left( C(\gamma)^2 - \normtwo{\bg}^2 \right)^{1/2}, \quad
A_i \triangleq \left( \normtwo{\bx_i}^2 - \frac{(\bx_i^\top \bg)^2}{\normtwo{\bg}^2} \right)^{1/2}.
\end{align*}
Then, for every index $i \in \sI$, the corresponding entry of $\widehatbbeta$ is zero: $\widehatbeta_i = 0$.
Consequently, the feature $\bx_i$ can be safely eliminated from $\bX$ before solving the full LASSO problem~\eqref{opt:ll}.
\end{theoremHigh}

In practice, if no prior solution $\widehatbbeta_{0}$ is available, we can initialize using the trivial solution $\widehatbbeta_{0}=\bzero$, which corresponds to the maximal regularization parameter  $\lambda_{0} = \lambda_{\max} \triangleq \norminf{\bX^\top \by}$. 
In this case, the SAFE test simplifies significantly. Specifically, the condition $\lambda > \max(S(\gamma, \bx_i), S(\gamma, -\bx_i))$ reduces to $\lambda > \rho_i \lambda_{\max}$, where
\begin{equation}\label{equation:safe_pracmax}
\rho_i \triangleq \frac{\normtwo{\by} \normtwo{\bx_i} + \abs{\by^\top \bx_i}}{\normtwo{\by} \normtwo{\bx_i} + \lambda_{\max}}.
\end{equation}

For scaled datasets, where the response and all features are normalized such that $\normtwo{\by} = 1$ and $\normtwo{\bx_i} = 1$ for every $i$, the quantity $\rho_i$ admits a clean geometric interpretation: 
\begin{equation}\label{equation:safe_pracmax2}
\rho_i = \frac{1 + \abs{\cos \alpha_i}}{1 + \max\limits_{1 \leq j \leq p} \abs{\cos \alpha_j}},
\end{equation}
where $\alpha_i$ denotes the angle between the $i$-th feature $\bx_i$ and the response vector $\by$. 
Thus, the SAFE test eliminates features based on their relative alignment with the response---specifically, how closely they align compared to the most aligned feature. In this setting, the SAFE rule closely resembles standard correlation-based feature selection methods \citep{ghaoui2010safe}.

\section{Optimality Condition of $\ell_1$-Minimization Problems}

As a recap, given measurements $\by$ and the assumption that the original signal $\bbeta$ is sparse or compressible, a natural approach to recover $\bbeta$ is to solve an optimization problem of the form
\begin{equation}\label{equation:optp0_in_opt}
\widehatbbeta = \argmin_{\bbeta} \normzero{\bbeta} \quad \text{s.t.} \quad \bbeta \in \mathcalB(\by),
\end{equation}  
where $\mathcalB(\by)$ ensures that the reconstructed signal $\widehatbbeta$ is consistent with the measurements $\by$. 
For example, if the measurements are exact and noise-free, we can define $\mathcalB(\by) = \{\bbeta \mid \bX\bbeta = \by\}$. 
In the presence of small bounded noise, we might instead use $\mathcalB(\by) = \{\bbeta \mid \normtwo{\bX\bbeta-\by} \leq \epsilon\}$. 
In both cases, \eqref{equation:optp0_in_opt} seeks the sparsest solution $\bbeta$ that satisfies the measurement constraints.

To make the problem more tractable, a common strategy is to replace the $\ell_0$-norm with its convex surrogate, the $\ell_1$-norm. This leads to the relaxed optimization problem:
\begin{equation}\label{equation:optp1_in_opt}
\widehatbbeta = \argmin_{\bbeta} \normone{\bbeta} \quad \text{s.t.} \quad \bbeta \in \mathcalB(\by).
\end{equation}  
Provided that $\mathcalB(\by)$ is convex, \eqref{equation:optp1_in_opt} becomes a convex optimization problem and can be solved efficiently.

The use of $\ell_1$-minimization to promote sparsity has a long history. It dates back at least to Beurling's work on Fourier transform extrapolation from partial observations \citep{beurling1938integrales}. In a different context, \citet{logan1965properties} showed in 1965 that a bandlimited signal can be perfectly recovered even when corrupted on a small interval (with later extensions provided by \citet{donoho1992signal}). His method involved finding the bandlimited signal closest to the observed data in the $\ell_1$-norm---further supporting the intuition that $\ell_1$ is well-suited for handling sparse errors.

Thus, there are multiple compelling reasons to expect that $\ell_1$-minimization yields accurate solutions for sparse signal recovery. Moreover, it provides a computationally efficient alternative to the intractable $\ell_0$-minimization problem. In this section, we present an optimality analysis of $\ell_1$-minimization, followed by a later discussion of algorithms for solving it in Chapter~\ref{chapter:spar_recov}.

\subsection{Necessary and Sufficient Condition}
Specifically, we now consider  the \textcolor{black}{system $\bX \in \real^{n \times p}$} and the $\ell_1$-minimization problem \eqref{opt:p1} (p.~\pageref{opt:p1}):
\begin{equation}\label{equation:optp1_optcond}
(\text{P}_1)
\quad 
\min_{\bbeta \in \real^p} \normone{\bbeta} 
\quad \text{s.t.}\quad 
\bX\bbeta = \by.
\end{equation}
A vector $\widehatbbeta$ is a solution to \eqref{opt:p1} if and only if:
$$
\normonebig{\widehatbbeta +\bn} \geq \normonebig{\widehatbbeta}, 
\quad \forall\,\bn \in \real^p \quad \text{s.t.} \quad \bn \in \nspace(\bX) \quad \text{and}\quad \by=\bX\widehatbbeta,
$$
where $\bn$ is in the nullspace of $\bX$, i.e., $\bX\bn = \bzero$.
Given that $\widehatbbeta$ is supported on $\sS \triangleq\supp(\widehatbbeta) = \{i \mid  \widehatbeta_i \neq 0\}$ and the inequality  $\absbig{\widehatbeta_i + n_i} \geq \absbig{\widehatbeta_i} + \sgn(\widehatbeta_i) n_i$, we have
$$
\normone{\widehatbbeta + \bn} = \sum_{i \in \sS} \abs{\widehatbeta_i + n_i} + \sum_{i \in \comple{\sS}} \abs{n_i}
\geq \sum_{i \in \sS} \abs{\widehatbeta_i} + \sum_{i \in \sS} \sgn(\widehatbeta_i) n_i + \sum_{i \in \comple{\sS} } \abs{n_i}.
$$
This analysis indicates the necessary and sufficient condition for \eqref{opt:p1}.
\begin{theoremHigh}[Necessary and sufficient condition for \eqref{opt:p1}]\label{theorem:nec_suf_p1}
Consider the optimization problem~\eqref{opt:p1}.
Then, $\widehatbbeta$ is a solution of \eqref{opt:p1} if and only if, for  all $\bn \in \nspace(\bX)$, it follows that
$$
\sum_{i \in \sS} \sgn(\widehatbeta_i) n_i \leq \sum_{i \in \comple{\sS}} \abs{n_i},
$$
where $\sS \triangleq\supp(\widehatbbeta) = \{i \mid  \widehatbeta_i \neq 0,\, i\in\{1,2,\ldots,p\}\}$.
\end{theoremHigh}
\begin{proof}[of Theorem~\ref{theorem:nec_suf_p1}]
\textbf{Necessity.} 
Suppose $\widehatbbeta$ is a solution of \eqref{opt:p1}. 
Choose $\gamma >0 $ sufficiently small so that the sign of each nonzero component $\widehatbeta_i$ is preserved in $\widehatbeta_i-\gamma n_i$.
Then for $i\in\{1,2,\ldots,p\}$,
$$
\abs{\widehatbeta_i - \gamma n_i} = 
\begin{cases} 
\widehatbeta_i - \gamma n_i = \widehatbeta_i - \gamma \sgn(\widehatbeta_i) n_i, & \text{if } \widehatbeta_i > 0; \\
-(\widehatbeta_i - \gamma n_i) = -\widehatbeta_i - \gamma \sgn(\widehatbeta_i) n_i, & \text{if } \widehatbeta_i < 0; \\
\gamma \abs{n_i}, & \text{otherwise},
\end{cases}
\qquad \forall\, i\in\{1,2,\ldots,p\}.
$$
Now assume, for contradiction, that there exists $\bn\in\nspace(\bX)$ such that $\sum_{i \in \sS} \sgn(\widehatbeta_i) n_i > \sum_{i \in \comple{\sS}} \abs{n_i}$.
This implies 
$$
\normone{\widehatbbeta - \gamma\bn} = \normonebig{\widehatbbeta} - \gamma \sum_{i \in \sS} \sgn(\widehatbeta_i) n_i + \gamma \sum_{i \in \comple{\sS}} \abs{n_i} < \normonebig{\widehatbbeta},
$$
leading to a contradiction. Hence, the stated inequality must hold.

\paragraph{Sufficiency.}
Suppose that there exists another feasible point $\bbeta'$ such that $\bX\bbeta' = \by$ and $\normone{\bbeta'} < \normonebig{\widehatbbeta}$.
Define $\bn \triangleq \widehatbbeta - \bbeta'$, which satisfies $\bX\bn = \bzero$ since $\bX\bbeta' = \bX\widehatbbeta$. This means $\bn \in \nspace(\bX)$.
Using the identity $\normonebig{\bbeta'} = \normonebig{\widehatbbeta - \bn}$, we analyze:
\begin{align*}
\normonebig{\bbeta'}=\normone{\widehatbbeta  - \bn} 
&= \sum_{i \in \sS} \abs{\widehatbeta_i - n_i} + \sum_{i \in \comple{\sS}} \abs{n_i}\\
&\geq \sum_{i \in \sS} \absbig{\widehatbeta_i} - \sum_{i \in \sS} \sgn(\widehatbeta_i)n_i + \sum_{i \in \comple{\sS}} \abs{n_i}
\geq \sum_{i \in \sS} \abs{\widehatbeta_i} =\normonebig{\widehatbbeta},
\end{align*}
where the first inequality follows from the fact that $\absbig{\widehatbeta_i - n_i} \geq \absbig{\widehatbeta_i} - \sgn(\widehatbeta_i)n_i$,
and the second inequality follows from the given condition that  $\sum_{i \in \sS} \sgn(\widehatbeta_i)n_i \leq \sum_{i \in \comple{\sS}} \abs{n_i}$.
Since the assumption was $\normone{\bbeta'} < \normonebig{\widehatbbeta}$, we obtain a contradiction.
This completes the proof.
\end{proof}

\subsection{Dual Certificate and Uniqueness\index{Dual certificate}}\label{section:dual_cert_ell1}

We introduced the dual of the $\ell_\infty$-minimization problem in Example~\ref{example:linf_min_dual}. 
Let $\bX\in\real^{n\times p}$.
The dual of the optimization problem of
\begin{subequations}
\begin{equation}
\min_{\bbeta\in \real^p} \norminf{\bbeta } \quad \text{s.t.} \quad \bX\bbeta = \by
\end{equation}
is
\begin{equation}\label{equation:dual_pinf_eq1}
\max_{\blambda \in \real^n}  \innerproduct{\by,\blambda} \quad \text{s.t.} \quad \normone{\bX^\top \blambda} \leq 1.
\end{equation}
\end{subequations}
The dual of the $\ell_1$-minimization problem can be derived analogously.
\begin{proposition}[Dual of $\ell_1$-minimization]\label{proposition:dual_p1}
Let $\bX\in\real^{n\times p}$.
The dual of the optimization problem~\eqref{opt:p1} 
\begin{subequations}
\begin{equation}\label{equation:ell1_dualprop}
(\text{P}_1)
\quad 
\min_{\bbeta\in \real^p} \normone{\bbeta } \quad \text{s.t.} \quad \bX\bbeta = \by
\end{equation}
is
\begin{equation}\label{equation:dual_p1_eq1}
\max_{\blambda \in \real^n} 
D(\blambda) 
= 
\begin{cases}
\innerproduct{\blambda, \by}, & \text{if } \norminf{\bX^\top \blambda} \leq 1; \\
-\infty, & \text{otherwise}.
\end{cases}
\end{equation}
Clearly, it is enough to maximize over the points $\blambda$ for which $D(\blambda) > -\infty$. Making this constraint explicit, the dual program to \eqref{equation:ell1_dualprop} is given by
\begin{equation}\label{equation:dual_p1_eq2}
\max_{\blambda \in \real^n} \innerproduct{\blambda, \by} \quad \text{s.t.}\quad \norminfbig{\bX^\top \blambda} \leq 1.
\end{equation}
\end{subequations}
\end{proposition}

\begin{proof}[of Proposition~\ref{proposition:dual_p1}]
The Lagrangian associated with \eqref{opt:p1}  is
$$
L(\bbeta, \blambda) = \normone{\bbeta} +\innerproduct{\blambda, \by - \bX\bbeta}.
$$
The Lagrange dual function is therefore
$$
D(\blambda) \triangleq \min_{\bbeta \in \real^p} \normone{\bbeta} 
-\innerproduct{\bX^\top \blambda, \bbeta} +\innerproduct{\blambda, \by}. 
$$
If $(\bX^\top \blambda)_i > 1$, then one can set $\beta_i$ arbitrarily large so that $D(\blambda) \to -\infty$. The same happens if $(\bX^\top \blambda)_i < -1$. If $\norminf{\bX^\top \blambda} \leq 1$, by \holders inequality (Theorem~\ref{theorem:holder-inequality})
\begin{align*}
\abs{(\bX^\top \blambda)^\top \bbeta} \leq \normone{\bbeta} \cdot \norminfbig{\bX^\top \blambda}
\leq \normone{\bbeta}, 
\end{align*}
so the Lagrangian is minimized by setting $\bbeta$ to zero and $D(\blambda) = \innerproduct{\blambda, \by}$. This completes the proof.
\end{proof}

By Theorem~\ref{theorem:slater_cond}, strong duality holds for this pair of primal and dual optimization problems provided the primal problem \eqref{equation:ell1_dualprop} is feasible.

Interestingly, the solution to the dual problem  of the $\ell_1$-minimization problem often provides valuable information about the support of the primal optimal solution.
Let the primal optimal solution $\widehatbbeta$ be supported on 
$$
\sS \triangleq\supp(\widehatbbeta) = \{i \mid  \widehatbeta_i \neq 0\}.
$$
Note that $ \bX\widehatbbeta = \bX_{\sS} \widehatbbeta_{\sS} $, where $ \widehatbbeta_{\sS} \triangleq \widehatbbeta[\sS] \in\real^{\abs{\sS}}$ denotes a subvector of $\widehatbbeta$ indexed by $\sS$, and  $ \bX_{\sS} \triangleq \bX[:,\sS] \in\real^{n\times \abs{\sS}}$ denotes a submatrix of $\bX$ consisting of the columns indexed by $ \sS $.

\begin{lemma}[Support of primal optimum using dual optimum]\label{lemma:dual_cer_p1_prem}
Suppose the primal problem~\eqref{equation:ell1_dualprop} is feasible.
Let $\widehatblambda$ be an optimal solution to the dual problem~\eqref{equation:dual_p1_eq2}. Then, for every primal optimal solution $\widehatbbeta$,
\begin{equation}
(\bX^\top \widehatblambda)_i = \sgn(\widehatbeta_i), 
\quad \text{for all } \widehatbeta_i \neq 0 \ (\text{i.e., } i \in \sS).
\end{equation}
\end{lemma}
\begin{proof}[of Lemma~\ref{lemma:dual_cer_p1_prem}]
Feasibility of the primal implies strong duality (by Theorem~\ref{theorem:slater_cond}). Hence,
\begin{align*}
\normonebig{\widehatbbeta} 
&= \by^\top \widehatblambda 
= (\bX\widehatbbeta)^\top \widehatblambda 
= (\widehatbbeta)^\top (\bX^\top \widehatblambda) 
= \sum_{i=1}^p (\bX^\top \widehatblambda)_i \widehatbeta_i. 
\end{align*}
On the other hand, by` \holders inequality, we have
$$
\sum_{i=1}^p (\bX^\top \widehatblambda)_i \widehatbeta_i 
\leq \norminfbig{\bX^\top \widehatblambda}\normonebig{\widehatbbeta}
\leq \normonebig{\widehatbbeta} ,
$$
with equality if and only if
$
(\bX^\top \widehatblambda)_i = \sgn(\widehatbeta_i)$, $\text{for all } \widehatbeta_i \neq 0. 
$
This completes the proof.
\end{proof}

In the study of sparse recovery or $\ell_1$-regularized optimization, a central question is whether an optimal solution to the $\ell_1$-minimization problem \eqref{opt:p1} is \textbf{unique}. 
While existence of solutions is guaranteed under mild conditions due to convexity, uniqueness often hinges on additional structural assumptions---particularly those related to the geometry of the data matrix $\bX$ and the support of the solution. 
The previous analysis established necessary optimality conditions, revealing that the optimal solution $\widehatbbeta$ 
satisfies certain sign and dual feasibility conditions involving the Lagrange multipliers. In particular, the vector $\bX^\top \blambda$ plays a crucial role in characterizing the support and sparsity pattern of $\widehatbbeta$. Building upon these insights, we now turn to the important question: under what conditions is this optimal solution uniquely determined? This leads us to the following key result, which provides sufficient conditions for the uniqueness of the solution to the $\ell_1$-minimization problem \eqref{opt:p1} based on the properties of the \textit{dual certificate}  and sheds lights on the \textit{restricted isometry property (RIP)} of the design matrix (Chapters~\ref{chapter:design} and \ref{chapter:recovery}).

Consider a certain sparse vector $\widehatbbeta \in \real^p$ with support $\sS$ such that $\bX\widehatbbeta = \by$. If there exists a vector $\blambda \in \real^n$ such that $\bX^\top \blambda$ is equal to the sign of $\widehatbbeta$ on $\sS$ and has magnitude smaller than one elsewhere, then $\blambda$ is feasible for the dual problem~\eqref{equation:dual_p1_eq1}, so by weak duality, $\normone{\bbeta}\geq \by^\top \blambda$ for any $\bbeta \in \real^p$ that is feasible for the primal problem \eqref{equation:ell1_dualprop}. 
By Lemma~\ref{lemma:dual_cer_p1_prem}, we then have
\begin{align}
\normone{\bbeta}
&\geq \by^\top \blambda 
= (\bX\widehatbbeta)^\top \blambda 
= (\widehatbbeta)^\top (\bX^\top \blambda) 
= \sum_{i=1}^p \widehatbeta_i \, \sgn(\widehatbeta_i) 
= \normonebig{\widehatbbeta}. 
\end{align}
Since $(\bX^\top \blambda)_i = \sgn(\widehatbeta_i)$ for all $\widehatbeta_i\neq 0$, geometrically, $\bX^\top \blambda$ is a subgradient of the $\ell_1$-norm at $\widehatbbeta$. The subgradient is \textbf{orthogonal} to the feasibility hyperplane given by $\bX\bbeta = \by$ (any vector $\bn$ within the hyperplane is the difference between two feasible vectors and therefore satisfies $\bX\bn = \bzero$). As a result, for any other feasible vector $\bbeta$
\begin{align}
\normone{\bbeta}
\geq \normonebig{\widehatbbeta} + (\bX^\top \blambda)^\top (\bbeta - \widehatbbeta) 
= \normonebig{\widehatbbeta} + \blambda^\top (\bX\bbeta - \bX\widehatbbeta)
= \normonebig{\widehatbbeta},
\end{align}
where the first inequality follows from the subgradient inequality (Definition~\ref{definition:subgrad}).

Alternatively, 
we have discussed the KKT conditions of convex optimization problems with equality constraints in Remark~\ref{remark:kkt_nutshell_cvx}, which is a result of  Slater's condition for convex optimization problems and the KKT conditions under strong duality (Theorems~\ref{theorem:slater_cond} and \ref{theorem:opt_cond_sd}).
Since $f(\bbeta) = \normone{\bbeta}$ is convex, Remark~\ref{remark:kkt_nutshell_cvx} and Equation~\eqref{equation:norm1_grad_lass} show that
$ \widehatbbeta $ is an optimal solution if $ \widehatbbeta $ is feasible and there exists $ \balpha= \bX^\top \blambda$, where $\blambda\in\real^n$ denotes the Lagrangian multipliers  (i.e., $\balpha$ lies in the row space of $\bX$, and $\balpha\perp \nspace(\bX)$)  with
\begin{equation}\label{equation:kkt_p1_bb}
\begin{cases}
\alpha_i = \sgn(\widehatbeta_i), & \widehatbeta_i \neq 0 \ (\text{i.e., } i \in \sS); \\
\abs{\alpha_i} \leq 1, & \widehatbeta_i = 0 \ (\text{i.e., } i \in \comple{\sS}),
\end{cases}
\end{equation}
where $\widehatbbeta$ is supported on $\sS \triangleq\supp(\widehatbbeta) = \{i \mid  \widehatbeta_i \neq 0\}$.

These two arguments show that the existence of a certain dual vector can be used to verify  that a given  primal feasible vector is  optimal, but they do not guarantee  uniqueness. 
It turns out that requiring  the magnitude of $\bX^\top \blambda$ to be strictly smaller than one on $\comple{\sS}$ is enough to guarantee it (as long as $\bX_{\sS}$ has full rank). In that case, we call the dual variable $\blambda$ a \textit{dual certificate} for the $\ell_1$-minimization problem.

\begin{theoremHigh}[Dual certificate for \eqref{opt:p1}]\label{theorem:dualcert_p1}
Consider the $\ell_1$-minimization problem \eqref{opt:p1}. 
Let $\widehatbbeta$ be an optimal solution of \eqref{opt:p1} such that  $\balpha=\bX^\top\blambda$ satisfies the following conditions:
Suppose further that
\begin{enumerate}[(i)]
\item $\alpha_i=(\bX^\top \blambda)_i = \sgn(\widehatbeta_i)$,  when $\widehatbeta_i \neq 0$ for all $i\in\{1,2,\ldots,p\}$.
\item $ \abs{\alpha_i} =\absbig{(\bX^\top \blambda)_i} < 1 $, when $ \widehatbeta_i = 0 $ for all $i\in\{1,2,\ldots,p\}$.
\item $ \bX_{\sS} $ has full column rank~\footnote{which is also implied by the RIP condition; see Chapter~\ref{chapter:design} and Lemma~\ref{lemma:dual_rec_linfty}.}, where $\sS \triangleq\supp(\widehatbbeta) = \{i \mid  \widehatbeta_i \neq 0\}$.
\end{enumerate}
\noindent Then $ \widehatbbeta $ is the \textbf{unique} solution of \eqref{opt:p1}. 
Any such $ \balpha $ or $ \blambda $ satisfying (i)--(iii) is called a dual certificate.
\end{theoremHigh}
\begin{proof}[of Theorem~\ref{theorem:dualcert_p1}]
Let $\bbeta \in \real^p$ be any feasible point, and define $\bn \triangleq \bbeta - \widehatbbeta$. 
Since $\bX_{\sS}$ is full-rank, we have $\bn_{\comple{\sS}} \neq \bzero$ unless $\bn = \bzero$, because otherwise $\bn_{\sS}$ would be a nonzero vector in the null space of $\bX_{\sS}$. 
Condition (ii) implies
$$
\normone{\bn_{\comple{\sS}}} > \balpha^\top \bn(\comple{\sS}),
$$
where $\bn(\comple{\sS})\in\real^p$ a projection that sets to zero all entries of a vector except the ones indexed by $\comple{\sS}$ (see Definition~\ref{definition:matlabnotation}). 
Using the fact that $\absbig{\widehatbeta_i + n_i} \geq \absbig{\widehatbeta_i} + \sgn(\widehatbeta_i)n_i$
and that $\widehatbbeta$ is supported on  $\sS$, the above inequality implies
\begin{align*}
\normone{\bbeta}
&= \normone{\widehatbbeta_{\sS} + \bn_{\sS}} + \normone{\bn_{\comple{\sS}}}  
> \normonebig{\widehatbbeta} + \balpha^\top \bn(\sS) + \balpha^\top \bn(\comple{\sS})  
= \normonebig{\widehatbbeta} + \blambda^\top \bX \bn
= \normonebig{\widehatbbeta}. 
\end{align*}

On the other hand, if $ \bn_{\comple{\sS}} = \bzero $, then since $ \bX_{\sS} $ has full column rank,
$$ 
\bX\bn = \bX_{\sS} \bn_{\sS} = \bzero 
\quad\implies\quad 
\bn_{\sS} = \bzero 
\quad\implies\quad 
\bn = \bzero.
$$
Therefore, for any $ \bn \in \nspace(\bX) $, $ \normonebig{\widehatbbeta+\bn} > \normonebig{\widehatbbeta} $ unless $ \bn = \bzero $.
This completes the proof.
\end{proof}

%

Note that ideally we would like to analyze the vector $\blambda$ satisfying this underdetermined system of $\abs{\sS}$ equations such that $\bX_\sS^\top \blambda$ has the smallest possible $\ell_\infty$-norm. Unfortunately, the solution to the optimization problem
\begin{equation}\label{equation:xslambda_inf}
\min_{\blambda} \norminf{\bX_\sS^\top \blambda}\quad \text{s.t.} \quad \bX_\sS^\top \blambda = \sign(\widehatbbeta_\sS)
\end{equation}
does not have a closed-form solution. However, the related least squares problem
\begin{equation}
\min_{\blambda} \normtwo{\blambda} \quad \text{s.t.} \quad \bX_\sS^\top \blambda = \sign(\widehatbbeta_\sS)  
\end{equation}
does have an explicit solution, so we can analyze it instead. Note that since $ \bX_{\sS} $ has full column rank ($ \bX_{\sS}^\top \bX_{\sS} $ invertible), {the solution is}
\begin{equation}\label{equation:l1_dualcert_lambda}
\blambda = \bX_\sS (\bX_\sS^\top \bX_\sS)^{-1} \sign(\widehatbbeta_\sS)  
\qquad \text{and} \qquad
\balpha\triangleq \bX^\top \blambda,
\end{equation}
which  is a high-dimensional LS problem (Corollary~\ref{corollary:high_dim_ls}) of the form $\bX_\sS^\top\blambda = \sign(\widehatbbeta_\sS)$.
Therefore, 
\begin{itemize}
\item if $ \abs{\alpha_i} \leq 1 $ for all $ i \in \comple{\sS} $, then $ \widehatbbeta $ is a solution;
\item if $ \abs{\alpha_i} < 1 $ for all $ i \in \comple{\sS} $, then $ \widehatbbeta $ is the unique solution.
\end{itemize}

For a given coefficient vector $\widehatbbeta\in\real^p$, the conditions in Theorem~\ref{theorem:dualcert_p1} can be restated  as follows: 
The index sets $\sS \triangleq \supp(\widehatbbeta)$ and $\sT \triangleq \comple{\sS}$ satisfy:
\begin{enumerate}[(i)]
\item $\nspace(\bX(\sS)) = \{\bzero\}$: this guarantees  that the linear system $\bX_{\sS} \bbeta_{\sS} = \by$ has a unique solution.
\item There exists a vector $\balpha\in\real^p$ such that $\balpha \in \cspace(\bX^\top)$, $\balpha_{\sS} = \sgn(\widehatbbeta_{\sS})$, and $\norminf{\balpha_\sT} < 1$: this ensures that any optimal solution must satisfy $\widehatbbeta_\sT = \bzero$.
\end{enumerate}
In fact, the \textit{restricted isometry property (RIP)} guarantees that such conditions hold for all sufficiently small supports $\sS$ and for arbitrary sign patterns; see Definition~\ref{definition:rip22} and Section~\ref{section:spar_rec_rip} for further details.

The argument in Theorem~\ref{theorem:dualcert_p1} can thus be naturally split into two parts.

\paragrapharrow{Condition (i) $\nspace(\bX(\sS)) = \{\bzero\}$ is necessary for the uniqueness.}

If $\bzero \neq \bn \in \nspace(\bX(\sS))\in\real^p$~\footnote{Once again, note that $\bX(\sS)\in\real^{n\times p}$ and $\bX_{\sS}=\bX[:,\sS]\in\real^{n\times \abs{\sS}}$; see Definition~\ref{definition:matlabnotation}.}, 
then all $\bbeta_\eta \triangleq  \widehatbbeta + \eta \widetildebn$ for any scalar $\eta$ is also optimal, where $\widetildebn_{\sS} = \bn_\sS$ and $\widetildebn_{\comple{\sS}}=\bzero$.
To see this, we have that $\bbeta_\eta$ is feasible since $\bX\bbeta_\eta = \bX\widehatbbeta = \by$.
Since $\widehatbbeta$ is the unique minimizer, we have $\normone{\bbeta_\eta} > \normonebig{\widehatbbeta}$. 
But for small $\eta$ around 0, namely $\abs{\eta} < \min_{i \in \sS} \absbig{\widehatbeta_i} / \norminf{\bn}$, it follows that $\sgn(\widehatbbeta + \eta \widetildebn) = \sgn(\widehatbbeta)$, whence we have 
$$
\normone{\bbeta_\eta} = \normone{\widehatbbeta_{\sS} + \eta\bn_\sS} 
= \innerproduct{ \sgn(\widehatbbeta_{\sS}),  (\widehatbbeta_{\sS} + \eta\bn_\sS) }
= \normonebig{\widehatbbeta} +   \innerproduct{\sgn(\widehatbbeta_{\sS}), \eta\bn_\sS}.
$$
This is a contradiction, because we can always choose a small $\eta \neq 0$ such that $\innerproductbig{\sgn(\widehatbbeta_{\sS}), \eta\bn_\sS}<0$.
Therefore, $\bX_{\sS}$ must have full column rank.

\paragrapharrow{Condition (ii) guarantees that $\widehatbbeta_\sT = \bzero$.}
To see this, we introduce the following optimization problem:
\begin{equation}\label{equation:dual_cer_nomrinf}
\min_{\balpha\in\real^p} \norminf{\balpha_\sT}
\quad \text{s.t.} \quad \balpha \in \cspace(\bX^\top), \; \balpha_{\sS} = \sgn(\widehatbbeta_{\sS}). 
\end{equation}
If the optimal objective value is strictly less than $1$, then there exists $\balpha$ obeying Condition (ii), so Condition (ii) is also necessary.

Let $\bv\in\real^p$ with $\bv_{\sS}=\sgn(\widehatbbeta_{\sS})$ and $\bv_{\sT}=\bzero$, and let $\bQ\in\real^{p\times r}$ contains the mutually orthonormal basis of $\nspace(\bX)$.
If $\bv \in \cspace(\bX^\top)$, set $\balpha \triangleq  \bv$, and the proof is complete.
Otherwise, let $\balpha \triangleq \bv + \bu$ for some $\bu\in\real^p$. Then the following problems are equivalent
$$
\balpha \in \cspace(\bX^\top) 
\quad\iff\quad \bQ^\top \balpha = \bzero 
\quad\iff\quad \bQ^\top \bu = -\bQ^\top \bv,
$$
where $\bu_{\sS} = \bzero$ since $\balpha_{\sS} = \sgn(\widehatbbeta_{\sS}) = \bv_{\sS}$. Since $\bv_\sT = \bzero$, we also have $\norminf{\balpha_\sT} = \norminf{\bu_\sT}$.
We can therefore translate \eqref{equation:dual_cer_nomrinf} and rewrite $\balpha \in \cspace(\bX^\top)$ as the following equivalent problem:
\begin{equation}\label{equation:dual_cer_nomrinf_equiv}
\min_{\bu\in\real^p} \norminf{\bu_\sT} 
\quad \text{s.t.} \quad \bQ^\top \bu = -\bQ^\top \bv,\; \bu_{\sS} = \bzero. 
\end{equation}
If the optimal objective value is strictly less than $1$, then Condition (ii) is necessary.

Problem \eqref{equation:dual_cer_nomrinf_equiv} is feasible and has a finite objective value. The dual of \eqref{equation:dual_cer_nomrinf_equiv} is
\begin{equation}\label{equation:dual_cer_nomrinf_equiv_dd}
\max_{\bgamma\in\real^p} \, (\bQ^\top \bv)^\top \bgamma 
\quad \text{s.t.} \quad \normone{(\bQ\bgamma)_\sT} \leq 1.
\end{equation}
By assumption, the strong duality of  problem \eqref{equation:dual_cer_nomrinf_equiv} holds.
Then if the optimal objective value of \eqref{equation:dual_cer_nomrinf_equiv_dd} is strictly less than $1$,  Condition (ii) is necessary.

\begin{lemma}\label{lemma:necess_cdii}
Suppose $\widehatbbeta$ is unique with support $\sS\triangleq \supp(\widehatbbeta)$ and let $\sT\triangleq \comple{\sS}$. 
Then the optimal objective of the following primal-dual problems is strictly less than 1:
\begin{align*}
\min_{\bu} \norminf{\bu_{\sT}} 
\quad &\text{s.t.} \quad \bQ^\top \bu = -\bQ^\top \bv,\; \bu_{\sS} = \bzero; \\
\max_{\bgamma} \, (\bQ^\top \bv)^\top \bgamma 
\quad &\text{s.t.} \quad \normone{(\bQ\bgamma)_\sT} \leq 1.
\end{align*}
\end{lemma}
\begin{proof}[of Lemma~\ref{lemma:necess_cdii}]
Since $\widehatbbeta$ is the unique solution of the $\ell_1$-minimization problem \eqref{opt:p1},   
$$
\normonebig{\widehatbbeta} < \normone{\widehatbbeta + \bn} 
=\normone{\widehatbbeta_{\sS}+\bn_{\sS}} + \normone{\bn_\sT}  , 
\quad \forall\,\bn \in \nspace(\bX) \setminus \{\bzero\}.
$$ 
Similar to the argument in the proof of Theorem~\ref{theorem:nec_suf_p1}, 
this shows  $\innerproduct{\bv_\sS, \bn_\sS} < \normone{\bn_\sT}$.
Therefore,
\begin{itemize}
\item if $\widehatbgamma  = \bzero$, then $\norminf{\widehatbu_\sT}= (\bQ^\top \bv)^\top \widehatbgamma = \bzero$.
\item if $\widehatbgamma \neq \bzero$, then $\bn \triangleq \bQ\widehatbgamma \in \nspace(\bX) \setminus \{\bzero\}$ satisfies
$$
\norminf{\widehatbu_\sT} = (\bQ^\top \bv)^\top \widehatbgamma 
= \innerproduct{\bv_\sS, \bn_\sS} < \normone{\bn_\sT}
= \normone{(\bQ\widehatbgamma)_\sT} \leq 1.
$$
\end{itemize}
In both cases, the optimal objective value is strictly less than  $1$. 
\end{proof}

\subsection{$\ell_1$-Minimization with Noise Measurements}

We now consider the case where the observations are contaminated by noise; that is, $\by = \bX\bbeta + \bepsilon$. 
In this setting, we employ the following two robust variants of the $\ell_1$-minimization problem, i.e.,  \eqref{opt:p1_penalize}  and \eqref{opt:p1_epsilon} (p.~\pageref{opt:p1_penalize} and \pageref{opt:p1_epsilon}):
\begin{subequations}
\begin{align}
\text{(P$_{1,\sigma}$)} \qquad &\min \normone{\bbeta} + \sigma\normtwo{\bX\bbeta - \by}^2  ;\\
\text{(P$_{1,\epsilon}$)}\qquad & \min \normone{\bbeta} \quad \text{s.t.} \quad \normtwo{\bX\bbeta - \by} \leq \epsilon .
\end{align}
\end{subequations}

As shown in Corollary~\ref{corollary:equigv_p1_epsi_pen_lag}, problems~\eqref{opt:p1_penalize} and~\eqref{opt:p1_epsilon} are equivalent in the sense that, for appropriate choices of $\sigma$ and $\epsilon$, they share the same set of solutions.
The following lemma establishes a stronger property: both the residual norm and the $\ell_1$-norm are constant over the entire solution set of each problem.
\begin{lemma}[Constants over the set of solutions in \eqref{opt:p1_penalize} and  \eqref{opt:p1_epsilon}]\label{lemma:const_conv_p1sigma}
Let $\sigma > 0$. 
Suppose the function $\normone{ \bbeta} + \sigma \normtwo{\bX\bbeta-\by}^2 $ is constant on a convex set $\sC$. 
Then both $\bX\bbeta-\by$ and $\normone{ \bbeta}$ are constants on $\sC$.

Let $\sQ_\sigma$ denote the set of solutions to problem~\eqref{opt:p1_penalize}.
Since $\normone{ \bbeta} + \sigma \normtwo{\bX\bbeta-\by}^2$ is constant over $\sQ_\sigma$, the above result shows that $\bX\bbeta-\by$ and $\normone{ \bbeta}$ are constants on $\sQ_\sigma$ in problem \eqref{opt:p1_penalize}.

Moreover, let  $\sQ_\epsilon$ denote the set of solutions to problem \eqref{opt:p1_epsilon}. Then both $\normtwo{\bX\bbeta-\by}$ and $\normone{ \bbeta}$ are constant on $\sQ_\epsilon$
\end{lemma}

\begin{proof}[of Lemma~\ref{lemma:const_conv_p1sigma}]
It suffices to consider the case where the convex set contains more than one point. 
Let $\bbeta_1,\bbeta_2\in \sC$ with $\bbeta_1\neq \bbeta_2$. 
Consider the line segment $\sL$ connecting $\bbeta_1$ and $\bbeta_2$. 
By  convexity of  $\sC$, we know $\sL \subset \sC$. Thus, $c \triangleq  \normone{ \bbeta} + \sigma \normtwo{\bX\bbeta-\by}^2$ is a constant on $\sL$. 
Suppose $\bX \bbeta_1 - \by \ne \bX \bbeta_2 - \by$. Then for any $0 < \lambda < 1$, we have
\begin{align*}
&\sigma \normtwo{\bX(\lambda \bbeta_1 + (1 - \lambda)\bbeta_2) - \by}^2 + \normonebig{(\lambda \bbeta_1 + (1 - \lambda)\bbeta_2)}  \\
&= \sigma \normtwo{\lambda(\bX \bbeta_1 - \by) + (1 - \lambda)(\bX \bbeta_2 - \by)}^2 + \normonebig{\lambda( \bbeta_1) + (1 - \lambda)( \bbeta_2)}  \\
&< \lambda\left(\sigma \normtwo{\bX \bbeta_1 - \by}^2 + \normonebig{ \bbeta_1}\right) 
+ (1 - \lambda)\left(\sigma \normtwo{\bX \bbeta_2 - \by}^2 + \normonebig{ \bbeta_2}\right)  \\
&= \lambda c + (1 - \lambda)c = c, 
\end{align*}
where the strict inequality follows from the strict convexity of $\sigma \normtwo{\cdot}^2$ and the convexity of $\normone{ \bbeta}$. This means that the points $\lambda \bbeta_1 + (1 - \lambda)\bbeta_2$ on $\sL$ attain a lower value than $c$, which leads to a contradiction. Therefore, we have $\bX \bbeta_1 - \by = \bX \bbeta_2 - \by$, which in turn implies that $\normone{ \bbeta_1} = \normone{ \bbeta_2}$.

For problem \eqref{opt:p1_epsilon}, if $\bzero \in \sQ_\epsilon$, then we have $\sQ_\epsilon = \{\bzero\}$, and the claim holds trivially.
Otherwise, suppose $\bzero \not\in \sQ_\epsilon$. 
Since the optimal objective $\normone{ \bbeta}$ is constant for all $\bbeta \in \sQ_\epsilon$, it remains to show that $\normtwo{\bX\bbeta-\by} = \epsilon$ for all $\bbeta \in \sQ_\epsilon$. 
Suppose on the contrary that there exists a nonzero $\widetildebbeta \in \sQ_\epsilon$ such that $\normtwobig{\bX \widetildebbeta - \by} < \epsilon$, we can find a non-empty ball $\sT$ centered at $\widetildebbeta$ with a sufficiently small radius $r > 0$ such that $\normtwo{\bX \balpha - \by} < \epsilon$ for all $\balpha \in \sT$. 
Let $\gamma \triangleq \min\left(\frac{r}{2\normtwo{ \widetildebbeta }}, \frac{1}{2}\right)$. 
We have $(1 - \gamma)\widetildebbeta \in \sT$ such that $(1 - \gamma)\widetildebbeta$ satisfies the constraint in problem \eqref{opt:p1_epsilon}, implying $\normonebig{(1 - \gamma) \widetildebbeta} < \normonebig{ \widetildebbeta}$, which leads to a contradiction and concludes the result.
\end{proof}

Using Lemma~\ref{lemma:const_conv_p1sigma}, we can establish a dual certificate for the two robust models under the same conditions as in Theorem~\ref{theorem:dualcert_p1} for the noiseless problem~\eqref{opt:p1}.
\begin{theoremHigh}[Dual certificate for \eqref{opt:p1_penalize} and \eqref{opt:p1_epsilon}]\label{theorem:dual_p1pen_epsi}
Suppose $\widehatbbeta$ is a solution to either \eqref{opt:p1_penalize} or \eqref{opt:p1_epsilon}. Then, $\widehatbbeta$ is the unique solution if and only if the same condition stated in Theorem~\ref{theorem:dualcert_p1} holds for $\widehatbbeta$.
\end{theoremHigh}
\begin{proof}[of Theorem~\ref{theorem:dual_p1pen_epsi}]
Since the results of Lemma \ref{lemma:const_conv_p1sigma} are identical for problems \eqref{opt:p1_penalize} and \eqref{opt:p1_epsilon}, we present the proof for problem \eqref{opt:p1_penalize} only.

By assumption, the solution set $\sQ_\sigma$ of \eqref{opt:p1_penalize} is nonempty, so we pick $\widetildebbeta \in \sQ_\sigma$. Let $\widetildeby  = \bX \widetildebbeta$, which is independent of the choice of $\widetildebbeta$ according to Lemma~\ref{lemma:const_conv_p1sigma}. 
Now consider the following noiseless $\ell_1$-minimization problem:
\begin{equation}\label{equation:exact_p1_pen_epsi_pv1}
\min_{\bbeta\in\real^p} \normone{ \bbeta} \quad \text{s.t.} \quad \bX \bbeta = \widetildeby, 
\end{equation}
and let $\sQ_q$ denote its solution set.

Now, we show that $\sQ_\sigma = \sQ_q$. Since $\bX \bbeta = \bX \widetildebbeta$ and $\normone{ \bbeta} = \normonebig{ \widetildebbeta}$ for all $\bbeta \in \sQ_\sigma$ 
and conversely any point $\bbeta$ obeying $\bX \bbeta = \bX \widetildebbeta$ and $\normone{ \bbeta} = \normonebig{ \widetildebbeta}$ must belong to $\sQ_\sigma$, it suffices to show that $\normone{ \bbeta} = \normonebig{ \widetildebbeta}$ for any $\bbeta \in \sQ_q$. 
Assuming this does not hold, then since problem \eqref{equation:exact_p1_pen_epsi_pv1} has $\widetildebbeta$ as a feasible solution and has a finite objective, we have a nonempty $\sQ_q$ and there exists $\bz \in \sQ_q$ satisfying $\normone{ \bz} < \normonebig{ \widetildebbeta}$. 
However, $\normtwo{\bX \bz - \by} = \normtwo{\widetildeby - \by} = \normtwobig{\bX \widetildebbeta - \by}$ and $\normone{ \bz} < \normonebig{ \widetildebbeta}$ mean that $\bz$ is a strictly better solution to problem \eqref{opt:p1_penalize} than $\widetildebbeta$, contradicting the assumption $\widetildebbeta \in \sQ_\sigma$.
Thus, $\sQ_\sigma = \sQ_q$.

Since $\sQ_\sigma = \sQ_q$, $\widetildebbeta$ is the unique solution to problem \eqref{opt:p1_penalize} if and only if it is the unique solution to problem \eqref{equation:exact_p1_pen_epsi_pv1}. Since problem \eqref{equation:exact_p1_pen_epsi_pv1} is in the same form of problem \eqref{opt:p1}, applying the part of Theorem~\ref{theorem:dualcert_p1} for problem \eqref{opt:p1},  we conclude that $\widetildebbeta$ is the unique solution to problem \eqref{opt:p1_penalize} if and only if the same condition as Theorem~\ref{theorem:dualcert_p1} holds. 
\end{proof}

\begin{theoremHigh}[Error bound of \eqref{opt:p1_epsilon} under measurement noise]\label{theorem:errbd_p1epsi}
Let $\bbeta^*$ be the original signal satisfying $\by = \bX\bbeta^* + \bepsilon$ with $\normtwo{\bepsilon} \leq \epsilon$. And let $\widehat{\bbeta}$ be a solution to the problem~\eqref{opt:p1_epsilon}:
$$
\text{(P$_{1,\epsilon}$)}\qquad
\min \normone{\bbeta} \quad \text{s.t.} \quad \normtwo{\bX\bbeta - \by} \leq \epsilon.
$$
Define the support of the true signal as $\sS \triangleq \supp(\bbeta^*)$ and its complement as $\sT \triangleq \comple{\sS}$.
Then the reconstruction error satisfies
$$
\normonebig{\widehat{\bbeta} - \bbeta^*} \leq C_3 \epsilon + C_4 \normonebig{\widehat{\bbeta}_\sT},
$$
where $C_3 \triangleq 2\sqrt{\abs{\sS}} \cdot \mathcalU(\sS)$, $C_4 \triangleq \normtwo{\bX} \sqrt{\abs{\sS}} \cdot \mathcalU(\sS)+1$, and 
$$
\mathcalU(\sS)
\triangleq \max_{\substack{\supp(\bu) = \sS, \\ \bu \neq \bzero}} \frac{\normtwo{\bu}}{\normtwo{\bX\bu}}.
~\footnote{$\mathcalU(\sS)$ is related to one side of the RIP bound; see Definition~\ref{definition:rip22}.} 
$$

\end{theoremHigh}

\begin{proof}[of Theorem~\ref{theorem:errbd_p1epsi}]
Observe that $\normonebig{\widehat{\bbeta} - \bbeta^*} = \normonebig{\widehat{\bbeta}_\sS - \bbeta^*_\sS} + \normonebig{\widehat{\bbeta}_\sT}$, where the first term satisfies
$$
\normone{\widehat{\bbeta}_\sS - \bbeta^*_\sS} 
\leq \sqrt{\abs{\sS}} \cdot \normtwo{\widehat{\bbeta}_\sS - \bbeta^*_\sS} 
\leq \sqrt{\abs{\sS}} \cdot \mathcalU(\sS) \cdot \normtwo{\bX_\sS (\widehat{\bbeta}_\sS - \bbeta^*_\sS)},
$$
where the first inequality follows from standard bounds on vector norms (Exercise~\ref{exercise:cauch_sc_l1l2}), and the last inequality follows from the definition of $\mathcalU(\sS)$.
Let $\bd \in\real^p$ such that $\bd_\sS = \widehat{\bbeta}_\sS$ and $\bd_\sT=\bzero$, that is $\supp(\bd)=\sS$. 
Therefore, by the triangle inequality,
$$
\normtwo{\bX_\sS (\widehat{\bbeta}_\sS - \bbeta^*_\sS)} 
= \normtwo{\bX(\bd - \bbeta^*)} 
\leq \normtwo{\bX(\bd - \widehat{\bbeta})} + {\normtwobig{\bX(\widehat{\bbeta} - \bbeta^*)}}.
$$
Since both $\widehat{\bbeta}$ and $\bbeta^*$ are feasible, we have ${\normtwobig{\bX(\widehat{\bbeta} - \bbeta^*)}}\leq 2\epsilon$.
The definition of the spectral norm and standard bounds on vector norms (Exercise~\ref{exercise:cauch_sc_l1l2}) also show that 
$$
\normtwo{\bX(\bd - \widehat{\bbeta})} \leq \normtwo{\bX} \normtwo{\bd - \widehat{\bbeta}} \leq \normtwo{\bX} \normone{\bd - \widehat{\bbeta}} = \normtwo{\bX} \normonebig{\widehat{\bbeta}_\sT}
$$
Combining the results shows 
$$
\normone{\widehat{\bbeta}_\sS - \bbeta^*_\sS} 
\leq 
\sqrt{\abs{\sS}} \cdot \mathcalU(\sS) \cdot (\normtwo{\bX} \normonebig{\widehat{\bbeta}_\sT} + 2\epsilon).
$$
This concludes the result after rearrangement.
\end{proof}

Suppose $\bbeta^*$ and $\balpha$ satisfy the same  condition as Theorem~\ref{theorem:dualcert_p1}. 
The vector then $\balpha=\bX^\top\blambda$ satisfies 
$$\balpha_\sS = \sgn(\bbeta_\sS^*)
\qquad\text{and}\qquad 
\norminf{\balpha_\sT} < 1,$$ 
where $\sS=\supp(\bbeta^*)$ and $\sT=\comple{\sS}$. 
Let $\widehatbbeta$ be a solution of \eqref{opt:p1}. Then
\begin{align*}
\normonebig{\widehat{\bbeta}_\sS} 
&\geq \innerproduct{\balpha_\sS, \widehat{\bbeta}_\sS};\\
\normonebig{\widehat{\bbeta}_\sT} 
&\leq (1 - \norminf{\balpha_\sT})^{-1} \left(\normonebig{\widehat{\bbeta}_\sT} 
- \innerproduct{\balpha_\sT, \widehat{\bbeta}_\sT}\right),
\end{align*}
where the second inequality follows from the triangle inequality 
$\normonebig{\widehat{\bbeta}_\sT} - \normonebig{\widehat{\bbeta}_\sT}\norminf{\balpha_\sT}
\leq  \big(\normonebig{\widehat{\bbeta}_\sT} - \innerproductbig{\balpha_\sT, \widehat{\bbeta}_\sT}\big)$.
Therefore,
$$
\normonebig{\widehat{\bbeta}_\sT} 
\leq (1 - \norminf{\balpha_\sT})^{-1} 
\left(\normonebig{\widehat{\bbeta}} - \innerproduct{\balpha, \widehat{\bbeta}}\right) 
\triangleq (1 - \norminf{\balpha_\sT})^{-1} D_1^{\balpha}(\widehat{\bbeta}, \bbeta^*),
$$
where
$D_1^{\balpha}(\widehat{\bbeta}, \bbeta^*) 
\triangleq \normonebig{\widehat{\bbeta}} - \normone{\bbeta^*} - \innerproductbig{\balpha, \widehat{\bbeta} - \bbeta^*}$
is the {Bregman distance} associated with $\normone{\cdot}$ (see Definition~\ref{definition:breg_dist}), and $\balpha\in\partial\normone{\bbeta^*}$.

Since $\balpha \in \cspace(\bX^\top)$ so $\balpha = \bX^\top \blambda$ for some vector $\blambda$ (i.e., the dual vector), 
using the optimality of $\widehat{\bbeta}$ shows 
\begin{align*}
D_1^{\balpha}(\widehat{\bbeta}, \bbeta^*) 
&\triangleq \normonebig{\widehat{\bbeta}} - \normone{\bbeta^*} - 
\innerproduct{\balpha, \widehat{\bbeta} - \bbeta^*}\\
&\leq 0+  \innerproduct{\bX^\top\blambda,  \bbeta^* - \widehat{\bbeta}}
\leq \normtwo{\blambda}\normtwo{\bX(\bbeta^* - \widehat{\bbeta})}
\leq 2\normtwo{\blambda} \epsilon.
\end{align*}
Therefore, we can bound the optimal solution in the complement set by
\begin{equation}
\normonebig{\widehat{\bbeta}_\sT} \leq 2(1 - \norminf{\balpha_\sT})^{-1} \normtwo{\blambda} \epsilon.
\end{equation}

The above analysis shows the following corollary from Theorem~\ref{theorem:errbd_p1epsi}
\begin{corollary}[Error bound of \eqref{opt:p1_epsilon}  under measurement noise]\label{corollary:errbd_p1epsi}
Suppose $\bbeta^*$ and $\balpha= \bX^\top \blambda$ satisfy the same  condition as Theorem~\ref{theorem:dualcert_p1}, and let $\bbeta^*$ be the original signal satisfying $\by = \bX\bbeta^* + \bepsilon$ with $\normtwo{\bepsilon} \leq \epsilon$. 
Let $\widehat{\bbeta}$ be a solution to the problem~\eqref{opt:p1_epsilon}:
$$
\text{(P$_{1,\epsilon}$)}\qquad
\min \normone{\bbeta} \quad \text{s.t.} \quad \normtwo{\bX\bbeta - \by} \leq \epsilon.
$$
Define $\sS \triangleq \supp(\bbeta^*)$ and $\sT \triangleq \comple{\sS}$.
Then,
$$
\normonebig{\widehat{\bbeta} - \bbeta^*} \leq C\epsilon,
$$
where
$$
C \triangleq 2\sqrt{\abs{\sS}} \cdot \mathcalU(\sS) 
+ \frac{2\normtwo{\blambda}(\normtwo{\bX} \sqrt{\abs{\sS}} \cdot \mathcalU(\sS)+1)}{1 - \norminf{\balpha_\sT}}.
$$
and
$$
\mathcalU(\sS)
\triangleq \max_{\substack{\supp(\bu) = \sS, \\ \bu \neq \bzero}} \frac{\normtwo{\bu}}{\normtwo{\bX\bu}}.
$$
\end{corollary}

Note that the condition $\norminf{\balpha_\sT} <1 $ is essential: as $\norminf{\balpha_\sT}\rightarrow 1$, the constant diverges, reflecting potential instability in recovery when the dual certificate is weak on the complement of the support.

A similar error bound can be derived for the penalized formulation~\eqref{opt:p1_penalize}:
$$
\min  \normone{\bbeta} + \sigma \normtwo{\bX\bbeta - \by}^2
$$ 
provided  $\sigma>0$ is chosen appropriately. This follows from the equivalence between problems~\eqref{opt:p1_epsilon} and~\eqref{opt:p1_penalize} established in Corollary~\ref{corollary:equigv_p1_epsi_pen_lag}; we omit the details here.

\section{Optimality Condition of $\ell_1$-Analysis Models*}

\newcommand{\tbbeta}{{\widebarbbeta}}

We now consider the following $\ell_1$-analysis formulations:
\begin{subequations}\label{equation:q1_analysis_opcd}
\begin{align}
(Q_{1})\qquad &\min_{\bbeta\in\real^p} \normone{\bPsi^\top \bbeta}, \quad \text{s.t.}\quad \bX \bbeta = \by; 
\label{equation:q1_opcd} \\
(Q_{1,\sigma}) \qquad &\min_{\bbeta\in\real^p}  \normone{\bPsi^\top \bbeta}  + \sigma \normtwo{\bX \bbeta - \by}^2; 
\label{equation:q1_penalize_opcd}  \\
(Q_{1,\epsilon})\qquad &\min_{\bbeta\in\real^p} \normone{\bPsi^\top \bbeta}, \quad \text{s.t.}\quad \normtwo{\bX\bbeta-\by} \leq \epsilon, 
\label{equation:q1_epsilon_opcd} 
\end{align}
\end{subequations}
where $\epsilon, \sigma$ are positive parameters,  $\bX\in\real^{n\times p}$, $\bbeta\in\real^p$, and  $\bPsi\in\real^{p\times q}$. 
Problem~\eqref{equation:q1_opcd} is appropriate when the observations are noise-free, i.e., $\bX\bbeta-\by=\bepsilon = \bzero$. If $\bPsi = \bI$, the identity matrix, these models reduce to the standard  $\ell_1$-synthesis models (i.e., \eqref{opt:p1}, p.~\pageref{opt:p1}, \eqref{opt:p1_penalize}, p.~\pageref{opt:p1_penalize}, and \eqref{opt:p1_epsilon}, p.~\pageref{opt:p1_epsilon}). 
We focus on the high-dimensional regime $n\leq p$ (see Lemma~\ref{lemma:errbd_noisem_q1ep2} or Theorem~\ref{theorem:robst_optq} for justification).

For certificates of these models, we consider the following condition.
\begin{condition}\label{condition:q1_cond1}
Let $\bPsi\in\real^{q\times p}$ and $\bX\in\real^{n\times p}$.
Given $\tbbeta \in \real^p$, index set $\sS = \supp(\bPsi^\top\tbbeta) \subset \{1,2,\ldots,q\}$  and $\sT\triangleq \comple{\sS}$ satisfy
\begin{enumerate}[(i)]
\item $\nspace(\bPsi(\sT)^\top) \cap \nspace(\bX) = \{\bzero\}$;
\item There exists $\balpha \in \real^q$ such that $\bPsi \balpha \in \cspace(\bX^\top)$, $\balpha_\sS = \sgn(\bPsi_\sS^\top\tbbeta)$, and $\norminf{\balpha_{\sT}} < 1$.
\end{enumerate}
\end{condition}
This condition is a natural extension of the dual certificate condition in Theorem~\ref{theorem:exact_q1} for the classical $\ell_1$-minimization problem \eqref{opt:p1}.
Note again that $\bPsi(\sT)\in\real^{q\times p}$  and $\bPsi_{\sT}=\bPsi[:,\sT]\in\real^{n\times \abs{\sT}}$; see Definition~\ref{definition:matlabnotation}.
When $\bPsi = \bI_p$, $\nspace(\bPsi(\sT)^\top) = \{\bbeta\in\real^p \mid\beta_i=0 \text{ for all } i\in\sT\}$. Then Condition~\ref{condition:q1_cond1} (i) reduces to requiring that $\bX_{\sS}$ has full column rank, exactly as in Theorem~\ref{theorem:exact_q1}.

More generally, Condition~\ref{condition:q1_cond1} (i) ensures that there is no nonzero perturbation $\Delta \bbeta$ satisfying both $\bPsi(\sT)^\top\tbbeta = \bPsi(\sT)^\top(\tbbeta + \Delta \bbeta)$ and $\bX\tbbeta = \bX(\tbbeta + \Delta \bbeta)$. Otherwise, there exists a nonempty interval $\sM = [\tbbeta - \gamma\Delta \bbeta, \tbbeta + \gamma\Delta \bbeta]$ for some sufficiently small $\gamma > 0$ so that $\bX \bbeta = \bX\tbbeta$ and $\normone{\bPsi^\top \bbeta}$ is constant for $\bbeta \in \sM$; 
consequently, $\tbbeta$ could not be the unique minimizer. 

Condition~\ref{condition:q1_cond1} (ii) asserts  the existence of a strictly complementary \textit{dual certificate} $\balpha$. 
To see this, recall the optimality condition for problem~\eqref{opt:q1}: $\bzero \in \bPsi \partial\normone{\bPsi^\top \bbeta} - \bX^\top\blambda$, where  $\blambda$ is the vector of Lagrangian multipliers (see Remark~\ref{remark:kkt_nutshell_cvx}).
In the synthesis case ($\bPsi=\bI$),  $\bX^\top\blambda$ is a subgradient of the $\ell_1$-norm at $\bbeta$ (see Theorem~\ref{theorem:dualcert_p1}); however, $\bX^\top\blambda$ is no longer a subgradient in this $\ell_1$-analysis case due to the transformation of $\bPsi$.
We can rewrite the condition as $\bzero = \bPsi \balpha - \bX^\top\blambda$, where $\balpha \in \partial\normone{\bPsi^\top \bbeta}$, which translates to $\balpha_\sS = \sgn(\bPsi_\sS^\top\bbeta)$ and $\norminf{\balpha_\sT} \leq 1$ by the subdifferential of the $\ell_1$-norm (Exercise~\ref{exercise:sub_norms}). 
The vector $\balpha$ thus certifies the optimality of $\bbeta=\tbbeta$. 
For uniqueness and robustness, we will later show that the strict inequality $\norminf{\balpha_\sS} < 1$ is essential.

Taken together, these observations confirm that Condition~\ref{condition:q1_cond1} is indeed a natural generalization of the dual certificate condition in Theorem~\ref{theorem:dualcert_p1} to the $\ell_1$-analysis framework.
We also introduce a slightly more general version:
\begin{condition}\label{condition:q1_cond2}
Let $\bPsi\in\real^{q\times p}$ and $\bX\in\real^{n\times p}$.
Given $\tbbeta \in \real^p$, let $\sS = \supp(\bPsi^\top\tbbeta)$. There exists a nonempty index set $\sT \textcolor{mydarkblue}{\subseteq} \comple{\sS}$ such that the index sets $\sS$, $\sT$ and $\sI \triangleq \comple{(\sS \cup \sT)}$ satisfy
\begin{enumerate}[(a)]
\item $\nspace(\bPsi(\sT)^\top) \cap \nspace(\bX) = \{\bzero\}$;
\item There exists $\balpha \in \real^q$ such that $\bPsi \balpha \in \cspace(\bX^\top)$, $\balpha_\sS = \sgn(\bPsi_\sS^\top\tbbeta)$, $\norminf{\balpha_\sT} < 1$, and $\norminf{\balpha_\sI} \leq 1$.
\end{enumerate}
\end{condition}

In Condition~\ref{condition:q1_cond2}, choosing a smaller $\sT$ relaxes requirement  (b) but enlarges the nullspace $\nspace(\bPsi(\sT)^\top)$, thereby tightening requirement  (a). 
Although this formulation allows greater flexibility in selecting $\sT$, we will show that Conditions~\ref{condition:q1_cond1} and~\ref{condition:q1_cond2} are in fact equivalent.

We also impose the following standard assumptions:

\begin{assumption}
Consider the problems in \eqref{equation:q1_analysis_opcd}, we assume that 
\begin{enumerate}
\item Matrix $\bX$ has full row rank.
\item Matrix $\bPsi$ has full row rank.
\item Maximum singular value $\sigma_{\max}(\bPsi) = 1$.
\end{enumerate}
\end{assumption}
Assumptions~1 and~2 prevent redundant equations or transform coefficients and are commonly adopted in compressed sensing and inverse problems. Assumption 3 is non-essential: any matrix $\bPsi$ can be normalized by dividing by $\sigma_{\max}(\bPsi)$, which simply rescales the $\ell_1$-penalty and leaves the solution set unchanged.

\subsection{Unique Recovery in Noise-Free Settings}

We begin by establishing a necessary and sufficient condition for the existence of a solution to problem~\eqref{opt:q1}.
\begin{theoremHigh}[Necessary and sufficient condition for \eqref{opt:q1}]\label{theorem:nec_suf_q1}
Consider the optimization problem~\eqref{opt:q1}.
A vector $\widehatbbeta$ is a solution of \eqref{opt:q1} if and only if, for  all $\bn \in \nspace(\bX)$, 
$$
\innerproduct{\bPsi_\sS^\top \bn, \sgn(\bPsi_\sS^\top \widehatbbeta)} 
\leq  \normonebig{\bPsi_{\comple{\sS}}^\top \bn},
$$
where $\sS \triangleq\supp(\bPsi^\top\widehatbbeta) $.
\end{theoremHigh}
\begin{proof}[of Theorem~\ref{theorem:nec_suf_q1}]
\textbf{Necessity.} 
Suppose $\widehatbbeta$ is a solution of \eqref{opt:q1}. 
Since the function $f(\eta) = \bX(\widehatbbeta + \eta\bn)$ is continuously differentiable, we can find $\eta$ small enough around 0 such that $\sgn(\bPsi^\top(\widehatbbeta + \eta\bn)) = \sgn(\bPsi^\top \widehatbbeta)$.
For small and nonzero $\eta$, we have
$$
\footnotesize
\begin{aligned}
\normonebig{\bPsi^\top\widehatbbeta}
&\leq  \normonebig{\bPsi^\top (\widehatbbeta + \eta\bn)} 
= \normonebig{\bPsi_\sS^\top (\widehatbbeta + \eta\bn)} + \normonebig{\bPsi_{\comple{\sS}}^\top (\widehatbbeta + \eta\bn)} \\
&\leq  \innerproduct{\bPsi_\sS^\top (\widehatbbeta + \eta \bn), \sgn(\bPsi_\sS^\top (\widehatbbeta + \eta\bn))} 
+ \normonebig{\bPsi_{\comple{\sS}}^\top \widehatbbeta} 
+ \normonebig{\eta\bPsi_{\comple{\sS}}^\top \bn} \\
&= \innerproduct{\bPsi_\sS^\top \widehatbbeta + \eta\bPsi_\sS^\top \bn, \sgn(\bPsi_\sS^\top \widehatbbeta)} 
+ \normonebig{\bPsi_{\comple{\sS}}^\top \widehatbbeta} 
+ \normonebig{\eta\bPsi_{\comple{\sS}}^\top \bn} \\
&= \innerproduct{\bPsi_\sS^\top \widehatbbeta, \sgn(\bPsi_\sS^\top \widehatbbeta)} + \eta\innerproduct{\bPsi_\sS^\top \bn, \sgn(\bPsi_\sS^\top \widehatbbeta)} 
+ \normonebig{\bPsi_{\comple{\sS}}^\top \widehatbbeta} 
+ \normonebig{\eta\bPsi_{\comple{\sS}}^\top \bn} \\
&= \normonebig{\bPsi^\top \widehatbbeta} + \eta \innerproduct{\bPsi_\sS^\top \bn, \sgn(\bPsi_\sS^\top \widehatbbeta)} +\normonebig{\eta\bPsi_{\comple{\sS}}^\top \bn}.
\end{aligned}
$$
Therefore, for any $\bn \in \nspace(\bX)$, we have 
$\innerproductbig{\bPsi_\sS^\top \bn, \sgn(\bPsi_\sS^\top \widehatbbeta)} 
\leq  \normonebig{\bPsi_{\comple{\sS}}^\top \bn}$.

\paragraph{Sufficiency.}
Suppose the stated inequality holds for all $\bn\in\nspace(\bX)$.
Assume, for contradiction, that there exists a feasible point $\bbeta'$ such that $\bX\bbeta' = \by$ and $\normone{\bPsi^\top\bbeta'} < \normonebig{\bPsi^\top\widehatbbeta}$.
Define $\bn \triangleq \widehatbbeta - \bbeta'$, which satisfies $\bX\bn = \bzero$ since $\bX\bbeta' = \bX\widehatbbeta$. This means $\bn \in \nspace(\bX)$.
Using the identity $\normonebig{\bbeta'} = \normonebig{\widehatbbeta - \bn}$, we consider the decomposition of the objective:
\begin{align*}
\normonebig{\bPsi^\top\bbeta'}
&=\normone{\bPsi^\top (\widehatbbeta  - \bn)} 
= \sum_{i \in \sS} \abs{ \bPsi^\top_i (\widehatbeta_i - n_i)} + \sum_{i \in \comple{\sS}} \abs{\bPsi^\top_i n_i}\\
&\geq \sum_{i \in \sS} \abs{\bPsi^\top_i\widehatbeta_i} - \sum_{i \in \sS} \sgn(\bPsi^\top_i\widehatbeta_i)\bPsi^\top_in_i + \sum_{i \in \comple{\sS}} \abs{\bPsi^\top_in_i}
\geq \sum_{i \in \sS} \abs{\bPsi^\top_i\widehatbeta_i} =\normonebig{\bPsi^\top\widehatbbeta},
\end{align*}
where the first inequality follows from the fact that $\absbig{b - n} \geq \absbig{b} - \sgn(b)n$,
and the second inequality follows from the given condition that  $\innerproductbig{\bPsi_\sS^\top \bn, \sgn(\bPsi_\sS^\top \widehatbbeta)} 
\leq  \normonebig{\bPsi_{\comple{\sS}}^\top \bn}$.
Since the assumption was $\normone{\bPsi^\top\bbeta'} < \normonebig{\bPsi^\top\widehatbbeta}$, we obtain a contradiction.
Hence, no such $\bbeta'$ exists, and  $\widehatbbeta$ is indeed optimal.
This completes the proof.
\end{proof}

The following theorem states that Conditions~\ref{condition:q1_cond1} and~\ref{condition:q1_cond2} are equivalent, and they are necessary and sufficient for the
unique solution $\widehatbbeta$ to problem \eqref{opt:q1}.
\begin{theoremHigh}[Dual certificate for \eqref{opt:q1} \citep{zhang2016one}]\label{theorem:exact_q1}
Under Assumption~1, let $\widehatbbeta$ be a solution to problem \eqref{opt:q1}. The following are equivalent:
\begin{enumerate}[(i)]
\item Solution $\widehatbbeta$ is unique.
\item Condition~\ref{condition:q1_cond1} holds for $\tbbeta = \widehatbbeta$.
\item Condition~\ref{condition:q1_cond2} holds for $\tbbeta = \widehatbbeta$.
\end{enumerate}
\end{theoremHigh}

\begin{proof}[of Theorem~\ref{theorem:exact_q1}]
This proof adapts arguments from  \citet{grasmair2011necessary, zhang2015necessary, zhang2016one}.
We establish the equivalence by proving the implications (iii) $\Rightarrow$ (i) $\Rightarrow$ (ii) $\Rightarrow$ (iii).

\paragraph{(iii) $\Rightarrow$ (i).} Consider any perturbation $\widehatbbeta + \bn$ where $\bn \in \nspace(\bX)\setminus\{\bzero\}$. 
Since by Condition~\ref{condition:q1_cond2},  $\balpha_\sS = \sgn(\bPsi_\sS^\top\widehatbbeta)$, $\norminf{\balpha_\sT} < 1$, and $\norminf{\balpha_\sI} \leq 1$, 
we can take a {subgradient} $\bg \in \partial\normonebig{\bPsi^\top \widehatbbeta}$ satisfying 
$\bg_\sS  = \balpha_\sS \equiv  \sgn(\bPsi_\sS^\top \widehatbbeta)$, 
$\bg_\sI = \balpha_\sI $ ($\norminf{\bg_\sI} =\norminf{\balpha_\sI}\leq 1$), 
and $\norminf{\bg_\sT} \leq 1$
such that $\innerproduct{\bg_\sT, \bPsi_\sT^\top \bn} = \normonebig{\bPsi_\sT^\top \bn}$  (which is similar to $\balpha$ while we do not require strict inequality here; actually $\norminf{\bg_\sT}=1$ is required to satisfy the latter equality). 
Then the subgradient inequality (Definition~\ref{definition:subgrad}) shows
\begin{align*}
\normonebig{\bPsi^\top (\widehatbbeta + \bn)}
&\geq \normonebig{\bPsi^\top \widehatbbeta} + \innerproduct{\bPsi \bg, \bn}
\stackrel{\dag}{=} \normonebig{\bPsi^\top \widehatbbeta} + \innerproduct{\bPsi \bg - \bPsi \balpha, \bn} \\
&= \normonebig{\bPsi^\top \widehatbbeta} + \innerproduct{\bg - \balpha, \bPsi^\top \bn}  
\stackrel{\ddag}{=} \normonebig{\bPsi^\top \widehatbbeta} + \innerproduct{\bg_\sT - \balpha_\sT, \bPsi_\sT^\top \bn} \\
&\geq \normonebig{\bPsi^\top \widehatbbeta} + \normonebig{\bPsi_\sT^\top \bn} (1 - \norminf{\balpha_\sT}), 
\end{align*}
where the equality $(\dag)$ follows from $\bPsi \balpha \in \cspace(\bX^\top) = \nspace(\bX)^\perp$ and $\bn \in \nspace(\bX)$, the equality ($\ddag$) follows from the definition of $\bg$, and the last inequality is an application of the \holders inequality $\innerproduct{\balpha_\sT,\bPsi_\sT^\top \bn} \leq \norminf{\balpha_\sT}\normone{\bPsi_\sT^\top \bn}$ and $\innerproduct{\bg_\sT, \bPsi_\sT^\top \bn} = \normonebig{\bPsi_\sT^\top \bn}$. 
Since $\bn \in \nspace(\bX)\setminus\{\bzero\}$ and $\nspace(\bPsi(\sT)^\top) \cap \nspace(\bX) = \{\bzero\}$, we have $\normonebig{\bPsi_\sT^\top \bn} > 0$. Together with the condition $\norminf{\balpha_\sT} < 1$, we have $\normonebig{\bPsi^\top (\widehatbbeta + \bn)}> \normonebig{\bPsi^\top \widehatbbeta}$ for every $\bn \in \nspace(\bX)\setminus\{\bzero\}$. 
This proves that $\widehatbbeta$ is the unique minimizer of problem \eqref{opt:q1}.

\paragraph{(i) $\Rightarrow$ (ii).} For every $\bn \in \nspace(\bX)\setminus\{\bzero\}$, we have $\bX(\widehatbbeta + \eta\bn) = \bX\widehatbbeta$.
Since the function $f(\eta) = \bX(\widehatbbeta + \eta\bn)$ is continuously differentiable, we can find $\eta$ small enough around 0 such that $\sgn(\bPsi^\top(\widehatbbeta + \eta\bn)) = \sgn(\bPsi^\top \widehatbbeta)$. Since $\widehatbbeta$ is the unique solution, for small and nonzero $\eta$, the proof of Theorem~\ref{theorem:nec_suf_q1} shows the following strict inequality
\begin{align*}
\normonebig{\bPsi^\top\widehatbbeta}
< \normonebig{\bPsi^\top \widehatbbeta} 
+ \eta \innerproduct{\bPsi_\sS^\top \bn,\sgn(\bPsi_\sS^\top \widehatbbeta)} +\normonebig{\eta\bPsi_{\comple{\sS}}^\top \bn}.
\end{align*}
Therefore, for any $\bn \in \nspace(\bX)\setminus\{\bzero\}$, we have
\begin{equation}\label{equ:exact_q1_pv1}
\innerproduct{\bPsi_\sS^\top \bn, \sgn(\bPsi_\sS^\top \widehatbbeta)} 
< \normonebig{\bPsi_{\comple{\sS}}^\top \bn}. 
\end{equation}
If the condition $\nspace(\bPsi(\comple{\sS})^\top) \cap \nspace(\bX) = \{\bzero\}$ does not hold, we can choose a nonzero vector $\bn \in \nspace(\bPsi(\comple{\sS})^\top) \cap \nspace(\bX)$. 
We also have $-\bn \in \nspace(\bPsi(\comple{\sS})^\top) \cap \nspace(\bX)$. 
By \eqref{equ:exact_q1_pv1}, we have $\innerproductbig{\bPsi_\sS^\top \bn, \sgn(\bPsi_\sS^\top \widehatbbeta)} < 0$ and $-\innerproductbig{\bPsi_\sS^\top \bn, \sgn(\bPsi_\sS^\top \widehatbbeta)} < 0$, which leads to a contradiction, and shows that $\nspace(\bPsi(\comple{\sS})^\top) \cap \nspace(\bX) = \{\bzero\}$, i.e., Condition~\ref{condition:q1_cond1} (i).

Next, we construct $\balpha$ satisfying  Condition~\ref{condition:q1_cond1} (ii).
Define  $\widetildebalpha\in\real^q$ by $\widetildebalpha_\sS = \sgn(\bPsi_\sS^\top \widehatbbeta)$ and $\widetildebalpha_{\comple{\sS}} = \bzero$. 
\begin{itemize}
\item If such $\widetildebalpha$ satisfies $\bPsi \widetildebalpha \in \cspace(\bX^\top)$, we are done.
\item If $\bPsi \widetildebalpha \notin \cspace(\bX^\top) \equiv \nspace(\bX)^\perp$, then we will construct a new vector to satisfy Condition~\ref{condition:q1_cond1} (ii). 
\end{itemize}
Let $\bQ\in\real^{p\times r}$ be a matrix whose columns form a basis for $\nspace(\bX)$. 
It follows that $\bgamma \triangleq \bQ^\top \bPsi \widetildebalpha\in\real^r$ must be a nonzero vector. 
Consider the optimization problem
\begin{equation}\label{equ:exact_q1_pv2}
\min_{\bu \in \real^q} \norminf{\bu} \quad \text{s.t.} \quad \bQ^\top \bPsi \bu = -\bgamma \text{ and } \bu_\sS = \bzero. 
\end{equation}
For any minimizer $\widehatbu$ of problem \eqref{equ:exact_q1_pv2}, we have $\bQ^\top \bPsi(\widetildebalpha + \widehatbu) = \bzero$ such that  $\bPsi(\widetildebalpha + \widehatbu)$ is orthogonal to the nullspace $\nspace(\bX)$, 
which implies $\bPsi(\widetildebalpha + \widehatbu) \in \nspace(\bX)^\perp \equiv \cspace(\bX^\top)$ and $(\widetildebalpha + \widehatbu)_\sS = \widetildebalpha_\sS = \sgn(\bPsi_\sS^\top \widehatbbeta)$. 
Thus, we shall show that the objective of problem \eqref{equ:exact_q1_pv2} is strictly less than 1 using  strong duality, thus completing for the proof of the part (ii) of Condition~\ref{condition:q1_cond1}. 
To this end, we rewrite problem \eqref{equ:exact_q1_pv2} in an equivalent form as:
\begin{equation}\label{equ:exact_q1_pv3}
\min_{\bu\in\real^q} \norminf{\bu_{\comple{\sS}}} \quad \text{s.t.} \quad \bQ^\top \bPsi_{\comple{\sS}} \bu_{\comple{\sS}} = -\bgamma, 
\end{equation}
whose Lagrange dual problem (i.e., the dual of an $\ell_\infty$-minimization problem; see Example~\ref{example:linf_min_dual}) is
\begin{equation}\label{equ:exact_q1_pv4}
\max_{\blambda\in\real^r} \innerproduct{\blambda, \bgamma} \quad \text{s.t.} \quad \normonebig{\bPsi_{\comple{\sS}}^\top \bQ \blambda} \leq 1. 
\end{equation}
Note that $\bQ \blambda \in \nspace(\bX)$, and 
$\abs{\innerproduct{\blambda, \bgamma}} = \absbig{\innerproductbig{\blambda, \bQ^\top \bPsi \widetildebalpha}} 
= \absbig{\innerproductbig{\bPsi_\sS^\top \bQ \blambda, \sgn(\bPsi_\sS^\top \widehatbbeta)}}$. 
Using \eqref{equ:exact_q1_pv1}, for any $\blambda$, we have
$$
\abs{\innerproduct{ \blambda, \bgamma}} =
\begin{cases}
\abs{\innerproduct{\bPsi_\sS^\top \bQ \blambda, \sgn(\bPsi_\sS^\top \widehatbbeta)}} = 0, & \text{if } \bQ \blambda = \bzero; \\
\abs{\innerproduct{\bPsi_\sS^\top \bQ \blambda, \sgn(\bPsi_\sS^\top \widehatbbeta)}} < \normone{\bPsi_{\comple{\sS}}^\top \bQ \blambda} \leq 1, & \text{otherwise}.
\end{cases}
$$
Hence, problem \eqref{equ:exact_q1_pv4} is feasible, and its objective value is strictly less than 1. By the strong duality property (which holds because the primal is feasible and the objective is convex and continuous), problems \eqref{equ:exact_q1_pv2} and \eqref{equ:exact_q1_pv3} attain their minimum, and the optimal value is also strictly less than 1. 
Let $\widehatbu$ be an optimal solution. Define $\balpha = \widetilde{\balpha} + \widehatbu$. Then:
\begin{itemize}
\item $\balpha_\sS = \widetildebalpha_\sS = \sign(\bPsi_\sS^\top \widehatbbeta)$,
\item $\norminf{\balpha_{\comple{\sS}}} = \norminf{\widehatbu_{\comple{\sS}}} < 1$,
\item $\bQ^\top \bPsi \balpha = \bQ^\top \bPsi (\widetilde{\balpha} + \widehatbu) = \bgamma - \bgamma = \bzero$, so $\bPsi \balpha \in \nspace(\bX)^\perp \equiv \cspace(\bX^\top)$.
\end{itemize}
Hence, $\balpha$ satisfies Condition~\ref{condition:q1_cond1} (ii).

\paragraph{(ii) $\Rightarrow$ (iii).} Let $\sT \triangleq \comple{\sS}$ and $\sI \triangleq \varnothing$.
Then Condition~\ref{condition:q1_cond1} directly implies Condition~\ref{condition:q1_cond2}, completing the cycle of implications.
\end{proof}

\subsection{Unique Recovery in Noisy Signal}

Similarly to Theorem~\ref{theorem:exact_q1}, in the presence of noise, Conditions~\ref{condition:q1_cond1} and~\ref{condition:q1_cond2} remain equivalent and provide necessary and sufficient conditions for the uniqueness of the solution $\widehatbbeta$ to either problem~\eqref{opt:q1_penalize} or problem~\eqref{opt:q1_epsilon}.
To establish this, we require the following lemma. Its proof follows analogously from that of Lemma~\ref{lemma:const_conv_p1sigma} and is omitted.
\begin{lemma}[Constants over the set of solutions in \eqref{opt:q1_penalize} and \eqref{opt:q1_epsilon}]\label{lemma:const_conv_q1sigma}
Let $\sigma > 0$. If the function $\normone{\bPsi^\top \bbeta} + \sigma \normtwo{\bX\bbeta-\by}^2 $ is constant on a convex set $\sC$, then both $\bX\bbeta-\by$ and $\normone{\bPsi^\top \bbeta}$ are constants on $\sC$.

Let $\sQ_\sigma$ denote the set of solutions to problem~\eqref{opt:q1_penalize}.
Since $\normone{\bPsi^\top \bbeta} + \sigma \normtwo{\bX\bbeta-\by}^2$ is constant over $\sQ_\sigma$, the above result shows that $\bX\bbeta-\by$ and $\normone{\bPsi^\top \bbeta}$ are constants on $\sQ_\sigma$ in problem \eqref{opt:q1_penalize}.

Moreover,
suppose that $\bPsi$ has full row rank in problem \eqref{opt:q1_epsilon}. 
Then both $\normtwo{\bX\bbeta-\by}$ and $\normone{\bPsi^\top \bbeta}$ are constant on $\sQ_\epsilon$, where $\sQ_\epsilon$ denotes the set of solutions of problem \eqref{opt:q1_epsilon}
\end{lemma}

\begin{theoremHigh}[Dual certificate for \eqref{opt:q1_penalize} and \eqref{opt:q1_epsilon} \citep{zhang2016one}]\label{thm:exact_q1_pen_epsi}
Under Assumption~1, let $\widehatbbeta$ be a solution to problem  \eqref{opt:q1_penalize}, or under Assumptions~1 and~2, let $\widehatbbeta$ be a solution to problem \eqref{opt:q1_epsilon}. 
Then the following statements are equivalent:
\begin{enumerate}[(i)]
\item Solution $\widehatbbeta$ is unique.
\item Condition~\ref{condition:q1_cond1} holds for $\tbbeta = \widehatbbeta$.
\item Condition~\ref{condition:q1_cond2} holds for $\tbbeta = \widehatbbeta$.
\end{enumerate}
\end{theoremHigh}
\begin{proof}[of Theorem~\ref{thm:exact_q1_pen_epsi}]
The proof closely mirrors that of Theorem~\ref{theorem:dual_p1pen_epsi} and relies on Lemma~\ref{lemma:const_conv_q1sigma}.
Since Lemma~\ref{lemma:const_conv_q1sigma} yields identical structural properties for the solution sets of both problems~\eqref{opt:q1_penalize} and~\eqref{opt:q1_epsilon}, we present the argument only for problem~\eqref{opt:q1_penalize}; the case of~\eqref{opt:q1_epsilon} follows analogously.

Let $\sQ_\sigma$ denote the (nonempty) solution set of problem~\eqref{opt:q1_penalize}, and fix any $\widehatbbeta \in \sQ_\sigma$. 
Let $\widehatby  = \bX \widehatbbeta$, which is independent of the choice of $\widehatbbeta$ according to Lemma~\ref{lemma:const_conv_q1sigma}. We introduce the following constrained optimization problem:
\begin{equation}\label{equ:exact_q1_pen_epsi_pv1}
\min_{\bbeta\in\real^p} \normone{\bPsi^\top \bbeta} \quad \text{s.t.} \quad \bX \bbeta = \widehatby, 
\end{equation}
and let $\sQ_q$ denote its solution set.

We claim that $\sQ_\sigma = \sQ_q$. Since $\bX \bbeta = \bX \widehatbbeta$ and $\normone{\bPsi^\top \bbeta} = \normonebig{\bPsi^\top \widehatbbeta}$ for all $\bbeta \in \sQ_\sigma$ 
and conversely any $\bbeta$ satisfying $\bX \bbeta = \bX \widehatbbeta$ and $\normone{\bPsi^\top \bbeta} = \normonebig{\bPsi^\top \widehatbbeta}$ belongs to $\sQ_\sigma$, it suffices to show that $\normone{\bPsi^\top \bbeta} = \normonebig{\bPsi^\top \widehatbbeta}$ for any $\bbeta \in \sQ_q$. 
Assuming, for contradiction, this does not hold, then since problem \eqref{equ:exact_q1_pen_epsi_pv1} has $\widehatbbeta$ as a feasible solution and has a finite objective, we have a nonempty $\sQ_q$ and there exists $\bz \in \sQ_q$ satisfying $\normone{\bPsi^\top \bz} < \normonebig{\bPsi^\top \widehatbbeta}$. But, $\normtwo{\bX \bz - \by} = \normtwo{\widehatby - \by} = \normtwobig{\bX \widehatbbeta - \by}$ and $\normone{\bPsi^\top \bz} < \normonebig{\bPsi^\top \widehatbbeta}$ mean that $\bz$ is a strictly better solution to problem \eqref{opt:q1_penalize} than $\widehatbbeta$, contradicting the assumption that $\widehatbbeta \in \sQ_\sigma$.
Hence, the norms must be equal, and $\sQ_\sigma = \sQ_q$.

Since $\sQ_\sigma = \sQ_q$, $\widehatbbeta$ is the unique solution to problem \eqref{opt:q1_penalize} if and only if it is the unique solution to problem \eqref{equ:exact_q1_pen_epsi_pv1}. 
But problem~\eqref{equ:exact_q1_pen_epsi_pv1} is of the same form as problem~\eqref{opt:q1}. 
Applying Theorem~\ref{theorem:exact_q1}, we conclude that uniqueness holds if and only if either Condition~\ref{condition:q1_cond1} or Condition~\ref{condition:q1_cond2} is satisfied (with $\widetildebbeta=\widehatbbeta$).
This completes the proof.
\end{proof}

\subsection{Stable Recovery in Noisy Signal}
To state the robustness result, define
\begin{equation}
\mathcalU(\sT) \triangleq \max_{\bbeta \in \nspace(\bPsi(\sT)^\top)\setminus\{\bzero\}} \frac{\normtwo{\bbeta}}{\normtwo{\bX \bbeta}}.
\end{equation}
Condition~\ref{condition:q1_cond1} (i) guarantees  that $0 < \mathcalU(\sT) < +\infty$. 
If $\bPsi = \bI$, 
any nonzero $\bbeta \in \nspace(\bPsi(\sT)^\top)\setminus\{\bzero\}$ is a sparse (nonzero) vector supported on $\sT^c$;
consequently, $\mathcalU(\sT)$ equals the reciprocal of the smallest singular value of the submatrix $\bX_{\comple{\sT}}$ (see \eqref{equation:svd_stre_bd}).

\begin{lemma}\label{lemma:errbd_noisem_q1ep}
Assume that vectors $\tbbeta$ and $\balpha$ satisfy Condition~\ref{condition:q1_cond1}. Let $\sS \triangleq \supp(\bPsi^\top \tbbeta)$ and $\sT \triangleq \comple{\sS}$. 
Then,
\begin{equation}\label{equation:errbd_noisem_q1ep_eq1}
\normone{\bPsi^\top \bbeta - \bPsi^\top \tbbeta} \leq C_3 \normtwo{\bX(\bbeta - \tbbeta)} + C_4 D_1^\balpha(\bbeta, \tbbeta), \quad \forall\, \bbeta,
\end{equation}
where $D_1^\balpha(\bbeta, \tbbeta) \triangleq \normone{\bPsi^\top \bbeta} - \normone{\bPsi^\top \tbbeta} - \innerproduct{\bPsi \balpha, \bbeta - \tbbeta}$ is the Bregman distance associated with the function $\normone{\bPsi^\top \cdot}$, with $\balpha \in \partial \normone{\bPsi^\top \tbbeta}$ (see Definition~\ref{definition:breg_dist}), the absolute constants $C_3, C_4$ are given by 
$$
C_3 = \mathcalU(\sT)\sqrt{\abs{\sS}}, 
\quad \text{and} \quad C_4 = \frac{1 + \cond(\bPsi)\normtwo{\bX}C_3}{1 - \norminf{\balpha_\sT}},
$$
where $\cond(\bPsi)=\sigma_{\max}(\bPsi)/\sigma_{\min}(\bPsi)$ denotes the condition number of $\bPsi$.
\end{lemma}
\begin{proof}[of Lemma~\ref{lemma:errbd_noisem_q1ep}]
This proof adapts the argument in \citet{zhang2016one} and proceeds in two steps.
The first part shows that for any $\bu \in \nspace(\bPsi(\sT)^\top)$,
\begin{subequations}
\begin{equation}\label{equation:errbd_noisem_q1ep_eq2}
\normone{\bPsi^\top \bbeta - \bPsi^\top \tbbeta} \leq \left(1 + \frac{C_3 \normtwo{\bX}}{\sigma_{\min}(\bPsi)}\right) \normone{\bPsi^\top(\bbeta - \bu)}+ C_3 \normtwo{\bX(\bbeta - \tbbeta)}.
\end{equation}
The second part shows that
\begin{equation}\label{equation:errbd_noisem_q1ep_eq3}
f(\bbeta) \triangleq \min \left\{ \normone{\bPsi^\top(\bbeta - \bu)} \mid \; \bu \in \nspace(\bPsi(\sT)^\top) \right\} \leq (1 - \norminf{\balpha_\sT})^{-1} D_1^\balpha(\bbeta, \tbbeta). 
\end{equation}
\end{subequations}
Combining the two parts and using the definition of $C_4$ prove the desired result.

\paragraph{Proof of \eqref{equation:errbd_noisem_q1ep_eq2}.} 
Let $\bu \in \nspace(\bPsi(\sT)^\top)$, which also means $\bPsi_\sT^\top \bu = \bzero$ and $\supp(\bPsi^\top\bu) = \sS$. 
By the triangle inequality, we get
\begin{equation}\label{equation:errbd_noisem_q1ep_pv1}
\normone{\bPsi^\top \bbeta - \bPsi^\top \tbbeta} \leq \normone{\bPsi^\top(\bbeta - \bu)} + \normone{\bPsi^\top(\bu - \tbbeta)}. 
\end{equation}
Since $\supp(\bPsi^\top\tbbeta)=\sS$, implying $\tbbeta \in \nspace(\bPsi(\sT)^\top)$, we have $\bu - \tbbeta \in \nspace(\bPsi(\sT)^\top)$ and thus $\normtwo{\bu - \tbbeta} \leq \mathcalU(\sT)\normtwo{\bX(\bu - \tbbeta)}$, where $\mathcalU(\sT) < +\infty$ follows from  Condition~\ref{condition:q1_cond1} (i). 
It can be shown that  $\supp(\bPsi^\top(\bu - \tbbeta)) = \sS$, whence we have 
\begin{equation}\label{equation:errbd_noisem_q1ep_pv2}
\begin{aligned}
\normone{\bPsi^\top(\bu - \tbbeta)} 
&\leq \sqrt{\abs{\sS}} \normtwo{\bPsi^\top(\bu - \tbbeta)} 
\leq \sqrt{\abs{\sS}} \normtwo{\bu - \tbbeta}  \\
&\leq \sqrt{\abs{\sS}} \, \mathcalU(\sT) \normtwo{\bX(\bu - \tbbeta)} 
= C_3 \normtwo{\bX(\bu - \tbbeta)}, 
\end{aligned}
\end{equation}
where the first inequality follows from the standard bounds on vector norms (Exerciser~\ref{exercise:cauch_sc_l1l2}), and  we have used the assumption $\sigma_{\max}(\bPsi) = 1$ (Assumption 3) and the definition $C_3 = \mathcalU(\sT)\sqrt{\abs{\sS}}$. 
Furthermore, using the singular value bound \eqref{equation:svd_stre_bd} and the standard bounds on vector norms (Exerciser~\ref{exercise:cauch_sc_l1l2}) again, we have 
\begin{equation}\label{equation:errbd_noisem_q1ep_pv3}
\begin{aligned}
\footnotesize
&\normtwo{\bX(\bu - \tbbeta)} 
\leq \normtwo{\bX(\bbeta - \bu)} + \normtwo{\bX(\bbeta - \tbbeta)} 
\leq \normtwo{\bX} \normtwo{\bbeta - \bu} + \normtwo{\bX(\bbeta - \tbbeta)}  \\
&\leq \frac{\normtwo{\bX} \normtwo{\bPsi^\top(\bbeta - \bu)}}{\sigma_{\min}(\bPsi)} + \normtwo{\bX(\bbeta - \tbbeta)}  
\leq \frac{\normtwo{\bX} \normone{\bPsi^\top(\bbeta - \bu)}}{\sigma_{\min}(\bPsi)} + \normtwo{\bX(\bbeta - \tbbeta)}.
\end{aligned}
\end{equation}
Therefore, we get \eqref{equation:errbd_noisem_q1ep_eq2} after combining \eqref{equation:errbd_noisem_q1ep_pv1}, \eqref{equation:errbd_noisem_q1ep_pv2}, and \eqref{equation:errbd_noisem_q1ep_pv3}.

\paragraph{Proof of \eqref{equation:errbd_noisem_q1ep_eq3}.} 
Since $\balpha_\sS = \sign(\bPsi_\sS^\top\tbbeta) $ and $\sS=\supp(\bPsi^\top\tbbeta)$, we have 
$\innerproduct{\bPsi \balpha, \tbbeta} = \normone{\bPsi^\top \tbbeta}$, so the Bregman distance simplifies to $D_1^\balpha(\bbeta, \tbbeta) = \normone{\bPsi^\top \bbeta} - \innerproduct{\bPsi \balpha, \bbeta}$. 
Thus, it suffices to show
$$
f(\bbeta) \leq (1 - \norminf{\balpha_\sT})^{-1} (\normonebig{\bPsi^\top \bbeta} - \innerproduct{\bPsi \balpha, \bbeta}).
$$
Since $\bu \in \nspace(\bPsi(\sT)^\top)$ is equivalent to $\bPsi_\sT^\top \bu = \bzero$, the Lagrangian of the minimization problem in \eqref{equation:errbd_noisem_q1ep_eq3} is
$$
L(\bu, \bnu) = \normone{\bPsi^\top(\bbeta - \bu)} + \innerproduct{\bnu, \bPsi_\sT^\top \bu} = \normone{\bPsi^\top(\bbeta - \bu)} + \innerproduct{\bPsi_\sT \bnu, \bu - \bbeta} 
+ \innerproduct{\bPsi_\sT \bnu, \bbeta}.
$$
Then, $f(\bbeta) = \min_\bu \max_\bnu L(\bu, \bnu)$. Following the minimax theorem, we derive that
\begin{align*}
f(\bbeta) 
&= \max_\bnu \min_\bu L(\bu, \bnu) 
= \max_\bnu \min_\bu \left\{ \normone{\bPsi^\top(\bbeta - \bu)} + \innerproduct{\bPsi_\sT \bnu, \bu }  \right\} \\
&\stackrel{\dag}{=}\max_\bnu \min_\bu \max_{\bw\in\sB_\infty}
\left\{ \innerproduct{\bw, {\bPsi^\top(\bbeta - \bu)}}  + \innerproduct{\bPsi_\sT \bnu, \bu} \right\}\\
&=\max_{\bz\in \cspace(\bPsi_\sT)} \min_\bu \max_{\bw\in\sB_\infty} 
\left\{ \innerproduct{\bPsi\bw, \bbeta} + \innerproduct{\bu, \bz-\bPsi\bw} \right\}\\
&\stackrel{\ddag}{=} \max_{\bz}  
\left\{\innerproduct{\bz, \bw} \mid  \bz\in \bPsi \cdot \sB_\infty \cap  \cspace(\bPsi_\sT)\right\}\\
&= \max_\bz \left\{ \innerproduct{c\bPsi \balpha + \bz, \bbeta } \mid \bz \in \bPsi \cdot \sB_\infty \cap \cspace(\bPsi_\sT) \right\} - \innerproduct{c\bPsi \balpha, \bbeta}, \; \forall\, c > 0 \\
&= c \max_\bz \left\{ \innerproduct{\bPsi \balpha + \bz, \bbeta} \mid \bz \in c^{-1} \bPsi \cdot \sB_\infty \cap \cspace(\bPsi_\sT) \right\} - c \innerproduct{ \bPsi \balpha, \bbeta}, \; \forall\, c > 0,
\end{align*}
where the equality ($\dag$) follows from the definition of the dual norm (see \eqref{equation:dual_norm_equa}), the equality ($\ddag$) follows by letting $\bz-\bPsi\bw=\bzero$ (otherwise the minimum over $\bu$ is $-\infty$),  $\sB_\infty \triangleq \{\bw \mid \norminf{\bw}\leq 1\}$ denotes the unit $\ell_\infty$-ball, and $\cspace(\bPsi_\sT)$ denotes the column space of $\bPsi_\sT$.
Let\ $c \triangleq (1 - \norminf{\balpha_\sT})^{-1}$
and $\sB_\sT \triangleq \{ \bb \in \real^q \mid \bb_\sS = \bzero \}$. 
Since $\balpha_\sS = \sign(\bPsi_\sS^\top \tbbeta)$ and  $\norminf{\balpha_\sT} < 1$ from Condition~\ref{condition:q1_cond1} (ii), we have $c < +\infty$ and get
$$
(\balpha + c^{-1} \sB_\infty \cap \sB_\sT) \subset \sB_\infty,
$$
from which we conclude
$$
(\bPsi \balpha + c^{-1} \bPsi \cdot \sB_\infty \cap \cspace(\bPsi_\sT)) \subset \bPsi \cdot \sB_\infty.
$$
Hence, for any $\bz \in c^{-1} \bPsi \cdot \sB_\infty \cap \cspace(\bPsi_\sT)$, it holds $\bPsi \balpha + \bz \in \bPsi \cdot \sB_\infty$, which by the \holders inequality implies $\innerproduct{\bPsi \balpha + \bz, \bbeta } \leq \normone{\bPsi^\top \bbeta}$. Therefore, $f(\bbeta) \leq c(\normone{\bPsi^\top \bbeta} - \innerproduct{\bPsi \balpha, \bbeta })$.
This completes the proof.
\end{proof}

\begin{lemma}\label{lemma:errbd_noisem_q1ep2}
Let $\bbeta^* \in \real^p$ be a fixed vector such that $\supp(\bPsi^\top \bbeta^*) = \sS$. 
Suppose there exists a subgradient $\balpha \in \partial \normone{\bPsi^\top \bbeta^*}$ satisfying  $\bPsi \balpha = \bX^\top \blambda \in \cspace(\bX^\top)$, 
where from $\bPsi \balpha = \bX^\top \blambda$ and the full-rankness of $\bX$ (with $n\leq p$), we have $\blambda = (\bX \bX^\top)^{-1} \bX \bPsi \balpha$.
Then, for every $\epsilon > 0$ and every data vector $\by$ satisfying $\normtwo{\bX \bbeta^* - \by} \leq \epsilon$, the following two statements hold:
\begin{enumerate}
\item Every minimizer $\bbeta_{\sigma}$ of problem \eqref{opt:q1_penalize} satisfies
$$
D_1^\balpha(\bbeta_{\sigma}, \bbeta^*) \leq \frac{(\epsilon + \sigma \normtwo{\blambda}/2)^2}{\sigma} \quad \text{and} \quad \normtwo{\bX \bbeta_{\sigma} - \by}\leq \epsilon + \sigma \normtwo{\blambda};
$$
\item Every minimizer $\bbeta_\epsilon$ of problem \eqref{opt:q1_epsilon} satisfies $D_1^\balpha(\bbeta_\epsilon, \bbeta^*) \leq 2\epsilon \normtwo{\blambda}$.
\end{enumerate}
\end{lemma}
\begin{proof}
See Theorem~3 of \citet{burger2004convergence} and Lemma 3.5 of \citet{grasmair2011necessary}.
\end{proof}

From the above lemmas, we now show that Condition~\ref{condition:q1_cond1} guarantees robustness of problems~\eqref{opt:q1_penalize} and~\eqref{opt:q1_epsilon} to arbitrary measurement noise in $\by$.
\begin{theoremHigh}[Error bound of \eqref{opt:q1_penalize} and \eqref{opt:q1_epsilon} under measurement noise \citep{zhang2016one}]\label{theorem:robst_optq}
Under Assumptions~1$\sim$3, let $\bbeta^* \in \real^p$ be a given signal, and let $\sS \triangleq \supp(\bPsi^\top \bbeta^*)$ and $\sT \triangleq \comple{\sS}$. For arbitrary noise $\bepsilon$, let $\by = \bX \bbeta^* + \bepsilon$ and $\epsilon = \normtwo{\bepsilon}$. 
If Condition~\ref{condition:q1_cond1} holds with $\tbbeta = \bbeta^*$, then
\begin{enumerate}[(i)]
\item For any $C_0 > 0$, there exists constant $C_1 > 0$ such that every minimizer $\bbeta_{\sigma}$ of problem \eqref{opt:q1_penalize} with parameter $\sigma = C_0\epsilon$ satisfies
$$
\normone{\bPsi^\top(\bbeta_{\sigma} - \bbeta^*)} \leq C_1\epsilon;
$$
\item Every minimizer $\bbeta_\epsilon$ of problem \eqref{opt:q1_epsilon} satisfies
$$
\normone{\bPsi^\top(\bbeta_\epsilon - \bbeta^*)} \leq C_2\epsilon.
$$
\end{enumerate}
The constants $C_1$ and $C_2$ are given explicitly as follows. Define
$$
\blambda = (\bX\bX^\top)^{-1}\bX\bPsi \balpha, \quad C_3 = \mathcalU(\sT)\sqrt{\abs{\sS}}, \quad \text{and} \quad C_4 = \frac{1 + \cond(\bPsi)\normtwo{\bX}C_3}{1 - \norminf{\balpha_\sT}},
$$
with which we have
\begin{align*}
C_1 &= 2C_3 + C_3C_0\normtwo{\blambda} + \frac{(1 + C_0\normtwo{\blambda}/2)^2 C_4}{C_0}, \\
C_2 &= 2C_3 + 2C_4\normtwo{\blambda}.
\end{align*}
\end{theoremHigh}

\begin{proof}[of Theorem~\ref{theorem:robst_optq}]
Firstly, we derive that
\begin{align*}
\normone{\bPsi^\top \bbeta_{\sigma} - \bbeta^*} &\leq C_3 \normtwo{\bX(\bbeta_{\sigma} - \bbeta^*)} + C_4 D_1^\balpha(\bbeta_{\sigma}, \bbeta^*)  \\
&\leq C_3 \normtwo{\bX \bbeta_{\sigma} - \by} + C_3 \normtwo{\bX \bbeta^* - \by} + C_4 D_1^\balpha(\bbeta_{\sigma}, \bbeta^*)  \\
&\leq C_3(\epsilon + \sigma \normtwo{\blambda}) + C_3 \epsilon + C_4 \frac{(\epsilon + \sigma \normtwo{\blambda}/2)^2}{\sigma}, 
\end{align*}
where the first and the third inequalities follow from Lemmas~\ref{lemma:errbd_noisem_q1ep} and~\ref{lemma:errbd_noisem_q1ep2}, respectively. Substituting $\sigma = C_0 \epsilon$  and simplifying yields the bound in part (i). 
Part (ii) follows analogously by applying the corresponding bound from Lemma~\ref{lemma:errbd_noisem_q1ep2} for the constrained formulation.
\end{proof}

Note that the condition  $\norminf{\balpha_\sT} < 1$ plays a crucial role: the constant $C_4$ is inversely proportional to $1 - \norminf{\balpha_\sT}$, so this gap directly controls stability.
Moreover, the numerator of $C_4$ also depends on  global properties of $\bPsi$ and $\bX$, 
whereas, $C_3$  depends only on the restricted geometry of the support $\sS$.

\begin{problemset}

\item \label{prob:1spar_ells} \textbf{Non-recovery for $s>1$.}
Let $\bX\in\real^{n\times p}$ with $n < p$. 
Show that there exists a $1$-sparse vector which is not a minimizer of 
$$
\min_{\bbeta\in\real^p}\norms{\bbeta}
\quad\text{s.t.}\quad  
\bX\bbeta=\by,
$$
where $s>1$.

\item \label{prob:sparsi_of_p1_n}\textbf{Linear independence of \eqref{opt:p1} under uniqueness.}
Let $\bX \in \real^{n\times p}$ be a measurement matrix. 
Assume that the minimizer $\widehatbbeta$ of the optimization problem
$$
(\text{P}_1)
\quad 
\min_{\bbeta \in \real^p} \normone{\bbeta} 
\quad \text{s.t.}\quad
\bX\bbeta = \by,
$$
is unique. Then the set of columns $\{\bx_i, i \in \supp(\widehatbbeta)\}$ is linearly independent. In particular,
$$
\normzerobig{\widehatbbeta} = \abs{\supp(\widehatbbeta)} \leq n.
$$

\item \label{prob:qnorm}\textbf{$\bS$-norm.}
Let $\norm{\cdot}$ be a valid vector norm on $\real^p$. And let $\bS\in\real^{n\times p}$ have full column rank. Show that $\norm{\cdot}_{\bS} \triangleq \norm{\bS\cdot}$ is also a valid norm satisfying  the three axioms of norms in Definition~\ref{definition:matrix-norm}.

\item \label{prob:push_through_ide} \textbf{Push-through identity.}
Let $\bX$ and $\bY$ be two $n \times p$ matrices. Show the following result:
\begin{equation}
\bX^\top(\bI_n + \bY\bX^\top)^{-1} = (\bI_p + \bX^\top\bY)^{-1}\bX^\top.
\end{equation}
Use the above result to show the following for any $n \times p$ matrix $\bZ$ and scalar $\lambda > 0$:
\begin{equation}
\bZ^\top(\lambda \bI_n + \bZ\bZ^\top)^{-1} = (\lambda \bI_p + \bZ^\top \bZ)^{-1}\bZ^\top.
\end{equation}
The push-through identity derives its name from the fact that we push in a matrix on the left and it comes out on the right. 
\textit{Hint: Premultiply and postmultiply the  above identities with  appropriate matrices. }

\item \textbf{Alternating least squares (ALS) \citep{lu2021numerical}.} Given a matrix $\bX\in\real^{n\times p}$, use the least squares (with or without regularization) models to find a matrix factorization $\bX\approx \bW\bZ$.

\item 
\textbf{Data least squares.}
While the OLS method accounts for errors in the response variable $\by$, the \textit{data least sqaures (DLS)} method considers errors in the predictor variables:
\begin{equation}
\bbeta^{DLS} = \mathop{\argmin}_{\bbeta, \widetildebX} \normfbig{\widetildebX}^2, \quad \text{s.t.}\quad \by\in\cspace(\bX+\widetildebX),
\end{equation}
where $\widetildebX$ represents a perturbation of $\bX$ (i.e., a noise in the predictor variables).
That is, $(\bX+\widetildebX) \bbeta^{DLS} = \by$, assuming the measured response $\by$ is noise-free.
The Lagrangian function and its gradient w.r.t. $\bbeta$ are, respectively, given by
$$
\begin{aligned}
	L(\bbeta, \widetildebX, \blambda) &= \trace(\widetildebX\widetildebX^\top) +\blambda^\top (\bX\bbeta+\widetildebX\bbeta-\by);\\
	\nabla_{\widetildebX} L(\bbeta, \widetildebX,\blambda) &= \widetildebX+\blambda\bbeta^\top = \bzero \quad\implies\quad \widetildebX=-\blambda\bbeta^\top,
\end{aligned}
$$
where $\blambda\in\real^n$  is a vector of Lagrange multipliers.
Show that the original problem is equivalent to 
$$
\mathop{\argmin}_{\bbeta}
\frac{(\bX\bbeta-\by)^\top (\bX\bbeta-\by)}{\bbeta^\top\bbeta} .
$$

\item \textbf{Total least squares.} Similar to  data least squares, the \textit{total least squares (TLS)} method considers errors in both the predictor variables and the response variables. The TLS problem can be formulated as:
\begin{equation}
\bbeta^{TLS} = \mathop{\argmin}_{\bbeta, \widetildebX, \widetildeby} \normf{[\widetildebX, \widetildeby]}^2, 
	\quad \text{s.t.}\quad (\by+\widetildeby)\in\cspace(\bX+\widetildebX), 
\end{equation}
where $\widetilde{\bX}$ and $\widetilde{\by}$ are perturbations in the predictor variables and the response variable, respectively.
Show that the problem can be equivalently stated as
\begin{equation}
\bbeta^{TLS} = \mathop{\argmin}_{\bz, \bD} \normf{\bD}^2, 
\quad \text{s.t.}\quad \bD\bz = -\bC\bz, 
\end{equation}
where $\bC\triangleq[\bX,\by]\in\real^{n\times (p+1)}$,  $\bD\triangleq[\widetildebX, \widetildeby]\in\real^{n\times (p+1)}$, and $\bz\triangleq\footnotesize\begin{bmatrix}
	\bbeta\\
	-1
\end{bmatrix}$.

\item Discuss the geometry of $\Theta_1$, $\Theta_2$ and $\Theta$ in the dual space of SAFE rule. Draw the figure. What do you observe?

\item Prove the result in \eqref{equation:safe_pracmax} and \eqref{equation:safe_pracmax2} rigorously.

\item \textbf{SAFE for LASSO with intercept.}
Show that the SAFE--LASSO theorem (Theorem~\ref{theorem:safe_lasso}) can be applied to the LASSO with intercept problem:
$$
 \min_{\bbeta,\bgamma} \frac{1}{2} \normtwo{\bX \bbeta + \bgamma - \by}^2 + \lambda \normone{\bbeta},
$$
where $\bgamma \in \real^n$ denotes the intercept term. 
\textit{Hint: Use a   transformation.}

\item \textbf{SAFE for elastic net.}
Show that the SAFE--LASSO theorem (Theorem~\ref{theorem:safe_lasso}) can be applied to the  \textit{elastic net} problem:
$$
\min_\bbeta \frac{1}{2} \normtwo{\bX \bbeta - \by}^2 + \lambda \normone{\bbeta} + \frac{1}{2} \sigma \normtwo{\bbeta}^2.
$$
\textit{Hint: Use a   transformation.}

\item 
Similar to the sequential DPP rule in Corollary~\ref{corollary:seq_dpp}, design an algorithm that leverages the SAFE rule to reduce memory usage for solving LASSO problems.

\item \textbf{Piecewise linear update of the fitted value.}
Consider the constrained LASSO problem \eqref{opt:lc} (p.~\pageref{opt:lc}):
$$
\mathopmin{\normone{\bbeta}\leq \Sigma}f(\bbeta)\triangleq \frac{1}{2}\normtwo{\by - \bX \bbeta}^2.
$$
Suppose  two constraint radii $\Sigma_1$ and $\Sigma_2$, lie within the same region of the solution path, where both the active set $\sS= \{i: \beta_i \neq 0\}$ and sign vector $\bs_\sS$ remain unchanged.
Here, the full sign vector $ \bs\in\real^p$ is defined component-wise as
$$
s_i =
\begin{cases}
	\sign(\beta_i), & \beta_i \neq 0; \\
	\in [-1,1], & \beta_i = 0,
\end{cases}
$$
Show that  the corresponding solutions $\bbeta_1$ and $\bbeta_2$ satisfy
$$
\bbeta_2 - \bbeta_1
= (\Sigma_2 - \Sigma_1)
\frac{(\bX_\sS^\top \bX_\sS)^{-1} \bs_\sS}
{\bs_\sS^\top (\bX_\sS^\top \bX_\sS)^{-1} \bs_\sS}.
$$
That is, the solution $\bbeta$ changes linearly with $\Sigma$
along the fixed direction
$$
\bv = \frac{(\bX_\sS^\top \bX_\sS)^{-1} \bs_\sS}
{\bs_\sS^\top (\bX_\sS^\top \bX_\sS)^{-1} \bs_\sS}.
$$

\item \textbf{Dual and primal optimum.} Does the converse of Lemma~\ref{lemma:dual_cer_p1_prem} hold?

\item  Prove that the dual of \eqref{equation:dual_cer_nomrinf_equiv} is \eqref{equation:dual_cer_nomrinf_equiv_dd}.

\item Derive the dual problem of $
\min_{\bbeta\in \real^p} \normtwo{\bbeta }$ subject to $ \bX\bbeta = \by
$

\item \label{prob:dantzig_dual_cert} \textbf{Dual certificate for Dantzig selector problem.}
Consider the following Dantzig selector (DS) problem \eqref{opt:p1_dantzig}:
\begin{equation}\label{equation:optdang_prob}
\begin{aligned}
\min_{\bbeta \in \real^p} \quad & \normone{\bbeta}
\quad \text{s.t.} \quad & \norminf{\bX^\top (\bX \bbeta - \by)} \leq D.
\end{aligned}
\end{equation}
\begin{enumerate}
\item Derive the dual problem of problem \eqref{equation:optdang_prob}.
\item Provide a dual certificate analysis for unique optimal solution of problem \eqref{equation:optdang_prob}.
\end{enumerate}
\textit{Hint: see Theorem~\ref{theorem:dualcert_p1} for \eqref{opt:p1}.}

\item \label{prob:group_bp} \textbf{Group basis pursuit (GBP).}
Consider the following group basis pursuit (GBP) problem:
\begin{equation}\label{equation:gbp_prob}
\begin{aligned}
\min_{\bbeta \in \real^{pT}} \quad & \sum_{t=1}^T \normtwo{\bbeta_t} 
\quad \text{s.t.} \quad  \bX \bbeta = \by,
\end{aligned}
\end{equation}
where $\bbeta = [\bbeta_1^\top, \bbeta_2^\top, \ldots, \bbeta_T^\top]^\top$, and $\forall\, t \in \{1, 2, \ldots, T\}$, $\bbeta_t \in \real^p$.

Assumptions: matrix $\bX$ has full row rank and the solution set $\sB$ of \eqref{equation:gbp_prob} is nonempty.
\begin{enumerate}
\item  Derive the dual problem of problem \eqref{equation:gbp_prob}.
\item Provide a dual certificate analysis for unique optimal solution of problem \eqref{equation:gbp_prob}.
\end{enumerate}
\textit{Hint: see Theorem~\ref{theorem:dualcert_p1} for \eqref{opt:p1}.}


\item \label{prob:nphard_ells} \textbf{NP-hardness of $\ell_s$-minimization for $0 < s < 1$ \citep{foucart2013invitation}.}
Given a matrix $\bX \in \real^{n\times p}$ and a vector $\by \in \real^n$, the $\ell_s$-minimization problem seeks a vector  $\bbeta \in \real^p$ that minimizes the $\ell_s$-quasinorm subject to the linear constraint  $\bX\bbeta = \by$. Assuming the NP-completeness of the \textit{partition problem}, which, given integers $x_1,  x_2, \ldots, x_k$, asks whether the index set $\{1,2,\ldots,k\}$ can be split into two disjoint subsets $\sS, \sT \subset \{1,2,\ldots,k\}$ such that $\sum_{i \in \sS} x_i = \sum_{j \in \sT} x_j$, prove that the $\ell_s$-minimization problem is NP-hard.
To establish this reduction, consider the following construction. Define the matrix $\bX$ and the vector $\by$ by
$$
\bX \triangleq 
\begin{bmatrix}
	x_1 & x_2 & \ldots & x_k & x_1 & x_2 & \ldots & x_k \\
	1 & 0 & \ldots & 0 & 1 & 0 & \ldots & 0 \\
	0 & 1 & \ldots & 0 & 0 & 1 & \ldots & 0 \\
	\vdots & \vdots & \ddots & \vdots & \vdots & \vdots & \ddots & \vdots \\
	0 & 0 & \ldots & 1 & 0 & 0 & \ldots & 1
\end{bmatrix}
\quad \text{and} \quad
\by = [0, 1, 1, \ldots, 1]^\top.
$$
\end{problemset}

\part{Design Matrices}

\newpage 
\chapter{Design Matrices and Properties}\label{chapter:design}
\begingroup
\hypersetup{
linkcolor=structurecolor,
linktoc=page,  
}
\minitoc \newpage
\endgroup

\lettrine{\color{caligraphcolor}I}
In sparse representation and \textit{compressed sensing} (also known as \textit{compressive sensing}, \textit{compressive sampling}, or \textit{sparse sampling}; see Section~\ref{section:sparse_opt_gen}), one must solve an underdetermined linear system $$
\bX\bbeta = \by,
$$ 
where $\bbeta$ is a sparse vector containing only a few nonzero entries, $\by$ is the vector of  \textit{measurements} (or \textit{reesponse}), and $\bX$ is the \textit{design matrix}  (or \textit{sensing matrix}). 
Sparsity has long been leveraged in signal processing and approximation theory---for tasks such as image compression \citep{pennebaker1992jpeg, taubman2002jpeg2000, rabbani2002overview, foucart2013invitation} and denoising \citep{donoho2002noising, joy2013denoising}---as well as in statistics and learning theory as a means to prevent overfitting \citep{vapnik2013nature, lu2021rigorous}. It also plays a central role in statistical estimation and model selection \citep{tibshirani1996regression, lu2021rigorous}, in modeling the human visual system \citep{olshausen1996emergence}, and especially in image processing, where the multiscale wavelet transform  yields nearly sparse representations for natural images \citep{mallat1999wavelet}.
In this chapter, we discuss the role of design matrices in solving sparse linear systems.

Transforming a signal into a new basis or frame---i.e., expressing it as $\by=\bX\bbeta$ so that $\bbeta$ represents the signal in that basis---can lead to a more concise representation. 
Such compression is valuable for reducing data storage and transmission costs, which can be substantial in many applications. Consequently, rather than transmitting the original high-dimensional signal, one may instead transmit only the (typically few) nonzero analysis coefficients from the basis or frame expansion. When the number of nonzero coefficients is small, we say the signal admits a sparse representation. 
Sparse signal models enable high compression rates, and in compressive sensing, the prior knowledge that a signal is sparse in some known basis or frame allows us to reconstruct the original signal from far fewer measurements than its ambient dimension would suggest. For sparse data, only the nonzero coefficients generally need to be stored or transmitted; all others can safely be assumed to be zero.

Mathematically, a signal $\bbeta$ is called $k$-sparse if it contains at most $k$ nonzero entries, i.e., $\normzero{\bbeta} \leq k$. We denote the set of all \textit{$k$-sparse signals} by
\begin{equation}
\sB_0[k] = \{ \bbeta : \normzero{\bbeta} \leq k \},
\end{equation}
~\footnote{Recall that  $\sB_0[k]=\sB_0[\bzero,  k]$ denotes the set of all $k$-sparse vectors (Definition~\ref{definition:open_closed_ball}). Although sometimes informally called the ``$\ell_0$-ball," the ``$\ell_0$-norm" is not a true norm.} where square brackets are used to emphasize that equality is allowed. 
Using this notation, problem \eqref{opt:s0_gen} (p.~\pageref{opt:s0_gen}) can be restated as:
$$
\min \; f(\bbeta) \quad \text{ s.t. }\quad\bbeta \in \sB_0[R] .
$$
In practice, many signals are not sparse in the standard basis but admit a sparse representation in some other basis $\bPhi$. In such cases, we still refer to $\bbeta$ as $k$-sparse, with the understanding that $\bbeta = \bPhi \balpha$  for some coefficient vector $\balpha$ satisfying $\normzero{\balpha} \leq k$.

In many applications, sparse optimization involves minimizing the objective function
$$
f(\bbeta) = \frac{1}{2}\normtwo{\by - \bX \bbeta}^2,
$$
where $\bX\in\real^{n\times p}$ is  the \textit{design matrix}, and $\by\in\real^n$ is the observed response vector.
Notably, when the design matrix is the identity, it preserves the geometric structure of the signal, thereby enabling exact recovery of sparse vectors. We will formalize this idea and establish precise conditions on $\bX$ that guarantee universal recovery---that is, the ability to uniquely reconstruct every $\bbeta \in \sB_0[k]$  from its measurements $\by=\bX\bbeta$.

However, a design matrix that can identify some sparse vectors does not necessarily ensure universal recovery. For example, consider a matrix $\bX \in \real^{n\times p}$ and two distinct $k$-sparse vectors $\bbeta_1, \bbeta_2 \in \sB_0[k]$ with  $\bbeta_1 \neq \bbeta_2$, yet $\bX \bbeta_1=\bX \bbeta_2=\by$. In this case, the measurements provide no way to distinguish between $\bbeta_1$ and $\bbeta_2$, making unique recovery impossible. Consequently, such a design matrix fails to support universal recovery. Moreover, this failure is not limited to a single pair of vectors: it typically implies that an entire subspace---or infinitely many pairs---of sparse vectors cannot be uniquely identified from their measurements.

Therefore, it is essential that the design matrix preserves the geometric structure of sparse vectors when projecting them from a high-dimensional space (of dimension $p$) into a lower-dimensional space (of dimension $n$), particularly in settings where $n \ll p$. The nullspace property and other properties in this chapter  formalize this requirement.

\subsection*{Basic Notations}
We now introduce some basic notation.
For any index set $\sS \subseteq \{1, 2,\ldots, p\}$, define the following sets:
\begin{subequations}\label{equation:cscsgamma_def}
\begin{align}
\sC[\sS] &\triangleq \{\bbeta \in \real^{p} \mid  \normone{\bbeta_{\comple{\sS}}} \leq \normone{\bbeta_{\sS}} \};\\
\sC[\sS; \gamma] &\triangleq \{\bbeta \in \real^{p} \mid  \normone{\bbeta_{\comple{\sS}}} \leq \gamma\normone{\bbeta_{\sS}} \};\\
\sC[k] &\triangleq \bigcup_{\sS: \abs{\sS} = k} \sC[\sS],
\end{align}
\end{subequations}
where $\comple{\sS}$ denotes  the complement of $\sS$.
Here, ``{square brackets}" are used (rather than parentheses) to indicate that equality is allowed in the defining inequalities.
Thus, $\sC[\sS]=\sC[\sS; 1]$ is the convex cone consisting of vectors whose $\ell_1$-mass outside $\sS$ does not exceed that inside $\sS$, i.e., vectors that concentrate most of their weight on the coordinates in $\sS$.
In contrast, $\sC[k]$ is the (non-convex) union of all such cones over index sets of size $k$; it captures vectors that concentrate most of their $\ell_1$-norm on some subset of $k$ coordinates.
Note that 
$$ \sB_0[k] \subset \sC[k], 
$$
since any $k$-sparse vector  places all its nonzero entries---and hence all its $\ell_1$-mass---within a set of at most $k$ coordinates.
Additionally, for a fixed $\sS$, we have 
$$
\sC[\sS; \gamma_1] \subseteq \sC[\sS; \gamma_2],
\quad\text{whenever }
\gamma_1\leq \gamma_2.
$$

In certain contexts, strict inequalities are required. In such cases, the above definitions can be restated using ``parentheses" to exclude equality:
\begin{subequations}\label{equation:cscsgamma_def2}
\begin{align}
\sC(\sS) &\triangleq \{\bbeta \in \real^{p} \mid  \normone{\bbeta_{\comple{\sS}}} < \normone{\bbeta_{\sS}} \};\\
\sC(\sS; \gamma) &\triangleq \{\bbeta \in \real^{p} \mid  \normone{\bbeta_{\comple{\sS}}} < \gamma\normone{\bbeta_{\sS}} \};\\
\sC(k) &\triangleq \bigcup_{\sS: \abs{\sS} = k} \sC(\sS),
\end{align}
\end{subequations}

\index{$k$-sparse approximation}
\index{$k$-term approximation}
\index{Optimal $k$-sparse approximation}
Before discussing these properties further, we introduce one more definition.
The \textit{$\ell_s$-error of best $k$-sparse approximation} (also called the \textit{$\ell_s$-distance to the set of $k$-sparse vectors}) is defined as:
\begin{equation}\label{equation:lsdist_spar}
\sigma_k(\bbeta)_s \triangleq \min_{\widehat{\bbeta} \in \sB_0[k]} \normsbig{\bbeta - \widehat{\bbeta}}.
\end{equation}

\section{Nullspace, Restricted Eigenvalue, and Other Design Properties}

We now introduce several desirable properties of the design matrix that are widely studied in the literature (see, e.g.,  \citet{candes2005decoding, cohen2009compressed, raskutti2010restricted, baraniuk2011introduction, jain2014iterative, hastle2015statistical, jain2017non}).

We begin with the nullspace property.
\begin{definition}[Nullspace property (NSP)\index{Nullspace property}]\label{definition:nullspace_prop}
A matrix  $\bX \in \real^{n\times p}$ is said to satisfy the \textit{nullspace property (NSP) or restricted nullspace property (RNP) over the set $\sS\subseteq\{1,2,\ldots,p\}$}, denoted  NSP($\sS$), if 
\begin{subequations}\label{equation:srnp_equ2}
\begin{equation}\label{equation:srnp_e1}
\qquad 	\nspace(\bX) \cap \sC[\sS] = \{\bzero\},
\end{equation}
where $\nspace(\bX) \triangleq \{\bbeta \in \real^{p}\mid \bX \bbeta = \bzero\}$ denotes the nullspace of $\bX$ (Definition~\ref{definition:nullspace}).
Equivalently, this condition can be expressed as:
\begin{equation}\label{equation:srnp_e2}
\text{or }\quad \normone{\bbeta_{\sS}} 
< 
\normone{\bbeta_{\comple{\sS}}}, 
\qquad \text{ for all } \bbeta \in \nspace(\bX)\setminus\{\bzero\};
\end{equation}
\begin{equation}\label{equation:srnp_e3}
\text{or }\quad 2\normone{\bbeta_{\sS}} 
< 
\normone{\bbeta}, 
\qquad \text{ for all } \bbeta \in \nspace(\bX)\setminus\{\bzero\}.
\end{equation}

On the other hand,  $\bX \in \real^{n\times p}$ is said to satisfy the \textit{NSP of order $k$}, denoted  NSP($k$), if $\nspace(\bX) \cap \sC[k] = \{\bzero\}$.
Note that, for a given $\bbeta \in \nspace(\bX)\setminus\{\bzero\}$, the
inequality  $\normone{\bbeta_{\sS}} 
< 
\normone{\bbeta_{\comple{\sS}}}$  holds for any index set $\sS\subset \{1,2,\ldots,p\}$ with $\abs{\sS}\leq k$ as soon as it holds for an index set of $k$ largest (in absolute value) entries of $\bbeta$.
Consequently, for NSP of order $k$, \eqref{equation:srnp_e2} or \eqref{equation:srnp_e3} can be equivalently stated as 
\begin{equation}
\normone{\bbeta} 
< 
2\sigma_k(\bbeta)_1, 
\qquad \text{ for all } \bbeta \in \nspace(\bX)\setminus\{\bzero\},
\end{equation}
where $\sigma_k(\bbeta)_s \triangleq \min_{\widehat{\bbeta} \in \sB_0[k]} \normsbig{\bbeta - \widehat{\bbeta}}$ is the $\ell_s$-distance to  $k$-sparse vectors.
\end{subequations} 
\end{definition}

While the NSP of order $k$ (i.e., NSP($k$)) is defined with respect to the union $\sC[k] = \bigcup_{\sS: \abs{\sS} = k} \sC[\sS]$,  the NSP over $\sS$ (i.e, NSP($\sS$)) applies to a fixed index set $\sS\subseteq \{1,2,\ldots,p\}$.
Therefore, NSP over $\sS$ is a less strict condition compared with NSP of order $k$ if $\abs{\sS}=k$.
In words, the NSP($\sS$) property holds when the only vector in the cone $\sC[\sS]$ that lies in the nullspace of $\bX$ is the zero vector.

Recall that the NSP of order $k$ states the following: for every nonzero vector $\bbeta\in\nspace(\bX)\setminus \{\bzero\}$ and for every index set $\sS\subset \{1,2,\ldots,p\}$ with $\abs{\sS}\leq k$, we have $\normone{\bbeta_{\sS}}< \normone{\bbeta_{\comple{\sS}}}$.
In other words, any nonzero vector in the nullspace of $\bX$ must be more ``spread out" outside any subset of at most $k$ coordinates than inside it.
This implies that no nonzero vector in $\nspace(\bX)$ can concentrate most of its $\ell_1$-mass on any $k$ coordinates. 
As a direct consequence, the nullspace cannot contain any nonzero $k$-sparse vector.

\begin{exercise}\label{exercise:nspi_2k}
	Let $\bX$ be a design matrix that satisfies the  NSP  of order $2k$. Show that for any two distinct $k$-sparse vectors $\bbeta_1, \bbeta_2 \in \sB_0[k]$ with $\bbeta_1 \neq \bbeta_2$, it holds that $\bX\bbeta_1 \neq \bX\bbeta_2$.
\end{exercise}

\paragrapharrow{Rescaled, reshuffled, and augmented  properties of NSP.}
Consider replacing the original measurement matrix $\bX$ with new measurement matrices $\widehatbX$ or $\widetildebX$,  defined as follows:
\begin{align*}
\widehatbX  &\triangleq \bP\bX, \quad \text{where } \bP \text{ is some invertible } n\times n \text{ matrix}; \\
\widetildebX &\triangleq \begin{bmatrix} \bX \\ \bY \end{bmatrix}, \quad \text{where } \bY \text{ is some } n' \times p  \text{ matrix}.
\end{align*}
When $\bP$ is a diagonal matrix, $\bP\bX$ corresponds to a rescaled version of $\bX$; when $\bP$ is a permutation matrix (Definition~\ref{definition:permutation-matrix}), $\bP\bX$ represents a reshuffled (i.e., row-permuted) version of  $\bX$.
We observe that 
\begin{equation}\label{equation:nsp_chg}
\nspace(\widehatbX)  = \nspace(\bX) 
\qquad\text{and}\qquad  
\nspace(\widetildebX)  \subseteq \nspace(\bX).
\end{equation}
Consequently, if the NSP condition holds for $\bX$, it also holds for both $\widehatbX$ and $\widetildebX$.

If a design matrix satisfies the NSP  of order $2k$, then no two distinct $k$-sparse vectors can be mapped to the same measurement vector (see Exercise~\ref{exercise:nspi_2k}). 
This injectivity on the set of $k$-sparse vectors is essential for guaranteeing global uniqueness in sparse recovery.
A stronger condition than NSP is the restricted eigenvalue (RE) property, which not only rules out sparse vectors in the nullspace but also ensures that sparse vectors retain a significant portion of their norm under the action of the measurement matrix.

\begin{definition}[Restricted eigenvalue (RE) property\index{Restricted eigenvalue property}]\label{definition:res_eig}
A matrix $\bX \in \real^{n\times p}$ is said to satisfy the \textit{restricted eigenvalue (RE) property with constant $\mu>0$} over a set $\sS$ if, for all $\bbeta \in \sS$, we have
$$
\mu\leq \frac{ \bbeta^\top\bX^\top\bX\bbeta}{\normtwo{\bbeta}^{2}}
\qquad\iff\qquad
\mu \cdot \normtwo{\bbeta}^{2} \leq  \normtwo{\bX\bbeta}^{2},
$$
where the right-hand side of the first inequality represents  the \textit{Rayleigh quotient}~\footnote{see, for example, \citet{lu2021numerical} for more details.} of the vector $\bbeta$ associated with the matrix $ \bX^\top\bX$, which also serves as the Hessian of the least squares loss function $f(\bbeta) = \frac{1}{2}\normtwo{\bX\bbeta-\by}^2$.
\end{definition}

This property implies that sparse vectors are not only excluded from the nullspace of $\bX$, but are also well-conditioned under the linear transformation defined by $\bX$: their Euclidean norm is preserved up to a factor of $\mu$.
In particular, if the RE property holds over the set of all $2k$-sparse vectors (i.e., $\sS=\sB_0[2k]$), then for any two $k$-sparse vectors $\bbeta_1,\bbeta_2\in\sB[k]$, their difference $\bbeta_1\bbeta_2$ is $2k$-sparse, and we have
\begin{equation}
\bbeta_1, \bbeta_2 \in \sB_0[k]
\quad\implies\quad
\normtwo{\bX(\bbeta_1 - \bbeta_2)}^{2} \geq \mu \cdot \normtwo{\bbeta_1 - \bbeta_2}^{2}.
\end{equation}
Consequently, the mapping $\bbeta\mapsto \bX\bbeta$ preserves distances between $k$-sparse vectors up to the constant $\mu$. 
This stability is crucial for both identifiability and robust recovery in high-dimensional statistics and compressed sensing.

\subsection{Nullspace Property-II}
The NSP condition is not the only condition that enforces the non-concentration of vectors in the nullspace.
In practice, the vectors we aim to recover---whether via basis pursuit or other reconstruction methods---are rarely exactly sparse. 
Exact sparsity is typically an idealization; more realistically, signals are only approximately sparse, meaning they are close (in some norm) to a $k$-sparse vector. 
In such cases, we would like the reconstruction error to be controlled by how well the true vector can be approximated by a $k$-sparse one. This desirable behavior is known as \textit{stability} with respect to sparsity defect.

To formalize this idea, we introduce an alternative formulation called the nullspace property-II (NSP$'$).

\begin{definition}[Nullspace property-II (NSP$'$)]\label{definition:nsp_ii}
A matrix $\bX\in\real^{n\times p}$ satisfies the \textit{NSP$'$} over a set $\sS \subset \{1,2,\ldots,p\}$ with constant $C$, denoted  NSP$'$($\sS, C$),  
if $\abs{\sS}=k$ and there exists a constant $C > 0$ such that
\begin{equation}\label{equation:nsp-ii}
\normtwo{\bbeta_\sS} \leq C \frac{\normone{\bbeta_{\comple{\sS}}}}{\sqrt{k}},
\quad \text{for all } \bbeta \in \nspace(\bX).
\end{equation}
It is said to satisfy the \textit{NSP$'$ of order $k$} with constant $C$, denoted  NSP$'$($k, C$), if it satisfies the  NSP$'$ with constant $C$ relative to any set $\sS \subseteq \{1,2,\ldots,p\}$ with $\abs{\sS} \leq k$.
\end{definition}

Note that the equality can be achieved in the definition of the  NSP$'$, while the equality is not achieved in the definition of the standard NSP.

Like the standard NSP, the NSP$'$ condition formalizes the idea that no nonzero vector in the nullspace of $\bX$ can be overly concentrated on a small subset of coordinates.
For example, suppose  $\bbeta$ is exactly $k$-sparse. 
Then there exists a  set $\sS$ with $\abs{\sS}\leq k$ such that $\normone{\bbeta_{\comple{\sS}}} = \bzero$.
Plugging this into \eqref{equation:nsp-ii} yields $\bbeta_\sS = \bzero$, which implies $\bbeta=\bzero$. 
Consequently, if  $\bX$ satisfies  NSP$'$,  the only $k$-sparse vector in $\nspace(\bX)$ is the zero vector $\bbeta = \bzero$.
That is, $\nspace(\bX)\cap \sB_0[k] = \{\bzero\}$, which is similar to the standard NSP condition.

\index{Instance-optimal}
To fully understand the implications of NSP$'$ in sparse recovery, we now consider how to evaluate the performance of reconstruction algorithms when the true signal $\bbeta$ is not necessarily sparse. 
Let $\Delta:\real^n\rightarrow \real^p$ denote a specific recovery map (e.g., a decoder). We are particularly interested in error guarantees of the form
\begin{equation}\label{equation:nsp-ii-3}
\normtwo{\Delta(\bX \bbeta) - \bbeta} \leq C \frac{\sigma_k(\bbeta)_1}{\sqrt{k}},
\quad \text{for all $\bbeta$},
\end{equation}
\footnote{This is also related to the instance optimality; Definition~\ref{definition:inst_opt}.}
where $\sigma_k(\bbeta)_1$ denotes the \textit{$\ell_1$-distance of best $k$-term approximation}:
$
\sigma_k(\bbeta)_1 \triangleq \min_{\widehat{\bbeta} \in \sB_0[k]} \normonebig{\bbeta - \widehat{\bbeta}},
$
i.e., the distance between a vector $\bbeta$ and the ball $\sB_0[k]$ under the $\ell_1$-norm.
This bound ensures two important properties:
\begin{itemize}
\item \textbf{Exact recovery} of all possible $k$-sparse solutions (i.e., $\bbeta' = \Delta (\bX\bbeta')$ if $\bbeta'\in\sB_0[k]$).
\item \textbf{Robustness} to non-sparse signals, where the reconstruction error scales with how well $\bbeta$ can be approximated by a $k$-sparse vector.
\end{itemize}
Such guarantees are called \textit{instance-optimal} \citep{cohen2009compressed, foucart2013invitation}, because the error bound adapts to each individual signal $\bbeta$, rather than holding only for a restricted class (e.g., strictly sparse vectors). These are also referred to as \textit{uniform guarantees}, since they hold uniformly over all possible $\bbeta$.

The choice of norms in \eqref{equation:nsp-ii-3} is not unique. One could measure reconstruction error using other $\ell_s$-norms. 
However, the choice of norm affects both the strength of achievable guarantees and the corresponding nullspace conditions. For example, replacing the $\ell_2$-$\ell_1$ structure in \eqref{equation:nsp-ii-3} with an $\ell_2$-$\ell_2$ bound would generally require significantly more measurements---and, in fact, such a guarantee is impossible without a prohibitively large number of samples; see Theorem~\ref{theorem:optimal_l2_p1_cohen} or \citet{cohen2009compressed}. 
Thus, the form in \eqref{equation:nsp-ii-3} represents the strongest feasible uniform guarantee for stable recovery from linear measurements.
Moreover, a specific $\ell_1$-$\ell_1$ bound will be discussed in Section~\ref{section:stable_nsp}.

The following result, adapted from \citet{cohen2009compressed}, shows that any recovery algorithm satisfying \eqref{equation:nsp-ii-3} necessarily requires the measurement matrix $\bX$ to satisfy NSP$'$ of order $2k$.
\begin{lemma}[Sparse recovery under NSP$'$ \citep{cohen2009compressed}]\label{lemms:spar_rec_nspii}
Let $\bX \in\real^{n\times p}$ be a measurement matrix and $\Delta:\real^n\rightarrow\real^p$ an arbitrary recovery map. 
If the pair $(\bX, \Delta)$ satisfies the error bound   \eqref{equation:nsp-ii-3}, then $\bX$ satisfies the NSP$'$ of order $2k$.
\end{lemma}
\begin{proof}[of Lemma~\ref{lemms:spar_rec_nspii}]
Let $ \bn \in \nspace(\bX) $, and let $ \sS $ be the index set corresponding to the $ 2k $ largest (in magnitude) entries of $ \bn $. 
Partition $ \sS $ into two disjoint subsets $ \sS_0 $ and $ \sS_1 $,each of size $k $. 
Define $ \bbeta \triangleq \bn(\sS_1) + \bn(\comple{\sS})$ and $ \bbeta' \triangleq -\bn(\sS_0) $, so that $ \bn = \bbeta - \bbeta' $.~\footnote{
Given $\bbeta\in\real^p$, note that $\bbeta(\sS)$ indicates it is an $\real^p$ vector where the entries corresponds to $\comple{\sS}$ are all zeros, and the entries correspond to $\sS$ are the same as $\bbeta$; see Definition~\ref{definition:matlabnotation}.}
By construction, $ \bbeta' \in \sB_0[k] $. 
Therefore, the instance-optimal guarantee \eqref{equation:nsp-ii-3} implies exact recovery: $ \bbeta' = \Delta(\bX \bbeta') $. Moreover, since $ \bn \in \nspace(\bX) $, we have
\begin{equation}
\bX \bn = \bX \left( \bbeta - \bbeta' \right) = \bzero
\qquad\implies\qquad 
\bX \bbeta' = \bX \bbeta.
\end{equation}
It follows that $ \bbeta' = \Delta(\bX \bbeta) $. 
Now consider the reconstruction error:
\begin{equation}
\normtwo{\bn_\sS} \leq \normtwo{\bn} = \normtwo{\bbeta - \bbeta'} = \normtwo{\bbeta - \Delta (\bX \bbeta)} \leq C \frac{\sigma_k(\bbeta)_1}{\sqrt{k}} = \sqrt{2} C \frac{\normone{\bn_{\comple{\sS}}}}{\sqrt{2k}},
\end{equation}
where the last inequality follows from \eqref{equation:nsp-ii-3}.
This completes the proof.
\end{proof}

\subsection{Stable Nullspace Property}\label{section:stable_nsp}

We have introduced the NSP$'$ property, which provides an $\ell_1$-$\ell_2$ bound on vectors in the nullspace.
Now we turn to an analogous condition based on an $\ell_1$-$\ell_1$ bound.
In particular, we will show that basis pursuit is stable under a slightly strengthened version of the nullspace property, known as the stable nullspace property (see Theorem~\ref{theorem:stablensp}).

\begin{definition}[Stable nullspace property\index{Stable NSP}]\label{definition:stable_nsp}
A matrix $\bX \in \real^{n\times p}$ is said to satisfy the \textit{stable nullspace property (stable NSP)}  over a set $\sS \subset \{1,2,\ldots,p\}$ with constant $0 < \rho < 1$, denoted  SNSP($\sS, \rho$),  if
$$
\normone{\bbeta_\sS} \leq \rho \normone{\bbeta_{\comple{\sS}}}, 
\quad \text{for all } \bbeta \in \nspace(\bX).
$$
It is said to satisfy the \textit{stable NSP of order $k$} with constant $0 < \rho < 1$, denoted  SNSP($k, \rho$), if it satisfies the stable NSP with constant $0 < \rho < 1$ for every subset $\sS \subset \{1,2,\ldots,p\}$ with $\abs{\sS} \leq k$.
\end{definition}

Note again that the equality can be achieved in the definition of the stable NSP, while the equality is not achieved in the definition of the standard NSP.

The following result  shows that, under the stable NSP relative to a set $\sS$, the $\ell_1$-distance between a vector $\bbeta \in \real^p$  and another vector $\balpha \in \real^p$ that yields the same measurements ($\bX\balpha = \bX\bbeta$) is controlled by the difference between their $\ell_1$-norms.
This fact plays a key role in establishing error bounds for basis pursuit under the stable NSP (see Theorem~\ref{theorem:stablensp}).

\begin{theoremHigh}[Standard consequence of stable NSP]\label{theorem:suff_nec_snsp}
A matrix $\bX \in \real^{n\times p}$ satisfies the stable NSP with constant $0 < \rho < 1$ over $\sS$ if and only if
\begin{equation}\label{equation:suff_nec_snsp}
\normone{\balpha - \bbeta} \leq \frac{1+\rho}{1-\rho} (\normone{\balpha} - \normone{\bbeta} + 2 \normone{\bbeta_{\comple{\sS}}})
\end{equation}
for all vectors $\balpha, \bbeta \in \real^p$ with $\bX\balpha = \bX\bbeta$.
\end{theoremHigh}
Observe that if $\sS$ consists of the indices of the $k$  largest (in magnitude) entries of  $\bbeta$, then $\normone{\bbeta_{\comple{\sS}}} = \sigma_k(\bbeta)_1$, the best $k$-term $\ell_1$ approximation error.
Consequently, if $\bX$ satisfies the stable NSP  of order $k$ with constant $0<\rho<1$, it holds that 
\begin{equation}\label{equation:suff_nec_snsp2}
\normone{\balpha - \bbeta} \leq \frac{1+\rho}{1-\rho} \big(\normone{\balpha} - \normone{\bbeta} + 2 \sigma_k(\bbeta)_1\big),
\quad \text{for all  $\balpha, \bbeta$ with $\bX\balpha = \bX\bbeta$.}
\end{equation}

Before proving Theorem~\ref{theorem:suff_nec_snsp}, we establish a useful auxiliary inequality.
\begin{lemma}\label{lemma:stabnsp_lem}
For any set $\sS \subset \{1,2,\ldots,p\}$ and any vectors $\bbeta, \balpha \in \real^p$, the following inequality holds:
$$
\normone{(\bbeta - \balpha)_{\comple{\sS}}} \leq \normone{\balpha} - \normone{\bbeta} + \normone{(\bbeta - \balpha)_{{\sS}}} + 2 \normone{\bbeta_{\comple{\sS}}}.
$$
\end{lemma}
\begin{proof}[of Lemma~\ref{lemma:stabnsp_lem}]
By the triangle inequality, we have 
\begin{align*}
\normone{\bbeta} 
&= \normone{\bbeta_{\comple{\sS}}} + \normone{\bbeta_\sS}\leq \normone{\bbeta_{\comple{\sS}}} + \normone{(\bbeta - \balpha)_\sS} + \normone{\balpha_\sS}; \\
\normone{(\bbeta - \balpha)_{\comple{\sS}}} 
&\leq \normone{\bbeta_{\comple{\sS}}} + \normone{\balpha_{\comple{\sS}}}.
\end{align*}
Taking the summation over the two inequalities yields the desired result.
\end{proof}
\begin{proof}[of Theorem~\ref{theorem:suff_nec_snsp}]
\textbf{$\Leftarrow$.} Suppose that the matrix $\bX$ satisfies \eqref{equation:suff_nec_snsp} for all vectors $\balpha, \bbeta \in \real^p$ with $\bX\balpha = \bX\bbeta$. 
Given a vector $\bu \in \nspace(\bX)$, since $\bX\bu(\comple{\sS}) = -\bX\bu(\sS)$, we can invoke \eqref{equation:suff_nec_snsp} with $\bbeta \triangleq -\bu(\sS)$ and $\balpha \triangleq \bu(\comple{\sS})$. 
Since $\bbeta(\comple{\sS})=\bzero$, we have 
$$
\begin{aligned}
&\normone{\bu} \leq \frac{1+\rho}{1-\rho} (\normone{\bu_{\comple{\sS}}} - \normone{\bu_\sS})
\qquad\implies\qquad 
\normone{\bu_\sS} \leq \rho\normone{\bu_{\comple{\sS}}}.
\end{aligned}
$$
This proves the stable NSP with constant $0 < \rho < 1$ over $\sS$.

\paragraph{$\Rightarrow$.}
Conversely, suppose that the matrix $\bX$ satisfies the stable NSP with constant $0 < \rho < 1$ over $\sS$. 
For $\balpha, \bbeta \in \real^p$ with $\bX\balpha = \bX\bbeta$, since $\bu \triangleq \balpha - \bbeta \in \nspace(\bX)$, the stable NSP yields 
$\normone{\bu_\sS} \leq \rho \normone{\bu_{\comple{\sS}}}$.
Applying Lemma~\ref{lemma:stabnsp_lem} to $\bbeta$ and $\balpha=\bbeta+\bu$, we obtain
$$
\normone{\bu_{\comple{\sS}}} \leq \normone{\balpha} - \normone{\bbeta} + \normone{\bu_\sS} + 2\normone{\bbeta_{\comple{\sS}}}.
$$
Substituting $\normone{\bu_\sS} \leq \rho \normone{\bu_{\comple{\sS}}}$ into this inequality yields
$$
\normone{\bu_{\comple{\sS}}} \leq \normone{\balpha} - \normone{\bbeta} + \rho \normone{\bu_{\comple{\sS}}} + 2\normone{\bbeta_{\comple{\sS}}}.
$$
Since $\rho < 1$, we can rearrange to get
$$
\normone{\bu_{\comple{\sS}}} \leq \frac{1}{1-\rho} (\normone{\balpha} - \normone{\bbeta} + 2\normone{\bbeta_{\comple{\sS}}}).
$$
Using $\normone{\bu_\sS} \leq \rho \normone{\bu_{\comple{\sS}}}$ once again, we conclude
$$
\normone{\bu} = \normone{\bu_{\comple{\sS}}} + \normone{\bu_\sS} \leq (1+\rho)\normone{\bu_{\comple{\sS}}} \leq \frac{1+\rho}{1-\rho} (\normone{\balpha} - \normone{\bbeta} + 2\normone{\bbeta_{\comple{\sS}}}).
$$
This completes the proof.
\end{proof}

\subsection{Robust Nullspace Property}
In realistic settings, it is impossible to measure a signal $\bbeta \in \real^p$ with infinite precision. 
This means that the observed measurement vector $\by \in \real^n$ only approximates the ideal measurement $\bX\bbeta \in \real^n$, satisfying
$$
\norm{\bX\bbeta - \by} \leq \epsilon
$$
for some $\epsilon \geq 0$ and  some norm $\norm{\cdot}$ on $\real^n$ (typically the $\ell_2$-norm). 
In such cases, a reconstruction algorithm should produce an estimate $\widehat{\bbeta} \in \real^p$ whose distance to the original vector $\bbeta^* \in \real^p$ is controlled by the measurement error $\epsilon \geq 0$. 
This desirable behavior is commonly referred to as \textit{robustness} to measurement error. 
We will show that if the standard basis pursuit problem \eqref{opt:p1} is replaced by the  convex optimization program \eqref{opt:p1_epsilon} (p.~\pageref{opt:p1_epsilon}):
\begin{equation}\label{opt:p1epsilon_rnsp}
\text{(P$_{1,\epsilon}$)}
\qquad\min_{\bbeta \in \real^p} \normone{\bbeta} \quad \text{s.t.} \quad \norm{\bX\bbeta - \by} \leq \epsilon,
\end{equation}
then robust recovery is guaranteed under a further strengthening of the nullspace property---namely, the robust nullspace property; see Theorem~\ref{theorem:robust_nsp} for further results.

\begin{definition}[$\ell_s$-robust nullspace property\index{Robust NSP}]\label{definition:ells_rob_nsp}
Let $s \geq 1$.
A matrix $\bX \in \real^{n\times p}$ is said to satisfy the \textit{$\ell_s$-robust nullspace property ($\ell_s$-robust NSP)} of order $k$ (with respect to the norm $\norm{\cdot}$ on $\real^n$) with constants $0 < \rho < 1$ and $\tau > 0$  over a set $\sS \subseteq \{1,2,\ldots,p\}$, denoted  RNSP($\ell_s,\sS, \rho, \tau$), if
\begin{subequations}
\begin{equation}\label{equation:ells_rob_nsp_e1}
\norms{\bbeta_\sS} \leq \frac{\rho}{k^{(1-\frac{1}{s})}} \normone{\bbeta_{\comple{\sS}}} + \tau \norm{\bX\bbeta}, \quad \text{for all } \bbeta \in \real^p.
\end{equation}
When $s=2$, $\rho=C$, and $\tau=0$, this reduces to the NSP$'$ (Definition~\ref{definition:nsp_ii}).
In the special case  $s=1$, the matrix $\bX \in \real^{n\times p}$ is simply said to satisfy the \textit{robust NSP (with respect to $\norm{\cdot}$) with constants $0 < \rho < 1$ and $\tau > 0$ over a set $\sS \subseteq \{1,2,\ldots,p\}$}, denoted  RNSP($\sS, \rho, \tau$), if
\begin{equation}\label{equation:ells_rob_nsp_e2}
\normone{\bbeta_\sS} \leq \rho \normone{\bbeta_{\comple{\sS}}} + \tau \norm{\bX\bbeta}, \quad \text{for all } \bbeta \in \real^p.
\end{equation}
It is said to satisfy the \textit{$\ell_s$-robust NSP} or \textit{robust NSP of order $k$} with constants $0 < \rho < 1$ and $\tau > 0$ if it satisfies the {$\ell_s$-robust NSP} or the robust NSP with constants $\rho, \tau$ relative to \textbf{any} set $\sS \subset \{1,2,\ldots,p\}$ with $\abs{\sS} \leq k$, respectively.
\end{subequations}
\end{definition}

Note that the robust NSP generalizes the stable NSP (Definition~\ref{definition:stable_nsp}): when $\bX\bbeta=\bzero$, the term $\tau\norm{\bX\bbeta} $ vanishes, and the robust NSP reduces exactly to the stable NSP.

\begin{remark}
For $1 \leq t \leq s$, the norm inequality $\normt{\bbeta_\sS}  \leq k^{(\frac{1}{t} - \frac{1}{s})} \norms{\bbeta_\sS}$  (see Problem~\ref{prob:cauch_sc_gen}) implies that the $\ell_s$-robust NSP with constants $0 < \rho < 1$ and $\tau > 0$ entails, for any $\sS \subset \{1,2,\ldots,p\}$ with $\abs{\sS}  \leq k$,
\begin{equation}\label{equation:ells_rob_nsp}
\normt{\bbeta_\sS} \leq \frac{\rho}{k^{(1-\frac{1}{t})}} \normone{\bbeta_{\comple{\sS}}} + \tau k^{(\frac{1}{t} - \frac{1}{s})} \norm{\bX\bbeta}, \quad \text{for all } \bbeta \in \real^p.
\end{equation}
Thus, for $1 \leq t \leq s$, the $\ell_s$-robust NSP implies the $\ell_t$-robust NSP with  the same constants,  up to a factor of  $k^{(\frac{1}{t} - \frac{1}{s})}$ in the measurement-dependent term. 
In particular, the $\ell_s$-robust NSP is a stronger condition than the $\ell_1$-robust NSP, since it controls the signal on $\sS$ in a stronger norm (e.g., $\ell_2$ or $\ell_\infty$) while still depending only on the $\ell_1$-norm of the tail and the measurement residual.

Note that the definition above does not require $\bbeta$ to lie in the nullspace of $\nspace(\bX)$. 
In fact, if $\bbeta \in \nspace(\bX)$, then the term $\tau \norm{\bX\bbeta}$ vanishes.
In this case, the robust NSP reduces exactly to the stable NSP given in Definition~\ref{definition:stable_nsp}. Thus, the robust NSP indeed generalizes the stable NSP.
\end{remark}

Similar to Theorem~\ref{theorem:suff_nec_snsp}, we now establish a stronger equivalence that holds for any fixed index set  $\sS\subseteq\{1,2,\ldots,p\}$.
\begin{theoremHigh}[Standard consequence  of robust NSP]\label{theorem:suff_nec_rnsp}
A matrix $\bX \in \real^{n\times p}$ satisfies the robust NSP with constants $0 < \rho < 1$ and $\tau > 0$ relative to a set $\sS$ if and only if
\begin{equation}\label{equation:suff_nec_rnsp}
\normone{\balpha - \bbeta} \leq \frac{1+\rho}{1-\rho} 
\big(\normone{\balpha} - \normone{\bbeta} + 2\normone{\bbeta_{\comple{\sS}}}\big) + \frac{2\tau}{1-\rho} \norm{\bX(\balpha - \bbeta)}
\end{equation}
for all vectors $\balpha, \bbeta \in \real^p$.
\end{theoremHigh}
\begin{proof}[of Theorem~\ref{theorem:suff_nec_rnsp}]
\textbf{$\Leftarrow$.}  The proof follows similarly from that of Theorem~\ref{theorem:suff_nec_snsp}. 
Assume that inequality~\eqref{equation:suff_nec_rnsp} holds for all $\balpha, \bbeta \in \real^p$. For any $\bu \in \real^p$, invoking \eqref{equation:suff_nec_rnsp} with   $\bbeta \triangleq -\bu(\sS)$ and $\balpha \triangleq \bu(\comple{\sS})$ yields
$$
\normone{\bu} \leq \frac{1+\rho}{1-\rho} (\normone{\bu_{\comple{\sS}}} - \normone{\bu_\sS}) + \frac{2\tau}{1-\rho} \norm{\bX\bu}
\qquad \implies\qquad 
\normone{\bu_\sS} \leq \rho \normone{\bu_{\comple{\sS}}} + \tau \norm{\bX\bu}.
$$
This proves the robust NSP with constants $0 < \rho < 1$ and $\tau > 0$ over $\sS$.

\paragraph{$\Rightarrow$.} 
Conversely, suppose that the matrix $\bX$ satisfies the robust NSP with constants $0 < \rho < 1$ and $\tau > 0$ over $\sS$. 
Let $\balpha, \bbeta \in \real^p$ be arbitrary,  and set $\bu \triangleq \balpha - \bbeta$. 
Then $\bu\in\real^p$, and the robust NSP and Lemma~\ref{lemma:stabnsp_lem} give
$$
\begin{aligned}
\normone{\bu_\sS} 
&\leq \rho \normone{\bu_{\comple{\sS}}} + \tau \norm{\bX\bu};\\
\normone{\bu_{\comple{\sS}}} 
&\leq \normone{\balpha} - \normone{\bbeta} +
\normone{\bu_{\sS}}+ 2\normone{\bbeta_{\comple{\sS}}}.
\end{aligned}
$$
Replacing the term $\normone{\bu_\sS} $ in the second inequality with the first inequality yields
$$
\normone{\bu_{\comple{\sS}}} \leq \frac{1}{1-\rho} (\normone{\balpha} - \normone{\bbeta} + 2\normone{\bbeta_{\comple{\sS}}} + \tau \norm{\bX\bu}).
$$
Using the robust NSP inequality  once again, we derive
$$
\begin{aligned}
	\normone{\bu} &= \normone{\bu_{\comple{\sS}}} + \normone{\bu_\sS} \leq (1+\rho)\normone{\bu_{\comple{\sS}}} + \tau \norm{\bX\bu} \\
	&\leq \frac{1+\rho}{1-\rho} (\normone{\balpha} - \normone{\bbeta} + 2\normone{\bbeta_{\comple{\sS}}}) + \frac{2\tau}{1-\rho} \norm{\bX\bu}.
\end{aligned}
$$
This completes the proof.
\end{proof}

Similarly, for the $\ell_s$-robust NSP, we obtain the following standard consequence.
This result strengthens the previous error bound by replacing the $\ell_1$-error estimate with an $\ell_t$-error estimate for $s \geq t \geq 1$. 

\begin{theoremHigh}[Standard consequence of $\ell_s$-robust NSP]\label{theorem:lsrob_nsp_gen}
Let $1 \leq t \leq s$. 
Suppose that the matrix $\bX \in \real^{n\times p}$ satisfies the $\ell_s$-robust NSP of order $k$ with constants $0 < \rho < 1$ and $\tau > 0$. Then, for all $\balpha, \bbeta \in \real^p$,
\begin{equation}\label{equation:ell2_suff_nec_rnsp}
\normt{\balpha - \bbeta} \leq \frac{C}{k^{(1-\frac{1}{t})}} (\normone{\balpha} - \normone{\bbeta} + 2\sigma_k(\bbeta)_1) + D k^{(\frac{1}{t} - \frac{1}{s})} \norm{\bX(\balpha - \bbeta)},
\end{equation}
where $C \triangleq (1+\rho)^2/(1-\rho)$ and $D \triangleq \tau(3+\rho)/(1-\rho)$.
\end{theoremHigh}
\begin{proof}[of Theorem~\ref{theorem:lsrob_nsp_gen}]
By \eqref{equation:ells_rob_nsp}, the $\ell_s$-robust NSP implies the $\ell_1$-robust and $\ell_t$-robust NSP ($t \leq s$) in the forms:
\begin{equation}\label{equation:lsrob_nsp_gen_e1}
	\normone{\bu_\sS} \leq \rho \normone{\bu_{\comple{\sS}}} + \tau k^{(1-\frac{1}{s})} \norm{\bX\bu},
\end{equation}
\begin{equation}\label{equation:lsrob_nsp_gen_e2}
	\normt{\bu_\sS} \leq \frac{\rho}{k^{(1-\frac{1}{t})}} \normone{\bu_{\comple{\sS}}} + \tau k^{(\frac{1}{t} - \frac{1}{s})} \norm{\bX\bu},
\end{equation}
for all $\bu \in \real^p$ and all $\sS \subset \{1,2,\ldots,p\}$ with $\abs{\sS} \leq k$. Thus, in view of \eqref{equation:lsrob_nsp_gen_e1}, applying Theorem~\ref{theorem:suff_nec_rnsp} with $\sS$ chosen as an index set of $k$ largest (in modulus) entries of $\balpha - \bbeta$ leads to 
\begin{equation}\label{equation:lsrob_nsp_gen_e3}
	\normone{\balpha - \bbeta} \leq \frac{1+\rho}{1-\rho} (\normone{\balpha} - \normone{\bbeta} + 2\sigma_k(\bbeta)_1) + \frac{2\tau}{1-\rho} k^{(1-\frac{1}{s})} \norm{\bX(\balpha - \bbeta)}.
\end{equation}
Then, choosing $\sS$ as an index set of $k$ largest (in modulus) entries of $\balpha-\bbeta$, Problem~\ref{prob:lsdist_spar} shows that
$$
\normt{\balpha - \bbeta} 
= \normt{(\balpha - \bbeta)_{\comple{\sS}}} + \normt{(\balpha - \bbeta)_\sS} 
\leq \frac{1}{k^{(1-\frac{1}{t})}} \normone{\balpha - \bbeta} + \normt{(\balpha - \bbeta)_\sS}.
$$
Now substitute the bounds from~\eqref{equation:lsrob_nsp_gen_e2} with $\bu=\balpha-\bbeta$  and~\eqref{equation:lsrob_nsp_gen_e3} into this inequality:
$$
\begin{aligned}
\normt{\balpha - \bbeta} 
&\stackrel{\eqref{equation:lsrob_nsp_gen_e2}}{\leq } 
\frac{1}{k^{(1-\frac{1}{t})}} \normone{\balpha - \bbeta} + \frac{\rho}{k^{(1-\frac{1}{t})}} \normone{(\balpha - \bbeta)_{\comple{\sS}}} + \tau k^{(\frac{1}{t} - \frac{1}{s})} \norm{\bX(\balpha - \bbeta)}\\
&\leq
\frac{1+\rho}{k^{(1-\frac{1}{t})}} \normone{\balpha - \bbeta}  + \tau k^{(\frac{1}{t} - \frac{1}{s})} \norm{\bX(\balpha - \bbeta)}\\
&\stackrel{\eqref{equation:lsrob_nsp_gen_e3}}{\leq }
C (\normone{\balpha} - \normone{\bbeta} + 2\sigma_k(\bbeta)_1) 
+ D\norm{\bX(\balpha - \bbeta)}.
\end{aligned}
$$
This completes the proof.
\end{proof}

\section{Spark, Mutual Coherence, and Uniqueness of Sparse Recovery}

A full row rank underdetermined linear system 
$$
\bX_{n\times p}\bbeta_{p\times 1} = \by_{n\times 1},
\quad\text{with}\quad  n < p,
$$
has infinitely many solutions. 
Suppose we seek the sparsest solution---that is, the one with the fewest nonzero entries. 
Can such a solution ever be unique? And how can we find the sparsest solution? This section briefly addresses the first question (see also Chapter~\ref{chapter:recovery} for further discussion); the second will be covered in Chapters~\ref{chapter:algouni} and~\ref{chapter:spar_recov}.
\subsection{Spark}
If we wish to be able to recover all sparse  signals $\bbeta$ from the measurements $\bX \bbeta$, 
then for any two distinct vectors $\bbeta, \bbeta' \in \sB_0[k] = \{ \bbeta \mid \normzero{\bbeta} \leq k \}$, we must have $\bX \bbeta \neq \bX \bbeta'$.
Otherwise, the measurements  $\by$ would not allow us to distinguish between $\bbeta$ and  $\bbeta'$. 
More formally,  if $\bX \bbeta = \bX \bbeta'$, then $\bX (\bbeta - \bbeta') = \bzero$, and the difference satisfies $\bbeta - \bbeta' \in \sB_0[2k]$. 
Therefore, $\bX$ uniquely represents all $\bbeta \in \sB_0[k]$ if and only if its nullspace $\nspace(\bX)$ contains no (nonzero) vector from $\sB_0[2k]$. 

This condition can be characterized in several equivalent ways. One of the most common is via the spark of the matrix $\bX$, a key concept introduced in \cite{donoho2003optimally} that plays a central role in analyzing the uniqueness of sparse solutions.

\begin{definition}[Spark \citep{gorodnitsky1997sparse, donoho2003optimally}\index{Spark}]
Given a matrix $\bX\in\real^{n\times p}$, its \textit{spark}, denoted $\sigma=\spark(\bX)$,  is defined as the smallest possible number such that there exists a subgroup of $\sigma$ columns from $\bX$ that are linearly dependent.
That is,
$$
\sigma=\min\{s  : \text{ there exists a subset of $s$ columns of $\bX$ that are linearly dependent}\}.
$$
It always holds that 
$$
\sigma = \spark(\bX) \in[2, \min\{n,p\}+1].
$$
\end{definition}

The spark is closely related to the \textit{Kruskal rank}, denoted  $\rank_k(\bX)$ or $k_{\bX}$, which is defined as the largest integer $r$ such that every set of $r$ columns of $\bX$ is linearly independent (see, e.g., \citet{lu2021numerical} for further details).
Recall that the (standard) rank of a matrix is the maximum number of linearly independent columns in $\bX$   (Definition~\ref{definition:rank}). 
While the definitions of rank and Kruskal rank (or spark) appear similar, computing the spark is significantly more difficult: it requires a combinatorial search over all possible subsets of columns of $\bX$.
By definition, any nonzero vector $\bbeta$ satisfying $\bX\bbeta=\bzero$ must have at least $\spark(\bX)$ nonzero entries. In other words,
\begin{equation}\label{equation:pro_spar1}
\normzero{\bbeta}\geq \spark(\bX)
\qquad \text{if} \qquad  
\bX\bbeta = \bzero, \text{ and } \bbeta\neq \bzero,
\end{equation}
since fewer than $\spark(\bX)$ columns cannot form a nontrivial linear combination equal to zero.

This property leads to a simple and powerful criterion for the uniqueness of sparse solutions to underdetermined systems.
\begin{theoremHigh}[Uniqueness under spark {\citep{gorodnitsky1997sparse, donoho2003optimally, bruckstein2009sparse}}\index{Unique recovery}]\label{theorem:opt_spar_l0}
If the linear system $\bX\bbeta = \by$ admits a solution $\bbeta$ satisfying $\normzero{\bbeta} < \spark(\bX)/2$, then this solution is necessarily the (unique) sparsest solution.
\end{theoremHigh}
\begin{proof}[of Theorem~\ref{theorem:opt_spar_l0}]
Let  $\widetilde{\bbeta}$ be any other solution such that $\bX\widetilde{\bbeta} = \by$. 
Then $\bX(\bbeta-\widetildebbeta)=\bzero$, so $\bbeta - \widetilde{\bbeta}$ lies in the nullspace of  $\bX$.
By \eqref{equation:pro_spar1}, 
\begin{equation}
\normzero{\bbeta} + \normzerobig{\widetilde{\bbeta}} \geq \normzerobig{\bbeta - \widetilde{\bbeta}} \geq \spark(\bX).
\end{equation}
Since we have a solution satisfying $\normzero{\bbeta} < \spark(\bX) / 2$, we conclude that any alternative solution $\widetilde{\bbeta}$ necessarily has more than $\spark(\bX) / 2$ nonzeros, proving that $\bbeta$ is the unique sparsest solution.
\end{proof}

The theorem can be put in a different way as follows.
\begin{corollary}\label{corollary:opt_spar_l0}
For any vector $ \by \in \real^n $, there exists at most one solution $ \bbeta \in \sB_0[k] $ such that $ \by = \bX\bbeta $ if and only if $ \spark(\bX) > 2k $.
\end{corollary}
\begin{proof}[of Corollary~\ref{corollary:opt_spar_l0}]
We first assume that, for every $ \by \in \real^n $, there is at most one solution $ \bbeta \in \sB_0[k] $ satisfying $ \by = \bX \bbeta $. Suppose, for contradiction, that $ \spark(\bX) \leq 2k $. 
By definition of spark, this implies that there exists a set of at most  $ 2k $ columns of $\bX$ that are linearly dependent.
Consequently, there exists a nonzero vector $ \bn \in \nspace(\bX) $ such that $ \bn \in \sB_0[2k] $. 
In this case, since $ \bn \in \sB_0[2k] $ we can write $ \bn = \bbeta - \bbeta' $, where $ \bbeta, \bbeta' \in \sB_0[k] $, whence we have that $ \bX (\bbeta - \bbeta') = \bzero $ and hence $ \bX \bbeta = \bX \bbeta' $. 
But this contradicts the assumption there exists at most one solution $ \bbeta \in \sB_0[k] $ such that $ \by = \bX \bbeta $. 
Hence, we must have that $ \spark(\bX) > 2k $.

Conversely, suppose $ \spark(\bX) > 2k $. Assume that for some $ \by\in\real^n $, there exist two vectors $ \bbeta, \bbeta' \in \sB_0[k] $ such that $ \by = \bX \bbeta = \bX \bbeta' $. 
Then  $ \bX (\bbeta - \bbeta') = \bzero $, so $ \bn \triangleq \bbeta - \bbeta' \in\nspace(\bx)$. 
Moreover, since both $\bbeta$ and $\bbeta'$ have at most $k$ nonzeros, their difference satisfies $\normzero{\bn}\leq 2k$.
However, because $ \spark(\bX) > 2k $, 
the nullspace of $\bX$ contains no nonzero vector with $\ell_0$-norm $\leq 2k$.
Therefore, $\bn=\bzero$, which implies $\bbeta=\bbeta'$. This proves uniqueness, completing the proof.
\end{proof}

Clearly, the spark of a matrix is highly informative: larger values of $\spark(\bX)$ guarantee uniqueness for sparser solutions. How large can the spark be? By definition,
$ 2 \leq \spark(\bX) \leq \min\{n,p\} + 1$.
For example, if $\bX\in\real^{n\times p}$ has entries drawn independently from a continuous distribution (e.g., Gaussian), then with probability one, every subset of $\min\{n,p\}$ columns is linearly independent. In particular, when $n<p$, we have $\spark(\bX)=n+1$ almost surely, since no more than $n$ columns can be linearly independent in $\real^n$. Thus, uniqueness is guaranteed for any solution with fewer than $(n+1)/2$ nonzeros.

\subsection{Mutual Coherence}\label{section:mutual_cohere}
While the spark, the standard NSP, and its robust variants all provide theoretical guarantees for the recovery of sparse signals (see Section~\ref{section:spar_ana_nsp} for details), 
verifying that a general matrix $\bX$ satisfies any of these properties is computationally intractable.
Specifically, each condition requires examining on the order of  $\binom{p}{k}$ submatrices, leading to combinatorial complexity. 
For example, Theorem~\ref{theorem:opt_spar_l0} shows that for any sparsest solution $\bbeta$ to the linear system $\bX\bbeta=\by$, it holds that $\normzero{\bbeta} \leq \spark(\bX)/2$.
However, computing the spark of a matrix is at least as difficult  as solving the $\ell_0$-minimization problem in \eqref{opt:p0} (p.~\pageref{opt:p0}).
In many practical settings, it is therefore preferable to use easily computable properties of $\bX$ that still yield meaningful recovery guarantees. One such property is the \textit{mutual coherence} of the matrix.

Mutual coherence provides a simple yet effective measure of how suitable a measurement matrix is for sparse recovery. In general, the smaller the mutual coherence, the better the performance of recovery algorithms.

\begin{definition}[Mutual coherence \citep{mallat1993matching, donoho2001uncertainty, bruckstein2009sparse}\index{Mutual coherence}]\label{definition:mutua_cohere}
The \textit{mutual coherence} (also known as \textit{pairwise coherence} or simply \textit{coherence}) of a given matrix $\bX\in\real^{n\times p}$ is defined as the largest absolute value of the normalized inner product between any two distinct columns of $\bX$. Denoting the $i$-th column in $\bX$ by $\bx_i$, mutual coherence is given by
\begin{equation}
\mu(\bX) = \max_{1 \leq i, j \leq p, i \neq j}
\frac{\abs{\innerproduct{\bx_i, \bx_j}}}{\normtwo{\bx_i} \normtwo{\bx_j}}
\in[0,1].
\end{equation}
\end{definition}

Mutual coherence characterizes the degree of linear dependence among the columns of $\bX$. It quantifies the maximum similarity between any pair of distinct columns: lower coherence implies that the columns are closer to being orthogonal. This is advantageous because it suggests that each column contributes unique information---a key requirement for reliable sparse representation and recovery.

For an orthogonal matrix  (Definition~\ref{definition:orthogn_mat}), the columns are mutually orthonormal, so the mutual coherence is exactly zero. In contrast, for underdetermined systems where $p>n$ (i.e., more columns than rows), $\mu(\bX)$ must be strictly positive. In such cases, we aim to minimize $\mu(\bX)$ to emulate the favorable properties of orthogonal matrices as closely as possible.

To illustrate the extremes, consider now a matrix $\bX\in\real^{n\times p}$ with $n\geq p$. In this case, $\mu(\bX)=0$ if and only if the columns of $\bX$ are orthonormal. In particular, when $\bX$ is square ($n=p$), this occurs precisely when $\bX$ is an orthogonal matrix.

In compressed sensing, however, we are primarily interested in the underdetermined setting where $n<p$. Here, perfect orthogonality among all columns is impossible, and there are fundamental limits on how small the coherence can be. Roughly speaking, small coherence implies that submatrices formed by a moderate number of columns are well-conditioned, which is crucial for stable sparse recovery.
It can be shown that the mutual coherence of any matrix $\bX\in\real^{n\times p}$ with $p>n$ always satisfies 
$$
\mu(\bX) \in \left[\sqrt{\frac{p-n}{n(p-1)}}, 1\right].
$$
The lower bound is known as the \textit{Welch bound} \citep{rosenfeld1997praise, strohmer2003grassmannian, welch2003lower}; a proof is provided in Theorem~\ref{theorem:low_mutucohe}. 
Note that when $p \gg n$, this bound simplifies approximately to  $\mu(\bX) \geq 1/\sqrt{n}$.

Equivalently, mutual coherence reflects---on an entry-wise level---how close the Gram matrix $\bG=\bX^\top\bX$ is to the $p$-dimensional identity matrix $\bI$.

\begin{lemma}[Interpretation of mutual coherence w.r.t. spectral norm]\label{lemma:int_coh_spec}
Let $\bX\in\real^{n\times p}$ be a matrix with unit-norm columns. 
Then 
\begin{equation}\label{equation:mutual_cohe_identity}
\normtwo{\bX^\top\bX - \bI} \leq \max_{i \in \{1,2,\ldots,p\}}\sum_{j \in \{1,2,\ldots,p\} \setminus \{i\}} \abs{\innerproduct{\bx_i, \bx_j}} ,
\end{equation}
where $\normtwo{\cdot}$ denotes the spectral norm of a matrix.
\end{lemma}
\begin{proof}[of Lemma~\ref{lemma:int_coh_spec}]
Let $\bH\triangleq \bX^\top\bX-\bI$. 
By  the Gershgorin circle theorem (Theorem~\ref{theorem:eign_disc}), every eigenvalue $\lambda$ of $\bH$ satisfies
$
\abs{\lambda - h_{ii}} \leq \sum_{j \neq i} \abs{h_{ij}}.
$
Since the diagonal entries of $\bH$ are zero  ($h_{ii} = 0$), this simplifies to
$$
\abs{\lambda} \leq \sum_{j \neq i} \abs{\innerproduct{ \bx_i, \bx_j}}.
$$
Taking the maximum over all $i \in  \{1, 2, \ldots, p\}$, we obtain
$$
\normtwo{\bH} = \max_{\lambda \in \Lambda(\bH)} \abs{\lambda} \leq \max_{i \in \{1, 2, \ldots, p\}} \sum_{j \neq i} \abs{\innerproduct{\bx_i, \bx_j}}.
$$
where $\Lambda(\bH)$ denotes the spectrum of $\bH$.
This completes the proof.
\end{proof}

A stronger condition than mutual coherence is the \textit{mutual incoherence condition}. 
This condition is defined with respect to any subset $\sS\subset\{1,2,\ldots,p\}$.
\begin{definition}[Mutual incoherence\index{Mutual incoherence}]\label{definition:mutual_incohe}
The \textit{mutual incoherence} (also known as the \textit{irrepresentable condition}) for a matrix $\bX\in\real^{n\times p}$ holds for a subset  $\sS\subset\{1,2,\ldots,p\}$  if there exists some $\gamma>0$ such that
\begin{equation}\label{equation:incoherence}
\max_{i\in\comple{\sS}} \normone{(\bX_\sS^\top \bX_\sS)^{-1} \bX_\sS^\top \bx_i} 
\leq 1 - \gamma.
\end{equation}
\end{definition}

To interpret this condition, observe that the submatrix $\bX_\sS \in \real^{n \times k}$ (where $\abs{\sS}=k$) consists of the columns of $\bX$ indexed by the support set $\sS$. For each index $i$ in the complement $\comple{\sS}$, the vector $(\bX_\sS^\top \bX_\sS)^{-1} \bX_\sS^\top \bx_i$ contains the coefficients obtained by regressing $\bx_i$ onto the columns of  $\bX_\sS$. The $\ell_1$-norm of this vector quantifies how well $\bx_i$ can be represented as a linear combination of the columns in $\bX_\sS$.
In the ideal scenario, every column $\bx_i$ with $i\in\comple{\sS}$ would be orthogonal to all columns in  $\bX_\sS$. In that case, the regression coefficients would all be zero, and we would have $\gamma=1$. Of course, in high-dimensional settings where $p\gg n$, such perfect orthogonality is impossible. Nevertheless, we may still hope for a form of ``near orthogonality," which is precisely what the mutual incoherence condition captures when $\gamma$ is close to 1.

It is worth noting that the term ``mutual incoherence" is sometimes used informally to refer to small pairwise correlations between columns---i.e., a small mutual coherence $\mu(\bX)$. However, the formal mutual incoherence condition in Definition~\ref{definition:mutual_incohe} is strictly stronger, as it involves collective behavior of groups of columns rather than just pairwise interactions.

Mutual coherence can sometimes be related to other fundamental properties such as spark, the NSP, and the restricted isometry property (RIP, as  a hindsight, see Definition~\ref{definition:rip22} for reference). For instance, spark and mutual coherence are connected via the Gershgorin circletheorem (Theorem~\ref{theorem:eign_disc}), as shown in the following proposition.

\begin{proposition}[Spark $\&$ mutual coherence]\label{proposition:spark_mutual}
For any matrix $\bX \in \real^{n\times p}$, it follows that
\begin{equation}
\textup{spark}(\bX) \geq 1 + \frac{1}{\mu(\bX)}.
\end{equation}
\end{proposition}
\begin{proof}[of Proposition~\ref{proposition:spark_mutual}]
First, normalize the columns of $\bX$ to unit  $\ell_2$-norm, obtaining a new matrix $\widetilde{\bX}$. 
This operation preserves both the spark and the mutual coherence. Let $\bG = \widetilde{\bX}^\top \widetilde{\bX}$ be the corresponding Gram matrix. Its entries satisfy:
$$
\{ g_{kk} = 1,\; \forall\, 1 \leq k \leq p \} \quad \text{and} \quad \{ \abs{g_{ij}} \leq \mu (\bX),\;  \forall\, 1 \leq i, j \leq p, \; i \neq j \}.
$$
Now consider any $l \times l$ principal submatrix of $\bG$, corresponding to a subset of $l$ columns from $\widetilde{\bX}$ (and computing their sub-Gram matrix). From the Gershgorin circle theorem (Theorem~\ref{theorem:eign_disc} or Theorem~\ref{theorem:pd_diag_domi}), if this minor is diagonally dominant (i.e., if $\sum_{j \neq i} \abs{g_{ij}} < \abs{g_{ii}}$ for every $i$), then this submatrix of $\bG$ is positive definite, and so those $l$ columns from $\widetilde{\bX}$ are linearly independent satisfying:
\begin{itemize}
\item $g_{ii} = 1$, $1 \leq i \leq l$;
\item $\abs{g_{ij}} \leq \mu(\bX)$, $1 \leq i, j \leq l$, $i \neq j$.
\end{itemize}
Since $l < \spark(\bX)$, the spark condition and the positive definiteness implies $(l - 1) \mu(\bX) < 1$ or, equivalently, $l < 1 + 1/\mu(\bX)$.
Therefore, we can choose $l$ such that $l < 1 + 1/\mu(\bX) \leq l+1 \leq \spark(\bX)$, yielding $\spark(\bX) \geq 1 + 1/\mu(\bX)$.
\end{proof}

By combining Theorem~\ref{theorem:opt_spar_l0} or Corollary~\ref{corollary:opt_spar_l0} with Proposition~\ref{proposition:spark_mutual}, we obtain the following condition on $\bX$ that guarantees uniqueness of sparse solutions.

\begin{theoremHigh}[Uniqueness under mutual coherence \citep{donoho2001uncertainty, donoho2003optimally}\index{Unique recovery}]\label{theorem:uniq_beta0_mutuocoh}
Consider the linear system $\bX\bbeta = \by$. 
If this system admits a solution $\bbeta$ satisfying
$$
\normzero{\bbeta} < \frac{1}{2} \left( 1 + \frac{1}{\mu(\bX)} \right),
$$
then  the solution $\bbeta$  is necessarily the (unique) sparsest solution.
Equivalently, if
$$
k < \frac{1}{2}\left(1 + \frac{1}{\mu(\bX)}\right),
$$
then for each measurement vector $\by \in \real^n$, there exists at most one signal $\bbeta \in \sB_0[k]$ such that $\by = \bX\bbeta$.
\end{theoremHigh}

This theorem once again partially  answers the following fundamental questions:
\begin{itemize}
\item  When can we guarantee uniqueness of the sparsest solution?
\item Can a candidate solution be verified as globally optimal for \eqref{opt:p0}?
\end{itemize} 
We have seen that any solution sufficiently sparse---specifically, sparser than the bounds above---is guaranteed to be unique among all solutions of comparable sparsity. Consequently, such a solution must also be the global minimizer of \eqref{opt:p0}. This shows that seeking sparse solutions is not merely a heuristic; under appropriate conditions, it leads to a well-posed inverse problem with strong theoretical guarantees.

Comparing Theorems~\ref{theorem:opt_spar_l0} and~\ref{theorem:uniq_beta0_mutuocoh}, we observe that they are structurally similar but rely on different assumptions. Theorem~\ref{theorem:opt_spar_l0}, which is based on the spark, is sharp and generally much stronger than Theorem~\ref{theorem:uniq_beta0_mutuocoh}, which depends only on mutual coherence---a quantity that provides a lower bound on the spark.

Indeed, mutual coherence cannot be smaller than approximately $1/\sqrt{n}$ (as implied by the Welch bound, i.e., $\mu(\bX) \in \left[\sqrt{\frac{p-n}{n(p-1)}}, 1\right]$), so the sparsity level guaranteed by Theorem~\ref{theorem:uniq_beta0_mutuocoh} never exceeds roughly $\sqrt{n}/2$.
That is, Theorem~\ref{theorem:uniq_beta0_mutuocoh}, together with the Welch bound, provides an upper bound on the level of sparsity $k$ that guarantees uniqueness using coherence: $k = \mathcalO\left(\sqrt{n}\right)$.
In contrast, the spark can be as large as $n$ (assuming $p>n$), in which case Theorem~\ref{theorem:opt_spar_l0} allows sparsity levels up to $n/2$---a significantly more permissive bound.

We now show that mutual coherence also implies the NSP condition.
\begin{proposition}[Coherence$\implies$NSP]\label{proposition:cohe2nspi}
Let $\bX\in\real^{n\times p}$.
Suppose that for some integer $k\in\{1,2,\ldots,p\}$, the mutual coherence satisfies the bound $\mu(\bX) < \frac{1}{3k}$.
Then $\bX$ satisfies the NSP  of order $k$.~\footnote{Consequently, by Corollary~\ref{corollary:exa_ell1_snsp}, $\ell_1$-minimization \eqref{opt:p1} recovers every $k$-sparse vector exactly.}
\end{proposition}
\begin{proof}[of Proposition~\ref{proposition:cohe2nspi}]
Without loss of generality, we may assume that all columns of $\bX$ are normalized to unit $\ell_2$-norm, i.e., $\normtwo{\bx_i} = 1$ for all $i$ (since mutual coherence and the NSP are invariant under column rescaling).
To simplify the analysis, suppose more generally that $\mu(\bX) < \frac{\delta}{k}$ for some $\delta>0$, and determine the smallest $\delta$ that guarantees the NSP. We will see that $\delta=1/3$ suffices.

Let $\sS\subseteq\{1,2,\ldots,p\}$ be any subset with $\abs{\sS}=k$, and let  $\bbeta \in \sC[\sS] \setminus \{\bzero\}$. 
To verify the NSP $\nspace(\bX)\cap \sC[\sS] = \{\bzero\}$, it suffices to show that $\normtwo{\bX\bbeta}^2 > 0$, and so we begin with the lower bound
$$
\normtwo{\bX\bbeta}^2 
= \normtwo{\bX_{\sS}\bbeta_{\sS} + \bX_{\comple{\sS}}\bbeta_{\comple{\sS}}}^2
\geq \normtwo{\bX_{\sS}\bbeta_{\sS}}^2 + 2\bbeta_{\sS}^\top\bX_{\sS}^\top\bX_{\comple{\sS}}\bbeta_{\comple{\sS}}.
$$
We now bound the cross term. Using the definition of mutual coherence and the triangle inequality,
\begin{align*}
2\abs{\bbeta_{\sS}^\top\bX_{\sS}^\top\bX_{\comple{\sS}}\bbeta_{\comple{\sS}}}
&\leq 2\absBig{\sum_{i \in \sS}\sum_{j \in \comple{\sS}} \abs{\beta_i}\abs{\beta_j}\innerproduct{\bx_i, \bx_j}} \\
&\stackrel{(i)}{\leq} 2\normone{\bbeta_{\sS}}\normone{\bbeta_{\comple{\sS}}}\mu(\bX) 
\stackrel{(ii)}{\leq} \frac{2\delta\normone{\bbeta_{\sS}}^2}{k} 
\stackrel{(iii)}{\leq} 2\delta\normtwo{\bbeta_{\sS}}^2,
\end{align*}
where the inequality (i) follows from the definition  of the mutual coherence; the inequality (ii) exploits the assumed bound on $\mu(\bX)$ combined with the fact that $\bbeta \in \sC[\sS]$; and the inequality (iii) follows from the fact that  $\normone{\bbeta_{\sS}} \leq \sqrt{k}\normtwo{\bbeta_{\sS}}$ (Exercise~\ref{exercise:cauch_sc_l1l2}) since the cardinality of $\sS$ is at most $k$.
Combining the two inequalities, we have
\begin{equation}\label{equation:cohe2nspi}
\normtwo{\bX\bbeta}^2 \geq \normtwo{\bX_{\sS}\bbeta_{\sS}}^2 - 2\delta\normtwo{\bbeta_{\sS}}^2.
\end{equation}
It remains to lower-bound $\normtwo{\bX_{\sS}\bbeta_{\sS}}^2$.
By Lemma~\ref{lemma:int_coh_spec}, we have 
$$
\normtwo{\bX_{\sS}^\top\bX_{\sS} - \bI} 
\leq \max_{i \in \sS}\sum_{j \in \sS \setminus \{i\}} 
\abs{\innerproduct{\bx_i, \bx_j}} 
\leq k\frac{\delta}{k} = \delta.
$$
This shows  all eigenvalues of $\bX_{\sS}^\top\bX_{\sS} - \bI$ lie in the interval $[-\delta, \delta]$ and $ -\delta\bI \preceq \bX_{\sS}^\top\bX_{\sS} - \bI \preceq \delta \bI$.
Therefore, $\normtwo{\bX_{\sS}\bbeta_{\sS}}^2 \geq (1-\delta)\normtwo{\bbeta_{\sS}}^2$, and combined with the bound \eqref{equation:cohe2nspi}, we conclude that $\normtwo{\bX\bbeta}^2 > (1-3\delta)\normtwo{\bbeta_{\sS}}^2$, so that $\delta = 1/3$ ensures strict positivity, which proves that $\bX$ satisfies the NSP of order $k$.
\end{proof}

Moreover, another direct application of the Gershgorin circle theorem (Theorem~\ref{theorem:eign_disc}) links mutual coherence to the RIP (see Definition~\ref{definition:rip22} for reference).
\begin{proposition}[Coherence$\implies$RIP]\label{proposition:cohe2rip}
If $\bX\in\real^{n\times p}$ has unit-norm columns and coherence $\mu = \mu(\bX)$, then $\bX$ satisfies the RIP of order $k$ with $\delta = (k - 1)\mu$ for all $k < 1/\mu$.
\end{proposition}
The proof follows the same reasoning as in Proposition~\ref{proposition:spark_mutual}, using diagonal dominance of sub-Gram matrices.

These results underscore the importance of small mutual coherence in compressed sensing  matrix design. Coherence bounds have been studied extensively for both deterministic and random constructions. For instance, there exist explicit deterministic matrices of size $n \times n^2$ that achieve the Welch lower bound $\mu(\bX) = 1/\sqrt{n}$, such as Gabor frames generated from the Alltop sequence \citep{herman2009high} and equiangular tight frames \citep{tropp2008conditioning}. Using such matrices, stable recovery of $k$-sparse signals is possible with $n = \mathcalO(k^2 \ln p)$ measurements.

In contrast, for random matrices with i.i.d. entries having zero mean and finite variance, the coherence concentrates around $\mu(\bX) = \sqrt{(2 \ln p)/n}$
in the asymptotic regime $(n,p\rightarrow \infty)$ \citep{cai2011limiting, cand2008near, donoho2006most}. Substituting this into the coherence-based sparsity bound yields
$k=\mathcalO(\sqrt{n/(\ln p)})$,
which---after rearrangement---again gives $n=\mathcalO(k^2\ln p)$. 
However, we will not delve into these details in this book,
while the RIP related results will be introduced in Chapters~\ref{chapter:spar_gauss} and \ref{chapter:ensur_rips}.

\subsection{$\ell_1$-Coherence}\label{section:ell1_cohere}
A more general notion than mutual coherence is the $\ell_1$-coherence, also known as the Babel function.
\begin{definition}[$\ell_1$-coherence \citep{foucart2013invitation}\index{$\ell_1$-coherence}]\label{definition:l1_coherence}
Let $\bX \in \real^{n \times p}$. The \textit{$\ell_1$-coherence} function (or \textit{Babel function}) $\mu_1$ of the matrix $\bX$ is defined for $1 \leq k \leq p - 1$ by
$$
\mu_1(\bX,k) 
=\max_{\substack{\sS \subset \{1,\ldots,p\} \\ \abs{\sS} = k}}
\max_{j\notin\sS}
\left\{ \sum_{i \in \sS} 
\frac{\abs{\innerproduct{\bx_i, \bx_j}}}{\normtwo{\bx_i}\normtwo{\bx_j}}
\right\}.
$$
\begin{subequations}
It is straightforward to verify that, for $1 \leq k \leq p - 1$,
\begin{equation}\label{equation:l1cohere_pro1}
\mu(\bX) \leq \mu_1(\bX,k) \leq k \cdot \mu(\bX),
\end{equation}
and more generally that, for $1 \leq k, t \leq p - 1$ with $k + t \leq p - 1$,
\begin{equation}\label{equation:l1cohere_pro2}
\max \{\mu_1(\bX,k), \mu_1(\bX,t)\} \leq \mu_1(\bX,k + t) \leq \mu_1(\bX,k) + \mu_1(\bX,t).
\end{equation}
\end{subequations}
\end{definition}

Mutual coherence captures the worst-case pairwise correlation between any two distinct columns of $\bX$, measuring only the largest absolute normalized inner product. In contrast, $\ell_1$-coherence quantifies the total interference that any set of $k$ columns can exert on another column. Specifically, it measures the maximum sum of absolute normalized inner products between a fixed column $\bx_j$ and any collection of $k$ other columns.

Note that by \eqref{equation:l1cohere_pro1}, we have $\mu_1(\bX,1)=\mu_1(\bX)$, so mutual coherence is simply the special case of $\ell_1$-coherence with $k=1$.

Both mutual coherence and the $\ell_1$-coherence function are \textit{invariant} under left multiplication by an orthogonal matrix $\bQ\in\real^{n\times n}$. 
Indeed, the columns of $\bQ\bX$ are  $\bQ\bx_1,\bQ\bx_2, \ldots, \bQ\bx_p$, and since orthogonal transformations preserve inner products,
$$
\innerproduct{\bQ\bx_i, \bQ\bx_j} = \innerproduct{\bx_i, \bx_j}. 
$$ 
Moreover, by the Cauchy--Schwarz inequality $\abs{\innerproduct{\bx_i, \bx_j}} \leq \normtwo{\bx_i} \normtwo{\bx_j}$, so the normalized inner product is always at most 1 in magnitude. 
Consequently, the mutual coherence satisfies
$\mu(\bX) \leq 1$.

As previously discussed, $\mu(\bX) = 0$ if and only if the columns of $\bX$ form an orthonormal system. In particular, when $\bX$ is square ($n=p$), this occurs precisely when $\bX$ is an orthogonal matrix.
However, in compressive sensing we are primarily interested in the underdetermined regime, where $n<p$. In this setting, perfect orthogonality among all columns is impossible, and there are fundamental limits on how small the mutual coherence can be. These limitations will be explored in Section~\ref{section:bd_cohe_l1}. 
For now, we simply note that small coherence implies that submatrices formed by a moderate number of columns are well-conditioned, which is essential for stable sparse recovery.

We have already shown in Proposition~\ref{proposition:cohe2rip} that mutual coherence implies the RIP (Definition~\ref{definition:rip22}) under mild conditions. An analogous result holds for $\ell_1$-coherence.
\begin{theoremHigh}[$\ell_1$-coherence and RIP]\label{theorem:l1co_and_rip}
Let $\bX \in \real^{n \times p}$ be a matrix with $\ell_2$-normalized columns (i.e., each column has unit norm), and let $1 \leq k \leq p$. 
Then, for every vector $\bbeta \in \real^p$ with $\bbeta\in\sB_0[k]$ (i.e., $\normzero{\bbeta}\leq k$),
$$
\big(1 - \mu_1(\bX,k - 1)\big) \normtwo{\bbeta}^2 
\leq \normtwo{\bX\bbeta}^2 \leq \big(1 + \mu_1(\bX,k - 1)\big) \normtwo{\bbeta}^2.
$$
Equivalently, for every index set $\sS \subseteq \{1,2,\ldots,p\}$ with $\abs{\sS} \leq k$, all eigenvalues of the Gram matrix $\bX_{\sS}^\top \bX_{\sS}$ lie in the interval 
$$
\lambda(\bX_{\sS}^\top \bX_{\sS})\in
\left[1 - \mu_1(\bX,k - 1), 1 + \mu_1(\bX,k - 1)\right].
$$ 
In particular, if $\mu_1(\bX,k - 1) < 1$,
then  $\bX$ satisfies the RIP of order $k$ with constant $\delta=\mu_1(\bX,k - 1)$ (see Definition~\ref{definition:rip22}), and  $\bX_{\sS}^\top \bX_{\sS}$ is nonsingular for all such $\sS$.
\end{theoremHigh}
\begin{proof}[of Theorem~\ref{theorem:l1co_and_rip}]
Let $\sS \subseteq \{1,2,\ldots,p\}$ with $\abs{\sS} \leq k$. 
Since $\bX_{\sS}^\top \bX_{\sS}$ is symmetric and positive semidefinite, it admits an orthonormal basis of eigenvectors with real, nonnegative eigenvalues  (Theorem~\ref{theorem:eigen_charac}). 
Denote its smallest and largest eigenvalues by  $\lambda_{\min}$ and $\lambda_{\max}$, respectively. 
For any $\bbeta \in \real^p$ supported on $\sS$, we have  $\bX\bbeta = \bX_{\sS} \bbeta_{\sS}$, thus  the maximum of
$
\normtwo{\bX\bbeta}^2 =\normtwo{\bX_{\sS} \bbeta_{\sS}}^2
$
over the set $\{\bbeta \in \real^p, \; \supp(\bbeta) \subseteq \sS, \; \normtwo{\bbeta} = 1\}$ is $\lambda_{\max}$ and that its minimum is $\lambda_{\min}$.~\footnote{Because $\max_{\bbeta\neq 0}\frac{\bbeta^\top\bA\bbeta}{\bbeta^\top\bbeta} = \lambda_{\max}$
and $\min_{\bbeta\neq 0}\frac{\bbeta^\top\bA\bbeta}{\bbeta^\top\bbeta} = \lambda_{\min}$
.}
Hence, the stated bounds on $\normtwo{\bX\bbeta}$ are equivalent to the eigenvalue inclusion claimed in the theorem.

Now, because the columns of $\bX$ are normalized, the diagonal entries of $\bX_{\sS}^\top\bX_{\sS}$ are all equal to 1.
Applying the Gershgorin circle theorem (Theorem~\ref{theorem:eign_disc}), each eigenvalue lies within at least one disc centered at 1 with radius
$$
r_i 
= \sum_{j \in \sS, j \neq i} \abs{(\bX_{\sS}^\top \bX_{\sS})_{ij}} 
= \sum_{j \in \sS, j \neq i} \abs{\innerproduct{\bx_j, \bx_i}} \leq \mu_1(\bX,k - 1), \quad i \in \sS,
$$
where the inequality follows from the definition of $\mu_1(\bX,k - 1)$ (since $\abs{\sS\setminus\{i\}}\leq k-1$).
Because all eigenvalues are real, they must lie in the interval $\left[1 - \mu_1(\bX,k - 1), 1 + \mu_1(\bX,k - 1)\right]$, as claimed.
\end{proof}

The above result shows that if $\mu_1(\bX,k - 1) < 1$, then $\bX_{\sS}^\top \bX_{\sS}$ is nonsingular for all such $\sS$. Expand writing here...

The result above shows that if $\mu_1(\bX,k - 1) < 1$, then the Gram matrix $\bX_{\sS}^\top \bX_{\sS}$ is nonsingular for every index set $\sS$ with $\abs{\sS}\leq k$. 
This guarantees that any $k$ columns of $\bX$ are linearly independent---a crucial property for sparse recovery. The following corollary strengthens this observation by providing a condition under which any $2k$ columns remain linearly independent, which is often required in the analysis of algorithms like orthogonal matching pursuit or in establishing uniqueness of $k$-sparse solutions.
\begin{corollary}\label{corollary:l1co_and_rip}
Let $\bX \in \real^{n \times p}$ have $\ell_2$-normalized columns, and let $k \geq 1$ be an integer.
If
$$
\mu_1(\bX,k) + \mu_1(\bX,k - 1) < 1,
$$
then for every index set $\sS \subseteq \{1,2,\ldots,p\}$ with $\abs{\sS} \leq 2k$, the matrix $\bX_{\sS}^\top \bX_{\sS}$ is invertible. 
Consequently,  $\bX_{\sS}$ has full column rank. 
In particular, this conclusion holds whenever
$$
\mu(\bX) < \frac{1}{2k - 1}.
$$
\end{corollary}`
\begin{proof}[of Corollary~\ref{corollary:l1co_and_rip}]
By inequality \eqref{equation:l1cohere_pro2}, we have $\mu_1(\bX,2k - 1) \leq \mu_1(\bX,k) + \mu_1(\bX,k - 1) < 1$. 
Now consider any index set $\sS \subseteq \{1,2,\ldots,p\}$ with $\abs{\sS} \leq 2k$.   Applying Theorem~\ref{theorem:l1co_and_rip} with sparsity level $2k$, we obtain that the smallest eigenvalue of $\bX_{\sS}^\top \bX_{\sS}$ satisfies $\lambda_{\min} \geq 1 - \mu_1(\bX,2k - 1) > 0$. 
Hence, $\bX_{\sS}^\top \bX_{\sS}$ is positive definite and therefore invertible.
To see that $\bX_{\sS}$ has full column rank, suppose  $\bX_{\sS} \bz= \bzero$ for some vector $\bz\in\real^{\abs{\sS}}$, which implies $\bX_{\sS}^\top\bX_{\sS} \bz= \bzero$.
Since  $\bX_{\sS}^\top \bX_{\sS}$ is invertible, it follows that  $\bz=\bzero$. This proves the first statement. 
Finally, using inequality~\eqref{equation:l1cohere_pro1}, we bound $\mu_1(\bX,k) + \mu_1(\bX,k - 1) \leq (2k - 1) \mu(\bX) < 1$ if $\mu(\bX) < 1/(2k - 1)$.
\end{proof}

\subsection{Bounds for Mutual Coherence and $\ell_1$-Coherence}\label{section:bd_cohe_l1}

In this subsection, we establish lower bounds for both the mutual coherence and the $\ell_1$-coherence function of a matrix $\bX \in \real^{n \times p}$
with $n < p$.

Matrices that achieve these lower bounds share a distinctive structure: their columns form equiangular tight frames, defined below.

\begin{definition}[Equiangular system, and tight frame system\index{Equiangular system}\index{Tight frame system}]\label{definition:equiv_tight}
A matrix $\bX\in\real^{n\times p}$ with $\ell_2$-normalized columns is called \textit{equiangular} if there is a exists  $C \geq 0$ such that
$$
\abs{\innerproduct{\bx_i, \bx_j}}= C \quad \text{for all } i \neq j,\; i, j \in \{1,2,\ldots,p\}.
$$
A matrix $\bX\in\real^{n\times p}$ is called a \textit{tight frame} if there exists a constant $\lambda > 0$ such that any of the following equivalent conditions holds:
\begin{enumerate}[(i)]
\item $\normtwo{\balpha}^2 = \lambda \sum_{j=1}^{p} \innerproduct{\balpha, \bx_j}^2$ for all $\balpha \in \real^n$,
\item $\balpha = \lambda \sum_{j=1}^{p} \innerproduct{\balpha, \bx_j} \bx_j$ for all $\balpha \in \real^n$,
\item $\bX \bX^\top = \frac{1}{\lambda} \bI_n$.
\end{enumerate}
As the name suggests, a system of $\ell_2$-normalized vectors is called an equiangular tight frame if it satisfies both properties above.
\end{definition}

Such systems are precisely those that attain the lower bound stated next, known as the \textit{Welch bound}.
\begin{theoremHigh}[Welch bound: lower bound of mutual coherence \citep{rosenfeld1997praise, strohmer2003grassmannian, welch2003lower}\index{Welch bound}]\label{theorem:low_mutucohe}
Let  $\bX \in \real^{n \times p}$ with $p>n$.
Then its mutual coherence satisfies
\begin{equation}\label{equation:low_mutucohe}
\mu(\bX) \geq \sqrt{\frac{p - n}{n(p - 1)}}.
\end{equation}
Equality holds if and only if the columns of  $\bX$ form an equiangular tight frame.
Note that when $p \gg n$, the lower bound is approximately $\mu(\bX) \geq 1/\sqrt{n}$.
\end{theoremHigh}
\begin{proof}[of Theorem~\ref{theorem:low_mutucohe}]
Without loss of generality, we assume that $\bX$ have  $\ell_2$-normalized columns.
Define the Gram matrix $\bG \triangleq \bX^\top \bX \in \real^{p \times p}$  and the matrix $\bH \triangleq \bX \bX^\top \in \real^{n \times n}$. 
Since the columns of $\bX$ are  $\ell_2$-normalized, we have
$
\trace(\bG) = \sum_{i=1}^{p} \normtwo{\bx_i}^2 = p$.
Recall that the Frobenius inner product
$
\innerproduct{\bA, \bB}_F \triangleq \trace(\bA\bB^\top) = \sum_{i,j=1}^{n} a_{ij} b_{ij}
$
induces the  Frobenius norm $\normf{\cdot}$  (Definition~\ref{definition:frobernius-in-svd}). 
By the Cauchy--Schwarz inequality,
\begin{equation}\label{equation:cohbound_eq1}
\trace(\bH) \leq \normf{\bH} \normf{\bI_n} = \sqrt{n} \sqrt{\trace(\bH \bH^\top)}.
\end{equation}
Using the cyclic invariance of the trace, we compute
\begin{align*}
\trace(\bH \bH^\top) 
&= \trace\big((\bX \bX^\top)(\bX \bX^\top)\big) = \trace(\bX^\top \bX \bX^\top \bX) 
= \trace(\bG \bG^\top) 
= \sum_{i,j=1}^{p,p} \abs{\innerproduct{\bx_i, \bx_j}}^2 \\
&= \sum_{i=1}^{p} \normtwo{\bx_i}^2 + \sum_{i,j=1, i \neq j}^{p,p} \abs{\innerproduct{\bx_i, \bx_j }}^2 
= p + \sum_{i,j=1, i \neq j}^{p,p} \abs{\innerproduct{ \bx_i, \bx_j}}^2.
\end{align*}
By the cyclic invariance of the trace again, it follows that $\trace(\bG) = \trace(\bH)$, whence we have 
\begin{equation}\label{equation:cohbound_eq3}
p^2 \leq n \left( p + \sum_{i,j=1, i \neq j}^{p,p} \abs{\innerproduct{\bx_i, \bx_j}}^2 \right).
\end{equation}
By the definition of mutual coherence,
\begin{equation}\label{equation:cohbound_eq2}
\abs{\innerproduct{\bx_i, \bx_j}}\leq \mu(\bX), 
\quad 
\text{ for all } 1 \leq i \neq j \leq p,
\end{equation}
we obtain
$
p^2 \leq n(p + (p^2 - p)\mu(\bX)^2),
$
which rearranges to the claimed bound \eqref{equation:low_mutucohe}.
Finally, equality holds if and only if equality holds in both \eqref{equation:cohbound_eq1} and \eqref{equation:cohbound_eq2}. Equality in \eqref{equation:cohbound_eq1} implies $\bH=\lambda\bI_n$ for some $\lambda>0$, i.e., $\bX$ is a tight frame. Equality in \eqref{equation:cohbound_eq2} implies that all off-diagonal inner products have the same magnitude, i.e., $\bX$ is equiangular. Hence, equality occurs precisely when $\bX$ forms an equiangular tight frame.
\end{proof}

The Welch bound can be extended to the $\ell_1$-coherence function for small values of its argument.
\begin{theoremHigh}[Welch bound: lower bound of $\ell_1$-coherence]\label{theorem:low_l1coherence}
Let  $\bX \in \real^{n \times p}$ with $p>n$. 
Then its $\ell_1$-coherence function satisfies 
\begin{equation}\label{equation:low_l1coherence}
\mu_1(\bX,k) \geq k \sqrt{\frac{p - n}{n(p - 1)}}, \quad \text{for} \quad k < \sqrt{p - 1}.
\end{equation}
Equality holds if and only if the columns  of  $\bX$ form an equiangular tight frame.
\end{theoremHigh}

The proof relies on the following auxiliary result.
\begin{lemma}\label{lemma:low_l1coherence}
Let $s < \sqrt{m}$, and let $(\alpha_1, \alpha_2, \ldots, \alpha_m)$ be a nonnegative, nonincreasing sequence satisfying
$$
\alpha_1 \geq \alpha_2 \geq \ldots \geq \alpha_m \geq 0 \quad \text{and} \quad \alpha_1^2 + \alpha_2^2 + \ldots + \alpha_m^2 \geq \frac{m}{s^2}.
$$
Then
$$
\alpha_1 + \alpha_2 + \ldots + \alpha_s \geq 1,
$$
with equality if and only if $\alpha_1 = \alpha_2 = \ldots = \alpha_m = 1/s$.
\end{lemma}
\begin{proof}[of Lemma~\ref{lemma:low_l1coherence}]
It suffices to prove the contrapositive:
$$
\left\{
\begin{array}{l}
\alpha_1 \geq \alpha_2 \geq \ldots \geq \alpha_m \geq 0 \\
\alpha_1 + \alpha_2 + \ldots + \alpha_s \leq 1
\end{array}
\right\} \implies \alpha_1^2 + \alpha_2^2 + \ldots + \alpha_m^2 \leq \frac{m}{s^2},
$$
with equality only when all $\alpha_i = 1/s$. 
This is equivalent to maximizing the convex function
$$
f(\alpha_1, \alpha_2, \ldots, \alpha_m) \triangleq \alpha_1^2 + \alpha_2^2 + \ldots + \alpha_m^2
$$
over the convex set
$$
\sS \triangleq \{ (\alpha_1, \ldots, \alpha_m) \in \real^m \mid  \alpha_1 \geq \ldots \geq \alpha_m \geq 0 \text{ and } \alpha_1 + \ldots + \alpha_s \leq 1 \}.
$$
Since $f$ is convex and $\sS$ is a compact convex polytope, the maximum of $f$ over $\sS$ is attained at an extreme point (vertex) of $\sS$.
The vertices of $\sS$ arise when enough of the inequality constraints become equalities. Up to symmetry and ordering, the relevant cases are:
\begin{itemize}
\item if $\alpha_1 = \alpha_2 = \ldots = \alpha_m = 0$, then $f(\alpha_1, \alpha_2, \ldots, \alpha_m) = 0$;
\item if $\alpha_1 + \ldots + \alpha_s = 1$ and $\alpha_1 = \ldots = \alpha_{\ell} > \alpha_{\ell+1} = \ldots = \alpha_m = 0$ for $1 \leq \ell \leq s$, then $\alpha_1 = \ldots = \alpha_{\ell} = 1/\ell$, and consequently $f(\alpha_1, \alpha_2, \ldots, \alpha_m) = 1/\ell$;
\item if $\alpha_1 + \ldots + \alpha_s = 1$ and $\alpha_1 = \ldots = \alpha_{\ell} > \alpha_{\ell+1} = \ldots = \alpha_m = 0$ for $s < \ell \leq m$, then $\alpha_1 = \ldots = \alpha_{\ell} = 1/s$, and consequently $f(\alpha_1, \alpha_2, \ldots, \alpha_m) = \ell/s^2$.
\end{itemize}
Since $s < \sqrt{m}$, we have  
$$
\max_{(\alpha_1, \ldots, \alpha_m) \in \sS} f(\alpha_1, \ldots, \alpha_m) = \max \left\{ \max_{1 \leq \ell \leq k} \frac{1}{\ell}, \max_{s < \ell \leq m} \frac{\ell}{s^2} \right\} = \max \left\{ 1, \frac{m}{s^2} \right\} = \frac{m}{s^2},
$$
with equality only when $\ell = m$ and $\alpha_1 = \alpha_2 = \ldots = \alpha_m = 1/s$.
\end{proof}

\begin{proof}[of Theorem~\ref{theorem:low_l1coherence}]
Without loss of generality, we assume that $\bX$ have  $\ell_2$-normalized columns.
From inequality \eqref{equation:cohbound_eq3}, we have
$
\sum_{i,j=1, i \neq j}^{p,p} \abs{\innerproduct{\bx_i, \bx_j}}^2 \geq \frac{p^2}{n} - p = \frac{p(p-n)}{n}
$.
This implies
$$
\max_{i \in \{1,2,\ldots,p\}} \sum_{j=1, j \neq i}^{p} \abs{\innerproduct{\bx_i, \bx_j}}^2 
\geq \frac{1}{p} \sum_{i,j=1, i \neq j}^{p,p} \abs{\innerproduct{ \bx_i, \bx_j}}^2 \geq \frac{p-n}{n}.
$$
For an index $i^\dag \in \{1,2,\ldots,p\}$ that achieves this maximum, and consider the $(p-1)$-dimensional vector of absolute inner products: $\big\{\abs{\innerproduct{\bx_{i^\dag}, \bx_j}}\big\}_{j=1, j \neq i^\dag}^{p}$.
Reorder these values in nonincreasing order as $\beta_1 \geq \beta_2 \geq \ldots \geq \beta_{p-1} \geq 0$, so that
$$
\beta_1^2 + \beta_2^2 + \ldots + \beta_{p-1}^2 \geq \frac{p-n}{n}.
$$
Invoking Lemma~\ref{lemma:low_l1coherence} with $m = p-1$ and $s=k$, and defining $\alpha_{\ell} \triangleq (\frac{1}{k}\sqrt{n(p-1)/(p-n)}) \beta_{\ell}$ gives $\alpha_1 + \alpha_2 + \ldots + \alpha_k \geq 1$. 
By the definition of $\ell_1$-coherence, it follows that
$$
\mu_1(\bX,k) \geq \beta_1 + \beta_2 + \ldots + \beta_k \geq k \sqrt{\frac{p-n}{n(p-1)}},
$$
which yields the desired bound~\eqref{equation:low_l1coherence}.

Finally, suppose equality holds in~\eqref{equation:low_l1coherence}. Then all intermediate inequalities must be tight. In particular:
\begin{itemize}
\item Equality in~\eqref{equation:cohbound_eq3} implies that $\bX$ is a tight frame (as in the proof of Theorem~\ref{theorem:low_mutucohe}).
\item Equality in Lemma~\ref{lemma:low_l1coherence} implies that all off-diagonal inner products involving column $i^\dag$ have the same magnitude: 
$\abs{\innerproduct{\bx_{i^\dag}, \bx_j }} = \sqrt{(p-n)/(n(p-1))}$ for all $j \in \{1,2,\ldots,p\}$, $j \neq i^\dag$.
\end{itemize}
This ensures that the choice of $i^\dag$ must be arbitrary in $\{1,2,\ldots,p\}$, it follows that all pairwise inner products between distinct columns have identical magnitude. Hence, $\bX$ is equiangular.
Conversely, the proof that equiangular tight frames yield equality in \eqref{equation:low_l1coherence} follows easily from the above discussions.
\end{proof}

\section{Restricted Isometry Property (RIP)}
The NSP condition is both necessary and sufficient for establishing error guarantees of the form
\begin{equation}
\normtwo{\Delta(\bX\bbeta) - \bbeta} \leq C \frac{\sigma_k(\bbeta)_1}{\sqrt{k}}, 
\quad \text{with }\quad
\sigma_k(\bbeta)_s \triangleq \min_{\widehat{\bbeta} \in \sB_0[k]} \normsbig{\bbeta - \widehat{\bbeta}}.
\end{equation}
However, these guarantees assume noise-free measurements. 
In practical settings---where measurements are corrupted by noise (e.g., due to sensor error or quantization)---stronger conditions are required to ensure stable and robust recovery.
To address this, \citet{candes2005decoding} introduced the restricted isometry property (RIP), a fundamental condition in compressive sensing that we briefly touched upon in the discussion of mutual coherence.

\begin{definition}[Restricted isometry property (RIP)\index{Restricted isometry property}]\label{definition:rip22}
A matrix $\bX \in \real^{n\times p}$ is said to satisfy the \textit{restricted isometry property (RIP) of order $k$ ($k<p$) with constant $\delta_{k} \in [0,1)$}, denoted  RIP$(k, \delta_k)$, if for $k$-sparse vector $\bbeta \in \sB_0[k] \triangleq \{\balpha \mid  \normzero{\balpha} \leq k\}$, we have
\begin{subequations}
\begin{equation}\label{equation:rippro22_e1}
(1-\delta_{k}) \cdot \normtwo{\bbeta}^{2} \leq  \normtwo{\bX\bbeta}^{2} \leq (1+\delta_{k}) \cdot \normtwo{\bbeta}^{2}.
\end{equation}
Using the spectral norm identity $\normtwo{\bX^\top\bX-\bI}
=\mathop{\max}_{\bu\in \real^p: \norm{\bu}_2=1} \abs{\bu^\top(\bX^\top\bX-\bI)\bu}
= \mathop{\max}_{\bu\in \real^p: \norm{\bu}_2=1} \abs{\normtwo{\bX\bu}^2-1}
$ by \eqref{equation:spectral_norm_eq2}, condition~\eqref{equation:rippro22_e1} can be equivalently restated as
\begin{equation}\label{equation:rippro22_e12}
\normtwo{\bX_{\sS}^\top\bX_{\sS} - \bI_{\abs{\sS}}} \leq \delta_k,
\quad \text{for all }\sS\subseteq \{1,2,\ldots,p\} \text{ with }\abs{\sS}\leq k.,
\end{equation}
\end{subequations}
where $\bX_{\sS}\triangleq\bX[:,\sS]\in\real^{n\times \abs{\sS}}$ denotes the columns of $\bX$ from the index set $\sS$.

\begin{subequations}
The smallest constant $\delta_k$ for which RIP$(k, \delta_k)$ holds is called the \textit{restricted isometry constant (RIC)} of order $k$. It can be defined as
\begin{align}
\delta_k 
&= \min \left\{ \delta\mid  (1 - \delta) \normtwo{\bv}^2 \leq 
\normtwo{\bX_{\sS} \bv}^2 \leq (1 + \delta) \normtwo{\bv}^2,\; 
\forall\, \abs{\sS} \leq k, 
\forall\, \bv \in \real^{\abs{\sS}} \right\} 
\label{equation:def_ric22}
\\
&= 
\max_{\sS\subseteq \{1,2,\ldots,p\}, \abs{\sS}\leq k} \normtwo{\bX_{\sS}^\top\bX_{\sS}-\bI_{\abs{\sS}}}. 
\label{equation:def_ric33}
\end{align}
\end{subequations}
\end{definition}

Note that either expression in~\eqref{equation:def_ric22} or~\eqref{equation:def_ric33} implies that all singular values of    $\bX_{\sS}$  lie in the interval 
$$
\sigma(\bX_{\sS}) 
\in [\sqrt{1-\delta_k}, \sqrt{1+\delta_k}]
\qquad 
\text{with} \qquad
\delta_k\leq 1.
$$ 
(see Theorem~\ref{theorem:sigbd_nearortho}).
Equivalently, the eigenvalues of $\bX_{\sS}^\top\bX_{\sS}$ lie in
$$
\lambda(\bX_{\sS}^\top\bX_{\sS})
\in 
[1-\delta_k, 1+\delta_k]
\qquad 
\text{with} \qquad
\delta_k\leq 1.
$$
\footnote{The quantity $\frac{\bbeta^\top\bX^\top\bX\bbeta}{\normtwo{\bbeta}^2}$ is known as the \textit{Rayleigh quotient} of the matrix $\bX^\top\bX$, which always lies between its smallest and largest eigenvalues as  $\bbeta$ varies.}
The strict inequality $\delta_k<1$ further implies that the smallest eigenvalue of $\bX_{\sS}^\top\bX_{\sS}$ is {positive}, which means that  $\bX_{\sS}$ has full column rank for any index set $\sS$ with $\abs{\sS}\leq k$, by the eigenvalue characterization theorem of positive definite matrices (Theorem~\ref{theorem:eigen_charac}).
Consequently, $\bX_{\sS}^\top\bX_{\sS}$ is {nonsingular}.

This relationship between eigenvalues (or singular values) and the RIP constant is instrumental in deriving key inequalities used later.
\begin{remark}[IHT inequalities]\label{remark:rip_iht}
From~\eqref{equation:rippro22_e12}, we obtain
\begin{equation}
\normtwo{\bX^\top \bX \bbeta - \bbeta} 
=\normtwo{\bX_{\sS}^\top \bX_{\sS} \bbeta_{\sS} - \bbeta_{\sS}} 
\leq \normtwo{\bX_{\sS}^\top \bX_{\sS} - \bI_{\abs{\sS}}}\normtwo{\bbeta_{\sS}}
\leq \delta_{k} \normtwo{\bbeta},
\end{equation}
where $\supp(\bbeta)=\sS$.
This is a standard bound used in iterative algorithms such as \textit{iterative hard-thresholding (IHT)}.
Alternatively, applying the triangle inequality yields:
\begin{align*}
\normtwo{\bX^\top \bX \bbeta} &\leq \normtwo{\bbeta} + \normtwo{\bX^\top \bX \bbeta - \bbeta} \leq (1 + \delta_{k}) \normtwo{\bbeta};\\
\normtwo{\bX^\top \bX \bbeta} &\geq \normtwo{\bbeta} - \normtwo{\bX^\top \bX \bbeta - \bbeta} \geq (1 - \delta_{k}) \normtwo{\bbeta}.
\end{align*}
Thus, we arrive at the key inequality:
\begin{equation}
(1 - \delta_{k}) \normtwo{\bbeta} \leq \normtwo{\bX^\top \bX \bbeta} \leq (1 + \delta_{k}) \normtwo{\bbeta}.
\end{equation}
\end{remark}

Informally, a matrix $\bX$ is said to \textit{possess the RIP} if $\delta_k$ is small for sufficiently large sparsity level $k$. 
Note that for a vector $\bbeta$ with support $\sS$, we have $\bX\bbeta = \bX_{\sS}\bbeta_{\sS}$. 
Thus, the RIP essentially requires that every submatrix of $\bX$ formed by selecting up to $k$ columns is well-conditioned---that is, it acts approximately like an isometry on sparse vectors.
In fact, for sparse optimization problems, the RIP of order $2k$ is particularly relevant (see, e.g., Theorem~\ref{theorem:ripnspI}).
Indeed, condition~\eqref{equation:rippro22_e1} implies that $\normtwo{\bX(\bbeta_1-\bbeta_2)}>0$
for any two distinct $k$-sparse vectors $\bbeta_1, \bbeta_2\in \sB_0[k]$. 
Consequently, their measurement vectors are also distinct: $\bX\bbeta_1\neq \bX\bbeta_2$.

For a concrete example of sparse recovery, suppose that $\bbeta^*$ is $k$-sparse. Then:
\begin{itemize}
\item We can uniquely  recover $\bbeta^*\in\sB_0[k]$ from the noiseless measurements $\bX\bbeta^*=\by$ provided that any $2k$ columns of $\bX$ are linearly independent. 
This is because if $\bX(\bbeta_1-\bbeta_2)=\bzero$  for two $k$-sparse vectors $\bbeta_1$ and $\bbeta_2$, then $\bbeta_1=\bbeta_2$, ensuring uniqueness of the solution.
\item Conversely, if $\bX$ contains a set of $2k$ linearly dependent columns, then unique recovery of $\bbeta^*$ is not guaranteed. In this case, there exists a nonzero vector $\bn = \bbeta_1 - \bbeta_2$ (with $\normzero{\bn} \leq 2k$) such that $\bX\bn = \bzero$, meaning $\bX\bbeta_1 = \bX\bbeta_2$ for two distinct $k$-sparse vectors.

\end{itemize}
This example motivates the need for a lower bound in the RIP of order $2k$.
If this lower bound fails---i.e., if $(1-\delta_{2k})$ is too small or zero---then it becomes possible to find two distinct $k$-sparse vectors $\bbeta_1$ and $\bbeta_2$ with disjoint supports such that $\bX\bbeta_1=\bX\bbeta_2$.
The lower bound in the RIP ensures that distinct sparse signals are mapped to well-separated measurement vectors, preventing ambiguity in recovery.

%

However, the RIP can be somewhat restrictive, as it requires the distortion of vector norms to be bounded symmetrically in the form $(1 \pm \delta_k)$ for $\delta_k \in [0,1)$.
By definition, if the matrix $\bX$ is orthogonal (i.e., $\bX\bX^\top=\bX^\top\bX=\bI$), then $\delta_k = 0$, since $\normtwo{ \bX \bbeta} = \normtwo{\bbeta}$ for all $\bbeta$.
Consequently, a nonzero RIC $\delta_k$ can be interpreted as a measure of how far $\bX$ deviates from orthogonality, but only when considering sparse linear combinations involving at most $k$ columns.

\index{Rayleigh quotient}
Recall that the expression ${\bv^\top\bA\bv}/{\normtwo{\bv}^2}$ is known as the \textit{Rayleigh quotient} of a nonzero vector $\bv$ with respect to a  matrix $\bA\in\real^{p\times p}$. This quantity always lies between the  largest and smallest eigenvalues of $\bA$. 
The ratio of the largest to the smallest eigenvalue---known as the condition number---characterizes how well-conditioned a matrix is \citep{lu2021numerical, lu2025practical}; a smaller condition number indicates better conditioning. In this context, the RIP condition implies two important properties:
\begin{itemize}
\item \textit{Norm preservation controlled by  $\delta_k$.} The inequality $(1 - \delta_k) \normtwo{\bbeta}^2 \leq\normtwo{\bX\bbeta}^2 \leq (1 + \delta_k) \normtwo{\bbeta}^2$ guarantees that $ \bX $ approximately preserves the Euclidean norm of every $k$-sparse vector $ \bbeta $, with the deviation governed by $\delta_k$.
\item \textit{Well-conditioned submatrices.} 
When $\delta_k$ is not too large, every submatrix of $\bX$ formed by selecting any $k$ columns is well-conditioned. In other words, such submatrices are close to being orthonormal. This property is essential for ensuring the stability and accuracy of sparse recovery algorithms.
\end{itemize}

\paragrapharrow{Rescaled, reshuffled, and augmented properties of RIP?}
We recall from~\eqref{equation:nsp_chg} that the NSP is invariant under certain transformations of the measurement matrix---specifically, when measurements are rescaled, reshuffled, or augmented. (Consequently, $k$-sparse recovery via $\ell_1$-minimization remains unaffected by such operations; see Chapter~\ref{chapter:recovery}.) 
However, these same operations can deteriorate the RICs.
Reshuffling measurements corresponds to replacing the measurement matrix $\bX \in \real^{n\times p}$ with  $\bP\bX$, where $\bP \in \real^{n\times n}$ is a permutation matrix (a special case of an orthogonal matrix). 
Since orthogonal transformations preserve singular values---i.e.,  $\sigma_i(\bQ\bX) = \sigma_i(\bX)$ for any orthogonal matrix $\bQ \in \real^{n\times n}$---this operation leaves the RICs unchanged. 

Rescaling measurements---i.e., replacing $\bX \in \real^{n\times p}$ with  $\bD\bX$, where $\bD \in \real^{n\times n}$ is a diagonal matrix---can also increase the RIC. This holds even for scalar rescaling: replacing $\bX$ by $\eta\bX$ for some $\eta\in\real$.
For instance, if $\bX \in \real^{n\times p}$ satisfies  $\delta_k(\bX) < 1/3$, then the scaled matrix $2\bX$ has RIC $\delta_k(2\bX) =5/6$. 

Augmenting a measurement, on the other hand, corresponds to appending a row to $\bX$, which can also increase the RIC.
For example, let $\bX \in \real^{n\times p}$ satisfy $\delta_k(\bX) < 1$, and choose any $\delta > \delta_k(\bX)$. Construct a new matrix $\widetildebX$ by appending the row vector $[0, \ldots, 0 ,\sqrt{1 + \delta}]$. 
Consider the 1-sparse vector $\bbeta \triangleq [0, \ldots, 0,  1]^\top$. 
Then $\normtwobig{\widetildebX\bbeta}^2 \geq 1 + \delta$,  implies that $\delta_1(\widetildebX) \geq \delta$, and hence $\delta_k(\widetildebX) > \delta_k(\bX)$ (by Exercise~\ref{exercise:order_delta_rip}).

\begin{exercise}
Given $ \bX \in \real^{n \times p} $, let $ \mu >0$ and $ \nu>0 $ be the largest and smallest constants such that
$$
\mu\cdot \normtwo{\bbeta}^2 \leq \normtwo{\bX\bbeta}^2 \leq \nu \cdot \normtwo{\bbeta}^2,
\quad \text{for all $ k $-sparse vectors $ \bbeta \in \real^p$.}
$$
Find the scaling factor $ \eta > 0 $ that minimizes $ \delta_k(\eta\b) $.
\textit{Hint: Think about $ (\nu - \mu)/(\nu + \mu) $.}
\end{exercise}

It is also important to note that in our definition of the RIP, we assume symmetric bounds around 1 (i.e., of the form $1\pm\delta_k$). 
However, this symmetry is adopted for notational convenience only. In practice, one may consider asymmetric bounds, leading to a useful generalization of the RIP---particularly relevant in settings where the design matrix is not tightly controlled, such as in gene expression analysis. This generalization is captured by the concepts of restricted strong convexity (RSC) and restricted smoothness (RSS).
\begin{definition}[Restricted strong convexity/smoothness property\index{Restricted strong convexity}\index{Restricted strong smoothness}]\label{definition:res_scss_mat}
A matrix $\bX \in \real^{n\times p}$ is said to satisfy the \textit{$\mu$-restricted strong convexity (RSC)} and \textit{$\nu$-restricted smoothness (RSS) property of order $k$} if, for all $\bbeta \in \sB_0[k]$, 
\begin{equation}\label{equation:res_scss_mat}
\mu \cdot \normtwo{\bbeta}^{2} \leq  \normtwo{\bX\bbeta}^{2} \leq \nu \cdot \normtwo{\bbeta}^{2}.
\end{equation}
\end{definition}

The key difference between RIP and RSC/RSS lies in their constants. 
RIP requires the lower and upper bounds to be symmetric about 1, i.e., $1 \pm \delta_{k}$, whereas RSC/RSS imposes no such symmetry. 
Readers may notice a close connection between Definition~\ref{definition:res_scss_mat} and Definition~\ref{definition:res_scss_func}, where restricted strong convexity and smoothness are defined for general functions. Indeed, Definition~\ref{definition:res_scss_func} can be viewed as a generalization of Definition~\ref{definition:res_scss_mat} to arbitrary (possibly non-quadratic) functions \citep{jain2014iterative, lu2025practical}. For twice-differentiable functions, both definitions essentially constrain the restricted eigenvalues of the Hessian matrix.

Importantly, given any matrix $\bX$ satisfying the asymmetric bounds in~\eqref{equation:res_scss_mat}, we can always rescale it to meet the symmetric RIP condition in~\eqref{equation:rippro22_e1}. Specifically, define
$$
\widetildebX\triangleq 
\sqrt{2/(\nu + \mu)}\bX.
$$
Then $\widetildebX$ satisfies the RIP of order $k$ with constant $\delta_{k} = (\nu - \mu)/(\nu + \mu)$.
While we do not prove this explicitly here, it is straightforward to verify. Moreover, all theorems in this book that assume the RIP actually remain valid as long as some scaling of $\bX$ satisfies the RIP. Therefore, without loss of generality, we may restrict our attention to the symmetric formulation in~\eqref{equation:rippro22_e1}.

Another useful property of the RIP is its monotonicity with respect to sparsity level:
If $\bX$ satisfies the RIP of order $k$ with constant $\delta_{k}$, then for any $k' < k$ it automatically satisfies the RIP of order $k'$ with constant $\delta_{k'} \leq \delta_{k}$. 
\begin{exercise}[Nondecreasing of RICs]\label{exercise:order_delta_rip}
Let $\bX\in\real^{n\times p}$ satisfy the RIP condition. Show that 
\begin{equation}
	\delta_1 \leq \delta_2\leq \ldots\leq \delta_p.
\end{equation}
\end{exercise}

Furthermore, if $\bX$
satisfies the RIP of order $k$ with a sufficiently small constant, then it will also automatically satisfy the RIP of order $\gamma k$ for certain $\gamma$, albeit with a somewhat worse constant.

\begin{exercise}[RIP variants \citep{needell2009cosamp}]
Suppose that $\bX$ satisfies the RIP of order $k$ with constant $\delta_{k}$. Let $\gamma$ be a positive integer. 
Show that $\bX$ satisfies the RIP of order $k' = \gamma \lfloor \frac{k}{2} \rfloor$ with constant $\delta_{k'} < \gamma \cdot \delta_{k}$, where $\lfloor \cdot \rfloor$ denotes the floor function.
\end{exercise}

This result is immediate for $\gamma = 1, 2$, but for $\gamma \geq 3$ (and $k \geq 4$),
it enables us to extend RIP guarantees from order $k$ to significantly higher sparsity levels. Note, however, that for the resulting bound to be meaningful (e.g., $\delta_{k'}<1$), the original constant $\delta_{k}$ must be sufficiently small.

\subsection{Restricted Orthogonality Property (ROP)}

Note that inequality \eqref{equation:rippro22_e1} implies that for any index set $\sS\subseteq \{1,2,\ldots,p\}$ with  $\abs{\sS}<k$, 
the eigenvalues  of the Gram matrix $\bX_{\sS}^\top\bX_{\sS}$ lie in interval $[1-\delta_k, 1+\delta_k]$.
In particular, when $k=1$, the constant $\delta_1$ reflects the deviation of the column norms of  $\bX$ from unity: $\delta_1=0$ if and only if every column of $\bX$ has unit Euclidean norm.
A related concept that captures the non-orthogonality between distinct columns (or disjoint subsets of columns) is the \textit{restricted orthogonality property (ROP)}.

To motivate this notion, we begin with a simple but essential observation, which leads naturally to the definition of the \textit{restricted orthogonality constant (ROC)}.
\begin{proposition}[Standard consequence of RIP]\label{proposition:conse_rip}
Let $\balpha, \bbeta \in \real^p$ be vectors such that $\normzero{\balpha} \leq s$ and $\normzero{\bbeta} \leq t$. If $\supp(\balpha) \cap \supp(\bbeta) = \varnothing$, then
$$
\abs{\innerproduct{\bX\balpha, \bX\bbeta}} \leq \delta_{s+t} \normtwo{\balpha} \normtwo{\bbeta}.
$$
\end{proposition}
\begin{proof}[of Proposition~\ref{proposition:conse_rip}]
Let $\sS \triangleq  \supp(\balpha) \cup \supp(\bbeta)$, and denote by $\balpha_{\sS}, \bbeta_{\sS} \in \real^{\abs{\sS}}$ the restrictions of  $\balpha, \bbeta \in \real^p$ to the coordinates in $\sS$. 
Since the supports of $\balpha$ and $\bbeta$ are disjoint,  $\innerproduct{\balpha_{\sS}, \bbeta_{\sS}} = 0$. 
We then compute:
$$
\begin{aligned}
\abs{\innerproduct{\bX\balpha, \bX\bbeta}} 
&= \abs{\innerproduct{\bX_{\sS}\balpha_{\sS}, \bX_{\sS}\bbeta_{\sS} } - \innerproduct{\balpha_{\sS}, \bbeta_{\sS}}}
= \abs{\innerproduct{(\bX_{\sS}^\top\bX_{\sS} - \bI)\balpha_{\sS}, \bbeta_{\sS}}}\\
&\leq \normtwo{\bX_{\sS}^\top\bX_{\sS} - \bI} \normtwo{\balpha_{\sS}} \normtwo{\bbeta_{\sS}} \leq \delta_{s+t} \normtwo{\balpha} \normtwo{\bbeta},
\end{aligned}
$$
where the last inequality follows from~\eqref{equation:def_ric33} and the facts that $\normtwo{\balpha_{\sS}} = \normtwo{\balpha}$ and $\normtwo{\bbeta_{\sS}} = \normtwo{\bbeta}$.
\end{proof}


This result motivates a complementary condition to the RIP, known as the restricted orthogonality property.

\begin{definition}[Restricted orthogonality  property (ROP)\index{Restricted orthogonality  property}]\label{definition:rop22}
A matrix $\bX\in\real^{n\times p}$ is said to satisfy the \textit{restricted orthogonality property (ROP)} of order ($k,k'$)  with constant $\theta_{k,k'}$ ($k+k'\leq p$), denoted  ROP$(k,k', \theta_{k,k'})$, if for all \textbf{disjoint} vectors $\bbeta\in \sB_0[k]$ and $\bbeta'\in\sB_0[k']$ (i.e., $\innerproduct{\bbeta,\bbeta'}=0$), we have 
\begin{subequations}
\begin{equation}\label{equation:def_roc22_e1}
\abs{\innerproduct{\bX\bbeta, \bX\bbeta'}} \leq \theta_{k,k'}\cdot \normtwo{\bbeta}\normtwo{\bbeta'}.
\end{equation}
The smallest such constant $\theta_{k,k'}$ is called the \textit{restricted orthogonality constant (ROC)}. 
It can be defined as the infimum of all parameters $\theta$ such that ROP$(k,k', \theta_{k,k'})$ holds, i.e., 
\begin{equation}\label{equation:def_roc22}
\small
\begin{aligned}
&\theta_{k,k'} = \min \sS;\\
\sS &\triangleq\left\{ \theta \mid
\abs{\innerproduct{\bX_{\sI}\bv_1, \bX_{\sI'}\bv_2}} \leq \theta\cdot \normtwo{\bv_1}\normtwo{\bv_2},
\forall \abs{\sI}\leq k, \abs{\sI'}\leq k',
\forall \bv_1 \in \real^{\abs{\sI}}, \bv_2 \in \real^{\abs{\sI'}}, \sI\cap \sI'=\varnothing
\right\}.
\end{aligned}
\end{equation}
Equivalently, the ROC admits the explicit characterization
\begin{equation}\label{equation:roc_def2}
\theta_{k,k'} 
= \max \Big\{\normtwo{\bX_{\sI}^\top\bX_{\sI'}}, \sI \cap \sI' = \varnothing, 
\abs{\sI} \leq k, \abs{\sI'} \leq k'\Big\}.
\end{equation}
\end{subequations}
Note that $\theta_{k,k'}\equiv \theta_{k',k}$ due to symmetry.
\end{definition}

The RIP and ROP constants are closely related. To see this, consider the setting of Definition~\ref{definition:rop22} and assume without loss of generality that $\normtwo{\bbeta}=\normtwo{\bbeta'}=1$ (since both properties are scale-invariant). 
Since $\sI$ and $\sI'$ are disjoint, there exist $\balpha\triangleq\bbeta+\bbeta'\in\sB_0[k+k']$ and $\balpha'\triangleq\bbeta-\bbeta'\in\sB_0[k+k']$---$(k+k')$-sparse vectors satisfying $\bX\balpha=\bX\bbeta+\bX\bbeta'$ and $\bX\balpha'=\bX\bbeta-\bX\bbeta'$---such that 
$$
\begin{aligned}
(1-\delta_{k+k'}) \cdot \normtwo{\bbeta+\bbeta'}^{2} \leq  \normtwo{\bX(\bbeta+\bbeta')}^{2} \leq (1+\delta_{k+k'}) \cdot \normtwo{\bbeta+\bbeta'}^{2};\\
(1-\delta_{k+k'}) \cdot \normtwo{\bbeta-\bbeta'}^{2} \leq  \normtwo{\bX(\bbeta-\bbeta')}^{2} \leq (1+\delta_{k+k'}) \cdot \normtwo{\bbeta-\bbeta'}^{2}.
\end{aligned}
$$
Since $\normtwo{\bbeta-\bbeta'}^{2} = \normtwo{\bbeta+\bbeta'}^{2}=2$, 
these inequalities simplify to:
\begin{equation}\label{equation:ric_roc_conn}
\begin{aligned}
2(1-\delta_{k+k'})  \leq  \normtwo{\bX(\bbeta+\bbeta')}^{2} \leq 2(1+\delta_{k+k'});\\
2(1-\delta_{k+k'})  \leq  \normtwo{\bX(\bbeta-\bbeta')}^{2} \leq 2(1+\delta_{k+k'}).
\end{aligned}
\end{equation}
Subtracting the two squared norms and using the identity
$$
\normtwo{\bX(\bbeta+\bbeta')}^2
-\normtwo{\bX(\bbeta-\bbeta')}^2
= 4\innerproduct{\bX\bbeta, \bX\bbeta'},
$$
one can further derive a quantitative link between $\theta_{k,k'}$
and $\delta_{k+k'}$, showing that small RIP constants imply small ROC constants; see Proposition~\ref{proposition:prop_rip_rop}.

\begin{exercise}[Nondecreasing of ROC]\label{exercise:nonde_roc}
Show that the ROC defined in \eqref{equation:def_roc22} satisfies   the following monotonicity properties:
\begin{equation}\label{equation:mont_roc22}
\begin{aligned}
\theta_{k,k'} &\leq \theta_{s,k'}, \quad \text{if}\quad k \leq s;\\
\theta_{k,k'} &\leq \theta_{k, s}, \quad \text{if}\quad k' \leq s.
\end{aligned}
\end{equation}
\end{exercise}

\subsection{Relationship between RIC, ROC, and other Constants}

In this subsection, we explore the connections among various constants associated with the design matrix $\bX$, including the RIC, mutual coherence, and the ROC.
\begin{proposition}[RIC and coherence]\label{proposition:ric_cohere}
Let $\bX\in\real^{n\times p}$ have $\ell_2$-normalized columns, i.e., $\normtwo{\bx_i} = 1$ for all $i \in \{1,2,\ldots,p\}$. Then, the RICs satisfy
$$
\delta_1 = 0, \quad \delta_2 = \mu(\bX), \quad \delta_k \leq \mu_1(\bX, k-1) \leq (k-1)\mu(\bX), \quad k \geq 2.
$$
\end{proposition}
\begin{proof}[of Proposition~\ref{proposition:ric_cohere}]
The $\ell_2$-normalization of the columns implies  that $\normtwo{\bX\be_i}^2 = \normtwo{\be_i}^2$  for each standard basis vector $\be_i$. 
Hence, the RIP condition holds exactly for 1-sparse vectors, yielding $\delta_1 = 0$.
For $\delta_2$, consider any pair of distinct indices $i\neq j$. The Gram matrix of the corresponding submatrix satisfies
$$
\delta_2 = \max_{1 \leq i \neq j \leq p} \normtwo{\bX[:,\{i,j\}]^\top\bX[:,\{i,j\}] - \bI_2}.
$$
The eigenvalues of this symmetric matrix  $\bX[:,\{i,j\}]^\top\bX[:,\{i,j\}] - \bI = \begin{bmatrixscript} 0 & \innerproduct{\bx_i, \bx_j} \\ \innerproduct{\bx_i, \bx_j} & 0 \end{bmatrixscript}$ are 
$\pm\abs{ \innerproduct{\bx_i, \bx_j}}$, and thus its spectral norm equals $\abs{\innerproduct{\bx_i, \bx_j}}$. Taking the maximum over $1 \leq i \neq j \leq p$ gives  $\delta_2 = \mu(\bX)$, the mutual coherence. 
The inequality $\delta_k \leq \mu_1(\bX,k-1) \leq (k-1)\mu(\bX)$ follows from Theorem~\ref{theorem:l1co_and_rip} and \eqref{equation:l1cohere_pro1}.
\end{proof}

We now present a comparison between restricted isometry constants and restricted orthogonality constants.
This result provides control over the principal angles between subspaces spanned by disjoint subsets of columns of $\bX$ of sizes $s$ and $t$, respectively. Such control is valuable because it enables guarantees of exact sparse recovery based solely on bounds involving the RICs; see Section~\ref{section:spar_rec_rip} for further details.
\begin{proposition}[RIC and ROC \citep{candes2005decoding, foucart2013invitation}]\label{proposition:prop_rip_rop}
For all positive integers $s,t$, we have that  
$$
\theta_{s,t}\leq \delta_{s+t} \leq 
\begin{cases}
&\theta_{s,t} + \max\{\delta_s, \delta_{t}\};\\
&\frac{1}{s+t}(s\delta_s + t\delta_t + 2\sqrt{st}\theta_{s,t}).
\end{cases}
$$
In the special case $t = s$, this yields
$$
\theta_{s,s} \leq \delta_{2s}  \leq \delta_s + \theta_{s,s}.
$$
\end{proposition}
\begin{proof}[of Proposition~\ref{proposition:prop_rip_rop}]
\textbf{The inequality $\theta_{s,t} \leq \delta_{s+t}$.}
The inequality  $\theta_{s,t} \leq \delta_{s+t}$ follows directly from definition. 
Alternatively, by homogeneity, it suffices to show that
\begin{equation}\label{equation:prop_rip_rop_pe1}
\abs{\innerproduct{\sum_{i \in \sI} \beta_i \bx_i, \sum_{j \in \sJ} \beta_{j} \bx_{j}}} \leq \delta_{s+t}
\end{equation}
whenever $\sI, \sJ\subseteq\{1,2,\ldots,p\}$ are disjoint ($\sI\cap\sJ=\varnothing$), $\abs{\sI} \leq s$, $\abs{\sJ} \leq t$, 
and $\bbeta\in\real^p$ satisfies $\normtwo{\bbeta_{\sI}} = \normtwo{\bbeta_{\sJ}}=1$, i.e., 
$\sum_{i \in \sI} \abs{\beta_i}^2 = \sum_{j \in \sJ} \abs{\beta_{j}}^2 = 1$.
To see this, \eqref{equation:ric_roc_conn} indicates that
$$
\begin{aligned}
2(1 - \delta_{s+t})&\leq \normtwo{\sum_{i \in \sI} \beta_i \bx_i + \sum_{j \in \sJ} \beta_{j} \bx_{j}}^2 \leq 2(1 + \delta_{s+t});\\
2(1 - \delta_{s+t})&\leq \normtwo{\sum_{i \in \sI} \beta_i \bx_i - \sum_{j \in \sJ} \beta_{j} \bx_{j}}^2 \leq 2(1 + \delta_{s+t}).
\end{aligned}
$$
Inequality \eqref{equation:prop_rip_rop_pe1} then follows from the parallelogram identity
$4\innerproduct{\ba, \bb} = {\normtwo{\ba + \bb}^2 - \normtwo{\ba - \bb}^2}$ for any vectors $\ba,\bb$.

\paragraph{The inequality $\delta_{s+t} \leq \theta_{s,t} + \max\{\delta_s, \delta_{t}\}$.}
For the reverse inequality, again by homogeneity, consider any $(s+t)$-sparse vector $\bbeta$ with $\normtwo{\bbeta}=1$.
Partition its support as $\sS = \sI \cup \sJ$, where $\abs{\sI} \leq s$, $\abs{\sJ} \leq t$, and define $ \nu \triangleq \sum_{i \in \sI} \abs{\beta_i}^2 $, so that $1 - \nu =\sum_{j \in \sJ} \abs{\beta_j}^2 $.
Then 
\begin{equation}\label{equation:prop_rip_rop_pe2}
\begin{aligned}
\abs{\normtwo{\bX\bbeta}^2 - 1}
&=
\abs{\normtwo{\bX_\sS\bbeta_\sS}^2 - \normtwo{\bbeta}^2}\\
&\leq \abs{\normtwo{\bX_{\sI} \bbeta_{\sI}}^2 - \nu}
+ \abs{\normtwo{\bX_{\sJ} \bbeta_{\sJ}}^2 - (1 - \nu)}
+ 2\abs{\innerproduct{\bX_{\sI} \bbeta_{\sI}, \bX_{\sJ} \bbeta_{\sJ}}} \\
&\leq \delta_s \nu + \delta_t (1 - \nu) + 2 \theta_{s,t} \sqrt{\nu(1 - \nu)}
\leq \max\{\delta_s, \delta_t\} (\nu + 1 - \nu) + \theta_{s,t} \\
&= \max\{\delta_s, \delta_t\} + \theta_{s,t},
\end{aligned}
\end{equation}
where the first inequality follows from the triangle inequality, and the second inequality follows from \eqref{equation:rippro22_e1} and \eqref{equation:def_roc22_e1}, the definitions of RICs and ROCs.

\paragraph{The inequality $\delta_{s+t} \leq \frac{1}{s+t}(s\delta_s + t\delta_t + 2\sqrt{st}\theta_{s,t})$.}
This inequality provides a sharper relationship between RICs and ROCs.
Once again, let $\bbeta\in\real^p$ be $(s+t)$-sparse with $\normtwo{\bbeta} = 1$. 
By \eqref{equation:prop_rip_rop_pe2}, we have 
\begin{equation}\label{equation:prop_rip_rop_pe3}
\abs{\normtwo{\bX\bbeta}^2 - \normtwo{\bbeta}^2} 
\leq \delta_s \nu + \delta_t (1-\nu) + 2\theta_{s,t} \sqrt{\nu(1 - \nu)} 
\triangleq f(\nu),
\end{equation}
where $f(\nu): [0,1]\rightarrow \real$.
It can be shown that there is a scalar $\nu^* \in [0,1]$ such that this function is nondecreasing on $[0,\nu^*]$ and then nonincreasing on $[\nu^*,1]$; see Problem~\ref{prob:nonin_nond_func}. Depending on the location of this $\nu^*$ with respect to $s/(s+t)$, the function $f$ is either nondecreasing on $[0,s/(s+t)]$ or nonincreasing on $[s/(s+t),1]$. By properly choosing the vector $\bbeta_{\sI}$, we can always assume that $\normtwo{\bbeta_{\sI}}^2$ is in one of these intervals. Indeed, if $\bbeta_{\sI}$ is made of $s$ smallest modulus components of $\bbeta$ while $\bbeta_{\sJ}$ is made of $t$ largest modulus components of $\bbeta$, then we have
$$
\frac{\normtwo{\bbeta_{\sI}}^2}{s} \leq \frac{\normtwo{\bbeta_{\sJ}}^2}{t} = \frac{1-\normtwo{\bbeta_{\sI}}^2}{t}, \quad \text{so that } \normtwo{\bbeta_{\sI}}^2 \leq \frac{s}{s+t}.
$$
Otherwise, if $\bbeta_{\sI}$ is made of $s$ largest modulus components of $\bbeta$, then we have $\normtwo{\bbeta_{\sI}}^2 \geq s/(s+t)$. 
Therefore, we obtain
$$
\abs{\normtwo{\bX\bbeta}^2 - \normtwo{\bbeta}^2} \leq f\left(\frac{s}{s+t}\right) = \delta_s \frac{s}{s+t} + \delta_t \frac{t}{s+t} + 2\theta_{s,t} \frac{\sqrt{st}}{s+t}.
$$
This completes the proof.
\end{proof}

We conclude this subsection by showing that RICs and ROCs of higher order can be bounded in terms of those of lower order.

\begin{proposition}[Different orders of RICs and ROCs]\label{proposition:diff_order_ricroc}
Let $k, s, t \geq 1$ be integers with $t \geq s$. 
Then
$$
\begin{aligned}
\theta_{t,k} &\leq \sqrt{\frac{t}{s}} \theta_{s,k},\\
\delta_t &\leq \frac{t-d}{s} \delta_{2s} + \frac{d}{s} \delta_s, \quad \text{where } d \triangleq \gcd(s,t).
\end{aligned}
$$
In the special case $t = cs$ for some integer $c\geq 1$, this simplifies to
$$
\delta_{cs} \leq c \cdot \delta_{2s}.
$$
\end{proposition}
\begin{proof}[of Proposition~\ref{proposition:diff_order_ricroc}]
Let   $\balpha \in \real^p$ be a $t$-sparse vector and  $\bbeta \in \real^p$ a $k$-sparse vector, with disjoint supports.
We aim to establish the two inequalities:
\begin{equation}\label{equation:diff_order_ricroc1}
\abs{\innerproduct{\bX\balpha, \bX\bbeta}} \leq \sqrt{\frac{t}{s}} \theta_{s,k} \normtwo{\balpha} \normtwo{\bbeta},
\end{equation}
\begin{equation}\label{equation:diff_order_ricroc2}
\abs{\normtwo{\bX\balpha}^2 - \normtwo{\balpha}^2} \leq \left(\frac{t-d}{s} \delta_{2s} + \frac{d}{s} \delta_s\right) \normtwo{\balpha}^2.
\end{equation}
Let $d$ be a common divisor of $s$ and $t$. We define integers $\mu, \nu$ such that
$$
s = \mu d, \qquad t = \nu d.
$$
Let $\sT = \{j_1, j_2, \ldots, j_t\}$ denote the support of $\balpha$. We consider the $\nu$ subsets $\sS_1, \sS_2, \ldots, \sS_\nu \subseteq \sT$ of size $s$ defined by
$$
\sS_i = \{j_{(i-1)d+1}, j_{(i-1)d+2}, \ldots, j_{(i-1)d+s}\},
$$
where indices are taken modulo $t$. 
With this construction, every index $j \in 
\sT$ appears in exactly $s/d = \mu$ of the subsets  $\sS_i$. 
Consequently, we can decompose $\balpha$ as
$$
\balpha = \frac{1}{\mu} \sum_{i=1}^\nu \balpha(\sS_i), \qquad \normtwo{\balpha}^2 = \frac{1}{\mu} \sum_{i=1}^\nu \normtwo{\balpha(\sS_i)}^2.
$$
Note again that $\balpha(\sS_i)\in\real^p$ (while $\balpha_{\sS_i}\in\real^s$; see Definition~\ref{definition:matlabnotation}).

\paragraph{Proof of \eqref{equation:diff_order_ricroc1}.}
Using the decomposition above and the definition of the ROC, we obtain
$$
\begin{aligned}
\innerproduct{\bX\balpha, \bX\bbeta} 
&\leq \frac{1}{\mu} \sum_{i=1}^\nu \innerproduct{\bX\balpha(\sS_i), \bX\bbeta} 
\leq \frac{1}{\mu} \sum_{i=1}^\nu \theta_{s,k} \normtwo{\balpha(\sS_i)} \normtwo{\bbeta}\\
&\leq \theta_{s,k} \frac{\sqrt{\nu}}{\mu} \sqrt{\left(\sum_{i=1}^\nu \normtwo{\balpha_{\sS_i}}^2\right)} \normtwo{\bbeta} = \theta_{s,k} \sqrt{\frac{\nu}{\mu}} \normtwo{\balpha} \normtwo{\bbeta},
\end{aligned}
$$
where the last inequality follows from the Cauchy--Schwarz inequality $\ba^\top\bb\leq \normtwo{\ba}\normtwo{\bb}$ by invoking with $\ba=\bone$ and $b_i = \normtwo{\balpha(\sS_i)}$ for all $i$.

\paragraph{Proof of \eqref{equation:diff_order_ricroc2}.}
We begin by writing
$$
\footnotesize
\begin{aligned}
&\normtwo{\bX\balpha}^2 - \normtwo{\balpha}^2 
= \abs{\innerproduct{\big(\bX^\top\bX - \bI\big)\balpha, \balpha}} 
\leq \frac{1}{\mu^2} \sum_{1 \leq i \leq \nu} \sum_{1 \leq j \leq \nu} 
\abs{\innerproduct{\big(\bX^\top\bX - \bI\big)\balpha(\sS_i), \balpha(\sS_j)}}\\
&= \frac{1}{\mu^2} \left( 
\sum_{1 \leq i \neq j \leq \nu} \abs{\innerproduct{\big(\bX^\top\bX - \bI\big)\balpha(\sS_i), \balpha(\sS_j)}}
+ \sum_{1 \leq i \leq \nu} \abs{\innerproduct{\big(\bX_{\sS_i}^\top\bX_{\sS_i} - \bI\big)\balpha_{\sS_i}, \balpha_{\sS_i}}} \right)\\
&\stackrel{\dag}{\leq} \frac{1}{\mu^2} \left( \sum_{1 \leq i \neq j \leq \nu} \delta_{2s} \normtwo{\balpha_{\sS_i}} \normtwo{\balpha_{\sS_j}} + \sum_{1 \leq i \leq \nu} \delta_s \normtwo{\balpha_{\sS_i}}^2 \right)
= \frac{\delta_{2s}}{\mu^2} \left( \sum_{1 \leq i \leq \nu} \normtwo{\balpha_{\sS_i}} \right)^2 - \frac{\delta_{2s} - \delta_s}{\mu^2} \sum_{1 \leq i \leq \nu} \normtwo{\balpha_{\sS_i}}^2\\
&\stackrel{\ddag}{\leq} \left( \frac{\delta_{2s} \nu}{\mu^2} - \frac{\delta_{2s} - \delta_s}{\mu^2} \right) \sum_{1 \leq i \leq \nu} \normtwo{\balpha_{\sS_i}}^2 
= \left( \frac{\nu}{\mu} \delta_{2s} - \frac{\delta_{2s} - \delta_s}{\mu} \right) \normtwo{\balpha}^2
= \left( \frac{t}{s} \delta_{2s} - \frac{1}{\mu} (\delta_{2s} - \delta_s) \right) \normtwo{\balpha}^2,
\end{aligned}
$$
where the inequality ($\dag$) follows from Proposition~\ref{proposition:conse_rip}, and the inequality $(\ddag)$
follows again from the Cauchy--Schwarz inequality $\ba^\top\bb\leq \normtwo{\ba}\normtwo{\bb}$ by invoking with $\ba=\bone$ and $b_i = \normtwo{\balpha_{\sS_i}}$ for all $i$.
To make the latter bound as tight as possible, we take $\mu$ as small as possible, i.e., we choose $d$ as the greatest common divisor of $s$ and $t$. This completes the proof.
\end{proof}

\subsection{Bounds for Restricted Isometry Property}
Just as with mutual coherence or $\ell_1$-coherence (see Section~\ref{section:bd_cohe_l1}), it is important to understand how small the $k$-th RIC of a matrix $\bX \in \real^{n \times p}$ can be. 
In the regime relevant to compressive sensing---namely, when $p \geq C n$---Theorem~\ref{theorem:dim_rip} below shows that the RIC must satisfy
$\delta_k \geq c \sqrt{k/n}$.
For $k = 2$, this implies the mutual coherence satisfies $\mu(\bX) \geq c' / \sqrt{n}$ by Proposition~\ref{proposition:cohe2rip}, which is reminiscent of the Welch bound of Theorem~\ref{theorem:low_mutucohe}.

\begin{theoremHigh}[Bounds for RICs \citep{krahmer2011new, foucart2013invitation}]\label{theorem:dim_rip}
Let $\bX \in \real^{n \times p}$ and $2 \leq k \leq p$. 
Suppose $p \geq C n$ and $\delta_k \leq \delta_*$,  where the constants $c$, $C$, and $\delta_*$ depend only on each other.
Then
\begin{equation}\label{equation:ric_bd_diss0}
n \geq c \frac{k}{\delta_k^2}.
\end{equation} 
For example, the choices $c = 1/162$, $C = 30$, and $\delta_* = 2/3$ are valid.
\end{theoremHigh}
\begin{proof}[of Theorem~\ref{theorem:dim_rip}]
We first observe that the statement cannot hold for  $k = 1$, since  $\delta_1 = 0$ whenever all columns of $\bX$ are normalized to unit  $\ell_2$-norm (see Proposition~\ref{proposition:ric_cohere}). 
Set $t \triangleq \lfloor k/2 \rfloor \geq 1$, and partition $\bX$ in blocks of size $n \times t$ (with possibly fewer columns in the last block):
$$
\bX = [\bX_1, \bX_2, \ldots, \bX_\nu], \quad p \leq \nu t.
$$
By Definitions~\ref{definition:rip22} and \ref{definition:rop22},  Exercise~\ref{exercise:order_delta_rip}, and Proposition~\ref{proposition:prop_rip_rop}, we have for all distinct  $i, j \in \{1,2,\ldots,\nu\}$, $i \neq j$,
$$
\normtwo{\bX_i^\top\bX_i - \bI} \leq \delta_t \leq \delta_k, 
\qquad \normtwo{\bX_i^\top\bX_j} \leq \theta_{t,t} \leq \delta_{2t} \leq \delta_k.
$$
Thus, all the eigenvalues of $\bX_i^\top\bX_i$ and all the singular values of $\bX_i^\top\bX_j$ satisfy
$$
1 - \delta_k \leq \lambda_h(\bX_i^\top\bX_i) \leq 1 + \delta_k, \qquad \sigma_h(\bX_i^\top\bX_j) \leq \delta_k, 
\qquad h=1,2,\ldots,t
$$
Define the Gram matrix and its dual:
$$
\bG \triangleq \bX^\top\bX = [\bX_i^\top\bX_j]_{1 \leq i,j \leq \nu } \in \real^{p \times p}, 
\qquad 
\bH \triangleq \bX\bX^\top \in \real^{n\times n}.
$$
\paragraph{Lower bound on $\trace(\bH)$.}
Using the cyclic invariance of the trace and the fact that each $\bX_i^\top\bX_i$ is symmetric, we obtain
\begin{equation}\label{equation:dim_rip1}
\trace(\bH) = \trace(\bG) = \sum_{i=1}^\nu \trace(\bX_i^\top\bX_i)
= \sum_{i=1}^\nu  \sum_{h=1}^t \lambda_h(\bX_i^\top\bX_i) \geq \nu t(1 - \delta_k).
\end{equation}
\paragraph{Upper bound via Frobenius norm.} 
Recall that for any matrices $\bA_1,\bA_2$, the Frobenius inner product satisfies $\innerproduct{\bA_1, \bA_2}_F = \trace(\bA_2^\top\bA_1)$, and by Cauchy--Schwarz,
$$
\trace(\bH)^2 = \innerproduct{\bI_n, \bH}_F^2 \leq \normf{\bI_n}^2 \normf{\bH}^2 
= n \trace(\bH^\top\bH).
$$
Now, using cyclic invariance again and the fact that $\bX\bX^\top=\sum_{i} \bX_i\bX_i^\top$, we have 
$$
\begin{aligned}
\trace(\bH^\top\bH) 
&= \trace(\bX\bX^\top\bX\bX^\top) = \trace(\bX^\top\bX\bX^\top\bX) = \trace(\bG\bG^\top)\\
&= \sum_{i=1}^\nu  \trace\left(\sum_{j=1}^\nu  \bX_i^\top\bX_j\bX_j^\top\bX_i\right)
\stackrel{\dag}{=} \sum_{1 \leq i \neq j \leq \nu } \sum_{h=1}^t \sigma_h(\bX_i^\top\bX_j)^2 + \sum_{i=1}^\nu  \sum_{h=1}^t \lambda_h(\bX_i^\top\bX_i)^2\\
&\leq \nu (\nu-1)t\delta_k^2 + \nu t(1+\delta_k)^2.
\end{aligned}
$$
where the equality $(\dag)$ follows from the fact that the trace of a symmetric matrix is equal to the sum of its eigenvalues, and the fact that $\sigma(\bA)^2 = \lambda(\bA^\top\bA)=\lambda(\bA\bA^\top)$ for any matrix $\bA$.
This implies the upper bound
\begin{equation}\label{equation:dim_rip2}
\trace(\bH)^2 \leq 
n \trace(\bH^\top\bH)
\leq  n\nu t\big((\nu-1)\delta_k^2 + (1+\delta_k)^2\big).
\end{equation}
Combining the bounds \eqref{equation:dim_rip1} and \eqref{equation:dim_rip2} yields the lower bound for the number of samples:
\begin{equation}\label{equation:dim_rip3}
n \geq \frac{\nu t(1-\delta_k)^2}{(\nu-1)\delta_k^2 + (1+\delta_k)^2}.
\end{equation}
If $(\nu-1)\delta_k^2 < (1+\delta_k)^2/5$, we would obtain, using $\delta_k \leq 2/3$,
$$
n > \frac{\nu t(1-\delta_k)^2}{6(1+\delta_k)^2/5} \geq \frac{5(1-\delta_k)^2}{6(1+\delta_k)^2} p \geq \frac{1}{30} p,
$$
which contradicts our assumption. We therefore have $(\nu-1)\delta_k^2 \geq (1+\delta_k)^2/5$, which yields, using $\delta_k \leq 2/3$ again and $k \leq 3t$,
$$
n \geq \frac{\nu t(1-\delta_k)^2}{6(\nu-1)\delta_k^2} \geq \frac{1}{54\delta_k^2} t \geq \frac{1}{162\delta_k^2} k.
$$
This completes the proof.
\end{proof}

We will see later in Chapter~\ref{chapter:ensur_rips} (e.g., Corollary~\ref{corollary:gaussben_m_bound_rip}) that certain random matrices $\bX \in \real^{n \times p}$ satisfy $\delta_k \leq \delta$ with high probability for any fixed $\delta \in(0,1)$, provided that
$$
n \geq C\delta^{-2}k\ln(ep/k).
$$
Thus, the lower bound~\eqref{equation:ric_bd_diss0} is optimal up to the logarithmic factor.

Theorem~\ref{theorem:dim_rip} examines measurement bounds when the RIP constant $\delta$ is fixed.
We can also ask a complementary question: How many measurements are necessary to satisfy the RIP, if we ignore the specific value of $\delta$ and focus only on the problem dimensions $n, p$, and $k$?
In this setting, we can derive a simple lower bound on $n$. To do so, we first establish a key lemma that will be used in the proof of the main result.

\begin{lemma}\label{lemma:mea_bd_rip}
Let $k$ and $p$ be integers satisfying $k < p/2$. 
Then there exists a set $\sS \subset \sB_0[k]$ such that 
\begin{enumerate}[(i)]
\item  $\normtwo{\bbeta} \leq \sqrt{k}$ for all	 $\bbeta \in \sS$.
\item $\normtwo{\bbeta - \balpha} \geq \sqrt{k/2}$  for  all distinct $\balpha,\bbeta  \in \sS$ with $\bbeta \neq \balpha$.
\item $\ln \abs{\sS} \geq \frac{k}{2} \ln \left( \frac{p}{k} \right)$.
\end{enumerate}
\end{lemma}
\begin{proof}[of Lemma~\ref{lemma:mea_bd_rip}]
Consider the set
$$
\sT\triangleq \{\bbeta \in \{0, +1, -1\}^p \mid \normzero{\bbeta}= k\}.
$$
By construction, every  $\bbeta \in \sT$ satisfies $\normtwo{\bbeta} = \sqrt{k}$. 
Thus, any subset of $\sT$ automatically satisfies condition (i).

Note that the cardinality of $\sT$ is $\abs{\sT} = \binom{p}{k} 2^k$.
Moreover, since  $\normzero{\bbeta - \balpha} \leq \normtwo{\bbeta - \balpha}^2$ (when $\balpha, \bbeta\in\sT$),
it follows that if  $\normtwo{\bbeta - \balpha}^2 \leq k/2$, then $\normzero{\bbeta - \balpha} \leq k/2$. 
Hence, for any fixed $\bbeta\in\sT$, the number of $\balpha\in\sT$ within squared $\ell_2$-distance $k/2$ satisfies
$$
\abs{\{\balpha \in \sT : \normtwo{\bbeta - \balpha}^2 \leq k/2\}} 
\leq 
\absbig{\{\balpha \in \sT \mid \normzero{\bbeta - \balpha} \leq k/2\}} 
\leq \binom{p}{k/2} 3^{k/2}.
$$
We now construct $\sS$ greedily: starting from $\sT$, we iteratively select vectors that are at least $\sqrt{k/2}$ apart in $\ell_2$-norm. After selecting $j$ points, at most $j \binom{p}{k/2} 3^{k/2}$ points in $\sT$ are excluded from further selection. 
Therefore, as long as
\begin{equation}\label{equa:lem_mea_bd_rip_e3}
\abs{\sS} \binom{p}{k/2} 3^{k/2} \leq \binom{p}{k} 2^k
\equiv \abs{\sT}.
\end{equation}
such a set $\sS$ of size $\abs{\sS}$ exists.
To verify this inequality, observe that
$$
\frac{\binom{p}{k}}{\binom{p}{k/2}} = \frac{(k/2)! (p - k/2)!}{k! (p - k)!} = \prod_{i=1}^{k/2} \frac{p - k + i}{k/2 + i} \geq \left(\frac{p}{k} - \frac{1}{2}\right)^{k/2},
$$
where the inequality holds because the function $f(i)\triangleq (p - k + i) / (k/2 + i)$ is decreasing in $i$~\footnote{Since $f'(i)=\frac{3k/2-p}{(k/2+i)^2}$ is negative under the assumption that $k<p/2$.}. 
Now set $\abs{\sS} = (p/k)^{k/2}$. Then
\begin{equation}\label{equa:lem_mea_bd_rip_e4}
	\abs{\sS} \left(\frac{3}{4}\right)^{k/2} = \left(\frac{3p}{4k}\right)^{k/2} 
	= \left(\frac{p}{k} - \frac{p}{4k}\right)^{k/2} \stackrel{\dag}{\leq} \left(\frac{p}{k} - \frac{1}{2}\right)^{k/2} \leq \frac{\binom{p}{k}}{\binom{p}{k/2}},
\end{equation}
where the inequality ($\dag$) follows from the assumption that $k < p/2$.
Combining \eqref{equa:lem_mea_bd_rip_e3} and \eqref{equa:lem_mea_bd_rip_e4}, it follows that \eqref{equa:lem_mea_bd_rip_e3} holds for $\abs{\sS} = (p/k)^{k/2}$, which establishes the lemma.
\end{proof}

Using this packing lemma, we now derive a fundamental lower bound on the number of measurements required for the RIP.

\begin{theoremHigh}[Measurement bound for RIP \citep{garnaev1984widths, baraniuk2011introduction}]\label{theorem:measur_bd_rip}
Let $\bX\in\real^{n\times p}$ satisfy the RIP of order $2k$ with constant $\delta \in (0, \frac{1}{2}]$
\footnote{Note that the restriction $\delta \leq \frac{1}{2}$ is made for convenience; minor adjustments to the constants allow the same argument to work for any $\delta_{\max} < 1$. }. 
Then,
\begin{equation}
n \geq C k \ln \left( \frac{p}{k} \right), \quad\text{with }C \triangleq \frac{1}{2 \ln (\sqrt{24} + 1)} \approx 0.28.
\end{equation}
\end{theoremHigh}
\begin{proof}[of Theorem~\ref{theorem:measur_bd_rip}]
Since $\bX$ satisfies the RIP or order $2k$, 
for any $\balpha,\bbeta\in\sS$ (with $\sS$ from Lemma~\ref{lemma:mea_bd_rip}), 
we have $\bbeta-\balpha\in\sB_0[2k]$, and thus
\begin{equation}
\normtwo{\bX\bbeta - \bX \balpha} \geq \sqrt{1 - \delta} \normtwo{\bbeta - \balpha} \geq \sqrt{k/4}, \quad \forall \, \balpha,\bbeta \in \sS.
\end{equation}
Similarly, using the upper RIP bound,
\begin{equation}
\normtwo{\bX\bbeta} \leq \sqrt{1 + \delta} \normtwo{\bbeta} \leq \sqrt{3k/2}, \quad \forall \, \bbeta \in \sS.
\end{equation}
The lower bound tells us that for any pair of points $\balpha,\bbeta \in \sS$, if we center balls of radius $\sqrt{k/4}/2 = \sqrt{k/16}$ at $\bX\bbeta$ and $\bX \balpha$, then these balls will be disjoint. In turn, the upper bound tells  
us that the entire set of balls is itself contained within a larger ball of radius $\sqrt{3k/2} + \sqrt{k/16}$. 
Let $\sB_2[r] = \{\bx \in \real^n \mid \normtwo{\bx} \leq r\}$ denote the $\ell_2$-ball. This implies that
\begin{equation}
\begin{aligned}
&&\text{Vol}\left(\sB_2\left[\sqrt{3k/2} + \sqrt{k/16}\right]\right) &\geq&&  \abs{\sS} \cdot \text{Vol}\left(\sB_2\left[\sqrt{k/16}\right]\right) \\
\implies&&\left(\sqrt{3k/2} + \sqrt{k/16}\right)^n &\geq&&  \abs{\sS} \cdot \left(\sqrt{k/16}\right)^n  \\
\implies&&(\sqrt{24} + 1)^n &\geq&& \abs{\sS}  \\
\implies& & n &\geq&& \frac{\ln  \abs{\sS}}{\ln (\sqrt{24} + 1)}.
\end{aligned}
\end{equation}
By the bound for $ \abs{\sS}$ from Lemma~\ref{lemma:mea_bd_rip}, this completes the proof.
\end{proof}

Theorem~\ref{theorem:dim_rip} provides a lower bound on $n$ that depends explicitly on the RIP constant $\delta_k$ and the sparsity level $k$, but does not depend on $p$ (the ambient dimension), except through the condition $p\geq Cn$.
Theorem~\ref{theorem:measur_bd_rip} gives a lower bound on $n$ that depends on $k$ and the logarithmic ratio $\ln(p/k)$, but assumes a fixed range for $\delta$ (e.g., $\delta \leq 1/2$) and does not explicitly involve $\delta$ in the bound.
Although the proofs above are somewhat indirect, a similar scaling in $p$ and $k$ can also be derived via the Gelfand width of the $\ell_1$-ball \citep{garnaev1984widths}.

\subsection{The RIP and Stability}

We will see later in this book that if a matrix $\bX$ satisfies the RIP, then this condition is sufficient to guarantee that a wide range of algorithms can successfully recover a sparse signal from noisy measurements; see Section~\ref{section:spar_rec_rip}.
Here, we examine whether the RIP is also necessary. It should be clear that the lower bound in the RIP is necessary if we wish to recover all $k$-sparse signals $\bbeta$ exactly from noiseless measurements $\bX\bbeta$, for the same reasons that the NSP is necessary (see Exercise~\ref{exercise:nspi_2k}). 
In fact, we can say even more about the necessity of the RIP by introducing the following notion of stability.

\begin{definition}[$C$-stable \citep{baraniuk2011introduction}]
Let $\bX\in\real^{n\times p}$ be a sensing matrix and $\Delta:\real^n\rightarrow \real^p$ a recovery algorithm. 
We say that the pair $(\bX, \Delta)$ is \textit{$C$-stable} if, for any $\bbeta \in \sB_0[k]$ and any $\bepsilon \in \real^n$, 
\begin{equation}
\normtwo{\Delta(\bX\bbeta + \bepsilon) - \bbeta} \leq C \normtwo{\bepsilon}.
\end{equation}
\end{definition}

This definition simply states that adding a small amount of noise to the measurements should not cause an arbitrarily large error in the recovered signal.
The following lemma shows that the existence of any decoding algorithm---even one that may be computationally impractical---that can stably recover signals from noisy measurements implies that $\bX$ must satisfy the lower bound of \eqref{equation:rippro22_e1}, with a constant determined by $C$.
\begin{lemma}[$C$-stable \citep{baraniuk2011introduction}\index{$C$-stable}]\label{lemma:cstable}
If the pair $(\bX, \Delta)$ is $C$-stable, then
\begin{equation}\label{equation:cstable}
\frac{1}{C} \normtwo{\bbeta} \leq \normtwo{\bX\bbeta}
\end{equation}
for all $\bbeta \in \sB_0[2k]$.
\end{lemma}
\begin{proof}[of Lemma~\ref{lemma:cstable}]
Let $\balpha, \bgamma \in \sB_0[k]$. Define
$$
\bepsilon_{\beta} \triangleq \frac{\bX(\balpha - \bgamma)}{2} 
\qquad \text{and} \qquad 
\bepsilon_{\alpha} = \frac{\bX(\bgamma - \balpha)}{2},
$$
and observe that
$$
\bX\bgamma + \bepsilon_{\beta} = \bX \balpha + \bepsilon_{\alpha} = \frac{\bX(\bgamma + \balpha)}{2}.
$$
Let $\widehatbbeta \triangleq \Delta(\bX\bgamma + \bepsilon_{\beta}) = \Delta(\bX \balpha + \bepsilon_{\alpha})$. 
By the triangle inequality and the definition of $C$-stability, 
$$
\begin{aligned}
\normtwo{\bgamma - \balpha} 
&= \normtwobig{\bgamma - \widehatbbeta + \widehatbbeta - \balpha} 
\leq \normtwobig{\bgamma - \widehatbbeta} + \normtwobig{\widehatbbeta - \balpha} \\
&\leq C \normtwo{\bepsilon_{\beta}} + C \normtwo{\bepsilon_{\alpha}} 
\leq  C \normtwo{\bX\bgamma - \bX \balpha}.
\end{aligned}
$$
Since this holds for any $\balpha,\bgamma \in \sB_0[k]$, it follows that for any $\bbeta\in\sB_0[2k]$ (by setting $\bbeta=\bgamma-\balpha$), we have $ \normtwo{\bbeta} \leq C\normtwo{\bX\bbeta}$, which completes the proof.
\end{proof}

Note that as $C \to 1$, the matrix $\bX$ must satisfy the lower bound of \eqref{equation:rippro22_e1} with $\delta_{k} \triangleq 1 - 1/C^2 \to 0$. Thus, to reduce the impact of noise on the reconstructed signal, we must choose $\bX$ so that it satisfies the RIP lower bound with a tighter (i.e., smaller) constant.

Since only the lower bound of RIP is necessary---and not the upper bound---we could avoid redesigning $\bX$ altogether by simply rescaling it. Indeed, if $\bX$ satisfies the RIP with $\delta_{2k} < 1$, then the rescaled matrix $\eta\bX$ will satisfy inequality \eqref{equation:cstable} for any desired constant $C$, provided $\eta$ is chosen appropriately. In settings where the noise magnitude is independent of $\bX$, this reasoning is valid: scaling $\bX$ increases the ``signal" component of the measurements without affecting the noise, thereby improving the signal-to-noise ratio. In principle, this allows the noise to become negligible relative to the signal.

However, in practice, we cannot arbitrarily scale $\bX$, since the noise does depend on $\bX$ in many real-world scenarios. 
For example, suppose the noise vector $\bepsilon$ arises from quantization using a finite dynamic-range quantizer with $B$ bits. If the original measurements lie in $[-R,R]$, and the quantizer is calibrated to this range, then rescaling $\bX$ by a factor $\eta$ causes the measurements to span $[-\eta R, \eta R]$. 
To accommodate this, the quantizer's dynamic range must also be scaled by $\eta$, resulting in quantization error $\eta\bepsilon$. Consequently, the reconstruction error does not decrease---it scales proportionally with $\eta$, offering no net improvement.

\subsection{The RIP, NSP, and Variants}\label{section:rip_imp_nsp}

In this subsection, we show that if a matrix satisfies the RIP, then it also satisfies various forms of the NSP. Consequently, the RIP is strictly stronger than the NSP.

We begin by illustrating how RIP implies the standard NSP. The following observation---which is an extension to the standard bounds on vector norms (Exercise~\ref{exercise:cauch_sc_l1l2})---is frequently used in our arguments, so we isolate it from the proof.
\begin{lemma}[Ordered entry norm inequality \citep{foucart2013invitation}]\label{lemma:ells_ellt}
Given $s > t > 0$, if $\bu \in \real^p$ and $\bv \in \real^m$ satisfy
\begin{equation}
\max_{i \in \{1,2,\ldots,p\}} \abs{u_i} \leq \min_{j \in \{1,2,\ldots, m\}} \abs{v_j},
\end{equation}
then
$$
\norm{\bu}_s \leq \frac{p^{1/s}}{m^{1/t}} \norm{\bv}_t.
$$
A special case occurs when $t=1$, $s=2$, and $p=m$, yielding
$$
\normtwo{\bu} \leq \frac{1}{\sqrt{p}} \normone{\bv}.
$$
\end{lemma}
\begin{proof}[of Lemma~\ref{lemma:ells_ellt}]
To prove this, observe that:
$$
\frac{\norm{\bu}_s}{p^{1/s}} = \left[ \frac{1}{p} \sum_{i=1}^p \abs{u_i}^s \right]^{1/s} \leq \max_{1 \leq i \leq p} \abs{u_i},
$$
$$
\frac{\norm{\bv}_t}{m^{1/t}} = \left[ \frac{1}{m} \sum_{j=1}^m \abs{v_j}^t \right]^{1/t} \geq \min_{1 \leq j \leq m} \abs{v_j}.
$$
This completes the proof.
\end{proof}

The following lemma establishes a key relationship between the $\ell_1$-norms of the largest and smallest entries in any vector from the null space of a matrix satisfying the RIP, showing that the energy of the significant components is strictly bounded by that of the less significant ones under suitable RIP conditions.
\begin{lemma}[RIP $\ell_1$-$\ell_1$ bound lemma]\label{lemma:rip2k_rho}
Let $\bX \in \real^{n \times p}$ satisfy the RIP of order $2k$ with constant $\delta_{2k}$.  
For any nonzero $\bbeta \in \nspace(\bX)$ and any index set $\sS$ consisting of the $k$ largest (in magnitude) coordinates of $\bbeta$, we have 
$$
\normone{\bbeta_{\sS}} \leq \rho \normone{\bbeta_{\comple{\sS}}} \quad \text{with} \quad \rho = \frac{\delta_{2k}}{1 - 2\delta_{2k}}.
$$
In particular, if $\delta_{2k} < \tfrac{1}{3}$, then $\rho < 1$.
\end{lemma}
\begin{proof}[of Lemma~\ref{lemma:rip2k_rho}]
Consider  $\bbeta \in \nspace(\bX) \setminus \{\bzero\}$, and let $\sS\triangleq \sS_0$ be the set of indices corresponding to the $k$ largest entries of  $\bbeta$ in absolute value.  
Partition $\comple{\sS}$ into disjoint blocks $\sS_1, \sS_2, \dots$,
each containing at most $k$ elements, such that the magnitudes of entries within each block are nonincreasing across blocks.
Since $\bbeta \in \nspace(\bX)$, We have the decomposition
$$
\bbeta = \bbeta(\sS) + \sum_{i \geq 1} \bbeta(\sS_i)
\qquad \implies\qquad 
\bzero = \bX\bbeta = \bX\bbeta(\sS) + \sum_{i \geq 1} \bX\bbeta(\sS_i),
$$
where $\bbeta(\sS)\in\real^p$ (see Definition~\ref{definition:matlabnotation}).
Taking the inner product with $\bX\bbeta(\sS)$ yields
$$
\bzero = \innerproduct{\bX\bbeta(\sS), \bX\bbeta}
= \normtwo{\bX\bbeta(\sS)}^2 + \sum_{i \geq 1} 
\innerproduct{\bX\bbeta(\sS), \bX\bbeta(\sS_i)}.
$$
Rearranging and taking absolute values gives
$$
\normtwo{\bX\bbeta(\sS)}^2
\leq \sum_{i \geq 1} \abs{\innerproduct{\bX\bbeta(\sS), \bX\bbeta(\sS_i)}}.
$$
From the RIP, for any two disjointly supported vectors $\bu,\bv$ with $\normzero{\bu} + \normzero{\bv} \leq 2k$, Proposition~\ref{proposition:conse_rip} shows that 
$$
\abs{\innerproduct{\bX \bu, \bX \bv}} \leq \delta_{2k} \, \normtwo{\bu} \, \normtwo{\bv}.
$$
Also, the RIP lower bound provides  $\normtwo{\bX \bbeta(\sS)}^2 
=\normtwo{\bX_\sS \bbeta_\sS}^2 
\geq (1 - \delta_{2k}) \normtwo{\bbeta_\sS}^2$.  
Applying these results, we obtain:
\begin{equation}\label{equation:tech_rip_ineq1}
(1 - \delta_{2k}) \normtwo{\bbeta_{\sS}}^2
\leq \delta_{2k} \, \normtwo{\bbeta_{\sS}} \sum_{i \geq 1} \normtwo{\bbeta_{\sS_i}}
\quad \implies\quad 
(1 - \delta_{2k}) \normtwo{\bbeta_{\sS}}
\leq \delta_{2k} \sum_{i \geq 1} \normtwo{\bbeta_{\sS_i}}.
\end{equation}

Because entries are nonincreasing across blocks, we have, for every $i\geq 1$,
$$
\norminf{\bbeta_{\sS_i}}
\leq \frac{1}{k}\sum_{v\in \sS_{i-1}} \abs{\beta_v} = \frac{1}{k}\normone{\bbeta_{\sS_{i-1}}}, 
\quad \forall\, i\geq 1.
$$
By the standard bounds on vector norms (Exercise~\ref{exercise:cauch_sc_l1l2}),
the $\ell_2$-norm of block $\sS_i$ satisfies:
$$
\normtwo{\bbeta_{\sS_i}} 
\leq \sqrt{k}\norminf{\bbeta_{\sS_i}}
\leq \sqrt{k}\cdot \frac{1}{k}\normone{\bbeta_{\sS_{i-1}}}
= \frac{1}{\sqrt{k}}\normone{\bbeta_{\sS_{i-1}}},
\quad \forall\, i\geq 1.
$$
(Note that the inequality can be directly obtained by Lemma~\ref{lemma:ells_ellt} such that  $\normtwo{\bbeta_{\sS_i}} \leq \frac{1}{\sqrt{k}} \normone{\bbeta_{\sS_{i-1}}}$.)
Summing these inequalities over $i\geq 1$ yields:
\begin{equation}\label{equation:tech_rip_ineq}
\sum_{i\geq 1}\normtwo{\bbeta_{\sS_i}}
\leq \frac{1}{\sqrt{k}}\sum_{i\geq 1}\normone{\bbeta_{\sS_{i-1}}}
= \frac{1}{\sqrt{k}}\sum_{i\ge0}\normone{\bbeta_{\sS_i}}
= \frac{1}{\sqrt{k}}\big(\normone{\bbeta_{\sS}} + \normone{\bbeta_{\comple{\sS}}}\big).
\end{equation}
~\footnote{Note that if you instead try to show the simpler-looking bound $\sum_i\normtwo{\bbeta_{\sS_i}} 
\leq \frac{1}{\sqrt{k}}\normone{\bbeta_{\comple{\sS}}}$ directly, we will fail in general---we need the extra $\normone{\bbeta_{\sS}}$ term as in \eqref{equation:tech_rip_ineq}.}
Use \eqref{equation:tech_rip_ineq} to bound the right-hand side of \eqref{equation:tech_rip_ineq1}, and use $\normone{\bbeta_{\sS}}\leq \sqrt{k}\normtwo{\bbeta_{\sS}}$ (Exercise~\ref{exercise:cauch_sc_l1l2}) to replace $\normone{\bbeta_{\sS}}$:
$$
\begin{aligned}
(1-\delta_{2k})\normtwo{\bbeta_{\sS}}
&\leq \delta_{2k}\cdot\frac{1}{\sqrt{k}}\big(\normone{\bbeta_{\sS}} + \normone{\bbeta_{\comple{\sS}}}\big)\\
&\leq \delta_{2k}\cdot\frac{1}{\sqrt{k}}\big(\sqrt{k}\normtwo{\bbeta_{\sS}} + \normone{\bbeta_{\comple{\sS}}}\big)
= \delta_{2k}\big(\normtwo{\bbeta_{\sS}} + \tfrac{1}{\sqrt{k}}\normone{\bbeta_{\comple{\sS}}}\big).
\end{aligned}
$$
Rearrange to collect the $\normtwo{\bbeta_{\sS}}$ terms and use $\frac{1}{\sqrt{k}}\normone{\bbeta_{\sS}}\leq\normtwo{\bbeta_{\sS}}$ to obtain:
$$
\frac{1}{\sqrt{k}}(1-2\delta_{2k})\normone{\bbeta_{\sS}}
\leq (1-2\delta_{2k})\normtwo{\bbeta_{\sS}}
\leq \frac{\delta_{2k}}{\sqrt{k}}\normone{\bbeta_{\comple{\sS}}}.
$$
Therefore, it holds that 
$$
\normone{\bbeta_{\sS}} \leq \frac{\delta_{2k}}{1-2\delta_{2k}}\,\normone{\bbeta_{\comple{\sS}}}.
$$
So we can take $\rho\triangleq\dfrac{\delta_{2k}}{1-2\delta_{2k}}$. In particular if $\delta_{2k}<\tfrac13$, then $1-2\delta_{2k}>1/3$ and $\rho<1$.
\end{proof}

\begin{remark}\label{remark:rip_roc_ct}
In $\normtwo{\bX \bbeta(\sS)}^2 
=\normtwo{\bX_\sS \bbeta_\sS}^2 
\geq (1 - \delta_{2k}) \normtwo{\bbeta_\sS}^2$, the vector $\bbeta_{\sS}$ is actually $k$-sparse, yet it was treated as $2k$-sparse.
Consequently, a tighter bound could have been used  $\normtwo{\bbeta_{\sS}}^2 \leq \normtwo{\bX\bbeta(\sS)}^2 / (1 - \delta_{k})$ (see Exercise~\ref{exercise:order_delta_rip}). 
The ROC $\theta_{k,k}<\delta_{2k}$ could also have been used instead of $\delta_{2k}$ (see Proposition~\ref{proposition:prop_rip_rop}). 
Using these refinements would lead to the improved sufficient condition $\delta_k + 2\theta_{k,k} < 1$ rather than the more conservative requirement $\delta_{2k}<1/3$.
\end{remark}

Using the lemma above, we now show that the RIP is a sufficient condition for the NSP (Definition~\ref{definition:nullspace_prop}) to hold.
\begin{theoremHigh}[RIP$\implies$NSP]\label{theorem:ripnspI}
If a matrix $\bX$ satisfies RIP($2k,\delta_{2k}$)  with $\delta_{2k} < 1/3$, then it satisfies the NSP of order $k$.
\end{theoremHigh}
\begin{proof}[of Lemma~\ref{theorem:ripnspI}]
By the technical lemma~\ref{lemma:rip2k_rho}: if $\bX$ satisfies RIP$(2k, \delta_{2k})$, then for any nonzero $\bbeta \in \nspace(\bX)$ and $\sT$ being the set of $k$ largest-magnitude entries of $\bbeta$,
$$
\normone{\bbeta_{\sT}} \leq \rho \normone{\bbeta_{\comple{\sT}}},
\quad \rho = \frac{\delta_{2k}}{1 - 2\delta_{2k}}.
$$
If $\delta_{2k} < \frac13$, then $\rho < 1$.

Let $\sS \subset \{1,2, \ldots, p\}$ with $\abs{\sS} = k$ be arbitrary, and let $\sT$ be the set of $k$ largest-magnitude entries of $\bbeta$.  
By construction of $\sT$:
$$
\normone{\bbeta_{\sS}} \leq \normone{\bbeta_{\sT}}, 
\quad
\normone{\bbeta_{\comple{\sS}}} \ge \normone{\bbeta_{\comple{\sT}}}.
$$
Therefore, 
$$
\normone{\bbeta_{\sS}} \leq \normone{\bbeta_{\sT}} \leq \rho \normone{\bbeta_{\comple{\sT}}} \leq \rho \normone{\bbeta_{\comple{\sS}}},
\quad 
\text{for {all} $\sS$ with $\abs{\sS}=k$}.
$$
This establishes the NSP of order $k$ for $\bX$.
\end{proof}

It is worth noting that the RIP($2k,\delta_{2k}$) condition imposes constraints on a very large number of submatrices---specifically, on $\binom{p}{2k}$ of them. In contrast to the pairwise incoherence condition, however, the RIP constant $\delta_{2k}$ itself does not inherently depend on $k$ (beyond the order of the property).

Results from random matrix theory show that various random ensembles satisfy the RIP with high probability, provided the number of measurements $n$ scales as $n \gtrsim k \ln \frac{e p}{k}$. For example, this holds for a standard Gaussian random matrix $\bX$ with i.i.d. entries drawn from $\normal(0, \frac{1}{n})$; see {Corollary~\ref{corollary:gaussben_m_bound_rip}} for details.
Thus, the RIP-based framework enables exact sparse recovery from significantly fewer measurements than what is required by pairwise incoherence---which  typically demands $n \gtrsim k^2 \ln p$ \citep{foucart2013invitation}. On the other hand, a major practical limitation of RIP is that, unlike pairwise incoherence, it is extremely difficult to verify directly, due to the combinatorial number $\binom{p}{2k}$ of submatrices that must be checked.

We next show that the RIP also implies NSP$'$.
We begin with a key lemma that will be used in the proof of Theorem~\ref{theorem:ripnsp}. This result holds for arbitrary vectors $\bbeta\in\real^p$, not only those in the null space of $\bX$. (Note that if $\bbeta\in\nspace(\bX)$, the argument could be significantly simplified.)
\begin{lemma}\label{lemma:ortho_ineq}
Suppose $\balpha$ and $\bbeta$ are orthogonal vectors. Then
$$
\normtwo{\balpha} + \normtwo{\bbeta} \leq \sqrt{2}\normtwo{\balpha+\bbeta}.
$$
\end{lemma}
\begin{proof}[of Lemma~\ref{lemma:ortho_ineq}]
Define the $2 \times 1$ vector $\bu = [\normtwo{\balpha}, \normtwo{\bbeta}]^\top$. 
By the standard bounds on vector norms (Exercise~\ref{exercise:cauch_sc_l1l2}), we have $\normone{\bu} \leq \sqrt{2}\normtwo{\bu}$. This gives 
$$
\normtwo{\balpha} + \normtwo{\bbeta} \leq \sqrt{2}\sqrt{\normtwo{\balpha}^2 + \normtwo{\bbeta}^2}.
$$
Since $\balpha$ and $\bbeta$ are orthogonal, $\normtwo{\balpha}^2 + \normtwo{\bbeta}^2 = \normtwo{\balpha+\bbeta}^2$, which completes the proof.
\end{proof}

The next lemma bounds the sum of the $\ell_2$-norms of successive blocks of $k$ largest entries (in magnitude) of a vector, in terms of the $\ell_1$-norm of its tail. It quantifies how the energy of a vector decays across increasingly less significant components.

\begin{lemma}[Partitioned ordered-set inequality]\label{lemma:ordered_set_ineq}
Let $\sS_0\subseteq \{1,2,\ldots, p\}$ be any index set with $\abs{\sS_0} \leq k$. 
For any vector $\bbeta \in \real^p$, define $\sS_1$ as the indices of the $k$ largest entries (in absolute value) of $\bbeta_{\comple{\sS_0}}$, 
$\sS_2$ as the indices of the next $k$ largest entries, and so on.  
Then
$$
\sum_{i \geq 2} \normtwo{\bbeta_{\sS_i}} 
\leq \frac{\normone{\bbeta_{\comple{\sS_0}}}}{\sqrt{k}}.
$$
\end{lemma}
\begin{proof}[of Lemma~\ref{lemma:ordered_set_ineq}]
By construction, the magnitudes of entries in $\bbeta_{\sS_i}$ are nonincreasing with $i$. In particular, for $i\geq 2$,
$$
\norminf{\bbeta_{\sS_i}} \leq \frac{\normone{\bbeta_{\sS_{i-1}}}}{k},
\quad \text{for $i \geq 2$},
$$
since each entry in $\sS_i$  no larger (in absolute value) than the average magnitude of the entries in $\sS_{i-1}$.
Using this result and applying standard bounds on vector norms (Exercise~\ref{exercise:cauch_sc_l1l2}), we obtain
$$
\sum_{i \geq 2} \normtwo{\bbeta_{\sS_i}}  
\leq \sqrt{k} \sum_{i \geq 2} \norminf{\bbeta_{\sS_i}} 
\leq \frac{1}{\sqrt{k}} \sum_{i \geq 1} \normone{\bbeta_{\sS_i}} 
= \frac{\normone{\bbeta_{\comple{\sS_0}}}}{\sqrt{k}}.
$$
This completes the proof.
\end{proof}

The following result, due to Emmanuel Candes, provides a crucial estimate for the norm of a vector restricted to a support set of size at most $2k$. It plays a central role in connecting the RIP to sparse recovery guarantees.
\begin{lemma}[Candes RIP $\ell_2$-$\ell_1$ bound lemma \citep{candes2008restricted}]\label{lemma:ripnsp_lem2}
Let $\bX\in\real^{n\times p}$ satisfy the RIP of order $2k$ (with $2k\leq p$), and let $\bbeta \in \real^p$, $\bbeta \neq \bzero$ be {arbitrary}. 
Let $\sS_0\subseteq \{1, 2, \ldots, p\}$ satisfy $\abs{\sS_0} \leq k$. 
Define $\sS_1$ as the indices of the $k$ largest entries (in magnitude) of  $\bbeta_{\comple{\sS_0}}$, and set $\sS \triangleq \sS_0 \cup \sS_1$. 
Then,
\begin{equation}
\normtwo{\bbeta_{\sS}} 
\leq \mu \frac{\normone{\bbeta_{\comple{\sS_0}}}}{\sqrt{k}} 
+ \nu \frac{\abs{\innerproduct{\bX_{\sS} \bbeta_{\sS}, \bX \bbeta}}}{\normtwo{\bbeta_{\sS}}},
\end{equation}
where $\mu = {\sqrt{2} \delta_{2k}}/{(1 - \delta_{2k})} $ 
and
$\nu = {1}/{(1 - \delta_{2k})}$.
\end{lemma}

\begin{proof}[of Lemma~\ref{lemma:ripnsp_lem2}]
Since $\bbeta_\sS \in \sB_0[2k]$, the lower RIP bound gives
\begin{equation}\label{equation:ripnsp_lem2_pv1}
(1-\delta_{2k})\normtwo{\bbeta_\sS}^2 \leq \normtwo{\bX_\sS \bbeta_\sS}^2
=
\innerproduct{\bX_\sS \bbeta_\sS, \bX \bbeta} 
- \innerproduct{\bX_\sS \bbeta_\sS, \sum_{i \geq 2} \bX_{\sS_i} \bbeta_{\sS_i}},
\end{equation}
where the sets $\sS_i$ for $i\geq 2$ are defined as in Lemma~\ref{lemma:ordered_set_ineq}, i.e., as $\sS_1$ is the index set corresponding to the $k$ largest entries of $\bbeta_{\comple{\sS_0}}$ in absolute value, define $\sS_2$ as the index set corresponding to the next $k$ largest entries, and so on. The above equality holds since $\bX_\sS \bbeta_\sS = \bX \bbeta - \sum_{i \geq 2} \bX_{\sS_i} \bbeta_{\sS_i}$.
We now bound the second term. Using linearity of the inner product,
\begin{align*}
\abs{\innerproduct{\bX_\sS \bbeta_\sS, \sum_{i \geq 2} \bX_{\sS_i} \bbeta_{\sS_i}}} 
&= \abs{\sum_{i \geq 2} \innerproduct{\bX_{\sS_0} \bbeta_{\sS_0}, \bX_{\sS_i} \bbeta_{\sS_i}} 
	+ \sum_{i \geq 2} \innerproduct{\bX_{\sS_1} \bbeta_{\sS_1}, \bX_{\sS_i} \bbeta_{\sS_i}}} \\
&\leq \sum_{i \geq 2} \abs{\innerproduct{\bX_{\sS_0} \bbeta_{\sS_0}, \bX_{\sS_i} \bbeta_{\sS_i}}} 
+ \sum_{i \geq 2} \abs{\innerproduct{\bX_{\sS_1} \bbeta_{\sS_1}, \bX_{\sS_i} \bbeta_{\sS_i}}} \\
&\stackrel{\dag}{\leq} \delta_{2k} \sum_{i \geq 2} \normtwo{\bbeta_{\sS_0}} \normtwo{\bbeta_{\sS_i}} + \delta_{2k} \sum_{i \geq 2} \normtwo{\bbeta_{\sS_1}} \normtwo{\bbeta_{\sS_i}} \\
&\stackrel{\ddag}{\leq} \sqrt{2}\delta_{2k}  \normtwo{\bbeta_\sS}\sum_{i \geq 2} \normtwo{\bbeta_{\sS_i}}
\leq \sqrt{2}\delta_{2k}  \normtwo{\bbeta_\sS}  \frac{\normone{\bbeta_{\comple{\sS_0}}}}{\sqrt{k}}.
\end{align*}
where the inequality ($\dag$) follows from the standard consequence of RIP (Proposition~\ref{proposition:conse_rip}), 
the inequality ($\ddag$) follows from the fact that  $\normtwo{\bbeta_{\sS_0}} + \normtwo{\bbeta_{\sS_1}} \leq \sqrt{2}\normtwo{\bbeta_\sS}$ by Lemma~\ref{lemma:ortho_ineq}, 
and the last inequality follows from Lemma~\ref{lemma:ordered_set_ineq}.
Substituting this bound back, we get 
\begin{align*}
(1-\delta_{2k})\normtwo{\bbeta_\sS}^2 
&\leq \abs{\innerproduct{\bX_\sS \bbeta_\sS, \bX \bbeta}} + 
\abs{\innerproduct{\bX_\sS \bbeta_\sS, \sum_{i \geq 2} \bX_{\sS_i} \bbeta_{\sS_i}}} \\
&\leq \abs{\innerproduct{\bX_\sS \bbeta_\sS, \bX \bbeta}} + \sqrt{2}\delta_{2k} \normtwo{\bbeta_\sS} \frac{\normone{\bbeta_{\comple{\sS_0}}}}{\sqrt{k}}.
\end{align*}
This proves the desired result after rearrangement. 
\end{proof}

A refined version of this lemma---tailored to bound $\normtwo{\balpha-\bbeta}$ with $\normone{\balpha} \leq \normone{\bbeta}$---is presented in Lemma~\ref{lemma:ell1_gen_diffbount} and will be used to establish noise-free signal recovery guarantees.

We now show that if a matrix satisfies the RIP, then it also satisfies the  NSP$'$ (Definition~\ref{definition:nsp_ii}). 
\begin{theoremHigh}[RIP$\implies$NSP$'$]\label{theorem:ripnsp}
Suppose that $\bX$ satisfies the RIP of order $2k$ with $\delta_{2k} < \sqrt{2} - 1$. Then $\bX$ satisfies the NSP$'$ of order $2k$ with constant
$$
C = \frac{\sqrt{2} \delta_{2k}}{1 - (1 + \sqrt{2}) \delta_{2k}}.
$$
\end{theoremHigh}
\begin{proof}[of Theorem~\ref{theorem:ripnsp}]
Recall that Lemma~\ref{lemma:ripnsp_lem2} holds for arbitrary vectors $\bbeta\in\real^p$.
To prove the theorem, we apply this lemma to the special case where $\bbeta\in\nspace(\bX)$, i.e., $\bX\bbeta=\bzero$.
Let $\sS\subseteq\{1,2,\ldots,p\}$ denote the index set corresponding to the $2k$ largest entries (in absolute value) of $\bbeta$. It suffices to show that
\begin{equation}\label{equation:ripnsp_eq0}
\normtwo{\bbeta_{\sS}} \leq C \frac{\normone{\bbeta_{\comple{\sS}}}}{\sqrt{k}}.
\end{equation}
To this end, let $\sS_0$ be the subset of $\sS$ containing the $k$ largest entries of $\bbeta$, and define $\sS_1=\sS\setminus \sS_0$ (so that $\abs{\sS_1}\leq k$). Applying Lemma~\ref{lemma:ripnsp_lem2} and noting that $\bX\bbeta=\bzero$, the second term in the lemma vanishes, yielding
\begin{equation}\label{equation:ripnsp_eq1}
\normtwo{\bbeta_{\sS}} \leq \mu \frac{\normone{\bbeta_{\comple{\sS_0}}}}{\sqrt{k}},
\end{equation}
where $\mu = {\sqrt{2} \delta_{2k}}/{(1 - \delta_{2k})} $.
Using the standard bounds on vector norms (Exercise~\ref{exercise:cauch_sc_l1l2}), we have
\begin{equation}\label{equation:ripnsp_eq2}
\normonebig{\bbeta_{\comple{\sS_0}}} = \normone{\bbeta_{\sS_1}} + \normonebig{\bbeta_{\comple{\sS}}} \leq \sqrt{k} \normtwo{\bbeta_{\sS_1}} + \normone{\bbeta_{\comple{\sS}}}.
\end{equation}
Combining \eqref{equation:ripnsp_eq1} and \eqref{equation:ripnsp_eq2} leads to 
$$
\normtwo{\bbeta_{\sS}} \leq \mu \left( \normtwo{\bbeta_{\sS_1}} + \frac{\normone{\bbeta_{\comple{\sS}}}}{\sqrt{k}} \right)
\quad\implies\quad 
(1 - \mu) \normtwo{\bbeta_{\sS}} \leq \mu \frac{\normone{\bbeta_{\comple{\sS}}}}{\sqrt{k}},
$$
where the implication follows since $\normtwo{\bbeta_{\sS_1}} \leq \normtwo{\bbeta_{\sS}}$.
The assumption $\delta_{2k} < \sqrt{2} - 1$ ensures that $\mu < 1$, and thus we may divide by $1 - \mu$ without changing the direction of the inequality to establish \eqref{equation:ripnsp_eq0} with constant
$
C \triangleq \frac{\mu}{1 - \mu} = \frac{\sqrt{2} \delta_{2k}}{1 - (1 + \sqrt{2}) \delta_{2k}}.
$
This completes the proof.
\end{proof}

The next argument establishes that the RIP implies the robust NSP (see Definition~\ref{definition:ells_rob_nsp}). 
This proof relies on the following auxiliary result, known as the \textit{square root lifting inequality}. 
It can be viewed as a refined counterpart of the standard vector norm inequality
$\normone{\bbeta} \leq \sqrt{p} \normtwo{\bbeta}$ (see Exercise~\ref{exercise:cauch_sc_l1l2}).
\begin{lemma}[Square root lifting inequality\index{Square root lifting inequality}]\label{lemma:sqrt_lift_ineq}
Let $\bbeta\in\real^p$ satisfy $\beta_1 \geq \beta_2 \geq \ldots \geq \beta_p \geq 0$. 
Then
$$
\sqrt{\beta_1^2 + \beta_2^2 + \ldots + \beta_p^2} \leq \frac{\beta_1 + \beta_2 + \ldots + \beta_p}{\sqrt{p}} + \frac{\sqrt{p}}{4} (\beta_1 - \beta_p).
$$
That is, $\normtwo{\bbeta} \leq \frac{\normone{\bbeta}}{\sqrt{p}} + \frac{\sqrt{p}}{4} (\beta_1 - \beta_p)$.
\end{lemma}
\begin{proof}[of Lemma~\ref{lemma:sqrt_lift_ineq}]
It suffices to prove the equivalent implication:
$$
\left\{
\begin{array}{l}
\beta_1 \geq \beta_2 \geq \ldots \geq \beta_p \geq 0 \\
\frac{\beta_1 + \beta_2 + \ldots + \beta_p}{\sqrt{p}} + \frac{\sqrt{p}}{4} \beta_1 \leq 1
\end{array}
\right\}
\implies \sqrt{\beta_1^2 + \beta_2^2 + \ldots + \beta_p^2} + \frac{\sqrt{p}}{4} \beta_p \leq 1.
$$
Define the convex function
$$
f(\beta_1, \beta_2, \ldots, \beta_p) \triangleq \sqrt{\beta_1^2 + \beta_2^2 + \ldots + \beta_p^2} + \frac{\sqrt{p}}{4} \beta_p
$$
and consider the convex polytope
$$
\sT \triangleq \left\{[\beta_1, \ldots, \beta_p] \in \real^p \mid  \beta_1 \geq \ldots \geq \beta_p \geq 0 \text{ and } \frac{\beta_1 + \beta_2 + \ldots + \beta_p}{\sqrt{p}} + \frac{\sqrt{p}}{4} \beta_1 \leq 1\right\}.
$$
Since $f$ is convex and $\sT$ is a compact convex set, the maximum of $f$ over $\sT$ is attained at an extreme point (vertex) of $\sT$. The vertices arise by setting $p$ of the $p+1$ defining inequalities to equality. We examine the relevant cases:
\begin{itemize}
\item If $\beta_1 = \ldots = \beta_p = 0$, then $f(\beta_1, \beta_2, \ldots, \beta_p) = 0$.
\item If $(\beta_1 + \beta_2 + \ldots + \beta_p)/\sqrt{p} + \sqrt{p} \beta_1/4 = 1$ and $\beta_1 = \ldots = \beta_k > \beta_{k+1} = \ldots = \beta_p = 0$ for some $1 \leq k \leq n-1$, then we have $\beta_1 = \ldots = \beta_k = 4\sqrt{p}/(4k+p)$, and consequently $f(\beta_1, \ldots, \beta_p) = 4\sqrt{p}/(4k+p) \leq 1$.
\item If $(\beta_1 + \beta_2 + \ldots + \beta_p)/\sqrt{p} + \sqrt{p} \beta_1/4 = 1$ and $\beta_1 = \ldots = \beta_p > 0$, then we have $\beta_1 = \ldots = \beta_p = 4/(5\sqrt{p})$, and consequently $f(\beta_1, \ldots, \beta_p) = 4/5 + 1/5 = 1$.
\end{itemize}
Thus, 
$\max_{[\beta_1, \ldots, \beta_p] \in \sT} f(\beta_1, \beta_2, \ldots, \beta_p) = 1$,
which proves the lemma.
\end{proof}

\begin{theoremHigh}[RIP$\implies$robust NSP]\label{theorem:rip2rnsp}
Let $\bX\in\real^{n\times p}$ satisfy the RIP condition of order $2k$ with the RIC
\begin{equation}\label{equation:stab_ell1_rip}
\delta_{2k} < \frac{77 - \sqrt{1337}}{82} \approx 0.4931.
\end{equation} 
Then  $\bX$ satisfies the $\ell_2$-robust NSP of order $k$ (w.r.t. $\norm{\cdot}=\normtwo{\cdot}$) with constants $0 < \rho < 1$ and $\tau > 0$ that depend only on $\delta_{2k}$.
\end{theoremHigh}

\begin{proof}[of Theorem~\ref{theorem:rip2rnsp}]
We need to find constants $0 < \rho < 1$ and $\tau > 0$ such that, for any $\bbeta \in \real^p$ and any $\sS \subseteq \{1,2,\ldots,p\}$ with $\abs{\sS} = k$,
\begin{equation}\label{equ:rip2rnsp_res}
\normtwo{\bbeta_\sS} \leq \frac{\rho}{\sqrt{k}} \normone{\bbeta} + \tau \normtwo{\bX\bbeta}.
\end{equation}
Given $\bbeta \in \real^p$, it suffices to consider an index set $\sS \triangleq \sS_0$ of $k$ largest entries of $\bbeta$ in modulus. 
As before, we partition $\comple{\sS}$ into disjoint blocks $\sS_1, \sS_2, \dots$,
each containing at most $k$ elements, such that the magnitudes of entries within each block are nonincreasing across blocks.
We begin by applying the RIP to the vector $\bbeta(\sS_0) + \bbeta(\sS_1)$, which is $2k$-sparse:
\begin{align}
(1 - \delta_{2k}) (\normtwo{\bbeta(\sS_0) + \bbeta(\sS_1)}^2) 
&\leq \normtwo{\bX(\bbeta(\sS_0) + \bbeta(\sS_1))}^2 
= \normtwobig{\bX\bbeta - \sum_{i \geq 2} \bX\bbeta(\sS_i)}^2 \nonumber \\
&= \normtwobig{\sum_{i \geq 2} \bX\bbeta(\sS_i)}^2 - 2 
\innerproduct{\bX\bbeta, \sum_{i \geq 2} \bX\bbeta(\sS_i)} + \normtwo{\bX\bbeta}^2. \label{equation:rip2robnsp_eq0}
\end{align}
Denote the sum of the last two terms by $\mu$.
For the first term, expand the squared norm:
\begin{align}
\normtwo{\sum_{i \geq 2} \bX\bbeta(\sS_i)}^2 
&= \innerproduct{\sum_{i \geq 2} \bX\bbeta(\sS_i), \sum_{j \geq 2} \bX\bbeta(\sS_j)} \nonumber \\
&= \sum_{i \geq 2} \normtwo{\bX\bbeta(\sS_i)}^2 + \sum_{i,j \geq 2, i \neq j}\innerproduct{\bX\bbeta(\sS_i), \bX\bbeta(\sS_j)} \nonumber \\
&\leq (1 + \delta_{2k}) \sum_{i \geq 2} \normtwo{\bbeta(\sS_i)}^2 + \sum_{i,j \geq 2, i \neq j} \delta_{2k} \normtwo{\bbeta(\sS_i)} \normtwo{\bbeta(\sS_j)} \nonumber \\
&= \sum_{i \geq 2} \normtwo{\bbeta(\sS_i)}^2 + \delta_{2k} \left( \sum_{i \geq 2} \normtwo{\bbeta(\sS_i)} \right)^2, \label{equation:rip2robnsp_eq1}
\end{align}
where the inequality follows from the RIP condition, and the standard consequence of RIP (Proposition~\ref{proposition:conse_rip}).

For each block $\sS_i$, let $\beta_i^-$ and $\beta_i^+$ denote the  smallest and largest absolute values of entries of $\bbeta$ on $\sS_i$.
Define  $\eta \triangleq \sum_{i \geq 2} \normone{\bbeta_{\sS_i}} = \normonebig{\bbeta_{\comple{\sS_0}}} - \normone{\bbeta_{\sS_1}}$. 
Then we can bound:
$$
\sum_{i \geq 2} \normtwo{\bbeta_{\sS_i}}^2 
= \sum_{i \geq 2} \sum_{j \in \sS_i} \abs{\beta_j}^2 
\leq \sum_{i \geq 2} \beta_2^+ \sum_{j \in \sS_i} \abs{\beta_j} 
= \beta_2^+ \eta.
$$
Second, by Lemma~\ref{lemma:sqrt_lift_ineq} and the monotonicity of the blocks (such that $\beta_i^- \geq \beta_{i+1}^+$), it follows that
$$
\sum_{i \geq 2} \normtwo{\bbeta_{\sS_i}} 
\leq \sum_{i \geq 2} \left( \frac{\normone{\bbeta_{\sS_i}}}{\sqrt{k}} + \frac{\sqrt{k}}{4} (\beta_i^+ - \beta_i^-) \right) 
\leq \frac{\eta}{\sqrt{k}} + \frac{\sqrt{k}}{4} \beta_2^+.
$$
Substituting the above two inequalities into \eqref{equation:rip2robnsp_eq1} yields
$$
\normtwo{\sum_{i \geq 2} \bX\bbeta(\sS_i)}^2 \leq \beta_2^+ \eta + \delta_{2k} \left( \frac{\eta}{\sqrt{k}} + \frac{\sqrt{k}}{4} \beta_2^+ \right)^2.
$$
Substituting the above inequality into \eqref{equation:rip2robnsp_eq0}, taking into account the fact that $\normtwo{\bbeta_{\sS_0} + \bbeta_{\sS_1}}^2 = \normtwo{\bbeta_{\sS_0}}^2 + \normtwo{\bbeta_{\sS_1}}^2 \geq \normtwo{\bbeta_{\sS_0}}^2 + k \beta_2^{+2}$, and using the inequality that 
$
\eta = \normonebig{\bbeta_{\comple{\sS_0}}} - \normone{\bbeta_{\sS_1}} 
\leq \normonebig{\bbeta_{\comple{\sS_0}}} - k \beta_2^+ 
= (1 - \nu) \normonebig{\bbeta_{\comple{\sS_0}}},
$ with  $\nu \triangleq k \beta_2^+ / \normonebig{\bbeta_{\comple{\sS_0}}}$
gives
\begin{equation}\label{equation:rip2robnsp_eq3}
\begin{aligned}
(1 &- \delta_{2k}) \normtwo{\bbeta_{\sS_0}}^2  
\leq - (1 - \delta_{2k}) k \beta_2^{+2} + \beta_2^+ \eta + \delta_{2k} \left( \frac{\eta}{\sqrt{k}} + \frac{\sqrt{k}}{4} \beta_2^+ \right)^2 + \mu\\
&\leq - (1 - \delta_{2k}) \frac{\nu^2}{k} \normonebig{\bbeta_{\comple{\sS_0}}}^2 
+ \frac{\nu(1-\nu)}{k} \normonebig{\bbeta_{\comple{\sS_0}}}^2 + \delta_{2k} \left( \frac{1-\nu}{\sqrt{k}} + \frac{\nu}{4\sqrt{k}} \right)^2 \normonebig{\bbeta_{\comple{\sS_0}}}^2 + \mu  \\
&= \left[ - (1 - \delta_{2k}) \nu^2 + \nu(1-\nu) + \delta_{2k} \left( 1 - \frac{3\nu}{4} \right)^2 \right] \frac{\normonebig{\bbeta_{\comple{\sS_0}}}^2}{k} + \mu.
\end{aligned}
\end{equation}
We need to prove that the quantity in square brackets divided by $(1 - \delta_{2k})$ is smaller than one, or equivalently that this quantity minus $(1 - \delta_{2k})$ is negative. This difference is the quadratic expression in $\nu$ given by
$$
f(\nu) \triangleq - \left( 2 - \frac{25}{16} \delta_{2k} \right) \nu^2 + \left( 1 - \frac{3}{2} \delta_{2k} \right) \nu - 1 + 2 \delta_{2k} \triangleq -a\nu^2 + b\nu - c.
$$
This expression is maximized at $\nu^* \triangleq b/(2a)$, so that
$$
f(\nu) \leq f(\nu^*) = \frac{b^2}{4a} - c = \frac{(1 - 3\delta_{2k}/2)^2}{4(2 - 25\delta_{2k}/16)} - 1 + 2\delta_{2k}.
$$
It is easy to show that $f(\nu^*) < 0$ if and only if $41\delta_{2k}^2 - 77\delta_{2k} + 28 > 0$. Thus $f(\nu) < 0$ holds as soon as $\delta_{2k}$ is smaller than the smallest root of $41x^2 - 77x + 28$. This is exactly the condition in \eqref{equation:stab_ell1_rip}. 
Under this condition, the quantity in square brackets of \eqref{equation:rip2robnsp_eq3} is at most $\rho^2 (1 - \delta_{2k})$ for some constant $0 < \rho < 1$ depending only on $\delta_{2k}$. 
To show the desired inequality \eqref{equ:rip2rnsp_res}, by \eqref{equation:rip2robnsp_eq3}, we have 
\begin{equation}\label{equation:rip2robnsp_eq4}
\normtwo{\bbeta_{\sS_0}}^2 \leq \frac{\rho^2 \normone{\bbeta_{\comple{\sS_0}}}^2}{k} + \frac{\mu}{1 - \delta_{2k}}.
\end{equation}
For the estimation of $\mu$, i.e., the sum of last two terms in the right-hand side  of \eqref{equation:rip2robnsp_eq0}, using the Cauchy--Schwarz and  triangle inequalities, we have
$$
\mu \leq 2 \normtwo{\bX\bbeta} \sum_{i \geq 2} \normtwo{\bX\bbeta(\sS_i)} + \normtwo{\bX\bbeta}^2 \leq 2 \sqrt{1 + \delta_{2k}} \normtwo{\bX\bbeta} \sum_{i \geq 2} \normtwo{\bbeta_{\sS_i}} + \normtwo{\bX\bbeta}^2.
$$
By Lemma~\ref{lemma:ells_ellt}, it follows that $\normtwo{\bbeta_{\sS_i}} \leq \normone{\bbeta_{\sS_{i-1}}} / \sqrt{k}$, whence we have 
$$
\sum_{i \geq 2} \normtwo{\bbeta_{\sS_i}} \leq \frac{\normone{\bbeta_{\comple{\sS_0}}}}{\sqrt{k}}.
$$
Substituting the above two inequalities into \eqref{equation:rip2robnsp_eq4}, we have
\begin{align*}
\normtwo{\bbeta_{\sS_0}}^2 &\leq \frac{\rho^2 \normone{\bbeta_{\comple{\sS_0}}}^2}{k} + 2 \frac{\sqrt{1 + \delta_{2k}}}{1 - \delta_{2k}} \normtwo{\bX\bbeta} \frac{\normone{\bbeta_{\comple{\sS_0}}}}{\sqrt{k}} + \frac{\normtwo{\bX\bbeta}^2}{1 - \delta_{2k}} \\
&\leq \left( \frac{\rho^2 \normone{\bbeta_{\comple{\sS_0}}}^2}{k} + 2 \frac{\sqrt{1 + \delta_{2k}}}{1 - \delta_{2k}} \normtwo{\bX\bbeta} \frac{\normone{\bbeta_{\comple{\sS_0}}}}{\sqrt{k}} + \frac{1 + \delta_{2k}}{\rho^2 (1 - \delta_{2k})^2} \normtwo{\bX\bbeta}^2 \right) \\
&= \left( \frac{\rho \normone{\bbeta_{\comple{\sS_0}}} }{\sqrt{k}} + \frac{\sqrt{1 + \delta_{2k}}}{\rho (1 - \delta_{2k})} \normtwo{\bX\bbeta} \right)^2.
\end{align*}
This completes the proof by letting $\tau \triangleq \sqrt{1 + \delta_{2k}} / (\rho (1 - \delta_{2k}))$. 
\end{proof}

\begin{problemset}
	
\item Prove that the set $\sC[k] \triangleq \bigcup_{\sS: \abs{\sS} = k} \sC[\sS]$ is non-convex.

\item \label{prob:lsdist_spar}\textbf{$\ell_s$-distance of $k$-sparse recovery.}
For any $s > t > 0$ and any $\bbeta \in \real^p$, show that
$$
\sigma_k(\bbeta)_s \leq \frac{c_{t,s}}{k^{(\frac{1}{t} - \frac{1}{s})}}\normt{\bbeta},
$$
where 
$
c_{t,s} \triangleq \left[ \left( \frac{t}{s} \right)^{\frac{t}{s}} \left( 1 - \frac{t}{s} \right)^{(1-\frac{t}{s})}  \right]^{\frac{1}{t}} \leq 1.
$
In the special case $t=1$ and $s=2$, the inequality becomes
$$
\sigma_k(\bbeta)_2 \leq \frac{1}{2\sqrt{k}} \normone{\bbeta}.
$$

\item \textbf{Mutual coherence of orthogonal pairs.} Define the {mutual coherence} between two orthogonal bases $ \bU = [\bu_1, \ldots, \bu_p]\in \real^p$ and $ \bV = [\bv_1, \bv_2, \ldots, \bv_p] \in \real^p$  as
$$
\mu(\bU, \bV) \triangleq  \max_{1 \leq i,j \leq p} \abs{\innerproduct{\bu_i, \bv_j}}.
$$
Prove that 
$$
1/\sqrt{p} \le \mu(\bU, \bV) \leq 1.
$$

\item Verify the equivalence of the three statements for a tight frame system in Definition~\ref{definition:equiv_tight}: (i)$\implies$(ii)$\implies$(iii)$\implies$(i).


\item \label{prob:bd_equivsys}\textbf{Bound on equiangular systems.}
Show that the cardinality $p$ of an equiangular system $(\bx_1, \bx_2, \ldots, \bx_p)$ of $\ell_2$-normalized vectors in $\real^n$ satisfies
$$
p \leq 
\frac{n(n+1)}{2}.
$$
Moreover, if equality holds, prove that the system $(\bx_1, \bx_2, \ldots, \bx_p)$ forms a tight frame for $\real^n$.

\item Prove the equivalence between \eqref{equation:def_ric22} and \eqref{equation:def_ric33}.

\item Prove the equivalence between \eqref{equation:def_roc22} and \eqref{equation:roc_def2}.

\item \textbf{Standard consequence of RIP.} For a specific instance of Proposition~\ref{proposition:conse_rip}, 
let $\bX$ satisfy the RIP condition of order $2k$ with constant $\delta_{2k}$, and let $\balpha$ and $\bbeta$ be vectors whose supports are disjoint and whose  support size is $\leq k$ (so the RIP applies to $\balpha\pm\bbeta$).
Prove directly that
$$
\abs{\innerproduct{\bX\balpha, \bX\bbeta}}
\leq 
\delta_{2k} \normtwo{\balpha}\normtwo{\bbeta}.
$$
\textit{Hint:
Do not use the result in Proposition~\ref{proposition:conse_rip}. 
Expand $\normtwo{\bX(\widetilde{\balpha}\pm\widetilde{\bbeta})}^2$ and subtract, where $\widetilde{\balpha}=\balpha/\normtwo{\balpha}$ and $\widetilde{\bbeta}=\bbeta/\normtwo{\bbeta}$; and then scale back.
}

\item \label{prob:cauch_sc_gen} \textbf{Vector norms.}
Let $\bbeta\in\real^p$. Show that 
$$\normt{\bbeta}  \leq p^{(\frac{1}{t} - \frac{1}{s})} \norms{\bbeta}, \quad \text{for }1 \leq t \leq s.$$

\item \label{prob:nonin_nond_func} Show that there is a scalar $\nu^*\in[0,1]$ such that the function defined in \eqref{equation:prop_rip_rop_pe3} is nondecreasing on $[0, \nu^*]$ and then nonincreasing on $[\nu^*, 1]$.

\item  \label{prob:rip_prop1} \textbf{Properties of RIP.} 
Let $\bX \in \real^{n\times p}$ satisfies RIP of order $k$ with constant $\delta_k$. 
Show that 
$$
\normtwo{\bX\bbeta} \leq \sqrt{1 + \delta_k} \left( \normtwo{\bbeta} + \frac{\normone{\bbeta}}{\sqrt{k}} \right),
\quad \text{for any $\bbeta\in\real^p$}.
$$

\item  \label{prob:rip_inv_boud} \textbf{Properties of RIP.} 
Let $\bX\in\real^{n\times p}$ satisfy the RIP of order $k$ with $\delta_k<1$, and let $\sS\subseteq\{1,2,  \ldots,p\}$ with $\abs{\sS}\leq k$.
Show that 
$$
\normtwo{(\bX_\sS^\top \bX_\sS)^{-1}} \leq \frac{1}{1 - \delta_k}.
$$
Also show that 
$$
\normtwo{\bX_\sS(\bX_\sS^\top \bX_\sS)^{-1}} \leq \sqrt{1 + \delta_k}/(1 - \delta_k).
$$
\textit{Hint: The spectral norm of a positive semidefinite or positive definite matrix is the largest eigenvalue of the matrix. Then examine the eigenvalue of the inverse of a matrix. Use the fact that the squared singular value of $\bX$ is the eigenvalue of $\bX^\top\bX$.}

\item \textbf{SSC.} A matrix $\bX\in\real^{n\times p}$ is said to satisfy the $\mu_k$-subset strong convexity (SSC) property and the $\nu_k$-subset smoothness property (SSS) of order $k$, if for all index set $\sS\subset\{1,2,\ldots,n\}$ with  $\abs{\sS}\leq k$, the following inequality holds for all $\bbeta\in\real^p$:
$$
\mu_k\normtwo{\bbeta}^2
\leq \normtwo{\bX^{\sS}\bbeta}^2
\leq \nu_k \normtwo{\bbeta}^2,
$$
where $\bX^\sS\in\real^{\abs{\sS}\times p}$ denotes the submatrix of $\bX$ formed by selecting the rows indexed by $\sS$.
Show that the SSC property implies that no subset $\sS$ with $\abs{\sS}\leq k$ and no pair of distinct vectors  $\bbeta_1,\bbeta_2\in\real^p$ can satisfy  $\bX^\sS\bbeta_1=\bX^\sS\bbeta_2$.
\end{problemset}

\newpage 
\chapter{Sparse Optimization Analysis}\label{chapter:recovery}
\begingroup
\hypersetup{
	linkcolor=structurecolor,
	linktoc=page,  
}
\minitoc \newpage
\endgroup

\lettrine{\color{caligraphcolor}I}
In many real-world applications---from signal processing and imaging to machine learning and genomics---the data we encounter often exhibit sparsity: that is, only a small number of components carry meaningful information. Exploiting this structural property allows us to recover signals or estimate models from limited or noisy observations with high accuracy and computational efficiency. This chapter delves into the theoretical foundations of sparse optimization, focusing on conditions under which sparse signals can be accurately recovered from linear measurements.

We begin in Section~\ref{section:spar_ana_nsp} by examining recovery guarantees based on the nullspace property (NSP), a fundamental condition that ensures exact recovery of sparse signals via convex optimization.  We analyze both noiseless and noisy settings, providing theoretical insights into stability and robustness.

Section~\ref{section:spar_ana_coherence} investigates recovery under coherence, a measure of similarity between columns of the measurement matrix. This approach offers intuitive understanding and practical tools for designing sensing matrices with favorable recovery properties. Moving forward, Section~\ref{section:spar_rec_rip} focuses on the restricted isometry property (RIP), which provides stronger and more versatile guarantees for sparse recovery. We present results for both noiseless and noisy scenarios, including instance-optimal error bounds and an analysis of the Dantzig selector---an alternative  to $\ell_1$-minimization.

The theoretical foundations of $\ell_1$-minimization leads naturally to the discussion of LASSO (least absolute shrinkage and selection operator) in Section~\ref{section:sparse_lin_reg_ana}, a widely used technique in statistics and machine learning for regression with sparsity-inducing regularization. We compare LASSO with robust linear regression methods and examine its behavior under different conditions.

To illustrate the relevance of the design matrix’s properties to sparse recovery, consider the following result:
\begin{theoremHigh}[Exact recovery of $\ell_0$ under RIP\index{Unique recovery}]\label{theorem:recov_rip_l0}
Suppose $\bX\in\real^{n\times p}$ is a design matrix satisfying the RIP condition  of order $2k$ with  $k \geq 1$ and restricted isometry constant $\delta_{2k} < 1$, and let $\sS \subset \{1,2,\ldots,p\} $ be such that $\abs{\sS} \leq k$. Let $\by \triangleq \bX \bbeta^*$ where $\bbeta^*$ is any vector supported on $\sS$ (i.e., $\sS=\supp(\bbeta^*)$). Then $\bbeta^*$ is the \textbf{unique} minimizer of problem \eqref{opt:p0} (p.~\pageref{opt:p0}):
$$
(\text{P}_0) \quad \min \normzero{\bbeta} \quad \text{s.t.}\quad \bX \bbeta = \by.
$$
\end{theoremHigh}
\begin{proof}[of Theorem~\ref{theorem:recov_rip_l0}]
Assume, for contradiction, that  $\by$ admits two distinct $k$-sparse representations: $\by = \bX \bbeta^* = \bX \bbeta'$, where $\bbeta^*$ and $\bbeta'$ are supported on sets $\sS$ and $\sS'$, respectively, with $\abs{\sS}, \abs{\sS'} \leq k$. 
Then
$ \bX (\bbeta^* - \bbeta') = \bzero. $
By construction, the difference $\bbeta^* - \bbeta'$ is supported on $\sS \cup \sS'$, which has cardinality at most $2k$. 
Applying the RIP condition (see Equation~\eqref{equation:rippro22_e1}) and using $\delta_{2k} < 1$, we obtain $\normtwo{\bbeta^* - \bbeta'}^2 = 0$, contradicting the assumption that the two representations are distinct.
\end{proof}

Note that this theorem is an abstract existence result: it establishes what is theoretically possible but does not provide an efficient algorithm to recover $\sS$ or the coefficients $\beta^*_i$ from $\by$ and $\{\bx_i\}_{i \in \{1,2,\ldots,p\}}$, other than exhaustive (brute-force) search---as discussed in Chapter~\ref{chapter:sparse_opt_cond}.
In contrast, as we show next, imposing slightly stronger conditions---such as a stricter bound on the restricted isometry constant $\delta_{2k}$ or alternative structural properties---ensures that the convex $\ell_1$-minimization program~\eqref{opt:p1} recovers $\bbeta^*$ exactly.

For example,
\citet{candes2005decoding} demonstrated that if $\bX \in \text{RIP}(2k, 0.33)$ (Definition~\ref{definition:rip22}), then the non-convex $\ell_0$-norm minimization problem~\eqref{opt:p0} and the convex $\ell_1$-minimization problem~\eqref{opt:p1} (p.~\pageref{opt:p1}) share identical solutions for all $k$-sparse vectors.
Moreover, in the presence of noise ($\by = \bX\bbeta + \bepsilon$), $\ell_1$-minimization stably approximates the sparsest near-solution, with a well-controlled stability constant (see Section~\ref{section:spar_rec_rip}).
Later, \citet{candes2008restricted} established a sharper sufficient condition: $\delta_{2k} < 0.4142$.
Under this bound, $\ell_1$-minimization exactly recovers any $k$-sparse signal in the noiseless case and provides robust estimates in the presence of noise. This chapter discusses further results related to this condition and extends the analysis to LASSO-type problems as well.

\section{Sparse Recovery under NSP}\label{section:spar_ana_nsp}

In Exercise~\ref{exercise:nspi_2k}, we demonstrated that the NSP of order $2k$ is sufficient for the unique recovery of distinct $k$-sparse signals, and that NSP$'$ is necessary when the recovery algorithm satisfies a specific bound as stated in \eqref{equation:nsp-ii-3} of Lemma~\ref{lemms:spar_rec_nspii}. In this section, we investigate the conditions under which the NSP is both sufficient and necessary for the unique recovery of sparse signals in the noiseless case. We further show that in noisy settings, the associated error bounds remain stable.

\subsection{Signal Recovery without Noise}

We begin by establishing the connection between the NSP (Definition~\ref{definition:nullspace_prop}) and the exact recovery of sparse vectors via basis pursuit (i.e., $\ell_1$-minimization).
\begin{theoremHigh}[Exact recovery of $\ell_1$ under NSP($\sS$)\index{Unique recovery}]\label{theorem:exa_ell1_srnp}
Let $\bX \in \real^{n\times p}$. 
Then, {every} vector $\bbeta^*\in\real^p$ supported on a set $\sS$ is the \textbf{unique} solution of \eqref{opt:p1} with $\by=\bX\bbeta^*$ if and only if $\bX$ satisfies the NSP condition relative to $\sS$, i.e., NSP($\sS$).
\end{theoremHigh}
\begin{proof}[of Theorem~\ref{theorem:exa_ell1_srnp}]
Fix an index set $\sS$, and assume that every vector $\bbeta^* \in \real^p$ supported on $\sS$ is the unique minimizer of $\normone{\balpha}$ subject to $\by=\bX\balpha$, where $\by=\bX\bbeta^*=\bX_{\sS}\bbeta^*_{\sS}$. 
Consider any nonzero $\bbeta \in \nspace(\bX) \setminus \{\bzero\}$. 
Then  $\bbeta(\sS)$ satisfies $\bX(-\bbeta(\comple{\sS})) = \bX\bbeta(\sS)$,~\footnote{Note again that $\bbeta(\sS)\in\real^p$ and $\bbeta_{\sS}\in\real^{\abs{\sS}}$; see Definition~\ref{definition:matlabnotation}.}  
and since $\bbeta\neq \bzero$, we have $-\bbeta(\comple{\sS}) \neq \bbeta(\sS)$. 
By the assumed uniqueness of the $\ell_1$ minimizer, it follows that $\normone{\bbeta_{\sS}} < \normone{\bbeta_{\comple{\sS}}}$, which is precisely the NSP($\sS$) condition.

Conversely, suppose the  NSP($\sS$) condition holds. 
Let  $\bbeta \in \real^p$ be supported on $\sS$, and let $\balpha \in \real^p$,  $\balpha \neq \bbeta$, satisfy $\bX\balpha = \bX\bbeta$.
Define $\bu \triangleq \bbeta - \balpha \in \nspace(\bX) \setminus \{\bzero\}$.
Then, using the triangle inequality and the NSP($\sS$) condition,
$$
\normone{\bbeta} \leq \normone{\bbeta_{\sS} - \balpha_{\sS}} + \normone{\balpha_{\sS}} = \normone{\bu_{\sS}} + \normone{\balpha_{\sS}} < \normone{\bu_{\comple{\sS}}} + \normone{\balpha_{\sS}} = \normone{- \balpha_{\comple{\sS}}} + \normone{\balpha_{\sS}} = \normone{\balpha}.
$$
Thus, $\bbeta$ uniquely minimizes the $\ell_1$-norm among all feasible solutions, completing the proof.
\end{proof}

Allowing the support set $\sS$ to vary over all subsets of size at most $k$, we immediately obtain the following corollary.
\begin{corollary}[Exact recovery of $\ell_1$ under  NSP of order $k$]\label{corollary:exa_ell1_snsp}
Let $\bX \in \real^{n\times p}$. 
Then, {every} $k$-sparse vector $\bbeta^* \in \real^p$ is the \textbf{unique} solution of \eqref{opt:p1} with $\by = \bX\bbeta^*$ if and only if $\bX$ satisfies the  NSP of order $k$.
\end{corollary}

This result shows that, whenever the NSP of order $k$ holds, the $\ell_1$-minimization problem \eqref{opt:p1} exactly solves the combinatorial $\ell_0$-minimization problem \eqref{opt:p0} for all $k$-sparse signals. 
Indeed, suppose $\bbeta^*$ is $k$-sparse and recovered uniquely via $\ell_1$-minimization from $\by=\bX\bbeta^*$. Let $\balpha$ be a solution to \eqref{opt:p0}; then $\normzero{\balpha}\leq \normzero{\bbeta^*}\leq k$, so $\balpha$ is also $k$-sparse. By uniqueness of the $\ell_1$ solution for all $k$-sparse vectors, we must have $\balpha=\bbeta^*$.
The following theorem formalizes this equivalence between the $\ell_0$ and $\ell_1$ problems under the NSP.

\begin{theoremHigh}[$\ell_0$ and $\ell_1$ equivalence under NSP]\label{theorem:l0_l1_equirnp}
Suppose that $\bbeta^* \in \real^p$ is the \textbf{unique} solution to the $\ell_0$ problem \eqref{opt:p0} and has support $\sS$. 
Then the basis pursuit problem  \eqref{opt:p1} has a \textbf{unique} solution equal to $\bbeta^*$ if and only if $\bX$ satisfies the NSP condition relative to $\sS$, i.e., NSP($\sS$).
\end{theoremHigh}
\begin{proof}[of Theorem~\ref{theorem:l0_l1_equirnp}]
Assume  $\bX$ satisfies  NSP($\sS$). 
Let $\widehatbbeta \in \real^p$ be the unique solution to problem \eqref{opt:p1}, and define the error vector $\Delta \triangleq \widehatbbeta - \bbeta^*$. Our goal is to show that $\Delta = \bzero$, and in order to do so, it suffices to show that $\Delta \in \nspace(\bX) \cap \sC(\sS)$. 
Since $\bbeta^*$ and $\widehatbbeta$ are optimal (and hence feasible) solutions to the $\ell_0$ and $\ell_1$ problems, respectively, they yield the same measurement $\bX\bbeta^* = \by = \bX\widehatbbeta$, showing that $\bX\Delta = \bzero$, and thus $\Delta\in\nspace(\bX)$. 
Moreover, since $\bbeta^*$ is  feasible for the $\ell_1$-minimization problem \eqref{opt:p1}, the unique optimality of $\widehatbbeta$ implies that $\normonebig{\widehatbbeta} < \normone{\bbeta^*} = \normone{\bbeta^*_{\sS}}+\normone{\bbeta^*_{\comple{\sS}}}$, whence we have
\begin{align*}
\normone{\bbeta^*_{\sS}} 
>  \normonebig{\widehatbbeta}
\geq \normonebig{\widehatbbeta_{\sS}}
= \normone{\bbeta^*_{\sS} + \Delta_{\sS}} + \normone{\Delta_{\comple{\sS}}}
\geq \normone{\bbeta^*_{\sS}} - \normone{\Delta_{\sS}} + \normone{\Delta_{\comple{\sS}}},
\end{align*}
where the  equality follows from the fact that $\widehatbbeta = \bbeta^* + \Delta$ and $\bbeta^*_{\comple{\sS}}=\bzero$, and the last inequality follows by the triangle inequality. Rearranging terms shows that $\Delta \in \sC(\sS)$. 
However, since $\Delta \nspace(\bX) \cap \sC(\sS)\neq \{\bzero\}$ would violate NSP($\sS$), we must have $\Delta=\bzero$. Thus $\widehatbbeta=\bbeta^*$.

Conversely, suppose the problem \eqref{opt:p1} has a unique solution equal to $\bbeta^*\in\real^p$ with support $\sS$. 
For any nonzero  $\bbeta \in \nspace(\bX) \setminus \{\bzero\}$, the vector $\bbeta(\sS)$ is the unique minimizer of $\normone{\balpha}$ subject to $\bX\balpha = \bX\bbeta(\sS)$. 
Since $\bX\bbeta=\bzero$ and $\bbeta \neq \bzero$, we have $\bX(-\bbeta(\comple{\sS})) = \bX\bbeta(\sS)$ and $-\bbeta(\comple{\sS}) \neq \bbeta(\sS)$. 
Hence, it follows that  $\normone{\bbeta_{\sS}} < \normone{\bbeta_{\comple{\sS}}}$, which establishes the NSP($\sS$) condition. 
This completes the proof.
\end{proof}

Note that it is also desirable that any reconstruction method be robust to certain transformations of the measurement system---such as rescaling, reordering, or augmenting rows of $\bX$. Basis pursuit enjoys such invariance properties, as shown in \eqref{equation:nsp_chg}.

To conclude this section, we extend the analysis to the more general  \textit{$\ell_s$-minimization} problem:
\begin{equation}\label{opt:ps}
\min_{\bbeta\in\real^p} \normone{\bbeta}  \qquad\text{s.t.}  \quad \by = \bX \bbeta, \tag{P$_s$} 
\end{equation}
where $0<s\leq 1$.
While the case $s=1$ is most relevant in practice due to convexity, the non-convex problems with $0<s<1$ are of theoretical interest, as they better approximate the $\ell_0$ objective.
Intuitively, smaller values of $s$ should yield solutions closer to the true $\ell_0$ minimizer. To formalize this, we introduce an $\ell_s$-analogue of the NSP.

The proof of the next result mirrors that of Theorem~\ref{theorem:exa_ell1_srnp} and is left as an exercise. It relies on the fact that the function $\bbeta\mapsto \norms{\bbeta}^s$ satisfies the triangle inequality for $0<s\leq 1$; see Problem~\ref{prob:ells_quasnorm_pow}.

\begin{theoremHigh}[Exact recovery of $\ell_s$]\label{theorem:ell_exact1}
Let $\bX \in \real^{n\times p}$ and $0 < s \leq 1$. 
Then every $k$-sparse vector $\bbeta^* \in \real^p$ is the \textbf{unique} solution of \eqref{opt:ps} with $\by = \bX\bbeta^*$ if and only if, for every index set $\sS \subseteq \{1,2,\ldots,p\}$ with $\abs{\sS}\leq k$,
\begin{equation}\label{equation:ell_exact1}
\norms{\balpha_{\sS}} < \norms{\balpha_{\comple{\sS}}}, \quad \text{for all } \balpha \in \nspace(\bX) \setminus \{\bzero\}.
\end{equation}
\end{theoremHigh}

The condition in \eqref{equation:ell_exact1} is actually the $\ell_s$-analogue of the standard NSP condition~\eqref{equation:srnp_e2} in Definition~\ref{definition:nullspace_prop}.

We now show that successful sparse recovery via $\ell_s$-minimization implies successful recovery via $\ell_t$-minimization whenever $0 < t < s \leq 1$---that is, using a smaller exponent yields at least as strong a recovery guarantee.

\begin{theoremHigh}[Exact recovery of $\ell_t<\ell_s$]\label{theorem:ell_exact2}
Let $\bX \in \real^{n \times p}$ and $0 < t < s \leq 1$. 
If every $k$-sparse vector $\bbeta^* \in \real^p$ is the \textbf{unique} solution of \eqref{opt:ps} with $\by = \bX\bbeta^*$, then every $k$-sparse vector $\bbeta^* \in \real^p$ is also the \textbf{unique} solution of (P$_t$) with $\by = \bX\bbeta^*$.
\end{theoremHigh}
\begin{proof}[of Theorem~\ref{theorem:ell_exact2}]
By  Theorem~\ref{theorem:ell_exact1},  it suffices to show the following: for any nonzero $\balpha \in \nspace(\bX) \setminus \{\bzero\}$, if $\sS\subseteq\{1,2,\ldots,p\}$ denotes a set of $k$ indices corresponding to the largest (in magnitude) entries of $\balpha$, then
\begin{equation}\label{theorem:ell_exact2_eq1}
\sum_{j \in \sS} \abs{\alpha_j}^t < \sum_{i \in \comple{\sS}} \abs{\alpha_i}^t,
\end{equation}
provided that the same inequality holds with  $s$ in place of $t$. 

Note that, since $\balpha\neq \bzero$ and $\sS$ contains the $k$ largest entries in absolute value, the complement $\balpha_{\comple\sS}$ cannot be identically zero unless $\balpha$ is $k$-sparse---but even in that case, the NSP condition would require $\balpha=\bzero$, a contradiction. Hence $\balpha_{\comple{\sS}}\neq \bzero$, and the denominator below is positive.

Dividing both sides of \eqref{theorem:ell_exact2_eq1} by $\sum_{i \in \comple{\sS}} \abs{\alpha_i}^t>0$,  the desired inequality is equivalent to
\begin{equation}\label{theorem:ell_exact2_eq2}
\sum_{j \in \sS}\frac{ \abs{\alpha_j}^t}{\sum_{i \in \comple{\sS}} \abs{\alpha_i}^t}
=
\sum_{j \in \sS}  \frac{1}{\sum_{i \in \comple{\sS}}(\abs{\alpha_i}/\abs{\alpha_j})^t} < 1.
\end{equation}
Now observe that  for each $i \in \comple{\sS}$ and $j \in \sS$, the ratio $\abs{\alpha_i}/\abs{\alpha_j} \leq 1$ by construction of $\sS$. 
This makes the left-hand side of \eqref{theorem:ell_exact2_eq2} a nondecreasing function of $0 < t \leq 1$. Hence, its value at $t < s$ does not exceed its value at $s$.
But by assumption and Theorem~\ref{theorem:ell_exact1}, the inequality holds strictly for $s$, i.e., the expression is $<1$ when $t=s$. Therefore, it is also $<1$ for $t<s$, which establishes \eqref{theorem:ell_exact2_eq1}. This completes the proof.
\end{proof}

\subsection{Signal Recovery with Noise}\label{section:noise_nsp}
The ability to perfectly reconstruct a sparse signal from noise-free measurements is a promising theoretical result. However, in most real-world applications, measurements are inevitably contaminated by some form of noise. For example, when processing data on a digital computer, we must represent signals using a finite number of bits, which introduces quantization error. Moreover, physical measurement systems are subject to various types of noise---such as thermal, electronic, or sensor noise---depending on the specific implementation and environment.

Perhaps somewhat surprisingly, it is possible to modify the standard $\ell_1$-minimization program
$$
\widehatbbeta = \argmin_{\bbeta} \normone{\bbeta} \qquad \text{s.t.} \quad \bbeta \in \mathcalB(\by).
$$
so that it stably recovers sparse (or approximately sparse) signals under a wide range of common noise models, e.g., model \eqref{opt:p1_epsilon}  (p.~\pageref{opt:p1_epsilon}) uses $\mathcalB(\by) = \{\bz \mid \normtwo{\bX \bz - \by} \leq \epsilon\}$.

The previous results establish exact recovery for the basis pursuit problem under various NSP conditions.
We then start by using the stable NSP (Definition~\ref{definition:stable_nsp}) or robust NSP (Definition~\ref{definition:ells_rob_nsp}) to show  that any solution lies in a neighborhood of the true vector, with the size of the neighborhood controlled by the best $k$-term approximation error.
\begin{theoremHigh}[Stability of $\ell_1$  under stable NSP]\label{theorem:stablensp}
Let $\bX \in \real^{n\times p}$ satisfy the {stable NSP} of order $k$ with constant $0 < \rho < 1$. 
Then, for any $\bbeta^* \in \real^p$, every solution  $\widehatbbeta$ of \eqref{opt:p1} with measurements $\by = \bX\bbeta^*$ satisfies the $\ell_1$-error bound:
\begin{equation}\label{equation:stablensp}
\normone{\bbeta^* - \widehatbbeta} \leq \frac{2(1+\rho)}{(1-\rho)} \sigma_k(\bbeta^*)_1,
\end{equation}
where $\sigma_k(\bbeta^*)_1 \triangleq \min_{\widehat{\bbeta} \in \sB_0[k]} \normonebig{\bbeta^* - \widehat{\bbeta}}$ denotes the $\ell_1$-error of the best $k$-term approximation.
\end{theoremHigh}
\begin{proof}[of Theorem~\ref{theorem:stablensp}]
Let $\sS\subseteq\{1,2,\ldots,p\}$ be a set of $k$ indices corresponding to the largest (in magnitude) entries of $\bbeta^*$.
Then, by definition, $\normone{\bbeta^*_{\comple{\sS}}} = \sigma_k(\bbeta^*)_1$.
Since $\widehatbbeta$ is a minimizer of \eqref{opt:p1}, it satisfies $\normonebig{\widehatbbeta} \leq \normone{\bbeta^*}$ and $\bX\widehatbbeta = \bX\bbeta^*$.
Applying Theorem~\ref{theorem:suff_nec_snsp} with  $\balpha = \widehatbbeta$ and $\bbeta=\bbeta^*$, we obtain
$$
\normone{\bbeta^* - \widehatbbeta} \leq \frac{2(1+\rho)}{(1-\rho)} \normone{\bbeta^*_{\comple{\sS}}} 
=   \frac{2(1+\rho)}{(1-\rho)} \sigma_k(\bbeta^*)_1,
$$
which is precisely the bound in \eqref{equation:stablensp}.
\end{proof}

In contrast to Theorem~\ref{theorem:exa_ell1_srnp} and Corollary~\ref{corollary:exa_ell1_snsp}, we can no longer guarantee uniqueness of the $\ell_1$-minimizer under the stable NSP alone---although non-uniqueness occurs only in degenerate cases. 
Nevertheless, even when multiple minimizers exist, the theorem ensures that every solution $\widehatbbeta$ of \eqref{opt:p1} with $\by=\bX\bbeta^*$ satisfies the error bound \eqref{equation:stablensp}.

The following theorem extends Theorem~\ref{theorem:stablensp}, recovering it as the special case $\epsilon = 0$ for the quadratically-constrained $\ell_1$-minimization problem \eqref{opt:p1_epsilon}.

\begin{theoremHigh}[Stability of \eqref{opt:p1_epsilon} under robust NSP]\label{theorem:robust_nsp}
Let $\bX \in \real^{n\times p}$ satisfy the {robust NSP} of order $k$ with constants $0 < \rho < 1$, $\tau > 0$, and $\norm{\cdot}=\normtwo{\cdot}$. 
Then, for any $\bbeta^* \in \real^p$, every solution $\widehatbbeta$ of \eqref{opt:p1_epsilon} with measurements $\by = \bX\bbeta^* + \bepsilon$ and noise satisfying $\normtwo{\bepsilon} \leq \epsilon$ obeys the $\ell_1$-error bound:
\begin{equation}\label{equation:robust_nsp}
\normone{\bbeta^* - \widehatbbeta} \leq \frac{2(1+\rho)}{(1-\rho)} \sigma_k(\bbeta^*)_1 + \frac{4\tau}{1-\rho} \epsilon.
\end{equation}
\end{theoremHigh}
\begin{proof}[of Theorem~\ref{theorem:robust_nsp}]
Let $\sS\subseteq\{1,2,\ldots,p\}$ be a set of $k$ indices corresponding to the largest (in magnitude) entries of $\bbeta^*$.
By definition, $\normone{\bbeta^*_{\comple{\sS}}} = \sigma_k(\bbeta^*)_1$.
Since $\widehatbbeta$ is a minimizer of \eqref{opt:p1_epsilon}, it satisfies $\normonebig{\widehatbbeta} \leq \normone{\bbeta^*}$ and $\normbig{\bX\widehatbbeta - \bX\bbeta^*}\leq 2\epsilon$.
Applying Theorem~\ref{theorem:suff_nec_rnsp} (standard consequence of robust NSP) with  $\balpha = \widehatbbeta$ and $\bbeta=\bbeta^*$, we obtain
$$
\normone{\bbeta^* - \widehatbbeta} \leq 
\frac{2(1+\rho)}{(1-\rho)} \normone{\bbeta^*_{\comple{\sS}}} + \frac{2\tau}{1-\rho} \normtwo{\bX(\widehatbbeta-\bbeta^*)}
\leq \frac{2(1+\rho)}{(1-\rho)} \sigma_k(\bbeta^*)_1 + \frac{4\tau}{1-\rho} \epsilon.
$$
which is precisely the bound in \eqref{equation:robust_nsp}.
\end{proof}

Note that when $\supp(\bbeta^*)=\sS$ with $\abs{\sS}=k$, i.e., $\bbeta^*$ is exactly $k$-sparse, the bound becomes 
$$
\normone{\bbeta^* - \widehatbbeta} \leq  \frac{4\tau}{1-\rho} \epsilon.
$$
Moreover, when  $\epsilon=0$, the result reduces to Corollary~\ref{corollary:exa_ell1_snsp}.

The following corollary strengthens the above robustness guarantee by providing error bounds not only in $\ell_1$, but in $\ell_t$  for all
 $1\leq t\leq 2$.

\begin{corollary}[Stability of \eqref{opt:p1_epsilon} under $\ell_2$-robust NSP]\label{corollary:ell2robust_nsp}
Let $\bX \in \real^{n\times p}$ satisfy the {$\ell_2$-robust NSP} of order $k$ with constants $0 < \rho < 1$, $\tau > 0$, and $\norm{\cdot}=\normtwo{\cdot}$. 
Then, for any $\bbeta^* \in \real^p$, every solution $\widehatbbeta$ of \eqref{opt:p1_epsilon} with measurements $\by = \bX\bbeta^* + \bepsilon$ and noise $\normtwo{\bepsilon} \leq \epsilon$ satisfies the $\ell_t$-error bound:
\begin{equation}\label{equation:ell2robust_nsp}
	\normt{\bbeta^* - \widehatbbeta} \leq \frac{C}{k^{(1-\frac{1}{t})}} \sigma_k(\bbeta^*)_1 + D k^{(\frac{1}{t} - \frac{1}{2})} \epsilon, \quad 1 \leq t \leq 2,
\end{equation}
where the constants $C, D > 0$ depend only on $\rho$ and $\tau$ (specifically, $C \triangleq 2(1+\rho)^2/(1-\rho)$ and $D \triangleq 2\tau(3+\rho)/(1-\rho)$).
In particular, the bounds for the extremal cases are:
\begin{align*}
\normone{\bbeta^* - \widehatbbeta} &\leq C \sigma_k(\bbeta^*)_1 + D \sqrt{k} \epsilon, \qquad t=1; \\
\normtwo{\bbeta^* - \widehatbbeta} &\leq \frac{C}{\sqrt{k}} \sigma_k(\bbeta^*)_1 + D \epsilon, \qquad t=2.
\end{align*}
\end{corollary}
\begin{proof}[of Corollary~\ref{corollary:ell2robust_nsp}]
The error bound \eqref{equation:ell2robust_nsp} follows from Theorem~\ref{theorem:lsrob_nsp_gen} (standard consequence of $\ell_s$-robust NSP) as follows: 
Let $\sS$ to be a set of $k$  indices corresponding to the largest (in magnitude) entries of $\bbeta^*$. 
Then $\normone{\bbeta^*_{\comple{\sS}}} = \sigma_k(\bbeta^*)_1$. 
The solution $\widehatbbeta$ of \eqref{opt:p1_epsilon} satisfies $\normonebig{\widehatbbeta} \leq \normone{\bbeta^*}$ and $\normtwobig{\bX\widehatbbeta - \bX\bbeta^*}\leq 2\epsilon$. 
Invoking Theorem~\ref{theorem:lsrob_nsp_gen}  with $\balpha = \widehatbbeta$ and $\bbeta = \bbeta^*$ can therefore prove the desired bound \eqref{equation:ell2robust_nsp}.
\end{proof}

We observe that the coefficient of $\sigma_k(\bbeta^*)_1$ is  constant when $t=1$ and decays  like $1/\sqrt{k}$ when $t=2$, 
whereas  the coefficient of $\epsilon$ grows  like $\sqrt{k}$ for $t=1$ but remains constant for $t=2$.
This reflects the well-known trade-off between sparsity-driven approximation error and noise sensitivity across different norms.

\section{Sparse Recovery under Coherence}\label{section:spar_ana_coherence}

In this section,  
we show that a small $\ell_1$-coherence (see Definition~\ref{definition:l1_coherence}) also guarantees the success of $\ell_1$-minimization. 
\begin{theoremHigh}\label{theorem:exact_rec_l1_l1cohe}
Let $\bX \in \real^{n\times p}$ be a matrix with $\ell_2$-normalized columns. If
\begin{equation}\label{equation:exact_rec_l1_l1cohe_eq1}
\mu_1(\bX,k) + \mu_1(\bX,k-1) < 1,
\end{equation}
then every $k$-sparse vector $\bbeta^* \in \real^p$ is exactly recovered from the measurements $\by = \bX \bbeta^*$ by solving the $\ell_1$-minimization problem  \eqref{opt:p1}.
\end{theoremHigh}
\begin{proof}[of Theorem~\ref{theorem:exact_rec_l1_l1cohe}]
By Corollary~\ref{corollary:exa_ell1_snsp}, it suffices to verify that $\bX$ satisfies the NSP condition of order $k$, i.e., that
\begin{equation}\label{equation:exact_rec_l1_l1cohe_eq2}
\normone{\balpha_{\sS}} < \normone{\balpha_{\comple{\sS}}}
\end{equation}
for every nonzero vector $\balpha \in \nspace(\bX) \setminus \{\bzero\}$ and every  index set $\sS \subseteq \{1,2,\ldots,p\}$ with $\abs{\sS} = k$. 
Since $\bX\balpha=\bzero$, we have $\sum_{j=1}^p \alpha_j \bx_j = \bzero$.
Taking the inner product of both sides with a column $\bx_i$ for some $i \in \sS$, and using the fact that $\normtwo{\bx_i}=1$, we obtain
$$
\alpha_i = \alpha_i \innerproduct{\bx_i, \bx_i}
= -\sum_{j=1, j \neq i}^p \alpha_j \innerproduct{ \bx_j, \bx_i} 
= -\sum_{h \in \comple{\sS}} \alpha_h \innerproduct{\bx_h, \bx_i} - \sum_{j \in \sS, j \neq i} \alpha_j \innerproduct{\bx_j, \bx_i},
$$
whence we have
$$
\abs{\alpha_i} \leq \sum_{h \in \comple{\sS}} \abs{\alpha_h} \abs{\innerproduct{\bx_h, \bx_i}} + \sum_{j \in \sS, j \neq i} \abs{\alpha_j} \abs{\innerproduct{\bx_j, \bx_i}}.
$$
Summing over all $i \in \sS$ and interchanging the order of summation gives
$$
\begin{aligned}
\normone{\balpha_{\sS}} 
&= \sum_{i \in \sS} \abs{\alpha_i} 
\leq \sum_{h \in \comple{\sS}} \abs{\alpha_h} 
\left(\sum_{i \in \sS} \abs{\innerproduct{\bx_h, \bx_i}} \right)
+ \sum_{j \in \sS} \abs{\alpha_j}  
\left(\sum_{i \in \sS, i \neq j} \abs{\innerproduct{\bx_j, \bx_i}}\right)\\
&\leq \sum_{h \in \comple{\sS}} \abs{\alpha_h} \cdot \mu_1(\bX,k) + \sum_{j \in \sS} \abs{\alpha_j} \cdot\mu_1(\bX,k-1) \\
&= \mu_1(\bX,k) \normone{\balpha_{\comple{\sS}}} + \mu_1(\bX,k-1) \normone{\balpha_{\sS}}.
\end{aligned}
$$
Rearranging terms, we obtain $\big(1 - \mu_1(\bX,k-1)\big) \normone{\balpha_{\sS}} \leq \mu_1(\bX,k) \normone{\balpha_{\comple{\sS}}}$.
Under the assumption \eqref{equation:exact_rec_l1_l1cohe_eq1}, we have $1 - \mu_1(\bX,k-1) > \mu_1(\bX,k)$, which implies  \eqref{equation:exact_rec_l1_l1cohe_eq2}. Thus, the NSP of order $k$ holds, and exact recovery via $\ell_1$-minimization follows from Corollary~\ref{corollary:exa_ell1_snsp}.
\end{proof}

\section{Sparse Recovery under RIP}\label{section:spar_rec_rip}

In this section, we establish the success of sparse recovery via basis pursuit for measurement matrices with small restricted isometry constants. We present several proofs of this result.
The first proof uses the condition $\delta_k + \theta_{k,k} + \theta_{k,2k} < 1$, which is derived from the dual certificate theorem (Theorem~\ref{theorem:dualcert_p1}).
The second proof is simpler and more intuitive. It shows that the condition $\delta_{2k} < 1/3$ is sufficient to guarantee exact recovery of all $k$-sparse vectors via $\ell_1$-minimization. 
The remaining proofs are more involved. We show that the weaker conditions $\delta_{2k} < 0.4142$ or $\delta_{2k} < 0.4931$ are actually sufficient to ensure stable and robust recovery of all $k$-sparse vectors via $\ell_1$-minimization.

\index{Dual certificate}
\subsection{Dual Certificate: Revisited}
As stated in Theorem~\ref{theorem:dualcert_p1}, a vector $\bbeta^*$, supported on $\sS=\supp(\bbeta^*)$, is the unique minimizer of \eqref{opt:p1} if the submatrix $\bX_{\sS}$ has full rank and there exists a vector $\blambda$ satisfying the following two properties:
\begin{enumerate}[(i)]
	\item $\innerproduct{\blambda, \bx_i} = \sign(\beta^*_i)$ for all $i \in \sS$;
	\item $\abs{\innerproduct{\blambda, \bx_i}} < 1$, for all $i \in \comple{\sS}$.
\end{enumerate}
This vector $\blambda$ is known as the \textit{dual certificate for problem  \eqref{opt:p1}}.
The next two lemmas demonstrate the existence of such a dual certificate under the restricted isometry property (RIP).
Lemma~\ref{lemma:dual_rec_l2} establishes that the inner products $\innerproduct{\blambda, \bx_i}$ for $i \in \comple{\sS}$ are suitably bounded, except possibly on a small ``exceptional set." Lemma~\ref{lemma:dual_rec_linfty} strengthens this result by eliminating the exceptional set entirely.

\begin{lemma}[Dual reconstruction property, $\ell_2$ version]\label{lemma:dual_rec_l2}
Suppose $\bX\in\real^{n\times p}$ is a design matrix that satisfies the RIP condition or order  $k \geq 1$ with constant $\delta_k < 1$, and let $\bbeta^*$ be a real vector supported on $\sS \subset \{1,2,\ldots,p\}$  with  $\abs{\sS} \leq k$. 
Then there exists a vector $\blambda \in \cspace(\bX)$ such that $\innerproduct{\blambda, \bx_i} = \sign(\beta^*_i)$ for all $i \in \sS$. Furthermore, there exists an ``exceptional set'' $\sT \subset \{1,2,\ldots,p\}$ which is disjoint from $\sS$, with cardinality at most
$ \abs{\sT} \leq k' $ ($k'\geq 1$),
such that
$$ 
\abs{\innerproduct{\blambda, \bx_i}} \leq \frac{\theta_{k,k'}}{(1 - \delta_k)\sqrt{k'}} \cdot \sqrt{k}, \quad \text{for all } i \in \comple{(\sS \cup \sT)}
$$
and
$$ 
\normtwo{\bX_{\sT}\blambda}
=\left( \sum_{i \in \sT} \abs{\innerproduct{\blambda, \bx_i}}^2 \right)^{1/2} 
\leq \frac{\theta_{k,k'}}{1 - \delta_k} \cdot \sqrt{k}.
$$
Additionally, $\normtwo{\blambda} \leq C \sqrt{k}$ for some constant $C > 0$  depending only upon $\delta_k$.
\end{lemma}
\begin{proof}[of Lemma~\ref{lemma:dual_rec_l2}]
The proof differs slightly but is based on  \citet{candes2005decoding}.
The RIP condition \eqref{equation:rippro22_e1} implies that the eigenvalue satisfies $\lambda(\bX_\sS^\top \bX_\sS) \in [ 1 - \delta_k, 1 + \delta_k]$.
Since  we assume $\delta_k < 1$, $\bX_\sS^\top \bX_\sS$ is invertible with
\begin{equation}\label{equation:spa_inv_eig}
\normtwo{(\bX_\sS^\top \bX_\sS)^{-1}} \leq \frac{1}{1 - \delta_k}
\quad\text{ and }\quad 
\normtwo{\bX_\sS(\bX_\sS^\top \bX_\sS)^{-1}} \leq \frac{\sqrt{1 + \delta_k}}{1 - \delta_k}
\end{equation}
by Problem~\ref{prob:rip_inv_boud}. 
Define $\blambda \in \cspace(\bX)$ by
$ \blambda \triangleq \bX_\sS (\bX_\sS^\top \bX_\sS)^{-1} \sign(\bbeta^*_\sS)$.
It is then clear that $\bX_\sS^\top \blambda = \sign(\bbeta^*_\sS)$, i.e., $\innerproduct{\blambda, \bx_i} = \sign(\beta^*_i)$ for all $i \in \sS$. 
Moreover, $\normtwo{\blambda} \leq C \cdot \normtwo{\sign(\beta^*_i)}$ with $C \triangleq \sqrt{1 + \delta_k}/(1 - \delta_k)$ by \eqref{equation:spa_inv_eig} and \eqref{equation:svd_stre_bd}. 
Now let  $\sS'\subseteq\{1,2,\ldots,p\}\setminus \sS$, i.e., $\sS'$ is any set in $ \{1,2,\ldots,p\}$ disjoint from $\sS$ with $\abs{\sS'} \leq k'$, and let $\bd_{\sS'} \triangleq \{d_i\}_{i \in \sS'} \in \real^{\abs{\sS'}}$ be arbitrary. 
Then the ROP condition \eqref{equation:def_roc22_e1} and \eqref{equation:spa_inv_eig} give
$$
\begin{aligned}
\abs{
\innerproduct{\bX_{\sS'}^\top \blambda, \bd_{\sS'}}
} 
&= 
\abs{
\innerproduct{\blambda, \bX_{\sS'} \bd_{\sS'}}
}
= \abs{\innerproduct{\sum_{i \in \sS} \big((\bX_{\sS}^\top \bX_\sS)^{-1} \sign(\bbeta^*_\sS)\big)_i \bx_i, \sum_{i \in \sS'} d_i \bx_i}}\\
&\leq 
\theta_{k,k'} \cdot \normtwo{(\bX_\sS^\top \bX_\sS)^{-1} \sign(\bbeta^*_\sS)} \cdot \normtwo{\bd_{\sS'}}
\leq \frac{\theta_{k,k'}}{1 - \delta_k} \normtwo{\sign(\bbeta^*_\sS)} \cdot \normtwo{\bd_{\sS'}}.
\end{aligned}
$$
Since $\bd_{\sS'}$ is arbitrary, it follows from \holders inequality by setting $\bd_{\sS'} \triangleq{\bX_{\sS'}^\top \blambda}/{\normtwo{\bX_{\sS'}^\top \blambda}}$ that
\begin{equation}\label{equation:dual_rec_l2_eq2}
\normtwo{\bX_{\sS'}^\top \blambda} 
=\left( \sum_{i \in \sS'} \abs{\innerproduct{\blambda, \bx_i}}^2 \right)^{1/2} 
\leq \frac{\theta_{k,k'}}{1 - \delta_k} \normtwo{\sign(\bbeta^*_\sS)}
\end{equation}
whenever $\sS' \subseteq \{1,2,\ldots,p\} \setminus \sS$ and $\abs{\sS'} \leq k'$.
Finally, define the exceptional set
$$
\sT \triangleq \left\{ i \in \{1,2,\ldots,p\} \setminus \sS \mid
\abs{\innerproduct{\blambda, \bx_i}} 
> \frac{\theta_{k,k'}}{(1 - \delta_k)\sqrt{k'}} \cdot \normtwo{\sign(\bbeta^*_\sS)} \right\},
$$
then $\abs{\sT}$ must obey $\abs{\sT} \leq k'$, since otherwise we could contradict \eqref{equation:dual_rec_l2_eq2} by taking a subset $\sS'$ of $\sT$ of cardinality $k'$.
\end{proof}

We now construct a dual vector with improved control over $\abs{\innerproduct{\blambda,\bx_i}}$ for indices $i$ outside the support $\sS$, by iteratively eliminating the exceptional set $\sT$ while keeping the values on $\sS$ fixed.
The following result and its proof are adapted from \citet{candes2005decoding}.

\begin{lemma}[Dual reconstruction property, $\ell_\infty$ version]\label{lemma:dual_rec_linfty}
Suppose $\bX\in\real^{n\times p}$ is a design matrix that satisfies the RIP and ROP condition with  $k \geq 1$ and constants such that  $\delta_k + \theta_{k,k} + \theta_{k,2k} < 1$, 
and let $\bbeta^*$ be a real vector supported on $\sS \subset \{1,2,\ldots,p\}$  with $\abs{\sS} \leq k$. 
Then there exists a vector $\blambda \in \cspace(\bX)$ such that $\innerproduct{\blambda, \bx_i} = \sign(\beta^*_i)$ for all $i \in \sS$. Furthermore, $\blambda$ satisfies
\begin{equation}
\abs{\innerproduct{\blambda, \bx_i}}
\leq \frac{\theta_{k,k}}{(1 - \delta_k - \theta_{k,2k})} <1,
\quad 
\text{for all $i \in \comple{\sS}$.}
\end{equation}
\end{lemma}
\begin{proof}[of Lemma~\ref{lemma:dual_rec_linfty}]
Set $\sS_0 \triangleq\sS $. 
By Lemma~\ref{lemma:dual_rec_l2}, there exists a vector $\blambda^\topone  \triangleq \bX_{\sS_0} (\bX_{\sS_0}^\top \bX_{\sS_0})^{-1} \sign(\bbeta^*_{\sS_0}) \in \cspace(\bX)$ and a set $\sS_1 \subseteq \{1,2,\ldots,p\}$ such that
$$
\begin{aligned}
&\sS_0 \cap \sS_1 = \varnothing \ \text{ with }\ \abs{\sS_1} \leq k; \\
&\innerproduct{\blambda^\topone, \bx_i} = \sign(\beta^*_i), \qquad \text{for all } i \in \sS_0; \\
&\abs{\innerproduct{\blambda^\topone, \bx_i}} \leq \frac{\theta_{k,k}}{(1 - \delta_k)} ,
\qquad \text{for all } i \in \comple{(\sS_0 \cup \sS_1)}; \\
&
\normtwo{\bX_{\sS_1}^\top\blambda^\toptone}
\leq \frac{\theta_{k,k}}{1 - \delta_k} \sqrt{k}; \\
&\normtwobig{\blambda^\topone} \leq C  \sqrt{k}.
\end{aligned}
$$
Since the assumption $\delta_k + \theta_{k,k} + \theta_{k,2k} < 1$ implies that $\delta_{2k} \leq \delta_k +\theta_{k,k}<1$ by Proposition~\ref{proposition:prop_rip_rop} such that $\bX$ satisfies the RIP condition of order $2k$.
We then iteratively construct a set $\sS_t$ with $\abs{\sS_t}\leq k$ for $t\geq 2$.
Let $\sI_t \triangleq \sS_t\cup \sS_0$, and for $i\in\sI_t$ define $\bb^\toptzero\in\real^{\abs{\sI_t}}$ with
$$
b_i^\toptzero
=
\begin{cases}
\innerproduct{\blambda^\toptzero,\bx_i}, \quad &i\in\sS_t;\\
0, \quad& i\in\sS_0.
\end{cases}
$$
Define further $\blambda^\toptone \triangleq  \bX_{\sI_t}(\bX_{\sI_t}^\top\bX_{\sI_t})^{-1}\bb^\toptzero$,
where $(\bX_{\sI_t}^\top\bX_{\sI_t})$ is invertible since $\bX$ satisfies the RIP condition of order $2k$. 
Applying Lemma~\ref{lemma:dual_rec_l2} iteratively gives a sequence of vectors 
$\blambda^\toptone  \triangleq \bX_{\sI_t} (\bX_{\sI_t}^\top \bX_{\sI_t})^{-1} \bb^\toptzero \in \cspace(\bX)$ and sets $\sS_{t+1} \subseteq \{1,2,\ldots,p\}$ for all $t \geq 1$ satisfying
\begin{subequations}
\begin{align}
&\sS_t \cap (\sS_0 \cup \sS_{t+1}) = \varnothing \ \text{ with }\  \abs{\sS_{t+1}} \leq k; \\
&\innerproduct{\blambda^\toptone, \bx_i} = \innerproduct{\blambda^\toptzero, \bx_i}, \qquad \text{for all } i \in \sS_t;
\label{equation:dual_rec_linfty_pva2}  \\
&\innerproduct{\blambda^\toptone, \bx_i} = 0, \qquad \text{for all } i \in \sS_0 ;\\
&\abs{\innerproduct{\blambda^\toptone, \bx_i}} 
\leq \frac{\theta_{k,k}}{1 - \delta_{k}} \left( \frac{\theta_{k,2k}}{1 - \delta_{2k}} \right)^t
\leq \frac{\theta_{k,k}}{1 - \delta_{2k}} \left( \frac{\theta_{k,2k}}{1 - \delta_{2k}} \right)^t, 
\quad \forall\, i \in \comple{(\sS_0 \cup \sS_t \cup \sS_{t+1})}; 
\label{equation:dual_rec_linfty_pva4} \\
&
\normtwo{\bX_{\sS_{t+1}}^\top\blambda^\toptone}
\leq \frac{\theta_{k,k}}{1 - \delta_{k}} \left( \frac{\theta_{k,2k}}{1 - \delta_{2k}} \right)^t \sqrt{k}
\leq \frac{\theta_{k,k}}{1 - \delta_{2k}} \left( \frac{\theta_{k,2k}}{1 - \delta_{2k}} \right)^t \sqrt{k}
; 
\label{equation:dual_rec_linfty_pva5}\\
&\abs{\blambda^\toptone} \leq \left( \frac{\theta_{k,2k}}{1 - \delta_{2k}} \right)^{t-1} C,
\end{align}
\end{subequations}
where we have used the facts that $\delta_k \leq \delta_{2k}$ (Exercise~\ref{exercise:order_delta_rip}).
The inequality \eqref{equation:dual_rec_linfty_pva4} follows since $\abs{\sI_t}\leq 2k$ (the ROC has $\theta_{k,2k}\equiv\theta_{2k,k}$, hence the parameter $\theta_{k,2k}$ instead of $\theta_{k,k}$)  and the iterative argument of $\normtwo{\bb^\toptzero}=\normtwobig{\bX_{\sS_t}^\top\blambda^\toptzero}$:
$$
\normtwo{\bX_{\sS_{t+1}}^\top\blambda^\toptone}
\leq 
\frac{\theta_{k,2k}}{1 - \delta_{2k}}  \normtwo{\bb^\toptzero}
=
\frac{\theta_{k,2k}}{1 - \delta_{2k}}  \normtwo{\bX_{\sS_t}^\top\blambda^\toptzero}.
$$
The inequality \eqref{equation:dual_rec_linfty_pva5} hence follows accordingly.

By hypothesis, we have $\theta_{k,2k}+\delta_{2k} \leq \theta_{k,2k} + \delta_k +\theta_{k,k}<1$ (Proposition~\ref{proposition:prop_rip_rop}) such that  $\frac{\theta_{k,2k}}{1 - \delta_{2k}} < 1$. Thus, if we set
$ \blambda \triangleq \sum_{t=1}^{\infty} (-1)^{t-1} \blambda^\toptzero $, then the series is absolutely convergent and, therefore, $\blambda$ is a well-defined vector in $\cspace(\bX)$. We now examine the coefficients
\begin{equation}\label{equation:dual_rec_linfty_pv1}
\innerproduct{\blambda, \bx_i} = \sum_{t=1}^{\infty} (-1)^{t-1} \innerproduct{\blambda^\toptzero, \bx_i},
\quad \text{for $i \in \{1,2,\ldots,p\}$}.
\end{equation}
On the one hand, for $i \in \sS_0$, it follows from the construction that $\innerproductbig{\blambda^\topone, \bx_i} = \sign(\beta^*_i)$ and $\innerproductbig{\blambda^\toptzero, \bx_i} = 0$ for all $t \geq 2$, and hence,
$$
\innerproduct{\blambda, \bx_i} = \sign(\beta^*_i), \quad \text{for all } i \in \sS_0.
$$
On the other hand, we examine $i$ with $i \in \comple{\sS_0}$. For example, let $i\in\sS_2$, it holds that 
$$
\innerproduct{\blambda, \bx_i} = \sum_{t=1}^{\infty} (-1)^{t-1} \innerproduct{\blambda^\toptzero, \bx_i}
=
\sum_{t\notin\{2,3\}}^{\infty} (-1)^{t-1} \innerproduct{\blambda^\toptzero, \bx_i},
$$
where the term $-\innerproductbig{\blambda^{(2)}, \bx_i} + \innerproductbig{\blambda^{(3)}, \bx_i}$ cancels out by \eqref{equation:dual_rec_linfty_pva2}.
And the other terms satisfy \eqref{equation:dual_rec_linfty_pva4}.
Therefore, $\innerproduct{\blambda, \bx_i}$ follows from a geometric series:
$$
\innerproduct{\blambda, \bx_i}
\leq 
\frac{\theta_{k,k}}{1-\delta_{2k}-\theta_{k,2k}},
\quad 
\text{for all $i \in \comple{\sS}$.}
$$
This completes the proof.
\end{proof}

Lemma~\ref{lemma:dual_rec_linfty} effectively solves the dual reconstruction problem: it guarantees the existence of a vector  $\blambda \in \cspace(\bX)$ satisfying both properties (i) and (ii) stated at the beginning of this section.
Combining this lemma with the dual certificate result in Theorem~\ref{theorem:dualcert_p1} immediately establishes the uniqueness of the solution to the $\ell_1$-minimization problem \eqref{opt:p1}.
Alternatively, the same conclusion can be reached directly using only the RIP. This leads to the following recovery guarantee.

\begin{theoremHigh}[Exact recovery of $\ell_1$ under RIP and ROP \citep{candes2005decoding}]\label{theorem:recov_rip_l1_t1}
Let $\bX\in\real^{n\times p}$ be a design matrix satisfying the  RIP and ROP conditions with  $k \geq 1$ and
\begin{equation}\label{equation:recov_rip_l1_t1_cond}
\delta_k + \theta_{k,k} + \theta_{k,2k} < 1,
\end{equation}
and let $\sS \subset \{1,2,\ldots,p\}$  be any index set with $\abs{\sS} \leq k$.
Let $\by \triangleq \bX \bbeta^*$ where $\bbeta^*$ is an arbitrary vector supported on $\sS$. 
Then $\bbeta^*$ is the \textbf{unique} minimizer of the optimization problem
$$
(\text{P}_1) \qquad \min \normone{\bbeta} \quad \text{s.t.}\quad \bX\bbeta = \by.
$$
\end{theoremHigh}
\begin{proof}[of Theorem~\ref{theorem:recov_rip_l1_t1}]
Standard results from convex analysis ensure that at least one minimizer $\widehatbbeta$ of  \eqref{opt:p1} exists. 
We need to prove that $\widehatbbeta = \bbeta^*$. 
Since $\bbeta^*$ is a feasible point of \eqref{opt:p1} and $\widehatbbeta$ is a minimizer of \eqref{opt:p1}, it follows that 
\begin{equation}\label{equation:recov_rip_l1_t1_e1}
\normonebig{\widehatbbeta} \leq \normone{\bbeta^*} = \sum_{i \in \sS} \abs{\beta^*_i}.
\end{equation}
By Lemma~\ref{lemma:dual_rec_linfty}, there exists a vector $\blambda$ satisfying the properties (i) and (ii) at the beginning of this Section. 
Then, it follows that 
\begin{equation}\label{equation:recov_rip_l1_t1_e2}
\small
\begin{aligned}
\normonebig{\widehatbbeta} 
&= \sum_{i \in \sS} \abs{\beta^*_i + (\widehatbeta_i - \beta^*_i)} + \sum_{i\in\comple{\sS}} \abs{\widehatbeta_i}
\geq \sum_{i \in \sS} \sign(\beta^*_i)\big(\beta^*_i + (\widehatbeta_i - \beta^*_i)\big) + \sum_{i\in\comple{\sS}} \widehatbeta_i 
{\abs{\innerproduct{\blambda, \bx_i}}} \\
&= \normone{\bbeta^*} 
+ \sum_{i \in \sS} (\widehatbeta_i - \beta^*_i) {\abs{\innerproduct{\blambda, \bx_i}}}
+ \sum_{i\in\comple{\sS}} \widehatbeta_i {\abs{\innerproduct{\blambda, \bx_i}}} 
= \normone{\bbeta^*} +  \innerproduct{\blambda, \bX\widehatbbeta - \bX\bbeta^*} \\
&= \normone{\bbeta^*} + \innerproduct{\blambda, \by - \by }
= \normone{\bbeta^*}.
\end{aligned}
\end{equation}
Combining \eqref{equation:recov_rip_l1_t1_e1} and \eqref{equation:recov_rip_l1_t1_e2} yields  $\normonebig{\widehatbbeta} = \normone{\bbeta^*}$.
Since $\abs{\innerproduct{\blambda, \bx_i}}$ is strictly less than 1 for all $i\in\comple{\sS}$, this in particular forces $\widehatbeta_i = 0$ for all $i\in\comple{\sS}$. 
To see this, if $\supp(\widehatbbeta)\neq \sS$, \eqref{equation:recov_rip_l1_t1_e2} also indicates that 
$$
0
=\sum_{i \in \sS} (\beta^*_i - \widehatbeta_i) \innerproduct{\blambda, \bx_i}- \sum_{i\in\comple{\sS}} \widehatbeta_i \innerproduct{\blambda, \bx_i}
=\normone{\bbeta^*} -  
\underbrace{\sum_{i\in\sS} \sign(\beta^*_i)\widehatbeta_i}_{\leq \sum_{i\in\sS} \abs{\widehatbeta_i}} 
-
\underbrace{\sum_{i\in\comple{\sS}} \widehatbeta_i \innerproduct{\blambda, \bx_i}}_{<\sum_{i\in\comple{\sS}} \abs{\widehatbeta_i}}
>0,
$$
leading to a contradiction with the fact that $\normonebig{\widehatbbeta} = \normone{\bbeta^*}$.
Thus, $\supp(\widehatbbeta)= \sS$ and $\by=\bX_{\sS}\widehatbbeta_{\sS} = \bX_{\sS} \bbeta^*_{\sS}$, 

Applying the RIP condition in \eqref{equation:rippro22_e1} and the fact  that $\delta_k < 1$, we conclude that $\widehatbbeta_{\sS}=\bbeta^*_{\sS}$. This proves the uniqueness of $\bbeta^*$.
\end{proof}

Note from Proposition~\ref{proposition:prop_rip_rop} that the hypothesis of \eqref{equation:recov_rip_l1_t1_cond} implies $\delta_{2k} < 1$, and is in turn implied by $\delta_k + \delta_{2k} + \delta_{3k} < 1$. Thus, the condition \eqref{equation:recov_rip_l1_t1_cond}  is roughly ``three times as strict'' as the condition required for Theorem~\ref{theorem:recov_rip_l0} under the $\ell_0$-minimization model ($\delta_{2k} < 1$).
By virtue of the previous discussion, we have the following simple result, which requires $\delta_{2k}<1/3$ and is also roughly ``three times as strict'' as the condition required for Theorem~\ref{theorem:recov_rip_l0}.

\begin{theoremHigh}[Exact recovery of $\ell_1$ under RIP]\label{theorem:ectg_el1_rip}
Let $\bX\in\real^{n\times p}$ satisfy the RIP condition of order $2k$ with RIC
\begin{equation}
\delta_{2k} < \frac{1}{3}.
\end{equation}
Then every $k$-sparse vector $\bbeta^* \in \real^p$ is the \textbf{unique} solution of
$$
\min_{\bbeta \in \real^p} \normone{\bbeta} \quad \text{s.t.} \quad \bX\bbeta = \bX\bbeta^*.
$$
\end{theoremHigh}

\begin{proof}[of Theorem~\ref{theorem:ectg_el1_rip}]
The proof follows simply from Corollary~\ref{corollary:exa_ell1_snsp} and Theorem~\ref{theorem:ripnspI}.
\end{proof}

It is instructive to go beyond exact recovery and establish stability and robustness---i.e., guarantees that persist when the measurements are noisy or the signal is only approximately sparse. The following exercise invites such a refinement.
\begin{exercise}\label{exercise:rob_rip033}
Refine the proof of Theorem~\ref{theorem:ectg_el1_rip} to establish stability and robustness of $k$-sparse recovery via basis pursuit under the condition $\delta_{2k} < 1/3$.
Specifically, show that the reconstruction error remains bounded in terms of the measurement noise level and the best $k$-term approximation error of the true signal.
\textit{Hint: See Theorem~\ref{theorem:stablensp}, or see Section~\ref{sectip:rip_noise} for further details.}
\end{exercise}

\subsection{Signal Recovery without Noise}

Exercise~\ref{exercise:rob_rip033} addressed stable recovery under the RIP of order $2k$ with constant $\delta_{2k}<1/3$.
We now relax this requirement to the weaker condition $\delta_{2k}<\sqrt{2}-1$.
To achieve this, we rely on the following general result, which refines the Candes RIP bound lemma (Lemma~\ref{lemma:ripnsp_lem2}).

\begin{lemma}[Candes RIP bound refinement lemma \citep{candes2008restricted}]\label{lemma:ell1_gen_diffbount}
Suppose that $\bX$ satisfies the RIP condition of order $2k$ with constant $\delta_{2k} < \sqrt{2} - 1$. Let $\balpha,\bbeta, \in \real^p$ be given, and define $\bd \triangleq \balpha - \bbeta$. 
Let $\sS_0$ denote the index set of the $k$ largest-magnitude entries of $\bbeta$, 
and let $\sS_1$ denote the index set of the $k$ largest-magnitude entries of $\bd_{\comple{\sS_0}}$. 
Set $\sS \triangleq  \sS_0 \cup \sS_1$. If $ \normonebig{\balpha} \leq \normone{\bbeta}$, then
\begin{equation}
\normtwo{\bd} 
\leq C_0 \frac{\sigma_k(\bbeta)_1}{\sqrt{k}} 
+ C_1 \frac{\abs{\innerproduct{\bX_{\sS} \bd_{\sS}, \bX \bd}}}{\normtwo{\bd_{\sS}}},
\end{equation}
where
$$
C_0 = 2 \frac{1 - (1 - \sqrt{2}) \delta_{2k}}{1 - (1 + \sqrt{2}) \delta_{2k}}
\qquad \text{and}\qquad 
C_1 = \frac{2}{1 - (1 + \sqrt{2}) \delta_{2k}}.
$$
\end{lemma}
\begin{proof}[of Lemma~\ref{lemma:ell1_gen_diffbount}]
Since $\bd = \bd(\sS) + \bd(\comple{\sS})$, the triangle inequality gives
\begin{equation}\label{equation:ell1_gen_diffbount_pv0}
\normtwo{\bd} \leq \normtwo{\bd_{\sS}} + \normtwo{\bd_{\comple{\sS}}}.
\end{equation}
We first bound $\normtwo{\bd_{\comple{\sS}}}$.
By the partitioned ordered sets inequality (Lemma~\ref{lemma:ordered_set_ineq}), we have
\begin{equation}\label{equation:ell1_gen_diffbount_pv1}
\normtwo{\bd_{\comple{\sS}}} 
= \normtwo{\sum_{i \geq 2} \bd_{\sS_i}} 
\leq \sum_{i \geq 2} \normtwo{\bd_{\sS_i}} 
\leq \frac{\normone{\bd_{\comple{\sS_0}}}}{\sqrt{k}},
\end{equation}
where the sets $\sS_i$ are defined as before, i.e., $\sS_1$ is the index set corresponding to the $k$ largest entries of $\bd_{\comple{\sS_0}}$ (in absolute value), $\sS_2$ as the index set corresponding to the next $k$ largest entries, and so on.

Next, we bound $\normone{\bd_{\comple{\sS_0}}}$. 
Since $ \normonebig{\bbeta+\bd}= \normonebig{\balpha} \leq \normone{\bbeta}$, the triangle inequality yields
\begin{align*}
\normone{\bbeta} 
&\geq \normone{\bbeta+\bd} 
= \normone{\bbeta_{\sS_0} + \bd_{\sS_0}} +\normone{ \bbeta_{\comple{\sS_0}} + \bd_{\comple{\sS_0}}} \nonumber \\
&\geq \normone{\bbeta_{\sS_0}} - \normone{\bd_{\sS_0}} + \normone{\bd_{\comple{\sS_0}}} - \normonebig{\bbeta_{\comple{\sS_0}}}.
\end{align*}
Rearranging terms and recalling that $\sigma_k(\bbeta)_1 = \normonebig{\bbeta_{\comple{\sS_0}}} = \normone{\bbeta - \bbeta(\sS_0)}$,  we obtain
\begin{align}
\normone{\bd_{\comple{\sS_0}}} 
&\leq \normone{\bbeta} - \normone{\bbeta_{\sS_0}} + \normone{\bd_{\sS_0}} + \normonebig{\bbeta_{\comple{\sS_0}}} \nonumber \\
&\leq \normone{\bbeta - \bbeta(\sS_0)} + \normone{\bd_{\sS_0}} + \normonebig{\bbeta_{\comple{\sS_0}}}
= \normone{\bd_{\sS_0}} + 2\sigma_k(\bbeta)_1. \label{equation:ell1_gen_diffbount_pv2}
\end{align}
Substituting \eqref{equation:ell1_gen_diffbount_pv2} into \eqref{equation:ell1_gen_diffbount_pv1} gives
\begin{equation}
\normtwo{\bd_{\comple{\sS}}} 
\leq \frac{\normone{\bd_{{\sS_0}}}+2\sigma_k(\bbeta)_1}{\sqrt{k}} 
\leq \normtwo{\bd_{\sS_0}} + \frac{2\sigma_k(\bbeta)_1}{\sqrt{k}},
\end{equation}
where the last inequality follows from the standard bounds on vector norms (Exercise~\ref{exercise:cauch_sc_l1l2}). 
By observing that $\normtwo{\bd_{\sS_0}} \leq \normtwo{\bd_{\sS}}$, substituting the above inequality into  \eqref{equation:ell1_gen_diffbount_pv0} yields
\begin{equation}\label{equation:ell1_gen_diffbount_pv3}
\normtwo{\bd} \leq 2\normtwo{\bd_{\sS}} + \frac{2\sigma_k(\bbeta)_1}{\sqrt{k}}.
\end{equation}

We now bound  $\normtwo{\bd_{\sS}}$. Combining Candes RIP bound lemma (Lemma~\ref{lemma:ripnsp_lem2}) with \eqref{equation:ell1_gen_diffbount_pv2} and again applying standard bounds on vector norms (Exercise~\ref{exercise:cauch_sc_l1l2}), we obtain
\begin{align}
\normtwo{\bd_{\sS}} 
&\leq \mu \frac{\normone{\bd_{\comple{\sS_0}}}}{\sqrt{k}} + \nu \frac{\abs{\innerproduct{\bX_{\sS} \bd_{\sS}, \bX \bd}}}{\normtwo{\bd_{\sS}}}  
\leq \mu \frac{\normone{\bd_{\sS_0}} + 2\sigma_k(\bbeta)_1}{\sqrt{k}} + \nu \frac{\abs{\innerproduct{\bX_{\sS} \bd_{\sS}, \bX \bd}}}{\normtwo{\bd_{\sS}}} \nonumber \\
&\leq \mu \normtwo{\bd_{\sS_0}} + 2\mu \frac{\sigma_k(\bbeta)_1}{\sqrt{k}} + \nu \frac{\abs{\innerproduct{\bX_{\sS} \bd_{\sS}, \bX \bd}}}{\normtwo{\bd_{\sS}}}, 
\end{align}
where $\mu, \nu$ are defined in Lemma~\ref{lemma:ripnsp_lem2}.
Since $\normtwo{\bd_{\sS_0}} \leq \normtwo{\bd_{\sS}}$, the above inequality becomes 
\begin{equation}
(1-\mu)\normtwo{\bd_{\sS}} \leq 2\mu \frac{\sigma_k(\bbeta)_1}{\sqrt{k}} + \nu \frac{\abs{\innerproduct{\bX_{\sS} \bd_{\sS}, \bX \bd}}}{\normtwo{\bd_{\sS}}}.
\end{equation}
The assumption that $\delta_{2k} < \sqrt{2}-1$ ensures that $\mu =\frac{\sqrt{2} \delta_{2k}}{1 - \delta_{2k}}  < 1$. Dividing by $(1-\mu)$ and substituting the resulting bound into \eqref{equation:ell1_gen_diffbount_pv3} yield
\begin{equation}
\normtwo{\bd} \leq \left(\frac{4\mu}{1-\mu} + 2\right) \frac{\sigma_k(\bbeta)_1}{\sqrt{k}} + \frac{2\nu}{1-\mu} \frac{\abs{\innerproduct{\bX_{\sS} \bd_{\sS}, \bX \bd}}}{\normtwo{\bd_{\sS}}}.
\end{equation}
Using the definitions of $\mu$ and $\nu$ yields the desired constants.
\end{proof}

Lemma~\ref{lemma:ell1_gen_diffbount} establishes an error bound for the class of $\ell_1$-minimization algorithms, assuming the measurement matrix $\bX$ satisfies the RIP. 
To derive concrete recovery guarantees, one must analyze how the constraint $\mathcalB(\by)$ influences the term $\abs{\innerproduct{\bX_{\sS} \bd_{\sS}, \bX \bd}}$. 
As an illustration, consider the noiseless setting where measurements are exact:  $\mathcalB(\by)=\{\balpha\mid \bX\balpha=\by\}$.
In this case, we obtain the following result.

\begin{theoremHigh}[Error bound of \eqref{opt:p1} under RIP \citep{candes2008restricted}]\label{theorem:error_p1}
Suppose that $\bX$ satisfies the RIP of order $2k$ with $\delta_{2k} < \sqrt{2}-1$, and let $\by = \bX\bbeta^*$ for some $\bbeta^*\in\real^p$. 
Then the solution $\widehatbbeta$ to \eqref{opt:p1} satisfies
\begin{equation}
\normtwobig{\widehatbbeta-\bbeta^*} \leq C_0 \frac{\sigma_k(\bbeta^*)_1}{\sqrt{k}},
\end{equation}
where $C_0$ is the constant defined in Lemma~\ref{lemma:ell1_gen_diffbount}.
\end{theoremHigh}
\begin{proof}[of Theorem~\ref{theorem:error_p1}]
Invoking Lemma~\ref{lemma:ell1_gen_diffbount} with $\balpha\triangleq\widehatbbeta$ and $\bbeta\triangleq \bbeta^*$ yields that, for $\bd = \widehatbbeta - \bbeta^*$,
\begin{equation}
\normtwo{\bd} 
\leq C_0 \frac{\sigma_k(\bbeta^*)_1}{\sqrt{k}} 
+ C_1 \frac{\abs{\innerproduct{\bX_{\sS} \bd_{\sS}, \bX \bd}}}{\normtwo{\bd_{\sS}}}.
\end{equation}
Furthermore, since $\bbeta^*, \widehatbbeta$ are feasible (i.e., $\bX\bbeta^*=\bX\widehatbbeta$), we also have that $\bX \bd = \bzero$.
Thus, the second term in the bound vanishes, yielding the claimed inequality. 
\end{proof}

Theorem~\ref{theorem:error_p1} is quite remarkable. In particular, if $\bbeta \in \sB_0[k] \triangleq \{\balpha\mid \normzero{\balpha} \leq k\}$, then---provided $\bX$ satisfies the RIP---we can recover any $k$-sparse vector $\bbeta^*$ exactly.
At first glance, this exact recovery might seem surprising, suggesting that the method could be highly sensitive to noise. However, as we will see next, Lemma~\ref{lemma:ell1_gen_diffbount} can also be used to show that the approach is in fact stable (see Theorem~\ref{theorem:error_p1epsilon}).

Finally, note that Theorem~\ref{theorem:error_p1} assumes the RIP holds for $\bX$. The argument can be easily adapted to instead assume the NSP, using the known implication that RIP implies NSP (see Section~\ref{section:rip_imp_nsp}).

\index{Instance optimality}
\subsection{Instance-Optimal Guarantees}
We now briefly return to the noise-free setting to examine \textit{instance-optimal guarantees} for recovering non-sparse signals more closely.
To begin, we introduce the concept of instance optimality.
\begin{definition}[Instance optimality\index{Instance optimality}]\label{definition:inst_opt}
Given $s \geq 1$, a pair of measurement matrix $\bX \in \real^{n \times p}$ and reconstruction map $\Delta : \real^n \to \real^p$ is called \textit{$\ell_s$-instance optimal of order $k$} with constant $C > 0$ if 
$$
\norms{\bbeta - \Delta(\bX\bbeta)} \leq C \sigma_k(\bbeta)_s,
\quad \text{for all $\bbeta \in \real^p$},
$$
where $\sigma_k(\bbeta)_s \triangleq \min_{\widehat{\bbeta} \in \sB_0[k]} \normsbig{\bbeta - \widehat{\bbeta}}$.
\end{definition}

In Theorem~\ref{theorem:error_p1}, we bounded the $\ell_2$-norm of the reconstruction error of problem~\eqref{opt:p1}
as
\begin{equation}
\normtwobig{\widehatbbeta-\bbeta^*} \leq C_0 \frac{\sigma_k(\bbeta^*)_1 }{\sqrt{k}}.
\end{equation}
This result can be generalized to other norms: for any $s \in [1, 2]$, one can bound the reconstruction error in the $\ell_s$-norm. 
For instance, with minor modifications to the same arguments, it can also be shown that
$$
\normonebig{\widehatbbeta - \bbeta^*} \leq C \sigma_k(\bbeta^*)_1;
$$
see Theorem~\ref{theorem:stab_ell1_rip} and \citet{candes2008restricted} for details.
This naturally raises the question: Can we replace the right-hand side with $\sigma_k(\bbeta^*)_2$, i.e., obtain a bound of the form
$$
\normtwobig{\widehatbbeta-\bbeta^*} \leq C \sigma_k(\bbeta^*)_2?
$$
Unfortunately, such an $\ell_2$-instance optimal guarantee comes at a steep cost in terms of the number of measurements, as quantified by the following theorem.

\begin{theoremHigh}[$\ell_2$-instance optimality \citep{cohen2009compressed}]\label{theorem:optimal_l2_p1_cohen}
Let $\bX\in\real^{n\times p}$ have full rank $n$ ($p\geq n$), and let $\Delta : \real^n \to \real^p$ be a recovery algorithm satisfying
\begin{equation}\label{equation:optimal_l2_p1_cohen}
\normtwo{\bbeta - \Delta(\bX\bbeta)} \leq C \sigma_k(\bbeta)_2, 
\qquad \text{for some $k \geq 1$},
\end{equation}
(i.e., $\ell_2$-instance optimal of order $k$). 
Then the number of measurements must satisfy $n > \left(1 - \sqrt{1 - 1/C^2}\right) p$.
\end{theoremHigh}
\begin{proof}[of Theorem~\ref{theorem:optimal_l2_p1_cohen}]
Let $\bd \in \real^p$ be any vector in the nullspace of $\bX$, i.e., $\bd\in\nspace(\bX)$. 
Decompose  $\bd$ as  $\bd = \bd(\sS) + \bd(\comple{\sS})$, where $\sS\subseteq\{1,2,\ldots,p\}$ is an arbitrary index set with $\abs{\sS} \leq k$. 
Set $\bbeta = \bd(\comple{\sS})$.
Since $\bd \in \nspace(\bX)$, we have  $\bX\bbeta = \bX_{\comple{\sS}} \bd_{\comple{\sS}} = \bX \bd - \bX_{\sS} \bd_{\sS} = -\bX_{\sS} \bd_{\sS}$. 
On the other hand, since $\bd({\sS}) \in \sB_0[k]$, $\sigma_k(\bd({\sS}))_2=0$; 
\eqref{equation:optimal_l2_p1_cohen} implies that $\Delta(\bX\bbeta)\equiv  \Delta(-\bX_{\sS} \bd_{\sS}) = -\bd_{\sS}$. 
Hence, $\normtwo{\bbeta - \Delta(\bX\bbeta)} = \normtwo{\bd_{\comple{\sS}} - (-\bd_{\sS})} = \normtwo{\bd}\leq C\sigma_k(\bbeta)$. 
Furthermore, by definition of $\sigma_k(\bbeta)_2$, we have $\sigma_k(\bbeta)_2 \leq \normtwo{\bbeta - \balpha}$ for all $\balpha \in \sB_0[k]$, including $\balpha = \bzero$. 
Thus, we obtain $\normtwo{\bd} \leq C \normtwo{\bd_{\comple{\sS}}}$ by setting $\balpha=\bzero$. 
Using the fact that  $\normtwo{\bd}^2 = \normtwo{\bd_{\sS}}^2 + \normtwo{\bd_{\comple{\sS}}}^2$, this yields
$$
\normtwo{\bd_{\sS}}^2 = \normtwo{\bd}^2 - \normtwo{\bd_{\comple{\sS}}}^2 \leq \normtwo{\bd}^2 - \frac{1}{C^2} \normtwo{\bd}^2 = \left(1 - \frac{1}{C^2}\right) \normtwo{\bd}^2.
$$
Since the choices of $\bd$ and $\sS$ are arbitrary whenever  $\bd \in \nspace(\bX)$ and  $\abs{\sS} \leq k$; in particular, let $\{\bu_i\}_{i=1}^{p-n}$ be an orthonormal basis for $\nspace(\bX)$, and define the vectors $\{\bd_j\}_{j=1}^p$ as follows:
$$
\bd_j = \sum_{i=1}^{p-n} u_{ij} \bu_i \equiv  \sum_{i=1}^{p-n}\innerproduct{\be_j, \bu_i} \bu_i = \sum_{i=1}^{p-n}(\bu_i\bu_i^\top)\be_j
\triangleq \bP_{\nspace} \be_j,
$$
where $\be_j$ denotes the $j$-th canonical  basis vector,  $u_{ij}$ denotes the $j$-th entry of $\bu_i$, 
and  $\bP_{\nspace}\triangleq \sum_{i=1}^{p-n}(\bu_i\bu_i^\top)$ denotes the orthogonal projection onto $\nspace(\bX)$ (see Definition~\ref{definition:orthogonal-projection-matrix}). 
Since $\normtwo{\bP_{\nspace} \be_j}^2 + \normtwo{\bP_{\nspace}^\perp \be_j}^2 = \normtwo{\be_j}^2 = 1$, we have that $\normtwo{\bd_j} \leq 1$. Thus, by setting $\sS = \{j\}$ for $\bd_j$, we observe that
$$
\abs{\sum_{i=1}^{p-n} (u_{ij})^2}^2 
= (d_{jj})^2 
\leq \left(1 - \frac{1}{C^2}\right) \normtwo{\bd_j}^2 
\leq 1 - \frac{1}{C^2}.
$$
Summing over $j = 1, 2, \ldots, p$, we obtain
$$
p \sqrt{1 - 1/C^2} 
\geq \sum_{j=1}^p \sum_{i=1}^{p-n} (u_{ij})^2 
= \sum_{i=1}^{p-n} \sum_{j=1}^p (u_{ij})^2 
= \sum_{i=1}^{p-n} \normtwo{\bu_i}^2 = p - n,
$$
which completes the proof after rearrangement.
\end{proof}

Thus, if we desire a bound of the form~\eqref{equation:optimal_l2_p1_cohen} that holds for all signals $\bbeta$ with a constant $C \approx 1$, then---regardless of the reconstruction algorithm---we must take nearly as many measurements as the ambient dimension, i.e., $n\approx p$.

\subsection{Signal Recovery with Noise}\label{sectip:rip_noise}

In Section~\ref{section:noise_nsp}, we discussed recovery guarantees under the stable or robust NSP, which ensure robust signal reconstruction in noisy settings.
The RIP also supports such stability and robustness guarantees and plays a central role in compressed sensing theory.
Significant effort has been devoted in the literature to identifying the largest allowable value of the RIC $\delta_{2k}$.
Specifically, \citet{candes2008restricted} showed that $\delta_{2k} < \sqrt{2} - 1 \approx 0.4142$ is sufficient to ensure exact recovery of any $k$-sparse signal. 
This bound was subsequently improved by \citet{foucart2009sparsest} to $\delta_{2k} < 2(3-\sqrt{2})/7\approx0.4531$, 
and later by \citet{foucart2010note} to $\delta_{2k} < 0.4652$. 
Further refinements include the result of \citet{cai2009shifting}, who established that $\delta_{2k} < 0.4721$ suffices when $k$ is a multiple of 4 or when $k$ is very large. 
The same authors (in \citet{foucart2010note}) later pushed this to $\delta_{2k} < 0.4734$ for large values of $k$. 
Most notably, \citet{mo2011new} derived the current best-known general bound, showing that exact recovery is guaranteed whenever $\delta_{2k} < (77 - \sqrt{1337})/{82} \approx 0.4931$ for arbitrary $k$.

We now turn to the analysis of reconstruction error under uniformly bounded noise. The following result provides a worst-case error bound for problem~\eqref{opt:p1_epsilon} and improves upon Theorem~\ref{theorem:ectg_el1_rip} by relaxing the sufficient condition from $\delta_{2k}<1/3$  to the less restrictive $\delta_{2k} < \sqrt{2}-1 \approx 0.4142$.

\begin{theoremHigh}[Error bound of \eqref{opt:p1_epsilon} under RIP \citep{candes2008restricted}]\label{theorem:error_p1epsilon}
Suppose the measurement matrix $\bX$ satisfies the RIP of order $2k$ with constant $\delta_{2k} < \sqrt{2}-1 \approx 0.4142$, and let $\by = \bX\bbeta^* + \bepsilon$ where $\normtwo{\bepsilon} \leq \epsilon$.
Then the solution $\widehatbbeta$ to \eqref{opt:p1_epsilon}  satisfies the error bound
\begin{equation}
\normtwobig{\widehatbbeta-\bbeta^*} \leq C_0 \frac{\sigma_k(\bbeta^*)_1}{\sqrt{k}} + C_2 \epsilon,
\end{equation}
where
$$
C_0 = 2 \frac{1 - (1 - \sqrt{2}) \delta_{2k}}{1 - (1 + \sqrt{2}) \delta_{2k}}
\qquad\text{and}\qquad 
C_2 = 4 \frac{\sqrt{1 + \delta_{2k}}}{1 - (1 + \sqrt{2}) \delta_{2k}}.
$$
\end{theoremHigh}
\begin{proof}[of Theorem~\ref{theorem:error_p1epsilon}]
Let $\bd\triangleq \widehatbbeta-\bbeta^*$. Our goal is to bound $\normtwo{\bd} \triangleq \normtwobig{\widehatbbeta-\bbeta^*}$. 
Since $\normtwo{\bepsilon} \leq \epsilon$, $\bbeta^* \in \mathcalB(\by)= \{\balpha \mid \normtwo{\bX \balpha - \by} \leq \epsilon\}$ is feasible, and $\widehatbbeta$ is the minimizer, we have $\normonebig{\widehatbbeta} \leq \normone{\bbeta^*}$. Thus we may invoke Lemma~\ref{lemma:ell1_gen_diffbount} with $\balpha \triangleq \widehatbbeta$ and $\bbeta\triangleq \bbeta^*$, and it remains to bound $\abs{\innerproduct{\bX_{\sS} \bd_{\sS}, \bX \bd}}$, where $\sS\subseteq\{1,2,\ldots,p\}$ is arbitrary with $\abs{\sS}\leq k$. To do this, using the triangle inequality, we have 
$$
\normtwo{\bX \bd} 
= \normtwobig{\bX (\widehatbbeta - \bbeta^*)}
= \normtwobig{\bX \widehatbbeta - \by + \by - \bX\bbeta^*}
\leq \normtwobig{\bX \widehatbbeta - \by} + \normtwo{\by - \bX\bbeta^*} \leq 2\epsilon,
$$
where the last inequality follows since $\widehatbbeta, \bbeta^* \in \mathcalB(\by)$. Combining this with the RIP and the Cauchy--Schwarz inequality, we obtain
$$
\abs{\innerproduct{\bX_{\sS} \bd_{\sS}, \bX \bd}} 
\leq \normtwo{\bX_{\sS} \bd_{\sS}} \normtwo{\bX \bd} 
\leq 2\epsilon \sqrt{1 + \delta_{2k}} \normtwo{\bd_{\sS}}.
$$
Thus,
$$
\normtwo{\bd} \leq C_0 \frac{\sigma_k(\bbeta^*)_1}{\sqrt{k}} + C_1 2\epsilon \sqrt{1 + \delta_{2k}} 
= C_0 \frac{\sigma_k(\bbeta^*)_1}{\sqrt{k}} + C_2 \epsilon.
$$
This completes the proof.
\end{proof}

Consider we have an \textit{oracle estimator}, which provides the exact support of the sparse vector $\bbeta^*$; that is, the $k$ indices corresponding to its nonzero entries, denoted by $\sS_0$.
In this case, a natural reconstruction strategy is to use the pseudo-inverse restricted to the known support:
\begin{align*}
\widehatbbeta_{\sS_0} &= \bX_{\sS_0}^+ \by = (\bX_{\sS_0}^\top \bX_{\sS_0})^{-1} \bX_{\sS_0}^\top \by, 
\qquad
\widehatbbeta_{\comple{\sS_0}} = \bzero,
\end{align*}
where we assume that  $\bX_{\sS_0}$ has full column rank,  so that there is a unique solution to the equation $\by = \bX_{\sS_0} \bbeta_{\sS_0}$. 
Under this assumption, and given noisy measurements $\by=\bX\bbeta^*+\bepsilon$, the reconstruction error becomes
$$
\normtwobig{\widehatbbeta-\bbeta^*} 
= \normtwo{(\bX_{\sS_0}^\top \bX_{\sS_0})^{-1} \bX_{\sS_0}^\top (\bX\bbeta^* + \bepsilon) - \bbeta^*}
= \normtwo{(\bX_{\sS_0}^\top \bX_{\sS_0})^{-1} \bX_{\sS_0}^\top \bepsilon}.
$$
Note that the  singulars value of $\bX_{\sS_0}$ lie in the range {$[1/\sqrt{1+\delta_{2k}}, 1/\sqrt{1-\delta_{2k}}]$} (see the discussion following Definition~\ref{definition:rip22}). 
Thus, if we consider the worst-case recovery error over all noise vectors $\bepsilon$ satisfying $\norminf{\bX^\top \bepsilon} \leq \epsilon$, then the reconstruction error obeys
\begin{equation}
\frac{\epsilon}{\sqrt{1+\delta_{2k}}} \leq \normtwobig{\widehatbbeta-\bbeta^*} \leq \frac{\epsilon}{\sqrt{1-\delta_{2k}}}.
\end{equation}
Thus, even with perfect knowledge of the true support (i.e., in the oracle setting), the best possible recovery guarantee cannot improve upon the bound in Theorem~\ref{theorem:error_p1epsilon} by more than a constant factor---highlighting the near-optimality of $\ell_1$-minimization under the RIP.

We can further refine the recovery guarantees by leveraging the connection between the RIP and the $\ell_2$-robust NSP. Specifically, using Corollary~\ref{corollary:ell2robust_nsp}, we obtain the following improved result.
\begin{theoremHigh}[Error bound of \eqref{opt:p1_epsilon} under RIP \citep{mo2011new}]\label{theorem:stab_ell1_rip}
Suppose that the measurement matrix  $\bX\in\real^{n\times p}$ satisfies the RIP  of order $2k$ with constant
\begin{equation}\label{equation:stab_ell1_rip2}
\delta_{2k} < \frac{77 - \sqrt{1337}}{82} \approx 0.4931,
\end{equation}
and let the observations be given by $\by = \bX\bbeta^* + \bepsilon$, where $\normtwo{\bepsilon} \leq \epsilon$.
Then any solution $\widehatbbeta$ of problem \eqref{opt:p1_epsilon}  satisfies the error bounds
$$
\normone{\bbeta^* - \widehatbbeta} \leq C \sigma_k(\bbeta^*)_1 + D \sqrt{k} \epsilon,
$$
$$
\normtwo{\bbeta^* - \widehatbbeta} \leq \frac{C}{\sqrt{k}} \sigma_k(\bbeta^*)_1 + D \epsilon,
$$
where the constants $C, D > 0$ depend only on $\delta_{2k}$.
\end{theoremHigh}
\begin{proof}[of Theorem~\ref{theorem:stab_ell1_rip}]
The RIP condition \eqref{equation:stab_ell1_rip2} implies that $\bX$  satisfies the  $\ell_2$-robust NSP of order $k$ (w.r.t. $\norm{\cdot}=\normtwo{\cdot}$) by Theorem~\ref{theorem:rip2rnsp}. 
The stated error bounds then follow directly from  Corollary~\ref{corollary:ell2robust_nsp}.
\end{proof}

In fact, interpolation arguments immediately yield corresponding $\ell_s$-error estimates for all $1 \leq s \leq 2$.

\paragrapharrow{Gaussian noise.}
We can also analyze the performance of these methods in the presence of Gaussian noise  \citep{haupt2006signal}.
To simplify the discussion, we restrict our attention to the case where the true signal $\bbeta^*$ is exactly $k$-sparse, i.e.,  $\bbeta^* \in \sB_0[k] \triangleq \{\balpha\mid \normzero{\balpha} \leq k\}$. 
In this setting, the best $k$-term approximation error vanishes:  $\sigma_k(\bbeta^*)_1 = 0$, and thus the error bounds in Theorem~\ref{theorem:error_p1epsilon} depend only on the noise vector $\bepsilon$.
Assume that the entries of $\bepsilon \in \real^n$ are i.i.d. according to a Gaussian distribution with mean zero and variance $\sigma^2$, i.e., $\epsilon_i\sim\normal(0, \sigma^2)$. 
Since the Gaussian distribution $\normal(0, \sigma^2)$ is itself sub-Gaussian (as a hindsight, see Chapter~\ref{chapter:ensur_rips} for more details), we may apply concentration results such as  Theorem~\ref{theorem:concen_SSG}. 
Specifically, there exists a constant $\kappa_0 > 0$ such that for any $e > 0$,
\begin{equation}
\Pr\left(\normtwo{\bepsilon}^2 \geq (1 + e) n\sigma^2\right) \leq \exp\left(-\kappa_0 ne^2 \right).
\end{equation}
Setting  $e = 3$,  we obtain
\begin{equation}
\Pr\left(\normtwo{\bepsilon} \geq 2 \sqrt{n}\sigma\right) \leq \exp\left(-9\kappa_0 n \right).
\end{equation} 
This leads to the following corollary for the special case of Gaussian noise.
\begin{corollary}[Error bound of \eqref{opt:p1_epsilon} under RIP and Gaussian noise]\label{corollary:errbd_p1epsi_rip_normal}
Suppose that the measurement matrix $\bX$ satisfies the RIP of order $2k$ with $\delta_{2k} < \sqrt{2}-1$. 
Assume further that $\bbeta^* \in \sB_0[k]$ and that the observations are given by $\by = \bX\bbeta^* + \bepsilon$, where the entries of $\bepsilon$ are i.i.d. $\normal(0, \sigma^2)$. 
Then when $\mathcalB(\by) = \{\balpha \mid \normtwo{\bX \balpha - \by} \leq 2 \sqrt{n} \sigma\}$, the solution $\widehatbbeta$ to   \eqref{opt:p1_epsilon} with $\epsilon=2 \sqrt{n} \sigma$ satisfies
\begin{equation}
\normtwobig{\widehatbbeta-\bbeta^*} \leq 8 \frac{\sqrt{1 + \delta_{2k}}}{1 - (1 + \sqrt{2}) \delta_{2k}} \sqrt{n} \sigma
\end{equation}
with probability at least $1 - \exp(-9\kappa_0 n)$.
\end{corollary}
The proof follows directly from Theorem~\ref{theorem:error_p1epsilon} by noting that  $\sigma_k(\bbeta^*)_1=0$ if $\bbeta^*\in\sB_0[k]$.
A similar result can also be derived for Theorem~\ref{theorem:stab_ell1_rip}, yielding analogous $\ell_1$- and $\ell_2$-error bounds under Gaussian noise with high probability.

\subsection{Signal Recovery with Dantzig Selector}
We now consider a slightly different noise model. Whereas Theorems~\ref{theorem:error_p1epsilon} and~\ref{theorem:stab_ell1_rip} assumed that the $\ell_2$-norm of the noise, $\normtwo{\bepsilon}$ is small, the theorem below analyzes the performance of a different recovery algorithm known as the \textit{Dantzig selector} under the assumption that $\norminf{\bX^\top \bepsilon}$ is small \citep{candes2007dantzig}; see Problem~\ref{prob:dantzig_dual_cert}. As we will see, this assumption leads to a particularly simple analysis of the algorithm’s performance in the presence of Gaussian noise.

\begin{theoremHigh}[Error bound of Dantzig selector under RIP\index{Dantzig selector}]\label{theorem:errbd_dantzig_rip}
Suppose that the measurement matrix $\bX$ satisfies the RIP of order $2k$ with $\delta_{2k} < \sqrt{2}-1$, and that we observe measurements of the form $\by = \bX\bbeta^* + \bepsilon$, where $\norminf{\bX^\top \bepsilon} \leq \epsilon$. 
Then when $\mathcalB(\by) = \{\balpha \mid \norminf{\bX^\top (\bX \balpha - \by)} \leq \epsilon\}$, the solution $\widehatbbeta$ to \eqref{opt:ps1} (p.~\pageref{opt:ps1}) satisfies
\begin{equation}
\normtwobig{\widehatbbeta-\bbeta^*} \leq C_0 \frac{\sigma_k(\bbeta^*)_1}{\sqrt{k}} + C_3 \sqrt{k} \epsilon,
\end{equation}
where
\begin{equation}
C_0 = 2 \frac{1 - (1 - \sqrt{2}) \delta_{2k}}{1 - (1 + \sqrt{2}) \delta_{2k}}
\qquad\text{and}\qquad 
C_3 = \frac{4 \sqrt{2}}{1 - (1 + \sqrt{2}) \delta_{2k}}.
\end{equation}
\end{theoremHigh}
\begin{proof}[of Theorem~\ref{theorem:errbd_dantzig_rip}]
The proof closely follows that of Theorem~\ref{theorem:error_p1epsilon}. 
Since $\norminf{\bX^\top \bepsilon} \leq \epsilon$, we again have that $\bbeta^* \in \mathcalB(\by)$, so $\normonebig{\widehatbbeta} \leq \normone{\bbeta^*}$.
Hence, Lemma~\ref{lemma:ell1_gen_diffbount} applies. 
Let $\bd=\widehatbbeta-\bbeta^*$.
As in the earlier proof, we aim to bound $\abs{\innerproduct{\bX_{\sS} \bd_{\sS}, \bX \bd}}$. 
First, by the triangle inequality,
$$
\norminfbig{\bX^\top \bX \bd} 
\leq \norminfbig{\bX^\top (\bX \widehatbbeta - \by)}
+ \norminfbig{\bX^\top (\by - \bX\bbeta^*)} 
\leq 2\epsilon,
$$
where the last inequality holds because both $\widehatbbeta$ and $\bbeta^*$ belong to $\mathcalB(\by)$. 
Applying the Cauchy--Schwarz inequality gives
$$
\abs{\innerproduct{\bX_{\sS} \bd_{\sS}, \bX \bd}} = \abs{\innerproduct{\bd_{\sS}, \bX_\sS^\top \bX \bd}} \leq 
\normtwo{\bd_{\sS}} 
\normtwo{\bX_\sS^\top \bX \bd}.
$$
Finally, since $\norminf{\bX^\top \bX \bd} \leq 2\epsilon$, each entry of $\bX^\top \bX \bd$ is bounded in absolute value by $2\epsilon$. 
Because $\abs{\sS}=\abs{\sS_0}+\abs{\sS_1}=2k$ from Lemma~\ref{lemma:ell1_gen_diffbount}, we obtain $\normtwo{\bX_\sS^\top \bX \bd} \leq \sqrt{2k} (2\epsilon)$. Combining these estimates yields
$$
\normtwo{\bd} \leq C_0 \frac{\sigma_k(\bbeta^*)_1}{\sqrt{k}} + C_1 2 \sqrt{2k} \epsilon 
= C_0 \frac{\sigma_k(\bbeta^*)_1}{\sqrt{k}} + C_3 \sqrt{k} \epsilon,
$$
This completes the proof.
\end{proof}

\paragrapharrow{Gaussian noise for Dantzig selector.}
We can similarly apply Theorem~\ref{theorem:errbd_dantzig_rip} to the case of Gaussian noise.  
Assume that the columns of $\bX$ are normalized to unit $\ell_2$-norm. Then each component of $\bX^\top \bepsilon$ is a Gaussian random variable with mean zero and variance $\sigma^2$.
Using standard Gaussian tail bounds (see Lemma~\ref{lemma:bd_cen_gaus}), we have
\begin{equation}
\Pr\left(\abs{(\bX^\top \bepsilon)_i} \geq t \sigma\right) 
\leq \exp\left(-\frac{t^2}{2}\right),
\quad \text{for $i = 1, 2, \ldots, p$}.
\end{equation}
Applying the union bound over all $p$ coordinates (Theorem~\ref{theorem:union_bound_proof}) and setting $t=2 \sqrt{\ln p}$, we obtain
\begin{equation}
\Pr\left(\norminfbig{\bX^\top \bepsilon} \geq 2 \sqrt{\ln p} \sigma\right) \leq p \exp(-2 \ln p) = \frac{1}{p}.
\end{equation}
Substituting this into Theorem~\ref{theorem:errbd_dantzig_rip}, and noting that  $\sigma_k(\bbeta^*)_1=0$ when  $\bbeta^*\in\sB_0[k]$, we arrive at the following corollary.

\begin{corollary}[Error bound of Dantzig selector under RIP and Gaussian noise \citep{candes2007dantzig}]\label{corollary:errbd_dantzig_rip_normal}
Suppose that the measurement matrix $\bX$ has unit-norm columns and satisfies the RIP of order $2k$ with $\delta_{2k} < \sqrt{2}-1$. Further assume that  $\bbeta^* \in \sB_0[k]$ and that the measurements are given by $\by = \bX\bbeta^* + \bepsilon$, where the entries of $\bepsilon$ are i.i.d. $\normal(0, \sigma^2)$. Then when $\mathcalB(\by) = \{\balpha \mid  \norminf{\bX^\top (\bX \balpha - \by)} \leq 2 \sqrt{\ln p} \sigma\}$, the solution $\widehatbbeta$ to \eqref{opt:ps1} satisfies
\begin{equation}
\normtwobig{\widehatbbeta-\bbeta^*} \leq  \frac{8 \sqrt{2}}{1 - (1 + \sqrt{2}) \delta_{2k}} \sqrt{k \ln p} \sigma,
\end{equation}
with probability at least $1 - \frac{1}{p}$.
\end{corollary}

Ignoring precise constants and the exact success probabilities (which we have not optimized), we observe that if $n = \mathcalO(k \ln p)$, then Corollaries~\ref{corollary:errbd_p1epsi_rip_normal} and~\ref{corollary:errbd_dantzig_rip_normal} yield essentially comparable error bounds. However, there is a subtle but important distinction: for fixed $n$ and $p$, the bound in Corollary~\ref{corollary:errbd_dantzig_rip_normal} adapts to the sparsity level $k$, providing a tighter guarantee when $k$ is small. In contrast, the bound in Corollary~\ref{corollary:errbd_p1epsi_rip_normal} does not improve as $k$ decreases. Thus, while both methods offer similar overall performance, the Dantzig selector may be preferable in settings where the true sparsity $k$ is significantly smaller than the ambient dimension. See \citet{candes2007dantzig} for a more detailed comparison of the relative merits of these approaches.

\index{LASSO}
\section{Sparse Linear Regression}\label{section:sparse_lin_reg_ana}

In many modern data problems---from genomics to finance---we encounter datasets where the number of features $p$ far exceeds the number of observations $n$. 
In such high-dimensional settings, classical linear regression breaks down: it becomes unstable, prone to overfitting, and often produces models that are difficult or impossible to interpret. Sparse linear regression addresses these challenges by seeking solutions in which only a small subset of coefficients are nonzero, effectively performing feature selection while fitting the model.

At its core, sparse linear regression builds on the familiar linear model
$$
\by = \bX\bbeta + \bepsilon,
$$
where $\by \in \real^n$ is the response vector, $\bX \in \real^{n \times p}$ contains the predictors, $\bbeta \in \real^p$ is the coefficient vector we wish to estimate, and $\bepsilon$ represents noise. The key assumption in the sparse setting is that most entries of $\bbeta$ are exactly zero---that is, only a few predictors truly influence the outcome. This assumption not only reflects scientific reality in many domains (e.g., only a handful of genes may affect a disease), but also leads to models that generalize better and are easier to interpret.

The most widely used method for achieving sparsity is the LASSO (least absolute shrinkage and selection operator), introduced by Robert Tibshirani in 1996 \citep{tibshirani1996regression}. The LASSO estimator solves the optimization problem
$$
\widehatbbeta_{\text{lasso}} = \arg\min_{\bbeta} \left\{ \normtwo{\by - \bX\bbeta}^2 + \lambda \normone{\bbeta} \right\},
$$
where $\lambda \geq 0$ is a tuning parameter that controls the degree of sparsity. The $\ell_1$ penalty $\normone{\bbeta} = \sum_{j=1}^p \abs{\beta_j}$ encourages many coefficients to be exactly zero---a property not shared by the $\ell_2$ penalty used in ridge regression. As $\lambda$ increases, more coefficients are shrunk to zero, yielding increasingly sparse models. In practice, $\lambda$ is typically selected via cross-validation to balance prediction accuracy and model simplicity.

Although LASSO is computationally efficient and performs well in many scenarios, it has limitations. For instance, when predictors are highly correlated, LASSO tends to select only one variable from a group, potentially omitting other relevant ones. It can also introduce bias by overshrinking large coefficients. Alternatives such as the \textit{Elastic Net}---which combines $\ell_1$ and $\ell_2$ penalties---help mitigate these issues, especially when features are grouped or highly correlated. At the other extreme, best subset selection directly minimizes prediction error subject to a constraint on the number of nonzero coefficients (i.e., $\normzero{\bbeta} \leq k$). While statistically appealing, this approach is combinatorially hard and generally infeasible for large $p$, though recent advances in optimization have renewed interest in it.

The appeal of sparse regression extends beyond computation. In applications like genome-wide association studies, researchers may measure expression levels for tens of thousands of genes but expect only a few dozen to be biologically relevant. A sparse model not only improves predictive performance but also highlights candidate genes for further experimental validation. Similarly, in signal processing, sparse representations underpin techniques like compressed sensing, where a signal can be accurately reconstructed from far fewer measurements than traditionally required---provided it is sparse in some basis.

From a theoretical perspective, consistent recovery of the true support (i.e., the set of nonzero coefficients) is possible even when  $p \gg n$, provided certain conditions on the design matrix $\bX$ hold. Notable among these are the restricted eigenvalue condition, which essentially require that relevant and irrelevant features are not too strongly correlated. Under such assumptions, methods like LASSO can provably recover the correct sparse structure with high probability as the sample size grows.

Of course, sparsity involves trade-offs. Imposing a sparse structure introduces bias into coefficient estimates, but this is often offset by a substantial reduction in variance---leading to lower overall prediction error. Moreover, while convex formulations like LASSO scale well to large problems, non-convex penalties (e.g., SCAD or MCP \citep{kim2008smoothly, zhang2010nearly}) can offer improved statistical properties at the cost of greater computational complexity and the risk of converging to local minima.

In summary, sparse linear regression transforms an ill-posed high-dimensional problem into a tractable one by embracing the principle of parsimony: the simplest explanation is often the best. This section has focused primarily on the theoretical foundations of the LASSO---particularly the conditions under which it consistently recovers the true sparse model. 
It is worth noting that LASSO is closely related to the broader framework of sparse recovery from compressed sensing and signal processing.
Thus, LASSO can be viewed as the statistical instantiation of $\ell_1$-based sparse recovery, adapted to settings where data are noisy and the primary goals are prediction and consistent support estimation rather than exact reconstruction. By automatically identifying and retaining only the most informative features, it delivers models that are not only accurate but also interpretable---a rare and valuable combination in the era of big data. Whether in biology, economics, engineering, or machine learning, sparse regression remains a foundational tool for turning overwhelming complexity into actionable insight.

\index{Robust linear regression}
\subsection{LASSO and Robust Linear Regression}
The problem of regression analysis is fundamental across many fields, including statistics, supervised machine learning, and optimization.
Minimizing the squared error often yields solutions that are highly sensitive to small changes in the data. 
To mitigate this sensitivity, various regularization methods have been proposed. Among the most well-known and widely cited are Tikhonov regularization ($\ell_2$ penalty) \citep{tikhonov1963solution} and the LASSO \citep{tibshirani1996regression}.

To reduce the computational complexity of solving the optimization problem \eqref{opt:p1} directly, consider a linear regression setting: given an observed response vector  $\by \in \real^n $ and a design matrix $\bX \in \real^{n\times p}$,  the goal is to find a coefficient vector $\bbeta \in \real^p$ such that
\begin{equation}
\widehat{y}_i = \beta_1 x_{i1} + \beta_2 x_{i2} + \cdots + \beta_p x_{ip}, \quad i = 1, 2, \ldots, n,
\end{equation}
or equivalently $\widehat{\by} = \sum_{i=1}^{p} \beta_i \bx_i = \bX \bbeta$.

\citet{tibshirani1996regression} introduced the LASSO   (least absolute shrinkage and selection operator) for linear regression. 
The core idea of LASSO is to minimize the sum of squared prediction errors subject to an upper bound $\Sigma$ on the $\ell_1$-norm of  the coefficient vector, i.e., model~\eqref{opt:lc} (p.~\pageref{opt:lc}):
\begin{equation}\label{equation:lc_ana}
\text{(L$_{C}$)} \qquad \min_{\bbeta} \normtwo{\by - \bX \bbeta}^2 \qquad \text{s.t.} \quad \normone{\bbeta} = \sum_{i=1}^{p} \abs{\beta_i} \leq \Sigma.
\end{equation}
The bound $\Sigma$ is a tuning parameter. 
When $\Sigma$ is sufficiently large, the constraint becomes inactive, and the solution reduces to ordinary least squares regression. However, for smaller values of $\Sigma$, some coefficients are driven exactly to zero, resulting in a sparse estimate of $\bbeta$. 
Selecting an appropriate $\Sigma$ thus enables automatic variable selection.

An equivalent formulation uses the Lagrangian form of the LASSO, i.e., model~\eqref{opt:ll} (p.~\pageref{opt:ll}):
\begin{equation}\label{equation:ll_ana}
\text{(L$_{L}$)} \qquad \min_{\bbeta} \left\{ \normtwo{\by - \bX \bbeta}^2 + \lambda \normone{\bbeta} \right\},
\end{equation}
where $\lambda \geq 0$ controls the strength of regularization and is related to $\Sigma$ (see Corollary~\ref{corollary:equigv_p1_epsi_pen_lag}).

The LASSO combines $\ell_1$-norm regularization with linear regression and exhibits two key properties:
\begin{enumerate}[(i)]
\item \textit{Contraction effect}. It shrinks the estimated coefficients toward zero, reducing model complexity.
\item \textit{Selection effect}. It automatically sets some coefficients exactly to zero, thereby performing feature selection and yielding a sparse solution.
\end{enumerate}
Consequently, LASSO serves as both a shrinkage and a variable selection method. It is closely related to soft-thresholding of wavelet coefficients and forward stagewise regression; see Section~\ref{section:lars}.

In practical applications, the observed design matrix $\bX$ may be corrupted by noise or even adversarial perturbations. To ensure robustness against such uncertainties, \textit{robust linear regression} considers the worst-case scenario by solving the following minimax problem:
\begin{equation}\label{equation:robustreg_prob}
\min_{\bbeta \in \real^p} \max_{\Delta \bX \in \sU} \normtwo{\by - (\bX + \Delta \bX) \bbeta},
\end{equation}
where $\Delta \bX = [\bdelta_1, \bdelta_2,\ldots, \bdelta_p]$ represents the perturbation applied to each column (feature) of $\bX$, and
\begin{equation}\label{equation:robustreg_prob_unc_set}
\sU \triangleq \big\{ [\bdelta_1,\bdelta_2, \ldots, \bdelta_p] \mid  \normtwo{\bdelta_i} \leq c_i, \ i = 1, 2,\ldots, p  \big\}
\end{equation}
is called the \textit{uncertainty set}, or the set of admissible perturbations of the matrix $\bX$.
This uncertainty set is featurewise uncoupled, meaning that perturbations to different features are constrained independently---unlike coupled uncertainty sets, which impose joint constraints across features.

Remarkably, when the uncertainty set is uncoupled and bounded in the $\ell_2$-norm per feature, the robust regression problem admits a tractable reformulation, as stated in the following theorem.

\begin{theoremHigh}[Robust regression \citep{xu2008robust}\index{Robust regression}]\label{theorem:rob_reg1}
The robust regression problem \eqref{equation:robustreg_prob} with uncertainty set $\sU$ given by \eqref{equation:robustreg_prob_unc_set} is equivalent to the following $\ell_1$-regularized regression problem:
\begin{equation}\label{equation:rob_reg1}
\min_{\bbeta \in \real^p} \left\{ \normtwo{\by - \bX \bbeta} + \sum_{i=1}^{p} c_i \abs{\beta_i} \right\}.
\end{equation}
\end{theoremHigh}
\begin{proof}[of Theorem~\ref{theorem:rob_reg1}]
Fix $\bbeta=\bbeta^*$. We show  that 
$$
\max_{\Delta \bX \in \sU} \normtwo{\by - (\bX + \Delta \bX) \bbeta^*} 
= \normtwo{\by - \bX \bbeta^*} + \sum_{i=1}^p c_i \abs{\beta_i^*}.
$$
By the triangle inequality, the left-hand side can be written as
\begin{equation}\label{equation:rob_reg1_eq1}
\small
\begin{aligned}
& \max_{\Delta \bX \in \sU} \normtwobig{\by - (\bX + \Delta \bX) \bbeta^*}\\
&= \max_{[\bdelta_1,\ldots, \bdelta_p] :  \normtwo{\bdelta_i} \leq c_i} \normtwobig{\by - \left( \bX + [\bdelta_1, \ldots, \bdelta_p] \right) \bbeta^*} 
= \max_{[\bdelta_1,\ldots, \bdelta_p] :  \normtwo{\bdelta_i} \leq c_i} \normtwo{ \by - \bX \bbeta^* - \sum_{i=1}^p \beta_i^* \bdelta_i} \\
&\leq \max_{[\bdelta_1,\ldots, \bdelta_p] :  \normtwo{\bdelta_i} \leq c_i} \left\{\normtwo{\by - \bX \bbeta^*} + \sum_{i=1}^p \normtwo{\beta_i^* \bdelta_i} \right\}
\leq \normtwo{\by - \bX \bbeta^*} + \sum_{i=1}^p \abs{\beta_i^*} c_i.
\end{aligned}
\end{equation}
To establish equality, define a unit vector $\bv $ as
$$
\bv \triangleq 
\begin{cases} 
{(\by - \bX \bbeta^*)}/{\normtwo{\by - \bX \bbeta^*}}, & \text{if } \bX \bbeta^* \neq \by; \\
\text{any unit norm vector}, & \text{otherwise},
\end{cases}
\qquad
\text{and}
\qquad 
\bdelta_i^* \triangleq -c_i \sgn(\beta_i^*) \bv, \,\forall\, i.
$$
Since $\normtwo{\bdelta_i^*} \leq c_i$, it follows that $\Delta \bX^* \triangleq [\bdelta_1^*, \bdelta_2^*, \ldots, \bdelta_p^*] \in \sU$, whence we have 
\begin{equation}\label{equation:rob_reg1_eq2}
\small
\begin{aligned}
& \max_{\Delta \bX \in \sU} \normtwo{\by - (\bX + \Delta \bX) \bbeta^*}
\geq \normtwo{\by - (\bX + \Delta \bX^*) \bbeta^*} \\
& = \normtwo{\by - \left( \bX + [\bdelta_1^*, \bdelta_2^*, \ldots, \bdelta_p^*] \right) \bbeta^*}
= \normtwobig{\left( \by - \bX \bbeta^* \right) - \sum_{i=1}^p \left( - \beta_i^* c_i \sgn(\beta_i^*) \bv \right)} \\
& = \normtwobig{\left( \by - \bX \bbeta^* \right) + \left( \sum_{i=1}^p c_i \abs{\beta_i^*} \right) \bv} 
= \normtwo{\by - \bX \bbeta^*} + \sum_{i=1}^p c_i \abs{\beta_i^*}.
\end{aligned}
\end{equation}
Combining \eqref{equation:rob_reg1_eq1} and \eqref{equation:rob_reg1_eq2} yields the desired equality for any fixed $\bbeta^*$
$$
\max_{\Delta \bX \in \sU} \normtwo{\by - (\bX + \Delta \bX) \bbeta^*} = \normtwo{\by - \bX \bbeta^*} + \sum_{i=1}^p c_i \abs{\beta_i^*}.
$$
Minimizing both sides over $\bbeta$ completes the proof.
\end{proof}

When the $\ell_1$-norm is replaced by an arbitrary norm $\norms{\cdot}$, the previous result generalizes as follows.
\begin{corollary}[Robust regression]\label{corollary:rob_reg2}
Let
\begin{equation}
\sU_a \triangleq \big\{ [\boldsymbol{\delta}_1, \ldots, \boldsymbol{\delta}_p] \mid \norms{\boldsymbol{\delta}_i} \leq c_i, \ i = 1,2, \ldots, p \big\}.
\end{equation}
Then the following robust regression problem 
$$
\min_{\bbeta \in \real^p} \max_{\Delta \bX \in \sU_a} \norm{\by - (\bX + \Delta \bX) \bbeta}_s,
$$ 
with uncertainty set $\sU_a$ is equivalent to the following arbitrary norm $\norms{\cdot}$-regularized regression problem:
\begin{equation}
\min_{\bbeta \in \real^p} \left\{ \norm{\by - \bX \bbeta}_s + \sum_{i=1}^{p} c_i \abs{\beta_i} \right\}.
\end{equation}
\end{corollary}

Thus, robust regression yields an objective that combines an $\ell_2$ (or more generally, $\ell_s$) loss with a weighted $\ell_1$ penalty. The following result shows that, in the case of the $\ell_2$ loss, this formulation is essentially equivalent to the standard regularized least squares model used in LASSO.
\begin{theoremHigh}[Equivalent forms \citep{xu2008robust}]\label{theorem:equiv_forms}
The optimization problems
\begin{equation}
\min_{\bbeta \in \real^p} \{\normtwo{\by - \bX \bbeta} + \lambda_1 \normone{\bbeta}\}
\qquad \text{and} \qquad
\min_{\bbeta \in \real^p} \{\normtwo{\by - \bX \bbeta}^{\textcolor{mylightbluetext}{2}} + \lambda_2 \normone{\bbeta}\}
\end{equation}
are equivalent in the sense that, for every choice of $\lambda_1>0$, there exists a corresponding $\lambda_2>0$ (and vice versa) such that both problems share the same set of minimizers.
\end{theoremHigh}
\begin{proof}
See Appendix of \citet{xu2008robust}.
\end{proof}

Now, suppose we set  $c_i = \lambda$ for all $i$ and assume that the columns $\bx_i$ of $\bX$ are orthonormal. Under these conditions, the robust regression problem \eqref{equation:rob_reg1} simplifies to
\begin{equation}
\min_{\bbeta \in \real^p} \left\{\normtwo{\by - \bX \bbeta} + \lambda \normone{\bbeta}\right\}.
\end{equation}
By Theorem~\ref{theorem:equiv_forms}, this problem is equivalent to the Lagrangian LASSO formulation \eqref{equation:ll_ana}. In other words, a robust linear regression problem---with featurewise $\ell_2$-bounded perturbations and orthonormal design---is equivalent to the LASSO problem.

Viewing the LASSO solution as the solution to a robust regression problem has two important implications:
\begin{enumerate}[(i)]
\item  \textit{Interpretability of regularization.} Robustness provides a concrete interpretation of the $\ell_1$ regularizer---it corresponds to protection against worst-case data perturbations, linking the mathematical penalty to a tangible physical property (noise resilience).
\item \textit{Analytical leverage.} More importantly, robustness is a powerful structural property that can be used to analyze the behavior of the solution. For instance, the sparsity of the LASSO estimator can be understood as a consequence of its robustness: under adversarial perturbations, the optimal strategy is to rely only on the most stable (i.e., nonzero) features. Moreover, by construction, a robust solution is optimal under the worst-case scenario, which offers strong guarantees in uncertain or noisy environments.
\end{enumerate}

\index{Prediction loss}
\index{Parameter estimation loss}
\subsection{Error Bounds of LASSO}
We conclude this section by presenting several theoretical results for the LASSO problem.
Given a LASSO estimate  $\widehat{\bbeta} \in \real^p$, its quality can be assessed in multiple ways, depending on the goals of the application.
In some settings, the primary concern is the predictive performance of $\widehat{\bbeta}$. 
In such cases, we evaluate it using the \textit{(in-sample) prediction loss}:
\begin{subequations}
\begin{equation}
\mathcalL_{\text{pred}}(\widehat{\bbeta}; \bbeta^*) = \frac{1}{n} \normtwo{\bX \widehat{\bbeta} - \bX \bbeta^*}^2,
\end{equation}
which measures the mean squared prediction error over the observed design matrix $\bX$.
In other applications---such as medical imaging, remote sensing, and compressed sensing---the main interest lies in accurately recovering the unknown parameter vector $\bbeta^*$ itself.
Here, it is more appropriate to consider the \textit{parameter estimation loss}, defined as the squared $\ell_2$-distance between the estimator and the true parameter:
\begin{equation}
\mathcalL_2(\widehat{\bbeta}; \bbeta^*) = \normtwobig{\widehat{\bbeta} - \bbeta^*}^2.
\end{equation}
Finally, when the goal is variable selection or support recovery, one may use a discrete loss that checks whether the estimated and true vectors share the same signed support:
\begin{equation}
\mathcalL_{\text{vs}}(\widehat{\bbeta}; \bbeta^*) =
\begin{cases}
0, & \text{if } \sign(\widehat{\beta}_i) = \sign(\beta_i^*) \text{ for all } i \in \{1,2, \ldots, p\}; \\
1, & \text{otherwise}.
\end{cases}
\end{equation}
\end{subequations}
This loss equals zero only when $\widehat{\bbeta}$ correctly identifies both the locations and signs of the nonzero entries of $\bbeta^*$.

For brevity, this subsection focuses on theoretical guarantees for the first two loss functions---prediction and parameter estimation. To establish these results, we require the following technical lemma.

\begin{lemma}\label{lemma:sparse_lem1}
Let $\balpha, \bbeta\in\real^p$  satisfy $\normone{\bbeta}\leq \normone{\balpha}$, where $\supp(\balpha) \triangleq\sS\subseteq \{1,2,\ldots,p\}$.
Denote $\be \triangleq \bbeta-\balpha$.  Then $\normone{\be_{\comple{\sS}}} \leq \normone{\be_{\sS}}$, i.e., $\be \in \sC[\sS; 1]$ (see definition in \eqref{equation:cscsgamma_def}), where $\comple{\sS}$ denotes the complement of $\sS$.
\end{lemma}
\begin{proof}[of Lemma~\ref{lemma:sparse_lem1}]
Let $\sT\triangleq\supp(\bbeta)$. The assumption  $\normone{\bbeta}\leq \normone{\balpha}$ indicates that
$$
\begin{aligned}
&\sum_{i\in\sT} \abs{\beta_i} \leq  \sum_{i\in\sS} \abs{\alpha_i}
\quad\implies\quad
\sum_{i\in\sT\setminus \sS} \abs{\beta_i}+ \sum_{i\in\sS \cap\sT} \abs{\beta_i}   \leq  \sum_{i\in\sS \setminus\sT} \abs{\alpha_i} + \sum_{i\in\sS \cap\sT} \abs{\alpha_i}\\
&\implies\quad
\sum_{i\in\sT\setminus \sS} \abs{\beta_i}    \leq  \sum_{i\in\sS \setminus\sT} \abs{\alpha_i} + \sum_{i\in\sS \cap\sT} \abs{\alpha_i - \beta_i},
\end{aligned}
$$
which obtains the desired result.
\end{proof}

With this intuition in place, we now present  a result that provides a bound on the parameter estimation error $\normtwobig{\widehat{\bbeta} - \bbeta^*}$ under the linear observation model $\by = \bX\bbeta^* + \bepsilon$, where $\bbeta^*$ is $k$-sparse and supported on a subset  $\sS\subseteq\{1,2,\ldots,p\}$ (i.e., $\abs{\sS}=k$).

\begin{theoremHigh}[Bounds of optimizer error for LASSO]\label{theorem:bound_lassl2error}
Assume the observations follow the linear model $\by = \bX\bbeta^* + \bepsilon$, where $\bX\in\real^{n\times p}$, and  $\bbeta^*$ is $k$-sparse with support $\sS\subseteq \{1, 2,\ldots, p\}$ (i.e., $\abs{\sS}=k$).
Suppose the design matrix $\bX$ satisfies the restricted eigenvalue (RE) property (Definition~\ref{definition:res_eig}) with constant $n\mu > 0$ over $\sC[\sS;3]$.~\footnote{Note that the first result only requires that the RE property holds over the subset $\sC[\sS;1]$.}
Then the following error bounds hold:
\begin{enumerate}[(i)]
\item \textit{Constrained LASSO.} Any estimate $\widehat{\bbeta}$ based on the constrained LASSO \eqref{opt:lc} with $\Sigma=\normone{\bbeta^*}$ satisfies the bound
\begin{equation}\label{equation:bound_lassl2error_11}
\normtwo{\widehat{\bbeta} - \bbeta^*} \leq \frac{4\sqrt{k}}{\mu n} \norminf{\bX^\top \bepsilon}.
\end{equation}

\item \textit{Lagrangian LASSO.} Given a regularization parameter $\lambda \geq {2} \norminf{\bX^\top \bepsilon} > 0$~\footnote{If we use the scaled Lagrangian LASSO loss in \eqref{opt:ll_scaled}, p.~\pageref{opt:ll_scaled}, then the requirement becomes $\lambda \geq \frac{2}{n} \norminf{\bX^\top \bepsilon} > 0$.}, any estimate $\widehat{\bbeta}$ from the regularized LASSO \eqref{opt:ll} satisfies the bound
\begin{equation}\label{equation:bound_lassl2error_22}
\normtwo{\widehat{\bbeta} - \bbeta^*} \leq \frac{3}{\mu} \sqrt{{k}} \frac{\lambda}{n}.
\end{equation}
\end{enumerate}
\end{theoremHigh}
\begin{proof}[of Theorem~\ref{theorem:bound_lassl2error}]
\textbf{(i) Constrained LASSO.} Since $\bbeta^*$ is feasible and $\widehat{\bbeta}$ is optimal, we have the inequality $\normtwobig{\by - \bX\widehat{\bbeta}}^2 \leq \normtwo{\by - \bX\bbeta^*}^2$. Defining the error vector $\widehatbe \triangleq \widehat{\bbeta} - \bbeta^*$, substituting in the relation $\by = \bX\bbeta^* + \bepsilon$ yields
\begin{equation}\label{equation:bound_lassl2error1}
\text{(Constrained LASSO)}: \qquad {\normtwo{\bX\widehatbe}^2} \leq {2}\bepsilon^\top \bX \widehatbe.
\end{equation}
Applying \holders inequality (Theorem~\ref{theorem:holder-inequality}) to the right-hand side yields the upper bound ${2} \abs{\bepsilon^\top \bX \widehatbe} \leq {2} \norminf{\bX^\top \bepsilon} \normone{\widehatbe}$. The inequality $\normonebig{\widehat{\bbeta}} \leq \Sigma = \normone{\bbeta^*}$ implies that $\widehatbe \in \sC[\sS; 1]$ (Lemma~\ref{lemma:sparse_lem1}):
\begin{equation}\label{equation:bound_lassl2error2}
\normonebig{\widehat{\bbeta}} \leq \normone{\bbeta^*} 
\;\;\implies\;\;
\normone{\widehatbe} = \normone{\widehatbe_{\sS}} + \normone{\widehatbe_{\comple{\sS}}} \leq 2 \normone{\widehatbe_{\sS}}
\stackrel{\dag}{\leq} 2 \sqrt{k} \normtwo{\widehatbe_{\sS}}\leq 2 \sqrt{k} \normtwo{\widehatbe},
\end{equation}
where the inequality ($\dag$) follows from the Cauchy--Schwarz inequality (Equation~\eqref{equation:vector_form_cauchyschwarz}):
\begin{equation}\label{equation:bound_lassl2error2_v2}
\normone{\widehatbe_{\sS}}=\abs{\widehatbe_{\sS}}^\top\bone \leq \normtwo{\widehatbe_{\sS}}\normtwo{\bone} = \sqrt{k}\normtwo{\widehatbe_{\sS}}
\leq  \sqrt{k}\normtwo{\widehatbe}.
\end{equation}
Note that \eqref{equation:bound_lassl2error2} holds under the condition $\normonebig{\widehat{\bbeta}} \leq \normone{\bbeta^*}$, while \eqref{equation:bound_lassl2error2_v2} is valid for any $\sS$ and $\widehatbe$.
On the other hand, applying the RE property to the left-hand side of the inequality \eqref{equation:bound_lassl2error1} yields $n \mu \normtwo{\widehatbe}^2 \leq  \normtwo{\bX \widehatbe}^2$. Combining these inequalities yields the desired result.

\paragraph{(ii) Lagrangian LASSO.} Define the function
\begin{equation}
G(\be) \triangleq \frac{1}{2} \normtwo{\by - \bX (\bbeta^* + \be)}^2 + \lambda \normone{\bbeta^* + \be}.
\end{equation}
Noting that $\widehatbe \triangleq \widehat{\bbeta} - \bbeta^*$ minimizes $G$ by construction, we have $G(\widehatbe) \leq G(\bzero)$.
Substituting in the relation $\by = \bX\bbeta^* + \bepsilon$ yields
\begin{equation}\label{equation:bound_lassl2error3}
\text{(Lagrangian LASSO)}:\qquad {\normtwo{\bX \widehatbe}^2} \leq {2} \bepsilon^\top \bX \widehatbe + 2\lambda \left\{ \normone{\bbeta^*} - \normone{\bbeta^* + \widehatbe} \right\}.
\end{equation}
Since $\bbeta_{\comple{\sS}}^* = \bzero$, we have $\normone{\bbeta^*} = \normone{\bbeta_{\sS}^*}$, and
$$
\normone{\bbeta^* + \widehatbe} = \normone{\bbeta_{\sS}^* + \widehatbe_{\sS}} + \normone{\widehatbe_{\comple{\sS}}} \geq \normone{\bbeta_{\sS}^*} - \normone{\widehatbe_{\sS}} + \normone{\widehatbe_{\comple{\sS}}}.
$$
Substituting these relations into the inequality \eqref{equation:bound_lassl2error3} and applying  \holders  yield
\begin{equation}\label{equation:bound_lassl2error5}
\begin{aligned}
{\normtwo{\bX \widehatbe}^2}
&\leq {2 }\bepsilon^\top \bX \widehatbe + 2\lambda \left\{ \normone{\widehatbe_{\sS}} - \normone{\widehatbe_{\comple{\sS}}} \right\} \\
&\leq {2 } \norminfbig{\bX^\top \bepsilon}\normone{\widehatbe} + 2\lambda \left\{ \normone{\widehatbe_{\sS}} - \normone{\widehatbe_{\comple{\sS}}} \right\}.
\end{aligned}
\end{equation}
Since ${2} \norminf{\bX^\top \bepsilon} \leq \lambda$ by assumption, using the fact that $\normone{\widehatbe_{\sS}} \leq \sqrt{k} \normtwo{\widehatbe}$ from \eqref{equation:bound_lassl2error2_v2}, we have:
\begin{equation}\label{equation:bound_lassl2error6}
0\leq {\normtwo{\bX \widehatbe}^2} \leq {\lambda} \left\{ \normone{\widehatbe_{\sS}} + \normone{\widehatbe_{\comple{\sS}}} \right\} + 2\lambda \left\{ \normone{\widehatbe_{\sS}} - \normone{\widehatbe_{\comple{\sS}}} \right\} \leq {3} \sqrt{k} \lambda \normtwo{\widehatbe}.
\end{equation}
The inequality \eqref{equation:bound_lassl2error6} also implies that $\normone{\widehatbe_{\comple{\sS}}} \leq 3 \normone{\widehatbe_{\sS}}$ such that $\widehatbe\in \sC[\sS; 3]$ (see definition in \eqref{equation:cscsgamma_def}).
Applying the $\mu$-RE condition to $\widehatbe$ ensures that $\mu \normtwo{\widehatbe}^2 \leq \frac{1}{n} \normtwo{\bX \widehatbe}^2$.
Combining this lower bound with the inequality \eqref{equation:bound_lassl2error6} yields
$
{\mu} \normtwo{\widehatbe}^2 \leq {3}  \sqrt{k} \frac{\lambda}{n} \normtwo{\widehatbe},
$
which implies the desired result.
\end{proof}

\begin{theoremHigh}[Bounds of prediction error for LASSO]\label{theorem:lassobound_predict_error}
Assume the observations follow the linear model $\by = \bX\bbeta^* + \bepsilon$, where $\bX\in\real^{n\times p}$.
Consider the Lagrangian LASSO \eqref{opt:ll} with a regularization parameter $\lambda \geq {2} \norminf{\bX^\top \bepsilon}$~\footnote{Again, if we use the scaled Lagrangian LASSO loss in \eqref{opt:ll_scaled}, then the requirement becomes $\lambda \geq \frac{2}{n} \norminf{\bX^\top \bepsilon} > 0$.}.
\begin{enumerate}[(i)]
\item If $\normone{\bbeta^*} < R_1$, then any optimal solution $\widehat{\bbeta}$ satisfies
\begin{equation}
{\normtwo{\bX (\widehat{\bbeta} - \bbeta^*)}^2} \leq 12 \, R_1 \, \lambda.
~\footnote{A refined bound is discussed in Problem~\ref{prob:refined_llbound}.}
\end{equation}
\item Suppose further that $\bbeta^*$ is $k$-sparse, supported on the subset $\sS\subseteq \{1, 2,\ldots, p\}$ (i.e., $\abs{\sS}=k$).
If the model matrix $\bX$ satisfies the restricted eigenvalue property with parameter $n\mu > 0$ over $\sC[\sS;3]$, then any optimal solution $\widehat{\bbeta}$ satisfies
\begin{equation}
{\normtwo{\bX (\widehat{\bbeta} - \bbeta^*)}^2} \leq 9\frac{ k \lambda^2}{n \mu}.
\end{equation}
\end{enumerate}
\end{theoremHigh}
\begin{proof}[of Theorem~\ref{theorem:lassobound_predict_error}]
The proof is a slight adaptation based on \citet{hastle2015statistical, lu2025practical}.
\paragraph{(i).} Once again, we define the error vector $\widehatbe \triangleq \widehat{\bbeta} - \bbeta^*$. Following the inequality \eqref{equation:bound_lassl2error3}, we have
\begin{align}
0 
&\leq{\normtwo{\bX \widehatbe}^2}  \leq {2} \norminfbig{\bX^\top \bepsilon} \normone{\widehatbe} + 2\lambda \big\{ \normone{\bbeta^*} - \normone{\bbeta^* + \widehatbe} \big\} \label{equation:lassobound_predict_error_ineq1}\\
&\overset{\dag}{\leq} 2\left\{ {\norminfbig{\bX^\top \bepsilon}} - \lambda \right\} \normone{\widehatbe} + 4 \lambda \normone{\bbeta^*}
\overset{\ddag}{\leq}  \lambda \Big\{ - \normone{\widehatbe} + 4 \normone{\bbeta^*} \Big\},
\end{align}
where the inequality $(\dag)$ follows from the triangular inequality
$-(\normone{\bbeta^*} - \normone{\widehatbe} ) \leq \normone{\bbeta^* - (-\widehatbe)}
$~\footnote{$\abs{\norm{\bbeta}-\norm{\balpha}} \leq \norm{\bbeta-\balpha}$ for any vector norm.}, and  the inequality $(\ddag)$ follows from the assumption that $\lambda \geq {2} \norminf{\bX^\top \bepsilon}$.
This indicates that $\normone{\widehatbe} \leq 4 \normone{\bbeta^*} \leq 4 R_1$. Combining these inequalities and applying the triangle inequality $\normone{\bbeta^*} - \normone{\bbeta^* + \widehatbe} \leq \normone{\bbeta^* -(\bbeta^* + \widehatbe)}=\normone{\widehatbe}$ in \eqref{equation:lassobound_predict_error_ineq1} yield that
\begin{equation}\label{equation:lassobound_predict_error_part1fin}
{\normtwo{\bX \widehatbe}^2} \leq 2\left\{ {\norminfbig{\bX^\top \bepsilon}} + \lambda \right\} \normone{\widehatbe} \leq 12 \lambda R_1,
\end{equation}
which establishes the desired result.

\paragraph{(ii).} Following \eqref{equation:bound_lassl2error6}, we have
$
{\normtwo{\bX \widehatbe}^2}  \leq 3  \sqrt{k}\lambda \normtwo{\widehatbe}.
$~\footnote{
Note that \citet{hastle2015statistical} use \eqref{equation:lassobound_predict_error_part1fin} and the fact $\normone{\widehatbe} = \normone{\widehatbe_{\sS}} + \normone{\widehatbe_{\comple{\sS}}} \leq 4 \normone{\widehatbe_{\sS}}
\leq 4 \sqrt{k} \normtwo{\widehatbe_{\sS}}\leq 4 \sqrt{k} \normtwo{\widehatbe}$ to obtain
${\normtwo{\bX \widehatbe}^2} \leq 12\sqrt{k}\normtwo{\widehatbe}$, which results in a looser bound.
}
By the proof of Theorem~\ref{theorem:bound_lassl2error}, the error vector $\widehatbe$ belongs to the cone $\sC[\sS;3]$, so that the $\mu$-RE condition guarantees that $\normtwo{\widehatbe}^2 \leq \frac{1}{n \mu} \normtwo{\bX \widehatbe}^2$. Combining these inequalities yields the desired result.
\end{proof}

\begin{remark}[Compatibility condition, $\ell_1$ estimation bound]
Another related condition is the \textit{compatibility condition} on the design matrix $\bX$. 
It is defined with respect to the true support set $\sS=\supp(\bbeta^*)$ and requires that, for some compatibility constant $\mu > 0$,
\begin{equation}\label{eq:compatibility}
\frac{1}{n}\normtwo{\bX \bu}^2 
\geq \frac{\mu}{k} \normone{\bu_\sS}^2, \quad \text{for all } \bu \in \sC[\sS; 3],
\text{ i.e., } \normone{\bu_{\comple{\sS}}} \leq 3\normone{\bu_\sS},
\end{equation}
where $k\triangleq \normzero{\bbeta^*}$.
Intuitively, this condition ensures that the truly active predictors are not too highly correlated.
Note that inequality~\eqref{equation:bound_lassl2error6} also implies
\begin{equation}\label{equation:bound_lassl2error6_impli2}
\normtwo{\bX\widehatbbeta - \bX\bbeta^*}^2 \leq 3\lambda \normone{\widehatbbeta_\sS - \bbeta^*_\sS}. 
\end{equation}
Since $\widehatbe=\widehatbbeta-\bbeta^*\in\sC[\sS;3]$ by \eqref{equation:bound_lassl2error6},
we can combine this with the compatibility condition to bound the parameter estimation error. Specifically, applying the lower bound from~\eqref{eq:compatibility} to the left-hand side of~\eqref{equation:bound_lassl2error6_impli2} yields
$$
\frac{n\mu}{k} \normone{\widehatbbeta_\sS - \bbeta^*_\sS}^2 
\leq \normtwo{\bX\widehatbbeta-\bX\bbeta^*}^2
\leq 3\lambda \normone{\widehatbbeta_\sS - \bbeta^*_\sS}.
$$
Dividing both sides by $\normonebig{\widehatbbeta_\sS - \bbeta^*_\sS}$, and using the fact that $\normonebig{\widehatbbeta_{\comple{\sS}}} \leq 3\normonebig{\widehatbbeta_\sS - \bbeta^*_\sS}$, which implies $\normonebig{\widehatbbeta - \bbeta^*} \leq 4\normonebig{\widehatbbeta_\sS - \bbeta^*_\sS}$, we get
\begin{equation}
\normone{\widehatbbeta - \bbeta^*} \leq \frac{k\lambda}{ n\mu}.
\end{equation}
This provides an $\ell_1$-bound for the parameter estimation loss. 
Compared to the $\ell_1$-bound in \eqref{equation:bound_lassl2error_22}, the dependence on the sparsity level is now linear in  $k$ rather than $\sqrt{k}$.
\end{remark}

\subsection{Gaussian Noise and Oracle Inequality}
We consider the linear model $\by=\bX\bbeta^*+\bepsilon$, where the noise vector $\bepsilon \in \real^n$ has i.i.d. Gaussian entries with mean zero and variance  $\sigma^2$. 
The design matrix $\bX \in \real^{n \times p}$ is assumed fixed, and each column satisfies $\max_{i=1,2,\ldots,p} \normtwo{\bx_i} \leq \sqrt{n}$.
(Note that this normalization can always be achieved by rescaling.)
Under these assumptions, the random vector $\bX^\top\bepsilon$ has independent Gaussian entries with mean zero and variance at most  $\max_{i=1,2,\ldots,p} \normtwo{\bx_i}^2 \sigma^2 \leq n\sigma^2$. 
Using standard Gaussian tail bounds (see Lemma~\ref{lemma:bd_cen_gaus}), we have
$$
\Pr\left(\abs{(\bX^\top \bepsilon)_i} \geq t \sqrt{n}\sigma\right) 
\leq \exp\left(-\frac{t^2}{2}\right),
\quad \text{for $i = 1, 2, \ldots, p$}.
$$
Applying the union bound over all $p$ coordinates (Theorem~\ref{theorem:union_bound_proof}) and setting the parameter $t=\sqrt{2(\ln(2p) + \xi)}$, we obtain
\begin{equation}
\Pr\left(\norminfbig{\bX^\top \bepsilon} \geq \sigma\sqrt{2n(\ln(2p) + \xi)} \right) 
\leq \frac{1}{2}\exp(-\xi), 
\quad \forall\, \xi>0.
\end{equation}
Thus, choosing $\lambda = 2\sigma \sqrt{2n(\ln(2p) + \xi)}= \mathcalO\big(\sqrt{n{\ln(p)}}\big)$ is valid for Theorem~\ref{theorem:lassobound_predict_error} with probability at least $1 - \exp(-\xi)/2$. 
Consequently, the prediction loss bounds take the form:
\begin{subequations}\label{equation:ll_bd_fastslow}
\begin{align}
\frac{\normtwo{\bX (\widehat{\bbeta} - \bbeta^*)}^2}{n} &\leq \mathcalO\big(R_1 \, \sqrt{\frac{\ln(p)}{n}}\big); \label{equation:lassobound_predict_error1} \\
\frac{\normtwo{\bX (\widehat{\bbeta} - \bbeta^*)}^2}{n} &\leq \mathcalO\big(\frac{ \abs{\sS} \, \ln(p)}{n}\big). \label{equation:lassobound_predict_error2}
\end{align}
\end{subequations}
The first bound \eqref{equation:lassobound_predict_error1}, which depends only on the $\ell_1$-norm radius $R_1$, is known as the ``slow rate," since the mean squared prediction error decays at the rate $1/\sqrt{n}$.
In contrast, the second bound \eqref{equation:lassobound_predict_error2} is referred to as the ``fast rate," exhibiting a faster $1/n$ decay.

Note that the fast rate relies on stronger assumptions: specifically, that $\bbeta^*$ is exactly $k$-sparse (i.e., hard sparsity), and---more critically---that the design matrix $\bX$ satisfies the RE condition.
In principle, good prediction performance should not require such structural conditions on $\bX$, suggesting that the RE assumption may be an artifact of the current proof technique rather than a fundamental requirement.

When analyzing the constrained form of the LASSO (problem \eqref{opt:lc}),
where $\Sigma \geq 0$ serves as the tuning parameter. 
Setting  $\Sigma = \normone{\bbeta^*}$ ensures that the true coefficient vector $\bbeta^*$ is feasible for \eqref{opt:lc}; under this choice, the same bounds in~\eqref{equation:ll_bd_fastslow} hold with probability at least  $1 - \exp(-\xi)/2$.

Even without assuming a linear model  $\by=\bX\bbeta^*+\bepsilon$, we can still derive a meaningful bound on the prediction loss that quantifies its excess risk relative to the best linear predictor.
Suppose instead that
$$
\by = f_0(\bX) + \bepsilon,
$$
for some unknown function $f_0 : \real^p \to \real$, where  $f_0(\bX) = [f_0(\bx^1), \ldots, f_0(\bx^n)]^\top \in \real^n$ denotes the row-wise application of $f_0$ to the rows of $\bX$. 
Then, for any solution $\widehatbbeta$ of the constrained LASSO problem~\eqref{opt:lc} and any  $\widebarbbeta$ satisfying  $\normone{\widebarbbeta} \leq \Sigma$, we have
$$
\normtwo{\by-\bX\widehatbbeta}^2
\leq \normtwo{\by-\bX\widebarbbeta}^2.
$$
Rearranging terms gives
$$
\begin{aligned}
\frac{1}{2}\normtwo{\bX\widehatbbeta - \bX\widebarbbeta}^2 
&\leq \innerproduct{\by - \bX\widebarbbeta, \bX\widehatbbeta - \bX\widebarbbeta} 
= \innerproduct{f_0(\bX) - \bX\widebarbbeta, \bX\widehatbbeta - \bX\widebarbbeta} + \innerproduct{\bepsilon, \bX\widehatbbeta - \bX\widebarbbeta}.
\end{aligned}
$$
Applying the polarization identity $\normtwo{\ba}^2 + \normtwo{\bb}^2 - \normtwo{\ba - \bb}^2 = 2\innerproduct{\ba, \bb }$ to the first inner product on the right-hand side yields
$$
\normtwo{\bX\widehatbbeta - \bX\widebarbbeta}^2 \leq \normtwo{\bX\widebarbbeta - f_0(\bX)}^2 + \normtwo{\bX\widehatbbeta - \bX\widebarbbeta}^2 - \normtwo{\bX\widebarbbeta - f_0(\bX)}^2 + 2\innerproduct{\bX^\top\bepsilon, \widehatbbeta - \widebarbbeta}.
$$
After rearranging, we obtain
\begin{align*}
	\normtwo{\bX\widehatbbeta - f_0(\bX)}^2 
	&\leq \normtwo{\bX\widebarbbeta - f_0(\bX)}^2 + 2\innerproduct{ \bX^\top\bepsilon, \widehatbbeta - \widebarbbeta} 
	\leq \normtwo{\bX\widebarbbeta - f_0(\bX)}^2 + 4\Sigma \norminf{\bX^\top\bepsilon}  \\
	&\leq \normtwo{\bX\widebarbbeta - f_0(\bX)}^2 + 4\Sigma\sigma\sqrt{{2n(\ln(2p) + \xi)}},
\end{align*}
where the last inequality holds with probability at least $1 - \exp(-\xi)/2$. 
Since this holds simultaneously for all $\widebarbbeta$ with  $\normone{\widebarbbeta} \leq \Sigma$ and  probability at least $1 - \exp(-\xi)/2$, we may take the infimum over such $\widebarbbeta$ on the right-hand side, yielding the \textit{oracle inequality}:
\begin{equation}\label{eq:oracle-inequality}
	\frac{1}{n}\normtwo{\bX\widehatbbeta - f_0(\bX)}^2 \leq \min_{\normone{\widebarbbeta} \leq \Sigma} \left( \frac{1}{n}\normtwo{\bX\widebarbbeta - f_0(\bX)}^2 \right) + 4\Sigma \sigma \sqrt{\frac{2(\ln(2p) + \xi)}{n}}, 
\end{equation}
with probability at least $1 - \exp(-\xi)/2$. 
This inequality states that, in terms of prediction loss, the constrained LASSO performs nearly as well as the best  $\ell_1$-constrained linear predictor of $f_0$.

Let $\widebarbbeta^{*}$ denote a minimizer of  $\frac{1}{n}\normtwo{\bX\widebarbbeta - f_0(\bX)}^2$ over $\normone{\widebarbbeta} \leq \Sigma$.~\footnote{The infimum is achieved by Weierstrass theorem (Theorem~\ref{theorem:weierstrass_them}) since we are minimizing a continuous function over a compact set.} 
Then, by the triangle inequality and~\eqref{eq:oracle-inequality}, we further obtain
\begin{equation}\label{eq:oracle-sharp}
	\frac{1}{n}\normtwo{\bX\widehatbbeta - \bX\widebarbbeta^{*}}^2 \leq 4\Sigma \sigma \sqrt{\frac{2(\ln(2p) + \xi)}{n}}, 
\end{equation}
with probability at least $1 - \exp(-\xi)/2$. 
This result is stronger than the oracle inequality: it shows that the LASSO's predictions are close to those of the best $\ell_1$-constrained linear predictor---even when that predictor itself is far from the true regression function $f_0$.

\begin{problemset}


\item \citep{candes2005decoding} Suppose $\bX$ is such that $\bX \bZ = \bzero$ and let $k \geq 1$ be a number obeying the hypothesis of Theorem~\ref{theorem:recov_rip_l1_t1}. 
Set $\by = \bZ\bbeta + \bz$, where $\bz$ is a real vector supported on a set of size at most $k$. Show that $\bbeta$ is the unique minimizer to
$$
(\text{P}_1') \qquad \min_{\bb \in \real^p} \normone{\by - \bZ \bb}.
$$

\item \label{prob:ells_quasnorm_pow} \textbf{Triangle inequality of $\ell_s$-quasinorm.}
For $0 < s < 1$, prove that the $s$-th power of the $\ell_s$-quasinorm satisfies the triangle inequality
$$
\norms{\balpha + \bbeta}^s \leq \norms{\balpha}^s + \norms{\bbeta}^s, \quad \balpha, \bbeta \in \real^p.
$$
Also show that
$$
\norms{\balpha_1 + \ldots + \balpha_k} \leq k^{\max(1/s - 1, 0)} (\norms{\balpha_1} + \ldots + \norms{\balpha_k}), \quad \balpha_1, \ldots, \balpha_k \in \real^p.
$$

\item \textbf{Recovery under $\ell_1$-coherence \citep{foucart2013invitation}.} Let $\nu < 1/2$, and suppose the matrix $\bX \in \real^{n\times p}$ satisfies
$$
\mu_1(\bX, k) \leq \nu.
$$
Prove that for any $\bbeta^* \in \real^p$ and observations from the linear model $\by = \bX\bbeta^* + \bepsilon$ with $\normtwo{\bepsilon} \leq \epsilon$, any minimizer $\widehatbbeta$ of $\normone{\bbeta}$ subject to $\normtwo{\bX\bbeta - \by} \leq \epsilon$ satisfies the $\ell_1$-error bound
$$
\normone{\bbeta^* - \widehatbbeta}\leq D_1 \sigma_k(\bbeta^*)_1 + D_2k\epsilon,
$$
for some positive constants $D_1$ and $D_2$ that depend only on $\nu$.

\item Use the argument in Theorem~\ref{theorem:rob_reg1} to prove Corollary~\ref{corollary:rob_reg2}.

\item Does the $\ell_1$-analysis models in \eqref{equation:q1_analysis_opcd} favor exact recovery? Discuss the conditions under which such recovery is possible.

\item \label{prob:ell1_re} \textbf{$\ell_1$-minimization with noise measurement under RE property \citep{hastle2015statistical}.}
Consider the $\ell_1$-minimization problem with noise measurement, \eqref{opt:p1_epsilon}:
\begin{equation}
\widehatbbeta = \arg\min_{\bbeta \in \real^p} \normone{\bbeta} \quad \text{s.t.} \quad \normtwo{\bX\bbeta-\by} \leq \epsilon, 
\end{equation}
where the constant $\epsilon > 0$ is a user-specified parameter.
\begin{enumerate}[(a)]
\item Suppose that $\epsilon$ is chosen such that $\bbeta^*$ is feasible (i.e., $\normtwo{\bX\bbeta^*-\by} \leq \epsilon$). Show that the error vector $\widehatbe = \widehatbbeta - \bbeta^*$ must satisfy the  constraint
$$
\normone{\widehatbe_{\comple{\sS}}} \leq \normone{\widehatbe_\sS}.
$$

\item Assuming the linear observation model $\by = \bX\bbeta^* + \bepsilon$, show that $\widehatbe$ satisfies the  inequality
$$
{\normtwo{\bX\widehatbe}^2} \leq 2 \left\{ {\norminf{\bX^\top \bepsilon}} \normone{\widehatbe} + \left({\epsilon^2} - {\normtwo{\bepsilon}^2} \right) \right\}.
$$

\item Assuming that   $\bX$ satisfies a $\mu$-RE property, use part (b) to derive a bound on the $\ell_2$-error $\normtwobig{\widehatbbeta - \bbeta^*}$.
\end{enumerate}

\item  \label{prob:lasso_uniq1}\textbf{LASSO uniqueness \citep{tibshirani2013lasso}.} 
For any $ \by, \bX $, and $ \lambda \geq 0 $, show that the LASSO problem~\eqref{opt:ll} has the following properties:
\begin{enumerate}[(a)]
\item Either the solution is unique, or there are uncountably infinitely many solutions.
\item All LASSO solutions yield the same fitted values, i.e., $ \bX\widehatbbeta $ is identical for every solution $\widehatbbeta$.
\item If $ \lambda > 0 $, then every LASSO solution $ \widehatbbeta $ has the same $ \ell_1 $-norm, $ \normonebig{\widehatbbeta} $.
\end{enumerate}
\textit{Hint: Consider the equivalence between the LASSO and $\ell_1$-minimization. Alternatively, argue by contradiction.}

\item \textbf{LASSO KKT.}
Consider the Lagrangian LASSO problem \eqref{opt:ll} with $\lambda>0$.
Show that  
\begin{subequations}
\begin{enumerate}[(a)]
\item $\widehatbbeta$ is a solution of \eqref{opt:ll} if and only if $\widehatbbeta$ satisfies~\eqref{equqtion:lassprob_kkt1} and~\eqref{equqtion:lassprob_kkt2} for some $\bnu\in\real^p$:
\begin{align}
\bX^\top(\by - \bX\widehatbbeta) &= \lambda \bnu, \label{equqtion:lassprob_kkt1} \\
\nu_i &\in 
\begin{cases}
\{\sgn(\widehatbeta_i)\} & \text{if } \widehatbeta_i \ne 0; \\
[-1,1] & \text{if } \widehatbeta_i = 0,
\end{cases}
\quad \text{for } i = 1,2,\dots p. \label{equqtion:lassprob_kkt2}
\end{align}
Here, $\bnu \in \real^p$ is  a subgradient of the function $f(\bbeta) = \normone{\bbeta}$ evaluated at $\bbeta = \widehatbbeta$ (Definition~\ref{definition:subgrad}). 
\item The subgradient $\bnu$ satisfying the above conditions is unique.
\item Define the equicorrelation set $\sS$ and equicorrelation signs $\bs$ by
\begin{align}
\sS &\triangleq \left\{ i \in \{1,2,\dots p\} \mid \abs{\bx_i^\top(\by - \bX\widehatbbeta)} = \lambda \right\}; \label{eq:equicorr_set}\\
\bs &\triangleq \sgn\big(\bX_{\sS}^\top(\by - \bX\widehatbbeta)\big).  \label{eq:equicorr_sign}
\end{align}
Prove that both $\sS$ and $\bs$ are uniquely determined (i.e., independent of the choice of LASSO solution).

\item Show  that $\widehatbbeta_{\comple{\sS}} = \bzero$ for any LASSO solution $\widehatbbeta$.
Consequently, the restriction of~\eqref{equqtion:lassprob_kkt1} to the indices in $\sS$ becomes
\begin{equation}
\bX_{\sS}^\top(\by - \bX_{\sS}\widehatbbeta_{\sS}) = \lambda \bs. \label{eq:ec_block}
\end{equation}
This means that $\lambda \bs$ lies in the row space of $\bX_{\sS}$, so $\lambda \bs = \bX_{\sS}^\top(\bX_{\sS}^\top)^+ \lambda \bs$. Using this fact, and rearranging~\eqref{eq:ec_block}, we get
$$
\bX_{\sS}^\top \bX_{\sS} \widehatbbeta_{\sS} = \bX_{\sS}^\top(\by - (\bX_{\sS}^\top)^+ \lambda \bs).
$$
\end{enumerate}
\end{subequations}

\textit{Hint: Use Problem~\ref{prob:lasso_uniq1}.}

\item \textbf{LASSO uniqueness.} For any $ \by, \bX $, and $ \lambda > 0 $, 
suppose that the submatrix $\bX_{\sS}$ corresponding to the equicorrelation set $\sS$ satisfies $ \nspace(\bX_{\sS}) = \{\bzero\} $, or equivalently  $ \rank(\bX_{\sS}) = \abs{\sS} $. 
Show that the Lagrangian LASSO solution is unique and is given by
$$
\widehatbbeta_{\comple{\sS}} = \bzero \quad \text{and} \quad \widehatbbeta_{\sS} = (\bX_{\sS}^\top \bX_{\sS})^{-1} (\bX_{\sS}^\top \by - \lambda \bs), 
$$
where $ \sS $ and $ \bs $ are the equicorrelation set and signs defined in~\eqref{eq:equicorr_set} and~\eqref{eq:equicorr_sign}. Note

\item \textbf{LASSO uniqueness.}
Assume that the matrix $ \bX $ has columns $ \bx_1,\bx_2, \ldots, \bx_p \in \real^n $ that are in general position. This means that for any integer $ k < \min\{n,p\} $, indices $ i_1, i_2,\ldots, i_{k+1} \in \{1,2,\ldots,p\} $, and signs $ s_1, s_2, \ldots, s_{k+1} \in \{-1,+1\} $, the affine span of the vectors $ s_1 \bx_{i_1},s_2 \bx_{i_2}, \ldots, s_{k+1} \bx_{i_{k+1}} $ 
does not contain any vector from the set $ \{\pm \bx_i \mid  i \ne i_1, \ldots, i_{k+1}\} $. 
Show that,  for any $ \by \in \real^n $ and $ \lambda > 0 $, the LASSO problem~\eqref{opt:ll} admits a unique solution.

\item  \label{prob:refined_llbound112} Consider the Lagrangian LASSO \eqref{opt:ll}, and let $\widehatbbeta$ denote any solution to it.
Show that for any coefficient vector $\bbeta \in \real^p$, the following inequality holds: 
$$
\normtwo{\bX\widehatbbeta - \bX\bbeta}^2 \leq 
2\innerproduct{\by - \bX\bbeta, \bX\widehatbbeta - \bX\bbeta} + 
2\lambda(\normone{\bbeta} - \normonebig{\widehatbbeta}),
$$
\textit{Hint: For any coefficient vector $\bbeta \in \real^p$, use the optimality of $\widehatbbeta$:
$
\frac{1}{2}\normtwo{\by - \bX\widehatbbeta}^2 + \lambda\normonebig{\widehatbbeta} \leq \frac{1}{2}\normtwo{\by - \bX\bbeta}^2 + \lambda\normone{\bbeta}.
$
}

\item \label{prob:refined_llbound} \textbf{Refined bound of prediction error for LASSO.}
Assume the observations follow the linear model $\by = \bX\bbeta^* + \bepsilon$, where $\bX\in\real^{n\times p}$.
Consider the Lagrangian LASSO \eqref{opt:ll} with a regularization parameter $\lambda \geq  \norminf{\bX^\top \bepsilon}$.
If $\normone{\bbeta^*} < R_1$, show that any optimal solution $\widehat{\bbeta}$ satisfies
$$
{\normtwo{\bX (\widehat{\bbeta} - \bbeta^*)}^2} \leq 4 \, R_1 \, \lambda.
$$
This result provides a refined bound compared to Theorem~\ref{theorem:lassobound_predict_error}.
Notably, it requires the less restrictive condition $\lambda \geq  \norminf{\bX^\top \bepsilon}$, rather than $\lambda \geq  2\norminf{\bX^\top \bepsilon}$.
\textit{Hint: Use Problem~\ref{prob:refined_llbound112} and   \holders inequality to bound the prediction loss.}

\end{problemset}

\newpage 
\chapter{Design of Sensing Matrices Using Gaussian  Ensembles}\label{chapter:spar_gauss}
\begingroup
\hypersetup{
linkcolor=structurecolor,
linktoc=page,  
}
\minitoc \newpage
\endgroup

\lettrine{\color{caligraphcolor}S}
Since properties such as the restricted isometry property (RIP), the restricted eigenvalue (RE) condition, nullspace property (NSP), and the restricted strong convexity (RSC) are crucial for guaranteed sparse recovery, it is important to study problem settings in which real-world or constructed data matrices actually satisfy these conditions.
Significant research has focused on explicitly constructing matrices that provably satisfy the RIP---for example, see \citet{baraniuk2008simple, bourgain2011explicit, jain2017non}.

To build intuition, we begin by introducing the necessary probabilistic tools related to Gaussian distributions.
We first analyze sparse recovery using random Gaussian matrices and then extend the results to sub-Gaussian random matrices in Chapter~\ref{chapter:ensur_rips}.

\section{Covering Numbers}\label{section:math_subgaus}

Before discussing Gaussian random variables and their application in ensuring the restricted isometry property, we introduce some essential mathematical tools concerning covering numbers.

A covering number refers to the minimum number of sets---typically balls in a metric space or substructures in a combinatorial setting---required to completely cover a given space or set.
Conversely, a packing number measures how many points can be placed within a set $\sS$ such that any two are at least $\varepsilon$ apart. This captures the ``spread" or ``separatedness" of the set.

Formally, we have the following definition.
\begin{definition}[Covering number,  packing number\index{Covering number}\index{Packing number}]\label{definitio:cov_pack_num}
Let $\sS$ be a subset of a metric space $(\sX,d)$, where $d(\cdot, \cdot)$ denotes the metric.
For $\varepsilon > 0$, the \textit{$\varepsilon$-covering number (or simply covering number)} of $\sS$, denoted  $\mathcalC(\sS,d,\varepsilon)$, is defined as the smallest integer $\mathcalC$ such that $\sS$ can be covered by $\mathcalC$ closed balls of radius $\varepsilon$, $\sB_d[\bx_i,\varepsilon] \triangleq \{\by \in \sX \mid  d(\bx_i, \by) \leq \varepsilon\}$, $\bx_i \in \sS$, $i = 1,2,\ldots,\mathcalC$. 
That is,
$$
\mathcalC\equiv \mathcalC(\sS,d,\varepsilon) = 
\min\left\{ k \in \naturalset \mid \exists\, \bx_1,\bx_2,\ldots,\bx_k\in \sX \text{ such that }  \sS \subset \bigcup_{i=1}^k \sB[\bx_i,\varepsilon] \right\}.
$$
Any such set  $\sC_{\varepsilon}\triangleq \{\bx_1,\bx_2,\ldots,\bx_\mathcalC\}$ is then called an \textit{$\varepsilon$-covering} of $\sS$.
Note that the centers $\bx_i$ do \textbf{not} necessarily have to lie in $\sS$.

For $\varepsilon > 0$, the \textit{$\varepsilon$-packing number (or simply packing number)} of $\sS$, denoted  $\mathcalP(\sS,d,\varepsilon)$,
is the largest number $\mathcalP$ of disjoint closed balls of radius $\varepsilon$ centered at points in $\sS$. Equivalently, it is the maximum number $\mathcalP$ of points 
$\bx_1, \bx_2, \ldots, \bx_\mathcalP$ such that $d(\bx_i,\bx_i) > \varepsilon$ for all $i\neq j$ (i.e., it is the maximum number of disjoint closed balls of radius $\varepsilon/2$ centered at points in $\sS$). 
Such a set $\sP_{\varepsilon}\triangleq \{\bx_1, \bx_2, \ldots, \bx_\mathcalP\}$ is called an \textit{$\varepsilon$-packing} or \textit{$\varepsilon$-separated} set.
\end{definition}

As $\varepsilon \rightarrow 0$, finer resolution is required, so the covering number $\mathcalC(\sS,d,\varepsilon)$ generally increases.

When  $\sX = \real^n$ is a vector space and  the metric is induced by a norm, i.e., $d(\ba, \bb) = \norm{\ba-\bb}$, we often write the following for brevity:
$$
\mathcalC(\sS,\norm{\cdot},\varepsilon)
\qquad \text{and}\qquad \mathcalP(\sS,\norm{\cdot},\varepsilon).
$$

\begin{subequations}
We now list several basic properties of covering numbers; packing numbers satisfy analogous properties. 
For arbitrary  sets $\sM,\sS \subset \sX$, the ``triangle inequality" holds:
\begin{equation}
\mathcalC(\sM \cup \sS, d, \varepsilon) \leq \mathcalC(\sM, d, \varepsilon) + \mathcalC(\sS, d, \varepsilon). 
\end{equation}
For any $\eta > 0$,
\begin{equation}
\mathcalC(\sS, \eta d, \varepsilon) = \mathcalC(\sS, d, \varepsilon/\eta). 
\end{equation}
If $\sX = \real^n$ and $d(\cdot,\cdot)$ is induced by a norm $\norm{\cdot}$, then 
\begin{equation}
\mathcalC(\eta \sS, d, \varepsilon) = \mathcalC(\sS, d, \eta^{-1}\varepsilon). 
\end{equation}
Moreover, if $d'(\cdot,\cdot)$ is another metric on $\sX$ satisfying  $d'(\ba, \bb) \leq d(\ba,\bb)$ for all $\ba,\bb \in \sX$, then
\begin{equation}
\mathcalC(\sS, d', \varepsilon) \leq \mathcalC(\sS, d, \varepsilon). 
\end{equation}
\end{subequations}

There is a fundamental relationship between covering and packing numbers:
\begin{lemma}[Key relationship between covering and packing numbers]\label{lemma:packing_covering}
Let $\sS$ be a subset of a metric space $(\sX,d)$, and let $\varepsilon > 0$. Then,
$$
\mathcalP(\sS, d, 2\varepsilon) \leq \mathcalC(\sS, d, \varepsilon) \leq \mathcalP(\sS, d, \varepsilon).
$$
\end{lemma}
\begin{proof}[of Lemma~\ref{lemma:packing_covering}]
Let $\{\bx_1,\bx_2,\ldots,\bx_\mathcalP\}$ be a $2\varepsilon$-packing set such that each $\bx_i\in \sS$ and $d(\bx_i, \bx_j)> 2\varepsilon$ for $i\neq j$,  
and let
$\{\by_1,\by_2, \ldots,\by_\mathcalC\}$ be an $\varepsilon$-covering set such that $\sS \subset \bigcup_{i=1}^\mathcalC \sB[\by_i,\varepsilon]$. 
Each $\bx_j$ must lie in at least one (closed) ball $\sB_d[\by_i,\varepsilon]$.
That is, we can assign to each point $\by_i$ a point $\bx_j$ with $d(\by_i,\bx_j) \leq \varepsilon$. 
Crucially, each ball $\sB_d[\by_i,\varepsilon]$ can contain \textbf{at most one} point from the $2\varepsilon$-packing set.
Indeed, if two distinct points $\bx_i, \bx_k$ ($i\neq k$) both belonged to $\sB_d[\by_i,\varepsilon]$, satisfying
$$
d(\bx_i, \by_j) \leq \varepsilon
\qquad \text{and}\qquad 
d(\bx_k, \by_j) \leq \varepsilon,
$$
then by the triangle inequality,
$$
d(\bx_i, \bx_k) 
\leq d(\bx_i, \by_j) + d(\bx_k, \by_j) 
\leq 2\varepsilon,
$$
contradicting the definition of a $2\varepsilon$-packing set. This proves the first inequality.

For the second inequality, let $\{\bx_1,\bx_2,\ldots,\bx_\mathcalP\}$ be an $\varepsilon$-packing of $\sS$. 
Since it's maximal, for any point $\bx\in\sS$, there exists some $\bx_i$ such that $d(\bx_i, \bx)\leq \varepsilon$. Otherwise,  we could add $\bx$ to the packing, contradicting maximality. Thus, $\sS\subset \bigcup_{i=1}^\mathcalP \sB[\by_i,\varepsilon]$. This means $\mathcalP$ is at least the covering number.
\end{proof}

The following result bounds the packing (and hence covering) number of a norm-induced unit sphere (norm-sphere) in finite dimensions.
\begin{theoremHigh}[Covering and packing numbers of a norm-sphere]\label{theorem:covnum_sphere}
Let $\norm{\cdot}$ be some norm on $\real^n$, and let $\sA$ be a subset of the unit ball $\sB \triangleq \sB_{\norm{\cdot}}[\bzero, 1] = \{\bx \in \real^n\mid \norm{\bx} \leq 1\}$. Then the packing and covering numbers of the $n$-dimensional sphere $\sB$ at scale $\varepsilon<1$ satisfy
\begin{equation}
\mathcalC(\sA, \norm{\cdot}, \varepsilon) 
\leq \mathcalP(\sA, \norm{\cdot}, \varepsilon) \leq \left(1 + \frac{2}{\varepsilon}\right)^n
\leq \left( \frac{3}{\varepsilon} \right)^n.
\end{equation}
\end{theoremHigh}
\begin{proof}[of Theorem~\ref{theorem:covnum_sphere}]
The first inequality follows from Lemma~\ref{lemma:packing_covering}. 
Let $\{\bx_1,\bx_2,\ldots,\bx_\mathcalP\} \subset \sA$ be a maximal $\varepsilon$-packing of $\sA$. 
Then the balls $\sB_{\norm{\cdot}}[\bx_\ell,\varepsilon/2]$ are pairwise disjoint, and they are all contained in the scaled  ball $(1 + \varepsilon/2)\sB$ with radius $(1 + \varepsilon/2)$. 
Comparing volumes yields
$$
\vol\left(\bigcup_{\ell=1}^\mathcalP \sB_{\norm{\cdot}}[\bx_\ell,\frac{\varepsilon}{2}]\right) 
= \mathcalP\, \vol\left(\frac{\varepsilon}{2}\sB\right) 
\leq \vol\left(\big(1 + \frac{\varepsilon}{2}\big)\sB\right).
$$
Since  the volume satisfies 
$$
\vol(\varepsilon\sB) = \varepsilon^n\, \vol(\sB)
\quad\text{ on $\real^n$},
$$ 
we obtain
$$
\mathcalP\big(\frac{\varepsilon}{2}\big)^n\, \vol(\sB) 
\leq \big(1 + \frac{\varepsilon}{2}\big)^n\, \vol(\sB)
\qquad\implies\qquad
\mathcalP \leq \big(1 + \frac{2}{\varepsilon}\big)^n.
$$
This completes the proof.
\end{proof}

Specifically, the following corollary provides a bound for the covering and packing numbers of the \textit{$n$-dimensional sphere}
$$
\sU^{n-1} \triangleq  \{\bx \in \real^n\mid \norm{\bx} = 1\},
$$
which contains the unit ($\ell_2$) norm vectors in $\real^n$.
Note the distinction  between the unit ball in $\real^n$ and the $n$-dimensional sphere, for the given norm $\norm{\cdot}: \real^n\rightarrow \real$:
\begin{subequations}
\begin{align}
\sB_{\norm{\cdot}}[\bzero, 1] &\triangleq \{\bx \in \real^n\mid \norm{\bx} \leq 1\};\\
\sU^{n-1} &\triangleq   \{\bx \in \real^n\mid \norm{\bx} = 1\}.
\end{align}
\end{subequations}
\begin{corollary}[Covering number of a sphere]\label{corollary:covnum_sphere}
The covering and packing numbers of the $n$-dimensional sphere $\sU^{n-1}$ at scale $\varepsilon<1$ satisfy
\begin{equation}\label{equ:covnum_sphere_eq1}
\mathcalC\left( \sU^{n-1}, \norm{\cdot}, \varepsilon \right)
\leq\mathcalP\left( \sU^{n-1}, \norm{\cdot}, \varepsilon \right) 
\leq \left( \frac{2 + \varepsilon}{\varepsilon} \right)^n \leq \left( \frac{3}{\varepsilon} \right)^n.
\end{equation}
\end{corollary}
We construct an $\varepsilon$-packing set $\sP_\varepsilon \subseteq \sU^{n-1}$ recursively as follows:
\begin{itemize}
\item We initialize $\sP_\varepsilon$ to the empty set.
\item We choose a point $\balpha \in \sU^{n-1}$ such that $\norm{\balpha - \bbeta} > \varepsilon$ for any $\bbeta \in \sP_\varepsilon$. We add $\balpha$ to $\sP_\varepsilon$ until there are no points in $\sU^{n-1}$ that are $\varepsilon$ away from any point in $\sP_\varepsilon$.
\end{itemize}
This procedure terminates after finitely many steps because $ \sU^{n-1}$ is compact; an infinite sequence would admit a convergent subsequence, contradicting the separation condition.

\section{Concentration of Gaussian Vectors}\label{section:conce_gausvec}

In low dimensions, the joint p.d.f. of i.i.d. standard Gaussian vectors is highly concentrated near the origin. However, this behavior changes dramatically as the dimension of the ambient space increases.
Specifically, consider an $n$-dimensional vector $\rvx=[\rx_1, \rx_2, \ldots,\rx_n]^\top$ whose entries are i.i.d. standard Gaussian random variables. 
\footnote{Note again that we use normal fonts of boldface lowercase letters to denote random vectors, and normal
fonts of boldface uppercase letters to denote random matrices. That is, $\rx, \rva, \rmX$ are random scalars,
vectors, or matrices; while $x, \ba, \bX$ are scalars, vectors, or matrices. In many cases, the two terms can be
used interchangeably; that is, $\rx = x$ denotes a realization of the variable.
}
The squared $\ell_2$-norm of $\rvx$, given by $\normtwo{\rvx}^2=\sum_{i=1}^n \rx_i^2$, follows a Chi-squared ($\chi^2$) distribution with $n$ degrees of freedom (see Definition~\ref{definition:chisquare_dist}). As $n$ grows, the p.d.f. of this $\chi^2$ random variable becomes increasingly concentrated around its mean, which equals $n$:
\begin{align}
\Exp[\normtwo{\rvx}^2]
&= \Exp \left[\sum_{i=1}^n \rx_i^2\right]
= \sum_{i=1}^n \Exp[\rx_i^2]
= n. 
\end{align}

Moreover, we present simple results concerning the expected values of the $\ell_1$-, $\ell_2$- and $\ell_\infty$-norms of standard Gaussian random vectors. The bound for the $\ell_\infty$-norm relies on properties of sub-Gaussian random variables discussed in Section~\ref{section:subgauss_rdv}.
\begin{theoremHigh}[Norms of standard Gaussian vectors]\label{theorem:norm_norm_vec}
Let $\rvx = [\rx_1, \rx_2, \ldots, \rx_n]^\top$ be a vector of (not necessarily independent) standard Gaussian random variables. Then,
\begin{enumerate}[(i)]
\item  $\Exp[\normone{\rvx}] = \sqrt{\frac{2}{\pi}} n$.

\item 
It holds that 
\begin{equation}\label{theorem:norm_norm_vec_eq1}
\Exp[\normtwo{\rvx}^2] = n
\qquad \text{and}\qquad 
\sqrt{\frac{2}{\pi}} \sqrt{n} \leq \Exp[\normtwo{\rvx}] \leq \sqrt{n}.
\end{equation}
If the entries of $\rvx$ are independent, then
\begin{equation}\label{theorem:norm_norm_vec_eq2}
\frac{n}{\sqrt{n+1}} \leq \Exp[\normtwo{\rvx}] \leq \sqrt{2} \frac{\Gamma((n+1)/2)}{\Gamma(n/2)} \leq \sqrt{n},
\end{equation}
and consequently $\Exp[\normtwo{\rvx}] \sim \sqrt{n}$ as $n \to \infty$.

\item It holds that
\begin{equation}\label{theorem:norm_norm_vec_eq3}
\Exp [\max_{i \in \{1,2,\ldots,n\}} \rx_i] \leq \sqrt{2 \ln(n)}  
\qquad \text{and} \qquad 
\Exp[\norminf{\rvx}] \leq \sqrt{2 \ln(2n)}.
\end{equation}
If the entries of $\rvx$ are independent then, for $n \geq 2$, then
\begin{equation}\label{theorem:norm_norm_vec_eq4}
\Exp[\norminf{\rvx}] \geq C \sqrt{\ln(n)}, \quad \text{with $C \approx 0.265$}.
\end{equation}
\end{enumerate}

\end{theoremHigh}
\begin{proof}[of Theorem~\ref{theorem:norm_norm_vec}]
\textbf{(i).} Using the density of a standard Gaussian random variable,
$$
\Exp[\abs{\rx_i}] 
= \frac{1}{\sqrt{2\pi}} \int_{-\infty}^\infty \abs{x} \exp(-\frac{x^2}{2}) dx 
= \sqrt{\frac{2}{\pi}} \int_0^\infty x \exp(-\frac{x^2}{2}) dx 
= \sqrt{\frac{2}{\pi}}.
$$
By linearity of expectation, $\Exp[\normone{\rvx}] = \sum_{i=1}^n \Exp[\abs{\rx_i}] = \sqrt{2/\pi} n$.

\paragraph{(ii).} 
Clearly, $\Exp[\normtwo{\rvx}^2] = \sum_{i=1}^n \Exp [\rx_i^2] = n$ for standard Gaussian random variables $\rx_i$.
Jensen's inequality for expectations of the square function gives $\Exp[\normtwo{\rvx}] \leq \sqrt{\Exp[\normtwo{\rvx}^2]} = \sqrt{n}$, while the Cauchy--Schwarz inequality for the inner product on $\real^n$ gives $\Exp[\normtwo{\rvx}] \geq \frac{1}{\sqrt{n}} \Exp[\normone{\rvx}] = \sqrt{2/\pi} \sqrt{n}$; see Exercise~\ref{exercise:cauch_sc_l1l2} for more details.

If the entries of $\rvx$ are independent, then $\normtwo{\rvx}^2$ follows the $\chi^2_{(n)}$ Chi-squared distribution (Definition~\ref{definition:chisquare_dist}). 
Therefore,
$$
\Exp[\normtwo{\rvx}] =
\int_0^\infty x^{\frac{1}{2}} \cdot \frac{x^{\frac{n}{2}-1}}{2^{\frac{n}{2}} \Gamma(\frac{n}{2})} e^{-\frac{x}{2}} dx
= \frac{2^{\frac{n}{2}+\frac{1}{2}}}{2^{\frac{n}{2}} \Gamma(\frac{n}{2})} \int_0^\infty t^{\frac{n}{2}-\frac{1}{2}} e^{-t} dt = \sqrt{2} \frac{\Gamma((n+1)/2)}{\Gamma(n/2)},
$$
where the last equality follows from  that the Gamma function is $\Gamma(x)=\int_{0}^{\infty} t^{x-1}\exp(-t)dt$.
Let $T_n \triangleq \Exp[\normtwo{\rvx}] \leq \sqrt{n}$, where the upper bound for Gaussian vector $\rvx$ of length $n$ was already shown above.
By the fact that $\Gamma(x+1)=x\Gamma(x)$ for $x>0$, we have 
$$
T_{n+1} T_n =2 \frac{\Gamma(n/2 + 1)}{\Gamma(n/2)} = n,
$$
which implies $T_n = n/T_{n+1} \geq n/\sqrt{n+1}$, since $T_{n+1}\leq \sqrt{n+1}$.
Therefore, $n/\sqrt{n+1}\leq T_n\leq \sqrt{n}$, and this proves the desired result.

\paragraph{(iii).}  
The bounds in \eqref{theorem:norm_norm_vec_eq3} follow from {Proposition~\ref{proposition:max_zero_sg}} (on maxima of sub-Gaussian variables). 
Indeed, Lemma~\ref{lemma:mom_stand_norm} shows that $\Exp [\exp(\theta \rx_i)] \leq \exp(\theta^2/2)$ for all $\theta\in\real$ so that the parameter $C$ is $C=1/2$ in Proposition~\ref{proposition:max_zero_sg} for Gaussian random variables.

For the lower bound under independence,  by Exercise~\ref{exercise:exp_as_int}
\begin{align*}
\Exp[\norminf{\rvx}] 
&= \int_0^\infty \Pr\left(\max_{i \in \{1,2,\ldots,n\}} \abs{\rx_i} > x\right) dx 
= \int_0^\infty \left(1 - \Pr\left(\max_{i \in \{1,2,\ldots,n\}} \abs{\rx_i} \leq x\right)\right) dx \\
&= \int_0^\infty \left(1 - \prod_{i=1}^n \Pr(\abs{\rx_i} \leq x)\right) dx 
\geq \int_0^\delta \big(1 - \big(\Pr(\abs{\rx} \leq x)\big)^n\big) dx \\
&\geq \delta \big(1 - \big(1 - \Pr(\abs{\rx} > \delta)\big)^n\big), 
\end{align*}
where $\delta >0$. For any $\delta>0$, we also have 
$$
\Pr(\abs{\rx} > \delta) = \sqrt{\frac{2}{\pi}} \int_\delta^\infty e^{-t^2/2} dt \geq \sqrt{\frac{2}{\pi}} \int_\delta^{2\delta} e^{-t^2/2} dt \geq \sqrt{\frac{2}{\pi}} \delta e^{-2\delta^2}.
$$
Choosing $\delta = \sqrt{\ln n / 2}$, for $n \geq 2$,
\begin{align*}
\Exp[\norminf{\rvx}] 
&\geq \sqrt{\frac{\ln n}{2}} \left(1 - \left(1 - \sqrt{\frac{\ln n}{\pi}} \frac{1}{n}\right)^n\right) 
\geq \sqrt{\frac{\ln n}{2}} \left(1 - \exp\left(-\sqrt{\frac{\ln n}{\pi}}\right)\right) ,
\end{align*}
where the second inequality follows since $(1-\frac{x}{n})^n \leq \exp(-x)$ if $n>x>0$.
For $n\geq 2$, the factor in parentheses is bounded below by a positive constant, yielding the claimed lower bound with $C = (1 - \exp(-\sqrt{\ln(2)/\pi})) / \sqrt{2} \approx 0.265$.
\end{proof}

We can further extend part (iii) of the previous theorem to bound the maximum squared $\ell_2$-norm over a sequence of standard Gaussian random vectors.
\begin{theoremHigh}[Maximum of norms of Gaussian vectors]\label{theorem:max_norm_normavec}
Let $\rvx_1, \rvx_2, \ldots, \rvx_n \in \real^p$ be a sequence of (not necessarily independent) standard Gaussian random vectors. Then, for any $\nu > 0$,
$$
\Exp \left[\max_{i \in \{1,2,\ldots,n\}} \normtwo{\rvx_i}^2\right] \leq 2(1 + \nu) \ln(n) + p(1 + \nu) \ln(1 + \nu^{-1}).
$$
Consequently,
$$
\Exp\left[\max_{i \in \{1,2,\ldots,n\}} \normtwo{\rvx_i}^2\right] \leq (\sqrt{2 \ln(n)} + \sqrt{p})^2.
$$
\end{theoremHigh}
\begin{proof}[of Theorem~\ref{theorem:max_norm_normavec}]
By the concavity of the logarithm and Jensen’s inequality,  for any $\theta > 0$,
\begin{align*}
\Exp \left[\max_{i } \normtwo{\rvx_i}^2\right]
&= \theta^{-1} \Exp\left[\ln \max_{i } \exp(\theta \normtwo{\rvx_i}^2)\right] 
\leq \theta^{-1} \ln \Exp \left[\max_{i } \exp(\theta \normtwo{\rvx_i}^2)\right] \\
&\leq \theta^{-1} \ln \left( \sum_{i=1}^{n} \Exp [\exp(\theta \normtwo{\rvx_i}^2)]\right),
\end{align*}
For each standard Gaussian random vector $\rvx=\rvx_i$, 
by the independence of the components of $\rvx$ and Lemma~\ref{lemma:mom_stand_norm}, we have for $\theta< 1/2$, 
\begin{align*}
\Exp [\exp(\theta \normtwo{\rvx}^2) ]
&= \Exp \left[\exp\left(\theta \sum_{i=1}^p \rx_i^2\right)\right] 
= \Exp \left[\prod_{i=1}^p \exp(\theta \rx_i^2)\right] 
= \prod_{i=1}^p \Exp [\exp(\theta \rx_i^2)]
\leq  (1 - 2\theta)^{-\frac{p}{2}},
\end{align*}
Therefore,
$$
\Exp \left[\max_{i } \normtwo{\rvx_i}^2 \right]
\leq \min_{0 < \theta < 1/2} \theta^{-1} \left(\ln n + \frac{p}{2} \ln \left((1 - 2\theta)^{-1}\right)\right).
$$
Choosing $\theta = (2 + 2\nu)^{-1}$ yields the first claim. Using the fact that $\ln(1 + \nu^{-1}) \leq \nu^{-1}$ and choosing $\nu = \sqrt{p / (2 \ln(n))}$, we obtain
\begin{align*}
\Exp \left[\max_{i } \normtwo{\rvx_i}^2\right] 
&\leq 2(1 + \nu) \ln(n) + p(1 + \nu^{-1})\\
&= 2 \ln(n) + 2 \sqrt{2p \ln(n)} + p = (\sqrt{2 \ln(n)} + \sqrt{p})^2.
\end{align*}
This completes the proof.
\end{proof}

The following result quantifies the variability of the squared $\ell_2$-norm of a standard Gaussian vector.
\begin{proposition}[Variance of the squared $\ell_2$-norm of a Gaussian vector]\label{proposition:var_l2norm_gvec}
Let $\rvx\in\real^n$ be an i.i.d. standard Gaussian random vector. The variance of $\normtwo{\rvx}^2$ is $2n$.
\end{proposition}
\begin{proof}[of Proposition~\ref{proposition:var_l2norm_gvec}]
Recall that $\Var[\normtwo{\rvx}^2] = \Exp[(\normtwo{\rvx}^2)^2] - (\Exp[\normtwo{\rvx}^2])^2$. 
We compute:
\begin{align*}
\Exp\left[\left(\normtwo{\rvx}^2\right)^2\right] 
&= \Exp\left[\left(\sum_{i=1}^n \rx_i^2\right)^2\right]
= \Exp\left[\sum_{i=1}^n \sum_{j=1}^n \rx_i^2 \rx_j^2\right]  
= \sum_{i=1}^n \sum_{j=1}^n \Exp[\rx_i^2 \rx_j^2]\\
&= \sum_{i=1}^n \Exp[\rx_i^4] + \sum_{i\neq j}^{n,n} \Exp[\rx_i^2]\Exp[\rx_j^2] 
\stackrel{\dag}{=} 3n + n(n-1)   
= n(n+2), 
\end{align*}
where the equality ($\dag$) follows since the 4-th moment of a standard Gaussian equals 3.
Thus,
$\Var[\normtwo{\rvx}^2] = n(n+2) - n^2 = 2n$ by Theorem~\ref{theorem:norm_norm_vec}.
\end{proof}

This result shows that as $n$ grows, the relative fluctuations of $ \normtwo{\rvx}$ around its mean  decay like $1/\sqrt{n}$. In other words, the squared norm becomes increasingly concentrated near $n$. The following theorem makes this concentration precise by providing a non-asymptotic tail bound.

\begin{theoremHigh}[Chebyshev tail bound for the $\ell_2$-norm of an i.i.d. standard Gaussian vector\index{Chebyshev tail bound}]\label{thm:cheby_tail_stdgausvec}
Let $\rvx\in\real^n$ be an i.i.d. standard Gaussian random vector. 
Then, for any $\varepsilon > 0$, 
\begin{equation}\label{equation:gauvec_chebyshev_bound}
\Pr\left(n(1 - \varepsilon) \leq \normtwo{\rvx}^2 \leq n(1 + \varepsilon)\right) \geq 1 - \frac{2}{n\varepsilon^2}. 
\end{equation}
\end{theoremHigh}
\begin{proof}[of Theorem~\ref{thm:cheby_tail_stdgausvec}]
Let $\ry\triangleq \normtwo{\rvx}^2$. Then Theorem~\ref{theorem:norm_norm_vec} shows $\ry-n = \ry-\Exp[\ry]$. 
Markov's inequality (Theorem~\ref{theorem:markov-inequality}) shows that 
$$
\Pr\big((\ry-\Exp[\ry])^2 \geq n^2\varepsilon^2 \big) \leq 
\frac{\Exp[(\ry-\Exp[\ry])^2]}{n^2\varepsilon^2} 
=\frac{\Var[\ry]}{n^2\varepsilon^2} =\frac{2n}{n^2\varepsilon}.
$$
Taking the complement yields the desired bound~\eqref{equation:gauvec_chebyshev_bound}.
\end{proof}

The bound in Theorem~\ref{thm:cheby_tail_stdgausvec} relies only on the variance to control the probability that the squared norm deviates from its mean. Consequently, it is significantly weaker than the following result, which leverages the favorable behavior of higher moments of the standard Gaussian distribution.
\begin{theoremHigh}[Chernoff tail bound for the $\ell_2$-norm of an i.i.d. standard Gaussian vector\index{Chernoff tail bound}]\label{thm:chernoff_tail_stdgausvec}
Let $\rvx\in\real^n$ be an i.i.d. standard Gaussian random vector. 
For any $\varepsilon \in (0,1)$, we have
\begin{equation}
\Pr\left(n(1 - \varepsilon) \leq \normtwo{\rvx}^2 \leq n(1 + \varepsilon)\right) \geq 1 - 2\exp\left(-\frac{n\varepsilon^2}{8}\right). \label{eq:chernoff_bound}
\end{equation}
\end{theoremHigh}

\begin{proof}[of Theorem~\ref{thm:chernoff_tail_stdgausvec}]
Let $\ry \triangleq \normtwo{\rvx}^2$. The desired inequality follows from the two one-sided bounds:
\begin{align*}
\Pr(\ry \geq n(1 + \varepsilon)) &\leq \exp\left(-\frac{n\varepsilon^2}{8}\right);  \\
\Pr(\ry \geq  n(1 - \varepsilon)) &\leq \exp\left(-\frac{n\varepsilon^2}{8}\right). 
\end{align*}
We prove the first; the second is analogous.
Fix $ \theta\in(0,1/2)$.  
By Markov's inequality (Theorem~\ref{theorem:markov-inequality}) and  the independence of the components of $\rvx$, we have for any $\tau>0$:
\begin{align*}
\Pr(\ry > \tau) 
&= \Pr(\exp(\theta y) > \exp(\theta \tau  )) 
\leq \exp(-\theta \tau)\,\Exp[\exp(\theta y)]  \\
&\leq \exp(-\theta \tau)\,\Exp\left[\exp\left(\theta \sum_{i=1}^n \rx_i^2\right)\right] 
= \exp(-\theta \tau)\prod_{i=1}^n \Exp[\exp(\theta\rx_i^2)] 
= \frac{\exp(-\theta \tau)}{(1 - 2\theta)^{n/2}}, 
\end{align*}
where the last equality is a consequence of the the moment lemma of standard Gaussian variables (Lemma~\ref{lemma:mom_stand_norm}).
Set $\tau \triangleq n(1 + \varepsilon)$ and
$\theta \triangleq \frac{1}{2} - \frac{1}{2(1 + \varepsilon)}$.
By taking the infimum over $\theta \in (0, 1/2)$, we have 
$$
\Pr(\ry > n(1 + \varepsilon)) 
\leq (1 + \varepsilon)^{n/2} \exp\left(-\frac{n\varepsilon}{2}\right)
= \exp\left(-\frac{n}{2}(\varepsilon - \ln(1 + \varepsilon))\right) 
\leq \exp\left(-\frac{n\varepsilon^2}{8}\right), 
$$
where the last equality follows from the fact that the function $g(x) \triangleq x - \frac{x^2}{4} - \ln(1 + x)$ is nonnegative between 0 and 1 (the derivative is nonnegative and $g(0) = 0$).
A similar argument gives the lower-tail bound, completing the proof.
\end{proof}

\paragrapharrow{Tail bound after projection.}
It is well known that the $\ell_2$-norm of the projection of i.i.d. standard Gaussian noise onto a fixed subspace concentrates around the square root of the subspace dimension. This is because, in any orthonormal basis of the subspace, the projected coordinates remain i.i.d. standard Gaussian.

\begin{lemma}[Projection of an i.i.d. Gaussian vector onto a subspace]\label{lem:proj_iidgaus}
Let $\mathcalV$ be a $k$-dimensional subspace of $\real^n$, and let $\rvx \in \real^n$ be a vector with i.i.d. standard Gaussian entries.
Then the squared norm of the orthogonal projection of  $\rvx$ onto $\mathcalV$, denoted $\normtwo{\mathcalP_{\mathcalV} \rvx}^2$, follows 
a $\chi^2$  distribution  with $k$ degrees of freedom. 
Equivalently, it has the same distribution as
\begin{equation}
\ry \triangleq \sum_{i=1}^k \ry_i^2 \label{eq:chi_squared_sum}
\end{equation}
where $\ry_1, \ry_2, \ldots, \ry_k$ are i.i.d. standard Gaussian random variables.
\end{lemma}

\begin{proof}[of Lemma~\ref{lem:proj_iidgaus}]
Let $\bQ \in \real^{n \times k}$ be a matrix whose columns form an orthonormal basis for $\mathcalV$ (i.e., $\bQ^\top\bQ=\bI_k$).
Then $\bQ\bQ^\top$ is the orthogonal projection matrix for the subspace $\mathcalV$ (Lemma~\ref{lemma:orthogo_genspa}). 
Thus, we have
\begin{align*}
\normtwo{\mathcalP_{\mathcalV} \rvx}^2 
&= \normtwo{\bQ\bQ^\top \rvx}^2 
= \rvx^\top \bQ\bQ^\top \bQ\bQ^\top \rvx 
= \rvx^\top \bQ\bQ^\top \rvx  
\triangleq \rvz^\top \rvz  
= \sum_{i=1}^k \rz_i^2, 
\end{align*}
where by Lemma~\ref{lemma:affine_mult_gauss}, the random vector $\rvz \triangleq \bQ^\top \rvx \in\real^k$ is Gaussian with mean zero and covariance matrix
\begin{align*}
\Sigma_{\rvz} 
= \bQ^\top \Cov[\rvx]\bQ 
= \bQ^\top \bI \bQ 
= \bI_k. 
\end{align*}
Therefore, the entries are independent standard Gaussians. 
\end{proof}

Because the projected norm squared is distributed as a $\chi^2_{(k)}$ random variable, we can directly apply Theorem~\ref{thm:chernoff_tail_stdgausvec} with $n$ replaced by $k$.

\begin{corollary}[Chernoff tail bound for the $\ell_2$-norm of projection an i.i.d. standard Gaussian vector]\label{corollary:chernf_proj_iidgaus}
Let $\mathcalV$ be a $k$-dimensional subspace of $\real^n$, and let $\rvx \in \real^n$ a vector of i.i.d. Gaussian entries. 
For any $\varepsilon \in (0,1)$, we have
$$
\Pr\left(\sqrt{k(1 - \varepsilon)} \leq \normtwo{\mathcalP_{\mathcalV} \rvx} \leq \sqrt{k(1 + \varepsilon)}\right)
\geq  1 - 2\exp\left(-\frac{k\varepsilon^2}{8}\right).
$$
\end{corollary}
This follows immediately from Lemma~\ref{lem:proj_iidgaus} and Theorem~\ref{thm:chernoff_tail_stdgausvec}, applied to the $k$-dimensional Gaussian vector $\bQ^\top\rvx$.

\section{Randomized Projections under Gaussian Random Matrices}
Previously, we showed that projections of Gaussian random vectors exhibit strong concentration properties.
We now turn to a key property of Gaussian random matrices.

In an introductory linear algebra course, one learns that dimensionality reduction via principal component analysis (PCA) involves projecting data onto low-dimensional subspaces that are optimal in the sense of preserving as much energy (i.e., variance) as possible. The principal directions correspond to the directions of maximum variation in the data and are obtained from the singular value decomposition (SVD) of the data matrix. However, computing the SVD can be computationally expensive---and may even be infeasible---when dealing with streaming data or when the projection must be applied in real time.

In such scenarios, a non-adaptive alternative to PCA is needed: one that selects the projection matrix before observing the data. A simple and effective approach is to use a random linear map, represented by a random matrix $\rmX\in\real^{n\times p}$ whose entries are sampled independently from the standard Gaussian distribution. The following lemma shows that applying such a matrix to any fixed unit vector yields a Gaussian random vector.

\begin{lemma}[Randomized projections under i.i.d. standard Gaussian matrix]\label{lemma:rdproj_iidgauss}
Let $\rmX\in\real^{n\times p}$ be a matrix with i.i.d. standard Gaussian entries. 
If $\bbeta \in \real^p$ is a deterministic vector with $\normtwo{\bbeta}=1$, then $\rmX\bbeta\in\real^n$ is an  i.i.d. standard  Gaussian random vector.
\end{lemma}
\begin{proof}[of Lemma~\ref{lemma:rdproj_iidgauss}]
For each $i\in\{1,2,\ldots,n\}$, the $i$-th entry of $\rmX\bbeta$ is given by $(\rmX\bbeta)_i=\innerproduct{\rvx^i, \bbeta}$, i.e., the inner product between $\bbeta$ and the $i$-th row $\rvx^i$ (viewed as a column vector in $\real^p$).
Since $\rvx^i\sim \normal(\bzero, \bI_p)$ and $\bbeta$ is deterministic, Lemma~\ref{lemma:affine_mult_gauss} implies that $\innerproduct{\rvx^i, \bbeta}$ is Gaussian.
Its mean is zero because $\Exp[\rvx^i]=\bzero$. The variance is
\begin{align}
\Var\left[\innerproduct{\rvx^i, \bbeta}\right]
&= \bbeta^\top \Cov[\rvx^i] \bbeta  
= \bbeta^\top \bI \bbeta 
= \normtwo{\bbeta}^2 
= 1.
\end{align}
Thus, each component of $\rmX\bbeta$ is a standard normal random variable. 
Moreover, the components are independent because each depends only on a distinct row of $\rmX$, and the rows are mutually independent by construction. Hence, $\rmX\bbeta\sim\normal(\bzero, \bI)$.
\end{proof}

A direct consequence of this result is a non-asymptotic concentration bound on the $\ell_2$-norm of $\rmX\bbeta$ for any fixed unit vector $\bbeta$.
\begin{corollary}[Chernoff tail bound for the $\ell_2$-norm of random Gaussian projection]\label{corollary:distor_gaus_mat}
Let $\rmX\in\real^{n\times p}$ be a  matrix with i.i.d. standard Gaussian entries. 
Then, for any deterministic $\bbeta \in \real^p$ with unit norm and any $\varepsilon \in (0,1)$, we have
\begin{equation}\label{equation:distor_gaus_mat}
\Pr\left(\sqrt{n(1 - \varepsilon)} \leq \normtwo{\rmX\bbeta} \leq \sqrt{n(1 + \varepsilon)} \right)
\geq 
1 - 2\exp\left(-\frac{n\varepsilon^2}{8}\right).
\end{equation}
\end{corollary}
This follows immediately from Lemma~\ref{lemma:rdproj_iidgauss}, which shows that $\rmX\bbeta$ is an $n$-dimensional i.i.d. standard Gaussian vector, together with Theorem~\ref{thm:chernoff_tail_stdgausvec} applied to its squared norm.

\section{Randomized Projections and Johnson--Lindenstrauss Lemma}
Dimensionality-reduction techniques are useful when they preserve the information we are interested in. 
In many applications, it is essential that the projection approximately preserves pairwise distances between data points. This property enables us to apply distance-based algorithms---such as nearest neighbors, clustering, or kernel methods---in the lower-dimensional space without significant loss of accuracy.

The following result, known as the Johnson--Lindenstrauss (JL) lemma, guarantees that random projections achieve precisely this: they preserve all pairwise distances among a finite set of points with high probability, even in a non-asymptotic setting. Remarkably, the required target dimension $n$ depends only logarithmically on the number of points $k$, and not at all on the ambient dimension $p$. The proof follows the classical argument of \citet{johnson1984extensions}.

\begin{theoremHigh}[Johnson--Lindenstrauss (JL) lemma\index{JL lemma}]\label{thm:jl_lemma}
Let $\rmX\in\real^{n\times p}$ be a random matrix with i.i.d. standard Gaussian entries. Let $\bbeta_1, \bbeta_2, \ldots, \bbeta_k \in \real^p$ be any fixed (deterministic) set of $k$ vectors. 
Then, for any $\delta \in (0,1)$ and any pair $i\neq j$,
\begin{equation}\label{eq:jl_inequality}
(1 - \delta)\normtwo{\bbeta_i - \bbeta_j}^2 
\leq \normtwo{\frac{1}{\sqrt{n}}\rmX(\bbeta_i - \bbeta_j)}^2 
\leq (1 + \delta)\normtwo{\bbeta_i - \bbeta_j}^2, 
\end{equation}
with probability at least $1-\frac{1}{k}$ as long as
\begin{equation}\label{eq:jl_dimension}
n \geq \frac{24 \ln(k)}{\delta^2}. 
\end{equation}
\end{theoremHigh}
\begin{proof}[of Theorem~\ref{thm:jl_lemma}]
We analyze the action of $\bX$ on the normalized difference vectors
$$
\balpha_{ij} \triangleq \frac{\bbeta_i - \bbeta_j}{\normtwo{\bbeta_i - \bbeta_j}},  
\quad 1 \leq i < j \leq k. 
$$
which satisfy $\normtwo{\balpha_{ij}}=1$ whenever $\bbeta_i \neq \bbeta_j$ (if $\bbeta_i = \bbeta_j$, the distance is zero and trivially preserved). 
Define the concentration event
$$
\mathcal{E}_{ij} = \left\{ n(1 - \delta) 
< \normtwo{\rmX\balpha_{ij}}^2 < n(1 + \delta) \right\}, 
\quad 1 \leq i < j \leq k. 
$$
By Corollary~\ref{corollary:distor_gaus_mat}, each event $\mathcal{E}_{ij}$ holds with high probability under condition~\eqref{eq:jl_dimension}:
$$
\Pr\left(\comple{\mathcal{E}_{ij}}\right) 
\leq 2\exp\left(-\frac{n\delta^2}{8}\right)
\leq \frac{2}{k^3}, 
\quad \text{if }n \geq \frac{24 \ln(k)}{\delta^2}.
$$
There are $\binom{k}{2} = k(k - 1)/2$ such pairs. Applying the union bound (Theorem~\ref{theorem:union_bound_proof}),
\begin{align*}
\Pr\left(\bigcup_{i,j} \comple{\mathcal{E}_{ij}}\right) 
\leq  \sum_{i,j} \Pr\left(\comple{\mathcal{E}_{ij}}\right)
\leq  \frac{k(k-1)}{2} \cdot \frac{2}{k^3} 
\leq  \frac{k-1}{k^2} = \frac{1}{k}.
\end{align*}
Thus, with probability at least $1-1/k$, inequality~\eqref{eq:jl_inequality} holds simultaneously for all pairs $i\neq j$. This completes the proof.
\end{proof}

Many algorithms---such as nearest neighbor search, clustering, and kernel methods---depend critically on the geometry of the data, particularly pairwise distances. The JL lemma guarantees that after projecting the data randomly into a space of dimension $n\ll p$, these distances remain nearly intact. Consequently, algorithmic behavior is preserved even though computations now occur in a dramatically smaller space.

In essence, the JL lemma decouples accuracy from ambient dimensionality: randomness enables extreme data compression while maintaining geometric fidelity. This makes randomized algorithms---for example, Gaussian sketching---simultaneously fast, scalable, simple to implement, and theoretically sound, a rare and powerful combination in algorithm design \citep{woodruff2014sketching, lu2021numerical, lu2021rigorous}.

\begin{figure}[h!]
\centering                      
\vspace{-0.15cm}                 
\subfigtopskip=2pt               
\subfigbottomskip=2pt            
\subfigcapskip=-5pt              
\includegraphics[width=0.98\linewidth]{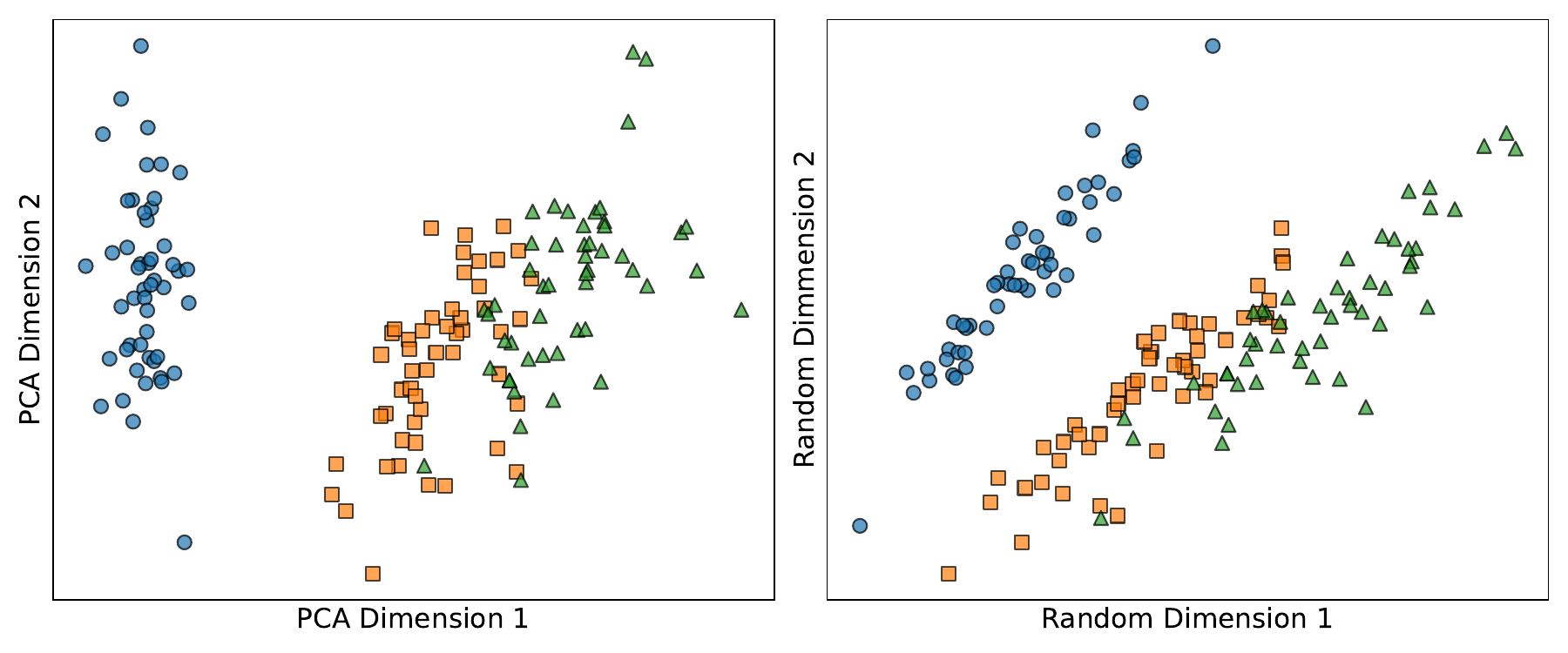}
\caption{
Approximate 2D projection of the Iris dataset using (left) two random Gaussian directions and (right) PCA. Each point represents an iris sample, colored by species.
}
\label{fig:pca_random_iris}
\end{figure}

\begin{example}[Dimensionality reduction via randomized projections] \label{example:random_projection_iris}
We consider the Iris dataset, introduced in \citet{fisher1936use}, which is a benchmark multivariate data set consisting of 150 samples from three iris species---Iris setosa, Iris versicolor, and Iris virginica (50 samples each). Each sample is described by four continuous features: sepal length, sepal width, petal length, and petal width (all in centimeters). The dataset is widely used in machine learning and statistics due to its clear class structure---I. setosa is linearly separable from the other two species, while I. versicolor and I. virginica exhibit partial overlap---making it ideal for illustrating classification, clustering, and dimensionality reduction techniques.
The objective is to project the data onto 2D for visualization. In Figure~\ref{fig:pca_random_iris}, we compare the result of randomly projecting the data by applying a $2 \times 4$ i.i.d. standard Gaussian matrix, with the result of PCA-based dimensionality reduction. 
In terms of keeping the different types of iris species, the randomized projection preserves the structure in the data as effectively as PCA.
\end{example}

This result also reveals a close connection between the concentration inequality~\eqref{equation:distor_gaus_mat} and the JL lemma.
Since~\eqref{equation:distor_gaus_mat} implies the RIP via Theorems~\ref{thm:sv_gaus_mat} and~\ref{theorem:rip_gaussian_uppdeltak}, in this sense, one may say that the JL lemma implies the RIP.

We now show a partial converse: under suitable conditions, a matrix satisfying the RIP can be turned into a JL embedding for a point set $\sS$ by randomizing the signs of its columns.
Indeed, without such sign randomization, the claim would be false. Randomizing the column signs ensures that the probability of the point set $\sS$ falling (or nearly falling) into the null space of $\bX\bD$ is exponentially small.

\begin{theoremHigh}\label{theorem:jl_radema}
Let $\sS \subset \real^p$ be a finite set of cardinality $\abs{\sS} = S$. Fix $\delta \in (0,1)$. 
Let $\rmX \in \real^{n \times p}$ be a matrix satisfying the RIP of order $2k$ with constant $\delta_{2k} \leq \delta/4$ for some integer  $k \geq 16\ln(4S/\varepsilon)$. 
Then, with probability at least $1 - \varepsilon$
$$
(1 - \delta)\normtwo{\bbeta}^2 
\leq \normtwo{\rmX\rmD\bbeta}^2 
\leq (1 + \delta)\normtwo{\bbeta}^2,
 \quad \text{for all } \bbeta \in \sS,
$$
where $\rmD\triangleq \diag(\bepsilon)\in\real^{p\times p}$, and 
$\bepsilon=[\epsilon_1, \epsilon_2, \ldots,\epsilon_p]^\top$ is  a Rademacher vector (i.e., each $\epsilon_i$ is an independent random variable taking values +1 and $-1$ with equal probability).
\end{theoremHigh}
\begin{proof}[of Theorem~\ref{theorem:jl_radema}]
Without loss of generality, assume that all vectors in $\sS$ are normalized: $\normtwo{\bbeta} = 1$. 
Fix any $\bbeta \in \sS$.  
We partition $\bbeta$ into blocks of size $k$ based on the magnitudes of its entries.
Specifically, let $\sT_1 \subset \{1,2,\ldots,p\}$ index the $k$ largest (in absolute value) entries of $\bbeta$, $\sT_2 \subset \{1,2,\ldots,p\} \setminus \sT_1$ index the next $k$ largest entries among the remaining coordinates, and so on. 
This yields a decomposition ${\rmX\rmD\bbeta}={\rmX\rmD \sum_i \bbeta(\sT_i)}$.
\footnote{Note again that $\bbeta(\sS)\in\real^p$ and $\bbeta_{\sS}\in\real^{\abs{\sS}}$; see Definition~\ref{definition:matlabnotation}.}
Now expand the squared norm:
\begin{align}
&\normtwo{\rmX\rmD\bbeta}^2 
= \normtwo{\rmX\rmD \sum_i \bbeta(\sT_i)}^2 \nonumber \\
&= \sum_i \normtwo{\rmX\rmD\bbeta(\sT_i)}^2 
+ 2\innerproduct{\rmX\rmD\bbeta(\sT_1), \rmX\rmD\bbeta(\comple{\sT_1})} 
+ \sum_{\substack{i,j \geq 2 \\ i \neq j}} 
\innerproduct{\rmX\rmD\bbeta(\sT_i), \rmX\rmD\bbeta(\sT_j)}  \label{equation:jl_radema_pv1}.
\end{align}

\paragraph{Bound the first term.}
Because  $\rmX$ satisfies the RIP of order $2k$ with constant $\delta_{2k} \leq \delta/4$, and because $\normtwobig{\rmD\bbeta(\sT_i)} = \normtwobig{\bbeta_{\sT_i}}$, the first term satisfies
$$
(1 - \delta/4)\normtwo{\bbeta}^2 
= (1 - \delta/4)\sum_i \normtwo{\bbeta_{\sT_i}}^2 
\leq \sum_i \normtwo{\rmX\rmD\bbeta(\sT_i)}^2 
\leq (1 + \delta/4)\normtwo{\bbeta}^2.
$$

\paragraph{Bound the second term.}
To bound the second term in~\eqref{equation:jl_radema_pv1}, we define the random variable
$$
\ry \triangleq 
\innerproduct{\rmX\rmD\bbeta(\sT_1), \rmX\rmD\bbeta(\comple{\sT_1})} 
= \innerproduct{\rvu, \bepsilon_{\comple{\sT_1}}} 
= \sum_{i \notin \sT_1} \epsilon_i \ru_i,
$$
where
$
\rvu 
\triangleq 
\diag(\bbeta_{\comple{\sT_1}})
\rmX_{\comple{\sT_1}}^\top 
\rmX_{\sT_1} 
\rmD_{\sT_1} \bbeta_{\sT_1}
\in \real^{\abs{\comple{\sT_1}}}$.
\footnote{Note that $\bD\bbeta = \diag(\bbeta)\bD$.}
Since $\rvu$ and $\bepsilon_{\comple{\sT_1}}$ are stochastically independent, Hoeffding's inequality for Rademacher sequences (Corollary~\ref{corollary:hoeffd_rademac}) shows that
\begin{equation}\label{equation:jl_radema_pv2_1}
\prob(\abs{\ry} \geq \mu) 
\leq 2\exp\left(-\frac{\mu^2}{2\normtwo{\rvu}^2}\right).
\end{equation}
We now bound $\normtwo{\rvu}$ above.
By the definition of dual norm~\eqref{equation:dual_norm_equa}, Cauchy--Schwarz ienquality, and the submultiplicativity of spectral norms (Theorem~\ref{theorem:submul_ortho_matnorm}), we obtain
\begin{align*}
\normtwo{\rvu} 
&= \max_{\normtwo{\ba} \leq 1} 
\innerproduct{\ba, \rvu}
= \max_{\normtwo{\ba} \leq 1} \sum_{i \geq 2} 
\innerproduct{\ba_{\sT_i}, 
\diag(\bbeta_{\sT_i})
\rmX_{\sT_i}^\top 
\rmX_{\sT_1}  \rmD_{\sT_1} 
\bbeta_{\sT_1}} \\
&\stackrel{\dag}{\leq} \max_{\normtwo{\ba} \leq 1} \sum_{i \geq 2} \normtwo{\rmX_{\sT_i}^\top 
\rmX_{\sT_1}} 
\normtwo{\ba_{\sT_i}} \norminf{\bbeta_{\sT_i}} \normtwo{\bbeta_{\sT_1}}
\stackrel{\ddag}{\leq} \frac{\delta_{2k}}{\sqrt{k}} \max_{\normtwo{\ba} \leq 1} 
\sum_{i \geq 2} \normtwo{\ba_{\sT_i}} \normtwo{\bbeta_{\sT_{i-1}}} \\
&\leq \frac{\delta_{2k}}{\sqrt{k}} \max_{\normtwo{\ba} \leq 1} 
\sum_{i \geq 2} \frac{1}{2} \left( \normtwo{\ba_{\sT_i}}^2 + \normtwo{\bbeta_{\sT_{i-1}}}^2 \right)
\leq \frac{\delta_{2k}}{\sqrt{k}} ,
\end{align*}
where the inequality ($\dag$) follows since $\normtwo{\diag(\bbeta)} = \norminf{\bbeta}$ (the largest singular value of the diagonal matrix $\diag(\bbeta)$ is $\norminf{\bbeta}$) and $\normtwo{\bD_{\sT_1}}=\norminf{\bepsilon} = 1$. 
The inequality ($\ddag$) follows from $\norminfbig{\bbeta_{\sT_i}} \leq k^{-1/2} \normtwobig{\bbeta_{\sT_{i-1}}}$ by the construction of the partitioning the ordered sets $\sT_1, \sT_2, \ldots$ (Lemma~\ref{lemma:ells_ellt}),  $\normtwo{\rmX_{\sT_i}^\top \rmX_{\sT_1}} \leq \delta_{2k}$ for $i \geq 2$ by the standard consequence of RIP (Proposition~\ref{proposition:conse_rip}),  and the fact that $\normtwo{\bbeta_{\sT_1}} \leq \normtwo{\bbeta} = 1$. 
Therefore, \eqref{equation:jl_radema_pv2_1} becomes
\begin{equation}\label{equation:jl_radema_pv2}
\prob(\abs{\ry} \geq \mu) 
\leq 2\exp\left(-\frac{\mu^2}{2\normtwo{\rvu}^2}\right)
\leq 2\exp\left(-\frac{k \mu^2 }{2\delta_{2k}^2}\right) 
\leq 2\exp\left(-\frac{8k\mu^2}{\delta^2}\right) .
\end{equation}

\paragraph{Bound the third term.}
Consider the remaining cross terms, we define the following random variable
$$
\rz \triangleq \sum_{\substack{i,j \geq 2 \\ i \neq j}} 
\innerproduct{\rmX\rmD \bbeta(\sT_i), \rmX\rmD \bbeta(\sT_j)} 
= \sum_{i,j \in \{1,2,\ldots,p\}} \epsilon_i \epsilon_j \ra_{ij}
= \bepsilon^\top \rmA \bepsilon ,
$$
where $\rmA \in \real^{p \times p}$ is a symmetric matrix with zero diagonal and off-diagonal entries:
$$
\ra_{ij} = 
\begin{cases}
\beta_i \rvx_i^\top \rvx_j \beta_j 
& \text{if } i, j \in \comple{\sT_1} \text{ and } i, j \text{ are contained in different blocks } \sT_\ell ; \\
0 
& \text{otherwise} .
\end{cases}
$$
The expression $\bepsilon^\top\rmA\bepsilon$ a homogeneous Rademacher chaos of order two.
To bound $\abs{\rz}$,  we apply Problem~\ref{prob:rade_chaos}, which requires estimates on $\normtwo{\rmA}$ and $\normf{\rmA}$.
Since $\rmA$ is symmetric, its spectral norm can be obtained by using the RIP and the block structure as before:
\begin{align*}
\normtwo{\rmA} 
&
= \max_{\normtwo{\ba} \leq 1} \innerproduct{\rmA\ba, \ba}
= \max_{\normtwo{\ba} \leq 1} \sum_{\substack{i,j \geq 2 \\ i \neq j}} 
\innerproduct{
\ba_{\sT_i}, 
\diag(\bbeta)_{\sT_i} 
\rmX_{ \sT_i}^\top \rmX_{ \sT_j} 
\diag(\bbeta)_{\sT_j} 
\ba_{\sT_j} 
} 
\\
&\leq \max_{\normtwo{\ba} \leq 1} 
\sum_{\substack{i,j \geq 2 \\ i \neq j}}
\normtwo{\ba_{\sT_i}} \normtwo{\ba_{\sT_j}} \norminf{\bbeta_{\sT_i}} \norminf{\bbeta_{\sT_j}}
\normtwo{\rmX_{ \sT_i}^\top \rmX_{ \sT_j}} \\
&\leq \delta_{2k} \max_{\normtwo{\ba} \leq 1} 
\sum_{\substack{i,j \geq 2 \\ i \neq j}} 
\normtwo{\ba_{\sT_i}} \normtwo{\ba_{\sT_j}} k^{-1/2} 
\normtwo{\bbeta_{\sT_{i-1}}} k^{-1/2} \normtwo{\bbeta_{\sT_{j-1}}} \\
&\leq \frac{\delta_{2k}}{k} \max_{\normtwo{\ba} \leq 1} 
\sum_{\substack{i,j \geq 2 \\ i \neq j}} 
\Bigg( \frac{\normtwobig{\bbeta_{\sT_{i-1}}}^2 + \normtwo{\ba_{\sT_i}}^2}{2} \Bigg) 
\Bigg( \frac{\normtwobig{\bbeta_{\sT_{j-1}}}^2 + \normtwo{\ba_{\sT_j}}^2}{2} \Bigg)
\leq \frac{\delta_{2k}}{k} 
\leq \frac{\delta}{4k} .
\end{align*}
For the Frobenius norm, observe that
\begin{align*}
\normf{\rmA}^2 
&= \sum_{\substack{s,t \geq 2 \\ s \neq t}} 
\sum_{i \in \sT_s} \sum_{j \in \sT_t} (\beta_i \rvx_i^\top \rvx_j \beta_j)^2 
= \sum_{\substack{s,t \geq 2 \\ s \neq t}} 
\sum_{i \in \sT_s} \beta_i^2 \rvx_i^\top 
\left(\sum_{j \in \sT_t} \beta_j \rvx_j \beta_j\rvx_j^\top\right) \rvx_i \\
&= \sum_{\substack{s,t \geq 2 \\ s \neq t}} 
\sum_{i \in \sT_s} \beta_i^2 \normtwo{\diag(\bbeta_{\sT_t}) \rmX_{ \sT_t}^\top \rvx_i}^2 
\leq \sum_{\substack{s,t \geq 2 \\ s \neq t}} \sum_{i \in \sT_s} \beta_i^2 \norminf{\bbeta_{\sT_t}}^2 \normtwo{\rmX_{ \sT_t}^\top \rvx_i}^2 \\
&\stackrel{\dag}{\leq} \delta_{2k}^2 \sum_{\substack{s,t \geq 2 \\ s \neq t}} \normtwo{\bbeta_{\sT_s}}^2 k^{-1} \normtwo{\bbeta_{\sT_t}}^2 
\leq \frac{\delta_{2k}^2}{k}
\leq \frac{\delta^2}{16 k},
\end{align*}
where the inequality ($\dag$) follows from the fact that  $\normtwo{\rmX_{ \sT_t}^\top \rvx_i} = \normtwo{\rmX_{ \sT_t}^\top \rvx_i} \leq \delta_{k+1} \leq \delta_{2k}$ by the standard consequence of RIP (Proposition~\ref{proposition:conse_rip}). 
Problem~\ref{prob:rade_chaos} shows that $\abs{\rz}$ satisfies
\begin{align*}
\prob(\abs{\rz} \geq \tau) 
&\leq 2\exp\left(-\min\left\{ \frac{3\tau^2}{128\,\normf{\rmA}^2}, \frac{\tau}{32\,\normtwo{\rmA}} \right\} \right) 
\leq 2\exp\left(-k\min\left\{ \frac{3\tau^2}{8\delta^2}, \frac{\tau}{8\delta} \right\} \right) .
\end{align*}
\paragraph{Combining bounds.}
Choosing $\mu \triangleq \delta/4$ and $\tau \triangleq \delta/2$, 
and plugging the three estimates  into \eqref{equation:jl_radema_pv1} shows that
\begin{equation}\label{equation:jl_radema_pv3}
(1 - \delta)\normtwo{\bbeta}^2 
\leq \normtwo{\rmX\rmD\bbeta}^2 
\leq (1 + \delta)\normtwo{\bbeta}^2,
\quad
\text{for any  $\bbeta \in \sS$}
\end{equation}
with probability at least
$$
1 - 2\exp(-k/2) - 2\exp\left(-k\min\{3/32, 1/16\}\right) \geq 1 - 4\exp(-k/16) .
$$
Finally, applying the union bound over all $S$ points in $\sS$, the inequality holds simultaneously for all $\bbeta\in\sS$ with probability at least
$1 - 4S \exp(-k/16) \geq 1 - \varepsilon$
under the condition $k \geq 16\ln(4S/\varepsilon)$. 
This completes the proof. 
\end{proof}

\begin{figure}[h!]
\centering                      
\vspace{-0.35cm}                 
\subfigtopskip=2pt               
\subfigbottomskip=2pt            
\subfigcapskip=-5pt              
\includegraphics[width=0.98\linewidth]{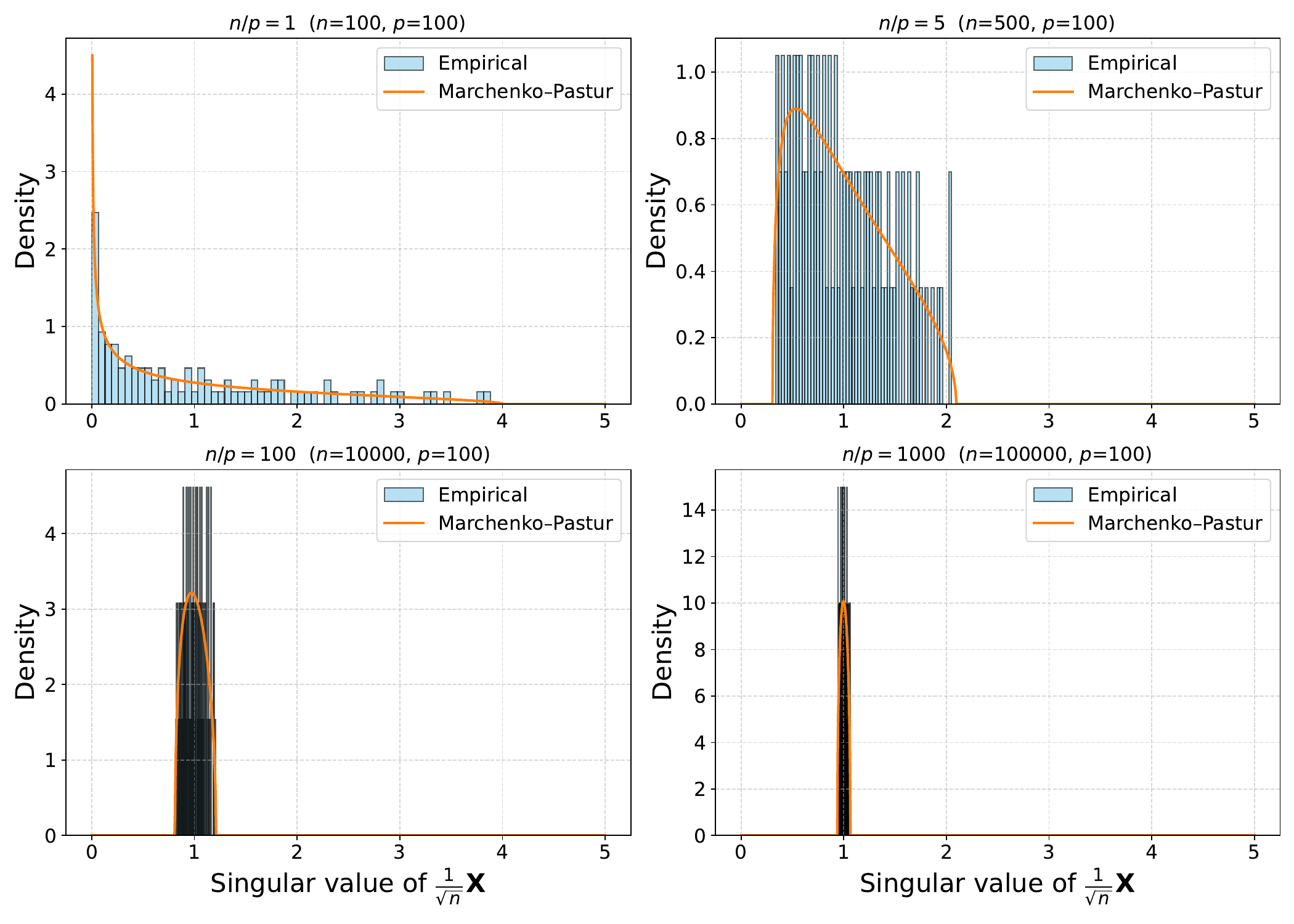}
\caption{Empirical distribution of the singular values of $n\times p$ matrices with i.i.d. standard Gaussian entries as the aspect ratio $n/p$ increases.
Providing also the theoretical Marchenko--Pastur distribution for each empirical distribution.}
\label{fig:singular_values_decay}
\end{figure}

\index{Gaussian matrices}
\section{Singular Values of Gaussian Random Matrices}
In this section, we analyze the singular values of random matrices whose entries are i.i.d. standard Gaussian random variables. 
Specifically, consider an $n\times p$ matrix $\bX$ with $n>p$. Numerical experiments show that as $n$ grows (with $p$ fixed), all $p$ singular values of $\bX$ concentrate tightly around  
$n$, as illustrated in Figure~\ref{fig:singular_values_decay}. Consequently,
\begin{equation}\label{eq:singular_approx}
	\bX \approx \bU (\sqrt{n} \bI) \bV^\top = \sqrt{n}  \bU \bV^\top,
\end{equation}
which means that $\bX$ behaves like a scaled orthogonal matrix. Geometrically, this implies that if we draw a fixed number of i.i.d. Gaussian vectors in increasingly high-dimensional ambient spaces, these vectors become nearly orthogonal with high probability as the dimension grows.

The following result provides a non-asymptotic probabilistic bound on the extreme singular values of such a matrix. The proof uses a covering number argument from~\cite{vershynin2010introduction}, which is flexible enough to extend to other sub-Gaussian distributions and settings.

\begin{theoremHigh}[Singular values of a Gaussian matrix]\label{thm:sv_gaus_mat}
Let $\rmX\in\real^{n \times k}$ be a random matrix with i.i.d. standard Gaussian entries, where  $n > k$. 
For any fixed $\varepsilon \in(0,1)$, the smallest and largest singular values of $\rmX$ satisfy
\begin{equation}\label{equ:sv_gaus_mat_eq1}
\sqrt{n} (1 - \varepsilon) \leq \sigma_k \leq \sigma_1 \leq \sqrt{n} (1 + \varepsilon)
\end{equation}
with probability at least $1 - 2 \left( \frac{12}{\varepsilon} \right)^k \exp \left( -\frac{n \varepsilon^2}{32} \right)  $, 
or
with probability at least $1 - 1/k$ as long as
\begin{equation}\label{equ:sv_gaus_mat_eq2}
n > 
\frac{32}{\varepsilon^2} \ln\left(2k\left(\frac{12}{\varepsilon}\right)^k\right).
\end{equation}
\end{theoremHigh}
The theorem  gives a concrete formula for how large $n$ needs to be (relative to $k$ and $\varepsilon$) to guarantee this behavior with probability at least $1-1/k$.
By the variational characterization of singular values (see~\eqref{equation:svd_stre_bd}), the bounds in~\eqref{equ:sv_gaus_mat_eq1} are equivalent to
\begin{equation}
\sqrt{n}(1 - \varepsilon) < \normtwo{\rmX\bbeta} < \sqrt{n}(1 + \varepsilon),
\quad \text{for all } \bbeta \in \sU^{k-1},
\end{equation}
where $\sU^{k-1} \triangleq  \{\balpha \in \real^k\mid \normtwo{\balpha} = 1\}$ is the $k$-dimensional sphere in $\real^k$, which contains the unit $\ell_2$-norm vectors in $\real^k$. 
Since this set is uncountably infinite, a direct union bound---such as the one used in the proof of the JL lemma---is not applicable. Instead, we use a net-based argument.

\begin{proof}[of Theorem~\ref{thm:sv_gaus_mat}]
Let $\varepsilon_1 \triangleq \varepsilon / 4$ and $\varepsilon_2 \triangleq \varepsilon / 2$. Consider the  $\varepsilon_1$-packing set  $\sP_{\varepsilon_1}$ of $\sU^{k-1}$ as that is constructed below  Corollary~\ref{corollary:covnum_sphere}. 
For any fixed $\bbeta\in\sU^{k-1}$, define the event
$$
\mathcalA_{\bbeta, \varepsilon_2} \triangleq 
\left\{ n(1 - \varepsilon_2) \normtwo{\bbeta}^2 \leq 
\normtwo{\rmX \bbeta}^2 \leq n(1 + \varepsilon_2) \normtwo{\bbeta}^2 \right\}.
$$
By Corollary~\ref{corollary:distor_gaus_mat}, for any fixed $\bbeta \in \sU^{k-1}$, $\Pr \left( \comple{\mathcalA_{\bbeta, \varepsilon_2}} \right) \leq 2 \exp \left( -n \varepsilon^2 / 32 \right)$.
Applying the union bound (Theorem~\ref{theorem:union_bound_proof}) over the finite set $\sP_{\varepsilon_1}$ gives
\begin{align*}
\Pr \left( \bigcup_{\bbeta \in \sP_{\varepsilon_1}} \comple{\mathcalA_{\bbeta, \varepsilon_2}} \right)
&\leq \sum_{\bbeta \in \sP_{\varepsilon_1}} \Pr \left( \comple{\mathcalA_{\bbeta, \varepsilon_2}} \right)
\leq \abs{\sP_{\varepsilon_1}} \cdot \Pr \left( \comple{\mathcalA_{\bbeta, \varepsilon_2}} \right)
\leq 2 \left( \frac{12}{\varepsilon} \right)^k \exp \left( -\frac{n \varepsilon^2}{32} \right)  
\leq \frac{1}{k},  
\end{align*}
where $\abs{\sP_{\varepsilon_1}}\leq \left( {12}/{\varepsilon} \right)^k$ by Corollary~\ref{corollary:covnum_sphere}, and  the last inequality follows from condition \eqref{equ:sv_gaus_mat_eq2}.

Thus, with probability at least $1-1/k$, the event $\mathcalA_{\bbeta, \varepsilon_2}$ holds simultaneously for all $\bbeta\in\sP_{\varepsilon_1}$.
To finish the proof, we need to show that if $\cup_{\bbeta \in \sP_{\varepsilon_1}} \comple{\mathcalA_{\bbeta, \varepsilon_2}}$ holds, then the bound holds for every element in $\sU^{k-1}$, not only for those in the packing set $\sP_{\varepsilon_1}$. 

For any arbitrary vector $\bbeta \in \sU^{k-1}$ on the sphere, there exists a vector in the $\varepsilon/4$-packing set $\balpha \in \sP_{\varepsilon_1}$ such that $\normtwo{\bbeta - \balpha} \leq \varepsilon_1=\varepsilon/4$. 
Assuming the good event $\comple{(\cup_{\balpha \in \sP_{\varepsilon_1}} \comple{\mathcalA_{\balpha,\varepsilon_2}})}$ holds (which has probability at least $1-1/k$), by the triangle inequality and  the definition of the spectral norm (Definition~\ref{definition:spectral_norm}), this implies
\begin{align*}
\normtwo{\rmX \bbeta}
&\leq \normtwo{\rmX \balpha} + \normtwo{\rmX (\bbeta - \balpha)}  
\leq \sqrt{n} \left(1 + \frac{\varepsilon}{2}\right) + \sigma_1 \normtwo{\bbeta - \balpha}  
\leq \sqrt{n} \left(1 + \frac{\varepsilon}{2}\right) + \frac{\sigma_1 \varepsilon}{4}, 
\end{align*}
where the second inequality follows because we assume $\cup_{\balpha \in \sP_{\varepsilon_1}} \comple{\mathcalA_{\balpha,\varepsilon_2}}$ holds.
The definition of the spectral norm also implies  $\sigma_1$ is the smallest upper bound on $\normtwo{\rmX \bbeta}$ for all $\bbeta$ on the sphere $\sU^{k-1}$, so the bound in the above inequality cannot be smaller:
$$
\sigma_1 \leq \sqrt{n} \left(1 + \frac{\varepsilon}{2}\right) + \frac{\sigma_1 \varepsilon}{4}.
$$
Rearranging yields
\begin{equation}\label{eq:sigma1_final}
\sigma_1 
\leq \sqrt{n} \left( \frac{1 + \varepsilon/2}{1 - \varepsilon/4} \right) 
= \sqrt{n} \left( 1 + \varepsilon - \frac{\varepsilon(1 - \varepsilon)}{4 - \varepsilon} \right) 
\leq \sqrt{n}(1 + \varepsilon). 
\end{equation}

For the lower bound, again under the same good event, $\comple{(\cup_{\balpha \in \sP_{\varepsilon_1}} \comple{\mathcalA_{\balpha,\varepsilon_2}})}$  (which has probability at least $1-1/k$), the triangle inequality also implies
\begin{align*}
\normtwo{\rmX \bbeta}
&\geq \normtwo{\rmX \balpha} - \normtwo{\rmX (\bbeta - \balpha)}
\geq \sqrt{n} \left(1 - \frac{\varepsilon}{2}\right) - \sigma_1 \normtwo{\bbeta - \balpha}\\
&\geq \sqrt{n} \left(1 - \frac{\varepsilon}{2}\right) - \frac{\varepsilon}{4} \sqrt{n}(1 + \varepsilon) 
\geq \sqrt{n}(1 - \varepsilon),
\end{align*}
where the penultimate inequality follows form \eqref{eq:sigma1_final}, and the last inequality follows since $\varepsilon<1$.
By the singular value property \eqref{equation:svd_stre_bd}, $\sigma_k$ is the largest lower bound on $\normtwo{\rmX \bbeta}$ for all $\bbeta$ on the sphere $\sU^{k-1}$, so $\sigma_k \geq \sqrt{n}(1 - \varepsilon)$ as long as $\comple{(\cup_{\balpha \in \sP_{\varepsilon_1}} \comple{\mathcalA_{\balpha,\varepsilon_2}})}$ holds.
This completes the proof.
\end{proof}

\section{Sparse Optimization with Gaussian Random Matrices}\label{section:spar_gauss}

With these tools in hand, we can easily establish the recovery guarantee of Gaussian random matrices in compressed sensing problems.

A common strategy to prove that compressed sensing succeeds for a given class of signals is to construct a dual certificate (see Theorem~\ref{theorem:dualcert_p1} and its application in Theorem~\ref{theorem:recov_rip_l1_t1} using RIP to find the certificate) and verify that it is valid for every signal in that class. We illustrate this approach using Gaussian measurement matrices.

\begin{theoremHigh}[Exact recovery of \eqref{opt:p1} under standard Gaussian matrix \citep{baraniuk2008simple}]\label{thm:recovery_l1_standnorm}
Let $\rmX \in \real^{n\times p}$ be a random matrix with i.i.d. standard Gaussian entries, and let the measurements $ \by=\rmX\bbeta^*\in \real^n$, where $\bbeta^*\in\real^p$ is $k$-sparse (i.e., has at most $k$ nonzero entries). 
Then $\bbeta^*$ is the unique solution to the $\ell_1$-minimization problem~\eqref{opt:p1} (p.~\pageref{opt:p1}) with probability at least $1 - \frac{1}{p}$, provided that the number of measurements satisfies
\begin{equation}
n \geq C k \ln p, 
\end{equation}
for some universal constant $C$.
\end{theoremHigh}
\begin{proof}[of Theorem~\ref{thm:recovery_l1_standnorm}]
The proof relies on the dual certificate condition for $\ell_1$-minimization (Theorem~\ref{theorem:dualcert_p1}).
According to that result, it suffices to show that for any support set $\sS\subseteq\{1,2,\ldots,p\}$ with $\abs{\sS}=k$ and any possible sign pattern $\bw \triangleq \sgn(\bbeta^*_\sS) \in \real^k$ of the nonzero entries of $\bbeta^*$, there exists a valid dual certificate $\blambda\in\real^n$ satisfying
$$
\rmX_\sS^\top \blambda = \bw. 
$$
We construct the candidate certificate as the least squares solution (see \eqref{equation:l1_dualcert_lambda}):
$\blambda_{\text{LS}} = \rmX_\sS (\rmX_\sS^\top \rmX_\sS)^{-1} \bw $.
To control $\blambda_\mathrm{LS}$, we resort to the bound on the singular values of a fixed $n \times k$ submatrix  by Theorem~\ref{thm:sv_gaus_mat}. Setting $\varepsilon \triangleq 0.5$, 
we denote by $\mathcal{E}$ the event that
$$
\mathcal{E}\triangleq 
\{0.5\sqrt{n} \leq \sigma_k \leq \sigma_1 \leq 1.5\sqrt{n}\},
$$
where
$$
\Pr(\mathcal{E}) 
\geq 1 - 2 \left( \frac{12}{\varepsilon} \right)^k \exp \left( -\frac{n \varepsilon^2}{32} \right) 
\approx 1 - \exp\left(-C' \frac{n}{k} \right) 
$$
for some absolute constant $C'>0$.
Conditioned on $\mathcal{E}$, the submatrix $\rmX_\sS$ has full rank and $\rmX_\sS^\top \rmX_\sS$ is invertible, 
so $\blambda_{\mathrm{LS}}$ is well-conditioned and guarantees $\rvx_i^\top\blambda_{\text{LS}} = \sign(\beta_i^*)$ for $i\in\sS$. 
To verify the second condition $\norminf{\bX_{\comple{\sS}}^\top\blambda_{\text{LS}}}<1$ (see Theorem~\ref{theorem:dualcert_p1}),  consider any index $i\in\comple{\sS}$, and let $\rvx_i$ denote the $i$-th column of $\rmX$. 
We need to bound $\abs{\rvx_i^\top \blambda_{\mathrm{LS}}}$ for all indices $i \in \comple{\sS}$.
Let $\rmU\bSigma\rmV^\top$ be the reduced SVD of $\rmX_\sS$. 
Conditioned on $\mathcal{E}$, we have
\begin{equation}\label{equ:recovery_l1_standnorm_pv1}
\normtwo{\blambda_{\mathrm{LS}}} 
= \normtwo{\rmV \bSigma^{-1} \rmU^\top \bw } 
\leq \frac{\normtwo{\bw}}{\sigma_k}
\leq 2 \sqrt{\frac{k}{n}}.
\end{equation}
For a fixed $i \in \comple{\sS}$ and a fixed vector $\bv \in \real^n$, $\rvx_i^\top \bv / \normtwo{\bv}$ is a standard Gaussian random variable by the definition of Gaussian random vectors and the linear transformation lemma  (Lemma~\ref{lemma:lin_tran_prob} or Lemma~\ref{lemma:affine_mult_gauss}), which, by Exercise~\ref{exercise:chernoff_gausvari}, implies
\begin{align*}
\Pr \left( \abs{\rvx_i^\top \bv} \geq 1 \right) 
&= \Pr \left( \frac{\abs{\rvx_i^\top \bv}}{\normtwo{\bv}} \geq \frac{1}{\normtwo{\bv}} \right)  
\leq 2 \exp \left( - \frac{1}{\normtwo{\bv}^2} / 2 \right) .
\end{align*}
Note that if $i \notin \sS$, then $\rvx_i$ and $\blambda_{\mathrm{LS}}$ are independent (they depend on different and hence independent entries of $\rmX$). This means that due to \eqref{equ:recovery_l1_standnorm_pv1},
$$
\Pr \left( \abs{\rvx_i^\top \blambda_{\mathrm{LS}}} \geq 1 \mid \mathcal{E} \right) 
= \Pr \left( \abs{\rvx_i^\top \bv} \geq 1  \right)
\leq 2 \exp \left( -\frac{n}{8k} \right),
\quad\text{for $ \normtwo{\bv} \leq 2 \sqrt{\frac{k}{n}}$}, \; i\notin \sS.
$$
As a result,
\begin{align*}
\Pr \left( \abs{\rvx_i^\top \blambda_{\mathrm{LS}}} \geq 1 \right) 
&\leq \Pr \left( \abs{\rvx_i^\top \blambda_{\mathrm{LS}}} \geq 1 \mid \mathcal{E} \right) + \Pr(\comple{\mathcal{E}}) 
\leq 2 \exp \left( -\frac{n}{8k} \right) + \exp \left( -C' \frac{n}{k} \right). 
\end{align*}
Applying the union bound over all $i \in \comple{\sS}$ (at most $p$ indices),
$$
\Pr \left( \bigcup_{i \in \comple{\sS}} \left\{ \abs{\rvx_i^\top \blambda_{\mathrm{LS}}} \geq 1 \right\} \right) 
\leq p \left( 2 \exp \left( -\frac{n}{8k} \right) + \exp \left( -C' \frac{n}{k} \right) \right). \label{eq:union_bound}
$$
Choosing a constant $C$ so that 
$n \geq C k \ln p$
ensures this probability is at most $\frac{1}{p}$. 
Hence, with probability at least $1-1/p$, the dual certificate conditions hold for all $k$-sparse vectors, which implies exact recovery via $\ell_1$-minimization.
\end{proof}

Rather than constructing a dual certificate, one can alternatively use the RIP to guarantee exact recovery. While verifying RIP  for a deterministic matrix is computationally intractable---requiring checks over all $\binom{p}{k}$ submatrices---it can be shown that random matrices satisfy RIP with high probability. The next two theorems establish this for i.i.d. Gaussian matrices.

\begin{theoremHigh}[RIP for Gaussian matrices]\label{theorem:rip_gaussian_uppdeltak}
Let $\rmX \in \real^{n \times p}$ have i.i.d. standard Gaussian entries. 
For any $\eta, \varepsilon \in (0,1)$, suppose
\begin{equation}\label{equation:rip_gaussian_uppdeltak_eq1}
n \geq 32\eta^{-2} \left( k \ln\left(\frac{ep}{k}\right) +  k \ln \left(\frac{12}{\eta}\right) + \ln(2/\varepsilon) \right). 
\end{equation}
Then, with probability at least $1 - \varepsilon$, the restricted isometry constant $\delta_k$ of the scaled matrix $\frac{1}{\sqrt{n}}\rmX$ satisfies
\begin{equation}\label{equation:rip_gaussian_uppdeltak_eq2}
\delta_k \leq \eta(2  +  \eta) . 
\end{equation}
\end{theoremHigh}
\begin{proof}[of Theorem~\ref{theorem:rip_gaussian_uppdeltak}]
Fix a support set $ \sS \subset \{1,2,\ldots,p\} $ with $ \abs{\sS}= k $. 
The submatrix  $ \rmX_\sS $ is an $ n \times k $ Gaussian matrix. 
The eigenvalues of $ \frac{1}{{n}} \rmX_\sS^\top \rmX_\sS - \bI $ lie in $ [\sigma_{\min}^2 - 1, \sigma_{\max}^2 - 1] $, where $ \sigma_{\min}$ and $ \sigma_{\max} $ are the smallest and largest singular values of $ \frac{1}{\sqrt{n}} \rmX_\sS $. 
Denote $ \widetilde{\rmX}_\sS \triangleq \frac{1}{\sqrt{n}} \rmX_\sS  $. 
For $\eta>0$, Theorem~\ref{thm:sv_gaus_mat} implies that
$$
\normtwo{\widetilde{\rmX}_\sS^\top \widetilde{\rmX}_\sS - \bI} 
\leq \max\left\{ (1 + \eta)^2 - 1, 1 - (1 - \eta)^2 \right\}
=2\eta + \eta^2,
$$
with probability at least $ 1 - 2 \left( \frac{12}{\eta} \right)^k \exp \left( -\frac{n \eta^2}{32} \right)   $. 
Taking the union bound over all $ \binom{p}{k} $ possible supports, by the definition of the restricted isometry constant (see Definition~\ref{definition:rip22}), $ \delta_k = \max_{\sS \subset \{1,2,\ldots,p\}, \abs{\sS} = k} \normtwo{\widetilde{\rmX}_\sS^\top \widetilde{\rmX}_\sS - \bI} $,
 we obtain
$$
\Pr\left( \delta_k > 2\eta + \eta^2 \right)
\leq 2 \binom{p}{k} \left( \frac{12}{\eta} \right)^k \exp \left( -\frac{n \eta^2}{32} \right) 
\leq 2 \left( \frac{ep}{k} \right)^k \left( \frac{12}{\eta} \right)^k \exp \left( -\frac{n \eta^2}{32} \right) ,
$$
where  the second inequality follows from the bound on the binomial coefficient (see Problem~\ref{prob:binom_ineq}). The last term is dominated by $ \varepsilon $ due to condition \eqref{equation:rip_gaussian_uppdeltak_eq1}. 
This completes the proof.
\end{proof}

Since $\delta_k \leq 2 \eta + \eta^2$, choosing  $\eta = 0.22$ yields $\delta_k \leq 0.4884< 0.4931$ (compare Theorem~\ref{theorem:stab_ell1_rip}  concerning $\ell_1$-minimization) under the condition $n \geq 662 \big( k \ln(ep/k) +  k \ln (12/\eta) + \ln(2/\varepsilon) \big) $ in this asymptotic regime.
\citet{foucart2013invitation} improves this bound to $n \geq 78.13 \left( k \ln(ep/k) + \ln(2/\varepsilon) \right)$ by improving the bound in Theorem~\ref{thm:sv_gaus_mat}; we omit the details here.

\begin{theoremHigh}[Exact recovery of \eqref{opt:p1} under Gaussian matrix]\label{theorem:rec_ell1_gauss}
Let $\rmX\in\real^{n \times p}$ be a random matrix with i.i.d. Gaussian entries. 
Suppose $k < p$ and $\varepsilon \in (0,1)$. 
If
\begin{equation*}
n \geq C_1k\ln(ep/k) +  C_2 k+ C_3\ln(2/\varepsilon)
\end{equation*}
for universal constants $C_1, C_2, C_3 > 0$. 
Then with probability at least $1 - \varepsilon$, every $k$-sparse vector $\bbeta$ is exactly recovered from the measurement $\rvy = \rmX\bbeta$ by solving the $\ell_1$-minimization \eqref{opt:p1}.
\end{theoremHigh}
\begin{proof}[of Theorem~\ref{theorem:rec_ell1_gauss}]
This follows by combining Theorem~\ref{theorem:rip_gaussian_uppdeltak} (which establishes RIP for Gaussian matrices) with Theorem~\ref{theorem:ectg_el1_rip} (which guarantees exact recovery under RIP). Note that recovery is unaffected by rescaling the measurement matrix, so the result holds for unnormalized Gaussian matrices as well.
\end{proof}

Once the RIP is established, recovery via $\ell_1$-minimization is not only exact for sparse signals, but also stable in the presence of sparsity defects and robust to additive noise in the measurements (see Theorem~\ref{theorem:error_p1epsilon} or Theorem~\ref{theorem:stab_ell1_rip}).

Additionally, these results are not limited to Gaussian matrices. They also hold for any random matrix ensemble satisfying the concentration inequality~\eqref{equation:concen_numsam_svd_eq1}, including matrices with independent isotropic sub-Gaussian rows (see Chapter~\ref{chapter:ensur_rips}).

\begin{problemset}
	
\item Show that any minimal $\varepsilon$-covering set must also be an $\varepsilon$-packing set.
Conversely, show that any maximal $2\varepsilon$-packing set must be an $\varepsilon$-covering set.
	
\item \label{prob:binom_ineq} \textbf{Exponential bounds on binomial coefficients.} 
For integers $p \geq k \geq 1$, prove that
$$
\left(\frac{p}{k}\right)^k \leq \binom{p}{k} \leq \left(\frac{ep}{k}\right)^k.
$$
\textit{Hint: Use the well-known Taylor expansion $\exp(k) = \sum_{i=0}^{\infty} k^i/(i!)$, find its lower bound, and  relate it to the expression for  $\binom{p}{k} $.
Alternatively, use \textit{Stirling's approximation}: $k! \geq (k/e)^k$.}

\item \textbf{Restricted eigenvalue property for Gaussian matrices \citep{raskutti2010restricted}.} 
Let  $\bX \in \real^{n \times p}$ be a random design matrix whose rows $\bx^{(i)} \in \real^p$ are drawn i.i.d. from $\normal(\bzero, \bSigma)$, where the covariance matrix $\bSigma$ is strictly positive definite. 
Show that the $\mu$-RE condition holds over the set $\sC[\sS; \gamma]$ with high probability, provided the sample size satisfies $n > C \abs{\sS} \log p$ for a sufficiently large universal constant $C$. 

\item \textbf{Mutual coherence for Gaussian matrices \citep{hastle2015statistical}.}
Let $\bX \in \real^{n \times p}$ be a random matrix with i.i.d. standard Gaussian entries.
Show that the {mutual incoherence condition (Definition~\ref{definition:mutual_incohe})} holds with high probability, provided the sample size satisfies $n > C k \log p$ for a sufficiently large numerical constant $C$. 
(It is known that the event $\mathcalE = \left\{ \sigma_{\min}\left( \frac{1}{\sqrt{n}}\bX_\sS^\top \right) \geq \frac{1}{2} \right\}$ holds with high probability; see Theorem~\ref{thm:sv_gaus_mat}.)
\begin{enumerate}[(i)]
\item 
Let $\bu\in\real^k$ and 
define $z_{i,\bu}\triangleq \bx_i^\top \bX_\sS (\bX_\sS^\top \bX_\sS)^{-1} \bu$. 
Show that 
$$
\gamma = 1 - \max_{i \in \comple{\sS}} \max_{\bu \in \{-1,+1\}^k} 
{z_{i,\bu}}.
$$

\item Using the event $\mathcalE$ defined above, show that there is a numerical constant $C_0 > 0$ such that
$$
\prob[z_{i,\bu} \geq \tau] 
\leq e^{-C_0 \frac{n \tau^2}{k}} + \prob[\comple{\mathcalE}],
 \quad \text{for any } \tau > 0,
$$
for each fixed index $i \in \comple{\sS}$ and vector $\bu \in \{-1,+1\}^k$.

\item Use part (ii) together with a union bound over all $i\in\comple{\sS}$ and $\bu \in \{-1,+1\}^k$ to complete the proof that the mutual incoherence condition holds.
\end{enumerate}

\item \label{prob:rd_colspace} \textbf{Randomized column-space approximation \citep{halko2011finding}.} 
Random projections enable efficient estimation of the column space of a low-rank matrix. Consider a data matrix $\bX \in \real^{n \times p}$ that is well-approximated by a rank-$r$ matrix, and examine the following algorithm:
\begin{enumerate}
\item Generate  a random matrix $\bA  \in\real^{p \times (r+k)}$ with i.i.d. standard Gaussian entries, where $k$ is a small oversampling parameter (e.g., 5).
\item Compute the sketch $\bB = \bX\bA \in \real^{n \times (r+k)}$.
\item Compute an orthonormal basis for the column space $\cspace(\bB)$, 
and  denote the resulting matrix by  $\widetildebU \in \real^{n \times (r+k)}$.
\end{enumerate}
Compare the column space of  $\widetildebU$ with that of $\bX$,
and contrast the computational cost of this randomized method with that of computing the exact column space of $\bX$ via a full SVD.
\textit{Hint:
See the core idea of JL lemma (Theorem~\ref{thm:jl_lemma}). 
Consider the matrix
\begin{align}
\bB &= \bX\bA 
= \bU \bSigma \bV^\top \bA
\triangleq \bU \bSigma \bC, 
\end{align}
where $\bX=\bU\bSigma\bV^\top$ is the reduced SVD of $\bX$, and $\bC\triangleq \bV^\top\bA\in\real^{r \times (r+k)}$.
If $\bX$ is exactly low rank,  $\bC$ is an i.i.d. standard Gaussian matrix since $\bV^\top \bV = \bI$ (see Lemma~\ref{lemma:affine_mult_gauss} on linear transformation of Gaussian vectors).
In that case the column space of $\bB$ is the same as that of $\bX$ because $\bC$ has full rank with high probability. When $\bX$ is only approximately low rank, then $\bC$ is a $\min\{n,p\} \times (r+k)$ i.i.d. standard Gaussian matrix. Surprisingly, for $p$ equal to a small integer, the product with $\bC$ conserves the subspace corresponding to the largest $r$ singular values with high probability, as long as the $(r+1)$-th singular values are sufficiently small. }

\item  \label{prob:rand_svd} \textbf{Randomized SVD.}
Given a matrix $\bX \in \real^{n \times p}$ that  is well-approximated by  a rank-$r$ matrix
\begin{enumerate}
\item Use the algorithm in Problem~\ref{prob:rd_colspace} to compute a semi-orthogonal matrix $\widetildebU \in \real^{n \times (r+k)}$ whose columns approximately span $\cspace(\bX)$.
\item Form the small matrix $\bY \in \real^{(r+k) \times p}$ defined by $\bY \triangleq \widetildebU^\top \bX$.
\item Compute the SVD of $\bY = \bU_\bY \bSigma_\bY \bV_\bY^\top$.
\item Output $\bU \triangleq (\widetildebU \bY)_{:,1:r}$, $\bSigma \triangleq (\bSigma_\bY)_{1:r,1:r}$ and $\bV \triangleq (\bV_\bY)_{:,1:r}$ as the SVD of $\bX$.
\end{enumerate}
Show that if $\bX \in \real^{n \times p}$ has exact rank $r$ and  
$\widetildebU$  spans its column space, then this algorithm recovers the exact SVD of $\bX$.
\textit{Note: Computing the SVD of $\bY$ costs  $\mathcalO(pr^2)$, so overall the complexity of the algorithm is governed by the second step which costs $\mathcalO(npr)$. 
This is dramatically faster than the classical SVD cost of  $\mathcalO(np \min\{n,p\})$, especially when $r\ll \min\{n,p\}$.}

\item \label{prob:rade_chaos}
\textbf{Homogeneous Rademacher chaos bound.}
Let $\bA \in \real^{p\times p}$ be a symmetric matrix with zero diagonal, and let $\bepsilon\in\{-1,+1\}^p$ be a Rademacher vector. 
Show that the homogeneous Rademacher chaos 
$
\rx = \bepsilon^\top \bA \bepsilon = \sum_{i \neq j} \epsilon_i \epsilon_j a_{ij}.
$
satisfies, 
\begin{align*}
\prob \left( \abs{\rx} \geq \tau \right)
&\leq 2 \exp \left( - \min \left\{ \frac{3\tau^2}{128 \normf{\bA}^2}, \frac{\tau}{32 \normtwo{\bA}} \right\} \right),
\quad\text{for $\tau > 0$}.
\end{align*}
\textit{Hint:
Use the moment generating function of $\rx = \bepsilon^\top \bA \bepsilon$.
Find a  decoupling inequality to bound $\Exp[\exp(\theta \rx)]$, then use convexity of $x \mapsto \exp(\theta x)$ and symmetry of $\bA$ to estimate the MGF of the quadratic form. 
}

\item Consider a $p \times k$ matrix $\bG$ whose entries are initialized as independent Rademacher random variables (i.e.,  $-1$ or $+1$ with equal probability)
and then scaled by $\sqrt{p}$. 
Discuss why the columns of $\bG$ will be (roughly) mutually orthogonal for large values of $p$ of the order of $10^6$. This trick is used frequently in machine learning for rapidly generating the random projection of an $n \times p$ data matrix $\bX$ as $\widetildebX = \bX\bG$.

\item \textbf{Stable NSP for Gaussian matrices.} 
Discuss under what conditions a random matrix $\rmX\in\real^{n\times p}$ with i.i.d. Gaussian entries satisfies the robust NSP. 
\textit{Hint: Use Theorem~\ref{theorem:rip_gaussian_uppdeltak} and Theorem~\ref{theorem:rip2rnsp}.}

\item \textbf{Near-orthogonality.} Two independent random vectors drawn uniformly from the unit sphere $\sU^{p-1}\subset\real^p$ are nearly orthogonal with high probability. Specifically, show that their inner product satisfies
$\abs{\innerproduct{\bu,\bv}}=\mathcalO(\sqrt{\frac{\ln(1/\varepsilon)}{p}})$ with probability at least $1 - \varepsilon$.

\item \label{prob:spe_err_ran_spe}
\textbf{Randomized low-rank approximation \citep{drineas2006fast}.}
Let $\bX\in\real^{n\times p}$, and let $\bC\in\real^{n\times m}$ contain $m$ columns of $\bX$ with $m<p$. 
Suppose further that the rank-$k$ truncated SVD of $\bC$ is $\bC\approx\bU_k\bSigma_k\bV_k^\top$ ($k\leq m$). 
Show that
\begin{subequations}\label{equation:spe_err_ran_spe}
\begin{align}
\normtwo{\bX - \bU_k \bU_k^\top \bX}^2 &\leq \normtwo{\bX - \bX_k}^2 + 2 \normtwo{\bX\bX^\top - \bC\bC^\top}; \label{equation:spe_err_ran_spe1}\\
\normf{\bX - \bU_k \bU_k^\top \bX}^2 &\leq \normf{\bX - \bX_k}^2 + 2 \sqrt{k} \normf{\bX \bX^\top - \bC \bC^\top}, \label{equation:spe_err_ran_spe2}
\end{align}
\end{subequations}
where $\bX_k$ denotes the best rank-$k$ approximation to $\bX$ obtained via its truncated SVD \citep{lu2022matrix}.

\end{problemset}

\newpage 
\chapter{Design of Sensing Matrices Using Sub-Gaussian  Ensembles}\label{chapter:ensur_rips}
\begingroup
\hypersetup{
linkcolor=structurecolor,
linktoc=page,  
}
\minitoc \newpage
\endgroup

\lettrine{\color{caligraphcolor}W}
While Gaussian random matrices offer strong theoretical guarantees for compressed sensing and sparse recovery, their practical use is often limited by computational and implementation constraints. In this chapter, we extend the analysis to sub-Gaussian ensembles---random matrices whose entries satisfy sub-Gaussian tail bounds---and demonstrate that many desirable properties of Gaussian matrices, such as the restricted isometry property and stable sparse recovery, continue to hold under much broader and more flexible distributions.

\section{Sub-Gaussian (SG) Random Variables}\label{section:subgauss_rdv}
Several common distributions, notably the Gaussian and Rademacher, are known to satisfy certain concentration-of-measure inequalities. 
We analyze this phenomenon from a more general perspective by considering the class of sub-Gaussian distributions \citep{buldygin2000metric, vershynin2010introduction}.

We begin by presenting the general definition of sub-Gaussian random variables and then specialize to the zero-mean case in the following section.
\begin{definition}[Sub-Gaussian random variables using concentration inequality\index{Sub-Gaussian variable}]\label{definition:subgau_concen_ineq}
A random variable $\rx$ is called \textit{sub-Gaussian} if there exist constants $\Phi, \xi > 0$ such that
\begin{equation}\label{equation:subgau_concen_ineq1}
\Pr(\abs{\rx} \geq \tau) \leq \Phi e^{-\xi \tau^2}, \quad \text{for all } \tau > 0.
\end{equation}
It is called \textit{sub-exponential} if
\begin{equation}
\Pr(\abs{\rx} \geq \tau) \leq \Phi e^{-\xi \tau}, \quad \text{for all } \tau > 0.
\end{equation}
\end{definition}

Clearly, a random variable $\rx$ is sub-Gaussian if and only if $\rx^2$ is sub-exponential. 
\footnote{Note again that we use normal fonts of boldface lowercase letters to denote random vectors, and normal
fonts of boldface uppercase letters to denote random matrices. That is, $\rx, \rva, \rmX$ are random scalars,
vectors, or matrices; while $x, \ba, \bX$ are scalars, vectors, or matrices. In many cases, the two terms can be
used interchangeably; that is, $\rx = x$ denotes a realization of the variable.
}
By Exercise~\ref{exercise:chernoff_gausvari} or Lemma~\ref{lemma:bd_cen_gaus},  a standard Gaussian random variable $\rx\sim \normal(0,1)$ is sub-Gaussian with $\Phi = 1$ and $\xi = 1/2$.
Moreover, Rademacher (i.e., symmetric Bernoulli) and bounded random variables are also sub-Gaussian.
According to {Theorem~\ref{theorem:hoeffding_ineq}} and Corollary~\ref{corollary:hoeffd_rademac}, sums of independent Rademacher variables are sub-Gaussian as well.

Lemmas~\ref{lemma:concent_res_fir} and \ref{lemma:concent_res_fir_inv} are converse to one another and establish a fundamental relationship between the tail probabilities and moments of a random variable---a connection essential for understanding the behavior of sub-Gaussian random variables. These lemmas build on the foundational notion of sub-Gaussianity by providing conditions under which concentration around the mean can be tightly controlled, thereby offering deeper insight into the probabilistic structure of such variables.

\begin{lemma}[Basic concentration lemma]\label{lemma:concent_res_fir}
Let $\rx$ be a random variable satisfying
$$
\norms{\rx}\triangleq
\big(\Exp[\abs{\rx}^s]\big)^{\frac{1}{s}} 
\leq \alpha \Phi^{\frac{1}{s}} s^{\frac{1}{\nu}},
\quad \text{for all } s_0 \leq s \leq s_1
$$
for some constants $\alpha, \Phi, \nu>0$ and $ s_1 > s_0 > 0$. Then
\begin{equation}
\Pr(\abs{\rx} \geq e^{\frac{1}{\nu}} \alpha \mu) 
\leq \Phi e^{-\frac{\mu^\nu}{\nu}},
\quad\text{for all $\mu \in [s_0^{1/\nu}, s_1^{1/\nu}]$}.
\end{equation}

\end{lemma}
\begin{proof}[of Lemma~\ref{lemma:concent_res_fir}]
By Markov's inequality (Theorem~\ref{theorem:markov-inequality}), for any $\kappa > 0$, we have 
$$
\Pr(\abs{\rx} \geq e^{\kappa} \alpha \mu) 
\leq \frac{\Exp[\abs{\rx}^s]}{(e^{\kappa} \alpha \mu)^s} 
\leq \Phi \left(\frac{\alpha s^{\frac{1}{\nu}}}{e^{\kappa} \alpha \mu}\right)^s.
$$
Choosing $s = \mu^\nu$ gives 
$\Pr(\abs{\rx} \geq e^{\kappa} \alpha \mu) 
\leq \Phi e^{-\kappa \mu^\nu}$, and  setting $\kappa = \frac{1}{\nu}$ yields the claim.
\end{proof}

Important special cases of Lemma~\ref{lemma:concent_res_fir} arise when $\nu = 1$  and $\nu= 2$. 
Indeed, 
if $(\Exp[\abs{\rx}^s])^{\frac{1}{s}} \leq \alpha \Phi^{\frac{1}{s}} s$ for all $s \geq 2$ then
\begin{equation}\label{equation:concent_res_fir_eq3}
\Pr(\abs{\rx} \geq e \alpha \mu) \leq \Phi e^{-\mu},
\quad \text{for all } \mu \geq 2.
\end{equation}
Similarly, if $(\Exp[\abs{\rx}^s])^{\frac{1}{s}} \leq \alpha \Phi^{\frac{1}{s}} \sqrt{s}$ for all $s \geq 2$, then
\begin{equation}\label{equation:concent_res_fir_eq2}
\Pr(\abs{\rx} \geq e^{1/2} \alpha \mu) \leq \Phi e^{-\mu^2/2},
\quad \text{for all } \mu \geq \sqrt{2}.
\end{equation}

The converse of Lemma~\ref{lemma:concent_res_fir} also holds.
\begin{lemma}[Basic moment bound lemma]\label{lemma:concent_res_fir_inv}
Suppose  a random variable $\rx$ satisfies, for some $\nu > 0$,
$$
\Pr(\abs{\rx} \geq e^{\frac{1}{\nu}} \alpha \mu) \leq \Phi e^{-\frac{\mu^\nu}{\nu}}, 
\quad \text{for all } \mu > 0.
$$
\begin{subequations}
Then, for all $s > 0$,
\begin{equation}\label{equation:concent_res_fir_inv_eq1}
\Exp[\abs{\rx}^s] \leq \Phi \alpha^s (e\nu)^{\frac{s}{\nu}} \Gamma\left(\frac{s}{\nu} + 1\right), 
\end{equation}
where $\Gamma(x)=\int_{0}^{\infty} t^{x-1}\exp(-t)dt$ is the \textit{Gamma function}.
As a consequence, for $s \geq 1$,
\begin{equation}\label{equation:concent_res_fir_inv_eq2}
\big(\Exp[\abs{\rx}^s]\big)^{\frac{1}{s}} \leq T_1 \alpha (T_{2,\nu} \Phi)^{\frac{1}{s}} s^{\frac{1}{\nu}},
\end{equation}
where $T_1 = e^{1/(2e)} \approx 1.2019$ and $T_{2,\nu} = \sqrt{\frac{2\pi}{\nu}} e^{\frac{\nu}{12}}$. 
\end{subequations}
\end{lemma}

\begin{proof}[of Lemma~\ref{lemma:concent_res_fir_inv}]
Using the absolute moments lemma (Lemma~\ref{lemma:mean_abs_mom}) and two changes of variables (i.e., $\tau\triangleq e^{\frac{1}{\nu}} \alpha \mu$ and $x=\frac{\mu^\nu}{\nu}$), we obtain
$$
\begin{aligned}
\footnotesize
\Exp[\abs{\rx}^s] 
&= s \int_0^\infty \Pr(\abs{\rx} > \tau) \tau^{s-1} d\tau
\stackrel{\tau\triangleq e^{\frac{1}{\nu}} \alpha \mu}{\longeq} s \alpha^s e^{\frac{s}{\nu}} \int_0^\infty 
\Pr(\abs{\rx}  \geq e^{\frac{1}{\nu}} \alpha \mu) \mu^{s-1} d\mu\\
&\leq s \alpha^s e^{\frac{s}{\nu}} \int_0^\infty \Phi e^{-\frac{\mu^\nu}{\nu}} \mu^{s-1} d\mu 
\stackrel{x=\frac{\mu^\nu}{\nu}}{\longeq} s \Phi \alpha^s e^{\frac{s}{\nu}} 
\int_0^\infty e^{-x} (x)^{\frac{s}{\nu} - 1} (\nu)^{\frac{s}{\nu} - 1} dx\\
&= \Phi \alpha^s (e\nu)^{\frac{s}{\nu}} \frac{s}{\nu} \Gamma\left(\frac{s}{\nu}\right),
\end{aligned}
$$
where the inequality follows from the assumption. 
and the last equality follows from the definition of Gamma function (see \eqref{equation:gamma_func}):
$
\Gamma(x) = \int_{0}^{\infty} t^{x-1} e^{-t} dt$ for  $x\geq 0$.
This proves \eqref{equation:concent_res_fir_inv_eq1} by the property of the Gamma function. 
Applying Stirling's formula---$\Gamma(x) = \sqrt{2\pi} x^{x-1/2} e^{-x} \exp\left(\frac{\theta(x)}{12x}\right)$ for positive $x$ and a correction term $0\leq \theta(x)\leq 1$---yields
$$
\Exp[\abs{\rx}^s] 
\leq \Phi \alpha^s (e\nu)^{\frac{s}{\nu}} \sqrt{2\pi} \left(\frac{s}{\nu}\right)^{\frac{s}{\nu} +\frac{ 1}{2}} e^{-\frac{s}{\nu}} e^{\frac{\nu}{12s}}
= \sqrt{2\pi} \Phi \alpha^s e^{\frac{\nu}{12s}} \left(s\right)^{\frac{s}{\nu} + \frac{ 1}{2}} \nu^{-\frac{ 1}{2}}.
$$
For $s \geq 1$, taking the $s$-th root gives
$$
\big(\Exp[\abs{\rx}^s]\big)^{\frac{1}{s}} 
\leq \left( \frac{\sqrt{2\pi}}{\sqrt{\nu} } e^{\frac{\nu}{12} } \Phi \right)^{\frac{1}{s}} \alpha s^{\frac{1}{2s}} s^{\frac{1}{\nu}}.
$$
Finally, the function $s^{{1}/{2s}}=\sqrt{s^{{1}/{s}}}$ attains its maximum at  $s = e$ (see Problem~\ref{prob:normineq}), so $s^{{1}/{2s}} \leq e^{{1}/{2e}}$. This completes the proof.
\end{proof}

\begin{remark}[Standard consequences for sub-Gaussian and sub-exponential variables]
Setting $\alpha \triangleq (2e\xi)^{-1/2}$ and $\nu \triangleq 2$, and applying Lemma~\ref{lemma:concent_res_fir_inv} with the change of variable $\tau = (2\xi) ^{-1/2}\mu$, we find that the moments of a sub-Gaussian variable $\rx$ satisfy
\begin{equation}\label{equation:moment_subgauss}
\big(\Exp[\abs{\rx}^s]\big)^{1/s} \leq \widetildeT\xi^{-1/2}\Phi^{1/s}s^{1/2}, 
\quad\forall\, s \geq 1,
\text{ with $\widetildeT \triangleq T_1T_{2,2}/\sqrt{2e}  \approx 1.0794$. }
\end{equation}
Similarly, setting $\alpha \triangleq (e\xi)^{-1}$ and $\nu \triangleq 1$, and using the change of variable $\tau=\xi^{-1}\mu$,
the moments of a sub-exponential variable $\rx$ satisfy
\begin{equation}\label{equation:moment_subexponen}
\big(\Exp[\abs{\rx}^s]\big)^{1/s} \leq \widehatT\xi^{-1}\Phi^{1/s}s, 
\quad \forall\, s \geq 1,
\text{ with $\widehatT \triangleq T_1T_{2,1}/e  \approx 1.2047$. }
\end{equation}
\end{remark}

We now present an equivalent characterization of  sub-Gaussian random variables.
\begin{theoremHigh}[Equivalent statement of sub-Gaussian]\label{theorem:stad_subgaus_txsqr}
Let $\rx$ be a random variable.
\begin{enumerate}[(i)]
\item If $\rx$ is sub-Gaussian, then there exist constants $c_1 > 0$ and $c_2 > 1$ such that $\Exp[\exp(c_1\rx^2)] \leq c_2$.
\item Conversely, if $\Exp[\exp(c_1\rx^2)] \leq c_2$ for some constants $c_1, c_2 > 0$, then $\rx$ is sub-Gaussian. More precisely, we have $\Pr(\abs{\rx} \geq \tau) \leq c_2\exp(-c_1 \tau^2)$ for all $\tau>0$.
\end{enumerate}
\end{theoremHigh}
\begin{proof}[of Theorem~\ref{theorem:stad_subgaus_txsqr}]
\textbf{(i).} 
Using the well-known Taylor expansion $\exp(z) = \sum_{i=0}^{\infty} \frac{z^i}{i!}$ for all $z\in\real$ around 0, we have 
$$
\Exp[\exp(c_1\rx^2)] 
= 1 + \sum_{i=1}^\infty \frac{c_1^i \Exp[\rx^{2i}]}{i!} 
$$
Applying the moment bound \eqref{equation:moment_subgauss} and Stirling's inequality---$i! \geq \sqrt{2\pi i} (i/e)^{i} \geq \sqrt{2\pi} (i/e)^{i}$ for $i \geq 1$---we obtain
$$
\Exp[\exp(c_1\rx^2)] 
= 1 + \sum_{i=1}^\infty \frac{c_1^i \Exp[\rx^{2i}]}{i!} 
\leq 
1 + \sum_{i=1}^\infty \frac{c_1^i \widetildeT^{2i}\xi^{-i}\Phi (2i)^i}{\sqrt{2\pi} (i/e)^{i}} 
=1 + \frac{\Phi}{\sqrt{2\pi}} \sum_{i=1}^\infty \left(\frac{c_1 \widetildeT^{2} (2e)}{\xi} \right)^i.
$$
If $\big({c_1 \widetildeT^{2} (2e)}/{\xi} \big) <1$, the series converges, and we may define $c_2$
accordingly: $c_2 \triangleq 1 + \frac{\Phi}{\sqrt{2\pi}} \sum_{i=1}^\infty \big({c_1 \widetildeT^{2} (2e)}/{\xi} \big)^i$.

\paragraph{(ii).} By Markov's inequality (Theorem~\ref{theorem:markov-inequality}),
$$
\Pr(\abs{\rx} \geq \tau) 
= \Pr\big(\exp(c_1\rx^2) \geq \exp(c_1 \tau^2)\big) 
\leq \Exp[\exp(c_1\rx^2)] e^{-c_1 \tau^2} 
\leq c_2 e^{-c_1 \tau^2}.
$$
This completes the proof. 
\end{proof}

We now consider the moment generating function of a zero-mean sub-Gaussian random variable.

\begin{theoremHigh}[Standard characterization of zero-mean sub-Gaussian variables]\label{theorem:stad_subgaus_tx}
Let $\rx$ be a random variable.
\begin{enumerate}[(i)]
\item If $\rx$ is sub-Gaussian with $\Exp[\rx] = 0$, then there exists a constant $C>0$ or a constant $c>0$ (depending only on $\Phi$ and $\xi$) such that
\begin{equation}\label{equation:stad_subgaus_tx}
\Exp[\exp(\theta \rx)] \leq \exp(\theta^2 \cdot C )=\exp(\theta^2 \cdot c^2  / 2), \quad \text{for all } \theta \in \real,
\end{equation}
where $C = c^2/2$.
\footnote{
Any constant $C$ satisfying \eqref{equation:stad_subgaus_tx} is called a \textit{sub-Gaussian parameter} of $\rx$. 
Typically, one chooses the smallest such $C$. 
The constant $c$ is closely related to the standard deviation; see Theorem~\ref{theorem:bk_subg_theo1}.
}
\item Conversely, if \eqref{equation:stad_subgaus_tx} holds for some $C>0$, 
then $\Exp[\rx] = 0$ and $\rx$ is sub-Gaussian with parameters $\Phi = 2$ and $\xi = \frac{1}{4C}$ in Definition~\ref{definition:subgau_concen_ineq}.
\end{enumerate}
\end{theoremHigh}

\begin{proof}[of Theorem~\ref{theorem:stad_subgaus_tx}]
\textbf{(i).} 
It suffices to consider $\theta \geq 0$, since the case $\theta < 0$ follows by replacing $\rx$ with $-\rx$, which preserves sub-Gaussianity and zero mean.
Using the well-known Taylor expansion $\exp(z) = \sum_{i=0}^{\infty} \frac{z^i}{i!}$ for all $z\in\real$ around 0, we obtain 
$$
\Exp[\exp(\theta \rx)] 
= 1 + \theta \Exp[\rx] + \sum_{i=2}^\infty \frac{\theta^i \Exp[ \rx^i]}{i!} 
= 1 + \sum_{i=2}^\infty \frac{\theta^i \Exp[\abs{\rx}^i]}{i!},
$$
where we used $\Exp[\rx]=0$. 
First, suppose $0 \leq \theta \leq \theta_0$, where  $\theta_0>0$ will be chosen later. 
Applying the moment bound  \eqref{equation:moment_subgauss} and Stirling's formula---$i! \geq \sqrt{2\pi i} (i/e)^{i}$ for $i \geq 1$---we get
\begin{align*}
\Exp[\exp(\theta \rx)] 
&\leq 1 + \Phi \sum_{i=2}^\infty \frac{\theta^i \widetildeT^i i^{i/2} \xi^{-i/2}}{i!} 
\leq 1 + \Phi \sum_{i=2}^\infty \frac{\theta^i \widetildeT^i i^{i/2} \xi^{-i/2}}{\sqrt{2\pi i} (i/e)^{i}}\\ 
&\leq 1 + \frac{\Phi (\widetildeT e \theta)^2}{\sqrt{2\pi }\xi} \sum_{i=0}^\infty \left( \widetildeT e \theta_0 \xi^{-1/2} \right)^i 
\leq 1 + \theta^2 \frac{\Phi (\widetildeT e)^2}{\sqrt{2\pi }\xi} \frac{1}{1 - \widetildeT e \theta_0 \xi^{-1/2}} \\ 
&= 1 + C_1 \theta^2 \leq \exp(C_1 \theta^2),
\end{align*}
where the third inequality holds since $i^{-(i+1)/2}<1$ for $i\geq 2$, and the penultimate inequality holds when $\widetildeT e \theta_0 \xi^{-1/2} < 1$ (using the inequality $\sum_{i=1}^{\infty} a^i\leq \frac{a}{1-a}$ for $0\leq a<1$). 
The last inequality is satisfied by setting
$$
\theta_0 \triangleq (2 \widetildeT e)^{-1} \sqrt{\xi},
$$
which gives the constant $C_1 = \sqrt{2} \Phi(\widetildeT e)^2 \xi^{-1} / \sqrt{\pi}$.

Now consider $\theta > \theta_0$. We aim to show $\Exp[\exp(\theta \rx - C_2 \theta^2)] \leq 1$ for a suitable $C_2>0$. 
Observe that
$$
\theta \rx - C_2 \theta^2 = -(\sqrt{C_2} \theta - \frac{\rx}{2\sqrt{C_2}})^2 + \frac{\rx^2}{4C_2} \leq \frac{\rx^2}{4C_2}.
$$
Invoking Theorem~\ref{theorem:stad_subgaus_txsqr} (i) with  $\widetilde{c}_1 > 0$ and $\widetilde{c}_2 \geq 1$, and choosing $C_2 = 1/(4\widetilde{c}_1)$, 
we have 
$$
\Exp[\exp(\theta \rx - C_2 \theta^2)] \leq \Exp[\exp(\widetilde{c}_1 \rx^2)] \leq \widetilde{c}_2.
$$
Define $\nu \triangleq \ln(\widetilde{c}_2) \theta_0^{-2}$. 
The above inequality and the condition $\theta>\theta_0$ show that 
\begin{align*}
\Exp[\exp(\theta \rx)] &\leq \widetilde{c}_2 \exp(C_2 \theta^2) = \widetilde{c}_2 \exp(-\nu \theta^2) \exp((C_2 + \nu) \theta^2) \\
&\leq \widetilde{c}_2 \exp(-\nu \theta_0^2) \exp((C_2 + \nu) \theta^2) 
\leq \exp((C_2 + \nu) \theta^2).
\end{align*}
Setting $C \triangleq \max\{C_1, C_2 + \nu\}$ completes the proof of part (i).

\paragraph{(ii).} 
Conversely, suppose \eqref{equation:stad_subgaus_tx} holds for some $C>0$. 
Let $\theta, \tau > 0$. 
By Markov's inequality,
$$
\Pr(\rx \geq \tau) = \Pr(\exp(\theta \rx) \geq \exp(\theta \tau)) 
\leq \frac{\Exp[\exp(\theta \rx)]  }{\exp(\theta \tau)}
\leq \exp(C \theta^2 - \theta \tau).
$$
Optimizing over $\theta$ gives the minimum at $\theta = \tau/(2C)$, yielding
$$
\Pr(\rx \geq \tau) \leq \exp(-\tau^2/(4C)).
$$
Applying the same argument to $-\rx$ gives 
$
\Pr(-\rx \geq \tau) \leq \exp(-\tau^2/(4C))$.
By the union bound, 
$$
\Pr(\abs{\rx} \geq \tau) \leq 2\exp(-\tau^2/(4C)),
$$ 
which shows that $\rx$ is sub-Gaussian with $\Phi=2$ and $\xi=1/(4C)$, as claimed.
To verify that $\Exp[\rx]=0$, note that for any real $\theta$,  $1 + \theta \rx \leq \exp(\theta \rx)$. 
Taking expectations and using \eqref{equation:stad_subgaus_tx}, we obtain for$\abs{\theta} < 1$,
$$ 
1 + \theta \Exp[\rx] \leq \Exp[\exp(\theta \rx)] \leq \exp(C \theta^2) \leq 1 + (C/2) \theta^2 + \mathcalO(\theta^4). 
$$
Letting $\theta \to 0$ forces  $\Exp[\rx] = 0$, completing the proof.
\end{proof}

\section{Strictly Sub-Gaussian (SSG) Random Variables}

As previously mentioned, when a sub-Gaussian random variable has zero mean, its moment generating function (MGF; see Definition~\ref{definition:momengenfunc}) satisfies a specific bound (see Theorem~\ref{theorem:stad_subgaus_tx}).
Moreover, key results concerning concentration inequalities and design properties---such as the restricted isometry property (RIP) for sub-Gaussian random matrices in Theorem~\ref{theorem:concen_numsam_rip}---assume that the underlying random variables are zero-mean sub-Gaussian.
We therefore adopt the MGF-based characterization to define this important subclass of sub-Gaussian random variables and explore how it enables fine-grained control over the behavior of design matrices.

\begin{definition}[(Zero-mean) sub-Gaussian random variables]\label{definition:sub_gaussian}
A random variable $\rx$ is called \textit{(zero-mean) sub-Gaussian (SG)}, denoted $\rx \sim \subnormal(c^2)$, if there exists a constant $c > 0$ such that~\footnote{More equivalent definitions are discussed in Problem~\ref{prob:sg_definition}.}
\begin{equation}\label{equation:subgau_mgf}
\Exp[\exp(\theta\rx )] \leq \exp(\theta^2 c^2  / 2), \quad \forall \, \theta \in \real.
\end{equation}
The quantity  $c^2$ is  referred to as the  \textit{variance proxy} of the sub-Gaussian distribution.
\end{definition}
This definition means that a random variable $\rx$ to is (zero-mean) sub-Gaussian if its Laplace transform (i.e., its MGF) is pointwise dominated by the MGF of a zero-mean Gaussian random variable with variance $c^2$.
Indeed, the function $\Exp[\exp(\theta\rx)]$ is the MGF of $\rx$, while the right-hand side of \eqref{equation:subgau_mgf} is precisely the MGF of a Gaussian random variable $\normal(0,c^2)$.
\footnote{Recall that the MGF of a Gaussian variable $\rx\sim\normal(\mu, \sigma^2)$ is $\Mom(\theta) = \exp(\theta\mu+\theta^2\sigma^2/2)$.}. 
Consequently, a (zero-mean) sub-Gaussian distribution is one whose tails decay at least as rapidly as those of a Gaussian distribution with variance $c^2$ (see also Theorem~\ref{theorem:bk_subg_theo1} for more insights).

Examples of such (zero-mean) sub-Gaussian distributions include the Gaussian distribution itself, the Rademacher distribution (symmetric Bernoulli variable) taking values $\pm 1$, and more generally, any zero-mean distribution with bounded support.
\begin{example}[Gaussian SG]
If $\rx \sim \normal(0, \sigma^2)$, i.e., $\rx$ is a zero-mean Gaussian random variable with variance $\sigma^2$, then $\rx \sim \subnormal(\sigma^2)$. Indeed, as mentioned above, the MGF of a Gaussian is given by $\Exp[\exp(\theta\rx)] = \exp( \theta^2 \sigma^2 / 2)$, and thus \eqref{equation:subgau_mgf} is trivially satisfied.
\end{example}

\begin{example}[Bounded SG and Rademacher]
Let $\rx$ be a zero-mean random variable satisfying $\abs{\rx} \leq B$ almost surely for some $B>0$. 
Then $\rx \sim \subnormal(B^2)$.
In particular, a Rademacher random variable (i.e., $\prob(\rx=\pm 1 )=1/2$) satisfies $\rx\sim\subnormal(1)$
\end{example}

A common way to characterize (zero-mean) sub-Gaussian random variables is through their moments. The following theorem, due to Buldygin and Kozachenko \citep{buldygin2000metric}, relates the first two moments to the variance proxy.

\begin{theoremHigh}[Buldygin--Kozachenko \citep{buldygin2000metric}]\label{theorem:bk_subg_theo1}
If $\rx \sim \subnormal(c^2)$, then,
\begin{equation}
\Exp[\rx] = 0 
\qquad \text{and}\qquad 
\Exp[\rx^2] \leq c^2.
\end{equation}
In other words, the variance of $\rx$ does not exceed its variance proxy $c^2$.
\end{theoremHigh}
\begin{proof}[of Theorem~\ref{theorem:bk_subg_theo1}]
Using the Taylor expansion $\exp(z) = \sum_{i=0}^{\infty} \frac{z^i}{i!}$ for the exponential function, for any $ \theta \in \real $,
$$
\sum_{k=0}^\infty \frac{\theta^k}{k!} \Exp[\rx^k] 
= \Exp[\exp(\theta \rx)] 
\leq \exp(c^2 \theta^2 / 2) 
= \sum_{k=0}^\infty \frac{c^{2 k }\theta^{2k}}{2^k k!}.
$$
Comparing terms up to second order as $\theta \rightarrow 0$, we obtain
$$
\Exp[\rx]\theta + \Exp[\rx^2]\frac{\theta^2}{2!} \leq \frac{c^2 \theta^2}{2} + o(\theta^2) \quad \text{as } \theta \to 0.
$$
Dividing both sides by $ \theta > 0 $ and letting $ \theta \to 0^+ $ yields $ \Exp[\rx] \leq 0 $. Dividing through by $ \theta < 0 $ and letting $ \theta \to 0^- $ gives $ \Exp[\rx] \geq 0 $. 
Hence, $ \Exp[\rx] = 0 $. (This was also shown in Theorem~\ref{theorem:stad_subgaus_tx} (ii).) 
With the mean established as zero, divide the inequality by $ \theta^2 $ and let $ \theta \to 0 $. 
This gives $ \Var[\rx] \leq c^2 $.
\end{proof}

Theorem~\ref{theorem:bk_subg_theo1} shows that if  $\rx \sim \subnormal(c^2)$, then $\Exp[\rx^2] \leq c^2$. 
In certain applications, it is useful to consider a more restrictive class of random variables for which this inequality becomes an equality.

\begin{definition}[Strictly sub-Gaussian random variables\index{Strictly sub-Gaussian variable}]\label{definition:strsub_gaussian}
A random variable $\rx$ is called \textit{strictly sub-Gaussian (SSG)} if it satisfies $\rx \sim \subnormal(\sigma^2)$ with $\sigma^2 = \Exp[\rx^2]$, denoted $\rx \sim \strsubnormal(\sigma^2)$, i.e., 
\begin{equation}
\Exp[\exp(\theta\rx)] \leq \exp( \theta^2 \sigma^2 / 2),  \text{ with }\sigma^2 = \Exp[\rx^2],  \quad \forall \, \theta \in \real.
\end{equation}
\end{definition}

\begin{example}[Strictly sub-Gaussian random variables]
The following are examples of SSG random variables:
\begin{itemize}
\item If $\rx \sim \normal(0, \sigma^2)$, then $\rx \sim \strsubnormal(\sigma^2)$.

\item If $\rx \sim \uniformdist(-1, 1)$, then $\Exp[\rx^2] = 1/3$, and one can verify that $\Exp[\exp(\theta \rx)] \leq \exp(\theta^2 / 6)$ for all $\theta \in \real$. Hence, $\rx \sim \strsubnormal(1/3)$.

\item Consider the symmetric three-point distribution defined by
\begin{equation}
\Pr(\rx=1) = \Pr(\rx=-1) = \frac{1-s}{2}, \quad \Pr(\rx=0) = s, \quad s \in [0,1).
\end{equation}
Its variance is $\Exp[\rx^2] = 1 - s$. This random variable is strictly sub-Gaussian if and only if $s \in [0, 2/3]$. For $s \in (2/3, 1)$, the tail behavior violates the strict sub-Gaussian condition, so $\rx$ is not SSG.
\end{itemize}
\end{example}

\citet{buldygin2000metric} provide an equivalent characterization of sub-Gaussian and strictly sub-Gaussian random variables in terms of their tail decay, which directly reflects their concentration-of-measure properties.
\begin{theoremHigh}[Buldygin--Kozachenko \citep{buldygin2000metric}]\label{theorem:bk_subg_theo}
A random variable $\rx$ is sub-Gaussian with variance proxy $c^2$, i.e., $\rx \sim \subnormal(c^2)$, if and only if there exist constants  $\tau_0 \geq 0$ and  $a \geq 0$ such that
\begin{equation}\label{equation:bk_subg_theo}
\Pr(\abs{\rx} \geq \tau) \leq 2 \exp\left(-\frac{\tau^2}{2a^2}\right), 
\quad \text{for all $\tau \geq \tau_0$}.
\end{equation}
Moreover, if $\rx \sim \strsubnormal(\sigma^2)$, then the above inequality holds    for all $\tau > 0$ with $a = \sigma$.
\end{theoremHigh}

The proof follows directly from Theorem~\ref{theorem:stad_subgaus_tx} (ii).
This matches the bound obtained by applying the Chernoff method directly to Gaussian variables (see Exercise~\ref{exercise:chernoff_gausvari}).
In fact, for a Gaussian random variable  $\rx\sim \normal(0, \sigma^2)$, Lemma~\ref{lemma:bd_cen_gaus} provides a slightly tighter tail bound:
\begin{equation}
\Pr(\abs{\rx} \geq \tau) \leq \exp\left(-\frac{\tau^2}{2\sigma^2}\right), 
\quad \text{for all $\tau \geq 0$}.
\end{equation}

Finally, sub-Gaussian distributions share a key stability property with Gaussian distributions: the sum of independent sub-Gaussian random variables is itself sub-Gaussian. 
This closure property is formalized in the following more general result.

\begin{proposition}\label{proposition:sgssg_vec}
Suppose that $\rvx = [\rx_1, \rx_2, \ldots, \rx_p]$, where each $\rx_i$ is i.i.d. with $\rx_i \sim \subnormal(c^2)$. 
Then,
\begin{equation}
\innerproduct{\rvx, \balpha} \sim \subnormal(c^2 \normtwo{\balpha}^2),
\quad \text{for any $\balpha \in \real^p$}.
\end{equation} 
Similarly, if each $\rx_i \sim \strsubnormal(\sigma^2)$, then  
\begin{equation}
\innerproduct{\rvx, \balpha} \sim \strsubnormal(\sigma^2 \normtwo{\balpha}^2), 
\quad\text{for any $\balpha \in \real^p$}.
\end{equation}
Moreover, it follows that 
$$
\Pr\left( \abs{\innerproduct{\rvx, \balpha}} \geq \tau \right) 
\leq 2 \exp\left(-\frac{\tau^2}{4C \normtwo{\balpha}^2}\right),
\quad \text{for all } \tau > 0.
$$
where $C={c^2}/{2}$ in the SG and $C={\sigma^2}/{2}$ in the SSG case.
\end{proposition}
\begin{proof}[of Proposition~\ref{proposition:sgssg_vec}]
Since the variables $\rx_i$ are i.i.d., the joint distribution factors and simplifies as:
\begin{align*}
\Exp\left[\exp\left(\theta \sum_{i=1}^p \alpha_i \rx_i\right)\right]
&= \Exp\left[\prod_{i=1}^p \exp(\theta \alpha_i \rx_i)\right]
= \prod_{i=1}^p \Exp\left[\exp(\theta \alpha_i \rx_i)\right] \\
&\leq \prod_{i=1}^p \exp(\theta^2\alpha_i^2 \cdot c^2   / 2)
= \exp\left(\left(\sum_{i=1}^p \alpha_i^2\right) \theta^2 c^2 / 2\right).
\end{align*}
Thus, $\innerproduct{\rvx, \balpha}$ is sub-Gaussian with variance proxy $c^2\normtwo{\balpha}^2$.
If the variables $\rx_i$ are strictly sub-Gaussian, the same argument applies with $c^2 = \sigma^2$, and one can verify that  $\Exp[\innerproduct{\rvx, \balpha }^2] = \sigma^2 \normtwo{\balpha}^2$, confirming strict sub-Gaussianity.
This proves the first part. 
The tail bound then follows directly from Theorem~\ref{theorem:stad_subgaus_tx} (ii).
\end{proof}

As a concrete example, let $\bepsilon = [\epsilon_1, \ldots, \epsilon_p]^\top$ be  a Rademacher sequence, and $\ry = \sum_{i=1}^p \alpha_i \epsilon_i$. 
Then
$$
\Exp [\exp(\theta \ry)] 
\leq \exp(\theta^2 \normtwo{\balpha}^2 / 2),
$$
confirming that $\ry$ is sub-Gaussian with variance proxy $\normtwo{\balpha}^2$.

The expected maximum of a finite collection of (possibly dependent) zero-mean sub-Gaussian random variables can be bounded as follows.
\begin{proposition}[Maximum of zero-mean sub-Gaussian variables]\label{proposition:max_zero_sg}
Let $\rvx = [\rx_1, \rx_2, \ldots, \rx_p]$, where each $\rx_i\sim \subnormal(c_i^2)$ (not necessarily independent)
satisfies $\Exp[\exp(\theta \rx_i)] \leq \exp(c_i^2 \theta^2/2)$, $i \in \{1,2,\ldots,p\}$. 
Define $C \triangleq \max_{i} c_i^2/2$. 
Then
\begin{equation}
\Exp [\max_{i \in\{1,2,\ldots,p\}} \rx_i] \leq \sqrt{4C \ln(p)},
\qquad\text{and}\qquad
\Exp [\max_{i \in \{1,2,\ldots,p\}} \abs{\rx_i}] \leq \sqrt{4C \ln(2p)}.
\end{equation}
And for any $\tau>0$, it holds that 
\begin{equation}
\prob\left(\max_{i \in\{1,2,\ldots,p\}} \rx_i \right) 
\leq p \exp(-\frac{\tau^2}{4C})
\quad\text{and}\quad
\prob\left(\max_{i \in\{1,2,\ldots,p\}} \abs{\rx_i} \right) 
\leq 2p \exp(-\frac{\tau^2}{4C}).
\end{equation}
\end{proposition}
\begin{proof}[of Proposition~\ref{proposition:max_zero_sg}]
The case $p = 1$ is trivial; assume $p \geq 2$. 
Let $\Theta > 0$ be a parameter to be chosen later. Using the monotonicity of the exponential and Jensen's inequality,
\begin{align*}
\Theta \Exp [\max_{i \in \{1,2,\ldots,p\}} \rx_i] 
&= \Exp [\ln \max_{i \in \{1,2,\ldots,p\}} \exp(\Theta \rx_i)] 
\leq \Exp \left[\ln \left( \sum_{i=1}^p \exp(\Theta \rx_i) \right)\right] \\
&\leq \ln \left( \sum_{i=1}^p \Exp [\exp(\Theta \rx_i)] \right) 
\leq \ln(p \exp(C \Theta^2)) = C \Theta^2 + \ln(p).
\end{align*}
Choosing $\Theta \triangleq \sqrt{C^{-1} \ln(p)}$ yields
$$
\sqrt{C^{-1} \ln(p)} \Exp [\max_{i \in \{1,2,\ldots,p\}} \rx_i] \leq \ln(p) + \ln(p),
$$
which implies $\Exp [\max_{i \in \{1,2,\ldots,p\}} \rx_i] \leq \sqrt{4C \ln(p)}$.

For the absolute maximum, we can rewrite the expectation by   $\Exp [\max_{i \in \{1,2,\ldots,p\}} \abs{\rx_i} ]
= \Exp [\max\{\rx_1, \ldots, \rx_p, -\rx_1, \ldots, -\rx_p\}]$. 
Applying the previous bound to this set of $2p$ variables gives the second expectation inequality.

The tail bounds follow from the union bound together with the sub-Gaussian tail inequality
\end{proof}

The bounds in this proposition are nearly sharp: for a sequence of independent standard Gaussian random variables, the expected maximum scales as $\sqrt{2\ln(p)}$, matching the above up to constant factors (see Theorem~\ref{theorem:norm_norm_vec} (iii)).

\section{Isotropic Sub-Gaussian  Random Vectors}

Sub-Gaussian distributions are closely connected to the concentration of measure phenomenon \citep{ledoux2001concentration}. For instance, combining Theorem~\ref{theorem:bk_subg_theo} and Proposition~\ref{proposition:sgssg_vec} yields deviation bounds for weighted sums of sub-Gaussian random variables.

More broadly, sub-Gaussian random matrices belong to a larger class of random matrices that exhibit strong concentration properties---specifically, they approximately preserve the Euclidean norm of vectors under linear transformation. This property is central to applications such as compressed sensing and high-dimensional statistics.

We begin by introducing key definitions.


\begin{definition}[Isotropic and sub-Gaussian random vectors\index{Isotropic vectors}]\label{definition:isotropic_subgvec}
Let $\rvx$ be a random vector in $\real^p$.
\begin{enumerate}[(i)]
\item If $\Exp[\innerproduct{\rvx,\bbeta}^2] = \normtwo{\bbeta}^2$ for all $\bbeta \in \real^p$, then $\rvx$ is called \textit{isotropic}.
\item If, for every unit vector $\bbeta \in \real^p$ (i.e., $\normtwo{\bbeta} = 1$), the scalar random variable $\innerproduct{\rvx,\bbeta}$ is sub-Gaussian with sub-Gaussian parameter $C\equiv c^2/2$ (see its Definition in Theorem~\ref{theorem:stad_subgaus_tx} or Definition~\ref{equation:subgau_mgf}) being independent of $\bbeta$ (and ideally independent of $p$); that is,
\begin{equation}\label{equation:subgaus_rdvec}
\Exp[\exp(\theta\innerproduct{\rvx,\bbeta})] 
\leq \exp(\theta^2 \cdot C), 
\quad \text{for all } \theta \in \real, \quad \normtwo{\bbeta} = 1,
\end{equation}
then $\rvx$ is called a \textit{sub-Gaussian random vector}.
\end{enumerate}
\end{definition}

Note that isotropic sub-Gaussian random vectors need not have independent entries.

The following proposition shows that a random vector with independent, mean-zero, unit-variance sub-Gaussian entries is both isotropic and sub-Gaussian; see Problem~\ref{prob:subg_vec} for further discussions.
\begin{proposition}[Isotropic and sub-Gaussian random vector]\label{proposition:iso_subg_vec}
Let $\rvx \in \real^p$ be a random vector whose entries are independent, mean-zero, sub-Gaussian random variables with variance $1$ nd a common sub-Gaussian parameter $C\equiv c^2/2$ as in \eqref{equation:subgaus_rdvec}. 
5Then $\rvx$ is an isotropic and sub-Gaussian random vector, with constants independent of the dimension $p$.
\end{proposition}

\begin{proof}[of Proposition~\ref{proposition:iso_subg_vec}]
Let $\bbeta \in \real^p$ with $\normtwo{\bbeta} = 1$. 
Since the the entries  $\rx_i$ of $\rvx$ are independent, zero-mean, and  have unit variance, 
\begin{equation*}
\Exp[\innerproduct{\rvx,\bbeta}^2] = \sum_{i,j=1}^{p} \beta_i\beta_{j} \Exp[\rx_i \rx_{j}] = \sum_{i=1}^{p}\beta_i^2 = \normtwo{\bbeta}^2.
\end{equation*}
Therefore, $\rvx$ is isotropic. Furthermore, according to Proposition~\ref{proposition:sgssg_vec} and Theorem~\ref{theorem:stad_subgaus_tx}, the random variable $\ry = \innerproduct{\rvx,\bbeta}$ is sub-Gaussian $\subnormal(c^2/2)$ (or with parameters $\Phi = 2$ and $\xi = 1/(4C\normtwo{\bbeta}^2) = 1/(4C)$). Hence, $\rvx$ is a sub-Gaussian random vector with parameters independent of $p$.
\end{proof}

Now consider a matrix $\rmX \in \real^{n \times p}$ whose entries are (unit-variance sub-Gaussian) random variables. 
 Such a matrix is referred to as a \textit{(sub-Gaussian) random matrix} or \textit{(sub-Gaussian) random matrix ensemble}.
More specifically, we focus on random matrices $\rmX \in \real^{n \times p}$ whose rows are independent, isotropic, and sub-Gaussian random vectors in $\real^p$:
\begin{equation*}
\rmX = \begin{bmatrix}
\rvx_1^\top \\
\rvx_2^\top \\
\vdots \\
\rvx_n^\top
\end{bmatrix}.
\end{equation*}
where each $\rvx_i\in\real^p $ satisfies the conditions in Definition~\ref{definition:subgau_concen_ineq}.

\begin{example}[Sub-Gaussian random matrices]\label{example:subgaus_rdmat}
Let $\rmX$ be an $n \times p$ random matrix.
\begin{enumerate}[(i)]
\item If the entries of $\rmX$ are independent Rademacher variables (i.e., taking values $\pm 1$ with equal probability), then $\rmX$ is called a \textit{Rademacher random matrix}.
\item If the entries of $\rmX$ are independent standard normal distributed random variables, then $\rmX$ is called a \textit{Gaussian random matrix}.
\item More generally, if all entries of $\rmX$ are independent mean-zero sub-Gaussian random variables of variance $1$ with the same constants $\Phi, \xi$ in the Definition~\ref{definition:subgau_concen_ineq} of sub-Gaussian random variables, that is,
\begin{equation}\label{equation:subgaus_rdmat_3}
\Pr(\abs{\rx_{ij}} \geq \tau) \leq \Phi e^{-\xi \tau^2}, \quad \text{for all } \tau > 0, \quad i \in \{1,2,\ldots,n\}, j \in \{1,2,\ldots,p\},
\end{equation}
then $\rmX$ is called a \textit{sub-Gaussian random matrix}.
Equivalently, one may require the moment-generating function bound
\begin{equation}
\Exp[\exp(\theta \rx_{ij})] 
\leq \exp(\theta^2 C), \quad \text{for all } \theta \in \real, \quad i \in \{1,2,\ldots,n\}, j \in \{1,2,\ldots,p\},
\end{equation}
for some constant $C\equiv c^2/2$ that is independent of $i,j$, and $p$ (see {Theorem~\ref{theorem:stad_subgaus_tx}}).
\end{enumerate}
Both Gaussian and Rademacher matrices are special cases of sub-Gaussian random matrices. Note that the entries of a sub-Gaussian matrix need not be identically distributed---only that they share uniform sub-Gaussian parameters.
\end{example}

A crucial property of random matrices with sub-Gaussian rows is their norm-preserving behavior, which underpins the restricted isometry property (RIP), which will be discussed in the sequel.
\begin{theoremHigh}[Norm preservation   under sub-Gaussian rows\index{Norm preservation}]\label{theorem:concen_under_sgrow}
Let $\rmX\in\real^{n \times p}$ be a random matrix with independent, isotropic, and sub-Gaussian rows with the same sub-Gaussian parameter $C\equiv c^2/2$ in \eqref{equation:subgaus_rdvec}. Then, 
$$
\Pr\left(\abs{n^{-1}\normtwo{\rmX\bbeta}^2 - \normtwo{\bbeta}^2} \geq \varepsilon\normtwo{\bbeta}^2\right) \leq 2\exp(-\widetildeC\varepsilon^2n),
\quad \text{for all $\bbeta \in \real^p$, all $\varepsilon \in (0,1)$},
$$
where $\widetildeC$ depends only on $C$.
\end{theoremHigh}
\begin{proof}[of Theorem~\ref{theorem:concen_under_sgrow}]
Let $\bbeta \in \real^p$. Without loss of generality, assume that $\normtwo{\bbeta} = 1$. 
Let $\rvx^{(1)}, \rvx^{(2)}, \ldots, \rvx^{(n)} \in \real^p$ denote the rows of $\rmX$,  and define the random variables
\begin{equation*}
\ry_i = \innerproduct{\rvx^{(i)},\bbeta}^2 - \normtwo{\bbeta}^2=\innerproduct{\rvx^{(i)},\bbeta}^2-1,
\quad i \in \{1,2,\ldots,n\}.
\end{equation*}
Since $\rvx^{(i)}$ is isotropic, we have $\Exp [\ry_i ]= 0$. 
Moreover, because $\innerproduct{\rvx^{(i)},\bbeta}$ is sub-Gaussian, $\ry_i$ is sub-exponential;
that is, there exist constants $\Phi, \xi$ (depending only on $C$) such that $\Pr(\abs{\ry_i} \geq \varepsilon) \leq \Phi \exp(-\xi \varepsilon)$. 
Observe that
\begin{equation*}
n^{-1}\normtwo{\rmX\bbeta}^2 - \normtwo{\bbeta}^2 
= \frac{1}{n}\sum_{i=1}^{n} \big(\innerproductbig{\rvx^{(i)},\bbeta}^2 - \normtwo{\bbeta}^2\big) = \frac{1}{n}\sum_{i=1}^{n} \ry_i.
\end{equation*}
Since the $\ry_i$ are independent, Bernstein’s inequality for sub-exponential random variables (Corollary~\ref{corollary:berns_subexp}) yields
\begin{align*}
\Pr\left(\abs{n^{-1}\sum_{i=1}^{n} \ry_i} \geq \varepsilon\right) 
&= \Pr\left(\abs{\sum_{i=1}^{n} \ry_i} \geq \varepsilon n\right) \\
&\leq 2\exp\left(-\frac{\xi^2\varepsilon^2n^2/2}{2\Phi n + \xi \varepsilon n}\right) 
\leq 2\exp\left(-\frac{\xi^2\varepsilon^2n}{4\Phi + 2\xi}\right),
\end{align*}
where the last inequality follows since  $\varepsilon \in (0,1)$. 
Therefore, the claim follows with $\widetildeC \triangleq \frac{\xi^2}{4\Phi + 2\xi}$.
\end{proof}

Finally, note that the scaled matrix $\widetilde{\rmX} \triangleq \frac{1}{\sqrt{n}}\rmX$, with $\rmX$ satisfying the assumptions of the previous result, satisfies
\begin{equation}\label{equation:norm_subg_ineq}
\Pr\left(\abs{\normtwo{\widetilde{\rmX}\bbeta}^2 - \normtwo{\bbeta}^2} 
\geq \varepsilon\normtwo{\bbeta}^2\right) \leq 2\exp(-\widetildeC\varepsilon^2n),
\quad \text{for all $\bbeta \in \real^p$, all $\varepsilon \in (0,1)$}.
\end{equation}
This concentration inequality serves as the foundation for proving the restricted isometry property of sub-Gaussian random matrices; 
see Proposition~\ref{proposition:concen_numsam_svd}.

\section{Concentration of Measure for SSG Random Variables}

In this section, we show that when the entries are strictly sub-Gaussian, the constant $\widetildeC$ in \eqref{equation:norm_subg_ineq} can be made explicit.
Moreover, for our purposes, it is particularly insightful to study the Euclidean norm of a vector composed of sub-Gaussian random variables. 
Specifically, suppose $\rvx=[\rx_1, \rx_2, \ldots, \rx_n]^\top$, where each $\rx_i$ is i.i.d. with $\rx_i \sim \subnormal(c)$. 
We are interested in how $\normtwo{\rvx}$ deviates from its expected value.

To establish the main result, we will use Markov's inequality for nonnegative random variables (Theorem~\ref{theorem:markov-inequality}) together with the following bound on the exponential moment of a sub-Gaussian random variable.

\begin{lemma}\label{lemma:subgau_bd}
Suppose $\rx \sim \subnormal(c^2)$. Then, for any $\lambda \in [0, 1)$,
$$
\Exp[\exp\left(\lambda \rx^2 / 2c^2\right)] \leq \frac{1}{\sqrt{1 - \lambda}},
$$
\end{lemma}
\begin{proof}[of Lemma~\ref{lemma:subgau_bd}]
The claim is trivial when $\lambda = 0$. 
Assume now that $\lambda \in (0, 1)$. 
Let $f(x)$ denote the probability density function of $\rx$. By the definition of sub-Gaussianity,
$$
\int_{-\infty}^{\infty} \exp(\theta x) f(x) \, dx 
\leq \exp(\theta^2 c^2 / 2), 
\quad \forall \, \theta \in \real.
$$
Multiplying both sides by $\exp(-\theta^2 c^2 / 2\lambda)$ gives
$$
\int_{-\infty}^{\infty} \exp\left(\theta x - \theta^2 c^2 / 2\lambda\right) f(x) \, dx 
\leq \exp(\theta^2 c^2 (\lambda - 1) / 2\lambda).
$$
Integrating both sides with respect to $\theta\in\real$ yields
$$
\int_{-\infty}^{\infty} \left(\int_{-\infty}^{\infty} 
\exp\left(\theta x - \theta^2 c^2 / 2\lambda\right) d\theta\right) f(x) \, dx
\leq \int_{-\infty}^{\infty} \exp(\theta^2 c^2 (\lambda - 1) / 2\lambda) \, d\theta.
$$
Evaluating the Gaussian integrals on both sides leads to
$$
\frac{1}{c} \sqrt{{2\pi\lambda}} \int_{-\infty}^{\infty} \exp\left(\lambda x^2 / 2c^2\right) f(x) \, dx \leq \frac{1}{c} \sqrt{\frac{2\pi\lambda}{1 - \lambda}},
$$
which simplifies to the desired inequality. This completes the proof.
\end{proof}

This lemma refines the properties of zero-mean sub-Gaussian variables and complements part (i) of Theorem~\ref{theorem:stad_subgaus_txsqr}.

We now present a key concentration result, which generalizes findings from \citet{dasgupta2003elementary} and extends the well-known behavior of Gaussian vectors (see Section~\ref{section:conce_gausvec}).
\begin{theoremHigh}[Concentration of SG variables \citep{dasgupta2003elementary, baraniuk2011introduction}]\label{theorem:concen_subgaus}
Let $\rvx = [\rx_1, \rx_2, \ldots, \rx_n]^\top$, where each $\rx_i$ is i.i.d. with $\rx_i \sim \subnormal(c^2)$ and $\Exp[\rx_i^2] = \sigma^2$. Then,
\begin{equation}\label{equation:concen_subgaus_res1}
\Exp[\normtwo{\rvx}^2] = n \sigma^2.
\end{equation}
Moreover, for any $\mu \in (0, 1)$ and for any $\nu \in [c^2/\sigma^2, \nu_{\max}]$, there exists a constant $\kappa \geq 4$ depending only on $\nu_{\max}$ and the ratio $\sigma^2 / c^2$ such that
\begin{equation}\label{equation:concen_subgaus_res2}
\Pr(\normtwo{\rvx}^2 \leq \mu n \sigma^2) \leq \exp\left(-n(1 - \mu)^2 / \kappa\right)
\end{equation}
and
\begin{equation}\label{equation:concen_subgaus_res3}
\Pr(\normtwo{\rvx}^2 \geq \nu n \sigma^2) \leq \exp\left(-n(\nu - 1)^2 / \kappa\right).
\end{equation}
\end{theoremHigh}

\begin{proof}[of Theorem~\ref{theorem:concen_subgaus}]
Equation \eqref{equation:concen_subgaus_res1} follows immediately from independence:
$$
\Exp[\normtwo{\rvx}^2] = \sum_{i=1}^{n} \Exp[\rx_i^2] = \sum_{i=1}^{n} \sigma^2 = n \sigma^2
$$
Since $\normtwo{\rvx}^2$ is nonnegative and the exponential function is monotone, we can apply Markov's inequality (Theorem~\ref{theorem:markov-inequality}):
\begin{align*}
\Pr(\normtwo{\rvx}^2 \geq \nu n \sigma^2) 
&= \Pr\left(\exp\big(\lambda \normtwo{\rvx}^2\big) \geq \exp\big(\lambda \nu n \sigma^2\big)\right) \\
&\leq \frac{\Exp\left[\exp\left(\lambda \normtwo{\rvx}^2\right)\right]}{\exp\left(\lambda \nu n \sigma^2\right)} 
= \frac{\prod_{i=1}^{n} \Exp\left[\exp\left(\lambda \rx_i^2\right)\right]}{\exp\left(\lambda \nu n \sigma^2\right)}.
\end{align*}
Since $\rx_i \sim \subnormal(c^2)$, for $\lambda\in[0,1/(2c^2))$, we have from Lemma~\ref{lemma:subgau_bd} that
\begin{align*}
&\Exp\left[\exp\left(\lambda \rx_i^2\right)\right] 
= \Exp\left[\exp\left(2c^2 \lambda \rx_i^2 / 2c^2\right)\right] \leq \frac{1}{\sqrt{1 - 2c^2 \lambda}}\\
\implies 
&\prod_{i=1}^{n} \Exp\left(\exp\left(\lambda \rx_i^2\right)\right) \leq \left(\frac{1}{1 - 2c^2 \lambda}\right)^{n/2}.
\end{align*}
Combining the preceding two results shows that 
$$
\Pr(\normtwo{\rvx}^2 \geq \nu n \sigma^2) \leq \left(\frac{\exp\left(-2\lambda \nu \sigma^2\right)}{1 - 2c^2 \lambda}\right)^{n/2}.
$$
Optimizing over $\lambda$ yields the choice
$$
\lambda = \frac{\nu \sigma^2 - c^2}{2c^2 \sigma^2 \nu},
$$
which lies in $[0,1/(2c^2)]$ provided  $\nu \in [c^2/\sigma^2, \nu_{\max}]$.
Substituting this optimal $\lambda$ gives
\begin{equation}\label{equ:concen_subgaus_e11}
\Pr(\normtwo{\rvx}^2 \geq \nu n \sigma^2) \leq \left(\nu \frac{\sigma^2}{c^2} \exp\left(1 - \nu \frac{\sigma^2}{c^2}\right)\right)^{n/2}.
\end{equation}
A similar calculation yields
\begin{equation}\label{equ:concen_subgaus_e1}
\Pr(\normtwo{\rvx}^2 \leq \mu n \sigma^2) \leq \left(\mu \frac{\sigma^2}{c^2} \exp\left(1 - \mu \frac{\sigma^2}{c^2}\right)\right)^{n/2}.
\end{equation}

To simply these bounds, define
$$
\kappa \triangleq
\max\left\{
4, 
2 \frac{(\nu_{\max} \sigma^2 / c^2 - 1)^2}{(\nu_{\max} \sigma^2 / c^2 - 1) - \ln(\nu_{\max} \sigma^2 / c^2)}
\right\}.
$$
Then, for all $\gamma \in [0, \nu_{\max} \sigma^2 / c^2]$, one has
$$
\ln(\gamma) \leq (\gamma - 1) - \frac{2(\gamma - 1)^2}{\kappa}
\quad\implies\quad 
\gamma \leq \exp\left((\gamma - 1) - \frac{2(\gamma - 1)^2}{\kappa}\right).
$$
Since $\sigma^2\leq c^2$ by Theorem~\ref{theorem:bk_subg_theo1},
setting $\gamma = \nu \sigma^2 / c^2\leq \nu\leq \nu_{\max}$  in \eqref{equ:concen_subgaus_e11} gives \eqref{equation:concen_subgaus_res3}.
Similarly, setting $\gamma = \mu \sigma^2 / c^2$  in \eqref{equ:concen_subgaus_e1} yields \eqref{equation:concen_subgaus_res2}. 
\end{proof}

This result shows that the squared norm of a sub-Gaussian vector concentrates sharply around its mean $n \sigma^2$, with exponentially decaying tail probabilities as $n$ grows. However, note that the upper-tail bound \eqref{equation:concen_subgaus_res3} only applies for $\nu \geq c^2 / \sigma^2 \geq 1$, as implied by Theorem~\ref{theorem:bk_subg_theo1}. Thus, for a general sub-Gaussian distribution, the concentration may not be arbitrarily tight.

Crucially, for strictly sub-Gaussian distributions---where $c^2 = \sigma^2$---this restriction disappears. In this favorable case, we obtain the following simplified and symmetric concentration inequality.

\begin{theoremHigh}[Concentration of SSG variables]\label{theorem:concen_SSG}
Let $\rvx = [\rx_1, \rx_2, \ldots, \rx_n]^\top$, where each $\rx_i$ is i.i.d. with $\rx_i \sim \strsubnormal(\sigma^2)$. Then,
\begin{equation}
\Exp[\normtwo{\rvx}^2] = n \sigma^2
\end{equation}
and for any $\varepsilon\in(0,1)$,
\begin{equation}\label{equation:concen_SSG_res2}
\Pr\left(\abs{\normtwo{\rvx}^2 - n \sigma^2} \geq \varepsilon n \sigma^2\right) \leq 2 \exp\left(-\frac{n \varepsilon^2}{\kappa}\right)
\end{equation}
with $\kappa = 2 / (1 - \ln(2)) \approx 6.52$.
\end{theoremHigh}
\begin{proof}[of Theorem~\ref{theorem:concen_SSG}]
Since each $\rx_i \sim \strsubnormal(\sigma^2)$, we have that $\rx_i \sim \subnormal(\sigma^2)$ and $\Exp[\rx_i^2] = \sigma^2$. 
Invoking Theorem~\ref{theorem:concen_subgaus} with $\mu = 1 - \varepsilon$ and $\nu = 1 + \varepsilon$  allows us to simplify and combine
the bounds in \eqref{equation:concen_subgaus_res2} and \eqref{equation:concen_subgaus_res3} to obtain \eqref{equation:concen_SSG_res2}. 
Since $\nu<2$ for $\varepsilon\in(0,1)$, we may take $\nu_{\max} = 2$.
Thus, $\kappa = 2 / (1 - \ln(2)) \approx 6.52$.
This completes the proof.
\end{proof}

As an immediate application, consider a noise vector $\bepsilon \in \real^n$ with i.i.d. entries $\epsilon_i\sim\normal(0, \sigma^2)$. Since the Gaussian distribution is strictly sub-Gaussian, Theorem~\ref{theorem:concen_SSG} implies that there exists a constant $c_0>0$ such that for any $e>0$,
\begin{equation}
\Pr\left(\normtwo{\bepsilon} \geq (1 + e) \sqrt{n} \sigma\right) \leq \exp\left(-c_0 ne^2 \right).
\end{equation}
In particular, setting $e=1$ gives
\begin{equation}
\Pr\left(\normtwo{\bepsilon} \geq 2 \sqrt{n} \sigma\right) \leq \exp\left(-c_0 n \right).
\end{equation}

Finally, Theorem~\ref{theorem:concen_SSG} also implies the following useful corollary, which shows that the transformation of a vector by a random matrix exhibits a concentration phenomenon.
The result is a special case of Theorem~\ref{theorem:concen_under_sgrow} where the constant is determined.

\begin{corollary}[Norm preservation for SSG matrices]\label{corollary:concen_SSG_vec}
Let  $\rmX\in\real^{n \times p}$ be a random matrix with i.i.d. entries $\rx_{ij} \sim \strsubnormal(1/n)$. 
For any fixed $\bbeta \in \real^p$, define $\rvy = \rmX \bbeta$. 
Then
\begin{equation}\label{equation:concen_SSG_vec_res1}
\Exp[\normtwo{\rvy}^2] = \normtwo{\bbeta}^2
\end{equation}
and for any $\varepsilon\in(0,1)$,
\begin{align}
\Pr\left(\abs{\normtwo{\rvy}^2 - \normtwo{\bbeta}^2} \geq \varepsilon \normtwo{\bbeta}^2\right) \leq 2 \exp\left(-\frac{n \varepsilon^2}{\kappa}\right) \label{equation:concen_SSG_vec_res2}\\
\implies 
\Pr\left( (1-\varepsilon) \normtwo{\bbeta}^2 \leq \normtwo{\rvy}^2 \leq  (1+\varepsilon) \normtwo{\bbeta}^2\right) 
\geq 1- 2 \exp\left(-\frac{n \varepsilon^2}{\kappa}\right), \label{equation:concen_SSG_vec_res3}
\end{align}
where $\kappa = 2 / (1 - \ln(2)) \approx 6.52$.
\end{corollary}
\begin{proof}[of Corollary~\ref{corollary:concen_SSG_vec}]
Let $\bx^{(i)}$ denote the $i$-th row of $\rmX$. Then $\ry_i = \innerproduct{\bx^{(i)}, \bbeta}$, and thus by Proposition~\ref{proposition:sgssg_vec}, $\ry_i \sim \strsubnormal(\normtwo{\bbeta}^2 / n)$. Applying Theorem~\ref{theorem:concen_SSG} to the $n$-dimensional random vector $\rvy$, we obtain \eqref{equation:concen_SSG_vec_res2}.
\end{proof}

\section{Design Properties under Sub-Gaussian Random Matrices}\label{section:desi_ssgsg}

We now return to again the question of how to construct matrices that satisfy the RIP. 
While it is possible to deterministically construct $n \times p$ matrices that satisfy the RIP of order $k$, early approaches relied on mutual coherence.
A matrix $\bX \in \real^{n \times p}$ with unit-norm columns is said to be $\mu$-coherent if  $\innerproduct{\bx_i, \bx_j} \leq \mu$ for all $i \neq j \in \{1,2,\ldots,p\}$ (see Definition~\ref{definition:mutua_cohere}). 
Such a  matrix  satisfies the RIP of order $k$ with constant $\delta = (k-1)\mu$ (see Proposition~\ref{proposition:cohe2rip}).

Deterministic constructions of $\mu$-coherent matrices with $\mu = \mathcalO\left(\frac{\ln p}{\sqrt{n \ln n}}\right)$ have been known since the work of \citet{kashin1975diameters}.
However, these constructions require $n = \widetilde{O}\left(\frac{k^2 \ln^3 p}{\delta^2}\right)$, which is quadratically larger than the number of rows needed by random constructions.
Other explicit constructions fare little better:  
\citet{devore2007deterministic} achieves $n = \mathcalO(k^2 \ln p)$, while \citet{indyk2008explicit} requires $n = \mathcalO(k p^\mu)$ for some constant $\mu>0$.
In many practical settings, such requirements on $n$ are prohibitively large.
The first significant improvement came from \citet{bourgain2011explicit}, who introduced deterministic combinatorial constructions satisfying the RIP with
$n = \widetilde{O}\left(\frac{k^{(2-\epsilon)}}{\delta^2}\right)$ for some constant $\epsilon > 0$. 
Although this reduces the quadratic dependence on $k$, the constructions remain complex and still fall short of the near-linear scaling achieved by random matrices.

Fortunately, these limitations can be overcome by randomizing the matrix construction.
A simple and effective approach is as follows: given dimensions $n$ and $p$, generate a random matrix $\rmX\in\real^{n\times p}$ by sampling its entries $\rx_{ij}$ independently from a suitable probability distribution.
One might first observe that if we only require  $\delta_{2k} > 0$  (i.e., mere injectivity on $2k$-sparse vectors), then setting $n=2k$ and drawing $\rmX$ from a continuous distribution---such as the Gaussian---suffices. With probability 1, every subset of $2k$ columns is linearly independent, so the RIP holds for some $\delta_{2k}\in(0,1)$.

However, in practice we typically need a specified RIP constant $\delta_{2k}$
(e.g., $\delta_{2k}<0.1$). To guarantee such a bound, we must control the singular values of all $\binom{p}{k}$ submatrices corresponding to $k$-sparse supports---a task that is computationally infeasible for realistic values of $p$ and $k$. Moreover, when $n=2k$, the resulting $\delta_{2k}$ is typically very close to 1, rendering the matrix unsuitable for applications like compressed sensing.

Therefore, we focus instead on constructing matrices that satisfy the RIP of order $2k$ with a prescribed constant $\delta_{2k}\in(0,1)$. Our approach follows the simple framework introduced by \citet{baraniuk2008simple} and later extended by \citet{mendelson2008uniform} to broader classes of random matrices.

In this section, we show how the concentration-of-measure properties of sub-Gaussian distributions yield a straightforward proof that random sub-Gaussian matrices satisfy the RIP. Specifically, we aim to prove that if an $n \times p$ matrix  $\rmX$ is drawn at random with sufficiently large $n$, then with high probability there exists a constant $\delta_k \in (0, 1)$ such that
\begin{equation}\label{equation:rip_prof}
(1 - \delta_k) \normtwo{\bbeta}^2 \leq \normtwo{\rmX \bbeta}^2 \leq (1 + \delta_k) \normtwo{\bbeta}^2
\end{equation}
holds for all $\bbeta \in\sB_0[k] \triangleq \{\balpha \mid  \normzero{\balpha} \leq k\}$;
that is, for all $k$-sparse vectors.
As is seen in the Gaussian case (see Theorem~\ref{theorem:rip_gaussian_uppdeltak}).

To ensure the RIP, we impose two conditions on the underlying distribution:
\begin{itemize}
\item \textit{Sub-Gaussianity.} We assume the entries of $\rmX$ are  sub-Gaussian.
This allows us to apply the clean concentration bound in Theorem~\ref{theorem:concen_under_sgrow} or Theorem~\ref{theorem:concen_subgaus}.
The argument can be adapted to strictly sub-Gaussian distributions using Corollary~\ref{corollary:concen_SSG_vec}, which yields explicit constants.

\item \textit{Approximate norm preservation.} We require that the matrix preserves the Euclidean norm in expectation. This is achieved by requiring isotropic or by setting
\begin{equation}
\Exp[\rx_{ij}^2] = \frac{1}{n},
\end{equation}
so that each entry has variance $1/n$. Under this normalization, for any fixed vector $\bbeta$, we have $\Exp[\normtwo{\rmX \bbeta}^2] = \normtwo{\bbeta}^2$, making the RIP a statement about concentration around the mean.
\end{itemize}

\subsection{RIP under Isotropic and Sub-Gaussian Random Ensembles}

We now show that a random matrix satisfying either the concentration inequality~\eqref{equation:norm_subg_ineq}---which holds when the rows are isotropic and sub-Gaussian---or inequality~\eqref{equation:concen_SSG_vec_res2}---which holds under strictly sub-Gaussian entries---also satisfies the restricted isometry property (RIP) of order $k$, provided that the number of rows scales at least as $k$ times a logarithmic factor. We begin by showing that any $k$-column submatrix of such a random matrix is well-conditioned, under an appropriate condition on its dimensions.
This is an extension to the singular values theorem of a matrix under Gaussian entries; see Theorem~\ref{thm:sv_gaus_mat}.

\begin{proposition}[RIP via tail bounds]\label{proposition:concen_numsam_svd}
Let $\sS \subset \{1,2,\ldots,p\}$ with $\abs{\sS} = k$. Suppose that an $n \times p$ random matrix $\rmX$ is drawn from a distribution satisfying the concentration inequality \eqref{equation:norm_subg_ineq} holds, that is,
\begin{equation}\label{equation:concen_numsam_svd_eq1}
\Pr\left(\abs{\normtwo{\rmX\bbeta}^2 - \normtwo{\bbeta}^2} > \tau\normtwo{\bbeta}^2\right) 
\leq 
2\exp(-\widetildeC\tau^2n), \quad \text{for all } \tau \in (0,1), \,\bbeta \in \real^p.
\footnote{Or the concentration inequality \eqref{equation:concen_SSG_vec_res2} holds with $\widetildeC=1/\kappa=(1-\ln(2))/2$.}
\end{equation}
Equivalently, Theorem~\ref{theorem:concen_under_sgrow} shows that $\sqrt{n}\rmX$ has
independent, isotropic, and sub-Gaussian rows with the same sub-Gaussian parameter $C$ in \eqref{equation:subgaus_rdvec}.
If
\begin{equation}\label{equation:concen_numsam_svd_eq2}
n \geq \frac{2}{3\widetildeC\delta^{2}}(7k + 2\ln(2/\varepsilon)),
\quad
\text{for $\varepsilon, \delta \in (0,1)$,}
\end{equation}
then, with probability at least $1 - \varepsilon$,
\begin{equation*}
\normtwo{\rmX_\sS^\top\rmX_\sS - \bI} \leq \delta.
\end{equation*}
\end{proposition}
\begin{proof}[of Proposition~\ref{proposition:concen_numsam_svd}]
By Theorem~\ref{theorem:covnum_sphere}, 
for $\rho \in (0,1/2)$, there exists a finite subset $\sA$ of the unit sphere $\sU = \{\bbeta \in \real^p \mid \supp(\bbeta) \subset \sS, \normtwo{\bbeta} = 1\}$ such that 
\begin{equation*}
\abs{\sA} \leq \left(1 + \frac{2}{\rho}\right)^k 
\qquad \text{and} \qquad 
\min_{\balpha \in \sA} \normtwo{\mathbf{z} - \balpha} \leq \rho \quad \text{for all } \mathbf{z} \in \sU.
\end{equation*}
That is, $\sA$ is the $\rho$-packing set of itself (Definition~\ref{definitio:cov_pack_num}).
For any  $\balpha \in \sA$, the concentration inequality \eqref{equation:concen_numsam_svd_eq1} gives, for $\tau  \in (0,1)$ depending on $\delta$ and $\rho$ to be determined later,
\begin{align*}
\Pr\left(\abs{\normtwo{\rmX\balpha}^2 - \normtwo{\balpha}^2} > \tau\normtwo{\balpha}^2 \right)
&\leq \sum_{\ba \in \sA} \Pr\left(\abs{\normtwo{\rmX\ba}^2 - \normtwo{\ba}^2} > \tau\normtwo{\ba}^2\right) \\
&\leq 2\,\abs{\sA}\exp(-\widetildeC\tau^2n)
\leq 2\left(1 + \frac{2}{\rho}\right)^k \exp(-\widetildeC\tau^2n).
\end{align*}
Since $\sA$ is a subset of the unit sphere $\sU$, the above inequality shows that the random matrix $\rmX$ satisfies
\begin{equation}\label{equation:concen_numsam_svd_pv1}
\abs{\normtwo{\rmX\balpha}^2 - \normtwo{\balpha}^2} \leq \tau, \quad \text{for all } \balpha \in \sA.
\end{equation}
with probability at least 
\begin{equation}\label{equation:concen_numsam_svd_pv2}
D\triangleq 1 - 2\left(1 + \frac{2}{\rho}\right)^k \exp(-\widetildeC\tau^2n).
\end{equation}
We now show that \eqref{equation:concen_numsam_svd_pv1} implies $\absbig{\normtwo{\rmX\bbeta}^2 - \normtwo{\bbeta}^2} \leq \delta$ for all $\bbeta \in \sU$, i.e., $\normtwo{\rmX_\sS^\top\rmX_\sS - \bI} \leq \delta$ once $\rho, \tau$ are chosen appropriately. Let $\rmH \triangleq \rmX_\sS^\top\rmX_\sS - \bI$. Then \eqref{equation:concen_numsam_svd_pv1} means that $\abs{\innerproduct{\rmH\balpha_{\sS},\balpha_{\sS}}} \leq \tau$ for all $\balpha \in \sA$ with probability exceeding $D$. 
Now consider a vector $\bbeta \in \sU$, for which we choose another vector $\balpha \in \sA$ satisfying $\normtwo{\bbeta - \balpha} \leq \rho < 1/2$. 
Using the triangle inequality and the Cauchy--Schwarz inequality (Theorem~\ref{theorem:cs_matvec}), we obtain
\begin{align*}
\abs{\innerproduct{\rmH\bbeta_{\sS},\bbeta_{\sS}}}
&= \abs{\innerproduct{\rmH\balpha_{\sS},\balpha_{\sS}} + \innerproduct{\rmH(\bbeta_{\sS} + \balpha_{\sS}),\bbeta_{\sS} - \balpha_{\sS}}}\\
& \leq \abs{\innerproduct{\rmH\balpha_{\sS},\balpha_{\sS}}} + \abs{\innerproduct{\rmH(\bbeta_{\sS} + \balpha_{\sS}),\bbeta_{\sS} - \balpha_{\sS}}} \\
&\leq \tau + \normtwo{\rmH}\normtwo{\bbeta_{\sS} + \balpha_{\sS}}\normtwo{\bbeta_{\sS} - \balpha_{\sS}} \leq \tau + 2\normtwo{\rmH}\rho.
\end{align*}
By the definition of the spectral norm (Definition~\ref{definition:spectral_norm}) and taking the supremum over all $\bbeta \in \sU$, we deduce that
\begin{equation*}
\normtwo{\rmH} \leq \tau + 2\normtwo{\rmH}\rho \quad \implies \quad \normtwo{\rmH} \leq \frac{\tau}{1 - 2\rho}.
\end{equation*}
Choosing $\tau \triangleq (1 - 2\rho)\delta\in(0,1)$ ensures $\normtwo{\rmH} \leq \delta$. By \eqref{equation:concen_numsam_svd_pv2} we conclude that
\begin{equation}\label{equation:concen_numsam_svd_pv3}
\Pr\left(\normtwo{\rmX_\sS^\top\rmX_\sS - \bI} \leq  \delta\right) 
\geq  1- 2\left(1 + \frac{2}{\rho}\right)^k \exp(-\widetildeC(1 - 2\rho)^2\delta^2n).
\end{equation}
Thus, the event $\normtwo{\rmX_\sS^\top\rmX_\sS - \bI} \leq \delta$ occurs with probability at least $1 - \varepsilon$ whenever
\begin{equation*}
n \geq \frac{1}{\widetildeC(1 - 2\rho)^2\delta^{2}}\left(k\ln(1 + 2/\rho) + \ln(2/\varepsilon)\right).
\end{equation*}
The choice $\rho \triangleq 0.2/(e^{3/2} - 1)<0.5$ makes $(\ln(1 + 2/\rho))/(1 - 2\rho)^2 \leq 14/3$ and $1/(1 - 2\rho)^2 \leq 4/3$.
Consequently, the required condition on $n$ is fulfilled if
\begin{equation*}
n \geq \frac{2}{3\widetildeC\delta^{2}}\left(7k + 2\ln(2/\varepsilon)\right),
\end{equation*}
which completes the proof.
\end{proof}

The above argument applies unchanged when replacing coordinate subspaces indexed by $\sS$ with arbitrary $k$-dimensional subspaces of $\real^p$. Moreover, the result depends only on the columns of $\rmX$ indexed by $\sS$; thus, it equally holds for an $n\times k$ sub-Gaussian random matrix $\rmZ$ (see Example~\ref{example:subgaus_rdmat}). In particular, for such a matrix,
\begin{equation*}
\normtwo{\frac{1}{n}\rmZ^\top\rmZ - \bI} \leq \delta
\end{equation*}
with probability at least $1 - \varepsilon$ provided that condition \eqref{equation:concen_numsam_svd_eq2} is satisfied.

We now present the main result of this section: the RIP for random matrices whose entries exhibit sub-Gaussian behavior.
\begin{theoremHigh}[RIP via tail bounds]\label{theorem:concen_numsam_rip}
Let   $\rmX\in\real^{n\times p}$ be a random matrix  drawn from a distribution satisfying the concentration inequality \eqref{equation:norm_subg_ineq}:
\begin{equation*}
\Pr\left(\abs{\normtwo{\rmX\bbeta}^2 - \normtwo{\bbeta}^2} > \tau\normtwo{\bbeta}^2\right) 
\leq 2\exp(-\widetildeC\tau^2n), 
\quad \text{for all $\tau \in (0,1)$ and $\bbeta \in \real^p$}.
\footnote{Or the concentration inequality \eqref{equation:concen_SSG_vec_res2} holds with $\widetildeC=1/\kappa=(1-\ln(2))/2$.}
\end{equation*}
Equivalently,  Theorem~\ref{theorem:concen_under_sgrow} shows that $\sqrt{n}\rmX$ has
independent, isotropic, and sub-Gaussian rows with the same sub-Gaussian parameter $C$ in \eqref{equation:subgaus_rdvec}.
Then the following statements hold: 
\begin{itemize}
\item If 
\begin{equation*}
n \geq 
\frac{2}{3\widetildeC\delta^{2}}\left[k \ln (e^9(p/k)^2) + 2\ln(2/\varepsilon)\right],
\quad \text{for $\delta, \varepsilon \in (0,1)$},
\end{equation*}
then with probability at least $1 - \varepsilon$ the restricted isometry constant $\delta_k$ of order $k$ of $\rmX$ satisfies $\delta_k \leq \delta$.

\item
Setting $\varepsilon = \exp(-3\widetildeC\delta^2 n/8)$ yields the condition
\begin{equation}\label{equation:concen_numsam_rip_cond}
n \geq \frac{8}{3\widetildeC\delta^2}k\ln\left(\frac{ep}{k}\right) + \frac{12k}{\widetildeC\delta^2},
\end{equation}
which guarantees that $\delta_k \leq \delta$ with probability at least $1 - 2\exp\left(-\delta^2 n/(2C)\right)$. 
\end{itemize}
\end{theoremHigh}
\begin{proof}[of Theorem~\ref{theorem:concen_numsam_rip}]
The conditioning of a single $k$-column submatrix $\rmX_\sS$ (with $\abs{\sS} = k$) was analyzed in Proposition~\ref{proposition:concen_numsam_svd}; we adopt the same notation here.
Recall from Definition~\ref{definition:rip22} that the restricted isometry constant of order $k$
is given by
$$
\delta_k = \max_{\sS \subset \{1,2,\ldots,p\}, \abs{\sS} = k} \normtwo{\rmX_\sS^\top\rmX_\sS - \bI}.
$$
Applying the union bound over all $\binom{p}{k}$ subsets $\sS \subset \{1,2,\ldots,p\}$ of size $k$, and using \eqref{equation:concen_numsam_svd_pv3}, we obtain
\begin{align*}
\Pr(\delta_k > \delta) 
&\leq \sum_{\sS \subset \{1,2,\ldots,p\}, \abs{\sS} = k} \Pr\left(\normtwo{\rmX_\sS^\top\rmX_\sS - \bI} \geq \delta\right) \\
&\leq 2\binom{p}{k}\left(1 + \frac{2}{\rho}\right)^k \exp(-\widetildeC(1 - 2\rho)^2\delta^2n) \\
&\leq 2\left(\frac{ep}{k}\right)^k \left(1 + \frac{2}{\rho}\right)^k \exp(-\widetildeC(1 - 2\rho)^2\delta^2n),
\end{align*}
where the last inequality follows from Problem~\ref{prob:binom_ineq}. 
This shows that the event $\delta_k\leq \delta$ holds with probability at least $1-\varepsilon$ provided that
$$
n \geq
\frac{1}{\widetildeC(1 - 2\rho)^2\delta^{2}}\left( k\ln(1 + 2/\rho) + k \ln (ep/k) + \ln(2/\varepsilon)\right).
$$
Once again, the choice $\rho \triangleq 0.2/(e^{3/2} - 1)<0.5$ yields that $\delta_k \leq \delta$ with probability at least $1 - \varepsilon$ provided
\begin{equation*}
n \geq 
\frac{1}{\widetildeC\delta^2}\left(\frac{4}{3}k\ln(ep/k) + \frac{14}{3}k + \frac{4}{3}\ln(2/\varepsilon)\right).
\end{equation*}
This  completes the proof.
\end{proof}

Since the entries of Gaussian and Rademacher random matrices are sub-Gaussian with variance 1 (see Example~\ref{example:subgaus_rdmat} and Problem~\ref{prob:rademacher_var}), it follows immediately that these matrices satisfy the RIP under condition~\eqref{equation:concen_numsam_rip_cond}.

\begin{corollary}\label{corollary:gaussben_m_bound_rip}
Let $\rmX\in\real^{n \times p}$ be a Gaussian or Rademacher random matrix. 
Then there exists a universal constant $C > 0$ such that the restricted isometry constant of $\frac{1}{\sqrt{n}}\rmX$ satisfies $\delta_k \leq \delta$ with probability at least $1 - \varepsilon$ provided
\begin{equation}\label{equation:gaussben_m_bound_rip}
n \geq C \delta^{-2} (k \ln(ep/k) + \frac{7}{2}k +  \ln(2/\varepsilon)).
\end{equation}
\end{corollary}
The normalization $\frac{1}{\sqrt{n}}\rmX$ is natural:
for any fixed  vector $\bbeta\in\real^p$, we have  $$
\Exp[\normtwobig{\frac{1}{\sqrt{n}}\rmX\bbeta}^2] = \normtwo{\bbeta}^2,
$$  
provided that the entries of $\rmX$ are independent, mean-zero, and have unit variance (as is the case for both Gaussian and Rademacher ensembles, or a $\strsubnormal(1/n)$ variable).
Consequently, the restricted isometry constant $\delta_k$ quantifies the maximal deviation of 
$\normtwobig{\frac{1}{\sqrt{n}}\rmX\bbeta}^2$ from its expectation, uniformly over all $k$-sparse vectors $\bbeta$.

For Gaussian matrices we will slightly improve on~\eqref{equation:gaussben_m_bound_rip} in the next subsection by making the constants explicit.

\subsection{RIP under SSG Random Ensembles}
We shall now show how the concentration of measure inequality in Corollary~\ref{corollary:concen_SSG_vec} (which tells us that the norm of a sub-Gaussian random vector strongly concentrates about its mean) can be used together with covering arguments to prove the RIP for sub-Gaussian random matrices. 
Similarly, our general approach will again be to construct nets of points in each $k$-dimensional subspace, apply \eqref{equation:concen_SSG_vec_res2} to all of these points through a union bound, and then extend the result from our finite set of points to all possible $k$-dimensional signals. 
The resulting  statement is often found in the literature.

\begin{theoremHigh}[RIP under SSG, Theorem 5.2 of \citet{baraniuk2008simple}]\label{theorem:ssg_sat_rip}
Fix $\delta \in (0, 1)$. 
Let $\rmX\in\real^{n\times p}$ be a random matrix whose entries $\rx_{ij}$ are i.i.d. with $\rx_{ij} \sim \strsubnormal(1/n)$. If
\begin{equation}\label{equation:ssg_sat_rip_res}
n \geq \kappa_1 k \ln\left(\frac{p}{k}\right),
\end{equation}
then $\rmX$ satisfies the RIP of order $k$ with the prescribed constant $\delta$, with probability exceeding $1 - 2\exp(-\kappa_2 n)$. 
Here, $\kappa_1 > 1$ is an arbitrary constant, and $\kappa_2 = \delta^2 / 2\kappa - 1/\kappa_1 - \ln(42e/\delta)$, where $\kappa = 2 / (1 - \ln(2)) \approx 6.52$.
\end{theoremHigh}
\begin{proof}[of Theorem~\ref{theorem:ssg_sat_rip}]
First, observe that it suffices to verify inequality~\eqref{equation:rip_prof} for vectors $\bbeta$ satisfying $\normtwo{\bbeta} = 1$, since $\rmX$ is linear.
Next, fix an index set $\sS \subset \{1, 2, \ldots, p\}$ with $\abs{\sS} = k$, and let $\sC_{\sS}$ denote the $k$-dimensional subspace spanned by the columns of $\rmX_{\sS}$. 
We construct a finite set of points $\sA_{\sS} \subseteq \sC_{\sS}$ such that  $\normtwo{\balpha} = 1$ for all $\balpha \in \sA_{\sS}$, and for every $\bbeta \in \sC_{\sS}$ with $\normtwo{\bbeta} = 1$, there exists some $\balpha \in \sA_{\sS}$ satisfying:
$$
\min_{\balpha \in \sA_{\sS}} \normtwo{\bbeta - \balpha} \leq \delta/14.
$$
By Theorem~\ref{theorem:covnum_sphere}, such a set $\sA_{\sS}$ can be chosen with cardinality at most $\abs{\sA_{\sS}} \leq (42/\delta)^k$. 
Repeating this construction for every possible index set $\sS$ of size $k$, we define the union
$$
\sA = \bigcup_{\sS : \abs{\sS} = k} \sA_{\sS}.
$$
There are $\binom{p}{k}$ such index sets $\sS$. 
Using Problem~\ref{prob:binom_ineq}, we can bound this number by $\binom{p}{k} \leq \left(\frac{ep}{k}\right)^k$.
Hence, $\abs{\sA} \leq (42ep/\delta k)^k$. 
Since the entries of $\rmX$ follow a SSG distribution, Corollary~\ref{corollary:concen_SSG_vec} implies inequality~\eqref{equation:concen_SSG_vec_res2}:
$$
\Pr\left(\abs{\normtwo{\rmX\bbeta}^2 - \normtwo{\bbeta}^2} \geq \varepsilon \normtwo{\bbeta}^2\right) \leq 2 \exp\left(-\frac{n \varepsilon^2}{\kappa}\right), 
\quad \text{for any $\bbeta\in\real^p$}.
$$ 
Applying the union bound over all $\balpha\in\sS$ with $\varepsilon = \delta/\sqrt{2}$, we find that with probability at least
\begin{equation}\label{equation:ssg_sat_rip_pr1}
\widetildep \triangleq 1 - 2\left(\frac{42ep}{\delta k}\right)^k \exp\left(-\frac{n\delta^2}{2\kappa}\right),
\end{equation}
the following holds for all $\balpha \in \sA$:
\begin{equation}\label{equation:ssg_sat_rip_pr11}
(1 - \delta/\sqrt{2})\normtwo{\balpha}^2 \leq \normtwo{\rmX \balpha}^2 \leq (1 + \delta/\sqrt{2})\normtwo{\balpha}^2, \quad \text{for all } \balpha \in \sA.
\end{equation}
Now, if $n$ satisfies condition~\eqref{equation:ssg_sat_rip_res}, then
$$
\ln\left(\frac{42ep}{\delta k}\right)^k \leq k\left(\ln\left(\frac{p}{k}\right) + \ln\left(\frac{42e}{\delta}\right)\right) \leq \frac{n}{\kappa_1} + n\ln\left(\frac{42e}{\delta}\right).
$$
Consequently, the probability $\widetildep$  in \eqref{equation:ssg_sat_rip_pr1} exceeds $1 - 2\exp(-\kappa_2 n)$, as required.

However, \eqref{equation:ssg_sat_rip_pr11} holds for all $\balpha\in\sS$. We now define $D$ as the smallest number such that
\begin{equation}\label{equation:ssg_sat_rip_pr2}
\normtwo{\rmX \bbeta} \leq \normtwo{\rmX}\normtwo{\bbeta} \leq \sqrt{1 + D}, \quad \text{for all } \bbeta \in \sB_0[k], \; \normtwo{\bbeta} = 1,
\end{equation}
where the first inequality follows from the definition of the spectral norm.
To prove the RIP condition, our goal is to show that $D \leq \delta$. For this, by the construction of the set $\sA$, for any $\bbeta \in \sB_0[k]$ with $\normtwo{\bbeta} = 1$, there exists a vector $\balpha \in \sA$ such that $\normtwo{\bbeta - \balpha} \leq \delta/14$ and  $\bbeta - \balpha \in \sB_0[k]$. 
Applying the triangle inequality,
$$
\normtwo{\rmX \bbeta} \leq \normtwo{\rmX \balpha} + \normtwo{\rmX(\bbeta - \balpha)} \leq \sqrt{1 + \delta/\sqrt{2}} + \sqrt{1 + D} \cdot \delta/14.
$$
Because $D$ is defined as the smallest constant satisfying~\eqref{equation:ssg_sat_rip_pr2}, it follows that $\sqrt{1 + D} \leq \sqrt{1 + \delta/\sqrt{2}} + \sqrt{1 + D} \cdot \delta/14$. 
Rearranging gives
$$
\sqrt{1 + D} \leq \frac{\sqrt{1 + \delta/\sqrt{2}}}{1 - \delta/14} \leq \sqrt{1 + \delta}.
$$
This establishes the upper bound of the RIP condition.
The lower bound follows similarly:
$$
\normtwo{\rmX \bbeta} \geq \normtwo{\rmX \balpha} - \normtwo{\rmX(\bbeta - \balpha)} \geq \sqrt{1 - \delta/\sqrt{2}} - \sqrt{1 + \delta} \cdot \delta/14 \geq \sqrt{1 - \delta},
$$
which completes the proof.
\end{proof}

Note that, in light of the measurement bounds for the RIP established in Theorem~\ref{theorem:measur_bd_rip}, condition~\eqref{equation:ssg_sat_rip_res} achieves the optimal number of measurements---up to a constant factor.

\subsection{Other Issues}
Using random matrices to construct $\rmX$ offers several additional benefits. To illustrate these, we again focus on the RIP condition.

First, random constructions yield democratic measurements: any sufficiently large subset of the measurements suffices for stable signal recovery \citep{davenport2009simple, laska2011democracy}. Consequently, systems based on random $\rmX$ are inherently robust to the loss or corruption of a small fraction of measurements.

Second---and perhaps more importantly---in practice we are often interested in signals that are sparse not in the standard basis, but with respect to some other orthonormal basis $\bQ$; see the next section. 
In such cases, what we actually require is that the product $\rmX \bQ$ satisfies the RIP. If $\rmX$ were constructed deterministically, we would need to tailor it specifically to the structure of $\bQ$. However, when $\rmX$ is drawn at random, this dependency disappears. For example, if $\rmX$ has i.i.d. Gaussian entries and $\bQ$ is any fixed orthonormal basis, then $\rmX \bQ$ also has i.i.d. Gaussian entries (see Lemma~\ref{lemma:affine_mult_gauss}). 
Therefore, as long as $n$ is sufficiently large, $\rmX \bQ$ satisfies the RIP with high probability---just as $\rmX$ does on its own. Although less immediate, analogous results hold for sub-Gaussian ensembles as well \citep{baraniuk2008simple}. This desirable feature is often referred to as \textit{universality}, and it represents a major advantage of random matrix constructions.

On the other hand, the work of \citet{agarwal2010fast} shows that the restricted strong convexity (RSC) and restricted strong smoothness (RSS) properties (Definition~\ref{definition:res_scss_mat}) hold whenever the rows of $\rmX$ are drawn independently from a sub-Gaussian distribution over $\real^p$ 
with a nonsingular covariance matrix. This result is particularly valuable because many real-world datasets---such as those arising in gene-expression analysis---can be reasonably modeled as sub-Gaussian random vectors. Hence, with high probability, such data matrices satisfy RSC/RSS, enabling the use of sparse recovery algorithms in practical applications.

If one is willing to accept a modest increase in the number of rows of $\rmX$, even more efficient constructions become available---particularly those that support fast matrix-vector multiplication. The seminal work of \citet{candes2005decoding} showed that selecting each row of a Fourier matrix independently with probability $\mathcalO\left({k \ln^4 p}/{p}\right)$ yields a matrix satisfying the RIP with high probability. This was later improved by \citet{haviv2017restricted}, who reduced the required sampling rate to $\mathcalO\left(k \ln^2 k {\ln(p)}/{p}\right)$. Crucially, multiplying such a structured matrix by a $k$-sparse vector takes only $\mathcalO(k \ln^2 p)$ time, compared to $\mathcalO(k^2 \ln p)$  time for a dense Gaussian matrix (under the typical scaling $n\sim k\ln(p/k)$).

\section{Sparse Optimization under Sub-Gaussian Ensembles}

Towards this end, observe that all the $\ell_1$-minimization results presented so far are deterministic: they hold simultaneously for all vectors $\bbeta$, provided the measurement matrix $\rmX$ satisfies the RIP
(see Section~\ref{section:spar_rec_rip}). 
This is a powerful theoretical property. However, as noted earlier, it is generally infeasible in practice to verify or guarantee that a given deterministic matrix satisfies the RIP. In fact, the only known constructions that provably satisfy the RIP do so with high probability---thanks to randomness. For instance, recall Theorem~\ref{theorem:concen_numsam_rip} or Theorem~\ref{theorem:ssg_sat_rip}, which illustrate how randomness enables slightly weaker---but still highly useful---probabilistic guarantees.

We now combine the results of this chapter to state the main theorem on sparse recovery via the $\ell_1$-minimization problem \eqref{opt:p1} (p.~\pageref{opt:p1}) from random measurements.

\begin{theoremHigh}[Recovery of \eqref{opt:p1} under Sub-Gaussian]\label{theorem:rec_ell1_subg}
Let $\rmX$ be an $n \times p$ random matrix with independent isotropic sub-Gaussian rows. Let $k < p$, $\varepsilon \in (0,1)$, and suppose
\begin{equation*}
n \geq C_1k\ln(ep/k) +  C_2 k+ C_3\ln(2/\varepsilon),
\end{equation*}
where the constants $C_1, C_2, C_3 > 0$ depend only on the sub-Gaussian parameters $\Phi, \xi$. 
Then, with probability at least $1 - \varepsilon$, every $k$-sparse vector $\bbeta$ is exactly recovered from $\rvy = \rmX\bbeta$ by solving the $\ell_1$-minimization problem \eqref{opt:p1}.
\end{theoremHigh}
\begin{proof}[of Theorem~\ref{theorem:rec_ell1_subg}]
The result follows by combining Theorem~\ref{theorem:concen_numsam_rip} (which ensures the RIP holds with high probability for sub-Gaussian matrices) and Theorem~\ref{theorem:ectg_el1_rip} (which guarantees exact recovery under the RIP). Note that exact sparse recovery is invariant under rescaling of the measurement matrix, so normalization does not affect the conclusion.
\end{proof}

\subsection{Sparsity in an Orthonormal Basis}
In many applications, sparsity is not expressed in the canonical basis but rather in some other orthonormal basis. That is, the signal of interest takes the form $\balpha = \bQ\bbeta$, where $\bQ\in\real^{p\times p}$ is an orthogonal matrix and $\bbeta\in\real^p$ is $k$-sparse. Taking measurements of $\balpha$ with a random matrix $\rmX$ gives
\begin{equation*}
	\rvy = \rmX\balpha = \rmX\bQ\bbeta.
\end{equation*}
To recover $\balpha$, it suffices to first recover the sparse coefficient vector $\bbeta$ and then compute $\balpha=\bQ\bbeta$. Thus, this generalized setting reduces to the standard compressive sensing problem with effective measurement matrix $\rmX' = \rmX\bQ$.

Consequently, we consider the model where $\rmX$ is a random $n\times p$ matrix and $\bQ$ is a fixed (deterministic) orthogonal matrix. Importantly, the analysis developed in the preceding sections extends naturally to this scenario.

\begin{corollary}\label{corollary:conce_orth_subg}
Let $\bQ \in \real^{p \times p}$ be a (fixed) orthogonal matrix. 
Suppose the $n \times p$ random matrix $\rmX$ satisfies the concentration inequality
\begin{equation}\label{equation:conce_orth_subg}
\Pr\left(\abs{\normtwo{\rmX\bbeta}^2 - \normtwo{\bbeta}^2} 
> \tau\normtwo{\bbeta}^2\right) \leq 2\exp(-\widetildeC \tau^2n), 
\quad \text{for all $\tau \in (0,1)$ and $\bbeta \in \real^p$},
\end{equation}
for some constant  $\widetildeC>0$.
Let $\delta, \varepsilon \in (0,1)$. Then the restricted isometry constant $\delta_k$ of $\rmX\bQ$ satisfies $\delta_k \leq \delta$ with probability at least $1 - \varepsilon$, provided that
\begin{equation*}
n 
\geq 
\frac{1}{\widetildeC\delta^2}\left(\frac{4}{3}k\ln(ep/k) + \frac{14}{3}k + \frac{4}{3}\ln(2/\varepsilon)\right).
\end{equation*}
\end{corollary}
\begin{proof}[of Corollary~\ref{corollary:conce_orth_subg}]
The key observation is that the concentration inequality \eqref{equation:conce_orth_subg} also holds for the matrix  $\rmX\bQ$. 
Indeed, for any $\bbeta \in \real^p$, define $\bbeta' \triangleq \bQ\bbeta$. Orthogonality of $\bQ$ shows
\begin{align*}
\Pr\left(\abs{\normtwo{\rmX\bQ\bbeta}^2 - \normtwo{\bbeta}^2} 
> \tau\normtwo{\bbeta}^2\right) 
&= \Pr\left(\abs{\normtwo{\rmX\bbeta'}^2 -\normtwo{\bQ^{-1}\bbeta'}^2} > \tau\normtwo{\bQ^{-1}\bbeta'}^2\right) \\
&= \Pr\left(\abs{\normtwo{\rmX\bbeta'}^2 - \normtwo{\bbeta'}^2} 
> \tau\normtwo{\bbeta'}^2\right) \leq 2\exp(-\widetildeC \tau^2n).
\end{align*}
Therefore, the claim follows directly from  Theorem~\ref{theorem:concen_numsam_rip}.
\end{proof}

This result implies that sparse recovery using sub-Gaussian random matrices is universal with respect to the choice of orthonormal basis in which the signal is sparse. Crucially, the orthogonal matrix $\bQ$ can be arbitrary---and notably, it need not be known at the encoding stage when measurements $\rvy=\rmX\bQ\bbeta$ are acquired. Knowledge of $\bQ $ is only required during decoding, when solving the $\ell_1$-minimization problem.

However, it is important to emphasize that this universality is probabilistic and basis-specific: it does not mean that a single fixed matrix $\bX$ works for all possible bases simultaneously. In fact, for any given deterministic $\bX$, one can construct an orthonormal basis $\bQ$ for which sparse recovery fails. The corollary only guarantees that, for any fixed $\bQ$, a randomly drawn $\rmX$ will enable successful recovery with high probability.

\subsection{Sparse Recovery with Strictly Sub-Gaussian}

To establish a robust recovery guarantee, we rely on the fact that strictly sub-Gaussian random matrices preserve the Euclidean norm of any fixed vector with high probability. 
Specifically, a slight modification of Corollary~\ref{corollary:concen_SSG_vec} shows that for any  $\bbeta \in \real^p$,
\begin{equation}\label{equation:instop_prob_ssg}
\Pr\left(\normtwo{\rmX\bbeta}^2 \geq 2 \normtwo{\bbeta}^2\right) \leq \exp(-\kappa_3 n), 
\quad \text{with $\kappa_3 \triangleq  1/\kappa$ and $\kappa = 2 / (1 - \ln(2)) \approx 6.52$.} 
\end{equation}
We now use this concentration property to derive the following result.

\begin{theoremHigh}[Error bound of \eqref{opt:p1} under SSG]\label{theorem:errbd_ssg_p1}
Let $\bbeta^* \in \real^p$ be a fixed vector, and assume $\delta_{2k} < \sqrt{2}-1$. 
Let $\rmX\in\real^{n \times p}$ be a strictly sub-Gaussian random matrix ($\rx_{ij}\sim\strsubnormal(1/n)$) satisfying $n \geq 2\kappa_1 k \ln(p/2k)$. 
Suppose we observe measurements $\by = \rmX\bbeta^*$.  Then with probability at least $1 - 2 \exp(-\kappa_2 n) - \exp(-\kappa_3 n)$, the solution $\widehatbbeta$ to problem \eqref{opt:p1} can be reduced  to the problem \eqref{opt:p1_epsilon} (p.~\pageref{opt:p1_epsilon}) with $\epsilon = 2 \sigma_k(\bbeta^*)_2$, and $\widehatbbeta$ satisfies
\begin{equation}
\normtwobig{\widehatbbeta-\bbeta^*} \leq \frac{8 \sqrt{1 + \delta_{2k}} - (1 + \sqrt{2}) \delta_{2k}+1}{1 - (1 + \sqrt{2}) \delta_{2k}} \sigma_k(\bbeta^*)_2.
\end{equation}
\end{theoremHigh}
\begin{proof}[of Theorem~\ref{theorem:errbd_ssg_p1}]
First, by Theorem~\ref{theorem:ssg_sat_rip}, the matrix $\rmX$ satisfies the RIP of order $2k$ with probability at least $1 - 2 \exp(-\kappa_2 n)$, provided  $n \geq 2\kappa_1 k \ln(p/2k)$. 
Next, let $\sS$ denote the index set corresponding to the $k$ entries of $\bbeta^*$ with largest magnitude, and decompose $\bbeta^*$ as $\bbeta^* = \bbeta^*(\sS) + \bbeta^*(\comple{\sS})$; note that $\normtwo{\bbeta^*_{\comple{\sS}}}=\sigma_k(\bbeta^*)_2$ is the best $k$-term approximation error. 
Since $\rmX$ is strictly sub-Gaussian, Proposition~\ref{proposition:sgssg_vec} implies that the vector $\rmX_{\comple{\sS}}\bbeta^*_{\comple{\sS}}$ is also strictly sub-Gaussian. 
Applying \eqref{equation:instop_prob_ssg}, we  obtain that with probability at least $1 - \exp(-\kappa_3 n)$, $\normtwo{\rmX_{\comple{\sS}}\bbeta^*_{\comple{\sS}}} \leq 2 \normtwo{\bbeta^*_{\comple{\sS}}} = 2 \sigma_k(\bbeta^*)_2$. 
Now observe that  $\bbeta^*(\sS) \in \sB_0[k]$, so we can write the measurement vector as $\rvy=\rmX\bbeta^* = \rmX_\sS\bbeta^*_\sS + \rmX_{\comple{\sS}}\bbeta^*_{\comple{\sS}} = \rmX_\sS\bbeta^*_\sS + \bepsilon$, where $\normtwo{\bepsilon}\leq \epsilon$. 
By the union bound, both the RIP condition (for order $2k$) and the above noise bound hold simultaneously with probability at least
 $1 - 2 \exp(-\kappa_2 n) - \exp(-\kappa_3 n)$. 
Under these conditions, we may apply Theorem~\ref{theorem:error_p1epsilon} to recover  $\bbeta^*(\sS)$ from noisy measurements
$\rvy=\rmX_\sS\bbeta^*_\sS+\bepsilon$.
Since $\bbeta^*(\sS)$ is exactly $k$-sparse, its best $k$-term $\ell_1$ approximation error is zero, i.e.,  $\sigma_k(\bbeta^*(\sS))_1 = 0$.
Therefore,
$$
\normtwobig{\widehatbbeta - \bbeta^*(\sS)} \leq C_2 \epsilon \leq 2 C_2 \sigma_k(\bbeta^*)_2,
$$
where $C_2 = 4 \frac{\sqrt{1 + \delta_{2k}}}{1 - (1 + \sqrt{2}) \delta_{2k}}$, as defined in Theorem~\ref{theorem:error_p1epsilon}.
Finally, by the triangle inequality,
$$
\normtwobig{\widehatbbeta-\bbeta^*} 
= \normtwobig{\widehatbbeta - \bbeta^*(\sS) + \bbeta^*(\sS) - \bbeta^*} \leq \normtwobig{\widehatbbeta - \bbeta^*(\sS)} 
+ \normtwo{\bbeta^*(\sS) - \bbeta^*} \leq (2C_2 + 1) \sigma_k(\bbeta^*)_2,
$$
which yields the claimed error bound after substituting the expression for $C_2$ and simplifying the numerator.
\end{proof}

Thus, although a deterministic guarantee of the form in~\eqref{equation:optimal_l2_p1_cohen} would require a prohibitively large number of measurements, it is possible to achieve comparable performance with high probability using far fewer measurements---significantly fewer than suggested by Theorem~\ref{theorem:optimal_l2_p1_cohen}.

We note, however, that the above result assumes the noise level parameter $\epsilon$ is chosen as $2\sigma_k(\bbeta^*)_2$, which requires partial knowledge of $\bbeta^*$. In practice, this limitation can be addressed using data-driven parameter selection methods such as cross-validation \citep{ward2009compressed}. Moreover, more refined analyses of $\ell_1$-minimization show that similar recovery guarantees can be obtained without oracle knowledge of $\sigma_k(\bbeta^*)_2$ \citep{wojtaszczyk2010stability}.

\begin{problemset}

\item \label{prob:normineq} Prove that the function $s^{1/s}$
attains its maximum at $s=e$, and therefore  $s^{1/s} \leq e^{1/e}$ for all $s>0$.
\textit{Hint: Take the logarithm of the function and compute its derivative.}

\item \label{prob:rademacher_var} Show that a Rademacher random variable has zero mean and unit variance.

\item State and prove a sub-Gaussian version of the JL lemma (Theorem~\ref{thm:jl_lemma}).

\item \label{prob:sg_definition} \textit{SG definition.}
Let $\rx$ be a centered random variable on $\real$. 
Show that each of the following statements implies the next (assume  $\sigma^2 > 0$  is a variance proxy in the first condition).
\begin{itemize}
\item \textit{Laplace transform.} For any $\theta \in \real$,
$
\Exp\left[\exp(\theta\rx)\right] \leq \exp\left(\frac{\sigma^2 \theta^2}{2}\right).
$

\item \textit{Tail estimation.} For any $t > 0$,
$
\max\left\{\prob(\rx \geq t),\, \prob(\rx \leq -t)\right\} \leq \exp\left(\frac{-t^2}{2\sigma^2}\right).
$

\item \textit{Moment condition.} For any positive integer $k$,
$
\Exp\left[\rx^{2k}\right] \leq k!(4\sigma^2)^k.
$

\item \textit{Orlicz condition.}
$
\Exp\left[\exp\left(\frac{\rx^2}{8\sigma^2}\right)\right] \leq 2.
$

\item \textit{Laplace transform.} For any $\tau \in \real$,
$
\Exp\left[\exp(\tau\rx)\right] \leq \exp\left(\frac{24\sigma^2 \tau^2}{2}\right).
$
\end{itemize}

\item   Let $\rx$ be a real centered random variable such that $\rx \in [a,b]$ almost surely. 
Show that
$$
\Exp[\exp(\theta\rx)] \leq \exp\left(\theta^2 \frac{(b-a)^2}{8}\right)
$$
for any $\theta \in\real$, and conclude that  $\rx$ is sub-Gaussian with variance proxy ${(b-a)^2}/{4}$.

\item \label{prob:subg_vec} \textbf{Sub-Gaussian vectors.}
A random vector $\rvx\in\real^p$ is said to be sub-Gaussian with variance proxy $c^2$, denoted $\rvx\sim\subnormal(c^2)$, if it is centered and for any $\ba\in\real^p$ with $\normtwo{\ba}=1$, the scalar random variable $\ba^\top\rvx$ is sub-Gaussian with variance proxy $c^2$ (Definition~\ref{definition:isotropic_subgvec}). 

Let $\rx_1, \rx_2, \ldots, \rx_p$ be independent $\subnormal(c^2)$ random variables. 
Show that the vector $\rvx=[\rx_1, \rx_2, \ldots, \rx_p]^\top$ is a sub-Gaussian vector with variance proxy $c^2$.

\item \textbf{Concentration of the sample mean.} Let $\rx_1, \rx_2, \ldots, \rx_p$ be independent $\subnormal(c^2)$ random variables. Show that 
$$
\prob\left(\frac{1}{p} \sum_{i=1}^{p} \rx_i \geq \tau\right)
\leq \exp\left(\frac{-\tau^2 p}{2c^2}\right).
$$

\item Let $\rx_1, \rx_2, \ldots, \rx_n\sim \subnormal(c^2)$ be a sequence of zero-mean sub-Gaussian random variables. 
For any $\delta > 0$, show that  it holds with probability at least $1 - \delta$,
$$
\max_{i=1,2\ldots,n} \abs{\rx_i} \leq c \sqrt{2 \ln(2n / \delta)}.
$$
\textit{Hint: Use Proposition~\ref{proposition:max_zero_sg}.}

\item \label{prob:subgaussian_max}
Let $ \rvx \sim \subnormal(c^2) $ be a $ p $-dimensional random vector. 
Show that 
$$
\Exp\left[ \max_{\normtwo{\ba} \leq 1} \ba^\top \rvx \right] = \Exp\left[ \max_{\normtwo{\ba} \leq 1} \abs{\ba^\top \rvx} \right] \leq 4c \sqrt{p}.
$$
Moreover, for any $ \tau > 0 $, show that
$$
\prob\left[ \max_{\normtwo{\ba} \leq 1} \abs{\ba^\top \rvx} > \tau \right] = \prob\left[ \max_{\normtwo{\ba} \leq 1} \ba^\top \rvx > \tau \right] \leq 6^p \exp\left( -\frac{\tau^2}{8c^2} \right).
$$
For any $ \delta > 0 $, taking $ \tau = \sqrt{8 \ln(6)c \sqrt{p}} + 2c \sqrt{2 \ln(1/\delta)} $, show that with probability $ 1 - \delta $, it holds that
$$
\max_{\normtwo{\ba} \leq 1} \ba^\top \rvx = \max_{\normtwo{\ba} \leq 1} \abs{\ba^\top \rvx} \leq 4c \sqrt{p} + 2c \sqrt{2 \ln(1/\delta)} = 4c \sqrt{p} \left( 1 + \sqrt{\frac{\ln(1/\delta)}{2p}} \right).
$$
\textit{Hint: Use Theorem~\ref{theorem:covnum_sphere} with radius $\varepsilon=1/2$.}


\item Let $\rx_1\sim\subnormal(c_1^2)$ and $\rx_2\sim \subnormal(c_2^2)$. 
Show that $\abs{\gamma}\rx_1\sim\subnormal(\gamma^2 c_1^2)$ and $\rx_1+\rx_2\sim \subnormal((c_1+c_2)^2)$.
Show that if $\rx_1$ and $\rx_2$ are independent, then $\rx_1+\rx_2\sim \subnormal(c_1^2+c_2^2)$.

\item Determine  the constant $C$ in Corollary~\ref{corollary:gaussben_m_bound_rip}.
\textit{Hint: Use Corollary~\ref{corollary:concen_SSG_vec}.}

\end{problemset}

\part{Algorithms}

\newpage 
\chapter{A Unified Algorithm for Sparse Optimization Problems}\label{chapter:algouni}
\begingroup
\hypersetup{
linkcolor=structurecolor,
linktoc=page,  
}
\minitoc \newpage
\endgroup

\lettrine{\color{caligraphcolor}W}
We will discuss algorithms for solving sparse regression and sparse recovery problems in greater detail in  Chapters~\ref{chapter:spar} and \ref{chapter:spar_recov}.
This Chapter primarily  introduces  a unified yet simple algorithm---\textit{primal-dual algorithm}---for solving the four core optimization problems covered in this book, namely the  Lagrangian LASSO, the  constrained LASSO, $\ell_1$-minimization, and $\ell_1$-minimization with noise measurements.

\section{Fenchel's Duality Theorem: Revisited}
Motivated by these considerations, we now examine the following convex optimization problem, which is central to this book:
\begin{equation}\label{opt:pg}
\min_{\bbeta \in \real^p} F(\bX\bbeta) + G(\bbeta), 
\tag{PG}
\end{equation}
with $\bX \in \real^{n\times p}$, and $F: \real^n \to (-\infty, \infty]$, $G: \real^p \to (-\infty, \infty]$ are proper convex functions.. 
Introducing an auxiliary variable $\balpha \triangleq \bX\bbeta$, we obtain the equivalent constrained formulation:
\begin{equation}\label{opt:pg_prime}
\min_{\bbeta \in \real^p, \balpha \in \real^n} F(\balpha) + G(\bbeta) \quad \text{s.t.}\quad \bX\bbeta - \balpha = \bzero.
\tag{PG$'$}
\end{equation}
The Lagrange dual function to this problem is given by
\begin{align}
D(\blambda) &= \min_{\bbeta, \balpha} \left\{ F(\balpha) + G(\bbeta) + \innerproduct{\bX^\top \blambda, \bbeta} - \innerproduct{\blambda, \balpha} \right\} \nonumber\\
&= -\max_{\balpha \in \real^n} \left\{ \innerproduct{\balpha, \blambda} - F(\balpha) \right\} 
- \max_{\bbeta \in \real^p} \left\{ \innerproduct{\bbeta, -\bX^\top \blambda } - G(\bbeta) \right\} \nonumber\\
&= -F^*(\blambda) - G^*(-\bX^\top \blambda), \label{equation:optdg_dev}
\end{align}
where $F^*$ and $G^*$ denote  the convex conjugate functions  of $F$ and $G$, respectively (Definition~\ref{definition:conjug_func}). Consequently, the dual problem of \eqref{opt:pg} is
\begin{equation}\label{opt:dg}
\max_{\blambda \in \real^n} D(\blambda) \equiv -F^*(\blambda) - G^*(-\bX^\top \blambda).
\tag{DG}
\end{equation}

Strictly speaking, \eqref{opt:dg} is not the Lagrange dual of the original unconstrained problem~\eqref{opt:pg}, since unconstrained problems do not naturally give rise to dual variables. However, because \eqref{opt:pg} and its constrained reformulation~\eqref{opt:pg_prime} are equivalent---and share the same optimal value---we refer to \eqref{opt:dg} as the dual of \eqref{opt:pg}.\footnote{In general, equivalent primal formulations may lead to different dual problems.}

A variant of Fenchel's duality theorem (Theorem~\ref{theorem:fenchel_dual_ori}), stated below, provides sufficient conditions under which strong duality holds between the primal problem~\eqref{opt:pg} and its dual~\eqref{opt:dg}.

\begin{corollary}[Fenchel's duality theorem; \citet{rockafellar2015convex}, Theorem 31.1\index{Fenchel's duality theorem}]\label{corollary:fenchel_dual}
Let $\bX \in \real^{n \times p}$, and let $F: \real^n \to (-\infty, \infty]$ and $G: \real^p \to (-\infty, \infty]$ be proper convex functions. 
Suppose that either $\domain(F) = \real^n$ or $\domain(G) = \real^p$, and that there exists some $\bbeta\in\real^p$ such that $\bX\bbeta \in \domain(F)$. 
If the  the optima in both  \eqref{opt:pg} and \eqref{opt:dg} are attained, then strong duality holds:
\begin{equation}
\min_{\bbeta \in \real^p} F(\bX\bbeta) + G(\bbeta) = \max_{\blambda \in \real^n} -F^*(\blambda) - G^*(-\bX^\top \blambda).
\end{equation}
\end{corollary}

Under the same conditions as in the theorem above,   strong duality implies the \textit{saddle-point property} \eqref{equation:saddl_point_gen} of the Lagrange function:
\begin{equation}\label{equation:saddl_point_gen_unit}
L((\bbeta^*, \balpha^*), \blambda) 
\leq L((\bbeta^*, \balpha^*), \blambda^*) 
\leq L((\bbeta, \balpha), \blambda^*) 
\quad \text{for all } \bbeta \in \real^p, \balpha\in\real^n,\blambda \in \real^n,
\end{equation}
where $(\bbeta^*, \balpha^*)$ and $\blambda^*$ denote primal and dual optimal points, respectively.
From \eqref{equation:optdg_dev}, the optimal value can also be expressed as the solution to the following min-max problem:
\begin{align}
&\min_{\bbeta, \balpha \in \real^p} \max_{\blambda \in \real^n} F(\balpha) + G(\bbeta) + \innerproduct{\bX^\top \blambda, \bbeta} - \innerproduct{\blambda, \balpha} \nonumber \\
&= \min_{\bbeta \in \real^p} \max_{\blambda \in \real^n} 
\left( \min_{\balpha \in \real^n} -\innerproduct{\blambda, \balpha} + F(\balpha) \right) + \innerproduct{\bX^\top \blambda, \bbeta} + G(\bbeta) \nonumber\\
&= \min_{\bbeta \in \real^p} \max_{\blambda \in \real^n} \innerproduct{\bX\bbeta, \blambda} + G(\bbeta) - F^*(\blambda), 
\label{equation:fenchel_dual_minmax}
\end{align}
where the last equality follows directly from the definition of the convex conjugate. 
The interchange of the minimum and maximum is justified because if  $((\bbeta^*, \balpha^*), \blambda^*)$ is a saddle point of the full Lagrangian  $L((\bbeta, \balpha), \blambda)$, then $(\bbeta^*, \blambda^*)$ is a saddle point of the reduced Lagrangian  $\mathcalL(\bbeta, \blambda) \triangleq \min_{\balpha} L((\bbeta, \balpha), \blambda)$. 
Thus, any primal-dual optimal pair $(\bbeta^*, \blambda^*)$ solves the saddle-point problem
\begin{equation}
\min_{\bbeta \in \real^p} \max_{\blambda \in \real^n} \innerproduct{\bX\bbeta, \blambda} + G(\bbeta) - F^*(\blambda).
\end{equation}
Moreover, the reduced Lagrangian $\mathcalL$ satisfies the saddle-point inequality:
\begin{equation}\label{equation:saddl_point_gen_unit22}
\mathcalL(\bbeta^*, \blambda) 
\leq \mathcalL(\bbeta^*, \blambda^*) 
\leq \mathcalL(\bbeta,  \blambda^*) 
\quad \text{for all } \bbeta \in \real^p,\blambda \in \real^n.
\end{equation}

\index{Primal-dual algorithm}
\section{Primal-Dual Algorithm}
We consider a general optimization framework for problems of the form
\begin{equation}\label{equation:gen_opt_prim_dual}
\min_{\bbeta \in \real^p} F(\bX\bbeta) + G(\bbeta), 
\end{equation}
where  $\bX \in \real^{n\times p}$, and  $F: \real^n \to (-\infty, \infty]$, $G: \real^p \to (-\infty, \infty]$ are proper closed and convex functions.
From \eqref{equation:optdg_dev}, the dual problem associated with \eqref{equation:gen_opt_prim_dual} is
\begin{equation}\label{equation:genDUAL_opt_prim_dual}
\max_{\blambda \in \real^n}  \left\{D(\blambda)\triangleq -F^*(\blambda) - G^*(-\bX^\top\blambda)\right\}.
\end{equation}
where  $F^*$ and $G^*$ denote the convex conjugate functions of $F$ and $G$, respectively (Definition~\ref{definition:conjug_func}).

By Fenchel's duality theorem (Corollary~\ref{corollary:fenchel_dual}), strong duality holds for the primal-dual pair \eqref{equation:gen_opt_prim_dual} and \eqref{equation:genDUAL_opt_prim_dual} under mild regularity conditions on $F$ and $G$. These conditions are always satisfied in the special cases of interest in this book.
Moreover, as shown in \eqref{equation:fenchel_dual_minmax}, jointly solving the primal and dual problems is equivalent to finding a saddle point of the Lagrangian:
\begin{equation}\label{equation:el1_pri_dual_saddlpint}
\min_{\bbeta \in \real^p} \max_{\blambda \in \real^n} \innerproduct{\bX\bbeta, \blambda} + G(\bbeta) - F^*(\blambda). 
\end{equation}

The algorithm described below relies on the proximal mappings of $F^*$ and $G$ (Definition~\ref{definition:projec_prox_opt}). 
For convenience, we introduce step-size parameters $\gamma > 0$ and $\phi>0$ and define
\begin{equation}
\prox_G(\gamma; \balpha) 
\triangleq \prox_{\gamma G}(\balpha) = \arg\min_{\bbeta \in \real^p} \left\{ \gamma G(\bbeta) + \frac{1}{2} \normtwo{\bbeta - \balpha}^2 \right\}, \quad \balpha \in \real^p, 
\end{equation}
and analogously,
\begin{equation}
\prox_{F^*}(\phi; \balpha) 
\triangleq \prox_{\phi F^*}(\balpha) = \arg\min_{\bbeta \in \real^n} \left\{ \phi F^*(\bbeta) + \frac{1}{2} \normtwo{\bbeta - \balpha}^2 \right\}, \quad \balpha \in \real^n.
\end{equation}

We assume that both $\prox_{F^*}(\phi; \balpha)$ and $\prox_G(\gamma; \balpha)$ can be evaluated efficiently---either in closed form or via simple linear-time algorithms. Note that, by the Moreau decomposition (Theorem~\ref{theorem:moreau_iden} or Theorem~\ref{theorem:ext_moreau_iden}), the proximal operator of $F^*$ can be computed easily once that of $F$ is known. Although the algorithm is well-defined for arbitrary convex $F$ and $G$, its practical efficiency hinges on the ability to compute these proximal maps repeatedly and cheaply.
Thus, the following algorithm works only when the proximal mappings are inexpensive.

\begin{algorithm}[h] 
\caption{Primal-Dual Method: the General Case \citep{chambolle2011first}}
\label{alg:prim_dual_method}
\begin{algorithmic}[1] 
\Require Functions $F(\bX\bbeta)$ and $G(\bbeta)$; 
\State {\bfseries input:}   $\xi \in [0,1]$, $\gamma,\phi > 0$ such that $\gamma\phi\normtwo{\bX} < 1$; 
\State \textbf{initialize:} $\bbeta^\topzero \in \real^p$, $\blambda^\topzero \in \real^n$, $\widetildebbeta^\topzero = \bbeta^\topzero$;
\For{$t=0, 1,2,\ldots$}
\State $\blambda^\toptone \gets \prox_{F^*}(\phi; \blambda^\toptzero + \phi \bX \widetildebbeta^\toptzero)$; \Comment{(PDM$_1$), dual proximal step}
\State $\bbeta^\toptone \gets \prox_G(\gamma; \bbeta^\toptzero - \gamma \bX^\top \blambda^\toptone)$; \Comment{(PDM$_2$), primal proximal step}
\State $\widetildebbeta^\toptone \gets \bbeta^\toptone + \xi (\bbeta^\toptone - \bbeta^\toptzero)$; \Comment{(PDM$_3$), extrapolation}
\State Exit if stopping criterion is met at $t = \bar{t}$;
\EndFor
\State {\bfseries return:}  Approximation $\widehatbbeta = \bbeta^{(\bar{t})}$ to the solution of primal problem \eqref{equation:gen_opt_prim_dual};
\State {\bfseries return:}  Approximation $\widehatblambda = \blambda^{(\bar{t})}$ to the solution of dual problem \eqref{equation:genDUAL_opt_prim_dual}.
\end{algorithmic} 
\end{algorithm}

Algorithm~\ref{alg:prim_dual_method} outlines the primal-dual method for solving \eqref{equation:gen_opt_prim_dual}.\footnote{Intuitively, the algorithm alternates between a proximal gradient ascent step in the dual variable $\blambda$ and a proximal gradient descent step in the primal variable $\bbeta$ (cf. Section~\ref{section:pgd}).}
The parameter $\xi$ controls the amount of extrapolation applied to the primal iterate $\bbeta$. Specifically:
\begin{itemize}
\item When $\xi=1$, we have $\widetildebbeta^\toptone = 2\bbeta^\toptone - \bbeta^\toptzero $ (full extrapolation).
\item When $\xi=0$, $\widetildebbeta^\toptone = \bbeta^\toptone $ (no extrapolation). The method reduces to the classical Arrow--Hurwicz algorithm, originally proposed in 1958 \citep{arrow1958studies}.
\end{itemize}
The choice of step sizes $\gamma$ and $\phi$ critically affects convergence; see Theorem~\ref{theorem:conv_primdual} for details.
In this work, we analyze Algorithm~\ref{alg:prim_dual_method} with $\xi = 1$. If either $F^*$ or $G$ is uniformly convex, further acceleration can be achieved by adaptively adjusting the parameters $\xi, \gamma,$ and $ \phi$ during the iterations; see \citet{chambolle2011first}.

\paragrapharrow{Variant.}
A symmetric variant of the algorithm is obtained by swapping the order of the primal and dual updates (i.e., $\blambda^\toptone$ and $\bbeta^\toptone$) and introducing an auxiliary dual variable  $\widetildeblambda^\toptzero$:
\begin{subequations}
\begin{align}
\bbeta^\toptone &\gets \prox_G(\gamma; \bbeta^\toptzero - \gamma \bX^\top \widetildeblambda^\toptzero); \\
\blambda^\toptone &\gets \prox_{F^*}(\phi; \blambda^\toptzero + \phi \bX \bbeta^\toptone); \\
\widetildeblambda^\toptone &\gets \blambda^\toptone + \xi (\blambda^\toptone - \blambda^\toptzero).
\end{align}
\end{subequations}

\paragrapharrow{Stopping criterion.}
A natural stopping criterion is based on the primal-dual gap (see \eqref{equation:cvx_prim_dual_gap}), which in this setting takes the form
$$
H(\bbeta, \blambda) = F(\bX\bbeta) + G(\bbeta) + F^*(\blambda) + G^*(-\bX^\top\blambda) \geq 0.
$$
At optimality $(\widehatbbeta, \widehatblambda)$, we have $H(\widehatbbeta, \widehatblambda) = 0$. 
Thus, one may terminate the algorithm at iteration $t$ when
$H(\bbeta^\toptzero, \blambda^\toptzero) \leq \varepsilon$ for a prescribed tolerance $\varepsilon > 0$.

However, in four of the key problems presented in the next section, the function $F$ may take the value $+\infty$, which can cause the primal-dual gap $H$ to be infinite during intermediate iterations---rendering it uninformative about the quality of the current iterate. In such cases, the algorithm can be modified (e.g., by restricting iterates to the effective domain of $F$) to avoid infinite values.

\section{Different Problems under the Primal-Dual Framework}\label{section:four_rep_prmdual}
Before proceeding with the convergence analysis of the algorithm, we illustrate how the primal-dual framework applies to several sparse optimization problems.

\subsection{Lagrangian LASSO} 
The Lagrangian LASSO (also known as  the $\ell_1$-regularized least squares) problem~\eqref{opt:ll} (p.~\pageref{opt:ll}) is given by
\begin{equation}\label{equation:ll_primdual}
\min_{\bbeta \in \real^p} \normone{\bbeta} + \frac{\eta}{2} \normtwo{\bX\bbeta - \by}^2, 
\end{equation}
with some regularization parameter $\eta > 0$. 
This formulation is equivalent to the standard Lagrangian LASSO~\eqref{opt:ll} after the reparameterization $\lambda = \eta^{-1}$. 
It fits the general form \eqref{equation:gen_opt_prim_dual} with $G(\bbeta) = \normone{\bbeta}$ and $F(\bX\bbeta) = \frac{\eta}{2} \normtwo{\bX\bbeta - \by}^2$.
Clearly, $F$ is not only convex but also smooth.
The convex conjugate of $F$ (see Example~\ref{example:self_conj}) is given by 
$$
F^*(\blambda) = \innerproduct{\by, \blambda} + \frac{1}{2\eta} \normtwo{\blambda}^2.
$$
The convex conjugate of $G$ (see Problem~\ref{prob:conju_norms}) is the indicator function of the unit $\ell_\infty$-ball:
\begin{equation}\label{equation:conj_indi}
G^*(\bnu) 
= \delta_{\sB_{\infty}[1]}(\bnu) =
\begin{cases}
0, & \text{if } \norminf{\bnu} \leq 1; \\
\infty, & \text{otherwise},
\end{cases} 
\end{equation}
where $\sB_{\infty}[1]$ denotes the closed unit-ball under the $\ell_\infty$-norm (Definition~\ref{definition:open_closed_ball}).
Consequently, the dual problem associated with \eqref{equation:ll_primdual} is
$$
\max_{\blambda \in \real^n} -\innerproduct{\by, \blambda} - \frac{1}{2\eta} \normtwo{\blambda}^2 \quad \text{s.t.}\quad \norminf{\bX^\top \blambda} \leq 1,
$$
and the corresponding saddle-point problem becomes
$$
\min_{\bbeta \in \real^p} \max_{\blambda \in \real^n} \innerproduct{\bX\bbeta - \by, \blambda} - \frac{1}{2\eta} \normtwo{\blambda}^2 + \normone{\bbeta}.
$$
A straightforward computation yields the proximal operator of $F$:
\begin{equation}\label{equation:primdual_ll_fprox}
\prox_F(\phi; \blambda) = \frac{\phi \eta}{\phi \eta+1} \by + \frac{1}{\phi \eta + 1} \blambda.
\end{equation}
By the extended Moreau decomposition (Theorem~\ref{theorem:ext_moreau_iden}), the proximal mapping of $F^*$ is given by
$$
\prox_{F^*}(\phi; \blambda) = \left(1 - \frac{\eta}{\eta + \phi}\right) \blambda - \frac{\phi^2}{\eta + \phi} \by.
$$
The proximal mapping of $\gamma G(\bbeta) = \gamma \normone{\bbeta}$ is the well-known  soft-thresholding operator (see Example~\ref{example:soft_thres}):
$$
\prox_G(\gamma; \bz) = \mathcalT_\gamma(\bz) = \left\{\mathcalT_\gamma(z_i)\right\}_{i=1}^p 
= [\abs{\bz} - \gamma \bone]_+ \hadaprod \sign(\bz),
$$
where $\hadaprod$ denotes element-wise multiplication and $[\cdot]_+$ applies the positive part component-wise.

Combining these components, the primal-dual algorithm for the Lagrangian LASSO \eqref{equation:ll_primdual} takes the following form:
\begin{subequations}
\begin{align}
\blambda^\toptone &\gets \left(1 - \frac{\eta}{\eta + \phi}\right) \blambda^\toptzero + \frac{\phi^2}{\eta + \phi} (\bX \widetildebbeta^\toptzero - \by); \\
\bbeta^\toptone &\gets \mathcalT_\gamma(\bbeta^\toptzero - \gamma \bX^\top \blambda^\toptone); \\
\widetildebbeta^\toptone &\gets \bbeta^\toptone + \xi (\bbeta^\toptone - \bbeta^\toptzero).
\end{align}
\end{subequations}

\subsection{Constrained LASSO}
The constrained LASSO (or the $\ell_1$-constrained least squares) problem~\eqref{opt:lc} (p.~\pageref{opt:lc}) is 
\begin{equation}\label{equation:lc_primdual}
\min_{\bbeta \in \real^p} \frac{1}{2} \normtwo{\bX\bbeta - \by}^2 \quad \text{s.t.}\quad \normone{\bbeta} \leq \Sigma.
\end{equation}
This can be cast in the form \eqref{equation:gen_opt_prim_dual} by setting $G(\bbeta) = \delta_{\sB_1[\Sigma]} (\bbeta)$ and $F(\bX\bbeta) = \frac{1}{2} \normtwo{\bX\bbeta - \by}^2$, where ${\sB_1[\Sigma]}$ denotes the $\ell_1$-ball with radius $\Sigma$.
The function $G$ is closed convex since the set $\sB_1[\Sigma]$ is closed and convex; see Exercises~\ref{exercise_convex_indica} and \ref{exercise_closed_indica}.
As before, $F$ is still continuous in this case. 
The convex conjugate of $F$ is given by 
$$
F^*(\blambda) = \innerproduct{\by, \blambda} + \frac{1}{2} \normtwo{\blambda}^2.
$$
For $G$, its conjugate is
$$
G^*(\bz) = \max_{\bbeta\in\real^p} \innerproduct{\bz, \bbeta} - G(\bbeta)
=
\max_{\bbeta\in{\sB_1[\Sigma]}} \innerproduct{\bz, \bbeta}
= \normone{\bbeta}\norminf{\bz},
$$
where the last equality follows from \holders inequality (Theorem~\ref{theorem:holder-inequality}).
Thus, the dual problem becomes
$$
\max_{\blambda \in \real^n} -\innerproduct{\by, \blambda} - \frac{1}{2} \normtwo{\blambda}^2 - \norminf{\bX^\top\blambda}.
$$
and the associated saddle-point problem is
$$
\min_{\bbeta \in \real^p} \max_{\blambda \in \real^n} 
\innerproduct{\bX\bbeta - \by, \blambda} - \frac{1}{2} \normtwo{\blambda}^2 + \delta_{\sB_1[\Sigma]} (\bbeta).
$$
From \eqref{equation:primdual_ll_fprox} with $\eta=1$, we obtain
$$
\prox_F(\phi; \blambda) = \frac{\phi}{\phi+1} \by + \frac{1}{\phi + 1} \blambda.
$$
By extended Moreau decomposition (Theorem~\ref{theorem:ext_moreau_iden}), the proximal mapping of $F^*$ is given by
$$
\prox_{F^*}(\phi; \blambda) = \left(1 - \frac{1}{1 + \phi}\right) \blambda - \frac{\phi^2}{1 + \phi} \by.
$$
The proximal mapping of $G$ is equivalent to the projection onto the $\ell_1$-ball $\sB_1[\Sigma]$ (see Definition~\ref{definition:projec_prox_opt}), which will be introduced in Section~\ref{section:alg_cons_lasso} for \textit{projected gradient descent methods}.

Putting everything together, the primal-dual method for the constrained LASSO \eqref{equation:lc_primdual} reads:
\begin{subequations}
\begin{align}
\blambda^\toptone &\gets \left(1 - \frac{1}{1 + \phi}\right) \blambda^\toptzero + \frac{\phi^2}{1 + \phi} (\bX \widetildebbeta^\toptzero - \by); \\
\bbeta^\toptone &\gets \project_{\sB_1[\Sigma]}(\bbeta^\toptzero - \gamma \bX^\top \blambda^\toptone); \\
\widetildebbeta^\toptone &\gets \bbeta^\toptone + \xi (\bbeta^\toptone - \bbeta^\toptzero).
\end{align}
\end{subequations}

\subsection{$\ell_1$-Minimization}
The $\ell_1$-minimization problem \eqref{opt:p1} (p.~\pageref{opt:p1}) is
\begin{equation}\label{equation:ell1_prim_dual}
\min_{\bbeta \in \real^p} \normone{\bbeta} \quad \text{s.t.}\quad \bX\bbeta = \by.
\end{equation}
This problem fits the general primal-dual form \eqref{equation:gen_opt_prim_dual} with $G(\bbeta) = \normone{\bbeta}$ and
$$
F(\bX\bbeta) = \delta_{\{\by\}}(\bX\bbeta) =
\begin{cases}
0, & \text{if } \bX\bbeta = \by; \\
\infty, & \text{if } \bX\bbeta \neq \by,
\end{cases}
$$
i.e., $F$ is the indicator function of the singleton set $\{\by\}$ (see Exercise~\ref{exercise_convex_indica} for its definition). 
Note that $F$ is trivially closed (Exercise~\ref{exercise_closed_indica}). 
The convex conjugate of $G$ is given by \eqref{equation:conj_indi}, and the conjugate of $F$ is given by
\begin{align*}
F^*(\blambda) &= \innerproduct{\blambda, \by}. 
\end{align*}
Because terms where the objective takes the value $-\infty$ can be ignored in maximization, we may explicitly enforce the constraint arising from the domain of $G^*$. Consequently, the dual problem \eqref{equation:genDUAL_opt_prim_dual} becomes
$$
\max_{\blambda \in \real^n} -\innerproduct{\by, \blambda} \quad \text{s.t.}\quad \norminf{\bX^\top\blambda} \leq 1.
$$
The associated saddle-point problem is
\begin{equation}
\min_{\bbeta \in \real^p} \max_{\blambda \in \real^n} \innerproduct{\bX\bbeta, \blambda} + \normone{\bbeta} - \innerproduct{\blambda, \by}.
\end{equation}

By the  extended Moreau decomposition (Theorem~\ref{theorem:ext_moreau_iden}), the proximal mapping of $F^*$ is therefore
$$
\prox_{F^*}(\phi; \blambda) = \blambda - \phi \prox_{\phi^{-1}F }(\blambda/\phi)= \blambda - \phi \by,
$$
where the proximal mapping of $F$ is the projection onto $\{\by\}$, that is, the constant map $\prox_F(\phi; \blambda) = \by$,  {for all } $\blambda \in \real^n$.

Combining these components, the primal-dual algorithm for the  $\ell_1$-minimization problem \eqref{equation:ell1_prim_dual} reads:
\begin{subequations}
\begin{align}
\blambda^\toptone &\gets \blambda^\toptzero + \phi (\bX \widetildebbeta^\toptzero - \by); \\
\bbeta^\toptone &\gets \mathcalT_\gamma(\bbeta^\toptzero - \gamma \bX^\top \blambda^\toptone); \\
\widetildebbeta^\toptone &\gets \bbeta^\toptone + \xi (\bbeta^\toptone - \bbeta^\toptzero).
\end{align}
\end{subequations}

\subsection{$\ell_1$-Minimization with Noise Measurements}
The quadratically-constrained $\ell_1$-minimization problem \eqref{opt:p1_epsilon} (p.~\pageref{opt:p1_epsilon}) is 
\begin{equation}\label{equation:p1epsi_primdual}
\min_{\bbeta \in \real^p} \normone{\bbeta} \quad \text{s.t.}\quad \normtwo{\bX\bbeta - \by} \leq \epsilon,
\end{equation}
where $\epsilon>0$ models the noise level.
This can be expressed in the form \eqref{equation:gen_opt_prim_dual} with $G(\bbeta) = \normone{\bbeta}$ and
$$
F(\bX\bbeta) = \delta_{\sB[\by, \epsilon]}(\bX\bbeta) =
\begin{cases}
0, & \text{if } \normtwo{\bX\bbeta - \by}\leq \epsilon; \\
\infty, & \text{otherwise}.
\end{cases}
$$
The function $F$ is closed because the set $\sB[\by, \epsilon]$ is closed (Exercise~\ref{exercise_closed_indica}). 
Its convex conjugate (see Problem~\ref{prob:conju_indic}) is the support function of the ball:
$$
F^*(\blambda) = \max_{\bz:\,\normtwo{\bz-\by} \leq \epsilon} \innerproduct{\bz, \blambda} =\innerproduct{\by, \blambda} + \epsilon \normtwo{\blambda}.
$$
The convex conjugate of $G$ remains as in \eqref{equation:conj_indi}. 
Therefore, the dual problem becomes
$$
\max_{\blambda \in \real^n} - \innerproduct{\by, \blambda} - \epsilon \normtwo{\blambda} \quad \text{s.t.}\quad \norminf{\bX^\top \blambda} \leq 1,
$$
and the corresponding saddle-point formulation is
$$
\min_{\bbeta \in \real^p} \max_{\blambda \in \real^n} \innerproduct{\bX\bbeta - \by, \blambda} - \epsilon \normtwo{\blambda} + \normone{\bbeta}.
$$
Example~\ref{example:proj_subs} shows that the proximal mapping of $F$ is the orthogonal projection onto the ball $\sB[\by, \epsilon]$,
\begin{align*}
\prox_F(\phi; \blambda) &= \argmin_{\bz \in \real^n:\,\normtwo{\bz - \by} \leq \epsilon} \normtwo{\bz - \blambda} 
=\by + \frac{\epsilon}{\max\{\normtwo{\blambda - \by}, \epsilon\}}(\blambda - \by)\\
&= 
\begin{cases}
\blambda & \text{if } \normtwo{\blambda - \by} \leq \epsilon; \\
\displaystyle \by + \frac{\epsilon}{\normtwo{\blambda - \by}} (\blambda - \by) & \text{otherwise}.
\end{cases}
\end{align*}
Applying the extended Moreau decomposition again yields the proximal mapping of $F^*$:
$$
\prox_{F^*}(\phi; \blambda) =
\begin{cases}
\bzero, & \text{if } \normtwo{\blambda - \phi \by} \leq \epsilon \phi; \\
\displaystyle \left(1 - \frac{\epsilon \phi}{\normtwo{\blambda - \phi \by}}\right)(\blambda - \phi \by), & \text{otherwise}.
\end{cases}
$$

Finally, the primal-dual algorithm for the quadratically constrained $\ell_1$-minimization problem \eqref{equation:p1epsi_primdual} is given by:
\begin{subequations}
\begin{align}
\blambda^\toptone &\gets \prox_{F^*}(\phi; \blambda^\toptzero + \phi \bX \widetildebbeta^\toptzero);\\
\bbeta^\toptone &\gets \mathcalT_\gamma(\bbeta^\toptzero - \gamma \bX^\top \blambda^\toptone); \\
\widetildebbeta^\toptone &\gets \bbeta^\toptone + \xi (\bbeta^\toptone - \bbeta^\toptzero).
\end{align}
\end{subequations}

\section{Optimality Condition of the Primal-Dual Algorithm}

The algorithm can be interpreted as a fixed point iteration:

\begin{theoremHigh}[Optmality condition of the primal-dual method]\label{theorem:opt_prim_dual_el1}
A point $(\widehatbbeta, \widehatblambda)$ is a fixed point of the iterations (PDM$_1$), (PDM$_2$), (PDM$_3$) (for any choice of $\xi$) if and only if $(\widehatbbeta, \widehatblambda)$ is a saddle point 
of \eqref{equation:el1_pri_dual_saddlpint}, that is, a primal-dual optimal point for the problems \eqref{equation:gen_opt_prim_dual} and \eqref{equation:genDUAL_opt_prim_dual}.
\end{theoremHigh}
\begin{proof}[of Theorem~\ref{theorem:opt_prim_dual_el1}]
From Algorithm~\ref{alg:prim_dual_method}, a fixed point must satisfy
\begin{align*}
\widehatblambda = \prox_{F^*}(\phi; \widehatblambda + \phi \bX \widehatbbeta)
\qquad\text{and}\qquad
\widehatbbeta  = \prox_G(\gamma; \widehatbbeta - \gamma \bX^\top \widehatblambda).
\end{align*}
Since $F^*$ and $G$ are proper closed and convex functions ($F^*$ is  proper closed and convex by Lemma~\ref{lemma:closedconv_conj} and Exercise~\ref{exercise:proper_conj}), 
it follows from the  Proximal Property-I (Lemma~\ref{lemma:prox_prop1}) that a fixed point $(\widehatbbeta, \widehatblambda)$ satisfies
\begin{align*}
\phi \bX \widehatbbeta \in  \phi \partial F^*(\widehatblambda)
\qquad\text{and}\qquad
- \gamma \bX^\top \widehatblambda \in \gamma \partial G(\widehatbbeta),
\end{align*}
where $\partial F^*$ and $\partial G$ denote the subdifferentials of $F^*$ and $G$, respectively (Definition~\ref{definition:subgrad}). 
Dividing both inclusions by the positive scalars $\phi$ and $\gamma$, we obtain the equivalent conditions
$$
\bzero \in -\bX\widehatbbeta + \partial F^*(\widehatblambda) 
\qquad\text{and}\qquad\bzero \in \bX^\top\widehatblambda + \partial G(\widehatbbeta).
$$
By the optimality condition under  subdifferential (Theorem~\ref{theorem:fetmat_opt}), these relations are equivalent to $\widehatbbeta$ being the minimum of the function $\bbeta \mapsto\innerproductbig{\bbeta, \bX^\top\widehatblambda} + G(\bbeta) - F^*(\widehatblambda)$ and $\widehatblambda$ being the maximum of the function $\blambda \mapsto \innerproductbig{\bX\widehatbbeta, \blambda} + G(\widehatbbeta) - F^*(\blambda)$. 
Together, this means that  $(\widehatbbeta, \widehatblambda)$ satisfies the saddle-point condition for the Lagrangian in~\eqref{equation:el1_pri_dual_saddlpint}. Hence, it is a saddle point of that problem.

Conversely, any saddle point of~\eqref{equation:el1_pri_dual_saddlpint} satisfies the above subdifferential inclusions, which---by the same proximal characterization---implies it is a fixed point of the primal-dual iterations. This completes the proof.
\end{proof}

Primal-dual algorithms  are often derived heuristically by splitting the primal and dual updates. This theorem validates that design:
The specific proximal steps for $G$ and $F^*$
are not arbitrary---they precisely encode the subdifferential optimality conditions required for a saddle point.

\section{Convergence of the Primal-Dual Algorithm}
We now discuss the convergence of the primal-dual method.
Recall that for the saddle-point problem~\eqref{equation:el1_pri_dual_saddlpint}, we define the (reduced) Lagrangian function
$$
\mathcalL(\bbeta, \blambda) \triangleq \innerproduct{\bX\bbeta, \blambda} + G(\bbeta) - F^*(\blambda),
$$
Note that this Lagrangian gives rise to the primal-dual gap, defined as
\begin{equation}\label{definition:primduialgap_re}
H(\bbeta, \blambda) \triangleq \max_{\blambda' \in \real^n} \mathcalL(\bbeta, \blambda') - \min_{\bbeta' \in \real^p} \mathcalL(\bbeta', \blambda). 
\end{equation}
It follows that $H(\bbeta, \blambda) \geq 0$ for all $\bbeta, \blambda$, and $H(\bbeta, \blambda) = 0$ if and only if $(\bbeta, \blambda)$ is a saddle point of $\mathcalL$; see the saddle-point characterization  in \eqref{equation:saddl_point_gen_unit}. 

We begin by establishing key properties of the primal-dual algorithm, which will later be used to prove its convergence.
To simplify notation, consider a sequence $(\ba^\toptzero)_{t \in \naturalset_0}$ (of scalars or vectors). 
We introduce the \textit{divided difference} operator:
\begin{align}
\Delta_\theta \ba^\toptzero &\triangleq \frac{\ba^\toptzero - \ba^\toptminus}{\theta};\\
\Delta_\theta \normtwo{\ba^\toptzero}^2 &\triangleq \frac{\normtwo{\ba^\toptzero}^2 - \normtwo{\ba^\toptminus}^2}{\theta}, \quad t \in \naturalset.
\end{align}
This quantity can be interpreted as a discrete derivative with step size $\theta$. The following identities closely mirror familiar rules from continuous calculus.

\begin{lemma}\label{lemma:disint_part}
Let $\ba, \ba^\toptzero \in \real^p$, $t \in \naturalset_0$. Then
\begin{equation}\label{equation:disint_part_eq1}
2\innerproduct{\Delta_\theta \ba^\toptzero, \ba^\toptzero - \ba} 
= \Delta_\theta \normtwo{\ba^\toptzero-\ba}^2 + \theta \normtwo{\Delta_\theta \ba^\toptzero}^2.
\end{equation}
Moreover, if $\bb^\toptzero$, $t \in \naturalset_0$, is another sequence of vectors, then the following discrete integration-by-parts formula holds for any $T \in \naturalset$,
\begin{equation}\label{equation:disint_part_eq2}
\theta \sum_{t=1}^T \left( \innerproduct{\Delta_\theta \ba^\toptzero, \bb^\toptzero} + \innerproduct{\ba^\toptminus, \Delta_\theta \bb^\toptzero} \right) = \innerproduct{\ba^{(T)}, \bb^{(T)} } - \innerproduct{\ba^\topzero, \bb^\topzero}. 
\end{equation}
\end{lemma}
\begin{proof}[of Lemma~\ref{lemma:disint_part}]
Let $\widetildeba^\toptzero \triangleq  \ba^\toptzero - \ba$. 
Then $\Delta_\theta \widetildeba^\toptzero = \Delta_\theta \ba^\toptzero$, and
\begin{align*}
2\theta \innerproduct{\Delta_\theta \ba^\toptzero, \ba^\toptzero - \ba}
&= 2\theta \innerproduct{\Delta_\theta \widetildeba^\toptzero, \widetildeba^\toptzero} = 2 \innerproduct{\widetildeba^\toptzero - \widetildeba^\toptminus, \widetildeba^\toptzero} \\
&= \innerproduct{\widetildeba^\toptzero - \widetildeba^\toptminus, \widetildeba^\toptzero - \widetildeba^\toptminus} + 
\innerproduct{\widetildeba^\toptzero - \widetildeba^\toptminus, \widetildeba^\toptzero + \widetildeba^\toptminus}\\
&=\normtwo{\widetildeba^\toptzero - \widetildeba^\toptminus}^2 + \theta \Delta_\theta\normtwo{ \widetildeba^\toptzero}^2.
\end{align*}
This proves \eqref{equation:disint_part_eq1}.

For the second part, we observe that
$$
\Delta_\theta \innerproduct{ \ba^\toptzero, \bb^\toptzero} 
= \frac{\innerproduct{\ba^\toptzero, \bb^\toptzero} - \innerproduct{\ba^\toptminus, \bb^\toptminus}}{\theta} 
= \innerproduct{\Delta_\theta \ba^\toptzero, \bb^\toptzero} + \innerproduct{\ba^\toptminus, \Delta_\theta \bb^\toptzero}.
$$
Telescoping the sum over $t = 1, 2, \ldots, T$ proves  the desired result.
\end{proof}

\begin{lemma}[Iterate inequality of  primal-dual algorithm]\label{lemma:primdual_ineq}
Let $\{\bbeta^\toptzero, \widetildebbeta^\toptzero, \blambda^\toptzero\}_{t \geq 0}$ be the sequence generated by \textnormal{(PDM$_1$), (PDM$_2$), (PDM$_3$)} in  Algorithm~\ref{alg:prim_dual_method}, and let $\bbeta \in \real^p$, $\blambda \in \real^n$ be arbitrary. Then, for $t \in \naturalset$,
\begin{align}
&\quad \frac{1}{2} \Delta_\phi \normtwo{\blambda - \blambda^\toptzero}^2 + \frac{1}{2} \Delta_\gamma \normtwo{\bbeta - \bbeta^\toptzero}^2 + \frac{\phi}{2} \normtwo{\Delta_\phi \blambda^\toptzero}^2 + \frac{\gamma}{2} \normtwo{\Delta_\gamma \bbeta^\toptzero}^2 \nonumber\\
&\leq \mathcalL(\bbeta, \blambda^\toptzero) - \mathcalL(\bbeta^\toptzero, \blambda) + \innerproduct{\bX(\bbeta^\toptzero - \widetildebbeta^\toptminus), \blambda - \blambda^\toptzero}. \label{equation:primdual_ineq}
\end{align}
Furthermore, assume that $\xi=1$ in Algorithm~\ref{alg:prim_dual_method}. 
Telescoping the sum of inequality \eqref{equation:primdual_ineq} over $t = 1, 2, \ldots, T$ gives
\begin{align}
&\frac{1}{2\phi} \left( \normtwo{\blambda - \blambda^{(T)}}^2 - \normtwo{\blambda - \blambda^\topzero}^2 \right) + \frac{1}{2\gamma} \left( \normtwo{\bbeta - \bbeta^{(T)}}^2 - \normtwo{\bbeta - \bbeta^\topzero}^2 \right) \nonumber\\
&+ \frac{1}{2\phi} \sum_{t=1}^T \normtwo{\blambda^\toptzero - \blambda^\toptminus}^2 + \frac{1}{2\gamma} \sum_{t=1}^T \normtwo{\bbeta^\toptzero - \bbeta^\toptminus}^2\nonumber \\
&\leq \sum_{t=1}^T \left( \mathcalL(\bbeta, \blambda^\toptzero) - \mathcalL(\bbeta^\toptzero, \blambda) \right) 
+ \gamma^2 \sum_{t=1}^T \innerproduct{\bX \Delta_\gamma^2 \bbeta^\toptzero, \blambda - \blambda^\toptzero}, \label{equation:primdual_ineq_eq2}
\end{align}
where $\Delta_\gamma^2\bbeta^\toptzero \triangleq \Delta_\gamma (\Delta_\gamma\bbeta^\toptzero)  $.
\end{lemma}
\begin{proof}[of Lemma~\ref{lemma:primdual_ineq}]
From (PDM$_1$), (PDM$_2$), and  Proximal Property-I (Lemma~\ref{lemma:prox_prop1}),
the iterates satisfy the relations (replacing $t+1$ by $t$)
$$
\blambda^\toptminus + \phi \bX \widetildebbeta^\toptminus - \blambda^\toptzero \in \phi \partial F^*(\blambda^\toptzero), \qquad
\bbeta^\toptminus - \gamma \bX^\top \blambda^\toptzero -\bbeta^\toptzero  \in \gamma \partial G(\bbeta^\toptzero),
$$
where $\partial F^*$ and $\partial G$ denote the subdifferentials of $F^*$ and $G$, respectively. 
Applying the subgradient inequality (Definition~\ref{definition:subgrad}), we obtain:
\begin{align*}
\phi F^*(\blambda) - \phi F^*(\blambda^\toptzero) 
&\geq \innerproduct{-\blambda^\toptzero + \blambda^\toptminus + \phi \bX \widetildebbeta^\toptminus, \blambda - \blambda^\toptzero};\\
\gamma G(\bbeta) - \gamma G(\bbeta^\toptzero) 
&\geq \innerproduct{-\bbeta^\toptzero + \bbeta^\toptminus - \gamma \bX^\top \blambda^\toptzero, \bbeta - \bbeta^\toptzero} .
\end{align*}
Dividing by $\phi$ and $\gamma$, respectively, and rewriting using the definition of the divided difference gives:
\begin{align*}
F^*(\blambda) - F^*(\blambda^\toptzero)
&\geq \innerproduct{\Delta_\phi \blambda^\toptzero, \blambda^\toptzero - \blambda }+ \innerproduct{\bX \widetildebbeta^\toptminus, \blambda - \blambda^\toptzero } ;\\
G(\bbeta) - G(\bbeta^\toptzero)
&\geq \innerproduct{\Delta_\gamma \bbeta^\toptzero, \bbeta^\toptzero - \bbeta} - \innerproduct{\bX(\bbeta - \bbeta^\toptzero), \blambda^\toptzero} .
\end{align*}
Summing these two inequalities and applying Lemma~\ref{lemma:disint_part} \eqref{equation:disint_part_eq1} yields:
\begin{align*}
\footnotesize
&
\left(\frac{1}{2} \Delta_\phi \normtwo{\blambda - \blambda^\toptzero}^2 
+ \frac{\phi}{2} \normtwo{\Delta_\phi \blambda^\toptzero}^2\right) 
+ \left( \frac{1}{2} \Delta_\gamma \normtwo{\bbeta - \bbeta^\toptzero}^2 
+ \frac{\gamma}{2} \normtwo{\Delta_\gamma \bbeta^\toptzero}^2\right) \\
&\leq F^*(\blambda) - F^*(\blambda^\toptzero) + G(\bbeta) - G(\bbeta^\toptzero) 
+ \innerproduct{\bX(\bbeta - \bbeta^\toptzero), \blambda^\toptzero} - \innerproduct{\bX \widetildebbeta^\toptminus, \blambda - \blambda^\toptzero} \\
&= 
\mathcalL(\bbeta, \blambda^\toptzero) - \mathcalL(\bbeta^\toptzero, \blambda)
+ \innerproduct{\bX(\bbeta^\toptzero - \widetildebbeta^\toptminus), \blambda - \blambda^\toptzero}.
\end{align*}
This proves \eqref{equation:primdual_ineq}.

Now suppose $\xi=1$. In this case, the extrapolation satisfies $\widetildebbeta^\toptminus =2\beta^\toptminus - \bbeta^\toptminusTWO$ so that
$$
\bbeta^\toptzero - \widetildebbeta^\toptminus = \bbeta^\toptzero - \bbeta^\toptminus - (\bbeta^\toptminus - \bbeta^\toptminusTWO) = \gamma^2 \Delta_\gamma (\Delta_\gamma\bbeta^\toptzero) \triangleq \gamma^2\Delta_\gamma^2\bbeta^\toptzero, \quad \text{for } t \geq 2.
$$
The formula holds for $t = 1$ if we set  $\bbeta^{(-1)} = \bbeta^\topzero$, which implies  $\widetildebbeta^\topzero = \bbeta^\topzero$ and $\Delta_\gamma \bbeta^\topzero = \bzero$.
Telescoping the sum of inequality \eqref{equation:primdual_ineq} over $t = 1, 2, \ldots, T$ proves \eqref{equation:primdual_ineq_eq2}.
\end{proof}

Inequality~\eqref{equation:primdual_ineq} suggests that the ideal choice would be $\widetildebbeta^\toptminus = \bbeta^\toptzero$, which would eliminate the last inner-product term. 
However, this would result in an implicit scheme, where each iteration requires solving a coupled system as difficult as the original problem.

Using the telescoped iterate inequality~\eqref{equation:primdual_ineq_eq2}, we now establish an upper bound on the cumulative primal-dual gap.

\begin{lemma}[Primal-dual gap inequality]\label{lemma:primdua_gap_ineq}
Let $\{\bbeta^\toptzero, \widetildebbeta^\toptzero, \blambda^\toptzero\}_{t \geq 0}$ be the sequence generated by {(PDM$_1$), (PDM$_2$), (PDM$_3$)} in  Algorithm~\ref{alg:prim_dual_method} with the parameter choice $\xi = 1$, and let $\bbeta \in \real^p$, $\blambda \in \real^n$ be arbitrary. Then, for any $T \in \naturalset$,
\begin{align}
&\sum_{t=1}^T \left( \mathcalL(\bbeta^\toptzero, \blambda) - \mathcalL(\bbeta, \blambda^\toptzero) \right)  + \frac{1 - \gamma\phi \normtwo{\bX}^2}{2\phi} \normtwo{\blambda - \blambda^{(T)}}^2  + \frac{1}{2\gamma} \normtwo{\bbeta - \bbeta^{(T)}}^2 \nonumber\\
&\quad
+ \frac{1 - \sqrt{\gamma\phi} \normtwo{\bX}}{2\phi} \sum_{t=1}^T \normtwo{\blambda^\toptzero - \blambda^\toptminus}^2
+ \frac{1 - \sqrt{\gamma\phi} \normtwo{\bX}}{2\gamma} \sum_{t=1}^{T-1} \normtwo{\bbeta^\toptzero - \bbeta^\toptminus}^2 
\nonumber\\
&\leq  \frac{1}{2\phi} \normtwo{\blambda - \blambda^\topzero}^2 + \frac{1}{2\gamma} \normtwo{\bbeta - \bbeta^\topzero}^2.  \label{equation:primdua_gap_ineq}
\end{align}
\end{lemma}
\begin{proof}[of Lemma~\ref{lemma:primdua_gap_ineq}]
We begin with the summed iterate inequality (Lemma~\ref{lemma:primdual_ineq} \eqref{equation:primdual_ineq_eq2}).
Using the discrete integration-by-parts formula Lemma~\ref{lemma:disint_part} \eqref{equation:disint_part_eq2}, together with the initialization $\Delta_\gamma \bbeta^\topzero = \bzero$ (by initializing $\bbeta^{(-1)} = \bbeta^\topzero$), and the definition $\Delta_\gamma^2\bbeta^\toptzero \triangleq \Delta_\gamma (\Delta_\gamma\bbeta^\toptzero)  $, 
we rewrite the second-order term as:
\begin{align}
\footnotesize
&\gamma^2 \sum_{t=1}^T \innerproduct{\bX \Delta_\gamma^2 \bbeta^\toptzero, \blambda - \blambda^\toptzero} 
= \gamma^2 \sum_{t=1}^T \innerproduct{\bX \Delta_\gamma \bbeta^\toptminus, \Delta_\gamma \blambda^\toptzero} 
+ \gamma \innerproduct{\bX \Delta_\gamma \bbeta^{(T)}, \blambda - \blambda^{(T)}} \nonumber\\
&\quad\quad =  \gamma \phi\sum_{t=1}^T \innerproduct{\bX \Delta_\gamma \bbeta^\toptminus, \Delta_\phi \blambda^\toptzero} 
+ \gamma \innerproduct{\Delta_\gamma \bbeta^{(T)}, \bX^\top(\blambda - \blambda^{(T)})}. \label{equation:primdua_gap_ineq_pv3}
\end{align}
For the first term on the right-hand side, by setting $\eta \triangleq  \sqrt{\gamma / \phi}$, we have
\begin{align}
\footnotesize
&\gamma \phi \innerproduct{\bX \Delta_\gamma \bbeta^\toptminus, \Delta_\phi \blambda^\toptzero}
\leq \gamma \phi \normtwo{\bX} \normtwo{\Delta_\gamma \bbeta^\toptminus}\normtwo{\Delta_\phi \blambda^\toptzero} \nonumber\\
&\quad\leq \frac{\gamma \phi \normtwo{\bX}}{2} \left( \eta \normtwo{\Delta_\gamma \bbeta^\toptminus}^2 + \eta^{-1} \normtwo{\Delta_\phi \blambda^\toptzero}^2 \right) \nonumber\\
&\quad= \frac{\phi \eta \normtwo{\bX}}{2\gamma} \normtwo{\bbeta^\toptminus - \bbeta^\toptminusTWO}^2 + \frac{\gamma \normtwo{\bX}}{2\eta \phi} \normtwo{\blambda^\toptzero - \blambda^\toptminus}^2 \nonumber\\
&\quad \stackrel{\eta =  \sqrt{\gamma / \phi}}{\longeq} \frac{\sqrt{\gamma \phi} \normtwo{\bX}}{2\gamma} \normtwo{\bbeta^\toptminus - \bbeta^\toptminusTWO}^2 + \frac{\sqrt{\gamma \phi} \normtwo{\bX}}{2\phi} \normtwo{\blambda^\toptzero - \blambda^\toptminus}^2, \label{equation:primdua_gap_ineq_pv2}
\end{align}
where the first inequality follows from the Cauchy--Schwarz inequality (Theorem~\ref{theorem:cs_matvec}), and the second inequality follows from the fact that $2ab \leq \eta a^2 + b^2/\eta$ for positive $a,b,\eta$.
Similarly, the second term satisfies
\begin{align*}
\footnotesize
\gamma \innerproduct{\Delta_\gamma \bbeta^{(T)}, \bX^\top(\blambda - \blambda^{(T)})}
&\leq \frac{\gamma}{2} \left( \normtwo{\Delta_\gamma \bbeta^{(T)}}^2 + \normtwo{\bX}^2 \normtwo{\blambda - \blambda^{(T)}}^2 \right) \\
&= \frac{1}{2\gamma} \normtwo{\bbeta^{(T)} - \bbeta^{(T-1)}}^2 + \frac{\gamma \phi \normtwo{\bX}^2}{2\phi} \normtwo{\blambda - \blambda^{(T)}}^2.
\end{align*}
Substituting these estimates into \eqref{equation:primdua_gap_ineq_pv3} and using that $\bbeta^{(-1)} = \bbeta^\topzero$, we obtain:
\begin{align*}
\footnotesize
&\gamma^2 \sum_{t=1}^T \innerproduct{\bX \Delta_\gamma^2 \bbeta^\toptzero, \blambda - \blambda^\toptzero} 
\leq \frac{\sqrt{ \gamma\phi} \normtwo{\bX}}{2\gamma} \sum_{t=1}^{T-1} \normtwo{\bbeta^\toptzero - \bbeta^\toptminus}^2 \\
&\quad + \frac{\sqrt{\gamma\phi } \normtwo{\bX}}{2\phi} \sum_{t=1}^T \normtwo{\blambda^\toptzero - \blambda^\toptminus}^2 
+ \frac{1}{2\gamma} \normtwo{\bbeta^{(T)} - \bbeta^{(T-1)}}^2 + \frac{\gamma \phi \normtwo{\bX}^2}{2\phi} \normtwo{\blambda - \blambda^{(T)}}^2.
\end{align*}
Plugging this bound back into the summed iterate inequality (Lemma~\ref{lemma:primdual_ineq} \eqref{equation:primdual_ineq_eq2}), and rearranging terms, yields the desired inequality~\eqref{equation:primdua_gap_ineq}.
\end{proof}

Using the results above, we now establish a descent lemma for the primal-dual algorithm.
\begin{lemma}[Decent lemma of primal-dual algorithm]\label{lemma:des_primdual}
Let $(\widehatbbeta, \widehatblambda)$ be a primal dual optimum, that is, a saddle point of the saddle-point problem \eqref{equation:el1_pri_dual_saddlpint}. 
Let $\{\bbeta^\toptzero, \widetildebbeta^\toptzero, \blambda^\toptzero\}_{t \geq 0}$ be the sequence generated by {(PDM$_1$), (PDM$_2$), (PDM$_3$)} in  Algorithm~\ref{alg:prim_dual_method} with the parameter choice $\xi = 1$ and $\gamma\phi \normtwo{\bX}^2 < 1$. Then, for any $T \in \naturalset$,
$$
\frac{1}{2\phi} \normtwo{\widehatblambda - \blambda^{(T)}}^2 + \frac{1}{2\gamma} \normtwo{\widehatbbeta - \bbeta^{(T)}}^2 \leq C \left( \frac{1}{2\phi} \normtwo{\widehatblambda - \blambda^\topzero}^2 + \frac{1}{2\gamma} \normtwo{\widehatbbeta - \bbeta^\topzero}^2 \right),
$$
where $C \triangleq (1 - \gamma\phi \normtwo{\bX}^2)^{-1}$. In particular, the iterates $\{\bbeta^\toptzero, \blambda^\toptzero\}$ remain bounded.
\end{lemma}

\begin{proof}[of Lemma~\ref{lemma:des_primdual}]
For a saddle point $(\widehatbbeta, \widehatblambda)$, the terms $\mathcalL(\widehatbbeta, \widehatblambda) - \mathcalL(\widehatbbeta, \blambda^\toptzero)$ are non-negative so that all terms on the left-hand side of \eqref{equation:primdua_gap_ineq} are positive. 
Consequently
$$
\frac{1 - \gamma\phi \normtwo{\bX}^2}{2\phi} \normtwo{\widehatblambda - \blambda^{(T)}}^2 +\frac{1}{2\gamma} \normtwo{\widehatbbeta - \bbeta^{(T)}}^2 
\leq 
\frac{1}{2\phi} \normtwo{\widehatblambda - \blambda^\topzero}^2 + \frac{1}{2\gamma} \normtwo{\widehatbbeta - \bbeta^\topzero}^2 .
$$
This completes the proof.
\end{proof}

The convergence of the primal-dual method is then established in the following theorem.
\begin{theoremHigh}[Convergence of primal-dual algorithm]\label{theorem:conv_primdual}
Let $\{\bbeta^\toptzero, \widetildebbeta^\toptzero, \blambda^\toptzero\}_{t \geq 0}$ be the sequence generated by {(PDM$_1$), (PDM$_2$), (PDM$_3$)} in Algorithm~\ref{alg:prim_dual_method} with the parameter choice $\xi = 1$ and $\gamma\phi \normtwo{\bX}^2 < 1$.
Suppose that the saddle-point problem \eqref{equation:el1_pri_dual_saddlpint} has a saddle point. 
Then the iterates  $\{(\bbeta^\toptzero, \blambda^\toptzero)\}$ converges to a saddle point $(\widehatbbeta, \widehatblambda)$ of \eqref{equation:el1_pri_dual_saddlpint} so that $\bbeta^\toptzero$ converges to a minimizer of problem \eqref{equation:gen_opt_prim_dual}.
\end{theoremHigh}
\begin{proof}[of Theorem~\ref{theorem:conv_primdual}]
By Lemma~\ref{lemma:des_primdual}, the sequence  $\{(\bbeta^\toptzero, \blambda^\toptzero)\}$ is bounded. Hence, there exists a convergent  subsequence $\{(\bbeta^{(t_k)}, \blambda^{(t_k)})\}_k$ such that 
$(\bbeta^{(t_k)}, \blambda^{(t_k)}) \to (\bbeta^\star, \blambda^\star)$ as $k \to \infty$. Choosing $(\bbeta, \blambda)$ to be a saddle point $(\widehatbbeta, \widehatblambda)$ in \eqref{equation:primdua_gap_ineq} makes all terms nonnegative. 
Since all terms on the left-hand side are nonnegative, we obtain in particular:
$$
\frac{1 - \sqrt{\gamma\phi} \normtwo{\bX}}{2\gamma} \sum_{t=1}^{T-1} \normtwo{\bbeta^\toptzero - \bbeta^\toptminus}^2 
\leq \frac{1}{2\gamma} \normtwo{\bbeta - \bbeta^\topzero}^2 + \frac{1}{2\phi} \normtwo{\blambda - \blambda^\topzero}^2.
$$
Because  the right-hand side is independent of $T$ and $\sqrt{\gamma\phi} \normtwo{\bX} < 1$, we conclude that the series $\sum_{t=1}^{\infty} \normtwobig{\bbeta^\toptzero - \bbeta^\toptminus}^2$ converges; 
hence,  $\normtwobig{\bbeta^\toptzero - \bbeta^\toptminus} \to 0$ as $t \to \infty$. 
An identical argument shows $\lim_{t \to \infty} \normtwobig{\blambda^\toptzero - \blambda^\toptminus} = 0$. 
In particular, the subsequence $(\bbeta^{(t_k - 1)}, \blambda^{(t_k - 1)})$ also converges to $(\bbeta^\star, \blambda^\star)$. 
It follows that $(\bbeta^\star, \blambda^\star)$ is a fixed point of the primal dual algorithm. By Theorem~\ref{theorem:opt_prim_dual_el1}, it is  primal-dual optimal (a saddle point).

Invoking  $(\bbeta, \blambda) = (\bbeta^\star, \blambda^\star)$ in \eqref{equation:primdual_ineq}, we have $\mathcalL(\bbeta^\toptzero, \blambda^\star) - \mathcalL(\bbeta^\star, \blambda^\toptzero) \geq 0$ by the primal-dual gap property. 
Similarly to the proof of \eqref{equation:primdual_ineq_eq2}, telescoping the sum of \eqref{equation:primdual_ineq} over $t = t_k$ to $t = T > t_k$ results in
\begin{align*}
&\frac{1}{2\phi} \left(\normtwo{\blambda^\star - \blambda^{(T)}}^2 - \normtwo{\blambda^\star - \blambda^{(t_k)}}^2\right) 
+ \frac{1}{2\gamma} \left(\normtwo{\bbeta^\star - \bbeta^{(T)}}^2 - \normtwo{\bbeta^\star - \bbeta^{(t_k)}}^2\right)  \\
&+ \frac{1}{2\phi} \sum_{t=t_k}^T \normtwo{\blambda^\toptzero - \blambda^\toptminus}^2 + \frac{1}{2\gamma} \sum_{t=t_k}^T \normtwo{\bbeta^\toptzero - \bbeta^\toptminus}^2 \nonumber 
\leq \gamma^2 \sum_{t=t_k}^T \innerproduct{\bX \Delta_\gamma^2 \bbeta^\toptzero, \blambda^\star - \blambda^\toptzero}\\
&= \gamma\phi \sum_{t=t_k}^T \innerproduct{\bX \Delta_\gamma \bbeta^\toptminus, \Delta_\phi \blambda^\toptzero} 
+ \gamma \innerproduct{\bX \Delta_\gamma \bbeta^{(T)}, \blambda^\star - \blambda^{(T)}} 
- \gamma \innerproduct{\bX \Delta_\gamma \bbeta^{(t_k )}, \blambda^\star - \blambda^{(t_k)}},
\end{align*}
where the last equality follows from the discrete integration-by-parts Lemma~\ref{lemma:disint_part} \eqref{equation:disint_part_eq2}.
Inequality \eqref{equation:primdua_gap_ineq_pv2} therefore implies
\begin{align*}
\footnotesize
&\frac{1}{2\phi} \normtwo{\blambda^\star - \blambda^{(T)}}^2 + \frac{1}{2\gamma} \normtwo{\bbeta^\star - \bbeta^{(T)}}^2 + \frac{1 - \sqrt{\gamma\phi} \normtwo{\bX}}{2\phi} \sum_{t=t_k}^T \normtwo{\blambda^\toptzero - \blambda^\toptminus}^2\\
&+ \frac{1 - \sqrt{\gamma\phi} \normtwo{\bX}}{2\gamma} \sum_{t=t_k}^{T-1} \normtwo{\bbeta^\toptzero - \bbeta^\toptminus}^2 \\
&+ \frac{1}{2\gamma} \left( \normtwo{\bbeta^{(T)} - \bbeta^{(T-1)}}^2 - \sqrt{\gamma\phi} \normtwo{\bX} \normtwo{\bbeta^{(t_k - 1)} - \bbeta^{(t_k - 2)}}^2 \right) \\
&+ \innerproduct{\bX(\bbeta^{(T)} - \bbeta^{(T-1)}), \blambda^\star - \blambda^{(T)}} - \innerproduct{\bX(\bbeta^{(t_k )} - \bbeta^{(t_k - 1)}), \blambda^\star - \blambda^{(t_k)} } \\
&\leq \frac{1}{2\phi} \normtwo{\blambda^\star - \blambda^{(t_k)}}^2 + \frac{1}{2\gamma} \normtwo{\bbeta^\star - \bbeta^{(t_k)}}^2.
\end{align*}
Since $\lim_{t \to \infty} \normtwobig{\bbeta^\toptzero - \bbeta^\toptminus} = \lim_{t \to \infty} \normtwobig{\blambda^\toptzero - \blambda^\toptminus} = 0$ and $\lim_{k \to \infty} \normtwobig{\bbeta^\star - \bbeta^{(t_k)}} = \lim_{k \to \infty} \normtwobig{\blambda^\star - \blambda^{(t_k)}} = 0$, it follows that $\lim_{T \to \infty} \normtwobig{\bbeta^\star - \bbeta^{(T)}} = \lim_{T \to \infty} \normtwobig{\blambda^\star - \blambda^{(T)}} = 0$. 
This completes the proof.
\end{proof}

\section{Rate of Convergence of the Primal-Dual Algorithm}
To analyze the rate of convergence of the primal-dual algorithm, we introduce the \textit{partial primal-dual gap}. 
For two sets $\sS_1 \subset \real^p$ and $\sS_2 \subset \real^n$, this quantity is defined as
\begin{equation}\label{definition:partial_primduialgap}
\mathcalP_{\sS_1,\sS_2}(\bbeta, \blambda) \triangleq \max_{\blambda' \in \sS_2} \mathcalL(\bbeta, \blambda') - \min_{\bbeta' \in \sS_1} \mathcalL(\bbeta', \blambda). 
\end{equation}
This is a variant of the standard primal-dual gap defined in \eqref{equation:cvx_prim_dual_gap}.
It is particularly convenient in our setting and is naturally motivated by the saddle-point problem~\eqref{equation:el1_pri_dual_saddlpint}.
Whenever the product set $\sS_1 \times \sS_2$ contains a saddle point $(\widehatbbeta, \widehatblambda)$, the partial gap satisfies $\mathcalP_{\sS_1,\sS_2}(\bbeta, \blambda) \geq 0$ for all $\bbeta, \blambda$, and $\mathcalP_{\sS_1,\sS_2}(\bbeta, \blambda) = 0$ if and only if $(\bbeta, \blambda)$ is itself a saddle point; see the saddle-point property \eqref{equation:saddl_point_gen_unit}. 
Thus, $\mathcalP_{\sS_1,\sS_2}(\bbeta, \blambda)$ serves as a natural measure of how far the pair $(\bbeta, \blambda)$  is from optimality in the saddle-point problem~\eqref{equation:el1_pri_dual_saddlpint}.

\begin{theoremHigh}[Rate of convergence of primal-dual algorithm]\label{theorem:conv_primdual_rate}
Assume the same conditions as in Theorem~\ref{theorem:conv_primdual}.
Define the averaged iterates $\widebarbbeta^{(T)} \triangleq T^{-1} \sum_{t=1}^T \bbeta^\toptzero$ and $\widebarblambda^{(T)} \triangleq T^{-1} \sum_{t=1}^T \blambda^\toptzero$. Let $\sS_1\subset\real^p, \sS_2\subset\real^n$ be bounded sets such that $\sS_1 \times \sS_2$ contains at least one saddle point  $(\widehatbbeta, \widehatblambda)$. Then,
\begin{equation}\label{equation:var_pridual_gap_ineq}
\mathcalP_{\sS_1,\sS_2}(\widebarbbeta^{(T)}, \widebarblambda^{(T)}) \leq \frac{D(\sS_1, \sS_2)}{T}, 
\end{equation}
where $D(\sS_1, \sS_2) 
\triangleq  (2\phi)^{-1} \max_{\blambda \in \sS_2} \normtwobig{\blambda - \widehatblambda}^2 
+ (2\gamma)^{-1} \max_{\bbeta \in \sS_1} \normtwobig{\bbeta - \widehatbbeta}^2$, 
and $\mathcalP_{\sS_1,\sS_2}$ denotes the partial primal-dual gap defined in \eqref{definition:partial_primduialgap}.
\end{theoremHigh}
\begin{proof}[of Theorem~\ref{theorem:conv_primdual_rate}]
By the convexity of  $G$ and the convex conjugate $F^*$, and using the primal-dual gap inequality \eqref{equation:primdua_gap_ineq}, we have for any $(\bbeta, \blambda)$,
\begin{align*}
\footnotesize
&\mathcalL(\widebarbbeta^{(T)}, \blambda) - \mathcalL(\bbeta, \widebarblambda^{(T)}) \\
&\quad = \left( \innerproduct{\bX\widebarbbeta^{(T)}, \blambda} + G(\widebarbbeta^{(T)}) - F^*(\blambda) \right) - \left( \innerproduct{\bX\bbeta, \widebarblambda^{(T)}}+ G(\bbeta) - F^*(\widebarblambda^{(T)}) \right) \\
&\quad \leq \frac{1}{T} \sum_{t=1}^T \left( \innerproduct{\bX\bbeta^\toptzero, \blambda} + G(\bbeta^\toptzero) - F^*(\blambda) \right) 
- \frac{1}{T} \sum_{t=1}^T \left( \innerproduct{\bX\bbeta, \blambda^\toptzero} + G(\bbeta) - F^*(\blambda^\toptzero) \right) \\
&\quad\leq \frac{1}{T} \left( \frac{1}{2\phi} \normtwo{\blambda - \blambda^\topzero}^2 
+
\frac{1}{2\gamma} \normtwo{\bbeta - \bbeta^\topzero}^2 
\right).
\end{align*}
Taking the supremum over all $(\bbeta, \blambda) \in \sS_1 \times \sS_2$ proves \eqref{equation:var_pridual_gap_ineq}.
\end{proof}

The inequality~\eqref{equation:var_pridual_gap_ineq} shows that the primal-dual algorithm achieves an $\mathcalO(1/T)$ convergence rate in terms of the partial primal-dual gap.
Note that this bound remains valid even when  $\sS_1\times \sS_2$ does not contain a saddle point; however, in that case the partial gap $\mathcalP_{\sS_1,\sS_2}$ may become negative and thus loses its interpretation as a meaningful error measure.

Importantly, inequality~\eqref{equation:var_pridual_gap_ineq} does not imply an $\mathcalO(1/T)$ rate for the individual sequences $\normtwobig{\bbeta^\toptzero - \widehatbbeta}$ or $\normtwobig{\blambda^\toptzero - \widehatblambda}$. Nevertheless, in practice the iterates often converge reasonably fast in these norms as well.

\begin{problemset}

\item \textbf{Hinge loss.}
Consider the one-dimensional function $f : \real \to \real$ defined by
$f(\beta) = \max\{1 - \beta, 0\}$.
Show that its convex conjugate satisfies, for any $\alpha \in \real$,
$$
f^*(\alpha) = \max_\beta \left[ \alpha\beta - \max\{1 - \beta, 0\} \right] 
=\alpha+ \delta_{[-1,0]}(\alpha).
$$

\item 
Let $f : \real \to \real$ be given by $f(\beta) = \frac{1}{s}\abs{\beta}^s$, where $s > 1$. 
Show that for any $\alpha \in \real$,
$$
f^*(\alpha) = \max_\beta \left\{ \beta \alpha - \frac{1}{s}\abs{\beta}^s \right\}
=
\frac{1}{t}\abs{\alpha}^t,
\quad \text{with $\frac{1}{s}+\frac{1}{t}=1$}.
$$

\item \textbf{Strictly convex quadratic.}
Let $ f : \real^p \to \real $ be defined by 
$$
f(\bbeta) = \frac{1}{2} \bbeta^\top \bX \bbeta + \by^\top \bbeta + z,
$$
where $ \bX$ is  strictly positive definite, $ \by \in \real^p $, and $ z \in \real $. Show that for any $ \balpha \in \real^p $,
\begin{align*}
f^*(\balpha) =
\frac{1}{2} (\balpha-\by)^\top\bX^{-1}(\balpha-\by)-z.
\end{align*}
What happens if $\bX$ is only positive semidefinite (i.e., not invertible)?

\item \label{prob:conju_indic}\textbf{Indicator function.} 
Let $ f = \indicatorS $ be the indicator function of a convex set $ \sS $ (see Exercise~\ref{exercise_convex_indica}). 
Show that its convex conjugate is given by
$$
f^*(\balpha) = \sup_{\bbeta \in \sS} \innerproduct{\bbeta, \balpha}
=\sigma_{\sS}(\balpha),
$$
where $\sigma_{\sS}(\balpha)=\max_{\bbeta\in\sS}\innerproduct{\balpha, \bbeta}$ denotes the support function of set $\sS$.

\item  \label{prob:conju_norms} \textbf{Norms.} Let $ f : \real^p \to \real $ be given by $ f(\bbeta) = \norm{\bbeta} $, where $\norm{\cdot}$ is an arbitrary norm. 
First, show that 
$$
f(\bbeta) = \sigma_{\sB_{\norm{\cdot}_*}[\bzero,1]}(\bbeta),
$$
where $\sigma_{\sS}(\balpha)=\max_{\bbeta\in\sS}\innerproduct{\balpha, \bbeta}$ denotes the support function of $\sS$, and $\norm{\cdot}_*$ denotes the dual norm. 
Then show that
$$
f^*(\balpha) = \delta_{\closure(\conv(\sB_{\norm{\cdot}_*}[\bzero,1]))}(\balpha),
$$
Since $ \sB_{\norm{\cdot}_*}[\bzero,1] $ is closed and convex, 
$ \closure(\conv(\sB_{\norm{\cdot}_*}[\bzero,1])) = \sB_{\norm{\cdot}_*}[\bzero,1] $, and therefore for any $ \balpha \in \real^p $,
$$
f^*(\balpha) = \delta_{\sB_{\norm{\cdot}_*}[\bzero,1]}(\balpha) =
\begin{cases}
0, & \norm{\balpha}_* \leq 1; \\
\infty, & \text{else}.
\end{cases}
$$

\item \label{prob:hing_svm}\textbf{Hinge loss and SVM.}\index{Support vector machine}\index{Hinge loss} 
Given data pairs $(\bx_1,y_1),(\bx_2,y_2),\ldots,(\bx_n,y_n) $, consider the soft-margin support vector machine (SVM) problem:
$$
\min_{\bbeta, b} \normtwo{\bbeta}^2 + \sum_{i=1}^{n}\max\{0, 1-y_i(\innerproduct{\bbeta, \bx_i} + b)\}.
$$
Discuss how the primal-dual algorithm can be applied to solve this optimization problem.

\item 
Describe how the primal-dual algorithm can be used to solve  the linear program problem:
$$
\min_{\bbeta} \innerproduct{\bbeta, \bc}
\quad\text{s.t.}\quad
\begin{cases}
\bX\bbeta&=\by; \\
\bbeta&\geq \bzero.
\end{cases}
$$

\item \citep{goldluecke2010convex} Potts-model over 2D grid can be written as following convex relaxation:
$$
\min_{\bbeta} \normone{\bD\bbeta} + \innerproduct{\bbeta, \bc}
\quad \text{s.t.}\quad 
0\leq \beta_i\leq 1, \quad \forall\, i, 
$$
where $\bD$ is a matrix encoding finite differences between neighboring pixels (e.g., horizontal and vertical gradients), $\normone{\bD\bbeta}$ approximates the sum of absolute differences between neighbors, and $\bc$ contains data fidelity  terms.
Discuss the application of the primal-dual algorithm to this model.
\end{problemset}

\newpage 
\chapter{Algorithms for Sparse Regression Problems}\label{chapter:spar}
\begingroup
\hypersetup{
linkcolor=structurecolor,
linktoc=page,  
}
\minitoc \newpage
\endgroup

\lettrine{\color{caligraphcolor}T}
This chapter presents a comprehensive overview of algorithms designed to solve sparse optimization problems, with a focus on the $\ell_0$-constrained and LASSO formulations. We begin by examining iterative methods for $\ell_0$-constrained problems, such as the iterative hard-thresholding (IHT) method, which directly enforces sparsity constraints. 
We then turn to constrained LASSO problems, discussing key approaches including projected gradient descent, proximal gradient methods, conditional gradient (Frank--Wolfe), and least angle regression (LARS). 
The chapter further explores Lagrangian formulations of the LASSO, covering advanced techniques like smoothing with Huber loss, FISTA, penalty methods, ADMM, and block coordinate descent. 
Finally, we extend the framework to other sparse models, including generalized LASSO, group LASSO, and sparse noise models. 
Together, these algorithms form the backbone of modern sparse recovery and feature selection methods in machine learning and signal processing.

\section{Algorithms for $\ell_0$-Constrainted Problems}\label{section:ell0_const_algo_IHT}

In Section~\ref{section:alg_cons_lasso}, we will introduce several algorithms for solving the constrained LASSO problem~\eqref{opt:lc} (p.~\pageref{opt:lc}),  including \textit{projected gradient descent}, \textit{proximal gradient method}, and \textit{conditional gradient method}. 
These methods are also applicable to general constrained optimization problems of the form $\min_{\bbeta}f(\bbeta)$ s.t. $\bbeta\in\sS$, provided that the feasible set $\sS$ is closed and convex.
However, in this section, we focus on sparse optimization problems where the constraint set is non-convex---specifically, $\ell_0$-constrained problems. Consequently, the aforementioned convex optimization techniques cannot be directly applied. Although there are some superficial similarities between the convex and non-convex settings, the theoretical guarantees and algorithmic behaviors differ significantly. To address this class of problems, we introduce the iterative hard-thresholding (IHT) algorithm.


In particular, we address the challenge of minimizing a general continuously differentiable objective function subject to a sparsity constraint, see \eqref{opt:s0_gen} (p.~\pageref{opt:s0_gen}):
\begin{equation}\label{equation:sparse_ellzero}
\text{(S$_0$)}\qquad 
\begin{aligned}
	\min_{\bbeta}  f(\bbeta) \quad \text{s.t.} \quad  \normzero{\bbeta} \leq k,
\end{aligned}
\end{equation}
where $f : \real^p \to \real$ is $L_b$-smooth (see Definition~\ref{definition:scss_func} and Theorem~\ref{theorem:equi_gradsch_smoo}), $k > 0$ is an integer satisfying $k<p$, and $\normzero{\bbeta}$ denotes the so-called ``$\ell_0$-norm" of $\bbeta$, which counts the number of nonzero components in $\bbeta$:
$$
\normzero{\bbeta} = \absbig{\{i \mid \beta_i \neq 0, \forall\, i\in\{1,2,\ldots,p\}\}}.
$$

It is important to note that the $\ell_0$-norm is not a true norm, as it fails to satisfy the homogeneity property:   for any scalar $\lambda \neq 0$ (Definition~\ref{definition:matrix-norm}), we have $\normzero{\lambda \bbeta} =  \normzero{\bbeta}$ (see Definition~\ref{definition:matrix-norm}). 
Moreover, we do not assume that $f$ is convex, and the constraint set  is inherently non-convex. 
As a result, much of the standard convex optimization theory does not apply to problem~\eqref{opt:s0_gen}. Nevertheless, this problem is of central importance in applications such as compressed sensing, motivating the development of specialized algorithms and adapted analytical tools.

We recall some useful notation. For any vector $\bbeta \in \real^p$, define its support set and its complement as
\begin{equation}
\sI_1(\bbeta)\triangleq\supp(\bbeta) \triangleq \{i \mid \beta_i \neq 0\}
\qquad\text{and}\qquad 
\sI_0(\bbeta) \triangleq \{i \mid \beta_i = 0\}.
\end{equation}
A vector $ \bbeta $ is called \textit{$ k $-sparse} if $ \abs{\supp(\bbeta)} \leq k $, or equivalently, if  $\bbeta\in \sB_0[k]$ (the closed $\ell_0$-ball of radius $k$ centered at the origin: $\sB_0[k]=\sB_0[\bzero, k]$ for brevity; see Definition~\ref{definition:open_closed_ball}). 
Using this notation, problem \eqref{opt:s0_gen} can be rewritten  as:
\begin{equation}
\min \; f(\bbeta) \qquad  \text{s.t.}\quad \bbeta \in \sB_0[k] .
\end{equation}
For a vector $\bbeta \in \real^p$ and $i \in \{1, 2, \ldots, p\}$,  let $[\bbeta]_i$ denote the \textit{$i$-th largest absolute value among its components}, so that 
\begin{equation}
	[\bbeta]_1 \geq [\bbeta]_2 \geq \ldots \geq [\bbeta]_p,
\end{equation}
where 
\begin{equation}
	[\bbeta]_1 = \max_{i=1,\ldots,p} \abs{\beta_i}
	\quad\text{and}\qquad 
	[\bbeta]_p = \min_{i=1,\ldots,p} \abs{\beta_i}.
\end{equation}

\subsection{Optimity Condition under $L$-Stationarity}\label{section:opt_l0_const}
Our goal now is to explore necessary optimality conditions for problem~\eqref{opt:s0_gen}. Because the $\ell_0$-constraint is non-convex, classical stationarity conditions from convex optimization do not generally apply. However, alternative notions of optimality---specifically tailored to the structure of sparsity constraints---can be formulated and analyzed.

We begin by extending the concept of stationarity to problem~\eqref{opt:s0_gen}, using the characterization of stationarity via the projection operator. Recall that a point $\bbeta^*$
is stationary for the problem of minimizing a continuously differentiable function over a closed and convex set $\sS$ if and only if
\begin{equation}\label{equation:lstat_proje}
	\bbeta^* = \projectS(\bbeta^* - \eta \nabla f(\bbeta^*)),
\end{equation}
where $\eta>0$ is an arbitrary positive scalar (see {Corollary~\ref{corollary:stat_point_uncons_convset_proj}}).
Although condition~\eqref{equation:lstat_proje} is expressed in terms of  $\eta$,  it is actually independent of the choice of $\eta$. 
However, this characterization does not extend directly to the case where $\sS = \sB_0[k]$, because the orthogonal projection onto non-convex sets such as $\sB_0[k]$ is not unique. Instead, the projection becomes a set-valued mapping defined by
$$
\project_{\sB_0[k]}(\bbeta) = \mathop{\argmin}_{\balpha} \{\normtwo{\balpha - \bbeta} \mid \balpha \in \sB_0[k]\}.
$$

By the Weierstrass theorem ({Theorem~\ref{theorem:weierstrass_them} \ref{weier2_prop_close_v3}}),
the existence of elements in $\project_{\sB_0[k]}(\bbeta)$ follows from the closedness of $\sB_0[k]$ and the coercivity  of the function $\normtwo{\balpha - \bbeta}$ (with respect to $\balpha$). 
To generalize condition \eqref{equation:lstat_proje} to the sparsity-constrained problem \eqref{opt:s0_gen}, we introduce the notion of ``$L$-stationarity.''

\begin{definition}[$L$-Stationarity\index{$L$-Stationarity}]\label{definition:l_stat}
A vector $\bbeta^* \in \sB_0[k]$ is called an \textit{$L$-stationary} point of problem \eqref{opt:s0_gen} if it satisfies 
\begin{equation}\label{equation:l_stationary}
\bbeta^* \in \project_{\sB_0[k]}\left(\bbeta^* - \frac{1}{L} \nabla f(\bbeta^*)\right).
\end{equation}
\end{definition}

Because the projection operator $\project_{\sB_0[k]}(\cdot)$ is set-valued, its output consists of all vectors obtained by keeping the
$k$ entries of  $\bbeta$ with the largest absolute values and setting the rest to zero. 
For example,
$$
\project_{\sB_0[2]}([3, 1, 1]^\top) = \big\{[3, 1, 0]^\top, [3, 0, 1]^\top\big\}.
$$

We now show that, under the assumption that $f $ is $L_b$-smooth, $L$-stationarity is a necessary condition for optimality whenever $L > L_b$.
Before proving this result, we first provide a more explicit characterization of $L$-stationarity.
\begin{lemma}[$L$-Stationarity]\label{lemma:l_stationary}
For any $L > 0$, a vector  $\bbeta^* \in \sB_0[k]$ satisfies $L$-stationarity if and only if
\begin{equation}\label{equation:l_stationarylem}
\abs{\frac{\partial f}{\partial \beta_i}(\bbeta^*)}
\begin{cases}
\leq L \cdot [\bbeta^*]_k, & \text{if } i \in \sI_0(\bbeta^*); \\
= 0, & \text{if } i \in \sI_1(\bbeta^*),
\end{cases}
\end{equation}
where $[\bbeta^*]_k$ denotes the $k$-th largest absolute value among the components of $\bbeta^*$.
\end{lemma}
\begin{proof}[of Lemma~\ref{lemma:l_stationary}]
Suppose that $\bbeta^*$ satisfies $L$-stationarity. 
Note that for any index $j \in \{1, 2, \ldots, p\}$, the $j$-th component of any vector in $\project_{\sB_0[k]}(\bbeta^* - \frac{1}{L} \nabla f(\bbeta^*))$ is either zero or equal to $x_j^* - \frac{1}{L} \frac{\partial f}{\partial x_j}(\bbeta^*)$. On the other hand, since $\bbeta^* \in \project_{\sB_0[k]}(\bbeta^* - \frac{1}{L} \nabla f(\bbeta^*))$, it follows that if $i \in \sI_1(\bbeta^*)$, then $\beta_i^* = \beta_i^* - \frac{1}{L} \frac{\partial f}{\partial \beta_i}(\bbeta^*)$, so that $\frac{\partial f}{\partial \beta_i}(\bbeta^*) = 0$. If $i \in \sI_0(\bbeta^*)$, then $\abs{\beta_i^* - \frac{1}{L} \frac{\partial f}{\partial \beta_i}(\bbeta^*)} \leq [\bbeta^*]_k$ and $\beta_i^* = 0$, which yields $\absbig{ \frac{\partial f}{\partial \beta_i}(\bbeta^*)} \leq L \cdot [\bbeta^*]_k$.
Thus, condition~\eqref{equation:l_stationarylem} holds.

Conversely, suppose that $\bbeta^*$ satisfies \eqref{equation:l_stationarylem}. If $\normzero{\bbeta^*} < k$, then $[\bbeta^*]_k = 0$, and by \eqref{equation:l_stationarylem} it follows that $\nabla f(\bbeta^*) = \bzero$; therefore, in this case, $\project_{\sB_0[k]}\left(\bbeta^* - \frac{1}{L} \nabla f(\bbeta^*)\right) = \project_{\sB_0[k]}(\bbeta^*)$ is the set $\{\bbeta^*\}$, which apparently satisfies the $L$-stationarity. 
On the other hand, if $\normzero{\bbeta^*} = k$, then $[\bbeta^*]_k \neq 0$ and $\abs{\sI_1(\bbeta^*)} = k$. By \eqref{equation:l_stationarylem},
$$
\abs{\beta_i^* - \frac{1}{L} \frac{\partial f}{\partial \beta_i}(\bbeta^*)}
\begin{cases}
\leq [\bbeta^*]_k, & i \in \sI_0(\bbeta^*), \\
= \abs{\beta_i^*}, & i \in \sI_1(\bbeta^*).
\end{cases}
$$
Therefore, the vector $\bbeta^* - \frac{1}{L} \nabla f(\bbeta^*)$ contains the $k$ components of $\bbeta^*$ with the largest absolute values and all other components are smaller or equal  to them in magnitude. This also implies the $L$-stationary and completes the proof.
\end{proof}

Clearly, condition~\eqref{equation:l_stationarylem} depends on the constant $L$; it becomes more restrictive as $L$ decreases.
This underscores a fundamental difference between the non-convex and convex settings. Before investigating the conditions under which $L$-stationarity constitutes a necessary optimality condition, we first establish a technical yet useful result.

\begin{lemma}[Descent lemma for sparse projection under SS]\label{lemma:des_lem_lstat}
Let  $f: \real^p\rightarrow \real$  be an $L_b$-smooth function, and let $L > L_b$. Then for any $\bbeta \in \sB_0[k]$ and any $\balpha \in \real^p$ satisfying
$
\balpha \in \project_{\sB_0[k]}\left( \bbeta - \frac{1}{L} \nabla f(\bbeta) \right),
$
we have
$$
f(\bbeta) - f(\balpha) \geq \frac{L - L_b}{2} \normtwo{\bbeta - \balpha}^2.
$$
\end{lemma}
\begin{proof}[of Lemma~\ref{lemma:des_lem_lstat}]
Since $\balpha \in \project_{\sB_0[k]}\left( \bbeta - \frac{1}{L} \nabla f(\bbeta) \right)$, 
it follows that
\begin{equation}\label{equation:des_lem_lstatpsi}
\begin{aligned}
\balpha 
&\in \mathop{\argmin}_{\bu \in \sB_0[k]} \normtwo{\bu - \left( \bbeta - \frac{1}{L} \nabla f(\bbeta) \right)}^2
=\mathop{\argmin}_{\bu \in \sB_0[k]}  \frac{L}{2} \normtwo{\bu - \left( \bbeta - \frac{1}{L} \nabla f(\bbeta) \right)}^2  
+C_{\bu}\\
&=\mathop{\argmin}_{\bu \in \sB_0[k]} 
\left\{
f(\bbeta) + \innerproduct{\nabla f(\bbeta), \bu - \bbeta} + \frac{L}{2} \normtwo{\bu - \bbeta}^2 \triangleq \Psi_L(\bu,\bbeta)
\right\} ,
\end{aligned}
\end{equation}
where $C_{\bu}\triangleq  f(\bbeta) -\frac{1}{2L} \normtwo{\nabla f(\bbeta)}^2$ denotes a constant with respect to $\bu$.
This implies that
$$
\Psi_L(\balpha, \bbeta) \leq \Psi_L(\bbeta, \bbeta) = f(\bbeta).
$$
By the $L_b$-smoothness of $f$ (Definition~\ref{definition:scss_func}), we also have 
$$
f(\bbeta) - f(\balpha) \geq f(\bbeta) - \Psi_{L_b}(\balpha, \bbeta)
=\Psi_L(\bbeta, \bbeta) - \Psi_{L_b}(\balpha, \bbeta) .
$$
Combining the preceding properties and the fact that $
\Psi_{L_b}(\balpha, \bbeta) = \Psi_L(\balpha, \bbeta) - \frac{L - L_b}{2} \normtwo{\bbeta - \balpha}^2
$ yields the desired result.
\end{proof}

We are now ready to demonstrate the necessity of  $L$-stationarity for optimality when  $L > L_b$ and $f$ is $L_b$-smooth.

\begin{theoremHigh}[Necessity under $L$-stationarity and smoothness]\label{theorem:necess_lstat}
Let  $f: \real^p\rightarrow \real$ be an $L_b$-smooth function (Definition~\ref{definition:scss_func}), and let $L > L_b$. 
If  $\bbeta^*$ is an optimal solution of problem~\eqref{opt:s0_gen}, then:
\begin{enumerate}[(i)]
\item $\bbeta^*$ is an $L$-stationary point.
\item The set $\project_{\sB_0[k]}\left(\bbeta^* - \frac{1}{L} \nabla f(\bbeta^*)\right)$ is a singleton.
\end{enumerate}
\end{theoremHigh}
\begin{proof}[of Theorem~\ref{theorem:necess_lstat}]
Suppose, for contradiction, that the projection set contains a vector $\balpha\neq \bbeta^*$, i.e.,
$\balpha \in \project_{\sB_0[k]}\left(\bbeta^* - \frac{1}{L} \nabla f(\bbeta^*)\right)$.
By Lemma~\ref{lemma:des_lem_lstat}, it follows that 
$
f(\bbeta^*) - f(\balpha) \geq \frac{L - L_b}{2} \normtwo{\bbeta^* - \balpha}^2,
$
contradicting the optimality of $\bbeta^*$. 
Thus, $\bbeta^*$ is the only vector in the set $\project_{\sB_0[k]}\left(\bbeta^* - \frac{1}{L} \nabla f(\bbeta^*)\right)$, and $\bbeta^*$ is $L$-stationary.
\end{proof}

In summary, we have shown that for any $L_b$-smooth objective function, $L$-stationarity with any $L > L_b$ is a necessary condition for optimality in problem~\eqref{opt:s0_gen}.

\subsection{The Iterative Hard-Thresholding (IHT) Method}
\begin{algorithm}[h] 
\caption{Iterative Hard-Thresholding (IHT) Method\index{Iterative hard-thresholding}}
\label{alg:pgd_iht}
\begin{algorithmic}[1] 
\Require Sparsity level $s$; 
\State {\bfseries input (a).} A differential function $f$; 
\State {\bfseries input (b).} $f(\bbeta) = \frac{1}{2n}\normtwo{\by - \bX \bbeta}$: design matrix $\bX\in\real^{n\times p}$ and response $\by\in\real^n$;
\State {\bfseries initialize:}   $\bbeta^\topone \in \sB_0[s]$;
\For{$t=1,2,\ldots$}
\State Pick a stepsize $\eta_t$;
\State $\balpha^\toptone \stackrel{(a)}{\leftarrow} \bbeta^\toptzero - \eta_t \nabla f(\bbeta^\toptzero) 
\quad \text{ or }\quad 
\stackrel{(b)}{\leftarrow} \bbeta^\toptzero - \eta_t  \frac{1}{n} \bX^\top (\bX\bbeta^\toptzero - \by)$;
\State $\bbeta^\toptone \in \mathcalP_{\sB_0[s]}(\balpha^\toptone)$; \Comment{Selecting the $s$ largest elements in magnitude of $\balpha^\toptone$}
\EndFor
\State \Return final $\bbeta\leftarrow \bbeta^\toptzero$;
\end{algorithmic} 
\end{algorithm}

We previously introduced optimality conditions for general $\ell_0$-constrained optimization problems in Section~\ref{section:opt_l0_const}:
\begin{equation}\label{equation:sparse_ellzero_algo}
\text{(S$_0$)}\qquad 
\begin{aligned}
	\min_{\bbeta}  f(\bbeta) \quad \text{s.t.} \quad  \normzero{\bbeta} \leq k,
\end{aligned}
\end{equation}
A natural approach to solving problem \eqref{opt:s0_gen} is to generalize the projected gradient descent (PGD) framework by alternating between a gradient descent step and a projection onto the non-convex sparsity set $\sB_0[k]=\{\balpha \mid \normzero{\balpha} \leq k\}$ (see Section~\ref{section:pgd} for further details). 
This yields the iterative scheme
\begin{subequations}
\begin{equation}
\bbeta^\toptone \in \project_{\sB_0[k]}\left(\bbeta^\toptzero - \frac{1}{L} \nabla f(\bbeta^\toptzero)\right), \quad t = 1, 2, \ldots.
\end{equation}
The method is known in the literature as the \textit{iterative hard-thresholding (IHT) method} (see Algorithm~\ref{alg:pgd_iht}), and we  adopt this terminology.
It applies both to general smooth functions $f$ (under the assumptions stated in  Lemma~\ref{lemma:conv_iht_smoot}) and to least squares objectives of the form $f(\bbeta) = \frac{1}{2n}\normtwo{\by - \bX \bbeta}$ (which satisfies certain design properties as previously discussed).
It can shown that the general step of the IHT method is equivalent to the relation
\begin{equation}
\bbeta^\toptone \in \mathop{\argmin}_{\bbeta \in \sB_0[k]} \Psi_L(\bbeta, \bbeta^\toptzero), \quad t = 1, 2, \ldots,
\end{equation}
\end{subequations}
where the surrogate function is defined as $\Psi_L(\bu,\bbeta)\triangleq f(\bbeta) +\innerproduct{\nabla f(\bbeta), \bu - \bbeta } + \frac{L}{2} \normtwo{\bu - \bbeta}^2$ as in \eqref{equation:des_lem_lstatpsi}.

When the constraint set $\sS$ is closed convex,  {projected gradient descent (see Algorithm~\ref{alg:pgd_gen})} can be interpreted as a fixed-point method for solving $
\bbeta^* =\mathcalP_{\sS}\big(\bbeta^* - \eta\nabla f(\bbeta^*)\big)
$ (Corollary~\ref{corollary:stat_point_uncons_convset_proj}).
Although $\sB_0[k]$ is non-convex, the IHT method can similarly be viewed as a fixed-point iteration aimed at enforcing the $L$-stationarity condition (Definition~\ref{definition:l_stat}). Specifically, under the assumption that $f$ is $L_b$-smooth and $L>L_b$, any limit point of the IHT iterates satisfies
$$
\bbeta^* =\mathcalP_{\sS}\big(\bbeta^* - \frac{1}{L}\nabla f(\bbeta^*)\big)
$$
which, by Theorem~\ref{theorem:necess_lstat}, is a necessary condition for optimality.

Several fundamental convergence properties of the IHT method are summarized in the following lemma.

\begin{lemma}[Properties of IHT under smoothness]\label{lemma:conv_iht_smoot}
Let  $f: \real^p\rightarrow \real$ be an $L_b$-smooth and lower-bounded function.
Suppose $\{\bbeta^\toptzero\}_{t > 0}$ is the sequence generated by the IHT method (Algorithm~\ref{alg:pgd_iht}(a)) with  constant stepsize $\eta\triangleq\eta_t\triangleq\frac{1}{L}$, where $L > L_b$, and a projection sparsity level $s=k$. Then,
\begin{enumerate}[(i)]
\item $f(\bbeta^\toptzero) - f(\bbeta^\toptone) \geq \frac{L - L_b}{2} \normtwobig{\bbeta^\toptzero - \bbeta^\toptone}^2$.
\item The sequence $\{f(\bbeta^\toptzero)\}_{t > 0}$ is nonincreasing.
\item $\normtwobig{\bbeta^\toptzero - \bbeta^\toptone} \to 0$ as $t\rightarrow \infty$.
\item For every $t = 1, 2, \ldots$, if $\bbeta^\toptzero \neq \bbeta^\toptone$, then $f(\bbeta^\toptone) < f(\bbeta^\toptzero)$.
\end{enumerate}
\end{lemma}
\begin{proof}[of Lemma~\ref{lemma:conv_iht_smoot}]
Part (i) follows directly  from Lemma~\ref{lemma:des_lem_lstat} by setting $\bbeta = \bbeta^\toptzero$ and $\balpha = \bbeta^\toptone$. 
Part (ii) is an immediate consequence of (i), since the right-hand side is nonnegative.
To prove (iii), observe that $\{f(\bbeta^\toptzero)\}_{t > 0}$ is nonincreasing and bounded below, hence convergent to some limit $B\in\real$. 
Therefore, $f(\bbeta^\toptzero)) - f(\bbeta^\toptone) \to B - B = 0$ as $t \to \infty$. By part (i) and the fact that $L > L_b$, it follows that $\normtwobig{\bbeta^\toptzero - \bbeta^\toptone} \to 0$. 
Finally, (iv) follows directly from (i): if $\bbeta^\toptzero\neq \bbeta^\toptone$, then the norm difference is positive, so the objective strictly decreases.
\end{proof}

As previously discussed, the IHT algorithm can be interpreted as a fixed-point method for enforcing the $L$-stationarity condition.
The following theorem establishes that every accumulation point of the sequence generated by IHT with a constant stepsize  $\frac{1}{L}$ is indeed an $L$-stationary point.

\begin{theoremHigh}[Convergence of IHT \citep{beck2014introduction}]\label{theorem:acc_iht_conv}
Assume the same conditions as in Lemma~\ref{lemma:conv_iht_smoot}.
Let $\{\bbeta^\toptzero\}_{t > 0}$ be the sequence generated by the IHT method (Algorithm~\ref{alg:pgd_iht}(a)) with  constant stepsize $\eta\triangleq\eta_t \triangleq \frac{1}{L}$, where $L > L_b$, and sparsity level $s=k$. 
Then, any accumulation point of $\{\bbeta^\toptzero\}_{t > 0}$ is an $L$-stationary point.
\end{theoremHigh}
\begin{proof}[of Theorem~\ref{theorem:acc_iht_conv}]
Let $\bbeta^*$ be an accumulation point of the sequence $\{\bbeta^\toptzero\}_{t > 0}$. 
Then there exists a subsequence $\{\bbeta^{(t_j)}\}_{j > 0}$ such that  $\bbeta^{(t_j)} \rightarrow \bbeta^*$ as $j\rightarrow \infty$. 
From  Lemma~\ref{lemma:conv_iht_smoot}, we have 
\begin{equation}\label{equation:acc_iht_conv1}
f(\bbeta^{(t_j)}) - f(\bbeta^{(t_j + 1)}) \geq \frac{L - L_b}{2} \normtwo{\bbeta^{(t_j)} - \bbeta^{(t_j + 1)}}^2.
\end{equation}
Since $\{f(\bbeta^{(t_j)})\}_{j > 0}$ and $\{f(\bbeta^{(t_j + 1)})\}_{j > 0}$, as nonincreasing and lower-bounded sequences,
they both converge to the same limit.
Consequently, $f(\bbeta^{(t_j)}) - f(\bbeta^{(t_j + 1)}) \to 0$ as $j \to \infty$, which combined with \eqref{equation:acc_iht_conv1} yields that
$
\bbeta^{(t_j + 1)} \to \bbeta^* \text{ as } j \to \infty.
$
Recall that for all $j > 0$,
$
\bbeta^{(t_j + 1)} \in \project_{\sB_0[k]}\left(\bbeta^{(t_j)} - \frac{1}{L} \nabla f(\bbeta^{(t_j)})\right).
$
We then consider the following two cases.

Let $i \in \sI_1(\bbeta^*)$ (i.e., in the support set; $\sI_1(\bbeta)\triangleq\supp(\bbeta) \triangleq \{i \mid \beta_i \neq 0\}$). By the convergence of $\bbeta^{(t_j)}$ and $\bbeta^{(t_j + 1)}$ to $\bbeta^*$, it follows that there exists an integer $J$ such that
$$
\beta_i^{(t_j)}, \beta_i^{(t_j+1)} \neq 0 \text{ for all } j > J,
$$
and therefore, for $j > J$,
$$
\beta_i^{(t_j+1)} = \beta_i^{(t_j)} - \frac{1}{L} \frac{\partial f}{\partial \beta_i}(\bbeta^{(t_j)}).
$$
Taking the limit as  $j\rightarrow \infty$ and using the continuity of $f$, we obtain that
$
\frac{\partial f}{\partial \beta_i}(\bbeta^*) = 0$ for $i \in \sI_1(\bbeta^*)$.

Now let $i \in \sI_0(\bbeta^*)$. If there exist an infinite number of indices $t_j$ for which $\beta_i^{(t_j+1)} \neq 0$, then as in the previous case, we obtain that $\beta_i^{(t_j+1)} = \beta_i^{(t_j)} - \frac{1}{L} \frac{\partial f}{\partial \beta_i}(\bbeta^{(t_j)})$ for these indices, implying (by taking the limit) that $\frac{\partial f}{\partial \beta_i}(\bbeta^*) = 0$. In particular, $\abs{\frac{\partial f}{\partial \beta_i}(\bbeta^*)} \leq L [\bbeta^*]_k$. On the other hand, if there exists an $M > 0$ such that for all $j > M$, $\beta_i^{(t_j+1)} = 0$, then
$$
\abs{\beta_i^{(t_j)} - \frac{1}{L} \frac{\partial f}{\partial \beta_i}(\bbeta^{(t_j)})} 
\leq \big[\bbeta^{(t_j)} - \frac{1}{L} \nabla f(\bbeta^{(t_j)}) \big]_k 
= \big[\bbeta^{(t_j + 1)}\big]_k.
$$
Taking the limit as $j\rightarrow \infty$, and using the continuity of the mapping $[\cdot]_k$ (which holds at points where the $k$-th and $(k+1)$-th largest magnitudes are distinct---a condition satisfied along the subsequence due to convergence to $\bbeta^*$), we obtain that
$$
\abs{\frac{\partial f}{\partial \beta_i}(\bbeta^*)} \leq L [\bbeta^*]_k,
\quad i \in \sI_0(\bbeta^*).
$$
Combining both cases, we conclude that $\bbeta^*$
satisfies the $L$-stationarity  by Lemma~\ref{lemma:l_stationary}.
\end{proof}

We now establish an important result: if the design matrix satisfies the restricted isometry property (RIP; see Definition~\ref{definition:rip22})  with suitable constants, then the IHT algorithm guarantees universal sparse recovery for the least squares objective: $f(\bbeta) = \frac{1}{2n}\normtwo{\by - \bX \bbeta}$.
\begin{theoremHigh}[Rate of convergence of IHT under RIP \citep{jain2014iterative, jain2017non}]\label{theorem:conv_iht_lasso}
Suppose $\frac{1}{\sqrt{n}}\bX \in \real^{n \times p}$ satisfies the RIP of order  $3k$ with constant $\delta_{3k} < \frac{1}{2}$. 
Let $\bbeta^* \in \sB_0[k] \subset \real^p$ be any $k$-sparse vector, and let $\by = \bX \bbeta^*$. Then the IHT algorithm (Algorithm~\ref{alg:pgd_iht}(b)), when executed with  constant stepsize $\eta_t = 1$ and sparsity level $s=k$, 
produces a sequence $\{\bbeta^\toptzero\}$ such that $\normtwobig{\bbeta^\toptzero - \bbeta^*} \leq \epsilon$ after at most $t = \mathcalO\left(\ln \frac{\normtwo{\bbeta^*}}{\epsilon}\right)$ iterations.
\end{theoremHigh}
\begin{proof}[of Theorem~\ref{theorem:conv_iht_lasso}]
Let $\sS^* \triangleq \text{supp}(\bbeta^*)$ and $\sS^t \triangleq \text{supp}(\bbeta^\toptzero)$. 
Define the index set $\sI^t \triangleq \sS^t \cup \sS^{t+1} \cup \sS^*$, which contains the supports of two consecutive iterates and the true support. 
Since each set has cardinality at most $k$, we have $\abs{\sI^t} \leq 3k$. 
Crucially, both error vectors $\bbeta^\toptzero - \bbeta^*$ and $\bbeta^\toptone - \bbeta^*$ are supported on $\sI^t$.

With stepsize $\eta_t = 1$, the gradient step in Algorithm~\ref{alg:pgd_iht}(b) yields $\balpha^\toptone = \bbeta^\toptzero - \frac{1}{n} \bX^\top (\bX \bbeta^\toptzero - \by)$.
Although $\sB_0[k]$ is non-convex, the projection $\bbeta^\toptone = \mathcalP_{\sB_0[k]}(\balpha^\toptone)$ satisfies the basic optimality property (Projection Property-O; see Lemma~\ref{lemma:proj_prop0}):
$$
\begin{aligned}
&\normtwo{\bbeta^\toptone - \balpha^\toptone}^2 \leq \normtwo{\bbeta^* - \balpha^\toptone}^2\\
&\implies \quad \normtwo{\bbeta_{\sI}^\toptone - \balpha_{\sI}^\toptone}^2 + \normtwo{\bbeta_{\comple{\sI}}^\toptone - \balpha_{\comple{\sI}}^\toptone}^2 \leq \normtwo{\bbeta_{\sI}^* - \balpha_{\sI}^\toptone}^2 + \normtwo{\bbeta_{\comple{\sI}}^* - \balpha_{\comple{\sI}}^\toptone}^2,
\end{aligned}
$$
where, for brevity, we denote $\sI \triangleq\sI^t$ such that $\bbeta_{\comple{\sI}}^\toptone = \bbeta_{\comple{\sI}}^* = \bzero$. 
Using the fact that $\by = \bX \bbeta^*$, and denoting $\widetildebX \triangleq \frac{1}{\sqrt{n}} \bX$ and $\widetildebX_{\sI} \triangleq\widetildebX[:,\sI]$, we obtain
$$
\begin{aligned}
\normtwo{\bbeta^\toptone - \bbeta^* } 
&=\normtwo{\bbeta_{\sI}^\toptone - \bbeta_{\sI}^*} 
\leq \normtwo{\bbeta_{\sI}^\toptone - \balpha_{\sI}^\toptone}
+\normtwo{\balpha_{\sI}^\toptone- \bbeta_{\sI}^*}
\leq 2\normtwo{\balpha_{\sI}^\toptone- \bbeta_{\sI}^*}\\
&=2\normtwo{
	\left( \bbeta_{\sI}^\toptzero - \widetildebX_{\sI}^\top \widetildebX (\bbeta^\toptzero - \bbeta^*) \right)
	- \bbeta_{\sI}^*
	}
\stackrel{\dag}{=}
2\normtwo{ (\bI - \widetildebX_{\sI}^\top \widetildebX_{\sI}) (\bbeta_{\sI}^\toptzero - \bbeta_{\sI}^*) }\\
&\leq 2\normtwo{\widetildebX_{\sI}^\top \widetildebX_{\sI}-\bI  } \normtwo{\bbeta_{\sI}^\toptzero - \bbeta_{\sI}^*}
\leq 2 \delta_{3k} \normtwo{\bbeta^\toptzero - \bbeta^*},
\end{aligned}
$$
where the equality ($\dag$) follows since  $\bbeta_{\comple{\sI}}^\toptzero = \bbeta_{\comple{\sI}}^* = \bzero$, 
and the last inequality follows from the RIP condition (Definition~\ref{definition:rip22} or Remark~\ref{remark:rip_iht}).
To achieve $\normtwobig{\bbeta^\toptzero - \bbeta^*} \leq \epsilon$, it suffices to take
$t= \mathcalO\left(\ln \frac{\normtwo{\bbeta^*}}{\epsilon}\right)$ with initialization $\bbeta^\topone=\bzero$, which completes the proof.
\end{proof}

The above theorem establishes the convergence rate of the IHT algorithm for the least squares objective.
To extend this result to a general smooth objective function, we first require the following lemma concerning the sparsity projection operator.
\begin{lemma}[Sparsity projection property]\label{lemma:spar_proj}
Let $\bbeta \in \real^p$, and let  $\widetildebbeta \triangleq \mathcal\project_{\sB_0[k]}(\bbeta)$ denote its best $k$-sparse approximation. 
Then, for any $k^*$-sparse vector $\bbeta^* \in \sB_0[k^*]\subseteq \real^p$  with $k^*\leq k\leq p$, the following inequality holds:
$$
\frac{\normtwobig{\widetildebbeta - \bbeta}^2}{p-k} \leq  \frac{\normtwo{\bbeta^* - \bbeta}^2}{p-k^*}
\qquad\implies\qquad 
\normtwobig{\widetildebbeta - \bbeta}^2 \leq \frac{p - k}{p - k^*} \normtwo{\bbeta^* - \bbeta}^2.
$$
\end{lemma}
\begin{proof}[of Lemma~\ref{lemma:spar_proj}]
Without loss of generality, assume the coordinates of $\bbeta$ are ordered so that  $\abs{\beta_1} \geq \abs{\beta_2} \geq \ldots \geq \abs{\beta_p}$. 
Since the projection operator $\mathcal\project_{\sB_0[k]}(\bbeta)$ operates by selecting the largest elements by magnitude, we have $\widetildebeta_1 = \beta_1, \widetildebeta_2 = \beta_2 \ldots, \widetildebeta_k = \beta_k$ and $\widetildebeta_{k+1} = \widetildebeta_{k+2} = \ldots = \widetildebeta_p = 0$.

Similarly, define $\widehatbbeta \triangleq \mathcalP_{\sB_0[k^*]}(\bbeta)$
so that   $\widehatbeta_1 = \beta_1, \widehatbeta_2 = \beta_2, \ldots, \widehatbeta_{k^*} = \beta_{k^*}$ and $\widehatbeta_{k^*+1} = \widehatbeta_{k^*+2} = \ldots = \widehatbeta_p = 0$. 
Since $k^*\leq k$, we compute 
\begin{equation}
\begin{aligned}
(p-k)&(p-k^*)\bigg(\frac{\normtwo{\widehatbbeta - \bbeta}^2}{p - k^*} - \frac{\normtwobig{\widetildebbeta - \bbeta}^2}{p - k} \bigg)
=(p-k) \sum_{i=k^*+1}^{p} \beta_i^2 - (p-k^*)\sum_{i=k+1}^{p} \beta_i^2\\
&=(p-k) \sum_{i=k^*+1}^{\textcolor{mylightbluetext}{k}} \beta_i^2 + \big((p-k)-(p-k^*)\big)\sum_{i=k+1}^{p} \beta_i^2\\
&\geq (p-k)(k-k^*)  \beta_{k}^2 +  (k^*-k)(p-k) \beta_{k+1}^2 \geq 0.\\
\end{aligned}
\end{equation}
Finally, by the optimality of the projection (the Projection Property-O; see Lemma~\ref{lemma:proj_prop0}), we have $ \normtwobig{\widehatbbeta - \bbeta} \leq \normtwo{\bbeta^* - \bbeta} $ for any $\bbeta^* \in \sB_0[k^*]$, which establishes the desired result when combining with the above inequality.
\end{proof}

For a more general objective function $f$ (not necessarily least squares), suppose there exists an optimal $k^*$-sparse parameter $\bbeta^*=\argmin_{\normzero{\bbeta}\leq k^*} f(\bbeta)$.
Our analysis combines the above projection property with the restricted strong convexity (RSC) and restricted smoothness (RSS) assumptions on $f$ to establish geometric convergence of IHT.

\begin{theoremHigh}[Rate of convergence of IHT under RSC \& RSS \citep{jain2014iterative}]\label{theorem:conv_iht_genfunc}
Assume $f$ satisfies  RSC and RSS with parameters  {$\mu_{2k+k^*}(f) = L_a>0$} and {$\nu_{2k+k^*}(f) = L_b>0$}, respectively (Definition~\ref{definition:res_scss_func} or Definition~\ref{definition:res_scss_mat}). 
Let Algorithm~\ref{alg:pgd_iht} be invoked with $f$, $s\triangleq k \geq 32 \left(\frac{L_b}{L_a}\right)^2 k^*$, and a constant stepsize $\eta_t = \frac{2}{3L_b}$. Also let $\bbeta^* = \argmin_{\normzero{\bbeta} \leq k^*} f(\bbeta)$. Then, the $T$-th iterate of Algorithm~\ref{alg:pgd_iht}(a), for $T = \mathcalO\left(\frac{L_b}{L_a} \cdot \ln\left(\frac{f(\bbeta^\topone)}{\epsilon}\right)\right)$ satisfies:
$ f(\bbeta^{(T)}) - f(\bbeta^*) \leq \epsilon.
$
\end{theoremHigh}

\begin{proof}[of Theorem~\ref{theorem:conv_iht_genfunc}]
Denote $\bg^\toptzero \triangleq \nabla f(\bbeta^\toptzero)$, recall that $\bbeta^\toptone = \mathcal\project_{\sB_0[k]}(\bbeta^\toptzero - \frac{\gamma}{L_b} \bg^\toptzero)$ where $\gamma = \frac{2}{3} < 1$. Let $\sS^t \triangleq \supp(\bbeta^\toptzero)$, $\sS^* \triangleq \supp(\bbeta^*)$, and $\sS^{t+1} \triangleq \supp(\bbeta^\toptone)$. Also, let $\sI \triangleq \sS^* \cup \sS^t \cup \sS^{t+1}$.
Given the RSS property (Definition~\ref{definition:res_scss_func}) and the fact that $\supp(\bbeta^\toptzero) \subseteq \sI$ and $\supp(\bbeta^\toptone) \subseteq \sI$ (implying that $\bbeta_{\comple{\sI}}^\toptone = \bbeta_{\comple{\sI}}^\toptzero = \bzero$), along with the fact that $\normtwo{\bu}^2 = \normtwo{\bu_{\sI}}^2 + \normtwo{\bu_{\comple{\sI}}}^2$ for any vector $\bu\in\real^p$, we have:
\begin{equation}\label{equation:conv_iht_genfunc1}
\begin{aligned}
&f(\bbeta^\toptone) - f(\bbeta^\toptzero) \leq \innerproduct{\bg^\toptzero, \bbeta^\toptone - \bbeta^\toptzero} + \frac{L_b}{2} \normtwo{\bbeta^\toptone - \bbeta^\toptzero}^2, \\
&= \frac{L_b}{2} \normtwo{\bbeta^\toptone_{\sI} - \bbeta^\toptzero_{\sI} + \frac{\gamma}{L_b} \cdot  \bg^\toptzero_{\sI}}^2 - \frac{\gamma^2}{2L_b} \normtwo{\bg^\toptzero_{\sI}}^2 + (1 - \gamma) \innerproduct{\bbeta^\toptone - \bbeta^\toptzero, \bg^\toptzero}.
\end{aligned}
\end{equation}
Since $\sS^t \setminus \sS^{t+1}$ and $\sS^{t+1}$ are disjoint, we have:
\begin{equation}\label{equation:conv_iht_genfunc2}
\small
\begin{aligned}
&\innerproduct{\bbeta^\toptone - \bbeta^\toptzero, \bg^\toptzero} 
= \innerproduct{\,\cancel{\bbeta^\toptone_{\sS^t \setminus \sS^{t+1}}} - \bbeta^\toptzero_{\sS^t \setminus \sS^{t+1}}, \bg^\toptzero_{\sS^t \setminus \sS^{t+1}}} 
+ \innerproduct{\bbeta^\toptone_{\sS^{t+1}} - \bbeta^\toptzero_{\sS^{t+1}}, \bg^\toptzero_{\sS^{t+1}}} \\
&\overset{\dag}{=} - \innerproduct{\bbeta^\toptzero_{\sS^t \setminus \sS^{t+1}}, \bg^\toptzero_{\sS^t \setminus \sS^{t+1}}} - \frac{\gamma}{L_b} \normtwo{\bg^\toptzero_{\sS^{t+1}}}^2 \\
&\overset{\ddag}{\leq} \frac{\gamma}{2L_b} \normtwo{\bg^\toptzero_{\sS^{t+1} \setminus \sS^t}}^2 - \frac{\gamma}{2L_b} \normtwo{\bg^\toptzero_{\sS^t \setminus \sS^{t+1}}}^2 - \frac{\gamma}{L_b} \normtwo{\bg^\toptzero_{\sS^{t+1}}}^2   \cancel{-\left( \frac{L_b}{2\gamma} \normtwo{\bbeta_{\sS^t \setminus \sS^{t+1}}}^2\right)} \\
&\overset{*}{=} - \frac{\gamma}{2L_b} \normtwo{\bg^\toptzero_{\sS^{t+1} \setminus \sS^t}}^2 - \frac{\gamma}{2L_b} \normtwo{\bg^\toptzero_{\sS^t \setminus \sS^{t+1}}}^2 - \frac{\gamma}{L_b} \normtwo{\bg^\toptzero_{\sS^t \cap \sS^{t+1}}}^2
\leq - \frac{\gamma}{2L_b} \normtwo{\bg^\toptzero_{\sS^t \cup \sS^{t+1}}}^2,
\end{aligned}
\end{equation}
where the equality $(\dag)$ follows from the gradient step, i.e., $\bbeta^\toptone_{\sS^{t+1}} = \bbeta^\toptzero_{\sS^{t+1}} - \frac{\gamma}{L_b} \bg^\toptzero_{\sS^{t+1}}$, the inequality $(\ddag)$ follows using the fact that $\bbeta^\toptone$ is obtained using hard-thresholding ($\sS^{t+1}$ contains the indices of largest $\abs{\sS^{t+1}}=\abs{\sS^{t}}$ magnitudes) and the fact that $\abs{\sS^t \setminus \sS^{t+1}} = \abs{\sS^{t+1} \setminus \sS^t}$, such that 
$
\normtwobig{\bbeta_{\sS^t \setminus \sS^{t+1}}^\toptzero - \frac{\gamma}{L_b} \bg_{\sS^t \setminus \sS^{t+1}}^\toptzero}^2 
\leq 
\normtwobig{\bbeta_{\sS^{t+1} \setminus \sS^t}^\toptone}^2 
= \frac{\gamma^2}{L_b^2} \normtwobig{\bg_{\sS^{t+1} \setminus \sS^t}^\toptzero}^2,
$
and the equality $(*)$ follows from $\normtwobig{\bg_{\sS^{t+1}}^\toptzero}^2 = \normtwobig{\bg_{\sS^{t+1} \setminus \sS^t}^\toptzero}^2 + \normtwobig{\bg_{\sS^t \cap \sS^{t+1}}^\toptzero}^2$.
Combining \eqref{equation:conv_iht_genfunc1} and \eqref{equation:conv_iht_genfunc2} yields
\begin{equation}\label{equation:conv_iht_genfunc3}
\begin{aligned}
&\quad f(\bbeta^\toptone) - f(\bbeta^\toptzero) \\
&\leq \frac{L_b}{2} \normtwo{\bbeta_{\sI}^\toptone - \bbeta_{\sI}^\toptzero + \frac{\gamma}{L_b} \bg_{\sI}^\toptzero}^2 - \frac{\gamma^2}{2L_b} \normtwo{\bg_{\sI}^\toptzero}^2 - \frac{\gamma(1-\gamma)}{2L_b} \normtwo{\bg_{\sS^t \cup \sS^{t+1}}^\toptzero}^2 \\
&= \frac{L_b}{2} \normtwo{\bbeta_{\sI}^\toptone - \bbeta_{\sI}^\toptzero + \frac{\gamma}{L_b} \bg_{\sI}^\toptzero}^2 
- \frac{\gamma^2}{2L_b} \normtwo{\bg_{\sI \setminus (\sS^t \cup \sS^*)}^\toptzero}^2 \\ 
&\qquad \quad - \frac{\gamma^2}{2L_b} \normtwo{\bg_{\sS^t \cup \sS^*}^\toptzero}^2 
- \frac{\gamma(1-\gamma)}{2L_b} \normtwo{\bg_{\sS^t \cup \sS^{t+1}}^\toptzero}^2.
\end{aligned}
\end{equation}
Next, we will upper bound the first two terms on the right-hand side  of the above inequality. Since $\sI \setminus (\sS^t \cup \sS^*) = \sS^{t+1} \setminus (\sS^t \cup \sS^*) \subseteq \sS^{t+1}$, we have $\bbeta_{\sI \setminus (\sS^t \cup \sS^*)}^\toptone = \bbeta_{\sI \setminus (\sS^t \cup \sS^*)}^\toptzero - \frac{\gamma}{L_b} \bg_{\sI \setminus (\sS^t \cup \sS^*)}^\toptzero$. 
However, since $\bbeta_{\sI \setminus \sS^t}^\toptzero = \bzero$, the preceding equality reduces to $\bbeta_{\sI \setminus (\sS^t \cup \sS^*)}^\toptone = -\frac{\gamma}{L_b} \bg_{\sI \setminus (\sS^t \cup \sS^*)}^\toptzero$. 
Let $\sT \subseteq \sS^t \setminus \sS^{t+1}$ such that $\abs{\sT} = \abs{\sS^{t+1} \setminus (\sS^t \cup \sS^*)}$. Such a choice is possible since $\abs{\sS^{t+1} \setminus (\sS^t \cup \sS^*)} = \abs{\sS^t \setminus \sS^{t+1}} - \abs{(\sS^{t+1} \cap \sS^*) \setminus \sS^t}$ (which itself is a consequence of the fact that $\abs{\sS^{t+1}} = \abs{\sS^t}$ such that $\abs{\sS^t \setminus \sS^{t+1}} = \abs{\sS^{t+1} \setminus \sS^t}$). Moreover, since $\bbeta^\toptone$ is obtained by hard-thresholding $\left(\bbeta^\toptzero - \frac{\gamma}{L_b} \bg^\toptzero\right)$, for any choice of $\sT$ made above, we have:
\begin{equation}
\frac{\gamma^2}{L_b^2} \normtwo{\bg_{\sS^{t+1} \setminus (\sS^t \cup \sS^*)}^\toptzero}^2 
= \normtwo{\bbeta_{\sS^{t+1} \setminus (\sS^t \cup \sS^*)}^\toptone}^2 \geq \normtwo{\bbeta_\sT^\toptzero - \frac{\gamma}{L_b} \bg_\sT^\toptzero}^2.
\end{equation}
Using the above equation,  the fact that $\bbeta_\sT^\toptone = \bzero$ (since $\sT \nsubseteq \sS^{t+1}$), and the fact that $\sI \setminus (\sS^t \cup \sS^*) = \sS^{t+1} \setminus (\sS^t \cup \sS^*)$, the first two terms of equality~\eqref{equation:conv_iht_genfunc3} becomes 
\begin{equation}\label{equation:conv_iht_genfunc4}
\small
\begin{aligned}
&\quad \frac{L_b}{2} \normtwo{\bbeta_{\sI}^\toptone - \bbeta_{\sI}^\toptzero + \frac{\gamma}{L_b} \bg_{\sI}^\toptzero}^2 - \frac{\gamma^2}{2L_b} \normtwo{\bg_{\sI \setminus (\sS^t \cup \sS^*)}^\toptzero}^2 \\
&\leq \frac{L_b}{2} \normtwo{\bbeta_{\sI}^\toptone - \bbeta_{\sI}^\toptzero + \frac{\gamma}{L_b} \bg_{\sI}^\toptzero}^2 - \frac{L_b}{2} \normtwo{\bbeta_\sT^\toptone - \bbeta_\sT^\toptzero + \frac{\gamma}{L_b} \bg_\sT^\toptzero}^2 \\
&= \frac{L_b}{2} \normtwo{\bbeta_{\sI \setminus \sT}^\toptone - \bbeta_{\sI \setminus \sT}^\toptzero + \frac{\gamma}{L_b} \bg_{\sI \setminus \sT}^\toptzero}^2.
\end{aligned}
\end{equation}
Since we construct the set $\sT$ such that $\sT \subseteq \sS^t \setminus \sS^{t+1}$ satisfying $\abs{\sT} = \abs{\sS^{t+1} \setminus (\sS^t \cup \sS^*)} = \abs{\sI \setminus (\sS^t \cup \sS^*)}$, we can bound the size of $\sI \setminus \sT$ as $\abs{\sI \setminus \sT} \leq \abs{\sS^{t}} + \abs{(\sS^t \setminus \sS^{t+1}) \setminus \sT} + \abs{\sS^*} \leq k + \abs{(\sS^{t+1} \cap \sS^*) \setminus \sS^t} + k^* \leq k + 2k^*$. Also, since $\sS^{t+1} \subseteq (\sI \setminus \sT)$, we have $\bbeta_{\sI \setminus \sT}^\toptone = \project_{\sB_0[k]} \left(\bbeta_{\sI \setminus \sT}^\toptzero - \frac{\gamma}{L_b} \bg_{\sI \setminus \sT}^\toptzero\right)$.

Given the above observation with \eqref{equation:conv_iht_genfunc4} and the assumption that $k \geq 32 \left(\frac{L_b}{L_a}\right)^2 k^*$, Lemma~\ref{lemma:spar_proj} is valid with $\widetildebbeta \triangleq \bbeta_{\sI \setminus \sT}^\toptone + \frac{\gamma}{L_b} \bg_{\sI \setminus \sT}^\toptzero$ and $\bbeta \triangleq\bbeta_{\sI \setminus \sT}^\toptzero $, whence we reduce the first two terms of equality~\eqref{equation:conv_iht_genfunc3} to
\begin{equation}\label{equation:conv_iht_genfunc5}
\small
\begin{aligned}
&\frac{L_b}{2} \normtwo{\bbeta_{\sI}^\toptone - \bbeta_{\sI}^\toptzero + \frac{\gamma}{L_b} \bg_{\sI}^\toptzero}^2 - \frac{\gamma^2}{2L_b} \normtwo{\bg_{\sI \setminus (\sS^t \cup \sS^*)}^\toptzero}^2 
\leq \frac{L_b}{2} \cdot \frac{\abs{\sI \setminus \sT} - k}{\abs{\sI \setminus \sT} - k^*} \normtwo{\bbeta_{\sI \setminus \sT}^* - \bbeta_{\sI \setminus \sT}^\toptzero + \frac{\gamma}{L_b} \bg_{\sI \setminus \sT}^\toptzero}^2 \\
&\overset{\dag}{\leq} \frac{L_b}{2} \cdot \frac{2k^*}{k + k^*} \normtwo{\bbeta_{\sI}^* - \bbeta_{\sI}^\toptzero + \frac{\gamma}{L_b} \bg_{\sI}^\toptzero}^2 
= \frac{2k^*}{k + k^*} \cdot \left(\gamma \innerproduct{\bbeta^* - \bbeta^\toptzero, \bg^\toptzero} + \frac{L_b}{2} \normtwo{\bbeta^* - \bbeta^\toptzero}^2 + \frac{\gamma^2}{2L_b} \normtwo{\bg_{\sI}^\toptzero}^2\right) \\
&\overset{\ddag}{\leq} \frac{2k^*}{k + k^*} \cdot \left(\gamma f(\bbeta^*) - \gamma f(\bbeta^\toptzero) + \frac{L_b - \gamma L_a}{2} \normtwo{\bbeta^* - \bbeta^\toptzero}^2 + \frac{\gamma^2}{2L_b} \normtwo{\bg_{\sI}^\toptzero}^2\right),
\end{aligned}
\end{equation}
where the inequality $(\dag)$ follows from $\abs{\sI \setminus \sT} \leq k + 2k^*$ as shown earlier and the observation that $\frac{z-a}{z-b}$ is a positive and increasing function on the interval $z \geq a$ if $a \geq b \geq 0$. Note that since we already have $\sS^{t+1} \subseteq (\sI \setminus \sT)$, we get $\abs{\sI \setminus \sT} \geq k$. The inequality $(\ddag)$ follows by the RSC property of functions (Definition~\ref{definition:res_scss_func}).

Plugging \eqref{equation:conv_iht_genfunc5} into the first two terms of equality~\eqref{equation:conv_iht_genfunc3}, and using the fact that $\sS^{t+1} \setminus (\sS^t \cup \sS^*) \subseteq (\sS^{t+1} \cup \sS^t )$, we get:
\begin{equation}
\small
\begin{aligned}
f(\bbeta^\toptone) - f(\bbeta^\toptzero) 
&\leq \frac{2k^*}{k + k^*} \cdot \left( \gamma f(\bbeta^*) - \gamma f(\bbeta^\toptzero) + \frac{L_b - \gamma L_a}{2} \normtwo{\bbeta^* - \bbeta^\toptzero}^2 + \frac{\gamma^2}{2L_b} \normtwo{\bg_\sI^\toptzero}^2 \right) \\
&\quad - \frac{\gamma^2}{2L_b} \normtwo{\bg_{\sS^t \cup \sS^*}^\toptzero}^2 - \frac{\gamma(1-\gamma)}{2L_b} \normtwo{\bg_{\sS^{t+1} \setminus (\sS^t \cup \sS^*)}^\toptzero}^2.
\end{aligned}
\end{equation}
Since $\gamma = 2/3$ and $k \geq 32 \left( \frac{L_b}{L_a} \right)^2 k^*$, so that we have $\frac{2k^*}{k + k^*} \leq \frac{L_a^2}{16L_b(L_b - \gamma L_a)}$. Since $L_b \geq L_a$, we also have $\frac{L_a^2}{16L_b(L_b - \gamma L_a)} \leq \frac{3}{16}$. Using these inequalities, rearrange the terms in \eqref{equation:conv_iht_genfunc5} above:
\begin{equation}
\begin{aligned}
f(\bbeta^\toptone) - f(\bbeta^\toptzero) &\leq \frac{2k^*}{k + k^*} \cdot \gamma \cdot \left( f(\bbeta^*) - f(\bbeta^\toptzero) \right) + \frac{L_a^2}{32L_b} \normtwo{\bbeta^* - \bbeta^\toptzero}^2 + \frac{1}{24L_b} \normtwo{\bg_\sI^\toptzero}^2 \\
&\quad - \frac{2}{9L_b} \normtwo{\bg_{\sS^t \cup \sS^*}^\toptzero}^2 - \frac{1}{9L_b} \normtwo{\bg_{\sS^{t+1} \setminus (\sS^t \cup \sS^*)}^\toptzero}^2.
\end{aligned}
\end{equation}
Decomposing $\normtwobig{\bg_\sI^\toptzero}^2 = \normtwobig{\bg_{\sS^t \cup \sS^*}^\toptzero}^2 + \normtwobig{\bg_{\sS^{t+1} \setminus (\sS^t \cup \sS^*)}^\toptzero}^2$ implies
\begin{equation}\label{equation:conv_iht_genfunc6}
\small
\begin{aligned}
f(\bbeta^\toptone) - f(\bbeta^\toptzero) &\leq \frac{2k^*}{k + k^*} \cdot \gamma \cdot \left( f(\bbeta^*) - f(\bbeta^\toptzero) \right) - \frac{1}{2L_b} \left( \frac{13}{36} \normtwo{\bg_{\sS^t \cup \sS^*}^\toptzero}^2 - \frac{L_a^2}{16} \normtwo{\bbeta^* - \bbeta^\toptzero}^2 \right) \\
&\quad - \frac{1}{2L_b} \cdot \left( \frac{4}{9} - \frac{1}{12} \right) \normtwo{\bg_{\sS^{t+1} \setminus (\sS^t \cup \sS^*)}^\toptzero}^2 \\
&\leq \frac{2k^*}{k + k^*} \cdot \gamma \cdot \left( f(\bbeta^*) - f(\bbeta^\toptzero) \right) - \frac{13}{72L_b} \left( \normtwo{\bg_{\sS^t \cup \sS^*}^\toptzero}^2 - \frac{L_a^2}{4} \normtwo{\bbeta^* - \bbeta^\toptzero}^2 \right). \\
\end{aligned}
\end{equation}
Using the RSC property, we have:
$$
\small
\begin{aligned}
f(\bbeta^\toptzero) - f(\bbeta^*) 
&\leq \innerproduct{\bg^\toptzero, \bbeta^\toptzero - \bbeta^*}- \frac{L_a}{2} \normtwo{\bbeta^* - \bbeta^\toptzero}^2 
= \innerproduct{\bg^\toptzero_{\sS^t \cup \sS^*}, \bbeta^\toptzero_{\sS^t \cup \sS^*} - \bbeta^*_{\sS^t \cup \sS^*}} - \frac{L_a}{2} \normtwo{\bbeta^* - \bbeta^\toptzero}^2 \\
&\leq \normtwo{\bg^\toptzero_{\sS^t \cup \sS^*}} \normtwo{\bbeta^\toptzero - \bbeta^*} - \frac{L_a}{2} \normtwo{\bbeta^* - \bbeta^\toptzero}^2.
\end{aligned}
$$
This shows 
$$
\small
\begin{aligned}
\normtwo{\bg_{\sS^t \cup \sS^*}^\toptzero}^2 
&- \frac{L_a^2}{4} \normtwo{\bbeta^* - \bbeta^\toptzero}^2
= \left( \normtwo{\bg_{\sS^t \cup \sS^*}^\toptzero} - \frac{L_a}{2} {\normtwo{\bbeta^* - \bbeta^\toptzero }} \right) \left( \normtwo{\bg_{\sS^t \cup \sS^*}^\toptzero} + \frac{L_a}{2} \normtwo{\bbeta^* - \bbeta^\toptzero} \right) \\
&\geq \frac{(f(\bbeta^\toptzero) - f(\bbeta^*))}{\normtwo{\bbeta^\toptzero - \bbeta^*}} \cdot \left( \normtwo{\bg_{\sS^t \cup \sS^*}^\toptzero} + \frac{L_a}{2} \normtwo{\bbeta^* - \bbeta^\toptzero} \right)
\geq \frac{L_a}{2} \cdot (f(\bbeta^\toptzero) - f(\bbeta^*)).
\end{aligned}
$$
Plugging this inequality into \eqref{equation:conv_iht_genfunc6} yields 
\begin{equation}\label{equation:conv_iht_genfunc7}
\small
\begin{aligned}
f(\bbeta^\toptone) - f(\bbeta^\toptzero) 
&\leq \frac{2k^*}{k + k^*} \cdot \gamma \cdot \left( f(\bbeta^*) - f(\bbeta^\toptzero) \right) - \frac{L_a}{12L_b} \left( f(\bbeta^\toptzero) - f(\bbeta^*) \right)\\
&=-\left(\frac{2k^*}{k + k^*} \cdot \gamma  + \frac{L_a}{12L_b}\right) \left(  f(\bbeta^\toptzero) - f(\bbeta^*)  \right) .
\end{aligned}
\end{equation}
The result  follows by observing that $\frac{2k^*}{k + k^*} \geq 0$.
\end{proof}

\section{Algorithms for Constrained LASSO Problems}\label{section:alg_cons_lasso}
Instead of enforcing an $\ell_0$-constraint---as in the IHT method---the original formulation of the constrained LASSO \eqref{opt:lc} (p.~\pageref{opt:lc}) is given by:
\begin{equation}\label{equation:loss_cons_lasso_scaled_in_algo} 
\mathopmin{\bbeta\in\real^p} f(\bbeta)\triangleq \frac{1}{2}\normtwo{\by - \bX \bbeta}^2
\quad \text{s.t.}\quad \{\normone{\bbeta}\leq \Sigma\} \equiv \{\bbeta \in \sB_{1}[\Sigma]\},
\end{equation}
where $\sB_{1}[\Sigma]\triangleq \{ \balpha \in \real^p \mid \normone{\balpha} \leq \Sigma \}$ denotes the closed $\ell_1$-ball of radius $\Sigma$.
Since $\sB_{1}[\Sigma]$ is convex, a natural approach to solving this problem is \textit{projected gradient descent (PGD)}, which alternates between a gradient descent step and a projection onto the feasible set $\sB_{1}[\Sigma]$.

\begin{algorithm}[h] 
\caption{Projected (Sub)Gradient Descent (PGD) Method\index{Projected gradient descent}}
\label{alg:pgd_gen}
\begin{algorithmic}[1] 
\Require A function $f(\bbeta)$ and a set $\sS$; 
\State {\bfseries initialize:}  $\bbeta^\topone$;
\For{$t=1,2,\ldots$}
\State Choose a stepsize $\eta_t$;
\State $\balpha^\toptone \leftarrow \bbeta^\toptzero - \eta_t \nabla f(\bbeta^\toptzero)$ or $\balpha^\toptone \leftarrow \bbeta^\toptzero - \eta_t f^\prime(\bbeta^\toptzero)$, where $ f^\prime(\bbeta^\toptzero)\in \partial f(\bbeta^\toptzero)$;
\State $\bbeta^\toptone \in \mathcalP_{\sS}(\balpha^\toptone)$;
\State Stop if a stopping criterion is satisfied at iteration $t=T$;
\EndFor
\State (Output Option 1) Output  $\bbeta_{\text{final}}\leftarrow \bbeta^{(T)}$;
\State (Output Option 2) Output  $\bbeta_{\text{avg}}\leftarrow \frac{1}{T}(\sum_{t=1}^{t}\bbeta^\toptzero)$ or $\sum_{t=1}^{T} \frac{2t}{T(T+1)} \bbeta^\toptzero$;
\State (Output Option 3) Output  $\bbeta_{\text{best}}\leftarrow \argmin_{t\in\{1,2,\ldots,T\}} f(\bbeta^\toptzero)$;
\end{algorithmic} 
\end{algorithm}

\subsection{Projected Gradient Descent (PGD)}\label{section:pgd}

Standard gradient descent is designed for unconstrained optimization problems. To handle constraints, the projected gradient descent (PGD) method---described in Algorithm~\ref{alg:pgd_gen}---solves problems of the form:
\begin{equation}\label{equation:p2_pgd}
\text{(C1)}: \qquad \min_{\bbeta\in\real^p } f(\bbeta) \quad \text{s.t.}\quad \bbeta\in\sS.
\end{equation}
Section~\ref{section:constr_convset} discusses optimality conditions for such constrained problems.
PGD is the most straightforward first-order method for solving (C1) and serves as a foundational building block in constrained optimization. Many advanced algorithms---including the proximal gradient method---are natural extensions of PGD; see Chapter~4 of \citet{lu2025practical} for further discussion.

When $\sS$ is closed and convex and the stepsize $\eta_t\equiv \eta$ is constant, a fixed point $\bbeta^*$ of the PGD update satisfies the stationarity condition:
$$
\bbeta^* =\mathcalP_{\sS}(\bbeta^* - \eta\nabla f(\bbeta^*))
$$
as stated in Corollary~\ref{corollary:stat_point_uncons_convset_proj}.
Because PGD closely resembles the proximal gradient method, we will analyze its convergence within the proximal framework in the following sections.

More precisely, the PGD update from iteration $t$ to $t+1$ can be equivalently written as:
\begin{tcolorbox}[colback=white,colframe=black]
\begin{minipage}{1\textwidth}
\begin{equation}\label{equation:pgd_decom_raw}
\small
\textbf{(PGD)}:\quad 
\begin{aligned}
\bbeta^\toptone 
&\leftarrow \mathcalP_{\sS}(\bbeta^\toptzero - \eta_t\nabla f(\bbeta^\toptzero)) \\
&= \mathop{\argmin}_{\bbeta\in\sS} \normtwo{\bbeta - \left(\bbeta^\toptzero - \eta_t \nabla f(\bbeta^\toptzero)\right)}^2\\
&=\mathop{\argmin}_{\bbeta\in\sS} 
\left\{
\frac{1}{2\eta_t} \normtwo{\bbeta-\bbeta^\toptzero}^2 + f(\bbeta^\toptzero)+ \innerproduct{\nabla f(\bbeta^\toptzero), \bbeta-\bbeta^\toptzero} 
\right\}.
\end{aligned}
\end{equation}
\end{minipage}
\end{tcolorbox}
\noindent
In the third line, the constant term $f(\bbeta^\toptzero)$ does not affect the minimizer and is included only for interpretability. 
This formulation reveals that each PGD step minimizes a local linear model of the objective $f$, regularized by a quadratic proximity term that penalizes large deviations from the current iterate $\bbeta^\toptzero$.

The linear term encourages progress in the direction of decreasing $f$, while the quadratic term ensures stability by keeping the next iterate close to the current one. This decomposition provides key intuition for designing more sophisticated first-order methods---such as the proximal gradient method, conditional gradient (Frank--Wolfe) method, and mirror descent---all of which balance local approximation with regularization tailored to the problem geometry.

\subsection*{Projection onto $\ell_1$-Ball}
To apply PGD to the constrained LASSO \eqref{opt:lc}, we must solve the following subproblem at each iteration: given a vector $\bv \in \real^p$, compute its Euclidean projection onto the $\ell_1$-ball of radius $\Sigma$ i.e., solve:
$$
\project_{\sB_{1}[\Sigma]}(\bv) 
\equiv 
\left\{
\min_{\bbeta\in\real^p} \frac{1}{2} \normtwo{\bbeta -\bv}^2 \quad \text{s.t.} \quad \normone{\bbeta} \leq \Sigma
\right\}.
$$
This is a convex optimization problem with a unique solution that can be computed efficiently \citep{duchi2008efficient}.

The projection $\bbeta^* = \project_{\sB_{1}[\Sigma]}(\bv) $ is obtained via \textit{soft-thresholding} (see Example~\ref{example:soft_thres}) and proceeds as follows:
\begin{enumerate}[(i)]
\item If $\normone{\bv} \leq \Sigma$, then $\bv$ already lies in the  in the $\ell_1$-ball, so  $\bbeta^* =\bv$.

\item Otherwise (i.e., $\normone{\bv} > \Sigma$), the projection is given by
$$
\beta^*_i = \sign(v_i) \cdot [\abs{v_i} - \lambda]_+, \quad \forall\, i\in\{1,2,\ldots,p\},
$$
where $[\cdot]_+ \triangleq  \max(\cdot, 0)$, and $\lambda \geq 0$ is chosen such that:
$$
\sum_{i=1}^p [\abs{v_i} - \lambda]_+ = \Sigma.
$$
In other words, $\lambda$ is the unique nonnegative scalar that ensures the soft-thresholded vector has $\ell_1$-norm exactly equal to $\Sigma$.
Intuitively, we uniformly shrink the magnitudes of all components of $\bv$ until the total $\ell_1$-norm reaches the constraint boundary $\Sigma$.
\end{enumerate}
Equivalently, we seek $\lambda$ such that $\normone{\mathcalT_{\lambda}(\bv)} = \Sigma$,
where the soft-thresholding operator $\mathcalT_{\lambda}(\bv)$ acts componentwise as $\big(\mathcalT_{\lambda}(\bv)\big)_i = \sign(v_i)[\abs{v_i} - \lambda]_+$.
This reduces to solving the scalar equation:
$\sum_{i=1}^p [\abs{v_i} - \lambda]_+ = \Sigma$.
Let $b_i \triangleq  \abs{v_i}$, and sort them in descending order: $b_{(1)} \geq b_{(2)} \geq \ldots \geq b_{(p)}$.
An efficient algorithm to compute $\lambda$  proceeds as follows:
\begin{tcolorbox}[colback=white,colframe=black]
\begin{minipage}{1\textwidth}
\begin{enumerate}
\item Compute $b_i = \abs{v_i}$ for all $i$.
\item Sort $\bb$ in descending order: $b_{(1)} \geq b_{(2)} \geq\ldots \geq b_{(p)}$.
\item Compute cumulative sums $S_t \triangleq \sum_{i=1}^t b_{(i)}$.
\item Find the largest $t \in \{1,2, \dots, p\}$ such that $b_{(t)} > \frac{S_t - \Sigma}{t}$.
This ensures $\lambda < b_{(t)}$.
\item Set 
\begin{equation}\label{equation:projec_ell1_st1}
\lambda \triangleq \frac{S_t - \Sigma}{t}.
\end{equation}
\item Apply soft-thresholding:
\begin{equation}\label{equation:projec_ell1_st2}
\beta^*_i = \sign(v_i) \cdot [\abs{v_i} - \lambda]_+.
\end{equation}
\end{enumerate}
\end{minipage}
\end{tcolorbox}
\noindent

\begin{example}
Let $\bv = [3, -2, 1]$, $\Sigma = 3$. 
Since $\normone{\bv} = 6 > 3$, projection is required.
The absolute values are $\bb = [3, 2, 1]$, already sorted in descending order.
\begin{itemize}
\item Try $t=1$: $S_1 = 3$, so $\lambda = (3 - 3)/1 = 0$, and $b_{(1)} = 3 > 0$, so $t=1$ is valid.
\item Try $t=2$: $S_2 = 3+2=5$, so $\lambda=(5-3)/2=1$. Check $b_{(2)}=2>1$ so  $t=2$ is valid.
\item Try $t=3$: $S_3= 3+2+1=6$, so $\lambda=(6-3)/3=1$. Check Check $b_{(3)}=1\ngtr 1$, so $t=3$ is not valid.
\end{itemize}
Thus, $\lambda = 1$ and we project:
\begin{align*}
\beta^*_1 &= \sign(3)(3 - 1) = 2; \\
\beta^*_2 &= \sign(-2)(2 - 1) = -1; \\
\beta^*_3 &= \sign(1)(1 - 1) = 0.
\end{align*}
The projection result is $\bbeta^* = [2, -1, 0]$ with $\normone{\bbeta^*} = 3 = \Sigma$.
\end{example}

\subsection*{Constrained LASSO with PGD}
Applying Algorithm~\ref{alg:pgd_gen} with the constraint set $\sS\triangleq \sB_1[\Sigma]$  to the least squares objective $f(\bbeta)\triangleq \frac{1}{2}\normtwo{\by - \bX \bbeta}^2$ yields a practical numerical method for solving the constrained LASSO problem.
Provided that we can project efficiently onto the $\ell_1$-ball, PGD behaves similarly to standard gradient descent applied to the unconstrained least squares problem. 
In particular, when the Gram matrix $\bX^\top\bX$  is well-conditioned, PGD enjoys a linear convergence rate---and this rate is fast when the condition number is small.

To formally establish this convergence guarantee, we introduce the proximal gradient method in the next subsection and clarify its relationship to projected gradient descent.

{
\begin{algorithm}[h] 
\caption{Proximal Gradient Method\index{Proximal gradient method}}
\label{alg:prox_gd_gen}
\begin{algorithmic}[1] 
\Require A function $f(\bbeta)$ and a closed convex function $g$ (usually non-smooth) satisfying (A1) and (A2) in {Theorem~\ref{theorem:prox_conv_ss_cvx}}; 
\State {\bfseries initialize:}  $\bbeta^\topone$;
\For{$t=1,2,\ldots$}
\State Choose a stepsize $\eta_t$;
\State $\balpha^\toptone \leftarrow \bbeta^\toptzero - \eta_t \nabla f(\bbeta^\toptzero)$; \Comment{(PROX$_1$)} 
\State $\bbeta^\toptone \leftarrow \prox_{\eta_t g}(\balpha^\toptone) \triangleq \mathcalT_{L_t}^{f,g}(\bbeta^\toptzero)$; \Comment{(PROX$_2$)} 
\State Stop if a stopping criterion is satisfied at iteration $t=T$;
\EndFor
\State (Output Option 1) Output  $\bbeta_{\text{final}}\leftarrow \bbeta^{(T)}$;
\State (Output Option 2) Output  $\bbeta_{\text{avg}}\leftarrow \frac{1}{T}(\sum_{t=1}^{t}\bbeta^\toptzero)$ or $\sum_{t=1}^{T} \frac{2t}{T(T+1)} \bbeta^\toptzero$;
\State (Output Option 3) Output  $\bbeta_{\text{best}}\leftarrow \argmin_{t\in\{1,2,\ldots,T\}} f(\bbeta^\toptzero)$;
\end{algorithmic} 
\end{algorithm}
}
\subsection{Proximal Gradient Method}\label{section:proxiGD_inClasso} 
\paragrapharrow{Proximal gradient method.}
In Section~\ref{section:proj_prox_sep}, we highlighted the close relationship between projection and proximal operators. Specifically, recall the identity 
$$
\prox_{\indicatorS}(\balpha) = \projectS(\balpha),
$$
where $\indicatorS$ denotes the indicator function of the set $\sS$.
This allows us to reinterpret the projection step in projected gradient descent as a proximal operation:
$$
\left\{\bbeta^\toptone \leftarrow \mathcalP_{\sS}(\balpha^\toptone)\right\}
\qquad\leadsto \qquad
\left\{\bbeta^\toptone \leftarrow \prox_{\indicatorS}(\balpha^\toptone)\right\}.
$$
The \textit{proximal gradient  method} (Algorithm~\ref{alg:prox_gd_gen}) generalizes this idea by replacing the indicator function with an arbitrary closed, convex, and possibly non-smooth function $g$:
\begin{equation}\label{equation:prox_obj1}
\left\{\bbeta^\toptone \leftarrow \prox_{\indicatorS}(\balpha^\toptone)\right\}
\qquad\leadsto \qquad
\left\{\bbeta^\toptone \leftarrow \prox_{\textcolor{mylightbluetext}{\eta_t g}}(\balpha^\toptone)\right\},
\end{equation}
where $\eta_t>0$ is the stepsize at iteration  $t$ (its role will become clear shortly)..
We assume that $g$ is closed and convex, which ensures that the proximal operator satisfies the standard properties  (Proximal Properties-O, I, II, III, IV; see Lemmas~\ref{lemma:prox_prop0}$\sim$\ref{lemma:prox_prop3}).
When $g\equiv \indicatorS$, the proximal gradient method reduces exactly to projected gradient descent.

Recall that in PGD, the update from iteration $t$ to $t+1$ can be written as:
\begin{center}
\framebox{
\begin{minipage}{0.95\textwidth}
\begin{equation}\label{equation:pgd_decom}
\small
\textbf{(PGD)}:\quad 
\begin{aligned}
\bbeta^\toptone 
&\leftarrow \mathop{\argmin}_{\bbeta\in\sS} \normtwo{\bbeta - \left(\bbeta^\toptzero - \eta_t \nabla f(\bbeta^\toptzero)\right)}^2\\
&=\mathop{\argmin}_{\bbeta\in\sS} 
\left\{
\frac{1}{2\eta_t} \normtwo{\bbeta-\bbeta^\toptzero}^2 + f(\bbeta^\toptzero)+ \innerproduct{\nabla f(\bbeta^\toptzero), \bbeta-\bbeta^\toptzero} 
\right\}.
\end{aligned}
\end{equation}
\end{minipage}
}
\end{center}
Here, the constant term $f(\bbeta^\toptzero)$ does not affect the minimizer and is included only for interpretability.
Thus, the PGD update minimizes a linear approximation of the smooth function $f$ around $\bbeta^\toptzero$, regularized by a quadratic proximity term. The linear term drives progress toward lower objective values, while the quadratic term keeps the new iterate close to the current one, ensuring stability.

The proximal gradient method extends this framework to unconstrained composite optimization by incorporating a non-smooth penalty $g$ directly into the objective:
\begin{subequations}\label{equation:prox_decom}
\begin{tcolorbox}[colback=white,colframe=black]
\begin{minipage}{1\textwidth}
\small
\begin{align}
\bbeta^\toptone 
&\leftarrow\mathop{\argmin}_{\bbeta\in\textcolor{mylightbluetext}{\real^p}} 
\left\{
\frac{1}{2\eta_t} \normtwo{\bbeta-\bbeta^\toptzero}^2 + f(\bbeta^\toptzero)+ \innerproduct{\nabla f(\bbeta^\toptzero), \bbeta-\bbeta^\toptzero }
+\textcolor{mylightbluetext}{g(\bbeta)}
\right\}\\
\textbf{(Prox)}:\qquad\,\,\,
&=\mathop{\argmin}_{\bbeta\in\textcolor{mylightbluetext}{\real^p}} \eta_t g(\bbeta) +\frac{1}{2}\normtwo{\bbeta-\left( \bbeta^\toptzero - \eta_t\nabla f(\bbeta^\toptzero) \right)}^2 \label{equation:prox_decom2}\\
&=\prox_{\eta_t g}\left(\bbeta^\toptzero - \eta_t\nabla f(\bbeta^\toptzero)\right)
\triangleq \mathcalT_{L_t}^{f,g}(\bbeta^\toptzero), \label{equation:prox_decom3}
\end{align}
\end{minipage}
\end{tcolorbox}
\noindent where $L_t \triangleq\frac{1}{\eta_t}$, and $\mathcalT_{L}^{f,g}(\bbeta)\triangleq\prox_{\frac{1}{L} g}\left(\bbeta - \frac{1}{L}\nabla f(\bbeta)\right)$ denotes the \textit{prox-grad operator}. 
We also define the \textit{gradient mapping} as
\begin{equation}
\mathcalG_{L}^{f,g}(\bbeta)\triangleq L\big(\bbeta-\mathcalT_{L}^{f,g}(\bbeta)\big).
\end{equation}
The update  \eqref{equation:prox_decom3} is computationally efficient whenever the proximal operator  $\prox_{\eta_t g}$ admits a closed-form solution or can be evaluated quickly (e.g., via soft-thresholding for $\ell_1$-regularization).
This formulation \eqref{equation:prox_decom3} is particularly well-suited for statistical learning problems that use penalized objectives---such as the Lagrangian (unconstrained) LASSO---rather than hard constraints of the form $g(\bbeta) \leq \Sigma$ (see Section~\ref{section:prox_gd_lag_lasso} for details).

Moreover, the update \eqref{equation:prox_decom3} can be rewritten as
\begin{equation}
\bbeta^\toptone 
\leftarrow 
\prox_{\eta_t g}\left(\bbeta^\toptzero - \eta_t\nabla f(\bbeta^\toptzero)\right)
\triangleq \mathcalT_{L_t}^{f,g}(\bbeta^\toptzero)
\triangleq \bbeta^\toptzero - \frac{1}{L_t} \mathcalG_{L_t}^{f,g}(\bbeta^\toptzero),
\end{equation}
\end{subequations}
which closely resembles a standard gradient descent (GD) step.
For this reason, $\mathcalG_{L_t}^{f,g}(\bbeta^\toptzero)$ is called the \textit{gradient mapping}. 
When the functions $f$ and $g$ are clear from  context,
we simplify notation by writing $\mathcalT_{L_t}(\cdot)$ and $\mathcalG_{L_t}(\cdot)$ instead of $\mathcalT_{L_t}^{f,g}(\cdot)$ and $\mathcalG_{L_t}^{f,g}(\cdot)$.

\index{Gradient mapping}\index{Prox-grad operator}

Equation~\eqref{equation:prox_decom} justifies the scaling by $\eta_t$ in \eqref{equation:prox_obj1}.
Altogether, the proximal gradient method is designed to minimize the composite objective
\begin{equation}\label{equation:prox_comp_prob}
\min \{F(\bbeta) \triangleq f(\bbeta)+g(\bbeta)\},
\end{equation}
where $f$ is smooth and $g$ is closed, convex, and possibly non-smooth.

Before presenting convergence results, we state the following descent lemma for the proximal gradient method.
\begin{lemma}[Descent lemma for proximal gradient method]\label{lemma:des_lem_prox}
Let $ f $ be a proper \textcolor{black}{closed} and $L_b$-smooth function, and let  $ g $ be a proper closed and convex function. 
Assume the following conditions hold:
\begin{itemize}
\item[(A1)] $ g: \real^p \rightarrow (-\infty, \infty] $ is proper, closed, and convex.
\item[(A2)] $ f: \real^p \rightarrow (-\infty, \infty] $ is proper, closed, and convex, $ \dom(f) $ is convex, $ \dom(g) \subseteq \textcolor{mylightbluetext}{\interior}(\dom(f)) $, and $ f $ is $ L_b $-smooth over $ \textcolor{mylightbluetext}{\interior}(\dom(f)) $.
\end{itemize} 
\noindent
Consider the composite function $ F \triangleq f + g $ and the prox-grad operator $ \mathcalT_{L} \triangleq \mathcalT_{L}^{f,g} $. 
Then, for any $ \bbeta \in \interior(\dom(f)) $ and $ L \in \left(\frac{L_b}{2}, \infty\right) $, the following inequality holds:
\begin{equation}\label{equation:des_lem_prox_res1}
F(\mathcalT_{L}(\bbeta)) -F(\bbeta) \leq \frac{\frac{L_b}{2} -L}{L^2} \normtwo{\mathcalG_{L}^{f,g}(\bbeta)}^2 \leq  0, 
\end{equation}
where $ \mathcalG_{L}^{f,g} : \interior(\dom(f)) \rightarrow \real^p $ is the operator defined by $\mathcalG_{L}^{f,g}(\bbeta)\triangleq L\big(\bbeta-\mathcalT_{L}(\bbeta)\big) $ for all $ \bbeta \in \interior(\dom(f)) $.

Moreover, for any $\balpha\in\real^p$, $\bbeta \in \interior(\dom(f))$, consider the \textit{proximal gap} $\mathcalD_f(\balpha, \bbeta) \triangleq f(\balpha) - f(\bbeta) - \innerproduct{\nabla f(\bbeta), \balpha - \bbeta }$~\footnote{When $f$ is convex, this is also known as the Bregman distance (Definition~\ref{definition:breg_dist}).} and  $L\triangleq L_b$. Then it holds that 
\begin{equation}\label{equation:des_lem_prox_res2}
F(\mathcalT_{L}(\bbeta)) - F(\balpha)  \leq   \frac{L}{2} \normtwo{\balpha-\bbeta}^2 -\frac{L}{2}\normtwo{\balpha - \mathcalT_{L}(\bbeta)}^2 - \mathcalD_f(\balpha, \bbeta).
\footnote{We can set $\balpha$ to the optimal point $\balpha\triangleq\bbeta^*$ to obtain convergence results.}
\end{equation}
When $\balpha=\bbeta$, \eqref{equation:des_lem_prox_res2} reduces to \eqref{equation:des_lem_prox_res1} by setting $L\triangleq L_b$.
\end{lemma}
\begin{proof}[of Lemma~\ref{lemma:des_lem_prox}]
\textbf{Equation~\eqref{equation:des_lem_prox_res1}.}
For brevity, denote $ \bbeta^+ \triangleq \mathcalT_{L}(\bbeta) $. By the smoothness of $f$ (Definition~\ref{definition:scss_func}),
$$
f(\bbeta^+) \leq f(\bbeta) + \innerproduct{\nabla f(\bbeta), \bbeta^+ - \bbeta} + \frac{L_b}{2} \normtwo{\bbeta - \bbeta^+}^2. 
$$
By the Proximal Property-I (Lemma~\ref{lemma:prox_prop1}), since $ \bbeta^+ = \prox_{\frac{1}{L}g}(\bbeta - \frac{1}{L}\nabla f(\bbeta)) $, we have
$$
\begin{aligned}
&\innerproduct{\bbeta - \frac{1}{L}\nabla f(\bbeta) - \bbeta^+, \bbeta - \bbeta^+} \leq \frac{1}{L} g(\bbeta) - \frac{1}{L} g(\bbeta^+), \quad \forall \, \bbeta\in\real^p\\
&\quad\implies\quad  g(\bbeta^+)    \leq  g(\bbeta)-L \normtwo{\bbeta^+ - \bbeta}^2  -\innerproduct{\nabla f(\bbeta), \bbeta^+ - \bbeta} .
\end{aligned}
$$
Combining the preceding two inequalities, we have 
$$
f(\bbeta^+) + g(\bbeta^+) \leq f(\bbeta) + g(\bbeta) + \left(-L + \frac{L_b}{2}\right) \normtwo{\bbeta^+ - \bbeta}^2.
$$
Thus, taking into account the definitions of $ \bbeta^+ $, $ \mathcalG_{L}^{f,g}(\bbeta) $, and the identities $ F(\bbeta) = f(\bbeta) + g(\bbeta) $, $ F(\bbeta^+) = f(\bbeta^+) + g(\bbeta^+) $, the desired result follows.

\paragraph{Equation~\eqref{equation:des_lem_prox_res2}.} Define the auxiliary function
$$
\varphi(\bu) \triangleq f(\bbeta) + \innerproduct{\nabla f(\bbeta), \bu - \bbeta} + g(\bu) + \frac{L}{2} \normtwo{\bu - \bbeta}^2. 
$$
This is precisely the first objective minimized in the proximal gradient update in \eqref{equation:prox_decom}.
Since $\varphi$ is an $L$-strongly convex function and $ \mathcalT_{L}(\bbeta) = \arg\min_{\bu \in \real^p} \varphi(\bu) $, it follows by Theorem~\ref{theorem:exi_close_sc}(ii) that
\begin{equation}\label{equation:des_lem_prox1}
\varphi(\balpha) - \varphi(\mathcalT_L(\bbeta)) 
\geq \frac{L}{2} \normtwo{\balpha - \mathcalT_L(\bbeta)}^2, 
\quad \text{for all }\balpha\in\real^p.
\end{equation}
This implies, by the smoothness of $f$, that
\begin{equation}\label{equation:des_lem_prox2}
\begin{aligned}
\varphi(\mathcalT_L(\bbeta)) &= f(\bbeta) + \innerproduct{\nabla f(\bbeta), \mathcalT_L(\bbeta) - \bbeta} + \frac{L}{2} \normtwo{\mathcalT_L(\bbeta) - \bbeta}^2 + g(\mathcalT_L(\bbeta)) \\
&\geq f(\mathcalT_L(\bbeta)) + g(\mathcalT_L(\bbeta)) = F(\mathcalT_L(\bbeta)).
\end{aligned}
\end{equation}
Thus, combining \eqref{equation:des_lem_prox1} and \eqref{equation:des_lem_prox2} yields that 
$$
\begin{aligned}
\underbrace{f(\bbeta) + \innerproduct{\nabla f(\bbeta), \balpha - \bbeta} + g(\balpha) + \frac{L}{2} \normtwo{\balpha - \bbeta}^2}_{=\varphi(\balpha)} - F(\mathcalT_L(\bbeta)) &\geq \frac{L}{2} \normtwo{\balpha - \mathcalT_L(\bbeta)}^2, \quad \forall\, \text{$\balpha \in \real^p $}.
\end{aligned}
$$
This proves the desired result.
\end{proof}

When $f$ is convex and strongly smooth, the proximal gradient method enjoys the following convergence guarantee.
\begin{theoremHigh}[Convergence of proximal gradient for convex and SS $f$: $\mathcalO(1/T)$ \citep{lu2025practical}]\label{theorem:prox_conv_ss_cvx}
Let $ f $ be a proper \textcolor{black}{closed} and $L_b$-smooth function, and let   $ g $ be a proper closed and convex function.
Assume the following conditions hold:
\begin{itemize}
\item[(A1)] $ g: \real^p \rightarrow (-\infty, \infty] $ is proper, closed, and convex.
\item[(A2)] $ f: \real^p \rightarrow (-\infty, \infty] $ is proper, closed, and convex, $ \dom(f) $ is convex, $ \dom(g) \subseteq \textcolor{mylightbluetext}{\interior}(\dom(f)) $, and $ f $ is $ L_b $-smooth over $ \textcolor{mylightbluetext}{\interior}(\dom(f)) $.
\end{itemize}
\noindent
Let $ \{\bbeta^\toptzero\}_{t > 0} $ be the sequence generated by the proximal gradient method (Algorithm~\ref{alg:prox_gd_gen}) applied to the composite problem $\min \{F(\bbeta) \triangleq f(\bbeta)+g(\bbeta)\}$  with a constant stepsize rule in which $ L_t  \triangleq L_b $ for all $ t > 0 $. 
Then, for any optimal solution $ \bbeta^*$ of $F\triangleq f+g$ and any integer $ T > 1 $, 
the following bound holds:
$$
F(\bbeta^{(T)}) - F(\bbeta^*) \leq \frac{ L_b \normtwo{\bbeta^\topone - \bbeta^*}^2}{2(T-1)}.
$$
Consequently, for any tolerance $\epsilon>0$, we have
$$
F(\bbeta^{(T)}) \leq  F(\bbeta^*) + \epsilon ,
\quad \text{for all $T\geq  \frac{ L_b \normtwo{\bbeta^\topone - \bbeta^*}^2}{2\epsilon} +1 $}.
$$
\end{theoremHigh}
\begin{proof}[of Theorem~\ref{theorem:prox_conv_ss_cvx}]
Since $f$ is convex, the proximal gap $\mathcalD_f$ qualifies a Bregman distance such that $\mathcalD_f(\bbeta^*, \bbeta^\toptzero)\geq 0$ (Remark~\ref{remark:bregnan_dist}).
For any $ t > 0 $, using the descent lemma in \eqref{equation:des_lem_prox_res2}, we obtain
$$
\begin{aligned}
\frac{2}{L_b} \left(  F(\bbeta^\toptone) - F(\bbeta^*) \right) 
&\leq   \normtwobig{\bbeta^* - \bbeta^\toptzero}^2 -\normtwobig{\bbeta^* - \bbeta^\toptone}^2 - \frac{2}{L_b} \mathcalD_f(\bbeta^*, \bbeta^\toptzero) \\
&\leq  \normtwobig{\bbeta^* - \bbeta^\toptzero}^2 - \normtwobig{\bbeta^* - \bbeta^\toptone}^2.
\end{aligned}
$$
Summing the above inequality over $ t = \{1, 2, \ldots, T-1\} $ yields a telescoping sum:
$$
\sum_{t=1}^{T-1} \left(  F(\bbeta^\toptone) - F(\bbeta^*) \right) \leq   \frac{ L_b}{2}\normtwo{\bbeta^* - \bbeta^\topone}^2 - \frac{ L_b}{2}\normtwo{\bbeta^* - \bbeta^{(T)}}^2
\leq  
\frac{ L_b}{2}\normtwo{\bbeta^* - \bbeta^\topone}^2.
$$
By the monotonicity of $ \{ F(\bbeta^\toptzero) \}_{t > 0} $ (Equation~\eqref{equation:des_lem_prox_res1} and Lemma~\ref{lemma:des_lem_prox}), we can conclude that
$$
\begin{aligned}
(T-1) \left( F(\bbeta^{(T)}) - F(\bbeta^*) \right) &\leq \sum_{t=1}^{T-1} \left( F(\bbeta^\toptone) - F(\bbeta^*) \right) \leq \frac{ L_b}{2} \normtwo{\bbeta^* - \bbeta^\topone}^2,
\end{aligned}
$$
which implies the desired result.
\end{proof}

\paragrapharrow{Constrained LASSO with proximal.}
Using an indicator function, the constrained LASSO problem in~\eqref{equation:loss_cons_lasso_scaled_in_algo} can be equivalently rewritten as the composite optimization problem:
\begin{equation}\label{equation:loss_cons_lasso_scaled_in_algo2} 
\min_{\bbeta\in\real^p}  F(\bbeta)\triangleq f(\bbeta)+g(\bbeta) \triangleq \frac{1}{2}\normtwo{\by - \bX \bbeta}^2 + \indicatorG_{\sB_{1}[\Sigma]}(\bbeta).
\end{equation}
Here, $\indicatorG_{\sB_{1}[\Sigma]}(\bbeta)$ denotes the indicator function of the scaled $\ell_1$-ball $\sB_1[\Sigma]$. 
$F(\bbeta)\triangleq f(\bbeta)+g(\bbeta)$ and $g(\bbeta)\triangleq \indicatorG_{\sB_{1}[\Sigma]}(\bbeta)$ is the indicator function.
Since this set is closed and convex, the indicator function $g$ is proper, closed, and convex (see Exercises~\ref{exercise_closed_indica} and~\ref{exercise_convex_indica}).

Applying the proximal gradient method (Algorithm~\ref{alg:prox_gd_gen}) to this formulation yields the following update rule at the $t$-th iteration:
\begin{subequations}\label{equation:prox_constrainedlso}
\begin{align}
\balpha^\toptone  &\leftarrow \bbeta^\toptzero - \eta_t \bX^\top (\bX\bbeta^\toptzero - \by);\\
\lambda^\toptone & \leftarrow  \text{Using the projection onto $\ell_1$-ball in \eqref{equation:projec_ell1_st1}};\\ 
\bbeta^\toptone &\leftarrow \prox_{\eta_t g}(\balpha^\toptone)= \mathcalT_{\eta_t\lambda^\toptone} (\balpha^\toptone)= \big[\absbig{\balpha^\toptone} - \eta_t \lambda^\toptone \bone\big]_+ \hadaprod \sign(\balpha^\toptone),
\end{align}
\end{subequations}
where $\mathcalT_\lambda(\cdot)$ denotes the soft-thresholding operator (Example~\ref{example:soft_thres}).
This is equivalent to the PGD update for constrained LASSO problem.
In words, the algorithm first takes a gradient descent step on the smooth loss $f$, then applies a proximal (shrinkage) step that enforces the $\ell_1$-norm constraint by projecting onto the feasible set. This shrinkage promotes sparsity in the iterates, which explains why the proximal gradient method is particularly well-suited for sparse recovery problems like the constrained LASSO.

Convergence of this scheme is guaranteed by Theorem~\ref{theorem:prox_conv_ss_cvx}. Specifically, since $f$ is $L_b$-smooth with $L_b=\normtwo{\bX^\top\bX}$ (Example~\ref{example:lipschitz_spar}),
a constant stepsize $\eta_t =1/L_b = \frac{1}{\normtwo{\bX^\top\bX}}$ ensures the $\mathcalO(1/T)$ convergence rate established in the theorem.

\begin{algorithm}[h] 
\caption{Conditional Gradient  Method\index{Conditional gradient  method}}
\label{alg:cond_gen}
\begin{algorithmic}[1] 
\Require A function $f(\bbeta)$ and a set $\sS$; 
\State {\bfseries initialize:}  $\bbeta^{(0)}$;
\For{$t=0,1,2,\ldots$}
\State Choose a stepsize $\eta_t\in[0,1]$;
\State $\widetildebbeta^\toptzero \leftarrow \mathop{\argmin}_{\bbeta\in \sS} \innerproduct{\nabla f(\bbeta^\toptzero), \bbeta}$; \Comment{linear optimization}
\State  $\bbeta^\toptone \leftarrow \bbeta^\toptzero + \eta_t (\widetildebbeta^\toptzero - \bbeta^\toptzero) = (1-\eta_t)\bbeta^\toptzero +  \eta_t \widetildebbeta^\toptzero $;  \Comment{``convex" update}
\State Stop if a stopping criterion is satisfied at iteration $t=T$;
\EndFor
\State \Return  $\bbeta_{\text{final}} \leftarrow \bbeta^{(T)}$;
\end{algorithmic} 
\end{algorithm}
\subsection{Conditional Gradient  (Frank--Wolfe) Method}\label{section:cond_gd}
\paragrapharrow{Conditional gradient  (CG, Frank--Wolfe) method}
We investigate  the \textit{conditional gradient method}, commonly known as the \textit{Frank--Wolfe (FW) algorithm}. 
This method offers an attractive alternative for solving constrained optimization problems, especially because it avoids the potentially expensive projection step required in PGD methods.
Consider again the following constrained optimization problem:
\begin{equation}\label{equation:cgd_prob}
\mathop{\min}_{\bbeta} f(\bbeta)\quad \text{s.t.}\quad \bbeta\in\sS.
\end{equation}
In the conditional gradient method, the next iterate is computed as a convex combination of the current iterate and a minimizer of a linear approximation of the objective function over the feasible set $\sS$ (see Algorithm~\ref{alg:cond_gen}).

Recall that in the PGD method, the update from the $t$-th to the $(t+1)$-th iteration is given by:
\begin{center}
\framebox{
\begin{minipage}{0.95\textwidth}
\begin{equation}\label{equation:pgd_decom_cg}
\small
\textbf{(PGD)}:\qquad 
\begin{aligned}
\bbeta^\toptone 
&\leftarrow \mathop{\argmin}_{\bbeta\in\sS} \normtwo{\bbeta - \left(\bbeta^\toptzero - \eta_t \nabla f(\bbeta^\toptzero)\right)}^2\\
&=\mathop{\argmin}_{\bbeta\in\sS} 
\left\{
\frac{1}{2} \normtwobig{\bbeta-\bbeta^\toptzero}^2+ \eta_t\innerproduct{\nabla f(\bbeta^\toptzero), \bbeta} 
\right\},
\end{aligned}
\end{equation}
\end{minipage}
}
\end{center}
where the quadratic term acts as a regularizer, ensuring the new iterate stays close to the  current point  $\bbeta^\toptzero$.

In contrast---while closely connected to the PGD update---the conditional gradient  method computes the next iterate as a convex combination of  a solution to the linear subproblem
$\mathop{\argmin}_{\bbeta\in\sS}\innerproductbig{\nabla f(\bbeta^\toptzero), \bbeta}$ and the current iterate  $\bbeta^\toptzero$.
This approach avoids quadratic regularization entirely. Nevertheless, the convex combination itself serves as a form of implicit regularization, preventing large deviations from the current iterate:
\begin{tcolorbox}[colback=white,colframe=black]
\begin{minipage}{1\textwidth}
\begin{equation}\label{equation:cgmethod}
\small
\textbf{(CG)}:\qquad 
\begin{aligned}
\bbeta^\toptone 
&\leftarrow (1-\eta_t)\bbeta^\toptzero + \eta_t\mathop{\argmin}_{\bbeta\in\sS}\innerproduct{\nabla f(\bbeta^\toptzero), \bbeta} .
\end{aligned}
\end{equation}
\end{minipage}
\end{tcolorbox}
Like PGD and proximal gradient method, the CG method achieves a convergence rate of $\mathcalO(T)$ for convex and smooth objective functions.
\begin{theoremHigh}[FW for convex and SS: $\mathcalO(1/T)$]\label{theorem:fw_under_smoo}
Let $f: \sS \rightarrow \real$ be a proper, differentiable, convex, and $L_b$-smooth function defined over the convex domain $\sS\subseteq\real^p$. 
Assume  that $f$ attains its global minimum at some point $\bbeta^* \in \sS$, and let  $D$ denote the diameter of $\sS$, defined as $D \triangleq  \max_{\bbeta, \bu \in \sS} \normtwo{\bbeta - \bu}$.
Let $ \{\bbeta^\toptzero\}_{t \geq 0} $
\footnote{For simplicity, we index the iterates starting from $t=0$ in the CG method.} 
be the sequence generated by the Frank--Wolfe method (Algorithm~\ref{alg:cond_gen}) for solving problem $\mathop{\min}_{\bbeta} f(\bbeta)$ subject to $\bbeta\in\sS$  with a dynamic stepsize rule in which $\eta_t = \frac{2}{t+2}$ at the $t$-th iteration. 
Then, for any $T\geq0$,
$$
f(\bbeta^{(T)}) - f(\bbeta^*) \leq \frac{2L_b D^2}{T+2}.
$$
\end{theoremHigh}

\begin{proof}[of Theorem~\ref{theorem:fw_under_smoo}]
By the smoothness, convexity, and the progress rule of conditional gradient method, we have
$$
\begin{aligned}
f(\bbeta^\toptone) 
&\leq f(\bbeta^\toptzero) + \nabla f(\bbeta^\toptzero)^\top (\bbeta^\toptone - \bbeta^\toptzero) + \frac{L_b}{2} \normtwo{ \bbeta^\toptone - \bbeta^\toptzero}^2\quad &&(\text{smoothness})\\
&=f(\bbeta^\toptzero) + \eta_t \nabla f(\bbeta^\toptzero)^\top (\widetildebbeta^\toptzero - \bbeta^\toptzero) + \frac{\eta_t^2 L_b}{2} \normtwo{ \widetildebbeta^\toptzero - \bbeta^\toptzero }^2\quad &&\text{(update rule)}\\
&\leq f(\bbeta^\toptzero) + \eta_t \nabla f(\bbeta^\toptzero)^\top (\bbeta^* - \bbeta^\toptzero) + \frac{\eta_t^2 L_b D^2}{2}. &&\text{(linear update rule)}
\end{aligned}
$$
By convexity of $f$, we also have:
$\nabla f(\bbeta^\toptzero)^\top (\bbeta^* - \bbeta^\toptzero) \leq f(\bbeta^*) - f(\bbeta^\toptzero)$ (see Theorem~\ref{theorem:conv_gradient_ineq});
Combining these two inequalities yields:
\begin{equation}\label{equation:fw_under_smoo1}
f(\bbeta^\toptone) - f(\bbeta^*) \leq (1 - \eta_t)(f(\bbeta^\toptzero) - f(\bbeta^*)) + \frac{\eta_t^2 L_b D^2}{2}.
\end{equation}
We now prove by induction that $f(\bbeta^\toptzero) - f(\bbeta^*) \leq \frac{2L_b D^2}{t+2}$ for all $t\geq 0$.

\paragraph{Base case $t = 0$.} When $t = 0$, we have $\eta_0 = \frac{2}{2} = 1$ and $\bbeta^{(0)}=\bzero$, and \eqref{equation:fw_under_smoo1} reduces to
$$
\begin{aligned}
f(\bbeta^\topone) - f(\bbeta^*) 
&\leq (1 - \eta_0)(f(\bbeta^{(0)}) - f(\bbeta^*)) + \frac{L_b}{2} \normtwo{\bbeta^\topone -\bbeta^*}^2 
=  \frac{L_b D^2}{2}
\leq \frac{2L_b D^2}{2}.
\end{aligned}
$$
\paragraph{Inductive step.} Assume the claim holds for some $t\geq 0$, i.e., 
$
f(\bbeta^\toptzero) - f(\bbeta^*) \leq \frac{2L_b D^2}{t+2}
$
holds for all integers up to $t$ and we show the claim for $t+1$. Using \eqref{equation:fw_under_smoo1} with $\eta_t=\frac{2}{t+2}$, we obtain:
\begin{align*}
f(\bbeta^\toptone) - f(\bbeta^*) &\leq \left(1 - \frac{2}{t+2}\right) \left(f(\bbeta^\toptzero) - f(\bbeta^*)\right) + \frac{2L_b D^2}{(t+2)^2} \\
&\leq \left(1 - \frac{2}{t+2}\right) \frac{2L_b D^2}{t+2} + \frac{2L_b D^2}{(t+2)^2} 
= 2L_b D^2 \frac{t+1}{(t+2)^2} \\
&= 2L_b D^2 \frac{t+1}{t+2} \cdot \frac{1}{t+2} 
\leq 2L_b D^2 \frac{t+2}{t+3} \cdot \frac{1}{t+2} 
= 2L_b D^2 \frac{1}{t+3}.
\end{align*}
Thus, the inequality also holds for the $t+1$ case, completing the induction.
\end{proof}

\begin{algorithm}[h] 
\caption{Conditional Gradient  Method}
\label{alg:conditional_gradient_lasso}
\begin{algorithmic}[1] 
\Require {Matrix $\bX$, vector $\by$, scalar $\Sigma > 0$, tolerance $\delta$};
\State {\bfseries goal:} 
Find $\bbeta$ such that $\frac{1}{2}\normtwo{\bX\bbeta - \by}^2$ is minimized subject to $\normone{\bbeta} \leq \Sigma$;
\State {\bfseries initialize:}  $\bbeta^{(0)} = \bzero$ (or any point in $\sS$);
\For{$t=0,1,2,\ldots$}

\State Compute gradient:
$\bg^\toptzero \leftarrow \nabla f(\bbeta^\toptzero) = \bX^\top(\bX\bbeta^\toptzero - \by)$;

\State Find coordinate with largest absolute gradient:
$i^* = \argmax_{i=1,2\ldots,p} \abs{g_i^\toptzero}$;

\State \algoalign{Compute linear minimizer:
	$\widetildebbeta^\toptzero \leftarrow -\sgn(g_{i^*}^\toptzero) \Sigma \be_{i^*}  $
	(i.e., $-\Sigma$ in coordinate $i^*$, zero elsewhere);}

\State Choose stepsize $\eta_t \in [0,1]$ (e.g., $\eta_t = \frac{2}{t+2}$);

\State  Update:
$\bbeta^\toptone \leftarrow (1 - \eta_t)\bbeta^\toptzero + \eta_t \widetildebbeta^\toptzero$;
\State Stop if a stopping criterion is satisfied at iteration $t=T$;
\EndFor
\State \Return  $\bbeta_{\text{final}} \leftarrow \bbeta^{(T)}$;
\end{algorithmic} 
\end{algorithm}

\paragrapharrow{Constrained LASSO with FW method.}
We apply the conditional gradient method to the constrained LASSO problem \eqref{opt:lc}; see Algorithm~\ref{alg:conditional_gradient_lasso}.
At the $t$-th iteration, the linear subproblem requires computing the gradient $\bg
\triangleq \nabla f(\bbeta^\toptzero)=\bX^\top(\bX\bbeta^\toptzero - \by)$. 
We then minimize the linear function  $\bg^\top \bbeta$ subject to the constraint $\normone{\bbeta} \leq \Sigma$.
By the duality between the $\ell_1$- and $\ell_\infty$-norms, the optimal value of this subproblem is
\begin{equation}\label{equation:lc_fw_subprob}
\widetildebbeta=\min_{\normone{\bbeta} \leq \Sigma} \bg^\top \bbeta = -\Sigma \norminf{\bg}.
\end{equation}
An optimal $\widetildebbeta$ is obtained by allocating  the entire $\ell_1$ budget $\Sigma$ to a single coordinate $j$ that achieves the maximum absolute gradient component, i.e., where $\abs{g_j} = \norminf{\bg}$, and choosing the sign opposite to $g_j$:
$$
\widetildebeta_i =
\begin{cases}
-\Sigma \sgn(g_j), & i = j, \\
0, & i \neq j.
\end{cases}
$$
This choice yields $\bg^\top \widetildebbeta = -\Sigma \norminf{\bg}$.
If multiple coordinates attain the maximum absolute value (i.e., there is a tie in $\norminf{\bg}$), the optimal solution is not unique. In such cases, one may distribute the $\ell_1$ budget among any subset of those coordinates, provided that each selected coordinate $i$ receives a sign opposite to $g_i$.
If $\bg = \bzero$, then every feasible $\bbeta$ (i.e., every $\bbeta$ with $\normone{\bbeta}\leq \Sigma$) is optimal, and the objective value is zero.

Since the feasible set is the  $\ell_1$-ball $\sB_1[\Sigma]$, its diameter (in the Euclidean norm) satisfies $D= \max_{\bbeta, \bu \in \sS} \normtwo{\bbeta - \bu}\leq 2\Sigma$, 
because any two points in the $\ell_1$-ball are at most $2\Sigma$ apart in $\ell_1$-norm, and $\normtwo{\cdot}\leq \normone{\cdot}$ (Exercise~\ref{exercise:cauch_sc_l1l2}). 
Thus, the method achieves an $\mathcalO(1/T)$ convergence rate, with constants depending on $\Sigma$.

\subsection{Least Angle Regression (LARS)\index{Least angle regression}}\label{section:lars}
\textit{Least angle regression (LARS)} is an efficient model-fitting algorithm that provides a unified framework for understanding several key regression techniques, including forward selection, forward stagewise regression, and the LASSO. It offers a geometrically intuitive approach to building linear models by iteratively adding predictors in a direction that maintains equal angles with their current correlations to the residual---hence the name least angle.

The algorithm proceeds sequentially, starting  with the predictor most correlated with the response. Rather than fully including this variable in the model right away, LARS gradually increases its coefficient, moving it continuously toward its least squares estimate. As the coefficient grows, the correlation between this predictor and the current residual decreases. The process pauses as soon as another predictor ``catches up"---that is, when its correlation with the residual becomes equal in magnitude. At this point, the second predictor joins the active set, and both coefficients are updated jointly, maintaining equal absolute correlation (i.e., equal angular distance) with the residual. This collaborative updating continues until all variables have been included, ultimately converging to the full least squares solution.

Before presenting the LARS algorithm, we first review two related methods: forward selection and forward stagewise regression, which serve as conceptual precursors to LARS \citep{efron2004least, weisberg2005applied}.

\subsection*{Forward Selection Algorithm\index{Forward selection algorithm}}

LARS is closely related to the classical model-selection method known as the \textit{forward selection algorithm} (a.k.a., \textit{forward stepwise regression}). 
Forward selection is a greedy algorithm that builds a model by sequentially adding the most informative predictor at each step.

Consider the linear model $ \by = \bX\bbeta $, where $ \by \in\real^n $  is the response vector and  $\bX\in\real^{n\times p}$ is the design matrix of $p$ predictors.
The algorithm begins with all coefficients set to zero and the initial residual defined as
$ \by^\topzero \triangleq \by $.
Among all candidate predictors in $\bX$, it selects the one---say, $ \bx_k $---that is most correlated with the current residual (equivalently, forms the smallest angle with it, or has the largest absolute cosine similarity).
The procedure then proceeds as follows:
The algorithm proceeds as follows:
\begin{enumerate}[(i)]
\item Compute the projection of $ \by^\topzero $ onto $ \bx_k $, yielding the fitted component $ \beta_k \bx_k $, and update the residual: $ \by^\topone = \by^\topzero - \beta_k \bx_k $. 
By construction, $ \by^\topone $ and $ \bx_k $ are orthogonal.
\item Treat  $ \by^\topone $ as the new response variable, remove $ \bx_k $ from the pool of available predictors, and repeat the process using the remaining features $\{\bx_i \mid i \neq k\}$.
\item Continue iterating until a stopping criterion is met---such as exhausting all predictors, achieving a sufficiently small residual, or reaching a pre-specified model size. The final estimate is $\bbeta = [\beta_1, \beta_2, \ldots, \beta_p]^\top$, where unused predictors have coefficient zero.
\end{enumerate}

\begin{figure}[h!]
\centering                      
\vspace{-0.35cm}               
\subfigtopskip=2pt             
\subfigbottomskip=2pt          
\subfigcapskip=-5pt              
\subfigure[Step 1.]{\label{fig:lars_forward1}
	\includegraphics[width=0.48\linewidth]{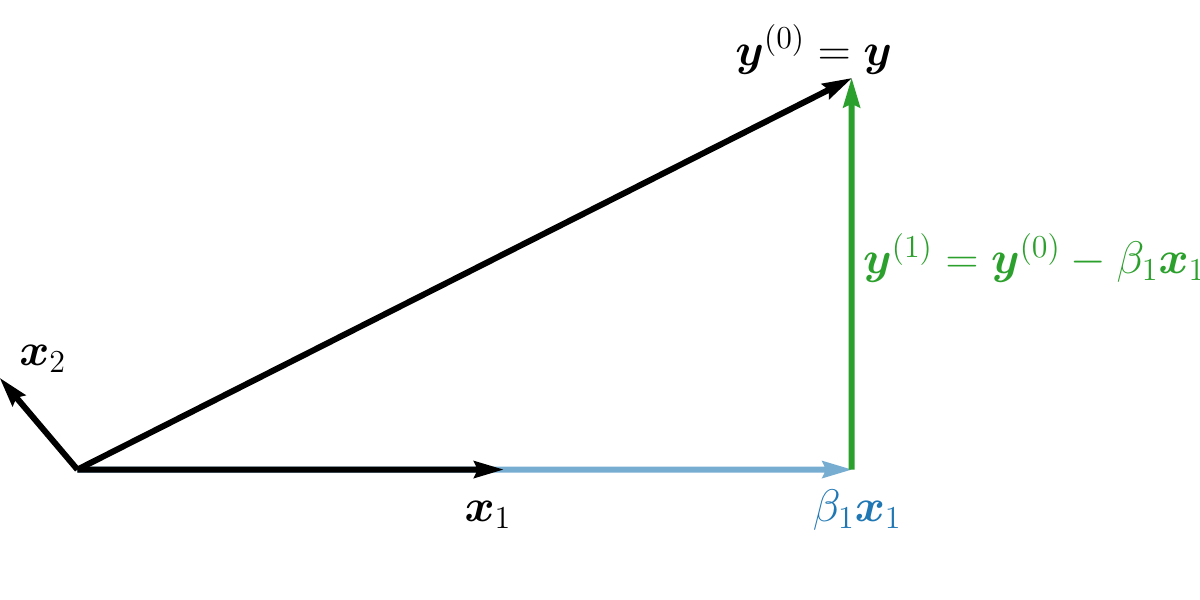}}
\subfigure[Step 2.]{\label{fig:lars_forward2}
	\includegraphics[width=0.48\linewidth]{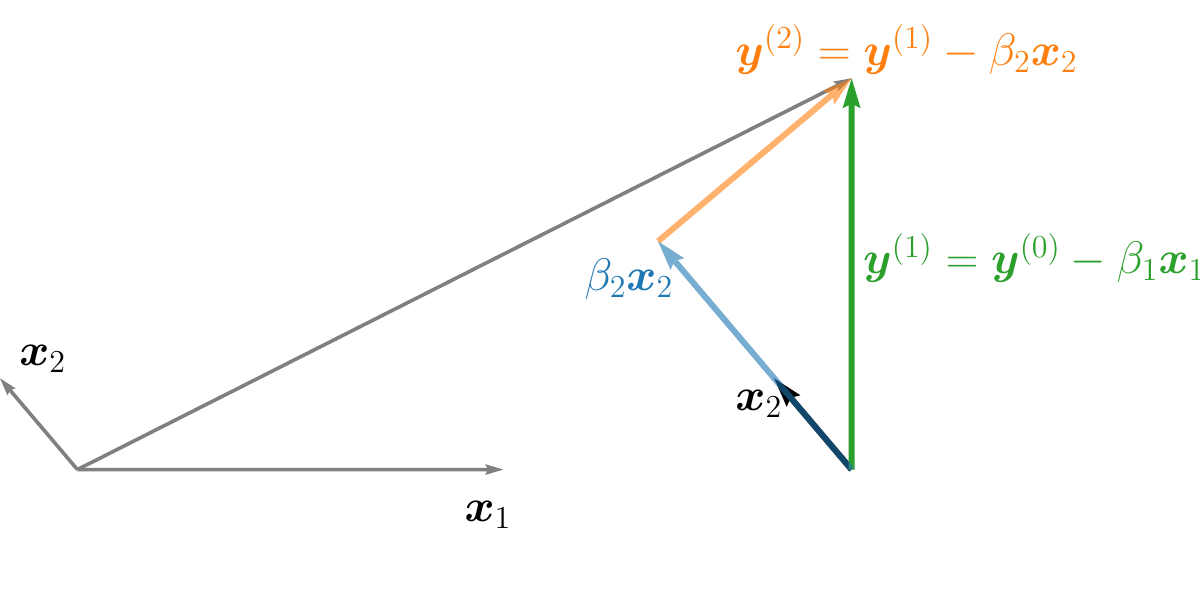}}
\caption{An example of the forward selection algorithm.}
\label{fig:forward_selection}
\end{figure}

To illustrate, consider a simple two-dimensional example (see Figure~\ref{fig:forward_selection}). In the first step, $ \bx_1 $ forms the smallest angle with $ \by^\topzero $. Its projection, $ \beta_1 \bx_1 $, serves as the initial approximation to $\by$, and the residual is $ \by^\topone = \by^\topzero - \beta_1\bx_1 $. In the next step, $ \by^\topone $ 
becomes the new target, and the only remaining predictor is $ \bx_2 $. After projecting onto $ \bx_2 $, the algorithm terminates, yielding the solution $ \bbeta = [\beta_1, \beta_2]^\top $.

Because forward selection performs only one update per selected variable, it is computationally efficient and fast. However, this simplicity comes at the cost of coarseness: once a variable is added, its coefficient is fixed, which can lead to suboptimal fits compared to more refined methods like LARS or LASSO.

\subsection*{Forward Stagewise Algorithm\index{Forward stagewise algorithm}}
Similar to the forward selection algorithm, the forward stagewise algorithm begins by identifying, among the set of feature variables $\bX$, the predictor $ \bx_k $ that is most correlated with the initial residual vector $ \by^\topzero \equiv \by$---that is, the one forming the smallest angle (or equivalently, having the largest absolute cosine similarity) with $\by$.
Let the stepsize be denoted by $ \eta >0 $. The algorithm proceeds as follows:
\begin{enumerate}[(i)]
\item Take a small step of size $\eta$ in the direction of $ \bx_k $: update the fit by adding $ \eta \bx_k $, and compute the new residual $ \by^\topone=\by^\topzero -\eta\bx_k $.

\item Treat $\by^{(1)}$ as the new response (residual), and retain all feature variables $\{\bx_i\}_{i=1}^p$---including those already used---as candidates for the next step. Repeat step (i), always selecting the predictor most correlated with the current residual, until a termination condition is met. Termination occurs when either the residual becomes sufficiently small or no further meaningful progress can be made.
\end{enumerate}

\begin{SCfigure}
\centering
\includegraphics[width=0.65\textwidth]{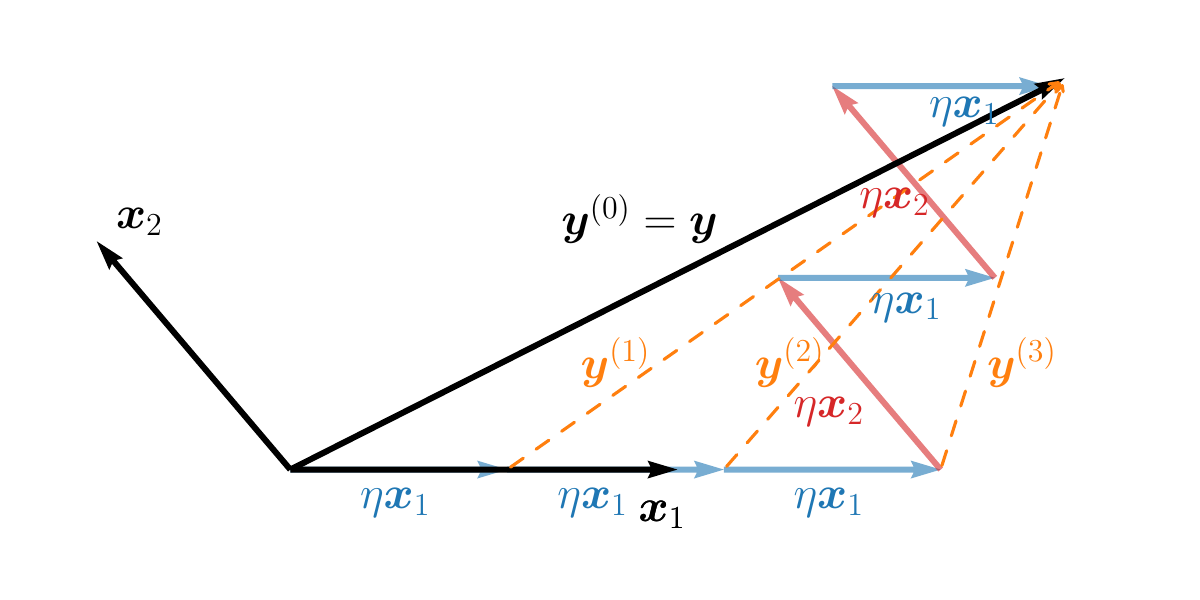}
\caption{Several steps in the forward stagewise algorithm for a two-dimensional case. 
Blue segments indicate movement in the direction of  $\bx_1$,  red segments indicate movement in the direction of $\bx_2$, and the orange vectors represent the residuals after each step. 
The staircase-like trajectory reflects the incremental, piecewise nature of the stagewise path.}
\label{fig:lars_forward3}
\end{SCfigure}

Figure~\ref{fig:lars_forward3} illustrates this process in a two-dimensional setting.
Suppose $\bx_1$ initially forms the smallest angle with $ \by^\topzero $. The algorithm takes repeated small steps in the direction of $\bx_1$, each of size $\eta\bx_1$, with the corresponding residuals shown as orange vectors. As the residual evolves, its correlation with $\bx_2$ increases. Once $ \bx_2 $ becomes more aligned with the residual than $ \bx_1 $ i.e., forms a smaller angle), the algorithm switches direction and begins taking steps along $ \bx_2 $ instead.

This procedure continues iteratively: at every step, all predictors remain eligible for selection, and any variable may be used multiple times. When the stepsize $\eta$ is chosen to be very small, the algorithm can closely approximate the optimal solution---but at the expense of increased computational cost due to the large number of iterations required.

The resulting path resembles a staircase, which is characteristic of the stagewise approach---hence the name forward stagewise algorithm.

\subsection*{LARS}

\textit{Least angle regression (LARS)} algorithm is a compromise between the two algorithms above: it is computationally more efficient than stagewise regression and more accurate than forward selection.

Let $\{\widehatby_0, \widehatby_1, \widehatby_2, \ldots, \widehatby_p\}$  denote the sequence of fitted value vectors generated by LARS as it progresses toward approximating the response vector $\by$. 
The algorithm is initialized with $\widehatby_0 = \bzero$ (corresponding to a model with only an intercept, which we assume has already been accounted for). 
The process terminates when one of the following conditions is met: 
there is no remaining residual to explain, all available predictors have been included, or the residual norm falls below a pre-specified tolerance.

In the one-dimensional case, LARS coincides exactly with forward selection. We compute 
$\widehatby_1 = \widehatby_0 + \widehatgamma_1 \bu_1$,
where $\bu_1$ is the unit vector in the direction of the selected predictor $\bx_1$. 
By choosing $\widehatgamma_1$ appropriately,  $\widehatby_1$ is aligned with $\widebarby_1$ (the orthogonal projection of the response variable $\by$ onto the one-dimensional subspace spanned by $\bx_1$), and the algorithm stops once the termination condition is satisfied.

\begin{figure}[h!]
\centering  
\vspace{-0.25cm}  
\subfigtopskip=2pt  
\subfigbottomskip=2pt  
\subfigcapskip=-5pt  
\subfigure[Two-dimensional case of LARS.]{\label{fig:lars2}
\includegraphics[width=0.47\linewidth]{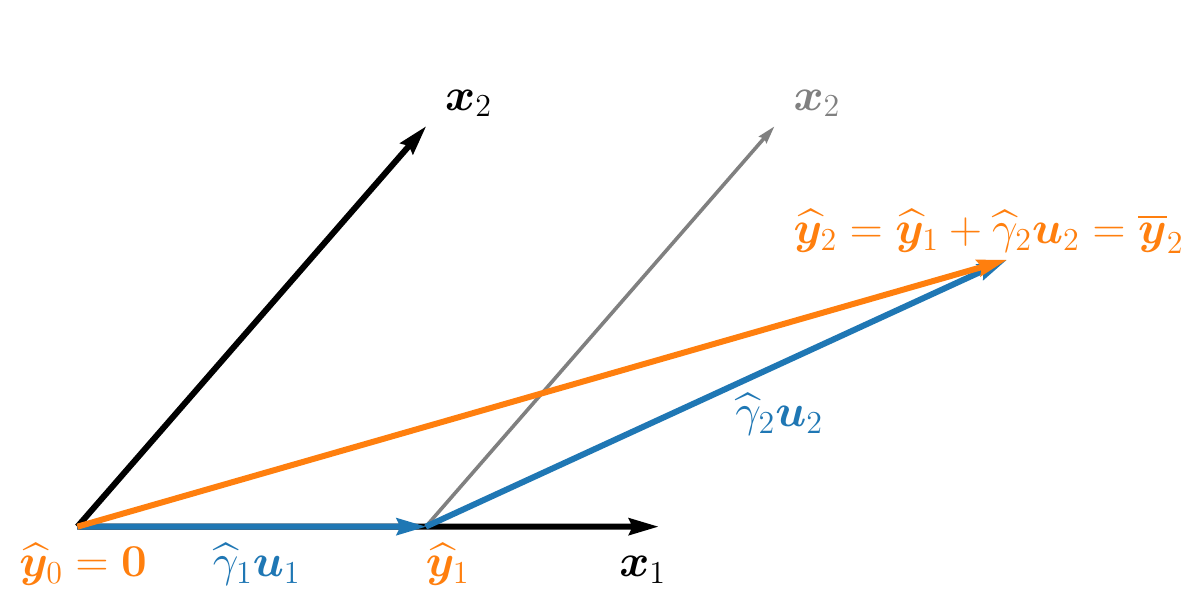}}
\quad 
\subfigure[Three-dimensional case of LARS.]{\label{fig:lars3}
\includegraphics[width=0.47\linewidth]{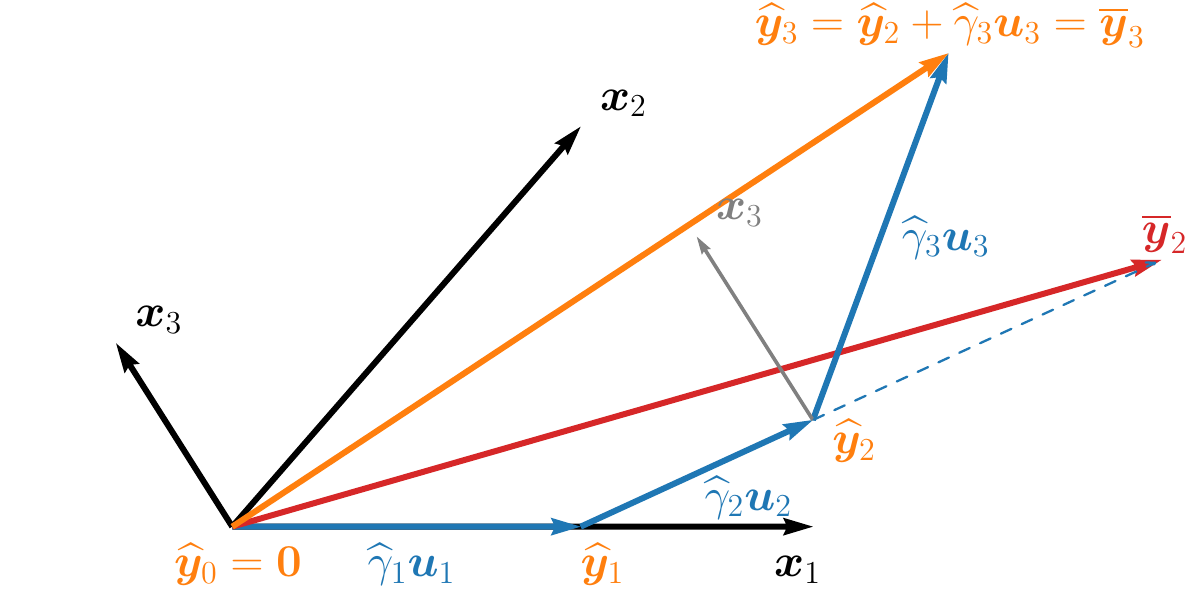}}
\caption{
Illustration of the LARS algorithm.
When $\bX$ is a two-dimensional matrix with two covariates (or features, predictors) $ \bx_1 $ and $ \bx_2 $, we find the component that is closest to $ \widebarby_2 $ (i.e., smallest angle, largest cosine value, and  $\widebarby_2$ is the projection of $\by$ onto $\cspace([\bx_1,\bx_2])$), say $ \bx_1 $. We move along the direction of $ \bx_1 $ by a stepsize $ \widehatgamma_1 \bu_1 $ such that the resulting residual $ \widebarby_2 - \widehatby_1 $ lies on the angle bisector between the feature vectors $ \bx_1 $ and $ \bx_2 $ (where $\bu_1$ is the unit vector along $\bx_1$). Let $ \widehatby_2 = \widehatby_1 + \widehatgamma_2 \bu_2 $, where $ \bu_2 $ is the unit vector along the angle bisector. By adjusting $ \widehatgamma_2 $, we can directly align $ \widehatby_2 $ with $ \widebarby_2 $ until the termination condition is met.
Note that $ \widehatby_2 = \widebarby_2 $ in the case $p=2$, but not for $p>2$, for example in Fig.~\ref{fig:lars3}.}
\label{fig:lars_all}
\end{figure}

In the two-dimensional case, suppose we have two predictors $\bx_1$ and $\bx_2$. 
Let  $\widebarby_2$ denote the projection of $\by$ onto the column space $\cspace([\bx_1,\bx_2])$.
LARS first identifies the predictor most correlated with the current residual $\widebarby_2$---in this case,  say $\bx_1$. 
Rather than moving all the way to $\widebarby_2$ along $\bx_1$, LARS takes a step of size $\widehatgamma_1$ in the direction of $\bu_1$ (unit vector along $\bx_1$), yielding
$$
\widehatby_1 = \widehatby_0 + \widehatgamma_1 \bu_1,
$$
chosen so that the new residual $\widebarby_2 - \widehatby_1$ lies on the angle bisector between
$\bx_1$ and $\bx_2$.
Next, the algorithm proceeds along this bisector direction. Let $\bu_2$
be the unit vector pointing along the bisector of the angle between $\bx_1$ and $\bx_2$. 
Then
$$
\widehatby_2 = \widehatby_1 + \widehatgamma_2 \bu_2.
$$
By adjusting $\widehatgamma_2$, we can move directly to $\widebarby_2$, completing the fit in two dimensions; see Figure~\ref{fig:lars2}.


For the general high-dimensional case, consider the $t$-th iteration of the algorithm:
\begin{enumerate}[(i)]
\item Move a step of size  $\widehatgamma_t$ in the equiangular (i.e., angle-bisecting) direction among the currently active predictors. Specifically, set
$$
\widehatby_t = \widehatby_{t-1} + \widehatgamma_t \bu_t,
$$
where $\bu_t$ is the unit vector in the equiangular direction defined by the active set.
This direction is chosen so that the updated residual remains equally correlated (in absolute value) with all active predictors. The next predictor to enter the active set---denoted $\bx_{t+1}$---is the one that becomes equally correlated with the residual at the end of this step. That is, the current residual $\widebarby_{t+1} - \widehatby_t$ lies on the angle bisector between $\bu_t$ and $\bx_{t+1}$.
\item Repeat step (i) until a termination condition is met.
\end{enumerate}

The key innovation of LARS occurs after the first predictor is selected. Unlike forward selection---which fully fits the selected variable before considering others---or stagewise regression---which takes tiny incremental steps---LARS moves in a direction that maintains equal angles (or equal absolute correlations) between the current residual and all active predictors. 
This equiangular direction ensures a balanced contribution from the selected variables. The algorithm continues along this path until another predictor ``catches up" in correlation with the residual, at which point it joins the active set, and a new equiangular direction is computed. This process yields a piecewise-linear solution path that is both efficient and geometrically elegant.

\paragrapharrow{Illustrating example.}
To further clarify the LARS algorithm, we make the following standard preprocessing assumptions:
\begin{itemize}
\item Each predictor $\bx_i$ is standardized to have mean zero and unit Euclidean norm.
\item The response vector $\by$ is centered so that it has mean zero.
\end{itemize}
These steps simplify computations and ensure all variables are on a comparable scale. Under this standardization, the (absolute) correlation between any predictor and the residual is proportional to their inner product, so angles can be interpreted directly via dot products.
Let $\widehatbbeta$ denote the current coefficient vector and $\widehatby$ be the current estimate.
LARS builds up estimates 
\begin{equation}
\widehatby = \bX\widehatbbeta
\end{equation}
in successive steps, each step adding one covariate to the model, so that after $t$ steps just $t$ of the coefficients $\widehatbeta_i$ are nonzero.  
Let $\bc(\widehatby)$ be the vector of \textit{current correlations}
\begin{equation}
\widehatbc \triangleq \bc(\widehatby) = \bX^\top(\by - \widehatby) \equiv \bX^\top(\by - \bX\widehatbbeta),
\end{equation}
such that $\widehatc_i = c_i(\widehatby)$ denotes the correlation between the $i$-th covariate $\bx_i$ and the current residual vector $\by - \bX\widehatbbeta$.

Figure~\ref{fig:lars2} illustrates the algorithm in the case of $p = 2$ predictors, $\bX = [\bx_1, \bx_2]$.
Because the predictors span a two-dimensional subspace, the current correlations depend only on the projection of $\by$ onto this space, denoted $\widebarby_2 = \project_{[\bx_1, \bx_2]}(\by)$. 
Thus,
\begin{equation}
	\bc(\widehatby) \triangleq  \bX^\top(\by - \widehatby) = \bX^\top(\widebarby_2 - \widehatby),
\end{equation}
where the $i$-th component $c_i(\widehatby)$ is precisely the cosine of the angle between $\bx_i$ and the residual $(\widebarby_2 - \widehatby)$ since the covariates $\bx_i$ have unit norms.

The algorithm starts at $\widehatby_0 = \bzero$. 
In Figure~\ref{fig:lars2} (and similarly in Figure~\ref{fig:lars3}), the initial residual  
$\widebarby_2 - \widehatby_0=\widebarby_2$ forms a smaller angle with $\bx_1$ than with $\bx_2$, i.e., $c_1(\widehatby_0) > c_2(\widehatby_0)$.
LARS then moves from  $\widehatby_0$ in the direction of $\bx_1$:
\begin{equation}
\widehatby_1 = \widehatby_0 + \widehat{\gamma}_1 \frac{\bx_1}{\normtwo{\bx_1}}
=\widehatby_0 + \widehat{\gamma}_1\bu_1.
\end{equation}
Classic forward selection algorithm would take $\widehat{\gamma}_1$ large enough to make $\widehatby_1$ equal $\widebarby_1\triangleq \project_{[\bx_1]}(\by)$, the orthogonal projection of $\by$ onto $\cspace([\bx_1])$;
while the forward stagewise algorithm  suggests  to choose $\widehat{\gamma}_1$ equal to some small value $\eta$, and then repeat the process many times.
LARS adopts an intermediate strategy: it selects $\widehat{\gamma}_1$ such that the new residual $\widebarby_2 - \widehatby_1$ becomes \textbf{equally} correlated with both $\bx_1$ and $\bx_2$.
Geometrically, this means the residual $\widebarby_2 - \widehatby_1$ lies along the angle bisector between $\bx_1$ and $\bx_2$, so that
$$
c_1(\widehatby_1) = c_2(\widehatby_1).
$$

Let $\bu_2 \triangleq \frac{\widebarby_2 - \widehatby_1}{\normtwo{\widebarby_2 - \widehatby_1}}$ denote the unit vector pointing along this bisector. The next LARS update is
\begin{equation}
	\widehatby_2 = \widehatby_1 + \widehat{\gamma}_2 \bu_2,
\end{equation}
with $\widehat{\gamma}_2$ chosen  so that $\widehatby_2 = \widebarby_2$ in the case $p = 2$ (Figure~\ref{fig:lars2}). 
With $p > 2$ covariates, however, $\widehat{\gamma}_2$ would be \textbf{smaller}, leading to another change of direction, as illustrated in {Figure~\ref{fig:lars3}} such that the new residual $ \widebarby_3 - \widehatby_2 $ lies on the angle bisector between $ \bu_2 $ and $ \bx_3 $ for some given $\bx_3$.

Subsequent LARS steps generalize this idea: once $k$ predictors are active, the algorithm proceeds in the equiangular direction---a unique direction that maintains equal absolute correlation with all active predictors. This direction generalizes the geometric notion of an angle bisector to higher dimensions, as exemplified by the vector $\bu_3$ in Figure~\ref{fig:lars3}.

\subsection*{Mathematical Formulas of LARS Algorithm}

With the geometric intuition behind the LARS algorithm established, we now turn to its precise mathematical formulation.
\paragrapharrow{Obtaining angle bisector with unit norm.}
Assume that the covariate vectors $\bx_1, \bx_2, \ldots, \bx_p$ are \emph{linearly independent}.
Suppose LARS has just completed step $t-1$, yielding the current estimate $\widehatby_{t-1}$, and is about to begin step $t$. 
Let $\sS$ denote the active set---the set of indices of predictors currently included in the model---with cardinality $\abs{\sS}=t$. Define the signed design matrix for the active variables as
\begin{subequations}\label{equation:ang_bisec}
\begin{equation}\label{equation:ang_bisec_eq0}
\bX_t \triangleq  [\ldots, s_i \bx_i, \ldots]_{i \in \sS},
\end{equation}
where the signs $s_i$ equal $\pm 1$.
\footnote{If the correlation of the residual vector with $\bx_i$ is positive, then $s_i=1$; otherwise $s_i=-1$; see \eqref{equation:lars_upd2} for more details. This signed covariates make the regression vector positive at the end of the algorithm, i.e., sign restriction; see Lemma~\ref{lemma:larscoef_reg}.} 
Next, define the Gram matrix and a scaling constant:
\begin{equation}\label{equation:ang_bisec_eq1}
\bG_t \triangleq \bX_t^\top \bX_t 
\qquad \text{and} \qquad 
\sigma_t \triangleq (\bone_t^\top \bG_t^{-1} \bone_t)^{-1/2},
\end{equation}
where $\bone_t$ is a vector of ones of length $\abs{\sS}=t$.
Then the vector
\begin{equation}\label{equation:ang_bisec_eq2}
\bu_t = \bX_t \bw_t, \quad \text{where } \bw_t \triangleq \sigma_t \bG_t^{-1} \bone_t,
\end{equation}
is a unit vector, satisfying
\begin{equation}\label{equation:ang_bisec_eq3}
\bX_t^\top \bu_t 
= \bX_t^\top \bX_t \bw_t
= \sigma_t \bone_t.
\end{equation}
\end{subequations}
Equality~\eqref{equation:ang_bisec_eq2} shows $\normtwo{\bu_t}=1$ such that $\bu_t$ is a unit-length vector; 
while \eqref{equation:ang_bisec_eq3} implies that the inner product (and hence the cosine of the angle) between $\bu_t$ and each signed predictor $s_i\bx_i$ ($i\in\sS$) is identical and equal to $\sigma_t>0$ (hence, the angles between $\bx_i$, $i\in\sS$ and $\bu_t$ are equal and less than $90^\circ$). Consequently, $\bu_t$ makes equal acute angles with all active predictors. 
Such a vector $\bu_t$ is known as the \textit{equiangular vector} or \textit{angle bisector} of the columns of $\bX_t$.

\paragrapharrow{The general step the LARS Algorithm.}
We can now describe the full LARS procedure.
Starting from $\widehatby_0 = \bzero$, the algorithm constructs the fitted values  $\widehatby$ through a sequence of steps, each typically larger than those in stagewise regression.
Suppose that $\widehatby_{t-1}$ is the current LARS estimate, lying in the column space of $\bX_{t-1}$ corresponding to an active set $\sS_-$ with $\abs{\sS_-} = \abs{\sS}-1=t-1$. Let the current correlation vector be
\begin{equation}\label{equation:lars_upd0}
\widehatbc = \bX^\top(\by- \widehatby_{t-1}),
\end{equation}
where $\by- \widehatby_{t-1}$ is the current residual. 
The general step consists of two parts:
(i) updating the estimate by moving in the equiangular direction within the current active subspace (i.e., column space of $\bX_t$), and
(ii) determining when a new predictor should enter the active set so as to update the active set to $\sS_+$ such that $\abs{\sS_+} = \abs{\sS}+1$.

The \emph{active set} $\sS\subseteq\{1,2,\ldots,p\}$ at step $t$ consists of all indices achieving the maximal absolute current correlation:
\begin{equation}\label{equation:lars_upd1}
\widehatC_t = \max_i \{\abs{\widehatc_i}\} 
\qquad \text{and} \qquad 
 \sS = \left\{ i \mid \abs{\widehatc_i} = \widehatC_t \right\}.
\end{equation}
For each $i\in\sS$, define the sign
\begin{equation}\label{equation:lars_upd2}
s_i \triangleq \sgn\{\widehatc_i\}, \quad \text{for } i \in \sS,
\end{equation}
and construct $\bX_t$, $\sigma_t$, and $\bu_t$ as in \eqref{equation:ang_bisec}. 
Additionally, compute the vector of directional correlations:
\begin{equation}\label{equation:lars_upd3}
\bz \triangleq \bX^\top \bu_t.
\end{equation}
The next LARS estimate is obtained by moving from $\widehatby_{t-1}$ in the direction $\bu_t$:
\begin{equation}\label{equation:lars_upd4}
\widehatby_t = \widehatby_{t-1} + \widehat{\gamma}_t \bu_t,
\end{equation}
where
\begin{equation}\label{equation:lars_upd5}
\widehat{\gamma}_t = \mathop{\text{min}^+}_{i \in \comple{\sS}} 
 \left\{ \frac{\widehatC_t - \widehatc_i}{\sigma_t - z_i}, \frac{\widehatC_t + \widehatc_i}{\sigma_t + z_i} \right\};
\end{equation}
``$\min^+$'' indicates that the minimum is taken over only positive components within each choice of $i$ in \eqref{equation:lars_upd5}.

To understand this choice  of the stepsize $\widehat{\gamma}_t$ \eqref{equation:lars_upd5}, consider a candidate update for any  $\gamma >0$. 
We define the potential next LARS estimate
\begin{equation}\label{equation:lars_upd6}
\by(\gamma) \triangleq  \widehatby_{t-1} + \gamma \bu_t,
\end{equation}
which yields the updated correlation for predictor $i$
\begin{equation}\label{equation:lars_upd7}
c_i(\gamma) = \bx_i^\top \big(\by - \by(\gamma)\big) = \widehatc_i - \gamma z_i, 
\quad \text{for all }i.
\end{equation}
For active indices $i \in \sS$, \eqref{equation:lars_upd1} along with the fact that $\bX_t^\top \bu_t = \sigma_t \bone_t$ in \eqref{equation:ang_bisec_eq3} yield
\begin{equation}\label{equation:lars_upd8}
\abs{c_i(\gamma)} = \widehatC_t - \gamma \sigma_t,
\quad \text{for all }i \in \sS,
\end{equation}
showing that all of the maximal absolute current correlations decline equally. 
For $i \in \comple{\sS}$, equating~\eqref{equation:lars_upd7} with~\eqref{equation:lars_upd8} shows that $c_i(\gamma)$ equals the maximal value at $\gamma = (\widehatC_t - \widehatc_i)/(\sigma_t - z_i)$. Likewise $-c_i(\gamma)$, the current correlation for the reversed covariate $-\bx_i$, achieves maximality at $(\widehatC_t + \widehatc_i)/(\sigma_t + z_i)$. 
Therefore, $\widehat{\gamma}_t$ in~\eqref{equation:lars_upd5} is the smallest positive value of $\gamma$ such that some new index $i_+$ joins the active set; $i_+$ is the minimizing index in~\eqref{equation:lars_upd5}, 
and the new active set $\sS_+$ is $\sS_+ = \sS \cup \{i_+\}$; the new maximum absolute correlation is $\widehatC_{t+1} = \widehatC_t - \widehat{\gamma}_t \sigma_t$ for $i\in\sS_+$.

\paragrapharrow{Geometric property of LARS  path.}
Let 
\begin{equation}\label{equation:larst_eq0}
\widehatC_t \bone_t \triangleq \bX_t^\top (\by-\widehatby_{t-1})
\end{equation}
denote the vector of correlations at the step $t$, where all active predictors share the same absolute correlation $\widehatC_t$ by construction.
Let $\widebarby_t$ denote the orthogonal projection of $\by$ onto the column space $\cspace(\bX_t)$. 
Since $\widehatby_{t-1} \in \cspace(\bX_{t-1})\subseteq \cspace(\bX_{t})$,  this projection can be written as
\begin{equation}\label{equation:larst_projet}
\widebarby_t = \widehatby_{t-1} + \bX_t \bG_t^{-1} \bX_t^\top (\by - \widehatby_{t-1}) = \widehatby_{t-1} + \frac{\widehatC_t}{\sigma_t} \bu_t,
\end{equation}
where $\widehatby_{t-1} + \bX_t \bG_t^{-1} \bX_t^\top (\by - \widehatby_{t-1})$ denotes the projection of the residual $\by - \widehatby_{t-1}$ onto $\cspace(\bX_t)$ (Lemma~\ref{lemma:projection-from-matrix}), 
and the last equality follows from \eqref{equation:ang_bisec_eq2} and the fact that the signed current correlations  all equal $\widehatC_t$ from \eqref{equation:larst_eq0}.
Since $\bu_t$ is a unit vector, \eqref{equation:larst_projet} says that $\widebarby_t - \widehatby_{t-1}$ has length
\begin{equation}\label{equation:larst_eq2}
\overline{\gamma}_t \triangleq \frac{\widehatC_t}{\sigma_t}
\quad\implies\quad 
\widebarby_t - \widehatby_{t-1} = \overline{\gamma}_t \bu_t.
\end{equation}

\begin{remark}[Geometric property of LARS  path]
Comparison with~\eqref{equation:lars_upd4}---the update of estimate $\widehatby_t 
= \widehatby_{t-1}  + \widehat{\gamma}_t \bu_t$---shows that the LARS estimate $\widehatby_t$ lies on the line from $\widehatby_{t-1}$ to $\widebarby_t$,
\begin{equation}\label{equation:larst_eq3}
\widehatby_t - \widehatby_{t-1} 
= \widehat{\gamma}_t \bu_t
= \frac{\widehat{\gamma}_t}{\overline{\gamma}_t} (\widebarby_t - \widehatby_{t-1}).
\end{equation}
By \eqref{equation:lars_upd5}, it is easy to see that $\widehat{\gamma}_t$ is always less than $\overline{\gamma}_t$, so that $\widehatby_t$ lies closer than $\widebarby_t$ to $\widehatby_{t-1}$. 
Therefore, the successive LARS estimates $\widehatby_t$ are always approaching but never reaching the OLS estimates $\widebarby_t$.

The only exception occurs at the final step. 
When the active set includes all $p$ predictors, i.e., $\sS = \{1,2,\ldots,p\}$, the stepsize formula \eqref{equation:lars_upd5} is no longer applicable (as there are no inactive variables). 
By convention, the algorithm takes $\widehat{\gamma}_p = \overline{\gamma}_p = \widehatC_p / \sigma_p$, making the current estimate $\widehatby_p = \widebarby_p$ and and the current coefficient {$\widehat{\bbeta}_p$} equal the OLS estimate for the full set of $p$ covariates.
\end{remark}

For notational clarity, we define: $\bc_t \triangleq \bX_t^\top (\by - \widehatby_{t-1})$ as the vector of correlations at step $t$, with $c_{ti}$ denoting  the $i$-th entry of $\bc_t$.
And we let $\widehatbbeta_t$ denote the full regression coefficient vector satisfying $\widehatby_t = \bX \widehatbbeta_t$, and let $\widehatbeta_{ti}$ denote the $i$-th entry of $\widehatbbeta_t$.
We now establish an important property of the LARS algorithm: new variables enter the model in the direction consistent with their current correlation with the residual.
\begin{lemma}[Sign restriction of LARS]\label{lemma:larscoef_reg}
Suppose LARS has completed step $t-1$, yielding the estimate $\widehatby_{t-1}$ and active set $\sS_t$ for step $t$, 
where the predictor  $\bx_t$ is the newest addition to the active set (i.e., $\sS_t=\sS_{t-1}\cup \{t\}$).
If $\bx_t$ is the only addition to the active set at the end of step $t-1$, then the coefficient vector $\bw_t = \sigma_t \bG_t^{-1} \bone_t$ for the equiangular vector $\bu_t = \bX_t \bw_t$, \eqref{equation:ang_bisec_eq2}, has its $t$-th component $w_{tt}$ agreeing in sign with the current correlation $c_{tt} = \bx_t^\top(\by - \widehatby_{t-1})$. Moreover, the regression vector $\widehatbbeta_t$ for $\widehatby_t = \bX \widehatbbeta_t$ has its $t$-th component $\widehatbeta_{tt}$ agreeing in sign with $c_{tt}$.
\end{lemma}

\begin{proof}[Lemma~\ref{lemma:larscoef_reg}]
The case $t=1$ is immediate,  since the first predictor enters in the direction of its correlation with $\by$. 
For $t\geq 2$, recall from \eqref{equation:larst_eq0} that
$\\bX_t^\top(\by - \widehatby_{t-1}) = \widehatC_t \bone_t$.
Using \eqref{equation:ang_bisec_eq2}, we can express
\begin{equation*}
\bw_t 
= \sigma_t \bG_t^{-1}\bone_t
= \sigma_t \widehatC_t^{-1} [(\bX_t^\top \bX_t)^{-1} \bX_t^\top (\by - \widehatby_{t-1})] 
\triangleq \sigma_t \widehatC_t^{-1} \widehatbw_t.
\end{equation*}
The term in square braces is the least squares coefficient vector in the regression of the current residual on $\bX_t$, and the term preceding it is positive.

Now observe that
\begin{equation}\label{eq:residual_decomp}
\bX_t^\top(\by - \widebarby_{t-1}) = [\bzero, \delta]^\top, \quad \text{with } \delta > 0,
\end{equation}
since $\bX_{t-1}^\top(\by - \widebarby_{t-1}) = \bzero_{t-1}$ by definition (because $\widebarby_{t-1}$ is the projection of $\by$ onto $\cspace(\bX_{t-1})$), and $c_t(\gamma) = \bx_t^\top(\by - \gamma \bu_{t-1})$ decreases more slowly in $\gamma$ than $c_i(\gamma)$ for $i \in \sS_{t-1}$:
\begin{equation*}
c_t(\gamma)
\begin{cases}
< c_i(\gamma), & \text{if } \gamma < \widehat{\gamma}_{t-1}; \\
= c_i(\gamma) = \widehatC_t, & \text{if } \gamma = \widehat{\gamma}_{t-1}; \\
> c_i(\gamma), & \text{if } \widehat{\gamma}_{t-1} < \gamma < \overline{\gamma}_{t-1}.
\end{cases}
\end{equation*}
Thus,
\begin{align}
\widehatbw_t
&= (\bX_t^\top \bX_t)^{-1} \bX_t^\top (\by - \widebarby_{t-1} + \widebarby_{t-1} - \widehatby_{t-1}) \nonumber \\
&= (\bX_t^\top \bX_t)^{-1} 
\begin{bmatrix} \bzero \\ 
\delta \end{bmatrix} 
+ (\bX_t^\top \bX_t)^{-1} \bX_t^\top [( \overline{\gamma}_{t-1} - \widehat{\gamma}_{t-1}) \bu_{t-1}], \label{eq:w_star_split}
\end{align}
where the last equality follows from \eqref{equation:larst_eq2} and \eqref{equation:larst_eq3}.
The $t$-th element of $\widehatbw_t$ is positive, because it is in the first term in~\eqref{eq:w_star_split} (since $\bX_t$ has full column rank and thus $\bX_t^\top \bX_t$ is positive definite), and in the second term it is 0 since $\bu_{t-1} =\bX_{t-1}\bw_{t-1}\in \cspace(\bX_{t-1})$.
This proves the first statement. 

For the second claim, note that the full coefficient vector updates as
\begin{equation}\label{eq:w_star_split33}
\widehatbeta_{tt} = \widehatbeta_{t-1,t} + \widehat{\gamma}_t w_{tt},
\end{equation}
where $\widehatbeta_{t-1,t} = 0$, since $\bx_t$  was not active before step $t$. 
Because both $\widehat{\gamma}_t>0$ and $w_{tt}$ share the sign of $c_{tt}$, it follows that  
$\widehatbeta_{tt}$ also shares that sign.
\end{proof}

Lemma~\ref{lemma:larscoef_reg} shows that new variables enter the LARS model in the ``correct" direction---i.e., with a coefficient sign matching their current correlation with the residual. This is a weaker but essential property compared to the stricter sign consistency required by the LASSO (see condition~\eqref{equation:larslass_res1}). As we shall see, this alignment is a key ingredient in the deep connection between LARS and the LASSO.

\subsection*{LARS for Constrained LASSO Problems}
LARS proceeds along the most natural compromise direction---the equiangular vector defined in \eqref{equation:ang_bisec_eq2}---whereas the constrained LASSO imposes additional restrictions on this equiangular strategy.

We now briefly show that the entire set of LASSO solutions can be generated by a minor modification of the standard LARS algorithm \eqref{equation:lars_upd0}$\sim$\eqref{equation:lars_upd5}.

\paragrapharrow{Piecewise linear update of the fitted values in LARS--LASSO.}
We begin by discussing the relationship between LARS and the LASSO.
Let $\widehatbbeta = \widehatbbeta(\Sigma)$  denote a solution to the constrained LASSO problem~\eqref{opt:lc}, and define the corresponding fitted values as $\widehatby = \widehatby(\Sigma) = \bX \widehatbbeta(\Sigma)$, where $\Sigma\geq 0$ is the constraint radius satisfying $\normone{\bbeta}\leq \Sigma$. 
Since both the loss function $f(\widehatbbeta) = \frac{1}{2}\normtwobig{\bX\widehatbbeta-\by}^2$ and the constraint $g(\widehatbbeta) = \normonebig{\widehatbbeta}$ are  convex functions of $\widehatbbeta$---and $f$ is strictly convex---it follows from standard convex analysis that $\widehatbbeta(\Sigma)$ and $\widehatby(\Sigma)$ are unique and continuous functions of the constraint radius $\Sigma$.

For a fixed $\Sigma$, define the active set as
$$
\sS \triangleq  \{ i \mid \widehatbeta_i(\Sigma) \neq 0 \}.
$$
Indeed, $\sS$ is also the active set that determines the equiangular direction $\bu_\sS$ in the LARS--LASSO computations.

We wish to characterize the trajectory of the LASSO solutions $\widehatbbeta(\Sigma)$ or equivalently of $\widehatby(\Sigma)$ as $\Sigma$ increases from 0 to its maximum effective value. Let $\sM$ be an open interval of the $\Sigma$ axis, with infimum $\Sigma_0 \triangleq \inf(\sM)$, within which the set $\sS$ of nonzero LASSO coefficients $\widehatbeta_i(\Sigma)$ remains constant:
$$
\sS \equiv \{i\mid \widehatbeta_i(\Phi) \neq 0\}, 
\quad \text{for all } \Phi \in \sM.
$$

\begin{lemma}[Piecewise linear update of the fitted values in LARS--LASSO]\label{lemma:larslasso_ud_est}
For $\Sigma \in \sM$, the LASSO estimates $\widehatby(\Sigma)$ satisfy
\begin{equation}\label{equation:haty_sigma_linear}
\widehatby(\Sigma) = \widehatby(\Sigma_0) + \sigma_\sS (\Sigma - \Sigma_0) \bu_\sS
\end{equation}
where $\bu_\sS$ is the equiangular vector $\bu_{\sS}=\bX_\sS \bw_\sS$, $\bw_\sS = \sigma_\sS \bG_\sS^{-1} \bone_\sS$.
Equivalently, in terms of the nonzero regression coefficients $\widehatbbeta_\sS(\Sigma)$,
\begin{equation}\label{equation:beta_sigma_linear}
\widehatbbeta_\sS(\Sigma) = \widehatbbeta_\sS(\Sigma_0) + \bS_\sS \sigma_\sS (\Sigma - \Sigma_0) \bw_\sS,
\end{equation}
where $\bS_\sS \triangleq  \diag(\bs_\sS)$ is the diagonal matrix with diagonal elements $s_i = \sign(\widehatbeta_i)$, $i \in \sS$, and $\bs_\sS$ is a vector containing these sign elements.~\footnote{
The matrix $\bS_\sS$ appears because, by definitions~\eqref{equation:ang_bisec_eq0} and~\eqref{equation:lars_upd2}, we have  $\widehatby(\Sigma) = \bX \widehatbbeta(\Sigma) = \bX_\sS \bS_\sS \widehatbbeta_\sS(\Sigma)$.
}
\end{lemma}
Note that~\eqref{equation:haty_sigma_linear} implies that the maximum absolute correlation $\widehatC(\Sigma)$ equals $\widehatC(\Sigma_0) - \sigma_\sS^2 (\Sigma - \Sigma_0)$ by \eqref{equation:ang_bisec_eq3}, so that $\widehatC(\Sigma)$ is a piecewise linear decreasing function of the LASSO parameter $\Sigma$. 

\begin{proof}[of Lemma~\ref{lemma:larslasso_ud_est}]
The lemma says that, for $\Sigma$ in $\sM$ where the active set $\sS$ and signs $\bs_{\sS} $
remain fixed, $\widehatby(\Sigma)$ moves linearly along the equiangular vector $\bu_\sS$ determined by $\sS$. 
Since $\widehatbbeta(\Sigma)$ satisfies the constrained LASSO problem \eqref{opt:lc} and has support $\sS$, it also minimizes
\begin{equation*}
f(\widehatbbeta_\sS) = \frac{1}{2}\normtwo{\by - \bX_\sS \bS_\sS \widehatbbeta_\sS}^2
\end{equation*}
subject to
\begin{equation}\label{eq:constraint}
\sum_{i \in \sS} s_i \widehatbeta_i = \Sigma 
\qquad \text{and} \qquad 
\sign(\widehatbeta_i) = s_i \quad \text{for } i \in \sS.
\end{equation}
The inequality constraint in the constraint LASSO problem \eqref{opt:lc} can be replaced by $\normonebig{\widehatbbeta} = \Sigma$ as long as $\Sigma$ is less than $\sum \absbig{\widehatbeta_i}$ for the full $p$-variable OLS solution $\widehatbbeta_p$.
Moreover, the fact that the minimizing point $\widehatbbeta_\sS(\Sigma)$ occurs strictly inside the simplex~\eqref{eq:constraint}, combined with the strict convexity of $f(\widehatbbeta_\sS)$, implies we can drop the second condition in~\eqref{eq:constraint} so that $\widehatbbeta_\sS(\Sigma)$ solves
\begin{equation}\label{eq:min_pieceproblem}
\min \quad \{ f(\widehatbbeta_\sS) \} \quad \text{s.t.} \quad \sum_{i \in \sS} s_i \widehatbeta_i = \Sigma.
\end{equation}
Introducing a Lagrange multiplier $\lambda$, \eqref{eq:min_pieceproblem} becomes
$$
\min \quad \frac{1}{2} \normtwo{\by - \bX_\sS \bS_\sS \widehatbbeta_\sS}^2 + \lambda \sum_{i \in \sS} s_i \widehatbeta_i.
$$
Differentiating with respect to $\bbeta_\sS$ and setting the gradient to zero yields
\begin{align}
&- \bS_\sS \bX_\sS^\top (\by - \bX_\sS \bS_\sS \widehatbbeta_\sS) + \lambda \bS_\sS \bone_\sS = \bzero \label{eq:larslasso_stat_condition010} \\
&\quad \implies\quad
\bbeta_\sS = (\bS_\sS\bX_\sS^\top \bX_\sS\bS_\sS)^{-1}(\bS_\sS\bX_\sS^\top \by - \lambda \bs_\sS). \label{eq:larslasso_stat_condition}
\end{align}

Now consider two values $\Sigma_1, \Sigma_2\in \sM$ with $\Sigma_0 < \Sigma_1 < \Sigma_2$, 
and let $\lambda_1 > \lambda_2$ be the corresponding Lagrange multipliers  (since larger $\Sigma$ implies smaller penalty weight).  
Subtracting the optimality conditions for $\Sigma_1$ and $\Sigma_2$ gives
$$
\bX_\sS^\top \bX_\sS \bS_\sS \big(\widehatbbeta_\sS(\Sigma_2) - \widehatbbeta_\sS(\Sigma_1)\big) = (\lambda_1 - \lambda_2) \bone_\sS,
$$
which implies
$$
\widehatbbeta_\sS(\Sigma_2) - \widehatbbeta_\sS(\Sigma_1) = (\lambda_1 - \lambda_2) \bS_\sS \bG_\sS^{-1} \bone_\sS.
$$
However, we can also prove $\bs_\sS^\top [\widehatbbeta_\sS(\Sigma_2) - \widehatbbeta_\sS(\Sigma_1)] = \Sigma_2 - \Sigma_1$ according to the LASSO definition:

\begin{mdframed}[hidealllines=\mdframehideline,backgroundcolor=\mdframecolor,frametitle={Relationship between $\widehatbbeta_\sS(\Sigma_2) - \widehatbbeta_\sS(\Sigma_1) $ and $ \Sigma_2 - \Sigma_1$:}]
For the active set, $\Sigma = \normone{\bbeta} = \sum_{i \in \sS} \abs{\beta_i} = \bs_\sS^\top \bbeta_\sS$.
Substitute \eqref{eq:larslasso_stat_condition}:
$$
\Sigma = \bs_\sS^\top (\bS_\sS\bX_\sS^\top \bX_\sS\bS_\sS)^{-1}(\bS_\sS\bX_\sS^\top \by - \lambda \bs_\sS).
$$
Define
$a \triangleq \bs_\sS^\top (\bS_\sS\bX_\sS^\top \bX_\sS\bS_\sS)^{-1} \bS_\sS\bX_\sS^\top \by$ 
and  $b \triangleq \bs_\sS^\top (\bS_\sS\bX_\sS^\top \bX_\sS\bS_\sS)^{-1} \bs_\sS$.
Then
$$
\Sigma = a - \lambda b
\quad \Longrightarrow \quad
\lambda = \frac{a - \Sigma}{b}.
$$
Thus, $\lambda$ is an affine (linear) function of $\Sigma$
while $\sS$ and $\bs_\sS$ remain fixed.
Differentiate \eqref{eq:larslasso_stat_condition} with respect to $\Sigma$:
$
\frac{\partial \bbeta_\sS}{\partial \Sigma}
= - (\bS_\sS\bX_\sS^\top \bX_\sS\bS_\sS)^{-1} \bs_\sS \cdot \frac{\partial \lambda}{\partial \Sigma}
$.
From $\lambda = (a - \Sigma)/b$, we get $\dfrac{\partial \lambda}{\partial \Sigma} = -\dfrac{1}{b}$.
Hence
$$
\frac{\partial\bbeta_\sS}{\partial \Sigma}
= \frac{(\bS_\sS\bX_\sS^\top \bX_\sS\bS_\sS)^{-1} \bs_\sS}
{\bs_\sS^\top (\bS_\sS\bX_\sS^\top \bX_\sS\bS_\sS)^{-1} \bs_\sS}
, \qquad
\frac{\partial\bbeta_{\comple{\sS}}}{\partial\Sigma} = \bzero.
$$
If $\Sigma_1$ and $\Sigma_2$ lie within the same region
where the active set $\sS$ and sign vector $\bs_\sS$ remain unchanged, then
\begin{equation}
\bbeta_2 - \bbeta_1
= (\Sigma_2 - \Sigma_1)
\frac{(\bS_\sS\bX_\sS^\top \bX_\sS\bS_\sS)^{-1} \bs_\sS}
{\bs_\sS^\top (\bS_\sS\bX_\sS^\top \bX_\sS\bS_\sS)^{-1} \bs_\sS}.
\end{equation}
That is, $\bbeta$ changes linearly with $\Sigma$
along the fixed direction
$
\bv \triangleq \frac{(\bX_\sS^\top \bX_\sS)^{-1} \bs_\sS}
{\bs_\sS^\top (\bX_\sS^\top \bX_\sS)^{-1} \bs_\sS}
$.
\end{mdframed}
Therefore,
$$
\Sigma_2 - \Sigma_1 = (\lambda_1 - \lambda_2) \bs_\sS^\top \bS_\sS \bG_\sS^{-1} \bone_\sS = (\lambda_1 - \lambda_2) \bone_\sS^\top \bG_\sS^{-1} \bone_\sS = (\lambda_1 - \lambda_2) \sigma_\sS^{-2},
$$
and hence
\begin{equation}\label{eq:beta_diff_final}
\widehatbbeta_\sS(\Sigma_2) - \widehatbbeta_\sS(\Sigma_1) 
= \bS_\sS \sigma_\sS^2 (\Sigma_2 - \Sigma_1) \bG_\sS^{-1} \bone_\sS 
= \bS_\sS \sigma_\sS (\Sigma_2 - \Sigma_1) \bw_\sS.
\end{equation}
Invoking $\Sigma_2 = \Sigma$ and $\Sigma_1 \to \Sigma_0$ gives~\eqref{equation:beta_sigma_linear} by the continuity of $\widehatbbeta(\Sigma)$, and implying the desired result \eqref{equation:haty_sigma_linear}. 
\end{proof}

Lemma~\ref{lemma:larslasso_ud_est}~\eqref{equation:haty_sigma_linear} shows that, for $\Sigma \in \sM$,  the constrained LASSO regression estimates $\widehatby(\Sigma)$ move in the LARS equiangular direction $\bu_\sS$ defined in \eqref{equation:ang_bisec_eq2}.
Moreover, the LASSO solution $\widehatbbeta(\Sigma)$ lies on the boundary of the diamond-shaped convex polytope (i.e., the ``$\ell_1$-norm" ball with radius $\Sigma$):
\begin{equation}
	\sB_1[\bzero, \Sigma] 
	= \left\{ \balpha \mid \sum \abs{\alpha_i} \leq \Sigma \right\},
\end{equation}
This set $\sB_1[\bzero, \Sigma]$  expands monotonically as $\Sigma$ increases. Lemma~\ref{lemma:larslasso_ud_est}~\eqref{equation:beta_sigma_linear} further shows  that, for $\Sigma \in \sM$, the regression coefficients $\widehatbbeta(\Sigma)$ move linearly along edge $\sS$ of the polytope, namely, the face where  $\beta_i = 0$ for all $i \notin \sS$.

\paragrapharrow{Sign restriction.}
Let $\widehatbbeta$ be a solution to the constrained LASSO problem \eqref{opt:lc}, and let $\widehatby = \bX \widehatbbeta$ denote the corresponding fitted values. 
It is straightforward to show that for any nonzero component $\widehatbeta_i$, its sign must match the sign $s_i$ of the current correlation $\widehatc_i = \bx_i^\top(\by - \widehatby)$,
\begin{equation}\label{equation:larslass_res1}
\sgn(\widehatbeta_i) = \sgn(\widehatc_i) = s_i.
\end{equation}

\begin{lemma}[Sign restriction of LARS--LASSO]\label{lamma:sign_corr_beta}
For any LASSO solution  $\widehatbbeta$, the current correlations satisfy
\begin{equation}\label{eq:c_i_sign}
\widehatc_i = \widehatC \cdot \sign(\widehatbeta_i), \quad \text{for } i \in \sS,
\end{equation}
where $\widehatc_i$ is the current correlation $\bx_i^\top(\by - \widehatby) = \bx_i^\top(\by - \bX \widehatbbeta)$. In particular, this implies that
\begin{equation}\label{eq:sign_equal}
\sign(\widehatbeta_i) = \sign(\widehatc_i), \quad \text{for } i \in \sS.
\end{equation}
\end{lemma}

\begin{proof}[of Lemma~\ref{lamma:sign_corr_beta}]
This follows immediately from the stationary condition \eqref{eq:larslasso_stat_condition010} by observing  that the $i$-th element of the left-hand side is $\widehatc_i$, and the right-hand side is $\lambda \cdot \sign(\widehatbeta_i)$ for $i \in \sS$. 
Moreover, $\lambda = \abs{\widehatc_i} = \widehatC$, which completes the proof.
\end{proof}

The standard LARS algorithm does not enforce the sign restriction \eqref{equation:larslass_res1}, so it cannot be used directly to compute LASSO solutions. However, it can be easily modified to incorporate this constraint.
\paragrapharrow{LASSO modification.}
Suppose we have just completed a LARS step, giving a new active set $\sS$ from the previous active set $\sS_-$, as in \eqref{equation:lars_upd1}, and that the current LARS estimate $\widehatby_{\sS_-}$ coincides with a valid LASSO solution $\widehatby = \bX \widehatbbeta$. 
At the beginning of this LARS step, we know (assume) that the sign restriction \eqref{equation:larslass_res1} holds.
However, this may no longer be true after moving further along the current equiangular direction.
To analyze this, define
\begin{equation}
\bw_\sS = \sigma_\sS \bG_\sS^{-1} \bone_\sS,
\end{equation}
a vector of length the size of $\abs{\sS}$.
Now construct the $p$-dimensional direction vector  
$\widehatbd$ by setting
$$
\widehatd_i=
\begin{cases}
s_i w_{\sS, i}, & i\in \sS;\\
0, & i\notin \sS.
\end{cases}
$$
Moving in the positive $\gamma$ direction along the LARS path \eqref{equation:lars_upd6}, 
the coefficient path is given by
\begin{equation}
\by(\gamma) = \bX \bbeta(\gamma), \quad \text{where } \beta_i(\gamma) \triangleq \widehatbeta_i + \gamma \widehatd_i, 
\text{ for $i \in \sS$},
\end{equation}
see also \eqref{eq:w_star_split33}.
Therefore, $\beta_i(\gamma)$ will change sign at
\begin{equation}
\gamma_i \triangleq -\widehatbeta_i / \widehatd_i;
\end{equation}
the first sign change happens at
\begin{equation}
\widetilde{\gamma} \triangleq  \min_{\gamma_i > 0} \{ \gamma_i \}, 
\end{equation}
say for covariate $\beta_{i_+}$.
If no such $\gamma_i$ exists (i.e., all $\gamma_i\leq 0$), we define  $\widetilde{\gamma}=\infty$.

If $\widetilde{\gamma}<\widehat{\gamma}_{\sS}$---where $\widehat{\gamma}_{\sS}$ is the stepsize to the next variable joining the active set, as defined in \eqref{equation:lars_upd5}---then the path $\beta_i(\gamma)$ ceases to be a valid LASSO solution for any $\gamma > \widetilde{\gamma}$ since the sign restriction~\eqref{equation:larslass_res1} would be violated: $\beta_i(\gamma)$  changes sign, but the corresponding correlation $c_i(\gamma)$ does not. 
This is because the continuous function $c_i(\gamma)$ cannot change sign within a single LARS step since $\abs{c_i(\gamma)} = \widehatC_{\sS} - \gamma \sigma_\sS > 0$ by \eqref{equation:lars_upd8}.
This leads to the following modification for computing LASSO paths:
\begin{remark}[LASSO modification]\label{remark:lass_modi}
If $\widetilde{\gamma} < \widehat{\gamma}_{\sS}$, stop the ongoing LARS step at $\gamma = \widetilde{\gamma}$ and remove the index $i_+$ (the one achieving $\widetildegamma$) from the calculation of the next equiangular direction. 
Specifically, update the fit and active set as
\begin{equation}\label{equation:larslasso_mod}
\widehatby_{\sS} = \widehatby_{\sS_-} + \widetilde{\gamma} \bu_\sS 
\qquad \text{and} \qquad 
\sS_+ = \sS - \{i_+\}
\end{equation}
instead of proceeding with the standard update \eqref{equation:lars_upd4}.
Otherwise, if  
$\widetilde{\gamma}\geq \widehat{\gamma}_{\sS}$, proceed with the usual LARS update \eqref{equation:lars_upd4}.
\end{remark}

In the original LARS algorithm, the active set $\sS$ grows monotonically. The LASSO modification, however, allows variables to be dropped from $\sS$ when their coefficients cross zero. 
The LASSO modification also assumes a ``one-at-a-time" condition
which means that at each step, at most one variable is either added to or removed from the active set. This condition typically holds for continuous (quantitative) data and can always be enforced by adding an infinitesimal random perturbation to the response vector $\by$ if necessary.
The complete procedure is summarized in Algorithm~\ref{alg:lars_lasso}.

\begin{theoremHigh}[LARS--LASSO relationship]\label{theorem:larslass_mod}
Under the above LASSO modification, and assuming the ``one-at-a-time'' condition, the LARS algorithm generates the entire LASSO solution path.
\end{theoremHigh}
\begin{proof}
See  \citet{efron2004least}.
\end{proof}

\begin{algorithm}[h] 
\caption{LARS algorithm with LASSO modification \citep{efron2004least}}
\label{alg:lars_lasso}
\begin{algorithmic}[1] 
\Require The data vector $\by \in \real^n$ and the input matrix $\bX \in \real^{n\times p}$. 

\State {\bfseries Initialization:}  $\sS_0 = \varnothing$, $\widehatby = \bzero$, and $\bX_{\sS_0} = \bX$;
\For{$t=1,2,\ldots$}

\State Compute the current correlation vector $\widehatbc_t = \bX^\top (\by - \widehatby_{t-1})$;

\State \algoalign{Update the active set $\sS_t = \sS_{t-1} \cup \left\{ i \mid \abs{\widehatc_{ti}} = C \right\}$ with $C = \max\left\{ \abs{\widehatc_{t1}}, \abs{\widehatc_{t2}}, \ldots \right\}$, and $\widehatc_{ti}$ being the $i$-th element of $\widehatbc_t$;}

\State Update the input matrix $\bX_{\sS_t} = [\ldots, s_i \bx_i, \ldots]_{i \in \sS_t}$, where $s_i = \sign(\widehatc_{ti})$;

\State \algoalign{Find the direction of the current minimum angle:
\begin{align*}
\bG_{\sS_t} &= \bX_{\sS_t}^\top \bX_{\sS_t} \in \real^{\abs{\sS_t} \times \abs{\sS_t}}, 
&&\sigma_{\sS_t} = (\bone_t^\top \bG_{\sS_t}^{-1} \bone_t)^{-1/2}, \\
\bw_{\sS_t} &= \sigma_{\sS_t}\bG_{\sS_t}^{-1}\bone_t \in\real^{\abs{\sS_t}}, 
&&\bu_t = \bX_{\sS_t} \bw_{\sS_t} \in \real^n.
\end{align*}}

\State Compute $\bz = \bX^\top \bu_t = [z_1, z_2, \ldots, z_p]^\top$ and estimate the coefficient vector
$$
\widehatbbeta_t = (\bX_{\sS_t}^\top \bX_{\sS_t})^{-1} \bX_{\sS_t}^\top \by
=
\bG_{\sS_t}^{-1}\bX_{\sS_t}^\top\by.
$$

\State 
\algoalign{
By \eqref{equation:lars_upd5}, compute
$$
\widehatgamma = \mathop{\text{min}^+}_{i \in \comple{\sS_t}} 
\left\{ \frac{C - \widehatc_{ti}}{\sigma_{\sS_t} - z_i}, \frac{C + \widehatc_{ti}}{\sigma_{\sS_t} + z_i} \right\}, 
\qquad
\widetildegamma = \min_{i \in \sS_t} \left\{ -\frac{z_i}{s_i w_{\sS_t,i}} \right\}^+,
$$
where $w_{\sS_t,i}$ is the $i$-th entry of $\bw_{\sS_t} = [w_1,w_2, \ldots]^\top$, and $(\cdot)^+$ denotes the positive minimum term. If there is no positive term then $ {\text{min}^+}\{\cdot\} = \infty$;
}

\If{$\widetildegamma < \widehatgamma$}
\State 
\algoalign{
The fitted vector $\widehatby_t$ and the active set $\sS_t$ are modified as follows:
$$
\widehatby_t = \widehatby_{t-1} + \widetildegamma \bu_t, 
\qquad
\sS_t = \sS_t \setminus \{i_+\},
$$
where the removed index $i_+$ is the index $i \in \sS_t$ such that $\widetildegamma$ is a minimum;
}
\Else
\State 
\algoalign{
The fitted vector $\widehatby_t$ and the active set $\sS_t$  are modified as follows:
$$
\widehatby_t = \widehatby_{t-1} + \widehatgamma \bu_t, 
\qquad
\sS_t = \sS_t \cup \{i_+\},
$$
where the added index $i_+$ is the index $i \in \comple{\sS_t}$ such that $\widehatgamma$ is a minimum;
}
\EndIf

\State Exit if some stopping criterion is satisfied;

\EndFor
\State \Return  $\bbeta = \widehatbbeta_t$;
\end{algorithmic} 
\end{algorithm}

\subsection*{LARS--LASSO Relationship}
We now briefly discuss Theorem~\ref{theorem:larslass_mod}, which establishes the equivalence between the modified LARS algorithm and the constrained LASSO problem. This equivalence stems from the fact that the sequence of coefficient estimates---the so-called solution path---generated by the modified LARS algorithm, as it proceeds step by step, coincides exactly with the set of solutions to the constrained LASSO problem as the constraint parameter $\Sigma$ increases from 0 to $\normonebig{\widehatbbeta^{OLS}}$, the $\ell_1$-norm of the ordinary least squares (OLS) solution.

The standard constrained LASSO problem is
$\min_{\bbeta} \frac{1}{2} \normtwo{\by - \bX\bbeta}^2 $ subvject to $ \normone{\bbeta} \leq \Sigma$, where $\Sigma$ is a nonnegative scalar that controls the total ``budget'' for the sum of the absolute values of the coefficients.
As $\Sigma$ increases from 0, the behavior of the solution $\widehatbbeta(\Sigma)$ evolves as follows:
\begin{itemize}
\item When $\Sigma = 0$, the only solution is $\bbeta = \bzero$.
\item As $\Sigma$ increases slightly, the optimal solution $\widehatbbeta(\Sigma)$ will have one or a few nonzero coefficients, chosen to reduce the least squares loss as much as possible per unit of ``$\ell_1$-cost."
\item As $\Sigma \to \infty$, the constraint becomes inactive, and $\widehatbbeta(\Sigma)$ approaches the OLS solution.
\end{itemize}
The {solution path} $\widehatbbeta(\Sigma)$ is the set of all optimal solutions as $\Sigma$ varies.

The equivalence between the modified LARS and constrained LASSO arises from the {optimality conditions (KKT conditions)} of the LASSO problem and how LARS satisfies them at each step.
At the optimal solution $\widehatbbeta$ for a given constraint $\Sigma$, the following KKT conditions (Theorem~\ref{theorem:opt_cond_sd}) must hold:
\begin{enumerate}
\item \textit{Subgradient condition.} The gradient of the least squares loss ($- \bX^\top(\by - \bX\widehatbbeta)$) must be in the subdifferential of $\lambda\normonebig{\widehatbbeta}$, where $\lambda$ is the Lagrangian dual variable. 
This translates to:
\begin{itemize}
\item For any $i$ where $\widehatbeta_i \neq 0$: $\absbig{\bx_i^\top(\by - \bX\widehatbbeta)} = \lambda$ (the absolute correlation between predictor $i$ and the residual is exactly $\lambda$).
\item For any $i$ where $\widehatbeta_i = 0$: $\absbig{\bx_i^\top(\by - \bX\widehatbbeta)} \leq \lambda$ (the absolute correlation is less than or equal to $\lambda$).
\end{itemize}
\item \textit{Complementary slackness.} $\lambda(\normonebig{\widehatbbeta} - \Sigma) = 0$. If $\normonebig{\widehatbbeta} < \Sigma$, then $\lambda = 0$ (no regularization). If $\lambda > 0$, then $\normonebig{\widehatbbeta} = \Sigma$ (the constraint is active).
\end{enumerate}
The modified LARS algorithm satisfies these KKT conditions by the following:
\begin{itemize}
\item \textit{Equal correlations.} At the start of each LARS step, the algorithm ensures that all predictors in the active set have \textbf{equal absolute correlation} with the current residual. This matches the KKT condition that nonzero coefficients must have $\abs{\bx_i^\top \br} = \lambda$ for some $\lambda$, where $\br$ denotes the residual.
\item \textit{Maximal correlation.} The algorithm also ensures that no inactive predictor has a higher absolute correlation with the residual than the active ones. This matches the KKT condition $\abs{\bx_i^\top \br} \leq \lambda$ for $i$ not in the active set.
\item \textit{Increasing $\Sigma$.} As LARS takes each step and moves the coefficients, the $\ell_1$-norm $\normone{\bbeta}$ increases. The point where LARS changes direction (when a new variable enters or an old one leaves) corresponds to a specific value of $\normone{\bbeta} = \Sigma_t$.
\item \textit{Path as $\Sigma$ increases.} The piecewise linear path generated by LARS, parameterized by the step or by $\normone{\bbeta}$, satisfies the KKT conditions for the constrained LASSO problem at every point along the path. The ``$\lambda$'' in the KKT conditions is implicitly determined by the common absolute correlation value of the active predictors at each step.
\end{itemize}

Importantly, neither the standard nor the modified LARS algorithm (Algorithm~\ref{alg:lars_lasso}) takes the constraint parameter $\Sigma$ as an explicit input.
Instead, it proceeds via geometric steps dictated by predictor--residual correlations. 
Nevertheless, at the end of each LARS step $t$, the resulting coefficient vector $\widehatbbeta_t$ is precisely the solution to the constrained LASSO problem with $\Sigma_t = \normonebig{\widehatbbeta_t}$:
$$
\min_{\bbeta\in\real^p} \frac{1}{2}\normtwo{\by - \bX\bbeta}^2 \quad\text{ s.t. }\quad
\normone{\bbeta} \leq \Sigma_t.
$$
In other words, the LARS path \emph{is} the LASSO solution path. Although $\Sigma$ is not an input to the algorithm, it serves as the natural parameter that indexes the path: each point on the LARS trajectory corresponds to a unique value of $\Sigma=\normone{\bbeta}$. Consequently, the full regularization path is obtained automatically, and one can subsequently select the optimal $\Sigma$---for example, via cross-validation.

\section{Algorithms for Lagrangian LASSO Problems}\label{section:algo_laglasso}
We introduce several algorithms for solving the constrained LASSO problem~\eqref{opt:lc}.
We now turn to the Lagrangian (or penalized) form of the LASSO problem. This formulation is more commonly used in machine learning due to its simplicity: it casts the problem as an unconstrained optimization with a single objective function, which is algorithmically easier to handle than a constrained problem like~\eqref{opt:lc}. Moreover, it integrates naturally into regularization frameworks.

In machine learning, regularization is a fundamental technique for preventing overfitting and promoting desirable model properties---such as sparsity; see Section~\ref{section:ls_regular}. 
The Lagrangian LASSO explicitly captures this trade-off between data fidelity (measured by the loss) and model complexity (penalized via the $\ell_1$-norm).

The standard form of the Lagrangian LASSO \eqref{opt:ll} (p.~\pageref{opt:ll}) is:
\begin{equation}\label{equation:lag_lasso_alg}
\min_{\bbeta\in \real^p} \quad F(\bbeta) \triangleq f(\bbeta) +g(\bbeta)\triangleq \frac{1}{2} \normtwo{\by - \bX \bbeta}^2 + \lambda \normone{\bbeta},
\end{equation}
where $g(\bbeta)\triangleq\lambda \normone{\bbeta}$.
\footnote{Note again that the Lagrangian LASSO is equivalent to the penalized $\ell_1$-minimization~\eqref{opt:p1_penalize} for sparse signal recovery problems.} 
Since both $f$ and $g$ are  convex (Exercise~\ref{exercise:conv_quad}), one might consider applying gradient descent as a straightforward approach.
However, the objective function $F(\bbeta)$ is non-smooth at points where any component of $\bbeta$ is zero (due to the nondifferentiability of the $\ell_1$-norm). 
Consequently, the classical gradient does not exist everywhere, and standard gradient descent cannot be applied directly.

To address this, subgradient methods can be employed. While $f$ is smooth---specifically, $\normtwo{\bX^\top\bX}$-smooth (see Example~\ref{example:lipschitz_spar})---and $g$ is Lipschitz continuous, their sum $F=f+g$ is neither smooth nor differentiable. As a result, subgradient methods applied to $F$ generally exhibit slow convergence and lack strong theoretical guarantees (e.g., no linear or even accelerated rates without additional structure).

\subsection{Smoothing with Huber Loss}
\begin{SCfigure}
\centering
\includegraphics[width=0.5\textwidth]{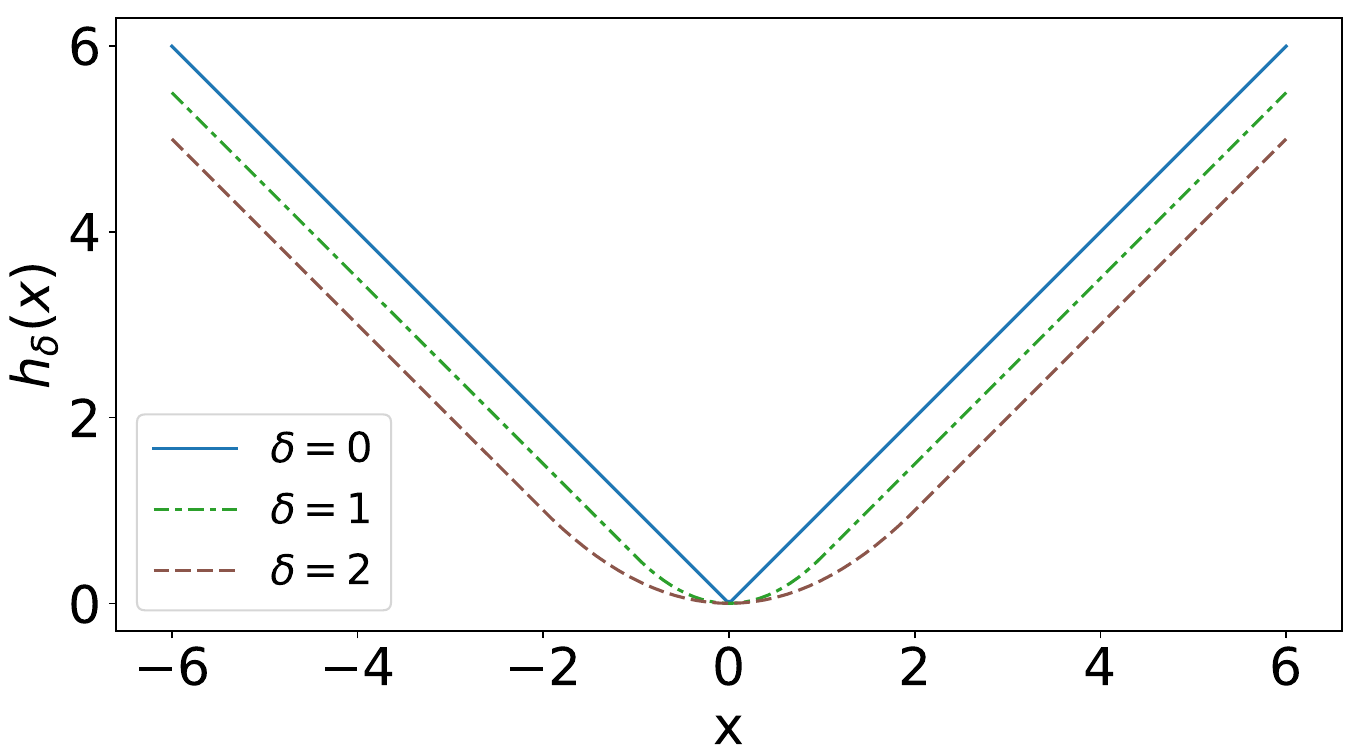}
\caption{Illustration of Huber loss function $h_\delta(x)$ for different values of $\delta$.}
\label{fig:huber}
\end{SCfigure}

An alternative approach to handling the non-smoothness of the LASSO objective is to recognize that the $\ell_1$-norm, $\normone{\bbeta}=\sum_{i=1}^{p}\abs{\beta_i}$, is a sum of absolute-value functions---each of which is nondifferentiable at zero. If we replace each absolute value with a smooth approximation, then standard gradient-based optimization methods become applicable.
One common choice for such an approximation is the \textit{Huber loss function}, defined as:
\begin{equation}
h_\delta(x) =
\begin{cases}
\frac{1}{2\delta}x^2, & \abs{x} < \delta; \\
\abs{x} - \frac{\delta}{2}, & \text{otherwise}.
\end{cases}
\end{equation}
When $\delta \to 0$, the smooth function $h_\delta(x)$ approaches the absolute value function $\abs{x}$. Figure~\ref{fig:huber} illustrates  the graph of $h_\delta(x)$ for various  values of $\delta$.

By substituting $\abs{\beta_i}$ with $h_\delta(\beta_i)$, we obtain a \textit{smoothed LASSO} problem:
\begin{equation}
\min_{\bbeta\in\real^p} \quad F_\delta(\bbeta) = \frac{1}{2}\normtwo{\bX\bbeta-\by}^2 + \lambda H_\delta(\bbeta),
\quad \text{with }H_\delta(\bbeta) \triangleq \sum_{i=1}^{p} h_\delta(\beta_i),
\end{equation}
where $\delta$ is the given smoothing parameter. 
The objective $F_\delta(\bbeta)$ is now continuously differentiable, and its gradient can be computed as:
\begin{equation} 
\nabla F_\delta(\bbeta) = \bX^\top(\bX\bbeta - \by) + \lambda \nabla H_\delta(\bbeta), 
\end{equation}
where $\nabla H_\delta(\bbeta)$ is defined component-wise by:
\begin{equation}
(\nabla H_\delta(\bbeta))_i =
\begin{cases}
\sign(\beta_i), & \abs{\beta_i} > \delta; \\
{\beta_i}/{\delta}, & \abs{\beta_i} \leq \delta,
\end{cases}
\quad i\in\{1,2,\ldots,p\}.
\end{equation}
Since both terms in $F_\delta(\bbeta)$ are convex and the gradient $\nabla F_\delta(\bbeta)$ is Lipschitz continuous, the function $F_\delta(\bbeta)$ is $L_b$-smooth with Lipschitz constant
$$
L_b = \normtwo{\bX^\top \bX} + \frac{\lambda}{\delta}
$$ 
(see Example~\ref{example:lipschitz_spar} and Exercise~\ref{exercise:sum_sc_conv}).
Therefore, by Problem~\ref{prob:pgd_smooth}, gradient descent with a constant stepsize $\eta=1/L_b $ is guaranteed to converge to the global minimum of $F_\delta$.

However, when $\delta$ is very small---i.e., when the smoothed problem closely approximates the original LASSO---the Lipschitz constant $L_b$ becomes large due to the $\lambda/\delta$ term. In this regime, the required stepsize $\eta=1/L_b$ becomes very small, which can lead to slow convergence in practice. Thus, careful tuning of $\delta$ (or the use of adaptive stepsize strategies) is essential when using this smoothing approach.

\subsection{Proximal Gradient and Minorization-Maximization Methods\index{Minorization-maximization method}}\label{section:prox_gd_lag_lasso}

Given that both functions $f$ and $g$   in the Lagrangian LASSO problem \eqref{equation:lag_lasso_alg} are convex, where $f$ is differentiable and $g$ is non-differentiable, problem \eqref{equation:lag_lasso_alg} can be considered  a composite model.
We have already introduced the proximal gradient method to address such composite problems  (Section~\ref{section:proxiGD_inClasso}).
Although alternative approaches exist---such as generalized conditional gradient and generalized mirror descent methods \citep{lu2025practical}---these typically require solving, at each iteration, a subproblem of the form  $\mathop{\argmin}_{\bbeta\in {\real^p}} \innerproductbig{\nabla f(\bbeta^\toptzero), \bbeta} + g(\bbeta)$. 
In the LASSO setting, this subproblem is essentially as difficult as the original problem~\eqref{equation:lag_lasso_alg}.
Therefore, the proximal gradient method remains the most practical and efficient choice.

Following Algorithm~\ref{alg:prox_gd_gen}, the update rule for the proximal gradient method applied to the Lagrangian LASSO problem~\eqref{equation:lag_lasso_alg} at the $t$-th iteration is:
\begin{subequations}\label{equation:prox_lagranglso}
\begin{align}
\balpha^\toptone  
&\leftarrow \bbeta^\toptzero - \eta_t \bX^\top (\bX\bbeta^\toptzero - \by);\\
\bbeta^\toptone 
&\leftarrow \prox_{\eta_t g}(\balpha^\toptone)
= \mathcalT_{\eta_t\lambda} (\balpha^\toptone)
= \big[\absbig{\balpha^\toptone} - \eta_t \lambda \bone\big]_+ \hadaprod \sign(\balpha^\toptone),
\end{align}
\end{subequations}
where $\mathcalT_\lambda(\cdot)$ denotes the soft-thresholding operator (Example~\ref{example:soft_thres}).
This two-step procedure first performs a standard gradient descent step on the smooth loss $f$, followed by a shrinkage (thresholding) operation that promotes sparsity in the iterates.
Convergence of this scheme is guaranteed by Theorem~\ref{theorem:prox_conv_ss_cvx} when using a constant stepsize $\eta_t = \frac{1}{\normtwo{\bX^\top\bX}}$ since $f$ is $\normtwo{\bX^\top\bX}$-smooth (Example~\ref{example:lipschitz_spar}).

Note that this update closely resembles the proximal gradient step for the constrained LASSO problem in~\eqref{equation:prox_constrainedlso}; however, in the constrained formulation, the regularization parameter $\lambda$ is not fixed but instead varies adaptively with the constraint level at each iteration.

\index{Iterative shrinkage-thresholding algorithm (ISTA)}
Finally, when the proximal gradient method is applied specifically to the Lagrangian LASSO problem, it is commonly referred to in the literature as the \textit{iterative shrinkage-thresholding algorithm (ISTA)}---a name that reflects how each iteration shrinks small coefficients toward zero, thereby encouraging sparse solutions.

\subsection*{Minorization-Maximization Method}
We further turn to a class of methods, known either as \textit{minorization-maximization} or \textit{majorization-minimization (MM)} algorithms.
These are particularly useful for nonconvex problems and belong to the family of \textit{auxiliary-variable methods}: they introduce an auxiliary function that locally upper-bounds (i.e., majorizes) the original objective, turning a difficult minimization into a sequence of simpler ones.

Although MM methods apply to both constrained and unconstrained settings, we describe them here for an unconstrained problem of the form $\min_{\bbeta \in \real^p} f(\bbeta)$, where $f : \real^p \to \real$ may be nonconvex.

\begin{definition}[Auxiliary function (majorizer)\index{Auxiliary function}\index{Majorizer}]\label{definition:aux_func}
A function $\Psi(\bbeta, \widetildebbeta)$ is called an \textit{auxiliary function} for $F(\bbeta)$ (or a \textit{majorizer} of $F$ at $\widetildebbeta$) if it satisfies:
\footnote{$\bbeta$ can be scalars, vectors, or matrices.}
\begin{equation}\label{equation:aux_func}
\Psi(\bbeta, \widetildebbeta) \geq F(\bbeta)
\quad 
\text{and}
\quad 
\Psi(\bbeta, \bbeta) = F(\bbeta).
\end{equation}
In other words, the auxiliary function $\Psi(\cdot, \widetildebbeta)$ is an upper bound of $F(\cdot)$, and the bound is tight at $\widetildebbeta=\bbeta$.
\end{definition}

\begin{lemma}[Nonincreasing in auxiliary functions]\label{lemma:noninmuaux}
If $ \Psi $ is an auxiliary function for $F$, then the sequence generated by
\begin{equation}\label{equation:aux_update}
\bbeta^\toptone = \mathop{\argmin}_{\bbeta} \, \Psi(\bbeta, \bbeta^\toptzero),
\quad \text{for } t = 0,1,2,\ldots
\end{equation}
satisfies  $F(\bbeta^\toptone)\leq F(\bbeta^\toptzero)$.
\end{lemma}
\begin{proof}[of Lemma~\ref{lemma:noninmuaux}]
By definition of the auxiliary function and the update rule: $ F(\bbeta^\toptone) \leq \Psi(\bbeta^\toptone, \bbeta^\toptzero) \leq \Psi(\bbeta^\toptzero, \bbeta^\toptzero) = F(\bbeta^\toptzero)$.
\end{proof}

An MM algorithm proceeds by initializing $\bbeta^\topzero$ and iteratively applying the update~\eqref{equation:aux_update}. 
The objective value decreases monotonically:
\begin{equation}\label{equation:mm_update}
\bbeta^\toptone 
= \argmin_{\bbeta \in \real^p} \Psi\left(\bbeta, \bbeta^\toptzero\right), \quad \text{for } t = 0,1,2,\ldots
\end{equation}
Note that $ F(\bbeta^\toptone) = F(\bbeta^\toptzero)$ only if $ \bbeta^\toptzero$ is a local minimum of $ \Psi(\bbeta, \bbeta^\toptzero)$ w.r.t. $\bbeta$. If the partial derivatives of $ F$ exist and are continuous in a small neighborhood of $ \bbeta^\toptzero$, this also implies that the gradient $ \nabla F(\bbeta^\toptzero) = \bzero$. Thus, by iterating the update in \eqref{equation:aux_update}, we obtain a sequence of estimates that converge to a local minimum $ \bbeta^* = \argmin_{\bbeta} F(\bbeta)$ of the objective function:
\begin{equation}
F(\bbeta^{(0)}) \geq 	 F(\bbeta^\topone) \geq  F(\bbeta^{(2)})\geq \ldots  \geq  F(\bbeta^\toptzero) \geq  F(\bbeta^\toptone)\geq \ldots \geq F(\bbeta^*).
\end{equation}
Definition~\ref{definition:aux_func} finds a majorizer $\Psi$ of $F$, and Lemma~\ref{lemma:noninmuaux} shows the minimization property in $\Psi$, hence the algorithm is often referred to as the \textit{majorization-minimization (MM) framework}.
The key advantage of the MM framework is that if the surrogate function $\Psi(\cdot, \bbeta^\toptzero)$ admits a closed-form minimizer or can be minimized efficiently, then each iteration becomes computationally tractable---even when direct minimization of $F$ is not.
Moreover, if $F$ is strictly convex, the MM algorithm converges to the global minimizer.

\paragrapharrow{Proximal gradient as an MM algorithm for Lagrangian LASSO problems.}

Recall from Sections~\ref{section:proxiGD_inClasso} and \ref{section:prox_gd_lag_lasso} that the proximal gradient method applies to objectives of the form $F = f + g$, where $f=\frac{1}{2}\normtwo{\by - \bX \bbeta}^2$ is convex and differentiable, and $g=\lambda \normone{\bbeta}$ is convex and  nondifferentiable. By applying the linear approximation theorem (Theorem~\ref{theorem:linear_approx}) to $f$, we obtain
\begin{align*}
F(\bbeta) 
&= f(\bbeta) + g(\bbeta)
= f(\widetildebbeta) + \innerproduct{\nabla f(\widetildebbeta), \widetildebbeta - \bbeta}
+ \frac{1}{2} \innerproduct{\widetildebbeta - \bbeta, \nabla^2 f(\bxi) (\widetildebbeta - \bbeta)} + g(\bbeta),
\end{align*}
where $\bxi = \gamma\bbeta + (1-\gamma)\widetildebbeta$ for some $\gamma \in [0,1]$. 
Since $f$ is $L_b$-smooth with $L_b\triangleq \normtwo{\bX^\top\bX}$ (Example~\ref{example:lipschitz_spar}), 
we have an upper bound on the Hessian, namely $\nabla^2 f(\bxi) \preceq L_b\bI_{p\times p}$, from which we obtain the inequality
$$
F(\bbeta) 
\leq 
f(\widetildebbeta) + \innerproduct{\nabla f(\widetildebbeta), \widetildebbeta - \bbeta} 
+ \frac{L_b}{2}\normtwo{\widetildebbeta - \bbeta}^2 + g(\bbeta)
\triangleq {\Psi(\bbeta,\widetildebbeta)},
$$
with equality holding when $\widetildebbeta = \bbeta$. 
At $\widetildebbeta = \bbeta^\toptzero$, the update at the $t$-th iteration becomes:
\begin{align}
\bbeta^\toptone 
&= \argmin_{\bbeta \in \real^p} \Psi\left(\bbeta, \bbeta^\toptzero\right)\\
&\implies \frac{\lambda}{L_b} \sign(\bbeta) + \bbeta  =\bbeta^\toptzero - \frac{1}{L_b} \bX^\top (\bX\bbeta^\toptzero - \by),
\end{align}
which is precisely the proximal gradient step with stepsize $\eta_t = {1}/{L_b}$.
Thus, the proximal gradient method for the Lagrangian LASSO can be interpreted as an instance of the MM framework, where the surrogate function $\Psi$ is constructed via a quadratic upper bound on the smooth part $f$.

\noindent
\begin{minipage}[t]{0.505\linewidth}
\begin{algorithm}[H]
\caption{FISTA-V1}
\label{alg:fistav1}
\begin{algorithmic}[1]
\State $\balpha^\topone = \bbeta^\topone \in \real^p$, $\gamma_1=\frac{1}{\zeta_1}=1$;
\For{$t=1,2,\ldots$}
\State Choose $\eta_{t}$ and $\zeta_t=\frac{1}{\gamma_t}$;
\State $\bbeta^\toptone \leftarrow \prox_{\eta_{t} g}(\balpha^\toptzero - \eta_{t} \nabla f(\balpha^\toptzero))$;
\State $\small\begin{aligned}
&\balpha^\toptone \leftarrow \bbeta^\toptone + \frac{\zeta_t-1}{\zeta_{t+1}} (\bbeta^\toptone - \bbeta^\toptzero);\\
&= \bbeta^\toptone + (\frac{\gamma_{t+1}}{\gamma_t}-\gamma_{t+1}) (\bbeta^\toptone - \bbeta^\toptzero);
\end{aligned}$
\EndFor
\State Output  $\bbeta_{\text{final}}\leftarrow \bbeta^{(T)}$;
\end{algorithmic}
\end{algorithm}
\end{minipage}%
\hfil 
\begin{minipage}[t]{0.495\linewidth}
\begin{algorithm}[H]
\caption{FISTA-V2}
\label{alg:fistav2}
\begin{algorithmic}[1]
\State $\bu^\topone =\balpha^\topone = \bbeta^\topone \in \real^p$, $\gamma_1=1$;
\For{$t=1,2,\ldots$}
\State Choose $\eta_{t}$ and $\gamma_t$;
\State $\bbeta^\toptone \leftarrow \prox_{\eta_{t} g}(\balpha^\toptzero - \eta_{t} \nabla f(\balpha^\toptzero))$;
\State $\bu^\toptone \leftarrow \bbeta^\toptzero + \frac{1}{\gamma_{t}}(\bbeta^\toptone - \bbeta^\toptzero)$;
\State $\balpha^\toptone \leftarrow (1 - \gamma_{t+1})\bbeta^\toptone + \gamma_{t+1} \bu^\toptone$;
\EndFor
\State Output  $\bbeta_{\text{final}}\leftarrow \bbeta^{(T)}$;
\end{algorithmic}
\end{algorithm}
\end{minipage}

\subsection{Fast Proximal Gradient Method---FISTA\index{FISTA}}\label{section:fista}
The \textit{fast proximal gradient} method, also known as the \textit{fast iterative shrinkage-thresholding algorithm (FISTA)} (Algorithm~\ref{alg:fistav2}) is an accelerated variant of the proximal gradient method introduced by \citet{beck2009fastgradient}. 
Recall that the standard proximal gradient method achieves a convergence rate of $\mathcalO\left({1}/{T}\right)$ for minimizing a sum of a smooth convex function and a non-smooth convex regularizer (Theorem~\ref{theorem:prox_conv_ss_cvx}). 
A natural question is whether this rate can be improved. FISTA addresses this by incorporating momentum, leading to a faster convergence $\mathcalO\left({1}/{T^2}\right)$.

FISTA operates in two steps at each iteration $t$:
\begin{itemize}
\item It constructs an extrapolated point $\balpha^\toptzero $ using information from the two most recent iterates.
\item It performs a proximal gradient step at this extrapolated point.
\end{itemize} 
Specifically, the update takes the form:
\begin{subequations}\label{equation:fista_intui}
\begin{align}
\balpha^\toptzero   &\leftarrow \bbeta^\toptzero + \frac{t-2}{t+1}(\bbeta^\toptzero - \bbeta^{(t-1)});\\
\bbeta^\toptone &\leftarrow \prox_{\eta_t g}(\balpha^\toptzero - \eta_t \nabla f(\balpha^\toptzero)).
\end{align}
\end{subequations}
The key distinction from the standard proximal gradient method is the  \textit{momentum term}, $\frac{t-2}{t+1}(\bbeta^\toptzero - \bbeta^{(t-1)})$, which leverages historical iterates to accelerate convergence.
Crucially, this acceleration incurs negligible additional computational cost per iteration, yet yields significantly faster empirical and theoretical performance.
If the stepsize $\eta_t$ is fixed and satisfies $\eta_t={1}/{L_b}$ (where $f$ is assumed to be $L_b$-smooth), then FISTA guarantees a convergence rate of  $\mathcalO\left({1}/{T^2}\right)$; see Theorem~\ref{theorem:fista_conv_ssf}. 
With this intuition, the complete FISTA algorithm is presented in Algorithm~\ref{alg:fistav1}, where we define $\zeta_t=\frac{1}{\gamma_t}$ for all $t>0$. 
When $\frac{1}{\gamma_t}=\zeta_t \triangleq \frac{t+1}{2}$, Algorithm~\ref{alg:fistav1} is equivalent to the intuitive update in~\eqref{equation:fista_intui}.

To facilitate theoretical analysis, an equivalent reformulation of FISTA is often used---shown in Algorithm~\ref{alg:fistav2}. While Algorithm~\ref{alg:fistav1} makes the momentum mechanism explicit and is thus more interpretable, Algorithm~\ref{alg:fistav2} structures the updates in a way that simplifies the convergence proof (see Theorem~\ref{theorem:fista_conv_ssf}).

\paragrapharrow{Convergence requirement.}
To achieve the accelerated $\mathcalO\left({1}/{T^2}\right)$  rate, the stepsizes $\eta_t$ 
and $\gamma_t = \frac{1}{\zeta_t}$ must satisfy certain conditions. Specifically, the following three requirements ensure the desired convergence behavior:
\begin{subequations}\label{equation:fista_require}
\small
\begin{align}
&f(\bbeta^\toptzero) \leq f(\balpha^\toptzero) + \innerproduct{\nabla f(\balpha^\toptzero), \bbeta^\toptzero - \balpha^\toptzero} + \frac{1}{2\eta_t} \normtwo{\bbeta^\toptzero - \balpha^\toptzero}^2, \text{ i.e., $\eta_t \leq \frac{1}{L_b}$};  \label{equation:fista_require1}\\
&\gamma_1 = 1, \quad \frac{(1 - \gamma_{t+1})\eta_{t+1}}{\gamma_{t+1}^2} \leq \frac{\eta_{t}}{\gamma_{t}^2}, \quad t > 1; \label{equation:fista_require2}\\
&\frac{\gamma_t^2}{\eta_t} = \mathcalO\left(\frac{1}{t^2}\right). \label{equation:fista_require3}
\end{align}
\end{subequations}
\paragrapharrow{Choice of $\gamma_t = \frac{1}{\zeta_t}$.}
A common and effective choice is to set $\eta_t = \frac{1}{L_b}$ and $\frac{1}{\zeta_t}=\gamma_t = \frac{2}{t+1}$ for $t>0$. 
This selection satisfies all conditions in~\eqref{equation:fista_require} and recovers the intuitive FISTA update in~\eqref{equation:fista_intui}.
However, this choice is not unique. Another valid sequence is defined recursively by:
$$
\small
\begin{aligned}
\gamma_1 = 1, \quad \frac{1}{\gamma_{t+1}} = \frac{1}{2} \left( 1 + \sqrt{1 + \frac{4}{\gamma_{t}^2}} \right)
\iff 
\zeta_{t+1} = \frac{1 + \sqrt{1 + 4\zeta_t^2}}{2},
\quad t>0.
\end{aligned}
$$
where $\{\gamma_t\} = \{\frac{1}{\zeta_t}\}$. 
This sequence also yields an $\mathcalO\left({1}/{T^2}\right)$ convergence rate.
The following lemma confirms that this recursive choice satisfies condition in~\eqref{equation:fista_require}.

\begin{lemma}\label{lemma:seq_gamma_zeta}
Let $\{\frac{1}{\gamma_t}\}_{t>0}=\{\zeta_t\}_{t> 0}$ be the sequence (see Algorithm~\ref{alg:fistav1}) defined by
$ \zeta_1 = 1, \; \zeta_{t+1} = \frac{1 + \sqrt{1 + 4\zeta_t^2}}{2},  t \geq 0. $
Then, $\zeta_t \geq \frac{t+1}{2}$ for all $t > 0$.
\end{lemma}
\begin{proof}[of Lemma~\ref{lemma:seq_gamma_zeta}]
We proceed by induction. 
Obviously, for $t = 1$, $\zeta_1 = 1 \geq \frac{1+1}{2}$, so the claim holds. 
Suppose that the claim holds for $t>1$ such that $\zeta_t \geq \frac{t+1}{2}$. 
We will prove that $\zeta_{t+1} \geq \frac{t+2}{2}$. By the recursive relation defining the sequence and the induction assumption,
$ \zeta_{t+1} = \frac{1 + \sqrt{1 + 4\zeta_t^2}}{2} \geq \frac{1 + \sqrt{1 + (t+1)^2}}{2} \geq \frac{1 + \sqrt{(t+1)^2}}{2} = \frac{t+2}{2}$. 
This completes the induction.
\end{proof}

\paragrapharrow{Monotone issue.}
The original FISTA algorithm is not a descent method---that is, the objective value
$F(\bbeta^\toptzero)$ may increase at some iterations.
To address this, we present a monotone (descent) variant of FISTA that requires only a minor modification to Algorithm~\ref{alg:fistav1}:
instead of always accepting the new iterate $\bbeta^\toptone$, we choose
$$
\bbeta^\toptone\in 
\argmin \left\{F(\bbeta)\mid \bbeta\in\{\bbeta^\toptzero, \bxi^\toptzero\}\right\}
$$
ensuring that $F(\bbeta^\toptone) \leq \min\{F(\bxi^\toptzero), F(\bbeta^\toptzero)\}$.
This modified procedure is given in Algorithm~\ref{alg:fistav_monotone}.
In practice, after computing the proximal update $\bxi^\toptzero$, we compare its objective value with that of the current iterate $\bbeta^\toptzero$. 
Only if $F(\bxi^\toptzero)\leq F(\bbeta^\toptzero)$ do we accept $\bxi^\toptzero$ as the next iterate; otherwise, we keep $\bbeta^\toptzero$ unchanged.
More compactly, Steps 6 to 10 of Algorithm~\ref{alg:fistav_monotone} can be equivalently expressed as:
$$
\balpha^\toptone = \bbeta^\toptone + \frac{\zeta_t}{\zeta_{t+1}} (\bxi^\toptzero - \bbeta^\toptone) + \frac{\zeta_t - 1}{\zeta_{t+1}} (\bbeta^\toptone - \bbeta^\toptzero).
$$

\noindent
\begin{minipage}[h]{1\linewidth}
\begin{algorithm}[H]
\caption{Monotine FISTA \citep{beck2009fastgradient}, Compare to Algorithm~\ref{alg:fistav1}}
\label{alg:fistav_monotone}
\begin{algorithmic}[1]
\Require A function $f(\bbeta)$ and a closed convex function $g$ (usually non-smooth) satisfying (A1) and (A2) in Theorem~\ref{theorem:fista_conv_ssf}; 
\State {\bfseries initialize:} $\balpha^\topone = \bbeta^\topone \in \real^p$, $\gamma_1=\frac{1}{\zeta_1}=1$;
\For{$t=1,2,\ldots$}
\State Choose $\eta_{t}$ and $\zeta_t=\frac{1}{\gamma_t}$;
\State $\bxi^\toptzero \leftarrow \prox_{\eta_{t} g}(\balpha^\toptzero - \eta_{t} \nabla f(\balpha^\toptzero))$;
\State Choose $\bbeta^\toptone$ such that $F(\bbeta^\toptone) \leq \min\{F(\bxi^\toptzero), F(\bbeta^\toptzero)\}$;
\If{$\bbeta^\toptone$ is $\bxi^\toptzero$}
\State $\small\begin{aligned}
\balpha^\toptone \leftarrow \bbeta^\toptone + \frac{\zeta_t-1}{\zeta_{t+1}} (\bxi^\toptzero - \bbeta^\toptzero)
\end{aligned}$ \Comment{Same as FISTA}
\ElsIf{$\bbeta^\toptone$ is $\bbeta^\toptzero$}
\State $\small\begin{aligned}
\balpha^\toptone \leftarrow \bbeta^\toptzero + \frac{\zeta_t}{\zeta_{t+1}} (\bxi^\toptzero - \bbeta^\toptzero);
\end{aligned}$
\EndIf
\State Stop if a stopping criterion is satisfied at iteration $t=T$;
\EndFor
\State \Return  $\bbeta_{\text{final}}\leftarrow \bbeta^{(T)}$;
\end{algorithmic}
\end{algorithm}
\end{minipage}%
\hfil 
\begin{minipage}[t]{1\linewidth}
\begin{algorithm}[H]
\caption{Line Search Algorithm at $t$-th Iteration}
\label{alg:fistav_line}
\begin{algorithmic}[1]
\State {\bfseries input:}  $\eta_t = \eta_{t-1} > 0$, $\rho < 1$. Reference point $\balpha^\toptzero$ and its gradient $\nabla f(\balpha^\toptzero)$;
\State Calculate candidate update $\bbeta^\toptone \leftarrow \prox_{\eta_t g}(\balpha^\toptzero - \eta_t \nabla f(\balpha^\toptzero))$;
\While {condition \eqref{equation:fista_require1} is not satisfied for $\bbeta^\toptone, \balpha^\toptzero$}
\State Reducing stepsize $\eta_t \gets \rho \eta_t$;
\State Recalculate $\bbeta^\toptone \leftarrow \prox_{\eta_t g}(\balpha^\toptzero - \eta_t \nabla f(\balpha^\toptzero))$;
\EndWhile
\State \Return update $\bbeta^\toptone$, stepsize $\eta_t$;
\end{algorithmic}
\end{algorithm}
\end{minipage}

\paragrapharrow{Line search methods.\index{Line search}}
In Algorithms~\ref{alg:fistav1} and \ref{alg:fistav2}, the stepsize must satisfy $\eta_t \leq \frac{1}{L_b}$ to ensure condition~\eqref{equation:fista_require1} holds. 
However, in many practical problems, the Lipschitz constant $L_b$ of  $\nabla f$ is unknown. 
In such cases, we can still satisfy condition~\eqref{equation:fista_require1} by using a backtracking line search to adaptively select $\eta_t$ at each iteration. Meanwhile, the sequence $\{\gamma_t\}$ can be chosen to simultaneously satisfy conditions~\eqref{equation:fista_require2} and~\eqref{equation:fista_require3}, thereby preserving the accelerated convergence rate of $\mathcalO\left({1}/{T^2}\right)$.

Algorithm~\ref{alg:fistav_line} implements such a line search strategy. At iteration $t$, it initializes $\eta_t$ with the previous stepsize $\eta_{t-1}$ and then decreases $\eta_t$ geometrically until condition~\eqref{equation:fista_require1} is met. Since this condition is guaranteed to hold for sufficiently small $\eta_t$ (due to the smoothness of $f$), the line search always terminates in finite steps.
Moreover, it is straightforward to verify that the other two conditions---\eqref{equation:fista_require2} and~\eqref{equation:fista_require3}---remain satisfied throughout the iterative process when $\{\gamma_t\}$ is updated appropriately (e.g., via the recursive rule in Lemma~\ref{lemma:seq_gamma_zeta}).

The following theorem establishes the convergence rate for FISTA with a fixed stepsize. The analysis for the adaptive (line search) variant follows analogously and is omitted for brevity.
\begin{theoremHigh}[FISTA for convex and SS $f$, $\mathcalO(1/T^2)$:  Compare to Theorem~\ref{theorem:prox_conv_ss_cvx}]\label{theorem:fista_conv_ssf}
Let $ f $ be a proper {closed} convex and $L_b$-smooth function, and let  $ g $ be a proper closed and convex function. 
Assume the following conditions hold:
\begin{itemize}
	\item[(A1)] $ g: \real^p \rightarrow (-\infty, \infty] $ is proper, closed, and convex.
	\item[(A2)] $ f: \real^p \rightarrow (-\infty, \infty] $ is proper, closed, and convex, $ \dom(f) $ is convex, $ \dom(g) \subseteq \textcolor{mylightbluetext}{\interior}(\dom(f)) $, and $ f $ is $ L_b $-smooth over $ \textcolor{mylightbluetext}{\interior}(\dom(f)) $.
\end{itemize}
\noindent
Let $ \{\bbeta^\toptzero\}_{t > 0} $ be the sequence generated by the FISTA method (Algorithm~\ref{alg:fistav2}) for solving problem $\min \{F(\bbeta) \triangleq f(\bbeta)+g(\bbeta)\}$  with a constant stepsize rule in which $\eta_t \triangleq \frac{1}{L_b}$ for all $ t > 0 $. Then, for any optimizer $ \bbeta^*$ of $F\triangleq f+g$ and $ T > 1 $, it follows that
$$
F(\bbeta^{(T)}) - F(\bbeta^*) \leq \frac{2L_b}{T^2} \normtwo{\bbeta^\topone - \bbeta^*}^2.~\footnote{If the Lipschitz constant $L_b$ is not known beforehand, line search Algorithm~\ref{alg:fistav_line} can be applied to achieve a same rate of convergence.}
$$
\end{theoremHigh}
\begin{proof}[of Theorem~\ref{theorem:fista_conv_ssf}]
Since $\bbeta^\toptone = \prox_{\eta_t g}(\balpha^\toptzero - \eta_t \nabla f(\balpha^\toptzero))$, by the Proximal Property-I (Lemma~\ref{lemma:prox_prop1}), we have
$$
-\bbeta^\toptone + \balpha^\toptzero - \eta_t \nabla f(\balpha^\toptzero) \in \eta_t \partial g(\bbeta^\toptone).
$$
By the subgradient inequality (Definition~\ref{definition:subgrad}), for any $\bbeta\in\real^p$, we have
$$
\eta_t g(\bbeta) \geq \eta_t g(\bbeta^\toptone) + \innerproduct{-\bbeta^\toptone + \balpha^\toptzero - \eta_t \nabla f(\balpha^\toptzero), \bbeta - \bbeta^\toptone}.
$$
Using the $L_b$-smoothness of $f$ (Definition~\ref{definition:scss_func}) and $\eta_t = \frac{1}{L_b}$, we obtain
$$
f(\bbeta^\toptone) \leq f(\balpha^\toptzero) + \innerproduct{\nabla f(\balpha^\toptzero), \bbeta^\toptone - \balpha^\toptzero} + \frac{1}{2\eta_t} \normtwobig{\bbeta^\toptone - \balpha^\toptzero}^2.
$$
Combining the above two inequalities, for any $\bbeta$, we have
\begin{equation}\label{equation:fista_conv_ssf1}
\small
\begin{aligned}
&F(\bbeta^\toptone) = f(\bbeta^\toptone) + g(\bbeta^\toptone)\\
&\leq g(\bbeta) + f(\balpha^\toptzero) + \innerproduct{\nabla f(\balpha^\toptzero), \bbeta - \balpha^\toptzero} 
+ \frac{1}{\eta_t} \innerproduct{\bbeta^\toptone - \balpha^\toptzero, \bbeta - \bbeta^\toptone} 
+ \frac{1}{2\eta_t}\normtwobig{\bbeta^\toptone - \balpha^\toptzero}^2 \\
&\stackrel{\dag}{\leq} g(\bbeta) + f(\bbeta) + \frac{1}{\eta_t} \innerproduct{\bbeta^\toptone - \balpha^\toptzero, \bbeta - \bbeta^\toptone} + \frac{1}{2\eta_t} \normtwobig{\bbeta^\toptone - \balpha^\toptzero}^2\\
&= F(\bbeta) + \frac{1}{\eta_t} \innerproduct{\bbeta^\toptone - \balpha^\toptzero, \bbeta - \bbeta^\toptone} + \frac{1}{2\eta_t} \normtwobig{\bbeta^\toptone - \balpha^\toptzero}^2,
\end{aligned}
\end{equation}
where the inequality ($\dag$) follows from the convexity of $f$.
Invoking \eqref{equation:fista_conv_ssf1} with $\bbeta \triangleq \bbeta^\toptzero$ and $\bbeta \triangleq \bbeta^*$, for which we multiply by $1 - \gamma_t$ and $\gamma_t$, respectively, and add them up to get
\begin{equation}\label{equation:fista_conv_ssf2}
\small
\begin{aligned}
&F(\bbeta^\toptone) - F(\bbeta^*) - (1 - \gamma_t)\big(F(\bbeta^\toptzero) - F(\bbeta^*)\big)\\
&\leq \frac{1}{\eta_t} \innerproduct{\bbeta^\toptone - \balpha^\toptzero, (1 - \gamma_t)\bbeta^\toptzero + \gamma_t \bbeta^* - \bbeta^\toptone} + \frac{1}{2\eta_t} \normtwobig{\bbeta^\toptone - \balpha^\toptzero}^2\\
&\stackrel{\dag}{=} \frac{1}{2\eta_t} \left\{\normtwo{\balpha^\toptzero - (1 - \gamma_t)\bbeta^\toptzero - \gamma_t \bbeta^*}^2 - \normtwo{\bbeta^\toptone - (1 - \gamma_t)\bbeta^\toptzero - \gamma_t \bbeta^*}^2\right\}\\
&= \frac{\gamma_t^2}{2\eta_t} \left\{\normtwo{\bu^\toptzero - \bbeta^*}^2 - \normtwo{\bu^\toptone - \bbeta^*}^2\right\},
\end{aligned}
\end{equation}
where the equality ($\dag$) follows from the fact that $2\innerproduct{\ba-\bb,\bc-\ba} +\normtwo{\ba-\bb}^2 = \normtwo{\bb-\bc}^2-\normtwo{\ba-\bc}^2$ for any $\ba,\bb,\bc\in\real^p$, and the last equality follows from the update rules:
$$
\begin{aligned}
\bu^\toptone &= \bbeta^\toptzero + \frac{1}{\gamma_{t}}(\bbeta^\toptone - \bbeta^\toptzero)
\qquad\text{and}\qquad
\balpha^\toptzero = (1 - \gamma_t)\bbeta^\toptzero + \gamma_t \bu^\toptzero.
\end{aligned}
$$
Recalling that $\eta_t$ and $\gamma_t$ are chosen such that
$
\frac{1 - \gamma_t}{\gamma_t^2} \eta_t \leq \frac{1}{\gamma_{t-1}^2} \eta_{t-1},
$
we obtain an iterate inequality for consecutive iterations
$$
\small
\begin{aligned}
&\frac{\eta_t}{\gamma_t^2} \big(F(\bbeta^\toptone) - F(\bbeta^*)\big) + \frac{1}{2} \normtwo{\bu^\toptone - \bbeta^*}^2 
\leq (1-\gamma_t)\frac{\eta_t}{\gamma_t^2} \big(F(\bbeta^\toptzero) - F(\bbeta^*)\big) + \frac{1}{2} \normtwo{\bu^\toptzero - \bbeta^*}^2\\
&\leq \frac{\eta_{t-1}}{\gamma_{t-1}^2} \big(F(\bbeta^\toptzero) - F(\bbeta^*)\big) + \frac{1}{2} \normtwo{\bu^\toptzero - \bbeta^*}^2
\leq \frac{\eta_1}{\gamma_1^2} \big(F(\bbeta^{(2)}) - F(\bbeta^*)\big) + \frac{1}{2} \normtwo{\bu^{(2)} - \bbeta^*}^2\\
&\leq \frac{(1 - \gamma_1)\eta_1}{\gamma_1^2} \big(F(\bbeta^\topone) - F(\bbeta^*)\big) + \frac{1}{2} \normtwo{\bu^\topone - \bbeta^*}^2
=\frac{1}{2} \normtwo{\bbeta^\topone - \bbeta^*}^2,
\end{aligned}
$$
where the last equality follows from the fact that  $\gamma_1 = 1$ and $\bu^\topone = \bbeta^\topone$.
Using Lemma~\ref{lemma:seq_gamma_zeta} concludes the result.
\end{proof}

\paragrapharrow{Applied to LASSO.}
Note that the FISTA method (Algorithm~\ref{alg:fistav1} and Lemma~\ref{lemma:seq_gamma_zeta}) is a slight acceleration of the standard proximal gradient method. When applied to the LASSO problem
$$
\min_{\bbeta\in\real^p}   \frac{1}{2} \normtwo{\bX\bbeta-\by}^2 + \lambda \normone{\bbeta},
$$
the FISTA updates take the following form:
\begin{subequations}
\begin{align}
\balpha^\toptone &\leftarrow \mathcalT_{{\eta_t\lambda}} \left( \bbeta^\toptzero - \eta_t  \bX^\top (\bX \bbeta^\toptzero - \by) \right);\\
\zeta_{t+1} &\leftarrow \frac{1 + \sqrt{1 + 4 \zeta_t^2}}{2};\\
\bbeta^\toptone &\leftarrow \balpha^\toptone + \left( \frac{\zeta_t - 1}{\zeta_{t+1}} \right) (\balpha^\toptone - \balpha^\toptzero),
\end{align}
\end{subequations}
where $\eta_t = 1/\normtwo{\bX^\top\bX}$ and $\zeta_1=1$.

\index{Penalty function}
\index{Basis pursuit}
\subsection{Penalty Function Method}\label{section:pen_equa_lass}

We have now introduced several approaches for solving the Lagrangian LASSO problem.
Note that when the regularization parameter $\lambda$ is small, the objective is dominated by the smooth loss term $f=\frac{1}{2}\normtwo{\bX\bbeta-\by}^2$. 
The proximal gradient method typically uses a fixed stepsize $\eta=1/L_b=1/\normtwo{\bX^\top\bX}$,which depends only on the data matrix $\bX$ and is independent of $\lambda$.
More critically, when $\lambda$ is very small, the proximal (soft-thresholding) step applies only very weak shrinkage. As a result, sparsity is barely enforced, and the iterates behave almost like those of standard gradient descent applied to an ill-conditioned least squares problem---leading to extremely slow convergence.
To address this issue, we turn to the penalty function method.
\subsection*{Quadratic Penalty Function Method\index{Quadratic penalty function}}
Consider optimization problems with equality constraints:
\begin{equation}\label{equation:pr_in_constchap_equ}
\begin{aligned}
& \mathopmin{\bbeta} & f(\bbeta)&, \\
& \text{s.t.} & g_i(\bbeta)& = 0, \,\, i \in \mathcalE,\,\,\text{with $\abs{\mathcalE}=k$}.\\
\end{aligned}
\end{equation} 
Among various possible penalty formulations, the quadratic penalty is the simplest and most widely used. It is defined as follows:
\begin{definition}[Quadratic penalty function for equality constraints]
For the equality-constrained optimization problem  \eqref{equation:pr_in_constchap_equ}, 
the quadratic penalty function is defined as
\begin{equation}
f_{\sigma}(\bbeta) \triangleq f(\bbeta) + \frac{1}{2} \sigma \sum_{i \in \mathcalE} g_i^2(\bbeta)
\triangleq f(\bbeta) + \frac{1}{2} \sigma \normtwo{\bg(\bbeta)}^2,
\end{equation}
where $\bg(\bbeta) : \real^p \rightarrow \real^k $ is the vector-valued constraint function whose $ i $-th component is $ g_i(\bbeta) $.
The term $\frac{1}{2} \sigma \normtwo{\bg(\bbeta)}^2$ is called the \textit{penalty term}, and $\sigma>0$ is the \textit{penalty parameter} (or \textit{penalty factor}).
\end{definition}

This penalty function assigns a high cost to points that violate the constraints. During optimization, iterates often lie outside the feasible region. As $\sigma$ increases, the penalty term dominates the objective, pushing the minimizer of $f_{\sigma}(\bbeta)$ closer to the feasible set.

Conversely, within the feasible region (i.e., where $\bg(\bbeta)=\bzero$), the penalty term vanishes, and minimizing $f_{\sigma}(\bbeta)$ is equivalent to minimizing the original objective $f(\bbeta)$. Thus, as $\sigma\rightarrow\infty$, the minimizers of $f_{\sigma}$ converge to the solution of the original constrained problem~\eqref{equation:pr_in_constchap_equ}, under mild regularity conditions.

\begin{algorithm}[h] 
\caption{Quadratic Penalty Function Method for Model \eqref{equation:pr_in_constchap_equ}}
\label{alg:quad_pen_eq}
\begin{algorithmic}[1] 
\Require A function $f(\bbeta)$ and a set of constraints $\{g_i(\bbeta)\}$; 
\State {\bfseries input:}  Choose the penalty factor growth coefficient $\rho > 1$;
\State {\bfseries input:}  Initialize $\bbeta^{(1)}$ and $\sigma_1 > 0$;
\For{$t=1,2,\ldots$}
\State With $\bbeta^\toptzero$ as the initial point, solve $\bbeta^\toptone = \mathop{\argmin}_{\bbeta} f_{\sigma_t}(\bbeta)$;  \Comment{(QPF$_1$)}
\State Update penalty parameter:  $\sigma_{t+1} = \rho \sigma_t$;\Comment{(QPF$_2$)}
\State Stop if a stopping criterion is satisfied at iteration $t=T$; \Comment{(QPF$_3$)}
\EndFor
\State {\bfseries return:} Final  $\bbeta^{(T)}$;
\end{algorithmic} 
\end{algorithm}

The quadratic penalty function method, as outlined in Algorithm~\ref{alg:quad_pen_eq}, proceeds as follows: it generates a sequence of penalty parameters $\{\sigma_t\}$ that increase exponentially (e.g., $\sigma_{t+1} = \rho \sigma_t$ with $\rho>1$), and for each $\sigma_t$, it approximately minimizes the penalized objective $f_{\sigma_t}(\bbeta)$.
The subproblem in Step (QPF$_1$)---an unconstrained optimization problem---can be solved using any iterative method suitable for smooth objectives, such as gradient descent, quasi-Newton methods, or Newton-type algorithms.

Importantly, the interpretation of argmin in Step (QPF$_1$) may vary depending on the problem and solver:
\begin{itemize}
\item $\bbeta^\toptone$ is the global minimizer of $f_{\sigma_t}(\bbeta)$;
\item $\bbeta^\toptone$ is a local minimizer of  $f_{\sigma_t}(\bbeta)$;
\item $\bbeta^\toptone$ is not necessarily a strict minimizer but satisfies the first-order optimality condition approximately, i.e.,  $\nabla f_{\sigma_t}(\bbeta^\toptone) \approx \bzero$.
\end{itemize}
Based on this, three practical considerations are essential when implementing Algorithm~\ref{alg:quad_pen_eq}:
\begin{itemize}
\item \textit{Choice of penalty parameter schedule.}
The rate at which $\sigma_t$ increases significantly affects performance. If $\sigma_t$ grows too rapidly, the subproblems become increasingly ill-conditioned and difficult to solve. Conversely, if $\sigma_t$ increases too slowly, the algorithm requires many outer iterations, leading to slow overall convergence. A more adaptive strategy is to adjust the growth of $\sigma_t$ based on the observed difficulty of solving the current subproblem: if the previous subproblem converged quickly, one may safely increase $\sigma_{t+1}$ more aggressively; otherwise, a more conservative update is advisable.

\item \textit{Risk of unboundedness for small $\sigma_t$.}  
When $\sigma_t$ is too small, the penalty term may be insufficient to enforce feasibility, and $f_{\sigma_t}(\boldsymbol{\beta})$ can become unbounded below. In such cases, the subproblem solver may diverge. If divergence is detected during Step (QPF$_1$), the algorithm should halt the current inner solve, increase $\sigma_t$, and retry.

\item \textit{Required accuracy of subproblem solutions.}  
To ensure convergence of the overall method, the subproblems must be solved with sufficient accuracy. In particular, the error in satisfying the first-order optimality conditions should tend to zero as $t \to \infty$.
\end{itemize}

\subsection*{Convergence Analysis}
We now examine the convergence properties of the quadratic penalty method for equality-constrained problems. For simplicity, we assume that for each  $\sigma_t>0$, the penalized function $f_{\sigma_t}(\bbeta)$ admits at least one minimizer. Note that this assumption may fail for certain problems---for instance, if the penalty is too weak to prevent the objective from decreasing indefinitely along infeasible directions. In such cases, the quadratic penalty method is not applicable, and we restrict our analysis to problems satisfying this basic regularity condition.

\begin{theoremHigh}[Penalty function method with optimal subproblem]\label{theorem:pen_cong_glo}
Let $\bbeta^\toptone$ be a global minimizer of $f_{\sigma_t}(\bbeta)$,
\footnote{A refinement of a suboptimal solution to the subproblem is discussed in Problem~\ref{prob:pen_cong_glo_full}.}
and suppose the penalty parameters $\{\sigma_t\}$ satisfy $\sigma_t\rightarrow \infty$ monotonically in Algorithm~\ref{alg:quad_pen_eq}. 
Then every limit point $\bbeta^*$ of the sequence $\{\bbeta^\toptzero\}_{t>0}$ is a global minimizer of the original constrained  problem \eqref{equation:pr_in_constchap_equ}.
\end{theoremHigh}

\begin{proof}[of Theorem~\ref{theorem:pen_cong_glo}]
Let $\widehatbbeta$ be a global minimum of the original problem \eqref{equation:pr_in_constchap_equ}, i.e.,
$ f(\widehatbbeta) \leq f(\bbeta)$, for all $ \bbeta$ satisfies $g_i(\bbeta) = 0$, $i \in \mathcalE$.
Since $\bbeta^\toptone$ minimizes $f_{\sigma_t}(\bbeta)$, we have $ f_{\sigma_t}(\bbeta^\toptone) \leq f_{\sigma_t}(\widehatbbeta), $
i.e.,
\begin{align}
&f(\bbeta^\toptone) + \frac{\sigma_t}{2} \sum_{i \in \mathcalE} g_i(\bbeta^\toptone)^2 \leq f(\widehatbbeta) + \frac{\sigma_t}{2} \sum_{i \in \mathcalE} g_i(\widehatbbeta)^2 \label{equation:pen_cong_glo1}\\
&\quad\implies\quad 
\sum_{i \in \mathcalE} g_i(\bbeta^\toptone)^2 \leq \frac{2}{\sigma_t} \big(f(\widehatbbeta) - f(\bbeta^\toptone)\big). \label{equation:pen_cong_glo2}
\end{align}
Let $\bbeta^*$ be any limit point of the sequence $\{\bbeta^\toptzero\}$.  
Taking the limit as $t \rightarrow\infty$ in \eqref{equation:pen_cong_glo2}, and using the continuity of  $f(\bbeta)$ and $g_i(\bbeta)$ 
together with $\sigma_t \rightarrow +\infty$, we obtain
$ \sum_{i \in \mathcalE} g_i(\bbeta^*)^2 = 0 $,
so $\bbeta^*$ is feasible for the original problem. 
Moreover, from~\eqref{equation:pen_cong_glo1}, we have
$ f(\bbeta^\toptone) \leq f(\widehatbbeta)$ for all $t$.
Passing to the limit gives $f(\bbeta^*) \leq f(\widehatbbeta)$. 
By the optimality of $\widehatbbeta$, it follows that $f(\bbeta^*) = f(\widehatbbeta)$, and thus $\bbeta^*$ is also a global minimizer.
\end{proof}

\begin{algorithm}[h]
\caption{Lagrangian LASSO Problem  via Penalty Function Method}
\label{alg:lass_penal}
\begin{algorithmic}[1]
\State {\bfseries input:}  Given initial value $\bbeta^\topone$, regularized parameter $\lambda$, initial parameter $\zeta_1\gg \lambda$, reducing factor $\gamma \in (0, 1)$.
\For{$t=1,2,\ldots$}
\State \algoalign{Using $\bbeta^\toptzero$ from last iteration as the initial guess, solve the problem $\bbeta^\toptone = \argmin_{\bbeta} \left\{ \frac{1}{2} \normtwo{\bX\bbeta-\by}^2 + \zeta_t \normone{\bbeta} \right\}$ by subgradient descent method or other unconstrained optimization method;}
\If{$\zeta_t = \lambda$}
\State Stop iteration, output $\bbeta^\toptone$;
\Else
\State Update/reducing the penalty factor $\zeta_{t+1} = \max\{\lambda, \gamma \zeta_t\}$;
\EndIf
\State Stop if a stopping criterion is satisfied at iteration $t=T$;
\EndFor
\State \Return Final $\bbeta^{(T)} $;
\end{algorithmic}
\end{algorithm}

\subsection*{Penalty Function Method Applied to Lagrangian LASSO}
Recall the Lagrangian LASSO problem \eqref{opt:ll}:
$$
\min_{\bbeta\in\real^p}   \frac{1}{2} \normtwo{\bX\bbeta-\by}^2 + \lambda \normone{\bbeta},
$$
where $\lambda > 0$ is a regularization parameter. Solving the LASSO problem can be viewed as an approximation to solving the following basis pursuit (BP) problem:
$$
\begin{aligned}
& \min_{\bbeta\in\real^p}   \normone{\bbeta}\quad \text{s.t.} \quad \bX\bbeta = \by,
\end{aligned}
$$
The BP problem is a non-smooth optimization problem with linear equality constraints. By applying a quadratic penalty function to enforce the constraint $\bX\bbeta = \by$, we obtain the penalized formulation:
$$
\min_{\bbeta\in\real^p}   \normone{\bbeta} + \frac{\sigma}{2} \normtwo{\bX\bbeta-\by}^2.
$$
Setting  $ \lambda \triangleq \frac{1}{\sigma}$  reveals a direct correspondence: the above penalized problem is precisely the LASSO problem with regularization parameter $\lambda$. This observation leads to two important insights:
\begin{itemize}
\item The solutions of the LASSO and BP problems are generally not identical. However, as $\lambda\rightarrow 0^+$ (equivalently, $\sigma\rightarrow \infty$), the solution of the LASSO problem converges to a solution of the BP problem---provided a feasible solution exists.

\item When $\lambda$ is very small (i.e., $\sigma$ is large), the quadratic term becomes dominant and highly sensitive to constraint violations, but the overall objective remains non-smooth due to the $\ell_1$-norm. Moreover, if one attempts to solve the BP problem directly via the penalty method with a fixed large $\sigma$, the resulting subproblem can become ill-conditioned, leading to slow convergence.

\end{itemize}
In the classical penalty function framework, the penalty parameter $\sigma$ is gradually increased toward infinity to enforce feasibility more strictly over successive iterations. Translating this into the LASSO context, this strategy corresponds to starting with a relatively large $\lambda$ (i.e., weak enforcement of the data-fitting term) and progressively decreasing $\lambda$ toward the desired value. 
This homotopy-like approach improves numerical stability and convergence behavior. The specific procedure is detailed in Algorithm~\ref{alg:lass_penal}.

\subsection{Alternating Direction Method of Multipliers (ADMM)\index{ADMM}}\label{section:laglass_admm}
\textit{Alternating direction methods of multipliers (ADMM)} is widely used to solve problems such as LASSO, sparse logistic regression, and basis pursuit. These applications benefit from ADMM's ability to efficiently handle sparsity constraints. Moreover, ADMM can also be applied to optimization problems arising in the training of support vector machines (SVMs; see Problem~\ref{prob:hing_svm}), particularly when dealing with large-scale datasets \citep{boyd2011distributed}.
We first provide a brief introduction to ADMM and then discuss its application to the LASSO problem.

\paragrapharrow{Fixed Augmented Lagrangian Method.\index{Augmented Lagrangian method}}
Consider the following convex optimization problem:
\begin{subequations}\label{equation:admm_subequs}
\begin{equation}\label{equation:admm_prob}
\mathopmin{\bbeta, \balpha} F(\bbeta)\triangleq f(\bbeta)+g(\balpha) 
\quad \text{s.t.}\quad \bX\bbeta+\bZ\balpha=\bff,
\end{equation}
where $\bX\in\real^{n\times p}$, $\bZ\in\real^{n\times q}$, $\bff\in\real^n$, and $f : \real^p \to (-\infty, \infty]$ and $g : \real^q \to (-\infty, \infty]$ are proper closed convex functions.
The Lagrangian function is 
\begin{equation}
L (\bbeta, \balpha, \bnu) = f(\bbeta)+g(\balpha) + \innerproduct{\bnu, \bX\bbeta+\bZ\balpha-\bff}.
\end{equation}
The dual objective function is obtained by minimizing the Lagrangian with respect to the primal variables $\bbeta, \balpha$ (see Section~\ref{section:sub_conjug}):
$$
Q(\bnu) = \min_{\bbeta, \balpha} L(\bbeta, \balpha, \bnu) = -f^*(-\bX^\top\bnu) - g^*(-\bZ^\top\bnu) - \innerproduct{\bnu, \bff}.
$$
Thus, the dual problem becomes
\begin{equation}\label{equation:admm_dualopt}
\min_{\bnu \in \real^n} f^*(-\bX^\top\bnu) + g^*(-\bZ^\top\bnu) + \innerproduct{\bnu, \bff}.
\end{equation}
\end{subequations}
If $f$ and $g$ are proper closed and convex, the dual problem can be solved using the \textit{proximal point method} (Problem~\ref{prob:proximal_point}, which is a consequence of the primal gradient method), due to the properness, closedness, and convexity of the conjugate functions (Lemma~\ref{lemma:closedconv_conj}, Exercise~\ref{exercise:proper_conj}).
Given a parameter $\rho>0$, the update at the $t$-th iteration is 
$$
\bnu^\toptone \leftarrow 
\argmin_{\bnu \in \real^n} 
\left\{G(\bnu)\triangleq f^*(-\bX^\top\bnu) + g^*(-\bZ^\top\bnu) + \innerproduct{\bnu, \bff} + \frac{1}{2\rho } \normtwobig{\bnu - \bnu^\toptzero}^2\right\}.
$$
By the first-order necessary and sufficient conditions of convex functions (Theorem~\ref{theorem:fetmat_opt}), it follows that 
$$
\bzero \in \partial G(\bnu^\toptone) =
-\bX\partial f^*(-\bX^\top\bnu^\toptone)  -\bZ\partial g^*(-\bZ^\top\bnu^\toptone) +  \bff + \frac{{\bnu^\toptone - \bnu^\toptzero}}{\rho }.
$$
From the definition of conjugate functions (Definition~\ref{definition:conjug_func}), we have
$$
\begin{aligned}
\min_{\bbeta\in\real^p} \innerproduct{\bX^\top\bnu^\toptone, \bbeta} + f(\bbeta)=  f^*(-\bX^\top\bnu^\toptone);\\
\min_{\balpha\in\real^q} \innerproduct{\bZ^\top\bnu^\toptone, \balpha} + g(\balpha)= g^*(-\bZ^\top\bnu^\toptone).
\end{aligned}
$$
By the conjugate subgradient theorem (Theorem~\ref{theorem:conju_subgra}), this is equivalent to updating $\bbeta$, $\balpha$, and $\bnu$ by
$$
\begin{aligned}
\bbeta^\toptone&\in \argmin_{\bbeta\in\real^p} \innerproduct{\bX^\top\bnu^\toptone, \bbeta} + f(\bbeta)=\partial f^*(-\bX^\top\bnu^\toptone);\\
\balpha^\toptone&\in \argmin_{\balpha\in\real^q} \innerproduct{\bZ^\top\bnu^\toptone, \balpha} + g(\balpha)=\partial g^*(-\bZ^\top\bnu^\toptone);\\
\bnu^\toptone &\leftarrow \bnu^\toptzero + \rho\left( \bX\bbeta^\toptone  + \bZ \balpha^\toptone  - \bff \right).
\end{aligned}
$$
Substituting the expression for $\bnu^\toptone$ into the first two equalities and applying  the first-order necessary and sufficient conditions of convex functions (Theorem~\ref{theorem:fetmat_opt}) again (since $f$ and $g$ are proper closed and convex), it follows that 
\begin{subequations}\label{equation:aug_ad_raw}
\begin{align}
\bzero&\in  \bX^\top\left(\bnu^\toptzero + \rho\big( \bX\bbeta^\toptone  + \bZ \balpha^\toptone  - \bff \big)\right) + \partial f(\bbeta^\toptone);\label{equation:aug_ad_raw1}\\
\bzero&\in  \bZ^\top\left(\bnu^\toptzero + \rho\big( \bX\bbeta^\toptone  + \bZ \balpha^\toptone  - \bff \big)\right) + \partial g(\balpha^\toptone);\label{equation:aug_ad_raw2}\\
\bnu^\toptone &\leftarrow \bnu^\toptzero + \rho\left( \bX\bbeta^\toptone  + \bZ \balpha^\toptone  - \bff \right).\label{equation:aug_ad_raw3}
\end{align}
\end{subequations}
Given $\bnu^\toptzero$, the quantities \eqref{equation:aug_ad_raw2} and \eqref{equation:aug_ad_raw3} correspond to the optimality conditions of the function
\begin{equation}\label{equation:admm_auglag}
L_\rho (\bbeta, \balpha, \bnu^\toptzero) = f(\bbeta)+g(\balpha) + \innerproduct{\bnu^\toptzero, \bX\bbeta+\bZ\balpha-\bff} + \frac{\rho}{2}\normtwo{\bX\bbeta+\bZ\balpha-\bff}^2,
\end{equation}
which is known as the \textit{augmented Lagrangian function} \citep{lu2025practical}.
When $\rho=0$, the augmented Lagrangian function reduces to the standard Lagrangian function; when $\rho>0$, it  acts as  a penalized version of the Lagrangian function.
The \textit{augmented Lagrangian method}  proceeds iteratively  (at the $t$-th iteration) as follows:
$$
\text{augmented Lagrangian:}
\quad 
\left\{
\begin{aligned}
(\bbeta^\toptone, \balpha^\toptone) 
&\in \mathop{\argmin}_{\bbeta, \balpha} L_\rho (\bbeta, \balpha, \bnu^\toptzero);\\
\bnu^\toptone
&\leftarrow\bnu^\toptzero + \rho(\bX\bbeta^\toptone+\bZ\balpha^\toptone -\bff ).
\end{aligned}
\right.
$$
However, unlike in standard augmented Lagrangian methods \citep{lu2025practical}, the parameter $\rho$ can remain fixed across iterations due to the convergence guarantees of the proximal point method (Problem~\ref{prob:proximal_point}).

\paragrapharrow{ADMM.}
One significant challenge in solving the problem is the coupling term between the variables $\bbeta$ and  $\balpha$, which appears as $\rho(\bbeta^\top\bX^\top\bZ\balpha)$.
ADMM addresses this issue by replacing the joint minimization over $(\bbeta,\balpha)$ with a single iteration of alternating minimization. 
Specifically, at the $t$-th iteration, ADMM updates $\bbeta, \balpha$, and the dual variable $\bnu$ from \eqref{equation:admm_auglag} sequentially as follows:
\begin{subequations}
\begin{tcolorbox}[colback=white,colframe=black]
\begin{minipage}{1\textwidth}
\begin{equation}\label{equation:upda_admm1}
\text{ADMM:}\quad
\left\{
\small
\begin{aligned}
\bbeta^\toptone&\in\mathop{\argmin}_{\bbeta} \left\{ f(\bbeta)+ \frac{\rho}{2}\normtwo{\bX\bbeta+\bZ\balpha^\toptzero -\bff +\frac{1}{\rho}\bnu^\toptzero}^2 \right\};\\
\balpha^\toptone&\in\mathop{\argmin}_{\balpha} \left\{ g(\balpha)+ \frac{\rho}{2}\normtwo{\bX\bbeta^\toptone+\bZ\balpha -\bff +\frac{1}{\rho}\bnu^\toptzero}^2 \right\};\\
\bnu^\toptone&\leftarrow\bnu^\toptzero + \rho(\bX\bbeta^\toptone+\bZ\balpha^\toptone -\bff ).
\end{aligned}
\right.
\end{equation}
\end{minipage}
\end{tcolorbox}
Let $\widetildebnu\triangleq\frac{1}{\rho}\bnu$, this scheme can be equivalently rewritten in terms of the scaled dual variable $\widetildebnu$:
\begin{tcolorbox}[colback=white,colframe=black]
\begin{minipage}{1\textwidth}
\begin{equation}\label{equation:admm_gen_up}
\text{ADMM:}\quad
\left\{
\small
\begin{aligned}
\bbeta^\toptone&\in\mathop{\argmin}_{\bbeta} \left\{ f(\bbeta)+ \frac{\rho}{2}\normtwo{\bX\bbeta+\bZ\balpha^\toptzero -\bff +\widetildebnu^\toptzero}^2 \right\};\\
\balpha^\toptone&\in\mathop{\argmin}_{\balpha} \left\{ g(\balpha)+ \frac{\rho}{2}\normtwo{\bX\bbeta^\toptone+\bZ\balpha -\bff +\widetildebnu^\toptzero}^2 \right\};\\
\widetildebnu^\toptone&\leftarrow\widetildebnu^\toptzero + (\bX\bbeta^\toptone+\bZ\balpha^\toptone -\bff ).
\end{aligned}
\right.
\end{equation}
\end{minipage}
\end{tcolorbox}
\end{subequations}

Note that the scaled-form ADMM update in \eqref{equation:admm_gen_up} can be derived from the following augmented Lagrangian expressed in terms of the scaled dual variable $\widetildebnu$:
\begin{equation}
L_\rho (\bbeta, \balpha, \widetildebnu) = f(\bbeta)+g(\balpha)  + \frac{\rho}{2}\normtwo{\bX\bbeta+\bZ\balpha-\bff +\widetildebnu}^2 - \frac{\rho}{2} \normtwo{\widetildebnu}^2.
\end{equation}

This iterative procedure enables the algorithm to converge to a solution that minimizes the objective function while satisfying the constraints. The main strength of ADMM lies in its ability to efficiently handle complex objectives and constraints---particularly in high-dimensional settings that are common in modern statistical learning and large-scale optimization problems.

\paragrapharrow{Stopping criteria.}
ADMM iterations continue until both the primal and dual residuals fall below prescribed tolerance levels. These residuals are defined as:
\begin{align*}
\br^\toptone &= \bbeta^\toptone - \balpha^\toptone \quad \text{(primal residual)}, \\
\bs^\toptone &= \rho(\balpha^\toptone - \balpha^\toptzero) \quad \text{(dual residual)}.
\end{align*}
The algorithm terminates when
$$
\normtwo{\br^\toptone} \leq \varepsilon_{\text{pri}}
\quad\text{and} \quad 
\normtwo{\bs^\toptone} \leq \varepsilon_{\text{dual}}.
$$
where  $\varepsilon_{\text{pri}}$ and  $\varepsilon_{\text{dual}}$ are user-specified tolerances.

\paragrapharrow{Convergence of ADMM.}
The following theorem
demonstrates an $\mathcalO(1/T)$ rate of convergence for the sequence generated by ADMM.
To state the result, we make the following assumptions.
For any $\bu \in \real^p$ and $\bv \in \real^q$, the optimal solution sets of the problems
$$
\min_{\bbeta \in \real^p} \left\{ f(\bbeta) + \frac{\rho}{2} \normtwo{\bD \bbeta}^2 + \innerproduct{\bu, \bbeta} \right\}
$$
and
$$
\min_{\balpha \in \real^q} \left\{ g(\balpha) + \frac{\rho}{2} \normtwo{\bZ \balpha}^2 + \innerproduct{\bv, \balpha} \right\}
$$
are nonempty.
Moreover, we assume there exist $\widetilde{\bbeta} \in \domain(f)$ and $\widetilde{\balpha} \in \domain(g)$ such that $\bD \widetilde{\bbeta} + \bZ \widetilde{\balpha} = \bff$.
These conditions ensure that the ADMM iterates are well-defined.
\begin{theoremHigh}[Convergence of ADMM, $\mathcalO(1/T)$]\label{theorem:ratconv_admm}
Suppose the above assumptions hold, 
and let $f : \real^p \to (-\infty, \infty]$ and $g : \real^q \to (-\infty, \infty]$ be proper closed convex functions.
Let $\big\{(\bbeta^\toptzero, \balpha^\toptzero)\big\}_{t > 0}$ be the sequence generated by ADMM for solving problem \eqref{equation:admm_prob}. Let $(\bbeta^*, \balpha^*)$ be an optimal solution of problem \eqref{equation:admm_prob} and $\bnu^*$  an optimal solution of the dual problem \eqref{equation:admm_dualopt}. 
Assume  $\zeta > 0$ is any constant satisfying $\zeta \geq  \normtwo{\bnu^*}$.
Then for all $T >0$,
\begin{align}
F(\overline{\bbeta}^{(T)}, \overline{\balpha}^{(T)}) - F(\bbeta^*, \balpha^*) &\leq \frac{ \normC{\balpha^* - \balpha^\topone}^2 + \frac{1}{\rho} \big(\zeta + \normtwobig{\bnu^\topone}\big)^2}{2T},
\end{align}
where $\bC \triangleq \rho \bZ^\top \bZ$, $\normC{\ba} = (\ba^\top\bC\ba)^{1/2}$, 
and the averaged iterates are defined as
$$
\overline{\bbeta}^{(T)} = \frac{1}{T} \sum_{t=1}^T \bbeta^\toptone
\qquad \text{and}\qquad 
\overline{\balpha}^{(T)} = \frac{1}{T} \sum_{t=1}^T \balpha^\toptone.
$$
\end{theoremHigh}
\begin{proof}[of Theorem~\ref{theorem:ratconv_admm}]
From the first-order optimality conditions for convex functions (Theorem~\ref{theorem:fetmat_opt}) and the ADMM update  rules \eqref{equation:upda_admm1}, it follows that $\bbeta^\toptone$ and $\balpha^\toptone$ satisfy
\begin{align}
- \rho \bX^\top \left( \bX\bbeta^\toptone + \bZ \balpha^\toptzero - \bff + \frac{1}{\rho} \bnu^\toptzero \right)
&\in \partial f(\bbeta^\toptone), \\
- \rho \bZ^\top \left( \bX\bbeta^\toptone + \bZ \balpha^\toptone - \bff + \frac{1}{\rho} \bnu^\toptzero \right) 
&\in \partial g(\balpha^\toptone).
\end{align}
Define the auxiliary variables
$$
\begin{aligned}
\widetildebbeta^\toptzero &\triangleq \bbeta^\toptone, 
\qquad
\widetildebalpha^\toptzero &\triangleq \balpha^\toptone, 
\qquad
\widetildebnu^\toptzero &\triangleq \bnu^\toptzero + \rho (\bX \bbeta^\toptone + \bZ \balpha^\toptzero - \bff).
\end{aligned}
$$
Using the above optimality conditions and the subgradient inequality (Definition~\ref{definition:subgrad}),  we obtain that for any $\bbeta \in \domain(f)$ and $\balpha \in \domain(g)$,
$$
\begin{aligned}
f(\bbeta) - f(\widetildebbeta^\toptzero) + \innerproduct{\rho \bX^\top \left( \bX \widetildebbeta^\toptzero + \bZ \balpha^\toptzero - \bff + \frac{1}{\rho} \bnu^\toptzero \right) , \bbeta - \widetildebbeta^\toptzero} &\geq 0, \\
g(\balpha) - g(\widetildebalpha^\toptzero) + \innerproduct{\rho \bZ^\top \left( \bX \widetildebbeta^\toptzero + \bZ \widetildebalpha^\toptzero - \bff + \frac{1}{\rho} \bnu^\toptzero \right), \balpha - \widetildebalpha^\toptzero} &\geq 0.
\end{aligned}
$$
By the definition of $\widetildebnu^\toptzero$,  these inequalities can be rewritten as
$$
\begin{aligned}
f(\bbeta) - f(\widetildebbeta^\toptzero) + \innerproduct{\bX^\top \widetildebnu^\toptzero , \bbeta - \widetildebbeta^\toptzero} 
&\geq 0, \\
g(\balpha) - g(\widetildebalpha^\toptzero) + \innerproduct{\bZ^\top \widetildebnu^\toptzero + \bC (\widetildebalpha^\toptzero - \balpha^\toptzero), \balpha - \widetildebalpha^\toptzero} 
&\geq 0,
\end{aligned}
$$
where $\bC = \rho \bZ^\top \bZ$.
Adding the two inequalities and using the dual update for $\bnu^\toptone$:
$
\bnu^\toptone - \bnu^\toptzero = \rho (\bX \widetildebbeta^\toptzero + \bZ \widetildebalpha^\toptzero - \bff),
$
we obtain that for any $\bbeta \in \domain(f)$, $\balpha \in \domain(g)$, and $\bnu \in \real^n$,
\begin{equation}\label{equation:admm_prof00}
\footnotesize
\begin{aligned}
F(\bbeta, \balpha) - F(\widetildebbeta^\toptzero, \widetildebalpha^\toptzero) + 
\innerproduct{\begin{bmatrix} 
\bbeta - \widetildebbeta^\toptzero \\ 
\balpha - \widetildebalpha^\toptzero \\ 
\bnu - \widetildebnu^\toptzero 
\end{bmatrix}, 
\begin{bmatrix} 
\bX^\top \widetildebnu^\toptzero \\ 
\bZ^\top \widetildebnu^\toptzero \\
- \bX \widetildebbeta^\toptzero - \bZ \widetildebalpha^\toptzero + \bff 
\end{bmatrix} 
- 
\begin{bmatrix} \bzero (\bbeta^\toptzero - \widetildebbeta^\toptzero) \\ 
\bC (\balpha^\toptzero - \widetildebalpha^\toptzero) \\ 
\frac{1}{\rho} \big(\bnu^\toptzero - \bnu^\toptone \big)
\end{bmatrix} }
\geq 0.
\end{aligned}
\end{equation}
 For any vectors $\ba,\bb,\bc,\bd$, and a positive semidefinite matrix $\bH$, we have:
$$
(\ba - \bb)^\top \bH (\bc - \bd) 
= 
\frac{1}{2} \left( \normH{\ba - \bd}^2 - \normH{\ba - \bc }^2 + \normH{\bb - \bc}^2 - \normH{\bb - \bd}^2 \right).
$$
Applying this with $\bH=\bC$ yields
\begin{equation}\label{equation:admm_prof2}
\begin{aligned}
(\balpha - \widetildebalpha^\toptzero)^\top \bC (\balpha^\toptzero - \widetildebalpha^\toptzero) 
&= \frac{1}{2} \normC{\balpha - \widetildebalpha^\toptzero }^2 - \frac{1}{2} \normC{\balpha - \balpha^\toptzero}^2 + \frac{1}{2} \normC{\widetildebalpha^\toptzero - \balpha^\toptzero}^2.
\end{aligned}
\end{equation}
Similarly, using the definition of $\widetildebnu^\toptzero$ and $\bnu^\toptone$, we derive
$$
\small
\begin{aligned}
&2 (\bnu - \widetildebnu^\toptzero)^\top (\bnu^\toptzero - \bnu^\toptone) 
= \normtwobig{\bnu - \bnu^\toptone}^2 - \normtwobig{\bnu - \bnu^\toptzero}^2 +\normtwobig{\widetildebnu^\toptzero - \bnu^\toptzero}^2 - \normtwobig{ \widetildebnu^\toptzero - \bnu^\toptone}^2\\
&= \normtwobig{\bnu - \bnu^\toptone}^2 - \normtwobig{\bnu - \bnu^\toptzero}^2 
+ \rho^2 \normtwo{\bX \widetildebbeta^\toptzero + \bZ \balpha^\toptzero - \bff}^2 
- \rho^2 \normtwo{\bZ (\balpha^\toptzero - \widetildebalpha^\toptzero)}^2.
\end{aligned}
$$
Therefore,
\begin{equation}\label{equation:admm_prof3}
\small
\begin{aligned}
\frac{1}{\rho} (\bnu - \widetildebnu^\toptzero)^\top (\bnu^\toptzero - \bnu^\toptone) \geq \frac{1}{2\rho} \left( \normtwobig{ \bnu - \bnu^\toptone}^2 - \normtwobig{ \bnu - \bnu^\toptzero}^2 \right) - \frac{\rho}{2} \normtwobig{\bZ (\balpha^\toptzero - \widetildebalpha^\toptzero)}^2.
\end{aligned}
\end{equation}
Now define the block-diagonal matrix and the stacked vectors:
$$
\bA \triangleq
\small
\begin{bmatrix}
\bzero & \bzero & \bzero \\
\bzero & \bC & \bzero \\
\bzero & \bzero & \frac{1}{\rho} \bI
\end{bmatrix}, 
\normalsize
\qquad
\bu \triangleq 
\small
\begin{bmatrix}
\bbeta \\
\balpha \\
\bnu
\end{bmatrix}, 
\normalsize
\qquad
\bu^\toptzero \triangleq 
\small
\begin{bmatrix}
\bbeta^\toptzero \\
\balpha^\toptzero \\
\bnu^\toptzero
\end{bmatrix}, 
\normalsize
\quad\text{and}\quad 
\widetildebu^\toptzero \triangleq 
\small
\begin{bmatrix}
\widetildebbeta^\toptzero \\
\widetildebalpha^\toptzero \\
\widetildebnu^\toptzero
\end{bmatrix}.
\normalsize
$$
Combining \eqref{equation:admm_prof2} and \eqref{equation:admm_prof3}, we obtain
$$
\left\langle
\small
\begin{bmatrix}
\bbeta - \widetildebbeta^\toptzero \\
\balpha - \widetildebalpha^\toptzero \\
\bnu - \widetildebnu^\toptzero
\end{bmatrix}, 
\begin{bmatrix}
\bzero (\bbeta^\toptzero - \widetildebbeta^\toptzero) \\
\bC (\balpha^\toptzero - \widetildebalpha^\toptzero) \\
\frac{1}{\rho} (\bnu^\toptzero - \bnu^\toptone)
\end{bmatrix} 
\normalsize
\right\rangle \geq \frac{1}{2} \normA{\bu - \bu^\toptone}^2 - \frac{1}{2} \normA{\bu - \bu^\toptzero}^2.
$$
Combining the preceding inequality with \eqref{equation:admm_prof00}, we obtain that for any $\bbeta \in \domain(f)$, $\balpha \in \domain(g)$, and $\bnu \in \real^n$,
\begin{equation}\label{equation:admm_prof4}
F(\bbeta, \balpha) - F(\widetildebbeta^\toptzero, \widetildebalpha^\toptzero) +\innerproduct{\bu - \widetildebu^\toptzero, \bG \widetildebu^\toptzero + \widetilde{\bff}} \geq \frac{1}{2} \normAbig{\bu - \bu^\toptone}^2 - \frac{1}{2} \normAbig{\bu - \bu^\toptzero}^2, 
\end{equation}
where
$$
\bG \triangleq 
\small
\begin{bmatrix}
\bzero & \bzero & \bX^\top \\
\bzero & \bzero & \bZ^\top \\
-\bX & -\bZ & \bzero
\end{bmatrix}
\normalsize
\qquad \text{and}\qquad 
\widetilde{\bff} \small
\triangleq 
\begin{bmatrix}
\bzero \\
\bzero \\
\bff
\end{bmatrix}.
$$
Note that
$$
\innerproduct{\bu - \widetildebu^\toptzero, \bG \widetildebu^\toptzero + \widetilde{\bff}}
= \innerproduct{\bu - \widetildebu^\toptzero, \bG (\widetildebu^\toptzero - \bu) + \bG \bu + \widetilde{\bff}}
= \innerproduct{\bu - \widetildebu^\toptzero, \bG \bu + \widetilde{\bff}},
$$
where the last equality follows from the fact that $\bG$ is skew-symmetric (i.e., $\bG^\top = -\bG$). 
Therefore, \eqref{equation:admm_prof4} can be rewritten as
$$
F(\bbeta, \balpha) - F(\widetildebbeta^\toptzero, \widetildebalpha^\toptzero) + 
\innerproduct{\bu - \widetildebu^\toptzero, \bG \bu + \widetilde{\bff} }
\geq \frac{1}{2} \normAbig{\bu - \bu^\toptone}^2 - \frac{1}{2} \normAbig{\bu - \bu^\toptzero}^2.
$$
Performing sum of the above inequality over $t = 1, 2, \ldots, T$ yields the inequality
$$
T\cdot  F(\bbeta, \balpha) - \sum_{t=1}^T F(\widetildebbeta^\toptzero, \widetildebalpha^\toptzero) + 
\innerproduct{T \cdot \bu - \sum_{t=1}^T \widetildebu^\toptzero, \bG \bu + \widetilde{\bff}}
\geq -\frac{1}{2} \normAbig{\bu - \bu^\topone}^2.
$$
Defining
$$
\overline{\bu}^{(T)} \triangleq \frac{1}{T} \sum_{t=1}^T \widetildebu^\toptzero, 
\qquad 
\overline{\bbeta}^{(T)} \triangleq \frac{1}{T} \sum_{t=1}^T \bbeta^\toptone,
\qquad 
\overline{\balpha}^{(T)} \triangleq \frac{1}{T} \sum_{t=1}^T \balpha^\toptone,
$$
and using the convexity of $F$, we obtain that
$$
\begin{aligned}
&F(\bbeta, \balpha) - F(\overline{\bbeta}^{(T)}, \overline{\balpha}^{(T)}) + \innerproduct{\bu - \overline{\bu}^{(T)}, \bG \bu + \widetilde{\bff}} + \frac{1}{2T} \normAbig{\bu - \bu^\topone}^2 \geq 0 \\
&\implies 
F(\bbeta, \balpha) - F(\overline{\bbeta}^{(T)}, \overline{\balpha}^{(T)}) + \innerproduct{\bu - \overline{\bu}^{(T)}, \bG \overline{\bu}^{(T)} + \widetilde{\bff}} + \frac{1}{2T} \normAbig{\bu - \bu^\topone}^2 \geq 0,
\end{aligned}
$$
where the implication follows again from  the skew-symmetry of $\bG$.
In other words, for any $\bbeta \in \domain(f)$ and $\balpha \in \domain(g)$,
\begin{equation}\label{equation:admm_prof5}
F(\overline{\bbeta}^{(T)}, \overline{\balpha}^{(T)}) - F(\bbeta, \balpha) + \innerproduct{\overline{\bu}^{(T)} - \bu, \bG \overline{\bu}^{(T)} + \widetilde{\bff}} \leq \frac{1}{2T} \normAbig{\bu - \bu^\topone}^2. 
\end{equation}

Let $(\bbeta^*, \balpha^*)$ be an optimal solution of problem \eqref{equation:admm_prob}. Then $\bX \bbeta^* + \bZ \balpha^* = \bff$. 
Let $\overline{\bnu}^{(T)} \triangleq \frac{1}{T} \sum_{t=1}^T \bnu^\toptzero$. Plugging $\bbeta = \bbeta^*$, $\balpha = \balpha^*$, and the expressions for $\overline{\bu}^{(T)}$, $\bu$, $\bu^\topone$, $\bG$, $\bZ$, $\widetilde{\bff}$ into \eqref{equation:admm_prof5}, we obtain
$$
\begin{aligned}
& F(\overline{\bbeta}^{(T)}, \overline{\balpha}^{(T)}) - F(\bbeta^*, \balpha^*) +\innerproduct{\overline{\bbeta}^{(T)} - \bbeta^*, \bX^\top \overline{\bnu}^{(T)}} + \innerproduct{\overline{\balpha}^{(T)} - \balpha^*, \bZ^\top \overline{\bnu}^{(T)}} \\
& \quad + \innerproduct{\overline{\bnu}^{(T)} - \bnu, -\bX \overline{\bbeta}^{(T)} - \bZ \overline{\balpha}^{(T)} + \bff}
\quad \leq\quad \frac{1}{2T} \left\{ \normCbig{\balpha^* - \balpha^\topone}^2 + \frac{1}{\rho} \normtwobig{ \bnu - \bnu^\topone}^2 \right\}.
\end{aligned}
$$
Cancelling terms and using the fact that $\bX \bbeta^* + \bZ \balpha^* = \bff$, the above inequality reduces to
$$
\begin{aligned}
&F(\overline{\bbeta}^{(T)}, \overline{\balpha}^{(T)}) - F(\bbeta^*, \balpha^*) + \innerproduct{\bnu, \bX \overline{\bbeta}^{(T)} + \bZ \overline{\balpha}^{(T)} - \bff} 
\leq \frac{ \normC{\balpha^* - \balpha^\topone}^2 + \frac{1}{\rho} \normtwobig{\bnu - \bnu^\topone }^2}{2T}.
\end{aligned}
$$
Taking $\bnu \in \sB[\bzero, \zeta]$ as the maximum of both sides obtains
$$
\begin{aligned}
&F(\overline{\bbeta}^{(T)}, \overline{\balpha}^{(T)}) - F(\bbeta^*, \balpha^*) + 
\zeta\normtwo{\bX \overline{\bbeta}^{(T)} + \bZ \overline{\balpha}^{(T)} - \bff}
\leq \frac{ \normC{\balpha^* - \balpha^\topone}^2 + \frac{1}{\rho} (\zeta + \normtwobig{\bnu^\topone})^2}{2T}.
\end{aligned}
$$
This completes the proof.
\end{proof}

In fact, by leveraging the strong duality between the primal problem~\eqref{equation:admm_prob} and its dual~\eqref{equation:admm_dualopt}, and under the assumption that $\zeta \geq 2 \normtwo{\bnu^*}$,  one can also establish the following bound on the feasibility residual:
$$
\normtwo{\bD \overline{\bbeta}^{(T)} + \bZ \overline{\balpha}^{(T)} - \bff} 
\leq \frac{ \normC{\balpha^* - \balpha^\topone}^2 + \frac{1}{\rho} \big(\zeta + \normtwobig{\bnu^\topone }\big)^2}{\zeta T},
$$
See \citet{beck2017first} for more details.

\paragrapharrow{ADMM applied to LASSO.}
Note that the Lagrangian LASSO problem \eqref{opt:ll} can be equivalently reformulated  using an auxiliary variable as follows:
$$
\min_{\bbeta\in\real^p} F(\bbeta)\triangleq f(\bbeta)+g(\balpha)
\triangleq
\frac{1}{2} \normtwo{\bX\bbeta-\by}^2 + \lambda \normone{\balpha} \quad \text{ s.t. }\bbeta=\balpha.
$$
This formulation fits the standard ADMM framework.
Given a smoothing parameter $\rho>0$, the updates at the $t$-th  iteration are:
$$
\begin{aligned}
\bbeta^\toptone
&\leftarrow \argmin_{\bbeta} \left\{ \frac{1}{2} \normtwo{\bX\bbeta-\by}^2 + \frac{\rho}{2} \normtwo{\bbeta - \balpha^\toptzero + \frac{1}{\rho} \bnu^\toptzero}^2 \right\}\\
&= \Big(\bX^\top \bX + \rho \bI\Big)^{-1} \Big(\bX^\top \by + \rho \balpha^\toptzero - \bnu^\toptzero\Big);\\
\balpha^\toptone
&\leftarrow \argmin_{\balpha} \left\{ \lambda \normone{\balpha} + \frac{\rho}{2} \normtwo{\bbeta^\toptone - \balpha + \frac{1}{\rho} \bnu^\toptzero}^2 \right\}\\
&= \prox_{(\lambda/\rho) \normone{\cdot}} \left( \bbeta^\toptone + \frac{1}{\rho} \bnu^\toptzero \right)
=\mathcalT_{(\lambda/\rho)}(\bbeta);\\
\bnu^\toptone
&\leftarrow \bnu^\toptzero +  \rho (\bbeta^\toptone - \balpha^\toptone),
\end{aligned}
$$
where $\mathcalT_\lambda(\cdot)$ denotes the soft-thresholding operator  (Example~\ref{example:soft_thres}).

Since $ \rho > 0 $, $ \bX^\top \bX + \rho \bI $ is always invertible. The update for $ \bbeta $  essentially solves  a ridge regression problem (an $ \ell_2 $-norm squared regularization least squares problem); while the update for $ \balpha $ involves  an $ \ell_1 $-norm proximal operator, which admits a closed-form solution via soft thresholding.
When solving the $ \bbeta $-subproblem, if a fixed penalty parameter $\rho$ is used throughout the iterations, one can precompute and cache a matrix factorization (e.g., Cholesky decomposition) of $ \bX^\top \bX + \rho \bI $; see, for example, \citet{lu2022matrix}. This significantly reduces the computational cost in subsequent iterations, as the expensive factorization need only be performed once.

\subsection{Block Coordinate Descent\index{Block coordinate descent}}

We briefly introduce how to solve the LASSO problem using the \textit{block coordinate descent (BCD)} method.
Because the  $\normone{\bbeta}$ term in the objective function is separable, each block variable corresponds to a single coordinate  of  $\bbeta$. 
For convenience, when updating the $i$-th block, we represent $\bbeta$ as:
$$
\bbeta =
\begin{bmatrix}
\beta_i \\
\bbeta_{-i}
\end{bmatrix},
$$
where $\beta_i$ denotes the $i$-th component of $\bbeta$, and $\bbeta_{-i}\in\real^{p-1}$ contains all components except the $i$-th.
Similarly, we partition the design matrix $\bX\in\real^{n\times p}$ as:
$$
\bX = \begin{bmatrix}
\bx_i & \bX_{-i}
\end{bmatrix},
$$
where $\bx_i\in\real^n$ is the $i$-th column of $\bX$, and $\bX_{-i}\in\real^{n\times {p-1}}$ consists of all columns except the  $i$-th.
This reordering does not change the problem---it merely isolates the variable being optimized.

During the update of the $i$-th coordinate, the LASSO objective becomes:
$$
\min_{\beta_i} \,
\lambda \abs{\beta_i} + \lambda \normone{\bbeta_{-i}} + \frac{1}{2} \normtwo{\bx_i \beta_i - (\by - \bX_{-i} \bbeta_{-i})}^2.
$$
Since $\bbeta_{-i}$ is fixed during this step, the term $\lambda \normone{\bbeta_{-i}}$ is constant and can be ignored for optimization purposes. Define the residual vector
$ \br_i \triangleq \by - \bX_{-i} \bbeta_{-i} $. 
Then the subproblem reduces to:
\begin{equation}\label{equation:bcd_lasso}
\min_{\beta_i} \, 
f_i(\beta_i) \triangleq \lambda \abs{\beta_i} + \frac{1}{2} \normtwo{\bx_i}^2 \beta_i^2 - \bx_i^\top \br_i \beta_i.
\end{equation}
The minimizer of \eqref{equation:bcd_lasso} has a closed-form solution given by soft thresholding:
\begin{equation}\label{equation:bcd_lasso_upd}
\beta_i^\toptzero \leftarrow \argmin_{\beta_i} f_i(\beta_i) =
\begin{cases}
\frac{\bx_i^\top \br_i - \lambda}{\normtwo{\bx_i}^2}, & \bx_i^\top \br_i > \lambda; \\
\frac{\bx_i^\top \br_i + \lambda}{\normtwo{\bx_i}^2}, & \bx_i^\top \br_i < -\lambda; \\
0, & \text{otherwise}.
\end{cases}
\end{equation}
Thus, the block coordinate descent algorithm for the LASSO problem proceeds by cyclically updating each coordinate using \eqref{equation:bcd_lasso_upd}. The full procedure is summarized in Algorithm~\ref{alg:bcd_lasso}.

\begin{algorithm}[H]
\caption{Block Coordinate Descent Method for LASSO}
\label{alg:bcd_lasso}
\begin{algorithmic}[1]
\State {\bfseries input:}   $ \bX, \by $, parameter $\lambda$;
\State {\bfseries initialize:}    $ \bbeta^\topone $;
\For{$t=1,2,\ldots$}
\For{$i = 1, 2, \ldots, p$}
\State Compute $\bbeta_{-i}, \br_i$ according to definitions;
\State Compute $ \beta_i^\toptzero $ by \eqref{equation:bcd_lasso_upd};
\EndFor
\State Stop if a stopping criterion is satisfied at iteration $t=T$;
\EndFor
\State {\bfseries return:} Final $\bbeta^{(T)}$;
\end{algorithmic}
\end{algorithm}

\index{Barzilai--Borwein (BB) rule}
\index{Bound-constrained quadratic program (BCQP)}
\subsection{Barzilai--Borwein Gradient Projection Algorithm}

Consider the standard \textit{bound-constrained quadratic program (BCQP)} problem
\begin{equation}\label{equation:bcqp_std}
\min_{\balpha} \left\{Q(\balpha)\triangleq\frac{1}{2}\balpha^\top\bY\balpha - \bz^\top \balpha  \right\}
\quad \text{s.t.} \quad \bm{\ell} \leq \balpha \leq  \bu.
\end{equation}
Projected gradient descent (PGD; see Algorithm~\ref{alg:pgd_gen}) methods provide a way of solving large-scale BCQP problems.
Let $\sS$ denote the feasible set of \eqref{equation:bcqp_std}:
\begin{equation}
\sS \triangleq \left\{ \balpha \in \real^p \mid \bm{\ell} \leq \balpha \leq \bu \right\}
\quad\implies\quad 
\project_{\sS}(\balpha) = \mathrm{mid}(\bm{\ell}, \balpha, \bu),
\end{equation}
where $\project_{\sS}$ denote the projection operator onto $\sS$, which  has a simple closed form; that is, $\project_{\sS}(\balpha)=\mathrm{mid}(\bm{\ell}, \balpha, \bu)$ is the vector whose $i$-th component is the median of the set $\{\ell_i, \alpha_i, u_i\}$. 
Given a current feasible iterate  $\balpha^\toptzero$, the projected gradient method computes the next iterate as
\begin{equation*}
\balpha^\toptone = \project_{\sS}(\balpha^\toptzero - \eta_t \bg^\toptzero),
\end{equation*}
where $\eta_t > 0$ is some stepsize and $\bg^\toptzero \triangleq \bY\balpha^\toptzero - \bz$ denotes the gradient; see Section~\ref{section:pgd}.

\citet{barzilai1988two} proposed a highly effective strategy for choosing the stepsize $\eta_t$, now known as the \textit{Barzilai--Borwein (BB) rule}. 
Two common variants are:
\begin{equation}
\eta_t^{BB1} = \frac{\bs^\toptminusTOP  \bs^\toptminus}{\bs^\toptminusTOP \bt^\toptminus}
\qquad \text{and}\qquad 
\eta_t^{BB2} = \frac{\bs^\toptminusTOP \bt^\toptminus}{\bt^\toptminusTOP \bt^\toptminus},
\end{equation}
where $\bs^\toptminus = \balpha^\toptzero - \balpha^\toptminus$ and $\bt^\toptminus = \bg^\toptzero - \bg^\toptminus$.

When combined with BB stepsizes, the projected gradient method is referred to as the \textit{Barzilai--Borwein gradient projection (BBGP)}  method. 
A notable feature of BBGP is that it is non-monotonic: the objective value $Q(\balpha^\toptzero)$   may increase on some iterations. Nevertheless, the method is simple, requires no matrix factorizations, and often exhibits fast practical convergence.

\paragrapharrow{Lagrangian LASSO with BB gradient projection.}
Interestingly, the Lagrangian LASSO problem \eqref{opt:ll}---originally an unconstrained convex optimization problem---can be reformulated as a BCQP. Specifically, consider
\begin{equation}\label{equation:laglasso_bb}
\min_{\bbeta} \left\{ \frac{1}{2} \normtwo{ \bX\bbeta - \by}^2 + \lambda \normone{\bbeta} \right\}, 
\quad
\bbeta \in \real^p, \bX \in \real^{n\times p}, \by \in \real^n.
\end{equation}
To transform this into a bound-constrained problem, we use the standard variable-splitting technique: for each component $\beta_i$, 
introduce nonnegative variables $u_i = (\beta_i)_+$ and $\ell_i = (-\beta_i)_+$ for all $i = 1,2,\ldots,p$, where $(\beta)_+ = \max(\beta,0)$ denotes the {positive-part operator}. Then, the coefficient vector $\bbeta$ can be split into a difference between nonnegative vectors:
\begin{equation}
\bbeta = \bu - \bm{\ell}, \quad \bu \geq \bzero, \quad \bm{\ell} \geq \bzero.
\end{equation}
And consequently, $\normone{\bbeta} = \bone_p^\top \bu + \bone_p^\top \bm{\ell}$, where $\bone_p = [1,1,\ldots,1]^\top$ is a $p\times 1$ summing vector consisting of $p$ ones. By setting $\bbeta = \bu - \bm{\ell}$, the Lagrangian LASSO problem \eqref{equation:laglasso_bb} can be rewritten as the following BCQP problem \citep{figueiredo2008gradient}:
\begin{equation}\label{equation:bcqp_lasso}
\min_{\balpha} \left\{ F(\balpha) \triangleq \bc^\top \balpha + \frac{1}{2} \balpha^\top \bY \balpha \right\} \quad \text{s.t.} \quad \balpha \geq \bzero,
\end{equation}
where 
\begin{equation}
\balpha = 
\begin{bmatrix}
\bu \\
\bm{\ell}
\end{bmatrix}, 
\qquad
\bc = \lambda \bone_{2p} + 
\begin{bmatrix}
-\bX^\top \by \\
\bX^\top \by
\end{bmatrix}, 
\qquad
\bY =
\begin{bmatrix}
\bX^\top \bX & -\bX^\top \bX \\
-\bX^\top \bX & \bX^\top \bX
\end{bmatrix}.
\end{equation}
To solve this BCQP problem~\eqref{equation:bcqp_lasso} with nonnegative constraint, \citet{figueiredo2008gradient} proposed the following two-step update:
\begin{subequations}
\begin{align}
\bw^\toptzero &\leftarrow \big(\balpha^\toptzero - \eta_t \nabla F(\balpha^\toptzero)\big)_+ ; \label{equation:bblasscv1}\\
\balpha^\toptone &\leftarrow  \balpha^\toptzero +\lambda_t (\bw^\toptzero - \balpha^\toptzero), \label{equation:bblasscv2}
\end{align}
\end{subequations}
where $\lambda_t\in[0,1]$.
Step \eqref{equation:bblasscv1} corresponds to a standard projected gradient descent step (projection onto the nonnegative orthant).
Step \eqref{equation:bblasscv2} introduces an extrapolation (or relaxation) similar in spirit to the conditional gradient method discussed in Section~\ref{section:cond_gd}.
The parameter $\lambda_t$ is typically chosen via a backtracking line search to ensure sufficient decrease in the objective $F$ \citep{bertsekas1997nonlinear}.

\section{Other Sparse Optimization Formulations}
\index{Generalized LASSO}

Apart from the standard LASSO problem, we also briefly introduce three important variants: the generalized LASSO, group LASSO, and sparse noise LASSO.

In the generalized LASSO, we assume that the regression coefficients exhibit a specific structural pattern---such as smoothness, piecewise constancy, or underlying trends. Rather than penalizing only the magnitude of the coefficients, this formulation penalizes their differences or higher-order variations. It is widely used in signal processing, image denoising, and time-series analysis, where such structural assumptions are natural.

The group LASSO is designed for settings in which predictors naturally form groups---for example, genes within biological pathways, dummy variables encoding a categorical feature, or multi-channel sensor measurements. In such cases, group LASSO selects or discards entire groups of variables simultaneously. This approach preserves model interpretability and respects the inherent grouping structure of the data.

The sparse noise LASSO (also known as robust LASSO) addresses scenarios where observations are corrupted by occasional large errors or outliers. Instead of assuming that all noise is small and Gaussian, it explicitly models sparse but significant (i.e., ``gross") errors in the data. This variant is particularly useful in real-world applications such as financial data analysis, sensor networks, or clinical measurements, where rare but severe anomalies can compromise data integrity.

\subsection{Generalized LASSO Problem\index{Generalized LASSO}}
The \textit{generalized LASSO problem} is defined as:
\begin{equation}\label{equation:gen_lasso}
\min_{\bbeta\in\real^p}  \,  \frac{1}{2} \normtwo{\bX\bbeta-\by}^2 +\lambda \normone{\bZ\bbeta},
\end{equation}
where $\bX\in\real^{n\times p}$, $\bZ\in\real^{l\times p}$, and $\by\in\real^n$.
In standard LASSO, the  $\normone{\bbeta}$ penalty promotes sparsity in  $\bbeta$. 
However, in many applications, $\bbeta$ itself may not be sparse, but becomes sparse after a linear transformation. 
A key example arises when $\bZ \in \real^{(p-1) \times p}$ is the first-order difference matrix:
$$
z_{ij} =
\begin{cases}
1, & j = i + 1, \\
-1, & j = i, \\
0, & \text{otherwise},
\end{cases}
$$
and $\bX = \bI$. 
In this case, the generalized LASSO reduces to:
$$
\min_{\bbeta\in\real^p}  \, \frac{1}{2} \normtwo{\bbeta - \by}^2 + \lambda \sum_{i=1}^{p-1} \abs{\beta_{i+1} - \beta_i},
$$
This formulation corresponds to the \textit{total variation (TV)} model used in image denoising. When $\bX=\bI$ and $\bZ$ is instead a second-order difference matrix, problem~\eqref{equation:gen_lasso} is known as \textit{trend filtering} \citep{rudin1992nonlinear, kim2009ell_1}.~\footnote{See {Problem~\ref{prob:denoise_rls}} for an illustration of its denoising effect.}

\paragrapharrow{Generalized LASSO with ADMM.}
In Section~\ref{section:laglass_admm}, we introduced the use of ADMM to solve the standard (Lagrangian) LASSO problem. Here, we extend that approach to the generalized LASSO.
To apply ADMM to problem~\eqref{equation:gen_lasso}, we introduce an auxiliary variable
$\balpha$ 
and enforce the constraint $\bZ\bbeta = \balpha$. 
This yields the equivalent constrained problem:
$$
\min_{\bbeta, \balpha} \quad \frac{1}{2} \normtwo{\bX\bbeta - \by}^2 + \lambda \normone{\balpha}
\quad \text{s.t.} \quad \bZ\bbeta - \balpha = \bzero.
$$
Let  $\bnu\in\real^l$ be the Lagrange multiplier associated with this constraint  (assuming $\bZ\in\real^{l\times p}$). 
The augmented Lagrangian is then:
$$
L_\rho(\bbeta, \balpha, \bnu) = \frac{1}{2} \normtwo{\bX\bbeta - \by}^2 + \lambda \normone{\balpha} + \bnu^\top (\bZ\bbeta - \balpha) + \frac{\rho}{2} \normtwo{\bZ\bbeta - \balpha}^2,
$$
where $\rho>0$ is the penalty parameter.
As in the standard LASSO case, the ADMM updates decouple the smooth quadratic term from the non-smooth $\ell_1$ penalty. Specifically, the $\bbeta$-update requires solving a linear system:
$$
\Big(\bX^\top \bX + \rho \bZ^\top \bZ\Big)\bbeta = \bX^\top \by + \rho \bZ^\top \Big( \balpha^\toptzero - \frac{1}{\rho}\bnu^\toptzero \Big),
$$
while the $\balpha$-update is performed via the proximal operator of the $\ell_1$-norm.
Thus, for a given $\rho>0$, the ADMM iterations (cf. Section~\ref{section:laglass_admm}) are:
\begin{subequations}
\begin{align}
\bbeta^\toptone &\leftarrow \Big(\bX^\top \bX + \rho \bZ^\top \bZ\Big)^{-1} \Big( \bX^\top \by + \rho \bZ^\top \Big( \balpha^\toptzero - \frac{1}{\rho}\bnu^\toptzero \Big) \Big);\\
\balpha^\toptone &\leftarrow \prox_{(\lambda/\rho) \normone{\cdot}} \Big( \bZ \bbeta^\toptone + \frac{1}{\rho}\bnu^\toptzero \Big);\\
\bnu^\toptone &\leftarrow \bnu^\toptzero +  \rho (\bZ\bbeta^\toptone - \balpha^\toptone).
\end{align}
\end{subequations}

\subsection{Group LASSO and Dual Polytope Projections\index{Group LASSO}}\label{section:group_las}
When the columns of $\bX$ can be naturally partitioned into disjoint subsets $\sG=\{g_1, g_2, \ldots, g_G\}$, 
where each subset corresponds to a coherent group of features (e.g., genes in a pathway or dummy variables from a categorical predictor), the group LASSO enforces sparsity at the group level rather than on individual coefficients \citep{yuan2006model}. Specifically, the groups cover all $p$ columns of $\bX$, i.e.,
$$
\bigcup_{\scriptsize g\in\sG}g=\{1,2,\ldots,p\},
$$
and typically the groups are assumed to be non-overlapping (though extensions exist for overlapping cases).
To achieve group-level sparsity, we apply an $\ell_1$ penalty to a vector in $\real^{\abs{\sG}}$, whose entries are the $\ell_2$-norms of the coefficient subvectors within each group. (Using the $\ell_2$-norm within groups would not induce group-wise selection on its own.) This strategy avoids penalizing individual coefficients directly and instead encourages entire groups of coefficients to be zero or nonzero together.

This approach is especially valuable when domain knowledge suggests that predictors form meaningful groups, enabling the model to select or discard whole sets of related features---enhancing both interpretability and statistical efficiency.

More formally, let $ \bX \in \real^{n\times p} $ and  $ \by \in \real^n $ denote the design matrix and response vector based on $n$ observations. 
The parameter vector $\bbeta\in\real^p$ is partitioned into $ G = \abs{\sG} $  blocks:
$$ 
\bbeta = [\bbeta_1^\top, \bbeta_2^\top, \ldots, \bbeta_G^\top]^\top \in \real^p,
$$
where $\bbeta_g\in\real^{p_g}$ corresponds to group $g$, with $\sum_{g=1}^G p_g = p$. 
The group LASSO assumes that only a few of the $\bbeta_g$ are nonzero.
The corresponding optimization problem is:
$$
\min_{\bbeta\in\real^p}  \, \frac{1}{2} \normtwo{\bX\bbeta - \by}^2 + \lambda \sum_{g=1}^G \sqrt{p_g} \normtwo{\bbeta_g},
$$
with feature groups $\bX\bbeta = \sum_{g=1}^G \bX_g \bbeta_g$, where $\bX_g \in \real^{n \times p_g}$ denoting the submatrix of $\bX$ containing the columns in group $g$.
The factor $\sqrt{p_g}$ compensates for varying group sizes, ensuring fair penalization across groups of different dimensions. The regularization parameter $\lambda>0$ controls the trade-off between data fidelity and group sparsity.
Note that the penalty uses the (unsquared) Euclidean norm $\normtwobig{\bbeta_g}$ for each group, making this a sum-of-norms regularization \citep{ohlsson2010segmentation}.

\subsection*{Group LASSO with ADMM}
Following the ADMM framework (see Section~\ref{section:laglass_admm}), we introduce auxiliary variables $ \balpha = [\balpha_1^\top, \balpha_2^\top, \ldots, \balpha_G^\top]^\top \in \real^p$ to decouple the smooth loss from the non-smooth group penalty. The problem is equivalently rewritten as:
$$
\min_{\bbeta, \balpha} \frac{1}{2} \normtwo{\bX\bbeta - \by}^2 + \lambda \sum_{g=1}^{G} \sqrt{p_g} \, \normtwo{\balpha_g} \quad \text{s.t.} \quad \bbeta_g = \balpha_g,\ \forall\, g.
$$
Let $ \bnu = [\bnu_1^\top, \bnu_2^\top, \ldots, \bnu_G^\top]^\top \in \real^p$
denote the scaled dual variables associated with these constraints. The scaled-form augmented Lagrangian is:
$$
L_{\rho}(\bbeta, \balpha, \bnu) = \frac{1}{2} \normtwo{\bX\bbeta - \by}^2 + \lambda \sum_{g=1}^{G} \sqrt{p_g} \, \normtwo{\balpha_g} 
+ \innerproduct{\bnu, \bbeta-\balpha}
+ \frac{\rho}{2} \sum_{g=1}^{G} \normtwo{\bbeta_g - \balpha_g}^2,
$$
where $\balpha_g$ are scaled Lagrange multipliers associated with the equality constraints $\bbeta_g = \balpha_g$, and $\rho > 0$ is a penalty parameter.

ADMM proceeds by iteratively updating $\bbeta$, $\balpha$, and $\bnu$; see \eqref{equation:admm_gen_up}.
The $\bbeta$-update is
$$
\bbeta^\toptone = \arg\min_{\bbeta} \frac{1}{2} \normtwo{\bX\bbeta-\by}^2 + \frac{\rho}{2} \normtwo{\bbeta - \balpha^\toptzero + \bnu^\toptzero}^2.
$$
This is a {ridge regression-type problem}, with a closed-form solution:
$$
\bbeta^\toptone = (\bX^\top \bX + \rho \bI)^{-1}(\bX^\top \by + \rho(\balpha^\toptzero - \bnu^\toptzero)).
$$
For large-scale problems, iterative solvers such as conjugate gradient or precomputed factorizations can be used to accelerate computation \citep{lu2025practical}.

The $\balpha$-update is
$$
\balpha^\toptone = \arg\min_\balpha \lambda \sum_{g=1}^G \sqrt{p_g} \normtwo{\balpha_g} + \frac{\rho}{2} \normtwo{\bbeta^\toptone - \balpha + \bnu^\toptzero}^2
$$
This decomposes {group-wise}. For each group $g\in\{1,2,\ldots,G\}$, the update is given by the proximal operator of the scaled $\ell_2$-norm
$$
\begin{aligned}
\balpha_g^\toptone 
&= \prox_{(\lambda \sqrt{p_g} / \rho) \normtwo{\cdot}} (\bbeta_g^\toptone + \bnu_g^\toptzero)\\
&=\begin{cases}
	\left(1 - \dfrac{\lambda \sqrt{p_g}}{\rho \normtwo{\bbeta_g^\toptone + \bnu_g^\toptzero}}\right)(\bbeta_g^\toptone + \bnu_g^\toptzero), & \text{if } \normtwo{\bbeta_g^\toptone + \bnu_g^\toptzero} > \dfrac{\lambda \sqrt{p_g}}{\rho}; \\
	0, & \text{otherwise},
\end{cases}
\end{aligned}
$$
which is known as the \textit{block soft-thresholding operator} (see Example~\ref{example:block_soft_thres}).

The scaled dual variable $\bnu$-update is
$$
\bnu^\toptone = \bnu^\toptzero + (\bbeta^\toptone - \balpha^\toptone)
$$

Combining these steps, the full ADMM iteration at step $t$ is:
\begin{subequations}
\begin{align}
\bbeta^\toptone &\leftarrow (\bX^\top \bX + \rho \bI)^{-1}(\bX^\top \by + \rho(\balpha^\toptzero - \bnu^\toptzero));  \\
\balpha_g^\toptone  &\leftarrow \prox_{(\lambda \sqrt{p_g} / \rho) \normtwo{\cdot}} (\bbeta_g^\toptone + \bnu_g^\toptzero), \quad \forall\, g; \\
\bnu^\toptone &\leftarrow \bnu^\toptzero + (\bbeta^\toptone - \balpha^\toptone). 
\end{align}
\end{subequations}

\subsection*{Group LASSO Dual}

Let again the design matrix be partitioned into $G$ groups, $\bX=[\bX_1,\bX_2,\ldots,\bX_g]$, where $\bX_g \in \real^{n \times p_g}$ and $p = \sum_{g=1}^G p_g$. 
The group LASSO problem is equivalently given by:
\begin{equation}\label{equation:group_lasso_opt}
\min_{\bbeta \in \real^p} \frac{1}{2} \normtwo{\by - \sum_{g=1}^G \bX_g \bbeta_g}^2 + \lambda \sum_{g=1}^G \sqrt{p_g} \normtwo{\bbeta_g}.
\end{equation}
Introducing the residual vector  $\br \triangleq  \by - \sum_{g=1}^G \bX_g \bbeta_g$, 
we reformulate the problem as:
\begin{equation}\label{equation:reformulated_group_lasso}
\min_{\bbeta, \br}  \frac{1}{2} \normtwo{\br}^2 + \lambda \sum_{g=1}^G \sqrt{p_g} \normtwo{\bbeta_g}
\quad \text{s.t.} \quad \br = \by - \sum_{g=1}^G \bX_g \bbeta_g.
\end{equation}
Let $\bnu \in \real^n$ denote the dual variable associated with the equality constraint. 
The Lagrangian of problem~\eqref{equation:reformulated_group_lasso} is:
\begin{equation}\label{equation:group_lasso_lagfunc}
L(\bbeta, \br, \bnu) = \frac{1}{2} \normtwo{\br}^2 + \lambda \sum_{g=1}^G \sqrt{p_g} \normtwo{\bbeta_g} 
+ \bnu^\top  \left(\by - \sum_{g=1}^G \bX_g \bbeta_g - \br\right).
\end{equation}
The dual function $D(\bnu)$ is obtained by minimizing the Lagrangian over $\bbeta$ and $\br$:
\begin{align*}
D(\bnu) &= \min_{\bbeta, \br} L(\bbeta, \br, \bnu) \\
&= \bnu^\top \by + \min_{\bbeta\in\real^p} \left( -\bnu^\top \sum_{g=1}^G \bX_g \bbeta_g + \lambda \sum_{g=1}^G \sqrt{p_g} \normtwo{\bbeta_g} \right) 
+ \min_{\br\in\real^n} \left( \frac{1}{2} \normtwo{\br}^2 - \bnu^\top \br \right).
\end{align*}
The last term matches the standard LASSO dual \eqref{equation:laglasso_dual2} and evaluates to $-\frac{1}{2}\normtwo{\bnu}^2$. 
We now analyze the second term.
Define
\begin{align*}
h(\bbeta) &= -\bnu^\top \sum_{g=1}^G \bX_g \bbeta_g + \lambda \sum_{g=1}^G \sqrt{p_g} \normtwo{\bbeta_g};\\
h_g(\bbeta_g) &= -\bnu^\top \bX_g \bbeta_g + \lambda \sqrt{p_g} \normtwo{\bbeta_g}, \qquad g = 1,2,\ldots,G.
\end{align*}
Since the groups are disjoint, and the minimization decouples across groups.
Then we can split $h(\bbeta)$ into a set of subproblems. Clearly $h_g(\bbeta_g)$ is convex but not smooth because it has a singular point at $\bzero$. 
Consider the subgradient of $h_g$,
$$
\partial h_g(\bbeta_g) = -\bX_g^\top \bnu + \lambda \sqrt{p_g} \bu_g, 
\quad \bu_g \in \partial \normtwo{\bbeta_g}, \qquad g = 1,2,\ldots,G,
$$
where $\bu_g$ is the subgradient of $\normtwo{\bbeta_g}$ (Exercise~\ref{exercise:sub_norms}):
\begin{equation}\label{equation:subgrad_beta_g}
\bu_g \in 
\begin{cases}
{\bbeta_g}/{\normtwo{\bbeta_g}}, & \text{if } \bbeta_g \neq \bzero; \\
\sB_2[\bzero ,1], & \text{if } \bbeta_g = \bzero,
\end{cases}
\end{equation}
and $\sB_2[\bzero ,1]$ denotes the unit $\ell_1$-ball.
The optimal solution $\widehatbbeta_g$ of $h_g$  satisfies
$$
\exists \widehatbu_g \in \partial \normtwobig{\widehatbbeta_g}, \quad -\bX_g^\top \bnu + \lambda \sqrt{p_g} \widehatbu_g = \bzero.
$$
If $\widehatbbeta_g = \bzero$, clearly $h_g(\widehatbbeta_g) = 0$. Otherwise, since $\lambda \sqrt{p_g} \widehatbu_g = \bX_g^\top \bnu$ and $\widehatbu_g = {\widehatbbeta_g}/{\normtwobig{\widehatbbeta_g}}$, we have
$$
h_g(\widehatbbeta_g) = -\lambda \sqrt{p_g} \frac{\widehatbbeta_g^\top}{\normtwobig{\widehatbbeta_g}} \bbeta_g + \lambda \sqrt{p_g} \normtwobig{\widehatbbeta_g} = 0.
$$
Therefore, this shows
$
\min_{\bbeta_g} h_g(\bbeta_g) = 0, \; g = 1,2,\ldots,G
$
and thus
$$
\min_\bbeta h(\bbeta) = \min_\bbeta \sum_{g=1}^G h_g(\bbeta_g) = \sum_{g=1}^G \min_{\bbeta_g} h_g(\bbeta_g) = 0.
$$
From~\eqref{equation:subgrad_beta_g}, it is easy to see $\normtwo{\bu_g}\leq 1$. Since $\lambda \sqrt{p_g} \widehatbu_g = \bX_g^\top \bnu$, we get a constraint on $\bnu$, i.e., $\bnu$ should satisfy:
\begin{equation}
\normtwo{\bX_g^\top \bnu}\leq \lambda \sqrt{p_g}, \quad g = 1,2,\ldots,G.
\end{equation}
Therefore, the dual function $D(\bnu)$ is:
$$
D(\bnu) = \bnu^\top \by - \frac{1}{2} \normtwo{\bnu}^2
=\frac{1}{2} \normtwo{\by}^2 - \frac{1}{2} \normtwo{\bnu - \by}^2.
$$
and the dual formulation of the group LASSO problem is:
\begin{equation}\label{equation:grouo_lasso_dual}
\max_{\bnu\in\real^n} D(\bnu) = \frac{1}{2} \normtwo{\by}^2 - \frac{1}{2} \normtwo{\bnu - \by}^2
\quad \text{s.t.} \quad \normtwo{\bX_g^\top \bnu} \leq \lambda \sqrt{p_g}, \;\forall \, g.
\end{equation}
As in the standard LASSO, a simple re-scaling of the dual variables $\bmu \triangleq \frac{\bnu}{\lambda}$ yields an equivalent formulation:
\begin{equation}\label{equation:grouo_lasso_dual_v2}
\max_{\bmu\in\real^n}  D(\bmu) = \frac{1}{2} \normtwo{\by}^2 - \frac{\lambda^2}{2} \normtwo{\bmu - \frac{\by}{\lambda}}^2
\quad \text{s.t.} \quad \normtwo{\bX_g^\top \bmu} \leq \sqrt{p_g}, \;\forall \, g.
\end{equation}

\paragrapharrow{Dual feasibility condition for sparsity.}
Clearly, problem~\eqref{equation:reformulated_group_lasso} is convex and its constraints are all affine. 
By Slater's condition, strong duality holds whenever the primal is feasible.
Let ($\widehatbbeta$, $\widehatbr$,  $\widehatbnu$) be optimal primal-dual solutions.
From the Lagrangian function \eqref{equation:group_lasso_lagfunc} and  KKT conditions (Theorem~\ref{theorem:opt_cond_sd}), we have
\begin{subequations}
\begin{align}
\bzero \in \partial_{\bbeta_g} L(\widehatbbeta, \widehatbr, \widehatbnu) 
&= - \bX_g^\top \widehatbnu + \lambda \sqrt{p_g} \bu_g,
\quad \text{where } \bu_g \in \partial \normtwobig{\widehatbbeta_g}, \quad \forall\,g; 
\label{equation:grouplasso_kkt_beta_group} \\
\nabla_{\br} L(\widehatbbeta, \widehatbr, \widehatbnu) 
&= \widehatbr -  \widehatbnu = \bzero;
\label{equation:grouplasso_kkt_z_group} \\
\nabla_\bnu L(\widehatbbeta, \widehatbr, \widehatbnu) 
&=  \by - \sum_{g=1}^G \bX_g \widehatbbeta_g - \widehatbr  = \bzero.
\label{equation:grouplasso_kkt_theta_group}
\end{align}
\end{subequations}
From~\eqref{equation:grouplasso_kkt_z_group} and~\eqref{equation:grouplasso_kkt_theta_group}, we have:
\begin{equation}
\by = \sum_{g=1}^G \bX_g \widehatbbeta_g +  \widehatbnu.
\end{equation}
From~\eqref{equation:grouplasso_kkt_beta_group}, we know there exists $\widehatbu_g \in \partial \normtwobig{\widehatbbeta_g}$ such that
\begin{equation}
\bX_g^\top \widehatbnu = \lambda\sqrt{p_g} \widehatbu_g.
\end{equation}
Then, for each group  $g = 1,2,\ldots,G$, the following holds:
$$
\bX_g^\top \widehatbnu \in 
\begin{cases}
\lambda\sqrt{p_g} \dfrac{\widehatbbeta_g}{\normtwobig{\widehatbbeta_g}}, & \text{if } \widehatbbeta_g \neq \bzero; \\
\lambda\sqrt{p_g} \sB[\bzero,1], & \text{if } \widehatbbeta_g = \bzero.
\end{cases}
$$
The dual problem is then important in the sense that it leads to a group-wise sparsity screening rule.
That is, the knowledge of $\widehatbnu$ allows us to identify the
zeros in $\widehatbbeta_g$ by checking the optimality condition:
\begin{equation}
\normtwo{\bX_g^\top \widehatbnu}< \lambda\sqrt{p_g}
\quad\implies\quad 
\widehatbbeta_g = \bzero.
\end{equation}
This condition is analogous to the standard LASSO dual feasibility condition \eqref{equation:lasso_stationarity_ori}, but enforces sparsity at the group level rather than on individual coefficients.
Note that this condition is sufficient but not necessary: it is possible to have  
$\widehatbbeta_g = \bzero$ even when $\normtwo{\bX_g^\top \widehatbnu}= \lambda\sqrt{p_g}$ (i.e., on the boundary of the dual feasible set).

\subsection*{Group LASSO with Dual Polytope Projections\index{Dual polytope projection}}
In Section~\ref{section:opcond_dpp}, we introduced dual polytope projection (DPP) rules to efficiently identify and eliminate inactive predictors---those with zero coefficients in the final solution---from the standard LASSO problem before solving it, by exploiting geometric properties in the dual space.
We now extend this idea to the group LASSO problem~\eqref{equation:group_lasso_opt}.

Consider the dual formulation~\eqref{equation:grouo_lasso_dual} or \eqref{equation:grouo_lasso_dual_v2}. 
It is easy to see that the dual optimal $\widehatbnu$ is the projection of ${\by}$ onto the feasible set. For each group $g$, the constraint $\normtwo{ \bX_g^\top \bnu} \leq \lambda \sqrt{p_g}$ confines $\bnu$ to an ellipsoid which is closed and convex. Therefore, the feasible set of the dual problem~\eqref{equation:grouo_lasso_dual} is the intersection of ellipsoids and thus closed and convex. Hence the dual solution $\widehatbnu(\lambda)$ under the penalty parameter $\lambda$  is also nonexpansive for the group LASSO problem.

Similarly to the standard LASSO problem, let $\lambda_{\max} \triangleq \max_g \normtwo{\bX_g^\top \by} / \sqrt{p_g}$, we can see that $\frac{\by}{\lambda_{\max}}$ is itself feasible, and $\lambda_{\max}$ is the largest parameter such that problem~\eqref{equation:group_lasso_opt} has a nonzero solution. That is, if $\lambda>\lambda_{\max}$, the primal solution $\widehatbbeta_g(\lambda)=\bzero$, $\forall\, g = 1, 2, \ldots, G$.
Clearly, the dual solution of \eqref{equation:grouo_lasso_dual} under $\lambda=\lambda_{\max}$ is $\widehatbnu(\lambda_{\max}) = {\by}$.

Building on the DPP theorem for the standard LASSO (Theorem~\ref{theorem:ddp}), we obtain the following result for the group LASSO.
\begin{theoremHigh}[Dual polytope projections (DPP) for group  LASSO]\label{theorem:ddp_grouplasso}
Suppose we are given a solution $\widehatbnu(\lambda')$ to the group LASSO dual~\eqref{equation:grouo_lasso_dual}  and a primal solution $\widehatbbeta(\lambda')$ for a specific $0<\lambda' \leq \lambda_{\max}$. 
Let $\lambda$ be a nonnegative value different from $\lambda'$. 
For each group $g\in\{1,2,\ldots,G\}$, if the following holds:
\begin{equation}\label{equation:group_lasso_dpp_thres}
\normtwo{\bX_g^\top \widehatbnu(\lambda')}
< \lambda'\sqrt{p_g} - \normf{\bX_g} \normtwo{\by} \abs{1 - \frac{\lambda'}{\lambda}}.
\end{equation}
then the $g$-th group coefficient vector satisfies $\widehatbbeta_g(\lambda) = \bzero$.
\end{theoremHigh}

\begin{subequations}
The above theorem shows that if  the following holds:
\begin{equation}\label{equation:group_lasso_dpp_thres2}
\normtwo{\bX_g^\top \big(\by-\sum_g\bX_g\widehatbbeta_g(\lambda')\big)}
< \lambda'\sqrt{p_g} - \normf{\bX_g} \normtwo{\by} \abs{1 - \frac{\lambda'}{\lambda}},
\end{equation}
then group $g$ is guaranteed to be inactive at $\lambda$, i.e.,  $\widehatbbeta_g(\lambda) = \bzero$.
In particular, when  $\lambda'=\lambda_{\max} \equiv  \max_g \normtwo{\bX_g^\top \by} / \sqrt{p_g}$, we have  
$\widehatbbeta_g(\lambda_{\max})=\bzero$ for all $g$, so the residual equals $\by$.
The above inequality becomes
\begin{equation}\label{equation:group_lasso_dpp_thres3}
\normtwo{\bX_g^\top {\by}} 
< \lambda_{\max}\sqrt{p_g} - \normf{\bX_g} \normtwo{\by} \left( \frac{\lambda_{\max}}{\lambda} - 1\right).
\end{equation}
\end{subequations}

\paragrapharrow{Sequential DPP for group LASSO.}
As in Corollary~\ref{corollary:seq_dpp}, we can derive a sequential DPP rule for the group LASSO.
Consider the group LASSO problem~\eqref{equation:group_lasso_opt} along a decreasing sequence of regularization parameters: $\lambda_{\max} = \lambda_0 > \lambda_1 > \ldots > \lambda_m$. For any integer $0 \leq t < m$, we have $\widehatbbeta_g(\lambda_{t+1}) = \bzero$ if $\widehatbbeta(\lambda_t)$ is known and the following holds:
\begin{equation}\label{equation:sdpp_grouplasso}
\normtwo{\bX_g^\top\big(\by - \sum_{g=1}^G \bX_g \widehatbbeta_g(\lambda_t)\big)} 
< \lambda_t\sqrt{p_g} - \normf{\bX_g} \normtwo{\by} \left( \frac{\lambda_t}{\lambda_{t+1}} - 1 \right).
\end{equation}
This sequential screening rule allows us to progressively eliminate irrelevant feature groups as $\lambda$ decreases, significantly reducing computational cost in path-following algorithms.

\subsection{Sparse Noise\index{Sparse noise}}
In the standard LASSO formulations \eqref{opt:ll} or \eqref{opt:lc} (p.~\pageref{opt:ll}), the residual (or noise) term $\by-\bX\bbeta$ is assumed to be dense.
In contrast, the following problem explicitly accounts for sparse noise:
\begin{equation}\label{equation:spar_noi_lass}
\min_{\bbeta \in \real^p} \left\{ F(\bbeta) \triangleq \normone{\bX \bbeta - \by} + \lambda \normone{\bbeta} \right\},
\end{equation}
where $\bX \in \real^{n\times p}$, $\by \in \real^n$, and $\lambda > 0$. 
We consider two optimization approaches to solve this problem:

\paragrapharrow{Subgradient descent.} 
Applying the subgradient method to~\eqref{equation:spar_noi_lass}, we use the subgradient of the $\ell_1$-norm given by $\partial \normone{\bbeta}\ni \sign(\bbeta)$ (see Exercise~\ref{exercise:sub_norms}). 
The update rule becomes:
$$
\bbeta^\toptone \leftarrow \bbeta^\toptzero - \eta_t \big(\bX^\top \sign(\bX \bbeta^\toptzero - \by) + \lambda \sign(\bbeta)\big).
$$
Here, $\eta_t$ denotes the stepsize. 
Following the guideline in Problem~\ref{prob:pgd_lipschitz} and the Lipschitzness of $\ell_1$-norm in Example~\ref{example:lipschitz_spar}, one may choose
$$
\eta_t = \frac{1}{\normtwo{F'(\bbeta^\toptzero)} \sqrt{t+1}},
$$ 
where $F'(\bbeta^\toptzero)$ is any subgradient of $F$ at $\bbeta^\toptzero$. 
Since the domain is unconstrained ($\sS=\real^p$), this method corresponds to an (unconstrained) subgradient descent algorithm.

\paragrapharrow{Proximal subgradient.}
We decompose the objective as $F(\bbeta)=f(\bbeta)+g(\bbeta)$, where $f(\bbeta) \triangleq \normone{\bX \bbeta - \by}$ and $g(\bbeta) \triangleq \lambda \normone{ \bbeta}$. 
The proximal (sub)gradient method then proceeds as:
$$
\bbeta^\toptone \leftarrow \prox_{\eta_t g}\big(\bbeta^\toptzero - \eta_t \bX^\top \sign(\bX \bbeta^\toptzero - \by)\big).
$$
Since $g(\bbeta) = \lambda \normone{\bbeta}$,
its proximal operator $\prox_{\eta_t g}$ is the well-known soft-thresholding operator. Specifically, by Example~\ref{example:soft_thres}, $\prox_{\eta_t g} = \mathcalT_{\lambda \eta_t}$,  and thus the general update rule becomes
$$
\bbeta^\toptone \leftarrow \mathcalT_{\lambda \eta_t}\big(\bbeta^\toptzero - \eta_t \bX^\top \sign(\bX \bbeta^\toptzero - \by)\big).
$$
A common choice for the stepsize is $\eta_t = \frac{1}{\normtwo{f'(\bbeta^\toptzero)} \sqrt{t+1}}$.

\begin{problemset}
\item Use the convergence of proximal gradient method in Theorem~\ref{theorem:prox_conv_ss_cvx} to prove the convergence of projected gradient descent method when applied to minimize a proper, closed, convex, and  $L_b$-smooth function $f$ over a closed convex set $\sS$.
What changes if we further assume that $f$ is also $L_a$-strongly convex?

\item \label{prob:pgd_smooth}\textbf{GD for convex and SS functions: $\mathcalO(1/T)$ \citep{lu2025practical}.}
Let $f:\real^p\rightarrow \real$ be a  differentiable, convex, and  $L_b$-smooth function defined on $\real^p$. 
Suppose $\{\bbeta^\toptzero\}_{t > 0}$ is the sequence generated by the gradient descent method (Algorithm~\ref{alg:pgd_gen} without projection) for solving   $\min_{\bbeta}f(\bbeta)$
with a constant stepsize $\eta = {1}/{L_b}$. 
Let $\bbeta^*$ be any minimizer  of $f$. 
Show that for any integer $T>0$,
$$
f(\bbeta^{(T)}) - f(\bbeta^*) \leq \frac{2L_b \normtwo{\bbeta^{(1)} -\bbeta^*}^2}{T-1}.
$$

\item \label{prob:pgd_lipschitz}\textbf{PGD for convex and Lipschitz functions: $\mathcalO(1/\sqrt{T})$ \citep{lu2025practical}}
Let $f:\sS\rightarrow \real$ be a proper, convex,  differentiable, and $L$-Lipschitz function defined over a closed, convex, and bounded set $\sS\subseteq\real^p$. 
Let $\bbeta^{(1)}, \bbeta^{(2)}, \ldots, \bbeta^{(T)}$ be the sequence of $T$ steps generated by the PGD method (Algorithm~\ref{alg:pgd_gen}) for solving  $\min_{\bbeta} f(\bbeta)$ subject to $\bbeta\in\sS$
with constant stepsize corresponding to  $\eta = \frac{R}{L\sqrt{T}}$, where $R$ is the upper bound on the distance $\normtwobig{\bbeta^{(1)}-\bbeta^{*}}$ from the initial point $\bbeta^{(1)}$ to some  optimal solution  $\bbeta^{*} \in \arg\min_{\bbeta\in\sS} f(\bbeta)$.
Show that for any $T\geq 1$,
$$
\small
\begin{aligned}
f\left(\frac{1}{T}\sum_{t=1}^{T}\bbeta^\toptzero\right) - f\left(\bbeta^{*}\right) \leqslant \frac{RL}{\sqrt{T}}.
\end{aligned}
$$
Show that this implies that for any integer $T$ satisfying $T\geq \frac{R^2L^2}{\epsilon^2}$, it follows that 
$$
\small
\begin{aligned}
f\left(\frac{1}{T}\sum_{t=1}^{T}\bbeta^\toptzero\right) - f\left(\bbeta^{*}\right) \leq \epsilon.
\end{aligned}
$$
Furthermore, show that if $R$ is unknown, a similar guarantee holds:
$$
f\left(\frac{1}{T}\sum_{t=1}^{T}\bbeta^\toptzero\right) - f\left(\bbeta^{*}\right) \leq \epsilon,
\quad \text{for $T=\mathcalO(\frac{1}{\epsilon^2})$ and $\eta = \frac{1}{\sqrt{T}}$}.
$$ 

\item  Discuss the \textit{Nesterov accelerated method} described  in Algorithm~\ref{alg:nesterov}. 
Using Theorem~\ref{theorem:fista_conv_ssf} as a guide, prove that this method achieves a convergence rate of $\mathcalO(1/T^2)$ for minimizing a convex and $L_b$-smooth function $f$,  when the stepsize is fixed as $\eta_t=1/L_b$ and the momentum parameter is set to $\gamma_t=2/(t+1)$.

\begin{minipage}[t]{0.98\linewidth}
\begin{algorithm}[H]
\caption{Nesterov Accelerated Method}
\label{alg:nesterov}
\begin{algorithmic}[1]
\State $\bbeta^\topone =\balpha^\topone = \widetildebbeta^\topone \in \real^p$, $\gamma_1=1$;
\For{$t=1,2,\ldots$}
\State Choose $\eta_{t}$ and $\gamma_{t+1}$;
\State  $\widetildebbeta^\toptone \leftarrow \prox_{\frac{\eta_t}{\gamma_t}g}\left(\widetildebbeta^\toptzero - \frac{\eta_t}{\gamma_t} \nabla f(\balpha^\toptzero)\right)$.
\State  $\bbeta^\toptone \leftarrow (1 - \gamma_t)\bbeta^\toptzero + \gamma_t \widetildebbeta^\toptone$.
\State  $\balpha^\toptone \leftarrow (1 - \gamma_{t+1})\bbeta^\toptone + \gamma_{t+1} \widetildebbeta^\toptone$.
\EndFor
\State Output  $\bbeta_{\text{final}}\leftarrow \bbeta^{(T)}$;
\end{algorithmic}
\end{algorithm}
\end{minipage}

\item \textbf{Proximal gradient method for nuclear norm regularization.}
Consider the composite objective function $F(\bX)=f(\bX)+g(\bX)$, where $f(\bX)=\normf{\bX-\bY}^2$ and $\bY\in\real^{n\times p}$  is a given matrix.
Suppose that $g(\bX) = \lambda \sum_i\sigma_i(\bX)$ is the nuclear norm (i.e., the sum of the singular values of $\bX$). 
Show that the proximal gradient update can be obtained by the following
procedure:
\begin{enumerate}[(a)]
\item Compute the singular value decomposition of the input matrix $\bY$, that is $\bY = \bU\bSigma\bV^\top$ where $\bSigma = \diag([\sigma_1(\bY), \sigma_2(\bY), \ldots])$ is a diagonal matrix of the singular values.
\item Apply the soft-thresholding operator (Example~\ref{example:soft_thres}) to compute the ``shrunken'' singular values
$$
\zeta_j \triangleq \mathcalT_{\lambda}(\sigma_j(\bY)), \quad \text{for } j = 1,2, \ldots
$$
\item Return the matrix $\widehat{\bY} = \bU \diag([\zeta_1,\zeta_2,\ldots]) \bV^\top$.
\end{enumerate}

\begin{minipage}[t]{0.98\linewidth}
\begin{algorithm}[H]
\caption{Proximal Point Method}
\label{alg:prox_point_gen}
\begin{algorithmic}[1] 
\Require A proper closed convex function $g$ ; 
\State {\bfseries input:}  Initialize $\bbeta^{(1)}$;
\For{$t=1,2,\ldots$}
\State Choose  a stepsize $\eta_t$;
\State $\bbeta^\toptone \leftarrow \prox_{\eta_t g}(\bbeta^\toptzero)$;
\EndFor
\State Output  $\bbeta_{\text{final}}\leftarrow \bbeta^{(T)}$;
\end{algorithmic} 
\end{algorithm}
\end{minipage}

\item \label{prob:proximal_point}\textbf{Proximal point method.}\index{Proximal point method}
Consider the unconstrained minimization problem
\begin{equation}\label{equation:prob_prox_point}
\min_{\bbeta \in \real^p} g(\bbeta),
\end{equation}
where $ g : \real^p \to (-\infty, \infty] $ is proper, closed, and convex. This is a special case of the composite problem~\eqref{equation:prox_comp_prob} with$ f \equiv 0 $. 
The proximal gradient update then reduces to
$$
\bbeta^\toptone \leftarrow \prox_{\eta_t g}(\bbeta^\toptzero).
$$
Setting  $ \eta_t = \eta>0$ for all $t $, yields the proximal point method, as described in Algorithm~\ref{alg:prox_point_gen}.
Assume that the optimal set of~\eqref{equation:prob_prox_point} is nonempty. 
Let $\bbeta^*$ be any optimal solution, and let $ \{ \bbeta^\toptzero \}_{t>0} $ be the sequence generated by the proximal point method with fixed stepsize $ \eta > 0 $. 
Show that 
\begin{enumerate}[(a)]
\item $ g(\bbeta^\toptzero) - g(\bbeta^*) \leq \frac{\normtwo{\bbeta^\topone - \bbeta^*}^2}{2\eta t} $  for $ t>0 $;
\item The sequence $ \{ \bbeta^\toptzero \}_{t>0} $ converges to some optimal point $\bbeta^*$ of $g$.
\end{enumerate}
\textit{Hint: Apply Theorem~\ref{theorem:prox_conv_ss_cvx} with $ f \equiv 0 $, which is a 0-smooth function.}

\item  Prove rigorously that the subproblem of the Frank--Wolfe (FW) method applied to the constrained LASSO problem~\eqref{equation:lc_fw_subprob} takes the stated form.

\item Building on the DPP rule for group LASSO (Theorem~\ref{theorem:ddp_grouplasso}) and the enhanced DPP rule for the standard LASSO (Corollary~\ref{corollary:enhanced_dpop}), derive an enhanced DPP screening rule for the group LASSO problem.

\item \label{prob:power_nmethod} \textbf{Power method \citep{lu2021numerical, lu2025practical}.}
The FW method  provides useful insight into the classical \textit{power method} for computing the dominant eigenpair of a matrix.
Consider the specific problem 
$$
\mathop{\min}_{\bbeta} f(\bbeta)=-\frac{1}{2} \bbeta^\top\bX\bbeta
\quad \text{s.t.}\quad 
\sS=\{\norm{\bbeta}_2\leq 1\},
$$
where $\bX$ is positive semidefinite. 
In the power method, one iteratively updates an estimate of the leading eigenvector and computes the associated Rayleigh quotient:
$$
\lambda^{{(t+1)}} \leftarrow \bbeta^{(t+1)\top}\bX\bbeta^\toptone.
$$
Discuss the connection between this procedure and the conditional gradient method applied to the above constrained problem.

\item \textbf{Failure of penalty function method.}
Consider the constrained optimization problem
$$
\begin{aligned}
& \min & -2x^2 + 9y^2 \quad
\text{s.t.} \quad  x = 1.
\end{aligned}
$$
By substituting the constraint  $x$, it is straightforward to verify that the unique optimal solution is $[1, 0]^\top$. Now consider the quadratic penalty function
$$
f_{\sigma}(x, y) = -2x^2 + 9y^2 + \frac{\sigma}{2}(x - 1)^2,
$$
where $\sigma>0$ is the penalty parameter.
Show that for any $\sigma \leq 4$, the penalty function is unbounded below.
This illustrates a key limitation of penalty function methods (Algorithm~\ref{alg:quad_pen_eq}): if the penalty parameter $\sigma$ is too small, the reduction in the original objective at infeasible points can outweigh the penalty incurred for violating the constraint. Consequently, the penalized problem may fail to have a minimizer or may yield solutions far from the true optimum. Therefore, the initial choice of $\sigma$ must be sufficiently large to ensure meaningful behavior.

\item \label{prob:pen_cong_glo_full}\textbf{Penalty function method with suboptimal subproblem \citep{lu2025practical}.}
Let $f(\bbeta)$ and $g_i(\bbeta)$, $i \in \mathcalE$ be continuously differentiable functions, and let $\sigma_t \rightarrow+\infty$ in Algorithm~\ref{alg:quad_pen_eq}. 
Suppose the solution $\bbeta^\toptone$ of the subproblem satisfies
$ \normtwobig{\nabla f_{\sigma_t}(\bbeta^\toptone)} \leq \varepsilon_t$ where $\varepsilon_t \rightarrow0$.
Furthermore, assume that  for any limit point $\bbeta^*$ of the sequence $\{\bbeta^\toptzero\}$, the gradients $\{ \nabla g_i(\bbeta^*) \}_{i \in \mathcalE}$ are linearly independent.
Show that $\bbeta^*$ is a KKT point  of the equality-constrained optimization problem  \eqref{equation:pr_in_constchap_equ}, and
$$
\lim_{t \rightarrow\infty} (\sigma_t g_i(\bbeta^\toptone)) = \lambda_i^*, \quad \forall i \in \mathcalE, 
$$
where $\lambda_i^*$ is the Lagrange multiplier corresponding to the constraint $g_i(\bbeta^*) = 0$.

\item 
\textbf{ADMM applied to matrix factorization \citep{lu2022matrix}.}
Consider the regularized matrix factorization problem
$$
\mathopmin{\bZ} \frac{1}{2}\normf{\bX-\bW\bZ}^2+r(\bZ)
$$
where $\bX\in\real^{n\times p}, \bW\in\real^{n\times k}$ are given, and $r(\cdot)$ is a (possibly nonsmooth) regularization function.
Introduce an auxiliary variable $\widetildebZ\in\real^{k\times p}$ to decouple the smooth least squares term from the regularizer. The problem becomes
\begin{equation}\label{equation:mf_admm_prob1}
\mathopmin{\bZ} \frac{1}{2}\normf{\bX-\bW\bZ}^2+r(\widetildebZ) 
\quad 
\text{s.t.}
\quad 
\bZ=\widetildebZ.
\end{equation}
Derive the  ADMM  update  steps for solving~\eqref{equation:mf_admm_prob1}.

\item Consider the unconstrained minimization of the quadratic function 
$
f(\bbeta) = \frac{1}{2} \bbeta^\top \bX \bbeta - \by^\top\bbeta + z,
$
where $ \bX \succ \bzero $ is a symmetric positive definite matrix, $ \by \in \real^p $, and $z\in\real$ is a constant.
\begin{enumerate}
\item Show that the optimal solution $ \bbeta^* $ exists and is unique, and express it explicitly in terms of $ \bX$ and $\by $.

\item Write the gradient descent update rule with a constant stepsize $\eta > 0$ for minimizing this function.

\end{enumerate}
\textit{Hint: Use Exercise~\ref{exercise:conv_quad}.}

\item \label{prob:denoise_rls} \textbf{Denoising via RLS.} 
Suppose we observe a noisy signal  $\by$  generated as $\by = \bbeta+\be$,  $\bbeta\in\real^p$  is the true (unknown) signal and
$\be$ is a noise vector. 
A naive reconstruction would solve
$
\min \normtwo{\bbeta-\by}^2,
$
whose solution is simply $\bbeta=\by$.
While mathematically optimal for this formulation, this estimate is useless in practice because it does not denoise the observation---it merely returns the noisy data.
To obtain a smoother, more plausible estimate, we incorporate a regularization term that penalizes abrupt changes between consecutive entries:
$
R(\bbeta) = \sum_{i=1}^{p-1} (\beta_i - \beta_{i+1})^2.
$
This encourages piecewise-smooth solutions.
Address the following:
\begin{itemize}
\item  Reformulate the denoising problem as a constrained least squares problem (e.g., by introducing auxiliary variables), and derive its closed-form solution.
\item Identify practical applications of this model. For instance, in financial time series analysis, the profit-and-loss (P\&L) of an asset typically evolves smoothly over consecutive days; thus, enforcing smoothness via $R(\bbeta)$ yields more realistic estimates than raw observations.
\end{itemize}

\item Discussion the relationship between the ADMM and the primal-dual method (see Chapter~\ref{chapter:algouni}).
In particular, explain how ADMM can be interpreted as a variant of primal-dual algorithms applied to problems with separable structure and linear constraints.
\end{problemset}

\newpage 
\chapter{Algorithms for Sparse Recovery  Problems}\label{chapter:spar_recov}
\begingroup
\hypersetup{
linkcolor=structurecolor,
linktoc=page,  
}
\minitoc \newpage
\endgroup

\lettrine{\color{caligraphcolor}I}
In Chapter~\ref{chapter:spar}, we explored a broad class of algorithms designed for sparse regression problems---tasks in which the goal is to estimate a sparse coefficient vector from noisy observations. We examined foundational methods such as iterative hard-thresholding (IHT) and proximal gradient descent, as well as more advanced techniques like the ADMM and coordinate descent. These algorithms are tailored to solve regularized optimization problems involving penalties such as the LASSO, group LASSO, and sparse noise models. They play a central role in high-dimensional statistical learning and have widespread applications in signal processing, machine learning, and econometrics.

Building on this foundation, we now turn to sparse recovery problems, which are closely related but differ in formulation and objective. While Chapter~\ref{chapter:spar} focused on regression settings---in which the response variable is continuous and modeled as a linear combination of predictors---this chapter addresses the broader challenge of reconstructing a sparse signal from incomplete or indirect measurements.

Here, we delve into specialized algorithms that leverage sparsity to recover signals from underdetermined linear systems, a core problem in compressed sensing, image reconstruction, and inverse problems. We begin with greedy methods such as matching pursuit (MP) and orthogonal matching pursuit (OMP), which iteratively select atoms from a dictionary to approximate the target signal. We then discuss thresholding-based approaches---including hard and soft thresholding---which are closely connected to the proximal operators introduced in Chapter~\ref{chapter:spar}.

Furthermore, we present iterative algorithms for $\ell_0$-minimization and its $\ell_1$-relaxation, such as the homotopy algorithm, iteratively reweighted least squares (IRLS), and Bregman iteration. These methods refine signal estimates through successive updates and extend the regularization and proximal optimization principles discussed earlier, adapting them specifically to the signal recovery context.

Thus, whereas Chapter~\ref{chapter:spar} emphasized estimation in regression frameworks, this chapter focuses on reconstruction in sensing and inverse problems---both domains relying heavily on sparsity assumptions and tools from convex and nonconvex optimization. The theoretical insights and algorithmic strategies developed in Chapters~\ref{chapter:recovery} and~\ref{chapter:spar} provide a natural foundation for understanding the more specialized recovery algorithms presented here. Together, these chapters offer a cohesive treatment of modern sparse modeling techniques across diverse application areas.

\section{Algorithms for $\ell_0$-Minimization}\label{section:algo_ell0}

\subsection{Greedy Algorithms\index{Greedy algorithm}}
We begin by introducing several greedy algorithms commonly used in compressed sensing.

The general approach to finding the exact sparse solution of an underdetermined linear system $\bX_{n\times p}\bbeta_{p\times 1} = \by_{n\times 1}$ (with $n \ll p$) under a sparsity constraint of $k$  involves solving the least squares  problem for every possible subset of $k$ columns of $\bX$: $\bX_{n\times k}\bbeta_{k\times 1} = \by$ (usually $n \gg k$) and then selecting the best solution among them. 
Because there are $\binom{p}{k}$ possible combinations of the overdetermined equations, this brute-force strategy is computationally prohibitive for even moderately large $p$ and $k$.

\textit{Greedy algorithms} circumvent this intractability by seeking a locally optimal choice at each step rather than a globally optimal solution \citep{foucart2013invitation}. Although they do not guarantee global optimality in all cases, greedy methods often yield either the exact sparse solution or a high-quality approximation in practice.

Typical greedy algorithms employ matching pursuit strategies, including:
\begin{enumerate}[(i)]
\item \textit{Basic matching pursuit (MP) method \citep{mallat1993matching}.} Rather than minimizing a global cost function, MP iteratively constructs a sparse approximation $\bbeta$ by expressing the observed signal $\by$ as a linear combination of a few selected columns---called atoms---from the dictionary matrix $\bX$. At each iteration, the atom most correlated with the current residual $\br=\by-\bX\bbeta$ is added to the active (support) set. If the residual norm decreases monotonically, convergence of the MP algorithm can be guaranteed. See Algorithm~\ref{alg:bmp}.

\item \textit{Orthogonal matching pursuit (OMP) \citep{pati1993orthogonal}.}
While standard MP does not enforce orthogonality between the residual and previously selected atoms, OMP improves upon this by ensuring that, at each iteration, the residual is orthogonal to all atoms in the current support set. This is achieved by solving a least squares problem over the selected indices. As a result, OMP guarantees optimality of the coefficients at each step, typically requires fewer iterations than MP, and exhibits greater numerical stability. The computational complexity of OMP is $\mathcalO(np)$ per iteration, and it can exactly recover any $k$-sparse signal provided $k \leq n / (2\ln p)$ (under suitable conditions on the measurement matrix). 
See Algorithm~\ref{alg:orth_match_pursuit}.

\item \textit{Regularized orthogonal matching pursuit (ROMP) \citep{needell2009uniform, needell2010signal}.}
Building on OMP, ROMP introduces a regularization step to enhance robustness. In each iteration, it first identifies a set of candidate atoms with significant correlation to the residual, then selects a subset of these candidates based on a regularization principle (e.g., grouping atoms with comparable correlation magnitudes). This structured selection helps ensure stable recovery under noise and satisfies theoretical guarantees similar to those of convex optimization methods.

\item \textit{Stagewise orthogonal matching pursuit (StOMP) \citep{donoho2007sparse}.}
StOMP is a faster, approximate variant of OMP designed for large-scale problems. Instead of selecting one atom per iteration, it adds multiple atoms whose correlations exceed a data-dependent threshold. This reduces the number of iterations at the expense of slightly lower approximation accuracy. With a per-iteration complexity of $\mathcalO(p)$, StOMP is well-suited for high-dimensional settings where speed is prioritized over exact recovery.

\item \textit{Compressive Sampling Matching Pursuit (CoSaMP) \citep{needell2009cosamp}.}  
The CoSaMP algorithm improves upon ROMP by incorporating both identification and pruning steps to maintain a fixed sparsity level throughout the iterations. It achieves higher reconstruction accuracy than OMP while maintaining a low computational complexity of $\mathcalO(p \ln^2 p)$. Under appropriate conditions on the measurement matrix, CoSaMP can recover any $k$-sparse signal with $k \leq n / \big(2\ln(1 + p/k)\big)$. See Algorithm~\ref{alg:cs_match_pursuit}.

\item \textit{Subspace basis pursuit (SP) \citep{dai2009subspace}.}
Algorithm~\ref{alg:sub_pursuit} presents the subspace pursuit method for compressed-sensing signal reconstruction. Unlike OMP, StOMP, and ROMP---which permanently add newly selected indices to the support set---SP maintains a fixed support size of $k$ (for a $k$-sparse signal) throughout the iterations. At each step, it updates the entire support set by merging previous estimates with new candidates and then pruning back to $k$ elements via least squares refinement. This dynamic updating allows SP to correct previously misidentified support elements, often leading to improved recovery performance.

\end{enumerate}

\begin{algorithm}
\caption{Matching Pursuit (MP)\index{Matching pursuit}}
\label{alg:bmp}
\begin{algorithmic}[1]
\Require Dictionary $\bX \in \real^{n\times p}$, measurement vector $\by \in  \real^n$, sparsity level $k$ (or tolerance $\varepsilon$);
\Ensure {Approximate the solution of \eqref{opt:p0}: $\argmin_{\bbeta} \normzero{\bbeta}$ subject to $\bX\bbeta = \by$;}
\State \textbf{initialize:} $\br^\topzero \gets \by$, $\sS \gets \emptyset$, $\bbeta \gets \bzero$;
\For{$t=0,1,2,\ldots$}
\State $i_t \gets \arg\max_{i \in \{1,\dots,p\}} \abs{\innerproduct{\br^\toptzero, \bx_i}}$; \Comment{(MP$_1$)} 
\State $\sS \gets \sS \cup \{i_t\}$; \Comment{(MP$_2$)} 
\State $\beta_{i_t} \gets \langle \br^\toptzero, \bx_{i_t} \rangle$;  \Comment{(MP$_3$)} 
\State $\br^\toptone \gets \br^\toptzero - \beta_{i_t} \bx_{i_t}$; \Comment{(MP$_4$)} 
\State Exit if $t\geq k$ or $\normtwo{\br^\toptzero} \leq \varepsilon$;
\EndFor
\State \Return $\bbeta$;
\end{algorithmic}
\end{algorithm}

\begin{algorithm}[H]
\caption{Compressive Sampling Matching Pursuit (CoSaMP)}
\label{alg:cs_match_pursuit}
\begin{algorithmic}[1]
\Require    Dictionary $ \bX\in\real^{n\times p}$, measurement vector $\by \in\real^n$, sparsity level $k$;
\State \textbf{initialize:}   $\bbeta^\topzero \in\sB_0[k]$ or simply $\bbeta^\topzero=\bzero$;
\State \textbf{define :} $\mathcalL_k(\bbeta) \triangleq \text{index set of $k$ largest entries of $\bbeta \in \real^p$ in modulus}$;
\State \textbf{define :} 
$\mathcalP_k(\bbeta) \triangleq \bbeta_{\mathcalL_k(\bbeta)}$;

\For{$t=0, 1,2,\ldots$}
\State $\sS^\toptone \gets \supp(\bbeta^\toptzero) \cup \mathcalL_{2k}\big(\bX^\top (\by-\bX\bbeta^\toptzero)\big)$;  \Comment{(CoSaMP$_1$)}
\State $\balpha^\toptone \gets \argmin_{\bbeta} \{\normtwo{\by-\bX\bbeta}, \supp(\bbeta) \subseteq \sS^\toptone\}$; \Comment{(CoSaMP$_2$)}
\State $\bbeta^\toptone \gets \mathcalP_k(\balpha^\toptone)$;  \Comment{(CoSaMP$_3$)}
\State Stop if a stopping criterion is satisfied at iteration $t=T$;
\EndFor
\State \Return Final $\bbeta^{(T)}$.
\end{algorithmic}
\end{algorithm}

\begin{algorithm}[ht]
\caption{Subspace Pursuit (SP) Algorithm\index{Subspace pursuit}}
\label{alg:sub_pursuit}
\begin{algorithmic}[1]
\Require Dictionary $\bX \in \real^{n\times p}$, measurement vector $\by \in  \real^n$, sparsity level $k$;
\Ensure {Approximate the solution of \eqref{opt:p0}: $\argmin_{\bbeta} \normzero{\bbeta}$ subject to $\bX\bbeta = \by$;}

\State \textbf{initialize:}
\begin{enumerate}[(i)]
\item $\sS^\topzero \gets \{\text{$k$ indexes corresponding to  largest-magnitude entries in $\bX^\top\by$}\}$;
\item $\br^\topzero \gets \by - \bX_{\sS^\topzero} \bX_{\sS^\topzero}^+ \by$, where  $\bA^+$ denotes the pseudo-inverse  of $\bA$;
\end{enumerate}
\For{$t=1,2,\ldots$}
\State \algoalign{$\sT^\toptzero \gets \sS^\toptminus \cup \{k$ indexes corresponding to  largest-magnitude entries  in $\bX^\top \br^\toptminus\}$;}

\State $\balpha^\toptzero \gets \argmin_{\bbeta} \{\normtwo{\by-\bX\bbeta}, \supp(\bbeta) \subseteq \sT^\toptzero\}$, i.e., compute $\balpha^\toptzero \gets \bX_{\sT^\toptzero}^+ \by$;

\State $\sS^\toptzero \gets \{k$ indexes corresponding to  largest-magnitude entries in $\balpha\}$;

\State $\bbeta^\toptzero \gets \argmin_{\bbeta} \{\normtwo{\by-\bX\bbeta}, \supp(\bbeta) \subseteq \sS^\toptzero\}$, i.e., compute $\bbeta^\toptzero \gets \bX_{\sS^\toptzero}^+ \by$;

\State $\br^\toptzero \gets \by - \bX_{\sS^\toptzero} \bX_{\sS^\toptzero}^+ \by$;

\State {Exit if} $\normtwo{\br^\toptzero} > \normtwo{\br^\toptminus}$, and let $\sS^\toptzero \gets \sS^\toptminus$;
\EndFor
\State \Return $\bbeta \gets \bbeta^\toptzero $.
\end{algorithmic}
\end{algorithm}

\subsection{Orthogonal Matching Pursuit Algorithm}
As mentioned previously, the \textit{orthogonal matching pursuit (OMP)} method is a greedy, stepwise least squares algorithm for sparse approximation.
At the $t$-th iteration, OMP adds one index to the current support estimate $\sS^\toptzero$
and updates the coefficient vector $\bbeta^\toptzero$ 
as the least squares solution supported on $\sS^\toptzero$ 
that best fits the measurements $\by$. A formal description of OMP is given in Algorithm~\ref{alg:orth_match_pursuit}, which we analyze in this subsection.

\begin{algorithm}[H]
\caption{Orthogonal Matching Pursuit (OMP)\index{Orthogonal matching pursuit}}
\label{alg:orth_match_pursuit}
\begin{algorithmic}[1]
\Require    Dictionary $ \bX\in\real^{n\times p}$, measurement vector $\by \in\real^n$;
\Ensure {Approximate the solution of \eqref{opt:p0}: $\argmin_{\bbeta} \normzero{\bbeta}$ subject to $\bX\bbeta = \by$;}

\State \textbf{initialize:}   $\sS^\topzero \gets \varnothing$, $\bbeta^\topzero \gets \bzero$;
\For{$t=0, 1,2,\ldots$}
\State $\sS^\toptone \gets \sS^\toptzero \cup 
\left\{i_{t+1} = \argmax_i \left\{\abs{\big(\bX^\top (\by-\bX\bbeta^\toptzero)\big)_i}, \forall\, i\right\}\right\}$;  \Comment{(OMP$_1$)} 
\State $\bbeta^\toptone \gets \argmin_{\bbeta} \{\normtwo{\by-\bX\bbeta}, \supp(\bbeta) \subseteq \sS^\toptone\}$; \Comment{(OMP$_2$)}  
\State Stop if a stopping criterion is satisfied at iteration $t=T$;
\EndFor
\State \Return Final $\bbeta^{(T)}$.
\end{algorithmic}
\end{algorithm}

\paragrapharrow{Stopping criteria \citep{tropp2010computational}.} 
A natural stopping condition for OMP is exact recovery:  $\bX\bbeta^\toptzero = \by$. 
However, in practice---due to measurement noise or numerical errors---it is more robust to stop when the residual energy falls below a tolerance $\varepsilon>0$, i.e., when
$$
\normtwo{\br^\toptzero}\triangleq  \normtwobig{\by - \bX\bbeta^\toptzero}\leq \varepsilon.
$$
Alternatively, one may stop when no column of $\bX$ exhibits significant correlation with the residual, i.e., when
$$
\norminf{\bX^\top\br^\toptzero} \leq \varepsilon.
$$
If an estimate $k$ of the true sparsity level is available, a third common stopping rule is to terminate after exactly $k$ iterations, yielding a $k$-sparse estimate $\bbeta^{(k)}$.
For example, if $\bX$ is a square orthogonal matrix, then OMP with this stopping criterion exactly recovers any $k$-sparse vector $\bbeta\in\real^p$ from $\by=\bX\bbeta$. In this case, the iterate $\bbeta^\toptzero$ coincides with the vector containing the $t$ largest (in magnitude) entries of $\bbeta$, with all other entries set to zero.

\paragrapharrow{Least squares.\index{Least squares}}
The update in Step (OMP$_2$) solves a least squares problem restricted to the current support $\sS^\toptone$:
\begin{equation}
\min_{\alpha} \normtwo{\by-\bX_{\sS^\toptone}\balpha}, \quad \balpha\in\real^{\abs{\sS^\toptone}.}
\end{equation} 
The solution satisfies   the \textit{normal equation} (see \eqref{equation:ne_firstkind}):
\begin{equation}\label{equation:ls_omp2}
\bX_{\sS^\toptone}^\top\bX_{\sS^\toptone} \bbeta_{\sS^\toptone} = \bX_{\sS^\toptone}^\top \by
\quad \implies \quad
\bX_{\sS^\toptone}^\top \big(\bX_{\sS^\toptone} \bbeta_{\sS^\toptone} - \by\big)=\bzero.
\end{equation}
This orthogonality condition---namely, that the residual is orthogonal to the column space of $\bX_{\sS^\toptone}$ explains the name \textit{orthogonal matching pursuit}.
This also shows that  the vector $\bX_{\sS^\toptone}\bbeta_{\sS^\toptone}$ is the {orthogonal projection} of $\by$ onto the column space of $\bX_{\sS^\toptone}$; see Section~\ref{section:ortho_proj_mat}.

Note that \eqref{equation:ls_omp2} also shows that as long as none of the stopping criteria are satisfied, the index $i_{t+1}$ selected in (OMP$_1$) does not belong to $\sS^\toptzero$, ensuring that the support grows by exactly one new element per iteration.

Step (OMP$_2$) is the most computationally expensive part of the OMP algorithm. However, it can be implemented efficiently by solving
\begin{equation}
\bbeta^\toptone \gets \argmin_{\bbeta} \normtwo{\by-\bX_{\sS^\toptone}\bbeta_{\sS^\toptone}},
\end{equation}
which is indeed a standard least squares problem. Computational efficiency can be further improved by maintaining a QR or Cholesky factorization of $\bX_{\sS^\toptone}$ and updating it incrementally as new columns are added \citep{lu2021rigorous}.

The selection of the index $i_{t+1}$ follows a greedy strategy: at each iteration, OMP chooses the atom most correlated with the current residual $\br^\toptzero=\by-\bX\bbeta^\toptzero$, i.e., the index maximizing $\abs{(\bX^\top (\by-\bX\bbeta^\toptzero))_i}$. The following proposition, applied with $\sS=\sS^\toptzero$ and $\bbeta_1=\bbeta^\toptzero$, provides insight into why this choice leads to a substantial reduction in the $\ell_2$-norm of the residual.

\begin{proposition}[Nonincreasing of iterations]\label{proposition:noninc_omp}
Let $\bX \in \real^{n\times p}$ be a matrix with $\ell_2$-normalized columns. Given a support set $\sS \subseteq \{1,2,\ldots,p\}$ and an index $i \in \{1,2,\ldots,p\}$, define
\begin{align*}
\bbeta_1 &= \argmin_{\bbeta} \{\normtwo{\by-\bX\bbeta}, \supp(\bbeta) \subseteq \sS\}, \\
\bbeta_2 &= \argmin_{\bbeta} \{\normtwo{\by-\bX\bbeta}, \supp(\bbeta) \subseteq \sS \cup \{i\}\}.
\end{align*}
Then the residual norms satisfy
$$
\normtwo{\by-\bX\bbeta_2}^2 \leq \normtwo{\by-\bX\bbeta_1}^2 - {\big(\bX^\top (\by-\bX\bbeta_1)\big)_i}^2.
$$
\end{proposition}
\begin{proof}[of Proposition~\ref{proposition:noninc_omp}]
Since any vector of the form $\bbeta_1 + \mu \be_i$ with $\mu \in \real$ is supported on $\sS \cup \{i\}$, the optimality of $\bbeta_2$ implies
$$
\normtwo{\by - \bX\bbeta_2}^2 \leq \min_{\mu \in \real} \normtwo{\by - \bX(\bbeta_1 + \mu \be_i)}^2.
$$
Expanding the right-hand side yields
\begin{align*}
\normtwo{\by - \bX(\bbeta_1 + \mu \be_i)}^2 
&= \normtwo{\by - \bX\bbeta_1}^2 + \mu^2 \normtwo{\bX\be_i}^2 - 2 \mu {\big(\bX^\top (\by - \bX\bbeta_1)\big)_i}.
\end{align*}
This is a quadratic function in $\mu$, minimized at $\mu^* = {(\bX^\top \big(\by - \bX\bbeta_1)\big)_i}$. This shows that
$$
\min_{\mu \in \real} \normtwo{\by - \bX(\bbeta_1 + \mu \be_i)}^2 
= \normtwo{\by - \bX\bbeta_1}^2 - {\big(\bX^\top (\by - \bX\bbeta_1)\big)_i}^2.
$$
Combining this with the earlier inequality completes the proof.
\end{proof}

More generally, the success of recovering $k$-sparse vectors using exactly $k$ iterations of the orthogonal matching pursuit algorithm is characterized by the following result.
\begin{theoremHigh}[Recovery under OMP]\label{theorem:rec_omp}
Let $\bX \in \real^{n\times p}$. 
Every nonzero vector $\bbeta^* \in \real^p$ supported on a set $\sS$ with $\abs{\sS}=k$ is exactly recovered from  $\by = \bX\bbeta^*$ in at most $k$ iterations of the OMP algorithm (Algorithm~\ref{alg:orth_match_pursuit}) if and only if the submatrix $\bX_\sS$ has full column rank, and
\begin{equation}\label{equation:rec_omp}
\max_{i \in \sS} \abs{(\bX^\top \bz)_i} > \max_{j \in \comple{\sS}} \abs{(\bX^\top \bz)_j},
\end{equation}
for all nonzero $\bz \in \{\bX\balpha, \supp(\balpha) \subseteq \sS\}$.
\end{theoremHigh}
\begin{proof}[of Theorem~\ref{theorem:rec_omp}]
Suppose that OMP recovers every vector supported on $\sS$ within $k = \abs{\sS}$ iterations.
First, note that if two distinct vectors supported on $\sS$ produced the same measurement $\by$, OMP could not distinguish between them---yet it is assumed to recover both exactly. Hence, $\bX_\sS$ must be injective, i.e., have full column rank; otherwise the difference between the two vectors will lie in the null space and thus the null space is not trivial.
Moreover, consider any nonzero $\by=\bX\bbeta^*$ with $\supp(\bbeta^*)=\sS$. Since OMP must select only indices from $\sS$ during its first iteration, no index $j\in\comple{\sS}$ can be chosen at this stage. This implies
$$
\max_{i \in \sS} \abs{(\bX^\top \by)_i} 
> \max_{j \in \comple{\sS}} \abs{(\bX^\top \by)_j}.
$$
Because every nonzero $\bz$ in the span of $\{\bx_i\mid i\in\sS\}$ can be written as $\bz=\bX\balpha$ for some $\balpha$ supported on $\sS$, the inequality \eqref{equation:rec_omp} holds for all such $\bz$.

Conversely, assume that $\bX_{\sS}$ has full column rank and that condition~\eqref{equation:rec_omp} holds. We prove by induction that, for all $0 \leq t \leq k$, the support estimate $\sS^\toptzero$ satisfies $\sS^\toptzero \subseteq \sS$ and $\abs{\sS^\toptzero} = t$.
This will imply $\sS^{(k)} = \sS$, hence $\bX\bbeta^{(k)}=\bX_{\sS}\bbeta^{(k)} = \by$ by (OMP$_2$), and finally $\bbeta^{(k)} = \bbeta^*$ by the full column rank of $\bX_\sS$. 

The base case $t = 0$ holds trivially since $\sS^{(0)} = \varnothing$. Now suppose the claim holds for some $t < k$. Define the residual at iteration $t$ as
$$
\bz^\toptzero \triangleq \by - \bX(\sS^\toptzero) \bbeta^\toptzero.
$$
\footnote{Note again that $\bX(\sS)\in\real^{n\times p}$ and $\bX_{\sS}\in\real^{n\times \abs{\sS}}$; see Definition~\ref{definition:matlabnotation}.} 
Since $\sS^\toptzero \subseteq \sS$ and $\bbeta^\toptzero$ is the least squares solution over $\sS^\toptzero$, it follows that $\bz^\toptzero \in \spn\{\bx_i \mid i \in \sS\}$, i.e., $\bz^\toptzero = \bX_{\sS} \boldsymbol{\alpha}$ for some $\boldsymbol{\alpha}$ supported on $\sS$ (since $\bbeta^*$ is supported on $\sS$ and $\bX_{\sS}$ has full column rank). By~\eqref{equation:rec_omp}, the index selected at step $t+1$,
$$
i_{t+1} = \arg\max_i \abs{(\bX^\top \bz^\toptzero)_i},
$$
must belong to $\sS$. Therefore, $\sS^\toptone = \sS^\toptzero \cup \{i_{t+1}\} \subseteq \sS$.

To show that $\abs{\sS^\toptone} = t+1$, observe that by the normal equations~\eqref{equation:ls_omp2}, the residual $\bz^\toptzero$ is orthogonal to all columns $\{\bx_i \mid i \in \sS^\toptzero\}$. Hence, $(\bX^\top \bz^\toptzero)_i = 0$ for all $i \in \sS^\toptzero$, which implies $i_{t+1} \notin \sS^\toptzero$. Thus, the support strictly grows by one new index at each step.

By induction, $\sS^{(k)} = \sS$. Then, by the least squares update (OMP$_2$), we have $\bX(\sS) \bbeta^{(k)} = \by$, and since $\bX_{\sS}$ has full column rank, it follows that $\bbeta^{(k)} = \bbeta^*$. This completes the proof.
\end{proof}

\begin{remark}[Exact recovery condition]
A more compact formulation of the necessary and sufficient condition in Theorem~\ref{theorem:rec_omp} is the \textit{exact recovery condition}:
\begin{equation}
\normone{\bX_\sS^+\bX_{\comple{\sS}}} < 1,
\end{equation}
where $\bX_\sS^+$ denotes the pseudo-inverse of $\bX_\sS$, and $\normone{\cdot} \equiv \norm{\cdot}_{1,1}$ denotes the induced matrix norm (Definition~\ref{definition:induced_norm_app}).
The existence of the pseudo-inverse $\bX_\sS^+ = (\bX_{\sS}^\top \bX_\sS)^{-1} \bX_{\sS}^\top$ is equivalent to  $\bX_\sS$ having  full column rank. 
Furthermore, condition \eqref{equation:rec_omp} is equivalent to
\begin{equation}
\norminf{\bX_{\sS}^\top \bX_\sS \balpha} > \norminf{\bX_{\comple{\sS}}^\top \bX_\sS \balpha}, \quad \text{for all } \balpha \in \real^k \setminus \{\bzero\}.
\end{equation}
Setting  $\bw \triangleq \bX_{\sS}^\top \bX_\sS \balpha$,  this becomes
$$
\norminf{\bw} 
> \norminf{\bX_{\comple{\sS}}^\top \bX_\sS (\bX_{\sS}^\top \bX_\sS)^{-1} \bw} 
= \norminf{\bX_{\comple{\sS}}^\top (\bX_\sS^+)^\top \bw}, \quad \text{for all } \bw \in \real^k \setminus \{\bzero\}.
$$
The inequality implies $\norminf{\bX_{\comple{\sS}}^\top (\bX_\sS^+)^\top} < 1$ by the matrix-vector product inequality \eqref{equation:induced_ineqy_intern}, which is equivalent to the condition $\normone{\bX_\sS^+ \bX_{\comple{\sS}}^\top} < 1$ by duality of induced norms (see \eqref{equation:mat_one_norm} and \eqref{equation:mat_inf_norm}).
\end{remark}

\paragrapharrow{CoSaMP.}
A key limitation of the OMP algorithm is that once an incorrect index is selected into the estimated support $\sS^\toptzero$, it remains in all subsequent supports $\sS^{(s)}$ for $s \geq t$. 
Consequently, if even a single wrong index is chosen during the first $k$ iterations, exact recovery of a $k$-sparse vector may fail---since OMP never removes indices from the support.

One possible remedy is to run more than $k$ iterations. However, this does not address the root issue: the lack of support correction. A more effective approach---when an estimate of the sparsity level $k$ is available---is to use algorithms that allow both addition and removal of indices from the support. One such method is the \textit{compressive sampling matching pursuit (CoSaMP)} algorithm, introduced earlier in this chapter \citep{needell2009cosamp}.

To describe CoSaMP, we first introduce two standard operators. For any vector $\bbeta \in \real^p$,  let
\begin{itemize}
\item $\mathcalL_k(\bbeta)$ denote index set of $k$ largest entries of $\bbeta \in \real^p$ in modulus.
\item $\mathcalP_k(\bbeta)=\bbeta_{\mathcalL_k(\bbeta)}$ denote the corresponding $k$-term approximation, i.e., the vector obtained by keeping only those $k$ entries and setting the rest to zero.
\end{itemize}
The mapping $\mathcalP_k \triangleq \mathcalP_{\sB_0[k]}$ is known as the \textit{hard-thresholding operator of order $k$}~\footnote{Note that $\mathcalP_{\sS}$ denotes the projection of the vector onto the set $\sS$; see Definition~\ref{definition:projec_prox_opt}.}.
Given the vector $\bbeta \in \real^p$, the operator $\mathcalP_k$ retains its $k$ largest-magnitude entries and zeros out the others. 
This operator is not always uniquely defined---for example, when multiple entries have identical magnitudes. To ensure well-definedness, we adopt a deterministic tie-breaking rule (e.g., lexicographic ordering of indices) to select $\mathcalL_k(\bbeta)$ from all admissible candidates.
Using these notations, the CoSaMP algorithm simply iteratively update using:
\begin{subequations}\label{equation:cosamp_eq}
\begin{align}
	\sS^\toptone &\gets \supp(\bbeta^\toptzero) \cup \mathcalL_{2k}\big(\bX^\top (\by-\bX\bbeta^\toptzero)\big);\\
	\balpha^\toptone &\gets \argmin_{\bbeta} \{\normtwo{\by-\bX\bbeta}, \supp(\bbeta) \subseteq \sS^\toptone\};\\
	\bbeta^\toptone &\gets \mathcalP_k(\balpha^\toptone).
\end{align}

\end{subequations}

In Sections~\ref{section:mutual_cohere}, \ref{section:ell1_cohere}, and \ref{section:spar_ana_coherence}, we show that the performance of sparse recovery algorithms improves when the coherence of the sensing matrix is small.
We now justify this claim in the context of the OMP algorithm.
\begin{theoremHigh}[OMP under $\ell_1$-coherence]\label{theorem:omp_recovery_mutocohe}
Let $\bX \in \real^{n \times p}$ be a matrix with $\ell_2$-normalized columns. If
\begin{equation}\label{eq:condition_omp}
\mu_1(\bX, k) + \mu_1(\bX, k - 1) < 1,
\end{equation}
then every $k$-sparse vector $\bbeta \in \real^p$ is exactly recovered from the measurement vector $\by = \bX\bbeta$ in at most  $k$ iterations of the OMP algorithm.
\end{theoremHigh}
\begin{proof}[of Theorem~\ref{theorem:omp_recovery_mutocohe}]
By Theorem~\ref{theorem:rec_omp}, it suffices to verify that for any support set $\sS \subseteq \{1,2,\ldots,p\}$ with  $\abs{\sS} = k$, the submatrix  $\bX_\sS$ has full column rank and that
\begin{equation}\label{eq:key_inequality}
\max_{i \in \sS} \abs{\innerproduct{\bz, \bx_i}} 
> 
\max_{j \in \comple{\sS}} \abs{\innerproduct{\bz, \bx_j}}, 
\quad\forall\,
\text{$\bz \in \{\bX\balpha, \supp(\balpha) \subseteq \sS\}\setminus\{\bzero\}$.}
\end{equation}
The full column rank of $\bX_\sS$ follows from Corollary~\ref{corollary:l1co_and_rip}.
Thus, it remains to establish inequality~\eqref{eq:key_inequality}.
Let  $\bz \triangleq \sum_{i \in \sS} \alpha_i \bx_i$ be a nonzero vector in the span of $\bX_{\sS}$, and let $q \in \sS$ be an index such that  $\abs{\alpha_q} = \max_{i \in \sS} \abs{\alpha_i} > 0$. 
For any $j \in \comple{\sS}$, by the definition of $\bz$ and the $\ell_1$-coherence, we obtain
$$
\abs{\innerproduct{\bz, \bx_j}}
= \abs{\sum_{i \in \sS} \alpha_i \innerproduct{\bx_i, \bx_j}} 
\leq \sum_{i \in \sS} \abs{\alpha_i} \abs{\innerproduct{\bx_i, \bx_j}} 
\leq \abs{\alpha_q} \mu_1(\bX, k), 
\quad \text{for $j \in \comple{\sS}$}.
$$
On the other hand, for the index $q\in\sS$, we have 
\begin{align*}
\abs{\innerproduct{\bz, \bx_q}}
&= \abs{ \sum_{i \in \sS} \alpha_i \innerproduct{ \bx_i, \bx_q}}
\geq \abs{\alpha_q} \abs{\innerproduct{\bx_q, \bx_q}} 
- \sum_{i \in \sS, i \neq q} \abs{\alpha_i} \abs{\innerproduct{\bx_i, \bx_q}}
\geq \abs{\alpha_q} - \abs{\alpha_q} \mu_1(\bX, k - 1).
\end{align*}
By assumption  \eqref{eq:condition_omp}, we have $1 - \mu_1(\bX, k - 1) > \mu_1(\bX, k)$, which implies  that  \eqref{eq:key_inequality} holds.
This completes the proof.
\end{proof}

\subsection{Thresholding-Based Algorithms}

If the sparsity level $k$ is very small, OMP is extremely fast, as its runtime depends essentially on the number of iterations---which equals $k$ when the algorithm succeeds. However, if $k$ is moderately large relative to the ambient dimension $p$, OMP can become computationally expensive.
However, the thresholding-based family of algorithms are  methods that do not require prior knowledge of the sparsity level $k$.

In this section, we describe additional algorithms that make use of the hard-thresholding operator $\mathcalL_k$ or $\project_k$; see its role in \eqref{equation:cosamp_eq}.
These methods form a distinct family, motivated by a simple yet powerful intuition: the action of the measurement matrix $\bX$ on a sparse vector can be approximately inverted by applying its transpose $\bX^\top$.

The simplest such method is the \textit{basic thresholding (BT)} algorithm. Given measurements $\by=\bX\bbeta\in\real^n$, where $\bbeta\in\real^p$ is $k$-sparse, BT estimates the support of $\bbeta$ as the set of indices corresponding to the $k$ largest (in absolute value) entries of the correlation vector $\bX^\top\by$. Once the support is identified, the algorithm computes the best least squares fit over that support. The procedure is summarized in Algorithm~\ref{alg:basic_threshold}.
Note that no initialization is required for BT, as the estimate is computed in a single step.

\begin{algorithm}[H]
\caption{Basic Thresholding (BT)\index{Basic thresholding}}
\label{alg:basic_threshold}
\begin{algorithmic}[1]
\Require    Dictionary $ \bX\in\real^{n\times p}$, measurement vector $\by \in\real^n$, sparsity level $k$;
\State \textbf{initialize:}   $\bbeta^\topzero \in\sB_0[k]$ or simply $\bbeta^\topzero=\bzero$;

\State $\sS^* \gets \mathcalL_k(\bX^\top \by)$; \Comment{(BT$_1$)}
\State $\bbeta^* \gets \argmin_{\bbeta} \{\normtwo{\by-\bX\bbeta}, \supp(\bbeta) \subseteq \sS^*\}$; \Comment{(BT$_2$)}

\State \Return Final $\bbeta\gets \bbeta^*$.
\end{algorithmic}
\end{algorithm}

A necessary and sufficient condition for exact recovery of a $k$-sparse vector via the BT algorithm closely resembles the condition for OMP (cf.~Theorem~\ref{theorem:rec_omp}  \eqref{equation:rec_omp}). It is stated below.

\begin{theoremHigh}[Recovery under BT\index{Exact recovery}]\label{theorem:rec_bt}
Let $\bbeta \in \real^p$ be a vector supported on a set $\sS$, and let $\by = \bX\bbeta$. 
Then $\bbeta$  is exactly recovered by the BT algorithm if and only if 
\begin{equation}\label{equation:rec_bt}
\min_{i \in \sS} \abs{(\bX^\top \by)_i} > \max_{j \in \comple{\sS}} \abs{(\bX^\top \by)_j}.
\end{equation}
\end{theoremHigh}
\begin{proof}[of Theorem~\ref{theorem:rec_bt}]
The algorithm recovers $\bbeta$ exactly if and only if the estimated support $\sS^*$ (computed in Step (BT$_1$) equals the true support $\sS$. This occurs precisely when every entry of $\bX^\top \by$ indexed by $\sS$ is strictly larger in magnitude than every entry indexed by the complement $\comple{\sS}$, which is exactly condition~\eqref{equation:rec_bt}. The result follows immediately.
\end{proof}

Although the recovery condition for BT is simple, it is rather restrictive.
A more sophisticated approach is the \textit{iterative hard-thresholding (IHT)} algorithm---an iterative method designed to solve the underdetermined linear system $\bX\bbeta=\by$ under the assumption that the solution $\bbeta$ is $k$-sparse.
\footnote{See Algorithm~\ref{alg:pgd_iht} for  its role in $\ell_0$-constraint problems.}

Rather than solving $\bX\bbeta=\by$ directly, IHT works with the normal equations of the least squares problem (see \eqref{equation:ne_firstkind}):
$$
\bX^\top \bX\bbeta = \bX^\top \by.
$$
This system can be rewritten as the fixed-point equation $\bbeta^\toptone = (\bI - \bX^\top \bX)\bbeta^\toptzero + \bX^\top \by$
Classical fixed-point iteration then suggests the update
$$
\bbeta^\toptone \leftarrow (\bI - \bX^\top \bX)\bbeta^\toptzero + \bX^\top \by.
$$
However, since we seek a $k$-sparse solution, we apply the hard-thresholding operator $\project_k$ at each step---keeping only the $k$ entries of largest magnitude and setting the rest to zero. This yields the IHT update rule.
Alternatively, it can be regarded as a standard IHT algorithm (Algorithm~\ref{alg:pgd_iht}) applied to the loss function $f(\bbeta)=\frac{1}{2}\normtwo{\bX\bbeta-\by}^2$ and the constraint set $\sS=\sB_0[k]$.
The resulting procedure is given in Algorithm~\ref{alg:basic_IHT}.

\begin{algorithm}[H]
\caption{Iterative Hard-Thresholding (IHT)\index{Iterative hard-thresholding}}
\label{alg:basic_IHT}
\begin{algorithmic}[1]
\Require    Dictionary $ \bX\in\real^{n\times p}$, measurement vector $\by \in\real^n$, sparsity level $k$;
\State \textbf{initialize:}   $\bbeta^\topzero \in\sB_0[k]$ or simply $\bbeta^\topzero=\bzero$;
\For{$t=0, 1,2,\ldots$}
\State $\bbeta^\toptone 
\leftarrow \mathcalP_k(\bbeta^\toptzero - \eta_t \bX^\top (\bX\bbeta^\toptzero-\by))$; 
\Comment{$\eta_t=1$ is a common choice.}
\State Stop if a stopping criterion is satisfied at iteration $t=T$;
\EndFor
\State \Return Final $\bbeta^{(T)}$.
\end{algorithmic}
\end{algorithm}

A key advantage of IHT is that it avoids computing orthogonal projections onto subspaces spanned by subsets of columns of $\bX$ (i.e., a least squares fit), which can be computationally expensive.
However, if we are willing to incur the cost of such projections---as in greedy methods like OMP---it is natural to ask: once a support has been selected, why not compute the best possible fit over that support?
This idea leads to the \textit{hard-thresholding pursuit (HTP)} algorithm, which combines the support selection of IHT with a least squares refinement step.

\begin{algorithm}[H]
\caption{Hard Thresholding Pursuit (HTP)}
\label{alg:hard_TP}
\begin{algorithmic}[1]
\Require    Dictionary $ \bX\in\real^{n\times p}$, measurement vectornal $\by \in\real^n$, sparsity level $k$;
\State \textbf{initialize:}   $\bbeta^\topzero \in\sB_0[k]$ or simply $\bbeta^\topzero=\bzero$;
\For{$t=0, 1,2,\ldots$}
\State $\sS^\toptone \leftarrow \mathcalL_k\big(\bbeta^\toptzero + \bX^\top (\by-\bX\bbeta^\toptzero)\big)$; \Comment{(HTP$_1$)}
\State $\bbeta^\toptone \leftarrow \argmin_{\bbeta} \{\normtwo{\by-\bX\bbeta}, \supp(\bbeta) \subseteq \sS^\toptone\}$; \Comment{(HTP$_2$)} 
\State Stop if a stopping criterion is satisfied at iteration $t=T$;
\EndFor
\State \Return Final $\bbeta^{(T)}$.
\end{algorithmic}
\end{algorithm}

To analyze the sparse recovery properties of these algorithms, we follow \citet{foucart2013invitation} and continue to use the $\ell_1$-coherence conditions.
Recovery guarantees under the RIP are discussed further in Problems~\ref{prob:iht_rip_convergence}$\sim$\ref{prob:iht_delta_3k_bound}.

\begin{theoremHigh}[Exact recovery of BT under $\ell_1$-coherence\index{Exact recovery}\index{$\ell_1$-coherence}]\label{theorem:bt_l1cohere}
Let $\bX \in \real^{n \times p}$ be a matrix with $\ell_2$-normalized columns and let $\bbeta^* \in \real^p$ be a $k$-sparse vector supported on a set $\sS\subseteq\{1,2,\ldots,p\}$ with $\abs{\sS}=k$. 
If
\begin{equation}\label{eq:condition_thresholding}
\mu_1(\bX,k) + \mu_1(\bX,k - 1) < \frac{\min_{i \in \sS} \abs{\beta_i^*}}{\max_{i \in \sS} \abs{\beta_i^*}},
\end{equation}
then the vector $\bbeta^* \in \real^p$ is exactly recovered from the measurement vector $\by = \bX\bbeta^*$ by the BT algorithm.
\end{theoremHigh}
\begin{proof}[of Theorem~\ref{theorem:bt_l1cohere}]
By  Theorem~\ref{theorem:rec_bt}, it suffices to show that for every $i \in \sS$ and every  $j \in \comple{\sS}$,
\begin{equation}\label{eq:key_ineq_thresholding}
\abs{\innerproduct{\bX\bbeta^*, \bx_i}} 
> 
\abs{\innerproduct{\bX\bbeta^*, \bx_j}}.
\end{equation}
Since $\bbeta^*$ is supported on $\sS$, 
we have  for any $i\in\sS$,
\begin{align*}
\abs{\innerproduct{\bX\bbeta^*, \bx_i}}
&= \abs{\sum_{q \in \sS} \beta_q^* \innerproduct{\bx_q, \bx_i}}
\geq \abs{\beta_i^*} 
- \sum_{q \in \sS, q \neq i} \abs{\beta_q^*} \abs{\innerproduct{\bx_q, \bx_i}}
\geq \abs{\beta_i^*} - \mu_1(\bX,k - 1) \max_{q \in \sS} \abs{\beta_q^*}.
\end{align*}
For any $j\in\comple{\sS}$, we obtain:
\begin{align*}
\abs{\innerproduct{\bX\bbeta^*, \bx_j}} 
=\abs{\innerproduct{\bX_\sS\bbeta^*_\sS, \bx_j}} 
= \abs{\sum_{q \in \sS} \beta_q^* \innerproduct{\bx_q, \bx_j}}
\leq \sum_{q \in \sS} \abs{\beta_q^*} \abs{\innerproduct{ \bx_q, \bx_j}} 
\leq \mu_1(\bX,k) \max_{q \in \sS} \abs{\beta_q^*}.
\end{align*}
Combining these two bounds and using assumption~\eqref{eq:condition_thresholding}, we find:
$$
\abs{\innerproduct{\bX\bbeta^*, \bx_i}} 
- \abs{\innerproduct{\bX\bbeta^*, \bx_j}} 
\geq \min_{q \in \sS} \abs{\beta_q^*} - \left( \mu_1(\bX,k) + \mu_1(\bX,k - 1) \right) \max_{i \in \sS} \abs{\beta_q^*} 
> 0.
$$
This establishes \eqref{eq:key_ineq_thresholding} and completes the proof.
\end{proof}

We now turn to the more sophisticated HTP algorithm. As with the OMP algorithm, we show that at most $k$ iterations suffice to recover any $k$-sparse vectors under a condition analogous to \eqref{eq:condition_omp}. 

\begin{theoremHigh}[Exact recovery of HTP under $\ell_1$-coherence]\label{theorem:htp_recovery}
Let $\bX \in \real^{n \times p}$ be a matrix with $\ell_2$-normalized columns. If
$$
2\mu_1(\bX,k) + \mu_1(\bX,k - 1) < 1,
$$
then every $k$-sparse vector $\bbeta^* \in \real^p$ supported on $\sS\subseteq \{1,2,\ldots,p\}$ with $\abs{\sS}=k$ is exactly recovered from the measurement vector $\by = \bX\bbeta^*$ in most $k$ iterations of the HTP algorithm.
\end{theoremHigh}
\begin{proof}[of Theorem~\ref{theorem:htp_recovery}]
Reorder the indices so that $i_1, i_2, \ldots, i_p$ such that
$$
\abs{\beta_{i_1}} \geq \abs{\beta_{i_2}} \geq \ldots \geq  \abs{\beta_{i_k}} > 
\abs{\beta_{i_{k+1}}} = \ldots = \abs{\beta_{i_p}} = 0.
$$
We prove by induction that for each iteration $t$ satisfying $0 \leq t \leq k - 1$, the set $\{i_1,i_2, \ldots, i_{t+1}\}$ is included in $\sS^\toptone$ defined by (HTP$_1$) with measurements $\by = \bX\bbeta^*$ as the set of largest absolute entries of
\begin{equation}\label{eq:z_n_plus_1}
\bz^\toptone 
\triangleq \bbeta^\toptzero + \bX^\top(\by - \bX\bbeta^\toptzero)
=\bbeta^\toptzero + \bX^\top\bX(\bbeta^* - \bbeta^\toptzero).
\end{equation}
This will imply that $\sS^{(k)} = \sS = \supp(\bbeta^*)$, and then step (HTP$_2$)---which computes the least squares solution on $\sS^{(k)}$---yields $\bbeta^{(k)} = \bbeta^*$.
To see this, it suffices to establish the key inequality:
\begin{equation}\label{eq:key_inequality_htp}
\min_{1 \leq q \leq t+1} \abs{z^\toptone_{i_q}} 
> 
\max_{j \in \comple{\sS}} \abs{z^\toptone_j}.
\end{equation}
That is, (HTP$_1$) selects the most significant coefficients.
We notice that, for every $i \in \{1,2,\ldots,p\}$,
\begin{align}
&z_i^\toptone
= \beta_i^\toptzero + \sum_{q=1}^p (\beta_q^* - \beta_q^\toptzero)\innerproduct{\bx_q, \bx_i}
= \beta_i^* + \sum_{q \neq i} (\beta_q^* - \beta_q^\toptzero)\innerproduct{\bx_q, \bx_i}\\
&\implies 
\abs{z_i^\toptone- \beta_i^*}
\leq \sum_{q \in \sS^\toptzero, q \neq i} \abs{\beta_q^* - \beta_q^\toptzero}
\abs{\innerproduct{\bx_q, \bx_i }} 
+ \sum_{q \in {\sS} \setminus \sS^\toptzero, q \neq i} 
\abs{\beta_q^*}\abs{\innerproduct{\bx_q, \bx_i}}.       \label{eq:bound_z_j}
\end{align}
The inequality holds since $\beta_q^\toptzero=0$ for $q\in\sS\setminus\sS^\toptzero$ by (HTP$_2$).
Since $t\leq k-1$ and $\abs{\sS}=k$,
the right-hand side of \eqref{eq:bound_z_j} can be divided into
\begin{equation}\label{eq:bound_z_j22}
\leq 
\begin{cases}
\sum_{q \in \sS^\toptzero, q \neq i} \abs{\beta_q^* - \beta_q^\toptzero}
\mu_1(\bX,k-1) 
+ \sum_{q \in {\sS} \setminus \sS^\toptzero} 
\abs{\beta_q^*}\mu_1(\bX,k-1),& i\in\sS^\toptzero;\\
\sum_{q \in \sS^\toptzero} \abs{\beta_q^* - \beta_q^\toptzero}
\mu_1(\bX,k)
+ \sum_{q \in {\sS} \setminus \sS^\toptzero, q \neq i} 
\abs{\beta_q^*}\mu_1(\bX,k-1),
& i\notin\sS^\toptzero.
\end{cases}
\end{equation}

We now prove \eqref{eq:key_inequality_htp} by induction.
Initially, for iteration $t = 0$, $\bbeta^\topzero)_{\sS^\topzero}=\bzero$ and $\sS^\topzero=\varnothing$. 
Substituting $\norminfbig{(\bbeta^* - \bbeta^\topzero)_{\sS^\topzero}} = 0$ into \eqref{eq:z_n_plus_1} and \eqref{eq:bound_z_j}, and using the fact that $2\mu_1(\bX,k) < 1$ implied by the assumption, we have
$$
\absbig{z_{i_1}^\topone} 
\geq (1 - \mu_1(\bX,k))\norminf{\bbeta^*} 
> \mu_1(\bX,k)\norminf{\bbeta^*} 
\geq \ab\absbig{z_j^\topone} ,
\quad \text{for all } j \in \comple{\sS}.
$$
The base case of \eqref{eq:key_inequality_htp} holds for $t = 0$. 
Assume that \eqref{eq:key_inequality_htp} holds for $t - 1$ with $t \geq 1$. 
This implies that $\{i_1, i_2, \ldots, i_t\} \subseteq \sS^\toptzero$. 
By Step (HTP$_2$), the residual $\by - \bX\bbeta^\toptzero$ is orthogonal to the space $\{\bX\bz \mid  \supp(\bz) \subseteq \sS^\toptzero\}$, i.e., that $\langle \by - \bX\bbeta^\toptzero, \bX\bz \rangle = \langle \bX^\top(\by - \bX\bbeta^\toptzero), \bz \rangle = 0$ for any $\bz \in \real^p$ supported on $\sS^\toptzero$. 
Substituting $\by = \bX\bbeta^*$ yields
$$
\left(\bX^\top\bX(\bbeta^* - \bbeta^\toptzero)\right)_{\sS^\toptzero} = \bzero.
$$
Therefore, for all  $i \in \sS^\toptzero$, we have $z_i^\toptone = \beta_i$ by the \eqref{eq:z_n_plus_1}. 
Using~\eqref{eq:bound_z_j} or \eqref{eq:bound_z_j22}, we obtain
$$
\abs{\beta_i^\toptzero - \beta_i^*} 
\leq \mu_1(\bX,k - 1)\norminf{(\bbeta^*- \bbeta^\toptzero)_{\sS^\toptzero}}
+ \mu_1(\bX,k - 1)\norminf{\bbeta^*_{{\sS} \setminus \sS^\toptzero}}, 
\quad\forall\, i \in \sS^\toptzero.
$$
Taking the maximum over $i \in \sS^\toptzero$ and rearranging terms:
\begin{equation}\label{eq:bound_x_minus_xn}
\norminf{(\bbeta^*- \bbeta^\toptzero)_{\sS^\toptzero}} \leq \frac{\mu_1(\bX,k - 1)}{1 - \mu_1(\bX,k - 1)}\norminf{\bbeta^*_{{\sS} \setminus \sS^\toptzero}}.
\end{equation}
Combining this with~\eqref{eq:z_n_plus_1} and~\eqref{eq:bound_z_j}, we obtain, for $1 \leq q \leq t+1$ and $j \in \comple{\sS}$,
\begin{align*}
\abs{z^\toptone_{i_q}} 
&\geq \left(1 - \frac{\mu_1(\bX,k)}{1 - \mu_1(\bX,k - 1)}\right)\abs{\beta^*_{i_{t+1}}}, \\
\abs{z^\toptone_j} 
&\leq \frac{\mu_1(\bX,k)}{1 - \mu_1(\bX,k - 1)} \abs{\beta^*_{i_{t+1}}}.
\end{align*}
By assumption, $\mu_1(\bX,k)/(1 - \mu_1(\bX,k - 1)) < 1/2$. 
Thus, \eqref{eq:key_inequality_htp} holds for iteration $t$, completing the induction.
\end{proof}

\section{Algorithms for $\ell_1$-Minimization}
We now turn to algorithms for $\ell_1$-minimization problems, which serve as a cornerstone of signal recovery in compressed sensing. As noted in model~\eqref{equation:prob_p1_equiv} (p.~\pageref{equation:prob_p1_equiv}), in the real-valued case, $\ell_1$-minimization can be reformulated as a linear program (LP).

Standard LP solvers for linear programs---based on interior-point methods or the classical simplex method---are readily available and widely used \citep{boyd2004convex, beck2014introduction}.
However, in practice, directly solving the LP formulation~\eqref{equation:prob_p1_equiv} is often challenging. In typical compressed sensing applications, the measurement matrix $\bX\in\real^{n\times p}$ is large (with $p \approx 10^6$) and dense, and the resulting LP is frequently ill-conditioned. Moreover, interior-point methods---which require factorizing matrices involving $\bX$, such as $\bX^\top\bX$ or augmented systems of size proportional to $n+p$---become computationally prohibitive at this scale. Consequently, general-purpose interior-point solvers are not practical for solving \eqref{equation:prob_p1_equiv} in large-scale compressed sensing problems.

While such general-purpose solvers are reliable and easy to use, they are not tailored to the specific structure of $\ell_1$-minimization problems. In fact, algorithms designed specifically for $\ell_1$-minimization can be significantly faster than generic linear programming approaches. This section introduces and analyzes several such specialized algorithms.

\subsection{Homotopy Algorithm\index{Homotopy method}}

We begin with the \textit{homotopy method}, which resembles OMP in spirit but---unlike OMP---is guaranteed to recover the exact solution of the $\ell_1$-minimization problem.
In topology, a homotopy describes a continuous deformation between two objects. The homotopy algorithm adopts this idea: it starts from a simple, easily computable solution and gradually deforms it---via a sequence of iterative updates---until it reaches the desired (and typically more complex) solution.
Thus, a crucial aspect of any homotopy method is the choice of an appropriate initial solution that is both simple and connected to the target solution through a well-behaved path.

\begin{subequations}
Specifically, the homotopy method  solves the $\ell_1$-minimization problem~\eqref{opt:p1} (p.~\pageref{opt:p1}):
\begin{equation}\label{opt:p1_homo}
(\text{P}_1)\qquad 
\min_{\bbeta\in\real^p} \normone{\bbeta} \quad \text{s.t.} \quad\bX\bbeta=\by,
\end{equation}
where $\bX \in \real^{n\times p}$ and $\by \in \real^n$. 
To solve this, for ${\lambda} > 0$, consider the Lagrangian form known as the LASSO (see \eqref{opt:ll}, p.~\pageref{opt:ll}):
\begin{equation}\label{opt:ll_homo}
F_{\lambda}(\bbeta) = \frac{1}{2}\normtwo{\bX\bbeta-\by}^2 + {\lambda} \normone{\bbeta}, \quad \bbeta \in \real^p,
\end{equation}
and let $\bbeta_{\lambda}$ denote a minimizer of $F_{\lambda}$. 
When ${\lambda} = {\lambda}_{\max}$ is sufficiently large, we have $\bbeta_{{\lambda}_{\max}} = \bzero$; see Section~\ref{section:optcd_el1_reg}. 
Moreover, as $\lambda\rightarrow 0$,  $ \bbeta_{\lambda} $ converges to $\bbeta^*$, where $\bbeta^*$ is a minimizer of \eqref{opt:p1}. 
A precise statement is given in the following theorem.
\end{subequations}

\begin{theoremHigh}[Recovery of $\ell_1$-minimization from LASSO]\label{theorem:recov_ell1_from_ll}
Assume that the system $\bX\bbeta=\by$ is consistent. 
If the minimizer $\bbeta^*$ of \eqref{opt:p1} is unique, then
$$
\lim_{{\lambda} \to 0} \bbeta_{\lambda} = \bbeta^*.
$$
More generally, if the minimizer of \eqref{opt:p1_homo} is not unique, then the set $\{\bbeta_{\lambda}\}_{\lambda>0}$ is bounded, and every accumulation point of $\bbeta_{\lambda}$ as $\lambda \rightarrow 0$ is a minimizer of \eqref{opt:p1}.
\end{theoremHigh}

\begin{proof}[of Theorem~\ref{theorem:recov_ell1_from_ll}]
To establish boundedness, observe that $F_{\lambda}(\bbeta_{\lambda}) \leq F_{\lambda}(\bzero) = \normtwo{\by}^2/2$. 
Let  $\{\lambda_i\}_{i \in \naturalset} \subset (0,\infty)$ be a sequence decreasing monotonically to $0$, and we write $\bbeta^i = \bbeta_{\lambda_i}$ for simplicity. 
Let $\bbeta^*$ be any minimizer of the $\ell_1$-minimization problem \eqref{opt:p1}. 
Then, since $\bX\bbeta^* = \by$ by \eqref{opt:p1}, we have 
\begin{equation}\label{equation:recov_ell1_from_ll_pv1}
\normone{\bbeta^i} \leq \frac{1}{\lambda_i} F_{\lambda_i}(\bbeta^i) \leq \frac{1}{\lambda_i} F_{\lambda_i}(\bbeta^*) = \normone{\bbeta^*}.
\end{equation}
Thus, all $\bbeta^i$ lie in the compact set
$$
\sS \triangleq \{\bbeta \in \real^p, \normone{\bbeta} \leq \normone{\bbeta^*}\}.
$$
Define  $F: \sS \to \real$ by $ F(\bbeta) = \normtwo{\bX\bbeta-\by}^2/2$. 
Since $\sS$ is compact, $F$ is continuous and attains its minimum on $\sS$ by the Weierstrass theorem (see Theorem~\ref{theorem:weierstrass_them}).
Moreover, denoting by $\widetilde{F}_{\lambda_i}$ the functions $F_{\lambda_i}$ restricted to $\sS$, we observe that the sequence $\widetilde{F}_{\lambda_i}$ is monotonically decreasing in $i$ that converges pointwise to $F$. 
By {Problem~\ref{prob:convergence_minimizers}}, any accumulation point of the sequence $\{\bbeta^i\}_i$ is a minimizer $\widehatbbeta$ of $F$, so $\bX\widehatbbeta = \by$. 
By the optimality of $\bbeta^*$ for the $\ell_1$-minimization problem, we have $\normonebig{\widehatbbeta} \geq \normone{\bbeta^*}$. 
While by definition of the set $\sS$ and \eqref{equation:recov_ell1_from_ll_pv1}, we also have $\normonebig{\widehatbbeta} = \lim_{i \to \infty} \normone{\bbeta^i} \leq \normone{\bbeta^*}$. 
Therefore,
$$
\normonebig{\widehatbbeta} = \normone{\bbeta^*},
$$
which implies that any accumulation point of $\{\bbeta^i\}_i$ is a minimizer of \eqref{opt:p1}. If the minimizer is unique, every subsequence of  $\{\bbeta^i\}_i$ converges to $\bbeta^*$, so the full sequence converges to $\bbeta^*$.
\end{proof}

\begin{remark}
Note that we use a converse statement to compute the Lagrangian LASSO solution using the penalty function method in Algorithm~\ref{alg:lass_penal}.
\end{remark}

Hence, $\bbeta_{\lambda_{\max}} = \bzero$  serves as the natural starting point for the homotopy algorithm when solving the $\ell_1$-minimization problem \eqref{opt:p1}.
The core idea of the homotopy method is to trace the solution path $\lambda \rightarrow \bbeta_{\lambda}$ from $\bbeta_{{\lambda}_{\max}} = \bzero$ down to $\bbeta^*$. As we will show below, this solution path is piecewise linear, and it suffices to compute only the breakpoints---the endpoints of the linear segments \citep{osborne2000new, efron2004least, foucart2013invitation}.

By the optimality condition for unconstrained convex problems (Theorem~\ref{theorem:fetmat_opt}), the minimizer of \eqref{opt:ll_homo} can be characterized using the subdifferential (Definition~\ref{definition:subgrad}). 
The subdifferential of $F_{\lambda}$ is given by
$$
\partial F_{\lambda}(\bbeta) = \bX^\top(\bX\bbeta - \by) + {\lambda} \partial \normone{\bbeta},
$$
where the subdifferential of the $\ell_1$-norm (Exercise~\ref{exercise:sub_norms}) is given by
$$
\partial \normone{\bbeta} = \{\bu \in \real^p \mid u_i \in \partial \abs{\beta_i}, i \in \{1,2,\ldots,p\}\}, 
\quad \text{where }
\partial \abs{x} = 
\begin{cases}
\{\sgn(x)\}, & \text{if } x \neq 0, \\
[-1, 1], & \text{if } x = 0.
\end{cases}
$$
A vector $\bbeta$ minimizes  $F_{\lambda}$ if and only if $\bzero \in \partial F_{\lambda}(\bbeta)$, see {Theorem~\ref{theorem:fetmat_opt}}. 
Using the expression above, this condition is equivalent to the following pair of relations for each coordinate $i\in\{1,2,\ldots,p\}$:
\begin{subequations}
\begin{align}
(\bX^\top(\bX\bbeta - \by))_i &= -{\lambda} \sgn(\beta_i), \quad &\text{if }& \beta_i \neq 0, \label{equ:ll_subd_neq0} \\
\abs{(\bX^\top(\bX\bbeta - \by))_i} &\leq {\lambda}, \quad &\text{if }& \beta_i = 0, \quad  \text{for }i \in \{1,2,\ldots,p\}. \label{equ:ll_subd_eq0}
\end{align}
\end{subequations}
The homotopy method begins with the initial solution $ \bbeta_{{{\lambda}}} = \bbeta^\topzero  = \bzero$. By condition \eqref{equ:ll_subd_eq0}, 
According to condition \eqref{equ:ll_subd_eq0}, the corresponding regularization parameter is chosen as ${\lambda} = {\lambda}^\topzero = \norminf{\bX^\top\by}$; see also Section~\ref{section:optcd_el1_reg} \eqref{equation:safe_lambdamax}.

In subsequent steps $t = 1, 2, \ldots$, the algorithm decreases ${\lambda}$, computes the corresponding minimizers  $\bbeta^\topone, \bbeta^\toptwo, \ldots$, and maintains an active (support) set $\sT_t$ at each step $t$. 
Let the residual vector at step $t$ be denoted by
$$
\br^\toptzero \triangleq \bX^\top(\bX\bbeta^\toptminus - \by).
$$

\paragraph{Step 1 of the homotopy method.} 
For step $1$, since $\bbeta^\topzero=\bzero$, define
\begin{equation}
i^\topone \triangleq \argmax_{i \in \{1,2,\ldots,p\}} \abs{(\bX^\top\by)_i} = \argmax_{i \in \{1,2,\ldots,p\}} \abs{r^\topone_i}.
\end{equation}
We assume throughout that this maximum is attained at a unique index $i$. (The case where the maximum is attained simultaneously at multiple indices---though theoretically possible---is extremely rare in practice and introduces technical complications that we omit here.)

Set the initial active set as $\sT_1 \triangleq \{i^\topone\}$. 
We now define the direction vector $\bd^\topone \in \real^p$ that describes the initial segment of the solution path:
$$
d^\topone_{i^\topone} \triangleq \normtwo{\bx_{i^\topone}}^{-2} \sgn((\bX^\top\by)_{i^\topone})
\qquad \text{and} \qquad 
d^\topone_i = 0, \quad \text{for } i \neq i^\topone.
$$
The first linear segment of the solution path is then
\begin{equation}\label{equation:homo_upp_ste1}
\bbeta = \bbeta(\eta) = \bbeta^\topzero + \eta \bd^\topone, \quad \eta \in [0, \eta^\topone],
\end{equation}
where  the upper bound $\eta^\topone$ will be determined shortly.
One can verify directly from the definition of $\bd^\topone$ that  condition \eqref{equ:ll_subd_neq0} holds for all $\bbeta = \bbeta(\eta)$ when we update the parameter  ${\lambda} = {\lambda}(\eta) = {\lambda}^\topzero - \eta$, $\eta \in [0, {\lambda}^\topzero]$, i.e.,:
$$
\big(\bX^\top(\bX\bbeta - \by)\big)_{i^\topone} 
= -{\lambda} \sgn(\beta_{i^\topone} ).
$$

The next breakpoint occurs at the largest  $\eta = \eta^\topone > 0$ such that condition \eqref{equ:ll_subd_eq0} continues to hold for all inactive coordinates $i\notin \sT_1$. 
Let $(a)_+ = \max\{a, 0\}$. 
This determines the upper bound $\eta^{(1)}$ in \eqref{equation:homo_upp_ste1}:
\begin{equation}\label{equation:homo_gammtopone}
\eta^\topone = \min_{i \notin \{i^\topone\}} \left\{ \left(\frac{{\lambda}^\topzero + r^\topone_i}{1 - (\bX^\top \bX \bd^\topone)_i}\right)_+, \left(\frac{{\lambda}^\topzero - r^\topone_i}{1 + (\bX^\top \bX \bd^\topone)_i}\right)_+ \right\}.
\end{equation}
Then the next solution is $\bbeta^\topone = \bbeta(\eta^\topone) = \eta^\topone \bd^\topone$, which minimizes $F_{\lambda}$ for ${\lambda} = {\lambda}^\topone \triangleq {\lambda}^\topzero - \eta^\topone$. 
By construction, this parameter ${\lambda}^\topone$ satisfies ${\lambda}^\topone = \norminf{\br^\toptwo}$. Let $i^\toptwo$ denote the index where the minimum in \eqref{equation:homo_gammtopone} is attained (where we again assume that the minimum is attained only at one index), and update the active set to $\sT_2 \triangleq \{i^\topone, i^\toptwo\}$.
This finish the step 1.

\paragraph{Step $t$ of the homotopy method.} 
At step $t$, we receive an active set $\sT_t$ from step $t-1$. 
The new direction  $\bd^\toptzero$ of the homotopy path is determined by solving the linear system
\begin{equation}\label{equation:homo_norm_equ}
\bX_{\sT_t}^\top \bX_{\sT_t} \bd_{\sT_t}^\toptzero = \sgn(\br_{\sT_t}^\toptzero).
\end{equation}
where $\bX_{\sT_t}$ denotes the submatrix of $\bX$ consisting of columns indexed by $\sT_t$, 
and $\br_{\sT_t}^\toptzero$ is the corresponding subvector of the residual.
This is a linear system  of size $\abs{\sT_t} \times \abs{\sT_t}$, with $\abs{\sT_t} \leq t$. 
For indices outside the active set, we set $d^\toptzero_i = 0$, for all  $i \notin \sT_t$. 
The next linear segment of the solution path is
\begin{equation}
\bbeta(\eta) = \bbeta^\toptminus + \eta \bd^\toptzero, \quad \eta \in [0, \eta^\toptzero].
\end{equation}
The maximal stepsize $\eta$ such that condition \eqref{equ:ll_subd_eq0} remains satisfied for all inactive coordinates is
\begin{equation}\label{equation:homo_stept_gamma}
\eta_+^\toptzero = \min_{i \notin \sT_t} \left\{ \left(\frac{{\lambda}^\toptminus + r^\toptzero_i}{1 - (\bX^\top \bX \bd^\toptzero)_i}\right)_+, \left(\frac{{\lambda}^\toptminus - r^\toptzero_i}{1 + (\bX^\top \bX \bd^\toptzero)_i}\right)_+ \right\}.
\end{equation}
The maximal stepsize such that condition \eqref{equ:ll_subd_neq0} continues to hold for all active coordinates (i.e., no component of $\bbeta$ changes sign or crosses zero prematurely) is
\begin{equation}\label{equation:homo_stept_gamma_neg}
\eta_-^\toptzero = \min_{i \in \sT_t} \left\{ \left(-\beta^\toptminus_i / d^\toptzero_i\right)_+ \right\}.
\end{equation}

The next breakpoint occurs at $\eta^\toptzero = \min\{\eta_+^\toptzero, \eta_-^\toptzero\}$ 
and the updated solution is $\bbeta^\toptzero = \bbeta(\eta^\toptzero)$:
\begin{itemize}
\item If $\eta^\toptzero=\eta_+^\toptzero$, then the minimizing index $i^\toptzero \notin \sT_t$ (from \eqref{equation:homo_stept_gamma}) is added to the active set: $\sT_{j+1} = \sT_t \cup \{i^\toptzero\}$.
\item  If $\eta^\toptzero = \eta_-^\toptzero$, then the corresponding index $i_-^\toptzero \in \sT_t$ is removed: $\sT_{j+1} = \sT_t \setminus \{i_-^\toptzero\}$.
\end{itemize}
This is just like the LASSO modification for LARS; see Remark~\ref{remark:lass_modi}.
Finally, update the regularization parameter as ${\lambda}^\toptzero = {\lambda}^\toptminus - \eta^\toptzero$. 
By construction, this satisfies ${\lambda}^\toptzero = \norminf{\br^\toptone}$.

The algorithm terminates when ${\lambda}^\toptzero = \norminf{\br^\toptone} = 0$, i.e., when the residual vanishes. At this point, the current iterate $\bbeta^\toptzero$ is the exact solution $\bbeta^*$ of the original $\ell_1$-minimization problem \eqref{opt:p1_homo}, and the algorithm outputs $\bbeta^*=\bbeta^\toptzero$.


\begin{theoremHigh}[Exact recovery of homotopy]\label{theorem:hm_conv}
Assume that the minimizer $\bbeta^*$ of the $\ell_1$-minimization problem \eqref{opt:p1} is unique. If, at each step of the algorithm, the minima in \eqref{equation:homo_stept_gamma} and \eqref{equation:homo_stept_gamma_neg} are attained at a single index $i$, then the homotopy algorithm described above terminates with output $\bbeta^*$.
\end{theoremHigh}
\begin{proof}[of Theorem~\ref{theorem:hm_conv}]
Following the description of the algorithm, it remains only to show that the algorithm terminates after finitely many steps. 
To this end, observe that the sign pattern $\sgn(\bbeta_{{\lambda}^{(a)}})$ is distinct at each iteration $a$. 
Indeed, suppose that for two parameters ${\lambda}^{(a)}$ and ${\lambda}^{(b)}$ with  $b > a$,  we had  $\sgn(\bbeta_{{\lambda}^{(a)}})_i = \sgn(\bbeta_{{\lambda}^{(b)}})_i$.
Let $0 \neq \sigma_i  \triangleq \sgn(\bbeta_{{\lambda}^{(a)}})_i = \sgn(\bbeta_{{\lambda}^{(b)}})_i $ for all indices $i$ in the support.
Then, by condition \eqref{equ:ll_subd_neq0}, we would have
$$
(\bX^\top \bX (\bbeta_{{\lambda}^{(a)}} - \bbeta_{{\lambda}^{(b)}}))_i = ({\lambda}^{(a)} - {\lambda}^{(b)}) \sigma_i,
$$
This implies that $\bbeta_{{\lambda}^{(b)}}$ lies on the same linear segment of the homotopy path that starts at $\bbeta_{{\lambda}^{(a)}}$. However, by construction, each $\bbeta_{{\lambda}^{(a)}}$ corresponds to a breakpoint---that is, an endpoint of a linear piece---and not an interior point of any segment.

Since the number of possible sign patterns in $\real^p$ is finite (at most $3^p$, accounting for positive, negative, and zero entries), the algorithm cannot generate an infinite sequence of distinct breakpoints. Therefore, it must terminate after finitely many steps.

At termination, we have $\lambda^\toptzero=0$, which implies $\bbeta_{{\lambda}^{(t)}}=\bbeta^*$ by Theorem~\ref{theorem:recov_ell1_from_ll}. Hence, the algorithm outputs the unique minimizer $\bbeta^*$.
\end{proof}

If the algorithm is terminated early at some iteration $t$, it clearly returns the minimizer of $F_{\lambda} = F_{{\lambda}^\toptzero}$, i.e., the solution to the Lagrangian LASSO problem with parameter $\lambda=\lambda^\toptzero$. 
Moreover, simple stopping criteria can be applied during the algorithm’s execution to solve the following constrained 
\begin{align}
\min \normone{\bbeta} &\quad \text{s.t.} \quad\normtwo{\bX\bbeta-\by} \leq \epsilon \\
\text{or} \quad  \normtwo{\bX\bbeta-\by} &\quad \text{s.t.} \quad\normone{\bbeta} \leq \Sigma.
\end{align}
The first formulation corresponds to $\ell_1$-minimization with noisy measurements, while the second is the constrained LASSO.

\subsection{Iteratively Reweighted Least Squares (IRLS)\index{Iteratively reweighted least squares}}

We now introduce the \textit{iteratively reweighted least squares (IRLS)} method. Although IRLS serves only as a proxy for $\ell_1$-minimization, its formulation is directly motivated by it, and under certain conditions, IRLS indeed recovers the $\ell_1$-minimizer. Assuming the stable nullspace property---a condition known to be equivalent to exact and stable sparse recovery via $\ell_1$-minimization---we will show that IRLS enjoys the same error guarantees as $\ell_1$-minimization. However, in general, the output of IRLS may differ from the true $\ell_1$-minimizer.

The algorithm originates from the elementary identity $\abs{x} = \frac{\abs{x}^2}{\abs{x}}$ for $x \neq 0$. 
This observation suggests that $\ell_1$-minimization can be reformulated as a weighted $\ell_2$-minimization. 
Specifically, let  $\bX \in \real^{n\times p}$ with $n \leq p$, and suppose $\bbeta^*$ is a minimizer of
\begin{equation}\label{opt:p1_irls}
\text{(P$_1$)}\quad
\min_{\bbeta \in \real^p} \normone{\bbeta} \quad \text{s.t.} \quad \bX\bbeta=\by.
\end{equation}
If $\beta_i^* \neq 0$ for all $i \in \{1,2,\ldots,p\}$, then $\bbeta^*$ also minimizes the weighted least squares problem:
\begin{equation}
\min_{\bbeta \in \real^p} \sum_{i=1}^p \frac{\beta_i^2 }{\absbig{\beta_i^*}} \quad \text{s.t.} \quad \bX\bbeta=\by.
\end{equation}
The advantage of this reformulation lies in the fact that minimizing the smooth quadratic function $\abs{x}^2$ is computationally more tractable than minimizing the non-smooth absolute value $\abs{x}$. However, two major drawbacks arise: first, the true minimizer $\bbeta^*$ is not known a priori; second, we cannot generally expect  $\beta_i^* \neq 0$ for all $i$, since we typically seek $k$-sparse solutions. In fact, by Problem~\ref{prob:sparsi_of_p1_n}, any unique $\ell_1$-minimizer is always at most $n$-sparse.

Despite these challenges, the above identity motivates an iterative strategy: at each step, we solve a weighted $\ell_2$-minimization problem, where the weights for the next iteration are updated based on the current solution.
Specifically, we can define $F(\bbeta, \bw) = \sum_{i}^{p}\beta_i^2w_i$ and update $\bbeta$ and $\bw$ iteratively:
\begin{align*}
\bbeta^\toptone&\leftarrow \argmin_{\bbeta\in\real^p} F(\bbeta, \bw^\toptzero);\\
\bw^\toptone&\leftarrow\argmin_{\bw\in\real^p} F(\bbeta^\toptone, \bw).
\end{align*}

For reasons of stability and convergence, a more robust ingredient---both in the formulation and analysis of the algorithm---is the following functional:
\begin{equation}\label{equation:p1_irls_form}
\mathcalF(\bbeta, \bw, \varepsilon) = \frac{1}{2} \left[ \sum_{i=1}^p \beta_i^2 w_i + \sum_{i=1}^p (\varepsilon^2 w_i + w_i^{-1}) \right],
\end{equation}
where $\bbeta \in \real^p$, $\varepsilon \geq 0$, and $\bw \in \real^p$ is a strictly positive weight vector ($w_i > 0$ for all $i \in \{1,2,\ldots,p\}$). 
The algorithm below makes use of the nonincreasing rearrangement of the current iterate $\bbeta^\toptzero\in \real^p$, denoted by
\begin{equation}\label{equation:nonincr_def}
[\bbeta^\toptzero]^\downarrow \in \real^p,
\end{equation} 
which satisfies $[\bbeta^\toptzero]^\downarrow_1 \geq [\bbeta^\toptzero]^\downarrow_2 \geq \ldots \geq [\bbeta^\toptzero]^\downarrow_p$.
The complete procedure is given in Algorithm~\ref{alg:ell1_irls}, which requires a sparsity parameter $k$.
\begin{algorithm}[H]
\caption{Iteratively Reweighted Least Squares (IRLS) \citep{daubechies2010iteratively}}
\label{alg:ell1_irls}
\begin{algorithmic}[1]
\Require    Dictionary $ \bX\in\real^{n\times p}$, measurement vector $\by \in\real^n$, $\eta > 0$, $k \in \{1,2,\ldots,p\}$;
\State \textbf{initialize:}  $\bw^\topzero = [1, 1, \ldots, 1]^\top \in \real^p$, $\varepsilon_0 = 1$;
\For{$t=1,2,\ldots$}
\State $\bbeta^\toptone \gets \arg \min_{\bbeta \in \real^p} \mathcalF(\bbeta, \bw^\toptzero, \varepsilon_t)\quad \text{s.t.} \quad \bX\bbeta = \by$; \Comment{(IRLS$_1$)}
\State $\varepsilon_{t+1} \gets \min\{\varepsilon_t, \eta [\bbeta^\toptone]^\downarrow_{k+1}\}$;  \Comment{(IRLS$_2$)}
\State $\bw^\toptone \gets \arg \min_{\bw > 0} \mathcalF(\bbeta^\toptone, \bw, \varepsilon_{t+1})$; \Comment{(IRLS$_3$)}
\State Stop if  a stopping criterion is met or $\varepsilon_{t+1}=0$; \Comment{(IRLS$_4$)}
\EndFor
\State \Return Final $\bbeta\gets \bbeta^\toptzero$.
\end{algorithmic}
\end{algorithm}

\paragrapharrow{Update for $\bbeta^\toptone$.}
Since $\bw^\toptzero$ and $\varepsilon_t$ are fixed in the minimization problem (IRLS$_1$), the second sum in the definition  \eqref{equation:p1_irls_form} of $\mathcalF$ is constant. 
Consequently, $\bbeta^\toptone$ minimizes a \textit{weighted least squares problem} of the form (see \eqref{opt:l2_weighted}, p.~\pageref{opt:l2_weighted}):
\begin{subequations}
\begin{equation}
\min_{\bbeta \in \real^p} \norm{\bbeta}_{2,\bw^\toptzero} = \left( \sum_{i=1}^p \beta_i^2 w_i^\toptzero \right)^{1/2}\quad \text{s.t.} \quad \bX\bbeta = \by.
\end{equation}
By \eqref{equation:weightls_sol1}, the minimizer $\bbeta^\toptone$ admits the explicit expression
\begin{equation}
\bbeta^\toptone = \bD_{t}^{-1/2} (\bX \bD_{t}^{-1/2})^+ \by,
\end{equation}
where $(\bX \bD_{t}^{-1/2})^+$ denotes the pseudo-inverse of $\bX \bD_{t}^{-1/2}$, and $\bD_{t} = \diag(\bw^\toptzero)$ 
is the diagonal matrix whose entries are given by the weight vector $\bw^\toptzero$. 
If $\bX$ has full row rank---which is typically the case in compressed sensing---then, due to the positivity of the weights $w^\toptzero>0$, the matrix $\bX \bD_{t}$ 
also has full row rank. 
In this setting, the pseudo-inverse simplifies, and the solution can be written as
\begin{equation}
\bbeta^\toptone = \bD_{t}^{-1} \bX^\top (\bX \bD_{t}^{-1} \bX^\top)^{-1} \by.
\end{equation}
Alternatively, defining $\bxi\triangleq  (\bX \bD_{t}^{-1} \bX^\top)^{-1} \by$, the solution may be expressed as
\begin{equation}\label{equation:irls_linsysdt}
\bbeta^\toptone = \bD_{t}^{-1} \bX^\top \bxi, \quad \text{where} \quad \bX \bD_{t}^{-1} \bX^\top \bxi = \by.
\end{equation}
Thus, computing $\bbeta^\toptone$ reduces to solving the linear system above for $\bxi$. 
For further background on least squares and weighted least squares problems, we refer the reader to \citet{lu2021rigorous}.
\end{subequations}

\paragrapharrow{Update for $\bw^\toptone$.}
The minimization of $\bw^\toptone$ in (IRLS$_3$) decouples across coordinates and can be performed explicitly:
\begin{equation}\label{equation:irls3_equiv}
w_i^\toptone = \frac{1}{\sqrt{(\beta_i^\toptone)^2 + \varepsilon_{t+1}^2}}, \quad i \in \{1,2,\ldots,p\}.
\end{equation}
This formula also illustrates the role of the smoothing parameter $\varepsilon_t$. 
In the naive formulation $w_i^\toptone = \absbig{\beta_i^\toptone}^{-1}$, the weight would become unbounded as $\beta_i^\toptone$ approaches zero. 
The introduction of $\varepsilon_{t+1}$ regularizes the weight $\bw^\toptone$; 
specifically, it ensures that 
$
\norminf{\bw^\toptone} \leq \varepsilon_{t+1}^{-2}.
$
Nevertheless, to recover the $\ell_1$-minimizer in the limit, we require that $\varepsilon_t\rightarrow 0$ as $t\rightarrow \infty$.
The update rule (IRLS$_2$) guarantees precisely this behavior: $\varepsilon_{t+1}= \min\{\varepsilon_t, \eta [\bbeta^\toptone]^\downarrow_{k+1}\}\leq \varepsilon_t$, so $\varepsilon_t$ is nonincreasing.
Moreover, if the iterates $\bbeta^\toptone$ converge to a $k$-sparse vector, then
$[\bbeta^\toptone]^\downarrow_{k+1}\rightarrow 0$, and consequently $\varepsilon_{t+1}\rightarrow 0$.
Thus, the sparsity parameter $k$ directly controls the target level of sparsity in the solution.

\begin{remark}
Alternative update rules for both $\varepsilon_{t+1}$ and the weights $\bw^\toptone$ are also possible. For instance, one may consider the truncated update
$
w_j^\toptone = \min\{\absbig{\beta_j^\toptone}^{-1}, \varepsilon_{t+1}^{-1}\}.
$
Convergence results for this variant have been established; see \citet{fornasier2011low} for details.
\end{remark}

To state the main convergence result for the IRLS algorithm, we introduce, for any $\varepsilon>0$, the auxiliary functional
\begin{equation}\label{equation:irls_auxit_epsi}
F_\varepsilon(\bbeta) \triangleq \sum_{i=1}^p \sqrt{\beta_i^2 + \varepsilon^2}.
\end{equation}
The proof of convergence relies on the following key properties of the IRLS iterates.
\begin{lemma}[Properties of the IRLS algorithm]\label{lemma:irls_bds_itera}
Let $\{\bbeta^\toptzero, \bw^\toptzero\}_{t\in \naturalset}$ be the sequence of iterates generated by the  IRLS Algorithm~\ref{alg:ell1_irls}. Then, for $t \in \naturalset$,
\begin{subequations}
\begin{equation}\label{equation:irls_bds_itera1}
\mathcalF(\bbeta^\toptzero, \bw^\toptzero, \varepsilon_t) = \sum_{i=1}^p \sqrt{(\beta_i^\toptzero)^2 + \varepsilon_t^2} = F_{\varepsilon_t}(\bbeta^\toptzero),
\end{equation}
and
\begin{align}
\mathcalF(\bbeta^\toptzero, \bw^\toptzero, \varepsilon_t) 
&\leq \mathcalF(\bbeta^\toptzero, \bw^\toptminus, \varepsilon_t) \leq \mathcalF(\bbeta^\toptzero, \bw^\toptminus, \varepsilon_{t-1}) \label{equation:irls_bds_itera2} \\
&\leq \mathcalF(\bbeta^\toptminus, \bw^\toptminus, \varepsilon_{t-1}). \label{equation:irls_bds_itera3}
\end{align}
Moreover, the sequence $\{\bbeta^\toptzero\}$ is bounded,
\begin{equation}\label{equation:irls_bds_itera4}
\normonebig{\bbeta^\toptzero} \leq \mathcalF(\bbeta^\topone, \bw^\topzero, \varepsilon_0) \triangleq B, \quad t \in \naturalset,
\end{equation}
and the weights $\bw^\toptzero$ are uniformly bounded away from zero:
\begin{equation}\label{equation:irls_bds_itera5}
w_i^\toptzero \geq B^{-1}, \quad i \in \{1,2,\ldots,p\}, \; t \in \naturalset.
\end{equation}
The iterates of the IRLS Algorithm~\ref{alg:ell1_irls} also satisfy
\begin{equation}\label{equation:irls_bds_itera6}
\sum_{t=1}^\infty \normtwo{\bbeta^\toptone - \bbeta^\toptzero}^2 \leq 2B^2,
\end{equation}
where $B$ is the constant defined in \eqref{equation:irls_bds_itera4}. Consequently, $\lim_{t \to \infty} (\bbeta^\toptone - \bbeta^\toptzero) = \bzero$.
\end{subequations}
\end{lemma}
\begin{proof}[of Lemma~\ref{lemma:irls_bds_itera}]
Equation \eqref{equation:irls_bds_itera1} follows directly from the definition of $\mathcalF$ in \eqref{equation:p1_irls_form} and the explicit update \eqref{equation:irls3_equiv} for the weights. 

The first inequality in \eqref{equation:irls_bds_itera2} holds because $\bw^\toptzero$ minimizes $\mathcalF(\bbeta^\toptzero,\cdot, \varepsilon_t)$ by construction  (IRLS$_3$).
The second inequality follows from the monotonicity $\varepsilon_{t+1} \leq \varepsilon_t$ and the positivity of  $\bw^\toptminus$. 
The final inequality  \eqref{equation:irls_bds_itera3} is a consequence of the optimality of $\bbeta^\toptzero$ in the weighted least squares subproblem  (IRLS$_1$).

From \eqref{equation:irls_bds_itera1} and the monotonicity in \eqref{equation:irls_bds_itera2}, we obtain
$$
\normonebig{\bbeta^\toptzero} \leq \sum_{i=1}^p \sqrt{(\beta_i^\toptzero)^2 + \varepsilon_t^2} = F_{\varepsilon_t}(\bbeta^\toptzero) = \mathcalF(\bbeta^\toptzero, \bw^\toptzero, \varepsilon_t) \leq \mathcalF(\bbeta^\topone, \bw^\topzero, \varepsilon_0) = B,
$$
which establishes \eqref{equation:irls_bds_itera4}. 
Using \eqref{equation:irls3_equiv}, we then have
$
(w_i^\toptzero)^{-1} = \sqrt{(\beta_i^\toptzero)^2 + \varepsilon_t^2} \leq B, \; i \in \{1,2,\ldots,p\},
$
proving \eqref{equation:irls_bds_itera5}.

To prove \eqref{equation:irls_bds_itera6}, observe that the monotonicity in \eqref{equation:irls_bds_itera2} implies
\begin{align*}
2 &\left( \mathcalF(\bbeta^\toptzero, \bw^\toptzero, \varepsilon_t) - \mathcalF(\bbeta^\toptone, \bw^\toptone, \varepsilon_{t+1}) \right) 
\geq 2 \left( \mathcalF(\bbeta^\toptzero, \bw^\toptzero, \varepsilon_t) - \mathcalF(\bbeta^\toptone, \bw^\toptzero, \varepsilon_t) \right) \\
&= \sum_{i=1}^p \left( (\beta_i^\toptzero)^2 - (\beta_i^\toptone)^2 \right) w_i^\toptzero 
= \innerproduct{\bbeta^\toptzero + \bbeta^\toptone, \bbeta^\toptzero - \bbeta^\toptone}_{\bw^\toptzero}\\
&\stackrel{\dag}{=} \innerproduct{\bbeta^\toptzero + \bbeta^\toptone, \bbeta^\toptzero - \bbeta^\toptone}_{\bw^\toptzero}  
- 2 \innerproduct{\bbeta^\toptone, \bbeta^\toptzero - \bbeta^\toptone}_{\bw^\toptzero} \\
\end{align*}
\begin{align*}
&=  \innerproduct{\bbeta^\toptzero - \bbeta^\toptone, \bbeta^\toptzero - \bbeta^\toptone }_{\bw^\toptzero}  
=\norm{\bbeta^\toptzero - \bbeta^\toptone}_{\bw^\toptzero}^2 = \sum_{i=1}^p w_i^\toptzero \abs{\beta_i^\toptzero - \beta_i^\toptone}^2 \\
&\geq B^{-1} \normtwobig{\bbeta^\toptzero - \bbeta^\toptone}^2,
\end{align*}
where  $\innerproductbig{\ba, \bb}_{\bw} = \sum_{i=1}^p a_i b_i w_i$ denotes the weighted inner product.
Since both $\bbeta^\toptzero$ and $\bbeta^\toptone$ satisfy $\bX\bbeta^\toptzero = \by = \bX\bbeta^\toptone$ by  (IRLS$_1$), we have $\bbeta^\toptzero - \bbeta^\toptone \in \nspace(\bX)$. 
Equality ($\dag$)  follows from the orthogonality property of weighted least squares solutions (Theorem~\ref{theorem:orthogonal_weight_lsii}): $ \innerproductbig{\bbeta^\toptone, \bbeta^\toptzero - \bbeta^\toptone}_{\bw^\toptzero} = 0$. 
And the last inequality follows from \eqref{equation:irls_bds_itera5}. Telescoping the sum of these inequalities over $t$ shows that
\begin{align*}
\sum_{t=1}^\infty \normtwobig{\bbeta^\toptzero - \bbeta^\toptone}^2 
&\leq 2B \sum_{t=1}^\infty \left( \mathcalF(\bbeta^\toptzero, \bw^\toptzero, \varepsilon_t) - \mathcalF(\bbeta^\toptone, \bw^\toptone, \varepsilon_{t+1}) \right) \\
&\leq 2B \mathcalF(\bbeta^\topone, \bw^\topone, \varepsilon_1) \leq 2B^2,
\end{align*}
This completes the proof.
\end{proof}

Note that \eqref{equation:irls_bds_itera3} shows that the objective value $\mathcalF$ decreases at each iteration, while \eqref{equation:irls_bds_itera6} guarantees that successive iterates become asymptotically stationary.

We further require a characterization of the minimizer  $\bbeta^{(\varepsilon)}$ of the smoothed objective $F_\varepsilon(\bbeta) = \sum_{i=1}^p \sqrt{\beta_i^2 + \varepsilon^2}$, defined in \eqref{equation:irls_auxit_epsi},
for the constrained optimization problem
\begin{equation}\label{equation:irls_auxit_epsi2}
\bbeta^{(\varepsilon)}\triangleq \arg\min_{\bbeta \in \real^p} F_\varepsilon(\bbeta)\quad \text{s.t.} \quad \bX\bbeta = \by.
\end{equation}
The minimizer $\bbeta^{(\varepsilon)}$  is \textbf{unique} due to the strict convexity of $F_\varepsilon$; see Problem~\ref{prob:irls_auxit_epsi2}.

\begin{lemma}[Optimality condition of minimization for $F_\varepsilon(\bbeta)$]\label{lemma:irls_lem3}
Let $\varepsilon > 0$, and let  $\bbeta \in \real^p$ satisfy $\bX\bbeta = \by$. Then $\bbeta = \bbeta^{(\varepsilon)}$ if and only if $\innerproduct{\bbeta, \bn}_{\bw_{\bbeta,\varepsilon}} = 0$ for all $\bn \in \nspace(\bX)$, where the weight vector $\bw_{\bbeta,\varepsilon}\in\real^p$ is defined componentwise by $(\bw_{\bbeta,\varepsilon})_i = (\beta_i^2 + \varepsilon^2)^{-1/2}$, $i=1,2,\ldots,p$.
\end{lemma}
\begin{proof}[of Lemma~\ref{lemma:irls_lem3}]
Suppose that $\bbeta = \bbeta^{(\varepsilon)}$ is the minimizer of \eqref{equation:irls_auxit_epsi2}. Let $\bn \in \nspace(\bX)$ be any vector in the nullspace of $\bX$, and define the function
$$g(\nu) = F_{\varepsilon}(\bbeta + \nu\bn) - F_{\varepsilon}(\bbeta), \quad \nu \in \real.
$$
Since $\bX(\bbeta + \nu\bn) = \by$ for all $\nu \in \real$,  the perturbed point $\bbeta+\nu\bn$ remains feasible.
By optimality of $\bbeta$, we have $g(\nu) \geq 0$ for all $\nu \in \real$, and clearly $g(0) = 0$. 
Hence, $g$ attains its minimum at $\nu=0$, so  $g'(0) = 0$. 
A direct computation yields
$$
g'(0) = \sum_{i=1}^p \frac{\beta_i n_i}{\sqrt{\beta_i^2 + \varepsilon^2}} = \innerproduct{\bbeta, \bn }_{\bw_{\bbeta,\varepsilon}} , 
$$
which implies $\innerproduct{ \bbeta, \bn }_{\bw_{\bbeta,\varepsilon}} = 0$ for all $\bn \in \nspace(\bX)$.

Conversely, assume that $\bbeta$ satisfies $\bX\bbeta = \by$ and $\innerproduct{\bbeta, \bn}_{\bw_{\bbeta,\varepsilon}} = 0$ for all $\bn \in \nspace(\bX)$. 
Consider the scalar function $f(u) \triangleq \sqrt{u^2 + \varepsilon^2}$, which is convex on $\real$. 
By the gradient inequality (Theorem~\ref{theorem:conv_gradient_ineq}), for any  $u, u_0 \in \real$, we have
$
\sqrt{u^2 + \varepsilon^2} \geq \sqrt{u_0^2 + \varepsilon^2} + \frac{u_0 {(u - u_0)}}{\sqrt{u_0^2 + \varepsilon^2}}. 
$
Therefore, for any $\bn \in \nspace(\bX)$, we have
$$ 
F_{\varepsilon}(\bbeta + \bn) \geq F_{\varepsilon}(\bbeta) + \sum_{i=1}^p \frac{\beta_i n_i}{\sqrt{\beta_i^2 + \varepsilon^2}} = F_{\varepsilon}(\bbeta) + \innerproduct{\bbeta, \bn}_{\bw_{\bbeta,\varepsilon}} = F_{\varepsilon}(\bbeta). 
$$
Since $\bn \in \nspace(\bX)$ is arbitrary, it follows that $\bbeta = \bbeta^{(\varepsilon)}$ is a minimizer of \eqref{equation:irls_auxit_epsi2}.
\end{proof}

We are now ready to establish the convergence of the IRLS algorithm.
The recovery theorem below is based on the notion of the stable NSP (Definition~\ref{definition:stable_nsp}) and closely mirrors analogous results for $\ell_1$-minimization. 
Recall that in {Sections~\ref{section:spar_gauss} and \ref{section:desi_ssgsg}}, it is shown directly that Gaussian and sub-Gaussian random matrices satisfy the RIP with high probability under suitable conditions. 
Moreover, by  {Theorem~\ref{theorem:rip2rnsp},  the RIP implies the NSP. Consequently, the various matrix ensembles discussed in this book also satisfy the stable NSP under appropriate assumptions.
\begin{theoremHigh}[Convergence of IRLS under stable NSP\index{Stable NSP}]\label{theorem:recv_irls_nsp}
Let $\bX \in \real^{n\times p}$ be a design matrix satisfying the stable NSP (Definition~\ref{definition:stable_nsp}) of order $k$ with parameter $\rho < 1$. 
Let $\bbeta^* \in \real^p$ and consider the measurements $\by = \bX\bbeta^*$. 
Consider the IRLS Algorithm~\ref{alg:ell1_irls} with sparsity parameters $k$ and smoothing parameter update rule using $\eta = 1/p$.
\footnote{The choice $\eta = 1/p$ is made to facilitate the theoretical analysis and ensure the validity of the  convergence result. In practice, however, larger values---such as $\eta = 1$---often yield faster convergence and better empirical performance. A rigorous convergence guarantee for such alternative choices is currently not available.}
Then the sequence of iterates $\{\bbeta^\toptzero\}_t$ converges to a vector $\widehatbbeta \in \real^p$ as $t \to \infty$, and the nonincreasing rearrangement $[\widehatbbeta]^\downarrow$ satisfies $[\widehatbbeta]^\downarrow_k = p \lim_{t \to \infty} \varepsilon_t$. 
Moreover, the following statements hold:
\begin{enumerate}[(i)]
\item If $\lim_{t \to \infty} \varepsilon_t = 0$, then $\widehatbbeta$ is $k$-sparse and solves the $\ell_1$-minimization problem \eqref{opt:p1}. 
Furthermore, if the true vector $\bbeta^*$ is itself $k$-sparse, then $\bbeta^* = \widehatbbeta$. More generally, we have the error bound
\begin{equation}\label{equation:recv_irls_nsp1}
\normone{\bbeta^* - \widehatbbeta} \leq \frac{2(1+\rho)}{1-\rho} \sigma_k(\bbeta^*)_1,
\end{equation}
where $\sigma_k(\bbeta^*)_1$ denotes the best $k$-term $\ell_1$-approximation error of $\bbeta^*$; see \eqref{equation:lsdist_spar}.

\item If $\varepsilon \triangleq \lim_{t \to \infty} \varepsilon_t > 0$, then $\widehatbbeta = \bbeta^{(\varepsilon)}$, i.e., the limit coincides with the minimizer of the smoothed functional $F_\varepsilon$ subject to $\bX \bbeta = \by$. 
In this case, if $\rho$ satisfies the tighter bound $\rho < 1 - \frac{2}{k+2}$, then for any $\widetildek < k - \frac{2\rho}{1-\rho}$,  it holds that 
\begin{equation}\label{equation:recv_irls_nsp2}
\normone{\bbeta^* - \widehatbbeta}  \leq C \sigma_{\widetildek}(\bbeta^*)_1 \quad \text{with} \quad C \triangleq \frac{2(1+\rho)}{1-\rho} \cdot \frac{k - \widetildek + 3/2}{k - \widetildek - \frac{2\rho}{1-\rho}}.
\end{equation}
In particular, if $\bbeta^*$ is $\widetildek$-sparse for some $\widetildek < k - \frac{2\rho}{1 - \rho}$, then $\lim_{t \to \infty} \varepsilon_t = 0$, so the case $\varepsilon > 0$ cannot occur.
\end{enumerate}
\end{theoremHigh}

\begin{proof}[of Theorem~\ref{theorem:recv_irls_nsp}]
Note that since  $0 \leq \varepsilon_{t+1} \leq \varepsilon_{t}$, the sequence $\{\varepsilon_t\}_{t \in \naturalset}$ always converges. 
Denote its limit by $\varepsilon \triangleq \lim_{t \to \infty} \varepsilon_t$.

\paragraph{(i) When $\lim_{t \to \infty} \varepsilon_t = 0$.} 
\textbf{Case (a).}
First, suppose that $\varepsilon_{t_0} = 0$ for some $t_0 \in \naturalset$. Then the algorithm terminates  by definition (the next (IRLS$_3$) step is not performed; see (IRLS$_4$)), and we can set $\bbeta^\toptzero = \bbeta^{(t_0)}$ for $t \geq t_0$ so that $\lim_{t \to \infty} \bbeta^\toptzero = \bbeta^{(t_0)} = \widehatbbeta$. By the definition of $\varepsilon_t$ in (IRLS$_2$), it follows that the nonincreasing rearrangement $[\bbeta^{(t_0)}]_{k+1}^\downarrow = 0$, which implies that $\widehatbbeta = \bbeta^{(t_0)}$ is $k$-sparse. 
Since $\bX$ satisfies the stable NSP of order $k$ (and hence the standard NSP of order $k$; see Definitions~\ref{definition:nullspace_prop} and \ref{definition:stable_nsp}),  Corollary~\ref{corollary:exa_ell1_snsp} guarantees  that $\widehatbbeta$ is the unique solution to the $\ell_1$-minimization problem~\eqref{opt:p1}.
Therefore, if the true vector $\bbeta^*$---where $\by=\bX\bbeta^*$---is $k$-sparse, then  $\bbeta$ is also the unique $\ell_1$-minimizer, i.e., $\bbeta^* = \widehatbbeta$. 
For a general $\bbeta^* \in \real^p$ satisfying $\by=\bX\bbeta^*$, not necessarily being $k$-sparse, Theorem~\ref{theorem:stablensp} yields the error bound
$$
\normone{\bbeta^* - \widehatbbeta} \leq \frac{2(1+\rho)}{1-\rho}\sigma_k(\bbeta^*)_1.
$$

\paragraph{Case (b).}
Now assume that $\varepsilon_t > 0$ for all $t \in \naturalset$, but still $\lim_{t \to \infty} \varepsilon_t = 0$. 
Then there exists an increasing sequence of indices $\{t_i\}$ such that $\varepsilon_{t_i} < \varepsilon_{t_{i-1}}$ for all $i \in \naturalset$. By the definition  of $\varepsilon_t$ in (IRLS$_2$) and the choice $\eta = 1/p$, this implies that the nonincreasing rearrangement of $\bbeta^{(t_i)}$ satisfies
$$
[\bbeta^{(t_i)}]_{k+1}^\downarrow < p\varepsilon_{t_{i-1}}, \quad i \in \naturalset.
$$
From property \eqref{equation:irls_bds_itera4}, the sequence $\{\bbeta^\toptzero\}$ is bounded. Hence, there exists a subsequence $\{t_{i_\ell}\}$ of $\{t_i\}$ such that $\{\bbeta^{(t_{i_\ell})}\}_\ell$ converges to some $\widehatbbeta$ satisfying $\bX\widehatbbeta = \by$.
Since $\norminfbig{[\widehatbbeta]^\downarrow - [\bbeta^{(t_{i_\ell})}]^\downarrow } \leq \norminfbig{\widehatbbeta - \bbeta^{(t_{i_\ell})} } $~\footnote{See Problem~\ref{prob:bound_nonincr}.}, it follows that $[\bbeta^{(t_{i_\ell})}]^\downarrow$ converges to $[\widehatbbeta]^\downarrow$. 
Therefore, we have
$$
[\widehatbbeta]_{k+1}^\downarrow = \lim_{\ell \to \infty} [\bbeta^{(t_{i_\ell})}]_{k+1}^\downarrow \leq \lim_{i \to \infty} p\varepsilon_{t_{i_\ell}} = 0,
$$
which shows that $\widehatbbeta$ is $k$-sparse. 
As before, the NSP of order $k$  implies that $\widehatbbeta$ is the unique $\ell_1$-minimizer. 
It remains to show that the full sequence $\{\bbeta^\toptzero\}$ converges to $\widehatbbeta$. Since $\bbeta^{(t_{i_\ell})} \to \widehatbbeta$ and $\varepsilon_{t_{i_\ell}} \to 0$ as $\ell \to \infty$, the identity \eqref{equation:irls_bds_itera1} implies that
$$
\lim_{\ell \to \infty} \mathcalF(\bbeta^{(t_{i_\ell})}, \bw^{(t_{i_\ell})}, \varepsilon_{t_{i_\ell}}) = \normonebig{\widehatbbeta}.
$$
By the monotonicity properties in \eqref{equation:irls_bds_itera2} and \eqref{equation:irls_bds_itera3}, 
we conclude that $\lim_{t \to \infty} \mathcalF(\bbeta^\toptzero, \bw^\toptzero, \varepsilon_t) = \normonebig{\widehatbbeta}$. 
Moreover, from~\eqref{equation:irls_bds_itera1} and the inequality $\sqrt{(\beta_i^\toptzero)^2+\varepsilon_t} - \varepsilon_t \leq \absbig{\beta_i^\toptzero}$, we obtain
$$
\mathcalF(\bbeta^\toptzero, \bw^\toptzero, \varepsilon_t) - p\varepsilon_t \leq \normonebig{\bbeta^\toptzero} \leq \mathcalF(\bbeta^\toptzero, \bw^\toptzero, \varepsilon_t),
$$
so that $\lim_{t \to \infty} \normonebig{\bbeta^\toptzero} = \normonebig{\widehatbbeta}$. 
Finally, applying the stable NSP and the standard consequence of stable NSP (Theorem~\ref{theorem:suff_nec_snsp}), we get
$$
\limsup_{t \to \infty} \normone{\bbeta^\toptzero - \widehatbbeta} \leq \frac{1+\rho}{1-\rho} \left( \lim_{t \to \infty} \normonebig{\bbeta^\toptzero} - \normonebig{\widehatbbeta} + 2\sigma_k(\widehatbbeta)_1 \right) = 0,
$$
which shows that $\bbeta^\toptzero \to \widehatbbeta$ since $\widehatbbeta$ is $k$-sparse. 
The error estimate~\eqref{equation:recv_irls_nsp1} follows similarly as above.

\paragraph{(ii) When $\varepsilon \triangleq \lim_{t \to \infty} \varepsilon_t > 0$.} We first show that $\bbeta^\toptzero \to \bbeta^{(\varepsilon)}$, where $\bbeta^{(\varepsilon)}$ is the minimizer of the smoothed functional $F_{\varepsilon}$ in problem \eqref{equation:irls_auxit_epsi2}. 
By property \eqref{equation:irls_bds_itera4}, the sequence $\{\bbeta^\toptzero\}$ is bounded thus has accumulation points. 
Let $\{\bbeta^{(t_i)}\}_i$ be a convergent subsequence with limit $\widehatbbeta$. We claim that $\widehatbbeta = \bbeta^{(\varepsilon)}$, which by uniqueness of $\bbeta^{(\varepsilon)}$ implies that every convergent subsequence converges to $\widehatbbeta$. This means in turn that $\widehatbbeta$ is the unique accumulation point of $\{\bbeta^\toptzero\}_t$ and therefore, $\bbeta^\toptzero$ converges to $\widehatbbeta$ as $t \to \infty$.

To see this,
as $t\to \infty$, by \eqref{equation:irls3_equiv}, it follows that $w_i^\toptzero = ((\beta_i^\toptzero)^2 + \varepsilon^2)^{-1/2} \leq \varepsilon^{-1}$, whence we have 
$$
\lim_{i \to \infty} w_i^{(t_i)} = \lim_{i \to \infty} ((\beta_i^{(t_i)})^2 + \varepsilon^2)^{-1/2} = (\bw_{\widehatbbeta,\varepsilon})_i \triangleq \widehatw_i, \quad i \in \{1,2,\ldots,p\},
$$
where  $(\bw_{\bbeta,\varepsilon})_i = (\beta_i^2 + \varepsilon^2)^{-1/2}$---the same notation as in Lemma~\ref{lemma:irls_lem3}. It follows from \eqref{equation:irls_bds_itera6} that $\{\bbeta^{(t_i+1)}\}_i$ also  converges to $\widehatbbeta$. 
By the orthogonality relation of weighted least squares from Theorem~\ref{theorem:orthogonal_weight_lsii} 
and the minimizing property (IRLS$_1$) of the weighted least squares optimization problem $\bbeta^{(t_i+1)}=\arg\min_{\bbeta} \mathcalF(\bbeta, \bw^{(t_i)}, \varepsilon_{t_i}) \text{ s.t. }  \bX\bbeta = \by$, we have, for every $\bn \in \nspace(\bX)$,
$$
\innerproduct{\widehatbbeta, \bn}_{\widehatbw} = \lim_{i \to \infty} \innerproduct{\bbeta^{(t_i+1)}, \bn}_{\bw^{(t_i)}} = 0.
$$
By Lemma~\ref{lemma:irls_lem3}, this orthogonality condition characterizes $\widehatbbeta = \bbeta^{(\varepsilon)}$.

We now establish the error bound \eqref{equation:recv_irls_nsp2}. For the true vector $\bbeta^* \in \real^p$ with $\bX\bbeta^*=\by$, the minimizing property of  $\bbeta^{(\varepsilon)}$ for $F_\varepsilon(\bbeta)$ yields
$$
\normonebig{\bbeta^{(\varepsilon)}} 
\leq F_\varepsilon(\bbeta^{(\varepsilon)}) 
\leq F_\varepsilon(\bbeta^*) 
= \sum_{i=1}^p \sqrt{\beta_i^2 + \varepsilon^2} \leq p\varepsilon + \sum_{i=1}^p \abs{\beta_i} 
= p\varepsilon + \normone{\bbeta^*}.
$$
Applying the stable NSP of order $k$ with $\widetildek \leq k$ (i.e., $\sigma_{k}(\bbeta^*)_1\leq \sigma_{\widetildek}(\bbeta^*)_1$) together with  the standard consequence of stable NSP (Theorem~\ref{theorem:suff_nec_snsp}), we obtain
\begin{equation}\label{equation:recv_irls_nsp_p1}
\normonebig{\bbeta^{(\varepsilon)} - \bbeta^*} \leq \frac{1+\rho}{1-\rho} \left( \normonebig{\bbeta^{(\varepsilon)}} - \normone{\bbeta^*} + 2\sigma_{\widetildek}(\bbeta^*)_1 \right) \leq \frac{1+\rho}{1-\rho} \left( p\varepsilon + 2\sigma_{\widetildek}(\bbeta^*)_1 \right).
\end{equation}
From the definition (IRLS$_2$) of $\varepsilon_t$ with $\eta = 1/p$ and  using the continuity of the nonincreasing rearrangement,
$\norminfbig{[\bbeta^\toptzero]^\downarrow  - [\bbeta^{(\varepsilon)}]^\downarrow } \leq \norminfbig{\bbeta^\toptzero - \bbeta^{(\varepsilon)} } $~\footnote{See Problem~\ref{prob:bound_nonincr}.},
it follows that
$$
p\varepsilon = \lim_{t \to \infty} p\varepsilon_t 
\leq \lim_{t \to \infty} [\bbeta^\toptzero]_{k+1}^\downarrow 
= [\bbeta^{(\varepsilon)}]_{k+1}^\downarrow,
$$
whence we have by multiplying the above inequality by $(k + 1 - \widetildek)$:
$$
\begin{aligned}
(k + 1 - \widetildek)p\varepsilon 
&\leq (k + 1 - \widetildek)[\bbeta^{(\varepsilon)}]_{k+1}^\downarrow \stackrel{\dag}{\leq} \normonebig{\bbeta^{(\varepsilon)} - \bbeta^*} + \sigma_{\widetildek}(\bbeta^*)_1\\
&\stackrel{\ddag}{\leq} \frac{1+\rho}{1-\rho} \left( p\varepsilon + 2\sigma_{\widetildek}(\bbeta^*)_1 \right) + \sigma_{\widetildek}(\bbeta^*)_1,
\end{aligned}
$$
where the inequality ($\dag$) follows from \ref{prob:bound_nonincr3} of Problem~\ref{prob:bound_nonincr}, and the inequality ($\ddag$) follows from \eqref{equation:recv_irls_nsp_p1}. 
Rearranging terms yields
$$
\left( k + 1 - \widetildek - \frac{1+\rho}{1-\rho} \right) \big(p\varepsilon + 2\sigma_{\widetildek}(\bbeta^*)_1\big) 
\leq 
2 \left(\frac{1}{2} + k + 1 - \widetildek\right)\sigma_{\widetildek}(\bbeta^*)_1.
$$
Under the assumption $k -\frac{2\rho}{1-\rho}  >  \widetildek$, we have  $k + 1 - \widetildek > (1+\rho)/(1-\rho)$, so the coefficient on the left is positive.
Therefore,
$$
p\varepsilon + 2\sigma_{\widetildek}(\bbeta^*)_1 \leq \frac{2(k - \widetildek) + 3}{(k - \widetildek) - \frac{2\rho}{1-\rho}} \sigma_{\widetildek}(\bbeta^*)_1.
$$
Substituting this into \eqref{equation:recv_irls_nsp_p1} yields the desired bound \eqref{equation:recv_irls_nsp2}.

Finally, suppose  that $\bbeta^*$ is $\widetildek$-sparse with $\widetildek \leq k$. Then by \eqref{equation:recv_irls_nsp2}, we have $\widehatbbeta = \bbeta^*$, so that  $\widehatbbeta$ is also $\widetildek$-sparse. Hence, $[\widehatbbeta]_{k+1}^\downarrow = 0$, and by the definition (IRLS$_2$) of $\varepsilon_t$, we have $\varepsilon = \lim_{t \to \infty} \varepsilon_t = 0$.
Thus, the case $\varepsilon>0$ cannot occur when $\bbeta^*$ is sufficiently sparse.
\end{proof}

The theorem establishes the stable recovery of the signal  $\widehatbbeta$, mirroring the result in Theorem~\ref{theorem:stablensp}.

We conclude this IRLS algorithm with an estimate of the rate of convergence.
However, this rate only becomes effective once the iterates are sufficiently close to their limit. No guarantees are provided for the initial phase of the algorithm, although practical experience suggests that this transient phase is typically short \citep{daubechies2010iteratively, foucart2013invitation}. The result below establishes linear convergence in the $\ell_1$-norm for the exactly sparse case.

\begin{theoremHigh}[Rate of convergence of IRLS under stable NSP]\label{theorem:recv_irls_nsp2}
Let $\bX \in \real^{n\times p}$ satisfy the stable NSP of order $k$ with constant $\rho < 1 - \frac{2}{k+2}$, where $\widetildek < k - \frac{2\rho}{1-\rho}$, and consider the measurements $\by=\bX\bbeta^*$, where  $\bbeta^* \in \real^p$ is $\widetildek$-sparse with support $\sS = \supp(\bbeta^*)$.
Let $0 < \nu < 1$ be such that
$$
\mu \triangleq \frac{\rho(1+\rho)}{1-\nu}\left(1 + \frac{1}{k+1-\widetildek}\right) < 1.
$$
Consider the sequence the sequence $\{\bbeta^\toptzero\}_t$ generated by the IRLS Algorithm~\ref{alg:ell1_irls} with parameters $k$ and smoothing update rule $\eta = 1/p$, 
and let $\widehatbbeta$ denote its limit. 
By Theorem~\ref{theorem:recv_irls_nsp}, we have $\widehatbbeta=\bbeta^*$. 
Suppose  $t_0 \in \naturalset$ satisfies
$$
\normone{\widehatbbeta - \bbeta^{(t_0)}} \leq R \triangleq \nu \min_{i \in \sS} \abs{\beta_i^*}.
$$
Then for all $t \geq t_0$,
\begin{equation}\label{equation:recv_irls_nsp2_eq1}
\normone{\widehatbbeta - \bbeta^\toptone} \leq \mu \normonebig{\widehatbbeta - \bbeta^\toptzero}.
\end{equation}
Consequently, the iterates converge linearly to  $\widehatbbeta$; 
specifically
$\normonebig{\bbeta^\toptzero - \widehatbbeta} \leq \mu^{(t-t_0)} \normonebig{\bbeta^{(t_0)} - \widehatbbeta}$ for all $t > t_0$.
\end{theoremHigh}

\begin{proof}[of Theorem~\ref{theorem:recv_irls_nsp2}]
The proof follows closely that of \citet{daubechies2010iteratively, foucart2013invitation}, with minor adaptations.
For each $t>0$,  denote the error vector $\bn^\toptzero \triangleq \bbeta^\toptzero - \widehatbbeta \in \nspace(\bX)$ (since both $\bbeta^\toptzero$ and $\widehatbbeta$ are feasible with $\by=\bX\bbeta^\toptzero=\bX\widehatbbeta$). 
By the minimizing property (IRLS$_1$) of $\bbeta^\toptone$ and the orthogonal condition for the minimizer of weighted least squares (Theorem~\ref{theorem:orthogonal_weight_lsii}),  we have
$$
0 = \innerproduct{\bbeta^\toptone, \bn^\toptone}_{\bw^\toptzero} = \innerproduct{\widehatbbeta + \bn^\toptone, \bn^\toptone}_{\bw^\toptzero}.
$$
Since $\supp (\widehatbbeta) = \supp (\bbeta^*) = \sS$, rearranging yields
\begin{equation}\label{equation:recv_irls_nsp2_pv1}
\sum_{i=1}^p \big(n_i^\toptone\big)^2 w_i^\toptzero 
= - \sum_{i \in \sS} \widehatbeta_i n_i^\toptone w_i^\toptzero
= - \sum_{i \in \sS} \frac{\widehatbeta_i n_i^\toptone}{\sqrt{(\beta_i^{(t)})^2 + \varepsilon_t^2}}  ,
\end{equation}
where the last equality follows from \eqref{equation:irls3_equiv}.

Now fix $t \geq t_0$, and denote $C_t \triangleq \normonebig{\bbeta^\toptzero-\widehatbbeta } =  \normonebig{\bn^\toptzero}\leq R$ (see arguments below). 
For any $i \in \sS$, we have $\absbig{n_i^\toptzero} \leq \normonebig{\bn^\toptzero} = C_t \leq \nu \absbig{\widehatbeta_i}$, which implies
\begin{equation}\label{equation:recv_irls_nsp2_pv2}
\frac{\absbig{\widehatbeta_i}}{\sqrt{(\beta_i^\toptzero)^2 + \varepsilon_t^2}} 
\leq \frac{\absbig{\widehatbeta_i}}{\absbig{\beta_i^\toptzero}} 
= \frac{\absbig{\widehatbeta_i}}{\absbig{\widehatbeta_i + n_i^\toptzero}} 
\leq \frac{1}{1-\nu}.
\end{equation}
By combining \eqref{equation:recv_irls_nsp2_pv1} and \eqref{equation:recv_irls_nsp2_pv2}, 
we obtain
$$
\sum_{i=1}^p \big(n_i^\toptone\big)^2 w_i^\toptzero 
\leq \frac{1}{1-\nu} \normone{\bn_{{\sS}}^\toptone} 
\leq \frac{\rho}{1-\nu} \normone{\bn_{\comple{\sS}}^\toptone},
$$
where the last inequality follows from  the definition of the stable NSP (Definition~\ref{definition:stable_nsp}).
Using the Cauchy--Schwarz inequality with  the above inequality  yields~\footnote{Since $\bx^\top\by\leq \normtwo{\bx}\normtwo{\by}$, it holds that 
$(\sum a_i)^2 \leq (\sum a_i^2 b_i )(\sum \frac{1}{b_i})$ if $a_i, b_i>0$ by setting 
$x_i\triangleq {\abs{a_i}}\sqrt{b_i}, y_i\triangleq \frac{1}{\sqrt{b_i}}$ for all $i$.}
\begin{align}
\normone{\bn_{\comple{\sS}}^\toptone}^2 
&\leq \left(\sum_{i \in \comple{\sS}} \big(n_i^\toptone\big)^2 w_i^\toptzero\right)
\left(\sum_{i \in \comple{\sS}} \sqrt{\big(\beta_i^\toptzero\big)^2 + \varepsilon_t^2}\right) \nonumber \\
&\leq \left(\sum_{i=1}^p \big(n_i^\toptone\big)^2  w_i^\toptzero\right) (\normonebig{\bbeta_{\comple{\sS}}^\toptzero} + p\varepsilon_t) 
\leq \frac{\rho}{1-\nu} \normone{\bn_{\comple{\sS}}^\toptone}(\normonebig{\bn^\toptzero} + p\varepsilon_t).
\label{equation:recv_irls_nsp2_pv3}
\end{align}

If $\bn_{\comple{\sS}}^\toptone = \bzero$, then $\bbeta_{\comple{\sS}}^\toptone = \bzero$, so that $\bbeta^\toptone$ is $\widetildek$-sparse, and the algorithm terminates by definition.  Since $\bbeta^\toptone - \widehatbbeta \in \nspace(\bX)$, and the NSP ensures that $\nspace(\bX)$ contains no nonzero $\widetildek$-sparse vectors, it follows that $\bbeta^\toptone = \widehatbbeta$.
Hence, $C_{t+1} = 0$, and \eqref{equation:recv_irls_nsp2_eq1} is trivially satisfied.

Otherwise, if $\bn_{\comple{\sS}}^\toptone \neq \bzero$, we may divide both sides of \eqref{equation:recv_irls_nsp2_pv3} by $\normonebig{\bn_{\comple{\sS}}^\toptone}$ to obtain
$$
\normone{\bn_{\comple{\sS}}^\toptone}\leq \frac{\rho}{1-\nu} \left(\normone{\bn^\toptzero} + p\varepsilon_t\right).
$$
Using the stable NSP again (Definition~\ref{definition:stable_nsp}), we bound the full error:
$$
\normone{\bn^\toptone} 
= \normone{\bn_\sS^\toptone} + \normone{\bn_{\comple{\sS}}^\toptone}
\leq (1+\rho) \normone{\bn_{\comple{\sS}}^\toptone}
\leq \frac{\rho(1+\rho)}{1-\nu} \left(\normone{\bn^\toptzero} + p\varepsilon_t\right).
$$
By the smoothing update rule (IRLS$_2$) for $\varepsilon_{t+1}$, the assumption of $\eta=1/p$, and \ref{prob:bound_nonincr3} of Problem~\ref{prob:bound_nonincr}, we have
$$
p\varepsilon_t 
\leq [\bbeta^\toptzero]_{k+1}^\downarrow
\leq \frac{1}{k+1-\widetildek} \left(\normone{\bbeta^\toptzero - \widehatbbeta} + \sigma_k(\widehatbbeta)_1\right) 
= \frac{\normonebig{\bn^\toptzero}}{k+1-\widetildek},
$$
where $\sigma_k(\widehatbbeta)_1 = 0$ since $\widehatbbeta$ is $\widetildek$-sparse.
Putting everything together,
$$
C_{t+1} = \normone{\bn^\toptone} \leq \frac{\rho(1+\rho)}{1-\nu} \left(1 + \frac{1}{k+1-\widetildek}\right) \normonebig{\bn^\toptzero} = \mu C_t.
$$
Because $\mu < 1$ by assumption, it follows that $C_{t+1} \leq R$, and thus the condition $C_{t} \leq R$ propagates for all $t\geq t_0$. 
Therefore, $C_{t+1} \leq \mu C_t$ for all $t \geq t_0$, completing the proof.
\end{proof}

Note that if $\rho$ is sufficiently small---specifically, if $\rho(1+\rho) < 2/3$---then for any $\widetildek \leq k-1$, there exists some $\nu > 0$ such that $\mu < 1$. 
Consequently $\bbeta^\toptzero$ converges linearly to $\widehatbbeta$ whenever $\bbeta^*$ is $(k-1)$-sparse.

\index{Augmented Lagrangian function}
\index{First-order augmented Lagrangian}
\subsection{First-Order Augmented Lagrangian Algorithm}

Similar to the penalty function method (Algorithm~\ref{alg:quad_pen_eq}), the $\ell_1$-minimization problem~\eqref{opt:p1} can be efficiently solved in its Lagrangian form (see e.g., \citet{boyd2004convex, figueiredo2008gradient, aybat2012first}):
\begin{equation}\label{equation:fal_lp_lag_form}
\min_{\bbeta \in \real^p} \left\{ \sigma \normone{\bbeta} + \frac{1}{2} \normtwo{\bX\bbeta - \by}^2 \right\},
\end{equation}
where the penalty parameter $\sigma>0$ is gradually decreased toward zero.

An \textit{augmented Lagrangian function (ALF)} for solving this problem is given by
\begin{equation}\label{eq:aug_lag_function}
L_{\sigma} (\bbeta, \blambda ) 
= \sigma \normone{\bbeta} - \sigma\blambda^\top (\bX\bbeta - \by) + \frac{1}{2} \normtwo{\bX\bbeta - \by}^2,
\end{equation}
where $\blambda$ is the vector of Lagrange multipliers associated with the equality constraints $\bX\bbeta = \by$.
In other words, the augmented Lagrangian function modifies the standard Lagrangian $L(\bbeta, \blambda) = \normone{\bbeta} - \blambda^\top (\bX\bbeta - \by)$ by adding a quadratic penalty term, which justifies the name ``augmented" Lagrangian \citep{lu2025practical}.

More generally,
the  key idea of the augmented Lagrangian method for solving an optimization problem with objective $f(\bbeta)$  and equality constraints $\bc(\bbeta) = [c_i(\bbeta)]_{i\in\mathcalE}$ is outlined  in Algorithm~\ref{alg:aug_lag_method}, where $\mathcalE$ denotes the index set of equality constraints.
We iteratively update the parameter $\bbeta^\toptzero$ and the multiplier $\blambda^\toptzero$ while progressively reducing the penalty parameter    $\sigma_t$ for each iteration $t$.
In the early iterations, it is common to fix the multiplier  (e.g., set  $ \blambda = \bzero $) and allow the penalty parameter $ \sigma $ to dominate, which helps drive the iterates toward feasibility and brings us closer to the true solution $ \bbeta_* $.
\begin{algorithm}[h] 
\caption{Augmented Lagrangian Method: the General Case \citep{lu2025practical}}
\label{alg:aug_lag_method}
\begin{algorithmic}[1] 
\Require A function $f(\bbeta)$ and a set of equality constraints $\{c_i(\bbeta)\}, i\in\mathcalE$; 
\State {\bfseries input:}   Choose initial point $ \bbeta^{(1)} $, multiplier $ \blambda^{(1)} $, penalty factor $ \sigma_1 > 0 $, penalty factor update constant $ \rho \in (0,1) $, constraint violation tolerance $ \varepsilon > 0 $, and precision requirement $ \gamma_t > 0 $.
\State {\bfseries define:} $L_\sigma (\bbeta, \blambda)=f(\bbeta)+ \blambda^\top\bc(\bbeta)+ \frac{1}{2}\sigma\bc(\bbeta)^\top\bc(\bbeta)$; 
\For{$t=1,2,\ldots$}
\State   \algoalign{With $ \bbeta^\toptzero $ as the initial point, solve $ \mathopmin{\bbeta}  L_{\sigma_t}(\bbeta, \blambda^\toptzero)$,
to obtain a solution $ \bbeta^\toptone $ satisfying the precision condition 
$\normtwobig{\nabla_{\bbeta} L_{\sigma_t}(\bbeta^\toptone, \blambda^\toptzero)} \leq \gamma_t$.
}
\If{$ \normtwobig{\bc(\bbeta^\toptone)} \leq \varepsilon $}
\State Return approximate solution $ \bbeta^\toptone, \blambda^\toptzero $, terminate iteration.
\EndIf
\State Update multipliers: $ \blambda^\toptone \gets \blambda^\toptzero +  \bc(\bbeta^\toptone) /\sigma_t$.
\State Update penalty factor: $ \sigma_{t+1} \gets \rho \sigma_t $.
\EndFor
\State \Return  final $\bbeta \gets \bbeta^\toptzero$;
\end{algorithmic} 
\end{algorithm}

As a refinement of the augmented Lagrangian method tailored to equality-constrained $\ell_1$-minimization problems, Algorithm~\ref{alg:fal} presents the \textit{first-order augmented Lagrangian (FAL)} algorithm for solving problem~\eqref{opt:p1}.
The subroutine \texttt{PGM} used in Step (FAL$_7$) is defined in Algorithm~\ref{alg:prox_gd_gen_fal} (see Section~\ref{section:proxiGD_inClasso}) \citep{aybat2012first}.~\footnote{The original paper \citet{aybat2012first} employs an accelerated proximal gradient method to solve the subproblem; for simplicity, we use the standard proximal gradient method here.}

The FAL algorithm solves the $\ell_1$-minimization problem by approximately solving a sequence of subproblems of the form
\begin{equation}\label{equation:fal_subproblem}
\begin{aligned}
\min_{\bbeta \in \real^p : \normone{\bbeta} \leq \eta_t} 
&\left\{ \sigma_t \normone{\bbeta} - \sigma_t \blambda^\toptzeroTOP (\bX\bbeta - \by) + \frac{1}{2} \normtwo{\bX\bbeta - \by}^2 \right\}\\
\equiv \min_{\bbeta \in \real^p : \normone{\bbeta} \leq \eta_t}
&P_t(\bbeta) 
\triangleq \sigma_t \normone{\bbeta} + \frac{1}{2} \normtwo{\bX\bbeta-\by - \sigma_t \blambda^\toptzero}^2
\end{aligned}
\end{equation}
for an appropriately chosen sequence of $\{\sigma_t, \blambda^\toptzero, \eta_t\}, t = 0, 1, \dots$
The FAL stopping criterion, denoted \texttt{FALstop}, is defined as
\begin{equation}\label{eq:fal_stop_noiseless}
\texttt{FALstop} =\norminf{\bbeta^\toptzero - \bbeta^\toptminus } \leq \gamma,
\end{equation}
where $\gamma$ is a prescribed tolerance.
Finally, the \texttt{PGM} subroutine computes an $\epsilon_t$-optimal solution to the subproblem~\eqref{equation:fal_subproblem} at iteration $t$.

\begin{algorithm}[ht]
\caption{\texttt{FAL}$\{(\sigma_t, \epsilon_t, \tau_t)_{t \in \integerset_+}, \widetildeeta\}$~\citet{aybat2012first}}
\label{alg:fal}
\begin{algorithmic}[1]
\Ensure Find $\bbeta = \argmin_{\bbeta} \left\{ \sigma_t \normone{\bbeta} - \sigma_t \blambda^\toptzeroTOP (\bX\bbeta - \by) + \frac{1}{2} \normtwo{\bX\bbeta - \by}^2 \right\}$;
\State \textbf{input:}  $L_b = \sigma_{\max}(\bX\bX^\top)$, $\rho_1,\rho_2,\rho_3\in(0,1)$;
\State \textbf{initialization:} $\bbeta^\topzero$, $\blambda^\topzero \gets \bzero$, $t\gets 0$;
\While{(\texttt{FALstop} not true)}
\State $g_t(\bbeta) \gets \sigma_t \normone{\bbeta}$, $f_t(\bbeta) \gets \frac{1}{2}\normtwo{\bX\bbeta-\by-\sigma_t \blambda^\toptzero}^2$;  \Comment{(FAL$_1$)}
\State $P_t(\bbeta) \gets f_t(\bbeta)+ g_t(\bbeta)$;\Comment{(FAL$_2$)}
\State $\eta_t \gets \widetildeeta + (\sigma_t / 2) \normtwo{\blambda^\toptzero}$; \Comment{(FAL$_3$)}
\State $\sS^\toptzero \gets \{\bbeta \in \real^p \mid \normonebig{\bbeta} \leq \eta_t\}$; \Comment{(FAL$_4$)}
\State $\ell_{t,\max} \gets \sigma_{\max}(\bX)(\eta_t + \normtwobig{\bbeta^\toptminus})\sqrt{2/\epsilon_t}$; 
\Comment{(FAL$_5$)}
\State \texttt{PGMstop} = $\{\ell \geq \ell_{t,\max}\}$ or $\{\exists\, \bg \in \partial P_t(\bzeta^{(\ell)}), \text{ with } \normtwo{\bg} \leq \tau_t\}$;
\Comment{(FAL$_6$)}
\State $\bbeta^\toptzero \gets \texttt{PGM}(g_t, f_t, L_b, \sS^\toptzero, \bbeta^\toptminus, \texttt{PGMstop})$, see Algorithm~\ref{alg:prox_gd_gen_fal};
\Comment{(FAL$_7$)}
\State $\blambda^\toptone \gets \blambda^\toptzero - (\bX\bbeta^\toptzero - \by)/\sigma_t$;\Comment{(FAL$_8$)}
\State \textcolor{black}{$\sigma_{t+1}\gets \rho_1 \sigma_t$, $\epsilon_{t+1}\gets \rho_2 \epsilon_t$, $\tau_{t+1}\gets \rho_3 \tau_t$};
\Comment{(FAL$_9$)}
\State $t \gets t + 1$;
\EndWhile
\State \textbf{output:} $\bbeta \gets \bbeta^\toptzero$
\end{algorithmic}
\end{algorithm}

\begin{algorithm}[h] 
\caption{Proximal Gradient Method: \texttt{PGM}$(g_t, f_t, L_b, \sS^\toptzero, \bbeta^\topzero, \texttt{PGMstop})$}
\label{alg:prox_gd_gen_fal}
\begin{algorithmic}[1] 
\Require A function $f(\bbeta)=f_t(\bbeta)$ and a closed convex function $g(\bbeta)=g_t(\bbeta)$ satisfying (A1) and (A2) in \textcolor{black}{Theorem~\ref{theorem:prox_conv_ss_cvx}}. $f(\bbeta)$ is $L_b$-smooth (see Example~\ref{example:lipschitz_spar});
\State {\bfseries initialize:}   $\bzeta^\topzero \gets \bbeta^\topzero$, $\ell \gets 0$;
\While{(\texttt{PGMstop} not true)}
\State Pick a stepsize $\eta_\ell = \frac{1}{L_\ell} = \frac{1}{L_b}$; 
\Comment{(PGM$_1$)}
\State $\balpha^{(\ell+1)} \gets \bzeta^{(\ell)} - \eta_\ell \nabla f(\bzeta^{(t)})$;
\Comment{(PGM$_2$)}
\State $\bzeta^{(\ell+1)} \gets \prox_{\eta_\ell g}(\balpha^{(\ell+1)}) \triangleq \mathcalT_{L_\ell}^{f,g}(\bzeta^\toptzero)$;
\Comment{(PGM$_3$)}
\EndWhile
\State \textbf{return} $\bzeta^\toptzero$;
\end{algorithmic} 
\end{algorithm}


More specifically, algorithm FAL takes as input a sequence $\{(\sigma_t, \epsilon_t, \tau_t)\}_{t \in \integerset_+}$, an initial point $\bbeta^\topzero$, and a bound $\widetildeeta$ on the $\ell_1$-norm of an optimal solution $\bbeta_*$ of the $\ell_1$-minimization problem. One such bound $\widetildeeta$ can be computed as follows. 
Let 
\begin{equation}\label{equation:fal_opt_l2}
\widetilde{\bbeta} \triangleq \argmin \{\normtwo{\bbeta} \mid \bX\bbeta = \by\} = \bX^\top(\bX\bX^\top)^{-1}\by
\end{equation} 
denote the solution to the corresponding $\ell_2$-minimization problem. 
By definition we have 
\begin{equation}\label{equation:fal_eta_ineq}
\normone{\bbeta_*} \leq \widetildeeta \triangleq \normonebig{\widetilde{\bbeta}}.
\end{equation}

We now describe each component of the algorithm in detail.
As mentioned above,  an augmented Lagrangian function for the $\ell_1$-minimization problem~\eqref{opt:p1} can be written as
\begin{equation}
P(\bbeta) \triangleq 
L_{\sigma} (\bbeta, \blambda ) 
= \sigma \normone{\bbeta} - \sigma \blambda^\top (\bX\bbeta-\by) + \frac{1}{2} \normtwo{\bX\bbeta-\by}^2,
\end{equation}
where $\sigma$ is the penalty parameter and $\blambda$ is a dual variable associated with the constraint $\bX\bbeta = \by$. 
According to Steps~(FAL$_2$)$\sim$(FAL$_5$) of Algorithm~\ref{alg:fal}, FAL approximately minimizes this augmented Lagrangian while keeping $\blambda$ fixed at its current value $ \blambda^\toptzero$. 
Specifically, it minimizes the function
\begin{align}
P_t(\bbeta) 
&= \sigma_t \normone{\bbeta} + \frac{1}{2} \normtwo{\bX\bbeta-\by - \sigma_t \blambda^\toptzero}^2 \label{equation:fal_ptbeta1}\\
&= \sigma_t \normone{\bbeta} - \sigma_t (\blambda^\toptzero)^\top (\bX\bbeta-\by) + \frac{1}{2} \normtwo{\bX\bbeta-\by}^2 + \frac{1}{2} \normtwo{\sigma_t \blambda^\toptzero}^2, \nonumber
\end{align}
over the constrained set $\sS^\toptzero = \{\bbeta \mid \normone{\bbeta} \leq \eta_t\}$ using Algorithm PGM (or PGD; see Algorithm~\ref{alg:pgd_gen}).

\paragrapharrow{The choice of $\eta_t$.}
When applying PGM, we must ensure that the minimizer of $P_t(\bbeta)$ over $\real^p$ lies within the feasible set $\sS^\toptzero $, i.e.,  $\sS^\toptzero \cap \argmin_{\bbeta \in \real^p} \{P_t(\bbeta) \} \neq \varnothing$. Let $\bbeta_*^\toptzero \in \argmin_{\bbeta \in \real^p} P_t(\bbeta) $. 
Since $P_t(\bbeta_*^\toptzero) \leq P_t(\bbeta_*)$, $\bX\bbeta_* = \by$ and $\normone{\bbeta_*} \leq \widetildeeta$ by \eqref{equation:fal_eta_ineq}, we have $\normonebig{\bbeta_*^\toptzero} \leq \eta_t$. 
Thus, $\bbeta_*^\toptzero \in \sS^\toptzero$, and this explains the choice of 
\begin{equation}\label{equation:fal_etat}
\eta_t= \widetildeeta + \frac{\sigma_t}{2} \normtwobig{\blambda^\toptzero}
\end{equation}
in algorithm FAL, which guarantees $\bbeta_*^\toptzero \in \sS^\toptzero$.

\paragrapharrow{Stopping criterion.}
We now discuss the stopping condition in Step~(FAL$_7$) of algorithm FAL.
The convergence rate of the proximal gradient method (see {Theorem~\ref{theorem:prox_conv_ss_cvx}}) ensures that
\begin{equation*}
P_t(\bzeta^\toplzero) \leq P_t(\bbeta_*^\toptzero) + \epsilon_t,
\quad \text{for} \quad \ell \geq   \frac{ L_b \normtwo{\bbeta^\toptminus - \bbeta_*^\toptzero}^2}{2\epsilon_t} +1,
\end{equation*}
where $\{\bzeta^\toplzero\}_{\ell \in \integerset_+}$ is the sequence of $\bzeta$-iterates when Algorithm PGM is applied to the $t$-th subproblem and $L_b$ denotes the Lipschitz constant of the gradient $\nabla f_t(\bbeta)$ (i.e., $L_b$-smooth; Definition~\ref{definition:scss_func}). 
Since $\normonebig{\bbeta_*^\toptzero} \leq \eta_t$ and $\normone{\cdot} \geq \normtwo{\cdot}$ (Exercise~\ref{exercise:cauch_sc_l1l2}),  the triangle inequality implies that
\begin{equation}
\normtwo{\bbeta_*^\toptzero - \bbeta^\toptminus}
\leq \normone{\bbeta_*^\toptzero} + \normtwo{\bbeta^\toptminus} \leq \eta_t + \normtwo{\bbeta^\toptminus}.
\end{equation}
Since $f_t(\bbeta) \triangleq \frac{1}{2}\normtwo{\bX\bbeta-\by-\sigma_t \blambda^\toptzero}^2$ and $\nabla f_t(\bbeta) = \bX^\top(\bX\bbeta-\by - \sigma_t\blambda^\toptzero)$, it follows that $L_b = \sigma_{\max}^2(\bX)$, where $\sigma_{\max}(\bX)$ denote the largest singular value of $\bX$. Thus, it follows that
\begin{equation*}
\frac{ L_b \normtwo{\bbeta^\toptminus - \bbeta_*^\toptzero}^2}{2\epsilon_t} +1
\leq 
\frac{\sigma_{\max}^2(\bX)}{2\epsilon_t}\left(\eta_t + \normtwo{\bbeta^\toptminus}\right)^2
\triangleq \ell_{t,\max}.
\end{equation*}
Therefore, running PGM for at least $\ell_{t,\max}$ iterations guarantees an $\epsilon_t$-optimal solution to the $t$-th subproblem.

Consequently, it follows that the stopping criterion \texttt{PGMstop} in Step~(FAL$_7$) of Algorithm~FAL ensures that the iterate $\bbeta^\toptzero$ satisfies one of the following two conditions:
\begin{subequations}
\begin{align}
& P_t(\bbeta^\toptzero) \leq \min_{\bbeta \in \real^p} P_t(\bbeta_*^\toptzero)  + \epsilon_t; \label{eq:apgstop_conditions1} \\
& \exists\, \bg \in \partial P_t(\bbeta^\toptzero) \text{ with } \normtwo{\bg} \leq \tau_t,
\label{eq:apgstop_conditions2}
\end{align}
\end{subequations}
where $\partial P_t(\bbeta^\toptzero)$ denotes the set of subgradients of the function $P_t$ at $\bbeta^\toptzero$, and $\bbeta_*^\toptzero \in \argmin_{\bbeta \in \real^p} P_t(\bbeta) $.

Finally, in Step~(FAL$_8$), the dual variables $\blambda^\toptzero$  are updated using the standard rule employed in augmented Lagrangian methods (cf. Algorithm~\ref{alg:aug_lag_method}).

\subsection*{Convergence of FAL}
We now present the main convergence results for algorithm FAL.
The following theorem establishes that every limit point of the iterate sequence $\{\bbeta^\toptzero\}_{t \in \integerset_+}$ generated by FAL is an optimal solution of the $\ell_1$-minimization problem~\eqref{opt:p1}.

To prove this result, we require the following lemmas.
\begin{lemma}[$\epsilon$-optimal gradient lemma]\label{lemma:smooth_prop_fal}
Let $ h : \real^p \to \real $ be a convex function, and suppose that $ f $ is $L_b$-smooth   (i.e., its gradient $\nabla h$ is Lipschitz continuous  with constant $L_b$, 
as stated in Theorem~\ref{theorem:equi_gradsch_smoo}). 
Fix $ \epsilon > 0 $, and let $ \widebarbbeta \in \real^p $ satisfy
$$
\sigma \normone{\widebarbbeta} + h(\widebarbbeta) - \big(\sigma \normonebig{ \widehatbbeta} + h(\widehatbbeta)\big) \leq \epsilon,
$$
where $ \widehatbbeta \in \argmin\{ \sigma \normone{\bbeta} + h(\bbeta) : \bbeta \in \real^p \} $ (i.e., $\widebarbbeta$ is $\epsilon$-optimal for this minimization problem). 
Then
\begin{equation}\label{equation:smooth_prop_fal}
\frac{1}{2L_b} \sum_{i: \abs{\nabla h_i(\widebarbbeta)} > \sigma} \big(\abs{\nabla h_i(\widebarbbeta)} - \sigma\big)^2 \leq \epsilon. 
\end{equation}
The bound \eqref{equation:smooth_prop_fal} implies $ \norminf{\nabla h(\widebarbbeta)} \leq \sqrt{2L_b\epsilon} + \sigma $.
\end{lemma}

\begin{proof}[of Lemma~\ref{lemma:smooth_prop_fal}]
Triangular inequality for $ \normone{\cdot} $ and $L_b$-smoothness (Definition~\ref{definition:scss_func}) implies that for all $ \balpha \in \real^p $
$$
\sigma \normone{\balpha} + h(\balpha) \leq \sigma \normone{\widebarbbeta} + h(\widebarbbeta) + \nabla h(\widebarbbeta)^\top(\balpha - \widebarbbeta) + \frac{L_b}{2} \normtwo{\balpha - \widebarbbeta}^2 + \sigma \normone{\balpha - \widebarbbeta}.
$$
Taking the minimum over $ \balpha \in \real^p $ on both sides yields
\begin{equation}\label{equation:smooth_prop_fal_pv1}
\sigma \normonebig{\widehatbbeta} + h(\widehatbbeta) \leq \sigma \normone{\widebarbbeta} + h(\widebarbbeta) + \min_{\balpha \in \real^p} \left\{ \nabla h(\widebarbbeta)^\top(\balpha - \widebarbbeta) + \frac{L_b}{2} \normtwo{\balpha - \widebarbbeta}^2 + \sigma \normone{\balpha - \widebarbbeta} \right\}. 
\end{equation}
Let $ \widebarbg \triangleq \nabla h(\widebarbbeta) $. 
The minimization problem on the right-hand side can be rewritten as
\begin{align*}
\widehatbalpha
&= \argmin_{\balpha \in \real^p} \left\{ \widebarbg^\top(\balpha - \widebarbbeta) + \frac{L_b}{2} \normtwo{\balpha - \widebarbbeta}^2 + \sigma \normone{\balpha - \widebarbbeta} \right\}  \\
&= \argmin_{\balpha \in \real^p} \left\{ \frac{1}{2} \normtwo{\balpha - \widebarbbeta - \frac{-\widebarbg}{L_b}}^2 + \frac{\sigma}{L_b} \normone{\balpha - \widebarbbeta} \right\}  \\
&\stackrel{\dag}{=} \widebarbbeta + \sign\left( \frac{- \widebarbg}{L_b} \right) \hadaprod \max\left\{\abs{\frac{- \widebarbg}{L_b}} - \frac{\sigma}{L_b}\bone, \bzero \right\} 
= \widebarbbeta -\frac{\sign(\widebarbg)}{L_b} \hadaprod \max\{ \abs{\widebarbg} - \sigma\bone, \bzero \} ,
\end{align*}
where the equality ($\dag$) follows from Example~\ref{example:soft_thres} (the proximal operator of the $\ell_1$-norm),
and all  vector operators such as $ \abs{\cdot} $, $ \sign(\cdot) $ and $ \max\{\cdot,\cdot\} $ are defined to operate component-wise. Substituting $ \widehatbalpha $ in \eqref{equation:smooth_prop_fal_pv1}, we get
\begin{align*}
&\min_{\balpha \in \real^p} \left\{ \widebarbg^\top(\balpha - \widebarbbeta) + \frac{L_b}{2} \normtwo{\balpha - \widebarbbeta}^2 + \sigma \normone{\balpha - \widebarbbeta} \right\} \nonumber \\
&= -\sum_i \frac{\abs{\overline{g}_i}}{L_b} \max\{ \abs{\overline{g}_i} - \sigma, 0 \} + \frac{1}{2L_b} \sum_i \max\{ \abs{\overline{g}_i} - \sigma, 0 \}^2 + \frac{\sigma}{L_b} \sum_i \max\{ \abs{\overline{g}_i} - \sigma, 0 \} \nonumber \\
&= \frac{1}{L_b} \sum_{i: \abs{\overline{g}_i} > \sigma} \left( -\abs{\overline{g}_i} + \frac{1}{2}(\abs{\overline{g}_i} - \sigma) + \sigma \right)(\abs{\overline{g}_i} - \sigma)
= -\frac{1}{2L_b} \sum_{i: \abs{\overline{g}_i} > \sigma} (\abs{\overline{g}_i} - \sigma)^2. 
\end{align*}
The bound \eqref{equation:smooth_prop_fal} follows from the assumption that $ \sigma \normone{\widebarbbeta} + h(\widebarbbeta) - (\sigma \normonebig{\widehatbbeta} + h(\widehatbbeta)) \leq \epsilon $. The bound \eqref{equation:smooth_prop_fal} clearly implies that $ \abs{\overline{g}_i} \leq \sqrt{2L_b\epsilon} + \sigma $ for all $ i $ by the standard bounds on vector norms (Exercise~\ref{exercise:cauch_sc_l1l2}), i.e., $ \norminf{\nabla h(\widebarbbeta)} \leq \sqrt{2L_b\epsilon} + \sigma $. 
\end{proof}

This lemma provides a quantitative bound on the gradient of $h$ at any $\epsilon$-optimal solution of the augmented Lagrangian subproblem~\eqref{equation:fal_ptbeta1}.
An immediate consequence is stated below.
\begin{lemma}\label{lemma:fal_aug_bd}
Suppose $ \bX \in \real^{n \times p} $ with $ n \leq p $ and full rank, 
and let $\blambda\in\real^n$ be a fixed vector. 
Define the function
$$
P(\bbeta) \triangleq \sigma \normone{\bbeta} + \frac{1}{2} \normtwo{\bX\bbeta - \by - \sigma\blambda}^2.
$$
Assume that $ \widebarbbeta $ is an $ \epsilon $-optimal solution to  $ \min_{\bbeta \in \real^p} P(\bbeta) $, i.e., $
0 \leq P(\widebarbbeta) - \min_{\bbeta \in \real^p} P(\bbeta) \leq \epsilon$.
Then the following bounds hold:
\begin{subequations}
\begin{align}
\norminf{\bX^\top(\bX\widebarbbeta - \by - \sigma\blambda)}
&\leq 
{\sqrt{2\epsilon}\, \sigma_{\max}(\bX) + \sigma};
\label{equation:fal_aug_bd_eq1} \\
\normtwo{\bX\widebarbbeta - \by - \sigma\blambda}
&\leq \frac{\sqrt{p}}{\sigma_{\min}(\bX)} \left( \sqrt{2\epsilon}\, \sigma_{\max}(\bX) + \sigma \right),
\label{equation:fal_aug_bd_eq2}
\end{align}
\end{subequations}
where $ \sigma_{\max}(\bX) $ and $ \sigma_{\min}(\bX) $ denote the largest and smallest nonzero singular values of $\bX$, respectively.
\end{lemma}
\begin{proof}[of Lemma~\ref{lemma:fal_aug_bd}]
Let $ h(\bbeta) \triangleq  \frac{1}{2} \normtwo{\bX\bbeta - \by - \sigma\blambda}^2 $. Then $ \nabla h(\bbeta) = \bX^\top(\bX\bbeta - \by - \sigma\blambda) $. For any $ \bbeta, \balpha \in \real^p $, we have
$$
\normtwo{\nabla h(\bbeta) - \nabla h(\balpha)} = \normtwo{\bX^\top(\bX\bbeta - \bX\balpha)} \leq \sigma_{\max}^2(\bX)\normtwo{\bbeta - \balpha},
$$
see the singular value inequality \eqref{equation:svd_stre_bd}. Thus, $ h : \real^p \to \real $ is a convex function and $h$ is $L_b$-smooth with constant $ L_b = \sigma_{\max}^2(\bX) $. 

Since $ \widebarbbeta $ is an $ \epsilon $-optimal solution to $ \min_{\bbeta \in \real^p} P(\bbeta) = \min_{\bbeta \in \real^p} \{ \sigma \normone{\bbeta} + h(\bbeta) \} $, Lemma~\ref{lemma:smooth_prop_fal} immediately implies \eqref{equation:fal_aug_bd_eq1}. 
Using the singular value inequality \eqref{equation:svd_stre_bd} and standard bounds on vector norms (Exercise~\ref{exercise:cauch_sc_l1l2}) again, the second bound follows from the fact that
$$
\normtwo{\bX\widebarbbeta - \by - \sigma\blambda} \leq \frac{\normtwo{\bX^\top(\bX\widebarbbeta - \by - \sigma\blambda)}}{\sigma_{\min}(\bX)} \leq \frac{\sqrt{p}}{\sigma_{\min}(\bX)} \norminf{\bX^\top(\bX\widebarbbeta - \by - \sigma\blambda)}.
$$
This completes the proof.
\end{proof}

Using the lemma above, we now establish the convergence of algorithm FAL.
\begin{theoremHigh}[Convergence of FAL \citep{aybat2012first}]\label{theorem:convergence_fal}
Assume that  $\bX\in\real^{n\times p}$ has full row rank. 
Fix an initial point $\bbeta^\topzero \in \real^p$, a constant $\widetildeeta > 0$ such that $\widetildeeta \geq \normone{\bbeta_*}$, and a sequence of parameters $\{(\sigma_t, \epsilon_t, \tau_t)\}_{t \in \integerset_+}$ satisfying the following conditions:
\begin{enumerate}[(i)]
\item penalty parameters, $\sigma_t \searrow 0$,
\item approximate optimality parameters, $\epsilon_t \searrow 0$ such that ${\epsilon_t}/{\sigma_t^2} \leq B_1$ for all $t \geq 1$,
\item approximate feasibility parameters, $\tau_t \searrow 0$ such that ${\tau_t}/{\sigma_t} \leq B_2$ for all $t \geq 1$, and ${\tau_t}/{\sigma_t} \to 0$ as $t \to \infty$.
\end{enumerate}
\noindent
Let $\mathcalX = \{\bbeta^\toptzero\}_{t \in \integerset_+}$ denote the sequence of iterates generated by algorithm FAL (Algorithm~\ref{alg:fal}) using these parameters. 
Then, $\mathcalX$ is a bounded sequence and any limit point $\widehatbbeta$ of $\mathcalX$ is an optimal solution of the $\ell_1$-minimization problem~\eqref{opt:p1}.
\end{theoremHigh}
\begin{proof}[of Theorem~\ref{theorem:convergence_fal}]
Since $\sigma_{\max}(\bX)$ is finite for all $t \geq 1$, it follows that the sequence $\mathcalX$ exists.

\paragraph{Boundedness.}
To show that $\mathcalX$ is bounded, we first establish a uniform bound on the sequence of dual multipliers $\{\blambda^\toptzero\}_{t \in \integerset_+}$, and then use \eqref{equation:fal_etat} to bound the iterates $\{\bbeta^\toptzero\}$. 
At iteration $t$, suppose Algorithm PGM terminates with an iterate  $\bbeta^\toptzero$ satisfying Condition~\eqref{eq:apgstop_conditions1}. 
Then Lemma~\ref{lemma:fal_aug_bd} applied to $P_t(\bbeta)  = \sigma_t \normone{\bbeta} + \normtwobig{\bX\bbeta-\by - \sigma_t\blambda^\toptzero}^2$, \eqref{equation:fal_ptbeta1},  guarantees that
\begin{equation}\label{equation:convergence_fal_pv1}
\norminf{\bX^\top(\bX\bbeta^\toptzero - \by - \sigma_t\blambda^\toptzero)}
\leq \sqrt{2\epsilon_t}\,\sigma_{\max}(\bX) + \sigma_t.
\end{equation}
Alternatively, if the iterate $\bbeta^\toptzero$ satisfies Condition~\eqref{eq:apgstop_conditions2}, then there exists $\bq^\toptzero \in \partial\normonebig{\bbeta^\toptzero}$ such that 
$\partial P_t(\bbeta^\toptzero) \ni\bg^\toptzero \triangleq \normtwobig{\sigma_t \bq^\toptzero + \bX^\top(\bX\bbeta^\toptzero - \by - \sigma_t\blambda^\toptzero)} \leq \tau_t$, whence by the triangle inequality and the standard bounds on vector norms (Exercise~\ref{exercise:cauch_sc_l1l2}), we have 
\begin{align}
\footnotesize
\norminf{\bX^\top(\bX\bbeta^\toptzero - \by - \sigma_t\blambda^\toptzero)} 
&\leq \norminf{\sigma_t \bq^\toptzero + \bX^\top(\bX\bbeta^\toptzero - \by - \sigma_t\blambda^\toptzero)} + \sigma_t\norminfbig{\bq^\toptzero} \nonumber \\
&\leq \tau_t + \sigma_t, \label{equation:convergence_fal_pv2}
\end{align}
where the last inequality follows from the fact that $\norminfbig{\bq^\toptzero} \leq 1$ for all $\bq^\toptzero \in \partial\normonebig{\bbeta^\toptzero}$; see Exercise~\ref{exercise:sub_norms}.

Since $\blambda^\topzero = \bzero$, $\blambda^\toptone = \blambda^\toptzero - {(\bX\bbeta^\toptzero - \by)}/{\sigma_t}$ for all $t \geq 1$ and $\bX$ has full row rank, it follows from \eqref{equation:svd_stre_bd}, \eqref{equation:convergence_fal_pv1},  \eqref{equation:convergence_fal_pv2}, and  the standard bounds on vector norms (Exercise~\ref{exercise:cauch_sc_l1l2}) that
\begin{align}
\normtwo{\blambda^\toptone} 
&\leq \frac{1}{\sigma_t\sigma_{\min}(\bX)} \normtwo{\bX^\top(\bX\bbeta^\toptzero - \by - \sigma_t\blambda^\toptzero)} \nonumber \\
&\leq \frac{\sqrt{p}}{\sigma_{\min}(\bX)} \left( \max\left\{ \sigma_{\max}(\bX)\sqrt{\frac{2\epsilon_t}{\sigma_t^2}},\; \frac{\tau_t}{\sigma_t} \right\} + 1 \right), \quad \forall\, t \geq 1.
\label{equation:falcv_theta_bound}
\end{align}
Combining the bounds $\frac{\epsilon_t}{\sigma_t^2} \leq B_1$, $\frac{\tau_t}{\sigma_t} \leq B_2$ (from Assumptions~(ii) and (iii)), and~\eqref{equation:falcv_theta_bound}, we have 
\begin{equation}\label{equation:fal_cv_uniform_dual_bound}
\normtwo{\blambda^\toptzero} \leq B_\blambda \triangleq \frac{\sqrt{p}}{\sigma_{\min}(\bX)} \left( \max\left\{ \sqrt{2B_1}\,\sigma_{\max}(\bX),\; B_2 \right\} + 1 \right), \quad \forall\, t > 1,
\end{equation}
which upper bounds $\normonebig{\bbeta^\toptzero}$ by
\begin{equation}\label{equation:_falcv_x_norm_bound}
\normone{\bbeta^\toptzero} \leq \eta_t =\widetildeeta + \frac{\sigma_t}{2}\normtwo{\blambda^\toptzero}^2 \leq B_b \triangleq \widetildeeta + \frac{1}{2}\sigma_1B_\blambda^2.
\end{equation}
Thus, the sequence $\mathcalX$  is bounded and therefore admits at least one limit point. 
Let $\widehatbbeta$ be any such limit point, and let $\sT \subset \integerset_+$ be an infinite subsequence such that $\lim_{t\in\sT} \bbeta^\toptzero = \widehatbbeta$.

\paragraph{Convergence under Condition~\eqref{eq:apgstop_conditions1}.}
Suppose that there exists a further subsequence $\sT_a \subset \sT$ such that for all $t\in\sT_a$, the iterate  $\bbeta^\toptzero$ satisfies Condition~\eqref{eq:apgstop_conditions1}. Then, for $t\in\sT_a$, we have that
\begin{align}
\normonebig{\bbeta^\toptzero} 
&\leq \frac{P_t(\bbeta^\toptzero)}{\sigma_t} 
\leq \frac{P_t(\bbeta_*^\toptzero) + \epsilon_t}{\sigma_t}
\leq \frac{P_t(\bbeta_*) + \epsilon_t}{\sigma_t} 
= \frac{P_t(\bbeta_*)}{\sigma_t} + \frac{\epsilon_t}{\sigma_t} 
\nonumber \\
&= \normone{\bbeta_*} + \frac{\sigma_t}{2}\normtwo{\blambda^\toptzero}^2 + \frac{\epsilon_t}{\sigma_t} \leq \normone{\bbeta_*} + \frac{1}{2}\sigma_t B_\blambda^2 + \sigma_t B_1,
\label{equation:fal_conv_cda_final_bound}
\end{align}
where the first inequality follows from the fact $f_t(\bbeta) \geq 0$, 
the second inequality follows from the stopping condition, 
the third inequality follows from the fact that $P_t(\bbeta_*^\toptzero) \leq P_t(\bbeta_*)$ ($\bbeta_*$ denotes the optimal solution of the $\ell_1$-minimization problem~\eqref{opt:p1}), 
the second equality follows from the fact that $\bX\bbeta_* = \by$, and the last inequality follows from the bounds $\normtwobig{\blambda^\toptzero} \leq B_\blambda$ and $\frac{\epsilon_t}{\sigma_t^2} \leq B_1$. 
Since $\sigma_t \searrow 0$, taking the limit as $t\in\sT_a\rightarrow \infty$, we obtain
\begin{equation}
\normonebig{\widehatbbeta} = \lim_{t\in\sT_a} \normonebig{\bbeta^\toptzero} \leq \normone{\bbeta_*} + \lim_{t\in\sT_a} \left\{ \sigma_t \left( \frac{1}{2} B_\blambda^2 + B_1 \right) \right\} = \normone{\bbeta_*}.
\end{equation}

Next, consider feasibility of the limit point $\widehatbbeta$. 
From the definition of $P_t$ and the stopping Condition~\eqref{eq:apgstop_conditions1}, 
\begin{align*}
\frac{1}{2}\normtwo{\bX\bbeta^\toptzero - \by - \sigma_t\blambda^\toptzero}^2 
&\leq P_t(\bbeta^\toptzero) \leq P_t(\bbeta_*^\toptzero) + \epsilon_t \leq P_t(\bbeta_*) + \epsilon_t \\
&\leq \sigma_t \normone{\bbeta_*} + \frac{1}{2} \normtwo{\sigma_t \blambda^\toptzero}^2 + \epsilon_t \leq \sigma_t \normone{\bbeta_*} + \frac{1}{2} (\sigma_t B_\blambda)^2 + \epsilon_t,
\end{align*}
where the first inequality follows from the fact $\sigma_t \normtwobig{\bbeta^\toptzero} \geq 0$, the third follows from the fact that $P_t(\bbeta_*^\toptzero) \leq P_t(\bbeta_*)$ (where $\sS^\toptzero \cap \argmin_{\bbeta \in \real^p} \{P_t(\bbeta) \}$), the fourth follows from the fact that $\bX\bbeta_* = \by$, and the last follows from the bound $\normtwobig{\blambda^\toptzero} \leq B_\blambda$. Taking the limit as $t\in\sT_a\rightarrow \infty$, we have
$$
\frac{1}{2} \normtwo{\bX\widehatbbeta - \by}^2 \leq 0,
$$
i.e. $\bX\widehatbbeta = \by$, the feasibility of $\widehatbbeta$. 
Since $\widehatbbeta$ is feasible, and $\normonebig{\widehatbbeta} \leq \normone{\bbeta_*}$, it follows that $\widehatbbeta$ is an optimal solution for the $\ell_1$-minimization problem~\eqref{opt:p1}.

\paragraph{Convergence under Condition~\eqref{eq:apgstop_conditions2}.}
Now, consider the complement case, i.e. there exists $T \in \sT$ such that for all $t\in\sT_b \triangleq \sT \cap \{t \geq T\}$, the iterate $\bbeta^\toptzero$  satisfies Condition~\eqref{eq:apgstop_conditions2}. For all $t\in\sT_b$, there exists $\bq^\toptzero \in \partial \normonebig{\bbeta^\toptzero}$ such that
\begin{equation}\label{equation:fal_cv_comppv1}
\partial P_t(\bbeta^\toptzero)
\ni
\normtwo{\sigma_t \bq^\toptzero + \bX^\top(\bX\bbeta^\toptzero - \by - \sigma_t\blambda^\toptzero)} \leq \tau_t.
\end{equation}
Then, for all $t\in\sT_b$, using \eqref{equation:svd_stre_bd} and the triangle inequality again:
\begin{align*}
&\normtwo{\bX\bbeta^\toptzero - \by} \leq \frac{1}{\sigma_{\min}(\bX)} \normtwo{\bX^\top(\bX\bbeta^\toptzero - \by)} \\
&\leq \frac{1}{\sigma_{\min}(\bX)} \left( \normtwo{\bX^\top(\bX\bbeta^\toptzero - \by - \sigma_t\blambda^\toptzero) + \sigma_t \bq^\toptzero} 
+ \normtwo{\sigma_t \bq^\toptzero} + \normtwo{\bX^\top(\sigma_t\blambda^\toptzero)} \right), \\
&\leq \frac{1}{\sigma_{\min}(\bX)} \left( \tau_t + \sigma_t \sqrt{p} + \sigma_{\max}(\bX) \sigma_t B_\blambda \right),
\end{align*}
where the last inequality follows from~\eqref{equation:fal_cv_comppv1}, the bound $\normtwobig{\blambda^\toptzero} \leq B_\blambda$, and fact that $\normtwo{\bq} \leq \sqrt{p}$ for any $\bq \in \partial \normonebig{\bbeta^\toptzero}$ (see Exercise~\ref{exercise:sub_norms}). 
Since $\tau_t\searrow 0$ and $\sigma_t\searrow 0$, taking the limit as $t\in\sT_b\rightarrow\infty$, we have $\normtwobig{\bX\widehatbbeta - \by} \leq 0$, implying $\bX\widehatbbeta = \by$, the feasibility of the limit $\widehatbbeta$.

For all $t\in\sT_b$, $\bq^\toptzero \in \partial \normonebig{\bbeta^\toptzero}$, therefore, $\norminfbig{\bq^\toptzero} \leq 1$. Hence, there exists a subsequence $\sT_b' \subset \sT_b$ such that $\lim_{t\in\sT_b'} \bq^\toptzero = \widehatbq$ exists. One can easily show that $\widehatbq \in \partial \normonebig{\widehatbbeta}$. Dividing both sides of~\eqref{equation:fal_cv_comppv1} by $\sigma_t$, and using the update rule (FAL$_8$) for $\blambda^\toptzero$, we have 
\begin{equation}\label{equation:fal_cv_comppv2}
\normtwo{\bq^\toptzero - \bX^\top\blambda^\toptone} \leq \frac{\tau_t}{\sigma_t},
\end{equation}
for all $t\in\sT_b'$. Since $\lim_{t\in\sT_b'} \bq^\toptzero = \widehatbq$, $\lim_{t\in\sT_b'} \frac{\tau_t}{\sigma_t} = 0$ (Assumption (iii)), and $\bX$ has full row rank, it follows that  $\lim_{t\in\sT_b'} \blambda^\toptone = \widehatblambda$ exists. Taking the limit of both sides of~\eqref{equation:fal_cv_comppv2} along $\sT_b'$, we have
\begin{equation}\label{equation:fal_cv_comppv3}
\widehatbq = \bX^\top\widehatblambda.
\end{equation}
Using~\eqref{equation:fal_cv_comppv3} together with that the fact that $\bX\widehatbbeta = \by$ and $\widehatbq \in \partial\normonebig{\widehatbbeta}$, it follows that the KKT conditions for optimality are satisfied at $\widehatbbeta$; thus, $\widehatbbeta$ is optimal for the $\ell_1$-minimization problem~\eqref{opt:p1} by Lemma~\ref{lemma:dual_cer_p1_prem}. 
This completes the proof.
\end{proof}

In compressed sensing, exact recovery occurs only when the optimization problem $\min\normone{\bbeta}$ subject to $\bX\bbeta = \by$
has a unique solution. Under the same conditions stated in Theorem~\ref{theorem:convergence_fal}, algorithm FAL converges to this unique solution.

\subsection*{Rate of Convergence of FAL}
We now analyze the finite-iteration performance of FAL, which leads to an explicit convergence rate given in Theorem~\ref{theorem:rate_fal}.

\begin{lemma}[Iterate property of FAL]\label{lemma:rate_fal_iterate}
Fix an initial point $\bbeta^\topzero \in \real^p$, a constant $\widetildeeta > 0$ such that $\widetildeeta \geq \normone{\bbeta_*}$, and  a parameter sequence  $\{(\sigma_t, \epsilon_t, \tau_t)\}_{t \in \integerset_+}$ satisfying all  conditions in Theorem~\ref{theorem:convergence_fal}. 
Additionally, assume that for all  $t \geq 1$, $\tau_t \leq c\,\epsilon_t$ for some $0 < c < 1$. Let $\{\bbeta^\toptzero\}_{t \in \integerset_+}$ denote the sequence of iterates generated by algorithm FAL (Algorithm~\ref{alg:fal}) with these parameters. 
Then, for all $t \geq 1$,
\begin{enumerate}[(i)]
\item $\normtwobig{\bX\bbeta^\toptzero - \by} \leq 2B_\blambda \sigma_t$.
\item $\normonebig{\bbeta^\toptzero} - \normone{\bbeta_*} \leq \max\left\{ {B_\blambda^2}/{2} + B_1 \cdot \max\{1,\,2cB_b\},\, \frac{1}{2}{(\frac{\sqrt{p}}{\sigma_{\min}(\bX)} + B_\blambda)^2} \right\} \sigma_t$.
\end{enumerate}
\noindent Note that  $B_\blambda = \frac{\sqrt{p}}{\sigma_{\min}(\bX)} \left( \max\left\{ \sqrt{2B_1}\,\sigma_{\max}(\bX),\, B_2 \right\} + 1 \right)$ and $B_b = \widetildeeta + \frac{1}{2}\sigma_1 B_\blambda^2$, as defined in \eqref{equation:fal_cv_uniform_dual_bound} and~\eqref{equation:_falcv_x_norm_bound}.
\end{lemma}

\begin{proof}[of Lemma~\ref{lemma:rate_fal_iterate}]
\textbf{(i) Feasibility error bound.}
From the proof of Theorem~\ref{theorem:convergence_fal}, we already have the uniform bounds $\normtwobig{\blambda^\toptzero} \leq B_\blambda$ and $\normonebig{\bbeta^\toptzero} \leq B_b$ by~\eqref{equation:fal_cv_uniform_dual_bound} and~\eqref{equation:_falcv_x_norm_bound}, respectively.
The bounds and the dual update in Step~(FAL$_8$) (i.e., $\blambda^\toptone = \blambda^\toptzero - (\bX\bbeta^\toptzero - \by)/\sigma_t$) of algorithm FAL implies that
$$
\normtwo{\bX\bbeta^\toptzero - \by} \leq \normtwo{\bX\bbeta^\toptzero - \by - \sigma_t\blambda^\toptzero} + \sigma_t\normtwo{\blambda^\toptzero} = \sigma_t\normtwo{\blambda^\toptone} + \sigma_t\normtwo{\blambda^\toptzero} \leq 2B_\blambda \sigma_t,
$$
which establishes part (i).

\paragraph{(ii) Objective suboptimality bound.}
The dual update in Step~(FAL$_8$) of algorithm FAL also implies that the augmented Lagrangian function satisfies
$$
P_t(\bbeta^\toptzero) 
= \sigma_t\normonebig{\bbeta^\toptzero} + \frac{1}{2}\normtwo{\bX\bbeta^\toptzero - \by - \sigma_t\blambda^\toptzero}^2 = \sigma_t \left( \normonebig{\bbeta^\toptzero} + \frac{\sigma_t}{2}\normtwo{\blambda^\toptone}^2 \right).
$$
For any $\bbeta_*^\toptzero \in \argmin_{\bbeta \in \real^p} P_t(\bbeta) $, and for all $t \geq 1$,
the above identity implies that
\begin{equation}\label{equation:falineq_lower_bound_bt}
\normonebig{\bbeta^\toptzero} \geq \frac{P_t(\bbeta^\toptzero_*)}{\sigma_t} - \frac{\sigma_t}{2}\normtwo{\blambda^\toptone}^2.
\end{equation}
To lower-bound  $P_t(\bbeta_*^\toptzero)$, consider the following primal-dual pair of the $\ell_1$-minimization problem (see Proposition~\ref{proposition:dual_p1}):
\begin{subequations}\label{equation:fal_basis_pursuit}
\begin{align}
\min_{\bbeta \in \real^p} \normone{\bbeta} \quad & \text{s.t.} \quad  \bX\bbeta = \by;
\label{equation:primal_fal_basis_pursuit} \\
\max_{\bw \in \real^n}  \by^\top \bw\quad & \text{s.t.} \quad \norminf{\bX^\top \bw} \leq 1.
\label{equation:dual_fal_basis_pursuit}
\end{align}
\end{subequations}
Let $\bw_* \in \real^n$ denote an optimal solution of the maximization problem in~\eqref{equation:dual_fal_basis_pursuit}. Next, consider the primal-dual pair of problems corresponding to the augmented Lagrangian  formulation for the $\ell_1$-minimization problem:
\begin{subequations}\label{equation:falpen_basis_pursuit}
\begin{align}
\min_{\bbeta \in \real^p} \sigma\normone{\bbeta} + \tfrac{1}{2}\normtwo{\bX\bbeta-\by - \sigma\blambda}^2;
\label{equation:primal_falpen_basis_pursuit}\\
\max_{\bw \in \real^n}  \sigma(\by + \sigma\blambda)^\top \bw - \tfrac{\sigma^2}{2}\normtwo{\bw}^2 \quad &  \text{s.t.}\quad\norminf{\bX^\top \bw} \leq 1.
\label{equation:dual_falpen_basis_pursuit}
\end{align}
\end{subequations}
Since $\bw_*$ is  also feasible for the maximization problem in~\eqref{equation:dual_falpen_basis_pursuit}, it follows that
\begin{align*}
P_t(\bbeta_*^\toptzero) 
&= \min_{\bbeta \in \real^p} \left( \sigma_t \normone{\bbeta} + \frac{1}{2} \normtwo{\bX\bbeta-\by - \sigma_t \blambda^\toptzero}^2 \right)  
\stackrel{\dag}{\geq} \sigma_t \left( \by^\top \bw_* + \sigma_t (\blambda^\toptzero)^\top \bw_* - \frac{\sigma_t}{2} \normtwo{\bw_*}^2 \right), \\
&\stackrel{\ddag}{\geq} \sigma_t \left( \normone{\bbeta_*} - \sigma_t \normtwo{\blambda^\toptzero} \normtwo{\bw_*} - \frac{\sigma_t}{2} \normtwo{\bw_*}^2 \right), 
\end{align*}
where the inequality ($\dag$) follows from weak duality for primal-dual problems~\eqref{equation:falpen_basis_pursuit}, and the inequality ($\ddag$) follows from strong duality for primal-dual problems~\eqref{equation:fal_basis_pursuit}---i.e. $\by^\top \bw_* = \normone{\bbeta_*}$ ($\bbeta_*$ denotes the optimal solution of the $\ell_1$-minimization problem \eqref{equation:primal_fal_basis_pursuit})---and the Cauchy--Schwarz inequality. 
Thus, the above inequality ($\ddag$) along with \eqref{equation:falineq_lower_bound_bt} implies that
$$
\normonebig{\bbeta^\toptzero} 
\geq \normone{\bbeta_*} - \frac{\sigma_t}{2} \left( 2\normtwo{\blambda^\toptzero} \normtwo{\bw_*} + \normtwo{\bw_*}^2 + \normtwo{\blambda^\toptone}^2 \right) 
\geq \normone{\bbeta_*} - \frac{\sigma_t}{2}\left( \frac{\sqrt{p}}{\sigma_{\min}(\bX)} + B_\blambda \right)^2 ,
$$
where the second inequality follows from the fact that $\sigma_{\min}(\bX)\normtwo{\bw_*}\leq \normtwo{\bX^\top \bw_*} \leq  	 \sqrt{p}\norminf{\bX^\top \bw_*} \leq \sqrt{p}$ (see Exercise~\ref{exercise:cauch_sc_l1l2} and \eqref{equation:svd_stre_bd}).

We now derive matching upper bounds under the two stopping conditions.
\paragraph{Iterate property under Condition~\eqref{eq:apgstop_conditions1}.}
Suppose the iterate $\bbeta^\toptzero$ satisfies the stopping Condition~\eqref{eq:apgstop_conditions1}. 
From inequality~\eqref{equation:fal_conv_cda_final_bound} in the proof of Theorem~\ref{theorem:convergence_fal}, we have
\begin{equation}
\normonebig{\bbeta^\toptzero} \leq \normone{\bbeta_*} + \frac{\sigma_t}{2} \normtwo{\blambda^\toptzero}^2 + \frac{\epsilon_t}{\sigma_t} \leq \normone{\bbeta_*} + \sigma_t \left( \frac{1}{2} B_\blambda^2 + B_1 \right), \label{eq:upper-bound-1}
\end{equation}
since  $\normtwobig{\blambda^\toptzero} \leq B_\blambda$ and ${\epsilon_t}/{\sigma_t^2} \leq B_1$.

\paragraph{Iterate property under Condition~\eqref{eq:apgstop_conditions2}.}
Suppose the iterate $\bbeta^\toptzero$ satisfies the stopping Condition~\eqref{eq:apgstop_conditions2}. Let $\bg^\toptzero \in \partial P_t(\bbeta^\toptzero)$ such that $\normtwo{\bg^\toptzero} \leq \tau_t$. Then the convexity of $P_t$ implies that
$$
P_t(\bbeta^\toptzero) - P_t(\bbeta_*) \leq -(\bg^\toptzero)^\top (\bbeta_* - \bbeta^\toptzero) 
\leq \normtwo{\bg^\toptzero} \normtwo{\bbeta_* - \bbeta^\toptzero} \leq \tau_t \normtwo{\bbeta_* - \bbeta^\toptzero}. 
$$
Now, an argument similar to the one that establishes the bound~\eqref{equation:fal_conv_cda_final_bound}, implies that
\begin{equation}\label{eq:upper-bound-2}
\normonebig{\bbeta^\toptzero} 
\leq \normone{\bbeta_*} + \frac{\sigma_t}{2} \normtwo{\blambda^\toptzero}^2 + \frac{\tau_t}{\sigma_t} \normtwo{\bbeta_* - \bbeta^\toptzero}  
\leq \normone{\bbeta_*} + \sigma_t \left( \frac{1}{2} B_\blambda^2 + 2 c B_b B_1 \right), 
\end{equation}
where we have used the fact that $\normtwobig{\blambda^\toptzero} \leq B_\blambda$, $\tau_t \leq c \epsilon_t$, and ${\epsilon_t}/{\sigma_t^2} \leq B_1$. The upper bound in (ii) follows from~\eqref{eq:upper-bound-1} and~\eqref{eq:upper-bound-2}. 
\end{proof}

Finally, we use the bounds established in the lemma above to derive a convergence rate for Algorithm FAL.
\begin{theoremHigh}[Rate of convergence of FAL \citep{aybat2012first}]\label{theorem:rate_fal}
Fix an $0 < \alpha < 1$. Then there exists---and one can explicitly construct---a sequence of parameters $\{ \sigma_t, \epsilon_t, \tau_t \}_{t \in \integerset_+}$ such that the iterates generated by algorithm FAL (Algorithm~\ref{alg:fal}) are $\epsilon$-feasible, i.e., $\normtwobig{\bX \bbeta^\toptzero - \by} \leq \epsilon$, 
and $\epsilon$-optimal, i.e., $\normonebig{\bbeta^\toptzero} - \normone{\bbeta_*} \leq \epsilon$, for all $t \geq T_{\text{FAL}}(\epsilon) = \mathcalO\left( \log_{\frac{1}{\alpha}} \left( \frac{1}{\epsilon} \right) \right)$. 
\end{theoremHigh}

\begin{proof}[of Theorem~\ref{theorem:rate_fal}]
We begin by rescaling the problem. Define $(\widetildebX, \widetildeby) \triangleq  \frac{1}{\sigma_{\max}(\bX)} (\bX, \by)$. 
Then for the rescaled problem, the smoothness parameter is $L_b = \sigma_{\max}(\widetildebX) = 1$, while the condition number remains unchanged: $\kappa(\widetildebX) = \kappa(\bX) = \sigma_{\max}(\bX) /\sigma_{\min}(\bX)$. 
We apply Algorithm FAL to the rescaled instance $(\widetildebX, \widetildeby)$.

Set the initial parameters as $\sigma_1 \triangleq 1$, $\epsilon_1 \triangleq 2$, 
and define the sequences recursively for $t\geq 1$ by
\begin{align*}
\sigma_{t+1} &= \alpha \cdot \sigma_t, \qquad 
\epsilon_{t+1} = \alpha^2 \cdot \epsilon_t,  \qquad 
\tau_t = \frac{\epsilon_t}{2\max\left\{1,\, \widetildeeta + \frac{9p}{2}\kappa(\bX)^2\right\}} . 
\end{align*}
With this choice, the constants appearing in Lemma~\ref{lemma:rate_fal_iterate} simplify as follows
$$
B_1 = \max_{t \geq 1} \left\{ \frac{\epsilon_t}{\sigma_t^2} \right\} = \frac{\epsilon_1}{\sigma_1^2} = 2, 
\qquad
B_2 = \max_{t \geq 1} \left\{ \frac{\tau_t}{\sigma_t} \right\} = \frac{\alpha^{t-1}}{\max\left\{1,\, \widetildeeta + \frac{9p}{2}\kappa(\bX)^2\right\}} \leq 1.
$$
Therefore, the uniform bounds on $\normtwobig{\blambda^\toptzero}$ and $\normonebig{\bbeta^\toptzero}$ are given by
\begin{align*}
B_\blambda &= \frac{\sqrt{p}}{\sigma_{\min}(\bX)} \left( \max\left\{ \sqrt{2B_1}\,\sigma_{\max}(\bX),\; B_2 \right\} + 1 \right) \leq 3\kappa(\bX)\sqrt{p}; \\
B_b &= \widetildeeta + \frac{1}{2} \sigma_1B_\blambda^2 \leq \widetildeeta + \frac{9p}{2}\kappa(\bX)^2.
\end{align*}
For the rescaled problem $(\widetildebX, \widetildeby)$,  
Lemma~\ref{lemma:rate_fal_iterate} guarantees that for all $t \geq 1$,
$$
\abs{\normonebig{\bbeta^\toptzero} - \normone{\bbeta_*}} \leq \max\left\{ \left( \frac{B_\blambda^2}{2} + B_1 \right),\; \frac{(\sqrt{p}\kappa(\bX) + B_\blambda)^2}{2} \right\} \sigma_1 \alpha^{t-1} \leq 8p\kappa(\bX)^2 \alpha^{t-1},
$$
where we use the fact that $\kappa(\bX) \geq 1$. Thus, $\normonebig{\bbeta^\toptzero} - \normone{\bbeta_*} \leq \epsilon$, for all $t \in \integerset_+$ such that
\begin{equation}\label{equation:falcvrate_eq1}
t \geq \log_{\frac{1}{\alpha}} \left( \frac{8p\kappa(\bX)^2}{\epsilon} \right) + 1. 
\end{equation}
From Lemma~\ref{lemma:rate_fal_iterate}, we also have that for all $t \geq 1$,
$$
\normtwo{\bX\bbeta^\toptzero - \by}
\leq 2B_\blambda \sigma_t 
= 2B_\blambda \sigma_1 \alpha^{t-1} \leq 6\sqrt{p}\kappa(\bX)\alpha^{t-1} .
$$
Thus, $\normtwobig{\bX\bbeta^\toptzero - \by} \leq \epsilon$ for all
\begin{equation}\label{equation:falcvrate_eq2}
t \geq \log_{\frac{1}{\alpha}}\left( \frac{6\sqrt{p}\kappa(\bX)}{\epsilon} \right) + 1. 
\end{equation}
Combining \eqref{equation:falcvrate_eq1} and~\eqref{equation:falcvrate_eq2} shows that for all $\epsilon > 0$, $T_{\text{FAL}}(\epsilon)$, the number of FAL iterations required to compute an $\epsilon$-feasible and $\epsilon$-optimal solution, is at most
\begin{equation}
T_{\text{FAL}}(\epsilon) \leq \left\lceil \log_{\frac{1}{\alpha}} \left( \frac{8p\kappa(\bX)^2}{\epsilon} \right) \right\rceil + 1. 
\end{equation}
This completes the proof.
\end{proof}

\index{Bregman iteration algorithm}
\index{Linearized Bregman iteration algorithm}
\subsection{Bregman Iteration Algorithms}

We present a unified algorithm for solving the four optimization problems discussed in this book (i.e., \eqref{opt:p1}, \eqref{opt:p1_epsilon}, \eqref{opt:ll}, and \eqref{opt:lc}) using the primal-dual framework introduced in Chapter~\ref{chapter:algouni}.
Among these models, the $\ell_1$-minimization formulation serves as a convex relaxation of the $\ell_0$-norm minimization problem, while the constrained LASSO model provides a linear prediction representation that is equivalent to the basis pursuit denoising model.

Now consider the general form of these four optimization problems (see Section~\ref{section:four_rep_prmdual}):
\begin{equation}\label{equation:breg_general_opt_model}
\widehatbbeta = \arg\min_{\bbeta} \left\{ J(\bbeta) + \sigma F(\bbeta) \right\},
\end{equation}
where $J : \sS \to \real$ and $H : \sS\to \real$ are nonnegative convex functions defined on a closed convex set $\sS\subseteq\real^p$. The function $J$ is generally non-smooth, whereas $F$ is differentiable.

A well-known iterative method for solving   problem~\eqref{equation:breg_general_opt_model} is the \textit{Bregman iteration}, which is based on the Bregman distance (Definition~\ref{definition:breg_dist}): 
\begin{equation}
D_\phi^{\bg}(\balpha, \bbeta) = \phi(\balpha) - \phi(\bbeta) - \innerproduct{\bg, \balpha - \bbeta},
\end{equation}
where $\phi(\balpha): \sS\subseteq \real^p\rightarrow \real$ is  a (possibly non-smooth) convex function, and $\bg \in \partial \phi(\bbeta)$ is a subgradient vector of $\phi$ at the point $\bbeta\in\real^p$.
Note that the Bregman distance is not a true metric in the classical sense, as it is generally asymmetric:  $D_\phi^{\bg}(\bbeta, \balpha) \neq D_\phi^{\bh}(\balpha, \bbeta)$.

Consider the first-order Taylor approximation of the non-smooth function $J(\bbeta)$ at the $t$-th iterate $\bbeta^\toptzero$, given by $J(\bbeta) = J(\bbeta^\toptzero) + \innerproductbig{ \bg^\toptzero, \bbeta - \bbeta^\toptzero}$, where $\bg^\toptzero\in\partial J(\bbeta^\toptzero)$. 
The approximation error is precisely measured by the Bregman distance:
\begin{equation}\label{equation:comp_bregman_error}
D_J^{\bg^\toptzero}(\bbeta, \bbeta^\toptzero) = J(\bbeta) - J(\bbeta^\toptzero) - \innerproduct{\bg^\toptzero, \bbeta - \bbeta^\toptzero}.
\end{equation}

\paragrapharrow{Bregman iteration algorithm.}
As early as 1965, Bregman proposed modifying the unconstrained problem~\eqref{equation:breg_general_opt_model} into the following iterative scheme~\citep{bregman1965method}:
\begin{align}
\bbeta^\toptone 
&\gets \arg\min_{\bbeta} \left\{ D_J^{\bg^\toptzero}(\bbeta, \bbeta^\toptzero) + \sigma F(\bbeta) \right\},
\label{equation:bregman_iter_update}\\
&\equiv  \arg\min_{\bbeta} \left\{ J(\bbeta) -\innerproduct{\bg^\toptzero, \bbeta - \bbeta^\toptzero} + \sigma F(\bbeta) \right\}.
\label{equation:bregman_iter_equiv}
\end{align}
This procedure is known as the \textit{Bregman iteration algorithm} or \textit{Bregman iterative regularization}.
Its core idea resembles that of the proximal gradient method, which itself builds upon the projected gradient approach (see Section~\ref{section:proxiGD_inClasso}).
In the following, we discuss the Bregman iteration algorithm and its generalizations.

At the $t$-th iteration, define the objective function of the Bregman subproblem as 
$$
L(\bbeta) = J(\bbeta) - \innerproductbig{\bg^\toptzero, \bbeta - \bbeta^\toptzero} + \sigma F(\bbeta).
$$ 
The first-order optimality condition $\bzero \in \partial L(\bbeta)$ implies $\bzero \in \partial J(\bbeta) - \bg^\toptzero + \sigma \nabla F(\bbeta)$. 
Hence, at the $(t+1)$-th iterative point $\bbeta^\toptone$, we obtain the update rule for the subgradient:
\begin{equation}\label{equation:bregman_update_g}
\bg^\toptone \gets \bg^\toptzero - \sigma \nabla F(\bbeta^\toptone)
\qquad \text{and}\qquad
\bg^\toptone \in \partial J(\bbeta^\toptone).
\end{equation}
Equations~\eqref{equation:bregman_iter_equiv} and~\eqref{equation:bregman_update_g} together constitute the \textit{Bregman iteration algorithm}, originally proposed by \citet{osher2005iterative} in the context of image processing.
The complete procedure is summarized in Algorithm~\ref{alg:breg_itera_gen}.

\begin{algorithm}[ht]
\caption{Bregman Iteration Algorithm \citep{osher2005iterative}}
\label{alg:breg_itera_gen}
\begin{algorithmic}[1]
\Ensure Approximate the solution of  $\min \left\{ J(\bbeta) + \sigma F(\bbeta) \right\}$;
\State \textbf{initialize:} $\bbeta^\topzero = \bzero$, $\bg^\topzero = \bzero$;
\For{$t=0, 1,2,\ldots$}        
\State $\bbeta^\toptone \gets \arg\min_{\bbeta} \left\{ D_J^{\bg^\toptzero}(\bbeta, \bbeta^\toptzero) + \sigma F(\bbeta) \right\}$; \Comment{(BI$_1$)}
\State $\bg^\toptone \gets \bg^\toptzero - \sigma \nabla F(\bbeta^\toptone)$;\Comment{(BI$_2$)}
\State {Exit if} $\normtwobig{\bbeta^\toptone - \bbeta^\toptzero} \leq \varepsilon$;
\EndFor

\State \Return $\bbeta \gets \bbeta^\toptzero$.
\end{algorithmic}
\end{algorithm}

We now recall the convergence properties of the Bregman iteration algorithm, stated below without proof.
\begin{theoremHigh}[Convergence of Bregman iteration algorithm \citep{osher2005iterative}]\label{theorem:conv_breg_gen_osher}
Assume that  $J$ and $F$ are convex functions, and that $F$ is differentiable. 
If a solution to Equation~\eqref{equation:bregman_iter_equiv} exists at each iteration, then the following hold:
\begin{enumerate}[(i)]
\item \textit{Monotonicity.} The sequence $\{F(\bbeta^\toptzero)\}$ is monotonically nonincreasing. Specifically,
$\sigma F(\bbeta^\toptone)
\leq \sigma F(\bbeta^\toptone) + D_{J}^{\bg^\toptzero } (\bbeta^\toptone, \bbeta^\toptzero)
\leq \sigma F(\bbeta^\toptzero)$.
\item The function $F$ will converge to the optimal solution $F(\bbeta^*)$, with the rate estimate $F(\bbeta^\toptzero) \leq F(\bbeta^*) + J(\bbeta^*)/t$.
\end{enumerate}
\end{theoremHigh}

\subsection*{Bregman Iteration Algorithm for Least Squares Objectives}
When $\sigma F(\bbeta) = \frac{1}{2}\normtwo{\bX\bbeta-\by}^2$, the Bregman iteration algorithm admits two equivalent formulations:
\begin{enumerate}[(1)]
\item Initialize $\bbeta^\topzero = \bzero$, $\bg^\topzero = \bzero$; and the iterations are
\begin{subequations}\label{equation:bregite_ls1}
\begin{align}
\bbeta^\toptone &\gets \arg\min_{\bbeta} \left\{ D_J^{\bg^\toptzero}(\bbeta, \bbeta^\toptzero) + \frac{1}{2} \normtwo{\bX\bbeta - \by}^2 \right\}; \label{equation:bregite_ls1_eq1} \\
\bg^\toptone &\gets \bg^\toptzero - \bX^\top(\bX\bbeta^\toptone - \by).
\label{equation:bregite_ls1_eq2} 
\end{align}
\end{subequations}
\item Initialize $\by^\topzero = \bzero$, $\bbeta^\topzero = \bzero$; and the iterations are
\begin{subequations}\label{equation:bregite_ls2}
\begin{align}
\by^\toptone &\gets \by + (\by^\toptzero - \bX\bbeta^\toptzero); 
\label{equation:bregite_ls2_eq1}\\
\bbeta^\toptone &\gets \arg\min_{\bbeta} \left\{ J(\bbeta) + \frac{1}{2} \normtwo{\bX\bbeta - \by^\toptone}^2 \right\}.
\label{equation:bregite_ls2_eq2}
\end{align}
\end{subequations}
\end{enumerate}
It has been shown in \citet{yin2008bregman} that these two versions are mathematically equivalent.

We now show that applying the Bregman iteration algorithm to the least-squares objective effectively approximates the $\ell_1$-minimization problem~\eqref{opt:p1}: 
\begin{equation}\label{opt:p1_breg}
(\text{P}_1)\qquad 
\min_{\bbeta\in\real^p} \mu\normone{\bbeta} \quad \text{s.t.} \quad\bX\bbeta=\by,
\end{equation}
where $\mu$ is a positive number, and we used $\mu=1$ previously.
We show that the Bregman iteration \eqref{equation:bregite_ls1} (or, equivalently, \eqref{equation:bregite_ls2})  generates a sequence
of solutions $\{\bbeta^\toptzero\}$ that converges to an optimal solution $\bbeta^*$ of the $\ell_1$-minimization problem \eqref{opt:p1_breg} in a
finite number of steps; that is, there exists a $T$ such that every $\bbeta^\toptzero$ for $t>T$ is a solution of \eqref{opt:p1_breg}. 

\begin{theoremHigh}[Approximate $\ell_1$-minimization problem \citep{yin2008bregman}]\label{theorem:equiv_ell1}
Consider the iterations in \eqref{equation:bregite_ls1} or \eqref{equation:bregite_ls2} for $J(\bbeta) = \mu \normone{\bbeta}$ and $\sigma F(\bbeta) = \frac{1}{2}\normtwo{\bX\bbeta-\by}^2$ in the Bregman iteration algorithm. Then:
\begin{enumerate}[(a)]
\item Suppose an iterate $\bbeta^\toptzero$ from \eqref{equation:bregite_ls1_eq1} satisfies $\bX\bbeta^\toptzero = \by$; then $\bbeta^\toptzero$ is a solution of the $\ell_1$-minimization problem \eqref{opt:p1_breg}.
\item There exists a finite integer $T < \infty$ such that every iterate $\bbeta^\toptzero$ with $t \geq T$ is an optimal solution of  the $\ell_1$-minimization problem \eqref{opt:p1_breg}.
\end{enumerate}
\end{theoremHigh}
\begin{proof}[of Theorem~\ref{theorem:equiv_ell1}]
\textbf{(a).}
Since $\bg^\toptzero \in \partial J(\bbeta^\toptzero)$,
for any $\bbeta$, by the nonnegativity of the Bregman distance (Remark~\ref{remark:bregnan_dist}), we have
\begin{align}
J(\bbeta^\toptzero) 
&\leq J(\bbeta) - \innerproduct{\bbeta - \bbeta^\toptzero, \bg^\toptzero} 
= J(\bbeta) - \innerproduct{\bbeta - \bbeta^\toptzero, \bX^\top(\by^\toptzero - \bX\bbeta^\toptzero)}  \\
&= J(\bbeta) - \innerproduct{\bX\bbeta - \bX\bbeta^\toptzero, \by^\toptzero - \bX\bbeta^\toptzero} 
= J(\bbeta) - \innerproduct{\bX\bbeta - \by, \by^\toptzero - \by},
\end{align}
where the second equality follows from the equivalence between formulations \eqref{equation:bregite_ls1} and \eqref{equation:bregite_ls2}; see Problem~\ref{prob:equiv_breg_lsobj} and \citet{yin2008bregman}.
Therefore, $\bbeta^\toptzero$ satisfies $J(\bbeta^\toptzero) \leq J(\bbeta)$ for any $\bbeta$ satisfying $\bX\bbeta = \by$; hence, $\bbeta^\toptzero$ is an optimal solution of the $\ell_1$-minimization problem \eqref{opt:p1_breg}.

\paragraph{(b).}
Let $(\sI_+^{(j)}, \sI_-^{(j)}, \sE^{(j)})$ be a partition of the index set $\{1, 2, \ldots, p\}$, and define
\begin{align}
\sU^{(j)} &\triangleq \sU(\sI_+^{(j)}, \sI_-^{(j)}, \sE^{(j)}) = \{ \bbeta \mid  \beta_i \geq 0,\ i \in \sI_+^{(j)};\ \beta_i \leq 0,\ i \in \sI_-^{(j)};\ \beta_i = 0,\ i \in \sE^{(j)} \}; \label{equation:breg_uuions} \\
\sH^{(j)} &\triangleq \min_{\bbeta} \left\{ \frac{1}{2} \normtwo{\bX\bbeta - \by}^2 : \bbeta \in \sU^{(j)} \right\}. 
\end{align}
There are only finitely many such partitions, $(\sI_+^{(j)}, \sI_-^{(j)}, \sE^{(j)})$, and the union of all possible $\sU^{(j)}$'s is $\mathcalH$, the entire space of $\bbeta\in\real^p$.

At iteration $t$, let $(\sI_+^\toptzero, \sI_-^\toptzero, \sE^\toptzero)$ be defined in terms of $\bg^\toptzero$ as follows:
\begin{equation}
\sI_+^\toptzero = \{ i \mid g_i^\toptzero = \mu \},\ 
\sI_-^\toptzero = \{ i \mid g_i^\toptzero = -\mu \},\ 
\sE^\toptzero = \{ i \mid g_i^\toptzero \in (-\mu, \mu) \}. \label{eq:3.21}
\end{equation}
In light of the definition \eqref{equation:breg_uuions} and the fact that $\bg^\toptzero \in \partial J(\bbeta^\toptzero) = \partial (\mu \normonebig{\bbeta^\toptzero})$ and $\mu\sign(\bbeta^\toptzero) \in \partial (\mu \normonebig{\bbeta^\toptzero}) $ by Exercise~\ref{exercise:sub_norms}, 
we have $\bbeta^\toptzero \in \sU^\toptzero$. 
To apply Theorem~\ref{theorem:conv_breg_gen_osher}, we let $\bbeta^*$ satisfy $F(\bbeta^*) = \frac{1}{2}\normtwobig{\bX\bbeta^* - \by}^2 = 0$. 
Using this $\bbeta^*$ in Theorem~\ref{theorem:conv_breg_gen_osher} (ii), we see that for each $j$ with $\sH^{(j)} > 0$, there is a sufficiently large $T_j$ such that $\bbeta^\toptzero$ is not in $\sU^{(j)}$ for $t \geq T_j$. Therefore, letting $T \triangleq \max_j\{T_j : \sH^{(j)} > 0\}$, we have $F(\bbeta^\toptzero) = 0$ for $t \geq T$. That is, $\bX\bbeta^\toptzero = \by$ for $t \geq T$.

Therefore, from \eqref{equation:bregite_ls1_eq2},  the subgradient stabilizes: $\bg^{(T-1)} = \bg^{(T)} = \ldots$, and   \eqref{equation:bregite_ls2_eq1} implies $\by^{(T)} = \by^{(T+1)} = \ldots$. 
Because the minimizers of Bregman iterations \eqref{equation:bregite_ls1_eq1} and \eqref{equation:bregite_ls2_eq2} are not necessarily unique, the iterates $\bbeta^\toptzero$ for $t \geq T$ are not necessarily the same. 
However, it follows from (a) that all $\bbeta^\toptzero$ for $t \geq T$ are optimal solutions of the $\ell_1$-minimization problem \eqref{opt:p1_breg}. 
\end{proof}

\subsection*{Linearized Bregman Iteration Algorithm}
The Bregman iteration algorithm provides an effective framework for solving structured optimization problems. However, at each iteration it requires minimizing the objective $D_J^{\bg^\toptzero}(\bbeta, \bbeta^\toptzero) + \sigma F(\bbeta)$, which can be computationally expensive when $F$ is nonlinear or non-quadratic.
To improve computational efficiency, \citet{yin2008bregman} proposed the \textit{linearized Bregman iteration algorithm}.

The key idea is to approximate the (possibly nonlinear) function $F(\bbeta)$ by its first-order Taylor expansion around the current iterate $\bbeta^\toptzero$ (Theorem~\ref{theorem:linear_approx}):
$$
F(\bbeta)\approx F(\bbeta^\toptzero) + \innerproductbig{\nabla F(\bbeta^\toptzero), \bbeta - \bbeta^\toptzero }.
$$ 
Then, the optimization problem~\eqref{equation:breg_general_opt_model} with $\sigma = 1$ becomes
\begin{equation}
\bbeta^\toptone 
\gets 
\arg\min_{\bbeta} \left\{ D_J^{\bg^\toptzero}(\bbeta, \bbeta^\toptzero) + F(\bbeta^\toptzero) 
+ \innerproduct{\nabla F(\bbeta^\toptzero), \bbeta - \bbeta^\toptzero} \right\}.
\end{equation}
Note that this linear approximation is only accurate near $\bbeta^\toptzero$, and the constant term $F(\bbeta^\toptzero)$ does not affect the minimizer with respect to $\bbeta$. However, to ensure stability and convergence---especially when $F$ is not strongly convex---a quadratic proximity term is added. This leads to the following regularized subproblem:
\begin{subequations}\label{equation:linearized_bregman_all}
\begin{align}
\bbeta^\toptone 
&\gets
\arg\min_{\bbeta} \left\{ D_J^{\bg^\toptzero}(\bbeta, \bbeta^\toptzero) + \innerproduct{\nabla F(\bbeta^\toptzero), \bbeta - \bbeta^\toptzero} + \frac{1}{2\eta} \normtwo{\bbeta - \bbeta^\toptzero}^2 \right\}
\label{equation:linearized_bregman_1}\\
&\equiv \arg\min_{\bbeta} \left\{ D_J^{\bg^\toptzero}(\bbeta, \bbeta^\toptzero) + \frac{1}{2\eta} \normtwo{\bbeta - (\bbeta^\toptzero - \eta \nabla F(\bbeta^\toptzero))}^2 \right\},
\label{equation:linearized_bregman_2}
\end{align}
\end{subequations}
where $\eta$ is positive and serves as the stepsize.
The equivalence follows  because~\eqref{equation:linearized_bregman_1} and~\eqref{equation:linearized_bregman_2} differ by only a constant term that is independent of $\bbeta$.
This update rule was used in  \citet{darbon2007fast} for solving the image deblurring/deconvolution problem $\min_{\bbeta} \{ TV(\bbeta) : \bX\bbeta = \by \}$.

\paragrapharrow{Connection to proximal gradient update and Bregman iteration update.}
The update \eqref{equation:linearized_bregman_2} is similar to the proximal gradient update \eqref{equation:prox_decom2} but is also different  by the use of regularization. 
While the proximal gradient update \eqref{equation:prox_decom2} minimizes $J$, \eqref{equation:linearized_bregman_2} minimizes the Bregman distance $D_J^{\bg^\toptzero}$ based on $J$. On the other hand, \eqref{equation:linearized_bregman_2} differs from \eqref{equation:bregite_ls1_eq1} by replacing the fidelity term $\normtwo{\bX\bbeta - \by}^2/2$ by the sum of its first-order approximation at $\bbeta^\toptzero$ and an $\ell_2$-proximity term at $\bbeta^\toptzero$ when $\sigma  F(\bbeta) = \frac{1}{2} \normtwo{\bX\bbeta - \by}^2$.

\subsection*{Linearized Bregman Iteration Applied to a Least Squares Objective}
In particular, when  $\sigma F(\bbeta) = \frac{1}{2} \normtwo{\bX\bbeta - \by}^2$, we have $\nabla F(\bbeta) = \bX^\top(\bX\bbeta - \by)$. 
Substituting this into \eqref{equation:linearized_bregman_2} yields
\begin{equation}\label{equation:linearized_bregman_3}
\bbeta^\toptone 
\gets 
\arg\min_{\bbeta} \left\{ D_J^{\bg^\toptzero}(\bbeta, \bbeta^\toptzero) 
+ \frac{1}{2\eta} \normtwo{\bbeta - \left( \bbeta^\toptzero - \eta \bX^\top(\bX\bbeta^\toptzero - \by) \right) }^2 \right\}.
\end{equation}
Consider the objective function in \eqref{equation:linearized_bregman_3}, which can be written as
$$
L(\bbeta) = J(\bbeta) - J(\bbeta^\toptzero) - \innerproduct{\bg^\toptzero, \bbeta - \bbeta^\toptzero} 
+ \frac{1}{2\eta} \normtwo{\bbeta - \left( \bbeta^\toptzero - \eta \bX^\top(\bX\bbeta^\toptzero - \by) \right)}^2.
$$
By the subdifferential optimality condition $\bzero \in \partial L(\bbeta)$, it follows that
\begin{align}
\bzero &\in \partial J(\bbeta) - \bg^\toptzero + \frac{1}{\eta} \left( \bbeta - \left( \bbeta^\toptzero - \eta \bX^\top(\bX\bbeta^\toptzero - \by) \right) \right)\\
\implies 
&\bg^\toptone - \bg^\toptzero 
=  - \bX^\top(\bX\bbeta^\toptzero - \by) - \frac{\bbeta^\toptone - \bbeta^\toptzero}{\eta},\label{equation:linbreg_gt1gt}
\end{align}
where $\bg^\toptone \in \partial J(\bbeta^\toptone)$.
Given the initialization  $\bbeta^\topzero=\bg^\topzero = \bzero$,
telescoping \eqref{equation:linbreg_gt1gt} over iterations 0 through $t$ gives the cumulative update for the subgradient:
\begin{equation}\label{equation:breg_g_update_sum}
\bg^\toptone= \sum_{i=1}^{t} \bX^\top(\by - \bX\bbeta^{(i)}) - \frac{\bbeta^\toptone}{\eta}
\triangleq \bu^\toptzero - \frac{\bbeta^\toptone}{\eta}
\end{equation}
where  $\bu^\toptzero \triangleq  \sum_{i=1}^{t} \bX^\top(\by - \bX\bbeta^{(i)})$.

From this, we derive two key iterative formulas. First, rearranging \eqref{equation:breg_g_update_sum} yields the update for $\bbeta$, and the definition of $\bu^\toptone$ gives the update for $\bu$:
\begin{subequations}\label{equation:breg_lin_beta_u}
\begin{align}
\bbeta^\toptone &\gets  \eta (\bu^\toptzero - \bg^\toptone);\label{equation:breg_lin_beta}\\
\bu^\toptone &\gets \bu^\toptzero + \bX^\top(\by - \bX\bbeta^\toptone).\label{equation:breg_lin_u_update}
\end{align}
\end{subequations}
This scheme constitutes the \textit{linearized Bregman iteration algorithm} for solving the composite optimization problem
\begin{equation}\label{equation:opt_lin_breg_gen}
\min_{\bbeta} \left\{ J(\bbeta) + \frac{1}{2} \normtwo{\bX\bbeta - \by}^2 \right\}.
\end{equation}
The algorithm is remarkably simple and efficient: each iteration involves only matrix-vector multiplications and a (possibly simple) subgradient or shrinkage operation---making it highly suitable for large-scale problems.

From Theorem~\ref{theorem:equiv_ell1}, the success of (linearized) Bregman iteration in approximating the $\ell_1$-minimization problem hinges on driving the residual to zero $\bX\bbeta^\toptzero = \by$, i.e., the least squares objective $\normtwobig{\bX\bbeta^\toptzero - \by} \searrow 0$.
The following result establishes convergence guarantees for the linearized variant.

\begin{theoremHigh}[Rate of convergence of linear Bregman iteration for $\ell_1$-minimization \citep{darbon2007fast, yin2008bregman}]\label{theorem:linbreg_xxtop}
Consider the linear Bregman iteration algorithm \eqref{equation:breg_lin_beta_u} applied to the least squares objective $\sigma F(\bbeta) = \frac{1}{2} \normtwo{\bX\bbeta - \by}^2$. 
\begin{enumerate}[(i)]
\item If $\eta \bX^\top \bX \prec \bI$ (i.e., $\bI-\eta \bX^\top \bX$ is strictly positive definite), then the residual norm $\normtwobig{\bX\bbeta^\toptzero - \by}$  decreases monotonically:
$$
\normtwo{\bX\bbeta^\toptone  - \by}^2 + \left( \frac{1}{\eta} - \normtwo{\bX^\top \bX} \right) \normtwo{\bbeta^\toptone  - \bbeta^\toptzero}^2 \leq \normtwo{\bX\bbeta^\toptzero - \by}^2.
$$

\item Suppose further that $J$ is strictly convex, and let $\bbeta^*$ be an optimal solution of $\min\{J(\bbeta) : \bX\bbeta = \by\}$. If $\bbeta^\toptzero \to \bw$, then
\begin{equation}
J(\bw) \leq J(\bbeta^*) + \frac{1}{\eta} \innerproduct{\bbeta^* - \bw, \bw}.
\end{equation}
\item If  $J$ is strictly convex and $\eta \bX\bX^\top \prec 2\bI$ (i.e., $\bI - \frac{\eta}{2} \bX\bX^\top $
is strictly positive definite), then the residual norm $\normtwobig{\bX\bbeta^\toptzero - \by}$ decays exponentially in $t$.
\end{enumerate}

\end{theoremHigh}
\begin{proof}[of Theorem~\ref{theorem:linbreg_xxtop}]
\textbf{(i).} 
Since the Bregman distance used in \eqref{equation:linearized_bregman_3} is nonnegative, we have
$$
\normtwo{\bbeta^\toptone  - \bbeta^\toptzero + \eta \bX^\top(\bX\bbeta^\toptzero - \by)}^2 
\leq \normtwo{\eta \bX^\top(\bX\bbeta^\toptzero - \by)}^2.
$$
Expanding both sides yields
$$
\normtwo{\bbeta^\toptone  - \bbeta^\toptzero}^2 + 2\eta \innerproduct{ \bbeta^\toptone  - \bbeta^\toptzero, \bX^\top(\bX\bbeta^\toptzero - \by)} \leq 0.
$$
After rearrangement, this shows that 
\begin{align*}
\footnotesize
\normtwo{\bbeta^\toptone  - \bbeta^\toptzero}^2 
+ \eta \normtwo{\bX\bbeta^\toptone  - \by}^2 - \eta \innerproduct{\bbeta^\toptone  - \bbeta^\toptzero, \bX^\top\bX(\bbeta^\toptone  - \bbeta^\toptzero) }
\leq \eta \normtwo{\bX\bbeta^\toptzero - \by}^2.
\end{align*}
Using \eqref{equation:svd_stre_bd}, we obtain
$$
\normtwo{\bX\bbeta^\toptone  - \by}^2 + \left( \frac{1}{\eta} \bI - \normtwo{\bX^\top \bX} \right) \normtwo{\bbeta^\toptone  - \bbeta^\toptzero}^2 \leq \normtwo{\bX\bbeta^\toptzero - \by}^2. 
$$

\paragraph{(ii).}
Let again $\bg^\toptzero \in\partial J(\bbeta^\toptzero)$. We have
$$
\bbeta^\toptone  = \arg\min_{\bbeta} \left\{ J(\bbeta) - J(\bbeta^\toptzero) - \innerproduct{\bbeta - \bbeta^\toptzero, \bg^\toptzero} + \frac{1}{2\eta} \normtwo{\bbeta - \bbeta^\toptzero + \eta \bX^\top(\bX\bbeta^\toptzero - \by)}^2 \right\}, 
$$
which, by \eqref{equation:breg_lin_beta_u}, becomes
\begin{align*}
\bbeta^\toptone  &= \arg\min_{\bbeta} \left\{ J(\bbeta) - J(\bbeta^\toptzero) - \innerproduct{\bbeta - \bbeta^\toptzero, \bu^\toptminus} \right. \\
&\quad\quad\quad\quad\quad\quad\quad\quad\left. + \frac{1}{\eta} \innerproduct{\bbeta - \bbeta^\toptzero, \bbeta^\toptzero} 
+ \frac{1}{2\eta} \normtwo{\bbeta - \bbeta^\toptzero + \eta \bX^\top(\bX\bbeta^\toptzero - \by)}^2 \right\}.
\end{align*}
By nonnegativity of the Bregman distance and the definition of $\bu^\toptminus \triangleq  \bX^\top\sum_{i=1}^{t-1} (\by - \bX\bbeta^{(i)})$, we have
\begin{align*}
J(\bbeta^\toptzero) 
&\leq J(\bbeta^*) - \innerproduct{\bbeta^* - \bbeta^\toptzero, \bg^\toptzero} 
= J(\bbeta^*) - \innerproduct{\bbeta^* - \bbeta^\toptzero, \bu^\toptminus} + \frac{1}{\eta} \innerproduct{\bbeta^* - \bbeta^\toptzero, \bbeta^\toptzero} \\
&= J(\bbeta^*) - \innerproduct{ \by - \bX\bbeta^\toptzero, \sum_{i=0}^{t-1} (\by - \bX\bbeta^{(i)})} 
+ \frac{1}{\eta} \innerproduct{\bbeta^* - \bbeta^\toptzero, \bbeta^\toptzero}.
\end{align*}
We will show that $\normtwobig{\bX\bbeta^{(i)} - \by}$ decays exponentially with $i$; then the middle term on the right-hand side above will vanish and the results follow.

\paragraph{(iii).}
By the strict convexity of $J$, there exists a symmetric positive definite operator $\bQ^\toptzero$ (Exercise~\ref{exercise:conv_quad}) such that $\bg^\toptone - \bg^\toptzero = \bQ^\toptzero (\bbeta^\toptone  - \bbeta^\toptzero)$, and \eqref{equation:linbreg_gt1gt} becomes 
$$
\bbeta^\toptone  - \bbeta^\toptzero = -\left( \bQ^\toptzero + \frac{1}{\eta} \bI \right)^{-1} \bX^\top(\bX\bbeta^\toptzero - \by),
$$
which implies 
$$
\bX\bbeta^\toptone  - \by = \left( \bI - \bX \left( \bQ^\toptzero + \frac{1}{\eta} \bI \right)^{-1} \bX^\top \right) (\bX\bbeta^\toptzero - \by). 
$$
To have exponential decay in $t$ of $\normtwobig{\bX\bbeta^\toptzero - \by}$, 
we require the spectral radius of the matrix $\big( \bI - \bX(\bQ^\toptzero + \frac{1}{\eta}\bI)^{-1}\bX^\top \big)$ to be strictly less than one. 
This holds if all eigenvalues of $\bX(\bQ^\toptzero + \frac{1}{\eta}\bI)^{-1}\bX^\top$ lie in $(0,2)$. 
The lower bound follows from the strict positive definiteness of $\bX(\bQ^\toptzero + \frac{1}{\eta}\bI)^{-1}\bX^\top$. 
To have the upper bound, we note that $\bX(\bQ^\toptzero + \frac{1}{\eta}\bI)^{-1}\bX^\top \prec \eta \bX\bX^\top$, which follows from the {Sherman--Morrison--Woodbury formula}; hence, it suffices to have $\eta \bX\bX^\top \prec 2\bI$. 
This completes the proof.
\end{proof}

Finally, suppose $J(\bbeta)\triangleq \mu \cdot j(\bbeta)$ for some base function $j$.
Note that the rate of exponential decay in $t$ of $\normtwobig{\bX\bbeta^\toptzero - \by}$ depends on the value of $\mu$ even when $\eta$ satisfies the condition $\bI - \frac{\eta}{2}\bX\bX^\top \succ 0$. Since $\bQ^\toptzero$ is linear in $\mu$ (assuming $\bbeta^\toptzero$ and $\bbeta^\toptone $ are fixed), we let $\bQ^\toptzero \triangleq \mu \widebarbQ^\toptzero$. When $\mu \gg \frac{1}{\eta}$, $(\bQ^\toptzero + \frac{1}{\eta}\bI)^{-1}$ is dominated by $\frac{1}{\mu} (\widebarbQ^{\toptzero})^{-1}$, so the minimum eigenvalue of $\bX(\bQ^\toptzero + \frac{1}{\eta}\bI)^{-1}\bX^\top$ diminishes linearly to 0 in $\mu$.
Consequently, while the residual still converges exponentially, the rate slows down as $\mu$ grows larger. In practice, this means we can control the speed of convergence by tuning the parameter $\mu$.

When $J(\bbeta) = \mu \normone{\bbeta}$, taking $\mu$ sufficiently large encourages the limit
 point $\bw=\lim_{t\rightarrow \infty}\bbeta^\toptzero$  to approach a solution of the constrained $\ell_1$-minimization problem: $\min\normone{\bbeta}$ subject to $\bX\bbeta = \by$. 
Although $J(\bbeta) = \mu \normone{\bbeta}$ is not strictly convex, it can be approximated arbitrarily well by a strictly convex perturbation (e.g., by adding a small $\ell_2$ term), as discussed in \citet{yin2008bregman}.

To summarize, for solving the $\ell_1$-minimization problem \eqref{opt:p1_breg}, which is traditionally formulated and solved as a constrained linear program, we show that a simple Bregman iterative scheme applied to its unconstrained Lagrangian relaxation $\mu\normone{\bbeta} + \frac{1}{2}\normtwo{\bX\bbeta-\by}^2$ yields an exact solution in a finite number of iterations. 
In practice, using a moderate value of the penalty parameter $\mu$ often leads to convergence in just a few iterations for most problem instances, making the method both efficient and easy to implement.

\begin{problemset}

\item \label{prob:iht_rip_convergence}\textbf{Convergence of IHT under RIP \citep{foucart2013invitation}.}
Let $\bX \in \real^{n \times p}$ satisfy the RIP of order $3k$ with constant 
$$
\delta_{3k} < \frac{1}{2}.
$$
For every $k$-sparse vector $\bbeta^* \in \real^p$, 
show that the sequence $\{\bbeta^\toptzero\}_{t\geq 0}$ defined by IHT (Algorithm~\ref{alg:basic_IHT}) with measurements $\by = \bX\bbeta^*$ converges to $\bbeta^*$.

\item \label{prob:iht_delta_3k_bound}\textbf{Stable recovery of IHT and HTP under RIP \citep{foucart2013invitation}.}
Let $\bX \in \real^{n \times p}$ satisfy the RIP of order $3k$ with constant 
$$
\delta_{3k} < \frac{1}{\sqrt{3}}.
$$
For $\bbeta^* \in \real^p$, $\bepsilon \in \real^n$, and $\sS \subseteq \{1,2,\ldots,p\}$ with $\abs{\sS} = k$,
show that the sequence $\{\bbeta^\toptzero\}_{t\geq 0}$ generated by IHT (Algorithm~\ref{alg:basic_IHT}) or by HTP (Algorithm~\ref{alg:hard_TP}) with $\by = \bX\bbeta^* + \bepsilon$ satisfies
$$
\normtwo{\bbeta^\toptzero_\sS - \bbeta_\sS^*}
\leq \mu^t \normtwo{\bbeta^{(0)}_\sS - \bbeta_\sS^*} 
+ \nu \normtwo{\bX_{\comple{\sS}}\bbeta_{\comple{\sS}} + \bepsilon},
\quad \forall \, t\geq 0,
$$
where $\mu = \sqrt{3}\,\delta_{3k} < 1$, $\nu \leq 2.18/(1 - \mu)$ for IHT, and $\mu = \sqrt{2\delta_{3k}^2/(1 - \delta_{3k}^2)} < 1$, $\nu \leq 5.15/(1 - \mu)$ for HTP.

\item \label{prob:convergence_minimizers}\textbf{Convergence.}
Let $F_t : \sS \subset \real^p \to \real$ be a decreasing sequence of continuous functions that converges pointwise to a function  $F : \sS \to \real$. 
That is, for every $\bbeta\in\sS$, $F_{t+1}(\bbeta) \leq F_t(\bbeta)$ and $\lim_{t \to \infty} F_t(\bbeta) = F(\bbeta)$. 
Assume that $F$ is continuous and coercive on $\sS$. 
Let $\bbeta^\toptzero$ denote a minimizer of $F_t$ over $\sS$. 
Show that every accumulation point of the sequence $\{\bbeta^\toptzero\}_{t \geq 0}$ is a minimizer of $F$. 
Moreover, if $F$ has a unique minimizer, then the entire sequence $\{\bbeta^\toptzero\}$ converges to that minimizer as $t \to \infty$.

\item \label{prob:irls_auxit_epsi2} Show that the function  $F_\varepsilon(\bbeta)$ in \eqref{equation:irls_auxit_epsi2} is strictly convex if $\varepsilon>0$.

\item \label{prob:bound_nonincr} Let $\balpha, \bbeta\in\real^p$, and denote their nonincreasing rearrangements by  $[\balpha]^\downarrow$ and $[\bbeta]^\downarrow $; see \eqref{equation:nonincr_def}. Show that 
\begin{enumerate}[(i)]
\item \label{prob:bound_nonincr1} $\norminf{[\balpha]^\downarrow - [\bbeta]^\downarrow } \leq \norminf{\balpha - \bbeta } $.
\item \label{prob:bound_nonincr2} $\abs{\sigma_k(\balpha)-\sigma_k(\bbeta)} \leq \normone{\balpha - \bbeta}$ for $k\in\{1,2,\ldots,p\}$.
\item \label{prob:bound_nonincr3} $(t-k) [\balpha]_t^\downarrow \leq \normone{\balpha-\bbeta} + \sigma_k(\bbeta)_1$ for any $t>k$.
\end{enumerate}

\item \label{prob:equiv_breg_lsobj} Show that the Bregman iteration algorithm applied to the least squares loss admits the two equivalent formulations given in \eqref{equation:bregite_ls1} and \eqref{equation:bregite_ls2}.

\item Discuss how the linearized Bregman iteration algorithm can be adapted to solve the constrained LASSO problem: $\min_{\bbeta}\frac{1}{2}\normtwo{\bX\bbeta-\by}^2$ subject to $\normone{\bbeta}\leq \Sigma$.


\end{problemset}

\newpage
\vskip 0.2in
\addcontentsline{toc}{chapter}{Bibliography}
\bibliography{bib}

\clearpage
\addcontentsline{toc}{chapter}{Index}
\printindex

\end{document}